\documentclass[a4paper, 10pt, openany]{book}

\usepackage[cp1251]{inputenc}

\usepackage{amsmath,amsthm,amssymb,tabularx,calc}

\usepackage{multicol}

\usepackage[cmtip,arrow,all,2cell]{xy}
\usepackage{pb-diagram,pb-xy}
\usepackage{makeidx}

\usepackage[final]{%graphics,
graphicx} \graphicspath{{pic/}}

\usepackage[english,russian]{babel}

\makeindex

\allowdisplaybreaks %разрешение переноса многострочных формул на новую страницу

\addtocounter{chapter}{-1}

 \makeatletter \@addtoreset {equation}{chapter}
 \renewcommand\theequation
 {\ifnum \c@chapter>\z@ \thechapter.\fi \@arabic\c@equation}
 \makeatother

\renewcommand{\thechapter}{\arabic{chapter}}

\renewcommand{\theequation}{\arabic{chapter}.\arabic{section}.\arabic{equation}}

\theoremstyle{plain}

\newtheorem{tm}{Теорема}[section]

\newtheorem{lm}{Лемма}[section]

\newtheorem{prop}{Предложение}[section]

\newtheorem{cor}{Следствие}[section]

\newtheorem{extm}{$\diamond$}[section]

\theoremstyle{definition}

\newtheorem{df}{Определение}[section]

\newtheorem{rem}{$\bold{!}$}[section]

\newtheorem{ex}{$\diamond$}[section]

\newtheorem{cex}{Контрпример}[section]

\newtheorem{exs}{$\diamond\diamond$}[section]

\newtheorem{er}{$\triangleright$}[section]

\newtheorem{ers}{$\triangleright\triangleright$}[section]

\theoremstyle{remark}

\newtheorem{rmk}{Замечание}

\newcommand{\beq}{\begin{equation}}
\newcommand{\eeq}{\end{equation}}
\newcommand{\btm}{\begin{tm}}
\newcommand{\etm}{\end{tm}}
\newcommand{\blm}{\begin{lm}}
\newcommand{\elm}{\end{lm}}
\newcommand{\bprop}{\begin{prop}}
\newcommand{\eprop}{\end{prop}}
\newcommand{\bcor}{\begin{cor}}
\newcommand{\ecor}{\end{cor}}
\newcommand{\bex}{\begin{ex}}
\newcommand{\eex}{\end{ex}}
\newcommand{\bexs}{\begin{exs}}
\newcommand{\eexs}{\end{exs}}
\newcommand{\bextm}{\begin{extm}}
\newcommand{\eextm}{\end{extm}}
\newcommand{\bcx}{\begin{cex}}
\newcommand{\ecx}{\end{cex}}
\newcommand{\bers}{\begin{ers}}
\newcommand{\eers}{\end{ers}}
\newcommand{\ber}{\begin{er}}
\newcommand{\eer}{\end{er}}
\newcommand{\bdf}{\begin{df}}
\newcommand{\edf}{\end{df}}
\newcommand{\brem}{\begin{rem}}
\newcommand{\erem}{\end{rem}}
\newcommand{\bpr}{\begin{proof}}
\newcommand{\epr}{\end{proof}}

  \newenvironment{list0}[1]
 {
 \begin{list}{}{
 \topsep=5pt %пробел сверху перед первой записью в списке
 \parsep=0pt  %пробел между абзацами в одной записи
 \itemsep=3pt %пробел между записями (дополнительно к \parsep)
 \itemindent=0pt %пробел слева от метки
 \labelwidth=20pt %ширина метки в записи
 \leftmargin=35pt %пробел слева перед первой записью в списке
 \rightmargin=10pt %пробел справа перед первой записью в списке
 \labelsep=10pt %пробел между меткой и записью
 }{#1}
 \end{list}
 }

\newcommand{\bit}{\begin{list0}}
\newcommand{\eit}{\end{list0}}

 \newenvironment{list1}[1]
 {
 \begin{list}{}{
 \itemsep=-1pt %пробел между записями (дополнительно к \parsep)
 \itemindent=20pt %пробел перед меткой
 \topsep=0pt %пробел сверху перед первой записью в списке
 \labelwidth=50pt %ширина метки в записи
 \leftmargin=50pt %пробел слева перед первой записью в списке
 \labelsep=20pt %пробел между меткой и записью
 }{#1}
 \end{list}
 }

\newcommand{\biti}{\begin{list1}}
\newcommand{\eiti}{\end{list1}}

  \newenvironment{list3}[1]
 {
 \begin{list}{}{
 \topsep=5pt %пробел сверху перед первой записью в списке
 \parsep=0pt  %пробел между абзацами в одной записи
 \itemsep=3pt %пробел между записями (дополнительно к \parsep)
 \itemindent=0pt %пробел слева от метки
 \labelwidth=20pt %ширина метки в записи
 \leftmargin=20pt %пробел слева перед первой записью в списке
 \rightmargin=0pt %пробел справа
 \labelsep=10pt %пробел между меткой и записью
 }{#1}
 \end{list}
 }

\newcommand{\biter}{\begin{list3}}
\newcommand{\eiter}{\end{list3}}

 \newenvironment{firstlist}[1]
 {
 \begin{list}{}{
 \topsep=10pt %пробел сверху перед первой записью в списке
 \labelwidth=50pt %ширина метки в записи
 \leftmargin=50pt %пробел слева перед первой записью в списке
 \labelsep=5pt %пробел между меткой и записью
 }{#1}
 \end{list}
 }

 \newenvironment{secondlist}[1]
 {
 \begin{list}{}{
 \topsep=10pt
 \itemindent=-20pt %пробел слева от метки
 \labelwidth=100pt %ширина метки в записи
 \leftmargin=100pt %пробел слева перед первой записью в списке
 \rightmargin=9pt  %пробел справа от записи
 }{#1}
 \end{list}
 }

\newcommand{\bfirstlist}{\begin{firstlist}}
\newcommand{\efirstlist}{\end{firstlist}}

\newcommand{\bsecondlist}{\begin{secondlist}}
\newcommand{\esecondlist}{\end{secondlist}}

 \def\smsize{\scriptsize %мелкий шрифт
 \parskip=1pt %дополнительный вертикальный пробел перед абзацем
 \baselineskip8pt plus 0.3pt minus 0.1pt
 \abovedisplayskip=1pt plus 1pt minus 1pt
 \belowdisplayskip=1pt plus 1pt minus 1pt
 \abovedisplayshortskip=1pt plus 1pt minus 1pt
 \belowdisplayshortskip=1pt plus 1pt minus 1pt
 \topsep=1pt
 }

\def \le {\leqslant}
\def \ge {\geqslant}

\def\({\hskip1pt(}
\def\){)\hskip1pt}
\def\>{\hskip1pt}

\def\e{\varepsilon}

\def\ph{\varphi}

\def\R{{\mathbb{R}}}

\def\tg{\operatorname{tg}}
\def\sgn{\operatorname{sgn}}
\def\arctg{\operatorname{arctg}}
\def\arcctg{\operatorname{arcctg}}
\def\ctg{\operatorname{ctg}}

\def\card{\operatorname{card}}

\def\ch{\operatorname{ch}}

\def\div{\text{\raise4.8pt\hbox{\rlap{\phantom{.}.\phantom{.}}}\raise2.4pt\hbox{\rlap{\phantom{.}.\phantom{.}}}\phantom{.}.\phantom{.}}}

\def \l { \left( }
\def \r { \right) }

\def \lll { \left\{ }
\def \rrr { \right\} }

\def \ml { \left| }
\def \mr { \right| }

\def \N {\mathbb{N}}

\def \Z {\mathbb{Z}}

\def \Q {\mathbb{Q}}

\def \d {\operatorname{\sf{d}}}

\def\Int{\operatorname{\sf{Int}}}

\def\diam{\operatorname{\sf{diam}}}

\def\notarrow{\operatorname{\kern-15pt\longrightarrow\kern-15pt{\not}}}

\def\D{\operatorname{\sf{D}}}

\def\convsupp{\operatorname{\sf{conv\, supp}}}

\def\Supp{\operatorname{\sf{S}}}

\def\Gen{\operatorname{\sf{Gen}}}

\def\VI{\text{\rotatebox{90}{$\le$}}}
\def\IA{\text{\rotatebox{90}{$\ge$}}}

\def\subto{\mathop{\hbox{$\subset$}}\limits_{\kern1pt\raise2pt\hbox{$\smash{\scriptstyle\to}$}}}

\def\subtt{\mathop{\hbox{$\subset$}\kern-4pt\raise-2.3pt\hbox{\resizebox{4pt}{3pt}{$\to$}}}}

\def\suptt{\mathop{\raise-2.3pt\hbox{\resizebox{4pt}{3pt}{$\gets$}}\kern-4pt\hbox{$\supset$}}}

\def\nfork{\mathop{\not\kern1pt\pitchfork}}

\def\V{\resizebox{8pt}{9pt}{$\vee$}}

\newcommand{\inorm}[1]{\left\|#1\right\|^{\int}}
\newcommand{\norm}[1]{\left\|#1\right\|}
\newcommand{\abs}[1]{\left\vert#1\right\vert}

\begin{document}

\title{\bf МАТЕМАТИЧЕСКИЙ АНАЛИЗ}

\author{\bf С.С.Акбаров}

\maketitle

%\smsize

\setcounter{tocdepth}{4} %параметр, объясняющий, насколько подробным
                         %должно быть оглавление

%\nmsize

Настоящий текст представляет собой черновик первого тома учебника по математическому анализу, который автор надеется опубликовать в обозримом будущем. Главную цель изложения автор видит в построении университетского курса анализа, как аксиоматической системы. Принципиально доказываются все (нетривиальные) формулируемые утверждения, за исключением  нескольких фактов общематематического значения (таких, как теорема Гёделя о неполноте или парадокс Банаха-Тарского), приводимых в тексте только для прояснения мотивировок, и никак не проявляющих себя в логической структуре курса.

По способу подачи материал делится на основной, излагаемый текстом в одну колонку, и иллюстративный, представленный двумя колонками. Разница между тем и другим состоит в том, что основной материал задуман, как логически последовательное изложение основных утверждений теории, в котором, в частности, не допускаются ссылки на утверждения, не доказанные на момент цитирования.

В иллюстративном материале, наоборот, приоритетом считается обеспечение читателя достаточным количеством примеров и упражнений для скорейшего привыкания к используемым в основном тексте понятиям и приемам, и, как следствие, уровень логической строгости здесь снижается. Однако, не намного, а только до той планки, на которой некоторые понятия позволяется упоминать существенно раньше, чем они будут формально определены (например, понятие последовательности впервые упоминается на с.\pageref{EX-posledovatelnosti}, хотя определяется только на с.\pageref{DEF:chislov-posledovatelnost}), а некоторым утверждениям позволяется быть сформулированными задолго до того, как они будут аккуратно доказаны в тексте (таковы, например, формулы для простейших геометрических величин в главе \ref{CH-definite-integral} -- площади правильной области на плоскости, объема тела вращения и т.п. -- которые мы по традиции приводим раньше, чем эти величины формально будут определены, и как следствие, доказательство этих формул переносится на несколько глав вперед). При таком подходе, в частности, все, что связано с {\it Исчислением}, включая определения элементарных (/стандартных) функций, описание формальных операций над ними и доказательство связи этих операций с дифференцированием и интегрированием, попадает в двухколоночный текст, поскольку идеологически превращается в иллюстративный материал.

При работе над текстом автор использовал многие идеи доказательств, а также некоторые задачи и упражнения из учебников и пособий, выходивших в разные годы в России, в частности:
 \bit{

 \item[1)] Г.~И.~Архипов, В.~А.~Садовничий, В.~Н.~Чубариков. Лекции по математическому анализу, М.: Высшая школа, 1999.

 \item[2)] Э.~Б.~Винберг. Курс алгебры. М.: Факториал-пресс, 2001.

 \item[3)] Г.~Грауэрт, И.~Либ, В.~Фишер. Дифференциальное и интегральное исчисление. М.: Мир, 1971.

 \item[4)] В.~А.~Зорич. Математический анализ, М.: Фазис, Т.1-2, 1997.

 \item[5)] А.~Н.~Колмогоров, С.~В.~Фомин, Элементы теории функций и функционального анализа. М.: Наука, 1972.

 \item[6)] Л.~Д.~Кудрявцев. Курс математического анализа, М.: Высшая школа, Т.1-2, 1981, Т.3, 1989.

 \item[7)] Л.~Д.~Кудрявцев, А.~Д.~Кутасов, В.~И.~Чехлов, М.~И.~Шабунин. Сборник задач по математическому анализу, М.: Физматлит, Т.1-3, 2003.

 \item[8)] У.~Рудин. Основы математического анализа. М.: Мир, 1976.

 \item[9)] М.~Спивак. Математический анализ на многообразиях. М.: Мир, 1968.

 \item[10)] Е.~Титчмарш. Теория функций. М.: Наука, 1980.

 }\eit\noindent
Автор будет признателен читателям за любые замечания и предложения по улучшению текста.

%\frontmatter

%\input{vvedenie.tex}
%\input{meropriy.tex}

%\mainmatter

\tableofcontents

\part{ФУНКЦИИ ОДНОГО ПЕРЕМЕННОГО}

\chapter{ЧИСЛА И МНОЖЕСТВА}
\label{ch-R&N}

Объектами изучения в математическом анализе являются числовые функции, поэтому
уместно начать учебник с вопроса: {\it что такое число?}

Для читателя это может быть новостью, но ответ на этот вопрос, при всей его
простодушной наивности, совсем не прост. Имеются довольно сложные теории,
объясняющие, как правильно следует определять натуральные, целые, рациональные
и вещественные числа, отталкиваясь, например, от понятия множества. Современный
взгляд на эту тему сформировался в конце 19 -- начале 20 века, и, вопреки
интуитивным ожиданиям, полученное наукой решение этой проблемы, удовлетворяющее
специалистов, оказалось настолько громоздким и хитроумным, что объяснить его
неспециалисту в наши дни стало совершенно невозможно.

Именно поэтому, как может заметить читатель, в учебниках по математическому
анализу {\it определение} понятию числа нигде не дается, а само понятие
вводится {\it аксиоматически} (с разной степенью строгости).\footnote{Мы не
принимаем в расчет <<определения>> типа {\it <<вещественным числом называется
бесконечная десятичная дробь>>}, потому что после этой фразы в воздухе повисает
вопрос {\it <<а что такое бесконечная десятичная дробь?>>}, или, {\it <<кстати,
я давно хотел кого-нибудь спросить, что такое вообще дробь?>>}} Это означает,
что, не объясняя точно, что такое число, Вам говорят, что числа образуют некое
множество, обладающее <<вот такими свойствами>> -- и далее приводится список
свойств чисел, из которых затем выводятся все остальные утверждения теории.
Таким образом, используется тот же прием, что и в геометрии, где понятиям точки
и прямой никаких определений, как известно, не дается, а перечисляются лишь
простейшие свойства этих объектов -- интуитивно понятные {\it аксиомы}, -- из
которых затем с помощью логических умозаключений выводятся следствия -- {\it
теоремы}, часто довольно сложные и неочевидные.

Такой подход в математике называется {\it аксиоматическим методом}, и этим
путем пойдем и мы.

\section{Множества и отображения}\label{SEC-mnozhestva}

Аксиоматические теории в математике делятся на два типа: те, что строятся
непосредственно из правил математической логики --- к ним относятся, например,
современная <<формальная>> арифметика и современная теория множеств, --- и те,
которые используют в своих конструкциях уже построенные аксиоматические теории
--- к ним принадлежит аксиоматическая теория вещественных чисел, которую, в
зависимости от вкуса, можно выводить либо из формальной арифметики (то есть
аксиоматической теории натуральных чисел), либо из (аксиоматической) теории
множеств. С точки зрения потребностей анализа второй путь предпочтительнее,
поскольку в этой науке постоянно приходится иметь дело с множествами,
операциями над ними, их отображениями, и поэтому мы выбираем его.

\paragraph{Об истории теории множеств.}

Читатель наверняка знаком с понятием множества по школе, поэтому мы
можем отложить на некоторое время обсуждение этого термина (см. ниже
пункт ``Наивная теория множеств'' на
с.\pageref{naivnaya-teoriya-mnozhestv}), напомнив лишь его
интуитивное определение: <<совокупность>> объектов, заданных
каким-нибудь правилом.

Эта часть математики родилась в 19 веке, и изначально представляла из себя
систему довольно нестрогих подходов и приемов, используемых математиками в
случаях, когда нужно поглядеть на математический объект, как на <<коллекцию>>
более мелких объектов (элементов), забывая (на время, или совсем) обо всех
остальных его свойствах. Несмотря на видимую непритязательность таких
намерений, в конце 19 -- начале 20 веков вдруг обнаружилось, что из-за
упомянутой нестрогости в тогдашней <<наивной>> теории множеств появляются
противоречия, получившие название {\it парадоксов} или {\it антиномий}.

\noindent\rule{160mm}{0.1pt}
\begin{multicols}{2}

\paragraph*{Антиномия Рассела.}\label{antinomiya-Russela} Самое,
по-видимому, яркое из таких противоречий было найдено в 1902 году
английским математиком Бертраном Расселом. Чтобы его описать, нужно
вспомнить простейшие обозначения и понятия теории множеств. Напомним,
что если объект $x$ принадлежит множеству $A$ (или, что то же самое,
$x$ является {\it элементом} $A$), то записывается это так:
$$
x\in A
$$
Наоборот, если $x$ не принадлежит $A$, то пишут
$$
x\notin A
$$

В теории множеств считается, что сами множества могут быть элементами
других множеств: если $A$ -- какое-то множество, то мы можем
построить какое-нибудь другое множество $B$, для которого $A$ будет
элементом:
$$
A\in B
$$
Наоборот, $B$ можно выбрать так, чтобы $A$ не было его элементом:
$$
A\notin B
$$

И вообще, если нам даны два множества $A$ и $B$, ничто не мешает нам задаться
вопросом, будет ли $A$ элементом $B$:
$$
A\in B \quad \text{или}\quad A\notin B\; ?
$$
При этом, понятное дело, из самого определения символа $\notin$
следует, что две эти возможности исключают друг друга: не бывает так,
чтобы $A\in B$ и одновременно $A\notin B$.

В частности, для произвольного данного множества $A$ мы можем
задуматься, не является ли оно элементом самого себя:
$$
A\in A \quad \text{или}\quad A\notin A\; ?
$$

Оставив в стороне вопрос, бывает ли вообще, чтобы $A\in A$, давайте
рассмотрим множество, которое мы обозначим буквой $R$, и которое
состоит из всевозможных таких множеств $A$, для которых верно
$A\notin A$. То есть, зададим $R$ правилом:
 \beq\label{Russel}
A\in R \;\Longleftrightarrow\; \text{$A$ -- множество, и $A\notin A$}
 \eeq
С точки зрения <<наивной>> теории множеств, $R$ будет неким
множеством, и поэтому мы также можем спросить себя, не будет ли оно
элементом самого себя:
$$
R\in R \quad \text{или}\quad R\notin R\; ?
$$

Так вот, парадокс Рассела состоит в том, что эти утверждения
эквивалентны:
 \begin{multline*}
R\in R \overset{\scriptsize\begin{matrix}\text{вспоминаем}\\
\text{определение
$R$} \\ \eqref{Russel} \\ \downarrow\end{matrix}}{\Longleftrightarrow} \big(\text{$R$ -- одно из тех $A$, для} \\
\text{которых выполняется $A\notin A$}\big)
\underset{\scriptsize\begin{matrix}\uparrow \\ \text{упрощаем}\\
\text{высказывание}
\end{matrix}}{\Longleftrightarrow} R\notin R
 \end{multline*}

На первый взгляд это кажется каким-то математическим фокусом, софизмом, вроде
тех, что мы приводили в ``Предисловии для преподавателей''. Но внимательно
проанализировав ход наших рассуждений, можно убедиться, что никаких формальных
ошибок здесь нет. Точнее сказать, их не будет, если искать их, опираясь на
представления о множествах и операциях над ними, сложившиеся ко времени, когда
появился этот парадокс.

\end{multicols}\noindent\rule[10pt]{160mm}{0.1pt}

Обнаружение антиномий, как можно догадаться, было весьма неприятным открытием
для математиков, поскольку поставило под сомнение само понятие истинности в
этой науке: как может математик быть уверенным, что доказанная им теорема
верна, если используемые им методы в соседней области так легко и наглядно
приводят к очевидному абсурду? Для разрешения возникшего кризиса {\it основания
математики} --- то есть та часть науки, которую, как набор рабочих
инструментов, явно или неявно используют в своих построениях все математики, а
именно, математическая логика и теория множеств (эти области, собственно
говоря, тогда только и стали оформляться как строгие математические дисциплины)
--- были подвергнуты безжалостному пересмотру.

Работа эта была начата в 20-х годах 20 века группой немецких математиков во
главе с Давидом Гильбертом, который поставил перед собой задачу формализовать
математику так, чтобы парадоксы вроде расселовского исчезли и больше уже не
появлялись в будущем. Через некоторое время его программа принесла определенные
положительные результаты: подходящим уточнением определений и введением новых
строгих правил для построения новых объектов удалось добиться устранения всех
накопленных в математике противоречий. Однако новым неприятным сюрпризом,
ставшим одним из главных итогов всей этой деятельности, явилась принципиальная
невозможность доказать, что подобные противоречия не возникнут в будущем.

Этот результат настолько важен, что было бы неприлично не привести его. Для
этого нам понадобятся два термина из математической логики. Аксиоматическую
теорию $\mathcal T$ логики называют {\it противоречивой}, если в ней можно
доказать {\it какое-нибудь} утверждение $P$ и одновременно его отрицание $\neg
P$. Один из простейших результатов этой науки гласит, что если такое возможно
(то есть, если $\mathcal T$ противоречива), то в $\mathcal T$ вообще {\it
любое} выводимое утверждение $P$ обладает этим свойством: вместе с $P$ выводимо
и $\neg P$.

Наоборот, аксиоматическая теория $\mathcal T$ называется {\it
непротиворечивой}, если в ней не бывает так, чтобы вместе с утверждением $P$
можно было вывести и его отрицание $\neg P$. Так вот, в самый разгар работ по
наведению порядка в основаниях математики никому неизвестным молодым
австрийским математиком Куртом Гёделем была доказана следующая замечательная
теорема.

\btm[\bf вторая теорема Гёделя о неполноте] Если какая-то аксиоматическая
теория $\mathcal T$ непротиворечива, включает в себя арифметику и имеет
конечный (или, более общая ситуация, так называемый перечислимый) набор аксиом,
то утверждение о её непротиворечивости невозможно доказать средствами самой
этой теории.
 \etm

Это утверждение означает, ни много, ни мало, что непротиворечивость арифметики,
теории множеств, а вместе с ними и математики в целом, доказать невозможно.
Действительно, для доказательства непротиворечивости любой более или менее
содержательной теории $\mathcal T$ (то есть, теории, включающей арифметику,
как, скажем, теория множеств) нужно, как выясняется, строить новую (более
широкую) теорию ${\mathcal T}'$, после чего немедленно встает вопрос о том, не
будет ли ${\mathcal T}'$ сама противоречива. Это важно, потому что если
${\mathcal T}'$ противоречива, то в ней вместе с утверждением о
непротиворечивости $\mathcal T$ становится выводимым и утверждение о
противоречивости $\mathcal T$ (и таким образом, теорема <<из ${\mathcal T}'$
следует непротиворечивость $\mathcal T$>> обесценивается противоположной
теоремой <<из ${\mathcal T}'$ следует противоречивость $\mathcal T$>>).

Как следствие, гарантий, что в будущем в математике не появится новых
противоречий, быть не может. С другой стороны, если такое противоречие будет
когда-нибудь найдено, то этому ни теорема Гёделя, ни какие-либо другие
результаты в математической логике не препятствуют, и означало бы это событие
лишь, что предложенная в 20 веке формализация оснований математики неудачна, и
следует искать другую.

Результаты Гёделя по существу похоронили надежды Гильберта и его
единомышленников, и единственное утешение в этом отношении для математиков с тех пор состоит в
том, что предпринятые логиками вслед многочисленные попытки найти такое
<<неустранимое>> противоречие в новой <<перестроенной>> математике также не
увенчались успехом, и в настоящее время здесь доминирует убеждение, что
подобные поиски безнадежны.

\paragraph{Соглашения и логическая символика.}

После описанных драматических событий читателю покажутся полной ерундой те проблемы,
из-за которых мы, собственно говоря, и начали разговор о множествах с этого
исторического экскурса. Дело в том, что еще одним неприятным итогом
исследований в этой области явилось то, что в своем новом виде теория множеств
превратилась в неожиданно громоздкую дисциплину, обучить которой быстро
невозможно --- для этого требуется минимум семестровый курс лекций по
математической логике. Именно поэтому аксиоматическая теория множеств обычно не
излагается в курсах анализа: приходится признать разумным ограничиваться
сведениями из <<наивной>> теории множеств. Человека знакомят на примерах со
стандартным набором соответствующих правил и приемов, после чего считается что
он будет способен воспринимать все, что касается множеств в курсе анализа.

Нам, естественно, ничего не остается, как сделать то же самое. В этом параграфе
мы изложим основы наивной теории множеств, на которой потом будем строить
теорию вещественных чисел. Неудобством такого подхода, правда, будет то, что мы
иногда будем приводить в качестве примеров объекты (такие как множество
вещественных чисел $\R$ или множество натуральных чисел $\N$), которые
формально будут определены только в следующих параграфах. Это нужно
исключительно для того, чтобы не обеднять набор иллюстраций излишней погоней за
строгостью: дело в том, что, если следовать формальным правилам, то все
содержательные примеры множеств и отображений (то есть числовые множества и
функции) следует приводить только когда числа уже определены, то есть гораздо
позже, чем вводятся собственно понятия множества и отображения.

Мы нарушаем это правило, следуя общей установке, что в иллюстрациях такое
допустимо (и, как и задумано, помещаем эти примеры, вместе с остальными, в
текст с двумя колонками). Тот же прием мы будем применять в других местах, в
частности, в примерах отображений, функций, и в главе
\ref{CH-definite-integral}, где будут обсуждаться приложения определенного
интеграла.

На протяжении всей книги мы будем пользоваться стандартными символами из
математической логики, смысл которых ясен из следующего перечисления:

\begin{samepage}
 \begin{align*}
 & \text{\sl логический символ} & & \text{\sl его смысл} & & \text{\sl название} \\
 & \neg & & \text{<<не>>} & & \text{логическое отрицание} \\
 & \& & & \text{<<и>>} & & \text{логическое пересечение} \\
 & \V & & \text{<<или>>} & & \text{логическое объединение} \\
 & \Longrightarrow & & \text{<<влечет за собой>>} & &\text{импликация} \\
 & \Longleftrightarrow & & \text{<<эквивалентно>>} & & \text{логическая эквивалентность} \\
 & \forall & & \text{<<для любого>>} & & \text{квантор всеобщности} \\
 & \exists & & \text{<<существует>>} & & \text{квантор существования}
 \end{align*}
\end{samepage}

\paragraph{<<Наивная>> теория множеств.}\label{naivnaya-teoriya-mnozhestv}

Взгляд теории множеств на мир заключается в том, что его объекты находятся друг
с другом в отношении принадлежности: если мы рассматриваем произвольные два
объекта $A$ и $B$, то непременно должно быть истинным или ложным (или независимым от используемой нами системы аксиом)
утверждение, звучащее так:
$$
\text{\it <<$A$ принадлежит $B$>>.}
$$
Если оно истинно, то пишут
$$
A\in B,
$$
а если ложно, то
$$
A\notin B.
$$
Чтобы отличать объекты теории множеств от тех, что рассматриваются в других
науках, а также для удобства построения речевых конструкций, вводятся два
термина: во-первых, сами объекты называются {\it множествами}, а, во-вторых,
если $A\in B$, то говорят, что $A$ является {\it элементом} $B$.

Про отношение принадлежности известно, что оно удовлетворяет некоему списку
условий, называемых {\it аксиомами теории множеств}. Эти аксиомы образуют
довольно сложную систему, обстоятельное изучение которой, как мы уже писали,
заняло бы слишком много времени. Мы поэтому не будем выписывать все эти
аксиомы\footnote{\label{Kelly}Заинтересованного читателя мы отсылаем за подробности к
учебнику Дж.Л.Келли по общей топологии (Келли Дж.Л. Общая топология. М.: Наука, 1981.), где (в Добавлении) в
компактной, хотя и несколько сухой, форме приводится полный список аксиом
теории множеств (их у Келли 9) и обсуждаются их следствия. Предупредим только,
что исходными объектами в ``строгой'' теории множеств являются так называемые
{\it классы}, под которыми понимаются образования, более общие, чем
множества.}, а приведем только три из них, наиболее важные: аксиомы объемности,
выделения и выбора. Третью из них (аксиому выбора) мы сможем сформулировать
только в конце этого параграфа (на
с.\pageref{AX:vybora}\label{axiomy-teorii-mnozhestv}), но первые две дадим
прямо сейчас:

\bit{\it

\item[]{\bf Аксиома объемности:} два множества $A$ и $B$ совпадают, если и
только если они <<имеют одинаковый набор элементов>>, то есть если справедливо
высказывание
$$
\forall x \quad (x\in A\;\Longleftrightarrow\; x\in B)
$$

\item[]{\bf Аксиома выделения:} если дано какое-то переменное
высказывание\footnote{В логике такие высказывания, зависящие от одной
переменной, называются {\it одноместными предикатами}.} $P(x)$ (куда в качестве
$x$ можно подставлять имена объектов теории множеств, то есть имена множеств),
то формула
$$
A=\{x:\, P(x)\}
$$
определяет некое множество $A$, про которое известно, что оно состоит из тех и
только тех объектов $x$, для которых высказывание $P(x)$ истинно.
 }\eit

Частным случаем конструкции, описываемой в аксиоме выделения считается запись
вида
$$
\{...;...;...;...\}
$$
обозначающая множество, элементы которого перечислены внутри фигурных
скобок через точку с запятой. Например, если даны три объекта $a$,
$b$, $c$ (неважно какой природы, но обычно под ними понимают какие-то
множества, потому что в наивной теории множеств считается, что других
объектов не бывает), то запись
$$
\{a;b;c\}
$$
обозначает множество, состоящее из трех объектов $a$, $b$, $c$:
$$
\{a;b;c\}=\{x:\; x=a\; \V\; x=b \; \V \; x=c\}
$$
Если же объектов два, $a$ и $b$, то под $\{a,b\}$ понимается
множество, состоящее из $a$ и $b$. Точно также, если имеется всего
один объект, $a$, то под $\{a\}$ понимается множество из одного
элемента $a$.

\noindent\rule{160mm}{0.1pt}
\begin{multicols}{2}

\bex {\bf Пустое множество.} Следующая формула определяет так
называемое {\it пустое множество}:
$$
\varnothing=\{x:\; x\ne x\}.
$$
В соответствии со сформулированным правилом, расшифровывается эта
запись так: элементами множества $\varnothing$ считаются те и только
те $x$ (объекты теории множеств), которые не равны самому себе.
Поскольку таких $x$ не бывает (высказывание $x=x$ считается верным
всегда), множество $\varnothing$ вообще не содержит никаких элементов
(отчего и называется пустым).
 \eex

\bex {\bf Множества, порождаемые пустым.} Не вполне очевидный факт (это наблюдение приписывается математику Джону фон Нейману)
состоит в том, что из пустого множества $\varnothing$ можно строить
разные другие множества, причем их будет бесконечно много. В качестве первого примера можно рассмотреть множество
$$
\{\varnothing\}=\{x:\; x=\varnothing\}
$$
(состоящее из таких $x$, которые совпадают с $\varnothing$, то есть содержащее
только один элемент: $x=\varnothing$). Это будет новое множество, не
совпадающее с $\varnothing$,
$$
\varnothing\ne\{\varnothing\}
$$
<<Наглядно>> объяснить это можно тем, что $\{\varnothing\}$ содержит один
элемент, а $\varnothing$ -- ни одного.

Затем можно рассмотреть множество
$$
\{\varnothing,\{\varnothing\}\}=\{x:\; x=\varnothing\;\V\;x=\{\varnothing\}\}.
$$
Его элементами будут два множества -- $\varnothing$ и
$\{\varnothing\}$ -- и это будет новое множество.

После этого можно определить множество из трех элементов
 \begin{multline*}
\{\varnothing,\{\varnothing\},\{\varnothing,\{\varnothing\}\}\}=\\=\{x:\;
x=\varnothing\;\V\;x=\{\varnothing\}\;\V\;x=\{\varnothing,\{\varnothing\}\}\},
 \end{multline*}
и так далее.

Помимо множеств, получающихся таким способом, можно рассмотреть
множества вида
$$
\{\{\varnothing\}\},\quad \{\{\{\varnothing\}\}\},\quad ... \quad,
$$
где каждая пара скобок $\{...\}$ обозначает множество, единственным
элементом которого является множество, лежащее внутри скобок. Эти
множества тоже будут отличаться друг от друга (и от множеств,
построенных выше), хотя все будут содержать по одному элементу.

В действительности, хоть это и странно звучит, пустое множество служит
объектом, из которого затем конструируются вообще все на свете числа
(натуральные, целые, рациональные, вещественные) в той теории чисел, которую
строят на фундаменте современной теории множеств (более того, при таком подходе
вообще все объекты в математике оказываются формально построенными из пустого множества
$\varnothing$, поскольку явно или неявно они конструируются из чисел). Нам это не
понадобится в дальнейшем, поэтому подробно останавливаться на этом мы не будем,
но, чтобы заинтриговать читателя, отметим, что в этой теории натуральные числа
(то есть числа, которые в большинстве учебников принято скупо определять как
``числа вида $1, 2, 3,...$'') определяются в точности, как та система множеств,
на которую мы обратили внимание здесь вначале:\footnote{Подробно об этом можно
прочитать в уже упоминавшейся нами в подстрочном примечании \ref{Kelly} книжке
Келли.}
 \begin{align*}
0&:=\varnothing,\\
1&:=\{\varnothing\}=0\cup\{0\},\\
2&:=\{\varnothing,\{\varnothing\}\}=1\cup\{1\},\\
3&:=\{\varnothing,\{\varnothing\},\{\varnothing,\{\varnothing\}\}\}
=2\cup\{2\},...
 \end{align*}
Множество определяемых таким способом чисел, начиная с 1, называется {\it
натуральным рядом}\label{DEF:N-1}, обозначается специальным символом $\N$, и
аккуратное определение ему (однако, как увидит читатель, в другом стиле,
опираясь на аксиоматически вводимую систему вещественных чисел) мы дадим ниже
на странице \pageref{DEF:N}.
 \eex

\bex {\bf Числовые множества.}\label{EX-chislovye-mnozhestva} Известные
читателю со школы числовые множества, такие как множество $\R$ всех
вещественных чисел, или множество $\N$ натуральных чисел, также можно приводить
в качестве примеров, однако их определение в рамках самой теории множеств, как
мы уже говорили, слишком сложно, чтобы его можно было дать в учебнике по
анализу (в следующем параграфе мы идем другим путем --- вводим $\R$
аксиоматически, а затем в \ref{SEC-chisla} \ref{natur-chisla-N} определяем
$\N$, исходя из аксиоматики $\R$). Поэтому мы лишь предлагаем читателю
поверить, что, говоря о множествах, вполне можно держать в воображении
привычные примеры числовых множеств -- и не только $\R$ или $\N$, но и отрезки
$[a;b]$, интервалы $(a,b)$, конечные множества чисел, и так далее.
 \eex

 \bex {\bf Геометрические объекты.} Точно также, геометрические объекты,
изучаемые в школе, такие как фигуры на плоскости
(треугольники, прямоугольники, круги), или объемные фигуры
в трехмерном пространстве, --- тоже можно представлять себе
как примеры множеств, хотя их аккуратное определение в
рамках теории множеств будет еще более сложным, и, как
следствие, мы его также не приводим.
 \eex

\end{multicols}\noindent\rule[10pt]{160mm}{0.1pt}

\paragraph{Пересечение, объединение и разность множеств.}

Если $A$ и $B$ --- два множества, то их {\it пересечением} называется множество
$A\cap B$, состоящее из элементов, принадлежащих обоим этим множествам
$$
A\cap B=\{x:\; x\in A\; \& \; x\in B \},
$$
а {\it объединением} -- множество $A\cup B$, состоящее из элементов,
принадлежащих хотя бы одному из них
$$
A\cup B=\{x:\; x\in A\; \V \; x\in B \}.
$$
Разность $A\setminus B$ множеств $A$ и $B$ определяется как множество,
состоящее из тех элементов множества $A$, которые не принадлежат $B$:
$$
A\setminus B=\{x:\; x\in A\; \& \; x\notin B \}
$$
Наглядно эти понятия обычно иллюстрируются следующими картинками:

\vglue50pt

Если каждый элемент множества $A$ является также элементом множества $B$, то
есть справедливо высказывание, записываемое на языке предикатов формулой
$$
\forall x \quad (x\in A\;\Longrightarrow\; x\in B),
$$
то говорят, что $A$ {\it содержится} в $B$ и записывают это так:
$$
A\subseteq B
$$
Соответствующая картинка выглядит так:

\vglue50pt

Понятно, что множества $A$ и $B$ совпадают тогда и только тогда, когда
$A\subseteq B$ и $B\subseteq A$:
 \beq\label{A-subseteq-B-&-B-subseteq-A}
A=B \quad\Longleftrightarrow\quad A\subseteq B\quad\&\quad B\subseteq A
 \eeq

\noindent\rule{160mm}{0.1pt}\begin{multicols}{2}

\bex Докажем равенство
 \beq\label{(A-setminus-B)-cap-C=(A-cap-C)-setminus-B}
(A\setminus B)\cap C=(A\cap C)\setminus B
 \eeq
Действительно,
 \begin{multline*}
x\in (A\setminus B)\cap C\quad\Longleftrightarrow\quad x\in A\setminus B\;\&\;
x\in C \quad\Longleftrightarrow \\ \Longleftrightarrow\quad x\in A\;\&\;
x\notin B\;\&\; x\in C \quad\Longleftrightarrow \\ \Longleftrightarrow\quad
x\in A\cap C\;\&\; x\notin B\quad\Longleftrightarrow\quad x\in (A\cap
C)\setminus B
 \end{multline*}
 \eex

\bers Докажите равенства:
 \beq\label{A-cup(B-setminus-C)=(A-cup-B)-setminus[(B-cap-C)-setminus-A]}
A\cup(B\setminus C)=(A\cup B)\setminus\Big[(B\cap C)\setminus A\Big]
 \eeq

\eers

\end{multicols}\noindent\rule[10pt]{160mm}{0.1pt}

\paragraph{Упорядоченные пары и декартово произведение множеств.}

Пусть $x$ и $y$ --- какие-нибудь два объекта (неважно какой природы,
но можно считать, что два множества, потому что в наивной теории
множеств считается, что других объектов не бывает).

 \bit{
\item[$\bullet$] {\it Упорядоченной парой} из $x$ и $y$ называется
множество
$$
(x,y):=\{\{x\};\{x;y\}\}
$$
(состоящее из двух элементов, первым из которых будет одноэлементное
множество $\{x\}$, а вторым --- двухэлементное множество $\{x;y\}$).
 }\eit

\bigskip

\centerline{\bf Свойства упорядоченных пар:}

 \bit{\it
\item[$1^\circ$.] Равенство $(x,y)=(u,v)$ выполняется только если
$x=u$ и $y=v$.

\item[$2^\circ$.] Равенство $(x,y)=(y,x)$ выполняется только если
$x=y$.
 }\eit
\bpr Заметим сразу, что $2^\circ$ следует из $1^\circ$: если
$1^\circ$ доказано, то равенство $(x,y)=(y,x)$ будет означать, что
$x=y$ и $y=x$, то есть что $x=y$.

Таким образом, остается доказать лишь $1^\circ$. Равенство
$(x,y)=(u,v)$, означающее
$$
\{\{x\};\{x;y\}\}=\{\{u\};\{u;v\}\},
$$
формально может выполняться только в двух ситуациях:
 \bit{
\item[---] либо если
$$
\begin{cases}\{x\}=\{u\} \\ \{x;y\}=\{u;v\}\end{cases}
$$
и тогда мы получаем
$$
\left\{\begin{array}{ccc}
\{x\}=\{u\} &\Longrightarrow & x=u \\
 & & \{x;y\}=\{u;v\}
\end{array}\right\}
\quad\Longrightarrow\quad \left\{\begin{array}{c}
x=u \\
y=v
\end{array}\right\}
$$

\item[---] либо если
$$
\begin{cases}\{x\}=\{u;v\} \\ \{x;y\}=\{u\} \end{cases};
$$
но такое невозможно, потому что множество с одним элементом $\{x\}$
не может быть равно множеству с двумя элементами $\{u;v\}$ (и точно
также $\{x;y\}$ не может быть равно $\{u\}$).
 }\eit
Значит остается один вариант: $x=u$, $y=v$.
 \epr

 \bit{
\item[$\bullet$] {\it Декартовым произведением} двух множеств $X$ и
$Y$ называется множество $X\times Y$ всех упорядоченных пар вида
$(x,y)$, где $x\in X$, $y\in Y$:
$$
X\times Y:=\{(x,y);\; x\in X,\; y\in Y\}
$$
 }\eit

\noindent\rule{160mm}{0.1pt}\begin{multicols}{2}

\bex {\bf Произведения конечных множеств.} Если множества $X$ и $Y$
конечны, то их декартово произведение $X\times Y$ также конечно, и
поэтому его элементы можно просто выписать. Например, если
$$
X=\{a;b;c\},\quad Y=\{\alpha;\beta\}
$$
то декартово произведение $X\times Y$ состоит из шести упорядоченных
пар, в которых на первом месте стоит какой-нибудь элемент из $X$, а
на втором --- из $Y$:
$$
(a,\alpha),\; (a,\beta),\; (b,\alpha),\; (b,\beta),\; (c,\alpha),\;
(c,\beta)
$$
Для таких случаев специально придуман способ изображения данных в
виде таблицы, заполняя которую невозможно какую-то пару забыть, или
выписать лишнюю:

 \begin{center}
   $\put(4,81){\line(1,-1){48}}$
\begin{tabular}{|c|c|c|}
  \hline
  &  &  \\
  $\begin{array}{cc}  & Y \\ X &  \end{array}$ & $\alpha$  & $\beta$ \\
  &  &  \\
  \hline
  &  &  \\
 $a$ & $(a,\alpha)$ &  $(a,\beta)$ \\
  &  &  \\
  \hline
  &  &  \\
 $b$ & $(b,\alpha)$ &  $(b,\beta)$ \\
  &  &  \\
  \hline
  &  &  \\
 $c$ & $(c,\alpha)$ &  $(c,\beta)$ \\
  &  &  \\
  \hline
\end{tabular}
\end{center}
 \eex

\bex {\bf Декартова плоскость.} Самый употребительный пример
декартова произведения --- декартова плоскость, знакомая читателю со
школы. Это плоскость с фиксированной на ней системой декартовых
координат, то есть двумя перпендикулярными прямыми, называемыми {\it
осями координат}, одна из которых считается горизонтальной (и
называется осью абсцисс), а другая --- вертикальной (и называется
осью ординат), и правилом, по которому первой координатой $x$ всякой
точки $A$ на плоскости считается ее ортогональная проекция на ось
абсцисс, а второй $y$ --- ортогональная проекция на ось ординат.

\vglue60pt

Декартова плоскость считается примером декартова произведения
$X\times Y$, потому что она представляет собой удобную наглядную
реализацию случая, когда $X$ и $Y$ совпадают с числовой прямой
$$
X=Y=\R
$$
Понятно, что в таком случае элементами декартова произведения будут
упорядоченные пары чисел
$$
(x,y),\quad x\in\R,\; y\in\R,
$$
которые удобно представлять себе как координаты точек на декартовой
плоскости, или как сами точки на плоскости.
 \eex

\bex {\bf Умножение на пустое множество.} Неожиданно важным с
методической точки зрения примером для будущего является декартово
произведение пустого множества $\varnothing$ на произвольное
множество $Y$. Поскольку, каким бы ни было $Y$, пары вида
$$
(x,y),\quad x\in \varnothing,\; y\in Y
$$
существовать не могут (так как не бывает таких $x$, чтобы $x\in
\varnothing$), декартово произведение $\varnothing\times Y$ состоит
из пустого набора элементов, то есть само является пустым множеством:
 \beq\label{varnothing-times-Y=-varnothing}
\varnothing\times Y=\varnothing
 \eeq
 \eex

\end{multicols}\noindent\rule[10pt]{160mm}{0.1pt}

\paragraph{Отображения.}

Помимо самих множеств важным понятием теории множеств являются так
называемые отображения. Определяются они так: если даны два множества
$X$ и $Y$, и дано логическое правило, которое можно обозначить,
скажем, буквой $F$, и которое каждому элементу $x$ множества $X$
ставит в соответствие некий (однозначно определенный для данного $x$)
элемент $F(x)$ множества $Y$, то такое правило $F$ называется {\it
отображением}\label{neformalnoe-opred-otobr} (множества $X$ в
множество $Y$), и коротко это записывается так:
$$
F:X\to Y
$$
Множество $X$ при этом называется {\it областью определения}, а $Y$
-- {\it областью значений} отображения $F$. Множество же
упорядоченных пар вида
$$
(x,F(x)),\qquad x\in X
$$
образует некое подмножество в декартовом произведении $X\times Y$,
называемое {\it графиком} отображения $F:X\to Y$.

Принятое нами изначально соглашение следовать наивной теории множеств (с соответствующим ей уровнем строгости) не позволяет вполне ясно уловить
небрежности в этих словах. В действительности, такое определение отображения нельзя считать строгим из-за некоторого произвола в понимании того, что следует считать <<логическим правилом>>. Чтобы обойти эту трудность в <<строгой>> теории множеств принят другой способ определять отображение, в котором не делается различий между самим отображением и его графиком: считается, что отображение $F:X\to Y$ есть подмножество в декартовом произведении $X\times Y$, удовлетворяющее некоему набору условий (и при таком подходе математик избавляется от нужды ссылаться в определении отображения на новый термин, <<логическое правило>>, который формально тоже надо определять). Точная же формулировка этого определения выглядит так:

 \bit{

\item[$\bullet$] {\it
Отображением}\index{отображение}\label{opredelenie-otobrazheniya}
$F:X\to Y$ из множества $X$ в множество $Y$ называется всякое
множество $F$ в декартовом произведении $X\times Y$, удовлетворяющее
следующим двум условиям:

 \bit{
\item[1)] {\bf определенность всюду на $X$:} для любого $x\in X$
имеется некоторый $y\in Y$ такой, что $(x;y)\in F$
$$
\forall x\in X\quad \exists y\in Y \quad (x,y)\in F
$$

\item[2)] {\bf однозначная зависимость от $x\in X$:} если для
некоторых $x\in X$ и $y,z\in Y$ выполняются включения $(x;y)\in F$ и
$(x;z)\in F$, то $y=z$:
$$
\forall x\in X\; \forall y,z\in Y \quad \Bigg[\;\Big(\; (x,y)\in F\;
\& \; (x,z)\in F\;\Big)\quad\Longrightarrow\quad y=z\;\Bigg]
$$
 }\eit
Если эти условия выполнены, то для всякого объекта $x\in X$
существует единственный объект $y\in Y$ такой, что пара $(x,y)$ лежит
в множестве $F$. Такой однозначно определяемый объект $y$ называется
{\it образом} объекта $x$ при отображении $F:X\to Y$ и обозначается
$F(x)$:
$$
(x,y)\in F\quad\Longleftrightarrow\quad F(x):=y
$$
 }\eit

\noindent Двумя важными характеристиками отображения являются его область
определения и область значений:

 \bit{
\item[$\bullet$] Для всякого отображения $f:X\to Y$\label{def-domain-range}
 \bit{
\item[---] множество $X$ называется {\it областью определения} отображения $f$,
для него полезно иметь специальное обозначение:
$$
\D(f):=X
$$

\item[---] множество $Y$ называется {\it областью значений} отображения $f$ (и
оно никак не обозначается),

\item[---] а подмножество в $Y$, состоящее из тех $y$, для которых существует
$x\in X$ такой, что $y=f(x)$ называется {\it образом} отображения $f$, и
обозначается оно так:
$$
\Supp(f):=\{y\in Y: \; \exists x\in X\quad y=f(x)\}
$$
  }\eit
Кроме того,
 \bit{\label{DF:obraz-mnozhestva}
\item[---] если $A\subseteq X$, то множество тех $y\in Y$, для которых
существует $x\in A$ такой, что $y=f(x)$ называется {\it образом множества $A$
при отображении $f$} и обозначается
$$
f(A):=\{y\in Y: \; \exists x\in A\quad y=f(x)\}
$$
(очевидно, $f(X)=\Supp(f)$),

\item[---] если $B\subseteq Y$, то множество тех $x\in X$, для которых $f(x)\in
B$ называется {\it прообразом множества $B$ при отображении $f$} и обозначается
$$
f^{-1}(B):=\{x\in X:\; f(x)\in B\}
$$
(очевидно, $f^{-1}(Y)=X=\D(f)$).
  }\eit

 }\eit\noindent

Среди отображений имеются три важных класса: инъективные, сюръективные и
биективные:

 \bit{
\item[$\bullet$]

Отображение $F:X\to Y$ называется
 \bit{

\item[---] {\it инъективным}, если образы любых двух различных
элементов $X$ также различны:
$$
\forall a,b\in X\qquad a\ne b\quad\Longrightarrow\quad F(a)\ne F(b)
$$

\item[---] {\it сюръективным}, если всякий элемент $y\in Y$ является
образом какого-нибудь элемента $x\in X$:
$$
\forall y\in Y\qquad \exists x\in X\qquad F(x)=y
$$

\item[---] {\it биективным}, если оно инъективно и сюръективно.
 }\eit
Если отображение $f:X\to Y$ биективно, то говорят также, что $f$ осуществляет
взаимно однозначное соответствие между $X$ и $Y$. В этом случае правило
$$
y=f(x)\quad\Longleftrightarrow\quad g(y)=x\qquad (x\in X,\; y\in Y)
$$
определяет так называемое {\it обратное отображение} $g:Y\to X$.

 }\eit

 \bit{
\item[$\bullet$] {\it Композицией отображений} $f:X\to Y$ и $g:Y\to Z$
называется отображение $g\circ f:X\to Z$, определяемое правилом
 \beq\label{DEF:kompozitsiya}
(g\circ f)(x)=g(f(x)),\qquad x\in X
 \eeq
 }\eit

\noindent\rule{160mm}{0.1pt}
\begin{multicols}{2}

\bex {\bf Отображения конечных множеств} явля\-ют\-ся, по-видимому,
простейшими примерами отображений. Рассмотрим в качестве иллюстрации
два множества
$$
X=\{A,B,C\},\quad Y=\{B,C,D,E\}
$$
где $A,B,C,D,E$ -- неповторяющиеся объекты теории множеств, неважно какие,
например,
$$
A=\varnothing,\; B=\{A\},\; C=\{B\},\; D=\{C\},\; E=\{D\}
$$
И определим отображение
$$
f:X\to Y
$$
правилом
$$
A\mapsto D,\quad B\mapsto B,\quad C\mapsto E
$$
то есть
$$
f(A)=D,\quad f(B)=B,\quad f(C)=E
$$
Это действительно будет отображение множества $X$ в множество $Y$, потому что
каждому элементу $X$ оно ставит в соответствие некий (даже явно указанный)
элемент $Y$.

Отображения конечных множеств удобно описывать в виде таблицы, и в нашем случае
эта таблица выглядит так:

\medskip

\vbox{

 \tabskip=0pt\offinterlineskip \halign to \hsize{
 \vrule#\tabskip=2pt plus3pt minus1pt
 &\strut\hfil\;#\hfil & \vrule#&\hfil#\hfil & \vrule#&\hfil#\hfil & \vrule#&\hfil#\hfil & \vrule#\tabskip=0pt\cr
 \noalign{\hrule} height2pt
 &\omit&&\omit&&\omit&&\omit&\cr
 & $x$ && $A$ & & $B$ & & $C$ &\cr
 \noalign{\hrule} height2pt
 &\omit&&\omit&&\omit&&\omit&\cr
 & $f(x)$ && $D$ & & $B$ && $E$ &\cr
 \noalign{\hrule}
 }

}

\medskip
 \noindent
Очевидно, это отображение инъективно, но не сюръективно. Нетрудно привести
также пример неинъективного отображения:

\medskip

\vbox{

 \tabskip=0pt\offinterlineskip \halign to \hsize{
 \vrule#\tabskip=2pt plus3pt minus1pt
 &\strut\hfil\;#\hfil & \vrule#&\hfil#\hfil & \vrule#&\hfil#\hfil & \vrule#&\hfil#\hfil & \vrule#\tabskip=0pt\cr
 \noalign{\hrule} height2pt
 &\omit&&\omit&&\omit&&\omit&\cr
 & $x$ && $A$ & & $B$ & & $C$ &\cr
 \noalign{\hrule} height2pt
 &\omit&&\omit&&\omit&&\omit&\cr
 & $f(x)$ && $D$ & & $D$ && $E$ &\cr
 \noalign{\hrule}
 }

}

\medskip
 \noindent
Оно также не будет сюръективно, так как для наших множеств $X$ и $Y$ вообще
никакое отображение $f:X\to Y$ не может быть сюръективным (из-за того, что
число элементов $Y$ больше, чем число элементов $X$). Чтобы построить
сюръективное отображение можно рассмотреть в качестве $Y$ множество с не более
чем тремя элементами, например,
$$
Y=\{C,D,E\}
$$
Тогда, например, сюръективное отображение $f:X\to Y$ можно будет определить
таблицей

\medskip
\vbox{
 \tabskip=0pt\offinterlineskip \halign to \hsize{
 \vrule#\tabskip=2pt plus3pt minus1pt
 &\strut\hfil\;#\hfil & \vrule#&\hfil#\hfil & \vrule#&\hfil#\hfil & \vrule#&\hfil#\hfil & \vrule#\tabskip=0pt\cr
 \noalign{\hrule} height2pt
 &\omit&&\omit&&\omit&&\omit&\cr
 & $x$ && $A$ & & $B$ & & $C$ &\cr
 \noalign{\hrule} height2pt
 &\omit&&\omit&&\omit&&\omit&\cr
 & $f(x)$ && $C$ & & $D$ && $E$ &\cr
 \noalign{\hrule}
 }
}\medskip

 \noindent
Между прочим, это отображение будет также инъективно, и поэтому биективно.
Соответствующее обратное отображение $g:Y\to X$ будет определяться <<обратной>>
таблицей:

\medskip
\vbox{
 \tabskip=0pt\offinterlineskip \halign to \hsize{
 \vrule#\tabskip=2pt plus3pt minus1pt
 &\strut\hfil\;#\hfil & \vrule#&\hfil#\hfil & \vrule#&\hfil#\hfil & \vrule#&\hfil#\hfil & \vrule#\tabskip=0pt\cr
 \noalign{\hrule} height2pt
 &\omit&&\omit&&\omit&&\omit&\cr
 & $y$ && $C$ & & $D$ & & $E$ &\cr
 \noalign{\hrule} height2pt
 &\omit&&\omit&&\omit&&\omit&\cr
 & $g(y)$ && $A$ & & $B$ && $C$ &\cr
 \noalign{\hrule}
 }
}\medskip

\eex

\bex {\bf Числовые функции.} Важным примером отображений являются
числовые функции. Мы пока не определили это понятие (это делается
ниже в \ref{SEC-numb-function}), поэтому упоминание о нем выглядит
преждевременно, однако, как и с числовыми множествами (пример на с.
\pageref{EX-chislovye-mnozhestva}), это тот случай, когда строгое
следование формальным правилам слишком утомительно, и, чтобы не
злоупотреблять терпением читателя, легче просто предложить ему
поверить, что знакомое ему со школы понятие функции есть частный
случай отображения, и его можно представлять себе как пример.

Скажем, формула
$$
f(x)=x^2
$$
определяет отображение (квадратичную функцию) множества $\R$ всех вещественных
чисел в себя:
$$
f:\R\to\R
$$
Это отображение не будет ни инъективным (потому что $f(-1)=1=f(1)$), ни
сюръективным (потому что ни при каком $x\in\R$ мы не получим $x^2=-1$), поэтому
у него нет обратного отображения.

Напротив, отображение
$$
f:\R\to\R,\qquad f(x)=x^3
$$
является биективным (то есть инъективным и сюръективным), и его обратное
отображение задается формулой
$$
g:\R\to\R,\qquad g(y)=\sqrt[3]{y}
$$
\eex

\bex\label{EX-posledovatelnosti} {\bf Последовательности.} Всякий список
объектов, занумерованных натуральными числами, называется {\it
последовательностью}. Например, десятичная запись числа
$$
476
$$
(как и любого другого натурального числа) есть
последовательность цифр (в данном случае длины 3). На
первом  месте в ней стоит 4, затем 7 и в конце 6. Если
записать это правило в виде таблицы

\medskip
\vbox{
 \tabskip=0pt\offinterlineskip \halign to \hsize{
 \vrule#\tabskip=2pt plus3pt minus1pt
 &\strut\hfil\;#\hfil & \vrule#&\hfil#\hfil & \vrule#&\hfil#\hfil & \vrule#\tabskip=0pt\cr
 \noalign{\hrule} height2pt
 &\omit&&\omit&&\omit&\cr
 & $1$ & & $2$ & & $3$ &\cr
 \noalign{\hrule} height2pt
 &\omit&&\omit&&\omit&\cr
 & $4$ & & $7$ && $6$ &\cr
 \noalign{\hrule}
 }
}\medskip

 \noindent
(где, например, 1 в верхней строчке указывает место, куда
следует записывать упоминаемую внизу цифру 4), то станет
понятно, что последовательность можно считать отображением
(в этом примере --- из множества чисел $\{1,2,3\}$ в
множество цифр $\{0,1,2,3,4,5,6,7,8,9\}$).

Если число целое, то соответствующая ему последовательность десятичных символов
конечна, но если оно нецелое, то последовательность может быть и бесконечной.
Например, число $\frac{15}{7}$ в десятичной записи имеет вид
$$
\frac{15}{7}=2,142857142857142857...
$$
Конечно, в этом случае таблицу полностью выписать невозможно, но правило можно
описать словами: на первом месте стоит 2, на втором -- запятая, потом
последовательно символы 1,4,2,8,5,7, и эта комбинация затем периодически
повторяется.

Последовательностями записываются не только рациональные числа (то есть не
только числа вида $\frac{m}{n}$, где $m,n$ -- целые). Скажем, числа $\sqrt{2}$
и $\pi$ записываются так:
$$
\sqrt{2}=1,4142135623730950488016887242097...
$$
$$
\pi=3,1415926535897932384626433832795...
$$
Сложность здесь только в том, что правила, описывающие эти последовательности,
довольно громоздки, поэтому мы их не будем приводить.

Понятно, что каждая бесконечная последовательность
$$
a_1,a_2,a_3,...
$$
также является отображением,
но только теперь ее областью определения будет множество $\N$ всех натуральных чисел:
$$
f:\N\to A,\qquad f(1)=a_1,\quad f(2)=a_2,\quad ...
$$
(здесь символом $A$ обозначено какое-нибудь множество, содержащее все элементы $\{a_n;\ n\in\N\}$).

Разумеется, последовательности совсем не обязательно должны состоять из чисел
--- бывают последовательности множеств, последовательности функций,
отображений, и вообще любых математических объектов.

Никто не запрещает также рассматривать последовательности нематематических
объектов, и люди обычно так и делают: в жизни последовательность действий
называется планом, последовательность символов алфавита (имеющая смысл) ---
словом, последовательность слов (также осмысленная) --- речью и т.д. \eex

\bex {\bf Движения в геометрии.} Мы опять же, не можем сейчас
определить это понятие, но одним из примеров отображений являются
движения. Как, наверное, помнит читатель, это преобразования
плоскости (или трехмерного пространства), сохраняющие расстояние.
Например, сдвиг на вектор, или центральная симметрия являются
движениями. \eex

\bex\label{pustoe-otobrazhenie} {\bf Пустое отображение.} Формальным
законам логики не противоречит конструкция отображения, определенного
на пустом множестве:
$$
F:\varnothing\to Y
$$
Под таковым естественно понимать подмножество в декартовом
произведении $\varnothing\times Y$, удовлетворяющее условиям 1) - 2)
на с.\pageref{opredelenie-otobrazheniya}. Поскольку, в силу
\eqref{varnothing-times-Y=-varnothing}, множество $\varnothing\times
Y$ пусто,
$$
\varnothing\times Y=\varnothing,
$$
в нем существует единственное подмножество $F\subseteq
\varnothing\times Y$, а именно
$$
F=\varnothing.
$$
Условия 1) - 2) для него выполняются вырожденным образом: поскольку
не существует $x\in X=\varnothing$, для таких $x$ выполняется все что
угодно, в том числе 1) и 2).

Такое отображение $F:\varnothing\to Y$ называется {\it пустым
отображением}\index{отображение!пустое}, причем оно единственно,
потому что $ F=\varnothing$, каким бы ни было $Y$. \eex

\end{multicols}\noindent\rule[10pt]{160mm}{0.1pt}

\paragraph{Семейства множеств и аксиома выбора.}

 \bit{
\item[$\bullet$] Пусть $A$ -- какое-то множество и каждому элементу $\alpha\in
A$ поставлено в соответствие множество $X_\alpha$. Такое соответствие
$\alpha\mapsto X_\alpha$ называется {\it семейством множеств}, и обозначается
оно так:
$$
\{X_\alpha;\; \alpha\in A\}
$$
Понятно, что семейство, то есть правило $\alpha\mapsto X_\alpha$,
есть просто по-другому обозначенное отображение  (с областью
определения $A$, но вместо $f(\alpha)$ мы пишем $X_\alpha$).

\item[$\bullet$] Если $\{X_\alpha;\; \alpha\in A\}$ -- семейство множеств, то
символы $\bigcup_{\alpha\in A} X_\alpha$ и $\bigcap_{\alpha\in A}^\infty
X_\alpha$ обозначают множества, называемые {\it объединением} и {\it
пересечением} множеств $X_\alpha$ и определяются так:
$$
\bigcup_{\alpha\in A} X_\alpha=\{x:\quad\exists \alpha\in A\quad x\in
X_\alpha\}\qquad\qquad \bigcap_{\alpha\in A} X_\alpha=\{x:\quad\forall
\alpha\in A\quad x\in X_\alpha\}
$$
 }\eit

Отметим следующие формулы:
\begin{align}
\Big(\bigcup_{\alpha\in A} X_\alpha\Big)\setminus Y &=\bigcup_{\alpha\in A}
\Big(X_\alpha\setminus Y\Big) & \Big(\bigcap_{\alpha\in A}
X_\alpha\Big)\setminus Y &=\bigcap_{\alpha\in A} \Big(X_\alpha\setminus
Y\Big) \label{UX/Y} \\
Y\setminus \Big(\bigcup_{\alpha\in A} X_\alpha\Big) &=\bigcap_{\alpha\in A}
\Big(Y\setminus X_\alpha\Big) & Y\setminus \Big(\bigcap_{\alpha\in A}
X_\alpha\Big)&=\bigcup_{\alpha\in A} \Big(Y\setminus X_\alpha\Big)\label{Y/UX}
\end{align}
 \bpr
Эти формулы доказываются одинаково. Покажем, например, как это
делается в случае с первой формулой в \eqref{Y/UX}:
 \begin{multline*}
y\in Y\setminus \Big(\bigcup_{\alpha\in A} X_\alpha\Big)
\quad\Longleftrightarrow\quad y\in Y\quad\&\quad \nexists \alpha\in A\quad y\in
X_\alpha \quad\Longleftrightarrow \quad
y\in Y\quad\&\quad \forall \alpha\in A\quad y\notin X_\alpha \quad\Longleftrightarrow\\
\Longleftrightarrow \quad
y\in Y\quad\&\quad \forall \alpha\in A\quad y\notin X_\alpha \quad\Longleftrightarrow\quad
 \forall \alpha\in A\quad y\in Y\quad\&\quad y\notin X_\alpha \quad\Longleftrightarrow
 \\ \Longleftrightarrow\quad
\forall \alpha\in A\quad y\in Y\setminus X_\alpha \quad\Longleftrightarrow\quad
y\in \bigcap_{\alpha\in A} \Big(Y\setminus X_\alpha\Big)
 \end{multline*}
 \epr

Теперь мы можем сформулировать обещанную на с.\pageref{axiomy-teorii-mnozhestv}
аксиому выбора. Интуитивно это утверждение выглядит довольно естественным,
поэтому может быть неожиданным, что оно оказывается независимым от остальных
аксиом теории множеств.

\bit{\it

\item[]{\bf Аксиома выбора:}\label{AX:vybora} для любого семейства непустых
множеств $\{X_\alpha;\; \alpha\in A\}$ найдется (хотя бы одно) отображение
$F:A\to\bigcup_{\alpha\in A}X_\alpha$, удовлетворяющее условию:
$$
\forall \alpha\in A \quad F(\alpha)\in X_\alpha
$$

}\eit

\noindent\rule{160mm}{0.1pt}\begin{multicols}{2}

 \paragraph*{Аксиома счетного выбора.}
 \biter{

\item[$\bullet$] Если в качестве индексного множества $A$ берется множество
$\N$ натуральных чисел (напомним, что наивно мы определили $\N$ на
с.\pageref{DEF:N-1}, а аккуратное определение пообещали дать на
с.\pageref{DEF:N}), то мы получаем частный случай семейства -- {\it
последовательность множеств} $\{X_n\}$.

 }\eiter

Объединение и пересечение последовательности множеств $\{X_n\}$
имеют специальные обозначения:
$$
\bigcup_{n=1}^\infty X_n:=\bigcup_{n\in\N} X_n=\{x:\quad\exists n\in\N\quad
x\in X_n\}
$$
$$
\bigcap_{n=1}^\infty X_n:=\bigcap_{n\in\N} X_n=\{x:\quad\forall n\in\N\quad
x\in X_n\}
$$

Для последовательностей множеств мы получаем частный случай аксиомы выбора, так называемую {\it аксиому счетного выбора}:

\biter{\it

\item[]{\bf Аксиома счетного выбора:}\label{AX:schentnogo-vybora} для любой последовательности непустых
множеств $\{X_n;\; n\in\N\}$ найдется (хотя бы одна) последовательность $\{x_n;\ n\in\N\}\subseteq\bigcup_{n=1}^\infty X_n$ такая, что для всякого $n\in\N$ объект $x_n$ является элементом множества $X_n$:
$$
\forall n\in\N \quad x_n\in X_n
$$

}\eiter

\paragraph*{Бесконечнократный выбор.}\label{beskonechnokratnyi-vybor}

Исторически в анализе утвердился стиль, в котором аксиома выбора (и даже ее
ослабленный вариант, аксиома счетного выбора) применяется почти исключительно в
комбинации с формулируемыми ниже теоремами  \ref{defin-induction} и
\ref{defin-induction-chast} об определениях полной и частной индукцией.
Появляющийся в результате прием {\it бесконечнократного выбора}, без ссылок на
мотивы, по которым его становится возможным использовать, сразу преподносится
слушателю, как естественный инструмент логики, наподобие употребления
логических связок или кванторов, а в счетном случае (то есть там, где
применяется аксиома счетного выбора) объявляется просто одним из способов
рассуждения по индукции (терминологически даже не делается различий между
использованием самих теорем \ref{defin-induction} и \ref{defin-induction-chast}
и их комбинаций с аксиомой выбора -- и тот, и другой прием обычно называется
{\it определением по индукции}; как следствие, например, способ определить
факториал $n!$ правилом \eqref{DEF:n!} подается как формально тот же самый
прием, что и, скажем, способ доказательства критерия бесконечности (теорема
\ref{TH:krit-besk-mnozh})).  Это затушевывает роль аксиомы выбора и этим
объясняется, почему у математиков так часто возникает впечатление, что эта
аксиома вовсе не применяется в классическом анализе, за исключением нескольких
экзотических ситуаций, вроде конструкции Витали неизмеримого по Лебегу
множества, знаменитых как раз разгоревшимися вокруг них спорами о законности
применяемых здесь методов.

Мы постараемся убедить читателя в ошибочности этого мнения на примере
нескольких первых утверждений, в доказательстве которых применяется аксиома
выбора:\label{beskonechnokratnyj-vybor} в критерии бесконечности множества
(теорема \ref{TH:krit-besk-mnozh}), затем в главе \ref{ch-x_n} при
доказательстве свойств $2^0$ и $3^0$ подпоследовательностей
(с.\pageref{podposledovatelnosti}) и наконец, при доказательстве теоремы
Больцано-Вейерштрасса \ref{Bol-Wei}. В них мы просто предложим для сравнения
два типа доказательства: классическое, основанное на бесконечнократном выборе,
и более формальное, с явными ссылками на аксиому выбора и теоремы об
определении индукцией.

В дальнейшем мы все же будем придерживаться традиционного стиля, бесспорным
достоинством которого является его наглядность. Мы полагаем, что обсужденных
примеров будет достаточно, чтобы показать читателю, что {\it всякий раз, когда
в доказательстве присутствует бесконечнократный выбор, под ним скрывается
аксиома выбора, обычно в комбинации с теоремами об определении индукцией}.
Внимательный человек, как мы надеемся, сможет после наших объяснений
восстановить при необходимости аккуратное доказательство, и мы считаем, этого
достаточно для осмысленного применения приема бесконечнократного выбора.

\paragraph*{Альтернатива: аксиома детерминированности.}

Заводя речь об аксиоме выбора, нельзя не упомянуть о спорах, которые она
породила среди математиков. Будучи предложена в 1904 году Эрнестом Цермело, как
инструмент, позволяющий вполне упорядочить произвольное множество, эта аксиома
сразу же стала объектом критики из-за своих многочисленных следствий, некоторые
из которых (как упомянутый выше пример Витали) противоречат интуитивным
ожиданиям и создают впечатление, что отказ от этой аксиомы или замена ее на
что-то более подходящее, позволили бы существенно упростить общую картину
математики.

Автору очевидно, что главную роль в этих спорах сыграл упоминавшийся выше
традиционно небрежный стиль доказательств в анализе, при котором
бесконечнократный выбор ошибочно объявляется стандартным приемом логики, либо
прямым следствием метода математической индукции, без ссылок на аксиому (хотя
бы  счетного) выбора. Во всяком случае, именно этим следует объяснять то, что
даже сейчас, спустя столетие, упреки в ``неконструктивности'' того или иного
доказательства все еще можно слышать от математиков.

Смысл эпитета ``неконструктивный'', которым принято выражать главную претензию
к аксиоме выбора, понять, по-видимому, невозможно. Но пытаться понять, чем эта
аксиома может не понравиться, можно, и эта задача даже становится легкой, если
перед глазами имеется альтернатива, способная заменить критикуемое решение
чем-то более предпочтительным.

Такая альтернатива была предложена намного позже, чем сама аксиома выбора, в
1962 году Яном Мычельским и Гуго Штейнгаузом, и называется она {\it аксиомой
детерминированности}. Чтобы объяснить, что это такое, обозначим символом
$\N^{\N}$ множество всевозможных бесконечных последовательностей натуральных
чисел:
$$
\{x_n;\ n\in\N\}\in\N^{\N}\quad\Leftrightarrow\quad \forall n\in\N\quad x_n\in\N
$$
Аккуратная формулировка аксиомы детерминированности выглядит довольно сложно:

\biter{\it

\item[]{\bf Аксиома детерминированности:}\label{AX:determinacy}
для любого множества $X\subseteq\N^{\N}$ выполняется одно из следующих двух утверждений:
    \biter{
\item[--] либо существует отображение $\sigma$, которое каждой конечной последовательности натуральных чисел четной длины, а также пустому множеству (рассматриваемому как последовательность нулевой длины) ставит в соответствие натуральное число
    $$
    \begin{cases}
    \varnothing\quad\mapsto\quad \sigma(\varnothing)\in\N \\
(x_1,...,x_{2k}) \quad \mapsto\quad \sigma(x_1,...,x_{2k})\in\N
\end{cases}
    $$
    таким образом, что для любой последовательности $\{b_n\}\in\N^{\N}$ последовательность $\{a_n\}$, определенная индуктивным правилом
    $$
    \begin{cases}
   a_1=\sigma(\varnothing) \\
   a_{k+1}=\sigma(a_1,b_1,...,a_k,b_k)
\end{cases}
    $$
    порождает последовательность
    $$
    (a_1,b_1,...,a_k,b_k,...)
    $$
    принадлежащую $X$;

\item[--] либо существует отображение $\tau$, которое каждой конечной последовательности натуральных чисел нечетной длины ставит в соответствие натуральное число
    $$
(x_1,...,x_{2k-1}) \quad \mapsto\quad \sigma(x_1,...,x_{2k-1})\in\N
    $$
    таким образом, что для любой последовательности $\{a_n\}\in\N^{\N}$ последовательность $\{b_n\}$, определенная индуктивным правилом
    $$
    \begin{cases}
   b_1=\sigma(a_1) \\
   b_{k+1}=\sigma(a_1,b_1,a_2,b_2,...,a_k)
\end{cases}
    $$
    порождает последовательность
    $$
    (a_1,b_1,...,a_k,b_k,...),
    $$
    не принадлежащую $X$.
    }\eiter
 }\eiter

Смысл этого громоздкого утверждения становится понятен, если его
интерпретировать, как гипотезу о возможных вариантах развития некоей
бесконечношаговой игры. Представим себе, что нам дано некое множество $X$
бесконечных последовательностей натуральных чисел. Зная описание этого
множества $X$, два игрока A и B ведут игру, правила которой выглядят так:
сначала игрок A записывает натуральное число $a_1$, затем игрок $B$ (глядя на
$a_1$) записывает число $b_1$, затем A (глядя на $a_1$ и $b_1$) записывает
число $a_2$, и так далее. После бесконечного числа ходов получается некая
бесконечная последовательноть
$$
(a_1,b_1,a_2,b_2,a_3,b_3,...)
$$
Если эта последовательность лежит в множестве $X$, то считается, что выиграл
игрок A, в противном случае победителем считается игрок B.

Аксиома детерминированности утверждает, что, {\it каким бы ни было исходное
множество $X$, либо у игрока A, либо у игрока B найдется выигрывающая
стратегия}. Для игрока A такой стратегией будет отображение
    $$
    \begin{cases}
\varnothing\quad\mapsto\quad \sigma(\varnothing)\in\N \\
(x_1,...,x_{2k}) \quad \mapsto\quad \sigma(x_1,...,x_{2k})\in\N
\end{cases}
    $$
позволяющее при любых ходах противника $b_1,b_2,...$ делать очередной ход по правилу
    $$
    \begin{cases}
   a_1=\sigma(\varnothing) \\
   a_{k+1}=\sigma(a_1,b_1,...,a_k,b_k)
\end{cases}
    $$
и при этом должно получиться, что итоговая последовательность
    $$
    (a_1,b_1,...,a_k,b_k,...)
    $$
лежит в $X$.

А для игрока B такой стратегией будет отображение
    $$
(x_1,...,x_{2k-1}) \quad \mapsto\quad \sigma(x_1,...,x_{2k-1})\in\N
    $$
позволяющее при любых ходах противника $a_1,a_2,...$ делать очередной ход по правилу
    $$
    \begin{cases}
   b_1=\sigma(a_1) \\
   b_{k+1}=\sigma(a_1,b_1,a_2,b_2,...,a_k)
\end{cases}
    $$
так, что итоге получится последовательность
    $$
    (a_1,b_1,...,a_k,b_k,...)
    $$
наоборот, не лежащая в $X$.

Удивительный факт состоит в том, что эта экзотическая гипотеза о существовании
выигрывающей стратегии в непонятно кому нужной формальной игре, как
оказывается, может составить конкуренцию на первый взгляд гораздо более простой
и осмысленной аксиоме выбора. Мы упомянем здесь только один результат, чтобы
заинтриговать (подготовленного) читателя:\footnote{Подробности можно прочитать
в книжке: В.Г.Кановей, Аксиома выбора и аксиома детерминированности. М.: Наука,
1984.} {\it если заменить аксиому выбора аксиомой детерминированности, то
всякое множество на прямой $\R$ станет измеримым по Лебегу}.

\end{multicols}\noindent\rule[10pt]{160mm}{0.1pt}

\section{Вещественные числа и числовые множества}\label{SEC-chisla}

\subsection{Вещественные числа $\R$}\label{SUBSEC-R}

\paragraph{Измерения по эталону и вещественные числа.}\label{izmereniya-po-etalonu}

Понятие вещественного числа является естественным развитием идеи {\sl измерения
по эталону}, уходящей корнями в глубокую древность и опирающейся на следующие
два фундаментальных допущения.
 \bit{
 \item[1.]
Первое из них --- {\sl возможность сравнивать с эталоном} --- заключается в
том, что {\it среди значений физической величины можно выбирать эталонные}, то
есть такие, которыми можно (в целых числах) оценивать все остальные значения
этой величины. Например, имея эталон длины -- метр -- можно длину $x$ любого
физического тела оценивать в метрах (то есть найти, сколько метров в этой длине
укладывается, а сколько -- нет):
$$
k\le x<k+1.
$$
При этом важно, что организовать методы измерения можно так, чтобы результаты
не зависели от попыток.
 \item[2.]
Второе же --- {\sl возможность измельчения эталона} --- предполагает, что {\it
любой эталон можно делить на равные части (точнее, на произвольное фиксированное число равных частей), которые также можно принимать за
эталоны}. Из этого принципа следует, что мы можем поделить наш исходный эталон
длины, скажем, на 10 равных частей, принять их длины за новые эталоны, и после
этого поглядеть, не войдут ли в оставшуюся часть длины $x-k$ (куда уже не
входит наш исходный эталон) новые, меньшие эталоны:
$$
\frac{l}{10}\le x-k<\frac{l+1}{10}.
$$
Таким образом, мы получим уже оценку в десятых долях эталона (то есть можно
будет сказать, сколько десятых долей эталона укладывается в величине, которую
мы измеряем):
$$
k+\frac{l}{10}\le x<k+\frac{l+1}{10}
$$
 }\eit
Поскольку (согласно второму допущению) измельчать эталон можно сколь угодно
долго, мы получаем, что физическую величину $x$ можно оценивать сверху и снизу
рациональными числами (количеством долей эталона)
$$
\frac{m}{n}\le x<\frac{m+1}{n},
$$
причем с любой точностью (то есть так, чтобы числа $\frac{m}{n}$ и
$\frac{m+1}{n}$ отличались друг от друга как угодно мало). В этом состоит
применяемый в естествознании {\it алгоритм измерения физической величины по
эталону}, и интерпретировать его удобнее всего, как приближение к некоему
идеальному значению, изначально имеющемуся у данной измеряемой величины $x$.

Разумеется, не всякая физическая величина допускает измерение по эталону так,
как мы это описали. Например, у векторных величин (таких, как скорость или
сила) оценивать при измерении необходимо не одну, а сразу несколько
характеристик (например, проекции данной величины на заранее выбранные
координатные оси), а у целочисленных величин (таких, как количество однотипных
предметов в данном объеме пространства) допущение о делимости эталонов не будет
справедливым. Более того, развитие физики показало, что даже для величин,
измерение которых традиционно организуется именно сравнением с эталоном (таких
как длина, масса или объем), деление эталона тоже не может быть бесконечным,
потому что, когда Вы доходите до планковских размеров, дальнейшее деление
утрачивает смысл (то есть попросту становится непонятно, что значит делить
дальше).

Тем не менее, вера в эти принципы (или, точнее сказать, невозможность на данном
этапе развития науки заменить их чем-то более эффективным) приводит к тому, что
физическая величина представляется такой, что ее можно с любой точностью
оценивать рациональными числами в эталонных единицах измерения. Само по себе
это еще не означает, что значения физической величины должны непременно
описываться вещественными числами, какими их изображают в математическом
анализе (дело в том, что в математике имеются альтернативные теории, например,
<<нестандартный анализ>>, где числа описываются иначе). Но при дополнительных
предположениях, что величины могут быть отрицательными, и что в качестве
эталона можно выбирать любое ненулевое значение физической величины,
математическая теория, формализующая идею измерения по эталону, оказывается
единственной, и ее аксиоматика выглядит следующим образом.

\paragraph{Аксиомы теории вещественных чисел.}

Итак, про вещественные числа говорится, что они образуют множество,
обозначаемое символом $\R$, и удовлетворяющее следующим двум группам свойств I
и II.

\bigskip

\centerline{\bf I. Арифметические операции над вещественными
числами}\label{AXIOMS}\nobreak
Для любых двух вещественных чисел $a$ и $b$ определены, и притом
единственным образом, два числа $a + b$ ({\it сумма}) и $a\cdot b$
({\it произведение}), причем выполняются следующие правила.
\nobreak
 \bit{
\item[A1.] {\bf Коммутативность сложения:}
для любых чисел $a$ и $b$ выполняется равенство
 \beq\label{a+b=b+a}
a+b=b+a
 \eeq
\item[A2.] {\bf Ассоциативность сложения:} для любых чисел $a$, $b$ и $c$
выполняется равенство
 \beq\label{(a+b)+c=a+(b+c)}
(a+b)+c=a+(b+c)
 \eeq
\item[A3.] {\bf Коммутативность умножения:} для любых чисел $a$ и $b$
выполняется равенство
 \beq\label{a-cdot-b=b-cdot-a}
a\cdot b=b\cdot a
 \eeq
\item[A4.] {\bf Ассоциативность умножения:} для любых чисел $a$ и $b$
выполняется равенство
\beq\label{(a-cdot-b)-cdot-c=a-cdot-(b-cdot-c)}
(a\cdot b)\cdot c=a\cdot (b\cdot c)
\eeq
\item[A5.] {\bf Дистрибутивность:} для любых чисел $a,b,c$
выполняется равенство
 \beq\label{(a+b)-cdot-c=a-cdot-c+b-cdot-c}
 (a+b)\cdot c=a\cdot c + b\cdot c
 \eeq
\item[A6.] {\bf Существование нуля:} существует число 0 такое, что
 \beq\label{a+0=a}
a+0=a
 \eeq
\item[A7.] {\bf Существование противоположного числа:} для любого числа $a$
существует такое число $-a$ (называемое {\it противоположным числом} к числу
$a$) что
 \beq\label{a+(-a)=0}
a+(-a)=0
 \eeq
\item[A8.] {\bf Существование единицы:} существует число $1\ne 0$ такое, что
для любого числа $a\ne 0$ справедливо равенство
 \beq\label{a-cdot-1=a}
a\cdot 1=a
 \eeq
\item[A9.] {\bf Существование обратного числа:} для любого числа $a\ne 0$
существует такое число $a^{-1}$ (называемое {\it обратным числом} к числу $a$),
что
 \beq\label{a-cdot-a^(-1)=1}
a\cdot a^{-1}=1
 \eeq
 }\eit

\bigskip

\centerline{\bf  II. Сравнение вещественных чисел}

Для любых вещественных чисел $a$ и $b$ имеет смысл высказывание <<$a$ меньше,
либо равно $b$>>, записываемое $a\le b$. Оно может быть верно или неверно, но
при этом выполняются следующие правила.

 \bit{
\item[A10.] {\bf Рефлексивность:} для любого числа $a$ справедливо
$a\le a$.

\item[A11.] {\bf Транзитивность:} если $a\le b$ и $b\le c$, то $a\le
c$.

\item[A12.] {\bf Антисимметричность:}
если $a\le b$ и $b\le a$, то $a=b$.

\item[A13.]
{\bf Линейность:} для любых чисел $a$ и $b$ верно одно из двух:
 \biti{
\item[--] либо $a\le b$;
\item[--] либо $b\le a$.
 }\eiti

\item[A14.] {\bf Монотонность сложения:}
если $a\le b$, то для любого $c$ справедливо $a+c\le b+c$.

\item[A15.] {\bf Сохранение знака $+$ при умножении:}
если $0\le a$ и $0\le b$, то $0\le a\cdot b$.

\item[A16.] {\bf Непрерывность вещественной прямой:}
пусть $X$ и $Y$ -- два множества вещественных чисел, причем для
любых чисел $x\in X$ и $y\in Y$ выполняется неравенство $x\le y$;
тогда существует хотя бы одно число $c$ такое, что для любых $x\in
X$ и $y\in Y$ выполняется неравенство $x\le c\le y$.

 }\eit

\bigskip

Из аксиом A1 - A16 следуют все остальные свойства вещественных
чисел (и, в частности, все теоремы математического анализа).

Помимо отношения <<меньше, либо равно>> $a\le b$ на множестве вещественных
чисел $\R$ вводятся еще три отношения, а именно,
 \bit{
\item[--] <<меньше>> $a<b$, означающее, что $a\le b$ и одновременно $a\ne b$;
\item[--] <<больше, либо равно>> $a\ge b$ означающее, что $b\le a$; \item[--]
<<больше>> $a>b$, означающее, что $b\le a$ и одновременно $b\ne a$.
 }\eit

Из аксиомы линейности A13 следует, что множество вещественных чисел $\R$ удобно
изображать в виде прямой, а сами вещественные числа -- в виде точек на ней.
Именно так и делают математики:

\vglue70pt \noindent

\noindent\rule{160mm}{0.1pt}
\begin{multicols}{2}

\begin{ex}\label{EX-!0}
Покажем, что из аксиом A6 (о существовании нуля) и A1 (коммутативности)
следует, что нуль единственен. Действительно, если $0$ и $\tilde{0}$ -- два
разных объекта со свойствами нуля
 \beq\label{a+0=a}
\forall a\quad a+0=a,
 \eeq
 \beq\label{a+tilde-0=a}
\forall a\quad a+\tilde{0}=a
 \eeq
то мы получили бы
$$
0=\eqref{a+0=a}=\tilde{0}+0=(A1)=0+\tilde{0}=\eqref{a+tilde-0=a}=\tilde{0}
$$
(в скобках мы указываем формулы или аксиомы, которые в данный момент используем
в своих рассуждениях).
\end{ex}

\bex\label{ER-0-1-1} Покажем, что для данного числа $a$ противоположное число
$-a$, существование которого постулируется в аксиоме A7, тоже должно быть
единственно. Действительно, если бы существовали два числа $x$ и $y$ со
свойствами $a+x=0$ и $a+y=0$, то мы получили бы
 \begin{multline*}
x=(A6)=x+0=x+(a+y)=(A2)=\\=(x+a)+y=(A1)=(a+x)+y=0+y=\\=(A1)=y+0=y
 \end{multline*}
\eex

\ber Докажите равенство:
 \beq\label{-0=0}
-0=0
 \eeq
\eer

\bex Покажем, что справедлива следующая эквивалентность, означающая, что
уравнения вида $a+x=b$ однозначно разрешимы в $\R$:
 \beq\label{a+x=b}
a+x=b\quad\Longleftrightarrow\quad  x=b+(-a)
 \eeq
 \eex
 \bpr
$$
a+x=b
$$
$$
\Updownarrow
$$
$$
(a+x)+(-a)=b+(-a)
$$
$$
\Updownarrow
$$
$$
(x+a)+(-a)=b+(-a)
$$
$$
\Updownarrow
$$
$$
x+(a+(-a))=b+(-a)
$$
$$
\Updownarrow
$$
$$
x+0=b+(-a)
$$
$$
\Updownarrow
$$
$$
x=b+(-a)
$$
\epr

\bex Справедливы следующие тождества:
 \begin{align}
& 0\cdot a=0 \label{0a=0} \\
& -(-a)=a \label{-(-a)=a} \\
& (-1)\cdot a=-a \label{(-1)a=-a} \\
& (-a)\cdot(-a)=a\cdot a \label{(-a)(-a)=aa}
 \end{align}
 \eex
\bpr 1. Докажем \eqref{0a=0}:
$$
0\cdot a=(0+0)\cdot a=0\cdot a+0\cdot a
$$
$$
\phantom{\text{\scriptsize \eqref{a+x=b}}}\quad\Downarrow\quad\text{\scriptsize
\eqref{a+x=b}}
$$
$$
0\cdot a+(-0\cdot a)=0\cdot a
$$
$$
\phantom{\text{\scriptsize (A7)}}\quad\Downarrow\quad\text{\scriptsize (A7)}
$$
$$
0=0\cdot a
$$

2. Для \eqref{-(-a)=a} доказательство выглядит так:
$$
a+(-a)=0
$$
$$
\Downarrow
$$
$$
(-a)+a=0
$$
$$
\Downarrow
$$
$$
\text{$a$ -- противоположное число для $-a$}
$$
$$
\phantom{{\scriptsize \begin{pmatrix}\text{у $-a$ не может}\\ \text{быть
двух}\\ \text{противоположных}\\
\text{чисел}\end{pmatrix}}}\quad\Downarrow\quad{\scriptsize
\begin{pmatrix}\text{у $-a$ не может}\\ \text{быть двух}\\
\text{противоположных}\\ \text{чисел}\end{pmatrix}}
$$
$$
a=-(-a)
$$

3. Для \eqref{(-1)a=-a} доказательство такое:
$$
(A7)
$$
$$
\Downarrow
$$
$$
1+(-1)=0
$$
$$
\Downarrow
$$
$$
(1+(-1))\cdot a=0\cdot a
$$
$$
\phantom{\text{\scriptsize
(A5),\eqref{0a=0}}}\quad\Downarrow\quad\text{\scriptsize (A5),\eqref{0a=0}}
$$
$$
1\cdot a+(-1)\cdot a=0
$$
$$
\Downarrow
$$
$$
a+(-1)\cdot a=0
$$
$$
\Downarrow
$$
$$
\text{$(-1)\cdot a$ -- противоположное число для $a$}
$$
$$
\Downarrow
$$
$$
(-1)\cdot a=-a
$$

4. Докажем \eqref{(-a)(-a)=aa}:
 \begin{multline*}
(-a)\cdot(-a)=\eqref{(-1)a=-a}=((-1)\cdot a)\cdot(-a)=\\=(A3)=(a\cdot
(-1))\cdot(-a)=(A4)=\\=a\cdot(
(-1)\cdot(-a))=\eqref{(-1)a=-a}=\\=a\cdot(-(-a))=\eqref{-(-a)=a}=a\cdot a
 \end{multline*}
 \epr

\begin{er}\label{ER-!1}
Докажите, по аналогии с примером \ref{EX-!0}, что единица 1, существование
которой постулируется в аксиоме A8, тоже единственна.
\end{er}

\ber Докажите равенство:
 \beq\label{1^(-1)=1}
1^{-1}=1
 \eeq
\eer

\begin{er}\label{ER-0-1-2}
По аналогии с примером \ref{ER-0-1-1} докажите, что для данного числа $a\ne 0$
обратное число $a^{-1}$ (о котором идет речь в аксиоме A9), должно быть
единственно (то есть, не может существовать двух чисел $x$ и $y$ со свойствами
$a\cdot x=1$ и $a\cdot y=1$).
\end{er}

\ber Докажите тождество:
 \beq\label{(a^(-1))^(-1)=a}
\left(a^{-1}\right)^{-1}=a
 \eeq
\eer

 \bex Справедливо тождество
  \beq\label{(-a)^(-1)=-a^(-1)}
(-a)^{-1}=-a^{-1}
  \eeq
 \eex
 \bpr
$$
(-a^{-1})\cdot(-a)=\eqref{(-a)(-a)=aa}=a^{-1}\cdot a=1
$$
$$
\phantom{{\scriptsize \begin{pmatrix}\text{в силу}\\ \text{упражнения \ref{ER-0-1-2},} \\ \text{у $-a$ не может}\\
\text{быть
двух}\\ \text{обратных}\\
\text{чисел}\end{pmatrix}}}\quad\Downarrow\quad{\scriptsize
\begin{pmatrix}\text{в силу}\\ \text{упражнения \ref{ER-0-1-2},} \\ \text{у $-a$ не может}\\ \text{быть двух}\\
\text{обратных}\\ \text{чисел}\end{pmatrix}}
$$
$$
(-a)^{-1}=-a^{-1}
$$
 \epr

\begin{er}\label{ER-0-1-3}
Докажите следующие утверждения:
 \begin{align}
&\label{x<y->x+a<y+a} x<y\quad\Longrightarrow\quad\forall a\quad x+a<y+a \\
&\label{xy=0=>x=0-vee-y=0} x\cdot y=0\quad\Longrightarrow\quad
x=0\quad\V\quad y=0 \\
&\label{x>0,y>0->x.y>0} x>0\quad\&\quad y>0\quad\Longrightarrow\quad x\cdot y>0
\\
&\label{x<0,y<0->x.y>0} x<0\quad\&\quad y<0\quad\Longrightarrow\quad x\cdot y>0
\\
&\label{x<0,y>0->x.y<0} x<0\quad\&\quad y>0\quad\Longrightarrow\quad x\cdot y<0
\\
&\label{1>0} 1>0 \\
&\label{-1<0} -1<0 \\
&\label{x<y=>-x>-y} x<y\quad\Longrightarrow\quad -x>-y \\
&\label{x>0=>x^(-1)>0} x>0\quad\Longrightarrow\quad x^{-1}>0 \\
&\label{x<0=>x^(-1)<0} x<0\quad\Longrightarrow\quad x^{-1}<0 \\
&\label{0<x<y=>x^(-1)>y^(-1)} 0<x<y\quad\Longrightarrow\quad x^{-1}>y^{-1} \\
&\label{x<y<0=>x^(-1)>y^(-1)} x<y<0\quad\Longrightarrow\quad x^{-1}>y^{-1} \\
&\label{x<y-a>0=>ax<ay} x<y\quad \&\quad a>0\quad \Longrightarrow\quad a\cdot
x<a\cdot y \\
&\label{x<y-a<0=>ax>ay} x<y\quad \&\quad a<0\quad \Longrightarrow\quad a\cdot
x>a\cdot y \\
&\label{a>1<=>0<a<1} a>1\quad\Longleftrightarrow\quad 0<a<1 \\
&\label{x>0-a>1=>ax>x} x>0\quad \&\quad a>1\quad \Longrightarrow\quad a\cdot
x>x \\
&\label{x>0-0<a<1=>ax<x} x>0\quad \&\quad 0<a<1\quad \Longrightarrow\quad
a\cdot x<x
 \end{align}

\end{er}

\end{multicols}\noindent\rule[10pt]{160mm}{0.1pt}

\paragraph{Собственные обозначения для чисел и числовые равенства.}
\label{SUBSEC-sobstv-obozn-chisla}

Буквы латинского алфавита -- a,b,c -- которыми мы обозначали числа в аксиомах
A1 -- A16, являются {\it переменными}. Это означает, что с их помощью
описываются общие свойства всех объектов данной теории (или части этих
объектов, если указана область изменения переменных).

Но помимо переменных, играющих в математическом языке роль <<имен
нарицательных>>, в нем используются также <<имена собственные>>, нужные для
обозначения <<конкретных объектов>> (в противоположность <<объектам вообще>>
или <<объектам из заданного множества>>). Для таких обозначений используется
обычно термин <<константа>>, но в математике его смысл, как правило, несколько
шире, чем тот, что нам нужен, поэтому мы будем использовать другой --- имена
собственные в математических формулах мы будем называть {\it собственными
обозначениями}.

Примерами собственных обозначений являются символы нуля 0 и единицы 1,
определенные аксиомами A6 и A8. Это следует из примера \eqref{EX-!0} и
упражнения \eqref{ER-!1}, где мы отмечали, что нуль и единица определены
однозначно. Неискушенному читателю может показаться неожиданным, но все
остальные собственные обозначения для чисел строятся из 0 и 1 с помощью
операций и отношений, описываемых аксиомами A1 -- A16 (с той только оговоркой,
что не всегда нужное тебе число, как, например, число $\pi$, можно определить
сразу -- иногда для этого нужна предварительная работа в виде системы
обозначений, определений и теорем).

\noindent\rule{160mm}{0.1pt}
\begin{multicols}{2}

\paragraph*{Числа $2,...,10$.}\label{chisla-2-10}

Например, число 2 можно сразу определить равенством
 \beq\label{DEF:2}
2:=1+1
 \eeq
(символ $:=$ здесь означает, что это соотношение следует понимать как
определение, то есть как декларацию, согласно которой символ слева от $:=$
будет использоваться в тексте всюду для обозначения объекта, лежащего справа от
этого знака).

Точно также, числа 3, 4, 5 можно определить равенствами
 \beq\label{0.2.1}
 \begin{split}
3&:=1+(1+1), \\ 4&:=1+(1+(1+1)), \\  5&:=1+(1+(1+(1+1)))
 \end{split}
 \eeq
Важное наблюдение здесь состоит в том, что одно и то же число можно определять
по-разному. Например, 3, 4, 5 можно иначе определить так:
 \beq\label{0.2.2}
  \begin{split}
 3&:=(1+1)+1, \\
 4&:=((1+1)+1)+1, \\
 5&:=(((1+1)+1)+1)+1
 \end{split}
 \eeq
Строго говоря, это будут другие определения, потому что записи в \eqref{0.2.1}
и \eqref{0.2.2} не одинаковы (дело в том, что, с точки зрения языка,
последовательность символов $1+(1+1)$ --- совсем не то же самое, что
последовательность символов $(1+1)+1$).

Тем не менее, определения \eqref{0.2.1} и \eqref{0.2.2} эквивалентны, потому
что из аксиом A1 -- A16 можно логическими средствами вывести, что числа,
определяемые формулами \eqref{0.2.1} --- те же самые, что числа, определяемые
формулами \eqref{0.2.2}. Действительно, эквивалентность определений для числа
$3$ получается подстановкой в тождество из аксиомы A2 вместо $a$, $b$, $c$,
символа $1$:
 \begin{multline*}
(1+1)+1=(a+b)+c\;\Bigg|_{\scriptsize\begin{matrix}a=1 \\ b=1 \\ c=1
\end{matrix}}=(A2)=\\
=a+(b+c)\;\Bigg|_{\scriptsize\begin{matrix}a=1 \\ b=1 \\ c=1
\end{matrix}} =1+(1+1)
 \end{multline*}
(символ $|$ означает подстановку, а фрагмент $=(A2)=$ есть ссылка на аксиому
A2, применяемую в этот момент). Эту запись можно сократить так:
$$
(1+1)+1=(A2)=1+(1+1)
$$
--- здесь читателю предлагается самому догадаться, какой должна быть
подстановка, чтобы можно было применить аксиому A2. Читатель теперь может
самостоятельно расшифровать сокращенную запись доказательства эквивалентности
определений числа 4:
 \begin{multline*}
((1+1)+1)+1=(A2)=(1+1)+(1+1)=\\=(A2)=1+(1+(1+1))
 \end{multline*}
Что же касается числа 5, то нам теперь достаточно сказать, что эквивалентность
определений для него доказывается по аналогии.

Итак, одни и те же собственные обозначения можно задавать формально по-разному,
нужно только следить, чтобы выписываемые тобой определения для одного и того же
символа были эквивалентны. Естественный (и самый распространенный) способ
вводить новые обозначения -- так называемый последовательный, в котором при
определении нового обозначения используются обозначения, введенные ранее.
Например, числа 3, 4, 5, 6, 7, 8, 9, 10 можно последовательно определить так:
 \begin{align*}
3&:=2+1, & 4&:=3+1, \\
5&:=4+1, & 6&:=5+1, \\
7&:=6+1, & 8&:=7+1, \\
9&:=8+1, & 10&:=9+1
 \end{align*}

\paragraph*{Числа вида $-a$.}

Когда определены числа $1$,...,$10$, их противоположные числа $-1$,...,$-10$
определять уже не нужно, поскольку они автоматически определяются аксиомой A7 о
существовании противоположного числа (дело в том, что, как мы отмечали в
упражнении \ref{ER-0-1-1}, из аксиомы A7 и других аксиом следует, что
противоположное число определяется единственным образом).

Вообще, всякий раз, когда мы дали собственное обозначение $a$ какому-нибудь
числу, аксиома A7 определяет его противоположное число и дает ему собственное
обозначение $-a$.

\paragraph*{Числа вида $a^{-1}$.}

Точно также, всякий раз, когда мы даем собственное обозначение $a$
какому-нибудь числу ($a\ne 0$), аксиома A9 дает собственное обозначение
обратному числу $a^{-1}$. В частности, обратные числа к $2,...,10$ имеют
собственные обозначения $2^{-1}$,...,$10^{-1}$.

\paragraph*{Дроби} определяются равенством
$$
\frac{a}{b}=a\cdot b^{-1}=b^{-1}\cdot a, \qquad (b\ne 0)
$$
Например,
$$
\frac{1}{2}:=2^{-1},\qquad \frac{2}{3}:=2\cdot 3^{-1},\qquad
\frac{3}{5}:=3\cdot 5^{-1}
$$

\ber Докажите тождества:
 \begin{align}
\frac{a}{b}+\frac{c}{d}&=\frac{a\cdot d+b\cdot c}{b\cdot d}
\label{(a/b)+(c/d)=(ad+bc)/(bd)} \\
\frac{a}{b}-\frac{c}{d}&=\frac{a\cdot d-b\cdot c}{b\cdot d}
\label{(a/b)-(c/d)=(ad-bc)/(bd)} \\
\frac{a}{b}\cdot\frac{c}{d}&=\frac{a\cdot c}{b\cdot d}
\label{(a/b)(c/d)=(ac)/(bd)}
 \end{align}
\eer

\end{multicols}\noindent\rule[10pt]{160mm}{0.1pt}

\subsection{Числовые множества}

Как мы уже говорили, числа также образуют множества. При этом, поскольку мы
строим теорию вещественных чисел, опираясь на теорию множеств, нам не нужно
давать строгое определение этому понятию. Здесь мы приведем простейшие примеры
числовых множеств и обсудим некоторые их свойства, которые понадобятся в
дальнейшем.

\noindent\rule{160mm}{0.1pt}
\begin{multicols}{2}

\bex {\bf Множества, определяемые явным перечислением своих элементов.} Это
самый простой пример числовых множеств. Например, запись
$$
\left\{-1;\frac{1}{2};3\right\}
$$
обозначает множество, состоящее из трех элементов: $-1$, $\frac{1}{2}$ и $3$.
\eex

\bex {\bf Числовые множества, определяемые отношениями порядка.} Отношения
$\le, \, <, \, \ge, \, >$ определяют на прямой $\R$ множества, известные как
интервал, отрезок и полуинтервалы:
 \begin{align*}
(a,b) & =\big\{ x\in \R : \quad a<x<b \big\} \\
[a,b] & =\big\{ x\in \R : \quad a\le x\le b \big\} \\
(a,b] & =\big\{ x\in \R : \quad a<x\le b \big\} \\
[a,b) & =\big\{ x\in \R : \quad a\le x<b \big\} \\
(a,+\infty) & =\big\{ x\in \R : \quad a<x \big\} \\
[a,+\infty) & =\big\{ x\in \R : \quad a\le x \big\} \\
(-\infty,b) & =\big\{ x\in \R : \quad x<b \big\} \\
(-\infty,b] & =\big\{ x\in \R : \quad x\le b \big\}
 \end{align*}
Точнее, в этом списке множество $(a,b)$ называется {\it интервалом},
$[a,b]$ -- {\it отрезком}, $(a,b]$ -- {\it полуинтервалом с замкнутым
правым концом}, $[a,b)$ -- полуинтервалом с замкнутым левым концом,
$(a,+\infty)$ -- открытым полуинтервалом с бесконечным правым концом,
$[a,+\infty)$ -- замкнутым полуинтервалом с бесконечным правым
концом, и так далее. \eex

\end{multicols}\noindent\rule[10pt]{160mm}{0.1pt}

\paragraph{Алгебраические операции над числовыми множествами и неравенства между ними.}

 \bit{
\item[$\bullet$] Если $X\subseteq\R$ --- числовое множество и $a\in\R$ --
число, то
 \bit{
\item[---] символом $-X$ обозначается множество, состоящее из чисел,
противоположных числам из $X$:
 \beq\label{protivop-mnozhestvo}
-X:=\{y\in\R:\;\exists x\in X\quad -x=y\}
 \eeq

\item[---] символами $X+a$ и $X-a$ обозначаются числовые множества, получаемые
сдвигом $X$ на $a$ вправо и влево:
 \begin{align}
& X+a=\{y\in\R:\;\exists x\in X\quad x+a=y\},\label{sdvigi-mnozhestv+} \\
& X-a=\{y\in\R:\;\exists x\in X\quad x-a=y\} \label{sdvigi-mnozhestv-}
 \end{align}

\item[---] символами $a\cdot X$ и $X\cdot a$ обозначается числовое множество,
получаемое умножением $X$ на $a$:
 \begin{align}
& a\cdot X=X\cdot a=\{y\in\R:\;\exists x\in X\quad a\cdot
x=y\},\label{umnozh-mnozhestva-na-chislo}
 \end{align}

\item[---] неравенство
 \beq\label{DF:a-le-X}
a\le X,
 \eeq
означает, что любой элемент множества $X$ не меньше числа $a$
 $$
\forall x\in X\quad a\le x;
 $$
точно так же формулы
$$
a< X,\quad X\le a,\quad X<a
$$
являются сокращенными записями следующих утверждений соответственно:
 $$
\forall x\in X\quad a<x,\qquad\forall x\in X\quad x\le a,\qquad\forall x\in
X\quad x<a
 $$

 }\eit

\item[$\bullet$] Если $X$ и $Y$ --- два числовых множества, то
 \bit{
\item[---] символами $X+Y$, $X-Y$, $X\cdot Y$, $\frac{X}{Y}$ обозначаются
множества, определяемые формулами
 \begin{align}
& X+Y=\{z\in\R:\;\exists x\in X\quad \exists y\in Y\quad x+y=z\}\label{DEF:X+Y}
 \\
& X-Y=\{z\in\R:\;\exists x\in X\quad \exists y\in Y\quad x-y=z\}\label{DEF:X-Y}\\
& X\cdot Y=\{z\in\R:\;\exists x\in X\quad \exists y\in Y\quad x\cdot
y=z\}\label{DEF:X-cdot-Y}
\\
& \frac{X}{Y}=\left\{z\in\R:\;\exists x\in X\quad \exists y\in Y\quad
\frac{x}{y}=z\right\}\label{DEF:X/Y}
 \end{align}

\item[---] неравенство
$$
X\le Y
$$
по определению означает, что любой элемент множества $X$ не превосходит любого
элемента множества $Y$:
$$
\forall x\in X\quad \forall y\in Y\qquad x\le y
$$
а неравенство
$$
X< Y
$$
-- что любой элемент множества $X$ меньше любого элемента множества $Y$:
$$
\forall x\in X\quad \forall y\in Y\qquad x<y
$$
 }\eit
 }\eit

\noindent\rule{160mm}{0.1pt}
\begin{multicols}{2}

\ber Проверьте справедливость равенств:
 \begin{align*}
& \left\{-1;\frac{1}{2};3\right\}+4=\left\{3;\frac{9}{2};7\right\} \\
& \left\{-1;\frac{1}{2};3\right\}-4=\left\{-5;-\frac{7}{2};-1\right\} \\
& (a;b)+c=(a+c;b+c) \\
& [a;b]+c=[a+c;b+c] \\
& (a;+\infty)+c=(a+c;+\infty) \\
& (-\infty;a)+c=(-\infty;a+c)
 \end{align*}
\eer

\end{multicols}\noindent\rule[10pt]{160mm}{0.1pt}

\paragraph{Минимум и максимум числового множества.}

 \bit{
\item[$\bullet$] Пусть $X$ --- числовое множество. Число $a$ называется {\it
минимумом} множества $X$, если $a$ принадлежит $X$, и все остальные элементы
$X$ не меньше $a$:
$$
a\in X\quad\&\quad a\le X \qquad(\text{то есть $a\in X$ и $\forall x\in X$
$a\le x$})
$$
В таких случаях пишут
$$
a=\min X
$$

\item[$\bullet$] Аналогично, число $b$ называется {\it максимумом} множества
$X$, если $b$ принадлежит $X$, и все остальные элементы $X$ не больше $b$:
$$
b\in X\quad\&\quad X\le b \qquad(\text{то есть $b\in X$ и $\forall x\in X$
$x\le b$})
$$
и в таких случаях пишут
$$
b=\max X
$$
 }\eit

\noindent\rule{160mm}{0.1pt}
\begin{multicols}{2}

\bexs Нетрудно сообразить, что минимум и максимум отрезка $X=[a;b]$
совпадают с его концами:
$$
a=\min[a;b],\qquad b=\max[a;b]
$$
Но при этом у интервала $X=(a;b)$ нет ни минимума, ни максимума:
 \beq\label{nexists-min(a;b)}
\nexists\min(a;b),\qquad \nexists\max(a;b)
 \eeq
<<Наглядное>> объяснение этому выглядит так: минимумом множества $(a;b)$ может
быть только число $a$, но оно не принадлежит $(a;b)$; и то же самое с
максимумом.

Разумеется, такие рассуждения нельзя считать строгим доказательством. Чтобы
провести аккуратные выкладки нам понадобится следующее утверждение:
 \beq\label{srednee-arifm}
 x<y\quad\Longrightarrow\quad x<\frac{x+y}{2}<y.
 \eeq
Для его доказательства в свою очередь понадобится вот что:
 \beq\label{0<a->0<a/2<a}
 0<a\quad\Longrightarrow\quad 0<\frac{a}{2}<a.
 \eeq
Это доказывается так: с одной стороны, из \eqref{x>0,y>0->x.y>0} получаем
 $$
0<a\quad\Longrightarrow\quad 0<\frac{a}{2}
 $$
отсюда, в силу \eqref{x<y->x+a<y+a} следует
 $$
\frac{a}{2}=0+\frac{a}{2}<\frac{a}{2}+\frac{a}{2}=a
 $$
и вместе это дает \eqref{0<a->0<a/2<a}.

Теперь из \eqref{0<a->0<a/2<a} следует \eqref{srednee-arifm}:
 \begin{align*}
 x &<y \\
   &\Downarrow \\
 0 &<y-x \\
   &\Downarrow \quad {\scriptsize \eqref{0<a->0<a/2<a}}\\
 0 <\frac{y-x}{2} &<y-x
 \end{align*}

Наконец, доказываем \eqref{nexists-min(a;b)}:
 \begin{align*}
 c&=\min(a;b) \\
   &\Downarrow \\
 c &\in (a;b) \\
   &\Downarrow \\
 a<& c<b \\
   &\Downarrow \quad {\scriptsize \eqref{srednee-arifm}}\\
 a<&\frac{a+c}{2}<c \\
   &\Downarrow \\
 \exists z=\frac{a+c}{2}&\in (a;b)\quad z<c \\
   &\Downarrow \\
 c&\ne\min(a;b)
 \end{align*}

И то же самое с максимумом.

\eexs

\end{multicols}\noindent\rule[10pt]{160mm}{0.1pt}

\bigskip

\centerline{\bf Свойства $\min$ и $\max$:}

\bit{\it

\item[$1^\circ$]\label{monotonnost-min-max-X} {\bf Монотонность:} если
множества $X$ и $Y$ имеют минимум (максимум), то
 \beq\label{X-subseteq-Y=>min(X)-ge-min(Y)}
X\subseteq Y\quad\Longrightarrow\quad \min X\ge \min Y\quad \Big(\max X\le\max
Y\Big)
 \eeq

\item[$2^\circ$]\label{min-X+a=min-X+a} {\bf Инвариантность относительно
сдвига:} если множество $X$ имеет минимум (максимум), то любой его сдвиг $X+a$
тоже имеет минимум (максимум), причем
 \beq\label{min(X+a)=min(X)+a}
\min (X+a)=(\min X)+a \qquad \Big(\max (X+a)=(\max X)+a\Big)
 \eeq

\item[$3^\circ$]\label{min-X=-max_X} {\bf Двойственность:} если множество $X$
имеет минимум (максимум), то противоположное множество $-X$ имеет максимум
(минимум), и
 \beq\label{min(-X)=-min(X)}
\max (-X)=-\min X \qquad \Big(\min (-X)=-\max X\Big)
 \eeq
 }\eit

\bpr Мы докажем первую половину каждого утверждения (вторая доказывается по
аналогии).

1. Пусть $X\subseteq Y$. Тогда
$$
a=\min X\quad\Longrightarrow\quad a\in X\subseteq Y \quad\Longrightarrow\quad
a\ge\min Y
$$

2.
 \begin{multline*}
p=\min X\quad\Longrightarrow\quad p\in X\quad\&\quad p\le X
\quad\Longrightarrow \\ \Longrightarrow\quad p+a\in X+a\quad\&\quad p+a\le X+a
\quad\Longrightarrow\quad p+a=\min(X+a)
 \end{multline*}

3.
 $$
p=\min X\quad\Longrightarrow\quad p\in X\quad\&\quad p\le X
\quad\Longrightarrow\quad -p\in -X\quad\&\quad -p\ge -X
\quad\Longrightarrow\quad -p=\max(-X)
 $$
\epr

\noindent\rule{160mm}{0.1pt}\begin{multicols}{2}

\paragraph{Минимум и максимум двух чисел.}

 \biter{
\item[$\bullet$] Из аксиомы A13 следует, что формулы
 \begin{align}
a\wedge b:=\begin{cases}a,& a\le b \\ b,& a>b\end{cases} \label{DEF:a-wedge-b}
\\
 a\vee b:=\begin{cases}b,& a\le b \\ a,& a>b\end{cases} \label{DEF:a-vee-b}
 \end{align}
 корректно определяют две операции над вещественными числами.
 }\eiter

\btm\label{TH:exists-min-a,b} Для любых двух чисел $a,b\in\R$ множество
$\{a,b\}$ обладает минимумом и максимумом, и справедливы формулы:
 \begin{align}
&\label{a-wedge-b=min-a-b} a\wedge b=\min\{a,b\} \\
&\label{a-vee-b=max-a-b} a\vee b=\max\{a,b\}
 \end{align}
 \etm
\bpr По аксиоме A13, выполняется одно из двух: либо $a\le b$, либо $a>b$.

В первом случае (когда $a\le b$) мы получим:
$$
a\wedge b=a=\min \{a,b\}.
$$
Здесь первое равенство выполняется по определению $a\wedge b$, а второе --
потому что $a\in\{a,b\}$ и $a\le\{a,b\}$. И точно так же
$$
a\vee b=b=\max\{a,b\},
$$
и первое равенство выполняется по определению $a\vee b$, а второе -- потому что
$b\in\{a,b\}$ и $\{a,b\}\le b$.

А во втором случае (когда $a>b$) мы получаем
$$
a\wedge b=b=\min\{a,b\},
$$
где второе равенство выполняется потому что $b\in\{a,b\}$ и $\{a,b\}\ge b$. И,
с другой стороны,
$$
a\vee b=a=\max\{a,b\},
$$
где второе равенство верно потому что $a\in\{a,b\}$ и $a\ge\{a,b\}$.
 \epr

\bex Очевидно,
 \begin{align*}
&0\wedge 1=0,&& 0\vee 1=1 \\
&0\wedge(-1)=-1,&& 0\vee(-1)=0.
 \end{align*}
 \eex

\ber Докажите, что выполняются следующие тождества:
 \begin{align*}
a\wedge b&=b\wedge a \\
a\vee b&=b\vee a \\
(a\wedge b)\wedge c&=a\wedge (b\wedge c) \\
(a\vee b)\vee c&=a\vee (b\vee c) \\
(a\wedge b)\vee c&=(a\vee c)\wedge (b\vee c) \\
(a\vee b)\wedge c&=(a\wedge c)\vee (b\wedge c) \\
 \end{align*}
\eer

\ber Докажите, что не существует числа $b\in\R$ такого, что равенство
$$
a\wedge b=a
$$
выполнялось бы для всех $a\in\R$. \eer

\ber Докажите, что не существует числа $b\in\R$ такого, что равенство
$$
a\vee b=a
$$
выполнялось бы для всех $a\in\R$. \eer

\ber Проверьте, какие из следующих тождеств верны:
 \begin{align*}
(a\wedge b)+c&=(a+c)\wedge (b+c) \\
(a+b)\wedge c&=(a\wedge c)+(b\wedge c) \\
(a\vee b)+c&=(a+c)\vee (b+c) \\
(a+b)\vee c&=(a\vee c)+(b\vee c)
 \end{align*}
\eer

\bprop\label{PROP:max(A-u-B)=max(A)-v-max(B)} Если числовые множества $A$ и $B$
оба обладают минимумом (максимумом), то их объединение $A\cup B$ тоже обладает
минимумом (максимумом), причем выполняется равенство:
 \begin{align}
&\min (A\cup B)=(\min A)\wedge(\min B)\\
\Big(&\max (A\cup B)=(\max A)\vee (\max B)
\Big)\label{max(A-u-B)=max(A)-v-max(B)}
 \end{align}
\eprop

\end{multicols}\noindent\rule[10pt]{160mm}{0.1pt}

\paragraph{Нижняя и верхняя грани числового множества.}
 \bit{
\item[$\bullet$] Говорят, что числовое множество $X$ {\it ограничено снизу},
если существует число $s$ такое, что для всякого $x\in X$ выполняется
неравенство $s\le x$. Коротко это условие можно записать с помощью неравенств,
в которых сравниваются числа и множества (мы определяли такие неравенства на
с.\pageref{DF:a-le-X}):
$$
\exists s\in\R\quad s\le X.
$$
В таких случаях принято говорить также, что число $s$ {\it ограничивает
множество $X$ снизу}.

\item[$\bullet$] Аналогично, говорят, что числовое множество $X$ {\it
ограничено сверху}, если существует число $t$ такое, что для всякого $x\in X$
выполняется неравенство $x\le t$. Это можно записать так:
$$
\exists t\in\R\quad X\le t.
$$
В таких случаях принято говорить также, что число $d$ {\it ограничивает
множество $X$ сверху}.

\item[$\bullet$] Если множество $X$ ограничено и снизу, и сверху, то говорят,
что множество $X$ ограничено.

\item[$\bullet$] {\it Точной нижней гранью} (или {\it точной нижней границей})
$\inf X$ числового множества $X$ называется

 \bit{

\item[---] символ $-\infty$, если $X$ не ограничено снизу; в этом случае точная
нижняя грань не является числом: запись
$$
\inf X=-\infty
$$
просто означает, что множество $X$ не ограничено снизу;

\item[---] наибольшее из чисел $s$, ограничивающих $X$ снизу, если $X$
ограничено снизу:
 \beq\label{DF:inf-X}
\inf X=\max\{s\in\R:\quad s\le X\}
 \eeq
запись
$$
\inf X>-\infty
$$
принято использовать в случаях, когда нужно коротко дать понять, что множество
$X$ ограничено снизу.
 }\eit

\item[$\bullet$] {\it Точной верхней гранью} (или {\it точной верхней
границей}) $\sup X$ числового множества $X$ называется

 \bit{

\item[---] символ $+\infty$, если $X$ не ограничено сверху; опять же в этом
случае точная нижняя грань не является числом: запись
$$
\inf X=+\infty
$$
просто означает, что множество $X$ не ограничено сверху;

\item[---] наименьшее из чисел $t$, ограничивающих $X$ сверху, если $X$
ограничено сверху:
 \beq\label{DF:sup-X}
\sup X=\min\{t\in\R:\quad  X\le t\}
 \eeq
запись
$$
\sup X<+\infty
$$
принято использовать в случаях, когда множество $X$ ограничено сверху.
 }\eit
 }\eit

\noindent\rule{160mm}{0.1pt}
\begin{multicols}{2}

\begin{ex}
Пусть $X=(a,b)$ --- интервал на вещественной прямой $\R$.

%\includegraphics[0,0][0,0]{10.eps}
%\picture{0pt}{0pt}{10.pcx}
\vglue70pt \noindent

Легко видеть, что  любое число $t\ge b$ ограничивает множество $X$ сверху. А
наименьшее из всех таких чисел, то есть число $B=b$ будет точной верхней гранью
множества $X$. С другой стороны, любое число $s\le a$ ограничивает $X$ снизу. А
наибольшее из всех таких чисел, то есть число $A=a$ будет точной нижней гранью
множества $X$. Таким образом,
$$
\inf (a,b)= a \qquad \sup (a,b)=b
$$
Точно также для отрезка и полуинтервалов
$$
\inf [a,b]=\inf (a,b]=\inf [a,b)= a
$$
$$\sup [a,b]=\sup (a,b]=\sup [a,b)=b
$$
\end{ex}

\end{multicols}\noindent\rule[10pt]{160mm}{0.1pt}

\btm Если множество $X$ обладает минимумом (максимумом), то он совпадает с
точной нижней (верхней) гранью этого множества:
$$
\exists \min X\quad\Longrightarrow\quad \min X=\inf X \qquad \Big(\exists \max
X\quad\Longrightarrow\quad \max X=\inf X\Big)
$$
 \etm

\begin{tm}[\bf о точной границе]\label{sup-inf}
Если $X$ -- непустое ограниченное снизу множество, то оно имеет точную нижнюю
грань. Аналогично, если $X$ -- непустое ограниченное сверху множество, то оно
имеет точную верхнюю грань.
\end{tm}

\begin{proof} Рассмотрим случай, когда $X$ -- непустое ограниченное сверху
множество (случай, когда $X$ ограничено снизу рассматривается аналогично).
Обозначим буквой $Y$ множество всех чисел, ограничивающих $X$ сверху:
$$
Y=\{y\in\R:\; X\le y\}=\{y\in\R:\; \forall x\in X\quad x\le y\}
$$
(поскольку $X$ ограничено сверху, хотя бы одно такое число $y$ существует, и
значит, $Y$ -- непустое множество). Тогда
$$
\forall x\in X \quad \forall y\in Y \quad x\le y
$$
и по аксиоме непрерывности вещественной прямой A16, найдется такое число $c$,
что
$$
\forall x\in X\quad \forall y\in Y \quad x\le c\le y
$$
Первое неравенство здесь означает, что $c$ ограничивает $X$ сверху, то есть,
что $c\in Y$. А второе -- что $c$ есть наименьшее число, лежащее в $Y$. То есть
$c=\min Y=\min\{y\in\R:\; X\le y\}=\sup X$.
\end{proof}

\bigskip

\centerline{\bf Свойства $\inf$ и $\sup$:}

\bit{\it

\item[$1^\circ$]\label{monotonnost-inf-sup} {\bf Монотонность:} если $X$ и $Y$
ограничены снизу (сверху), то
$$
X\subseteq Y\quad\Longrightarrow\quad \inf X\ge \inf Y\quad \Big(\sup X\le\sup
Y\Big).
$$

\item[$2^\circ$]\label{inf-X+a=inf-X+a} {\bf Инвариантность относительно
сдвига:} если множество $X$ ограничено снизу (сверху), то любой его сдвиг $X+a$
тоже ограничен снизу (сверху), причем
 $$
\inf (X+a)=(\inf X)+a \qquad \Big(\sup (X+a)=(\sup X)+a\Big)
 $$

\item[$3^\circ$]\label{sup-X=-sup_X} {\bf Двойственность:} если множество $X$
ограничено снизу (сверху), то противоположное множество $-X$ ограничено сверху
(снизу), и
 $$
\sup (-X)=-\inf X \qquad \Big(\sup (-X)=-\inf X\Big)
 $$
}\eit

\bpr Здесь мы также докажем только первую половину каждого утверждения
(предложив читателю самостоятельно по аналогии доказать вторую).

1. В первом случае получаем цепочку:
\begin{align*}
X &\subseteq Y \\
  &\Downarrow \\
\forall c\in\R\quad c\le Y\;\; &\Longrightarrow\;\; c\le X \\
  &\Downarrow \\
\{c\in\R:\;\; c\le Y\} &\subseteq \{c\in\R:\;\; c\le X\} \\
  &\Downarrow \qquad (\text{\scriptsize свойство $1^\circ$ на с. \pageref{monotonnost-min-max-X} }) \\
\inf Y=\max\{c\in\R:\;\; c\le Y\} &\le \max\{c\in\R:\;\; c\le X\}=\inf X
\end{align*}

2. Обозначим
$$
B=\{b\in\R:\;\; b\le X\}, \qquad C=\{c\in\R:\;\; c\le X+a\}
$$
Тогда получаем цепочку
\begin{align*}
b\in B \quad\Longleftrightarrow\quad b\le X \quad&\Longleftrightarrow\quad
b+a\le X+a \quad\Longleftrightarrow\quad b+a\in C \\
  &\Downarrow \\
 B+a &=C \\
  &\Downarrow \\
\Big(\inf X\Big)+a=\Big(\max B\Big)+a =& (\text{\scriptsize свойство $2^\circ$
на с. \pageref{min-X+a=min-X+a}}) =\max(B+a)=\max C=\inf(X+a)
\end{align*}

3. Обозначим
$$
B=\{b\in\R:\;\; b\le X\}, \qquad C=\{c\in\R:\;\; -X \le c\}
$$
Тогда получаем цепочку
\begin{align*}
c\in C \quad\Longleftrightarrow\quad -X\le c \quad&\Longleftrightarrow\quad
-c\le X \quad\Longleftrightarrow\quad -c\in B \quad\Longleftrightarrow\quad c\in -B\\
  &\Downarrow \\
 -B &=C \\
  &\Downarrow \\
-\inf X=-\max B =& (\text{\scriptsize свойство $3^\circ$ на с.
\pageref{min-X=-max_X}}) =\min\Big(-B\Big)=\min C=\sup(-X)
\end{align*}
\epr

\subsection{Модуль}

{\it Модулем} (или {\it абсолютной величиной}) вещественного числа $x\in \R$
называется число $|x|$, определяемое формулой
 \beq\label{|x|=cases}
|x|=\begin{cases}x, & \text{если}\,\, x\ge 0
 \\ x, & \text{если}\,\, x<0 \end{cases}
 \eeq
Ниже на с.\pageref{|x|-ex-func} мы отметим, что отображение $x\mapsto|x|$
является функцией. Свойства этой функции так часто используются в анализе, что
их изучению мы посвятим целый раздел.

\btm Справедливо тождество
 \beq\label{|x|=max-x;x}
|x|=\max\{-x;x\}=(-x)\vee x,\qquad x\in\R
 \eeq
\etm
 \bpr
Здесь нужно рассмотреть три случая:
 \bit{
\item[1)] если $x>0$, то $|x|=x$ и, с другой стороны, $\max\{-x;x\}=x$, поэтому
$|x|=\max\{-x;x\}$;

\item[2)] если $x=0$, то $|x|=0=x$ и, с другой стороны, $\max\{-x;x\}=0=x$,
поэтому $|x|=\max\{-x;x\}$;

\item[3)] если $x<0$, то $|x|=-x$ и, с другой стороны, $\max\{-x;x\}=-x$,
поэтому $|x|=\max\{-x;x\}$.
 }\eit
Мы получаем, что, \eqref{|x|=max-x;x} выполняется независимо от знака $x$.
 \epr

\paragraph{Алгебраические тождества с модулем.}
 \btm Для любых $x,y\in\R$, $n\in\N$ справедливы следующие равенства:
 \begin{align}
 |x\cdot y|&=|x|\cdot |y| &
\label{module-4^0}\\
 |x^n| &= |x|^n & \label{module-5^0}
 \end{align}
 \etm
 \bpr
1. При доказательстве $6^\circ$ употребляются те же рассуждения, что и в случае
с $2^\circ$: если оба знака у $x$ и $y$ положительны, получаем
$$
\begin{cases} x\ge 0\\ y\ge 0\end{cases}
\quad\Longrightarrow\quad \underbrace{|x\cdot y|}_{\scriptsize\begin{matrix}\|
\\ x\cdot y\end{matrix}}=\underbrace{|x|}_{\scriptsize\begin{matrix}\|
\\ x\end{matrix}}\cdot \underbrace{|y|}_{\scriptsize\begin{matrix}\|
\\ y\end{matrix}}
$$
Во-вторых, если оба знака отрицательны,
$$
\begin{cases} x\le 0\\ y\le 0\end{cases}
\quad\Longrightarrow\quad \underbrace{|x\cdot y|}_{\scriptsize\begin{matrix}\|
\\ x\cdot y\end{matrix}}=\underbrace{|x|}_{\scriptsize\begin{matrix}\|
\\ -x\end{matrix}}\cdot \underbrace{|y|}_{\scriptsize\begin{matrix}\|
\\ -y\end{matrix}}
$$
Наконец, если знаки разные, например, $x\ge 0$, а $y\le 0$, то
$$
\begin{cases} x\ge 0\\ y\le 0\end{cases}
\quad\Longrightarrow\quad\underbrace{|x\cdot y|}_{\scriptsize\begin{matrix}\|
\\ -x\cdot y\end{matrix}}=\underbrace{|x|}_{\scriptsize\begin{matrix}\|
\\ x\end{matrix}}\cdot \underbrace{|y|}_{\scriptsize\begin{matrix}\|
\\ -y\end{matrix}}
$$

2. Свойство $7^\circ$ получается индукцией из $6^\circ$.
 \epr

\paragraph{Алгебраические неравенства с модулем.}
 \btm\label{TH:nerav-s-modulem} Для любых $x,y\in\R$ справедливы следующие неравенства:
 \begin{align}
|x|&\ge 0,\quad\text{причем}\quad
|x|=0\;\Longleftrightarrow\; x=0 &  \label{module-1^0}\\
|x|-|y|\le |x&+y|\le |x|+|y| & \label{module-2^0}\\
-|x-y|\le |x|&-|y|\le |x-y| & \label{module-3^0}
 \end{align}
 \etm
\bpr 1. Если $x\ge 0$, то $|x|=x\ge 0$. Если же $x<0$, то $|x|=-x>0$. Далее,
если $x=0$, то $|x|=\max\{0;0\}=0$. И наоборот, если $|x|=0$, то либо
$|x|=x=0$, либо $|x|=-x=0$, то есть опять $x=0$.

2. Докажем сначала правую половину \eqref{module-2^0}:
 \beq\label{|x+y|-le-|x|+|y|}
 |x+y|\le |x|+|y|
 \eeq
Во-первых, если оба знака у $x$ и $y$ положительны, получаем
$$
\begin{cases} x\ge 0\\ y\ge 0\end{cases}
\quad\Longrightarrow\quad \underbrace{|x+y|}_{\scriptsize\begin{matrix}\|
\\ x+y\end{matrix}}\le \underbrace{|x|}_{\scriptsize\begin{matrix}\|
\\ x\end{matrix}}+\underbrace{|y|}_{\scriptsize\begin{matrix}\|
\\ y\end{matrix}}
$$
Во-вторых, если оба знака отрицательны,
$$
\begin{cases} x\le 0\\ y\le 0\end{cases}
\quad\Longrightarrow\quad \underbrace{|x+y|}_{\scriptsize\begin{matrix}\|
\\ -(x+y)\end{matrix}}\le \underbrace{|x|}_{\scriptsize\begin{matrix}\|
\\ -x\end{matrix}}+\underbrace{|y|}_{\scriptsize\begin{matrix}\|
\\ -y\end{matrix}}
$$
Наконец, если знаки разные, например, $x\ge 0$, а $y\le 0$, то
$$
\begin{cases} x\ge 0\\ y\le 0\end{cases}
\quad\Longrightarrow\quad
\begin{cases} x+y=|x|-|y|\le |x|+|y|\\ -x-y=-|x|+|y|\le |x|+|y|\end{cases}
\quad\Longrightarrow\quad |x+y|=\max\{x+y;-x-y\}\le |x|+|y|
$$
Из \eqref{|x+y|-le-|x|+|y|} теперь получаем левую половину \eqref{module-2^0}:
$$
|x|=|(x+y)-y|\le |x+y|+|-y|=|x+y|+|y|\quad\Longrightarrow\quad |x|-|y|\le |x+y|
$$

3. Докажем правую половину \eqref{module-3^0}:
$$
|x|=|(x-y)+y|\le\eqref{|x+y|-le-|x|+|y|}\le |x-y|+|y| \quad\Longrightarrow\quad
|x|-|y|\le |x-y|
$$
После этого становится очевидной левая половина:
$$
|y|-|x|\le |y-x|\quad\Longrightarrow\quad -|y|+|x|\ge -|y-x|=-|x-y|
$$
\epr

\paragraph{Решение неравенств с модулем.}

 \btm Справедливы следующие правила решения неравенств с модулем:
 \begin{align}
 |x|<\varepsilon \;\Longleftrightarrow\; x\in
(-\varepsilon,&\varepsilon) \;\Longleftrightarrow\; |x|\in (-\varepsilon,
\varepsilon)  & \label{|x|<e}\\
 |x|>\e \;\Longleftrightarrow\; x\in (-\infty; -\e)\cup
&(\e;+\infty) \;\Longleftrightarrow\; |x|\in (-\infty;-\e)\cup (\e;+\infty) &
\label{|x|>e}
 \end{align}
 \etm
 \bpr
6. Если $x\ge 0$, то $|x|=x$, и цепочка \eqref{|x|<e} становится очевидной:
$$
\underbrace{|x|<\varepsilon}_{\scriptsize\begin{matrix}\Updownarrow
\\ 0\le x<\e \end{matrix}} \;\Longleftrightarrow\; x\in (-\varepsilon,\varepsilon)
\;\Longleftrightarrow\; \underbrace{|x|}_{\scriptsize\begin{matrix}\|
\\ x \end{matrix}}\in (-\varepsilon, \varepsilon)
$$
Если же $x\le 0$, то $|x|=-x$, и получается
$$
\underbrace{|x|<\varepsilon}_{\scriptsize\begin{matrix}\Updownarrow
\\ 0\le -x<\e \\ \Updownarrow \\ 0\ge x>-\e \end{matrix}} \;\Longleftrightarrow\; x\in
(-\varepsilon,\varepsilon) \;\Longleftrightarrow\;
\underbrace{|x|}_{\scriptsize\begin{matrix}\|
\\ -x \end{matrix}}\in (-\varepsilon, \varepsilon)
$$

7. То же самое с $3^\circ$: если $x\ge 0$, то $|x|=x$, и цепочка \eqref{|x|>e}
приобретает вид,
$$
\underbrace{|x|>\varepsilon}_{\scriptsize\begin{matrix}\Updownarrow
\\ x>\e \end{matrix}} \;\Longleftrightarrow\; x\in (-\infty;\varepsilon)\cup
(\varepsilon;+\infty) \;\Longleftrightarrow\;
\underbrace{|x|}_{\scriptsize\begin{matrix}\|
\\ x \end{matrix}}\in (-\infty;\varepsilon)\cup (\varepsilon;+\infty)
$$
а если $x\le 0$, то $|x|=-x$, и получается
$$
\underbrace{|x|>\varepsilon}_{\scriptsize\begin{matrix}\Updownarrow
\\ -x>\e \\ \Updownarrow \\ x<-\e \end{matrix}} \;\Longleftrightarrow\; x\in
(-\infty;\varepsilon)\cup (\varepsilon;+\infty) \;\Longleftrightarrow\;
\underbrace{|x|}_{\scriptsize\begin{matrix}\|
\\ -x \end{matrix}}\in (-\infty;\varepsilon)\cup (\varepsilon;+\infty)
$$
 \epr

\section{Натуральные числа и конечные множества}

\subsection{Натуральные числа $\N$ и целые неотрицательные числа $\Z_+$}\label{natur-chisla-N}

\paragraph{Определение натуральных чисел $\N$.}

В школьном курсе математики объясняется, что натуральные
числа -- это числа вида $1,2,3,...$ Однако, в этих словах
имеется некоторое лукавство, потому что такое
<<определение>> нельзя считать математически строгим.
Действительно, почему, например, число 5 является
натуральным? Здесь про него ничего не сказано, и поэтому
формально у нас нет возможности проверить, будет оно
натуральным, или нет. То же самое можно сказать про числа
$\frac{1}{5}$, $\sqrt{2}$, $\pi$ и любые другие.

В этом разделе мы объясним, как определяются натуральные числа, исходя из
аксиом $A1 - A16$. Все начинается с так называемых индуктивных множеств:
 \bit{
 \item[$\bullet$]
Множество вещественных чисел $X$ называется {\it
индуктивным}, если
 \bit{
\item[1)] $X$ содержит единицу (о которой говорится в аксиоме A8):
$$
  1\in X
$$
\item[2)] вместе с каждым числом $x\in X$ ему принадлежит также
число $x+1$:
$$
  \forall x\in X\qquad x+1\in X
$$
 }\eit
 }\eit

\noindent\rule{160mm}{0.1pt}
\begin{multicols}{2}

\bexs Само множество $\R$ является индуктивным. Полуинтервалы вида $(a,
\infty)$ где $a<1$, или $[a, \infty)$ где $a\le 1$ также являются индуктивными
множествами. Наоборот, скажем, множества $(1,\infty)$, $[2,\infty)$, $[0,5]$ не
являются индуктивными. \eexs

\end{multicols}\noindent\rule[10pt]{160mm}{0.1pt}

Теперь главное определение:
 \bit{\label{DEF:N}
 \item[$\bullet$] Наименьшее индуктивное множество обозначается
символом $\N$ и называется {\it множеством натуральных чисел}, или {\it
натуральным рядом}. Слово ``наименьшее'' здесь означает, что выполняются два
условия:
 \bit{
\item[$(i)$] $\mathbb{N}$ является индуктивным множеством;
\item[$(ii)$] если $X$ -- какое-нибудь другое
индуктивное  множество, то $\mathbb{N}$ содержится в $X$.
 }\eit
 }\eit

\begin{proof}[Доказательство существования и единственности.]\label{dokazatelstvo-dlya-N}

1. Покажем сначала, что такое множество $\mathbb{N}$ действительно существует.
Обозначим через $\N$ пересечение всех индуктивных множеств $X$:
$$
  \N=\bigcap_{\scriptsize\begin{matrix}X - \text{индуктивное}\\ \text{множество}\end{matrix}} X
$$
Иными словами,
$$
y\in \N \quad \Leftrightarrow \quad  y \,\, \text{принадлежит любому
индуктивному множеству} \, X
$$
Тогда мы получим:
 \bit{
\item[1)] $1\in\N$, потому что 1 лежит в любом множестве $X$, и

\item[2)] если $x\in\N$, то есть $x\in\N$ для всякого индуктивного множества
$X$, то $x+1\in X$ (поскольку $X$ индуктивно); значит, $x+1\in X$ для всякого
индуктивного множества $X$, поэтому $x+1\in\N$.
 }\eit\noindent
Таким образом, $\N$ -- индуктивное множество. С другой стороны, $\N$ содержится
в любом множестве $X$, то есть, в любом другом индуктивном множестве. Значит,
$\N$ как раз и есть множество натуральных чисел.

2. Убедимся, что множество $\mathbb{N}$ определяется однозначно. Действительно,
если бы нам было дано какое-то другое множество $\widetilde{\mathbb{N}}$ с теми
же свойствами, то мы получили бы, что, во-первых,
$$
\mathbb{N}\subseteq\widetilde{\mathbb{N}}
$$
(потому что $\mathbb{N}$ содержится в любом индуктивном множестве, в частности,
в $\widetilde{\mathbb{N}}$), и, во-вторых,
$$
\widetilde{\mathbb{N}}\subseteq\mathbb{N}
$$
(потому что $\widetilde{\mathbb{N}}$ содержится в любом индуктивном множестве,
в частности, в $\mathbb{N}$). Вместе это означает (в силу
\eqref{A-subseteq-B-&-B-subseteq-A}), что
$$
\mathbb{N}=\widetilde{\mathbb{N}}.
$$
\end{proof}

\paragraph{Принцип математической индукции.}

\begin{tm}[\bf принцип математической индукции в абстрактной форме]
\label{induction-abs} Пусть $E$ -- какое-нибудь подмножество в
$\mathbb{N}$, обладающее свойствами:
 \bit{
\item[$(a)$] $1\in E$
\item[$(b)$] $\forall n\in E \quad n+1 \in E$
 }\eit
Тогда $E=\mathbb{N}$.
\end{tm}
\begin{proof} С одной стороны, нам дано, что
$E\subseteq \mathbb{N}$. С другой же стороны, свойства $(a),(b)$
означают, что $E$ является индуктивным множеством. Поэтому, в силу
условия $(ii)$ определения 3.2, $\mathbb{N}$ содержится в $E$:
$\mathbb{N}\subseteq E$. Итак,
$$
E\subseteq \mathbb{N} \quad \& \quad  \mathbb{N}\subseteq E
$$
то есть, $E=\mathbb{N}$. \end{proof}

На практике удобно пользоваться следующей переформулировкой принципа
математической индукции.

\begin{tm}[\bf принцип математической индукции]\label{induction}
Пусть дана последовательность утверждений
\begin{equation}
\varPhi_n \qquad\qquad n\in \mathbb{N} \label{0.3.1}
\end{equation}
со следующими свойствами:
 \bit{
\item[$(a')$] первое утверждение $\varPhi_1$ истинно;

\item[$(b')$] если при каком-нибудь $n\in\N$ истинно утверждение $\varPhi_n$, то
истинно и утверждение $\varPhi_{n+1}$.
 }\eit\noindent
Тогда все утверждения $\varPhi_n$ истинны (для всех $n\in \mathbb{N}$).
\end{tm}
\begin{proof} Обозначим через $E$ множество всех
таких чисел $n\in \mathbb{N}$, для которых утверждение $\varPhi_n$ истинно. Тогда
условие $(a')$ будет эквивалентно условию $(a)$ теоремы \ref{induction-abs}, а
условие $(b')$ будет эквивалентно условию $(b)$. Значит, $E=\mathbb{N}$, то
есть $\varPhi_n$ истинно для каждого $n\in \mathbb{N}$.
\end{proof}

\paragraph{Свойства натуральных чисел.}

Первое следствие из принципа математической индукции выглядит так:

\btm\label{zamknutost-N-otn-+-i-umn} Множество $\N$ натуральных чисел замкнуто
относительно сложения и умножения: если $m,n\in\N$, то $m+n\in\N$ и $m\cdot
n\in\N$.
 \etm
\bpr 1. Зафиксируем $m\in\N$ и обозначим через $E$ множество чисел $n\in\N$
таких, что $m+n\in\N$:
$$
E=\{n\in\N:\quad m+n\in\N\}
$$
Нам нужно убедиться, что $E=\N$. Это делается индукцией. Сначала замечаем, что
$1\in E$, потому что из $m\in\N$ следует $m+1\in\N$. Затем берем какое-нибудь
$k\in E$, то есть
$$
m+k\in\N
$$
и замечаем, что тогда $k+1\in E$, потому что
$$
m+(k+1)=(m+k)+1\in\N
$$
По принципу индукции все это означает, что $E=\N$.

2. Опять фиксируем $m\in\N$ и обозначаем через $E$ множество чисел $n\in\N$
таких, что $m\cdot n\in\N$:
$$
E=\{n\in\N:\quad m\cdot n\in\N\}
$$
Нам нужно убедиться, что $E=\N$. Это тоже делается индукцией. Сначала замечаем,
что $1\in E$, потому что из $m\in\N$ следует $m\cdot 1=m\in\N$. Затем берем
какое-нибудь $k\in E$, то есть
$$
m\cdot k\in\N
$$
и замечаем, что тогда $k+1\in E$, потому что из только что доказанной
замкнутости $\N$ относительно сложения следует
$$
m\cdot k\in\N \quad\Longrightarrow\quad m\cdot (k+1)=(m\cdot k)+m\in\N
$$
Снова по принципу индукции $E=\N$.
 \epr

Ниже нам понадобятся еще несколько утверждений об $\N$.

 \bit{
\item[$\bullet$] {\it Дискретным интервалом} с концами $m\in\N$ и $n\in\N$
называется множество
$$
\{m,...,n\}:=\{k\in\N:\; m\le k\le n\}
$$
(очевидно, если $m>n$, то $\{m,...,n\}=\varnothing$).

\item[$\bullet$] {\it Начальным дискретным интервалом}
\label{DEF:nach-interval-v-N} называется дискретный интервал, у которого левый
конец равен 1:
$$
\{1,...,n\}=\{k\in\N:\; k\le n\}
$$
 }\eit

\bigskip

\centerline{\bf Свойства натуральных чисел}
 \bit{\it
 \nopagebreak[4]
\item[$1^\circ$.] Для всякого $n\in\N$ выполняется $n\ge 1$. Поэтому,
 \beq\label{min-N=1}
\min\N=1
 \eeq

\item[$2^\circ$.] Если $n\in\N$ и $n\ne 1$, то $n-1\in\N$.

\item[$3^\circ$.] Для любого $n\in\N$ множество чисел $x\in\N$, больших $n-1$,
имеет наименьший элемент, а именно $n$:
 \beq\label{min-X_n=n}
\min\{x\in\N:\; n-1<x\}=n
 \eeq

\item[$4^\circ$.]\label{minimum-v-N} Каждое непустое подмножество
$M\subseteq\N$ обладает минимальным элементом.

\item[$5^\circ$.]\label{N:m<n=>m-le-n-1} Для любых $m,n\in\N$ справедливы
импликации:
 \beq\label{m<n+1=>m-le-n}
m<n+1\qquad\Longrightarrow\qquad m\le n
 \eeq
 \beq\label{m<n=>m+1-le-n}
m<n\qquad\Longrightarrow\qquad m+1\le n
 \eeq

\item[$6^\circ$.]\label{m<n->n-m-in-N} Если $m,n\in\N$ и $m<n$, то $n-m\in\N$.

\item[$7^\circ$.]\label{konechnaya-induktsiya} Пусть $I$ -- подмножество в
начальном интервале $\{1,...,n\}$, удовлетворяющее условиям
 \bit{
\item[$(a)$] $1\in I$

\item[$(b)$] $\forall k<n \qquad \Big(k\in I \;\Longrightarrow\; k+1 \in
I\Big)$
 }\eit
тогда $I=\{1,...,n\}$.

\item[$8^\circ$.]\label{osn-sv-nachalnyh-intervalov} Пусть $E$ -- подмножество
в $\N$, удовлетворяющее условиям
 \bit{
\item[$(a)$] $1\in E$

\item[$(b)$] $\forall k\in E\qquad \{1,...,k\}\subseteq E$
 }\eit
Тогда
 \bit{
\item[---] либо $E=\N$,

\item[---] либо $E$ -- начальный интервал:
 \beq
\exists n\in E\qquad \{1,...,n\}=E
 \eeq

 }\eit

}\eit

\bpr Все эти утверждения доказываются математической индукцией.

1. Формула $n\ge 1$ верна для любого $n\in\N$, потому что, во-первых, она верна
при $n=1$: $1\ge 1$. И, во-вторых, если она верна при некотором $n=m$, $m\ge
1$, то и при $n=m+1$ она тоже верна: $m+1\ge 1+1\ge 1$.

2. Рассмотрим множество
 $$
E=\{1\}\cup (\N+1)=\{x\in\R:\quad x=1\;\V\; \exists n\in\N\quad x=n+1\}
 $$
и заметим, что $E=\N$. Действительно, во-первых, $1\in E$, а во-вторых, если
$m\in E$, то либо $m=1$, и тогда $m+1\in\N+1\subseteq E$, либо $m=n+1$ для
некоторого $n\in\N$, и тогда $m+1=(n+1)+1\in\N+1\subseteq E$. Теперь из $E=\N$
получаем: если $n\in\N=E=\{1\}\cup (\N+1)$ и $n\ne 1$, то $n\in\N+1$, то есть
$n=k+1$ для некоторого $k\in\N$.

3. Обозначим $X_n=\{x\in\N:\; n-1<x\}$ и заметим, что
 \beq\label{X_n+1=X_n+1}
\underbrace{X_{n+1}}_{\scriptsize\begin{matrix}\| \\
\{y\in\N:\;n<y\}\end{matrix}}=\underbrace{X_n+1}_{\scriptsize\begin{matrix}\| \\
\{y\in\R:\;\exists x\in X_n\quad y=x+1\}
\end{matrix}}
 \eeq
Действительно, с одной стороны, $X_n+1\subseteq X_{n+1}$, потому что
 \begin{multline*}
y\in X_n+1\quad\Longrightarrow\quad \exists x\in X_n\quad y=x+1
\quad\Longrightarrow\quad \exists x\in\N\quad n-1<x \;\&\; y=x+1
\quad\Longrightarrow\\ \Longrightarrow\quad \exists x\in\N\quad n<x+1=y
\quad\Longrightarrow\quad y\in X_{n+1}
\end{multline*}
а, с другой -- $X_{n+1}\subseteq X_n+1$, потому что
 \begin{multline*}
y\in X_{n+1}\quad\Longrightarrow\quad y\in\N\;\&\;
\underbrace{n<y}_{\scriptsize\begin{matrix}\Downarrow \\ y>\min\N=1
\\ \Downarrow \\ 2^\circ \\ \Downarrow \\ x:=y-1\in\N
\end{matrix}} \quad\Longrightarrow\quad \exists x\in\N\quad y=x+1\;\&\;n<x+1
\quad\Longrightarrow \\ \Longrightarrow\quad \exists x\in\N\quad
y=x+1\;\&\;n-1<x \quad\Longrightarrow\quad \exists x\in X_n\quad y=x+1
\quad\Longrightarrow\quad y\in X_n+1
 \end{multline*}

Теперь доказываем \eqref{min-X_n=n} по индукции. Во-первых, убеждаемся, что эта
формула верна при $n=1$. Действительно,
$$
\underbrace{x\in\N\stackrel{1^\circ}{\quad\Longrightarrow\quad} x\ge
1>0\quad\Longrightarrow\quad x\in X_1}
$$
$$
\Downarrow
$$
$$
\N\subseteq X_1
$$
$$
\Downarrow
$$
$$
X_1=\N
$$
$$
\Downarrow
$$
$$
\min X_1=\min\N=1
$$
Далее предполагаем, что \eqref{min-X_n=n} верна при $n=m$:
 \beq\label{min-X_m=m}
\min X_m=m
 \eeq
и докажем, что тогда она верна при $n=m+1$:
 $$
\min X_{m+1}=\eqref{X_n+1=X_n+1}=\min (X_m+1)=(\min
X_m)+1=\eqref{min-X_m=m}=m+1
 $$

4. Пусть $M\subseteq\N$ и $M\ne\varnothing$. Если $1\in M$, то доказывать
нечего, потому что тогда $1=\min M$. Поэтому будем считать, что $1\notin M$.
Тогда множество $E=\N\setminus M$, наоборот, должно содержать 1:
$$
1\in E
$$
Покажем, что существует такое $n\in\N$, что
 \beq\label{1,...,n-subseteq-E-&-n+1-notin-E}
\{1,...,n\}\subseteq E\quad\&\quad n+1\notin E
 \eeq
Действительно, если бы это было не так, то есть при любом $n\in E$ мы получали
бы, что из $\{1,...,n\}\subseteq E$ следует $n+1\in E$, то это означало бы, что
$E=\N$, то есть $M=\varnothing$.

Теперь, возвращаясь от $E$ к $M$, условие
\eqref{1,...,n-subseteq-E-&-n+1-notin-E} можно переформулировать так:
$$
\{1,...,n\}\cap M=\varnothing\quad\&\quad n+1\in M
$$
Это и означает, что $n+1=\min M$.

5. Пусть $m,n\in\N$ и $m<n+1$, то есть $m-1<n$. Тогда $n\in\{x\in\N:\; m-1<x\}$
и по свойству $3^\circ$ получаем $n\ge \min\{x\in\N:\; m-1<x\}=m$. Это
доказывает импликацию \eqref{m<n+1=>m-le-n}. А из нее следует
\eqref{m<n=>m+1-le-n}:
$$
m<n\quad\overset{\eqref{x<y->x+a<y+a}}{\Longrightarrow}\quad
m+1<n+1\quad\overset{\eqref{m<n+1=>m-le-n}}{\Longrightarrow}\quad m+1\le n
$$

6. Пусть $m\in\N$. Рассмотрим множество
$$
E=\{1,...,m\}\cup(\N+m)
$$
и заметим, что $E=\N$. Действительно, во-первых, $1\in E$, а во-вторых, если
$k\in E$, то возможны три случая, в каждом из которых получается $k+1\in E$:
 \bit{
\item[1)] $k<m=(m-1)+1\quad\Longrightarrow\quad$ (применяем $5^\circ$)
$\quad\Longrightarrow\quad k\le m-1\quad\Longrightarrow\quad k+1\le m
\quad\Longrightarrow\quad k+1\in E$;

\item[2)] $k=m\quad\Longrightarrow\quad k+1=1+m\in(\N+m)\subseteq E$;

\item[3)] $k>m\quad\Longrightarrow\quad k\notin
\{1,...,m\}\quad\Longrightarrow\quad k\in(\N+m)\quad\Longrightarrow\quad
k=n+m\quad (n\in\N)\quad\Longrightarrow\quad k+1=(n+1)+m\quad (n+1\in\N)
\quad\Longrightarrow\quad k+1\in(\N+m)\subseteq E$;
 }\eit
Таким образом, получается что $\N=E=\{1,...,m\}\cup(\N+m)$, и поэтому если
$n\in\N$, $n>m$, то, поскольку $n\notin \{1,...,m\}$, имеем $n\in\N+m$, то есть
$n=m+k$, $k\in\N$, или, иными словами, $n-m=k\in\N$.

7. Рассмотрим множество $M=\N\setminus I$. Тогда получаем цепочку
$$
I\subseteq \{1,...,n\} \quad\Longrightarrow\quad n+1\notin
I\quad\Longrightarrow\quad n+1\in M \quad\Longrightarrow\quad
n+1\ge\underbrace{\min M}_{\scriptsize\begin{matrix}\uparrow \\
\text{существует}\\ \text{в силу $4^\circ$} \end{matrix}}
$$
Нам нужно убедиться, что $\min M=n+1$: тогда получится, что $\{1,...,n\}\cap
M=\varnothing$, то есть $\{1,...,n\}\subseteq I$ и значит $I=\{1,...,n\}$.

Предположим, что это не так: $\min M\ne n+1$. Тогда $\min M<n+1$, и по свойству
$5^\circ$ $\min M\le n$. С другой стороны, в силу $(a)$, $1\in I$, поэтому
$1\notin M$, и значит $\min M>1$, то есть по свойству $5^\circ$, $\min M-1\ge
1$, и
$$
2\le \min M\le n
$$
Если теперь положить $k=\min M-1$, то получается вот что:
$$
1\le k<n\quad\&\quad k\notin M \quad\Longrightarrow\quad k\in I\quad\&\quad k<n
\stackrel{(b)}{\quad\Longrightarrow\quad} k+1=\underbrace{\min M\in
I}_{\scriptsize\begin{matrix}\uparrow \\
\text{невозможно,}\\ \text{потому что} \\ M\cap I=\varnothing \end{matrix}}
$$

8. Если $E\ne\N$, то найдется $x\in\N$ такое что $x\notin E$. То есть множество
$$
M=\{x\in\N:\quad x\notin E\}
$$
непусто. Значит, по уже доказанному свойству $4^\circ$, $M$ имеет минимальный
элемент:
$$
m=\min M
$$
По условию (a), $1\in E$, и поэтому $1\notin M$. С другой стороны, по уже
доказанному свойству $1^\circ$, $1\le m$. Значит, $1<m$. Поэтому по свойству
$6^\circ$, число $n=m-1$ тоже лежит в $\N$:
$$
n=m-1\in\N
$$
Тогда, во-первых, $n\notin M$ (потому что иначе число $m=n+1$ не было бы
минимумом для $M$), и значит $n\in E$, откуда в силу (b),
 \beq\label{(1,...,n)-subseteq-E}
 \{1,...,n\}\subseteq E
 \eeq
А, во-вторых, для $x\in\N$ справедлива импликация
 \beq\label{x>n=>x-notin-E}
x>n\quad\Longrightarrow\quad x\notin E,
 \eeq
потому что иначе мы получили бы
 $$
\underbrace{
\begin{matrix}
n<x \\
{\scriptsize\eqref{m<n=>m+1-le-n}}\quad\Downarrow\quad\phantom{\scriptsize\eqref{m<n=>m+1-le-n}}\\
m=n+1\le x
\end{matrix}
\qquad
\begin{matrix}
x\in E \\
\phantom{\scriptsize\text{\rm (b)}}\quad\Downarrow\quad{\scriptsize\text{\rm (b)}} \\
\{1,...,x\}\subseteq E
\end{matrix}}
 $$
 $$
 \Downarrow
 $$
 $$
m\in E
 $$
а это невозможно, поскольку $m\in M$. Утверждения \eqref{(1,...,n)-subseteq-E}
и \eqref{x>n=>x-notin-E} вместе означают, что $\{1,...,n\}=E$.
 \epr

\paragraph{Целые неотрицательные числа $\Z_+$.}

 \biter{
\item[$\bullet$] Множество целых неотрицательных чисел, обозначаемое символом
$\Z_+$ (и называемое иногда {\it расширенным натуральным рядом}), определяется
как объединение множества $\N$ натуральных чисел и нуля $\{0\}$:
 \beq\label{DF:Z_+}
\Z_+=\{0\}\cup \N
 \eeq
 }\eiter

Для этого множества, как и для $\N$ справедлив принцип математической индукции:

\begin{tm}[\bf принцип математической индукции в абстрактной форме для $\Z_+$]
\label{induction-abs-Z_+} Пусть $E$ -- какое-нибудь подмножество в $\Z_+$,
обладающее свойствами:
 \bit{
\item[$(a)$] $0\in E$

\item[$(b)$] $\forall n\in E \quad n+1 \in E$
 }\eit
Тогда $E=\Z_+$.
\end{tm}

\begin{tm}[\bf принцип математической индукции для $\Z_+$]\label{induction-Z_+}
Пусть дана последовательность утверждений
\begin{equation}
\varPhi_n \qquad\qquad n\in \Z_+ \label{0.3.1-Z}
\end{equation}
со следующими свойствами:
 \bit{
\item[$(a')$] первое утверждение $\varPhi_0$ верно;

\item[$(b')$] если при каком-нибудь $n\in\N$ верно утверждение $\varPhi_n$, то
верно и утверждение $\varPhi_{n+1}$.
 }\eit\noindent
Тогда все утверждения $\varPhi_n$ справедливы (для всех $n\in\Z_+$).
\end{tm}

\paragraph{Определения по индукции.}

Читатель мог заметить, что среди примеров собственных обозначений для чисел на
странице \pageref{chisla-2-10} не было чисел, больших 10. Мы могли бы, конечно,
по аналогии последовательно определить числа 11, 12, или присвоить собственные
обозначения вообще любому конечному набору натуральных чисел, но чтобы дать
собственные обозначения {\it всем} натуральным числам (которых бесконечно
много), нужно использовать так называемый прием определения по индукции. О
собственных обозначениях для чисел из $\N$ (точнее, об их десятичной записи) мы
поговорим ниже на странице \pageref{SUBSEC-desyat-zapis-N}, а здесь мы приведем
две теоремы, служащие обоснованием приема определения по индукции и
проиллюстрируем их примерами. Первая из них используется в случаях, когда нужно
определить бесконечную последовательность элементов:

\begin{tm}[\bf об определениях полной индукцией]\label{defin-induction}
Пусть $A$ -- произвольное множество (необязательно числовое) и дано отображение
$G$, которое {\bf каждой} конечной последовательности элементов $a_1,...,a_n\in
A$, $n\in\N$, ставит в соответствие определенный элемент $G(a_1,...,a_n)\in A$.
Тогда для всякого начального элемента $a_1\in A$ отображение $G$ однозначно
определяет бесконечную последовательность $\{a_n\}\subseteq A$, содержащую
$a_1$ в качестве первого элемента и удовлетворяющую условию
 \beq\label{ind-pravilo-N}
\forall n\in\N\qquad G(a_1,...,a_n)=a_{n+1}
 \eeq
\end{tm}
 \bpr
Зафиксируем $a_1\in A$ и обозначим через $E$ подмножество в $\N$, состоящее из таких
$n\in\N$, для которых существует конечная последовательность $a_1,...,a_n\in A$
со следующим свойством:
 \beq\label{ind-pravilo}
\forall k=2,...,n \qquad G(a_1,...,a_{k-1})=a_k
 \eeq
Тогда, во-первых, $1\in E$, потому что для $n=1$ нужная последовательность
должна состоять только из одного числа, и таким числом будет $a_1$. И,
во-вторых, если $n\in E$, то это значит, что существует последовательность
$a_1,...,a_n\in A$ для которой выполняется \eqref{ind-pravilo}. Тогда можно
положить $a_{n+1}=G(a_1,...,a_n)$, и мы получим что $n+1\in E$.

Таким образом, по принципу математической индукции \eqref{induction-abs},
$E=\N$. Отсюда следует, что для всякого $n\in \N$ можно выбрать $a_n$ по такому
правилу: выбираем какую-нибудь конечную последовательность $\{a_1,...,a_n\}$,
подчиненную условию \eqref{ind-pravilo} (такая последовательность существует
поскольку $n\in E$), и в качестве $a_n$ берем ее последний элемент. Нам нужно
только убедиться, что число $a_n$ определяется таким алгоритмом однозначно.
Действительно, если предположить, что существует какое-то другое число $a_n'$,
которое можно определить таким же образом, то есть для которого можно построить
конечную последовательность $\{a_1',...,a_n'\}$ так, чтобы
 $$
a_1'=a_1\qquad\&\qquad \forall k=1,...,n-1 \qquad G(a_1',...,a_k')=a_{k+1}'
 $$
то, рассмотрев множество индексов
$$
I=\{k=1,...,n:\quad \forall i\le k\quad a_i'=a_i\}
$$
мы получим, что, во-первых, $1\in I$ (потому что $a_1'=a_1$), и, во-вторых,
$$
k<n\qquad\&\qquad k\in I
$$
$$
\Downarrow
$$
$$
k<n\qquad\&\qquad \forall i\le k\quad a_i'=a_i
$$
$$
\Downarrow
$$
$$
a_{k+1}'=G(a_1',...,a_k')=G(a_1,...,a_k)=a_{k+1}
$$
$$
\Downarrow
$$
$$
k+1\le n\qquad\&\qquad \forall i\le k+1\quad a_i'=a_i
$$
$$
\Downarrow
$$
$$
k+1\in I
$$
Отсюда, по принципу конечной индукции (свойство $7^\circ$ на
с.\pageref{konechnaya-induktsiya}), $I=\{1,...,n\}$, то есть, в частности,
$a_n'=a_n$.
 \epr

Вторая теорема бывает нужна, когда нет гарантии, что определяемая тобой
последовательность будет бесконечной:

\begin{tm}[\bf об определениях частной индукцией]\label{defin-induction-chast}
Пусть $A$ -- произвольное множество (необязательно числовое) и дано отображение
$G$, которое {\bf некоторым} конечным последовательностям элементов
$a_1,...,a_n\in A$, $n\in\N$, ставит в соответствие определенные элементы
$G(a_1,...,a_n)\in A$. Тогда для всякого начального элемента $a_1\in A$ это
правило однозначно определяет (конечную или бесконечную) последовательность
$\{a_n\}\subseteq A$, содержащую $a_1$ в качестве первого элемента и
удовлетворяющую условию
 \beq\label{ind-pravilo-N-chast}
G(a_1,...,a_n)=a_{n+1},\qquad n<N
 \eeq
где величина $N$ зависит от множества $A$, отображения $G$ и начального
элемента $a_1$, и определяется условием:
 \beq\label{dlina-ind-posledov}
N=\sup\Big\{n\in\N:\quad \exists \{a_2,...,a_n\}\subseteq A\quad \forall
k\in\{2,...,n\} \qquad G(a_1,...,a_{k-1})=a_k\Big\}
 \eeq
Если $N$ конечна, то она представляет собой длину последовательности $\{a_n\}$.
\end{tm}
 \bpr
Зафиксируем $a_1$ и обозначим через $E$ подмножество в $\N$, от которого
берется верхняя грань в формуле \eqref{dlina-ind-posledov}
$$
N=\sup E,
$$
то есть состоящее из таких $n\in\N$, для которых существует конечная
последовательность $a_2,...,a_n\in A$ со следующим свойством:
 \beq\label{ind-pravilo-chast}
\forall k=2,...,n \qquad G(a_1,...,a_{k-1})=a_k
 \eeq
Это множество $E$ будет удовлетворять условиям (a) и (b) свойства $8^\circ$ на
с.\pageref{osn-sv-nachalnyh-intervalov}:
 \bit{
\item[$(a)$] $1\in E$ (потому что для $n=1$ нужная последовательность должна
состоять только из одного числа, и таким числом будет $a_1$), и

\item[$(b)$] $\forall n\in E$ $\{1,...,n\}\subseteq E$ (потому что если $n\in
E$, то есть существуют $a_1,...,a_n\in A$ со свойством
\eqref{ind-pravilo-chast}, то для всякого $k\le n$ последовательность
$a_2,...,a_k$ будет обладать теми же свойствами, только с заменой $k$ на
какое-нибудь $i$, а после этого $n$ на $k$).
 }\eit
По свойству $8^\circ$ на с.\pageref{osn-sv-nachalnyh-intervalov} отсюда
следует, что либо $E=\N$, либо $E=\{1,...,N\}$ для некоторого $N\in\N$. В обоих
случаях мы получаем, что каждому $n\in E$ можно поставить в соответствие
некоторый элемент $a_n$ по такому правилу: выбираем какую-нибудь конечную
последовательность $\{a_2,...,a_n\}$, подчиненную условию \eqref{ind-pravilo}
(такая последовательность существует поскольку $n\in E$), и в качестве $a_n$
берем ее последний элемент. Здесь нужно только убедиться, что элемент $a_n$
определяется таким алгоритмом однозначно. Это делается в точности так же, как
при доказательстве теоремы \ref{defin-induction}.
 \epr

\noindent\rule{160mm}{0.1pt}
\begin{multicols}{2}

\paragraph*{Степени с натуральным показателем.}\label{opr-stepeni}

 \biter{
\item[$\bullet$] {\it Степенью} $a^n$ числа $a\in\R$ с показателем $n\in\N$,
как известно, называется число
$$
a^n=\underbrace{a\cdot a\cdot...\cdot a}_{\text{$n$ множителей}}
$$
 }\eiter
Теорема \ref{defin-induction} позволяет придать точный смысл этим словам: чтобы
определить $a^n$ для любого $n\in \N$, нужно положить
 \beq\label{a^1=a}
a^1=a
 \eeq
а затем сказать, что если $a^n$ уже определено, то $a^{n+1}$ определяется
формулой:
$$
a^{n+1}=a^n\cdot a
$$
Коротко такое определение по индукции удобно записывать так:
$$
a^1=a,\quad a^{n+1}=a^n\cdot a \quad (a\in\R,\ n\in\N)
$$
В дальнейшем нам будет удобно отдельно положить
$$
a^0=1\qquad (a\in\R),
$$
и тогда индуктивное определение для $a^n$ можно будет переписать в виде
 \beq\label{def:a^n}
a^0=1,\quad a^{n+1}=a^n\cdot a \quad (a\in\R,\ n\in\Z_+)
 \eeq

\brem Это определение, между прочим, предполагает, что $0^0$ равно единице:
 \beq\label{0^0=1}
0^0=1
 \eeq
\erem

\paragraph*{Факториал $n!$ и двойной факториал $n!!$.}

 \biter{
\item[$\bullet$] {\it Факториал} числа $n\in\Z_+$ определяется по индукции
следующим правилом:
 \beq\label{DEF:n!}
\kern-40pt 0!=1, \qquad (n+1)!=(n+1)\cdot n!
 \eeq

\item[$\bullet$] {\it Двойной факториал} числа $n\in\Z_+$ определяется по
индукции следующими правилами:
 \biter{
\item[---] если число $n$ четно, то есть имеет вид $n=2m$, $m\in\Z_+$, то
 \beq\label{DEF:(2m)!!}
\kern-40pt 0!!=1, \quad (2m+2)!!=(2m+2)\cdot (2m)!!
 \eeq
\item[---] если число $n$ нечетно, то есть имеет вид $n=2m-1$, $m\in\N$, то
 \beq\label{DEF:(2m-1)!!}
\kern-40pt 1!!=1, \quad (2m+1)!!=(2m+1)\cdot (2m-1)!!
 \eeq
 }\eiter
 }\eiter

\paragraph*{Число сочетаний $C_n^k$.}
Для любых $n,k\in\Z_+$, $k\le n$ число
 \beq\label{def-C_n^k}
C_n^k=\frac{n!}{(n-k)!\cdot k!},\qquad (0\le k\le n)
 \eeq
называется {\it числом сочетаний из $n$ по $k$}.

 \bprop Справедливо следующее тождество
 \beq\label{C_n^k+C_n^k-1=C_n+1^k}
C_n^k+C_n^{k-1}=C_{n+1}^k,
 \eeq
называемое {\rm тождеством Паскаля}.
 \eprop
\bpr Это доказывается простой проверкой:
$$
C_n^k+C_n^{k-1}=C_{n+1}^k
$$
$$
\Updownarrow
$$
$$
\frac{n!}{(n-k)!k!}+\frac{n!}{(n-k+1)!(k-1)!}=\frac{(n+1)!}{(n+1-k)!k!}
$$
 \begin{tabbing}
\hspace{25ex}\=\hspace{3ex}\=     \kill
 \> $\Updownarrow$ \> {\scriptsize (делим
на $n!$)}
 \end{tabbing}
$$
\frac{1}{(n-k)!k!}+\frac{1}{(n-k+1)!(k-1)!}=\frac{n+1}{(n+1-k)!k!}
$$
 \begin{tabbing}
\hspace{25ex}\=\hspace{3ex}\=     \kill
 \> $\Updownarrow$ \> {\scriptsize (умножаем на $k!$)}
 \end{tabbing}
$$
\frac{1}{(n-k)!}+\frac{k}{(n-k+1)!}=\frac{n+1}{(n+1-k)!}
$$
 \begin{tabbing}
\hspace{25ex}\=\hspace{3ex}\=     \kill
 \> $\Updownarrow$ \> {\scriptsize (умножаем на $(n+1-k)!$)}
 \end{tabbing}
$$
n-k+1+k=n+1
$$
\epr

\end{multicols}\noindent\rule[10pt]{160mm}{0.1pt}

\paragraph{Доказательство формул по индукции.}

Здесь мы покажем, как с помощью математической индукции доказываются разные
формулы. Те, что мы докажем, понадобятся нам в дальнейшем.

\noindent\rule{160mm}{0.1pt}
\begin{multicols}{2}

\bex {\bf Неравенство Бернулли.} {\it Для любых $n\in \mathbb{N}$ и
$\alpha > -1$ выполняется неравенство}
\begin{equation}
(1+\alpha)^n \ge 1+ n\cdot \alpha \label{nerav-Bernoulli}
\end{equation}
\begin{proof} Действительно, при $n=1$ получается
$$
(1+\alpha)^1 \ge 1+ 1\cdot \alpha,
$$
то есть в этом случае неравенство верно. Предположим, что оно
верно для какого-нибудь $n=k$, то есть, что справедливо
\begin{equation}
(1+\alpha)^k \ge 1+ k\cdot \alpha \label{0.3.3}
\end{equation}
Тогда для $n=k+1$ мы получаем
 \begin{multline*}
(1+\alpha)^n= (1+\alpha)^{k+1} = (1+\alpha)\cdot (1+\alpha)^k \ge\\
\ge (\text{применяем \eqref{0.3.3}}) \ge (1+\alpha)\cdot (1+ k\cdot
\alpha)=\\=1 + (k+1)\cdot \alpha + k\cdot \alpha^2 \ge 1 + (k+1)\cdot \alpha =
1 + n\cdot \alpha
\end{multline*}

Мы получили, что из того, что \eqref{nerav-Bernoulli} выполняется для какого-то
$n=k$ автоматически следует, что оно выполняется для $n=k+1$. Таким образом,
\eqref{nerav-Bernoulli} обладает свойствами $(a')$ и $(b')$ теоремы
\ref{induction}, и поэтому оно выполняется для любых $n\in \mathbb{N}$.
\end{proof} \eex

 \bex При $n\in\Z_+$ справедливо неравенство
 \beq\label{n!-ge-2^(n-1)}
n!\ge 2^{n-1}
 \eeq
 \eex
 \bpr
1. Докажем сначала это неравенство для $n\in\N$. При $n=1$ это верно:
$$
\underbrace{1!}_{1}\le \underbrace{2^0}_{1}
$$
Предположим, что это верное при $n=k$:
 $$
k!\ge 2^{k-1}
 $$
Тогда при $n=k+1$ получаем:
$$
k\in\N
$$
$$
\Downarrow
$$
$$
k\ge 1
$$
$$
\Downarrow
$$
$$
n=k+1\ge 2
$$
$$
\Downarrow
$$
 \begin{multline*}
n!=(k+1)!=\underbrace{(k+1)}_{\scriptsize\begin{matrix}\text{\rotatebox{90}{$\le$}}\\
2 \end{matrix}}\cdot
\underbrace{k!}_{\scriptsize\begin{matrix}\text{\rotatebox{90}{$\le$}}\\
2^{k-1} \end{matrix}}\ge\\ \ge 2 \cdot 2^{k-1}=2^k=2^{n-1}
 \end{multline*}

2. После того, как \eqref{n!-ge-2^(n-1)} доказано для $n\in\N$, остается
заметить, что при $n=0$ оно тоже верно:
$$
0!=1\ge\frac{1}{2}=2^{-1}.
$$
Следовательно, оно верно для $n\in\Z_+$.
 \epr

 \bex При $n,m\in\Z_+$ справедливо неравенство
 \beq\label{(n+m)!-ge-n!(n+1)^m}
(n+m)!\ge n!\cdot (n+1)^m
 \eeq
 \eex
 \bpr
Проведем индукцию по $m\in\Z_+$. При $m=0$ это неравенство верно:
$$
(n+0)!=n!\ge n!\cdot (n+1)^0
$$
Предположим, что оно верно при $m=k$:
 $$
(n+k)!\ge n!\cdot (n+1)^k
 $$
Тогда при $m=k+1$ получаем:
 \begin{multline*}
(n+m)!=(n+k+1)!=\\=\underbrace{(n+k+1)}_{\scriptsize\begin{matrix}\text{\rotatebox{90}{$\le$}}\\
n+1 \end{matrix}}\cdot \underbrace{(n+k)!}_{\scriptsize\begin{matrix}\text{\rotatebox{90}{$\le$}}\\
n!\cdot (n+1)^k \end{matrix}}\ge (n+1)\cdot n!\cdot (n+1)^k=\\=n!\cdot
(n+1)^{k+1}=n!\cdot (n+1)^m
 \end{multline*}

 \epr

\bex {\bf Сумма отрезка геометрической прогрессии.} {\it Для всякого
вещественного числа $q\ne 1$ и любых чисел $m,n\in\Z_+$, $m\le n$, справедливы
формулы:}
 \begin{align}
&\sum_{i=0}^n q^i =\frac{1-q^{n+1}}{1-q}, \label{Geom-progr-N}\\
&\sum_{i=m}^n q^i =\frac{q^m-q^{n+1}}{1-q} \label{Geom-progr-m<n}
 \end{align}
 \eex
\bpr 1. Докажем сначала \eqref{Geom-progr-N}. При $n=0$ формула верна:
$$
\sum_{i=0}^0 q^i=q^0=1=\frac{1-q^1}{1-q}.
$$
Предположим, что формула верна при $n=k$:
 \beq\label{ind-dlya-geom-progr}
\sum_{i=0}^k q^i =\frac{1-q^{k+1}}{1-q},
 \eeq
Тогда при $n=k+1$ получаем:
 \begin{multline*}
\sum_{i=0}^{k+1} q^i =\sum_{i=0}^k q^i+q^{k+1}=\eqref{ind-dlya-geom-progr}=\\=
\frac{1-q^{k+1}}{1-q}+q^{k+1}=\frac{1-q^{k+1}+q^{k+1}-q^{k+2}}{1-q}=\\=
\frac{1-q^{k+2}}{1-q},
 \end{multline*}
что и есть формула \eqref{Geom-progr-N} при $n=k+1$.

2. После того, как \eqref{Geom-progr-N} доказано, формула
\eqref{Geom-progr-m<n} получается вычислением:
 \begin{multline*}
\sum_{i=m}^n q^i =\sum_{i=0}^n q^i -\sum_{i=0}^{m-1} q^i =\\=
\frac{1-q^{n+1}}{1-q}-\frac{1-q^m}{1-q}=\frac{q^m-q^{n+1}}{1-q}
 \end{multline*}
\epr

\bex {\bf Бином Ньютона.} {\it Для любых $a,b\in\R$ и $n\in\N$
справедливо равенство}
 \beq\label{binom-Newtona}
(a+b)^n=\sum_{k=0}^n C^k_n \cdot a^{n-k} \cdot b^k
 \eeq
где $C_n^k$ -- число сочетаний, определенное выше формулой \eqref{def-C_n^k}.
 \bpr Сначала проверяем ее при $n=1$:
$$
\underbrace{(a+b)^1}_{\scriptsize\begin{matrix}\| \\ a+b \end{matrix}} =
\underbrace{\sum_{k=0}^1 C^k_1 \cdot a^{1-k} \cdot
b^k}_{\scriptsize\begin{matrix}\| \\ C^0_1 \cdot a^1 \cdot b^0+C_1^1\cdot
a^0\cdot b^1 \end{matrix}}
$$
Затем предполагаем, что она верна при $n=m$
 \beq\label{binom-Newtona-m}
(a+b)^m=\sum_{k=0}^m C^k_m \cdot a^{m-k} \cdot b^k
 \eeq
и проверяем, что тогда она будет верна при $n=m+1$:
 \begin{multline*}
(a+b)^{m+1}=(a+b)\cdot (a+b)^m
=\eqref{binom-Newtona-m}=\\=(a+b)\cdot\sum_{k=0}^m C^k_m \cdot a^{m-k} \cdot
b^k=\\=\sum_{k=0}^m C^k_m \cdot a^{m+1-k} \cdot b^k+\underbrace{\sum_{k=0}^m
C^k_m \cdot a^{m-k} \cdot b^{k+1}}_{\scriptsize\text{заменяем $k$ на
$i-1$}}=\\= \sum_{k=0}^m C^k_m \cdot a^{m+1-k} \cdot
b^k+\underbrace{\sum_{i=1}^{m+1} C^{i-1}_m \cdot a^{m+1-i} \cdot
b^i}_{\scriptsize\text{заменяем $i$ на $k$}}=\\= \sum_{k=0}^m C^k_m \cdot
a^{m+1-k} \cdot b^k+\sum_{k=1}^{m+1} C^{k-1}_m \cdot a^{m+1-k} \cdot b^k=\\=
a^{m+1}+\sum_{k=1}^m C^k_m \cdot a^{m+1-k} \cdot b^k +\\+\sum_{k=1}^m C^{k-1}_m
\cdot a^{m+1-k} \cdot b^k+b^{m+1}=\\= a^{m+1}+\sum_{k=1}^m
\underbrace{(C^k_m+C^{k-1}_m)}_{\scriptsize\begin{matrix}\| \\ \phantom{,} C_{n+1}^k, \\
\text{по тождеству}\\ \text{Паскаля \eqref{C_n^k+C_n^k-1=C_n+1^k}}
\end{matrix}} \cdot a^{m+1-k} \cdot b^k +b^{m+1}=\\=a^{m+1}+
\sum_{k=1}^m C^k_{m+1} \cdot a^{m+1-k} \cdot b^k+b^{m+1}= \\=\sum_{k=0}^{m+1}
C^k_{m+1} \cdot a^{m+1-k} \cdot b^k
 \end{multline*}
\epr \eex

\ber{\bf Сумма отрезка арифметической прогрессии.} Докажите  методом
математической индукции формулу:
 \beq\label{arifm-progressija}
\sum_{i=1}^n i=\frac{n\cdot (n+1)}{2}
 \eeq
\eer

\ber Докажите следующие тождества для двойного факториала:
 \begin{align}
& n!!=n\cdot (n-2)!!, && n\ge 2 \label{n!!=n(n-2)!!} \\
& n!=n!!\cdot (n-1)!!, && n\in\N \label{n!=n!!(n-1)!!} \\
& (2m)!!=2^m\cdot m!, && m\in\Z_+ \label{chetnyi-dvoinoi-faktorial}\\
& (2m-1)!!=\frac{(2m-1)!}{2^{m-1}\cdot (m-1)!}, && m\in\N.
\label{nechetnyi-dvoinoi-faktorial}
 \end{align}
\eer

\begin{ers}
Докажите методом математической индукции:

 \biter{
\item[1)] $1^2+3^2+...+(2n-1)^2=\frac{n(4n^2-1)}{3}$

\item[2)] $\frac{1}{n+1}+\frac{1}{n+2}+...+\frac{1}{3n+1}>1$

\item[3)] $\sum_{k=1}^n k (3k-1)=n^2(n+1)$

\item[4)] $\frac{1}{n+1}+\frac{1}{n+2}+...+\frac{1}{2n}>\frac{11}{24}$

\item[5)] $\sum_{k=1}^n \frac{1}{(4k-3)\cdot (4k+1)}=\frac{n}{4n+1}$

\item[6)] $\prod_{k=1}^n \left( 1-\frac{1}{(k+1)^2} \right)=\frac{n+2}{2n+2}$

\item[7)] $1\cdot 2^2+2\cdot 3^2+...+(n-1)\cdot n^2= \frac{n\cdot (n^2-1)\cdot
(3n+2)}{12} \quad (n\ge 2)$,

\item[8)] $\frac{n}{2}<\sum_{k=0}^{2n-1} \frac{1}{k} < n$.

 }\eiter
\end{ers}

\end{multicols}\noindent\rule[10pt]{160mm}{0.1pt}

\paragraph{Индуктивная сумма $\sum_{k=1}^n a_k$ и индуктивное произведение $\prod_{k=1}^n a_k$}
элементов последовательности $\{a_n\}$ определяются индуктивными правилами так:
$$
\sum_{k=1}^1 a_k=a_m,\qquad \sum_{k=1}^{n+1} a_k=\left(\sum_{k=1}^n
a_k\right)+a_{n+1}
$$
$$
\prod_{k=1}^1 a_k=a_m,\qquad \prod_{k=1}^{n+1} a_k=\left(\prod_{k=1}^n
a_k\right)\cdot a_{n+1}
$$
Часто приходится рассматривать также суммы и произведения по интервалам
натуральных чисел $\{m,...,n\}$. Читатель может по аналогии дать индуктивные
определения этим понятиям, мы же воспользуемся уже записанными формулами, и
положим
$$
\sum_{k=m}^n a_k=\sum_{k=1}^{n} a_k-\sum_{k=1}^{m-1} a_k
$$
$$
\prod_{k=m}^n a_k=\prod_{k=1}^{n} a_k\Big/\prod_{k=1}^{m-1} a_k\qquad
(\text{если $a_k\ne 0$}).
$$

\noindent\rule{160mm}{0.1pt}\begin{multicols}{2}

\ber Докажите индукцией следующие тождества:
 \begin{align}
&  a^n=\prod_{k=1}^n a, \qquad a\in\R,\ n\in\N \\
&  n!=\prod_{k=1}^n k, \qquad  n\in\N \label{n!=prod} \\
& (2m)!!=\prod_{k=1}^m(2k), \qquad m\in\Z_+ \\
& (2m-1)!!=\prod_{k=1}^m(2k-1), \qquad m\in \N
 \end{align}
\eer

\ber Покажите индукцией, что для любой последовательности чисел $x_1,...,x_n$
выполняются соотношения:
 \begin{align}
& \left|\sum_{k=1}^n x_k\right|\le\sum_{k=1}^n |x_k|,
\label{ind-nerav-s-modulem} \\
& \left|\prod_{k=1}^n x_k\right|=\prod_{k=1}^n |x_k|
 \end{align}
\eer

\end{multicols}\noindent\rule[10pt]{160mm}{0.1pt}

\paragraph{Принцип Архимеда.}

\begin{tm}\label{N-unbounded}
Множество натуральных чисел $\mathbb{N}$ не ограничено сверху.
\end{tm}
\begin{proof} Предположим, что наоборот,
$\mathbb{N}$ ограничено сверху. Тогда по теореме \ref{sup-inf}, $\mathbb{N}$
должно иметь точную верхнюю грань:
\begin{equation}
\exists \sup \mathbb{N} = B \in \R \label{0.4.1}
\end{equation}
Поскольку $B$ -- точная верхняя грань (то есть наименьшее из чисел,
ограничивающих $\mathbb{N}$ сверху), число $B-1$ не может ограничивать
$\mathbb{N}$ сверху. Значит, существует какое-то $n\in \mathbb{N}$, большее чем
$B-1$: $n>B-1$. То есть
$$
B<n+1
$$
Это означает, что $B$ не может ограничивать $\mathbb{N}$ сверху (потому что $B$
меньше некоторого числа $n+1\in \mathbb{N}$). Таким образом, мы получили
противоречие с \eqref{0.4.1}. Оно означает, что наше предположение о том, что
$\mathbb{N}$ ограничено сверху неверно.
\end{proof}

Из теоремы \ref{N-unbounded} следует важный вывод:

\begin{tm}[\bf принцип Архимеда]\label{Archimed-principle}
\index{принцип!Архимеда} Для любого вещественного числа $C\in \R$ найдется
натуральное число $n\in \mathbb{N}$ такое, что $n>C$.
\end{tm}
\begin{proof} Зафиксируем число $C\in
\R$, и предположим противное, то есть что для него не существует такого $n\in
\mathbb{N}$, чтобы $n>C$. Тогда для всякого $n\in \mathbb{N}$ мы получаем $n\le
C$. Иными словами, $C$ ограничивает $\mathbb{N}$ сверху. Этого не может быть по
теореме \ref{N-unbounded}. Значит, наше предположение неверно.
\end{proof}

\subsection{Конечные множества}

 \biter{\label{DEF:finite-sets}
\item[$\bullet$] Множество $A$ (необязательно числовое) мы будем называть {\it
конечным}, если можно построить биективное отображение $A$ на какой-нибудь
отрезок $\{1,...,n\}=\{k\in\N:\ 1\le k\le n\}$ натурального ряда $\N$, то есть
существует отображение
$$
f:\{1,...,n\}\to A
$$
такое что выполняются два условия:
 \bit{
\item[1)] $\forall k,l\in \{1,...,n\}\quad \Big(k\ne l\quad\Rightarrow\quad
f(k)\ne f(l)\Big)$;

\item[2)] $\forall a\in A\quad \exists k\in \{1,...,n\}\quad f(k)=a$.
 }\eit

Пустое множество $\varnothing$ также считается конечным, потому что пустое
отображение
$$
F:\varnothing\to\varnothing
$$
о котором мы говорили в примере \ref{pustoe-otobrazhenie}, можно считать
биекцией между $\varnothing$ и тем же $\varnothing$, но рассматриваемым как
пустой отрезок натурального ряда:
$$
\varnothing=\{k\in\N:\ 1\le k\le 0\}
$$
 }\eiter

\noindent\rule{160mm}{0.1pt}\begin{multicols}{2}

\paragraph{Принцип Дирихле и его следствия.}

\blm[\bf принцип Дирихле] Ни при каком $n\in\N$ отображение
$f:\{1,...,n+1\}\to\{1,...,n\}$ не может быть инъективным. \elm
 \bpr Проведем индукцию по $n$. При $n=1$ мы получаем отображение
 $$
f:\{1,2\}\to\{1\},
 $$
которое переводит два разных числа 1 и 2 в одно и то же число 1, и поэтому оно
не будет инъективным. Предположим, что наша лемма верна при некотором $n=k$,
где $k\in\N$: никакое отображение
$$
f:\{1,...,k+1\}\to\{1,...,k\}
$$
не может быть инъективным. Покажем, что тогда это будет верно и для $n=k+1$.
Пусть
$$
f:\{1,...,k+2\}\to\{1,...,k+1\}
$$
-- произвольное отображение. Рассмотрим несколько случаев.

1. Предположим сначала, что
 \beq\label{f(k+2)=k+1}
f(k+2)=k+1
 \eeq
Тогда:
 \biter{
\item[---] либо
$$
f\Big(\{1,...,k+1\}\Big)\subseteq\{1,...,k\}
$$
то есть $f$  переводит множество $\{1,...,k+1\}$ в множество $\{1,...,k\}$, и
тогда по предположению индукции,
$$
\exists i\ne j\in \{1,...,k+1\}\quad f(i)=f(j)
$$
\item[---] либо
$$
f\Big(\{1,...,k+1\}\Big)\nsubseteq\{1,...,k\}
$$
то есть
$$
\exists i\in \{1,...,k+1\}\quad
f(i)=\underbrace{k+1}_{\scriptsize\begin{matrix}\|
\\ \eqref{f(k+2)=k+1}\\ \| \\ f(k+2)\end{matrix}},
$$
и мы опять получаем, что $f$ не инъективно:
$$
i\ne k+2\quad\&\quad f(i)=f(k+2).
$$
 }\eiter

2. Предположим, что наоборот,
 \beq\label{f(k+2)-ne-k+1}
f(k+2)\ne k+1
 \eeq
Тогда
$$
f(k+2)\in\{1,...,k\}
$$
Рассмотрим отображение
 \begin{multline*}
\sigma:\{1,...,k+1\}\to\{1,...,k+1\}, \\
\sigma(j)=\begin{cases}k+1,& j=f(k+2) \\ f(k+2),& j=k+1 \\ j,& j\notin
\{f(k+2),k+1\}
\end{cases},
 \end{multline*}
и пусть $g$ -- композиция отображений $f$ и $\sigma$:
 \begin{multline*}
g=\sigma\circ f:\{1,...,k+2\}\to\{1,...,k+1\}, \\
g(i)=\sigma\big(f(i)\big), \quad i\in\{1,...,k+2\}
 \end{multline*}
Тогда
 \beq\label{g(k+2)=k+1}
g(k+2)=\sigma\big(f(k+2)\big)=k+1,
 \eeq
и мы получаем два случая:
 \biter{
\item[---] либо
$$
g\Big(\{1,...,k+1\}\Big)\subseteq\{1,...,k\}
$$
то есть $g$  переводит множество $\{1,...,k+1\}$ в множество $\{1,...,k\}$, и
тогда по предположению индукции,
$$
\exists i\ne j\in \{1,...,k+1\}\quad
\underset{\scriptsize\begin{matrix}\|\\
\sigma\big(f(i)\big)\end{matrix}}{g(i)}=
\underset{\scriptsize\begin{matrix}\|\\
\sigma\big(f(j)\big)\end{matrix}}{g(j)}
$$
$$
\phantom{\scriptsize\text{$\sigma$ --
биекция}}\quad\Downarrow\quad{\scriptsize\text{$\sigma$ -- биекция}}
$$
$$
\exists i\ne j\in \{1,...,k+1\}\quad f(i)=f(j)
$$

\item[---] либо
$$
g\Big(\{1,...,k+1\}\Big)\nsubseteq\{1,...,k\}
$$
то есть
$$
\exists i\in \{1,...,k+1\}\quad
g(i)=\underbrace{k+1}_{\scriptsize\begin{matrix}\|
\\ \eqref{g(k+2)=k+1}\\ \| \\ g(k+2)\end{matrix}},
$$
$$
\Downarrow
$$
$$
\exists i\ne k+2\quad\&\quad
\underset{\scriptsize\begin{matrix}\|\\
\sigma\big(f(i)\big)\end{matrix}}{g(i)}=
\underset{\scriptsize\begin{matrix}\|\\
\sigma\big(f(k+2)\big)\end{matrix}}{g(k+2)}
$$
$$
\phantom{\scriptsize\text{$\sigma$ --
биекция}}\quad\Downarrow\quad{\scriptsize\text{$\sigma$ -- биекция}}
$$
$$
\exists i\ne k+2\quad\&\quad f(i)=f(k+2),
$$
и мы опять получаем, что $f$ не инъективно.
 }\eiter

 \epr

\blm\label{LM:(1,...,k)->(1,...,l)-injec} Если $k,l\in\N$ и существует инъекция
$$
f:\{1,...,k\}\to\{1,...,l\},
$$
то $k\le l$. \elm
 \bpr
Предположим, что $k>l$. Тогда по свойству \eqref{m<n=>m+1-le-n},
$$
l+1\le k
$$
$$
\Downarrow
$$
$$
\{1,...,l+1\}\subseteq \{1,...,k\}
$$
$$
\Downarrow
$$
$$
\exists \ \text{инъекция} \ f|_{\{1,...,l+1\}}:\{1,...,l+1\}\to \{1,...,l\}
$$
 \epr

\blm\label{LM:(1,...,k)->(1,...,l)-bijec} Если $k,l\in\N$ и существует биекция
$$
f:\{1,...,k\}\to\{1,...,l\},
$$
то $k=l$. \elm
 \bpr
Поскольку $f$ является инъекцией, мы по лемме
\ref{LM:(1,...,k)->(1,...,l)-injec} сразу получаем:
$$
k\le l
$$
Обратное отображение
$$
f^{-1}:\{1,...,l\}\to\{1,...,k\},
$$
тоже будет инъекцией, поэтому опять по лемме
\ref{LM:(1,...,k)->(1,...,l)-injec} мы получаем:
$$
l\le k
$$
Вместе то и другое означает, что $k=l$.
 \epr

\paragraph{Мощность конечного множества.}

\btm\label{TH:chislo-elem-konech-mnozh} Если $X$ -- конечное множество, то
число $n\in\Z_+$, для которого можно построить биекцию $f:\{1,...,n\}\to X$,
определяется однозначно.
 \etm

 \biter{
\item[$\bullet$] Если множество $X$ конечно, то число $n$, для которого можно
построить биекцию $f:\{1,...,n\}\to X$ (оно будет единственно по теореме
\ref{TH:chislo-elem-konech-mnozh}), называется {\it мощностью} (или {\it числом
элементов}) множества $X$ и обозначается символом $\card X$:
 \begin{multline}\label{DEF:card(X)}
\card X=n\quad\Longleftrightarrow \\
\Longleftrightarrow\quad \exists \
\text{биекция}\ f:\{1,...,n\}\to X
 \end{multline}
Тот факт, что $X$ конечно, коротко записывают формулой
 $$
\card X<\infty
 $$
Наоборот, запись
 $$
\card X=\infty
 $$
означает, что множество $X$ не является конечным (такие множества называются бесконечными, и мы поговорим о них ниже на с.\pageref{DEF:besk-mnozhestva}).
 }\eiter

 \bpr[Доказательство теоремы \ref{TH:chislo-elem-konech-mnozh}]
Пусть $m\in\Z_+$ и $n\in\Z_+$ -- два числа, такие, что существуют биекции
$$
f:\{1,...,m\}\to X,\qquad g:\{1,...,n\}\to X
$$
Тогда отображение
$$
g^{-1}\circ f:\{1,...,m\}\to \{1,...,n\}
$$
будет тоже биекцией, и поэтому по лемме \ref{LM:(1,...,k)->(1,...,l)-bijec},
$m=n$.
 \epr

\paragraph{Конечные множества в $\R$.}

\btm\label{TH:max-finite-set} Всякое конечное множество $X$ в $\R$ обладает
минимумом и максимумом. \etm
 \bpr
Покажем, что $X$ имеет минимумом (существование максимума доказывается по
аналогии). Проведем индукцию по мощности $\card X$ (числу элементов) множества
$X$.

1. При $\card X=1$ мы получаем, что $X$ состоит всего из одного элемента $a$, и
ясно, что $a$ будет минимумом $X$:
$$
a=\min X
$$

2. Предположим, что мы доказали наше утверждение для всех множеств $X$ мощности
$\card X=k$, где $k\in\N$.

3. Пусть $X$ -- множество мощности $\card X=k+1$. Рассмотрим какую-нибудь
биекцию $f:\{1,...,n\}\to X$ и обозначим $K=\{f(i);\; i\in\{1,...,k\}\}$.
Тогда:
$$
X=K\cup\{f(k+1)\},
$$
Рассмотрим число
$$
a=(\min K)\wedge f(k+1)
$$
(напомним, что мы определили операцию $\wedge$ формулой \eqref{DEF:a-wedge-b}).
Это число $a$ будет минимумом множества $X$, потому что, во-первых,
$$
a=(\min K)\wedge f(k+1)=\min\{\min K, f(k+1)\}
$$
$$
\Downarrow
$$
$$
a=\min K\quad \V\quad a=f(k+1)
$$
$$
\Downarrow
$$
$$
a\in K\quad \V\quad a\in\{f(k+1)\}
$$
$$
\Downarrow
$$
$$
a\in K\cup \{f(k+1)\}=X
$$
И, во-вторых,
$$
a=(\min K)\wedge f(k+1)=\min\{\min K, f(k+1)\}
$$
$$
\Downarrow
$$
$$
a\le\min K\quad \&\quad a\le f(k+1)
$$
$$
\Downarrow
$$
$$
a\le K\quad \&\quad a\le f(k+1)
$$
$$
\Downarrow
$$
$$
a\le K\cup \{f(k+1)\}=X
$$
 \epr

\paragraph{Конечные множества в $\N$.}

\btm\label{TH:ogran-mnozh-v-N-konechno} Множество $A$ в $\N$ конечно тогда и
только тогда, когда оно ограничено:
$$
\card A<\infty\quad\Longleftrightarrow\quad \exists N\in\N\quad
A\subseteq\{1,...,N\}
$$
\etm
 \bpr
Первая половина этого утверждение следует из теоремы \ref{TH:max-finite-set}:
если множество $A$ конечно, то оно либо пусто, и тогда ограничено, либо
непусто, и тогда оно является образом $\Supp(f)$ некоторого инъективного
отображения $f:\{1,...,n\}\to\N$. Отсюда по теореме \ref{TH:max-finite-set}
получаем, что $A$ имеет максимум, и поэтому ограничено сверху. С другой
стороны, оно всегда ограничено снизу единицей 1, и значит оно просто
ограничено.

Докажем, что, наоборот, если $A$ ограничено, то оно конечно. Опять, если $A$
пусто, то доказывать ничего не нужно, поэтому мы будем считать, что
$A\ne\varnothing$. Тогда, положив
$$
a_1=\min A,
$$
можно определить отображение $G$, которое некоторым последовательностям
$\{a_1,...,a_k\}\subseteq A$ будет ставить в соответствие число
$$
G(a_1,...,a_k)=\min A\setminus\{a_1,...,a_k\}
$$
(если множество $A\setminus\{a_1,...,a_k\}$ непусто, то такое число существует
по свойству $4^\circ$ на с.\pageref{minimum-v-N}). По теореме
\ref{defin-induction-chast} об определениях частной индукцией, отображение $G$
определяет некую последовательность $\{a_n\}\subseteq A$, удовлетворяющую
условию
$$
G(a_1,...,a_n)=a_{n+1},
$$
для всех $n$, меньших длины $N$ этой последовательности.

1. Докажем, что она строго возрастает:
 \beq\label{a_n<a_(n+1)-konech-mnozh}
\forall n<N\qquad a_n<a_{n+1}
 \eeq
Действительно, во-первых,
 $$
\{1,...,a_{n-1}\}\subseteq \{1,...,a_{n-1},a_n\}
 $$
 $$
 \Downarrow
 $$
 $$
A\setminus \{1,...,a_{n-1}\}\supseteq A\setminus\{1,...,a_{n-1},a_n\}
 $$
 $$
\phantom{\scriptsize
\eqref{X-subseteq-Y=>min(X)-ge-min(Y)}}\quad\Downarrow\quad{\scriptsize
\eqref{X-subseteq-Y=>min(X)-ge-min(Y)}}
 $$
 \begin{multline*}
a_n=\min A\setminus \{1,...,a_{n-1}\}\le
\\ \le \min A\setminus\{1,...,a_{n-1},a_n\}=a_{n+1}
 \end{multline*}
И, во-вторых, если бы оказалось, что $a_n=a_{n+1}$, то мы получили бы
$$
a_n=a_{n+1}=\min A\setminus\{1,...,a_n\}
$$
 $$
 \Downarrow
 $$
$$
a_n\in A\setminus\{1,...,a_n\},
$$
что невозможно.

2. Из \eqref{a_n<a_(n+1)-konech-mnozh} следует условие
 \beq\label{n-le-a_(n+1)-konech-mnozh}
\forall n\qquad n\le a_n
 \eeq
Его можно доказать индукцией. При $n=1$ оно очевидно:
$$
1\le a_1.
$$
А если оно верно при $n=k$
$$
k\le a_k,
$$
то из этого следует
$$
k\le a_k<a_{k+1}
$$
 $$
\phantom{\scriptsize \eqref{m<n=>m+1-le-n}}\quad\Downarrow\quad{\scriptsize
\eqref{m<n=>m+1-le-n}}
 $$
$$
k+1\le a_{k+1}
$$
То есть \eqref{n-le-a_(n+1)-konech-mnozh} будет верно и при $n=k+1$.

3. Заметим теперь, что поскольку $A$ ограничено, существует $c\in\R$ такое что
$$
A\le c
$$
По принципу Архимеда (теорема \ref{Archimed-principle}), должно существовать
$m\in\N$ такое что
$$
A\le c<m
$$
Отсюда следует, что последовательность $\{a_n\}$ конечна. Действительно, если
бы она была бесконечной, то мы получили бы, что ее элемент с номером $m$ лежит
в $A$
$$
a_m\in A
$$
А с другой стороны, это невозможно, потому что
$$
A<m\le a_m
$$

4. Итак, мы имеем конечную последовательность $\{a_n\}\subseteq A$, то есть
отображение
$$
f:\{1,...,N\}\to A,\qquad f(n)=a_n
$$
где $N\in\N$ -- длина последовательности. Это отображение инъективно, потому
что
$$
m<n$$
 $$
 \Downarrow
 $$
$$
f(m)=a_m<a_n=f(n)
$$
 $$
 \Downarrow
 $$
$$
f(m)\ne f(n)
$$
Нам остается проверить, что оно сюръективно. Действительно, если бы оказалось,
что
$$
A\setminus\{a_1,...,a_N\}\ne\varnothing,
$$
то положив
$$
a_{N+1}=\min A\setminus\{a_1,...,a_N\}
$$
мы получили бы, что $N$ не может быть длиной последовательности $\{a_n\}$,
определяемой отображением $G$ по теореме \ref{defin-induction-chast}:
 \begin{multline*}
\exists \{a_2,...,a_{N+1}\}\subseteq A\\ \forall k\in\{2,...,N+1\} \qquad
G(a_1,...,a_{k-1})=a_k
 \end{multline*}
$$
\Downarrow
$$
 \begin{multline*}
N+1\le \sup\Big\{n\in\N:\quad \exists \{a_2,...,a_n\}\subseteq A\\
\forall k\in\{2,...,n\} \qquad G(a_1,...,a_{k-1})=a_k\Big\}
 \end{multline*}
$$
\Downarrow
$$
 \begin{multline*}
N\ne \sup\Big\{n\in\N:\quad \exists \{a_2,...,a_n\}\subseteq A\\
\forall k\in\{2,...,n\} \qquad G(a_1,...,a_{k-1})=a_k\Big\}
 \end{multline*}
 \epr

\end{multicols}\noindent\rule[10pt]{160mm}{0.1pt}

 \bigskip
\centerline{\bf Свойства конечных множеств:}\label{svoistva-konechnyh-mnozhestv}
 \bit{\it

\item[$1^\circ$.] Если $f:X\to Y$ -- инъективное отображение и $Y$ -- конечное
множество, то $X$ -- тоже конечное множество.

\item[$2^\circ$.] Если $f:X\to Y$ -- сюръективное отображение и $X$ -- конечное
множество, то $Y$ -- тоже конечное множество.

\item[$3^\circ$.]\label{subset-finite-set} Всякое подмножество $X$ любого
конечного множества $Y$ также является конечным множеством.

\item[$4^\circ$.] Пересечение $X\cap Y$ любых двух конечных множеств $X$ и $Y$
тоже является конечным множеством.

\item[$5^\circ$.] Объединение $X\cup Y$ любых двух конечных множеств $X$ и $Y$
тоже является конечным множеством.

 }\eit

\bpr

1. Пусть $f:X\to Y$ -- инъективное отображение. Если $Y$ -- конечное множество,
то существует биекция $g:\{1,...,n\}\to Y$. Рассмотрим обратное отображение
$g^{-1}:Y\to \{1,...,n\}$ и композицию $h=g^{-1}\circ f$.
$$
\begin{diagram}
\node[2]{Y}\arrow{se,t}{g^{-1}}\\
\node{X} \arrow[2]{e,b}{h}\arrow{ne,t}{f} \node[2]{\{1,...,n\}}
\end{diagram}
$$
Множество $A=h(X)$ содержится в дискретном интервале $\{1,...,n\}$, то есть
является ограниченным подмножеством в $\N$. Значит, в силу примера
\ref{TH:ogran-mnozh-v-N-konechno}, оно конечно, то есть существует биекция
$\alpha:\{1,...,k\}\to A$. Поскольку отображение $h:X\to Y$ инъективно, оно
будет биекцией между $X$ и $A=h(X)$. Значит существует обратное отображение
$h^{-1}:A\to X$. Рассмотрим композицию $\beta=h^{-1}\circ\alpha$:
$$
\begin{diagram}
\node[2]{\{1,...,k\}}\arrow{se,t}{\alpha}\arrow{sw,t}{\beta} \\
\node{X}  \node[2]{A}\arrow[2]{w,b}{h^{-1}}
\end{diagram}
$$
Здесь $\alpha$ и $h^{-1}$ -- биекции, значит $\beta$ -- тоже биекция. Поэтому
$X$ -- конечное множество.

2. Пусть $f:X\to Y$ -- сюръективное отображение. Если $X$ -- конечное
множество, то существует биекция $g:\{1,...,n\}\to X$. Композиция этих
отображений $h=f\circ g$ должна быть сюръекцией:
$$
\begin{diagram}
\node[2]{X}\arrow{se,t}{f}\\
\node{\{1,...,n\}} \arrow[2]{e,b}{h}\arrow{ne,t}{g} \node[2]{Y}
\end{diagram}
$$
Рассмотрим множество
$$
A=\Big\{j\in\{1,...,n\}:\ j=\min h^{-1}\big(h(j)\big) \Big\}
$$
Оно является подмножеством в дискретном интервале $\{1,...,n\}$, то есть
ограниченным множеством в $\N$. Значит, в силу примера
\ref{TH:ogran-mnozh-v-N-konechno}, оно конечно, то есть существует биекция
$\alpha:\{1,...,m\}\to A$. Покажем, что композиция $\beta=h|_A\circ\alpha$ тоже
является биекцией:
$$
\begin{diagram}
\node{A} \arrow[2]{e,t}{h|_A} \node[2]{Y}\\
\node[2]{\{1,...,m\}}\arrow{nw,b}{\alpha}\arrow{ne,b}{\beta}
\end{diagram}
$$
Сначала проверим сюръективность:
$$
y\in Y
$$
$$
\phantom{\scriptsize \text{$h$ -- сюръективное
отображение}}\quad\Downarrow\quad{\scriptsize \text{$h$ -- сюръективное
отображение}}
$$
$$
\exists k\in\{1,...,n\}\quad h(k)=y
$$
$$
\Downarrow
$$
$$
h^{-1}(y)=\Big\{k\in\{1,...,n\}:\; h(k)=y\Big\}\ne\varnothing
$$
$$
\Downarrow
$$
$$
\exists j=\min h^{-1}(y)
$$
$$
\phantom{\scriptsize\begin{pmatrix} h^{-1}(y)\owns j\\ \Downarrow \\
y=h(j)\\ \Downarrow \\ h^{-1}(y)=h^{-1}(h(j))\\ \Downarrow \\ j=\min
h^{-1}(y)=\min h^{-1}(h(j))\\ \Downarrow \\ j\in A
\end{pmatrix}}\quad\Downarrow\quad{\scriptsize\begin{pmatrix} h^{-1}(y)\owns j\\ \Downarrow \\
y=h(j)\\ \Downarrow \\ h^{-1}(y)=h^{-1}(h(j))\\ \Downarrow \\ j=\min
h^{-1}(y)=\min h^{-1}(h(j))\\ \Downarrow \\ j\in A
\end{pmatrix}}
$$
$$
\exists i\in \{1,...,m\}\quad \alpha(i)=j
$$
$$
\Downarrow
$$
$$
\exists i\in \{1,...,m\}\quad \beta(i)=h(\alpha(i))=h(j)=y
$$
Теперь инъективность. Пусть $i,j\in\{1,...,m\}$ и $\beta(i)=\beta(j)$.
Обозначим $y=\beta(i)=\beta(j)$, тогда
$$
h(\alpha(i))=\beta(i)=y=\beta(j)=h(\alpha(j))
$$
$$
\Downarrow
$$
$$
\alpha(i)\in h^{-1}(y)\owns \alpha(j)
$$
$$
\phantom{\scriptsize \text{$\alpha(i),\alpha(j)\in
A$}}\quad\Downarrow\quad{\scriptsize \text{$\alpha(i),\alpha(j)\in A$}}
$$
$$
\alpha(i)=\min h^{-1}(y)=\alpha(j)
$$
$$
\phantom{\scriptsize \text{$\alpha$ -- инъективное
отображение}}\quad\Downarrow\quad{\scriptsize \text{$\alpha$ -- инъективное
отображение}}
$$
$$
i=j
$$

3. Если $X$ -- подмножество в конечном множестве $Y$, то тождественное
отображение
$$
f:X\to Y,\qquad f(x)=x
$$
будет инъекцией, поэтому в силу уже доказанного свойства $1^\circ$, $X$ должно
быть конечным.

4. Пусть $X$ и $Y$ -- конечные множества. Их пересечение $X\cap Y$ будет
подмножеством в $X$ (и в $Y$), поэтому по только что доказанному свойству
$3^\circ$, $X\cap Y$ должно быть конечным.

5. Пусть $X$ и $Y$ -- конечные множества. Тогда существуют биекции
$f:\{1,...,m\}\to X$ и $g:\{1,...,n\}\to Y$. Рассмотрим отображение
$$
h:\{1,...,m+n\}\to X\cup Y,\qquad h(k)=\begin{cases}f(k), & k\le m
\\ g(k-m), & k>m \end{cases}
$$
Это будет сюръекция, потому что
$$
y\in X\cup Y\quad\Longrightarrow\quad \left[\begin{matrix}y\in
X\quad\Longrightarrow\quad \exists i\in\{1,...,m\} \quad h(i)=y \\
y\in Y\quad\Longrightarrow\quad \exists j\in\{1,...,m\} \quad h(j)=y
\end{matrix}\right]\quad\Longrightarrow\quad \exists k\in\{1,...,m+m\} \quad h(k)=y
$$
Поэтому в силу уже доказанного свойства $2^\circ$, множество $Y$ должно быть
конечно.
 \epr

\subsection{Бесконечные множества}

 \bit{
\item[$\bullet$] Множество $X$ (необязательно числовое) называется {\it бесконечным}\label{DEF:besk-mnozhestva}, если оно не является конечным (в смысле
определения на с.\pageref{DEF:finite-sets}), то есть если не существует биективного отображения $X$ на какой-нибудь
отрезок $\{1,...,n\}=\{k\in\N:\ 1\le k\le n\}$ натурального ряда $\N$.
 }\eit

Напомним, что в примере \ref{EX-posledovatelnosti} выше мы вводили понятие (конечной и бесконечной) последовательности. Теперь, когда множество натуральных чисел $\N$ строго определено, понятие последовательности также можно считать строго определенным. В частности, бесконечная последовательность -- это произвольное отображение $f:\N\to A$ в произвольное множество $A$. Значения такого отображения принято обозначать буквами с нижним индексом:
$$
f(1)=a_1,\quad f(2)=a_2,\quad  f(3)=a_3,\quad ... \qquad (n\in\N)
$$

 \bit{
\item[$\bullet$] Последовательность $\{a_n\}$ называется {\it инъективной}, если ее элементы не повторяются:
$$
 i\ne j\quad\Longrightarrow\quad a_i\ne a_j
$$ }\eit

Следующая теорема называется {\it критерием бесконечности множества}:

\btm\label{TH:krit-besk-mnozh} Множество $X$ бесконечно в том и только в том случае, если в нем найдется инъективная бесконечная последовательность элементов:
$$
\exists \{a_n;\ n\in\N\}\subseteq X:\qquad \forall i,j\in\N\qquad \Big(\ i\ne j\quad\Longrightarrow\quad a_i\ne a_j \ \Big)
$$
\etm
\bpr
1. Достаточность. Пусть в $X$ имеется инъективная бесконечная последовательность, то есть существует инъективное отображение $g:\N\to X$. Если бы $X$ было конечно, то есть если бы существовало биективное отображение $f:\{1,...,n\}\to X$ для некоторого $n\in\N$, мы получили бы, что композиция инъективных отображений $g$ и $f^{-1}$
$$
f^{-1}\circ g:\N\to \{1,...,n\}
$$
тоже является инъективным отображением. Но тогда по свойству $1^\circ$ на с.\pageref{svoistva-konechnyh-mnozhestv} множество $\N$ должно быть конечно, и значит по теореме \ref{TH:ogran-mnozh-v-N-konechno} ограничено. Поскольку это не так, наше предположение о конечности $X$ ложно, и это множество должно быть бесконечно.

2. Необходимость.  Пусть $X$ бесконечно, покажем, что тогда в нем существует инъективная бесконечная последовательность. Здесь наступает момент, когда нам нужно вспомнить о своем старом обещании. На с.\pageref{beskonechnokratnyj-vybor} мы ссылались на теорему \ref{TH:krit-besk-mnozh}, как на пример утверждения, в доказательстве которого используется аксиома выбора (сформулированная нами на с.\pageref{AX:vybora}). Мы говорили также, что обычно доказательство этого и подобных утверждений строится таким образом, что присутствие аксиомы выбора в нем затушевывается, и используемый при этом прием называется {\it бесконечнократным выбором}. Чтобы читатель смог составить себе представление, о чем здесь идет речь, мы пообещали там же привести для сравнения два доказательства, традиционное, в котором аксиома выбора присутствует неявно, и аккуратное, содержащее прямую ссылку на нее. В обоих доказательствах части, в которых обсуждается достаточность, обычно совпадают (если достаточность не считается очевидной), поэтому эту часть мы привели сразу и без комментариев. Различаются только части, посвященные необходимости, и вот как они выглядят в двух вариантах.

\biter{

\item[а)] {\it Доказательство приемом бесконечнократного выбора}\footnote{Здесь приведено доказательство из учебника: А.Н.Колмогоров, С.В.Фомин, Элементы теории функций и функционального анализа. М.: Наука, 1972 (с.21).}. Выберем в $X$ произвольный элемент $a_1$. Поскольку $X$ бесконечно, в нем найдется элемент $a_2$, отличный от $a_1$, затем найдется элемент $a_3$, отличный от $a_1$ и от $a_2$ и так далее. Продолжая эттот процесс (который не может оборваться из-за ``нехватки'' элементов, ибо $X$ бесконечно), мы получаем инъективную последовательность
$$
a_1, a_2, a_3, ...
$$
На этом доказательство заканчивается. Как мы говорили на с.\pageref{beskonechnokratnyj-vybor}, удобство такого стиля состоит в его наглядности. Но неудобством, как может заметить читатель, является то, что остается совершенно непонятно, как можно за конечное время проделать процедуру выбора бесконечное число раз. Тем не менее, прием бесконечнократного выбора вполне применим, и чтобы это увидеть, нужно научиться переводить предлагаемые им конструкции на язык, в котором процедура выбора бесконечное число раз за конечное время отсутствует. Вот как это делается в нашей конкретной ситуации.

\item[б)] {\it Доказательство с явной ссылкой на аксиому выбора и теорему \ref{defin-induction-chast} об определениях частной индукцией}. Поскольку множество $X$ бесконечно, никакая инъективная конечная последовательность его элементов $\{a_1,...,a_n\}$ не может исчерпывать все $X$:
$$
X\setminus \{a_1,...,a_n\}\ne\varnothing
$$
Рассмотрим семейство непустых множеств $X\setminus \{a_1,...,a_n\}$, индексируемое инъективными конечными последовательностями $\{a_1,...,a_n\}$ элементов $X$. По аксиоме выбора найдется отображение $G$, которое каждой конечной инъективной последовательности $\{a_1,...,a_n\}\subseteq X$ ставит в соответствие некоторый элемент
 \beq\label{krit-beskonechnosti}
G(a_1,...,a_n)\in X\setminus \{a_1,...,a_n\}
 \eeq
Зафиксируем теперь какой-нибудь начальный элемент $a_1\in X$. Тогда по теореме \ref{defin-induction-chast} найдется (конечная или бесконечная) последовательность $\{a_n\}\subseteq A$ такая, что при всех $n$, меньших некоторого $N$ (которое либо конечно, либо бесконечно), выполняется
$$
G(a_1,...,a_n)=a_{n+1}
$$
Из условия
$$
a_{n+1}=G(a_1,...,a_n)\in X\setminus \{a_1,...,a_n\}
$$
следует, что такая последовательность обязательно инъективна, и нам остается только убедиться, что величина
$$
N=\sup\Big\{n\in\N:\quad \exists \{a_2,...,a_n\}\subseteq A\quad \forall
k\in\{2,...,n\} \qquad G(a_1,...,a_{k-1})=a_k\Big\}
$$
бесконечна. Действительно, если бы она была конечна, то есть $N\in\N$, это бы означало, что полученную нами инъективную последовательность  $\{a_1,...,a_N\}$ нельзя доопределить, то есть что не существует элемента
$$
a_{N+1}=G(a_1,...,a_N)\in X\setminus \{a_1,...,a_N\}
$$
Это противоречит условию \eqref{krit-beskonechnosti}, которое выполняется для любой инъективной последовательности $\{a_1,...,a_n\}$.
 }\eiter

\epr

\section{Целые числа, делимость и рациональные числа}

\subsection{Целые числа $\Z$}\label{SEC-tselye-chisla}

\paragraph{Определение множества целых чисел.}

 \bit{
\item[$\bullet$] Число $x\in\R$ называется {\it целым}, если оно
 \bit{

\item[--] либо натуральное,
$$
x\in\N
$$
\item[--] либо равно нулю
$$
x=0
$$
\item[--] либо противоположно натуральному
$$
x\in-\N
$$
(то есть $x$ имеет вид $x=-n$, где $n\in\N$).
 }\eit
В соответствии с этим, множество целых чисел, обозначаемое символом $\Z$,
описывается формулой
 \beq\label{DF:Z}
\Z=(-\N)\cup \{0\}\cup \N
 \eeq
 }\eit

\btm\label{Z-zamknuto-otn-+-i-} Множество целых чисел $\Z$ замкнуто
относительно сложения, вычитания и умножения: если $m,n\in\Z$, то $m+n\in\Z$,
$m-n\in\Z$ и $m\cdot n\in\Z$.
 \etm
\bpr 1. Пусть $m,n\in\Z$. Покажем сначала, что $m+n\in\Z$. Для этого придется
рассмотреть несколько случаев.

Во-первых, заметим, что если какое-то из этих чисел равно нулю, например,
$m=0$, то их сумма очевидно, принадлежит $\Z$: $m+n=n\in\Z$. Поэтому будем
считать, что оба они отличны от нуля: $m\ne 0$, $n\ne 0$.

Далее рассматриваем три случая:
 \bit{
\item[---] если $m>0$ и $n>0$, то $m,n\in\N$, и по теореме
\ref{zamknutost-N-otn-+-i-umn} получаем $m+n\in\N\subseteq\Z$;

\item[---] если $m<0$ и $n<0$, то $m,n\in -\N$, то есть $-m\in\N$ и $-n\in\N$ и
опять по теореме \ref{zamknutost-N-otn-+-i-umn} получаем $-m-n\in\N$, откуда
$m+n\in-\N\subseteq\Z$;

\item[---] если $m<0$ и $n>0$, то $m\in-\N$, $n\in\N$, то есть $-m=k\in\N$ и
$n\in\N$; здесь нужно рассмотреть три частных случая:
 \bit{

\item[---] если $k<n$, то по свойству $7^\circ$ на с. \pageref{m<n->n-m-in-N},
$n+m=n-k\in\N\subseteq\Z$;

\item[---] если $k=n$, то $n+m=n-k=0\in\Z$;

\item[---] если $k>n$, то по свойству $7^\circ$ на с. \pageref{m<n->n-m-in-N},
$-(n+m)=k-n\in\N$, откуда $n+m\in-\N\subseteq\Z$;

 }\eit
 }\eit

2. Когда доказано, что $\Z$ замкнуто относительно сложения, замкнутость
относительно вычитания становится очевидной, потому что из самого определения
$\Z$ ясно, что если $n\in\Z$, то $-n\in\Z$. Поэтому при $m,n\in\Z$ получаем
$-n\in\Z$ и можно применить уже доказанную замкнутость относительно сложения:
$m-n=m+(-n)\in\Z$.

3. Замкнутость относительно умножения доказывается так же, как и относительно
сложения. Пусть $m,n\in\Z$. Чтобы убедиться, что $m\cdot n\in\Z$, нужно
рассмотреть несколько случаев.

Во-первых, опять же нужно заметить, что если какое-то из этих чисел равно нулю,
например, $m=0$, то их сумма очевидно, принадлежит $\Z$: $m\cdot n=0\in\Z$.
Поэтому будем считать, что оба они отличны от нуля: $m\ne 0$, $n\ne 0$.

Далее рассматриваем три случая:
 \bit{
\item[---] если $m>0$ и $n>0$, то $m,n\in\N$, и по теореме
\ref{zamknutost-N-otn-+-i-umn} получаем $m\cdot n\in\N\subseteq\Z$;

\item[---] если $m<0$ и $n<0$, то $m,n\in -\N$, то есть $-m\in\N$ и $-n\in\N$ и
опять по теореме \ref{zamknutost-N-otn-+-i-umn} получаем $m\cdot n=(-m)\cdot
(-n)\in\N\subseteq\Z$;

\item[---] если $m<0$ и $n>0$, то $m\in-\N$, $n\in\N$, то есть $-m\in\N$ и
$n\in\N$; опять по теореме \ref{zamknutost-N-otn-+-i-umn} получаем $-m\cdot
n=(-m)\cdot n\in -\N$, откуда $m\cdot n\in-\N\subseteq\Z$.
 }\eit
\epr

\btm\label{TH:m<n=>m-le-n-1} Для любых $m,n\in\Z$
 \beq\label{m<n=>m-le-n-1}
m<n+1\quad\Longrightarrow\quad m\le n
 \eeq
 \etm
 \bpr
Здесь нужно рассмотреть несколько случаев.

1. Если $m\in\N$, то $1\le m$, поэтому $0\le m-1<n$, и, поскольку
$n\in\Z=(-\N)\cup \{0\}\cup \N$, такое возможно только если $n\in\N$. Мы
получаем, что оба числа $m$ и $n$ лежат в $\N$, а этот случай уже рассмотрен в
свойстве $5^\circ$ на с.\pageref{N:m<n=>m-le-n-1}.

2. Если $m=0$, то $0<n+1$, то есть $-1<n$, поэтому либо $n=0$, и тогда $m=0\le
0=n$, либо $n\in\N$, и тогда $n\ge 1>0$, откуда $m=0\le n$.

3. Если $m<0$, то $m\in-\N$, и приходится рассматривать несколько частных
случаев.
 \bit{

\item[---] если $n\in\N$, то мы получаем $m\le 0\le n$.

\item[---] если $n=0$, то $m\le 0=n$.

\item[---] если $n<0$, то $n\in-\N$, и мы получаем
$$
m<n+1\quad\Longrightarrow\quad -m>-n-1\quad\Longrightarrow\quad -n<-m+1;
$$
поскольку при этом $-m\in\N$ и $-n\in\N$, мы опять попадаем в ситуацию,
описанную в свойстве $5^\circ$ на с.\pageref{N:m<n=>m-le-n-1}, и получаем, что
$-n\le -m$. То есть $m\le n$.
 }\eit
 \epr

\bcor\label{COR:n-le-m<n+1} Если $m,n\in\Z$, причем $n\le m<n+1$, то $m=n$.
\ecor

\bpr Из \eqref{m<n=>m-le-n-1} получаем: $m<n+1$ $\Longrightarrow$ $m\le n$, и
вместе с неравенством $n\le m$ это означает, что $m=n$. \epr

\btm\label{max-i-min-Z} Любое непустое ограниченное снизу (сверху) множество
целых чисел $E\subseteq\Z$ имеет минимальный (максимальный) элемент.
$$
\exists \min E \qquad (\max E)
$$
причем
 \beq\label{E<n=>max-E<n}
n<E\quad\Longrightarrow\quad n<\min E \qquad (E<n\quad\Longrightarrow\quad \max
E<n)
 \eeq
 \etm
\bpr Пусть $E$ ограничено сверху. Тогда по принципу Архимеда
\ref{Archimed-principle}, найдется натуральное число $n\in\N$ ограничивающее
$E$ сверху. Теперь получаем:
 \begin{align*}
  E & <n \\
    & \Downarrow  \\
 -n &<-E \\
 & \Downarrow \\
  0 & <-E+n \qquad (\text{причем $-E+n\subseteq \Z$, по теореме \ref{Z-zamknuto-otn-+-i-}}) \\
 & \Downarrow \qquad\text{\scriptsize \eqref{DF:Z}} \\
 -E+n & \subseteq\N \\
 & \Downarrow \qquad (\text{\scriptsize свойство $5^\circ$ на с. \pageref{minimum-v-N}}) \\
  \exists & \min(-E+n)=\min(-E)+n>0 \\
 & \Downarrow \qquad(\text{\scriptsize свойство $2^\circ (i)$ на с. \pageref{min-X+a=min-X+a}}) \\
  \exists & \min(-E)>-n \\
 & \Downarrow\qquad (\text{\scriptsize свойство $2^\circ (ii)$ на с. \pageref{min-X+a=min-X+a}}) \\
  \exists & \max E=-\min(-E)<n
 \end{align*}

 \epr

\noindent\rule{160mm}{0.1pt}\begin{multicols}{2}

\paragraph{Степени с целым показателем.}\label{stepeni-n-in-Z}

Напомним, что на странице \pageref{opr-stepeni} мы уже определили степень $a^n$
с показателем $n\in\Z_+$ индуктивной формулой \eqref{def:a^n}:
$$
a^0=1,\quad a^{n+1}=a^n\cdot a,\qquad n\in\Z_+\qquad{\scriptsize
\eqref{def:a^n}}
$$
(при этом $a$ может быть любым). Для целых отрицательных показателей степень
определяется формулой
 \beq\label{DF:a^n-n<0}
a^{-n}:=\big(\kern-10pt\underbrace{a^{-1}}_{\scriptsize\begin{matrix}\text{обратный}\\
\text{элемент} \\ \text{для $a$} \\ \text{из аксиомы} \\
A9\end{matrix}}\kern-10pt\big)^n=\l\frac{1}{a}\r^n,\qquad n\in\N,
 \eeq
(в этом случае $a$ должно быть ненулевым).

 \biter{
\item[$\bullet$] Отображение $(a,n)\mapsto a^n$ называется {\it степенным
отображением} с целым показателем. Его свойства удобно собрать в следующей
теореме.
 }\eiter

\btm[\bf о степенном отображении]\label{TH-o-step-otobr-malaya} Отображение
$$
(a,n)\mapsto a^n,
$$
обладает следующими свойствами:
 \bit{
\item[$Z_0$:] оно определено в следующих двух ситуациях:
\label{usloviya-sushestv-a^n}
 \biter{
\item[---] при $a\ne 0$ и $n\in\Z$,

\item[---] при $a=0$ и $n\in\Z_+$,
 }\eiter\noindent
или, что то же самое, в следующих двух:
 \biter{
\item[---] при $n\in\Z_+$ и тогда $a\in\R$,

\item[---] при $n\in-\N$ и тогда $a\ne 0$;
 }\eiter\noindent

\item[$Z_1$:] следующие две группы тождеств выполняются всякий раз, когда  обе
части тождества определены:
 \biter{

\item[---] {\bf показательные законы}\index{закон!показательный}:
 \begin{align}
&a^0=1,\label{a^0-Z} \\
&a^{-n}=\frac{1}{a^n}, \label{a^(-n)} \\
&a^{m+n}=a^m\cdot a^n;\label{a^(m+n)}
 \end{align}

\item[---] {\bf степенные законы}\index{закон!степенной}:
 \begin{align}
&1^n=1,\label{1^n} \\
&\left(\frac{1}{a}\right)^n=\frac{1}{a^n},\label{(1/a)^n} \\
&(a\cdot b)^n=a^n\cdot b^n \label{(ab)^n}
 \end{align}

 }\eiter\noindent

\item[$Z_2$:] тождество, называемое {\bf накопительным
законом}\index{закон!накопительный},
 \begin{align}\label{nakop-zakon-Z}
&(a^m)^n=a^{m\cdot n}
 \end{align}
выполняется для следующих значений переменных:
 \biter{
\item[---] при $a\ne 0$ и $m,n\in\Z$,

\item[---] при $a=0$ и $m,n\in\Z_+$;
 }\eiter

\item[$Z_3$:] для {\bf нулевого основания} степень, в случаях, когда она
определена, описывается формулой
 \beq\label{0^n}
0^n=\begin{cases}1,& n=0 \\ 0,& n>0\end{cases}
 \eeq
а для {\bf положительных оснований} выполняется следующее {\bf условие
сохранения знака}:
 \beq\label{x^n>0}
\kern-30pt a>0\quad\Rightarrow\quad a^n>0\quad (n\in\Z)
 \eeq

\item[$Z_4$:] для положительных оснований выполняются следующие {\bf условия
монотонности}:
 \biter{
\item[---] если $n>0$, то возведение в степень $n$ сохраняет знак неравенства:
 \beq\label{monot-a^n-n>0}
\kern-40pt 0<x<y \quad\Longrightarrow\quad x^n<y^n
 \eeq
\item[---] если $n<0$, то возведение в степень $n$ меняет знак неравенства:
 \beq\label{monot-a^n-n<0}
\kern-40pt 0<x<y \quad\Longrightarrow\quad x^n>y^n
 \eeq

\item[---] если $a>1$, то потенцирование с основанием $a$ сохраняет знак
неравенства:
 \beq\label{monot-a^n-a>1}
\kern-40pt m<n \quad\Longrightarrow\quad a^m<a^n
 \eeq

\item[---] если $0<a<1$, то потенцирование с основанием $a$ меняет знак
неравенства:
 \beq\label{monot-a^n-0<a<1}
\kern-40pt m<n \quad\Longrightarrow\quad a^m>a^n
 \eeq
 }\eiter
 }\eit
 \etm
\bpr

1. Начнем с показательных законов.  Формула \eqref{a^0-Z} есть просто часть
определения \eqref{def:a^n} степени $a^n$ с неотрицательным показателем.
Формула \eqref{a^(-n)} доказывается рассмотрением трех случаев: если $n\in\N$,
то
$$
a^{-n}=\eqref{DF:a^n-n<0}=\frac{1}{a^n}
$$
если $n=0$, то
 \begin{multline*}
a^{-n}=a^0=\eqref{a^0-Z}=1=\\=\frac{1}{1}=\eqref{a^0-Z}=\frac{1}{a^0}=\frac{1}{a^n}
 \end{multline*}
а если $n\in-\N$, то есть $n=-k$, $k\in\N$, то поскольку \eqref{a^(-n)} уже
доказано для положительных степеней, должно выполняться $a^{-k}=\frac{1}{a^k}$,
или, что то же самое, $\frac{1}{a^{-k}}=a^k$, поэтому
$$
a^{-n}=a^k=\frac{1}{a^{-k}}=\frac{1}{a^n}
$$

Для доказательства \eqref{a^(m+n)} нам понадобится следующее вспомогательное
тождество:
 \beq\label{a^(m+1)=a^m-a-m-in-Z}
a^{m+1}=a^m\cdot a,\qquad m\in\Z
 \eeq
При $m\in\Z_+$ оно является просто частью определения \eqref{def:a^n}, поэтому
для его доказательства нам нужно лишь проверить, что оно верно при $m\in-\N$,
то есть при $m=-k$, где $k\in\N$:
 $$
a^{1-k}=a^{-k}\cdot a,\qquad k\in\N
 $$
Если сделать замену $k-1=n$, $n\in\Z_+$, то мы получим эквивалентную формулу
 $$
a^{-n}=a^{-n-1}\cdot a,\qquad n\in\Z_+
 $$
которую можно переписать так:
 $$
\frac{a^{-n}}{a}=a^{-n-1},\qquad n\in\Z_+
 $$
и доказывается это цепочкой:
 \begin{multline*}
a^{-n-1}=\eqref{a^(-n)}=\frac{1}{a^{n+1}}=\eqref{def:a^n}=\frac{1}{a^n\cdot
a}=\\=\eqref{(a/b)(c/d)=(ac)/(bd)}=\frac{1}{a^n}\cdot\frac{1}{a}=\eqref{a^(-n)}=
\frac{a^{-n}}{a}
 \end{multline*}

Докажем теперь \eqref{a^(m+n)}. Зафиксируем $m\in\Z$. Для $n\in\Z_+$ это
доказывается индукцией:
 \begin{align*}
& a^{m+0}=a^m=a^m\cdot a=a^m\cdot a^0
 \\
& a^{m+n+1}\overset{\eqref{a^(m+1)=a^m-a-m-in-Z}}{=}\kern-5pt
\underbrace{a^{m+n}}_{\scriptsize
\begin{matrix}
\text{\rotatebox{90}{$=$}}\\ \phantom{,} a^m\cdot a^n, \\ \text{посылка}\\
\text{индукции}\end{matrix}}\kern-5pt\cdot a=a^m\cdot a^n\cdot a=
\\
&\qquad \overset{\eqref{a^(m+1)=a^m-a-m-in-Z}}{=}a^m\cdot a^{n+1}
 \end{align*}
Для $n\in-\N$, то есть $n=-k$, где $k\in\N$, формула \eqref{a^(m+n)}
переписывается так:
$$
a^{m-k}=a^m\cdot a^{-k}=\eqref{a^(-n)}=\frac{a^m}{a^k}
$$
то есть так:
$$
a^{m-k}\cdot a^k=a^m
$$
Заменив $m-k=l$, мы получим:
$$
a^l\cdot a^k=a^{k+l},\qquad l\in\Z,\ k\in\N
$$
Это просто иначе записанное равенство \eqref{a^(m+n)}, в котором одна из
степеней неотрицательна, и для этого случая мы его уже доказали.

2. Проверим степенные законы. Они все доказываются сначала для неотрицательных
степеней индукцией, а затем для отрицательных степеней выводятся из уже
доказанных формул. Например, формулу \eqref{1^n} нужно доказать сначала для
$n\in\Z_+$ индукцией:
 \begin{align*}
&1^0= \eqref{def:a^n}=1, \\
&1^{n+1}=\eqref{def:a^n}=\kern-5pt \underbrace{1^n}_{\scriptsize
\begin{matrix}
\text{\rotatebox{90}{$=$}}\\ \phantom{,}1, \\ \text{посылка}\\
\text{индукции}\end{matrix}}\kern-5pt\cdot 1=1\cdot 1=1
 \end{align*}
А отсюда уже становится ясно, что она верна и при $n\in-\N$, то есть при
$n=-k$, $k\in\N$:
 \begin{multline*}
1^n=1^{-k}=\eqref{DF:a^n-n<0}=\left(1^{-1}\right)^k=\\=\eqref{1^(-1)=1}=
\underbrace{1^k=1}_{\scriptsize \begin{matrix}\text{для $k\ge 0$}\\
\text{уже}\\ \text{доказано}\end{matrix}}
 \end{multline*}

Точно так же формула \eqref{(1/a)^n} доказывается сначала для $n\in\Z_+$
индукцией:
 \begin{align*}
&\l\frac{1}{a}\r^0\overset{\eqref{def:a^n}}{=}1=\frac{1}{1}\overset{\eqref{def:a^n}}{=}\frac{1}{a^0},
 \\
& \l\frac{1}{a}\r^{n+1}\overset{\eqref{def:a^n}}{=}\kern-5pt
\underbrace{\l\frac{1}{a}\r^n}_{\scriptsize
\begin{matrix}
\text{\rotatebox{90}{$=$}}\\ \phantom{,}\frac{1}{a^n}, \\ \text{посылка}\\
\text{индукции}\end{matrix}}\kern-5pt\cdot \frac{1}{a}=\frac{1}{a^n}\cdot
\frac{1}{a}\overset{\eqref{(a/b)(c/d)=(ac)/(bd)}}{=}\\
&\qquad=\frac{1}{a^n\cdot a} \overset{\eqref{def:a^n}}{=}\frac{1}{a^{n+1}}
 \end{align*}
А затем для $n=-k$, $k\in\N$ мы получаем:
 \begin{multline*}
\l\frac{1}{a}\r^n=\l\frac{1}{a}\r^{-k}=\eqref{a^(-n)}=\\=\underbrace{\frac{1}{\l\frac{1}{a}\r^k}=
\frac{1}{\frac{1}{a^k}}}_{\scriptsize\begin{matrix}\text{здесь применяется}\\ \eqref{(1/a)^n},\\
\text{уже доказанное}
\\ \text{для $k\ge 0$}\end{matrix}}=\eqref{a^(-n)}=\frac{1}{a^{-k}}=\frac{1}{a^n}
 \end{multline*}

Ну и формула \eqref{(ab)^n} тоже доказывается сначала для $n\in\Z_+$ индукцией:
 \begin{align*}
& (a\cdot b)^0 \overset{\eqref{def:a^n}}{=} 1=1\cdot 1
\overset{\eqref{def:a^n}}{=}a^0\cdot b^0,
 \\
& (a\cdot b)^{n+1}\overset{\eqref{def:a^n}}{=}\kern-5pt \underbrace{(a\cdot
b)^n}_{\scriptsize
\begin{matrix}
\text{\rotatebox{90}{$=$}}\\ \phantom{,} a^n\cdot b^n, \\ \text{посылка}\\
\text{индукции}\end{matrix}}\kern-5pt\cdot (a\cdot b)=a^n\cdot b^n\cdot a\cdot
b=\\ &\qquad=a^n\cdot a\cdot b^n\cdot b
\overset{\eqref{def:a^n}}{=}a^{n+1}\cdot b^{n+1}
 \end{align*}
А затем для $n=-k$, $k\in\N$, получается как следствие из уже доказанного:
 \begin{multline*}
(a\cdot b)^n= (a\cdot b)^{-k} =\eqref{a^(-n)}=\\=\underbrace{\frac{1}{(a\cdot
b)^k}=
\frac{1}{a^k\cdot b^k}}_{\scriptsize\begin{matrix}\text{здесь применяется}\\ \eqref{(ab)^n},\\
\text{уже доказанное}
\\ \text{для $k\ge 0$}\end{matrix}}\overset{\eqref{(a/b)(c/d)=(ac)/(bd)}}{=}
\frac{1}{a^k}\cdot \frac{1}{b^k}=\\=\eqref{a^(-n)}=a^{-k}\cdot b^{-k}=a^n\cdot
b^n
 \end{multline*}

3. Докажем накопительный закон \eqref{nakop-zakon-Z}. Зафиксируем $m\in\Z$ и
проведем сначала индукцию по $n\in\Z_+$.
 \begin{align*}
& (a^m)^0 \overset{\eqref{def:a^n}}{=}
1\overset{\eqref{def:a^n}}{=}a^0=a^{0\cdot m},
 \\
& (a^m)^{n+1}\overset{\eqref{def:a^n}}{=}\kern-5pt
\underbrace{(a^m)^n}_{\scriptsize
\begin{matrix}
\text{\rotatebox{90}{$=$}}\\ \phantom{,} a^{m\cdot n}, \\ \text{посылка}\\
\text{индукции}\end{matrix}}\kern-5pt\cdot a^m=a^{m\cdot n}\cdot a^m=\\
&\qquad\overset{\eqref{def:a^n}}{=}a^{m\cdot n+m} =a^{m\cdot(n+1)}
 \end{align*}
Это доказывает формулу \eqref{nakop-zakon-Z} для случая $m\in\Z$ и $n\in\Z_+$.
Если же $n\in-\N$, то есть $n=k$, $k\in\N$, то мы получаем:
 \begin{multline*}
(a^m)^n=(a^m)^{-k}\overset{\eqref{a^(-n)}}{=}
\overbrace{\frac{1}{(a^m)^k}=\frac{1}{a^{m\cdot
k}}}^{\scriptsize\begin{matrix}\text{здесь применяется}\\ \eqref{nakop-zakon-Z},\\
\text{уже доказанное}
\\ \text{для $k\ge 0$}\end{matrix}}\overset{\eqref{a^(-n)}}{=}\\=a^{-m\cdot k}=a^{m\cdot(-k)}=a^{m\cdot n}
 \end{multline*}

4. Формулу \eqref{0^n} при $n=0$ можно считать следствием \eqref{a^0-Z} (или
можно сказать что раньше мы отмечали это утверждение как равенство
\eqref{0^0=1}). Затем для $n=1$ ее можно вывести из \eqref{def:a^n} (или, опять
же, можно заявить, что мы уже этот факт отмечали в \eqref{a^1=a}):
$$
0^1=0^0\cdot 0=1\cdot 0=0
$$
А для остальных $n\in\N$ она доказывается индукцией:
 \begin{align*}
& 0^{n+1}\overset{\eqref{def:a^n}}{=}\kern-5pt \underbrace{0^n}_{\scriptsize
\begin{matrix}
\text{\rotatebox{90}{$=$}}\\ \phantom{,} 0, \\ \text{посылка}\\
\text{индукции}\end{matrix}}\kern-5pt\cdot 0=0
 \end{align*}
Условие сохранения знака \eqref{x^n>0} для $n\in\Z_+$ доказывается индукцией:
 \begin{align*}
& a^0 \overset{\eqref{def:a^n}}{=} 1>0,
 \\
& a^{n+1}\overset{\eqref{def:a^n}}{=}\kern-5pt \underbrace{a^n}_{\scriptsize
\begin{matrix}
\text{\rotatebox{90}{$<$}}\\ \phantom{,} 0, \\ \text{посылка}\\
\text{индукции}\end{matrix}}\kern-5pt\cdot \underbrace{a}_{\scriptsize
\begin{matrix}
\text{\rotatebox{90}{$<$}}\\ 0\end{matrix}}\overset{\eqref{x>0,y>0->x.y>0}}{>}0
 \end{align*}
А после этого для $n\in-\N$, то есть для $n=k$, $k\in\N$, мы получаем:
 \begin{align*}
 a^k&>0 && {\scriptsize\text{(уже доказано)}} \\
 & \Downarrow && {\scriptsize \eqref{x>0=>x^(-1)>0}} \\
a^n=a^{-k}\overset{\eqref{a^(-n)}}{=}\frac{1}{a^k} & >0 &&
 \end{align*}

5. Импликация \eqref{monot-a^n-n>0} доказывается индукцией по $n\in\N$. При
$n=1$ она становится {\it тавтологией}\index{тавтология} (то есть в ней
следствие является частью посылки, и поэтому утверждение будет очевидно
настолько, что в обычной речи это звучало бы глупо):
$$
0<x<y \quad\Longrightarrow\quad x<y
$$
Предположим, что мы доказали \eqref{monot-a^n-n>0} для какого-то $n=k$
 $$
0<x<y\quad\Longrightarrow\quad x^k<y^k
 $$
Тогда для $n=k+1$ мы получим:
 $$
\underbrace{\begin{matrix}
x^k<y^k,\\
{\scriptsize\eqref{x<y-a>0=>ax<ay}}\Downarrow\phantom{\scriptsize\eqref{x<y-a>0=>ax<ay}} \\
x\cdot x^k<x\cdot y^k,
\end{matrix}
\qquad
\begin{matrix}
x<y \\
\phantom{\scriptsize\eqref{x<y-a>0=>ax<ay}}\Downarrow{\scriptsize\eqref{x<y-a>0=>ax<ay}} \\
x\cdot y^k<y\cdot y^k
\end{matrix}}
 $$
 $$
 \Downarrow
 $$
 $$
x^n=x^{k+1}=x\cdot x^k<x\cdot y^k<y\cdot y^k=y^{k+1}=y^n
 $$
Когда \eqref{monot-a^n-n>0} доказано, \eqref{monot-a^n-n<0} становится его
следствием: если $n<0$, то есть $n=-k$, $k>0$, то
$$
0<x<y
$$
$$
\phantom{\text{\scriptsize\eqref{monot-a^n-n>0}}}\quad\Downarrow\quad\text{\scriptsize\eqref{monot-a^n-n>0}}
$$
$$
x^k<y^k
$$
$$
\phantom{\text{\scriptsize
\eqref{0<x<y=>x^(-1)>y^(-1)}}}\quad\Downarrow\quad\text{\scriptsize
\eqref{0<x<y=>x^(-1)>y^(-1)}}
$$
$$
x^n=x^{-k}=\frac{1}{x^k}>\frac{1}{y^k}=y^{-k}=y^n
$$

6. Импликацию \eqref{monot-a^n-a>1} можно доказать, например, с помощью
следующего вспомогательного утверждения:
 \beq\label{a>1=>a^n>1}
a>1\quad\Rightarrow\quad \forall n\in\N\quad a^n>1
 \eeq
Оно доказывается индукцией: подставив $n=1$ мы получим тавтологию,
$$
a>1\quad\Rightarrow\quad a^1=a>1
$$
а если оно верно при каком-то $n=k$
$$
a>1\quad\Rightarrow\quad a^k>1
$$
то при $n=k+1$ мы получим
$$
a^n=a^{k+1}=\underbrace{a^k}_{\scriptsize\begin{matrix}\text{\rotatebox{90}{$<$}}\\
0\end{matrix}}\cdot \underbrace{a}_{\scriptsize\begin{matrix}\text{\rotatebox{90}{$<$}}\\
1\end{matrix}}\overset{\eqref{x>0-a>1=>ax>x}}{>}a^k>1
$$
Когда \eqref{a>1=>a^n>1} доказано, \eqref{monot-a^n-a>1} становится его
следствием -- при $a>1$ получаем:
$$
m<n
$$
$$
\Downarrow
$$
$$
n-m>0
$$
$$
\phantom{\scriptsize\eqref{a>1=>a^n>1}}\Downarrow{\scriptsize\eqref{a>1=>a^n>1}}
$$
$$
a^{n-m}>1
$$
$$
\Downarrow
$$
$$
a^{n-m}-1>0
$$
$$
\Downarrow
$$
$$
a^n-a^m=a^m(a^{n-m}-1)\overset{\eqref{x>0,y>0->x.y>0}}{>}0
$$
$$
\Downarrow
$$
$$
a^m<a^n
$$

Наконец, из \eqref{monot-a^n-a>1} выводится \eqref{monot-a^n-0<a<1}:
$$
\underbrace{
\begin{matrix}
m<n \\
\Downarrow \\
-m>-n
\end{matrix}
\quad
\begin{matrix}
0<a<1 \\
\phantom{\scriptsize\eqref{a>1<=>0<a<1}}\Downarrow{\scriptsize\eqref{a>1<=>0<a<1}} \\
\frac{1}{a}>1
\end{matrix}}
$$
$$
\Downarrow
$$
$$
a^m=\l\frac{1}{a}\r^{-m}\overset{\eqref{monot-a^n-a>1}}{>}\l\frac{1}{a}\r^{-n}=a^n
$$

 \epr

\end{multicols}\noindent\rule[10pt]{160mm}{0.1pt}

\paragraph{Целая и дробная части числа.}\label{subsec-tselaya-chast}

 \bit{
\item[$\bullet$] {\it Целой частью}\index{целая часть} вещественного числа
$x\in\R$ называется целое число $[x]\in\Z$, определяемое правилом
 \beq\label{tselaya-chast-chisla}
[x]=\max \{n\in\Z:\; n\le x\},
 \eeq
то есть, максимум множества $E=\{n\in\Z:\; n\le x\}$ целых чисел, не
превосходящих $x$ (по теореме \ref{max-i-min-Z}, такой максимум существует).
 }\eit

 \bigskip
\centerline{\bf Свойства целой части:}
 \bit{\it
\item[$1^\circ$.] Целая часть числа $x\in[0;1)$ равна нулю:
 \beq\label{tsel-chast-0-1}
x\in[0;1)\qquad\Longrightarrow\qquad [x]=0
 \eeq

\item[$2^\circ$.] Для любого $x\in\R$ справедливо двойное неравенство
 \beq\label{opr-tsel-chasti}
[x]\le x< [x]+1
 \eeq

\item[$3^\circ$.] При сдвиге на целое число целая часть сдвигается на это
число:
 \beq\label{[x+m]=[x]+m}
[x+m]=[x]+m,\qquad x\in\R,\ m\in\Z
 \eeq

\item[$4^\circ$.] Тождественность на целых числах: если $x\in\Z$, то $[x]=x$.

\item[$5^\circ$.] Монотонность: если $x,y\in\R$, то
 \beq\label{x<y=>tch-x<tch-y}
x\le y \qquad\Longrightarrow\qquad [x]\le [y]
 \eeq

\item[$6^\circ$.] Сохранение двойных неравенств с целыми концами: если
$m,n\in\Z$ и $x\in\R$, то
 \beq\label{m<x<n=>m<tch-x<n}
m\le x<n \qquad\Longrightarrow\qquad  m\le [x]<n
 \eeq

 }\eit

\bpr

1. Если $x\in[0;1)$, то с одной стороны,
$$
0\in\{n\in\Z:\; n\le x\}\quad\Longrightarrow\quad 0\le \max\{n\in\Z:\; n\le
x\}=[x]
$$
а с другой --
$$
\{n\in\Z:\; n\le x\}\le
x<1\quad\overset{\eqref{E<n=>max-E<n}}{\Longrightarrow}\quad
[x]=\max\{n\in\Z:\; n\le x\}<1
$$
То есть $0\le [x]<1$. По следствию \ref{COR:n-le-m<n+1} это означает, что
$[x]=0$.

2. Обозначим $E=\{n\in\Z:\; n\le x\}$. Тогда, во-первых,
$$
\{n\in\Z:\; n\le x\}\le x\quad\Longrightarrow\quad [x]=\max\{n\in\Z:\; n\le
x\}\le x
$$
и, во-вторых, если бы $[x]+1\le x$, то мы получили бы цепочку
$$
[x]+1\le x\quad\Longrightarrow\quad [x]+1\in E\quad\Longrightarrow\quad
[x]+1\le \max E=[x]
$$
Поскольку последнее невозможно, это означает, что $x<[x]+1$.

3. Если $x\in\Z$, то $x\in\{n\in\Z:\; n\le x\}\le x$, и поэтому
$\max\{n\in\Z:\; n\le x\}=x$.

4. При $x\in\R$, $m\in\Z$ получаем:
 \begin{multline*}
[x+m]=\max\{n\in\Z:\ n\le x+m\}=\max\{n\in\Z:\ n-m\le x\}=\\=\max\{k+m:\
k\in\Z\ \& \ k\le x\}=\max\{k\in\Z:\ k\le x\}+m
 \end{multline*}

5. Если $x\le y$, то $\{n\in\Z:\; n\le x\}\subseteq\{n\in\Z:\; n\le y\}$, и
поэтому $[x]=\max\{n\in\Z:\; n\le x\}\le\max\{n\in\Z:\; n\le y\}=[y]$.

6. Если $m\le x$, то по уже доказанным свойствам $2^\circ$ и $3^\circ$,
$m=[m]\le [x]$. А если $x<n$, то для всякого $k\in\Z$ условие $k\le x$
автоматически влечет условие $k<n$ а это, в силу \eqref{m<n=>m-le-n-1}, влечет
условие $k\le n-1$. Отсюда
$$
[x]=\max\{k\in\Z:\; k\le x\}\le \max\{k\in\Z:\; k\le n-1\}=n-1<n
$$
 \epr

 \bit{
\item[$\bullet$] {\it Дробной частью}\index{дробная часть} $\{x\}$
вещественного числа $x\in\R$ называется разность между $x$ и его целой частью
$[x]$:
 \beq\label{drobnaya-chast-chisla}
\{x\}:=x-[x]
 \eeq
 }\eit

 \bigskip
\centerline{\bf Свойства дробной части:}
 \bit{
\item[$1^\circ$.] Дробная часть числа $x\in[0;1)$ совпадает с $x$:
 \beq\label{drob-chast-0-1}
x\in[0;1)\qquad\Longrightarrow\qquad \{x\}=x
 \eeq

\item[$2^\circ$.] Для любого $x\in\R$ справедливо двойное неравенство
 \beq\label{drob-chast<1}
0\le\{x\}< 1
 \eeq

\item[$3^\circ$.] При сдвиге на целое число дробная часть не меняется:
 \beq\label{drob-chast(x+m)}
\{x+m\}=\{x\},\qquad x\in\R,\ m\in\Z
 \eeq
 }\eit
\bpr 1. Если $x\in[0;1)$, то в силу \eqref{tsel-chast-0-1}, $[x]=0$, поэтому
$\{x\}=x-0=x$.

2. Неравенство \eqref{drob-chast<1} следует сразу из неравенства
\eqref{opr-tsel-chasti}:
$$
[x]\le x<[x]+1\quad\Longrightarrow\quad 0\le x-[x]<1
$$

3. А тождество \eqref{drob-chast(x+m)} -- из \eqref{[x+m]=[x]+m}:
$$
\{x+m\}=(x+m)-[x+m]=(x+m)-([x]+m)=x-[x]=\{x\}
$$
 \epr

\noindent\rule{160mm}{0.1pt}\begin{multicols}{2}

\paragraph{Десятичная запись целых чисел.}
\label{SUBSEC-desyat-zapis-N} Как известно, целые неотрицательные числа
$n\in\Z_+$ (мы определили это множество выше формулой \eqref{DF:Z_+}) можно
записывать цифрами $0$, $1$, $2$, $3$, $4$, $5$, $6$, $7$, $8$, $9$. Например,
записи 487 и 1204 означают числа
 \begin{align*}
& 487=7+8\cdot 10+4\cdot 10^2, \\
& 1204=4+0\cdot 10+2\cdot 10^2+1\cdot 10^3.
 \end{align*}
Тот факт, что любое число $n\in\Z_+$ можно записать таким образом, и что такая
запись однозначно определяет число (то есть не может случиться, чтобы разные
числа записывались одинаково), следует из принципа Архимеда
\ref{Archimed-principle}, принципа математической индукции и свойств
отображения $x\mapsto[x]$. Здесь мы объясним, почему это так.

\btm\label{TH-desyatichnaya-zapis-v-N} Для всякого числа $x\in\N$ найдется
число $n\in\Z_+$ и конечная последовательность чисел
$a_0,a_1,...,a_n\in\{0;1;2;...;9\}$ такая, что
 \begin{multline}\label{desyatichnaya-zapis-chisla}
x=a_0+a_1\cdot 10+a_2\cdot 10^2+...+a_n\cdot 10^n=\\=\sum_{k=0}^n a_k\cdot 10^k
 \end{multline}
Числа $n$ и $a_i$ в этой формуле можно выбрать так, чтобы старший коэффициент
был ненулевым
 \beq\label{0---a_n-ne-0}
a_n\ne 0,
 \eeq
и тогда представление $x$ в виде \eqref{desyatichnaya-zapis-chisla} становится
единственным.
 \etm

 \biter{
\item[$\bullet$] Формула \eqref{desyatichnaya-zapis-chisla} называется {\it
десятичным представлением} числа $x\in\N$, а последовательность цифр
$\{a_n,a_{n-1},...,a_1,a_0\}$ -- {\it последовательностью коэффициентов} этого
представления. Коротко эту формулу принято записывать в виде
$$
x=a_n...a_1a_0
$$
и такая запись называется {\it десятичной записью} числа $x\in\N$. Например,
запись
$$
x=1825
$$
расшифровывается как равенство
$$
x=1\cdot 10^3+8\cdot 10^2+2\cdot 10^1+5\cdot 10^1
$$

\item[$\bullet$] Если $x\in-\N$, то по теореме \ref{TH-desyatichnaya-zapis-v-N}
число $-x\in\N$ представимо в виде
$$
-x=\sum_{k=0}^n a_k\cdot 10^k
$$
и для этого используется запись
$$
x=-a_n...a_1a_0
$$
называемая {\it десятичной записью} числа $x\in-\N$. Например, запись
$$
x=-247
$$
расшифровывается как равенство
$$
-x=2\cdot 10^2+4\cdot 10^1+7\cdot 10^1
$$
или, что эквивалентно, как равенство
$$
x=-2\cdot 10^2-4\cdot 10^1-7\cdot 10^1
$$

 }\eiter

Для доказательства теоремы \ref{TH-desyatichnaya-zapis-v-N} нам понадобятся три
леммы.

\blm\label{LM-x<10^n} Для всякого $x\in\R$ найдется натуральное число
$n\in\Z_+$ такое, что
 \beq\label{x<10^n}
x<10^{n+1}
 \eeq
 \elm
\bpr Рассмотрим число $\frac{x-1}{9}$. По принципу Архимеда
\ref{Archimed-principle}, для него найдется $n\in\N$ такое, что
$\frac{x-1}{9}\le n$. Отсюда получаем следствия:
$$
\frac{x-1}{9}<n+1
$$
$$
\Downarrow
$$
$$
x-1<(n+1)\cdot 9
$$
$$
\Downarrow
$$
$$
x<1+(n+1)\cdot9 \kern-30pt
\underset{\scriptsize\begin{matrix}\uparrow \\
\text{здесь мы применяем} \\
\text{неравенство Бернулли} \\ \eqref{nerav-Bernoulli}
\end{matrix}}{\le}\kern-30pt (1+9)^{n+1}=10^{n+1}
$$
 \epr

\blm\label{lemma-dlya-10} Если $x\in\R$ и выполняется неравенство
 $$
0\le x<10^{n+1},
 $$
то найдется число $a\in \{0;1;2;...;9\}$ такое, что
 \beq\label{a-cdot-10^n-1-le-x<a+1-cdot-10^n}
a\cdot 10^n\le x<(a+1)\cdot 10^n
 \eeq
 \elm
\bpr Обозначим через $a$ целую часть числа $\frac{x}{10^n}$:
$$
a=\left[\frac{x}{10^n}\right]
$$
Тогда, во-первых, из свойства целой части сохранять неравенства с целыми
концами \eqref{m<x<n=>m<tch-x<n} следует, что $a\in \{0;1;2;...;9\}$,
$$
0\le x<10^{n+1}
$$
$$
\Downarrow
$$
$$
0\le \frac{x}{10^n}<10
$$
$$
\phantom{{\scriptsize
\eqref{m<x<n=>m<tch-x<n}}}\quad\Downarrow\quad{\scriptsize
\eqref{m<x<n=>m<tch-x<n}}
$$
$$
0\le \underbrace{\left[\frac{x}{10^n}\right]}_{a} <10
$$
И, во-вторых, по определению целой части \eqref{opr-tsel-chasti},
$$
\underbrace{\left[\frac{x}{10^n}\right]}_{a}\le \frac{x}{10^n}<
\underbrace{\left[\frac{x}{10^n}\right]+1}_{a+1}
$$
$$
\Downarrow
$$
$$
a\cdot10^n\le x<(a+1)\cdot 10^n
$$
 \epr

\blm\label{LM-desyat-zap-starshij-koeff} Если число $x$ представлено в виде
\eqref{desyatichnaya-zapis-chisla}, то старший коэффициент $a_n$ в этой формуле
вычисляется по формуле
 \beq\label{koeff-desyatichnoi-zapisi-chisla}
a_n=\left[\frac{x}{10^n}\right]
 \eeq
 \elm

\bpr Здесь применяется формула \eqref{Geom-progr-N} для суммы первых членов
геометрической прогрессии:
 \begin{multline*}
\sum_{k=0}^{n-1}
\underbrace{a_k}_{\scriptsize\begin{matrix}\text{\rotatebox{90}{$\ge$}}
\\ 9\end{matrix}}\cdot 10^k\le \sum_{k=0}^{n-1} 9\cdot 10^k=
9\cdot \sum_{k=0}^{n-1} 10^k=\\=\eqref{Geom-progr-N}=9\cdot\frac{10^n-1}{10-1}=
9\cdot\frac{10^n-1}{9}=\\=10^n-1<10^n
 \end{multline*}
$$
\Downarrow
$$
$$
0\le\frac{1}{10^n}\cdot\sum_{k=0}^{n-1} a_k\cdot 10^k<1
$$
$$
\Downarrow
$$
 \begin{multline*}
\frac{x}{10^n}=\frac{1}{10^n}\cdot\left(\sum_{k=0}^{n-1} a_k\cdot 10^k+a_n\cdot
10^n\right)=\\= \underbrace{\frac{1}{10^n}\cdot\sum_{k=0}^{n-1} a_k\cdot
10^k}_{\scriptsize\begin{matrix} \text{\rotatebox{-90}{$\in$}} \\
[0,1)\end{matrix}} +\underbrace{a_n}_{\scriptsize\begin{matrix}
\text{\rotatebox{-90}{$\in$}} \\ \Z \end{matrix}}
 \end{multline*}
$$
\Downarrow
$$
$$
\left[\frac{x}{10^n}\right]=a_n
$$
\epr

\bpr[Доказательство теоремы \ref{TH-desyatichnaya-zapis-v-N}] Пусть $x\in\N$.
Применяя  лемму \ref{LM-x<10^n}, можно выбрать $n\in\Z_+$ так, чтобы
выполнялось \eqref{x<10^n}, то есть, \eqref{0-le-x<10^n+1}:
 \beq\label{0-le-x<10^n+1}
0\le x<10^{n+1},
 \eeq
1. Покажем, что при выполнении этого неравенства $x$ можно представить в виде
\eqref{desyatichnaya-zapis-chisla}.

Это доказывается индукцией по $n$.

1) При $n=0$ условие \eqref{0-le-x<10^n+1} превращается в $0\le x<10$, и
формула \eqref{desyatichnaya-zapis-chisla} становится верна, если положить
$a=x$.

2) Предположим, что наше утверждение верно при $n=m-1$:
 \beq\label{0-le-x<10^m-predpol-induktsii}
0\le x<10^m\quad\Rightarrow\quad x=\sum_{k=0}^{m-1} a_k\cdot 10^k
 \eeq
и покажем, что тогда оно верно и при $n=m$:
 $$\label{0-le-x<10^m+1-sledstvie-induktsii}
0\le x<10^{m+1}\quad\Rightarrow\quad x=\sum_{k=0}^m a_k\cdot 10^k
 $$
Действительно,
$$
0\le x<10^{m+1}
$$
$$
\phantom{{\scriptsize \text{\scriptsize лемма
\ref{lemma-dlya-10}}}}\quad\Downarrow\quad{\scriptsize \text{\scriptsize лемма
\ref{lemma-dlya-10}}}
$$
 \begin{multline*}
\exists a_m\in\{0;1;2;...;9\}:\\ a_m\cdot 10^m\le x <(a_m+1)\cdot 10^m
 \end{multline*}
$$
\phantom{{\scriptsize \text{\scriptsize (вычитаем $a_m\cdot
10^m$)}}}\quad\Downarrow\quad{\scriptsize \text{\scriptsize (вычитаем $a_m\cdot
10^m$)}}
$$
$$
0\le x-a_m\cdot 10^m< 10^m
$$
 $$
\phantom{{\scriptsize \begin{pmatrix} \text{применяем посылку}\\
\text{индукции: \eqref{0-le-x<10^m-predpol-induktsii}}
\end{pmatrix}}}
\quad\Downarrow\quad{\scriptsize \begin{pmatrix} \text{применяем посылку}\\
\text{индукции: \eqref{0-le-x<10^m-predpol-induktsii}}
\end{pmatrix}}
 $$
$$
x-a_m\cdot 10^m  =\sum_{k=0}^{m-1} a_k\cdot 10^k
$$
 $$
\phantom{{\scriptsize \begin{pmatrix} \text{переносим $a_m\cdot 10^m$}\\
\text{в правую часть}
\end{pmatrix}}}
\quad\Downarrow\quad{\scriptsize \begin{pmatrix} \text{переносим $a_m\cdot 10^m$}\\
\text{в правую часть}
\end{pmatrix}}
 $$
$$
x =\sum_{k=0}^{m-1} a_k\cdot 10^k+a_m\cdot 10^m=\sum_{k=0}^{m} a_k\cdot 10^k
$$

2. Мы доказали, что $x$ представимо в виде \eqref{desyatichnaya-zapis-chisla}.
Теперь убедимся, что такое представление единственно, с точностью до добавления
или исключения нулевых слагаемых (то есть слагаемых со старшими коэффициентами
$a_n=0$). Пусть нам даны два разных представления:
$$
x=\sum_{k=0}^{n} a_k\cdot 10^k \qquad\text{и}\qquad x=\sum_{k=0}^{l} b_k\cdot
10^k
$$
Добавив слагаемые с нулевыми коэффициентами (то есть, взяв дополнительные
$a_i=0$ или $b_i=0$), можно добиться, чтобы количество слагаемых в этих суммах
стало одинаковым:
 \beq\label{odnoznachnost'-desyat-zapisi}
\sum_{k=0}^n a_k\cdot 10^k=\sum_{k=0}^n b_k\cdot 10^k
 \eeq
Нам нужно проверить, что тогда
 \beq\label{a_k=b_k-desyat-zapis}
a_k=b_k,\qquad k=0,...,n.
 \eeq
Это также делается индукцией по $n$.

1) Если $n=0$, то \eqref{odnoznachnost'-desyat-zapisi} сразу превращается в
\eqref{a_k=b_k-desyat-zapis}: $a_0=b_0$.

2) Предполагаем, что это верно при $n=m$,
 \begin{multline}\label{posylka-ind-dlya-odnoz-desyat-razl}
\sum_{k=0}^m a_k\cdot 10^k=\sum_{k=0}^m b_k\cdot 10^k \quad\Rightarrow\\
\Rightarrow\quad\forall k=0,...,m\quad a_k=b_k
 \end{multline}
тогда для $n=m+1$ получаем:
 \beq
\sum_{k=0}^{m+1} a_k\cdot 10^k =\sum_{k=0}^{m+1} b_k\cdot 10^k
\label{desyat-zap-m+1}
 \eeq
 $$
\phantom{\text{\scriptsize (лемма \ref{LM-desyat-zap-starshij-koeff})}}
\quad\Downarrow\quad\text{\scriptsize (лемма
\ref{LM-desyat-zap-starshij-koeff})}
 $$
 $$
a_{m+1}=\left[\sum_{k=0}^{m+1} a_k\cdot 10^k\right] = \left[\sum_{k=0}^{m+1}
b_k\cdot 10^k\right]=b_{m+1}
 $$
 $$
\phantom{{\scriptsize \begin{pmatrix} \text{старшие слагаемые}\\
\text{в \eqref{desyat-zap-m+1} равны,}\\ \text{и их можно}\\ \text{отбросить}
\end{pmatrix}}}
\quad\Downarrow\quad{\scriptsize \begin{pmatrix} \text{старшие слагаемые}\\
\text{в \eqref{desyat-zap-m+1} равны,}\\ \text{и их можно}\\ \text{отбросить}
\end{pmatrix}}
 $$
 $$
\sum_{k=0}^m a_k\cdot 10^k =\sum_{k=0}^m b_k\cdot 10^k
 $$
 $$
\phantom{{\scriptsize \begin{pmatrix} \text{применяем посылку}\\
\text{индукции: \eqref{posylka-ind-dlya-odnoz-desyat-razl}}
\end{pmatrix}}}
\quad\Downarrow\quad{\scriptsize \begin{pmatrix} \text{применяем посылку}\\
\text{индукции: \eqref{posylka-ind-dlya-odnoz-desyat-razl}}
\end{pmatrix}}
 $$
 $$
 \forall k=0,...,m\quad a_k =b_k
 $$
 $$
 {\scriptsize \begin{pmatrix}\text{а при $k=m+1$ равенство}\\ \text{$a_k=b_k$ уже доказано.}\end{pmatrix}}
 $$

3. Остается убедиться, что если $x\ne 0$, то $n$ можно выбрать так, чтобы
$a_n\ne 0$. Для этого можно рассмотреть множество $E$, состоящее из тех
$n\in\Z_+$, для которых справедливо представление
\eqref{desyatichnaya-zapis-chisla} и найти его минимум $n=\min E$. Тогда
$a_n\ne 0$, потому что иначе мы получили бы, что старшее слагаемое в
\eqref{desyatichnaya-zapis-chisla} можно выбросить,
$$
x=\sum_{k=0}^{n-1} a_k\cdot 10^k+ \underbrace{a_n\cdot
10^n}_{0}=\sum_{k=0}^{n-1} a_k\cdot 10^k
$$
и это означало бы, что $n-1\in E$ (то есть $n$ не может быть минимумом $E$).
 \epr

\end{multicols}\noindent\rule[10pt]{160mm}{0.1pt}

\subsection{Деление с остатком и делимость}

\paragraph{Деление с остатком.}

\btm[\bf о делении с остатком]\label{TH:delenie-s-ostatkom} Для любых двух
чисел $a\in\R$ и $b>0$ существуют и единственны числа $q\in\Z$ и $r\in[0;b)$
такие, что
 \beq\label{a=b-cdot-q+r}
 a=b\cdot q+r,
 \eeq
Если вдобавок $a\in\Z$ и $b\in\N$, то $r\in\Z_+$.
 \etm
 \bpr
Числа $q$ и $r$ можно определить формулами
 $$
q=\left[\frac{a}{b}\right],\qquad r=b\cdot\left\{\frac{a}{b}\right\}
 $$
Тогда по определению целой части $q\in\Z$, а по определению дробной части,
$\frac{r}{b}\in[0;1)$, и значит $r\in[0,b)$. Нам остается доказать, что такие
числа $q$ и $r$ единственны. Предположим, что существуют два разложения одного
числа $a$:
$$
a=b\cdot q+r=b\cdot\tilde{q}+\tilde{r}
$$
Тогда
$$
q+\frac{r}{b}=\tilde{q}+\frac{\tilde{r}}{b}
$$
Из этого следует, что $q=\tilde{q}$, потому что если бы это было не так, то
есть например $q<\tilde{q}$, то, по свойству \eqref{m<n=>m-le-n-1}, мы получили
бы, что $q+1\le\tilde{q}$, и поскольку $\frac{r}{b}<1$,
$$
q+\frac{r}{b}<q+1\le\tilde{q}\le\tilde{q}+\frac{\tilde{r}}{b}
$$
А с другой стороны, равенство $\frac{r}{b}=\tilde{q}-q+\frac{\tilde{r}}{b}$
означает, что $\frac{r}{b}$ и $\tilde{q}-q+\frac{\tilde{r}}{b}$ отличаются на
целое число $\tilde{q}-q$, поэтому
$$
\frac{r}{b}=\left\{\frac{r}{b}\right\}=\left\{\tilde{q}-q+\frac{\tilde{r}}{b}\right\}=
\eqref{drob-chast(x+m)}=\frac{\tilde{r}}{b},
$$
и отсюда $r=\tilde{r}$.

Остается заметить, что поскольку $q\in\Z$, то при $a\in\Z$ и $b\in\N$ число
$r=a-b\cdot q$ должно быть целым.
 \epr

\noindent\rule{160mm}{0.1pt}\begin{multicols}{2}

\paragraph{Делимость в множестве $\N$ натуральных чисел.}

 \biter{

\item[$\bullet$] Говорят что число $a\in\N$ {\it делится} на число $b\in\N$, и
записывают это формулой
$$
a\div b,
$$
если найдется число $m\in\N$ такое, что
$$
a=m\cdot b.
$$
В этом случае
 \biter{
\item[---] число $a$ называется {\it кратным} числа $b$,

\item[---] число $b$ называется {\it делителем} числа $a$.
 }\eiter

\item[$\bullet$] Если $A\subseteq\N$ и $B\subseteq\N$ -- множества, то записи
 $$
A\div b,\quad a\div B,\quad A\div B
 $$
означают соответственно, что $a\div b$ выполняется для любого $a\in A$, для
любого $b\in B$, и для любых $a\in A$ и $b\in B$.
 }\eiter

\bex Следующие утверждения можно проверить, например, с помощью калькулятора:
 $$
6\div 2, \quad 123\div 3, \quad 430095\div 541
 $$
\eex

 \bigskip
\centerline{\bf Свойства отношения делимости:}
 \bit{\it
\item[$1^\circ$.] Всегда $x\div 1$ и $x\div x$.

\item[$2^\circ$.]\label{x:y=>ax:y} Если $x\div y$, то $(a\cdot x)\div y$ для
любого $a\in\N$.

\item[$3^\circ$.] Если $x\div y$ и $y\div z$, то $x\div z$.

\item[$4^\circ$.]\label{x:a,y:a=>x-y:a} Если $x\div a$ и $y\div a$, то
$(x+y)\div a$, а если вдобавок $x>y$, то $(x-y)\div a$.

\item[$5^\circ$.]\label{x:y=>x>y} Если $x\div y$, то $x\ge y$.

 }\eit

Доказательство этих свойств мы оставляем читателю.

\bers Докажите методом математической индукции\footnote{В первых двух примерах
нужно использовать бином Ньютона.}:
 \biter{
\item[1)] $n^5-n\div 5$ при любом $n\in\N$;

\item[2)] $n^7-n\div 7$ при любом $n\in\N$;

\item[3)] $n\cdot (2n^2-3n+1)\div 6$ при любом $n\in\N$;

\item[4)] $11^{n+1}+12^{2n-1}\div 133$ при любом $n\in\N$;

\item[5)] $5\cdot 2^{3n-2}+3^{3n-1}\div 19$ при любом $n\in\N$.

 }\eiter
 \eers

\paragraph{Наибольший общий делитель и наименьшее общее кратное.}

 Пусть $A\subseteq\N$.
 \biter{
\item[$\bullet$] Число $x\in\N$ называется {\it общим кратным} множества $A$,
если оно кратно любому числу $a\in A$ (то есть делится на любое $a\in A$):
$$
\forall a\in A\qquad x\div a,
$$
коротко это записывается формулой
$$
x\div A.
$$
Множество всех общих кратных для $A$ обозначается $A^\triangledown$:
 \begin{multline}\label{DEF:A^triangledown}
A^\triangledown:=\{x\in\N: \ x\div A\}=\\=\{x\in\N: \ \forall a\in A\quad x\div
a\}
 \end{multline}
Очевидно,
 \beq\label{A^triangledown-div-A}
A^\triangledown\div A
 \eeq

\item[$\bullet$] {\it Hаименьшее общее кратное} множества $A$ (то есть
наименьшее число среди всех $x\in A^\triangledown$) обозначается
$\triangledown(A)$:
 \begin{multline}\label{DEF:NOK}
\triangledown(A):=\min A^\triangledown=\\=\min\{x\in\N: \ \forall a\in A\quad
x\div a\}
 \end{multline}
В частном случае, когда $A$ состоит из двух элементов $a$ и $b$, наименьшее
общее кратное множества $A$ называется {\it наименьшим общим кратным чисел $a$
и $b$} и обозначается
 \begin{multline}\label{DEF:a-V-b}
a\,\triangledown\, b:=\triangledown\{a,b\}=\\=\min\{x\in\N:\ x\div a\ \& \
x\div b \}
 \end{multline}
  }\eiter

 \bprop\label{PROP:(A-cup-B)^triangledown=A^triangledown-cap-B^triangledown}
Следующее соотношение справедливо для любых двух множеств $A\subseteq\N$ и
$B\subseteq\N$:
 \beq\label{(A-cup-B)^triangledown=A^triangledown-cap-B^triangledown}
(A\cup B)^\triangledown =A^\triangledown\cap B^\triangledown
 \eeq
 \eprop
 \bpr
$$
x\in (A\cup B)^\triangledown
$$
$$
\Updownarrow
$$
$$
x\div (A\cup B)
$$
$$
\Updownarrow
$$
$$
x\div A\quad \&\quad x\div B
$$
$$
\Updownarrow
$$
$$
x\in A^\triangledown \quad \&\quad x\in B^\triangledown
$$
$$
\Updownarrow
$$
$$
x\in A^\triangledown\cap B^\triangledown
$$
 \epr

\brem Если в формуле
\eqref{(A-cup-B)^triangledown=A^triangledown-cap-B^triangledown} поменять
местами $\cup$ и $\cap$, она перестает быть верной:
$$
(A\cap B)^\triangledown\ne A^\triangledown\cup B^\triangledown
$$
Например, при $A=\{2\}$ и $B=\{3\}$ мы получаем
$$
(\{2\}\cap\{3\})^\triangledown=\varnothing^\triangledown=\N,
$$
но
$$
\{2\}^\triangledown\cup\{3\}^\triangledown=2\N\cup 3\N
$$
Это не одно и то же, потому что, например,
$$
5\in\N\quad\&\quad 5\notin2\N\cup 3\N
$$
 \erem

\bprop Если $A$ конечно (необязательно непусто), то $A^\triangledown$ непусто,
и поэтому существует наименьшее общее кратное $\triangledown(A)$. \eprop
 \bpr
Если $A$ пусто, то число $1$ будет его общим кратным, поскольку при $a\in A$
условие $1\div a$ выполняется тривиально (из-за отсутствия таких $a$). Если же
$A$ непусто, то, поскольку оно конечно, его можно представить в виде конечной
последовательности (с неповторяющимися элементами)
$$
A=\{a_1,...,a_n\}
$$
а затем рассмотреть произведение всех чисел, входящих в $A$
$$
x=\prod_{a\in A} a=\prod_{i=1}^n a_i
$$
Это число делится на любой элемент $a\in A$, то есть является общим кратным для
множества $A$:
$$
x\in A^\triangledown
$$
Значит, множество $A^\triangledown$ непусто. Как любое непустое подмножество в
$\N$, оно, в силу свойства $5^\circ$ на с.\pageref{minimum-v-N}, должно иметь
минимальный элемент $\min A^\triangledown=\triangledown(A)$.
 \epr

Снова пусть $A\subseteq\N$.
 \biter{
\item[$\bullet$] Число $x\in\N$ называется {\it общим делителем} множества $A$,
если оно является делителем для любого числа $a\in A$ (то есть любое $a\in A$
делится на $x$):
$$
\forall a\in A\qquad a\div x,
$$
коротко это записывается формулой
$$
A\div x.
$$
Множество всех общих делителей для $A$ обозначается $A^\vartriangle$:
 \begin{multline}\label{DEF:A^vartriangledown}
A^\vartriangle:=\{x\in\N: \ A\div x\}=\\=\{x\in\N: \ \forall a\in A\quad a\div
x\}
 \end{multline}
Очевидно,
 \beq\label{A-div-A^vartriangle}
A\div A^\vartriangle
 \eeq

\item[$\bullet$] {\it Наибольший общий делитель} множества $A$ (то есть
наибольшее число среди всех $x\in A^\vartriangle$) обозначается
$\vartriangle(A)$:
 \begin{multline}\label{DEF:NOD}
\vartriangle(A)=\max A^\vartriangle=\\=\max\{x\in\N: \ \forall a\in A\quad
a\div x\}
 \end{multline}
В частном случае, когда $A$ состоит из двух элементов $a$ и $b$, наибольший
общий делитель множества $A$ называется {\it наибольшим общим делителем чисел
$a$ и $b$} и обозначается
 \beq\label{DEF:a-wedge-b}
a\vartriangle b:=\max\{x\in\N:\ a\div x\ \& \ b\div x \}
 \eeq
  }\eiter

 \bprop\label{PROP:(A-cup-B)^vartriangle=A^vartriangle-cap-B^vartriangle}
Следующее соотношение справедливо для любых двух множеств $A\subseteq\N$ и
$B\subseteq\N$:
 \beq\label{(A-cup-B)^vartriangle=A^vartriangle-cap-B^vartriangle}
(A\cup B)^\vartriangle=A^\vartriangle\cap B^\vartriangle.
 \eeq
 \eprop
 \bpr Это доказывается по аналогии с предложением
\ref{PROP:(A-cup-B)^triangledown=A^triangledown-cap-B^triangledown}:
$$
x\in (A\cup B)^\vartriangle
$$
$$
\Updownarrow
$$
$$
(A\cup B)\div x
$$
$$
\Updownarrow
$$
$$
A\div x\quad \&\quad B\div x
$$
$$
\Updownarrow
$$
$$
x\in A^\vartriangle \quad \&\quad x\in B^\vartriangle
$$
$$
\Updownarrow
$$
$$
x\in A^\vartriangle\cap B^\vartriangle
$$
 \epr

\brem Если в формуле
\eqref{(A-cup-B)^vartriangle=A^vartriangle-cap-B^vartriangle} поменять местами
$\cup$ и $\cap$, то она перестает быть верной:
$$
(A\cap B)^\vartriangle\ne A^\vartriangle\cup B^\vartriangle
$$
Например, при $A=\{2\}$ и $B=\{3\}$ мы получаем
$$
(\{2\}\cap\{3\})^\vartriangle=\varnothing^\vartriangle=\N,
$$
но
$$
\{2\}^\vartriangle\cup\{3\}^\vartriangle=\{1,2\}\cup\{1,3\}=\{1\}.
$$
\erem

\bprop Если $A$ непусто (необязательно конечно), то $A^\vartriangle$ непусто и
конечно, и поэтому существует наибольший общий делитель $\vartriangle(A)$.
\eprop
 \bpr
Число 1, очевидно, является общим делителем для множества $A$:
$$
A\div 1
$$
Это значит, что $1\in A^\vartriangle$, и поэтому $A^\vartriangle$ непусто. С
другой стороны, если взять какое-нибудь $a\in A$ (поскольку $A$ непусто, такое
$a$ найдется), то мы получим:
 $$
A^\vartriangle\subseteq \{a\}^\vartriangle\overset{\scriptsize\begin{matrix}\text{$5^\circ$ на}\\
\text{с.\pageref{x:y=>x>y}}\\ \downarrow\end{matrix}}{\subseteq} \{1,...,a\}
 $$
То есть $A^\vartriangle$ содержится в конечном множестве $\{1,...,a\}$. Значит,
по свойству $1^\circ$ на с.\pageref{subset-finite-set}, $A^\vartriangle$ само
конечно. После этого применяется теорема \ref{TH:max-finite-set}: поскольку
множество $A^\vartriangle$ конечно, оно имеет максимальный элемент.
 \epr

\btm Следующие соотношения справедливы для любого непустого конечного
$A\subseteq\N$:
 \begin{align}
\triangledown(A)\ge\max & A\ge\min A\ge\vartriangle(A), \label{triangledown(A)>max-A}\\
A^\triangledown\div\triangledown(A)\div &A\div\vartriangle(A)\div A^\vartriangle,
\label{A^triangledown-div-triangledown-A} \\
\vartriangle(A)&=\triangledown(A^\vartriangle), \label{vartriangle(A)=triangledown(A^vartriangle)} \\
\triangledown(A)&=\vartriangle(A^\triangledown).
\label{triangledown(A)=vartriangle(A^triangledown)}
 \end{align}
\etm
 \bpr
1. В цепочке \eqref{triangledown(A)>max-A} нужно проверить два крайних
неравенства.
$$
\triangledown(A)\ge\max A,\qquad \min A\ge\vartriangle(A)
$$
Они оба следуют из свойства $5^\circ$ на с.\pageref{x:y=>x>y}. Например,
первое:
$$
x\in A^\triangledown
$$
$$
\Downarrow
$$
$$
x\div A
$$
$$
\phantom{{\scriptsize\text{$5^\circ$ на
с.\pageref{x:y=>x>y}}}}\quad\Downarrow\quad{\scriptsize\text{$5^\circ$ на
с.\pageref{x:y=>x>y}}}
$$
$$
x\ge A
$$
$$
\Downarrow
$$
$$
x\ge\max A
$$
Это верно для любого $x\in A^\triangledown$, значит
$$
\triangledown(A)=\min A^\triangledown=\min_{x\in A^\triangledown}x \ge\max A,
$$
И то же самое со вторым:
$$
x\in A^\vartriangle
$$
$$
\Downarrow
$$
$$
A\div x
$$
$$
\phantom{{\scriptsize\text{$5^\circ$ на
с.\pageref{x:y=>x>y}}}}\quad\Downarrow\quad{\scriptsize\text{$5^\circ$ на
с.\pageref{x:y=>x>y}}}
$$
$$
A\ge x
$$
$$
\Downarrow
$$
$$
\min A\ge x
$$
Это верно для любого $x\in A^\vartriangle$, значит
$$
\min A\ge\max_{x\in A^\vartriangle}x=\max A^\vartriangle=\vartriangle(A)
$$

2. В цепочке \eqref{A^triangledown-div-triangledown-A} середина очевидна:
 \beq\label{triangledown(A)-div-A-div-vartriangle(A)}
\triangledown(A)\div A\div\vartriangle(A)
 \eeq
-- с одной стороны,
$$
\underbrace{A^\triangledown\div
A}_{\eqref{A^triangledown-div-A}}\quad\Rightarrow\quad \triangledown(A)=\min
A^\triangledown\div A
$$
а с другой,
$$ \underbrace{A\div A^\vartriangle}_{\eqref{A-div-A^vartriangle}}
\quad\Rightarrow\quad A\div\max A^\vartriangle=\vartriangle(A)
$$
Поэтому нужно доказать только два крайних соотношения:
$$
A^\triangledown\div\triangledown(A),\qquad \vartriangle(A)\div A^\vartriangle.
$$
Мы будем их доказывать вместе с соотношениями
\eqref{vartriangle(A)=triangledown(A^vartriangle)} и
\eqref{triangledown(A)=vartriangle(A^triangledown)}, причем в следующей
последовательности:
 \beq\label{A^triangledown-div-triangledown(A)}
A^\triangledown\div\triangledown(A)
 \eeq
$$
\Downarrow
$$
$$
\vartriangle(A)=\triangledown(A^\vartriangle)
$$
$$
\Downarrow
$$
 \beq\label{vartriangle(A)-div-A^vartriangle}
\vartriangle(A)\div A^\vartriangle
 \eeq
$$
\Downarrow
$$
 $$
\triangledown(A)=\vartriangle(A^\triangledown)
 $$

3. Докажем первое звено в этой цепочке, то есть формулу
\eqref{A^triangledown-div-triangledown(A)}. Пусть $x\in A^\triangledown$, то
есть $x\div A$. Нам нужно показать, что $x\div\triangledown(A)$. Разделим $x$
на $\triangledown(A)$ с остатком: по теореме \ref{TH:delenie-s-ostatkom},
существуют $q\in\Z$ и $r\in\Z$ такие, что
$$
x=\triangledown(A)\cdot q+r,\qquad 0\le r<\triangledown(A)
$$
Если $r=0$, то мы сразу получаем, что $x\div\triangledown(A)$. Если же $r>0$,
то для всякого $a\in A$ мы получаем цепочку:
$$
\phantom{{\scriptsize \text{(уже доказано)}}}\quad\triangledown(A)\div
A\quad{\scriptsize \text{(уже доказано)}}
$$
$$
\Downarrow
$$
$$
\underbrace{x\div a,\quad \triangledown(A)\cdot q\div a,\quad
x>\triangledown(A)}
$$
$$
\phantom{{\scriptsize \text{$5^\circ$ на
с.\pageref{x:a,y:a=>x-y:a}}}}\quad\Downarrow\quad{\scriptsize \text{$4^\circ$
на с.\pageref{x:a,y:a=>x-y:a}}}
$$
$$
r=x-\triangledown(A)\cdot q\div a
$$
Это верно для любого $a\in A$, значит $r$ -- общее кратное для $A$:
$$
r\div A
$$
или, что то же самое,
$$
r\in A^\triangledown
$$
С другой стороны,
$$
r<\triangledown(A)=\min A^\triangledown
$$
Это противоречие означает, что наше предположение, что $r>0$, было неверно.

4. Докажем \eqref{vartriangle(A)=triangledown(A^vartriangle)}, то есть второе
звено в нашей цепочке. Для этого расшифруем уже доказанную формулу
\eqref{A^triangledown-div-triangledown(A)} так:
 \beq\label{A^triangledown-div-triangledown(A)-*}
x\div A\quad\Rightarrow\quad x\div \triangledown(A)
 \eeq
Тогда из \eqref{A-div-A^vartriangle} мы получаем цепочку:
 $$
A\div A^\vartriangle
 $$
$$
\Downarrow
$$
 $$
\forall a\in A\quad a\div A^\vartriangle
 $$
$$
\phantom{\scriptsize\eqref{A^triangledown-div-triangledown(A)-*}}
\quad\Downarrow\quad{\scriptsize\eqref{A^triangledown-div-triangledown(A)-*}}
$$
 $$
\forall a\in A\quad a\div \triangledown(A^\vartriangle)
 $$
$$
\Downarrow
$$
 $$
A\div \triangledown(A^\vartriangle)
 $$
$$
\Downarrow
$$
 \beq\label{triangledown(A^vartriangle)-in-A^vartriangle}
\triangledown(A^\vartriangle)\in A^\vartriangle
 \eeq
С другой стороны, из \eqref{triangledown(A)-div-A-div-vartriangle(A)} получаем:
$$
\triangledown(A^\vartriangle)\div A^\vartriangle
$$
$$
\Downarrow
$$
$$
\triangledown(A^\vartriangle)\ge A^\vartriangle
$$
Вместе с \eqref{triangledown(A^vartriangle)-in-A^vartriangle} это означает, что
$$
\triangledown(A^\vartriangle)=\max A^\vartriangle=\vartriangle(A),
$$
то есть как раз справедливо \eqref{vartriangle(A)=triangledown(A^vartriangle)}.

5. Теперь формула \eqref{vartriangle(A)-div-A^vartriangle}, то есть третье
звено в нашей цепочке, получается в одну строчку:
 $$
\vartriangle(A)\overset{\eqref{vartriangle(A)=triangledown(A^vartriangle)}}{=}
\triangledown(A^\vartriangle)\overset{\eqref{triangledown(A)-div-A-div-vartriangle(A)}}{\div}
A^\vartriangle
 $$

6. Остается доказать последнее звено цепочки -- формулу
\eqref{triangledown(A)=vartriangle(A^triangledown)}. Здесь последовательность
рассуждений напоминает пункт 4. Сначала расшифруем уже доказанную формулу
\eqref{vartriangle(A)-div-A^vartriangle} так:
 \beq\label{vartriangle(A)-div-A^vartriangle-*}
A\div x\quad\Rightarrow\quad \vartriangle(A)\div x
 \eeq
Тогда из \eqref{A^triangledown-div-A} мы получаем цепочку:
 $$
A^\triangledown\div A
 $$
$$
\Downarrow
$$
 $$
\forall a\in A\quad A^\triangledown\div a
 $$
$$
\phantom{\scriptsize\eqref{vartriangle(A)-div-A^vartriangle-*}}
\quad\Downarrow\quad{\scriptsize\eqref{vartriangle(A)-div-A^vartriangle-*}}
$$
 $$
\forall a\in A\quad \vartriangle(A^\triangledown)\div a
 $$
$$
\Downarrow
$$
 $$
\vartriangle(A^\triangledown)\div A
 $$
$$
\Downarrow
$$
 \beq\label{vartriangle(A^triangledown)-in-A^triangledown}
\vartriangle(A^\triangledown)\in A^\triangledown
 \eeq
С другой стороны, из \eqref{triangledown(A)-div-A-div-vartriangle(A)} получаем:
$$
A^\triangledown\div \vartriangle(A^\triangledown)
$$
$$
\Downarrow
$$
$$
A^\triangledown\ge \vartriangle(A^\triangledown)
$$
Вместе с \eqref{vartriangle(A^triangledown)-in-A^triangledown} это означает,
что
$$
\vartriangle(A^\triangledown)=\min A^\triangledown=\triangledown(A),
$$
то есть как раз справедливо
\eqref{triangledown(A)=vartriangle(A^triangledown)}.
 \epr

\btm Для любых $a,b\in\N$ справедливы импликации
 \begin{align}
&\left\{\begin{matrix}x\div a \\ x\div b
\end{matrix}\right\}\quad\Rightarrow\quad x\div (a\,\triangledown\, b) \label{x:a,x:b=>x:a-v-b}
\\
&\left\{\begin{matrix}a\div x \\ b\div x \end{matrix}
\right\}\quad\Rightarrow\quad (a\vartriangle b)\div x \label{a:x,b:x=>a-l-b:x}
\end{align}
и равенство
 \beq\label{a-v-b.a-l-b=a.b}
(a\,\triangledown\, b)\cdot (a\vartriangle b)=a\cdot b
 \eeq
\etm
 \bpr
1. Импликации следуют из соотношений \eqref{A^triangledown-div-triangledown-A}:
$$
A^\triangledown\div\triangledown(A)
$$
$$
\Downarrow
$$
$$
\{a;b\}^\triangledown\div\triangledown\{a;b\}
$$
$$
\Downarrow
$$
$$
x\in\{a;b\}^\triangledown\quad\Rightarrow\quad
x\div\triangledown\{a;b\}=(a\,\triangledown\, b)
$$
$$
\Downarrow
$$
$$
x\div\{a;b\}\quad\Rightarrow\quad x\div\triangledown\{a;b\}=(a\,\triangledown\,
b)
$$
$$
\Downarrow
$$
$$
\eqref{x:a,x:b=>x:a-v-b}
$$
И точно так же,
$$
\vartriangle(A)\div A^\vartriangle
$$
$$
\Downarrow
$$
$$
\vartriangle\{a;b\}\div \{a;b\}^\vartriangle
$$
$$
\Downarrow
$$
$$
x\in\{a;b\}^\vartriangle\quad\Rightarrow\quad \vartriangle\{a;b\}\div x
$$
$$
\Downarrow
$$
$$
x\div\{a;b\}\quad\Rightarrow\quad (a\vartriangle b)=\vartriangle\{a;b\}\div x
$$
$$
\Downarrow
$$
$$
\eqref{a:x,b:x=>a-l-b:x}
$$

2. Чтобы доказать \eqref{a-v-b.a-l-b=a.b}, обозначим
$$
x=\frac{a\cdot b}{a\,\triangledown\, b}
$$
и заметим, что $x\in\N$, потому что
$$
\underbrace{(a\cdot b)\div a,\quad (a\cdot b)\div b}
$$
$$
\phantom{\scriptsize \eqref{x:a,x:b=>x:a-v-b}}\quad\Downarrow\quad{\scriptsize
\eqref{x:a,x:b=>x:a-v-b}}
$$
$$
(a\cdot b)\div (a\,\triangledown\, b)
$$
Нам нужно убедиться, что
$$
x=a\vartriangle b
$$
Для этого отдельно докажем два неравенства
$$
a\vartriangle b\ge x, \qquad x\ge a\vartriangle b
$$
Первое из них доказывается такой цепочкой:
$$
\underbrace{a=x\cdot\frac{a\,\triangledown\, b}{b},\qquad
b=x\cdot\frac{a\,\triangledown\, b}{a}}
$$
$$
\Downarrow
$$
$$
\underbrace{a\div x,\qquad b\div x}
$$
$$
\phantom{\scriptsize \eqref{a:x,b:x=>a-l-b:x}}\quad\Downarrow\quad{\scriptsize
\eqref{a:x,b:x=>a-l-b:x}}
$$
$$
(a\vartriangle b)\div x
$$
$$
\Downarrow
$$
$$
a\vartriangle b\ge x
$$
А второе так:
$$
\underbrace{\frac{a\cdot b}{a\vartriangle b}=a\cdot\frac{b}{a\vartriangle
b},\qquad \frac{a\cdot b}{a\vartriangle b}=\frac{a}{a\vartriangle b}\cdot b}
$$
$$
\Downarrow
$$
$$
\underbrace{\frac{a\cdot b}{a\vartriangle b}\div a,\qquad \frac{a\cdot
b}{a\vartriangle b}\div b}
$$
$$
\phantom{\scriptsize \eqref{x:a,x:b=>x:a-v-b}}\quad\Downarrow\quad{\scriptsize
\eqref{x:a,x:b=>x:a-v-b}}
$$
$$
\frac{a\cdot b}{a\vartriangle b}\div (a\,\triangledown\, b)
$$
$$
\Downarrow
$$
$$
\frac{a\cdot b}{a\vartriangle b}\ge a\,\triangledown\, b
$$
$$
\Downarrow
$$
$$
x=\frac{a\cdot b}{a\,\triangledown\, b}\ge a\vartriangle b
$$
 \epr

\paragraph{Основная теорема арифметики.}

 \biter{
\item[$\bullet$] Число $a\in\N$, $a>1$, называется
 \biter{
\item[---] {\it простым}, если оно делится только на $1$ и на $a$:
$$
a^\vartriangle=\{1;a\}
$$

\item[---] {\it составным}, если кроме $1$ и $a$ у него есть еще какие-то
делители:
$$
a^\vartriangle\ne\{1;a\}
$$
 }\eiter
 }\eiter

\bex\label{EX:2-prostoe} Число $2$ является простым, потому что
$$
2\div x
$$
$$
\Downarrow
$$
$$
1\le x\le 2
$$
$$
\Downarrow
$$
$$
x\in\{1;2\}.
$$
То есть
$$
2^\vartriangle=\{1;2\}
$$
 \eex

 \biter{
\item[$\bullet$] Числа $a,b\in\N$ называются {\it взаимно простыми}, если они
имеют только один общий делитель, а именно, число $1$:
$$
\{a,b\}^\vartriangle=\{1\}
$$
Понятно, что это равносильно тому, что их наибольший общий делитель равен
единице:
$$
a\vartriangle b=1.
$$
Это в свою очередь означает, что
 \beq\label{a-v-b=a.b}
a\,\triangledown\, b=a\cdot b,
 \eeq
поскольку
$$
(a\,\triangledown\, b)\cdot \underbrace{(a\vartriangle
b)}_{1}\overset{\eqref{a-v-b.a-l-b=a.b}}{=}a\cdot b
$$
 }\eiter

\bex\label{EX:prostye=>vzaimno-prostye} Если $a$ и $b$ -- два различных простых
числа, то они взаимно просты, потому что
$$
a^\vartriangle=\{1;a\},\quad b^\vartriangle=\{1;b\}
$$
$$
\Downarrow
$$
 \begin{multline*}
\{a,b\}^\vartriangle=\eqref{(A-cup-B)^vartriangle=A^vartriangle-cap-B^vartriangle}=
\{a\}^\vartriangle\cap\{b\}^\vartriangle=\\=\{1;a\}\cap\{1;b\}=\{1\}
 \end{multline*}
 \eex

\bprop\label{PROP:xa:b=>x:b} Если $a$ и $b$ -- взаимно простые числа, то для
любого $x\in\N$
$$
(x\cdot a)\div b\qquad\Longleftrightarrow\qquad  x\div b
$$
\eprop
 \bpr
В обратную сторону это следует из $2^\circ$ на с.\pageref{x:y=>ax:y}:
$$
(x\cdot a)\div b\qquad\Longleftarrow\qquad  x\div b
$$
Покажем, что она верна в прямую сторону:
$$
(x\cdot a)\div b\qquad\Longrightarrow\qquad  x\div b
$$
Для этого нужно добавить к $(x\cdot a)\div b$ очевидное соотношение $(x\cdot
a)\div a$:
$$
(x\cdot a)\div b
$$
$$
\Downarrow
$$
$$
\underbrace{(x\cdot a)\div a\quad\&\quad  (x\cdot a)\div b}
$$
$$
\phantom{\scriptsize \eqref{x:a,x:b=>x:a-v-b}}\quad\Downarrow\quad{\scriptsize
\eqref{x:a,x:b=>x:a-v-b}}
$$
$$
(x\cdot a)\div a\,\triangledown\, b\overset{\eqref{a-v-b=a.b}}{=}a\cdot b
$$
$$
\Downarrow
$$
$$
\frac{x}{b}=\frac{x\cdot a}{a\cdot b}\in\N
$$
$$
\Downarrow
$$
$$
x\div b
$$
 \epr

\bprop\label{PROP:a_1...a_n:b=>a_i:b} Если произведение простых чисел
$$
a_1\cdot a_2\cdot...\cdot a_n
$$
делится на простое число $b$, то $b$ содержится в последовательности
$a_1,...,a_n$:
$$
\exists i\in\{1,...,n\}\qquad a_i=b
$$
 \eprop
\bpr Это доказывается индукцией по $n$.

1. При $n=1$ утверждение становится тривиальным:
$$
a_1\div b
$$
$$
\Downarrow
$$
$$
b\in \{a_1\}^\vartriangle=\{1;a_1\}
$$
$$
\phantom{\scriptsize (b>1)}\quad\Downarrow\quad{\scriptsize (b>1)}
$$
$$
b=a_1
$$

2. Предположим, что оно верно при $n=k$:
 \beq\label{predp-ind:a_1...a_k:b=>a_i=b}
a_1\cdot a_2\cdot...\cdot a_k\div b\quad\Rightarrow\quad \Big(\exists
i=1,...,k\qquad a_i=b \Big)
 \eeq

3. Покажем, что тогда оно будет верно и при $n=k+1$. Пусть
$$
a_1\cdot a_2\cdot...\cdot a_k\cdot a_{k+1}\div b
$$
Если $b=a_{k+1}$, то наше утверждение доказано. Если же $b\ne a_{k+1}$, то
обозначив $x=a_1\cdot a_2\cdot...\cdot a_k$, мы получим:
$$
x\cdot a_{k+1}\div b,
$$
причем $a_{k+1}$ и $b$ простые, не совпадающие числа, и значит, в силу примера
\ref{EX:prostye=>vzaimno-prostye}, взаимно простые. Поэтому по предложению
\ref{PROP:xa:b=>x:b}, $x$ должно делиться на $b$:
$$
x=a_1\cdot a_2\cdot...\cdot a_k\div b
$$
Но тогда по предположению индукции \eqref{predp-ind:a_1...a_k:b=>a_i=b},
$b=a_i$ для некоторого $i=1,...,k$. То есть в любом случае мы получаем $b=a_i$
для некоторого $i=1,...,k,k+1$.
 \epr

\btm[\bf основная теорема арифметики] Каждое число $x\in\N$, $x>1$, можно
единственным образом разложить в произведение
 \beq\label{razlozh-na-prostye-mnozhiteli}
x=p_1^{n_1}\cdot p_2^{n_2}\cdot...\cdot p_k^{n_k},
 \eeq
в котором $k$ и $\{n_1,n_2,...,n_k\}$ -- натуральные числа, а
$\{p_1,p_2,...,p_k\}$ -- простые, расположенные в порядке возрастания:
$$
p_1<p_2<...<p_k
$$
 \etm
 \bpr
Проведем индукцию по $x\in\N$, $x\ge 2$.

1. Пусть $x=2$. Тогда в виде \eqref{razlozh-na-prostye-mnozhiteli} его можно
представить так:
$$
2=2^1
$$
(то есть в данном случае $k=1$, $p_1=2$, $n_1=1$). Единственность такого
разложения следует из того, что, как мы убедились в примере \ref{EX:2-prostoe},
число 2 простое: если
 \beq\label{2=q_1^(m_1)-q_l^(m_l)}
2=q_1^{m_1}\cdot q_2^{m_2}\cdot...\cdot q_l^{m_l}
 \eeq
-- какое-нибудь другое разложение, то мы во-первых получим, что все числа $q_i$
должны быть одинаковыми и равными двойке, потому что
$$
2\div q_i
$$
$$
\Downarrow
$$
$$
q_i\in 2^\vartriangle=\{1;2\}
$$
$$
\phantom{\scriptsize (q_i>1)}\quad\Downarrow\quad{\scriptsize (q_i>1)}
$$
$$
q_i=2
$$
То есть разложение \eqref{2=q_1^(m_1)-q_l^(m_l)} должно иметь вид
$$
2=2^m
$$
И, во-вторых, здесь $m$ должно быть равно 1, потому что иначе, то есть при
$m>1$, мы получили бы
$$
2=2^1\overset{\eqref{monot-a^n-a>1}}{<}2^m
$$

2. Предположим далее, что наше утверждение верно при всех $x$, меньших
некоторого $N$,
$$
x<N.
$$
Покажем тогда, что оно верно и при $x=N$.

Здесь придется рассмотреть два случая. Во-первых, если $x=N$ -- простое число,
то наше утверждение для него верно по тем же причинам, по которым оно было
верно для уже рассмотренного случая $x=2$: те же рассуждения, но с заменой 2 на
$N$ приводят к тому же результату.

Во-вторых, если $x=N$ -- составное число, то есть
$$
N=K\cdot L\qquad (K,L\in\N,\ K,L>1)
$$
то
$$
K<N,\qquad L<N
$$
поэтому по предположению индукции числа $K$ и $L$ раскладываются по формуле
\eqref{razlozh-na-prostye-mnozhiteli}. Отсюда можно сделать вывод, что их
произведение $N=K\cdot L$ тоже  раскладывается по формуле
\eqref{razlozh-na-prostye-mnozhiteli}.

Остается только объяснить, почему для такого $N$ разложение будет единственно.
Предположим, что этих разложений имеется два:
 \beq\label{edinstv-razlozh-na-prostye-mnozh}
N=p_1^{n_1}\cdot ...\cdot p_k^{n_k}=q_1^{m_1}\cdot ...\cdot q_l^{m_l}
 \eeq
Поглядим на это с точки зрения предложения \ref{PROP:a_1...a_n:b=>a_i:b}: число
$N$ делится на простое число $q_1$, а, с другой стороны, $N$ является
произведением простых чисел $p_1,...,p_k$ (с разными степенями). Значит, среди
чисел $p_1,...,p_k$ какое-то должно быть равно $q_1$:
$$
\exists i\in\{1,...,k\}\qquad p_i=q_1
$$
По той же причине среди чисел $q_1,...,q_l$ какое-то должно быть равно $p_1$:
$$
\exists j\in\{1,...,l\}\qquad q_j=p_1
$$
Поскольку числа $p_1,...,p_k$ упорядочены по возрастанию, получаем
$$
p_1\le p_i=q_1
$$
С другой стороны, числа $q_1,...,q_l$ тоже упорядочены по возрастанию, и
поэтому
$$
q_1\le q_j=p_1
$$
Вместе это означает, что
$$
p_1=q_1
$$
откуда следует, что $i=1=j$.

Теперь разложение \eqref{edinstv-razlozh-na-prostye-mnozh} можно записать так:
$$
N=p_1^{n_1}\cdot p_2^{n_2}\cdot...\cdot
p_k^{n_k}=\underbrace{p_1^{m_1}}_{\scriptsize\begin{matrix} \text{$q_1$}\\
\text{заменено} \\ \text{на $p_1$}\end{matrix}}\cdot q_2^{n_2}\cdot ...\cdot
q_l^{m_l}
$$
Отсюда можно вывести, что $n_1=m_1$. Действительно, если бы оказалось, что
$n_1>m_1$, то поделив на $p_1^{m_1}$, мы получили бы
$$
\frac{N}{p_1^{m_1}}=p_1^{n_1-m_1}\cdot p_2^{n_2}\cdot...\cdot
p_k^{n_k}=q_2^{m_2}\cdot ...\cdot q_l^{m_l}
$$
и оказалось бы, что это число, будучи произведением простых чисел $q_2,...,q_l$
(с разными степенями), делится на простое число $p_1$, которое не принадлежит
последовательности  $q_2,...,q_l$ (потому что $p_1=q_1$ меньше любого из этих
чисел). И точно так же, если бы $n_1<m_1$, то поделив на $p_1^{n_1}$, мы
получили бы
$$
\frac{N}{p_1^{n_1}}=p_2^{n_2}\cdot...\cdot p_k^{n_k}=p_1^{m_1-n_1}\cdot
q_2^{m_2}\cdot ...\cdot q_l^{m_l}
$$
то есть произведение простых чисел $p_2,...,p_k$ (с разными степенями),
делилось бы на простое число $p_1$, которое не принадлежит последовательности
$p_2,...,p_k$ (потому что $p_1$ меньше любого из этих чисел).

Итак, $n_1=m_1$, и поэтому разложение \eqref{edinstv-razlozh-na-prostye-mnozh}
можно переписать так:
$$
N=p_1^{n_1}\cdot p_2^{n_2}\cdot...\cdot p_k^{n_k}=\underbrace{p_1^{n_1}}_{\scriptsize\begin{matrix} \text{$q_1$ и $m_1$}\\
\text{заменены} \\ \text{на $p_1$ и $n_1$}\end{matrix}}\cdot q_2^{n_2}\cdot
...\cdot q_l^{m_l}
$$
Поделим это на $p_1^{n_1}$:
$$
\frac{N}{p_1^{n_1}}=p_2^{n_2}\cdot...\cdot p_k^{n_k}=q_2^{m_2}\cdot ...\cdot
q_l^{m_l}
$$
Из $p_1>1$ и $n_1>0$ следует, что $p_1^{n_1}>1$, поэтому
$$
\frac{N}{p_1^{n_1}}<N
$$
и по предположению индукции, для числа $\frac{N}{p_1^{n_1}}$ разложение
\eqref{razlozh-na-prostye-mnozhiteli} должно быть единственно. То есть,
$$
k=l\quad\&\quad \Big(\forall i=2,...,k\quad p_i=q_i\quad\&\quad n_i=m_i\Big)
$$
С другой стороны, равенства $p_1=q_1$ и $n_1=m_1$ мы уже доказали раньше. То
есть в формуле \eqref{edinstv-razlozh-na-prostye-mnozh} числа справа и слева от
последнего равенства совпадают, и это значит, что разложение $N$ на множители
однозначно.
 \epr

\end{multicols}\noindent\rule[10pt]{160mm}{0.1pt}

\subsection{Рациональные числа $\Q$}\label{SEC-ratsionalnye-chisla}

\paragraph{Определение и свойства рациональных чисел.}

 \bit{
\item[$\bullet$] Число $x\in\R$ называется {\it рациональным}, если оно
представимо в виде
$$
x=\frac{a}{b}, \qquad a\in\Z,\quad b\in\N
$$
В соответствии с этим, множество рациональных чисел, обозначаемое символом
$\Q$, описывается формулой
 \beq\label{DF:Q}
\Q=\frac{\Z}{\N}
 \eeq
 }\eit

 \bigskip
\centerline{\bf Свойства множества $\Q$:}
 \bit{\it

\item[$1^\circ$.] Сумма $x+y$ двух рациональных чисел $x$ и $y$, также является
рациональным числом:
 \beq
x,y\in \Q\qquad\Longrightarrow\qquad x+y\in \Q
 \eeq

\item[$2^\circ$.] Число $-x$, противоположное рациональному числу $x$, также
является рациональным числом:
 \beq
x\in \Q\qquad\Longrightarrow\qquad -x\in \Q
 \eeq

\item[$3^\circ$.] Произведение $x\cdot y$ двух рациональных чисел $x$ и $y$,
также является рациональным числом:
 \beq
x,y\in \Q\quad\Longrightarrow\quad x\cdot y\in \Q
 \eeq

\item[$4^\circ$.] Число $x^{-1}$, обратное рациональному числу  $x\ne 0$, также
является рациональным числом:
 \beq
x\in \Q\qquad\Longrightarrow\qquad x^{-1}\in \Q
 \eeq

 }\eit

\btm В любом интервале $(a,b)$ на вещественной прямой $\R$ найдутся
рациональные числа. \etm
 \bpr
Пользуясь принципом Архимеда (теорема \ref{Archimed-principle}), подберем
число $n\in\N$ такое, что
$$
n>\frac{1}{b-a}
$$
и положим
$$
m=[a\cdot n]+1
$$
Тогда
$$
a<\frac{m}{n}<b
$$
Первое неравенство здесь доказывается напрямую:
$$
a\cdot n\overset{\eqref{opr-tsel-chasti}}{<}[a\cdot
n]+1=m\quad\Longrightarrow\quad a<\frac{m}{n}
$$
А правое -- от противного: если бы оказалось, что $b\le\frac{m}{n}$, то мы
получили бы
 \begin{multline*}
[a\cdot n]\overset{\eqref{opr-tsel-chasti}}{<}a\cdot n
\quad\Longrightarrow\quad \frac{[a\cdot n]}{n}\le
a<b\le\frac{m}{n}=\frac{[a\cdot n]+1}{n}=\frac{[a\cdot n]}{n}+\frac{1}{n}
\quad\Longrightarrow\\ \Longrightarrow\quad b-a\le\frac{1}{n}
\quad\Longrightarrow\quad n\le \frac{1}{b-a},
 \end{multline*}
и последнее противоречит выбору $n$.
 \epr

\noindent\rule{160mm}{0.1pt}\begin{multicols}{2}

\paragraph{Несократимые дроби.}

\btm\label{TH:nesokratimoe-predstavlenie} Всякое положительное рациональное
число $x\in\Q$ представимо в виде дроби, в которой числитель и знаменатель
взаимно просты:
 \beq
x=\frac{p}{q}, \qquad p,q\in\N,\qquad p\vartriangle q=1
 \eeq
\etm
 \bit{
\item[$\bullet$] Такие дроби называются {\it несократимыми}.
 }\eit
\bpr Поскольку $x$ рационально, оно имеет вид $x=\frac{a}{b}$, где $a\in\Z$ и
$b\in\N$. С другой стороны, $x$ положительно, поэтому $a$ тоже должно быть
положительно, и значит $a\in\N$. Поскольку $a\div (a\vartriangle b)$ и $b\div
(a\vartriangle b)$, числа
$$
p=\frac{a}{a\vartriangle b},\qquad q=\frac{b}{a\vartriangle b}
$$
должны быть натуральными. Понятно, что
$$
x=\frac{a}{b}=\frac{\frac{a}{a\vartriangle b}}{\frac{b}{a\vartriangle
b}}=\frac{p}{q}
$$
Остается проверить, что $p\vartriangle q=1$. Это делается так:
 $$
x\in\{p;q\}^\vartriangle
 $$
 $$
 \Downarrow
 $$
 $$
\left\{\begin{matrix} p\div x\\
q\div x\end{matrix}\right\}
 $$
 $$
 \Downarrow
 $$
 $$
\left\{\begin{matrix}
a=p\cdot(a\vartriangle b) \div x\cdot(a\vartriangle b)\\
b=q\cdot(a\vartriangle b)\div x\cdot(a\vartriangle b)
\end{matrix}\right\}
 $$
 $$
 \Downarrow
 $$
 $$
x\cdot(a\vartriangle b)\in\{a;b\}^\vartriangle
 $$
 $$
 \Downarrow
 $$
 $$
x\cdot(a\vartriangle b)\le \max\{a;b\}^\vartriangle=(a\vartriangle
b)
 $$
 $$
 \Downarrow
 $$
 $$
x\le 1
 $$
Из этой цепочки видно, что
$$
p\vartriangle q=\max\{p;q\}^\vartriangle=1
$$
\epr

\paragraph{Существование иррациональных чисел.}

 \biter{
\item[$\bullet$] Число $x\in\R$ называется {\it иррациональным}, если оно не
является рациональным:
$$
x\notin\Q
$$
 }\eiter\noindent
Следующий пример показывает, что такие числа в самом деле существуют.

\bex\label{EX:sqrt(2)-irrat} Существует число $c\in[1;2]$, удовлетворяющее
условию\footnote{Читатель, конечно, узнает в числе $c$ из примера
\ref{EX:sqrt(2)-irrat} число $\sqrt{2}$. Мы не вставляем это в формулировку,
потому что понятие корня будет определено нами только в
\ref{SUBSEC-elem-funktsii} главы \ref{ch-ELEM-FUNCTIONS}.}
 \beq\label{c^2=2}
c^2=2
 \eeq
Это число иррационально:
 $$
c\notin\Q
 $$
 \eex
 \bpr
1. Рассмотрим два множества:
$$
X=\{x\in\R: \quad x>0\quad\&\quad x^2\le 2\}
$$
$$
Y=\{y\in\R: \quad y>0\quad\&\quad y^2\ge 2\}
$$
Заметим, что $1\in X$ и $2\in Y$, поэтому $X$ и $Y$ непустые. Докажем
неравенство для множеств
 $$
X\le Y
 $$
Действительно, если бы для каких-то $x\in X$ и $y\in Y$ выполнялось обратное
неравенство
$$
x>y,
$$
то, поскольку $y>0$, по свойству степенного отображения \eqref{monot-a^n-n>0},
мы получили бы
$$
x^2>y^2
$$
при том что из условий $x\in X$ и $y\in Y$ следует
$$
x^2\le 2\le y^2.
$$

Из неравенства $X\le Y$, в силу аксиомы непрерывности A16, следует, что
найдется число $c\in\R$, лежащее между $X$ и $Y$:
$$
X\le c\le Y
$$
Заметим сразу, что, поскольку $1\in X$ и $2\in Y$, число $c$ должно
удовлетворять условию
 \beq\label{1-le-c-le-2}
1\le c\le 2
 \eeq

2. Покажем, что
$$
c^2=2
$$
Предположим, что это не так, например, пусть
 \beq\label{sqrt(2)-hyp-1}
c^2<2
 \eeq
Тогда число $\e=2-c^2$ будет положительно:
$$
\e=2-c^2>0
$$
Отсюда сразу следует, что
$$
c+\frac{\e}{3 c}>c\ge X
$$
и следовательно число $c+\frac{\e}{3 c}$ не лежит в $X$. Но с другой стороны,
 $$
 \underbrace{
\begin{matrix}
{\scriptsize\eqref{1-le-c-le-2}} \\
\Downarrow \\
1\le c^2 \\
\Downarrow \\
-c^2\le-1 \\
\Downarrow \\
\e=2-c^2\le 2-1=1
\end{matrix}
\qquad
\begin{matrix}
{\scriptsize\eqref{1-le-c-le-2}} \\
 \Downarrow \\
 0<\frac{1}{c}\le 1 \\
 \Downarrow \\
 0<\frac{1}{c^2}\le 1
\end{matrix}}
 $$
 $$
 \Downarrow
 $$
 \begin{multline*}
\left(c+\frac{\e}{3 c}\right)^2=c^2+2\cdot\frac{c\cdot \e}{3\cdot
c}+\left(\frac{\e}{3
c}\right)^2=\\=c^2+\frac{2}{3}\cdot\e+\underbrace{\frac{1}{9}}_{\scriptsize\begin{matrix}\text{\rotatebox{-90}{$\le$}}\\
\frac{1}{3}\end{matrix}}\cdot\underbrace{\frac{1}{c^2}}_{\scriptsize\begin{matrix}\text{\rotatebox{-90}{$\le$}}\\
1\end{matrix}}\cdot\underbrace{\e^2}_{\scriptsize\begin{matrix}\text{\rotatebox{-90}{$\le$}}\\
\phantom{,}\e,\\ \text{так как}\\ 0<\e\le 1\end{matrix}} \le
c^2+\frac{2}{3}\cdot\e+\frac{1}{3}\cdot\e=\\= c^2+\e=c^2+(2-c^2)=2
 \end{multline*}
То есть
$$
c+\frac{\e}{3 c}>c>0\quad\&\quad  \left(c+\frac{\e}{3 c}\right)^2\le 2
$$
и это значит, что число $c+\frac{\e}{3 c}$ наоборот лежит в $X$. Это
противоречие означает, что наше предположение \eqref{sqrt(2)-hyp-1} было
неверным.

Предположим, наоборот, что
 \beq\label{sqrt(2)-hyp-2}
c^2>2
 \eeq
Тогда число $\e=c^2-2$ будет положительно
$$
\e=c^2-2>0
$$
Отсюда сразу следует, что
$$
c-\frac{\e}{3 c}<c\le Y
$$
и поэтому число $c-\frac{\e}{3 c}$ не лежит в $Y$. Но с другой стороны,
во-первых,
$$
c-\frac{\e}{3 c}=c-\frac{c^2-2}{3c}=\frac{2c^2+2}{3c}>0
$$
И, во-вторых,
 \begin{multline*}
\left(c-\frac{\e}{3 c}\right)^2=c^2-\frac{2\cdot \e}{3}+\left(\frac{\e}{3
c}\right)^2=\\= c^2-\e+\underbrace{\frac{\e}{3}+\left(\frac{\e}{3
c}\right)^2}_{\scriptsize\begin{matrix}\text{\rotatebox{90}{$<$}}\\
0\end{matrix}}>c^2-\e=c^2-(c^2-2)=2
 \end{multline*}
То есть мы получаем
$$
c-\frac{\e}{3 c}>0\quad\&\quad  \left(c-\frac{\e}{3 c}\right)^2\ge 2
$$
и это значит, что число $c-\frac{\e}{3 c}$ наоборот лежит в $Y$. Это
противоречие говорит о том, что наше предположение \eqref{sqrt(2)-hyp-2} тоже
было неверным.

3. Нам осталось доказать, что число $c$ иррационально. Предположим, что оно
рационально:
$$
c\in\Q
$$
Тогда, поскольку $c>0$, по теореме \ref{TH:nesokratimoe-predstavlenie} это
число можно представить в виде несократимой дроби
$$
c=\frac{m}{n}, \qquad m\vartriangle n=1
$$
Мы получаем цепочку:
$$
\frac{m^2}{n^2}=c^2=2
$$
$$
\Downarrow
$$
$$
m^2=2\cdot n^2
$$
$$
\Downarrow
$$
$$
m^2\div 2
$$
$$
\Downarrow
$$
$$
 \begin{matrix}
\text{в разложении \eqref{razlozh-na-prostye-mnozhiteli} числа $m^2$}\\
\text{имеется множитель $2^s$, $s\in\N$}
 \end{matrix}
$$
$$
\Downarrow
$$
$$
 \begin{matrix}
\text{в разложении \eqref{razlozh-na-prostye-mnozhiteli} числа $m$}\\
\text{имеется множитель $2^t$, $t\in\N$}
 \end{matrix}
$$
$$
\Downarrow
$$
$$
m\div 2
$$
$$
\Downarrow
$$
$$
m=2\cdot k,\qquad k>1
$$
$$
\Downarrow
$$
$$
4\cdot k^2=m^2=2\cdot n^2
$$
$$
\Downarrow
$$
$$
2\cdot k^2=n^2
$$
$$
\Downarrow
$$
$$
n^2\div 2
$$
$$
\Downarrow
$$
$$
 \begin{matrix}
\text{по тем же причинам, что и для $m$,}\\
\text{в разложении \eqref{razlozh-na-prostye-mnozhiteli} числа $n$}\\
\text{должен быть множитель $2^r$, $r\in\N$}
 \end{matrix}
$$
$$
\Downarrow
$$
$$
n\div 2
$$
Мы получаем, что оба числа $m$ и $n$ делятся на 2. Это противоречит тому, что у
них наибольший общий делитель равен 1. Это означает, что наше предположение о
рациональности $c$ было ошибочным.
 \epr

\paragraph{Рациональные числа с нечетным знаменателем
$\frac{\Z}{2\N-1}$.}\label{Z/2N-1}

Напомним, что выше формулой \eqref{umnozh-mnozhestva-na-chislo} мы определили
умножение числового множества $X$ на число $a$:
$$
a\cdot X=\{y\in\R:\;\exists x\in X\quad a\cdot x=y\}
$$
формулой \eqref{sdvigi-mnozhestv-} -- сдвиг множества
$$
X-a=\{y\in\R:\;\exists x\in X\quad x-a=y\},
$$
а формулой \eqref{DEF:X/Y} -- частное числовых множеств:
$$
\frac{X}{Y}=\left\{z\in\R:\;\exists x\in X\quad \exists y\in Y\quad
\frac{x}{y}=z\right\}
$$
Этого достаточно, чтобы сообразить, что понимается под множеством
$\frac{\Z}{2\N-1}$, о котором мы поведем речь здесь, но перед тем, как
приступить собственно к изучению его свойств, мы рассмотрим несколько
подготовительных примеров.

\bex Множество четных натуральных чисел $\{2,4,6,8,...\}=\{2n;\ n\in\N\}$
удобно записывать в виде
$$
2\N.
$$
Точно также множество нечетных натуральных чисел $\{1,3,5,7,...\}=\{2n-1;\
n\in\N\}$ удобно записывать в виде
$$
2\N-1.
$$
 \eex

\bex Выше мы уже приводили формулу \eqref{DF:Q}
$$
\Q=\frac{\Z}{\N},
$$
которую можно рассматривать, как определение множества рациональных чисел $\Q$.
 \eex

\bex\label{EX:Q=Z/2N} Непривычному взгляду не сразу будет очевидной
справедливость равенства
 \beq\label{Q=Z/2N}
\Q=\frac{\Z}{2\N}
 \eeq
 \eex
 \bpr Ясно, что выполняется включение
 $$
\frac{\Z}{2\N}\subseteq\Q
 $$
Покажем, что верно и обратное включение:
 $$
\Q\subseteq\frac{\Z}{2\N}
 $$
Действительно, если $r=\frac{p}{q}\in\Q$, $p\in\Z$, $q\in\N$, то положив
$m=2p$, мы получим:
$$
r=\frac{p}{q}=\frac{m}{2q},\quad m\in\Z,\ q\in\N
$$
То есть $r\in\frac{\Z}{2\N}$.
 \epr

\ber Докажите равенства:
 \begin{align*}
&\frac{2\Z}{2\N}=\Q
 \\
&\frac{\N}{\N}=\{r\in\Q:\quad r>0\}
 \end{align*}
 \eer

\bex После примера \ref{EX:Q=Z/2N} можно было бы ожидать, что если заменить в
\eqref{Q=Z/2N} четные числа $2\N$ на нечетные $2\N-1$, то равенство останется
верным. Но это не так:
$$
\Q\ne\frac{\Z}{2\N-1}
$$
 \eex
\bpr Число $\frac{1}{2}$ лежит в $\Q$, но не в $\frac{\Z}{2\N-1}$ (потому что
его невозможно представить в виде $\frac{\text{целое}}{\text{нечетное}}$). \epr

 \biter{
\item[$\bullet$] Число $x\in\Q$ называется {\it рациональным числом с нечетным
знаменателем}, если его можно представить в виде дроби, в которой числитель -
целое число, а знаменатель -- нечетное натуральное число:
 $$
x=\frac{m}{2n-1},\qquad m\in\Z,\ n\in\N
 $$
Нетрудно сообразить, что это равносильно тому, что, будучи представлено в виде
несократимой дроби, $x$ имеет нечетный знаменатель. Множество всех рациональных
чисел с нечетным знаменателем описывается формулой
$$
\frac{\Z}{2\N-1}
$$
 }\eiter

 \bigskip
\centerline{\bf Свойства множества $\frac{\Z}{2\N-1}$:}
 \biter{\it

\item[$1^\circ$.] Сумма $x+y$ двух рациональных чисел  $x$ и $y$ с нечетным
знаменателем, также является рациональным с нечетным знаменателем:
 $$
x,y\in \frac{\Z}{2\N-1}\quad\Longrightarrow\quad x+y\in \frac{\Z}{2\N-1}
 $$

\item[$2^\circ$.] Число $-x$, противоположное рациональному числу $x$ с
нечетным знаменателем, также является рациональным с нечетным знаменателем:
 $$
x\in \frac{\Z}{2\N-1}\quad\Longrightarrow\quad -x\in \frac{\Z}{2\N-1}
 $$

\item[$3^\circ$.] Произведение $x\cdot y$ двух рациональных чисел $x$ и $y$ с
нечетным знаменателем, также является рациональным с нечетным знаменателем:
 $$
x,y\in \frac{\Z}{2\N-1}\quad\Longrightarrow\quad x\cdot y\in \frac{\Z}{2\N-1}
 $$
 }\eiter

\brem Можно заметить, что для рациональных чисел с четным знаменателем, которые
можно определить как те, которые будучи представлены в виде несократимой дроби
имеют четный знаменатель (и, возможно, знак минус перед дробью), утверждение
$2^\circ$ в этом списке не будет справедливо: равенство
$$
\frac{1}{2}+\frac{1}{2}=1
$$
означает, что множество таких чисел не замкнуто относительно сложения. Однако,
утверждения $1^\circ$ и  $3^\circ$ для этого класса все же будут выполняться.
 \erem

 \biter{

\item[$\bullet$] Рациональное число $x$ с нечетным знаменателем называется

 \biter{
\item[---] {\it четным}, если его можно представить в виде дроби с четным
числителем и нечетным знаменателем:
$$
x=\frac{2m}{2n-1},\qquad m\in\Z,\quad n\in\N
$$
это равносильно тому, что в разложении $x$ на несократимую дробь числитель
является четным числом; очевидно, что множество четных чисел с нечетным
знаменателем описывается формулой:
$$
\frac{2\Z}{2\N-1}
$$

\item[---] {\it нечетным}, если его нельзя представить в виде дроби с четным
числителем и нечетным знаменателем; это равносильно тому, что при любом
представлении в виде дроби с нечетным знаменателем, числитель тоже будет
нечетным:
$$
x=\frac{2m-1}{2n-1},\qquad m\in\Z,\quad n\in\N
$$
в частности, в разложении $x$ на несократимую дробь числитель должен быть
нечетным числом; множество нечетных чисел с нечетным знаменателем описывается
формулой:
$$
\frac{2\Z-1}{2\N-1}
$$
 }\eiter
 }\eiter

Доказательство следующей теоремы мы оставляем читателю:

\btm Множества $\frac{2\Z}{2\N-1}$  и $\frac{2\Z-1}{2\N-1}$ не пересекаются и в
объединении дают множество $\frac{\Z}{2\N-1}$:
 \begin{align*}
&\frac{2\Z}{2\N-1}\cap \frac{2\Z-1}{2\N-1}=\varnothing
\\
&\frac{2\Z}{2\N-1}\cup \frac{2\Z-1}{2\N-1}=\frac{\Z}{2\N-1}
 \end{align*}
\etm

 \biter{
\item[$\bullet$] {\it Четностью} рационального числа $x$ с нечетным
знаменателем называется число $\sgn_2 x$, определяемое правилом:
 \beq\label{sgn_2}
\sgn_2 x=\begin{cases} +1,& \text{$x$ -- четное}\\ -1,& \text{$x$ -- нечетное}
\end{cases}
 \eeq
 }\eiter

Напомним, что выше (см.\eqref{DEF:a-vee-b}) мы условились записью $a\vee b$
обозначать операцию взятия максимума двух чисел:
$$
a\vee b=\max\{a,b\}
$$
Такая запись позволяет придать ``алгебраический'' вид последнему из тождеств в
следующем списке ($x,y\in\frac{\Z}{2\N-1}$):
 \begin{align}
& \sgn_2 x=\frac{1}{\sgn_2 x}\label{sgn_2=1/sgn_2} \\
& \sgn_2(-x)=\sgn_2 x \label{sgn_2(-x)=sgn_2(x)} \\
& \sgn_2(x+y)=\sgn_2 x\cdot\sgn_2 y \label{sgn_2(x+y)=sgn_2(x)-sgn_2(y)} \\
& \sgn_2(x\cdot y)=\sgn_2 x\ \vee\ \sgn_2 y
\label{sgn_2(xy)=max(sgn_2(x),sgn_2(y))}
 \end{align}
 \bpr
Здесь нужно просто рассмотреть несколько случаев. Для двух последних тождеств
рассуждения выглядят так:

1) если $x,y\in\frac{2\Z}{2\N-1}$, то $x+y\in\frac{2\Z}{2\N-1}$ и $x\cdot
y\in\frac{2\Z}{2\N-1}$, поэтому
$$
\underbrace{\sgn_2(x+y)}_{\scriptsize\begin{matrix}\text{\rotatebox{90}{$=$}}\\
1\end{matrix}}=\underbrace{\sgn_2x}_{\scriptsize\begin{matrix}\text{\rotatebox{90}{$=$}}\\
1\end{matrix}}\cdot\underbrace{\sgn_2y}_{\scriptsize\begin{matrix}\text{\rotatebox{90}{$=$}}\\
1\end{matrix}},
$$
$$
\underbrace{\sgn_2(x\cdot y)}_{\scriptsize\begin{matrix}\text{\rotatebox{90}{$=$}}\\
1\end{matrix}}=\underbrace{\sgn_2x}_{\scriptsize\begin{matrix}\text{\rotatebox{90}{$=$}}\\
1\end{matrix}}\ \vee \ \underbrace{\sgn_2y}_{\scriptsize\begin{matrix}\text{\rotatebox{90}{$=$}}\\
1\end{matrix}},
$$

2) если $x\in\frac{2\Z}{2\N-1}$, $y\in\frac{2\Z-1}{2\N-1}$, то
$x+y\in\frac{2\Z-1}{2\N-1}$ и $x\cdot y\in\frac{2\Z}{2\N-1}$, поэтому
$$
\underbrace{\sgn_2(x+y)}_{\scriptsize\begin{matrix}\text{\rotatebox{90}{$=$}}\\
-1\end{matrix}}=\underbrace{\sgn_2x}_{\scriptsize\begin{matrix}\text{\rotatebox{90}{$=$}}\\
1\end{matrix}}\cdot\underbrace{\sgn_2y}_{\scriptsize\begin{matrix}\text{\rotatebox{90}{$=$}}\\
-1\end{matrix}},
$$
$$
\underbrace{\sgn_2(x\cdot y)}_{\scriptsize\begin{matrix}\text{\rotatebox{90}{$=$}}\\
1\end{matrix}}=\underbrace{\sgn_2x}_{\scriptsize\begin{matrix}\text{\rotatebox{90}{$=$}}\\
1\end{matrix}}\ \vee \ \underbrace{\sgn_2y}_{\scriptsize\begin{matrix}\text{\rotatebox{90}{$=$}}\\
-1\end{matrix}},
$$

3) если $x,y\in\frac{2\Z-1}{2\N-1}$, то $x+y\in\frac{2\Z}{2\N-1}$ и $x\cdot
y\in\frac{2\Z-1}{2\N-1}$, поэтому
$$
\underbrace{\sgn_2(x+y)}_{\scriptsize\begin{matrix}\text{\rotatebox{90}{$=$}}\\
1\end{matrix}}=\underbrace{\sgn_2x}_{\scriptsize\begin{matrix}\text{\rotatebox{90}{$=$}}\\
-1\end{matrix}}\cdot\underbrace{\sgn_2y}_{\scriptsize\begin{matrix}\text{\rotatebox{90}{$=$}}\\
-1\end{matrix}},
$$
$$
\underbrace{\sgn_2(x\cdot y)}_{\scriptsize\begin{matrix}\text{\rotatebox{90}{$=$}}\\
-1\end{matrix}}=\underbrace{\sgn_2x}_{\scriptsize\begin{matrix}\text{\rotatebox{90}{$=$}}\\
-1\end{matrix}}\ \vee \ \underbrace{\sgn_2y}_{\scriptsize\begin{matrix}\text{\rotatebox{90}{$=$}}\\
-1\end{matrix}}.
$$
 \epr

\end{multicols}\noindent\rule[10pt]{160mm}{0.1pt}

\chapter{ПРЕДЕЛ ПОСЛЕДОВАТЕЛЬНОСТИ}
\label{ch-x_n}

\section{Числовые последовательности и их пределы}

В примере \ref{EX-posledovatelnosti} выше мы уже вводили понятие (бесконечной)
последовательности, хотя в тот момент это нарушало логику изложения, поскольку
определение было дано раньше, чем был определен натуральный ряд $\N$. Теперь
самое время вернуться к последовательностям снова.

\subsection{Числовые последовательности}

 \bit{\label{DEF:chislov-posledovatelnost}
\item[$\bullet$] Пусть каждому натуральному числу $n\in\N$ по какому-нибудь
правилу поставлено в соответствие вещественное число $x_n$ (единственное для
данного $n$). Тогда такое соответствие (семейство $\{ x_n \}$) называется {\it
(бесконечной) числовой последовательностью}, или просто {\it
последовательностью}.
 }\eit

\noindent  Таким образом, формально числовая последовательность -- просто
отображение с областью определения $\N$ и областью значений $\R$ (об этом мы
тоже говорили в примере \ref{EX-posledovatelnosti}).

\paragraph{Способы описания числовой последовательности.}

В математическом анализе имеется только три способа определить числовую
последовательность. Мы перечислим их в следующих примерах.

\noindent\rule{160mm}{0.1pt}\begin{multicols}{2}

\bex {\bf Последовательности, заданные явно.}\label{EX:posl-zadan-yavno} Первый
способ -- написать формулу вида
 \beq\label{zadanie-posledov-yavno}
x_n=f(n),
 \eeq
в которой $f$ -- заранее определенная функция. Например, следующие формулы {\it
явно задают последовательности}, поскольку функции в правых частях были уже
нами определены:
$$
x_n=\frac{1}{n},\qquad x_n=n,\qquad x_n=(-1)^n.
$$
Естественно, под $f$ в формуле \eqref{zadanie-posledov-yavno} можно также
понимать какую-нибудь функцию, составленную из уже определенных, с помощью явно
описанных операций, например, операции композиции,
$$
x_n=2^{n^2}
$$
или алгебраических операций,
$$
x_n=n^2-n
$$
или операции индуктивного суммирования,
$$
x_n=\sum_{k=1}^n\frac{1}{k}
$$
или каких-то других операций (например, дифференцирования и интегрирования
которые мы опишем ниже в главах \ref{ch-f'(x)} и \ref{CH-definite-integral}).
Во всех этих случаях считается, что формула задает последовательность явно.
\eex

\bex{\bf Последовательности, заданные рекуррентно.}\label{EX:rekurr-posledov}
Второй способ -- задать последовательность индуктивно с помощью теоремы
\ref{defin-induction} (об определениях полной индукцией). В таких случаях говорят также, что последовательность
задана {\it рекуррентно}. Например, правило
$$
x_1=2,\quad x_{n+1}=\frac{1}{2}\left(x_n+\frac{1}{x_n}\right)
$$
индуктивно задает некую числовую последовательность (ниже в примере \ref{x_(n+1)=1/2(x_n+1/x_n)} мы будем изучать ее свойства).

Другой пример -- знаменитая {\it последовательность Фибоначчи}, задаваемая
правилом
 \beq\label{DEF:Fibonacci}
x_1=x_2=1,\quad x_{n+2}=x_{n+1}+x_n,
 \eeq
(ниже в главе \ref{CH-step-ryady} мы докажем неочевидный факт, что эту
последовательность можно задать и явной формулой
\eqref{formula-dlya-Fibonacci}). Отметим, что в этом примере ссылка на теорему \ref{defin-induction}, из
которой мы выводили саму возможность определений по индукции не вполне
корректна, потому что в качестве начальных данных в определении
\eqref{DEF:Fibonacci} используются два числа, $x_1$ и $x_2$ (а не одно, $x_1$).
Для обоснования такого индуктивного приема нужно заменить теорему
\ref{defin-induction} на формально другое утверждение, в котором отображение
$G$ будет действовать на последовательности элементов $a_1,...,a_n\in A$ {\it
длиной не меньше $2$}. Поскольку пример Фибоначчи используется в нашем учебнике
исключительно как иллюстрация, мы предоставляем читателю самостоятельно
продумать эти детали.
 \eex

\bex{\bf Последовательности, заданные с помощью аксиомы
выбора.}\label{EX:posl-zadan-axiomoi-vybora} Последний способ --
воспользоваться аксиомой счетного выбора (мы ее формулировали на с.\pageref{AX:schentnogo-vybora}). В
соответствии с ней, для всякой последовательности непустых множеств $\{X_n;\ n\in\N\}$
можно подобрать последовательность $\{x_n;\ n\in\N\}\subseteq \bigcup_{n\in\N}X_n$ такую, что для всякого $n\in\N$ объект $x_n$ является элементом множества $X_n$:
$$
\forall n\in\N \quad x_n\in X_n.
$$
Этот способ может показаться слишком абстрактным, но, как правило, именно он
(за редкими исключениями) используется при доказательстве общих утверждений
там, где бывает нужно построить последовательность с заданными свойствами.
Ниже на примере доказательств свойств подпоследовательностей $2^0$ и $3^0$ на
с.\pageref{podposledovatelnosti} и при доказательстве теоремы
Больцано-Вейерштрасса \ref{Bol-Wei} мы подробно обсуждаем этот эффект в связи с замечаниями о приеме бесконечнократного выбора, сделанными нами на с.\pageref{beskonechnokratnyi-vybor}.
 \eex

\end{multicols}\noindent\rule[10pt]{160mm}{0.1pt}

\paragraph{Способы изображения числовой последовательности.}

Множество натуральных чисел $\N$, по которому бегает аргумент $n$
последовательности $\{x_n\}$, обладает довольно редко встречающимся в
математике свойством {\it перечислимости}. Это значит, что можно указать
алгоритм, который сначала выдаст первый элемент этого множества, $n=1$, а затем
для каждого $n\in \N$ сгенерирует следующий за ним элемент $n+1$, и таким
образом будут перечислены все элементы $\N$. Из этого следует, что
последовательность $\{x_n\}$ можно представлять себе в виде таблицы, в которой
первую строчку (множество значений аргумента) занимают натуральные числа, а
вторая строчка будет занята соответствующими вещественными числами (значениями
последовательности):

\medskip
\vbox{
 \tabskip=0pt\offinterlineskip \halign to \hsize{
 \vrule#\tabskip=2pt plus3pt minus1pt
 &\strut\hfil\;#\hfil & \vrule#&\hfil#\hfil & \vrule#&\hfil#\hfil &
 \vrule#&\hfil#\hfil & \vrule#&\hfil#\hfil & \vrule#&\hfil#\hfil
 & \vrule#&\hfil#\hfil & \vrule#&\hfil#\hfil & \vrule#\tabskip=0pt\cr
 \noalign{\hrule} height2pt
 &\omit&&\omit&&\omit&&\omit&&\omit&&\omit&&\omit&&\omit&\cr
 & $1$ & & $2$ & & $3$ &  & $4$ & & $5$ & & $6$ & & $7$ & & ... &\cr
 \noalign{\hrule} height2pt
 &\omit&&\omit&&\omit&&\omit&&\omit&&\omit&&\omit&&\omit&\cr
 & $x_1$ & & $x_2$ && $x_3$ && $x_4$ && $x_5$ && $x_6$ && $x_7$ && ... &\cr
 \noalign{\hrule}
 }
}\medskip

\noindent Поскольку ячеек у такой таблицы должно быть бесконечно много, всю ее
нарисовать, конечно, невозможно, однако в качестве зрительного образа, который
можно держать в голове, когда речь заходит о последовательностях, такое
представление бывает полезно.

Но оно фокусирует внимание наблюдателя на первых
элементах последовательности и бывает удобно только там, где важно ее
``начало'', а на остающийся ``хвост'' можно не обращать внимание. Нас же,
наоборот, будет главным образом интересовать этот ``хвост'': первые значения
как раз не будут нам особенно нужны (причем, каким бы длинным ни был этот
начальный отрезок последовательности), а важной для нас будет лишь
закономерность, позволяющая понять, к чему стремится данная последовательность
(и стремится ли она к чему-то).

Пытаться уловить эту закономерность (разумеется, не на произвольных примерах, а
на специально подобранных, обучающих) гораздо удобнее, пользуясь другим
зрительным образом, в котором последовательность изображается системой
направленных дуг, символизирующих переходы от $x_n$ к $x_{n+1}$:

%\picture{0pt}{0pt}{55.pcx}
\vglue100pt \noindent Это можно представлять себе как движение ``блохи,
прыгающей по прямой'': сначала эта ``блоха находилась в точке $x_1$'', затем
``прыгнула в точку $x_2$'', и так далее. Понятно, что и здесь всю систему дуг
нарисовать никогда не бывает возможно (потому что их тоже бесконечно много),
однако как зрительный образ такое представление оказывается удивительно
полезным.

\noindent\rule{160mm}{0.1pt}\begin{multicols}{2}

\bex $x_n=\frac{1}{n}$

%\picture{0pt}{0pt}{56.pcx}
\vglue100pt

\eex

\bex $x_n=n$

%\picture{0pt}{0pt}{57.pcx}
\vglue100pt

\eex

\bex $x_n=\frac{(-1)^n}{n}$

%\picture{0pt}{0pt}{58.pcx}
\vglue110pt \eex

\bex $x_n=(-1)^n\cdot n$

%\picture{0pt}{0pt}{59.pcx}
\vglue110pt \eex

\bex $x_n=(-1)^n$

%\picture{0pt}{0pt}{60.pcx}
\vglue100pt \eex

\bex $x_n=3$  (постоянная последовательность)

%\picture{0pt}{0pt}{61.pcx}
\vglue100pt

\eex

\end{multicols}\noindent\rule[10pt]{160mm}{0.1pt}

\subsection{``Почти все $n\in\N$''}

Пусть $P$ обозначает какое-нибудь переменное высказывание о натуральных числах
(например, ``$n>10$'', или ``$n$ четное'', и т.д.).

 \bit{
 \item[$\bullet$]
Говорят, что {\it почти все числа $n\in \mathbb{N}$ обладают свойством
$P$}\index{``почти все''}, если только конечное (возможно, пустое) множество
чисел $n\in \mathbb{N}$ не обладает этим свойством:
 \beq\label{pochti-vse-1}
\card\{n\in\N:\quad \neg P\}<\infty
 \eeq
Это эквивалентно тому, что свойство $P$ выполняется для всех $n\in\N$, начиная
с некоторого номера:
 \beq\label{pochti-vse-2}
\exists N\in\N\qquad \forall n\ge N\qquad P
 \eeq
 }\eit

\bpr[Доказательство эквивалентности] Пусть $P$ -- какое-то переменное
высказывание о натуральных числах. Обозначим буквой $A$ множество тех чисел
$n\in\N$, для которых $P$ ложно, а через $B$ множество тех чисел $n\in\N$, для
которых $P$ истинно:
$$
A=\{n\in\N: \quad \neg P\},\qquad B=\{n\in\N: \quad P\}
$$
Тогда
$$
\eqref{pochti-vse-1}
$$
$$
\Updownarrow
$$
$$
\card A<\infty
$$
$$
\phantom{\scriptsize\text{теорема
\ref{TH:ogran-mnozh-v-N-konechno}}}\quad\Updownarrow\quad{\scriptsize\text{теорема
\ref{TH:ogran-mnozh-v-N-konechno}}}
$$
$$
\exists N\in\N\qquad A\subseteq\{1,...,N\}
$$
$$
\Updownarrow
$$
$$
\exists N\in\N\qquad \{n\in\N: \quad n>N\}\subseteq B
$$
$$
\Updownarrow
$$
$$
\eqref{pochti-vse-2}
$$

 \epr

\noindent\rule{160mm}{0.1pt}\begin{multicols}{2}

\begin{ex}
Ясно, что почти все числа $n\in \mathbb{N}$ удовлетворяют
неравенству
$$
n \ge 5
$$
потому что чисел $n\in \mathbb{N}$, не удовлетворяющих этому
неравенству имеется только конечное множество (а именно $n\in \{1,
2, 3, 4\}$).
\end{ex}

\begin{ex}
Наоборот, неверно, что почти все числа $n\in \mathbb{N}$
удовлетворяют неравенству
$$
n < 5
$$
потому что чисел $n\in \mathbb{N}$, не удовлетворяющих этому
неравенству бесконечно много (а именно, $n= 5,6,7,8,...$).
\end{ex}

\begin{er}
Сообразим, будет ли неравенство
$$
(-1)^n >0
$$
выполняется для почти всех $n\in \mathbb{N}$?

Понятно, что это не так, потому что оно выполняется только для
четных чисел ($n= 2,4,6,8,...$), а для нечетных чисел ($n=
1,3,5,7,...$) -- которых бесконечно много -- оно неверно.
\end{er}

\begin{er}
Верно ли, что обратное неравенство
$$
(-1)^n \le 0
$$
выполняется для почти всех $n\in \mathbb{N}$?

Это, конечно, тоже не так, потому что этому неравенству
удовлетворяют только нечетные числа ($n= 1,3,5,7,...$), а для
четных чисел ($n= 2,4,6,8,...$) -- которых бесконечно много -- оно
неверно.
\end{er}
 \end{multicols}\noindent\rule[10pt]{160mm}{0.1pt}

\begin{tm}[\bf Архимеда]\label{Archimed-theo}
\index{теорема!Архимеда} Каково бы ни было число  $C\in \R$, почти все числа
$n\in \mathbb{N}$ удовлетворяют неравенству
\begin{equation}
n>C \label{1.1.1}
\end{equation}
\end{tm}
\begin{proof} Это следует из принципа Архимеда
(теорема \ref{Archimed-principle}). Выберем какое-нибудь $N\in
\mathbb{N}$ так, чтобы $N>C$. Тогда все числа $n$, начиная с числа
$N$, удовлетворяют неравенству \eqref{1.1.1}, потому что $n\ge N>C$.
\end{proof}

\noindent\rule{160mm}{0.1pt}\begin{multicols}{2}

\ber
 Следующая ``абстрактная'' задача может считаться тестом
 на понимание термина ``почти все''.
 Пусть $M$ какое-нибудь множество натуральных чисел ($M\subseteq \N$),
 рассмотрим пять возможных ситуаций:

 \biter{
\item[(A)] все числа $n\in\N$ лежат в множестве $M$ (то есть $M=\N$);
\item[(B)] конечный набор чисел $n\in\N$ не лежит в множестве $M$,
а остальные числа $n\in\N$ (которых -- бесконечное множество)
лежат в $M$;
\item[(C)] бесконечный набор чисел $n\in\N$ не лежит в множестве $M$,
и одновременно бесконечный набор чисел $n\in\N$ лежит в $M$ (такое
бывает, например, когда $M$ -- множество четных чисел);
\item[(D)] конечный набор чисел $n\in\N$ лежит в множестве $M$,
а остальные числа $n\in\N$ (которых -- бесконечное множество) не
лежат в $M$;
\item[(E)] все числа $n\in\N$ не лежат в множестве $M$  (то есть $M=\varnothing$).
 }\eiter

Ответьте на вопросы:
 \biter{
 \item[1)] в каких из этих случаев почти все числа $n\in\N$ лежат в множестве $M$?
 \item[2)] в каких случаях почти все числа $n\in\N$ не лежат в множестве $M$?
 \item[3)] в каких случаях не верно ни то, ни другое?
 }\eiter

Ответы:
 \biter{
 \item[1)] (A) и (B) (чисел, лежащих в $M$ больше, чем остальных);
 \item[2)] (D) и (E) (чисел, не лежащих в $M$ больше, чем остальных);
 \item[3)] (C).
 }\eiter
 \eer

\end{multicols}\noindent\rule[10pt]{160mm}{0.1pt}

\bit{ \item[$\bullet$] Говорят, что {\it почти все элементы последовательности
$\{ x_n \}$ принадлежат множеству $E$}, если для почти всех номеров $n\in
\mathbb{N}$ соответствующие числа $\{ x_n \}$ принадлежат множеству $E$.
 }\eit

\noindent\rule{160mm}{0.1pt}\begin{multicols}{2}

\begin{ex}
Проверим, что почти все элементы последовательности $\{ x_n =\frac{1}{n} \}$
лежат в интервале $\l 0; \frac{1}{4}\r$.

%\picture{0pt}{0pt}{62.pcx}
\vglue100pt \noindent
 Действительно,

\begin{gather*}
x_n \in \l 0; \frac{1}{4}\r \quad \Leftrightarrow \quad \frac{1}{n}
\in \l 0; \frac{1}{4} \r \quad \Leftrightarrow \\
\Leftrightarrow \quad
\begin{cases}
 \phantom{.}\overset{\phantom{.}}{\frac{1}{n}<\frac{1}{4}}
 \\
 {\tiny \boxed{\frac{1}{n}>0}}
\end{cases}
 \put(-20,-34){\vector(-1,3){5}}\put(-85,-47){
 \boxed{\tiny \begin{matrix}\text{это неравенство
выполняется автоматически для любого $n\in \mathbb{N}$,} \\
\text{поэтому его можно отбросить}\end{matrix}}} \quad
\Leftrightarrow \quad \frac{1}{n} < \frac{1}{4} \quad \Leftrightarrow
\quad n>4
 \end{gather*}
а последнее неравенство в этой цепочке выпол\-няет\-ся для почти всех
$n\in \mathbb{N}$, по теореме Архимеда \ref{Archimed-theo}.
\end{ex}

\begin{ex}
Верно ли, что почти все элементы последовательности $\{ x_n =\frac{1}{n} \}$
лежат в интервале $\l \frac{1}{4};\infty\r$?

%\picture{0pt}{0pt}{63.pcx}
\vglue100pt \noindent Нет, конечно, потому что
$$
x_n \in \l \frac{1}{4};\infty\r \quad \Leftrightarrow \quad
\frac{1}{n} > \frac{1}{4} \quad \Leftrightarrow \quad n<4
$$
а последнее неравенство выпол\-няется {\bf не для почти всех} чисел
$n\in \mathbb{N}$, а только для $n\in\{ 1,2,3,$ $...,9 \}$.
\end{ex}

\begin{ex}
Верно ли, что почти все элементы последовательности $\{ x_n =(-1)^n \}$ лежат в
интервале $(0;2)$?

%\picture{0pt}{0pt}{64.pcx}
\vglue100pt \noindent Чтобы это понять, разберемся сначала для
каких номеров $n\in \mathbb{N}$ числа $x_n$ лежат в интервале
$(0;2)$:
 \begin{gather*}
x_n \in (0;2) \quad \Leftrightarrow \quad (-1)^n \in (0;2) \quad
\Leftrightarrow \\ \Leftrightarrow  \quad
\begin{cases}
\phantom{.}\underset{\phantom{\cdot}}{(-1)^n>0}
\\
\boxed{(-1)^n<2}
\end{cases}
 \put(-20,-34){\vector(-1,3){5}}\put(-110,-47){
 \boxed{\tiny \begin{matrix}
 \text{это неравенство выполняется автоматически для любого
 $n\in \mathbb{N}$,} \\ \text{поэтому его можно отбросить}\end{matrix}}}
\quad \Leftrightarrow \quad (-1)^n>0
 \end{gather*}
Последнее неравенство в этой цепочке выполняется {\bf не для почти
всех} чисел $n\in \mathbb{N}$ а только для четных: $n\in \{
2,4,6,... \}$. Значит, соотношение $x_n \in (0;2)$ тоже
выполняется {\bf не для почти всех} чисел $n\in \mathbb{N}$.
\end{ex}

\begin{ex}
 \biter{
\item[(a)] Найдите какие-нибудь два числа $\alpha$ и $\beta$ такие, чтобы почти
все элементы последовательности $x_n= \frac{n-1}{n}$ лежали в интервале
$(\alpha ; \beta)$.

\item[(b)] Найдите все такие числа $\alpha$ и $\beta$, чтобы почти все элементы
последовательности $x_n= \frac{n-1}{n}$ лежали в интервале $(\alpha ; \beta)$.
 }\eiter

Нарисуем на прямой первые несколько элементов нашей последовательности

%\picture{0pt}{0pt}{65.pcx}
\vglue100pt \noindent

Из рисунка видно, что, например, вне интервала $(\frac{1}{2} ; 2)$
лежит только конечный набор чисел $x_n$. Отсюда следует

Ответ для пункта (a): можно взять, например, $\alpha=\frac{1}{2}$ и $\beta=2$.

Более того, если взять какие-нибудь другие числа  $\alpha$ и
$\beta$ так чтобы $\alpha < 1 \le \beta$, то начиная с какого-то
номера все $x_n$ будут лежать в интервале $(\alpha ; \beta)$.

%\picture{0pt}{0pt}{66.pcx}
\vglue100pt \noindent

При этом ясно, что если $\alpha \ge 1$ или $\beta < 1$, то
утверждение, что почти все числа $x_n$ лежат в интервале $(\alpha
; \beta)$ будет неверно. Поэтому ответ для второй половины задачи
выглядит так

Ответ для пункта (b): $\alpha < 1 \le \beta$.
\end{ex}

\begin{ex}
 \biter{
\item[(a)] Найдите какие-нибудь два числа $\alpha$ и $\beta$ такие, чтобы почти
все элементы последовательности $x_n= (-1)^n$ лежали в интервале $(\alpha ;
\beta)$. \item[(b)] Найдите все такие числа $\alpha$ и $\beta$, чтобы почти все
элементы последовательности $x_n= (-1)^n$ лежали в интервале $(\alpha ;
\beta)$.
 }\eiter

Снова нарисуем на прямой первые несколько элементов нашей последовательности

%\picture{0pt}{0pt}{67.pcx}
\vglue100pt \noindent

Из рисунка видно, что, например, вне интервала $(-2 ; 2)$ нет
вообще никаких чисел $x_n$. Отсюда следует

Ответ для пункта (a): можно взять, например, $\alpha=-2$ и $\beta=2$.

И вообще, если взять какие-нибудь другие числа  $\alpha$ и $\beta$
так чтобы $\alpha < -1$ и $\beta > 1$, то все $x_n$ будут лежать в
интервале $(\alpha ; \beta)$.

%\picture{0pt}{0pt}{68.pcx}
\vglue100pt \noindent

Если же $\alpha \ge 1$ или $\beta \le 1$, то утверждение, что
почти все числа $x_n$ лежат в интервале $(\alpha ; \beta)$ будет
неверно. Поэтому ответ для второй половины задачи выглядит так

Ответ для пункта (b): $\alpha < -1$ и $\beta > 1$.
\end{ex}

\begin{ers}
 Ничего не доказывая, но пользуясь интуицией, сообразите,
 какими должны быть величины $\alpha, \beta$, чтобы почти все
 элементы последовательности $x_n$ лежали в интервале  $(\alpha;\beta)$?
 \biter{
\item[1)] $x_n=(-1)^n\cdot n$ \\
{\smsize (Ответ:$\alpha=-\infty,\,\beta=\infty$)}

\item[2)] $x_n=\frac{(-1)^n}{n}+\frac{1+(-1)^n}{2}$ \\
{\smsize (Ответ: $\alpha<0,\, \beta>1$)}

\item[3)] $x_n=2+\frac{1+(-1)^n}{n}$ \\
{\smsize (Ответ: $\alpha<2,\,\beta>2$)}

\item[4)] $x_n=2-\frac{1+(-1)^n}{n}$ \\
{\smsize (Ответ: $\alpha<2,\, \beta>2$)}

\item[5)] $x_n=n^{(-1)^n}$ \\
{\smsize (Ответ: $\alpha<0,\, \beta=\infty$)}

\item[6)] $x_n=n\cdot (1+(-1)^n)$ \\
{\smsize (Ответ: $\alpha<0,\, \beta=\infty$)}

\item[7)] $x_n=\frac{2+(-1)^n}{2}-\frac{1}{n}$ \\
{\smsize (Ответ: $\alpha<\frac{1}{2}, \, \beta>\frac{3}{2}$)}

\item[8)] $x_n=(-1)^n\cdot \l 1-\frac{1}{n} \r$ \\
{\smsize (Ответ: $\alpha\le -1, \, \beta\ge 1$)}

\item[9)] $x_n=(-1)^n\cdot \l 1+\frac{1}{n} \r$ \\
{\smsize (Ответ: $\alpha< -1, \, \beta> 1$)}

\item[10)] $x_n=(-1)^{n-1}\cdot \l 2+\frac{3}{n} \r$ \\
{\smsize (Ответ: $\alpha<-2, \, \beta>2$)}

 }\eiter
\end{ers}
\end{multicols}\noindent\rule[10pt]{160mm}{0.1pt}

\subsection{Предел числовой последовательности}\label{SEC:predel-posledov}

\paragraph{Окрестность точки.}

\bit{

\item[$\bullet$] {\it Окрестностью} точки $c\in \R$ называется всякий интервал
с центром в точке $c$, то есть всякий интервал вида $(c-\varepsilon ; c+
\varepsilon)$. Число $\varepsilon$ при этом называется {\it радиусом}
окрестности $(c-\varepsilon ; c+ \varepsilon)$.

}\eit

\bprop\label{PROP:x-in-(a,b)=>U(x)-sub-(a,b)} Всякая точка $x$ произвольного интервала $(a,b)$ содержится в нем вместе с некоторой своей окрестностью:
$$
x\in(a,b)\qquad\Longrightarrow\qquad \exists\e>0: \quad (x-\varepsilon ; x+ \varepsilon)\subseteq (a,b)
$$
\eprop
\bpr Положим
$$
\e=\min\{b-x;x-a\}.
$$
Тогда
 \begin{multline*}
\begin{cases}
\e\le b-x \\
\e\le x-a
\end{cases}
\quad\Longrightarrow\quad
\begin{cases}
\e\le b-x  \\
-\e\ge a-x
\end{cases}
\quad\Longrightarrow\quad
\begin{cases}
x+\e\le b  \\
x-\e\ge a
\end{cases}
\quad\Longrightarrow \\ \Longrightarrow\quad a\le x-\e<x<x+\e\le b
\quad\Longrightarrow\quad x\in(x-\e,x+\e)\subseteq (a,b)
 \end{multline*}

\epr

\paragraph{Конечный предел последовательности.}

\bit{

\item[$\bullet$] Точка $c\in \R$ называется {\it пределом последовательности}
\index{предел!последовательности!конечный} $\{ x_n \}$, если любая ее
окрестность $(c-\varepsilon ; c+ \varepsilon)$ содержит почти все элементы
последовательности $\{ x_n \}$. В этом случае пишут
$$
\lim_{n\to \infty} x_n = c
$$
или
$$
x_n \underset{n\to \infty}{\longrightarrow} c
$$
и говорят, что {\it последовательность $\{ x_n \}$ стремится к числу $c$}.

Если последовательность $\{ x_n \}$ имеет предел (то есть стремится к какому-то
числу $c$), то говорят, что {\it последовательность $\{ x_n \}$ сходится}
\index{последовательность!сходящаяся}. Если же такого числа $c$ не существует,
то говорят, что {\it последовательность $\{ x_n \}$ расходится}
\index{последовательность!расходящаяся}.
 }\eit

\noindent\rule{160mm}{0.1pt}\begin{multicols}{2}

\begin{ex}
Покажем, что
$$
\lim_{n\to \infty} \frac{n-1}{n} = 1
$$
то есть что точка $c=1$ является пределом последовательности $\{
x_n=\frac{n-1}{n} \}$.

Действительно, возьмем какую-нибудь окрестность $(1-\varepsilon ; 1+
\varepsilon)$ точки $c=1$,

%\picture{0pt}{0pt}{69.pcx}
\vglue100pt \noindent и поймем, что означает, что $x_n$ лежит в
$(1-\varepsilon ; 1+ \varepsilon)$:

\begin{gather*}
x_n \in (1-\varepsilon ; 1+ \varepsilon) \quad \Leftrightarrow \quad
1-\varepsilon < x_n <1+\varepsilon \quad \Leftrightarrow \\
\Leftrightarrow \quad
\begin{cases} { x_n>1-\varepsilon}\\{x_n <1+\varepsilon}\end{cases}
\quad \Leftrightarrow \quad
\begin{cases} {\frac{n-1}{n}>1-\varepsilon}\\{\frac{n-1}{n}
<1+\varepsilon}\end{cases}
 \quad \Leftrightarrow \\ \Leftrightarrow \quad
\begin{cases} {\frac{n-1}{n}-1>-\varepsilon}\\{\frac{n-1}{n}-
1<\varepsilon}\end{cases}
\quad \Leftrightarrow \quad
\begin{cases} { -\frac{1}{n}>-\varepsilon }\\{-
\frac{1}{n}<\varepsilon}\end{cases}
 \quad \Leftrightarrow \\ \Leftrightarrow \quad
\begin{cases}
\phantom{A}\underset{\phantom{\cdot}}{\frac{1}{n}<\varepsilon}
\\
{\tiny \boxed{-\frac{1}{n}<\varepsilon}}
\end{cases}\put(-20,-35){\vector(-1,3){5}}
\put(-85,-50){\boxed{\tiny \begin{matrix}\text{поскольку $\varepsilon>0$,}\\
\text{это неравенство выполняется автоматически для любого $n\in
\mathbb{N}$,}\\
\text{значит его можно отбросить}\end{matrix}}} \quad \Leftrightarrow
\quad \frac{1}{n} < \varepsilon \quad \Leftrightarrow \quad
n>\frac{1}{\varepsilon}
\end{gather*}

По теореме Архимеда \ref{Archimed-theo}, последнее неравенство
выполняется для почти всех $n\in \mathbb{N}$. Это означает, что
почти все числа $x_n$ содержатся в интервале $(1-\varepsilon ;
1+\varepsilon)$.

Итак, мы получили, что любая окрестность $(1-\varepsilon ; 1+\varepsilon)$
точки $c$ содержит почти все элементы последовательности $\{ x_n=\frac{n-1}{n}
\}$. Это как раз означает, что точка $c=1$ является пределом последовательности
$\{ x_n=\frac{n-1}{n} \}$.
\end{ex}

\begin{ex}
Покажем, что
$$
1\ne \lim_{n\to \infty} (-1)^n
$$
то есть что число $c=1$ не является пределом последовательности $\{ x_n=(-1)^n
\}$.

Возьмем $\varepsilon =1$ и рассмотрим соответствующую окрестность
$(0;2)$ точки $c=1$. Тогда
\begin{gather*}
x_n \in (0;2) \quad \Leftrightarrow \\ \Leftrightarrow  \quad
\begin{cases}
\phantom{A}\underset{\phantom{\cdot}}{0< (-1)^n}
\\
\boxed{(-1)^n<2}
\end{cases}\put(-20,-33){\vector(-1,3){5}}\put(-110,-45)
{\boxed{\tiny\begin{matrix} \text{это неравенство
выполняется автоматически для любого $n\in \mathbb{N}$,} \\
\text{поэтому его можно отбросить}\end{matrix}}} \quad
\Leftrightarrow \quad 0< (-1)^n
\end{gather*}
Последнее неравенство выполняется только для четных $n\in
\mathbb{N}$.

Мы получили, что бесконечное множество элементов последовательности $\{
x_n=(-1)^n \}$ не принадлежит окрестности $(0;2)$ точки $c=1$. Таким образом,
будет неверно утверждать, что почти все $\{ x_n=(-1)^n \}$ содержатся в
окрестности $(0;2)$ точки $c=1$. Это означает, что $c=1$ не может быть пределом
последовательности $\{ x_n=(-1)^n \}$.
\end{ex}

\begin{ex}
Покажем, что
$$
\lim_{n\to \infty} \frac{2n+3}{n+2} = 2
$$
то есть что точка $c=2$ является пределом последовательности $\{
x_n=\frac{2n+3}{n+2} \}$.

Возьмем какую-нибудь окрестность $(2-\varepsilon ; 2+ \varepsilon)$
точки $c=2$

%\picture{0pt}{0pt}{70.pcx}
\vglue100pt \noindent Тогда

\begin{gather*}
x_n \in (2-\varepsilon ; 2+\varepsilon ) \quad \Leftrightarrow \quad
2-\varepsilon < x_n <2+\varepsilon
 \quad \Leftrightarrow \\ \Leftrightarrow\quad
\begin{cases} {x_n>2-\varepsilon}\\{x_n<2+\varepsilon}\end{cases}
\quad \Leftrightarrow \quad
\begin{cases} {\frac{2n+3}{n+2}>2-
\varepsilon}\\{\frac{2n+3}{n+2}<2+\varepsilon}\end{cases}
 \quad \Leftrightarrow \\ \Leftrightarrow \quad
\begin{cases} {\frac{2n+3}{n+2}-2>-\varepsilon}\\{\frac{2n+3}{n+2}-
2<\varepsilon}\end{cases}
\quad \Leftrightarrow \quad
\begin{cases} {-\frac{1}{n+2}>-\varepsilon}\\{-\frac{1}{n+2}<\varepsilon
}\end{cases}
 \quad \Leftrightarrow \\ \Leftrightarrow  \quad
\begin{cases}
\phantom{...}\underset{\phantom{.}}{\frac{1}{n+2}<\varepsilon}
\\
{\tiny \boxed{-\frac{1}{n+2}<\varepsilon}}
\end{cases}\put(-20,-35){\vector(-1,3){5}}
\put(-100,-52){
\boxed{\tiny \begin{matrix}\text{поскольку $\varepsilon>0$,} \\
\text{это неравенство выполняется автоматически для любого $n\in
\mathbb{N}$,} \\
\text{значит его можно отбросить}\end{matrix}}}
\quad \Leftrightarrow
\quad \frac{1}{n+2} <\varepsilon
\quad \Leftrightarrow \\
\Leftrightarrow \quad n+2 >\frac{1}{\varepsilon} \quad
\Leftrightarrow \quad n >\frac{1}{\varepsilon}-2
\end{gather*}
По теореме Архимеда \ref{Archimed-theo}, последнее неравенство
выполняется для почти всех $n\in \mathbb{N}$. Это означает, что
почти все числа $x_n$ содержатся в интервале $(2-\varepsilon ; 2+
\varepsilon)$.

Итак, мы получили, что любая окрестность $(2-\varepsilon ; 2+ \varepsilon)$
точки $c=2$ содержит почти все элементы последовательности $\{
x_n=\frac{2n+3}{n+2} \}$. Это означает, что точка $c=2$ является пределом
последовательности $\{ x_n=\frac{2n+3}{n+2} \}$.
\end{ex}

\begin{ex}
Покажем, что
$$
\lim_{n\to \infty} \frac{(-1)^n}{n} = 0
$$
то есть что точка $c=0$ является пределом последовательности $\{
x_n=\frac{(-1)^n}{n} \}$.

Возьмем какую-нибудь окрестность $(-\varepsilon;\varepsilon)$ точки
$c=0$

%\picture{0pt}{0pt}{71.pcx}
\vglue140pt \noindent и покажем, что почти все числа
$x_n=\frac{(-1)^n}{n}$ содержатся в интервале
$(-\varepsilon;\varepsilon)$. Действительно,

\begin{multline*}
x_n \in (-\varepsilon;\varepsilon) \quad \Leftrightarrow \quad
(\text{свойство} \, 2^0 \, \text{пункта 0.6}) \quad \Leftrightarrow
\\ \Leftrightarrow \quad |x_n|< \varepsilon \quad \Leftrightarrow \quad
\left|\frac{(-1)^n}{n}\right|< \varepsilon \quad \Leftrightarrow \\
\Leftrightarrow \quad (\text{свойства} \, 5^0 ,\, 6^0 \, \text{пункта
0.6}) \quad \Leftrightarrow \\ \Leftrightarrow \quad
\frac{|-1|^n}{n}< \varepsilon \quad \Leftrightarrow \quad
\frac{1}{n}< \varepsilon \quad \Leftrightarrow \quad
n>\frac{1}{\varepsilon} \end{multline*}

По теореме Архимеда \ref{Archimed-theo}, последнее неравенство выполняется для
почти всех $n\in \mathbb{N}$. Это означает, что почти все числа $x_n$
содержатся в окрестности $(-\varepsilon;\varepsilon)$ точки $c=0$. Таким
образом, точка $c=0$ действительно является пределом последовательности $\{
x_n=\frac{2n+3}{n+2} \}$.
\end{ex}

\begin{ex}
Покажем, что
$$
0\ne \lim_{n\to \infty} n^{(-1)^n}
$$
то есть что точка $c=0$ не является пределом последовательности $\{
x_n=n^{(-1)^n} \}$.

Возьмем интервал $(-1;1)$, являющийся окрестностью точки $c=0$ (при
$\varepsilon=1$).

%\picture{0pt}{0pt}{72.pcx}
\vglue100pt \noindent и посмотрим, что будет означать условие $x_n
\in (-1 ; 1)$.

Для этого нужно будет рассмотреть отдельно четные и нечетные
номера $n$. Тогда условие $x_n \in (-1 ; 1)$ распадется в
совокупность двух систем:
 \begin{gather*}
x_n \in (-1 ; 1) \quad \Leftrightarrow \quad
\begin{cases} {-1 < x_n}\\{x_n <1}\end{cases} \quad \Leftrightarrow \\
 \put(5,50){\boxed{\tiny \begin{matrix}\text{поскольку всегда
$n^{(-1)^n}>0$,} \\ \text{это неравенство выполняется автоматически
для любого $n\in \mathbb{N}$,} \\ \text{значит его можно
отбросить}\end{matrix}}}\put(75,38){\vector(-1,-3){5}}
 \Leftrightarrow \quad
\begin{cases}
\boxed{-1< n^{(-1)^n} }
\\
\phantom{..}\overset{\phantom{.}}{n^{(-1)^n} <1}
\end{cases}
\quad \Leftrightarrow \quad n^{(-1)^n} <1
 \quad \Leftrightarrow \\ \Leftrightarrow \quad
 \left[\begin{matrix} \left\{
\begin{matrix}
\text{$n$ -- четное: $n=2k, \, k\in \N$}
\\
\text{тогда $(-1)^n=(-1)^{2k}=1$}
\\
\text{и при подстановке получается:}
\\
{(2k)^{1} <1}
\end{matrix}
\right\}
\\
\left\{
\begin{matrix}
\text{$n$ -- нечетное: $n=2k-1, \, k\in \N$}
\\
\text{тогда $(-1)^n=(-1)^{2k}=-1$}
\\
\text{и при подстановке получается:}
\\
{(2k-1)^{-1} <1}
\end{matrix}
\right\}
\end{matrix} \right]
 \quad \Leftrightarrow \\ \Leftrightarrow \quad
 \left[\begin{matrix} \left\{
\begin{matrix}
n=2k
\\
k\in \N
\\
2k<1
\end{matrix}
\right\}
\\
\left\{
\begin{matrix}
n=2k-1
\\
k\in \N
\\
\frac{1}{2k-1} <1
\end{matrix}
\right\}
\end{matrix} \right]
\quad \Leftrightarrow \\
\Leftrightarrow \quad
\put(-10,93){\phantom{.}}\put(-10,73){\boxed{\tiny
\begin{matrix}\text{эта система не имеет решений,} \\ \text{поэтому
ее можно выбросить из совокупности}\end{matrix}}}
 \put(60,65){\vector(-1,-3){5}}
 \left[\begin{matrix} \boxed{ \left\{
\begin{matrix} n=2k
\\
k\in \N
\\
k<\frac{1}{2}
\end{matrix}
\right\} } \\
\\
\left\{
\begin{matrix}
n=2k-1
\\
k\in \N
\\
2k-1>1
\end{matrix}
\right\}
\end{matrix} \right]
\quad \Leftrightarrow \\ \Leftrightarrow \quad  \left\{
\begin{matrix}
n=2k-1
\\
k\in \N
\\
k>1
\end{matrix}
\right\}  \quad \Leftrightarrow \quad \left\{
\begin{matrix}
n=2k-1
\\
k\in \N
\\
k>1
\end{matrix}
\right\}
\end{gather*}
Мы видим, что условие $x_n \in (-1 ; 1)$ выполняется только для индексов
$n=3,5,7...$. Таким образом, бесконечное множество чисел $x_n$ (с четными
индексами $n$) не содержатся в окрестности $(-1 ; 1)$ точки $c=0$. Значит эта
точка не может быть пределом последовательности
 $\{ x_n=n^{(-1)^n} \}$.
\end{ex}

\end{multicols}\noindent\rule[10pt]{160mm}{0.1pt}

\begin{tm}[\bf о единственности предела]
Если
$$
x_n \underset{n\to \infty}{\longrightarrow} a
$$
и
$$
x_n \underset{n\to \infty}{\longrightarrow} b
$$
то
$$
a=b
$$
\end{tm}
\begin{proof} Предположим, что $a\ne b$, например
$a<b$. Возьмем $\varepsilon=\frac{b-a}{2}$, тогда мы получим, что почти все $x_n$ лежат в двух непересекающихся интервалах -- в $(a-\e;a+\e)$ и в $(b-\e;b+\e)$:
$$
\underbrace{a-\varepsilon < \kern-7pt\overset{\scriptsize\begin{matrix}x_n\\ \phantom{\tiny\begin{matrix}n\\ \downarrow\\ \infty \end{matrix}}\ \downarrow \ {\tiny\begin{matrix}n\\ \downarrow\\ \infty \end{matrix}}\end{matrix}}{a}\kern-7pt < a+\varepsilon}_{\scriptsize\begin{matrix}\Downarrow\\ x_n\in(a-\e,a+\e)\\ \text{для почти всех $n\in\N$}\end{matrix}}=\underbrace{b-\varepsilon < \kern-7pt\overset{\scriptsize\begin{matrix}x_n\\ \phantom{\tiny\begin{matrix}n\\ \downarrow\\ \infty \end{matrix}}\ \downarrow \ {\tiny\begin{matrix}n\\ \downarrow\\ \infty \end{matrix}}\end{matrix}}{b}\kern-7pt <
b+\varepsilon}_{\scriptsize\begin{matrix}\Downarrow\\ x_n\in(b-\e,b+\e)\\ \text{для почти всех $n\in\N$}\end{matrix}}
$$
Но это невозможно: если почти все $x_n$ лежат в $(a-\e;a+\e)$, то вне $(a-\e;a+\e)$ должен лежать лишь конечный набор чисел $x_n$. В частности, в  $(b-\e;b+\e)$ тоже  должен лежать лишь конечный набор чисел $x_n$.
\end{proof}

\noindent\rule{160mm}{0.1pt}\begin{multicols}{2}

\begin{ers}
Доказать что
 \biter{
\item[1)] $\lim\limits_{n\to \infty} \frac{1}{n^3}=0$
\item[2)] $\lim\limits_{n\to \infty} \frac{n-3}{2n+5}=\frac{1}{2}$
\item[3)] $\lim\limits_{n\to \infty} \frac{n^2+1}{n^2+2}=1$
\item[4)] $\lim\limits_{n\to \infty} \frac{1-n}{n+1}=-1$
\item[5)] $\lim\limits_{n\to \infty} \frac{1}{1+2^n}=0$
\item[6)] $1\ne \lim\limits_{n\to \infty} (-1)^n$
\item[7)] $0\ne \lim\limits_{n\to \infty} (-1)^n \cdot n$
 }\eiter
\end{ers}

\end{multicols}\noindent\rule[10pt]{160mm}{0.1pt}

\paragraph{Бесконечный предел последовательности.}

 \bit{
\item[$\bullet$] Говорят, что числовая последовательность $\{ x_n \}$
\index{предел!последовательности!бесконечный}
 \bit{
\item[--] {\it стремится к} $+\infty$
$$
x_n\underset{n\to \infty}{\longrightarrow} +\infty
$$
если любой интервал $(E;+\infty)$ содержит почти все элементы $\{
x_n \}$.

\item[--] {\it стремится к} $-\infty$
$$
x_n\underset{n\to \infty}{\longrightarrow} -\infty
$$
если любой интервал $(-\infty;E)$ содержит почти все элементы $\{
x_n \}$.

\item[--] {\it стремится к} $\infty$
$$
x_n\underset{n\to \infty}{\longrightarrow} \infty
$$
если всякое множество $(-\infty;-E)\cup (E;+\infty)$ содержит почти все
элементы $\{ x_n \}$; нетрудно проверить, что это равносильно тому, что
последовательность из модулей $\{ |x_n| \}$ стремится к $+\infty$:
$$
|x_n|\underset{n\to \infty}{\longrightarrow} +\infty
$$
 }\eit
 }\eit

\noindent\rule{160mm}{0.1pt}\begin{multicols}{2}

\begin{ex}
Покажем, что
$$
\lim_{n\to \infty} n^2=+\infty
$$
Возьмем какой-нибудь интервал $(E; +\infty)$

%\picture{0pt}{0pt}{73.pcx}
\vglue100pt \noindent Тогда
$$
x_n \in (E; +\infty) \quad \Leftrightarrow \quad x_n>E \quad
\Leftrightarrow \quad n^2>E
$$
Теперь надо рассмотреть два случая:

1) если $E<0$, то неравенство $n^2>E$ выполняется для всех $n\in
\mathbb{N}$;

2) если $E\ge 0$, то неравенство $n^2>E$ эквивалентно неравенству
$n>\sqrt{E}$, а оно выполняется для почти всех $n\in \mathbb{N}$
по теореме Архимеда \ref{Archimed-theo}.

Таким образом, мы получаем, что в любом случае почти все числа $x_n$ содержатся
в интервале $(E; +\infty)$. Поскольку это верно для любого интервала $(E;
+\infty)$, последовательность $x_n=n^2$ действительно стремится к $+\infty$.
\end{ex}

\begin{ex}
Покажем, что последовательность $\{ x_n=n^{(-1)^n} \}$ не стремится к
бесконечности:
$$
n^{(-1)^n}\underset{n\to \infty}{\notarrow} \infty
$$

Возьмем множество $(-\infty; -1)\cup (1;+\infty)$

%\picture{0pt}{0pt}{74.pcx}
\vglue100pt \noindent и рассмотрим нечетные номера $n$, то есть
числа вида $n=2k-1$, где $k\in \mathbb{N}$. Тогда
$(-1)^n=(-1)^{2k-1}=-1$, и
 \begin{multline*}
x_n \in (-\infty; -1)\cup (1;+\infty) \quad \Leftrightarrow \quad
|x_n|>1 \quad \Leftrightarrow \\ \Leftrightarrow \quad n^{(-1)^n} >1
 \quad \Leftrightarrow \\ \Leftrightarrow \quad
 {\smsize\begin{pmatrix}
\text{применяем равенства}\\
n=2k-1,\\
(-1)^n=(-1)^{2k-1}=-1
\end{pmatrix}}
 \quad \Leftrightarrow \\ \Leftrightarrow  \quad
 (2k-1)^{-1} >1 \quad \Leftrightarrow \quad 2k-1<1
 \quad \Leftrightarrow \\ \Leftrightarrow  \quad k<1
 \end{multline*}
Последнее неравенство не выполняется ни при каком $k=1,2,3...$ Это означает,
что все числа $x_{2k-1}$ не содержатся в множестве $(-\infty; -1)\cup
(1;+\infty)$. Таким образом, бесконечное множество чисел $x_n$ не содержатся в
``окрестности бесконечности'' $(-\infty; -1)\cup (1;+\infty)$. Значит,
бесконечность не может быть пределом последовательности
 $\{ x_n=n^{(-1)^n} \}$.
\end{ex}

\begin{ers}
Доказать что
 \biter{
\item[1)] $\lim\limits_{n\to \infty} \frac{4n+1}{3n-2}\ne+\infty$

\item[2)] $\lim\limits_{n\to \infty} \frac{n^2}{1-2n}=-\infty$
 }\eiter
\end{ers}

\end{multicols}\noindent\rule[10pt]{160mm}{0.1pt}

\subsection{Бесконечно малые и бесконечно большие
последовательности} \label{->0-&-->infty}

\bit{ \item[$\bullet$] Последовательность  $\{ x_n \}$ называется {\it
бесконечно малой} \index{последовательность!бесконечно малая}, если выполняются
следующие равносильные условия:
 \bit{
\item[(i)]  $\{ x_n \}$ стремится к нулю:
$$
x_n \underset{n\to \infty}{\longrightarrow} 0
$$
\item[(ii)]  $\{ x_n \}$ по модулю стремится к нулю:
$$
|x_n| \underset{n\to \infty}{\longrightarrow} 0
$$
 }\eit
 }\eit
\begin{proof} Здесь необходимо проверить, что
условия $(i)$ и $(ii)$ действительно равносильны. Для этого перепишем свойство
модуля $2^0$ пункта 0.6, подставив вместо $x$ последовательность $\{ x_n \}$:
$$
x_n\in (-\varepsilon;+\varepsilon) \quad \Longleftrightarrow \quad
|x_n|\in (-\varepsilon;+\varepsilon)
$$
После этого получается следующая цепочка равносильностей:
$$
x_n \underset{n\to \infty}{\longrightarrow} 0
$$
$$
\Updownarrow
$$
$$
\text{ для всякого $\varepsilon>0$ соотношение $x_n\in
(-\varepsilon;+\varepsilon)$ выполняется для почти всех $n\in
\mathbb{N}$ }
$$
$$
\Updownarrow
$$
$$
\text{ для всякого $\varepsilon>0$ соотношение $|x_n|\in
(-\varepsilon;+\varepsilon)$ выполняется для почти всех $n\in
\mathbb{N}$ }
$$
$$
\Updownarrow
$$
$$
|x_n| \underset{n\to \infty}{\longrightarrow} 0 \qquad
$$ \end{proof}

\bit{ \item[$\bullet$] Последовательность  $\{ x_n \}$ называется {\it
бесконечно большой} \index{последовательность!бесконечно большая}, если
выполняются следующие равносильные условия:
 \bit{
\item[(i)] $\{ x_n \}$ стремится к бесконечности
$$
x_n \underset{n\to \infty}{\longrightarrow} \infty
$$
\item[(ii)] $\{ x_n \}$ по модулю стремится к бесконечности
$$
|x_n| \underset{n\to \infty}{\longrightarrow} \infty
$$
 }\eit
 }\eit
\begin{proof}[Доказательство равносильности условий $(i)$ и $(ii)$.]
$$
x_n \underset{n\to \infty}{\longrightarrow} \infty
$$
$$
\Updownarrow
$$
$$
\text{ для всякого $\varepsilon>0$ соотношение $x_n\in
(-\infty;-E)\cup (E:+\infty)$ выполняется для почти всех $n\in
\mathbb{N}$ }
$$
$$
\phantom{\scriptsize \eqref{|x|>e}}\quad\Updownarrow\quad{\scriptsize \eqref{|x|>e}}
$$
$$
\text{ для всякого $\varepsilon>0$ соотношение $|x_n|\in
(-\infty;-E)\cup (E:+\infty)$ выполняется для почти всех $n\in
\mathbb{N}$ }
$$
$$
\Updownarrow
$$
$$
|x_n| \underset{n\to \infty}{\longrightarrow} \infty \qquad
$$ \end{proof}

\begin{tm}[\bf о связи между бесконечно малыми и бесконечно
большими последова\-тельностями]\label{0<->infty}

Пусть $x_n \ne 0$. Тогда последовательность $\{ x_n \}$ бесконечно малая тогда и только тогда, когда последовательность $\{ \frac{1}{x_n} \}$ бесконечно большая.
\end{tm}
\begin{proof}
$$
\text{$\{x_n\}$ -- бесконечно малая}
$$
$$
\Updownarrow
$$
$$
x_n\underset{n\to\infty}{\longrightarrow}0
$$
$$
\Updownarrow
$$
$$
\forall \e>0 \qquad |x_n|<\e \quad\text{-- верно для почти всех $n\in\N$}
$$
$$
\Updownarrow
$$
$$
\forall \e>0 \qquad \left|\frac{1}{x_n}\right|>\frac{1}{\e} \quad\text{-- верно для почти всех $n\in\N$}
$$
$$
\Updownarrow
$$
$$
\forall E>0 \qquad \left|\frac{1}{x_n}\right|>E \quad\text{-- верно для почти всех $n\in\N$}
$$
$$
\Updownarrow
$$
$$
\frac{1}{x_n}\underset{n\to\infty}{\longrightarrow}\infty
$$
$$
\Updownarrow
$$
$$
\text{$\left\{\frac{1}{x_n}\right\}$ -- бесконечно большая.}
$$
\end{proof}

\begin{tm}[\bf о связи между бесконечно малыми последовательностями
и конечными пределами]\label{x->a<=>x-a->0}

Последовательность $\{ x_n \}$ стремится к числу $a$
$$
x_n \underset{n\to \infty}{\longrightarrow} a
$$
тогда и только тогда, когда она имеет вид
$$
x_n=a+\alpha_n
$$
где $\{ \alpha_n \}$ -- некоторая бесконечно малая последовательность.
\end{tm}
\begin{proof} Положим
$$
\alpha_n=x_n-a
$$
и заметим сразу, что

\begin{multline*}
x_n \in (a-\varepsilon; a+\varepsilon) \quad \Longleftrightarrow
\quad a-\varepsilon<x_n< a+\varepsilon \quad \Longleftrightarrow
\quad
\begin{cases}{a-\varepsilon<x_n}\\{x_n< a+\varepsilon}\end{cases}
\quad \Longleftrightarrow \quad \\
\quad \Longleftrightarrow \quad
\begin{cases}{-\varepsilon<x_n-a}\\{x_n-a< \varepsilon}\end{cases}
\quad \Longleftrightarrow \quad -\varepsilon<x_n-a< \varepsilon \quad
\Longleftrightarrow \quad -\varepsilon<\alpha_n< \varepsilon \quad
\Longleftrightarrow \quad \alpha_n\in (-\varepsilon;\varepsilon)
\end{multline*} поэтому окрестность $(a-\varepsilon;
a+\varepsilon)$ точки $a$ содержит почти все элементы последовательности $x_n$
в том и только в том случае, если окрестность $(-\varepsilon; \varepsilon)$
точки $0$ содержит почти все элементы последовательности $\alpha_n$.

Отсюда следует, что если $x_n \underset{n\to \infty}{\longrightarrow} a$ (то
есть любая окрестность $(a-\varepsilon; a+\varepsilon)$ точки $a$ содержит
почти все элементы последовательности $x_n$), то это равносильно тому, что
$\alpha_n \underset{n\to \infty}{\longrightarrow} 0$ (то есть тому, что любая
окрестность $(-\varepsilon; \varepsilon)$ точки $0$ содержит почти все элементы
последовательности $\alpha_n$).
\end{proof}

\bigskip

\centerline{\bf Свойства бесконечно малых последовательностей}\label{svoistva-besk-malyh-posl}

 \bit{
\item[$1^\circ$] Если $\{ \alpha_n \}$ -- бесконечно малая последовательность,
то для любого числа $C\in \R$ последова\-тельность $\{ C\cdot \alpha_n \}$ --
тоже бесконечно малая.

\item[$2^\circ$] Если $\{ \alpha_n \}$ и $\{ \beta_n \}$ -- бесконечно малые
последовательности, то $\{ \alpha_n + \beta_n \}$ и $\{ \alpha_n - \beta_n \}$
-- тоже бесконечно малые последовательности.

\item[$3^\circ$] Если $\{ \alpha_n \}$ и $\{ \beta_n \}$ -- бесконечно малые
последовательности, то $\{ \alpha_n \cdot \beta_n \}$ -- тоже бесконечно малая
последовательность.

}\eit

\begin{proof}

1. Пусть $\{ \alpha_n \}$ -- бесконечно малая последовательность. Тогда если $C=0$, то
$$
C\cdot\alpha_n=0\underset{n\to\infty}{\longrightarrow}0
$$
и поэтому $\{ C\cdot \alpha_n \}$ -- тоже бесконечно малая. Рассмотрим
случай, когда $C\ne 0$. Тогда:
$$
\text{$\{ \alpha_n \}$ -- бесконечно малая}
$$
$$
\Updownarrow
$$
$$
\alpha_n\underset{n\to\infty}{\longrightarrow}0
$$
$$
\Updownarrow
$$
$$
\forall\e>o \qquad |\alpha_n|<\e \quad \text{-- верно для почти всех $n\in\N$}
$$
$$
\Updownarrow
$$
$$
\forall\e>o \qquad |\alpha_n|<\frac{\e}{|C|} \quad \text{-- верно для почти всех $n\in\N$}
$$
$$
\Updownarrow
$$
$$
\forall\e>o \qquad |C\cdot\alpha_n|<\e \quad \text{-- верно для почти всех $n\in\N$}
$$
$$
\Updownarrow
$$
$$
\text{$\{C\cdot \alpha_n \}$ -- бесконечно малая}
$$

2. Пусть $\{ \alpha_n \}$ и $\{ \beta_n \}$ -- бесконечно малые
последовательности, покажем что тогда $\{ \alpha_n + \beta_n \}$ -- тоже
бесконечно малая последовательность:
$$
\text{$\{ \alpha_n \}$ и $\{ \beta_n \}$ -- бесконечно малые}
$$
$$
\Downarrow
$$
$$
\alpha_n\underset{n\to\infty}{\longrightarrow}0,\qquad \beta_n\underset{n\to\infty}{\longrightarrow}0
$$
$$
\Downarrow
$$
$$
\forall\e>o \qquad |\alpha_n|<\e\quad\&\quad |\beta_n|<\e \quad \text{-- верно для почти всех $n\in\N$}
$$
$$
\Downarrow
$$
$$
\forall\e>o \qquad |\alpha_n|<\frac{\e}{2}\quad\&\quad |\beta_n|<\frac{\e}{2} \quad \text{-- верно для почти всех $n\in\N$}
$$
$$
\Downarrow
$$
$$
\forall\e>o \qquad |\alpha_n+\beta_n|\le |\alpha_n|+|\beta_n|<\frac{\e}{2}+\frac{\e}{2} \quad \text{-- верно для почти всех $n\in\N$}
$$
$$
\Downarrow
$$
$$
\text{$\{ \alpha_n + \beta_n \}$ -- бесконечно малая}
$$
Точно так же рассматривается случай $\{ \alpha_n - \beta_n \}$.

3. Пусть снова $\{ \alpha_n \}$ и $\{ \beta_n \}$ -- бесконечно малые
последовательности, покажем что $\{ \alpha_n \cdot \beta_n \}$ -- тоже
бесконечно малая последовательность:
$$
\text{$\{ \alpha_n \}$ и $\{ \beta_n \}$ -- бесконечно малые}
$$
$$
\Downarrow
$$
$$
\alpha_n\underset{n\to\infty}{\longrightarrow}0,\qquad \beta_n\underset{n\to\infty}{\longrightarrow}0
$$
$$
\Downarrow
$$
$$
\forall\e>o \qquad |\alpha_n|<\e\quad\&\quad |\beta_n|<\e \quad \text{-- верно для почти всех $n\in\N$}
$$
$$
\Downarrow
$$
$$
\forall\e>o \qquad |\alpha_n|<\sqrt{\e}\quad\&\quad |\beta_n|<\sqrt{\e} \quad \text{-- верно для почти всех $n\in\N$}
$$
$$
\Downarrow
$$
$$
\forall\e>o \qquad |\alpha_n\cdot\beta_n|=|\alpha_n|\cdot|\beta_n|<\sqrt{\e}\cdot\sqrt{\e}=\e \quad \text{-- верно для почти всех $n\in\N$}
$$
$$
\Downarrow
$$
$$
\text{$\{ \alpha_n\cdot\beta_n \}$ -- бесконечно малая.}
$$
 \end{proof}

\bigskip

\centerline{\bf Свойства бесконечно больших последовательностей}
 \bit{
\item[$1^\circ$] Если $C\ne 0$ и $y_n \underset{n\to
\infty}{\longrightarrow} \infty$, то $C\cdot y_n \underset{n\to
\infty}{\longrightarrow} \infty$.

\item[$2^\circ$] Если $x_n \underset{n\to \infty}{\longrightarrow}
C\ne \infty$, а $y_n \underset{n\to \infty}{\longrightarrow} \infty$,
то $x_n+y_n \underset{n\to \infty}{\longrightarrow} \infty$.

\item[$3^\circ$] Если $x_n \underset{n\to \infty}{\longrightarrow}
C\ne 0$, а $y_n \underset{n\to \infty}{\longrightarrow} \infty$, то
$x_n\cdot y_n \underset{n\to \infty}{\longrightarrow} \infty$.

\item[$4^\circ$]\label{arifm-oper-s-predelami-posl} Если $x_n
\underset{n\to \infty}{\longrightarrow} C\ne \infty$, а $y_n
\underset{n\to \infty}{\longrightarrow} \infty$, то $\frac{y_n}{x_n}
\underset{n\to \infty}{\longrightarrow} \infty$.

\item[$5^\circ$]\label{x_n->8,y_n>x_n=>y_n->8} Если $x_n \underset{n\to
\infty}{\longrightarrow} \infty$ и $|x_n|\le y_n$, то $y_n \underset{n\to
\infty}{\longrightarrow} \infty$.

}\eit

\begin{proof} Мы докажем только второе из этих свойств (остальные доказываются аналогично):
$$
\underbrace{
\begin{matrix}
x_n \underset{n\to \infty}{\longrightarrow} C\ne \infty &  y_n \underset{n\to \infty}{\longrightarrow} \infty \\
\Downarrow & \Downarrow \\
|x_n|-|C|\le\eqref{module-3^0}\le\kern-20pt\underbrace{|x_n-C|<1}_{\scriptsize \begin{matrix}\uparrow\\ \text{верно для почти всех $n\in\N$}\end{matrix}} &  \forall D>0 \underbrace{|y_n|>D+|C|+1}_{\scriptsize \begin{matrix}\uparrow\\ \text{верно для почти всех $n\in\N$}\end{matrix}} \\
\Downarrow &  \\
\underbrace{|x_n|<|C|+1}_{\scriptsize \begin{matrix}\uparrow\\ \text{верно для почти всех $n\in\N$}\end{matrix}} & \\
\end{matrix}}_{}
$$
$$
\Downarrow
$$
$$
\forall D>0\qquad |y_n+x_n|\ge\eqref{module-2^0}\ge \kern-10pt
\overbrace{|y_n|}^{\scriptsize \begin{matrix}D+|C|+1\\
\text{\rotatebox{90}{$>$}}  \end{matrix}}\kern-10pt\underbrace{-|x_n|}_{\scriptsize \begin{matrix}
\text{\rotatebox{90}{$<$}} \\ -|C|-1 \end{matrix}}>(D+|C|+1)-|C|-1=D\quad\text{-- верно для почти всех $n\in\N$}
$$
$$
\Downarrow
$$
$$
x_n+y_n \underset{n\to \infty}{\longrightarrow} \infty
$$
\end{proof}

\noindent\rule{160mm}{0.1pt}\begin{multicols}{2}

\bex Для любого $\alpha\ne 0$ последовательность $\{\alpha\cdot n\}_{n\in\N}$
является бесконечно большой, точнее:
 \beq\label{n-alpha->infty}
\boxed{\quad \lim_{n\to\infty}n\cdot\alpha=\begin{cases}+\infty,& \text{если
$\alpha>0$} \\ -\infty,& \text{если $\alpha<0$}
\end{cases}
\quad}
 \eeq
\eex
 \bpr
1. Пусть сначала $\alpha>0$. Тогда для любого $E>0$ условие
$$
n\cdot\alpha>E
$$
эквивалентно условию
$$
n>\frac{E}{\alpha}
$$
которое по теореме Архимеда \ref{Archimed-theo} выполняется для почти всех $n$.
Это означает, что
$$
n\cdot\alpha\underset{n\to\infty}{\longrightarrow}+\infty
$$

2. Если же $\alpha<0$, то для любого $E>0$ условие
$$
n\cdot\alpha<-E
$$
эквивалентно условию
$$
n>-\frac{E}{\alpha}
$$
и опять по теореме Архимеда \ref{Archimed-theo} это выполняется для почти всех
$n$. Это означает, что
$$
n\cdot\alpha\underset{n\to\infty}{\longrightarrow}-\infty
$$
 \epr

\bex Последовательность $\{n^k\}_{n\in\N}$ является
 \biter{
\item[--] бесконечно малой, если $k<0$;

\item[--] бесконечно большой, если $k>0$.
 }\eiter\noindent
Точнее,
 \beq\label{k>0=>n^k->infty}
\boxed{\quad \lim_{n\to \infty} n^k=
\begin{cases}0,& \text{\rm если}\quad k<0\\
 +\infty, & \text{\rm
если} \quad k>0\end{cases} \quad}
 \eeq
\eex
\begin{proof}

1. Пусть $k>0$. Применяя неравенство Бернулли \eqref{nerav-Bernoulli}, мы
получим
 \begin{multline*}
n^k=(1+(n-1))^k \overset{\eqref{nerav-Bernoulli}}{\ge}\\ \ge
1+k\cdot(n-1)=1-k+k\cdot
n\overset{\eqref{n-alpha->infty}}{\underset{n\to\infty}{\longrightarrow}}+\infty
 \end{multline*}
$$
\phantom{\text{\scriptsize $5^\circ$ на
с.\pageref{x_n->8,y_n>x_n=>y_n->8}}}\quad\Downarrow\quad\text{\scriptsize
$5^\circ$ на с.\pageref{x_n->8,y_n>x_n=>y_n->8}}
$$
$$
n^k\underset{n\to\infty}{\longrightarrow}\infty
$$

2. Пусть $k<0$. Тогда $-k>0$, и, по уже доказанному,
$$
n^{-k}\underset{n\to\infty}{\longrightarrow}\infty
$$
$$
\phantom{\text{\scriptsize теорема
\ref{0<->infty}}}\quad\Downarrow\quad\text{\scriptsize теорема \ref{0<->infty}}
$$
$$
n^k=\frac{1}{n^{-k}}\underset{n\to\infty}{\longrightarrow}0
$$
\end{proof}

\bex Последовательность $a^n$ является
 \biter{
\item[--] бесконечно малой, если $|a|<1$;

\item[--] бесконечно большой, если $|a|>1$.
 }\eiter
То есть,
 \beq\label{a>1=>a^n->infty}
 \boxed{\quad
\lim_{n\to \infty} a^n=
\begin{cases} \infty,& \text{если $|a|>1$} \\
0, &\text{если $|a|<1$} \end{cases} \quad }
 \eeq
\eex
\begin{proof}

1. Если $|a|>1$, то это означает либо $a>1$, либо $a<-1$. Рассмотрим каждый из
этих случаев отдельно. Если $a>1$, то применяя неравенство Бернулли
\eqref{nerav-Bernoulli}, мы получим
$$
|a|^n=a^n=(1+(a-1))^n\overset{\eqref{nerav-Bernoulli}}{\ge}
1+n\cdot(a-1)\underset{n\to\infty}{\longrightarrow}\infty
$$
$$
\phantom{\text{\scriptsize $5^\circ$ на
с.\pageref{x_n->8,y_n>x_n=>y_n->8}}}\quad\Downarrow\quad\text{\scriptsize
$5^\circ$ на с.\pageref{x_n->8,y_n>x_n=>y_n->8}}
$$
$$
a^n\underset{n\to\infty}{\longrightarrow}\infty
$$
Если же $a<-1$, то $|a|=-a>1$, и по уже доказанному,
$$
|a|^n=(-a)^n\underset{n\to\infty}{\longrightarrow}\infty
$$

2. Пусть наоборот, $|a|<1$. Тогда $\left|\frac{1}{a}\right|=\frac{1}{|a|}>1$, и
мы получаем
$$
\frac{1}{|a|^n}=\left|\frac{1}{a}\right|^n\overset{\scriptsize
\begin{matrix}\text{уже}\\ \text{доказано} \\ \phantom{.} \end{matrix}}{\underset{n\to\infty}{\longrightarrow}}\infty
$$
$$
\phantom{\text{\scriptsize теорема
\ref{0<->infty}}}\quad\Downarrow\quad\text{\scriptsize теорема \ref{0<->infty}}
$$
$$
|a|^n\underset{n\to\infty}{\longrightarrow}0
$$

\end{proof}

\begin{er}
Какие из последовательностей являются бесконечно малыми, а какие -- бесконечно
большими?
 \biter{
 \item[1)] $x_n=\frac{1}{n}$
 \item[2)] $x_n=\frac{1}{n^2}$
 \item[3)] $x_n=n$
 \item[4)] $x_n=(-1)^n\cdot n$
 \item[5)] $x_n=\frac{(-1)^n}{n}$
 \item[7)] $x_n=n^2$
 }\eiter
\end{er}

\end{multicols}\noindent\rule[10pt]{160mm}{0.1pt}

\subsection{Арифметические операции с пределами}
\label{lim-ariphm}

Следующие свойства пределов связаны с арифметическими операциями.

{\it
 \bit{\label{lim-ariphm-1}
\item[$1^0.$] Если $x_n \underset{n\to \infty}{\longrightarrow} a$, то для
всякого числа $C\in \R$ справедливо $C\cdot x_n \underset{n\to
\infty}{\longrightarrow} C\cdot a$.

\item[$2^0.$] Если $x_n \underset{n\to \infty}{\longrightarrow} a$ и
$y_n \underset{n\to \infty}{\longrightarrow} b$, то $x_n+y_n
\underset{n\to \infty}{\longrightarrow} a+b$ и $x_n-y_n
\underset{n\to \infty}{\longrightarrow} a-b$.

\item[$3^0.$]\label{x_n->x,y_n->y=>x_n-cdot-y_n->x-cdot-y} Если $x_n
\underset{n\to \infty}{\longrightarrow} a$ и $y_n \underset{n\to
\infty}{\longrightarrow} b$, то $x_n\cdot y_n \underset{n\to
\infty}{\longrightarrow} a\cdot b$.

\item[$4^0.$] Если $x_n \underset{n\to \infty}{\longrightarrow} a$ и
$y_n \underset{n\to \infty}{\longrightarrow} b$, причем $a\ne 0$ и
$x_n\ne 0$ (при любом $n$), то $\frac{y_n}{x_n} \underset{n\to
\infty}{\longrightarrow} \frac{b}{a}$.
 }\eit
}

\begin{proof}[Доказательство свойств $1^0$-$3^0$]

1. Свойство $1^0$ доказывается такой логической цепочкой:
$$
x_n \underset{n\to\infty}{\longrightarrow} a
$$
$$
\phantom{\scriptsize\text{теорема \ref{x->a<=>x-a->0}}}\quad\Downarrow\quad{\scriptsize\text{теорема \ref{x->a<=>x-a->0}}}
$$
$$
\kern5pt x_n=a+\kern-5pt\underbrace{\alpha_n}_{\scriptsize\begin{matrix}\phantom{\tiny \begin{matrix}n\\ \downarrow\\ \infty \end{matrix}}\ \downarrow\ {\tiny \begin{matrix}n\\ \downarrow\\ \infty \end{matrix}}\\ 0 \end{matrix}}
$$
$$
\Downarrow
$$
$$
\kern83pt C\cdot x_n=C\cdot a+\kern-83pt\underbrace{C\cdot \alpha_n}_{\scriptsize\begin{matrix}\phantom{\tiny \begin{matrix}n\\ \downarrow\\ \infty \end{matrix}}
\ \phantom{\text{(свойство $1^0$ на с.\pageref{svoistva-besk-malyh-posl})}}\quad \downarrow\ {\tiny \begin{matrix}n\\ \downarrow\\ \infty \end{matrix}}\quad \text{(свойство $1^0$ на с.\pageref{svoistva-besk-malyh-posl})}\\ 0 \end{matrix}}
$$
$$
\phantom{\scriptsize\text{теорема \ref{x->a<=>x-a->0}}}\quad\Downarrow\quad{\scriptsize\text{теорема \ref{x->a<=>x-a->0}}}
$$
$$
C\cdot x_n \underset{n\to\infty}{\longrightarrow} C\cdot a
$$

2. Свойство $2^0$:
$$
x_n \underset{n\to\infty}{\longrightarrow} a,\qquad y_n \underset{n\to\infty}{\longrightarrow} b
$$
$$
\phantom{\scriptsize\text{теорема \ref{x->a<=>x-a->0}}}\quad\Downarrow\quad{\scriptsize\text{теорема \ref{x->a<=>x-a->0}}}
$$
$$
 x_n=a+\kern-5pt\underbrace{\alpha_n}_{\scriptsize\begin{matrix}\phantom{\tiny \begin{matrix}n\\ \downarrow\\ \infty \end{matrix}}\ \downarrow\ {\tiny \begin{matrix}n\\ \downarrow\\ \infty \end{matrix}}\\ 0 \end{matrix}},\qquad
 y_n=b+\kern-5pt\underbrace{\beta_n}_{\scriptsize\begin{matrix}\phantom{\tiny \begin{matrix}n\\ \downarrow\\ \infty \end{matrix}}\ \downarrow\ {\tiny \begin{matrix}n\\ \downarrow\\ \infty \end{matrix}}\\ 0 \end{matrix}}
$$
$$
\Downarrow
$$
$$
\kern78pt x_n+y_n=a+b+\kern-78pt\underbrace{\alpha_n+\beta_n}_{\scriptsize\begin{matrix}\phantom{\tiny \begin{matrix}n\\ \downarrow\\ \infty \end{matrix}}
\ \phantom{\text{(свойство $2^0$ на с.\pageref{svoistva-besk-malyh-posl})}}\quad \downarrow\ {\tiny \begin{matrix}n\\ \downarrow\\ \infty \end{matrix}}\quad \text{(свойство $2^0$ на с.\pageref{svoistva-besk-malyh-posl})}\\ 0 \end{matrix}}
$$
$$
\phantom{\scriptsize\text{теорема \ref{x->a<=>x-a->0}}}\quad\Downarrow\quad{\scriptsize\text{теорема \ref{x->a<=>x-a->0}}}
$$
$$
x_n+y_n \underset{n\to\infty}{\longrightarrow} a+b
$$
И точно так же для $x_n-y_n$.

3. Свойство $3^0$:
$$
x_n \underset{n\to\infty}{\longrightarrow} a,\qquad y_n \underset{n\to\infty}{\longrightarrow} b
$$
$$
\phantom{\scriptsize\text{теорема \ref{x->a<=>x-a->0}}}\quad\Downarrow\quad{\scriptsize\text{теорема \ref{x->a<=>x-a->0}}}
$$
$$
 x_n=a+\kern-5pt\underbrace{\alpha_n}_{\scriptsize\begin{matrix}\phantom{\tiny \begin{matrix}n\\ \downarrow\\ \infty \end{matrix}}\ \downarrow\ {\tiny \begin{matrix}n\\ \downarrow\\ \infty \end{matrix}}\\ 0 \end{matrix}},\qquad
 y_n=b+\kern-5pt\underbrace{\beta_n}_{\scriptsize\begin{matrix}\phantom{\tiny \begin{matrix}n\\ \downarrow\\ \infty \end{matrix}}\ \downarrow\ {\tiny \begin{matrix}n\\ \downarrow\\ \infty \end{matrix}}\\ 0 \end{matrix}}
$$
$$
\Downarrow
$$
$$
\kern78pt x_n\cdot y_n=a\cdot b+\kern-45pt\underbrace{\underbrace{\alpha_n\cdot b}_{\scriptsize\begin{matrix}\phantom{\tiny \begin{matrix}n\\ \downarrow\\ \infty \end{matrix}}
\ \downarrow\ {\tiny \begin{matrix}n\\ \downarrow\\ \infty \end{matrix}}\\ 0 \\ \begin{pmatrix}\text{свойство $1^0$}\\ \text{на с.\pageref{svoistva-besk-malyh-posl}}\end{pmatrix}\end{matrix}}
\kern-12pt +\kern-12pt\overbrace{a\cdot\beta_n}^{\scriptsize\begin{matrix}
\begin{pmatrix}\text{свойство $1^0$}\\ \text{на с.\pageref{svoistva-besk-malyh-posl}}\end{pmatrix}\\ 0\\
\phantom{\tiny \begin{matrix}n\\ \downarrow\\ \infty \end{matrix}}
\ \uparrow\ {\tiny \begin{matrix}\infty\\ \uparrow\\ n \end{matrix}}\end{matrix}}
\kern-12pt +\kern-12pt\underbrace{\alpha_n\cdot\beta_n}_{\scriptsize\begin{matrix}\phantom{\tiny \begin{matrix}n\\ \downarrow\\ \infty \end{matrix}}
\ \downarrow\ {\tiny \begin{matrix}n\\ \downarrow\\ \infty \end{matrix}}\\ 0 \\ \begin{pmatrix}\text{свойство $3^0$}\\ \text{на с.\pageref{svoistva-besk-malyh-posl}}\end{pmatrix}\end{matrix}}}_{\scriptsize\begin{matrix}\phantom{\tiny \begin{matrix}n\\ \downarrow\\ \infty \end{matrix}}
\ \phantom{\text{(свойство $2^0$ на с.\pageref{svoistva-besk-malyh-posl})}}\quad \downarrow\ {\tiny \begin{matrix}n\\ \downarrow\\ \infty \end{matrix}}\quad \text{(свойство $2^0$ на с.\pageref{svoistva-besk-malyh-posl})}\\ 0 \end{matrix}}
$$
$$
\phantom{\scriptsize\text{теорема \ref{x->a<=>x-a->0}}}\quad\Downarrow\quad{\scriptsize\text{теорема \ref{x->a<=>x-a->0}}}
$$
$$
x_n\cdot y_n \underset{n\to\infty}{\longrightarrow} a\cdot b
$$
\epr

Для доказательства последнего свойства нам понадобится следующая

\begin{lm}\label{1/x} Если последовательность
$\{ x_n \}$  стремится к числу $a$, причем $x_n\ne 0$ и $a\ne 0$, то
последовательность $\{ \frac{1}{x_n} \}$ стремится к числу $\frac{1}{a}$:
$$
x_n \underset{n\to \infty}{\longrightarrow} a \quad (x_n\ne 0
\quad \& \quad a\ne 0) \quad \Longrightarrow \quad \frac{1}{x_n}
\underset{n\to \infty}{\longrightarrow} \frac{1}{a}
$$
\end{lm}
\begin{proof} По условию, $a\ne 0$, значит либо
$a<0$, либо $a>0$.

1. Рассмотрим сначала случай, когда $a>0$. Тогда в силу \eqref{x>0=>x^(-1)>0}, $\frac{1}{a}>0$. Выберем произвольное
$\varepsilon$ такое, что
$$
0< \varepsilon< \frac{1}{a}
$$
Тогда
$$
\e\cdot a>0>-\e\cdot a
$$
$$
\Downarrow
$$
$$
1+\e\cdot a>1>1-\e\cdot a
$$
$$
\Downarrow
$$
$$
\frac{1}{1+\e\cdot a}<1<\frac{1}{1-\e\cdot a}
$$
$$
\Downarrow
$$
$$
\frac{a}{1+\e\cdot a}<a<\frac{a}{1-\e\cdot a}
$$
$$
\phantom{{\scriptsize\text{предложение \ref{PROP:x-in-(a,b)=>U(x)-sub-(a,b)}}}}\quad\Downarrow\quad{\scriptsize\text{предложение \ref{PROP:x-in-(a,b)=>U(x)-sub-(a,b)}}}
$$
$$
\exists \delta>0:\qquad \frac{a}{1+\e\cdot a}<a-\delta<a<a+\delta<\frac{a}{1-\e\cdot a}
$$
Запомним это число $\delta$ и заметим теперь цепочку:
$$
x_n\underset{n\to\infty}{\longrightarrow}a
$$
$$
\Downarrow
$$
$$
x_n\in(a-\delta,a+\delta)\quad\text{-- выполняется для почти всех $n\in\N$}
$$
$$
\Downarrow
$$
$$
\frac{a}{1+\e\cdot a}<a-\delta<x_n<a+\delta<\frac{a}{1-\e\cdot a} \quad\text{-- выполняется для почти всех $n\in\N$}
$$
$$
\Downarrow
$$
$$
\frac{a}{1+\e\cdot a}<x_n<\frac{a}{1-\e\cdot a} \quad\text{-- выполняется для почти всех $n\in\N$}
$$
$$
\Downarrow
$$
$$
\frac{1+\e\cdot a}{a}>\frac{1}{x_n}>\frac{1-\e\cdot a}{a} \quad\text{-- выполняется для почти всех $n\in\N$}
$$
$$
\Downarrow
$$
$$
\frac{1}{a}+\e>\frac{1}{x_n}>\frac{1}{a}-\e \quad\text{-- выполняется для почти всех $n\in\N$}
$$
$$
\Downarrow
$$
$$
\frac{1}{x_n}\in\l\frac{1}{a}-\e,\frac{1}{a}+\e\r \quad\text{-- выполняется для почти всех $n\in\N$}
$$

Мы получили, что любая ``маленькая'' окрестность $\left(\frac{1}{a} - \varepsilon;
\frac{1}{a} + \varepsilon \right)$ точки $\frac{1}{a}$ (радиуса $\varepsilon <
\frac{1}{a}$) содержит почти все элементы последовательности $\frac{1}{x_n}$.
Отсюда следует, что любая ``большая'' окрестность $\left(\frac{1}{a}
- \varepsilon;\frac{1}{a} + \varepsilon \right)$ точки $\frac{1}{a}$ (радиуса
$\varepsilon \ge \frac{1}{a}$) тоже содержит почти все элементы последовательности
$\frac{1}{x_n}$ (потому что в ней можно выбрать какую-нибудь ``маленькую''
окрестность, которая будет содержать почти все $\frac{1}{x_n}$). Вместе это означает, что $\frac{1}{x_n} \underset{n\to \infty}{\longrightarrow}
\frac{1}{a}$.

2. Мы доказали нашу лемму для случая, когда $a>0$. Рассмотрим теперь случай
$a<0$. Тогда будет по уже доказанному свойству $1^0$ с.\pageref{lim-ariphm-1} последовательность $y_n=-x_n$ стремиться к числу $-a>0$. Значит, в силу уже доказанного в нашей лемме,
$$
\frac{1}{y_n} \underset{n\to \infty}{\longrightarrow} \frac{1}{-a}
$$
откуда опять по  свойству $1^0$ с.\pageref{lim-ariphm-1},
$$
\frac{1}{x_n}= -\frac{1}{y_n} \underset{n\to
\infty}{\longrightarrow} -\frac{1}{-a}= \frac{1}{a}
$$
\end{proof}

\bpr[Доказательство свойства $4^0$.]
$$
\underbrace{
\begin{matrix}
\overset{\scriptsize\begin{matrix}0\\ \text{\rotatebox{90}{$\ne$}}\end{matrix}}{x_n} \underset{n\to\infty}{\longrightarrow} \overset{\scriptsize\begin{matrix}0\\ \text{\rotatebox{90}{$\ne$}}\end{matrix}}{a},& y_n \underset{n\to\infty}{\longrightarrow} b \\
{\scriptsize\text{лемма \ref{1/x}}}\ \Downarrow \ \phantom{\scriptsize\text{лемма \ref{1/x}}} & \\
\frac{1}{x_n} \underset{n\to\infty}{\longrightarrow} \frac{1}{a} &
\end{matrix}
}
$$
$$
\phantom{\scriptsize\text{свойство $3^0$ на с.\pageref{lim-ariphm-1}}}\ \Downarrow \ {\scriptsize\text{свойство $3^0$ на с.\pageref{lim-ariphm-1}}}
$$
$$
\frac{y_n}{x_n}=\frac{1}{x_n}\cdot y_n \underset{n\to\infty}{\longrightarrow} \frac{1}{a}\cdot b=\frac{b}{a}
$$
\end{proof}

\subsection{Ограниченные последовательности}
\label{bounded-seq}

Последовательность  $\{ x_n \}$ называется {\it
ограниченной}\index{последовательность!ограниченной}, если множество чисел
$X=\{ x_n \}$ ограничено (сверху и снизу), то есть существуют такие числа
$C,D$, что для любого $n\in \mathbb{N}$ выполняется
$$
C\le x_n\le D
$$

%\picture{0pt}{0pt}{75.pcx}
\vglue100pt \noindent

\noindent\rule{160mm}{0.1pt}\begin{multicols}{2}

\begin{ex}
Проверьте, что следующие последовательности ограничены:
$$
x_n= (-1)^n, \, x_n= \frac{n+1}{n},\, x_n=\left\{\frac{n^2+1}{n}\right\}.
$$
Какие из них сходятся (то есть имеют конечный предел)?
\end{ex}

\end{multicols}\noindent\rule[10pt]{160mm}{0.1pt}

\begin{tm}[\bf об ограниченности сходящейся последовательности]
Всякая сходящаяся последовательность ограничена.
\end{tm}

\bigskip

\centerline{\bf Свойства ограниченных последовательностей}
 \bit{
\item[$1^{\circ\circ}$] Если $\{ x_n \}$ -- ограниченная последовательность, а
$\{ y_n \}$ -- бесконечно малая последовательность, то $\{ x_n \cdot y_n \}$ --
бесконечно малая последовательность.

\item[$2^{\circ\circ}$] Если $\{ x_n \}$ -- ограниченная последовательность, а
$\{ y_n \}$  -- бесконечно большая последовательность, то $\{ \frac{x_n}{y_n}
\}$ -- бесконечно малая последовательность.

\item[$3^{\circ\circ}$] Если $\{ x_n \}$ -- ограниченная последовательность, а
$\{ y_n \}$ -- бесконечно большая последовательность, то $\{ x_n+y_n \}$ --
бесконечно большая последовательность.

}\eit

\noindent\rule{160mm}{0.1pt}\begin{multicols}{2}
\begin{ex}
Найти предел
$$
 \lim\limits_{n\to \infty} \frac{(-1)^n}{n}
$$
Поскольку последовательность в числителе ($x_n=(-1)^n$) ограничена, а
последовательность в знаменателе ($y_n=n$) бесконечно большая, по свойству
$2^{00}$ дробь $\frac{x_n}{y_n}=\frac{(-1)^n}{n}$ должна быть бесконечно малой,
то есть
$$
 \lim\limits_{n\to \infty} \frac{(-1)^n}{n}=0
$$
\end{ex}

\begin{ers}
Найдите пределы
 \begin{align*}
& \lim\limits_{n\to \infty} \frac{\left\{\frac{n^3}{n+1}\right\}}{n},\\
& \lim\limits_{n\to \infty} \frac{n+5}{n+\left\{n+\frac{1}{2}\right\}},\\
& \lim\limits_{n\to \infty}
\frac{n+\left\{\frac{n^2+n+1}{n+2}\right\}}{n-\left\{\frac{n^3+5}{n+5}\right\}}.
 \end{align*}
\end{ers}
\end{multicols}\noindent\rule[10pt]{160mm}{0.1pt}

\subsection{Вычисление пределов}

Теперь, когда доказаны главные свойства сходящихся последовательностей, мы,
наконец, можем научиться вычислять простейшие пределы. Для этого полезно
переформулировать утверждения из \ref{SEC:predel-posledov}\ref{->0-&-->infty} и
\ref{SEC:predel-posledov}\ref{lim-ariphm} в виде следующей серии формул.

Это, во-первых, ``арифметические формулы'':
\begin{align}
\lim_{n\to \infty} (x_n+y_n) &=\lim_{n\to \infty} x_n + \lim_{n\to
\infty} y_n \label{1.5.1}
 \\
\lim_{n\to \infty} (x_n-y_n) &=\lim_{n\to \infty} x_n - \lim_{n\to
\infty} y_n \label{1.5.2}
 \\
\lim_{n\to \infty} (C\cdot x_n) &=C\cdot \lim_{n\to \infty} x_n
\label{1.5.3}
 \\
\lim_{n\to \infty} (x_n\cdot y_n) &=\lim_{n\to \infty} x_n \cdot
\lim_{n\to \infty} y_n \label{1.5.4}
 \\
\lim_{n\to \infty} \frac{y_n}{x_n} &= \frac{\lim_{n\to \infty}
y_n}{\lim_{n\to \infty} x_n} \label{1.5.5}
\end{align}

И, во-вторых -- ``эвристические'' (здесь $C$ обозначает ненулевое
число: $0\ne C\ne \infty$):
\begin{align}
\frac{C}{\infty} &=\, 0 & \frac{C}{0} &=\,\infty, \label{1.5.6}
 \\
\frac{\infty}{C} &=\, \infty & \frac{0}{C} &=\,0 \label{1.5.7}
 \\
C+\infty &=\, \infty & C\cdot \infty  &=\,\infty \label{1.5.8}
 \\
\infty+\infty &=\, \infty & \infty\cdot \infty  &=\,\infty
\label{1.5.8-1}
 \\
\frac{\infty}{\infty} &=\; ? & \frac{0}{0} &=\; ? \label{1.5.9}
 \\
\infty-\infty &=\; ? & 0\cdot \infty &=\; ? \label{1.5.10}
\end{align}\noindent
Здесь вопросительные знаки в последних четырех формулах означают, что в
рассматриваемых случаях сказать определенно, чему равен предел, невозможно, и
требуется дополнительное исследование. Эти ситуации объединяются общим термином
{\it неопределенность}\index{неопределенность}, и что он означает легче всего
понять из следующих примеров.

\noindent\rule{160mm}{0.1pt}\begin{multicols}{2}

\begin{ex}
Вычислим предел
$$
\lim_{n\to \infty} \frac{n-1}{2n+3}
$$
Если попробовать применить формулы \eqref{1.5.1}--\eqref{1.5.5}
сразу, то мы увидим, что ответа не получается:
\begin{multline*}
\lim_{n\to \infty}
\frac{n-1}{2n+3}
 \overset{\tiny\begin{matrix}\eqref{1.5.4}\\ \downarrow\end{matrix}}{=}
 \frac{\lim\limits_{n\to \infty}(n-1)}{\lim\limits_{n\to \infty} (2n+3)}
 \overset{\tiny\begin{matrix}\eqref{1.5.3}\\ \downarrow\end{matrix}}{=}
 \\= \frac{\lim\limits_{n\to \infty} n-1}{2\cdot
\lim\limits_{n\to \infty} n+3}= \frac{\infty -1}{2\cdot \infty
+3}
  \overset{\tiny\begin{matrix}\eqref{1.5.8}\\ \downarrow\end{matrix}}{=}
\frac{\infty}{\infty} \overset{\tiny\begin{matrix}\eqref{1.5.9}\\
\downarrow\end{matrix}}{=}?
 \end{multline*}
Появившийся вопросительный знак означает, что таким способом найти предел
невозможно, и это тот случай, когда принято говорить, что возникает
неопределенность. Это в свою очередь означает, что нужно искать какой-то другой
способ (при котором этой неопределенности не возникнет).

Для этого обычно бывает достаточно преобразовать подходящим образом выражение
под пределом. В нашем случае можно поделить числитель и знаменатель на $n$:
 \begin{multline*}
\lim\limits_{n\to \infty}
\kern-20pt\underbrace{\frac{n-1}{2n+3}}_{\smsize
\begin{matrix} \text{делим числитель}
\\ \text{и знаменатель на $n$} \end{matrix}} \kern-20pt=
\lim\limits_{n\to \infty} \frac{1-\frac{1}{n}}{2+\frac{3}{n}}=
\frac{\lim\limits_{n\to \infty} (1-\frac{1}{n})} {\lim\limits_{n\to
\infty} (2+\frac{3}{n})}=\\= \frac{1-\lim\limits_{n\to \infty}
\frac{1}{n}} {2+\lim\limits_{n\to \infty} \frac{3}{n}}= \frac{1-
\lim\limits_{n\to \infty} \frac{1}{n}} {2+3\cdot \lim\limits_{n\to
\infty} \frac{1}{n}}= \frac{1-0} {2+3\cdot 0}= \frac{1}{2}
\end{multline*} Итак, со второй попытки мы получили ответ:
$\frac{1}{2}$.

На этом обсуждение этого примера можно закончить, но перед этим мы
хотели бы показать, как можно сократить количество писанины в
подобных вычислениях, одновременно, сделав их более наглядными.
Коротко и просто нашу первую попытку можно записать так:
$$
\lim_{n\to \infty} \frac{n-1}{2n+3}=\frac{\infty-1}{2\cdot\infty+3}=
\frac{\infty}{\infty}=?
$$
А вторую -- так:
$$
\lim_{n\to \infty} \kern-20pt\underbrace{\frac{n-1}{2n+3}}_{\smsize
\begin{matrix} \text{делим числитель}
\\ \text{и знаменатель на $n$} \end{matrix}} \kern-20pt=\lim_{n\to \infty}
 \frac{{1-
 {\smsize
 \boxed{\frac{1}{n}}\put(3,12){\vector(2,1){10}\put(3,3){0}}
 }
 }}{2+
 {\smsize
 \boxed{\frac{3}{n}}\put(2,-10){\vector(2,-1){10}\put(3,-9){0}}
 }
 }=
\frac{1-0}{2+0}= \frac{1}{2}
$$
Ниже мы постараемся всюду пользоваться такой ``наглядной'' записью,
поэтому будет хорошо, если читатель уже сейчас обратит на нее
внимание.
\end{ex}

\begin{ex}
Вычислим предел
$$
\lim\limits_{n\to \infty} \frac{3n^2+5n+7}{8n^2-4n+1}
$$
Здесь также первая попытка ничего не дает, потому что
неопределенность возникает на этот раз в знаменателе:
$$
\lim_{n\to \infty} \frac{3n^2+5n+7}{8n^2-4n+1}=
\frac{3\cdot\infty^2+5\cdot\infty+7}{8\cdot\infty^2-4\cdot\infty+1}=
\frac{\infty}{?}
$$
Поэтому мы снова преобразуем дробь, прежде чем пользоваться правилами
\eqref{1.5.1}--\eqref{1.5.5}:
 \begin{multline*}
\lim_{n\to \infty}
\kern-5pt\underbrace{\frac{3n^2+5n+7}{8n^2-4n+1}}_{\smsize
\begin{matrix} \text{делим числитель}
\\ \text{и знаменатель на $n^2$} \end{matrix}} \kern-5pt
 =
\lim\limits_{n\to \infty} \frac{3+
 {\smsize
 \boxed{\frac{5}{n}}\put(2,15){\vector(1,1){10}\put(3,10){0}}
 }
 +
 {\smsize
 \boxed{\frac{7}{n^2}}\put(2,15){\vector(1,1){10}\put(3,10){0}}
 }
 }{8-
 {\smsize
 \boxed{\frac{4}{n}}\put(2,-10){\vector(1,-1){10}\put(2,-15){0}}
 }
  +
  {\smsize
 \boxed{\frac{1}{n^2}}\put(2,-10){\vector(1,-1){10}\put(2,-15){0}}
 }
 }= \\ =
\frac{3+0+0}{8-0-0}= \frac{3}{8} \end{multline*}
\end{ex}

\begin{ex}
Вычислим предел
$$
\lim\limits_{n\to \infty} \frac{5n^2-n+1}{n^3+n^2+2}
$$
Здесь у нас уже достаточно опыта, чтобы сообразить сразу, что первая
попытка ``подставить в лоб'' снова приведет к неопределенности,
поэтому мы можем сразу перейти к преобразованию дроби:
 \begin{multline*}
 \lim_{n\to \infty} \kern-15pt
\underbrace{\frac{5n^2-n+1}{n^3+n^2+2}}_{\smsize
 \begin{matrix} \text{делим числитель} \\
\text{и знаменатель на} \\ \text{максимальную степень $n$,}
\\ \text{то есть на $n^3$} \end{matrix}}\kern-15pt=
 \lim_{n\to \infty}
 \frac{
 {\smsize
 \boxed{\frac{5}{n}}\put(2,15){\vector(1,1){10}\put(3,10){0}}
 }
 -
 {\smsize
 \boxed{\frac{1}{n^2}}\put(2,15){\vector(1,1){10}\put(3,10){0}}
 }
 +
 {\smsize
 \boxed{\frac{1}{n^3}}\put(2,15){\vector(1,1){10}\put(3,10){0}}
 }
 }
 {1+
 {\smsize
 \boxed{\frac{1}{n}}\put(2,-10){\vector(1,-1){10}\put(2,-15){0}}
 }
 +
 {\smsize
 \boxed{\frac{2}{n^3}}\put(2,-10){\vector(1,-1){10}\put(2,-15){0}}
 }
 }= \\=
 \frac{0-0+0} {1+0+0}=0 \end{multline*}
\end{ex}

\begin{ex} В следующем примере читатель уже сразу должен
сообразить, что без преобразований ничего не получится, поэтому о
``подстановке в лоб'' мы уже и не заикаемся:
 \begin{multline*}
 \lim_{n\to \infty} \kern-18pt
\underbrace{\frac{n^4+n^2+5}{n^3-3n+2}}_{\smsize
 \begin{matrix} \text{делим числитель} \\
 \text{и знаменатель на}
\\ \text{максимальную степень $n$,}
\\ \text{то есть на $n^4$} \end{matrix}}\kern-18pt
 =
 \lim_{n\to \infty}
 \frac{1+
 {\smsize
 \boxed{\frac{1}{n^2}}\put(2,15){\vector(1,1){10}\put(3,10){0}}
 }
 +
 {\smsize
 \boxed{\frac{5}{n^4}}\put(2,15){\vector(1,1){10}\put(3,10){0}}
 }
 }
 {
 {\smsize
 \boxed{\frac{1}{n}}\put(2,-10){\vector(1,-1){10}\put(2,-15){0}}
 }
 -
 {\smsize
 \boxed{\frac{3}{n^3}}\put(2,-10){\vector(1,-1){10}\put(2,-15){0}}
 }
 +
 {\smsize
 \boxed{\frac{2}{n^4}}\put(2,-10){\vector(1,-1){10}\put(2,-15){0}}
 }
 }=\\=
 \frac{1+0+0}{0-0+0}= \frac{1}{0}
 \underset{\tiny \begin{matrix}\uparrow \\ \eqref{1.5.6}\end{matrix}}{=}
 \infty
 \end{multline*}
\end{ex}

\begin{ex}
 \begin{multline*}
\lim\limits_{n\to \infty}
 \overbrace{\frac{(n+5)^3-n(n+7)^2}{n^2}}^{\smsize
 \text{раскрываем скобки}} =\\= \lim\limits_{n\to \infty}
\Big\{n^3+15n^2+75n+125-\\-n(n^2+14n+49)\Big\}:\Big\{n^2\Big\}=\\=
\lim\limits_{n\to \infty} \frac{15n^2+75n+125-14n^2+49n}{n^2}= \\=
\lim\limits_{n\to \infty} \frac{n^2+26n+125}{n^2}=\\=
\lim\limits_{n\to \infty}
 \Bigg(1+
 {\smsize
 \boxed{\frac{26}{n}}\put(2,-10){\vector(1,-1){10}\put(2,-15){0}}
 }
 +
 {\smsize
 \boxed{\frac{125}{n^2}}\put(2,-10){\vector(1,-1){10}\put(2,-15){0}}
 } \;
 \Bigg)=1+0+0=1 \end{multline*}
\end{ex}

\bex Докажем формулу:
 \beq\label{C_n^k/n^l->?}
 \boxed{\quad
\frac{C_n^k}{n^l}\underset{n\to\infty}{\longrightarrow}\begin{cases}0,& k<l\\
\frac{1}{k!},& k=l \\ \infty,& k>l\end{cases}\quad}
 \eeq
\eex \bpr Начнем со случая $k=l$:
 \begin{multline*}
\frac{C_n^k}{n^k}=\frac{n!}{n^k\cdot
k!\cdot(n-k)!}=\\=\frac{\overbrace{n\cdot(n-1)\cdot...\cdot(n-k+1)}^{\text{$k$
множителей}}}{\underbrace{n\cdot n\cdot...\cdot n}_{\text{$k$ множителей}}\cdot
k!}=\\
=\frac{\Big( 1-\boxed{\text{\scriptsize
$\frac{1}{n}$}}\put(2,15){\vector(1,1){10}\put(3,10){0}}\
\Big)\cdot...\cdot\Big( 1-\boxed{\text{\scriptsize
$\frac{k-1}{n}$}}\put(2,15){\vector(1,1){10}\put(3,10){0}}\ \Big)}{k!}
\underset{n\to\infty}{\longrightarrow}\frac{1}{k!}
 \end{multline*}
После этого при $k<l$ мы получим:
$$
\frac{C_n^k}{n^l}=
\boxed{\frac{1}{n^{l-k}}}\put(2,-12){\vector(1,-1){10}\put(2,-15){0}} \cdot
\boxed{\frac{C_n^k}{n^k}}\put(2,-12){\vector(1,-1){10}\put(2,-15){$\frac{1}{k!}$}}
\underset{n\to\infty}{\longrightarrow}0
$$
А при $k>l$ мы получим:
$$
\frac{C_n^k}{n^l}=
\boxed{n^{k-l}}\put(-10,-6){\vector(1,-1){10}\put(-2,-17){$\infty$}} \cdot
\boxed{\frac{C_n^k}{n^k}}\put(2,-12){\vector(1,-1){10}\put(2,-15){$\frac{1}{k!}$}}
\underset{n\to\infty}{\longrightarrow}\infty
$$
 \epr

\begin{ers}
Вычислите пределы
 \biter{
\item[1)] $\lim\limits_{n\to \infty} \frac{2-n}{3n+4}$

\item[2)] $\lim\limits_{n\to \infty} \frac{n^2-3}{3n-4n^2}$

\item[3)] $\lim\limits_{n\to \infty} \frac{3+ 2n+ n^2}{4-9n+4n^3}$

\item[4)] $\lim\limits_{n\to \infty}
\frac{n^2+1}{2n+1}-\frac{3n^2+1}{6n+1}$

\item[5)] $\lim\limits_{n\to \infty} \frac{n^2+2n+1}{n+5}$

\item[6)] $\lim\limits_{n\to \infty} \frac{n^3-n^2+1}{5-n-n^2}$
 }\eiter
\end{ers}

\end{multicols}\noindent\rule[10pt]{160mm}{0.1pt}

\section{Классические теоремы о последовательностях}

\subsection{Предельный переход в неравенствах}

\begin{tm}[\bf о предельном переходе в неравенстве с конечным пределом]\label{x_n<_y_n}\index{теорема!о предельном переходе в неравенстве}\footnote{Этот
результат используется ниже в разных ситуациях, в частности, в главе
\ref{ch-cont-f(x)} при доказательстве теоремы \ref{sign-pres} о сохранении
знака непрерывной функцией} Если последовательности $\{ x_n \}$, $\{ y_n \}$
связаны неравенством
\begin{equation}
x_n\le y_n \label{2.1.1}\end{equation} и сходятся
$$
x_n \underset{n\to \infty}{\longrightarrow} a \qquad y_n
\underset{n\to \infty}{\longrightarrow} b
$$
то их пределы связаны тем же неравенством
$$
a\le b
$$
\end{tm}\begin{proof} Предположим, что $b<a$. Тогда можно
взять $\varepsilon=\frac{a-b}{2}$. Из $x_n \underset{n\to
\infty}{\longrightarrow} a$ будет следовать, что почти все
числа $x_n$ лежат в окрестности $(a-\varepsilon;
a+\varepsilon)$. То есть, для почти всех номеров $n\in
\mathbb{N}$ выполняется неравенство
$$
a-\varepsilon<x_n<a+\varepsilon
$$
С другой стороны, из $y_n \underset{n\to
\infty}{\longrightarrow} a$ будет следовать, что почти все
числа $y_n$ лежат в окрестности $(b-\varepsilon;
b+\varepsilon)$. То есть, для почти всех номеров $n\in
\mathbb{N}$ выполняется неравенство
$$
b-\varepsilon<y_n<b+\varepsilon
$$
Таким образом, у нас получается, что для почти всех номеров
$n\in \mathbb{N}$ выполняются неравенства
$$
y_n<b+\varepsilon=a-\varepsilon<x_n
$$
То есть, для почти всех $n\in \mathbb{N}$
$$
y_n<x_n
$$
Это противоречит условию \eqref{2.2.1}. \end{proof}

Аналогично доказывается

\begin{tm}[\bf о предельном переходе в неравенстве с бесконечным пределом]\label{x_n<_y_n,x_n->infty} Пусть последовательности $\{ x_n \}$, $\{ y_n \}$
связаны неравенством
$$
x_n\le y_n
$$
тогда
 \bit{
 \item[---] если
$x_n \underset{n\to \infty}{\longrightarrow}+\infty,$
то
$y_n \underset{n\to \infty}{\longrightarrow}+\infty;$

 \item[---] если
$y_n \underset{n\to \infty}{\longrightarrow}-\infty,$
то
$x_n \underset{n\to \infty}{\longrightarrow}-\infty.$

 }\eit

\end{tm}

\begin{tm}[\bf о двух милиционерах]\label{milit}\index{теорема!о двух
милиционерах}\footnote{Эта теорема также используется
довольно часто, например, в  при
доказательстве теоремы Больцано-Вейерштрасса \ref{Bol-Wei},
при доказательстве непрерывности
модуля, синуса и косинуса (предложение \ref{PROP:sin-cos-nepreryvny}) и при доказательстве замечательных
пределов (теоремы \ref{sin_x/x} и \ref{II-lim}).}

Пусть последовательности $\{ x_n \}$, $\{ y_n \}$, $\{ z_n \}$ связаны
неравенствами
$$
x_n\le y_n\le z_n
$$
причем две крайние из них сходятся к одному пределу
$$
x_n \underset{n\to \infty}{\longrightarrow} C \qquad z_n
\underset{n\to \infty}{\longrightarrow} C
$$
Тогда средняя последовательность тоже стремится к этому пределу:
$$
y_n \underset{n\to \infty}{\longrightarrow} C
$$
\end{tm}\begin{proof} Возьмем произвольную окрестность
$(C-\varepsilon;C+\varepsilon)$ точки $C$. Тогда

1) поскольку $x_n \underset{n\to \infty}{\longrightarrow}
C$, почти все числа $x_n$ содержатся в окрестности
$(C-\varepsilon;C+\varepsilon)$; поэтому для почти всех
номеров $n\in \mathbb{N}$ выполняется неравенство
$$
C-\varepsilon<x_n
$$
то есть для почти всех номеров $n\in \mathbb{N}$
выполняется
$$
C-\varepsilon<x_n\le y_n
$$

2) поскольку $z_n \underset{n\to \infty}{\longrightarrow}
C$, почти все числа $z_n$ содержатся в окрестности
$(C-\varepsilon;C+\varepsilon)$; поэтому для почти всех
номеров $n\in \mathbb{N}$ выполняется неравенство
$$
z_n<C+\varepsilon
$$
то есть для почти всех номеров $n\in \mathbb{N}$
выполняется
$$
y_n\le z_n<C+\varepsilon
$$

Итак, мы получили, что для почти всех номеров $n\in
\mathbb{N}$ выполняются неравенства
$$
C-\varepsilon< y_n \qquad y_n<C+\varepsilon
$$
Это означает, что почти все числа $y_n$ содержатся в
окрестности $(C-\varepsilon;C+\varepsilon)$ точки $C$.
Поскольку это верно для любой окрестности, мы получаем, что
$y_n \underset{n\to \infty}{\longrightarrow} C$,
\end{proof}

\noindent\rule{160mm}{0.1pt}\begin{multicols}{2}

\bex Докажем формулу:
 \beq\label{n^k/a^n->?}
\frac{n^k}{a^n}\underset{n\to\infty}{\longrightarrow}\begin{cases}0,& |a|>1\\
\infty,& |a|<1\end{cases}\qquad (k\in\Z)
 \eeq
\eex \bpr Очевидно, здесь достаточно рассмотреть случай, когда $a>0$.

1. Если $k<0$ и $a>1$, то все очевидно:
$$
\frac{n^k}{a^n}=\underbrace{\frac{1}{n^{-k}}}_{\scriptsize\begin{matrix}\downarrow\\
0\end{matrix}}\cdot\underbrace{\frac{1}{a^n}}_{\scriptsize\begin{matrix}\downarrow\\
0\end{matrix}}\underset{n\to\infty}{\longrightarrow}0
$$

2. Пусть $k\ge 0$ и $a>1$. Тогда $a=1+\e$, $\e>0$, и при $n\ge k+1$ мы получаем
цепочку:
$$
a^n=(1+\e)^n=\sum_{i=0}^nC_n^i\e^i\ge C_n^{k+1}\e^{k+1}
$$
$$
\Downarrow
$$
$$
0\le\frac{n^k}{a^n}\le \frac{n^k}{C_n^{k+1}\e^{k+1}}=
\underbrace{\frac{1}{\e^{k+1}}}_{\text{константа}}\cdot\kern-15pt\underbrace{\frac{n^k}{C_n^{k+1}}}_{\scriptsize\begin{matrix}
\phantom{ \eqref{C_n^k/n^l->?}}\ \downarrow\ \eqref{C_n^k/n^l->?}\\
0\end{matrix}}\kern-15pt\underset{n\to\infty}{\longrightarrow}0
$$
$$
\phantom{\text{\scriptsize теорема \ref{milit}}}\ \Downarrow\ \text{\scriptsize
теорема \ref{milit}}
$$
$$
\frac{n^k}{a^n}\underset{n\to\infty}{\longrightarrow}0
$$

3. Если $k\le 0$ и $0<a<1$, то, по уже доказанному,
$$
\frac{n^{-k}}{(a^{-1})^n}\underset{n\to\infty}{\longrightarrow}0
$$
$$
\Downarrow
$$
$$
\frac{n^k}{a^n}=\l\frac{n^{-k}}{(a^{-1})^n}\r^{-1}\underset{n\to\infty}{\longrightarrow}\infty
$$

4. Наконец, при $k>0$ и $0<a<1$ все опять очевидно:
$$
\frac{n^k}{a^n}=\underbrace{n^k}_{\scriptsize\begin{matrix}\downarrow\\
\infty\end{matrix}}\cdot\underbrace{\l\frac{1}{a}\r^n}_{\scriptsize\begin{matrix}\downarrow\\
\infty\end{matrix}}\underset{n\to\infty}{\longrightarrow}\infty
$$
 \epr

\end{multicols}\noindent\rule[10pt]{160mm}{0.1pt}

\subsection{Монотонные последовательности и
теорема Вейерштрасса}

Последовательность  $\{ x_n \}$ называется
 \bit{
\item[--] {\it
возрастающей}\index{последовательность!возрастающая}, если
$$
x_1< x_2< x_3< ...< x_n < x_{n+1}< ...
$$

%\picture{0pt}{0pt}{76.pcx}

\vglue100pt \noindent

\item[--] {\it
убывающей}\index{последовательность!убывающая}, если
$$
x_1> x_2> x_3> ...> x_n > x_{n+1}> ...
$$

%\picture{0pt}{0pt}{76.pcx}

\vglue100pt \noindent

\item[--] {\it неубывающей}\index{последовательность!неубывающая}, если
$$
x_1\le x_2\le x_3\le ...\le  x_n\le x_{n+1}\le ...
$$

%\picture{0pt}{0pt}{76.pcx}

\vglue100pt \noindent

\item[--] {\it невозрастающей}\index{последовательность!невозрастающая}, если
$$
x_1\ge x_2\ge x_3\ge ...\ge  x_n\ge x_{n+1}\ge ...
$$

%\picture{0pt}{0pt}{77.pcx}

\vglue100pt \noindent

\item[--] {\it монотонной}\index{последовательность!монотонная}, если она
неубывающая или невозрастающая;
 }\eit

 \bigskip

\centerline{\bf Свойства монотонных последвательностей}

 \bit{\it
\item[$1^\circ$]\label{vozrastanie-posl=>neubyvanie-posl} Если
последовательность $x_n$ возрастает, то она неубывает.

\item[$2^\circ$] Если последовательность $x_n$ убывает, то она невозрастает.
 }\eit

\noindent\rule{160mm}{0.1pt}\begin{multicols}{2}

\begin{ex}
Последовательность $x_n=\frac{1}{n}$ монотонна (точнее, убывает) и ограничена;

%\picture{0pt}{0pt}{78.pcx}

\vglue100pt

\end{ex}

\begin{ex} Последовательность $x_n=n$ монотонна
(точнее, возрастает) и неограничена;

%\picture{0pt}{0pt}{79.pcx}

\vglue100pt

\end{ex}

\begin{ex} Последовательность $x_n=\frac{(-1)^n}{n}$ немонотонна, но
ограничена;

%\picture{0pt}{0pt}{80.pcx}

\vglue100pt
\end{ex}

\end{multicols}\noindent\rule[10pt]{160mm}{0.1pt}

\begin{tm}[\bf Вейерштрасса о монотоных ограниченных последовательностях]
\label{Wei-I}\index{теорема!Вейерштрасса!о
последовательностях}\footnote{Этот результат используется
далее в главе \ref{ch-th-x_n} при доказательстве теоремы о
вложенных отрезках \ref{I_1>I_2>...} и в главе
\ref{ch-lim_f(x)} при доказательстве второго замечательного
предела (лемма \ref{e}}

Всякая монотонная ограниченная последовательность $\{ x_n \}$ сходится (то есть
имеет конечный предел), причем
 \bit{
\item[---] если $\{ x_n \}$ неубывающая последовательность, то
$$
  \lim_{n\to\infty} x_n=\sup_{n\in \mathbb{N}} x_n
$$
\item[---] если $\{ x_n \}$ невозрастающая последовательность, то
$$
  \lim_{n\to\infty} x_n=\inf_{n\in \mathbb{N}} x_n
$$
 }\eit
\end{tm}\begin{proof} Если последовательность $\{ x_n \}$
монотонна, то это значит, что она или неубывающая, или невозрастающая. Докажем
теорему для случая, когда $\{ x_n \}$ -- неубывающая (случай, когда она
невозрастающая рассматривается аналогично):
\begin{equation}
C\le x_1\le x_2\le x_3\le ...\le  x_n\le x_{n+1}\le ...\le
D \label{2.2.1}\end{equation} Рассмотрим множество $X$,
состоящее из элементов $\{ x_n \}$. Из цепочки неравенств
\eqref{2.2.1} следует, что оно будет ограничено сверху и
непусто. Значит, по теореме \ref{sup-inf} (о точной
границе), оно имеет точную верхнюю грань $B$:
$$
B=\sup \{ x_n ; n\in \mathbb{N}\}
$$
Покажем, что $\{ x_n \}$ стремится к $B$. Возьмем
окрестность $(B-\varepsilon ,B+\varepsilon)$ точки $B$.
Поскольку $B$ -- точная верхняя грань для $\{ x_n \}$,
число $B-\varepsilon < B$ не может быть верхней гранью для
$\{ x_n \}$. Это значит, что существует такое $k$, что
$$
B-\varepsilon< x_k
$$
При этом для любого $n\ge k$ мы получим
$$
B-\varepsilon < x_k \le x_n
$$
то есть все числа  $\{ x_n \}$, начиная с номера $k$, лежат
в интервале $(B-\varepsilon ,B+\varepsilon)$.

%\picture{0pt}{0pt}{81.pcx}

\vglue100pt \noindent Итак, получается, что любая
окрестность $(B-\varepsilon ,B+\varepsilon)$ точки $B$
содержит почти все точки последова\-тельности $\{ x_n \}$.
Значит
$$
\lim_{n\to \infty} x_n = B \qquad $$ \end{proof}

\noindent\rule{160mm}{0.1pt}\begin{multicols}{2}

\begin{ex}\label{x_(n+1)=1/2(x_n+1/x_n)}
Покажем, что последовательность, заданная рекуррентно
$$
x_1=2, \quad
x_{n+1}=\frac{1}{2}\left(x_n+\frac{1}{x_n}\right)
$$
сходится, и найдем ее предел.

1. Вычислим сначала несколько первых элементов последовательности $\{ x_n \}$,
и нарисуем картинку:

%\picture{0pt}{0pt}{83.pcx}

\vglue100pt

2. Заметим, что все числа $\{ x_n \}$ положительны:
$$
x_n > 0
$$
(это можно отдельно доказать индукцией).

3. Из картинки видно, что на первых нескольких элементах $\{ x_n \}$ наша
последовательность невозрастает. Попробуем доказать, что она невозрастает
всегда:
\begin{equation}
x_{n+1}\le x_n \label{2.2.6}\end{equation} Поймем сначала,
что означает это неравенство:
 \begin{gather*}
x_{n+1}\le x_n \quad \Leftrightarrow \quad
\overbrace{\frac{1}{2}\left(x_n+\frac{1}{x_n}\right)\le
x_n}^{\smsize \begin{matrix}\text{умножаем на $2x_n$;}\\
\text{поскольку $x_n>0$,}\\ \text{знак неравенства не меняется}\\
\downarrow \end{matrix}}\quad \Leftrightarrow \\
\Leftrightarrow \quad x_n^2+1\le 2 x_n^2 \quad
\Leftrightarrow \quad x_n^2\ge 1 \quad \Leftrightarrow
\quad x_n\ge 1
 \end{gather*}
 Теперь докажем, последнее неравенство:
\begin{equation}
x_n\ge 1 \label{2.2.7}\end{equation} Это делается
математической индукцией.

a) Сначала проверяем наше неравенство при $n=1$:
$$
x_1=2\ge 1
$$

b) Предполагаем, что оно верно при $n=k$:
\begin{equation}
x_k\ge 1 \label{2.2.8}\end{equation}

c) Доказываем, что тогда оно будет верно при $n=k+1$.
\begin{equation}
x_{k+1}\ge 1 \label{2.2.9}\end{equation} Действительно,
 \begin{gather*}
x_{k+1}\ge 1 \quad \Leftrightarrow \quad
 \overbrace{\frac{1}{2}\left(x_k+\frac{1}{x_k}\right)\ge 1}^
 {\smsize \begin{matrix}\text{умножаем на $2x_k$,}\\
\text{и, поскольку $x_n>0$,}\\ \text{знак неравенства не меняется}\\
\downarrow\end{matrix}}\quad \Leftrightarrow \\
\Leftrightarrow \quad x_k^2+1\ge 2 x_k \quad
\Leftrightarrow \quad x_k^2-2 x_k+1\ge 0 \quad \Leftrightarrow \\
\Leftrightarrow \quad (x_k-1)^2\ge 0
 \end{gather*}
 а последнее неравенство выполняется
всегда (потому что квадрат любого числа неотрицателен).

Итак, мы доказали \eqref{2.2.6}, а значит и равносильное
ему условие \eqref{2.2.7}. Эти неравенства означают, что
$\{ x_n \}$ невозрастает и ограничена сверху. Значит, по
теореме Вейерштрасса \ref{Wei-I} она имеет предел.
Обозначим его буквой $c$:
$$
\lim_{n\to \infty} x_n=c
$$
и заметим, что $c\ge 0$, поскольку $x_n\ge 0$. Тогда из
формулы
$$
x_{n+1}=\frac{1}{2}\left(x_n+\frac{1}{x_n}\right)
$$
получаем
$$
c=\lim_{n\to \infty} x_{n+1}=\lim_{n\to
\infty}\frac{1}{2}\left(x_n+\frac{1}{x_n}\right) =
\frac{1}{2}\left(с+\frac{1}{c}\right)
$$
То есть $c$ удовлетворяет уравнению
$$
c=\frac{1}{2}\left(с+\frac{1}{c}\right)
$$
решая которое  мы получаем
$$
2c^2=c^2+1 \quad \Leftrightarrow \quad с^2=1
$$
откуда (с учетом $c\ge 0$) имеем $c=1$.

Ответ: $\lim\limits_{n\to \infty} x_n =1$
\end{ex}

\begin{ers}
Докажите, что последовательность, заданная рекуррентно сходится, и найдите ее
предел:
 \biter{
\item[1)] $x_1=3, \, x_{n+1}=\frac{5x_n}{3x_n+1}$;

\item[2)] $x_1=1, \, x_{n+1}=\frac{4x_n}{2+x_n^2}$;

\item[3)] $x_1=-1, \, x_{n+1}=\frac{x_n(x_n+4)}{2}$;

\item[4)] $x_1=2, \, x_{n+1}=\frac{x_n(6-x_n)}{3}$.
 }\eiter
\end{ers}

\end{multicols}\noindent\rule[10pt]{160mm}{0.1pt}

\subsection{Теорема о вложенных отрезках}

Пусть дана последовательность отрезков $[a_n; b_n]$ такая, что каждый
последующий содержится в предыдущем:
$$
[a_1;b_1]\supseteq [a_2;b_2]\supseteq [a_3;b_3]\supseteq
... \supseteq [a_n;b_n]\supseteq ...
$$
Понятно, что в этом случае концы отрезков связаны
неравенствами
\begin{equation}
a_1\le a_2\le a_3 \le ... \le a_n \le ... \le b_n\le ...\le
b_3\le b_2 \le b_1 \label{2.3.1}\end{equation} Пусть, кроме
того, длины этих отрезков стремятся к нулю:
\begin{equation}
b_n-a_n \underset{n\to \infty}{\longrightarrow} 0 \label{2.3.2}\end{equation}
Такая последовательность отрезков называется {\it последовательностью вложенных
отрезков}.

\begin{tm}[\bf о вложенных отрезках]\label{I_1>I_2>...}
\index{теорема!о вложенных отрезках}\footnote{Этот результат используется в
главе \ref{ch-th-x_n} при доказательстве теоремы Больцано-Вейерштрасса
\ref{Bol-Wei} и в главе \ref{ch-th-cont-f(x)} при доказательстве теоремы Коши о
промежуточном значении \ref{Cauchy-I}.} Для любой последовательности вложенных
отрезков $[a_n; b_n]$ существует единственная точка $C$, принадлежащая всем
отрезкам этой последовательности,
\begin{equation}
C\in [a_n;b_n] \label{2.3.3}\end{equation} При этом
\begin{equation}\lim_{n\to \infty} a_n=C=\lim_{n\to \infty} b_n \label{2.3.4}\end{equation}\end{tm}\begin{proof} Неравенства \eqref{2.3.1} означают,
что последовательности $\{ a_n \}$ и $\{ b_n \}$ монотоны и ограничены. Значит,
они сходятся (то есть имеют пределы). Обозначим буквами $A$ и $B$ эти пределы:
$$
A=\lim_{n\to \infty} a_n \qquad B=\lim_{n\to \infty} b_n
$$
Покажем, что эти числа совпадают. Действительно,
$$
B-A=\lim_{n\to \infty} b_n-\lim_{n\to \infty} a_n=
{\smsize\begin{pmatrix}\text{используем свойство}\\
\text{$2^0$ из \ref{->0-&-->infty}} \end{pmatrix}}
=\lim_{n\to \infty} (b_n-a_n)=\left(\text{используем
\eqref{2.3.2}}\right)=0
$$
то есть
$$
A=B
$$
Итак, мы получили, что можно взять $C=A=B$, и тогда будет
выполняться \eqref{2.3.4}. По теореме Вейерштрасса
\ref{Wei-I}, получаем
$$
a_1\le a_2\le a_3 \le ... \le a_n \le ...\le \sup_{n\in
\mathbb{N}} a_n=A=C=B=\inf_{n\in \mathbb{N}} b_n  \le...
\le b_n\le ...\le b_3\le b_2 \le b_1
$$
и это означает, что $C$ принадлежит всем интервалам
$[a_n;b_n]$. Наконец, такая точка $C$ (обладающая свойством
\eqref{2.3.3}), единственна, потому что если бы
существовала какая-то другая точка $C'$ с тем же самым
свойством \eqref{2.3.3}, отличющаяся от  $C$, например,
большая, чем $C$, то мы получили бы, цепочку неравенств
$$
a_1\le a_2\le a_3 \le ... \le a_n \le ...\le C<C' \le...
\le b_n\le ...\le b_3\le b_2 \le b_1
$$
из которых следовало бы, что расстояние между концами
любого отрезка $[a_n;b_n]$ не меньше чем фиксированное
число $\Delta=C'-C$
$$
 b_n - a_n \ge C'-C=\Delta>0
$$
и поэтому $b_n - a_n$ не может стремиться к нулю.
\end{proof}

\subsection{Подпоследовательности и теорема
Больцано-Вейерштрасса}\label{SEC-BOL-WEI}

Пусть нам дана какая-то последовательность вещественных чисел
$$
x_n
$$
и пусть кроме того дана бесконечно большая последовательность натуральных
чисел:
$$
n_k \underset{k\to \infty}{\longrightarrow} +\infty
$$
Тогда можно рассмотреть последовательность $y_k=x_{n_k}$ (зависящую от индекса
$k$). Она называется {\it подпоследовательностью} числовой последовательности
$x_n$.

\noindent\rule{160mm}{0.1pt}\begin{multicols}{2}

\begin{ex}
Пусть, скажем, нам дана последовательность
$$
x_n=\frac{1}{n}
$$
Если взять
$$
n_k=k^2
$$
то мы получим подпоследовательность
$$
y_k=x_{n_k}=\frac{1}{n_k}=\frac{1}{k^2}
$$
\end{ex}

\begin{ex}
Пусть
$$
x_n=(-1)^n
$$
Если взять
$$
n_k=2k
$$
то мы получим подпоследовательность
$$
y_k=x_{n_k}=(-1)^{n_k}=(-1)^{2k}=1
$$
\end{ex}

\begin{ers}
Какую нужно взять последовательность индексов $\{ n_k \}$, чтобы
последовательность $\{ y_k \}$ была подпоследовательностью последовательности
$\{ x_n \}$?
 \biter{
\item[1)] $x_n=\frac{1}{n}, \, y_k=\frac{1}{k^3}$

\item[2)] $x_n=(-1)^n, \, y_k=-1$
 }\eiter
\end{ers}

\begin{ers}
Является ли данная последовательность $\{ y_k \}$ подпоследовательностью
последовательности $\{ x_n \}$? (Если да, то указать соответствующую
последовательность индексов $\{ n_k \}$.)
 \biter{
\item[1)] $x_n=n, \, y_k=k^2$

\item[2)] $x_n=n^2, \, y_k=k$

\item[3)] $x_n=n, \, y_k=2$
 }\eiter
\end{ers}

\end{multicols}\noindent\rule[10pt]{160mm}{0.1pt}

\bigskip

\centerline{\bf Свойства подпоследовательностей}

 \bit{\it
 \label{podposledovatelnosti}
\item[$1^\circ$.] Если последовательность $\{ x_n \}$ стремится к числу $a$, то
любая ее подпоследовательность $\{ x_{n_k}\}$ тоже стремится к числу $a$:
$$
 x_n\underset{n\to \infty}{\longrightarrow} a
 \quad \Longrightarrow\quad
  x_{n_k}\underset{k\to \infty}{\longrightarrow} a
$$
\item[$2^\circ$.] Если последовательность $\{ x_n \}$ не стремится к числу $a$,
то найдется такая подпоследовательность $\{ x_{n_k}\}$ и такое число
$\varepsilon>0$, что $\{ x_{n_k}\}$ отстоит от $a$ больше чем на $\e$:
$$
x_n\underset{n\to \infty}{\notarrow} a
 \quad\Longrightarrow\quad
  \exists\e>0 \quad \exists x_{n_k}\quad
\forall k\quad |x_{n_k}-a|>\varepsilon
$$
\item[$3^\circ$.] Если последовательность $\{ x_n \}$ неограничена, то она
содержит некоторую бесконечно большую подпоследовательность:
$$
\sup_{n\in \mathbb{N}} |x_n|=+\infty
 \quad\Longrightarrow\quad
 \exists x_{n_k}\underset{k\to \infty}{\longrightarrow}\infty
$$

\item[$4^\circ$.] Если бесконечный набор элементов последовательности $\{x_n\}$
содержится в множестве $M$, то $\{x_n\}$ содержит некоторую
подпоследовательность, полностью содержащуюся в $M$:
$$
 \exists x_{n_k}\subseteq M
$$
 }\eit

\begin{proof} 1. Докажем сначала $1^\circ$. Возьмем
произвольную окрестность $(a-\varepsilon;a+\varepsilon)$
точки $a$. Поскольку $x_n\underset{n\to
\infty}{\longrightarrow} a$, почти все точки $\{ x_n \}$
должны лежать в $(a-\varepsilon;a+\varepsilon)$. Значит,
точки с индексами $n_k$ тоже почти все лежат в
$(a-\varepsilon;a+\varepsilon)$. Поскольку это верно для
любой окрестности $(a-\varepsilon;a+\varepsilon)$ мы
получаем $x_{n_k}\underset{n\to \infty}{\longrightarrow}
a$.

2. Следуя обещанию на с.\pageref{beskonechnokratnyj-vybor}, мы дадим два доказательства свойства $2^0$: традиционное, использующее прием бесконечнократного выбора, и аккуратное, со ссылкой на аксиому (в данном случае счетного) выбора.

\biter{

\item[а)] {\it Доказательство приемом бесконечнократного выбора}. Если $x_n\underset{n\to \infty}{\notarrow} a$, то это
означает, что найдется такая окрестность $(a-\e, a+\e)$ точки $a$, вне
которой лежит бесконечное число элементов $x_n$. То есть множество
$$
N=\{n\in\N:\quad x_n\notin (a-\e, a+\e)\}
$$
бесконечно.

1) Выберем произвольное $n_1\in N$.

2) Поскольку $N$ бесконечно, найдется $n_2\in N$ такое что $n_2\ge 2$. Зафиксируем это $n_2$.

3) Поскольку $N$ бесконечно, найдется $n_3\in N$ такое что $n_3\ge 3$. Зафиксируем это $n_3$.

И так далее. Поступая подобным образом, мы получим последовательность индексов $\{n_k\}$ такую, что
$$
n_k\ge k,\qquad k\in\N
$$
поэтому по теореме \ref{x_n<_y_n,x_n->infty} о предельном переходе в неравенстве, $n_k\underset{k\to \infty}{\longrightarrow}+\infty$, то есть $\{x_{n_k}\}$ -- подпоследовательность в $\{x_n\}$.

\item[б)] {\it Доказательство с явной ссылкой на аксиому выбора}.
Еще раз, если $x_n\underset{n\to \infty}{\notarrow} a$, то найдется такая окрестность $(a-\e, a+\e)$ точки $a$, вне
которой лежит бесконечное число элементов $x_n$. То есть множество
$$
N=\{n\in\N:\quad x_n\notin (a-\e, a+\e)\}
$$
бесконечно. По теореме \ref{TH:ogran-mnozh-v-N-konechno}, оно не ограничено, то есть для всякого $k\in\N$
$$
N_k=\{n\in\N:\quad n\le k\}\ne \varnothing
$$
По аксиоме счетного выбора (с.\pageref{AX:schentnogo-vybora}), найдется последовательность
$$
n_k\in N_k
$$
Для нее мы получаем во-первых:
$$
x_{n_k}\notin (a-\e, a+\e), \qquad k\in\N
$$
то есть
$$
|x_{n_k}-a|> \varepsilon, \qquad k\in\N
$$
И, во-вторых,
$$
n_k\ge k,\qquad k\in\N,
$$
откуда по теореме \ref{x_n<_y_n,x_n->infty} о предельном переходе в неравенстве, $n_k\underset{k\to \infty}{\longrightarrow}+\infty$, то есть $\{x_{n_k}\}$
-- подпоследовательность в $\{x_n\}$.
 }\eiter

3.\label{PROOF:3^0-podposled} Как и в предыдущем случае, мы даем два доказательства свойства $3^0$: традиционное и с явной ссылкой на аксиому (счетного) выбора.

\biter{

\item[а)] {\it Доказательство приемом бесконечнократного выбора}. Пусть $\sup_{n\in \mathbb{N}}
|x_n|=+\infty$, построим подпоследовательность $x_{n_k}\underset{k\to\infty}{\longrightarrow}\infty$.

1) Сначала выбирается индекс $n_1$. Это делается так:
поскольку $\{ x_n \}$ неограничена, найдется такой номер
$n_1$, что
$$
 \left| x_{n_1}\right|\ge 1
$$
Этот номер $n_1$ фиксируется.

2) Затем выбирается индекс $n_2$: поскольку $\{ x_n \}$ неограничена, найдется такой
номер $n_2$, что
$$
 \left| x_{n_2}\right|\ge 2
$$
Этот номер $n_2$ фиксируется.

$k$). Когда индексы $n_1,n_2,...,n_{k-1}$ уже выбраны,
выбирается индекс $n_k$: поскольку $\{
x_n \}$ неограничена, найдется такой номер $n_k$, что
$$
 \left| x_{n_k}\right|\ge k
$$
Этот номер $n_k$ фиксируется.

И так далее. В результате получается последовательность $\{ x_{n_k}\}$ со
свойством
$$
\forall k\in \mathbb{N}\qquad \left| x_{n_k}\right|\ge k
$$
из которого и следует $x_{n_k}\underset{k\to\infty}{\longrightarrow}\infty$.

\item[б)] {\it Доказательство с явной ссылкой на аксиому выбора}.
Пусть $\sup_{n\in \mathbb{N}}
|x_n|=+\infty$. Это значит, что
$$
\forall k\in\N\qquad \{n\in\N:\quad |x_n|\ge k\}\ne\varnothing
$$
По аксиоме счетного выбора (с.\pageref{AX:schentnogo-vybora}), найдется последовательность $\{n_k\}$ такая, что
$$
n_k\in\{n\in\N:\quad |x_n|\ge k\},\qquad k\in\N
$$
то есть
$$
|x_{n_k}|\ge k,\qquad k\in\N.
$$
Отсюда следует $x_{n_k}\underset{k\to\infty}{\longrightarrow}\infty$.
}\eiter

4. Свойство $4^\circ$ доказывается тем же приемом, только номера $n_k$ здесь нужно
выбирать так, чтобы выполнялись не одно, а два условия:
$$
n_{k+1}>n_k \qquad \& \qquad x_{n_{k+1}}\in M.
$$
\end{proof}

Можно заметить, что иногда получается, что, хотя данная последовательность $\{
x_n \}$ расходится (то есть не имеет предела), тем не менее, выбранная
подпоследовательность $\{ y_k \}$ сходится.

\noindent\rule{160mm}{0.1pt}\begin{multicols}{2}

\begin{ers}
Для данной последовательности $\{ x_n \}$ укажите какую-нибудь сходящуюся
подпоследовательность $\{ x_{n_k}\}$, если она существует.
 \biter{
\item[1)] $x_n=(-1)^n\cdot n$

\item[2)] $x_n=\frac{(-1)^n}{n}+\frac{1+(-1)^n}{2}$

\item[3)] $x_n=n^{(-1)^n}$

\item[4)] $x_n=n\cdot (1+(-1)^n)$

\item[5)] $x_n=\frac{2+(-1)^n}{2}-\frac{1}{n}$

\item[6)] $x_n=(-1)^{n-1}\cdot \l 2+\frac{3}{n}\r$
 }\eiter
\end{ers}

\end{multicols}\noindent\rule[10pt]{160mm}{0.1pt}

\begin{tm}[\bf Больцано-Вейерштрасса]\label{Bol-Wei}
\index{теорема!Больцано-Вейерштрасса}\footnote{Этот результат
используется в главе \ref{ch-th-x_n} при доказательстве критерия Коши
сходимости числовой последовательности (теорема \ref{Cauchy-crit-seq}), а также
в главе \ref{ch-cont-f(x)} при доказательстве теорем Вейерштрасса об
ограниченности и об экстремумах (теоремы \ref{Wei-II} и \ref{Wei-III}) и
теоремы Кантора о равномерной ограниченности (теорема \ref{Kantor}) } Всякая
ограниченная последовательность $\{ x_n \}$ имеет сходящуюся
подпоследовательность $\{ x_{n_k}\}$.
\end{tm}

Для доказательства условимся {\it длиной отрезка} $I=[a,b]$ называть число
$$
\diam I=b-a
$$
Нам удобно начать со следующей леммы:

\blm\label{LM:Bol-Wei}
Если отрезок $I$ содержит бесконечно много элементов последовательности $\{x_n\}$,
$$
\card\{n\in\N:\ x_n\in I\}=\infty
$$
то в $I$ содержится отрезок $J$ в два раза меньшей длины,
$$
\diam J=\frac{\diam I}{2},
$$
также содержащий бесконечно много элементов $\{x_n\}$:
$$
\card\{n\in\N:\ x_n\in J\}=\infty
$$
\elm
\bpr
Это очевидно: разделим отрезок $I$ на два равных (по длине) отрезка $J$ и $K$. По крайней мере в одном из них должно содержаться бесконечное число элементов  $\{x_n\}$ (потому что иначе получилось бы, что $J$ и $K$ содержат конечный набор элементов $\{x_n\}$, и значит их объединение $I=J\cup K$ тоже содержит конечный набор элементов.
\epr

\begin{proof}[Доказательство теоремы \ref{Bol-Wei}.] Мы обещали на с.\pageref{beskonechnokratnyj-vybor} дать два доказательства этого утверждения: традиционное, приемом бесконечнократного выбора, и современное, в котором аксиома выбора явно упоминается.

\biter{

\item[а)] {\it Доказательство приемом бесконечнократного выбора}.

1) Последовательность $\{ x_n \}$ ограничена, значит она лежит в некотором отрезке $I_1=[a_1;b_1]$. Положим $n_1=1$, и заметим на будущее, что
$x_{n_1}\in I_1=[a_1;b_1]$.

2) По лемме \ref{LM:Bol-Wei}, в $I_1$ содержится некий отрезок $I_2=[a_2;b_2]$ вдвое меньшей длины,
$$
\diam I_2=\frac{\diam I_1}{2},
$$
также содержащий бесконечно много чисел $\{ x_n \}$. Выбираем какой-нибудь номер $n_2\ge 2$, такой что $x_{n_2}\in I_2=[a_2;b_2]$ (такой номер $n_2$ найдется, потому что $x_n\in [a_2;b_2]$ для бесконечного количества номеров $n$).

3) Опять по лемме \ref{LM:Bol-Wei}, в $I_2$ содержится некий отрезок $I_3=[a_3;b_3]$ вдвое меньшей длины,
$$
\diam I_3=\frac{\diam I_2}{2},
$$
также содержащий бесконечно много чисел $\{ x_n \}$. Выбираем какой-нибудь номер $n_3\ge 3$, такой что $x_{n_3}\in I_3=[a_3;b_3]$ (такой номер $n_3$ найдется, потому что $x_n\in [a_3;b_3]$ для бесконечного количества номеров $n$).

Продолжая эту процедуру, мы получим последовательность отрезков $I_k=[a_k;b_k]$ и чисел $x_{n_k}$ со следующими свойствами:
\begin{equation}\label{2.4.1}
[a_1;b_1]\supseteq [a_2;b_2]\supseteq [a_3;b_3]\supseteq ... \supseteq
[a_k;b_k]\supseteq ... \qquad b_k-a_k \underset{k\to \infty}{\longrightarrow} 0
\qquad x_{n_k}\in [a_k;b_k]
\end{equation}
Первые два свойства
означают, что $[a_k;b_k]$ -- последовательность вложенных отрезков, значит, по
теореме \ref{I_1>I_2>...} о вложенных отрезках, существует число $C$ такое что
$$
\lim_{k\to \infty} a_k=C=\lim_{k\to \infty} b_k
$$
Последнее же свойство в цепочке \eqref{2.4.1} означает, что
$$
a_k\le x_{n_k}\le b_k
$$
и поэтому мы можем применить теорему о двух милиционерах
\ref{milit}, из которой следует, что
$$
\lim_{k\to \infty} x_{n_k}=C
$$
То есть последовательность $x_{n_k}$ имеет предел.

\item[б)] {\it Доказательство с явной ссылкой на аксиому выбора}. Опять заметим сначала, что поскольку последовательность $\{ x_n \}$ ограничена, она лежит в некотором отрезке $I=[a;b]$. Пусть ${\mathcal I}$ обозначает множество отрезков $J$, содержащихся в $I$, и содержащих бесконечный набор элементов последовательности $\{x_n\}$
    $$
    J\in {\mathcal I}\qquad\Longleftrightarrow\qquad J\subseteq I\quad \&\quad \card\{n\in\N:\ x_n\in J\}=\infty
    $$
    (множество ${\mathcal I}$ непусто, потому что, например, $I\in{\mathcal I}$). Пусть теперь  $I_1,I_2,...,I_k$ -- конечная последовательность отрезков, удовлетворяющая таким условиям:
    \beq\label{Bol-Wei-1}
\{I_1,...I_k\}\subseteq {\mathcal I},\qquad I_1\supset I_2\supset ...\supset I_k,\qquad \forall i<k \quad \diam I_{i+1}=\frac{\diam I_i}{2},
    \eeq
    Тогда по лемме \ref{LM:Bol-Wei}, найдется отрезок $J$ со свойствами:
    \beq\label{Bol-Wei-2}
J\in{\mathcal I},\qquad    J\subset I_k,\qquad \diam J=\frac{\diam I_k}{2}.
    \eeq
    Это наблюдение можно переформулировать таким хитрым образом: для любой последовательности $I_1,I_2,...,I_k$, удовлетворяющей условиям \eqref{Bol-Wei-1}, множество ${\mathcal J}(I_1,I_2,...,I_k)$ всех отрезков $J$, удовлетворяющих условиям \eqref{Bol-Wei-2}, непусто:
    $$
    {\mathcal J}(I_1,I_2,...,I_k)=\left\{J\in{\mathcal I}:\quad  J\subset I_k,\quad \diam J=\frac{\diam I_k}{2}\right\}\ne\varnothing
    $$
    Теперь мы применяем аксиому выбора: должно существовать отображение $(I_1,I_2,...,I_k)\mapsto G(I_1,I_2,...,I_k)$, которое каждой последовательности $(I_1,I_2,...,I_k)$ отрезков, удовлетворяющих условиям \eqref{Bol-Wei-1}, ставит в соответствие отрезок $G(I_1,I_2,...,I_k)$, принадлежащий ${\mathcal J}(I_1,I_2,...,I_k)$, то есть удовлетворяющий условиям
    \beq\label{Bol-Wei-3}
G(I_1,I_2,...,I_k)\in{\mathcal I},\qquad   G(I_1,I_2,...,I_k)\subset I_k,\qquad \diam G(I_1,I_2,...,I_k)=\frac{\diam I_k}{2}.
    \eeq
После этого мы применяем теорему \ref{defin-induction-chast} об определениях частной индукцией: если положить $I_1=I$, то это определит последовательность $\{I_k\}\in {\mathcal I}$, содержащую $I_1$ в качестве первого элемента и
удовлетворяющую условию
 $$
G(I_1,...,I_k)=I_{k+1},\qquad k<N
 $$
где величина $N$  определяется условием:
 $$
N=\sup\Big\{k\in\N:\quad \exists \{I_2,...,I_k\}\subseteq {\mathcal I}\quad \forall
i\in\{2,...,k\} \qquad G(I_1,...,I_{i-1})=I_i\Big\}
 $$
Но здесь верхняя грань равна бесконечности, потому что наше отображение $G$ было определено на всех конечных последовательностях со свойствами \eqref{Bol-Wei-1}. Поэтому последовательность $\{I_k\}$, определяемая теоремой \ref{defin-induction-chast} должна быть бесконечной.

Теперь каждый отрезок $\{I_k\}$ содержит бесконечный набор элементов $\{x_n\}$, поэтому
$$
\{n\ge k:\ x_n\in I_k\}\ne\varnothing
$$
Отсюда по аксиоме счетного выбора (с.\pageref{AX:schentnogo-vybora}) следует, что существует последовательность номеров
$$
n_k\in \{n\ge k:\ x_n\in I_k\}
$$
то есть такая, что
$$
n_k\ge k\quad\&\quad x_{n_k}\in I_k
$$
Обозначив теперь через $a_k$ и $b_k$ концы отрезков $I_k]$, мы получим в точности условия \eqref{2.4.1}, и дальше рассуждения повторяют пункт а).
}\eiter

\end{proof}

\subsection{Критерий Коши сходимости
последовательности}

\begin{tm}[\bf критерий Коши сходимости числовой последовательности]
\label{Cauchy-crit-seq}\index{критерий!Коши!сходимости
последовательности}\footnote{Этот результат понадобится нам только в главе 18
при доказательстве критерия Коши сходимости числового ряда (теорема 4.1)}

Числовая последовательность $\{ x_n \}$ тогда и только тогда сходится, когда
она удовлетворяет следующим двум эквивалентным условиям:
 \bit{
\item[(i)] для любых двух бесконечно больших
последовательностей индексов $\{ p_i \},\, \{ q_i
\}\subseteq \mathbb{N}$
\begin{equation}
p_i \underset{i\to \infty}{\longrightarrow}\infty, \quad q_i \underset{i\to
\infty}{\longrightarrow}\infty \label{2.5.1}\end{equation} соответствующие
подпоследовательности $\{ x_{p_i}\}$ и $\{ x_{q_i}\}$ стремятся друг к другу:
\begin{equation}
x_{p_i}-x_{q_i}\underset{i\to \infty}{\longrightarrow} 0
\label{2.5.2}\end{equation}

\item[(ii)] для любой последовательности $\{ l_i
\}\subseteq \mathbb{N}$ и любой бесконечно большой последовательности $\{ k_i
\}\subseteq \mathbb{N}$
\begin{equation}
   k_i \underset{i\to \infty}{\longrightarrow}\infty
\label{2.5.3}\end{equation} выполняется соотношение
\begin{equation}
x_{k_i+l_i}-x_{k_i}\underset{i\to \infty}{\longrightarrow}
0 \label{2.5.4}\end{equation}
 }\eit
\end{tm}

Для доказательства нам понадобится две леммы.

\begin{lm}\label{lm-1-for-Cauchy}
Если последовательность $\{ x_n \}$ бесконечно большая
$$
x_n\underset{n\to \infty}{\longrightarrow}\infty
$$
то она содержит некоторую подпоследовательность $\{ x_{n_k}\}$ со свойством
\begin{equation}
x_{n_k}-x_{n_{k-1}}\underset{k\to \infty}{\longrightarrow}\infty
\label{2.5.5}\end{equation}\end{lm}

\begin{proof} Подпоследовательность $\{x_{n_k}\}$ строится бесконечнократным выбором (традиционно этот прием называется ``построением по индукции'').

1. Индекс $n_1$ выбирается произвольным, например, $n_1=1$.

2. Затем выбирается индекс $n_2$. Для этого произносится
фраза: поскольку $\{ x_n \}$ неограничена, найдется такой
номер $n_2$, что
$$
 \left| x_{n_2}- x_{n_1}\right|\ge 2
$$
Этот номер $n_2$ фиксируется.

$k$. Когда индексы $n_1,n_2,...,n_{k-1}$ уже выбраны,
выбирается индекс $n_k$. Опять говорится так: поскольку $\{
x_n \}$ неограничена, найдется такой номер $n_k$, что
$$
 \left| x_{n_k} -x_{n_{k-1}}\right|\ge k
$$
Этот номер $n_k$ фиксируется.

И так далее. В результате получается последовательность $\{ x_{n_k}\}$ со
свойством
$$
\forall k\in \mathbb{N}\qquad \left|
x_{n_k}-x_{n_{k-1}}\right|\ge k
$$
из которого и следует \eqref{2.5.5}. \end{proof}

\begin{lm}\label{lm-2-for-Cauchy}
Если последовательность $\{ x_n \}$ обладает свойством $(i)$ теоремы
\ref{Cauchy-crit-seq}, то она ограничена.
\end{lm}\begin{proof} Предположим, что это не так, то
есть что $\{ x_n \}$ -- неограниченная последовательность, обладающая свойством
$(i)$. Тогда по свойству $3^0$ предыдущего параграфа, нашлась бы
подпоследовательность $\{ x_{n_k}\}$ такая что
$$
x_{n_k}\underset{k\to \infty}{\longrightarrow}\infty
$$
а по лемме \ref{lm-1-for-Cauchy} это означало бы, что можно выбрать
подпоследовательность $\{ x_{n_{k_i}}\}$ последовательности $\{ x_{n_k}\}$
такую что
$$
x_{n_{k_i}}-x_{n_{k_{i-1}}}\underset{k\to
\infty}{\longrightarrow}\infty
$$
После этого мы бы могли взять две последовательности
$$
  p_i=n_{k_i}, \quad q_i=n_{k_{i-1}}
$$
и у нас бы получилось, что
$$
x_{p_i}-x_{q_i}= x_{n_{k_i}}-x_{n_{k_{i-1}}}\underset{k\to
\infty}{\longrightarrow}\infty
$$
Но это невозможно, потому что $\{ x_n \}$ обладает свойством $(i)$. Значит,
наше исходное предположение (о том, что $\{ x_n \}$ -- неограниченная
последовательность) неверно. То есть $\{ x_n \}$ обязательно ограничена.
\end{proof}

\begin{proof}[Доказательство теоремы \ref{Cauchy-crit-seq}]
Убедимся сначала что условия $(i)$ и $(ii)$ действительно
эквивалентны.

1. Импликация $(i)\Rightarrow (ii)$ очевидна: если выполняется $(i)$ то есть
любые две подпоследовательности $\{ x_{p_i}\}$ и $\{ x_{q_i}\}$
последовательности $\{ x_n \}$ стремятся друг к другу
$$
x_{p_i}-x_{q_i}\underset{i\to \infty}{\longrightarrow} 0,
$$
то автоматически выполняется и $(ii)$, потому что для любой последовательности
$\{ l_i \}\subseteq \mathbb{N}$ и любой бесконечно большой последовательности
$\{ k_i \}\subseteq \mathbb{N}$ мы можем положить $p_i=k_i+l_i$ и $q_i=k_i$, и
тогда будет выполняться соотношение
$$
x_{k_i+l_i}-x_{k_i}=x_{p_i}-x_{q_i}\underset{i\to
\infty}{\longrightarrow} 0
$$

2. Докажем импликацию $(ii)\Rightarrow (i)$. Пусть выполняется $(ii)$, то есть
для любой последовательности $\{ l_i \}\subseteq \mathbb{N}$ и любой бесконечно
большой последовательности $\{ k_i \}\subseteq \mathbb{N}$ выполняется
соотношение \eqref{2.5.4}. Возьмем какие-нибудь две бесконечно большие
последовательности индексов
$$
p_i \underset{i\to \infty}{\longrightarrow}\infty, \quad
q_i \underset{i\to \infty}{\longrightarrow}\infty
$$
Положим
$$
k_i=\min\{ p_i; q_i \}, \quad l_i=\max\{ p_i; q_i \}-k_i
$$
Тогда
$$
k_i+l_i=\max\{ p_i; q_i \}
$$
и поэтому
$$
x_{p_i}-x_{q_i}= \lll
\begin{array}{c}
{ x_{\max\{ p_i; q_i \}}-x_{\min\{ p_i; q_i \}}, \quad
\text{если}\,\, p_i>q_i }\\{ x_{\min\{ p_i; q_i
\}}-x_{\max\{ p_i; q_i \}}, \quad \text{если}\,\, p_i<q_i
}\\
{ 0, \quad \text{если}\,\, p_i=q_i }\end{array}\rrr= \lll
\begin{array}{c}{
x_{ k_i+l_i}-x_{k_i}, \quad \text{если}\,\, p_i>q_i }\\{
x_{k_i}-x_{k_i+l_i}, \quad \text{если}\,\, p_i<q_i
}\\
{ 0, \quad \text{если}\,\, p_i=q_i }\end{array}\rrr
$$
В любом случае получается
$$
\left|x_{p_i}-x_{q_i}\right|=\left|x_{k_i+l_i}-x_{k_i}\right|
$$
поэтому
$$
\left|x_{p_i}-x_{q_i}\right|=\left|x_{k_i+l_i}-x_{k_i}\right|
\underset{i\to \infty}{\longrightarrow} 0
$$
и значит,
$$
x_{p_i}-x_{q_i}\underset{i\to \infty}{\longrightarrow} 0
$$
Последовательности $\{ p_i \}$ и $\{ q_i \}$ здесь
выбирались с самого начала произвольными. Это значит, что
выполняется $(i)$. Таким образом, мы доказали, что из
$(ii)$ следует $(i)$.

3. Докажем теперь, что если последовательность $\{ x_n \}$ сходится, то есть
существует конечный предел
$$
  \lim_{n\to \infty} x_n=a
$$
то $\{ x_n \}$ обязательно должна обладать свойством $(i)$. Возьмем любые две
бесконечно большие последовательности индексов $\{ p_i \},\, \{ q_i \}\subseteq
\mathbb{N}$
$$
p_i \underset{i\to \infty}{\longrightarrow}\infty, \quad
q_i \underset{i\to \infty}{\longrightarrow}\infty
$$
Соответствующие подпоследовательности $\{ x_{p_i}\}$ и $\{ x_{q_i}\}$ тоже
стремятся к числу $a$ по свойству $1^0$ предыдущего параграфа:
$$
x_{p_i}\underset{i\to \infty}{\longrightarrow} a, \quad
x_{q_i}\underset{i\to \infty}{\longrightarrow} a
$$
поэтому их разность должна стремиться к нулю:
$$
x_{p_i}- x_{q_i}\underset{i\to \infty}{\longrightarrow}
a-a=0
$$
Это значит, что $\{ x_n \}$ действительно должна обладать
свойством $(i)$.

4. Теперь докажем, что наоборот, если $\{ x_n \}$ обладает свойством $(i)$, то
она обязательно сходится. Действительно, по лемме \ref{lm-2-for-Cauchy}, $\{
x_n \}$ ограничена. Значит, по теореме Больцано-Вейерштрасса \ref{Bol-Wei}, $\{
x_n \}$ содержит некоторую сходящуюся подпоследовательность $\{ x_{n_i}\}$:
\begin{equation}
x_{n_i}\underset{i\to \infty}{\longrightarrow} a \label{2.5.6}\end{equation}
Теперь возьмем следующие две последовательности индексов
$$
p_i=i, \quad q_i=n_i
$$
Поскольку $\{ x_n \}$ обладает свойством $(i)$, должно
выполняться соотношение
\begin{equation}
x_i-x_{n_i}=x_{p_i}-x_{q_i}\underset{i\to
\infty}{\longrightarrow} 0 \label{2.5.7}\end{equation} Из
\eqref{2.5.6} и \eqref{2.5.7} получаем
$$
x_i=\l x_i-x_{n_i}\r +x_{n_i}\underset{i\to
\infty}{\longrightarrow} 0+a=a
$$
То есть последовательность  $\{ x_n \}$ стремится к некоторому конечному
пределу $a$. \end{proof}

\noindent\rule{160mm}{0.1pt}\begin{multicols}{2}

\begin{ers}
Проверьте с помощью критерия Коши сходимость
последовательностей:
 \biter{
\item[1)] $x_n=3^{(-1)^n}$;
 }\eiter
\end{ers}

\end{multicols}\noindent\rule[10pt]{160mm}{0.1pt}

\subsection{Теорема Штольца}

\begin{tm}[\bf Штольц]\label{Stoltz}\index{теорема!Штольца}\footnote{Эта теорема используется при
доказательстве правила Лопиталя \ref{Lopital_x->a_infty/infty} для раскрытия
неопределенностей типа $\frac{\infty}{\infty}$.} Пусть $\{ y_n \}$ и $\{ z_n
\}$ --  две последовательности, причем
 \bit{
\item[(i)]  $\{ y_n \}$ строго монотонна:
\begin{equation}
  y_1<y_2<...<y_n<...
\qquad \text{или}\qquad
  y_1>y_2>...>y_n>...
\label{2.6.1}\end{equation}

\item[(ii)]  $\{ y_n \}$ бесконечно большая:
\begin{equation}\label{2.6.2}
  y_n\underset{n\to\infty}{\longrightarrow}\infty
\end{equation}

\item[(iii)] последовательность $\frac{z_n-z_{n-1}}{y_n-y_{n-1}}$ имеет
конечный предел:
\begin{equation}\lim_{n\to\infty}\frac{z_n-z_{n-1}}{y_n-y_{n-1}}=A \label{2.6.3}\end{equation}
 }\eit
Тогда последовательность $\{ \frac{z_n}{y_n}\}$ тоже имеет конечный предел, и
он совпадает с $A$:
\begin{equation}\label{2.6.4}
\lim_{n\to\infty}\frac{z_n}{y_n}= \lim_{n\to\infty}\frac{z_n-z_{n-1}}{y_n-y_{n-1}}=A
\end{equation}
\end{tm}
\begin{proof}
Заметим сразу, что достаточно рассмотреть случай, когда $y_n$ возрастает:
$$
  y_1<y_2<...<y_n<...,
$$
поскольку случай убывающей последовательности
$$
  y_1>y_2>...>y_n>...
$$
сводится к случаю возрастающей заменой $\tilde{y}_n=-y_n, \, \l
\tilde{y}_1<\tilde{y}_2<...<\tilde{y}_n<... \r$. Заодно можно считать, что все
$y_n$ положительны:
\begin{equation}\label{2.6.5}
  0<y_1<y_2<...<y_n<...
\end{equation}
Из условия $(iii)$ следует, что
$\frac{z_n-z_{n-1}}{y_n-y_{n-1}}$ представимо в виде
\begin{equation}\frac{z_n-z_{n-1}}{y_n-y_{n-1}}=A+\alpha_n, \label{2.6.6}\end{equation}
где
\begin{equation}\alpha_n \underset{n\to\infty}{\longrightarrow} 0 \label{2.6.7}\end{equation}
Зафиксируем $\varepsilon>0$, и подберем $N$ так чтобы
\begin{equation}
  \forall n>N\quad |\alpha_n|<\frac{\varepsilon}{2}\label{2.6.8}\end{equation}
(это можно сделать, потому что $|\alpha_n|
\underset{n\to\infty}{\longrightarrow} 0$).

Из \eqref{2.6.2} следует, что существует $M>N$ такое, что
$$
  \forall n>M\quad y_n>
\frac{2}{\varepsilon}\cdot \Big|z_N-A\cdot y_N\Big|$$ то есть
\begin{equation}
  \forall n>M\quad \Big| z_N-A\cdot y_N\Big|<\frac{\varepsilon}{2}\cdot y_n
\label{2.6.9}
\end{equation}
Зафиксируем эти $N$ и $M$. Из \eqref{2.6.6} получаем: для всякого $n>M>N$
$$
z_n-z_{n-1}=(\alpha_n+A)\cdot (y_n-y_{n-1})=\alpha_n\cdot(y_n-y_{n-1})+A\cdot
(y_n- y_{n-1})
$$
$$
\Downarrow
$$
$$
\begin{cases}
z_n-z_{n-1}=\alpha_n\cdot (y_n-y_{n-1})+A\cdot (y_n-y_{n-1})\\
z_{n-1}-z_{n-2}=\alpha_{n-1}\cdot (y_{n-1}-y_{n-2})+A\cdot (y_{n-1}-y_{n-2})\\
...\\
z_{N+1}-z_N=\alpha_{N+1}\cdot(y_{N+1}-y_N)+A\cdot (y_{N+1}- y_N)
\end{cases}
$$
$$
\phantom{\text{\scriptsize (складываем
уравнения)}}\quad\Downarrow\quad\text{\scriptsize (складываем уравнения)}
$$
 \begin{multline*}
\overbrace{z_n-z_{n-1}+...+z_{N+1}-z_N}^{\scriptsize\begin{matrix}\text{после сокращений остаются}\\
\text{только первое и последнее слагаемое}\end{matrix}}=\alpha_n\cdot
(y_n-y_{n-1})+...+\alpha_{N+1}\cdot(y_{N+1}-y_N)+\\+ A\cdot
(\underbrace{y_n-y_{n-1}+...+y_{N+1}- y_N}_{\scriptsize\begin{matrix}\text{после сокращений остаются}\\
\text{только первое и последнее слагаемое}\end{matrix}})
 \end{multline*}
$$
\Downarrow
$$
$$
z_n-z_N=\alpha_n\cdot (y_n-y_{n-1})+...+\alpha_{N+1}\cdot(y_{N+1}-y_N)+
A\cdot(y_n-y_N)
$$
$$
\Downarrow
$$
$$
z_n-A\cdot y_n=\alpha_n\cdot (y_n-y_{n-1})+...+\alpha_{N+1}\cdot(y_{N+1}-y_N)+
z_N-A\cdot y_N
$$
$$
\Downarrow
$$
\begin{align*}
\Big| z_n-A\cdot y_n\Big| &= \Big| \alpha_n (y_n-y_{n-1})+...+
\alpha_{N+1} (y_{N+1}-y_N)+z_{N}-A\cdot y_{N}\Big|\le \\
 &\le
 \kern-20pt\underbrace{|\alpha_n|}_{\scriptsize\begin{matrix}
 \phantom{\ \eqref{2.6.8}}\text{\rotatebox{90}{$>$}}\ \eqref{2.6.8} \\ \frac{\e}{2}
 \end{matrix}}\kern-20pt \cdot
\Big| y_n-y_{n-1}\Big|+...+
\kern-18pt\underbrace{|\alpha_{N+1}|}_{\scriptsize\begin{matrix}
 \phantom{\ \eqref{2.6.8}}\text{\rotatebox{90}{$>$}}\ \eqref{2.6.8} \\ \frac{\e}{2}
 \end{matrix}}\kern-18pt\cdot \Big| y_{N+1}-y_N\Big|+\Big|
z_{N}-A\cdot y_{N}\Big|<  \\
 & < \frac{\varepsilon}{2}\cdot
\Big|\kern-10pt\underbrace{y_n-y_{n-1}}_{\scriptsize\begin{matrix} \phantom{\
\eqref{2.6.1}} \text{\rotatebox{90}{$<$}}\ \eqref{2.6.1} \\ 0
 \end{matrix}}\kern-10pt\Big|+...+
\frac{\varepsilon}{2}\cdot\Big|\kern-10pt\underbrace{y_{N+1}-y_N}_{\scriptsize\begin{matrix}
\phantom{\ \eqref{2.6.1}} \text{\rotatebox{90}{$<$}}\ \eqref{2.6.1} \\ 0
 \end{matrix}}\kern-10pt\Big|+ \Big| z_{N}-A\cdot y_{N}\Big|=
\\
& = \frac{\varepsilon}{2}\cdot\Big( y_n-y_{n-1}\Big)+...+
\frac{\varepsilon}{2}\cdot\Big(
y_{N+1}-y_N\Big)+ \Big| z_{N}-A\cdot y_{N}\Big|=\\
 &= \frac{\varepsilon}{2}\cdot\Big(
\underbrace{y_n-y_{n-1}+...+y_{N+1}-y_N}_{\smsize
\begin{matrix}\text{после сокращений остаются}\\
\text{только первое и последнее слагаемое}\end{matrix}}\Big)+ \Big|
z_{N}-A\cdot y_{N}\Big|=\\
 &= \frac{\varepsilon}{2}\cdot\Big(
\underbrace{y_n-\kern-22pt\overbrace{y_N}^{\scriptsize\begin{matrix}0\\
\phantom{\ \eqref{2.6.5}} \text{\rotatebox{90}{$>$}}\ \eqref{2.6.5}
 \end{matrix}}\kern-22pt}_{\scriptsize\begin{matrix}
\text{\rotatebox{90}{$>$}}\\ y_n
 \end{matrix}} \Big)+
\underbrace{\Big| z_{N}-A\cdot y_{N}\Big|}_{\scriptsize\begin{matrix}
\phantom{\ \eqref{2.6.9}} \text{\rotatebox{90}{$>$}}\ \eqref{2.6.9} \\
\frac{\varepsilon}{2}\cdot y_n
 \end{matrix}}<
\frac{\varepsilon}{2}\cdot y_n+ \frac{\varepsilon}{2}\cdot y_n=
\varepsilon\cdot y_n
\end{align*}
$$
\Downarrow
$$
$$
\Big| z_n-A\cdot y_n\Big|< \varepsilon\cdot y_n
$$
$$
\Downarrow
$$
$$
\Big| \frac{z_n}{y_n}-A\Big|< \e
$$
У нас получилось, что для всякого $\varepsilon>0$ найдется $M$ такое, что
$$
\forall n>M\quad \Big| \frac{z_n}{y_n}-A\Big|< \varepsilon
$$
То есть
$$
\forall n>M \quad -\varepsilon< \frac{z_n}{y_n}-A< \varepsilon
$$
Иными словами, для любого $\varepsilon>0$ почти все числа
$\frac{z_n}{y_n}-A$ лежат в интервале $(-\varepsilon,
\varepsilon)$. Значит,
$$
\frac{z_n}{y_n}-A\underset{n\to\infty}{\longrightarrow} 0
$$
то есть,
$$
\frac{z_n}{y_n}\underset{n\to\infty}{\longrightarrow} A
\qquad $$
\end{proof}

\noindent\rule{160mm}{0.1pt}\begin{multicols}{2}

\section{Приложения предела последовательности}

\subsection{Десятичная запись вещественных чисел.}

\blm\label{LM:desyat-zapis-x<1} Для всякого вещественного числа $x\in(0,1)$
последовательность
 \beq\label{x_n=(x10^n)/10^n}
x_N=\frac{[x\cdot 10^N]}{10^N}
 \eeq
обладает следующими свойствами:
 \biter{

\item[(i)] для всякого $N\in\N$ выполняется неравенство
 \beq\label{0<x-x_n<1/10^n}
0\le x-x_N<\frac{1}{10^N}
 \eeq
из которого следует, что $x_N$ сходится к $x$:
 \beq
x_N\underset{N\to\infty}{\longrightarrow}x
 \eeq

\item[(ii)] последовательность $\{a_{-k};\ k\in\N\}$, определенная правилом
 \begin{align}
\kern-15pt & a_{-1}=x_1\cdot 10=[x\cdot 10],\label{desyat-zapis-a_(-1)}\\
\kern-15pt & a_{-k}=(x_k-x_{k-1})\cdot 10^k\quad(k>1)
\label{desyat-zapis-a_(-k)}
 \end{align}
состоит из целых чисел в интервале $\{0,1,...,9\}$:
 \beq\label{0<a_(-n)<9}
\forall k\in\N\quad a_{-k}\in\{0,1,...,9\}
 \eeq

\item[(iii)] справедливо тождество
 \beq\label{x_n=sum-a_(-k)10^(-k)}
x_N=\sum_{k=1}^N a_{-k}\cdot 10^{-k}
 \eeq
 }\eiter
  \elm
\bpr

1. Начнем с \eqref{0<x-x_n<1/10^n}:
$$
\{x\cdot 10^N\}\in[0,1)\qquad{\scriptsize \eqref{drob-chast<1}}
$$
$$
\Downarrow
$$
 \begin{multline*}
x-x_N=x-\frac{[x\cdot 10^N]}{10^N}=\\=\frac{x\cdot10^N}{10^N}-\frac{[x\cdot
10^N]}{10^N}=\frac{x\cdot10^N-[x\cdot
10^N]}{10^N}=\\=\eqref{drobnaya-chast-chisla}=\frac{\{x\cdot
10^N\}}{10^N}\in\left[0,\frac{1}{10^N}\right)
 \end{multline*}

2. Докажем \eqref{0<a_(-n)<9}. Сначала заметим, что
 $$
\forall k\in\N\qquad x_k\cdot 10^k=[x\cdot 10^k]\in\Z
 $$
$$
\Downarrow
$$
$$
a_{-k}=\left\{\begin{matrix}x_1, & k=1\\ (x_k-x_{k-1})\cdot 10^k, &
k>1\end{matrix}\right\}\in\Z
$$
Теперь для $k=1$ получаем:
$$
0<x<1
$$
$$
\Downarrow
$$
$$
0<x\cdot 10<10
$$
$$
\phantom{\scriptsize\eqref{m<x<n=>m<tch-x<n}}\qquad\Downarrow\qquad{\scriptsize\eqref{m<x<n=>m<tch-x<n}}
$$
$$
0\le [x\cdot 10]<10
$$
$$
\Downarrow
$$
$$
a_{-1}=[x\cdot 10]\in\{0,1,...,9\}
$$
Затем для $k>1$ получаем:
$$
\phantom{\scriptsize \eqref{0<x-x_n<1/10^n}}\qquad0\le
x-x_{k-1}<\frac{1}{10^{k-1}}\qquad{\scriptsize \eqref{0<x-x_n<1/10^n}}
$$
$$
\Downarrow
$$
$$
0\le (x-x_{k-1})\cdot 10^k<10
$$
$$
\phantom{\scriptsize\eqref{m<x<n=>m<tch-x<n}}\qquad\Downarrow\qquad{\scriptsize\eqref{m<x<n=>m<tch-x<n}}
$$
$$
0\le \underbrace{[(x-x_{k-1})\cdot
10^k]}_{\scriptsize\begin{matrix}\text{\rotatebox{90}{$=$}}\\ [x\cdot
10^k-x_{k-1}\cdot 10^k] \\ \text{\rotatebox{90}{$=$}}\\ [x\cdot
10^k]-x_{k-1}\cdot 10^k\\ \text{\rotatebox{90}{$=$}}\\ \l\frac{[x\cdot
10^k]}{10^k}-x_{k-1}\r\cdot 10^k \\ \text{\rotatebox{90}{$=$}}\\
\l x_k-x_{k-1}\r\cdot 10^k\\ \text{\rotatebox{90}{$=$}}\\
a_{-k}
\end{matrix}}<10
$$
$$
\Downarrow
$$
$$
a_{-k}\in\{0,1,...,9\}
$$

3. Формула \eqref{x_n=sum-a_(-k)10^(-k)} доказывается индукцией. При $N=1$ она
принимает вид:
$$
x_1=a_{-1}\cdot 10^{-1},
$$
и это эквивалентно \eqref{desyat-zapis-a_(-1)}. Далее, если
\eqref{x_n=sum-a_(-k)10^(-k)} верна для $N=m$,
$$
\sum_{k=1}^m a_{-k}\cdot 10^{-k}=x_m,
$$
то для $N=m+1$ мы получаем:
 \begin{multline*}
\sum_{k=1}^{m+1} a_{-k}\cdot 10^{-k}=\\=\overbrace{\sum_{k=1}^m a_{-k}\cdot
10^{-k}}^{x_m}+ a_{-m-1}\cdot 10^{-m-1}=\\= x_m+a_{-m-1}\cdot
10^{-m-1}=\eqref{desyat-zapis-a_(-k)}=\\=x_m+\big((x_{m+1}-x_m)\cdot
10^{m+1}\big)\cdot 10^{-m-1}=\\=x_m+(x_{m+1}-x_m)=x_{m+1}
 \end{multline*}
\epr

Из леммы \ref{LM:desyat-zapis-x<1} и теоремы \ref{TH-desyatichnaya-zapis-v-N}
следует

\btm\label{TH-desyatichnaya-zapis-v-R} Всякое число $x>0$ представимо в виде
 \beq\label{desyatichnaya-zapis-chisla-x>0}
x=\sum_{k=0}^n a_k\cdot 10^k+\lim_{N\to\infty}\sum_{k=1}^N a_{-k}\cdot 10^{-k}
 \eeq
где $n\in\Z_+$, $a_i\in\{0,1,...,9\}$. \etm
 \bpr
Разложим $x$ в сумму целой и дробной части:
$$
x=[x]+\{x\}
$$
Если $[x]=0$, то в \eqref{desyatichnaya-zapis-chisla-x>0} нужно положить $n=0$,
$a_0=0$, если же $[x]\ne 0$, то по теореме \ref{TH-desyatichnaya-zapis-v-N}
$[x]$ можно представить в виде
$$
[x]=\sum_{k=0}^n a_k\cdot 10^k\qquad (a_k\in\{0,1,...,9\})
$$
То же самое с дробной частью: если $\{x\}=0$, то в
\eqref{desyatichnaya-zapis-chisla-x>0} нужно положить $a_{-k}=0$ при $k\in\N$,
а если $\{x\}\ne 0$, то  $\{x\}$ можно по лемме \ref{LM:desyat-zapis-x<1}
записать в виде
$$
\{x\}=\lim_{N\to\infty}\sum_{k=1}^N a_{-k}\cdot 10^{-k}\qquad
(a_{-k}\in\{0,1,...,9\}).
$$
 \epr

 \biter{
\item[$\bullet$] Формула \eqref{desyatichnaya-zapis-chisla-x>0} называется {\it
десятичным представлением} числа $x>0$, а последовательность
$\{a_n,a_{n-1},...,a_1,a_0,a_{-1},...\}$ -- {\it последовательностью
коэффициентов} этого представления. Если начиная с номера $N+1$ все числа
$a_{-k}$ ($k\ge N+1$) равны нулю, то справедливо равенство
$$
x=\sum_{k=0}^n a_k\cdot 10^k+\sum_{k=1}^N a_{-k}\cdot 10^{-k}
$$
которое коротко принято записывать так:
$$
x=a_n...a_1a_0,a_{-1}...a_{-N}
$$
(запятая отделяет цифры с индексами $0$ и $-1$), и это называется {\it
десятичной записью} числа $x$.
 }\eiter

\bex Запись
$$
x=45,28
$$
расшифровывается как равенство
$$
x=4\cdot 10^1+5\cdot 10^0+2\cdot 10^{-1}+8\cdot 10^{-2}
$$
\eex

 \biter{
\item[$\bullet$]
Помимо этого, запись вида
$$
x=a_n...a_1a_0,a_{-1}...a_{-N}...
$$
(в которой перечислены все числа $a_k$ до номера $a_{-N}$ включительно, а после
них стоит многоточие), называется {\it приближенной десятичной записью} числа
$x$ и означает, что $x$ отличается от числа $a_n...a_1a_0,a_{-1}...a_{-N}$ не
больше чем на $10^{-N}$:
$$
0\le x-a_n...a_1a_0,a_{-1}...a_{-N}<10^{-N}
$$
 }\eiter

 \bex Запись
$$
x=102,67...
$$
означает, что справедлива оценка
$$
0\le x-102,67<10^{-2}
$$
 \eex

 \biter{
\item[$\bullet$] Если $x<0$, то по теореме \ref{TH-desyatichnaya-zapis-v-R}
число $-x>0$ представимо в виде
$$
-x=\sum_{k=0}^n a_k\cdot 10^k+\sum_{k=1}^N a_{-k}\cdot 10^{-k}
$$
Опять же, если здесь начиная с номера $N+1$ все числа $a_{-k}$ ($k\ge N+1$)
равны нулю, то справедливо равенство
$$
-x=\sum_{k=0}^n a_k\cdot 10^k+\sum_{k=1}^N a_{-k}\cdot 10^{-k}
$$
и его принято коротко записывать так:
$$
x=-a_n...a_1a_0,a_{-1}...a_{-N}
$$
(запятая отделяет цифры с индексами $0$ и $-1$). Это называется {\it десятичной
записью} числа $x<0$.
 }\eiter

 \bex Запись
$$
x=-34,72
$$
расшифровывается как равенство
$$
-x=3\cdot 10^1+4\cdot 10^0+7\cdot 10^{-1}+2\cdot 10^{-2}
$$
 \eex

 \biter{
\item[$\bullet$]
Помимо этого, запись вида
$$
x=-a_n...a_1a_0,a_{-1}...a_{-N}...
$$
(в которой перечислены все числа $a_k$ до номера $a_{-N}$ включительно, а после
них стоит многоточие), называется {\it приближенной десятичной записью} числа
$x<0$ и означает, что $-x$ отличается от числа $a_n...a_1a_0,a_{-1}...a_{-N}$
не больше чем на $10^{-N}$:
$$
0\le -x-a_n...a_1a_0,a_{-1}...a_{-N}<10^{-N}
$$
Понятно, что это равносильно оценке
$$
10^{-N}<x+a_n...a_1a_0,a_{-1}...a_{-N}\le 0
$$
 }\eiter

 \bex Запись
$$
x=-2,31...
$$
означает оценку
$$
10^{-2}<x+2,31\le 0
$$
 \eex

\subsection{Число Непера $e$.}\label{SEC:e}

Из теоремы Вейерштрасса \ref{Wei-I} выводится следующее утверждение, неожиданно
оказывающееся важным в анализе:

 \bit{\label{e}
\item[$\bullet$] Последовательность $\l 1+\frac{1}{n}\r^n$ имеет конечный
предел, который называется {\it числом Непера} и обозначается буквой $e$:
\begin{equation}\label{5.9.1}
e=\lim_{n\to \infty}\left(1+\frac{1}{n}\right)^n
\end{equation}
 }\eit
\begin{proof} Рассмотрим сначала вспомогательную
последовательность $y_n=(1+\frac{1}{n})^{n+1}$ и докажем что она сходится. Для
$n\ge 2$ мы получаем:
 \begin{multline*}
\frac{y_{n-1}}{y_n}=\frac{(1+\frac{1}{n-1})^n}{(1+\frac{1}{n})^{n+1}}=
\frac{(\frac{n}{n-1})^n}{(\frac{n+1}{n})^{n+1}}=\\=
\frac{n^{2n}}{(n^2-1)^n}\cdot \frac{n}{n+1}=
\left(\frac{n^2}{n^2-1}\right)^n\cdot \frac{n}{n+1}=\\=
\left(1+\frac{1}{n^2-1}\right)^n\cdot \frac{n}{n+1}\ge \left({\smsize
\begin{array}{c} \text{неравенство}\\ \text{Бернулли}\\
\eqref{nerav-Bernoulli}\end{array}}\right)\ge\\ \ge
\left(1+\frac{n}{n^2-1}\right)\cdot \frac{n}{n+1}> \\ >
\left(1+\frac{n}{n^2}\right)\cdot \frac{n}{n+1}=
\left(1+\frac{1}{n}\right)\cdot \frac{n}{n+1}=1
\end{multline*} То есть последовательность $y_n$
монотонно убывает:
$$
y_{n-1}> y_n
$$
С другой стороны, все числа $y_n$ положительны, и поэтому она ограничена:
$$
0<...<y_n<y_{n-1}<...<y_3<y_2<y_1
$$
По теореме Вейерштрасса \ref{Wei-I}, это означает, что последовательность $y_n$
сходится. Обозначим ее предел буквой $e$:
 \beq\label{(1+1/n)^n+1}
e=\lim_{n\to \infty} y_n=\lim_{n\to \infty}\l 1+\frac{1}{n}\r^{n+1}
 \eeq
Теперь сходимость последовательности $\l 1+\frac{1}{n}\r^n$ следует из правила
вычисления предела дроби (свойство $4^\circ$ на с.
\pageref{arifm-oper-s-predelami-posl}):
$$
\l 1+\frac{1}{n}\r^n=\frac{\l 1+\frac{1}{n}\r^{n+1}}{1+\frac{1}{n}}= \frac{
\boxed{y_n}\put(-5,9){\vector(1,2){5}\put(2,16){$e$}}
}{1+\boxed{\frac{1}{n}}\put(-5,-12){\vector(1,-2){5}\put(2,-16){0}}
}\underset{n\to \infty}{\longrightarrow}\frac{e}{1}=e \qquad
$$
Остается заметить, что в этом изложении мы, формально говоря, определили число
$e$ иначе, чем обещали: вместо \eqref{5.9.1}, мы определили $e$ по формуле
\eqref{(1+1/n)^n+1}, а затем \eqref{5.9.1} вывели как следствие. Но, поскольку
из наших рассуждений следует, что пределы в \eqref{(1+1/n)^n+1} и \eqref{5.9.1}
совпадают, неважно, какую из этих формул считать определением, а какую
следствием.
\end{proof}

\btm\label{TH:e=sum-1/k!} Число Непера $e$ удовлетворяет соотношению:
 \beq\label{e=sum-1/k!}
e=\lim_{n\to\infty}\sum_{k=0}^n\frac{1}{k!}
 \eeq
причем для любого $n\in\N$ справедливо двойное неравенство:
 \beq\label{0<e-sum-1/k!<1/n!n}
0<e-\sum_{k=0}^n\frac{1}{k!}\le \frac{1}{n!\cdot n}
 \eeq
или, что то же самое двойное неравенство
 \beq\label{0<e-sum-1/k!<1/n!n-1}
\sum_{k=0}^n\frac{1}{k!}<e\le \sum_{k=0}^n\frac{1}{k!}+\frac{1}{n!\cdot n}
 \eeq
\etm
 \bpr
Обозначим
$$
E_n=\sum_{k=0}^n\frac{1}{k!}
$$
и докажем сначала равенство \eqref{e=sum-1/k!}.

1. Поскольку
 \begin{multline*}
E_n=\sum_{k=0}^n\frac{1}{k!}<\sum_{k=0}^n\frac{1}{k!}+\frac{1}{(n+1)!}=\\=\sum_{k=0}^{n+1}\frac{1}{k!}=E_{n+1},
 \end{multline*}
эта последовательность возрастает.

С другой стороны, из \eqref{n!-ge-2^(n-1)} получаем:
$$
\forall k\in\Z_+\qquad k!\ge 2^{k-1}
$$
$$
\Downarrow
$$
$$
\forall k\in\Z_+\qquad \frac{1}{k!}\le \frac{1}{2^{k-1}}
$$
$$
\Downarrow
$$
 \begin{multline*}
E_n=\sum_{k=0}^n\frac{1}{k!}\le
\sum_{k=0}^n\frac{1}{2^{k-1}}=2\cdot\sum_{k=0}^n\left(\frac{1}{2}\right)^k=\\=\eqref{Geom-progr-N}=
2\cdot\frac{1-\left(\frac{1}{2}\right)^{n+1}}{1-\frac{1}{2}}=
2\cdot\frac{1-\left(\frac{1}{2}\right)^{n+1}}{\frac{1}{2}}=\\=
4\cdot\left(1-\left(\frac{1}{2}\right)^{n+1}\right)=4-\frac{1}{2^{n-1}}<4,
 \end{multline*}
и значит она ограничена сверху.

По теореме Вейерштрасса \ref{Wei-I} получаем, что она имеет предел. Обозначим
его буквой $E$:
$$
E_n\underset{n\to\infty}{\longrightarrow} E
$$
Наша задача -- проверить, что $E=e$. С одной стороны,
 \begin{multline*}
\l 1+\frac{1}{n}\r^n=\eqref{binom-Newtona}=\sum_{k=0}^n C_n^k\cdot
\frac{1}{n^k}=\\=\sum_{k=0}^n \frac{n!}{k!\cdot(n-k)!}\cdot \frac{1}{n^k}=\\=
\sum_{k=0}^n \frac{1}{k!}\cdot \frac{n!}{(n-k)!}\cdot
\left(\frac{1}{n}\right)^k=\\=\sum_{k=0}^n\frac{1}{k!}\cdot \underbrace{n\cdot
...\cdot(n-k+1)}_{\text{$k$ множителей}}\cdot
\underbrace{\frac{1}{n}\cdot...\cdot\frac{1}{n}}_{\text{$k$
множителей}}=\\=\sum_{k=0}^n\frac{1}{k!}\cdot
\underbrace{\frac{n}{n}\cdot\frac{n-1}{n}\cdot...\cdot\frac{n-k+1}{n}}_{\text{$k$
множителей}}=\\=\sum_{k=0}^n\frac{1}{k!}\cdot \underbrace{1\cdot\l
1-\frac{1}{n}\r\cdot...\cdot\l 1-\frac{k-1}{n}\r}_{\text{$k$ множителей, и все
$\le 1$}}\le\\ \le \sum_{k=0}^n\frac{1}{k!}=E_n
 \end{multline*}
$$
\Downarrow
$$
$$
e=\lim_{n\to\infty}\l 1+\frac{1}{n}\r^n\le\lim_{n\to\infty}E_n=E
$$
А с другой -- при любом $m\le n$
 \begin{multline*}
\overbrace{\sum_{k=0}^m\frac{1}{k!}\cdot 1\cdot\l
1-\frac{1}{n}\r\cdot...\cdot\l 1-\frac{k-1}{n}\r}^{\begin{matrix}
E_m=\sum_{k=0}^m\frac{1}{k!}\\ \phantom{\scriptsize
\text{$(n\to\infty)$}}\uparrow {\scriptsize \text{$(n\to\infty)$}}
\end{matrix}}\le\\ \le \sum_{k=0}^n\frac{1}{k!}\cdot 1\cdot\l
1-\frac{1}{n}\r\cdot...\cdot\l 1-\frac{k-1}{n}\r=\\=\kern-10pt\underbrace{\l
1+\frac{1}{n}\r^n}_{\begin{matrix}\phantom{\scriptsize \text{$(n\to\infty)$}}\downarrow {\scriptsize \text{$(n\to\infty)$}} \\
e \end{matrix}}
 \end{multline*}
$$
\Downarrow
$$
$$
\kern50pt  E_m\le e\qquad (m\in\N)
$$
$$
\Downarrow
$$
$$
E=\lim_{m\to\infty} E_m\le e
$$

2. Докажем теперь \eqref{0<e-sum-1/k!<1/n!n}. Первое неравенство здесь следует
из возрастания последовательности $E_n$:
$$
e=\lim_{m\to\infty} E_m=\sup_{m\in\Z_+} E_m\ge E_{n+1}>E_n
$$
$$
\Downarrow
$$
$$
e-E_n>0
$$
А второе -- из неравенства \eqref{(n+m)!-ge-n!(n+1)^m}: для любых $n,p\in\N$ мы
получим
$$
\phantom{{\scriptsize\eqref{(n+m)!-ge-n!(n+1)^m}}}\qquad(n+m)!\ge n!\cdot
(n+1)^m\qquad{\scriptsize\eqref{(n+m)!-ge-n!(n+1)^m}}
$$
$$
\Downarrow
$$
$$
\frac{1}{(n+m)!}\le\frac{1}{n!\cdot (n+1)^m}
$$
$$
\Downarrow
$$
 \begin{multline*}
E_{n+p}-E_n=\sum_{k=0}^{n+p}\frac{1}{k!}-\sum_{k=0}^n\frac{1}{k!}=\sum_{k=n+1}^{n+p}\frac{1}{k!}=\\=
{\scriptsize\left|\begin{matrix}m=k-n \\ k=m+n \\ 1\le m\le p
\end{matrix}\right|}=\sum_{m=1}^{p}\frac{1}{(n+m)!}\le\\ \le
\sum_{m=1}^{p}\frac{1}{n!\cdot(n+1)^m}=
\frac{1}{n!}\cdot\sum_{m=1}^{p}\left(\frac{1}{n+1}\right)^m=\\=\eqref{Geom-progr-m<n}=
\frac{1}{n!}\cdot
\frac{\frac{1}{n+1}-\left(\frac{1}{n+1}\right)^{p+1}}{1-\frac{1}{n+1}}=\\=
\frac{1}{n!}\cdot \frac{1-\left(\frac{1}{n+1}\right)^p}{n+1-1}=
\frac{1}{n!}\cdot \frac{1-\left(\frac{1}{n+1}\right)^p}{n}
 \end{multline*}
$$
\Downarrow
$$
$$
\underbrace{E_{n+p}-E_n}_{\scriptsize\begin{matrix}(p\to\infty)\downarrow\phantom{(p\to\infty)}
\\ e-E_n\end{matrix}}\le \underbrace{\frac{1}{n!}\cdot
\frac{1-\left(\frac{1}{n+1}\right)^p}{n}}_{\scriptsize\begin{matrix}\phantom{(p\to\infty)}\downarrow
(p\to\infty)
\\ \frac{1}{n!\cdot n}\end{matrix}}
$$
$$
\Downarrow
$$
$$
e-E_n\le\frac{1}{n!\cdot n}
$$
 \epr

\brem Двойное неравенство \eqref{0<e-sum-1/k!<1/n!n-1} позволяет оценивать
число $e$. Для нахождения двух первых знаков после запятой можно взять $n=7$:
$$
\sum_{k=0}^7\frac{1}{k!}<e\le \sum_{k=0}^7\frac{1}{k!}+\frac{1}{7!\cdot 7}
$$
Значения величин справа и слева можно приближенно вычислить (на бумаге, или с
помощью калькулятора), получив оценки сверху и снизу, из которых затем
выводится оценка для $e$: с одной стороны,
$$
\begin{matrix}
\frac{1}{8!\cdot 8}=0,000198412698... \\
 \Downarrow \\
 0,00019<\frac{1}{8!\cdot 8}<0,0002
\end{matrix}
$$
с другой --
 $$
\begin{matrix}
\sum_{k=0}^8\frac{1}{k!}=2,71827877... \\
\Downarrow \\
2,718<\sum_{k=0}^8\frac{1}{k!}<2,719
\end{matrix}
 $$
И вместе это дает:
 \begin{multline*}
2,718<\sum_{k=0}^8\frac{1}{k!}<e< \sum_{k=0}^8\frac{1}{k!}+\frac{1}{8!\cdot
8}<\\<2,719+0,0002=2,7192
 \end{multline*}
$$
\Downarrow
$$
 \beq\label{otsenka-e}
\kern32pt\boxed{\quad 2,71<e<2,72\quad}
 \eeq
$$
\Downarrow
$$
 \beq\label{otsenka-e-1}
\kern32pt\boxed{\quad e=2,71...\quad}
 \eeq
\erem

\btm Число $e$ иррационально:
$$
e\notin\Q
$$
 \etm
 \bpr
Предположим, что оно рационально, то есть имеет вид
$$
e=\frac{m}{n}
$$
для некоторых $m,n\in\N$. Из оценки \eqref{otsenka-e} следует, что $n\ge 2$
(потому что при $n=1$ мы получили бы, что $e$ целое). Запомним это. Из двойного
неравенства \eqref{0<e-sum-1/k!<1/n!n} мы получаем:
$$
0<\frac{m}{n}-\sum_{k=n}^n\frac{1}{k!}\le \frac{1}{n!\cdot n}
$$
$$
\Downarrow
$$
$$
0<\frac{m\cdot n!}{n}-\sum_{k=n}^n\frac{n!}{k!}\le \frac{1}{n}
$$
$$
\Downarrow
$$
 \begin{multline*}
0<\underbrace{\underbrace{m\cdot (n-1)!}_{\text{целое}}-\sum_{k=n}^n
\underbrace{n\cdot(n-1)\cdot...\cdot(n-k)}_{\text{целое}}}_{\text{целое}}\le\\
\le \frac{1}{n}\le\frac{1}{2}
 \end{multline*}
Мы получаем, что в полуинтервале $\left(0;\frac{1}{2}\right]$ лежит какое-то
целое число. Понятно, что это невозможно, и значит наше предположение (о
рациональности $e$) было неверным.
 \epr

\end{multicols}\noindent\rule[10pt]{160mm}{0.1pt}

\chapter{ФУНКЦИЯ, ЕЕ НЕПРЕРЫВНОСТЬ И ПРЕДЕЛ}\label{ch-cont-f(x)}

\section{Числовые функции}\label{SEC-numb-function}

\subsection{Определение функции и примеры}

Понятие функции формализует в математике идею взаимной зависимости физических
величин. Идея же состоит в следующем.

Как мы уже говорили на с.\pageref{SEC-chisla}, физические величины, такие как
расстояние, масса, скорость, и так далее, поддаются количественной оценке в
числах. Люди находят и устанавливают правила, позволяющие измерять эти величины
(сравнивая их с эталонами), то есть выражать их числами: 200 метров, 1,8
килограмма, и т.д. Собственно говоря, в этом состоит смысл самого понятия
физической величины: параметр, существование которого предполагается
теоретически, должен быть снабжен системой правил, позволяющих его измерять и
сравнивать значения в различных ситуациях -- только в этом случае он
приобретает статус физической величины.

При этом очень важное наблюдение, с которого и начинается
вся наука, состоит в том, что при фиксированных правилах измерений
обнаруживается, что разные величины могут быть связаны друг с другом довольно
жесткими соотношениями (законами природы): если Вам известны результаты
измерений какой-то величины, то это, как правило, означает, что Вы имеете
возможность вычислить какие-то другие величины в рассматриваемой Вами ситуации
(без необходимости их измерять).

\noindent\rule{160mm}{0.1pt}\begin{multicols}{2}

\bex\label{EX:l=pi.D} Например, по результатам измерения диаметра окружности
$D$, можно с достаточной точностью судить о том, какой будет ее длина $l$,
 \beq\label{l-D-ex-func}
l=\pi D,
 \eeq
или площадь ограничиваемого ею круга:
 \beq\label{S-D-ex-func}
S=\frac{\pi}{4} D^2.
 \eeq
Важно, что Вам не нужно измерять $l$ и $S$, если Вы измерили $D$, потому что
зависимости \eqref{l-D-ex-func} и \eqref{S-D-ex-func} позволяют вычислять $l$ и
$S$ по данному $D$. \eex

\bex\label{EX:V=rho.m} Точно также, зная (с достаточной точностью) массу тела
$m$ и плотность $\rho$ материала, из которого оно изготовлено, Вы можете (опять
же, достаточно точно) определить объем $V$ этого тела по формуле
$$
V=\rho\cdot m.
$$
Вам для этого не нужно ставить экспериментов, требующих времени и денег (если,
скажем, изучаемый Вами предмет довольно громоздок) -- все легко и быстро
вычисляется на бумаге (или в уме). \eex

\end{multicols}\noindent\rule[10pt]{160mm}{0.1pt}

В этих примерах одна величина выражается через другую с помощью логических
правил, позволяющих проводить {\it вычисления}, то есть находить нужные
значения, возможно приближенные, пользуясь заранее разработанными алгоритмами.
Инженер, привыкший пользоваться подобными вычислениями на практике, обычно
только так и понимает функциональную зависимость. Но для исследователя, в
частности, математика, такое толкование понятия функции было бы
 \bit{
\item[--] слишком неудобным,
потому что,  если определять функциональные зависимости, как правила, для которых алгоритмы вычисления к
настоящему времени установлены, то придется сразу же в этом
определении описывать и сами эти алгоритмы вычислений, и это само по себе
утяжеляет определение настолько, что оправданность такого подхода никогда не
будет очевидной,

\item[--] и слишком узким, потому что такое определение никогда нельзя будет считать
окончательным, или хотя бы годным на обозримое будущее, поскольку со временем,
и весьма скоро, его неизбежно придется корректировать, по мере того, как будут
обнаруживаться новые зависимости между величинами, для которых будут
придумываться новые алгоритмы вычислений -- наблюдая за развитием науки, нелепо полагать, что человечество
когда-нибудь достигнет такого состояния, когда зависимости между различными
физическими величинами, уже открытыми и еще не найденными, будут описываться
зафиксированными к тому моменту алгоритмами вычислений.
 }\eit
Эти соображения приводят математиков к следующему самому широкому определению
функции, как просто правила, которое можно описать средствами языка:

 \bit{
\item[$\bullet$] {\it Числовой функцией} называется произвольное отображение
$f:X\to Y$, у которого область определения $X$ и область значений $Y$ оба
являются числовыми множествами:
$$
X\subseteq \R,\quad Y\subseteq \R
$$
Как и у любого отображения (см. определения на с.\pageref{def-domain-range}),
область определения $X$ функции $f$ мы обозначаем $\D(f)$.
 }\eit\noindent
Интуитивно яснее, хотя и менее точно, звучит следующая переформулировка этого
определения, скопированная с {\sl неформального} определения отображения,
данного нами на с. \pageref{neformalnoe-opred-otobr}:

 \bit{
\item[$\bullet$] Пусть даны два числовых множества $X$ и $Y$, и пусть задано
некое правило, которое {\sl каждому} элементу $x\in X$ ставит в соответствие
{\sl единственный} (зависящий от $x$) элемент $y=f(x)\in Y$. Такое правило
обозначается $f:X\to Y$ или просто $f$ и называется {\it числовой функцией} с
{\it областью определения} $X$ и {\it множеством значений} $Y$.
 }\eit

Напомним, что {\sl формальное} определение на с.
\pageref{opredelenie-otobrazheniya} не различает отображение и его график. В
соответствии с ним, мы можем сказать, что числовая функция $f:X\to Y$
($X,Y\subseteq\R$) представляет из себя подмножество в декартовом произведении
$X\times Y$, удовлетворяющее условиям 1) - 2) на
с.\pageref{opredelenie-otobrazheniya}. Обычно все же такое различие проводится:
человеческому воображению удобно представлять себе функцию как правило, а не
как подмножество, поэтому мы будем функцию $f:X\to Y$ считать правилом, а к ее
графику относиться как к множеству в декартовом произведении $X\times Y$,
состоящему из упорядоченных пар вида $(x;f(y))$, $x\in X$. Поскольку при этом
оба множества $X$ и $Y$ содержатся в $\R$, их декартово произведение $X\times
Y$ можно считать подмножеством в декартовой плоскости
$$
X\times Y\subseteq\R\times\R.
$$
Поэтому каждую пару $(x;f(y))$, $x\in X$, можно считать точкой на декартовой
плоскости
$$
(x;f(y))\in X\times Y\subseteq\R\times\R,\qquad x\in X,
$$
а {\it график функции}\index{график!функции} $f:X\to Y$ --- подмножеством на
декартовой плоскости $\R\times\R$ (состоящим из точек $(x;y)$ удовлетворяющих
условию $f(x)=y$).

\noindent\rule{160mm}{0.1pt}\begin{multicols}{2}

\bex{\bf Степенная функция с целым показателем.} Напомним, что на странице
\pageref{stepeni-n-in-Z} мы определили степени с целым показателем. Степенной
функцией с показателем $n\in\Z$ называется функция
 $$
f(x)=x^n
 $$
При $n\ge 0$ она определена на всей числовой прямой $\R$, а при $n<0$ -- только
на множестве $\R\setminus\{0\}=(-\infty;0)\cup(0;+\infty)$.
 \eex

\bex {\bf Модуль числа.}\label{|x|-ex-func} Выше формулой \eqref{|x|=cases} мы
определили модуль числа $x\in\R$:
$$
|x|=\begin{cases}x, & \text{если}\,\, x\ge 0
 \\ x, & \text{если}\,\, x<0 \end{cases}
$$
Понятно, что правило $x\mapsto|x|$ будет функцией, определенной всюду на
числовой прямой $\R$. \eex

\bex {\bf Целая часть числа.} Формулой \eqref{tselaya-chast-chisla} выше мы
определили целую часть числа $x\in\R$:
 $$
[x]=\max \{n\in\Z:\; n\le x\},
 $$
Правило $x\mapsto[x]$ будет функцией, определенной всюду на числовой прямой
$\R$. Ее график интересен тем, что имеет бесконечно много разрывов (что такое
разрыв мы объясним ниже на с. \pageref{SEC-nepreryvnost-na-mnozhestve}):

\vglue120pt

\eex

\bex {\bf Дробная часть числа.} Формулой \eqref{drobnaya-chast-chisla} мы
определили дробную часть числа $x\in\R$:
 $$
\{x\}:=x-[x]
 $$
Это тоже будет функция, определенная на $\R$, и ее график тоже имеет
бесконечное число разрывов:

\vglue80pt

\eex

\bex {\bf Сигнум} числа $x\in\R$ определяется правилом
 \beq\label{DEF:sgn}
\sgn x=\begin{cases}1, & \text{если}\,\, x>0
 \\
 0, & \text{если}\,\,x=0
 \\
 -1, & \text{если}\,\, x<0
 \end{cases}
 \eeq
График этой функции имеет разрыв только в одной точке $x=0$:

\vglue80pt

\eex

 \bex\label{EX:func-Dirichlet}
{\bf Функция Дирихле.} Функция $D$ на $\R$, определенная правилом
 \beq\label{func-Dirichle}
D(x)=\begin{cases}1,& x\in\Q \\ 0,& x\notin\Q \end{cases}
 \eeq
называется {\it функцией Дирихле}. В главе \ref{CH-functional-sequen} мы
докажем формулу \eqref{vyrazhenie-Dirichlet-cherez-cos}, выражающую эту функцию
как двойной предел стандартных функций. \eex

\bex {\bf Пустая функция.} В список примеров полезно включить {\it пустую
функцию}, то есть пустое отображение, определенное нами в
\ref{pustoe-otobrazhenie}, и рассматриваемое как отображение из $\varnothing$ в
$\R$. \eex

\end{multicols}\noindent\rule[10pt]{160mm}{0.1pt}

\subsection{Функции с симметриями}\label{func-s-simmetr}

\paragraph{Четные и нечетные функции.}

 \bit{

\item[$\bullet$] Функция $f$ называется

 \biter{
\item[--] {\it четной}, если
 \biter{
\item[(1)] для всякой точки $x$ из области определения $D(f)$ функции $f$ точка
$-x$ также лежит в области определения,
$$
\forall \, x\in D(f) \quad -x\in D(f)
$$
\item[(2)] выполняется тождество
$$
f(-x)=f(x)
$$
 }\eiter
Эти условия означают, что график функции $f$ должен быть симметричен
относительно оси ординат:

\vglue120pt

\item[--] {\it нечетной}, если
 \biter{
\item[(1)] для всякой точки $x$ из области определения $D(f)$ функции $f$ точка
$-x$ также лежит в области определения,
$$
\forall \, x\in D(f) \quad -x\in D(f)
$$
\item[(2)] выполняется тождество
$$
f(-x)=-f(x)
$$
 }\eiter
Эти условия означают, что график функции $f$ должен быть симметричен
относительно начала координат:

\vglue120pt
 }\eiter
 }\eit

\noindent\rule{160mm}{0.1pt}\begin{multicols}{2}

\bex\label{EX:(-x)^(2n)=x^(2n)} Степенная функция с четным показателем является
четной:
 \beq\label{(-x)^(2n)=x^(2n)}
(-x)^{2n}=x^{2n},\qquad n\in\Z
 \eeq
 \eex
\bpr Заметим, что это достаточно доказать для $n\ge 0$, потому что случай
отрицательных степеней получается отсюда в качестве следствия:
$$
(-x)^{-2n}=\frac{1}{(-x)^{2n}}=\frac{1}{x^{2n}}=x^{-2n}
$$
Для $n\ge 0$ доказательство проводится по индукции.

1. При $n=0$ утверждение очевидно:
$$
(-x)^0=1=x^0
$$

2. Предположим, что мы доказали это для какого-то $n=k$:
 \beq\label{(-x)^2k=x^2k}
(-x)^{2k}=x^{2k}
 \eeq

3. Тогда для $n=k+1$ мы получаем:
 \begin{multline*}
(-x)^{2(k+1)}=(-x)^{2k+2}=\\=\underbrace{(-x)^{2k}}_{\scriptsize\begin{matrix}
\|\\ \phantom{,}x^{2k}, \\ \text{в силу \eqref{(-x)^2k=x^2k}}\end{matrix}}
\cdot\underbrace{(-x)^2}_{\scriptsize\begin{matrix} \|\\ \phantom{,}x^2, \\
\text{в силу \eqref{(-a)(-a)=aa}}\end{matrix}}= x^{2k}\cdot
x^2=\\=x^{2k+2}=x^{2(k+1)}
 \end{multline*}
\epr

\bex\label{EX:(-x)^(2n-1)=-x^(2n-1)} Степенная функция с нечетным показателем
является нечетной:
 \beq\label{(-x)^(2n-1)=-x^(2n-1)}
(-x)^{2n-1}=-x^{2n-1},\qquad n\in\Z
 \eeq
 \eex
\bpr Как и в предыдущем примере, здесь достаточно доказать утверждение для
$n\ge 0$, потому что для $n=-k$, $k\in\Z_+$, это выводится как следствие:
$$
(-x)^{-2k-1}=\frac{1}{(-x)^{2k+1}}=-\frac{1}{x^{2k+1}}=-x^{-2k-1}
$$
Для $n\ge 0$ доказательство проводится по индукции.

1. При $n=0$ тождество \eqref{(-x)^(2n-1)=-x^(2n-1)} превращается в уже
доказанное тождество \eqref{(-a)^(-1)=-a^(-1)}:
$$
(-x)^{-1}=-x^{-1}
$$

2. Предположим, что мы доказали это для какого-то $n=k$:
 \beq\label{(-x)^(2k-1)=x^(2k-1)}
(-x)^{2k-1}=-x^{2k-1}
 \eeq

3. Тогда для $n=k+1$ мы получаем:
 \begin{multline*}
(-x)^{2(k+1)-1}=(-x)^{2k+1}=\\=\underbrace{(-x)^{2k-1}}_{\scriptsize\begin{matrix}
\|\\ \phantom{,}-x^{2k-1}, \\ \text{в силу
\eqref{(-x)^(2k-1)=x^(2k-1)}}\end{matrix}}
\cdot\underbrace{(-x)^2}_{\scriptsize\begin{matrix} \|\\ \phantom{,}x^2, \\
\text{в силу \eqref{(-a)(-a)=aa}}\end{matrix}}=-x^{2k-1}\cdot
x^2=\\=-x^{2k+1}=x^{2(k+1)-1}
 \end{multline*}
\epr

\bex Модуль является четной функцией:
 \beq
|-x|=|x|
 \eeq
\eex \bpr Если $x>0$, то
$$
|\underbrace{-x}_{\scriptsize\begin{matrix}\text{\rotatebox{90}{$>$}}\\
0\end{matrix}}|=-(-x)=x=|\underbrace{x}_{\scriptsize\begin{matrix}\text{\rotatebox{90}{$<$}}\\
0\end{matrix}}|
$$
Если $x=0$, то
$$
|-x|=|-0|=|0|=|x|
$$
Если $x<0$, то
$$
|\underbrace{-x}_{\scriptsize\begin{matrix}\text{\rotatebox{90}{$<$}}\\
0\end{matrix}}|=-x=|\underbrace{x}_{\scriptsize\begin{matrix}\text{\rotatebox{90}{$>$}}\\
0\end{matrix}}|
$$
\epr

\brem На всякий случай полезно заметить, что функции не обязаны быть только
четными или нечетными. Например, функция
$$
f(x)=1+x
$$
не будет ни четной ни нечетной. Докажите это. \erem

\end{multicols}\noindent\rule[10pt]{160mm}{0.1pt}

\paragraph{Периодические функции.}

 \bit{
\item[$\bullet$] Число $T$ называется {\it периодом} числовой функции $f$, если
 \biter{
\item[(a)] оно положительно:
$$
T>0,
$$

\item[(b)] область определения $\D(f)$ функции $f$ инвариантна относительно
сдвигов на $T$ вправо и влево:
$$
\D(f)-T=\D(f)=\D(f)+T
$$
(то есть для любого $x\in\D(f)$ числа $x+T$ и $x-T$ тоже лежат в $\D(f)$),

\item[(c)] справедливо двойное тождество:
 \beq\label{f(x+T)=f(x)}
f(x-T)=f(x)=f(x+T),\qquad x\in\D(f)
 \eeq
(здесь из-за условия (b) достаточно выполнения какого-нибудь одного из этих
двух равенств, например $f(x)=f(x+T)$).
 }\eiter

\item[$\bullet$] Число $h$ называется  {\it полупериодом} числовой функции $f$,
если число $2h$ является периодом функции $f$.

\item[$\bullet$] Числовая функция $f$ называется
 \biter{

\item[---] {\it периодической}, если она обладает каким-нибудь периодом $T$,

\item[---] {\it $T$-периодической}, если число $T$ является ее периодом,

 }\eiter
 }\eit

\noindent\rule{160mm}{0.1pt}\begin{multicols}{2}

 \bex{\bf Непериодические функции.} Понятно, что не всякая функция будет
периодической. Например, функция
 $$
f(x)=x,\qquad x\in\R
 $$
-- не периодическая, потому что для нее тождество \eqref{f(x+T)=f(x)}
эквивалентно тождеству
$$
x+T=x,\qquad x\in\R
$$
которое выполняется только если $T=0$. То есть периода (числа $T>0$,
удовлетворяющего \eqref{f(x+T)=f(x)}) у этой функции не существует.
 \eex

\bex{\bf Константы.} Всякую функцию, тождественно равную какому-нибудь числу
$C\in\R$,
 $$
f(x)=C,\qquad x\in\R
 $$
можно считать периодической, потому что для нее любое число $T>0$ будет
периодом:
$$
f(x+T)=C=f(x),\qquad x\in\R
$$
 \eex

\bex{\bf Дробная часть.} Функция $x\mapsto\{x\}$, сопоставляющая числу $x$ его
дробную часть $\{x\}$ (мы определили ее формулой
\eqref{drobnaya-chast-chisla}), является периодической, потому что, например,
число 1 является для нее периодом:
$$
\{x+1\}=\{x\},\qquad x\in\R
$$
Можно заметить, что вообще любое натуральное число $n\in\N$ является периодом
для $x\mapsto\{x\}$:
$$
\{x+n\}=\{x\},\qquad x\in\R
$$
Других же периодов здесь нет:
$$
\Big\{T>0: \quad \forall x\in\R\quad \{x+n\}=\{x\}\ \Big\}=\N
$$
Как следствие, эта функция обладает важным свойством, которое встречается не у
всякой периодической функции (например, у констант его нет): у нее есть
наименьший период -- число 1.
 \eex

\bex{\bf Функция Дирихле} $x\mapsto D(x)$, определенная в главе \ref{ch-R&N}
формулой \eqref{func-Dirichle}, также будет периодической, потому что любое
положительное рациональное число $r\in\Q$ является для нее периодом
 $$
D(x+r)=D(x),\qquad x\in\R
 $$
Однако среди таких чисел $r$ не существует наименьшего. Поэтому $D$, будучи
периодической, не обладает наименьшим периодом.
 \eex

\end{multicols}\noindent\rule[10pt]{160mm}{0.1pt}

\btm Если функция $f$ является периодической и обладает наименьшим периодом
$M$, то любой ее период $T$ имеет вид
$$
T=n\cdot M
$$
где $n\in\N$.
 \etm
 \bpr
Разделим $T$ на $M$ с остатком по теореме \ref{TH:delenie-s-ostatkom}:
$$
T=n\cdot M+r,\qquad 0\le r<M
$$
Если бы оказалось, что $r>0$, то мы получили бы, что $r=T-n\cdot M$ -- период
функции $f$, меньший $M$. То есть $M$ в этом случае не могло бы быть наименьшим
периодом. Значит $r=0$, и это то, что нам нужно.
 \epr

\btm\label{TH:period-f-bez-naim-perioda} Если функция $f$ является
периодической, но не обладает наименьшим периодом $M$, то
 \bit{
\item[(i)] нижняя грань всех ее периодов равна нулю:
 \beq\label{inf-T>0=0}
\inf\Big\{T>0:\quad \forall x\in\D(f)\quad f(x+T)=f(x)\Big\}=0
 \eeq

\item[(ii)] на любом непустом интервале $I=(a,b)$ функция $f$ принимает все
свои значения:
 \beq
f\Big(I\cap\D(f)\Big)=f\Big(\D(f)\Big)
 \eeq
(иными словами, каковы бы ни были $a$ и $b$, $a<b$, для любой точки $x\in\D(f)$
найдется точка $t\in (a,b)\cap\D(f)$ такая, что $f(t)=f(x)$).
 }\eit
\etm
 \bpr
1. Обозначим буквой $P$ множество всех периодов для $f$, а его нижнюю грань
буквой $\alpha$:
$$
P=\Big\{T>0:\quad \forall x\in\D(f)\quad f(x+T)=f(x)\Big\},\qquad \alpha=\inf P
$$
Если $\alpha\notin P$, то для всякого $\beta>\alpha$ найдется $T\in P$ такое,
что
$$
\alpha<T<\beta
$$
Возьмем произвольное $\e>0$. Для него существует какое-то $T_1\in P$ со
свойством
$$
\alpha<T_1<\alpha+\e
$$
Поскольку $T_1>\alpha$, можно подобрать $T_2\in P$ такое, что
$$
\alpha<T_2<T_1<\alpha+\e
$$
Числа $T_1$ и $T_2$ являются периодами для $f$. Поэтому их разность $T=T_1-T_2$
(которая больше нуля, потому что $T_1>T_2$) тоже будет периодом для $f$:
$$
f(x+T)=f(x+T_1-T_2)=f(x+T_1)=f(x)
$$
С другой стороны, $T_1$ и $T_2$ лежат в интервале $(\alpha,\alpha+\e)$ длины
$\e$, поэтому их разность должна быть меньше $\e$:
$$
T=T_1-T_2<\e
$$
Мы получили, что для всякого $\e>0$ найдется период $T\in P$, меньший $\e$.
Значит, $\alpha=\inf P=0$

2. Рассмотрим сначала интервал $I=(-\e,\e)$, $\e>0$. В силу уже доказанной
формулы \eqref{inf-T>0=0}, найдется период $T<\e$. Возьмем произвольную точку
$x\in\D(f)$, и обозначим через $n$ целую часть числа $\frac{x}{T}$:
$$
n=\left[\frac{x}{T}\right]
$$
Тогда точка $t=x-nT$ будет обладать нужными свойствами:
$$
t\in(-\e,\e),\qquad t\in\D(f),\qquad f(t)=f(x)
$$
Действительно, последние два условия --  $t\in\D(f)$ и $f(t)=f(x)$ --
выполняются потому что $T$ -- период. А первое следует из свойств целой части:
$$
n=\left[\frac{x}{T}\right]
$$
$$
\phantom{\text{\scriptsize\eqref{opr-tsel-chasti}}}\Downarrow\quad
\text{\scriptsize\eqref{opr-tsel-chasti}}
$$
$$
n\le \frac{x}{T}<n+1<n+\frac{\e}{T}
$$
$$
\Downarrow
$$
$$
nT\le x<nT+\e
$$
$$
\Downarrow
$$
$$
0\le \underbrace{x-nT}_{t}<\e
$$
$$
\Downarrow
$$
$$
t\in(-\e,\e)
$$

3. Пусть, наконец, $I=(a,b)$ -- произвольный интервал. Положим
$\e=\frac{b-a}{3}$ и, пользуясь формулой \eqref{inf-T>0=0}, подберем период $T$
так, чтобы
$$
\frac{\e}{T}>1
$$
Тогда $a+\e<b-\e$, поэтому $\frac{a+\e}{T}<\frac{b-\e}{T}$, причем расстояние
между точками $\frac{a+\e}{T}$ и $\frac{b-\e}{T}$ больше единицы:
$$
\frac{b-\e}{T}-\frac{a+\e}{T}=\frac{(b-\e)-(a+\e)}{T}=\frac{(b-a)-2\e}{T}=\frac{3\e-2\e}{T}=\frac{\e}{T}>1
$$
Значит, между ними лежит хотя бы одно целое число $n\in\Z$:
$$
\frac{a+\e}{T}<n<\frac{b-\e}{T}
$$
Умножив на $T$ мы получим:
$$
a+\e<nT<b-\e
$$
$$
\Downarrow
$$
$$
a<nT-\e\quad\&\quad nT+\e<b
$$
$$
\Downarrow
$$
$$
(-\e;\e)+nT=(nT-\e;nT+\e)\subseteq(a;b)
$$
То есть $(-\e;\e)$ -- такой интервал, сдвиг которого на число $nT$, кратное
периоду $T$, попадает в интервал $(a;b)$. Мы уже доказали, что для любой точки
$x\in\D(f)$ в интервале $(-\e;\e)$ найдется точка $t$, на которой функция $f$
определена и принимает то же значение:
$$
f(t)=f(x)
$$
Положив $s=t+nT$, мы получим точку из интервала $(nT-\e;nT+\e)\subseteq(a;b)$ с
теми же свойствами:
$$
f(s)=f(t+nT)=f(t)=f(x)
$$
 \epr

\subsection{Свойства функций, связанные с отношением
порядка}

\paragraph{Ограниченные функции и точная грань функции на множестве.} \label{SEC-tochnaya-gran-funktsii}

Пусть $f$ -- числовая функция и $E$ -- подмножество в ее области определения:
$$
E\subseteq \D(f)
$$
Напомним (см. определение на с.\pageref{DF:obraz-mnozhestva}), что {\it
образом} этого множества под действием функции $f$ называется множество
$$
f(E)=\Big\{ f(x);\quad x\in E \Big\}
$$

 \bit{
\item[$\bullet$] Точная нижняя грань множества $f(E)$ называется {\it точной
нижней гранью функции $f$ на множестве $E$} и обозначается $\inf_{x\in E}
f(x)$:
 \beq
\inf_{x\in E} f(x)=\inf f(E)=\inf \Big\{ f(x);\quad x\in E \Big\}
 \eeq
Если функция $f$ {\it ограничена снизу} на множестве $E$, то есть существует
число $A$ такое, что
$$
\forall x\in E \quad A\le f(x)
$$
(это эквивалентно тому, что множество $f(E)$ ограничено снизу), то точная
нижняя грань $f$ на $E$ совпадает с наибольшим из таких чисел $A$:
$$
\inf_{x\in E} f(x)=\max\{A\in\R:\quad \forall x\in E \quad A\le f(x) \}.
$$
В таких случаях $\inf_{x\in E} f(x)$ само является числом (а не символом
$-\infty$, как тоже иногда бывает), и это коротко записывается формулой:
$$
\inf_{x\in E} f(x)>-\infty
$$
В противном случае говорят, что функция $f$ {\it не ограничена снизу} на
множестве $E$, и записывается это формулой:
$$
\inf_{x\in E} f(x)=-\infty
$$

\item[$\bullet$] Точная верхняя грань множества $f(E)$ называется {\it точной
верхней гранью функции $f$ на множестве $E$} и обозначается $\sup_{x\in E}
f(x)$:
 \beq
\sup_{x\in E} f(x)=\sup f(E)=\sup \Big\{ f(x);\quad x\in E \Big\}
 \eeq
Если функция $f$ {\it ограничена сверху} на множестве $E$, то есть существует
число $B$ такое, что
$$
\forall x\in E \quad f(x)\le B
$$
(это эквивалентно тому, что множество $f(E)$ ограничено сверху), то точная
верхняя грань $f$ на $E$ совпадает с наименьшим из таких чисел $B$:
$$
\sup_{x\in E} f(x)=\min\{B\in\R:\quad \forall x\in E \quad f(x)\le B \}.
$$
В таких случаях $\sup_{x\in E} f(x)$ является числом (а не символом $+\infty$),
и это коротко записывается формулой:
$$
\sup_{x\in E} f(x)<+\infty
$$
В противном случае говорят, что функция $f$ {\it не ограничена сверху} на
множестве $E$, и записывается это формулой:
$$
\sup_{x\in E} f(x)=+\infty
$$

\item[$\bullet$] Числовая функция $f$ называется {\it ограниченной} на
множестве $E\subseteq\D(f)$, если она ограничена на $E$ сверху и снизу, то есть
существуют числа $A$ и $B$, такие что:
$$
\forall x\in E \quad A\le f(x)\le B
$$

%\picture{0pt}{0pt}{52.pcx}
\vglue40pt
 }\eit

 \bigskip
\centerline{\bf Свойства точной верхней и точной нижней грани функции:}
 \bit{
\item[$1^\circ$] {\bf Точная грань от константы:}
$$
\Big(\forall x\in E\quad f(x)=C\Big)\quad\Longrightarrow\quad \inf\limits_{x\in
E} f(x)=C=\sup\limits_{x\in E} f(x)
$$

\item[$2^\circ$] {\bf Монотонность:}
$$
\Big(\forall x\in E\quad f(x)\le g(x)\Big) \quad \Longrightarrow \quad
\inf_{x\in E} f(x)\le \inf_{x\in E} g(x) \quad \& \quad \sup_{x\in E} f(x)\le
\sup_{x\in E} g(x)
$$

\item[$3^\circ$] {\bf Антиоднородность:}
 \beq\label{antiodnorodnost-inf-i-sup}
 \begin{split}
& \inf_{x\in E} \Big( -f(x)\Big)= -\sup_{x\in E} f(x)
 \\
& \sup_{x\in E} \Big( -f(x)\Big)=-\inf_{x\in E} f(x)
 \end{split}
 \eeq

\item[$4^\circ$] {\bf Полуаддитивность:}
 \beq \label{poluaddit-sup}
 \begin{split}
 &
\inf_{x\in E} \Big( f(x)+g(x)\Big)\ge
 \inf_{x\in E} f(x)+\inf_{x\in E} g(x)
 \\
 &
\sup_{x\in E} \Big( f(x)+g(x)\Big)\le
 \sup_{x\in E} f(x)+\sup_{x\in E} g(x)
 \end{split}
 \eeq

\item[$5^\circ$] {\bf Полумультипликативность:} если $f(x),g(x)\ge 0$, то
 \beq \label{polumultiplik-sup}
  \begin{split}
& \inf_{x\in E} \Big( f(x)\cdot g(x)\Big)\ge
 \inf_{x\in E}f(x)\cdot \inf_{x\in E} g(x)
 \\
& \sup_{x\in E} \Big( f(x)\cdot g(x)\Big)\le
 \sup_{x\in E} f(x)\cdot \sup_{x\in E} g(x)
  \end{split}
 \eeq

\item[$6^\circ$] {\bf Связь с модулем:}\footnote{Это тождество понадобится нам
в \ref{SEC-svojstva-opredelyonnogo-integrala} главы
\ref{CH-definite-integral}.}
 \begin{equation} \label{sup_x,y in X |f(x)-f(y)| sup_x,y in X (f(x)-f(y))}
 \sup_{x,y\in E} \Big|f(x)-f(y)\Big|=\sup_{x,y\in E} \Big(f(x)-f(y)\Big)
 \end{equation}
 }\eit

\begin{proof}
В каждом пункте (за исключением $6^\circ$) мы докажем только первую половину
утверждения (вторая доказывается по аналогии).

1.
 \begin{multline*}
\Big(\forall x\in E\quad f(x)=C\Big)\quad\Longrightarrow\quad \{A:\;\; \forall
x\in E\quad A\le f(x)\}=\{A:\;\; A\le C\}=(-\infty; C]\quad\Longrightarrow \\
\Longrightarrow\quad \inf_{x\in E}f(x)=\max\{A:\;\; \forall x\in E\quad A\le
f(x)\}=\max(-\infty; C]=C
 \end{multline*}

2. Пусть $\forall x\in E\quad f(x)\le g(x)$. Тогда
 \begin{align*}
\{C:\;\;\forall x\in E\;\; C\le f(x)\}&\subseteq \{C:\;\;\forall x\in E\;\;
C\le
g(x)\} \\
&\Downarrow\qquad (\text{\scriptsize свойство $1^\circ$ на стр.
\pageref{svoistvo-1-inf-sup}})
\\
\inf_{x\in E}f(x)=\max\{C:\;\;\forall x\in E\;\; C\le f(x)\}&\le
\max\{C:\;\;\forall x\in E\;\; C\le g(x)\}=\inf_{x\in E}g(x)
 \end{align*}

3.
$$
\inf_{x\in E}\Big(-f(x)\Big)=\inf\Big(-f(E)\Big)=(\text{\scriptsize свойство
$3^\circ$ на стр. \pageref{inf-X+a=inf-X+a}})=-\sup f(E)=-\sup_{x\in E} f(x)
$$
.

4.
 \begin{multline*}
\begin{cases} \forall x\in E\quad f(x)\ge \inf\limits_{x\in E} f(x) \\
\forall x\in E\quad g(x)\ge \inf\limits_{x\in E} g(x)
\end{cases}\quad\Longrightarrow\quad
\forall x\in E\quad f(x)+g(x)\ge \inf\limits_{x\in E} f(x)+\inf\limits_{x\in E}
g(x) \quad\Longrightarrow \\ \Longrightarrow \quad \inf\limits_{x\in
E}\Big(f(x)+g(x)\Big)\ge \inf\limits_{x\in E} f(x)+\inf\limits_{x\in E}g(x)
 \end{multline*}

5.
 \begin{multline*}
\begin{cases} \forall x\in E\quad f(x)\ge \inf\limits_{x\in E} f(x) \ge 0\\
\forall x\in E\quad g(x)\ge \inf\limits_{x\in E} g(x)\ge 0
\end{cases}\quad\Longrightarrow\quad
\forall x\in E\quad f(x)\cdot g(x)\ge \inf\limits_{x\in E} f(x)\cdot
\inf\limits_{x\in E} g(x) \quad\Longrightarrow \\ \Longrightarrow \quad
\inf\limits_{x\in E}\Big(f(x)\cdot g(x)\Big)\ge \inf\limits_{x\in E} f(x)\cdot
\inf\limits_{x\in E}g(x)
 \end{multline*}

6. Во-первых, нужно заметить неравенство:
$$
 \sup_{x,y\in E} \Big(f(x)-f(y)\Big)\le \sup_{x,y\in E} \Big|f(x)-f(y)\Big|
$$
Оно следует из свойства монотонности:
$$
 \forall x,y\in E \qquad \Big(f(x)-f(y)\Big)\le \Big|f(x)-f(y)\Big|
 \quad \Longrightarrow \quad
 \sup_{x,y\in E} \Big(f(x)-f(y)\Big)\le \sup_{x,y\in E} \Big|f(x)-f(y)\Big|
$$
После этого нужно убедиться, что это неравенство не может быть строгим.
Предположим, что
$$
 \sup_{x,y\in E} \Big(f(x)-f(y)\Big)<\sup_{x,y\in E} \Big|f(x)-f(y)\Big|
$$
то есть что для каких-то элементов $a,b\in E$ выполняется
$$
 \forall x,y\in E\quad \Big(f(x)-f(y)\Big)<\Big|f(a)-f(b)\Big|
$$
В частности, взяв $x=a,y=b$, мы получили бы
$$
f(a)-f(b)<\Big|f(a)-f(b)\Big|
$$
а, взяв $x=b,y=a$, получили бы
$$
f(b)-f(a)<\Big|f(a)-f(b)\Big|
$$
Вместе эти два неравенства выполняться не могут, потому что
$\Big|f(a)-f(b)\Big|$ равен либо $f(a)-f(b)$, либо $f(b)-f(a)$.
\end{proof}

\paragraph{Монотонные функции}

 \bit{
\item[$\bullet$] Числовая функция $f$ называется
 \bit{
\item[--] {\it возрастающей} на множестве $E\subseteq \D(f)$, если для любых
точек $x,y\in E$ таких что $x<y$ выполняется неравенство $f(x)<f(y)$;

%\picture{0pt}{0pt}{46.pcx}
\vglue100pt

\item[--] {\it убывающей} на множестве $E\subseteq \D(f)$, если для любых точек
$x,y\in E$ таких что $x < y$ выполняется неравенство $f(x) > f(y)$;

%\picture{0pt}{0pt}{47.pcx}
\vglue100pt

\item[--] {\it невозрастающей} на множестве $E\subseteq \D(f)$, если для любых
точек $x,y\in E$ таких что $x < y$ выполняется неравенство $f(x) \ge f(y)$;

%\picture{0pt}{0pt}{48.pcx}
\vglue100pt

\item[--] {\it неубывающей} на множестве $E\subseteq \D(f)$, если для любых
точек $x,y\in E$ таких что $x < y$ выполняется неравенство $f(x) \le f(y)$;

%\picture{0pt}{0pt}{49.pcx}
\vglue100pt

\item[--] {\it монотонной} на множестве $E\subseteq \D(f)$, если она
невозрастающая или неубывающая;

\item[--] {\it строго монотонной} на множестве $E\subseteq \D(f)$, если она
возрастающая или убывающая.

 }\eit
 }\eit

Доказательство следующих свойств мы оставляем читателю:

 \bigskip

\centerline{\bf Свойства монотонных функций}

 \bit{\it
\item[$1^\circ$]\label{vozrastanie=>neubyvanie} Если функция $f$ возрастает на
множестве $E$, то она неубывает
 на множестве $E$.

\item[$2^\circ$] Если функция $f$ убывает  на множестве $E$, то она
невозрастает  на множестве $E$.

\item[$3^\circ$] Если функция $f$ неубывает на множестве $E$, то на всяком
ограниченном снизу (сверху) подмножестве $M\subseteq E$ таком, что $\inf M\in
E$ ($\sup M\in E$) функция $f$ ограничена снизу (сверху) и удовлетворяет
неравенству
 \begin{align}\label{f(inf-E)<inf-f(E)}
& f(\inf M)\le \inf_{x\in M} f(x) && \Big( \sup_{x\in M} f(x)\le f(\sup M)\Big)
 \end{align}

\item[$4^\circ$] Если функция $f$ невозрастает на множестве $E$, то на всяком
ограниченном сверху (снизу) подмножестве $M\subseteq E$ таком, что $\sup M\in
E$ ($\inf M\in E$) функция $f$ ограничена снизу (сверху) и удовлетворяет
неравенству
 \begin{align}\label{f(sup-E)<inf-f(E)}
& f(\sup M)\le \inf_{x\in M} f(x) && \Big( \sup_{x\in M} f(x)\le f(\inf M)\Big)
 \end{align}

\item[$5^\circ$] Композиция $g\circ f$ функций $f$ и $g$
 \bit{
\item[--] возрастает (неубывает), если обе функции $f$ и $g$ возрастают
(неубывают), либо обе убывают (невозрастают),

\item[--] убывает (невозрастает), если одна из функций $f$ и $g$ возрастает
(неубывает), а другая убывает (невозрастает).
 }\eit
 }\eit

\noindent\rule{160mm}{0.1pt}\begin{multicols}{2}

\bex\label{EX:x^(2n)-vozrastaet-x>0} Степенная функция $x\mapsto x^{2n}$,
$n\in\N$,
 \biter{
\item[---] убывает на полуинтервале $(-\infty;0]$,

\item[---] возрастает на полуинтервале $[0;+\infty)$.
 }\eiter
Как следствие,
 \beq\label{x-ne-0=>x^(2n)>0}
\boxed{\quad x\ne 0\quad\Longrightarrow\quad x^{2n}>0 \quad}
 \eeq
\eex
 \bpr
Возрастание на полуинтервале $[0;+\infty)$ следует из теоремы
\ref{TH-o-step-otobr-malaya}: в силу \eqref{monot-a^n-n>0}, функция $x\mapsto
x^{2n}$ возрастает на интервале $(0;+\infty)$:
$$
0<x<y\quad\Rightarrow\quad x^{2n}<y^{2n}
$$
А случай $x=0$ получается из условия сохранения знака \eqref{x^n>0}:
$$
0<y\quad\Rightarrow\quad 0^{2n}\overset{\eqref{0^n}}{=}0<y^{2n}
$$
С другой стороны, в силу \eqref{(-x)^(2n)=x^(2n)}, наша функция четная:
$$
(-x)^{2n}=x^{2n}
$$
Отсюда следует, что она убывает на полуинтервале $(-\infty;0]$:
$$
x<y\le 0
$$
$$
\Downarrow
$$
$$
0\le-y<-x
$$
$$
\Downarrow
$$
$$
\underbrace{(-y)^{2n}}_{y^{2n}}<\underbrace{(-x)^{2n}}_{x^{2n}}
$$
$$
\Downarrow
$$
$$
y^{2n}<x^{2n}
$$
Для доказательства \eqref{x-ne-0=>x^(2n)>0} нужно просто рассмотреть два
случая: $x>0$ и $x<0$.
 \epr

\bex\label{EX:x^(2n-1)-vozrastaet} Степенная функция $x\mapsto x^{2n-1}$,
$n\in\N$, возрастает всюду на своей области определения (то есть на прямой
$\R$). Как следствие,
 \beq\label{x>0=>x^(2n-1)>0}
\boxed{\quad\begin{split}& x>0\quad\Longrightarrow\quad x^{2n-1}>0 \\
& x<0\quad\Longrightarrow\quad x^{2n-1}<0 \end{split}\quad}
 \eeq
\eex
 \bpr
Возрастание на полуинтервале $[0;+\infty)$ здесь доказывается так же, как в
предыдущем примере \ref{EX:x^(2n)-vozrastaet-x>0}. Поскольку вдобавок, в силу
\eqref{(-x)^(2n-1)=-x^(2n-1)}, наша функция нечетная,
$$
(-x)^{2n-1}=-x^{2n-1},
$$
она должна возрастать и на полуинтервале $(-\infty;0]$:
$$
x<y\le 0
$$
$$
\Downarrow
$$
$$
0\le-y<-x
$$
$$
\Downarrow
$$
$$
\underbrace{(-y)^{2n-1}}_{-y^{2n-1}}<\underbrace{(-x)^{2n-1}}_{-x^{2n-1}}
$$
$$
\Downarrow
$$
$$
-y^{2n-1}<-x^{2n-1}
$$
$$
\Downarrow
$$
$$
x^{2n-1}<y^{2n-1}
$$
Остается заметить, что значения функции на полуинтервале $(-\infty;0]$ меньше
значений на $[0;+\infty)$. Это следует из того, что эти интервалы имеют общую
точку $0$:
$$
x<0<y
$$
$$
\Downarrow
$$
$$
\underbrace{\underbrace{x^{2n-1}<0}_{\scriptsize\begin{matrix}\text{потому что}\\
x,0\in(-\infty;0]\end{matrix}} \quad \&\quad \underbrace{0<y^{2n-1}}_{\scriptsize\begin{matrix}\text{потому что}\\
0,y\in[0;+\infty)\end{matrix}}}
$$
$$
\Downarrow
$$
 \beq\label{x^(2n-1)<0<y^(2n-1)}
\kern40pt x^{2n-1}<0<y^{2n-1}
 \eeq
Формулы \eqref{x>0=>x^(2n-1)>0} можно вывести из возрастания функции $x\mapsto
x^{2n-1}$, а можно объявить их следствием уже доказанного неравенства
\eqref{x^(2n-1)<0<y^(2n-1)}.
 \epr

\bex\label{EX:monotonnost-x^(-2n)} Степенная функция $x\mapsto x^{-2n}$,
$n\in\N$,
 \biter{
\item[---] возрастает на интервале $(-\infty;0)$,

\item[---] убывает на интервале $(0;+\infty)$.
 }\eiter
При этом,
 \beq\label{x-ne-0=>x^(-2n)>0}
\boxed{\quad x\ne 0\quad\Longrightarrow\quad x^{-2n}>0 \quad}
 \eeq
\eex
 \bpr
Убывание на интервале $(0;+\infty)$ следует из условия монотонности
\eqref{monot-a^n-n<0} теоремы \ref{TH-o-step-otobr-malaya}. Из него и условия
четности \eqref{(-x)^(2n)=x^(2n)} выводится возрастание на интервале
$(-\infty;0)$ по аналогии с тем, как это делалось в примере
\ref{EX:x^(2n)-vozrastaet-x>0}. Свойства \eqref{x-ne-0=>x^(-2n)>0} следуют из
из условия сохранения знака \eqref{x^n>0} теоремы \ref{TH-o-step-otobr-malaya}
и условия четности \eqref{(-x)^(2n)=x^(2n)}.
 \epr

\bex\label{EX:monotonnost-x^(-(2n-1))} Степенная функция $x\mapsto
x^{-(2n-1)}$, $n\in\N$,
 \biter{
\item[---] убывает на интервале $(-\infty;0)$,

\item[---] убывает на интервале $(0;+\infty)$.
 }\eiter
При этом,
 \beq\label{x>0=>x^(-(2n-1))>0}
\boxed{\quad\begin{split}& x>0\quad\Longrightarrow\quad x^{-(2n-1)}>0 \\
& x<0\quad\Longrightarrow\quad x^{-(2n-1)}<0 \end{split}\quad}
 \eeq
\eex
 \bpr
Как и в предыдущем примере убывание на интервале $(0;+\infty)$ следует из
условия монотонности \eqref{monot-a^n-n<0} теоремы
\ref{TH-o-step-otobr-malaya}. Из него и условия нечетности
\eqref{(-x)^(2n-1)=-x^(2n-1)} выводится убывание на интервале $(-\infty;0)$ по
аналогии с тем, как это делалось в примере \ref{EX:x^(2n-1)-vozrastaet}.
Свойства \eqref{x>0=>x^(-(2n-1))>0} следуют из из условия сохранения знака
\eqref{x^n>0} теоремы \ref{TH-o-step-otobr-malaya} и условия нечетности
\eqref{(-x)^(2n-1)=-x^(2n-1)}.
 \epr

\begin{er}
Нарисуйте график какой-нибудь функции $f:\R\to \R$ со следующими свойствами:
 \biter{
\item[1)] на интервале $(-\infty;-1)$ функция $f$ убывает и ограничена;
\item[2)] на интервале $(-1;0)$ функция $f$ возрастает и ограничена; \item[3)]
на интервале $(0;1)$ функция $f$ убывает, ограничена снизу, но не ограничена
сверху; \item[4)] на интервале $(1;+\infty)$ функция $f$ возрастает и
ограничена;
 }\eiter
Понятно, что эта задача имеет неединственное решение. В качестве возможного
ответа мы предлагаем следующую картинку:

%\picture{0pt}{0pt}{53.pcx}
\vglue100pt
\end{er}

\begin{er}
Нарисуйте график какой-нибудь функции $f:\R\to \R$ со следующими свойствами:
 \biter{
\item[1)] на интервале $(-\infty;-2)$ функция $f$ возрастает и ограничена;
\item[2)] на интервале $(-2;0)$ функция $f$ убывает, ограничена сверху, но не
ограничена снизу; \item[3)] на интервале $(0;1)$ функция $f$ возрастает и не
ограничена ни сверху, ни снизу; \item[4)] на интервале $(1;+\infty)$ функция
$f$ возрастает и ограничена снизу.
 }\eiter
Один из возможных ответов выглядит так:

%\picture{0pt}{0pt}{54.pcx}
\vglue100pt
\end{er}

\end{multicols}\noindent\rule[10pt]{160mm}{0.1pt}

\section{Непрерывные функции}\label{SEC-nepr-funct}

\subsection{Что такое непрерывная функция?}\label{SEC-nepreryvnost-functsii}

Пусть нам дано уравнение, связывающее две физические величины,
$$
  y=f(x),
$$
и мы говорим, что $y$ непрерывно зависит от $x$. Что это означает? Интуитивный
ответ очевиден: непрерывная функция -- это такая, у которой график --
непрерывная линия.

\noindent\rule{160mm}{0.1pt}\begin{multicols}{2}

\bex Например, функция
$$
f(x)=x^2
$$
непрерывна, потому что у нее график -- непрерывная линия

%\picture{50pt}{0pt}{x^2.pcx}

\vglue180pt \eex

\bex И наоборот, функция {\it сигнум}\label{razryvnost-signuma}
$$
\sgn x=\begin{cases}1, & \text{если}\,\, x>0
 \\
 0, & \text{если}\,\,x=0
 \\
 -1, & \text{если}\,\, x<0
 \end{cases}
$$
разрывна, потому что у нее график имеет разрыв:

%\picture{100pt}{0pt}{sgn.pcx}

\vglue120pt

\eex

\end{multicols}\noindent\rule[10pt]{160mm}{0.1pt}

Однако, это соображение не может быть точным определением, потому что оно
объясняет смысл одного нового понятия (непрерывная функция) с помощью другого
нового понятия (непрерывная линия).

Как же быть? Оказывается, имеется очень простой способ объяснить, что такое
непрерывная функция, опираясь на известное нам уже (то есть, старое для нас)
понятие предела последовательности. Это делается с помощью следующей серии
определений.

 \bit{
\item[$\bullet$] Функция $f$, определенная на множестве $E\subseteq\R$,
называется
 \biter{
\item[--] {\it непрерывной в точке $a\in E$ на множестве
$E$}\index{функция!непрерывная!в точке на множестве}, если $f$ определена на
множестве $E$, и для всякой последовательности точек $x_n \in E$, сходящейся к
точке $a$,
$$
x_n\underset{n\to \infty}{\longrightarrow} a \qquad \Big(
x_n\in E \Big)
$$
соответствующая последовательность $f(x_n)$ значений функции $f$ стремится к
$f(a)$:
$$
f(x_n)\underset{n\to \infty}{\longrightarrow} f(a)
$$

%\picture{0pt}{0pt}{p5.pcx}

\vglue100pt \noindent (при этом точка $a$ называется {\it
точкой непрерывности функции $f$ на множестве $E$});

\item[--] {\it разрывной в точке $a$ на множестве $E$}, если $f$ определена на
множестве $E$, и найдется последовательность
$$
x_n\underset{n\to \infty}{\longrightarrow} a \qquad \Big(
x_n\in E \Big)
$$
такая, что $f(x_n)$ не стремится к $f(a)$
$$
f(x_n)\underset{n\to \infty}{\notarrow} f(a)
$$

%\picture{0pt}{0pt}{p6.pcx}

\vglue100pt \noindent (при этом точка $a$ называется {\it
точкой разрыва функции $f$ на множестве
$E$})\index{разрыв}\index{точка разрыва}
 }\eiter
 }\eit

 \bit{
\item[$\bullet$] Функция $f$, определенная на множестве $E$, называется
 \bit{
\item[--] {\it непрерывной на множестве $E$}\index{функция!непрерывная!на
множестве}, если она непрерывна в любой точке $a\in E$ на множестве $E$, то
есть для всякой последовательности точек $x_n \in E$, сходящейся к любой точке
$a\in E$,
$$
x_n\underset{n\to \infty}{\longrightarrow} a \qquad \Big(
x_n, a\in E \Big)
$$
соответствующая последовательность $f(x_n)$ значений функции $f$ стремится к
$f(a)$:
$$
f(x_n)\underset{n\to \infty}{\longrightarrow} f(a)
$$
\item[--] {\it разрывной на множестве $E$}, если она разрывна в некоторой точке
$a\in E$ на множестве $E$, то есть существует такая точка $a\in E$ и такая
последовательности точек $x_n \in E$, сходящаяся к $a$,
$$
x_n\underset{n\to \infty}{\longrightarrow} a \qquad \Big(
x_n, a\in E \Big)
$$
что соответствующая последовательность $f(x_n)$ значений функции $f$ не
стремится к $f(a)$:
$$
f(x_n)\underset{n\to \infty}{\notarrow} f(a)
$$
 }\eit

\item[$\bullet$] Функция $f$ называется {\it непрерывной} (соответственно, {\it
разрывной}), если она непрерывна (соответственно, разрывна) на своей области
определения $\D(f)$.
 }\eit

\noindent\rule{160mm}{0.1pt}\begin{multicols}{2}

\bex Функция $f(x)=x^2$ будет непрерывна на множестве $E=\R$, потому что (по
свойству $3^0$ сходящихся последовательностей на
с.\pageref{x_n->x,y_n->y=>x_n-cdot-y_n->x-cdot-y}) из условия
$x_n\underset{n\to\infty}{\longrightarrow} a$  следует условие
$f(x_n)=(x_n)^2\underset{n\to\infty}{\longrightarrow} a^2=f(a)$.\eex

\bex\label{EX:razryvnost-sgn} Наоборот, функция сигнум $f(x)=\sgn x$,
определенная выше формулой \eqref{DEF:sgn}, разрывна на множестве $E=\R$,
точнее, в точке $x=0$, потому что если взять последовательность
$$
x_n=\frac{1}{n}\underset{n\to\infty}{\longrightarrow} 0
$$
то окажется, что
$$
f(x_n)=1\underset{n\to\infty}{\longrightarrow} 1\ne 0=f(0)
$$
Однако, в любой другой точке она непрерывна. Например, в любой точке $a>0$:
если взять $\varepsilon=a>0$ и рассмотреть окрестность
$(0;2a)=(a-\varepsilon;a+\varepsilon)$, то для всякой последовательности
$$
x_n\underset{n\to \infty}{\longrightarrow} a
$$
мы получим, что почти все числа $x_n$ лежат в интервале
$(0;2a)$, значит, для почти всех $n\in \mathbb{N}$
выполняется неравенство
$$
x_n>0
$$
поэтому для почти всех $n\in \mathbb{N}$
$$
f(x_n)=\sgn x_n=1
$$
и следовательно,
$$
f(x_n)\underset{n\to \infty}{\longrightarrow} 1=f(a)
$$
Аналогично доказывается, что $f(x)=\sgn x$ непрерывна в
любой точке $a<0$.
\end{ex}

\ber Покажите, что функция
$$
f(x)=\sgn^2 x=\begin{cases}
 1, & \text{если}\,\, x\ne 0
 \\
 0, & \text{если}\,\,x=0
\end{cases}
$$
также разрывна в точке $0$, но непрерывна в любой точке $a\ne 0$. \eer

\bex\label{EX:drob-chast-razryvna} Покажем, что функция $f(x)=\{x\}$ (дробная
часть числа) разрывна в точке 1. Действительно, для последовательности
$$
x_n=1-\frac{1}{n}\underset{n\to \infty}{\longrightarrow} 1
$$
мы получим
$$
\{x_n\}=\eqref{drob-chast-0-1}=x_n\underset{n\to \infty}{\longrightarrow} 1,
$$
то есть
$$
\{x_n\}\underset{n\to \infty}{\notarrow} 0=\{1\}.
$$
В силу периодичности функции $f(x)=\{x\}$, отсюда следует, что она разрывна в
любой точке $a\in\Z$. С другой стороны, нетрудно убедиться, что она непрерывна
в любой точке $a\notin\Z$. \eex

\ber В каких точках непрерывна функция $f(x)=[x]$ (целая часть числа)? \eer

\begin{er}
Для следующих функций найдите точки разрыва (если они есть) и докажите, что в
этих точках функции действительно разрывны:
$$
f(x)=\begin{cases} 1-x, & \text{если}\; x>0
 \\ -1-x, & \text{если}\; x<0\\
 0, & \text{если}\; x=0
\end{cases},
$$
$$
g(x)=\begin{cases}
 0, & \text{если}\; x<0
 \\
 x^2, & \text{если}\;0\le x\le 1\\
 0, & \text{если}\; x>1
\end{cases},
$$
$$
h(x)=\begin{cases}
 \frac{\pi}{2}, & \text{если}\,\, x\in (-\infty;-1)
 \\
 \arcsin x, & \text{если}\,\,x\in [-1;1]
\end{cases}
$$
\end{er}

\bex{\bf Непрерывность модуля} Функция $f(x)=|x|$ непрерывна. \eex
\begin{proof}
Для любой точки $a\in \R$ и любой последовательности $x_n\underset{n\to
\infty}{\longrightarrow} a$ получаем цепочку следствий:
 \begin{gather*}
x_n \underset{n\to \infty}{\longrightarrow} a \\
 \Downarrow \put(20,0){\text{\smsize
 $\begin{pmatrix}\text{применяем}\\ \text{теорему
 \ref{x->a<=>x-a->0}}\end{pmatrix}$}}\\
x_n-a \underset{n\to \infty}{\longrightarrow} 0 \\
\Downarrow \put(20,0){\text{\smsize
 $\begin{pmatrix}\text{вспоминаем}\\ \text{определение}
 \\ \text{бесконечно малой}\\
\text{посл-ти в \ref{->0-&-->infty}}\end{pmatrix}$}}\\
|x_n-a| \underset{n\to \infty}{\longrightarrow} 0  \\
\Downarrow \put(20,0){\text{\smsize
 $\begin{pmatrix}\text{применяем}\\
 \text{свойство модуля $4^0$}\\
 \text{из $\S 6$ главы \ref{ch-R&N}}\end{pmatrix}$}}\\
\underbrace{-|x_n-a|}_{\scriptsize\begin{matrix}\downarrow\\ 0 \end{matrix}}\le
|x_n| -
|a|\le \underbrace{|x_n-a|}_{\scriptsize\begin{matrix}\downarrow\\ 0 \end{matrix}}  \\
\Downarrow \put(20,0){\text{\smsize
  $\begin{pmatrix}\text{применяем}\\ \text{теорему о двух}\\
  \text{милиционерах
\ref{milit}}\end{pmatrix}$}}\\
|x_n|-|a| \underset{n\to \infty}{\longrightarrow} 0
 \end{gather*}
Мы получили, что если $x_n\underset{n\to \infty}{\longrightarrow} a$, то
$|x_n|\underset{n\to \infty}{\longrightarrow} |a|$. Это означает, что функция
$f(x)=|x|$ непрерывна в любой точке $a\in \R$. Таким образом,  $f(x)=|x|$ --
непрерывная функция.
\end{proof}

\bex{\bf Непрерывность линейной функции.} Функция вида
\beq\label{DEF:lineynaya-functsiya}
f(x)=k x+b,\qquad k,b\in\R
\eeq
называется линейной. Покажем, что она непрерывна.
 \eex
\begin{proof} Пусть
$$
x_n\underset{n\to \infty}{\longrightarrow} a
$$
тогда, в силу арифметических свойств пределов последовательностей ($1^0, 3^0$
из \ref{lim-ariphm}), имеем
$$
f(x_n)=k x_n+b\underset{n\to \infty}{\longrightarrow} k a+b=f(a)
$$
то есть $f(x)$ непрерывна в точке $a$ \end{proof}

 \bex\label{nepr-x^n} {\bf Непрерывность степенной функции.}
Функция $f(x)=x^k$ непрерывна при любом $k\in\Z$. В частности, функция
$x\mapsto\frac{1}{x}$ непрерывна.
 \eex
\begin{proof}
Напомним, что при $k\ge 0$ степень $x^k$ определялась индуктивным соотношением
\eqref{def:a^n}, а при $k<0$ -- формулой \eqref{DF:a^n-n<0}. Рассмотрим три
случая.

1. При $k=0$ эта функция будет по определению константой
$$
x^0=1,
$$
поэтому в этом случае она непрерывна тривиальным образом.

2. Пусть $k\in\N$, тогда функция $f(x)=x^k$ определена на всей числовой прямой
$\R$. Зафиксируем точку $a\in\R$ и рассмотрим произвольную последовательность
$$
x_n\underset{n\to \infty}{\longrightarrow} a
$$
Тогда, в силу арифметических свойств пределов последовательностей ($4^0$ на с.
\pageref{lim-ariphm}), имеем
$$
f(x_n)=(x_n)^k\underset{n\to \infty}{\longrightarrow} a^k=f(a)
$$
то есть $f(x)$ непрерывна в точке $a$. Поскольку точка $a\in D(f)$ выбиралась
произвольно, мы получаем, что функция $f(x)=x^k$ непрерывна.

2. Пусть $k=-m$, $m\in \N$, тогда функция $f(x)=x^{-m}$ определена на множестве
$(-\infty;0)\cup (0;+\infty)$. Зафиксируем точку $a\in (-\infty;0)\cup
(0;+\infty)$ и рассмотрим произвольную последовательность
$$
x_n\underset{n\to \infty}{\longrightarrow} a
$$
тогда, в силу арифметических свойств пределов последовательностей ($4^0$ на с.
\pageref{lim-ariphm}), имеем
$$
(x_n)^m\underset{n\to \infty}{\longrightarrow}  a^m
$$
а затем по лемме \ref{1/x} главы \ref{ch-x_n},
$$
f(x_n)=\frac{1}{(x_n)^m}\underset{n\to
\infty}{\longrightarrow}\frac{1}{a^m}=f(a)
$$
то есть $f(x)$ непрерывна в точке $a$. Поскольку точка $a\in D(f)$ выбиралась
произвольно, мы получаем, что функция $f(x)=x^k$ непрерывна.
 \end{proof}

\end{multicols}\noindent\rule[10pt]{160mm}{0.1pt}

Интуитивный смысл непрерывности в том, что график функции $f$ непрерывной на
каком-нибудь интервале $(\alpha;\beta)$ является непрерывной кривой.

\noindent\rule{160mm}{0.1pt}\begin{multicols}{2}

\begin{er} По графику функции

%\picture{0pt}{0pt}{53.pcx}

\vglue50pt

\noindent определите, в каких точках она является непрерывной, и будет ли она
непрерывной
 \biter{
\item[1)] на интервале $(-\infty; 1)$?

\item[2)] на интервале $(-\infty;0)$?

\item[3)] на интервале $(1;+\infty)$?

\item[4)] на интервале $(0;1)$?
 }\eiter
\end{er}

\end{multicols}\noindent\rule[10pt]{160mm}{0.1pt}

\subsection{Непрерывность и монотонные последовательности}

\begin{tm}\label{TH:nepr-monot} Для всякой функции $f$, определенной на
множестве $M$, следующие условия эквивалентны:
 \bit{
\item[(i)] $f$ непрерывна в точке $a$ на $M$;

\item[(ii)] для любой строго монотонной последовательности
$x_n\underset{n\to\infty}{\longrightarrow} a$
 \begin{equation}\label{f(x_n)->f(a)}
f(x_n)\underset{n\to\infty}{\longrightarrow} f(a)
 \end{equation}
 }\eit
 \end{tm}
\begin{proof} Ясно, что из $(i)$ следует $(ii)$, потому что
если $f$ непрерывна в $a$, то $f(x_n)\underset{n\to\infty}{\longrightarrow}
f(a)$ для любой последовательности $x_n\underset{n\to\infty}{\longrightarrow}
a$ (необязательно монотонной). Докажем, что наоборот если $(i)$ не выполняется,
то не выполняется и $(ii)$. Пусть $f$ -- разрывная функция на $M$, то есть
существует последовательность $x_n\in M$ такие что
$$
x_n\underset{n\to\infty}{\longrightarrow} a\qquad \&\qquad f(x_n)\underset{n\to
\infty}{\notarrow} f(a)
$$
План наших действий состоит в том, чтобы, последовательно переходя от $x_n$ к
ее подпоследовательностям, построить строго монотонную последовательность с
теми же свойствами.

1. Прежде всего, по свойству подпоследовательностей $2^\circ$ на
с.\pageref{podposledovatelnosti}, второе условие означает, что существует
$\e>0$ и последовательность натуральных чисел
$n_k\underset{k\to\infty}{\longrightarrow}\infty$, такие что
$$
\forall k\in\N\qquad f(x_{n_k})\notin (f(a)-\e;f(a)+\e)
$$
Обозначив $y_k=x_{n_k}$, мы получим последовательность с такими свойствами:
$$
y_k\underset{k\to\infty}{\longrightarrow} a\qquad \&\qquad \forall k\in\N\qquad
f(y_k)\notin (f(a)-\e;f(a)+\e)
$$

2. Заметим далее, что числа $y_k$ не могут совпадать с $a$
$$
\forall k\in\N\qquad y_k\ne a
$$
потому что иначе мы получили бы, что для некоторых $k$ числа $f(y_k)=f(a)$
лежат в интервале $(f(a)-\e;f(a)+\e)$. Отсюда следует, что либо бесконечный
набор этих чисел лежит левее $a$, либо бесконечный набор этих чисел лежит
правее $a$. Будем для определенности считать, что выполняется первое (второй
случай рассматривается по аналогии):
$$
y_k<a.
$$
Тогда по свойству $4^\circ$ на с.\pageref{podposledovatelnosti}, найдется
последовательность $k_i\underset{i\to\infty}{\longrightarrow}\infty$ такая, что
$$
\forall i\in\N\qquad y_{k_i}<a
$$
Обозначив $z_i=y_{k_i}$, мы получим последовательность с такими свойствами:
$$
z_i\underset{i\to\infty}{\longrightarrow} a\qquad \&\qquad \Big(\forall
i\in\N\qquad z_i<a\Big) \qquad \&\qquad\Big(\forall i\in\N\qquad  f(z_i)\notin
(f(a)-\e;f(a)+\e)\Big)
$$

3. После этого мы индуктивно определим две последовательности $i_m$ и $t_m$.
Сначала полагаем
$$
i_1=1,\qquad t_1=z_1
$$
Затем, если для $m=1,...,l$ числа $i_m$ и $t_m$ определены, мы замечаем вот
что. Поскольку $z_i<a$, $z_i\underset{i\to\infty}{\longrightarrow} a$ и
$t_l<a$, найдется такое $i_{l+1}$, что $z_{i_{l+1}}$ лежит в интервале
$(t_l;a)$:
$$
z_{i_{l+1}}>t_l
$$
Выбираем такое $i_{l+1}$ и полагаем
$$
t_{l+1}=z_{i_{l+1}}
$$
Полученная последовательность $t_m=z_{m+1}$ обладает нужными нам свойствами.
Во-первых, как подпоследовательность $z_i$, она стремится к $a$:
$$
t_m=z_{i_m}\underset{m\to\infty}{\longrightarrow} a
$$
Во-вторых, последовательность значений $f$ на ней не заходит в интервал
$(f(a)-\e;f(a)+\e)$:
$$
\forall m\in\N\qquad  f(t_m)=f(z_{i_m})\notin (f(a)-\e;f(a)+\e)
$$
И, в-третьих, она строго монотонна:
$$
t_m<z_{i_{m+1}}=t_{m+1}
$$
\end{proof}

\subsection{Операции над непрерывными функциями}\label{SUBSEC:oper-nad-nepr-func}

Из непрерывных функций можно конструировать новые
непрерывные функции с помощью алгебраических операций и
операции композиции.

\paragraph{Арифметические операции с непрерывными функциями.}

\bit{

\item[$\bullet$] Если $f$ и $g$ -- две функции, определенные на множестве $X$,
то их {\it суммой, разностью и произведением} называются функции на $X$,
обозначаемые $f+g$, $f-g$ и $f\cdot g$, и определенные формулами
\begin{align}
& (f+g)(x):=f(x)+g(x), && x\in X \label{DEF:f+g}\\
& (f-g)(x):=f(x)-g(x), && x\in X \label{DEF:f-g}\\
& (f\cdot g)(x):=f(x)\cdot g(x), && x\in X \label{DEF:f-cdot-g}
\end{align}

\item[$\bullet$] Если $f$ и $g$ -- две функции, определенные на множестве $X$,
причем $f$ нигде не обращается в нуль на $X$, то {\it отношением} этих функций
$\frac{g}{f}$ называется функция на $X$, определенная формулой
\begin{align}
& \frac{g}{f}(x):=\frac{g(x)}{f(x)}, && x\in X \label{DEF:g/f}
\end{align}

\item[$\bullet$] Если $f$ -- функция, определенная на множестве $X$, а
$\lambda$ -- число, то их {\it произведением} называется функция $\lambda\cdot
f$ на $X$, определенная формулой
\begin{align}
& (\lambda\cdot f)(x):=\lambda\cdot f(x), && x\in X \label{DEF:lambda-f}
\end{align}

}\eit

\begin{tm}[\bf об арифметических операциях с непрерывными функциями]
\label{cont-alg} Если функции $f$ и $g$ непрерывны (на области определения), то
непрерывны (на области определения) и функции
$$
f+g, \quad f-g, \quad \lambda\cdot f, \quad f\cdot g, \quad \frac{g}{f}
$$
(где $\lambda\in\R$ -- произвольное число).
\end{tm}\begin{proof} Пусть точка $a$ и
последовательность $\{ x_n \}$ таковы, что обе функции $f$ и $g$ в них
определены, причем
$$
x_n\underset{n\to \infty}{\longrightarrow} a
$$
Тогда, поскольку функции $f$ и $g$ непрерывны в точке $x=a$, мы получаем
$$
f(x_n)\underset{n\to \infty}{\longrightarrow} f(a) \qquad
g(x_n)\underset{n\to \infty}{\longrightarrow} g(a)
$$
Поэтому, в силу свойств $1^0,2^0,4^0$ из \ref{lim-ariphm},
мы получаем
$$
f(x_n)+g(x_n)\underset{n\to \infty}{\longrightarrow}
f(a)+g(a) \qquad f(x_n)-g(x_n)\underset{n\to
\infty}{\longrightarrow} f(a)-g(a) \qquad f(x_n)\cdot
g(x_n)\underset{n\to \infty}{\longrightarrow} f(a)\cdot
g(a)
$$
Это верно для произвольной точки $a$ и последовательности $x_n\underset{n\to
\infty}{\longrightarrow} a$, поэтому функции $f+g$, $f-g$, $f\cdot g$
непрерывны. Аналогично доказывается непрерывность для функции $\lambda\cdot f$.

Если, кроме того, дополнительно потребовать, чтобы $f(a)\ne
0$ и $f(x_n)\ne 0$, то по свойству $5^0$ из
\ref{lim-ariphm} мы получим
$$
\frac{g(x_n)}{f(x_n)}\underset{n\to
\infty}{\longrightarrow}\frac{g(a)}{f(a)}
$$
и поскольку это будет верно для любой точки $a\in D(\frac{g}{f})$ и любой
последовательности $x_n\in D(\frac{g}{f}),$ $x_n\underset{n\to
\infty}{\longrightarrow} a$, мы получим, что функция $\frac{g}{f}$ тоже должна
быть непрерывна.
\end{proof}

\paragraph{Непрерывность композиции.}\label{SUBSEC-cont-of-comp}

Вспомним, что на с.\pageref{DEF:kompozitsiya}  мы определили понятие композиции
двух отображений. В частном случае, когда отображения являются числовыми
функциями, их композиция также становится числовой функцией, и мы можем
повторить это определение, выбросив из него слово ``отображение'':

 \bit{
\item[$\bullet$] Пусть $g$ и $f$ -- две функции. Тогда функция $g\circ f$,
определенная правилом
 \beq\label{DEF:kompozitsiya-funktsij}
(g\circ f)(x)=g(f(x)),
 \eeq
называется {\it композицией} функций $g$ и $f$. Эту функцию называют также {\it
сложной функцией}, составленной из функций $f$ и $g$.
 }\eit

\noindent\rule{160mm}{0.1pt}\begin{multicols}{2}

\begin{ex}
Пусть $g(y)=\frac{1}{y}$ и $f(x)=x+1$. Тогда
$$
h(x)=g(y)\Big|_{y=f(x)}=\frac{1}{y}\Big|_{y=x+1}=\frac{1}{x+1}
$$
Если же взять композицию этих функций в обратном порядке, то мы получим другую
функцию:
$$
r(y)=f(x)\Big|_{x=g(y)}=x+1\Big|_{x=\frac{1}{y}}=\frac{1}{y}+1
$$
\end{ex}

\begin{ex}
Представить функции в виде композиции двух функций из списка на
с.\pageref{|x|-ex-func}:
$$
G(x)=|x+1|, \quad H(x)=\frac{1}{\{x\}}
$$
Решение:
$$
G(x)=|x+1|=|y|\Big|_{y=x+1},
$$
$$
H(x)=\frac{1}{\{x\}}=\frac{1}{y}\Big|_{y=\{x\}}
$$
\end{ex}

\begin{ers}
Составьте композиции из следующих функций в разном порядке (то есть найдите
$g(y)\Big|_{y=f(x)}$ и $f(x)\Big|_{x=g(y)}$), и проверьте, будут ли эти функции
одинаковы:
 \biter{
\item[1)] $g(y)=y^2$ и $f(x)=D(x)-\frac{1}{2}$.

\item[2)] $g(y)=\sgn y$ и $f(x)=|x|$.

 }\eiter
\end{ers}

\begin{ers}
Представьте функции в виде композиции двух функций из списка на
с.\pageref{|x|-ex-func}:
 \biter{
\item[1)] $[\frac{x}{2}]$,

\item[2)] $\frac{[x]}{2}$.
 }\eiter
\end{ers}

\end{multicols}\noindent\rule[10pt]{160mm}{0.1pt}

\begin{tm}[\bf о композиции непрерывных функций]\label{cont-composition}
Если функции $f$ и $g$ непрерывны (на области определения), то их композиция
$h(x)=g(f(x))$ тоже непрерывна (на области определения).
\end{tm}\begin{proof} Пусть точка $a$ и
последовательность $\{ x_n \}$ лежат в области определения функции
$h(x)=g(f(x))$, причем
$$
x_n\underset{n\to \infty}{\longrightarrow} a
$$
Тогда, поскольку $f(x)$ непрерывна в точке $x=a$,
$$
f(x_n)\underset{n\to \infty}{\longrightarrow} f(a)
$$
Отсюда, поскольку $g(y)$ непрерывна в точке $y=f(a)$,
$$
h(x_n)=g(f(x_n))\underset{n\to \infty}{\longrightarrow}
g(f(a))
$$
Мы получили, что если $x_n\underset{n\to
\infty}{\longrightarrow} a$, то
$h(x_n)=g(f(x_n))\underset{n\to \infty}{\longrightarrow}
g(f(a))$. Это означает, что $h(x)=g(f(x))$ непрерывна в
точке $x=a$. Поскольку точка $a$ с самого начала выбиралась
произвольной, получаем, что $h(x)=g(f(x))$ непрерывна в
произвольной точке, то есть всюду на области определения.
\end{proof}

\section{Классические теоремы о непрерывных
функциях}\label{SEC-th-o-nepr-func}

Непрерывные функции обладают несколькими важными
свойствами, которые называются {\it классическими теоремами
о непрерывных функциях}, потому что были доказаны в 19 веке
известными мировыми математиками. Эти теоремы нужны для
доказательства других важных фактов математического
анализа, в частности, теоремы Ролля \ref{Roll} (необходимой
для доказательства теорем Лагранжа, Коши и Тейлора из главы
\ref{ch-th-smooth-f(x)}), а также теоремы об
интегрируемости непрерывной функции из главы 14.

\subsection{Теорема о сохранении знака}

\begin{tm}[\bf о сохранении знака непрерывной функцией]
\index{теорема!о сохранении знака непрерывной
функцией}\label{sign-pres}\footnote{Эта теорема
используется далее в главе \ref{ch-graph-f(x)} при
доказательстве теоремы \ref{tochki-peregiba} о точках
перегиба} Пусть функция $f$ определена в некоторой
окрестности точки $c$ и непрерывна в точке $c$, причем
$f(c)\ne 0$. Тогда существует окрестность
$(c-\delta;c+\delta)$ точки $c$ такая, что всюду в этой
окрестности функция $f$ имеет тот же знак, что и $f(c)$:
$$
  \forall x\in (c-\delta;c+\delta) \quad \sgn f(x)=\sgn f(c)
$$
\end{tm}\begin{proof} Пусть для определенности $f(c)>0$,
покажем что для некоторого $\delta>0$ выполняется
соотношение
\begin{equation}\forall x\in (c-\delta;c+\delta) \quad f(x)>0 \label{4.1.1}\end{equation}

%\picture{0pt}{0pt}{101.pcx}

\vglue140pt \noindent Доказательство проводится методом от
противного. Предположим, что (1.1) не выполняется ни для
какого $\delta>0$, то есть что
$$
\forall \delta>0 \quad \exists x\in (c-\delta;c+\delta)
\quad f(x)\le 0
$$
Тогда, в частности, оно не выполняется для чисел
$\delta_n=\frac{1}{n}$:
$$
\forall n\in \mathbb{N}\quad \exists x_n\in \l
c-\frac{1}{n};c+\frac{1}{n}\r \quad f(x_n)\le 0
$$
Мы получаем, что имеется последовательность
\begin{equation}
x_n\underset{n\to \infty}{\longrightarrow} c
\label{4.1.2}\end{equation} для которой
\begin{equation}
f(x_n)\le 0 \label{4.1.3}\end{equation} Поскольку $f(x)$
непрерывна в точке $c$, из \eqref{4.1.2} следует
$$
f(x_n)\underset{n\to \infty}{\longrightarrow} f(c)
$$
а по теореме \ref{x_n<_y_n} (о предельном переходе в
неравенстве), из \eqref{4.1.3} следует
$$
f(c)\le 0
$$
Это противоречит предположению, что $f(c)>0$. Значит,
\eqref{4.1.1} все-таки должно выполняться для какого-то
$\delta>0$. \end{proof}

\subsection{Теорема Коши о промежуточном значении.}

\begin{lm}[\bf о промежуточном значении]\label{lm-Cauchy}
Пусть функция $\ph$ непрерывна на отрезке $[a;b]$, и на
концах его принимает значения разных знаков:
$$
  \ph(a)<0<\ph(b) \qquad (\text{или}\,\, \ph(a)>0>\ph(b) )
$$
Тогда существует точка $c\in (a;b)$ в которой $\ph(x)$
равна нулю:
$$
  \ph(c)=0
$$
\end{lm}\begin{proof} Пусть для определенности
$$
\ph(a)<0<\ph(b)
$$

%\picture{0pt}{0pt}{102.pcx}

\vglue140pt \noindent Разделим отрезок $[a;b]$ пополам.
Если значение функции в середине отрезка равно нулю, то
теорема доказана. Если же нет, то выберем из двух половинок
ту, у которой на концах функция $\ph$ принимает значения
разных знаков, и обозначим ее $[a_1;b_1]$:
$$
\ph(a_1)<0<\ph(b_1), \qquad [a_1;b_1]\subseteq [a;b]
$$
Проделаем то же самое с отрезком $[a_1;b_1]$: разделим его
пополам, выберем ту половинку, на концах которой $\ph(x)$
меняет знак, и обозначим ее $[a_2;b_2]$:
$$
\ph(a_2)<0<\ph(b_2), \qquad [a_2;b_2]\subseteq [a_1;b_1]
$$
Затем поделим отрезок $[a_2;b_2]$, и так далее. В
результате у нас получится цепочка вложенных отрезков, на
концах которых функция меняет знак:
$$
[a;b]\supset [a_1;b_1]\supset [a_2;b_2]\supset ... \qquad
b_n-a_n=\frac{b-a}{2^n}\underset{n\to
\infty}{\longrightarrow} 0, \qquad \ph(a_n)<0<\ph(b_n)
$$
По теореме \ref{I_1>I_2>...} о вложенных отрезках,
существует единственная точка $c$, принадлежащая всем
отрезкам $[a_n;b_n]$, причем
$$
\lim_{n\to \infty} a_n=c=\lim_{n\to \infty} b_n
$$
Отсюда следует, что с одной стороны, поскольку
$\ph(a_n)<0$,
$$
\ph(c)=\lim_{n\to \infty}\ph(a_n)\le 0
$$
а с другой, поскольку $\ph(b_n)>0$,
$$
\ph(c)=\lim_{n\to \infty}\ph(b_n)\ge 0
$$
И следовательно, $\ph(c)=0$. \end{proof}

\begin{tm}[\bf Коши о промежуточном значении]\label{Cauchy-I}
\index{теорема!Коши!о промежуточном значении} \footnote{Эта теорема
используется далее в $\S 6$ главы 14 при доказательстве теоремы о среднем для
неопределенного интеграла.} Пусть функция $f$ непрерывна на отрезке $[a;b]$, и
пусть $C$ -- произвольное число, лежащее между $f(a)$ и $f(b)$:
$$
  f(a)<C<f(b) \qquad \Big(\text{или}\,\, f(a)>C>f(b) \Big)
$$
Тогда на интервале $(a;b)$ найдется точка $c$ такая, что
$$
  f(c)=C
$$
\end{tm}\begin{proof} Рассмотрим функцию
$$
  \ph(x)=f(x)-C
$$
Она будет непрерывна на отрезке $[a;b]$ (как разность
непрерывных функций) и принимает на концах этого отрезка
значения разных знаков
$$
 \ph(a)=f(a)-C<0<f(b)-C=\ph(b) \qquad
\Big(\text{или}\,\, \ph(a)=f(a)-C>0>f(b)-C=\ph(b) \Big)
$$
Поэтому, по лемме \ref{lm-Cauchy}, существует точка $c\in
[a;b]$ такая, что $\ph(c)=0$. Отсюда $f(c)=\ph(c)+C=C$.
\end{proof}

\subsection{Теоремы Вейерштрасса}

\paragraph{Теорема Вейерштрасса об ограниченности.}

Напомним, что функция $f$ называется
 \bit{
\item[--] {\it ограниченной сверху} на множестве $E$, если
существует число $B$ такое, что
$$
\forall x\in E \quad f(x)\le B
$$

\item[--] {\it ограниченной снизу} на множестве $E$, если
существует число $A$ такое, что
$$
\forall x\in E \quad A\le f(x)
$$

\item[--] {\it ограниченной} на множестве $E$, если она
ограничена на $E$ сверху и снизу, то есть существуют такие
числа $A$ и $B$ что:
$$
\forall x\in E \quad A\le f(x)\le B
$$

 }\eit

\begin{tm}[\bf Вейерштрасса об ограниченности]\label{Wei-II}\index{теорема!Вейерштрасса!об ограниченности}\footnote{Эта
теорема используется в следующем параграфе при доказательстве теоремы
Вейерштрасса об экстремумах \ref{Wei-III}} Если функция $f$ определена и
непрерывна на отрезке $[a;b]$, то она ограничена на этом отрезке.
\end{tm}\begin{proof} Предположим обратное, то есть что
$f(x)$ не ограничена на отрезке $[a;b]$. Значит $f(x)$ не
ограничена сверху или снизу. Пусть для определенности
$f(x)$ не ограничена сверху:
$$
\forall B\in \R\quad \exists x\in [a;b] \quad f(x)> B
$$
Тогда для всякого $n\in \mathbb{N}$ можно подобрать точку
$x_n \in [a;b]$ такую, что
\begin{equation}
f(x_n)> n \label{4.3.1}\end{equation} Последовательность точек $\{ x_n \}$
принадлежит отрезку $[a;b]$ и значит ограничена. Следовательно, по теореме
Больцано-Вейерштрасса \ref{Bol-Wei},  $\{ x_n \}$ содержит сходящуюся
подпоследовательность $\{ x_{n_k}\}$:
\begin{equation}
x_{n_k}\underset{k\to \infty}{\longrightarrow} c \in [a;b]
\label{4.3.2}\end{equation} Из формулы \eqref{4.3.1}
получаем
$$
f(x_{n_k})>n_k\underset{k\to \infty}{\longrightarrow}
+\infty
$$
поэтому
\begin{equation}
f(x_{n_k})\underset{k\to \infty}{\longrightarrow} +\infty
\label{4.3.3}\end{equation} Но с другой стороны, функция $f$ непрерывна в любой точке отрезка $[a;b]$, и в
частности в точке $c\in [a;b]$, значит из \eqref{4.3.2}
следует
\begin{equation}
f(x_{n_k})\underset{k\to \infty}{\longrightarrow} f(c)\ne
+\infty \label{4.3.4}\end{equation} Мы получили
противоречие между \eqref{4.3.3} и \eqref{4.3.4}, которое
означает, что наше исходное предположение было неверно.
Значит, функция $f$ все-таки ограничена. \end{proof}

\paragraph{Теорема Вейерштрасса об экстремумах.}

Пусть функция $f$ определена на множестве $E\subseteq \R$. Точка
$x_{\max}\in E$ называется {\it точкой максимума} функции $f$ на множестве
$E$, если
$$
\forall x\in E \quad f(x)\le f(x_{\max})
$$
Число $f(x_{\max})$ при этом называется {\it максимумом}
функции $f$ на множестве $E$ и обозначается
$$
f(x_{\max})=\max_{x\in E} f(x)
$$

Наоборот, точка $x_{\min}\in E$ называется {\it точкой
минимума} функции $f$ на множестве $E$, если
$$
\forall x\in E \quad f(x_{\min})\le f(x)
$$
Число $f(x_{\min})$ при этом называется {\it минимумом}
функции $f$ на множестве $E$ и обозначается
$$
f(x_{\min})=\min_{x\in E} f(x)
$$

{\it Экстремумом} функции $f$ на множестве $E$
называется минимум или максимум  функции $f$ на
множестве $E$.

\begin{rmk}
Максимум и минимум функции -- не то же самое, что точная
верхняя и точная нижняя грань.
\end{rmk}

\noindent\rule{160mm}{0.1pt}\begin{multicols}{2}

\begin{ex}
Пусть $E=[-1;1]$ и $f(x)=x^2$.

%\picture{0pt}{0pt}{103.pcx}

\vglue120pt \noindent Тогда на множестве $E$ у функции $f$ имеется одна точка минимума и две точки максимума,
причем минимум совпадает с точной нижней гранью, а максимум
-- с точной верхней гранью:
$$
\min_{x\in [-1;1]} x^2=x^2\Big|_{x=0}=0=\inf_{x\in [-1;1]}
x^2
$$
$$
\max_{x\in [-1;1]}
x^2=x^2\Big|_{x=-1}=x^2\Big|_{x=1}=1=\sup_{x\in [-1;1]} x^2
$$
\end{ex}

\begin{ex}
Пусть $E=(-1;1)$ и $f(x)=x^2$.

%\picture{0pt}{0pt}{104.pcx}

\vglue120pt \noindent Тогда на множестве $E$ у функции $f$ имеется одна точка минимума и {\it нет точек
максимума}, причем точная нижняя грань (существует и)
совпадает с минимумом, а точная верхняя грань существует,
но не совпадает с максимумом, которого на множестве
$E=(-1;1)$ нет (потому что нет точек максимума):
$$
\min_{x\in (-1;1)} x^2=x^2\Big|_{x=0}=0=\inf_{x\in (-1;1)}
x^2
$$
$$
\sup_{x\in (-1;1)} x^2=1 \qquad \nexists \max_{x\in (-1;1)}
x^2
$$
\end{ex}

\begin{ex}
Пусть $E=\R$ и $f(x)=\{x\}$.

%\picture{0pt}{0pt}{44.pcx}

\vglue120pt \noindent Тогда на множестве $E$ у функции $f$ {\it нет точек
минимума}, хотя точная нижняя грань существует:
$$
\inf_{x\in \R}\{x\}=0\qquad \nexists \min_{x\in \R}\{x\}
$$
При этом максимум на этом множестве есть:
$$
\sup_{x\in \R}\{x\}=\max_{x\in \R}\{x\}=1
$$
\end{ex}

\end{multicols}\noindent\rule[10pt]{160mm}{0.1pt}

Итак, экстремумы могут существовать, а могут и не
существовать.

\begin{tm}[\bf Вейерштрасса об экстремумах]\label{Wei-III}
\index{теорема!Вейерштрасса!об экстремумах} \footnote{Эта теорема используется
далее в главе \ref{ch-th-smooth-f(x)} при доказательстве теоремы Ролля, а также
в $\S 6$ главы 14 при доказательстве интегральной теоремы о среднем, и в $\S 7$
главы 14 при доказательстве теоремы об интегрируемости непрерывной функции.}
Если функция $f$ определена и непрерывна на отрезке $[a;b]$, то она имеет
минимум и максимум на этом отрезке.
\end{tm}\begin{proof} Докажем существование максимума. По
теореме \ref{Wei-II}, функция $f$ ограничена сверху на
$[a;b]$:
$$
\exists C\quad \forall x\in [a;b] \quad f(x)\le C
$$
поэтому по теореме \ref{sup-inf} о точной границе, у
функции $f$ существует точная верхняя грань на этом
множестве:
$$
\exists \sup_{x\in [a;b]} f(x)=M
$$
Это значит, что во-первых, все значения функции $f$ на
отрезке $[a;b]$ не превосходят $M$
\begin{equation}\forall x\in [a;b] \quad f(x)\le M \label{4.4.1}\end{equation}
и во-вторых если взять любое число $\alpha <M$, то
обязательно оно окажется меньше, чем какое-то значение
$f(x)$ на отрезке $[a;b]$
$$
\forall \alpha <M \quad \exists x\in [a;b] \quad \alpha
<f(x)
$$

Возьмем в качестве $\alpha$ числа вида $M-\frac{1}{n}$. Тогда у нас получится
целая последовательность $\{ x_n \}$:
\begin{equation}\forall n\in \mathbb{N}\quad \exists x_n\in [a;b] \quad
M-\frac{1}{n} <f(x_n) \label{4.4.2}\end{equation} Поскольку последовательность
$\{ x_n \}$ лежит в отрезке $[a;b]$, она ограничена. Значит по теореме
Больцано-Вейерштрасса \ref{Bol-Wei}, из нее можно выбрать сходящуюся
подпоследовательность:
$$
x_{n_k}\underset{k\to \infty}{\longrightarrow} c \qquad (c
\in [a;b] )
$$
Поскольку $f(x)$ непрерывна в любой точке отрезка $[a;b]$
и, в частности, в точке $c\in [a;b]$, мы получаем
$$
f(x_{n_k})\underset{k\to \infty}{\longrightarrow} f(c)
$$
С другой стороны, из \eqref{4.4.1} и \eqref{4.4.2} следует
$$
M-\frac{1}{n_k} <f(x_{n_k})\le M
$$
откуда
$$
M=\lim_{k\to\infty}\left( M-\frac{1}{n_k}\right) \le
\lim_{k\to\infty} f(x_{n_k})=f(c)\le M
$$
то есть
$$
f(c)=M
$$
Еще раз вспомнив \eqref{4.4.1}, получим
$$
\forall x\in [a;b] \quad f(x)\le M=f(c)
$$
Это означает, что $c\in [a;b]$ является точкой максимума
функции $f$ на множестве $[a;b]$. \end{proof}

\subsection{Теорема Кантора о равномерной
непрерывности}\label{SEC-teorema-Kantora}

Говорят, что функция $f$ {\it равномерно непрерывна на множестве} $E$ если
для любых двух последовательностей аргументов $\alpha_n \in E$ и $\beta_n\in E$
стремящихся друг к другу, соответствующие последовательности значений $\{
f(\alpha_n) \}$ и $\{ f(\beta_n) \}$ тоже стремятся друг к другу:
 \begin{equation}
 \forall \alpha_n,\beta_n\in E \quad
 \alpha_n - \beta_n \underset{n\to \infty}{\longrightarrow} 0
 \quad \Longrightarrow \quad
 f(\alpha_n) - f(\beta_n) \underset{n\to \infty}{\longrightarrow} 0
 \label{ravnomernaya nepreryvnost po Heine}
 \end{equation}

\noindent\rule{160mm}{0.1pt}\begin{multicols}{2}

\begin{ex}
Линейная функция $f(x)=k x+b$ равномерно непрерывна на множестве $E=\R$, потому
что если
$$
\alpha_n - \beta_n \underset{n\to \infty}{\longrightarrow}
0
$$
то
 \begin{multline*}
f(\alpha_n) - f(\beta_n) =(k\alpha_n+b) - (k\beta_n+b)=\\= k\alpha_n -
k\beta_n=k(\alpha_n - \beta_n) \underset{n\to \infty}{\longrightarrow} 0
 \end{multline*}
\end{ex}

\begin{ex}
Квадратичная функция $f(x)=x^2$ не является равномерно непрерывной на множестве
$E=\R$, потому что если взять, например,
$$
\alpha_n=n+\frac{1}{n}, \quad \beta_n=n
$$
то мы получим, что
$$
\alpha_n - \beta_n =\frac{1}{n}\underset{n\to
\infty}{\longrightarrow} 0
$$
но при этом
 \begin{multline*}
f(\alpha_n) - f(\beta_n) =\left(\alpha_n \right)^2 - \left(\beta_n \right)^2=
\l n+\frac{1}{n}\r^2 - n^2=\\=2n\cdot \frac{1}{n}+\frac{1}{n^2}=
2+\frac{1}{n^2}\underset{n\to \infty}{\longrightarrow} 2\ne 0
 \end{multline*}
\end{ex}

\end{multicols}\noindent\rule[10pt]{160mm}{0.1pt}

\begin{tm}[\bf Кантора о равномерной непрерывности]\label{Kantor}\index{теорема!Кантора}\footnote{Эта теорема используется далее в
$\S 7$ главы 14 при доказательстве теоремы об
интегрируемости непрерывной функции.} Если функция $f$
определена и непрерывна на отрезке $[a;b]$, то она
равномерно непрерывна на этом отрезке.
\end{tm}\begin{proof} Доказательство проводится от
противного. Предположим, что $f(x)$ не является равномерно непрерывной на
$[a;b]$. Тогда существуют последовательности аргументов $\alpha_n \in [a;b]$ и
$\beta_n\in [a;b]$ которые стремятся друг к другу
$$
\alpha_n - \beta_n \underset{n\to \infty}{\longrightarrow}
0
$$
но при этом последовательности значений не стремятся друг к другу:
$$
f(\alpha_n) - f(\beta_n) \underset{n\to \infty}{\notarrow} 0
$$
Рассмотрим последовательность
$$
\gamma_n=f(\alpha_n) - f(\beta_n)
$$
То, что она не стремится к нулю
$$
\gamma_n \underset{n\to \infty}{\notarrow} 0
$$
означает, что найдется некоторая окрестность нуля $(-\varepsilon;\varepsilon)$
которая не содержит бесконечное количество чисел $\{ \gamma_n \}$. Значит можно
выбрать подпоследовательность $\{ \gamma_{n_k}\}$, которая вообще не будет
заходить в окрестность $(-\varepsilon;\varepsilon)$:
\begin{equation}\forall k \quad |\gamma_{n_k}|=|f(\alpha_n) - f(\beta_n)|\ge
\varepsilon \label{4.5.1}\end{equation} Рассмотрим теперь индексы $\{ n_k \}$.
Им соответствует некоторая подпоследовательность $\{ \alpha_{n_k}\}$
последовательности $\{ \alpha_n \}$. Поскольку она ограничена (из-за того, что
$\alpha_{n_k}\in [a;b]$), по теореме Больцано-Вейерштрасса \ref{Bol-Wei} мы
можем выбрать у нее некоторую сходящуюся подпоследовательность $\{
\alpha_{n_{k_j}}\}$:
$$
\alpha_{n_{k_j}}\underset{j\to \infty}{\longrightarrow} c
\qquad (c\in [a;b])
$$
Тогда мы получим, что соответствующая последовательность $\{ \beta_{n_{k_j}}\}$
стремится к тому же пределу $c$, потому что
$$
\beta_{n_{k_j}}=(\beta_{n_{k_j}}-\alpha_{n_{k_j}}) +
\alpha_{n_{k_j}}\underset{j\to \infty}{\longrightarrow}
0+c=c
$$
Отсюда следует, что
$$
f(\alpha_{n_{k_j}})- f(\beta_{n_{k_j}}) \underset{j\to
\infty}{\longrightarrow} c-c=0
$$
а это противоречит условию \eqref{4.5.1}, из которого
видно, что расстояние между $f(\alpha_{n_{k_j}})$ и
$f(\beta_{n_{k_j}})$ не может быть меньше $\varepsilon$.

Итак, мы предположили, что $f(x)$ не является равномерно
непрерывной на $[a;b]$ и получили противоречие. Это
означает, что предположение неверно, и функция $f$
все-таки равномерно непрерывна на $[a;b]$. \end{proof}

\subsection{Теоремы о монотонных функциях.}

\paragraph{Монотонные функции на отрезке.}

\btm[\bf о монотонной функции на отрезке]
\label{Teor-ob-obratnoi-funktsii-na-R^1} Пусть функция $f$ определена и строго
монотонна на отрезке $I=[a,b]$, и пусть $J$ -- отрезок с концами $f(a)$ и
$f(b)$:
$$
J=\begin{cases}[f(a);f(b)],& \text{если $f$ возрастает} \\
[f(b);f(a)],& \text{если $f$ убывает.}
\end{cases}
$$
Тогда:
 \bit{
\item[(i)] если $f$ биективно отображает отрезок $I$ на отрезок $J$, то
 \bit{
\item[(a)] она непрерывна на $I$, и

\item[(b)] ее обратная функция $f^{-1}:J\to I$ также непрерывна и строго
монотонна.
 }\eit
\item[(ii)] если $f$ непрерывна на $I$, то
 \bit{
\item[(a)] она биективно отображает отрезок $I$ на отрезок $J$, и

\item[(b)] ее обратная функция $f^{-1}:J\to I$ также непрерывна и строго
монотонна.
 }\eit
 }\eit\noindent
При этом характер монотонности наследуется обратной функцией:
 \bit{
\item[---] если $f:I\to J$ возрастает, то и $f^{-1}:J\to I$ возрастает,

\item[---] если $f:I\to J$ убывает, то и $f^{-1}:J\to I$ убывает.
 }\eit
\etm
 \bpr
Будем считать для определенности, что функция $f$ строго возрастает:
 \beq\label{f-vozrastaet-*}
s<t\qquad\Longrightarrow\qquad f(s)<f(t)
 \eeq
(случай убывания рассматривается аналогично).

Докажем (i). Пусть $f$ биективно отображает отрезок $I$ на отрезок $J$. Тогда у
нее имеется обратная функция $f^{-1}:J\to I$.

1. Прежде всего заметим, что $f^{-1}:J\to I$ также будет возрастать.
Действительно, если бы это было не так, то для некоторых $y_1,y_2\in J$ мы
получили бы
$$
y_1<y_2\qquad \&\qquad f^{-1}(y_1)\ge f^{-1}(y_2)
$$
Тогда можно было бы обозначить $x_1=f^{-1}(y_1)$ и $x_2=f^{-1}(y_2)$, и
получилось бы:
$$
\underbrace{f(x_1)}_{y_1}<\underbrace{f(x_2)}_{y_2}\qquad \&\qquad
\underbrace{x_1}_{f^{-1}(y_1)}\ge \underbrace{x_2}_{f^{-1}(y_2)}
$$
Это противоречит утверждению $1^\circ$ на с.\pageref{vozrastanie=>neubyvanie}:
поскольку $f$ возрастает, она должна неубывать, значит условие $x_1\ge x_2$
должно влечь за собой условие $f(x_1)\ge f(x_2)$, а у нас получается наоборот,
$f(x_1)<f(x_2)$.

2. Покажем, что функция $f:I\to J$ непрерывна. Пусть
$x_n\underset{n\to\infty}{\longrightarrow}x$, покажем что
$f(x_n)\underset{n\to\infty}{\longrightarrow}f(x)$. По теореме
\ref{TH:nepr-monot}, последовательность $x_n$ можно считать строго монотонной.
Нам придется рассмотреть два случая.
 \bit{
\item[A.] Пусть сначала $x_n$ возрастает:
$$
x_1<x_2<...<x_n<...<x
$$
Тогда предел $x$ этой последовательности равен ее точной верхней грани
 \beq\label{lim-x_n=sup-x_n}
\lim_{n\to\infty} x_n=\sup_{n\in\N} x_n=x
 \eeq
и то же для предела последовательности $f(x_n)$ поскольку, в силу
\eqref{f-vozrastaet-*}, она тоже возрастает:
 \beq\label{lim-f(x_n)=sup-f(x_n)}
\lim_{n\to\infty} f(x_n)=\sup_{n\in\N} f(x_n)
 \eeq
Обозначив $y_n=f(x_n)$, мы получим
$$
\sup_{n\in\N} y_n=\sup_{n\in\N} f(x_n)\le \eqref{f(inf-E)<inf-f(E)}\le
f\big(\sup_{n\in\N} x_n\big)=\eqref{lim-x_n=sup-x_n}=f(x)
$$
$$
\Downarrow
$$
$$
\sup_{n\in\N} y_n\le f(x)
$$
$$
\phantom{\text{\scriptsize ($f^{-1}$ возрастает, как уже
доказано)}}\quad\Downarrow\quad\text{\scriptsize ($f^{-1}$ возрастает, как уже
доказано)}
$$
 \beq\label{f^(-1)(sup-y_n)<x}
f^{-1}\Big(\sup_{n\in\N} y_n\Big)\le f^{-1}\big(f(x)\big)=x
 \eeq
А с другой стороны, опять, поскольку $f^{-1}$ возрастает,
$$
x=\sup_{n\in\N} x_n=\sup_{n\in\N} f^{-1}\big(f(x_n)\big)=\sup_{n\in\N}
f^{-1}(y_n)\le\eqref{f(inf-E)<inf-f(E)}\le f^{-1}\Big(\sup_{n\in\N} y_n\Big)
$$
$$
\Downarrow
$$
 \beq\label{x<f^(-1)(sup-y_n)}
 x\le f^{-1}\Big(\sup_{n\in\N} y_n\Big)
 \eeq
Неравенства \eqref{f^(-1)(sup-y_n)<x} и \eqref{x<f^(-1)(sup-y_n)} вместе дают
  $$
  x=f^{-1}\Big(\sup_{n\in\N} y_n\Big)
  $$
То есть
  $$
f(x)=\sup_{n\in\N} y_n=\sup_{n\in\N}
f(x_n)=\eqref{lim-f(x_n)=sup-f(x_n)}=\lim_{n\to\infty} f(x_n)
  $$

\item[B.] Пусть теперь наоборот, $x_n$ убывает:
$$
x<...<x_n<...<x_2<x_1
$$
Тогда предел $x$ этой последовательности равен ее точной нижней грани
 \beq\label{lim-x_n=inf-x_n}
\lim_{n\to\infty} x_n=\inf_{n\in\N} x_n=x
 \eeq
и то же для предела последовательности $f(x_n)$ поскольку, в силу
\eqref{f-vozrastaet-*}, она тоже возрастает:
 \beq\label{lim-f(x_n)=inf-f(x_n)}
\lim_{n\to\infty} f(x_n)=\inf_{n\in\N} f(x_n)
 \eeq
Обозначив $y_n=f(x_n)$, мы получим
$$
\inf_{n\in\N} y_n=\inf_{n\in\N} f(x_n)\ge \eqref{f(inf-E)<inf-f(E)}\ge
f\big(\inf_{n\in\N} x_n\big)=\eqref{lim-x_n=inf-x_n}=f(x)
$$
$$
\Downarrow
$$
$$
\sup_{n\in\N} y_n\ge f(x)
$$
$$
\phantom{\text{\scriptsize ($f^{-1}$ возрастает, как уже
доказано)}}\quad\Downarrow\quad\text{\scriptsize ($f^{-1}$ возрастает, как уже
доказано)}
$$
 \beq\label{f^(-1)(inf-y_n)>x}
f^{-1}\Big(\inf_{n\in\N} y_n\Big)\ge f^{-1}\big(f(x)\big)=x
 \eeq
А с другой стороны, опять, поскольку $f^{-1}$ возрастает,
$$
x=\inf_{n\in\N} x_n=\inf_{n\in\N} f^{-1}\big(f(x_n)\big)=\inf_{n\in\N}
f^{-1}(y_n)\ge\eqref{f(inf-E)<inf-f(E)}\ge f^{-1}\Big(\inf_{n\in\N} y_n\Big)
$$
$$
\Downarrow
$$
 \beq\label{x>f^(-1)(inf-y_n)}
 x\ge f^{-1}\Big(\inf_{n\in\N} y_n\Big)
 \eeq
Неравенства \eqref{f^(-1)(inf-y_n)>x} и \eqref{x>f^(-1)(inf-y_n)} вместе дают
  $$
  x=f^{-1}\Big(\inf_{n\in\N} y_n\Big)
  $$
То есть
  $$
f(x)=\inf_{n\in\N} y_n=\inf_{n\in\N}
f(x_n)=\eqref{lim-f(x_n)=inf-f(x_n)}=\lim_{n\to\infty} f(x_n)
  $$

 }\eit

3. Мы показали, что если функция $f:I\to J$ возрастает и биективна, то она
автоматически непрерывна. Поскольку обратная функция $f^{-1}:J\to I$ тоже
биективна и, как мы уже убедились, тоже возрастает, мы получаем, что она тоже
должна быть непрерывной.

Докажем (ii). Пусть $f$ непрерывна на $I$.

1. Заметим, что $f$ отображает $I$ в $J$. Это следует из того, что функция $f$
неубывает ($1^\circ$ на с.\pageref{vozrastanie=>neubyvanie}): если $x\in
I=[a,b]$, то
$$
a\le x\le b\quad\Longrightarrow\quad f(a)\le f(x)\le
f(b)\quad\Longrightarrow\quad f(x)\in J=[f(a),f(b)]
$$

2. Проверим, что $f$ инъективно отображает $I$ в $J$. Для любых $x,y\in I$
получаем:
 \beq\label{inj-f-dlya-teor-ob-obr-func}
x\ne y\quad\Longrightarrow\quad\left[\begin{matrix}
x<y\quad\Longrightarrow\quad f(x)<f(y) \\ y<x\quad\Longrightarrow\quad
f(y)<f(x)
\end{matrix}\right]\quad\Longrightarrow\quad\left[\begin{matrix} f(x)<f(y) \\
f(y)<f(x)\end{matrix}\right]\quad\Longrightarrow\quad f(x)\ne f(y)
 \eeq

3. Покажем далее, что $f$ сюръективно отображает $I$ в $J$, то есть что
$f(I)=J$. Для любого $C\in J=[f(a),f(b)]$ получаем:
 \bit{
\item{---} либо $C$ лежит на левом конце отрезка $J=[f(a),f(b)]$, то есть
$C=f(a)$, и тогда точка $x=a\in I=[a,b]$ будет прообразом для $C$: $C=f(x)$;

\item{---} либо $C$ лежит на правом конце отрезка $J=[f(a),f(b)]$, то есть
$C=f(b)$, и тогда точка $x=b\in I=[a,b]$ будет прообразом для $C$: $C=f(x)$;

\item{---} либо $C$ лежит внутри отрезка $J=[f(a),f(b)]$, то есть
$f(a)<C<f(b)$, и тогда по теореме Коши о промежуточном значении \ref{Cauchy-I}
найдется точка $x\in (a,b)$  такая, что $C=f(x)$.
 }\eit

4. Инъективность и сюръективность отображения $f:I\to J$ означают, что оно
является биекцией, и поэтому правило
$$
f(x)=y\qquad\Longleftrightarrow\qquad x=f^{-1}(y)
$$
определяет обратное к нему отображение $f^{-1}:J\to I$.

5. Проверим, что функция $f^{-1}:J\to I$ строго монотонна. Пусть $y_1<y_2$ --
две точки из $J$. Обозначим $x_1=f^{-1}(y_1)$ и $x_2=f^{-1}(y_2)$. Поскольку
отображение $f^{-1}$ биективно, $x_1$ не может быть равно $x_2$. Значит, либо
$x_1<x_2$, либо $x_1>x_2$. Но второе, в силу строгого возрастания $f$, влечет
за собой неравенство $y_1>y_2$:
$$
x_1>x_2\qquad\Longrightarrow\qquad y_1=f(x_1)>f(x_2)=y_2,
$$
что противоречит исходному выбору $y_1<y_2$. Значит, $x_1>x_2$ тоже невозможно.
Мы получаем, что возможно только неравенство $x_1<x_2$, и это означает, что
функция $f^{-1}$ строго возрастает:
$$
y_1>y_2\qquad\Longrightarrow\qquad x_1=f^{-1}(y_1)>f^{-1}(y_2)=x_2.
$$

6. Осталось убедиться, что функция $f^{-1}:J\to I$ непрерывна. Предположим, что
это не так, то есть что существует последовательность $y_n\in J=[f(a),f(b)]$
такая что
 \beq\label{predp-razr-f^-1}
y_n\underset{n\to\infty}{\longrightarrow} y\qquad\&\qquad
f^{-1}(y_n)\underset{n\to\infty}{\kern-12pt\longrightarrow\kern-15pt{\not}}
f^{-1}(y)
 \eeq
Обозначим $x_n=f^{-1}(y_n)$ и $x=f^{-1}(y)$, тогда
$$
x_n=f^{-1}(y_n)\underset{n\to\infty}{\kern-12pt\longrightarrow\kern-15pt{\not}}
f^{-1}(y)=x
$$
и по свойству $2^\circ$ на с.\pageref{podposledovatelnosti}, это означает, что
существует подпоследовательность $x_{n_k}$, лежащая вне некоторой окрестности
$(x-\e;x+\e)$ точки $x$:
 \beq\label{x_n_k-notin-(x-e;x+e)}
x_{n_k}\notin(x-\e;x+\e),\qquad \e>0
 \eeq
С другой стороны, последовательность $x_{n_k}$ лежит в отрезке $[a,b]$, поэтому
по теореме Больцано-Вейерштрасса \ref{Bol-Wei}, у нее должна быть сходящаяся
подпоследовательность:
$$
x_{n_{k_i}}\underset{n\to\infty}{\longrightarrow} t \qquad (t\in I=[a,b])
$$
Из \eqref{x_n_k-notin-(x-e;x+e)} следует, что $t\notin (x-\e;x+\e)$, поэтому
$t\ne x$. Но, с другой стороны,
$$
x_{n_{k_i}}\underset{n\to\infty}{\longrightarrow} t
$$
$$
\Downarrow
$$
$$
y_{n_{k_i}}=f(x_{n_{k_i}})\underset{n\to\infty}{\longrightarrow} f(t)
$$
$$
\Downarrow
$$
$$
f(t)=\lim_{i\to\infty}y_{n_{k_i}}=(\text{свойство $1^\circ$ на с.
\pageref{podposledovatelnosti}})=\lim_{n\to\infty}y_n=y
$$
$$
\Downarrow
$$
$$
t=f^{-1}(f(t))=f^{-1}(y)=x
$$
То есть $t\ne x$ и одновременно $t=x$. Это противоречие означает, что наше
предположение \eqref{predp-razr-f^-1} о разрывности функции $f^{-1}$ было
неверно. \epr

\paragraph{Монотонные функции на интервале.}

\btm[\bf о монотонной функции на интервале]
\label{Teor-ob-obratnoi-funktsii-na-int-v-R^1} Пусть функция $f$ определена и
строго монотонна на интервале $I=(a,b)$, и пусть $J$ -- интервал с концами
$\inf_{x\in I}f(x)$ и $\sup_{x\in I}f(x)$:
$$
J=\Big(\inf_{x\in I}f(x);\sup_{x\in I}f(x)\Big)
$$
Тогда:
 \bit{
\item[(i)] если $f$ биективно отображает интервал $I$ на интервал $J$, то
 \bit{
\item[(a)] она непрерывна на $I$, и

\item[(b)] ее обратная функция $f^{-1}:J\to I$ также непрерывна и строго
монотонна.
 }\eit
\item[(ii)] если $f$ непрерывна на $I$, то
 \bit{
\item[(a)] она биективно отображает интервал $I$ на интервал $J$, и

\item[(b)] ее обратная функция $f^{-1}:J\to I$ также непрерывна и строго
монотонна.
 }\eit
 }\eit\noindent
При этом характер монотонности наследуется обратной функцией:
 \bit{
\item[---] если $f:I\to J$ возрастает, то и $f^{-1}:J\to I$ возрастает,

\item[---] если $f:I\to J$ убывает, то и $f^{-1}:J\to I$ убывает.
 }\eit
 \etm

\brem В условиях теоремы из равенства $f^{-1}(J)=I$ следуют две формулы,
которые, несмотря на свою очевидность, оказываются полезными:\footnote{Мы
используем формулы \eqref{inf(y-in-J)f^(-1)(y)=inf-I} ниже в предложении
\ref{LM:x^(1/2n-1)-nepr}.}
 \begin{align}\label{inf(y-in-J)f^(-1)(y)=inf-I}
&\inf_{y\in J}f^{-1}(y)=\inf I &&\sup_{y\in J}f^{-1}(y)=\sup I
 \end{align}
\erem

 \bpr
Как и в доказательстве предыдущей теоремы будем считать, что $f$ строго
возрастает.

Докажем (i). Пусть $f$ биективно отображает интервал $I$ на интервал $J$. Тогда
у нее имеется обратная функция $f^{-1}:J\to I$.

1. Прежде всего, как и в теореме \ref{Teor-ob-obratnoi-funktsii-na-R^1}, без
труда доказывается, что обратная функция $f^{-1}:J\to I$ также возрастает.

2. Выберем отрезок $[\alpha,\beta]\subseteq I$, и покажем, что $f$ биективно
отображает его на отрезок $[f(\alpha),f(\beta)]\subseteq J$. Во-первых, из-за
монотонности, $f$ действительно отображает $[\alpha,\beta]$ в
$[f(\alpha),f(\beta)]$:
$$
x\in[\alpha,\beta]\quad\Longrightarrow\quad \alpha\le x\le\beta
\quad\Longrightarrow\quad f(\alpha)\le f(x)\le f(\beta)
$$
Во-вторых, если $C\in [f(\alpha),f(\beta)]$, то, поскольку $f:I\to J$
биективно, должно существовать $c\in I$ такое, что
$$
f(c)=C
$$
Нам нужно показать, что это $c$ лежит в отрезке $[\alpha,\beta]$. Если бы это
было не так, то мы получили бы:
 \bit{
\item[---] либо $c<\alpha$, и тогда $C=f(c)<f(\alpha)$, а это невозможно,
потому что $C\in [f(\alpha),f(\beta)]$,

\item[---] либо $\beta<c$, и тогда $f(\beta)<f(c)=C$, и это тоже невозможно,
потому что $C\in [f(\alpha),f(\beta)]$.
 }\eit\noindent
Таким образом, что $c\in[\alpha,\beta]$, и это доказывает сюръективность
отображения $f:[\alpha,\beta]\to [f(\alpha),f(\beta)]$. С другой стороны,
будучи строго монотонным, оно инъективно, поэтому является биекцией.

3. По теореме \ref{Teor-ob-obratnoi-funktsii-na-R^1}, функция
$f:[\alpha,\beta]\to [f(\alpha),f(\beta)]$ должна быть непрерывной. Это верно
для любого отрезка $[\alpha,\beta]\subseteq I$, поэтому наша исходная функция
$f:I\to J$ непрерывна на интервале $I$.

4. Для доказательства (i) остается только убедиться, что обратная функция
$f^{-1}:J\to I$ тоже непрерывна на $J$. Для этого возьмем произвольный отрезок
$[A,B]\in J$ ($A<B$) и положим
$$
\alpha=f^{-1}(A),\qquad \beta=f^{-1}(B)
$$
Поскольку, как мы уже заметили, $f^{-1}$ возрастает, получаем
$$
\alpha<\beta
$$
и, опять мы это уже показали, функция $f$ будет биективно (и строго монотонно)
отображать отрезок $[\alpha,\beta]$ на отрезок $[A,B]= [f(\alpha),f(\beta)]$.
Значит, по теореме \ref{Teor-ob-obratnoi-funktsii-na-R^1}, обратная функция
$f^{-1}:[A,B]\to [\alpha,\beta]$ тоже должна быть непрерывной. Это верно для
любого отрезка $[A,B]\in J$, поэтому функция $f^{-1}:J\to I$ должна быть
непрерывна.

Докажем (ii). Пусть $f$ непрерывна на $I$.

1. Если $x\in I=(a,b)$, то выбрав $\alpha$ и $\beta$ так, чтобы
$x\in[\alpha,\beta]\subseteq (a,b)=I$. Тогда:
$$
a<\alpha<x<\beta<b
$$
$$
\Downarrow
$$
$$
\inf_{x\in I}f(t)\le f(\alpha)<f(x)<f(\beta)\le \sup_{x\in I}f(t)
$$
$$
\Downarrow
$$
$$
f(x)\in \Big(\inf_{x\in I}f(t),\sup_{x\in I}f(t)\Big)=J
$$
То есть $f(I)\subseteq J$.

2. Инъективность отображения $f:I\to J$ доказывается той же цепочкой
\eqref{inj-f-dlya-teor-ob-obr-func}, что и для случая, когда $I$ -- отрезок.

3. Сюръективность отображения $f:I\to J$. Пусть $y\in J=\left(\inf_{t\in
I}f(t),\sup_{t\in I}f(t)\right)$. Тогда:
 \begin{multline*}
\inf_{t\in I}f(t)<y<\sup_{t\in I}f(t)\qquad\Longrightarrow\qquad \exists
t_1,t_2\in I:\qquad f(t_1)<y<f(t_2)\qquad\Longrightarrow\\
\Longrightarrow\qquad y\in [f(t_1);f(t_2)]\qquad\underset{\text{теорема
\ref{Teor-ob-obratnoi-funktsii-na-R^1}}}{\Longrightarrow}\qquad \exists
x\in[t_1;t_2]\quad y=f(x)
 \end{multline*}

4. Снова инъективность и сюръективность отображения $f:I\to J$ означают, что
оно является биекцией, и поэтому правило
$$
f(x)=y\qquad\Longleftrightarrow\qquad x=f^{-1}(y)
$$
определяет обратное к нему отображение $f^{-1}:J\to I$.

5. Строгая монотонность функции $f^{-1}:J\to I$ доказывается так же, как и в
случае, когда $I$ -- отрезок.

6. Остается проверить непрерывность функции $f^{-1}:J\to I$. Пусть $y_1<y_2$ --
произвольные точки из $J$. Рассмотрим точки $x_1=f^{-1}(y_1)$ и
$x_2=f^{-1}(y_2)$ из $I$. Поскольку функция $f^{-1}:J\to I$ строго возрастает,
должно выполняться неравенство $x_1<x_2$. Рассмотрим отрезок $[x_1,x_2]$. На
нем функция $f$ непрерывна и строго возрастает, поэтому, как мы уже доказали,
$f$ должна непрерывно и биективно отображать $[x_1,x_2]$ на
$[f(x_1),f(x_2)]=[y_1,y_2]$, а обратное отображение $f^{-1}$ будет непрерывно
отображать $[y_1,y_2]$ в $[x_1,x_2]$.

Из этого можно сделать вывод, что функция $f^{-1}:J\to I$ непрерывна на каждом
отрезке $[y_1,y_2]$ из интервала $J$. Это автоматически означает, что $f^{-1}$
непрерывна на $J$.
 \epr

\brem\label{REM:Teor-ob-obratnoi-funktsii-na-poluint-v-R^1} Утверждения,
аналогичные теоремам \ref{Teor-ob-obratnoi-funktsii-na-R^1} и
\ref{Teor-ob-obratnoi-funktsii-na-int-v-R^1} справедливы и для полуинтервалов.
Мы предлагаем читателю самостоятельно их сформулировать и доказать. \erem

\section{Предел функции}\label{sec-pred-funct}

В примере \ref{EX:drob-chast-razryvna} выше мы заметили, что функция
$f(x)=\{x\}$ (дробная часть числа) разрывна в точке $x=1$ (да и вообще во всех
целых точках $a\in\Z$), то есть для некоторых последовательностей аргументов
$x_n$ выполняются соотношения
$$
x_n\underset{n\to\infty}{\longrightarrow} 1,\qquad
f(x_n)\underset{n\to\infty}{\notarrow} f(1)
$$
Если приглядеться к этой функции, то можно заметить, что в поведении
последовательности $f(x_n)$ имеется следующая простая закономерность: если
$x_n$ стремится к 1 слева, то $f(x_n)$ стремится к 1,
 \beq\label{lim_(x->1-0)-drob-chast-x}
x_n\underset{n\to\infty}{\longrightarrow} 1,\qquad
x_n<1\qquad\Longrightarrow\qquad f(x_n)\underset{n\to\infty}{\longrightarrow}
1,
 \eeq
а если $x_n$ стремится к 1 справа, то $f(x_n)$ стремится к 0:
 \beq\label{lim_(x->1+0)-drob-chast-x}
x_n\underset{n\to\infty}{\longrightarrow} 1,\qquad
1<x_n\qquad\Longrightarrow\qquad f(x_n)\underset{n\to\infty}{\longrightarrow}
0.
 \eeq
Это наблюдение оказывается очень полезным, потому что позволяет просто и
наглядно объяснить, как ведет себя наша функция вблизи своих точек разрыва:
если сказать, что при приближении аргумента к точке разрыва слева функция
стремится к единице, а при приближении аргумента справа -- к нулю, то изобразив
это графически картинкой

\vglue80pt

\noindent мы даем человеческому сознанию наглядный образ, с помощью которого
все вопросы о пределах последовательностей вида $f(x_n)$ решаются без труда.

Условие \eqref{lim_(x->1-0)-drob-chast-x} коротко записывается формулой
$$
\lim_{x\to 1-0}f(x)=1,
$$
и в таких случаях говорят, что функция $f$ имеет предел слева в точке 1. А
условие \eqref{lim_(x->1+0)-drob-chast-x} записывается формулой
$$
\lim_{x\to 1+0}f(x)=0,
$$
и в таких случаях говорят, что функция $f$ имеет предел справа в точке 1,
равный 0.

В этом параграфе мы подробно обсудим понятие предела функции. Ниже в главе
\ref{ch-f'(x)} будет дано его важное применение: с его помощью определяется
понятие {\it производной}, играющие центральную роль во всем математическом
анализе.

\subsection{Определение предела функции}\label{subsec-opred-pred-funct}

Когда в $\S 2$ главы \ref{ch-th-x_n} мы определяли предел последовательности,
мы приводили два определения: отдельно для конечного и бесконечного предела. В
случае с пределом функции ситуация многократно усложняется из-за того, что
аргумент функции $x$ и ее значение $f(x)$ могут стремиться в различных
ситуациях к конечным величинам или бесконечности. В соответствии с этим
различаются четыре основных вида пределов функции, причем каждый из них имеет
два дополнительных ``подвида'', потому что аргумент $x$ предполагается иногда
стремящимся к $a$ с одной стороны.

Представление об этой картине дает следующая таблица:

\medskip

\vbox{\tabskip=0pt\offinterlineskip \halign to \hsize{ \vrule#\tabskip=2pt
plus3pt minus1pt & \strut\hfil\;#\hfil & \vrule#&\hfil#\hfil &
\vrule#&\hfil#\hfil & \vrule#\tabskip=0pt\cr \noalign{\hrule}
height2pt&\omit&&\omit&&\omit&\cr \noalign{\hrule}
height2pt&\omit&&\omit&&\omit&\cr & \omit && \omit && \omit & \cr & \omit &&
Функция стремится && Функция стремится & \cr & ТИПЫ ПРЕДЕЛОВ ФУНКЦИИ && к
конечной величине: && к бесконечности: & \cr height2pt&\omit&&\omit&&\omit&\cr
& \omit && $f(x)\to A$ && $f(x)\to \infty$ & \cr & \omit && \omit && \omit &
\cr height2pt&\omit&&\omit&&\omit&\cr \noalign{\hrule}
height2pt&\omit&&\omit&&\omit&\cr \noalign{\hrule}
height2pt&\omit&&\omit&&\omit& \cr & Аргумент стремится & & Конечный предел & &
Бесконечный предел &\cr & к конечной величине: & & в точке: & & в точке: &\cr &
$x\to a$ & & $\lim\limits_{x\to a} f(x)=A$ & & $\lim\limits_{x\to a}
f(x)=\infty$ &\cr height2pt&\omit&&\omit&&\omit&\cr \noalign{\hrule} &
возможно, с одной стороны: & & \omit & & \omit &\cr & \omit & & Конечный предел
& & Бесконечный предел &\cr & слева & & в точке $a$ слева: & & в точке $a$
слева: &\cr & $x\to a-0$ & & $\lim\limits_{x\to a-0} f(x)=A$ & &
$\lim\limits_{x\to a+0} f(x)=\infty$ & \cr & \omit & & Конечный предел & &
Бесконечный предел & \cr & или справа && в точке $a$ справа: && в точке $a$
справа: & \cr & $x\to a+0$ && $\lim\limits_{x\to a+0} f(x)=A$ &&
$\lim\limits_{x\to a+0} f(x)=\infty$ &\cr height2pt&\omit&&\omit&&\omit&\cr
\noalign{\hrule} height2pt&\omit&&\omit&&\omit&\cr \noalign{\hrule} & Аргумент
стремится & & Конечный предел & & Бесконечный предел &\cr & к бесконечности: &
& в бесконечности: & & в бесконечности: &\cr & $x\to \infty$ & &
$\lim\limits_{x\to \infty} f(x)=A$ & & $\lim\limits_{x\to \infty} f(x)=\infty$
&\cr height2pt&\omit&&\omit&&\omit&\cr \noalign{\hrule} & возможно, в одну
сторону: & & \omit & & \omit &\cr & \omit & & Конечный предел & & Бесконечный
предел &\cr & влево & & в минус бесконечности: & & в минус бесконечности: &\cr
& $x\to -\infty$ & & $\lim\limits_{x\to -\infty} f(x)=A$ & & $\lim\limits_{x\to
-\infty} f(x)=\infty$ & \cr & \omit & & Конечный предел & & Бесконечный предел
&\cr & или вправо & & в плюс бесконечности: & & в плюс бесконечности: &\cr &
$x\to +\infty$ && $\lim\limits_{x\to +\infty} f(x)=A$ && $\lim\limits_{x\to
+\infty} f(x)=\infty$ &\cr height2pt&\omit&&\omit&&\omit&\cr \noalign{\hrule}
height2pt&\omit&&\omit&&\omit&\cr \noalign{\hrule} } }

\medskip

К этим 12 типам пределов следует добавть еще 12, соответствующих случаям, когда
функция $f$ стремится к бесконечности с определенным знаком:

\medskip

\vbox{\tabskip=0pt\offinterlineskip \halign to \hsize{ \vrule#\tabskip=2pt
plus3pt minus1pt & \strut\hfil\;#\hfil & \vrule#&\hfil#\hfil &
\vrule#&\hfil#\hfil & \vrule#\tabskip=0pt\cr \noalign{\hrule}
height2pt&\omit&&\omit&&\omit&\cr \noalign{\hrule}
height2pt&\omit&&\omit&&\omit&\cr & \omit && \omit && \omit & \cr & \omit &&
Функция стремится && Функция стремится & \cr & ТИПЫ ПРЕДЕЛОВ ФУНКЦИИ && к минус
бесконечности: && к плюс бесконечности: & \cr height2pt&\omit&&\omit&&\omit&\cr
& \omit && $f(x)\to -\infty$ && $f(x)\to +\infty$ & \cr & \omit && \omit &&
\omit & \cr height2pt&\omit&&\omit&&\omit&\cr \noalign{\hrule}
height2pt&\omit&&\omit&&\omit&\cr \noalign{\hrule}
height2pt&\omit&&\omit&&\omit& \cr & Аргумент стремится & & Бесконечный предел
& & Бесконечный предел &\cr & к конечной величине: & & в точке: & & в точке:
&\cr & $x\to a$ & & $\lim\limits_{x\to a} f(x)=-\infty$ & & $\lim\limits_{x\to
a} f(x)=+\infty$ &\cr height2pt&\omit&&\omit&&\omit&\cr \noalign{\hrule} &
возможно, с одной стороны: & & \omit & & \omit &\cr & \omit & & Бесконечный
предел & & Бесконечный предел &\cr & слева & & в точке $a$ слева: & & в точке
$a$ слева: &\cr & $x\to a-0$ & & $\lim\limits_{x\to a-0} f(x)=-\infty$ & &
$\lim\limits_{x\to a+0} f(x)=+\infty$ & \cr & \omit & & Бесконечный предел & &
Бесконечный предел & \cr & или справа && в точке $a$ справа: && в точке $a$
справа: & \cr & $x\to a+0$ && $\lim\limits_{x\to a+0} f(x)=-\infty$ &&
$\lim\limits_{x\to a+0} f(x)=+\infty$ &\cr height2pt&\omit&&\omit&&\omit&\cr
\noalign{\hrule} height2pt&\omit&&\omit&&\omit&\cr \noalign{\hrule} & Аргумент
стремится & & Бесконечный предел & & Бесконечный предел &\cr & к бесконечности:
& & в бесконечности: & & в бесконечности: &\cr & $x\to \infty$ & &
$\lim\limits_{x\to \infty} f(x)=-\infty$ & & $\lim\limits_{x\to \infty}
f(x)=+\infty$ &\cr height2pt&\omit&&\omit&&\omit&\cr \noalign{\hrule} &
возможно, в одну сторону: & & \omit & & \omit &\cr & \omit & & Бесконечный
предел & & Бесконечный предел &\cr & влево & & в минус бесконечности: & & в
минус бесконечности: &\cr & $x\to -\infty$ & & $\lim\limits_{x\to -\infty}
f(x)=-\infty$ & & $\lim\limits_{x\to -\infty} f(x)=+\infty$ & \cr & \omit & &
Бесконечный предел & & Бесконечный предел &\cr & или вправо & & в плюс
бесконечности: & & в плюс бесконечности: &\cr & $x\to +\infty$ &&
$\lim\limits_{x\to +\infty} f(x)=-\infty$ && $\lim\limits_{x\to +\infty}
f(x)=+\infty$ &\cr height2pt&\omit&&\omit&&\omit&\cr \noalign{\hrule}
height2pt&\omit&&\omit&&\omit&\cr \noalign{\hrule} } }

\medskip

Таким образом, чтобы определить предел функции, нужно формально сформулировать
24 разных определения. Чтобы не утомлять читателя такими вразумлениями, обычно
поступают иначе. Говорят так: пусть $a$ и $A$ обозначают числа, или символы
бесконечности, возможно с определенным знаком
$$
a, A \in \R \quad \text{или}\quad a,A=\infty, \, -\infty, \, +\infty
$$
и пусть для всякого такого $a$
 \bit{
\item[---] {\it проколотой окрестностью} величины $a$ называется всякое
множество вида
$$
U=\left\{ \begin{matrix} (a-\delta,a)\cup (a,a+\delta), & \text{где}\, \delta>0
&
\Big(\text{если $a$ -- число}\Big) \\
(-\infty,-E)\cup (E,+\infty), & \text{где}\, E>0  &
\Big(\text{если}\, a=\infty \Big);\\
(-\infty,-E), & \text{где}\, E>0 &
\Big(\text{если}\, a=-\infty \Big); \\
(E,+\infty), & \text{где}\, E>0 & \Big(\text{если}\, a=+\infty \Big);
\end{matrix}\right\}
$$
\item[---] {\it левой полуокрестностью} $a$ называется всякий интервал вида
$$
U=(a-\delta,a), \quad \text{где}\, \delta>0 \quad \Big(\text{если $a$ --
число}\Big)
$$
\item[---] {\it правой полуокрестностью} $a$ называется всякий интервал вида
$$
U=(a,a+\delta), \quad \text{где}\, \delta>0 \quad \Big(\text{если $a$ --
число}\Big)
$$
 }\eit
Далее следует определение, впервые сформулированное немецким математиком
Генрихом Гейне:\footnote{Генрих Гейне (1821-1881) -- известный немецкий
математик, поэт и публицист (1797-1856).}

\bit{ \item[$\bullet$]\label{Opred-predela-funktsii} Величина $A$ называется
{\it пределом функции}\index{предел!функции!по Гейне} $f(x)$ при $x$
стремящимся к $a$
$$
A=\lim_{x\to a} f(x) \quad \l f(x)\underset{x\to a}{\longrightarrow} A \r
$$
если
 \bit{
\item[1)] $f(x)$ определена в некоторой выколотой окрестности $U$ величины $a$
и

\item[2)] для любой последовательности $x_n\in U$, стремящейся к $a$
$$
  x_n\underset{n\to\infty}{\longrightarrow} a \quad \Big( x_n\ne a \Big)
$$
соответствующая последовательность значений $f(x_n)$ стремится к $A$:
$$
  f(x_n)\underset{n\to\infty}{\longrightarrow} A
$$
 }\eit
}\eit

Отдельно определяются односторонние пределы:

\bit{ \item[$\bullet$] Величина $A$ называется {\it пределом слева (справа)
функции} $f(x)$ в точке $a$
$$
  A=\lim_{x\to a-0} f(x) \quad \l A=\lim_{x\to a+0} f(x) \r
$$
если
 \bit{
\item[1)] $f(x)$ определена в некоторой левой (правой) полуокрестности $U$
точки $a$ и

\item[2)] для любой последовательности $x_n\in U$, стремящейся к $a$
$$
  x_n\underset{n\to\infty}{\longrightarrow} a \quad \Big( x_n\ne a \Big)
$$
соответствующая последовательность значений $f(x_n)$ стремится к $A$:
$$
  f(x_n)\underset{n\to\infty}{\longrightarrow} A
$$
 }\eit
}\eit

\noindent\rule{160mm}{0.1pt}\begin{multicols}{2}

Перепишем эти общие определения для некоторых конкретных ситуаций.

\begin{ex}[\bf конечный предел функции в точке по Гейне]
Применительно к случаю конечного предела функции в точке наше определение можно
сформулировать следующим образом.

Число $A$ называется {\it пределом} функции $f$ в точке $a\in \R$, если
 \biter{
\item[1)] функция $f$ определена на некотором множестве вида $(a-\delta; a)\cup
(a; a+\delta)$, и \item[2)] для всякой последовательности аргументов $\{ x_n
\}$, стремящейся к точке $a$, но не попадающей в $a$
$$
x_n \underset{n\to \infty}{\longrightarrow} a \qquad x_n \in (a-\delta; a)\cup
(a; a+\delta)
$$
соответствующая последовательность значений функции $\{ f(x_n) \}$ стремится к
числу $A$:
$$
f(x_n) \underset{n\to \infty}{\longrightarrow} A
$$
 }\eiter

%\picture{0pt}{0pt}{86.pcx}

\vglue140pt \noindent В этом случае пишут
$$
f(x)\underset{x\to a}{\longrightarrow} A
$$
или
$$
\lim_{x\to a} f(x)=A
$$
\end{ex}

\begin{ex}[\bf конечный предел функции в точке слева по Гейне]
Число $A$ называется {\it пределом слева} функции $f$ в точке $a$, если
 \biter{
\item[1)] функция $f$ определена на некотором интервале вида $(a-\delta; a)$, и
\item[2)] для всякой последовательности аргументов $\{ x_n \}$, стремящейся
слева к точке $a$
$$
x_n \underset{n\to \infty}{\longrightarrow} a \qquad x_n \in (a-\delta; a)
$$
соответствующая последовательность значений функции $\{ f(x_n) \}$ стремится к
числу $A$:
$$
f(x_n) \underset{n\to \infty}{\longrightarrow} A
$$
 }\eiter

%\picture{0pt}{0pt}{84.pcx}

\vglue140pt \noindent В этом случае пишут
$$
f(x)\underset{x\to a-0}{\longrightarrow} A
$$
или
$$
\lim_{x\to a-0} f(x)=A
$$
\end{ex}

\begin{ex}[\bf бесконечный предел функции в точке по Гейне]
Говорят, что функция $f$ {\it стремится к} $\infty$) в точке $a$, если
 \biter{
\item[1)] функция $f$ определена на некотором множестве вида $(a-\delta; a)\cup
(a;a+\delta)$, и \item[2)] для всякой последовательности аргументов $\{ x_n
\}$, стремящейся к точке $a$, но не попадающей в $a$
$$
x_n \underset{n\to \infty}{\longrightarrow} a \qquad x_n \in (a-\delta; a)\cup
(a;a+\delta)
$$
соответствующая последовательность значений функции $\{ f(x_n) \}$ стремится к
$\infty$:
$$
f(x_n) \underset{n\to \infty}{\longrightarrow}\infty
$$
 }\eiter
%\picture{0pt}{0pt}{89.pcx}

\vglue100pt \noindent В этом случае пишут
$$
f(x)\underset{x\to a}{\longrightarrow} +\infty \qquad (\text{или $-\infty$, или
$\infty$})
$$
или
$$
\lim_{x\to a} f(x)=+\infty \qquad (\text{или $-\infty$, или $\infty$})
$$
\end{ex}

\begin{ex}[\bf конечный предел функции в бесконечности по Гейне]
Число $A$ называется {\it пределом функции $f$ при $x\to +\infty$}, если
 \biter{
\item[1)] функция $f$ определена на некотором интервале вида $(E;+\infty)$, и
\item[2)] для всякой последовательности аргументов $\{ x_n \}\subset
(\beta;+\infty)$, стремящейся к $+\infty$
$$
x_n \underset{n\to \infty}{\longrightarrow} +\infty
$$
соответствующая последовательность значений функции $\{ f(x_n) \}$ стремится к
числу $A$:
$$
f(x_n) \underset{n\to \infty}{\longrightarrow} A
$$
 }\eiter

%\picture{0pt}{0pt}{90.pcx}

\vglue140pt \noindent В этом случае пишут
$$
f(x)\underset{x\to +\infty}{\longrightarrow} A
$$
или
$$
\lim_{x\to +\infty} f(x)=A
$$
\end{ex}

\begin{er}
Сформулируйте определения предела функции для остальных ситуаций, описанных в
двух таблицах этого параграфа, и приведите графические иллюстрации.
\end{er}

\begin{er}
Нарисуйте график какой-нибудь функции $f:\R\to \R$, удовлетворяющей следующим
условиям (всем сразу):
 \biter{
\item[1)] $\lim\limits_{x\to -\infty} f(x)=-2$; \item[2)] $\lim\limits_{x\to
-1} f(x)=0$; \item[3)] $\lim\limits_{x\to 1-0} f(x)=-\infty$; \item[4)]
$\lim\limits_{x\to 1+0} f(x)=2$; \item[5)] $\lim\limits_{x\to +\infty}
f(x)=+\infty$;
 }\eiter
\end{er}

\end{multicols}\noindent\rule[10pt]{160mm}{0.1pt}

\subsection{Нахождение предела функции по определению}

\noindent\rule{160mm}{0.1pt}\begin{multicols}{2}

\begin{ex}
Покажем что
$$
\lim_{x\to 0}\frac{1}{x}=\infty
$$
Действительно, если взять последовательность аргументов $\{ x_n \}$,
стремящуюся к нулю, но не попадающую в ноль
$$
x_n \underset{n\to \infty}{\longrightarrow}  0 \qquad x_n \ne 0
$$
то мы получим
$$
f(x_n)=\frac{1}{x_n}\underset{n\to \infty}{\longrightarrow}\infty
$$
Это верно для всякой последовательности $x_n \underset{n\to
\infty}{\longrightarrow} 0\quad x_n \ne 0 $, поэтому $\lim\limits_{x\to
0}\frac{1}{x}=\infty$.
\end{ex}

\begin{ex}
Покажем что
$$
\lim_{x\to \infty}\frac{x+2}{2x-1}=\frac{1}{2}
$$
Действительно, если взять последовательность аргументов $\{ x_n \}$,
стремящуюся к бесконечности
$$
x_n \underset{n\to \infty}{\longrightarrow}  \infty
$$
то мы получим
$$
f(x_n)=\frac{x_n+2}{2x_n-1}=\frac{1+\frac{2}{x_n}}{2-\frac{1}{x_n}}\underset{n\to
\infty}{\longrightarrow}\frac{1}{2}
$$
Это верно для всякой последовательности $x_n \underset{n\to
\infty}{\longrightarrow}\infty$, поэтому $\lim\limits_{x\to
0}\frac{x+2}{2x-1}=\frac{1}{2}$.
\end{ex}

\begin{ex}
Покажем что
$$
\lim_{x\to 0-0}\frac{1}{x}=-\infty \qquad \lim_{x\to 0+0}\frac{1}{x}=+\infty
$$
Действительно, если взять последовательность аргументов $\{ x_n \}$,
стремящуюся к 0 слева,
$$
x_n \underset{n\to \infty}{\longrightarrow}  0 \qquad x_n < 0
$$
то мы получим
$$
\lim_{n\to \infty}\frac{1}{x_n}={\smsize (\text{применяем теорему
\ref{0<->infty}})}=-\infty
$$
Это верно для всякой последовательности $x_n \underset{n\to
\infty}{\longrightarrow} 0\quad x_n < 0 $, поэтому $\lim\limits_{x\to
0-0}\frac{1}{x}=-\infty$.

Аналогично, если
$$
x_n \underset{n\to \infty}{\longrightarrow} 0 \qquad x_n > 0
$$
то мы получим
$$
\lim_{n\to \infty}\frac{1}{x_n}= {\smsize (\text{применяем теорему
\ref{0<->infty}})}=+\infty
$$
Это верно для всякой последовательности $x_n \underset{n\to
\infty}{\longrightarrow} 0\quad x_n > 0 $, поэтому $\lim\limits_{x\to
0+0}\frac{1}{x}=+\infty$.
\end{ex}

\begin{ex}
Покажем что
$$
\lim_{x\to +\infty}\frac{x^2+1}{|x|(x+1)}=1 \qquad \lim_{x\to
-\infty}\frac{x^2+1}{|x|(x+1)}=-1
$$
Действительно, если взять последовательность аргументов $\{ x_n \}$,
стремящуюся к $+\infty$,
$$
x_n \underset{n\to \infty}{\longrightarrow}  +\infty
$$
то мы получим
 \begin{multline*}
\lim_{n\to \infty} f(x_n)=\lim_{n\to \infty}\frac{x_n^2+1}{|x_n|(x_n+1)}=\\=
 {\smsize \begin{pmatrix} |x_n|=x_n, \\
 \text{потому что}\\ \text{почти все}\,\, x_n>0
 \end{pmatrix}}=
  \underbrace{\lim_{n\to \infty}\frac{x_n^2+1}{x_n (x_n+1)}}_
  {\smsize \begin{matrix}\uparrow \\
  \text{делим числитель}\\ \text{и знаменатель на}
 \, \, x_n^2 \end{matrix}}=\\=
 \lim_{n\to \infty}\frac{1+\frac{1}{x_n^2}}{1+\frac{1}{x_n}}
 ={\smsize
 \begin{pmatrix}\text{применяем} \\ \text{теорему \ref{0<->infty}}
 \end{pmatrix}}=
 \frac{1+0}{1+0}=1
\end{multline*} Это верно для всякой последовательности
$x_n \underset{n\to \infty}{\longrightarrow} +\infty$, поэтому
$\lim\limits_{x\to +\infty}\frac{x^2+1}{|x|(x+1)}=1$.

Аналогично, если
$$
x_n \underset{n\to \infty}{\longrightarrow} -\infty
$$
то мы получим
 \begin{multline*}
\lim_{n\to \infty} f(x_n)= \lim_{n\to \infty}\frac{x_n^2+1}{|x_n|(x_n+1)}=\\=
 {\smsize \begin{pmatrix} |x_n|=-x_n, \\
 \text{потому что}\\ \text{почти все}\,\, x_n<0 \end{pmatrix}}=
 \underbrace{\lim_{n\to \infty}\frac{x_n^2+1}{-x_n (x_n+1)}}_{\smsize
 \begin{matrix}\uparrow \\ \text{делим числитель}\\ \text{и знаменатель на}
 \; x_n^2 \end{matrix}}=\\=
\lim_{n\to \infty}\frac{1+\frac{1}{x_n^2}}{-1-\frac{1}{x_n}}
 ={\smsize\begin{pmatrix}\text{применяем}\\
 \text{теорему \ref{0<->infty}}\end{pmatrix}}=
 \frac{1+0}{-1-0}=-1
 \end{multline*}
 Это верно для всякой последовательности $x_n \underset{n\to
\infty}{\longrightarrow} -\infty$, поэтому $\lim\limits_{x\to
-\infty}\frac{x^2+1}{|x|(x+1)}=1$.
\end{ex}

\begin{ex}
Покажем что
$$
\lim_{x\to +\infty}\frac{1}{x^3}=0
$$
Действительно, если взять последовательность аргументов $\{ x_n \}$,
стремящуюся к $\infty$,
$$
x_n \underset{n\to \infty}{\longrightarrow}  \infty
$$
то мы получим
$$
\lim_{n\to \infty} f(x_n)=\lim_{n\to \infty}\left(\frac{1}{x_n}
 \right)^3
={\smsize \begin{pmatrix}\text{применяем} \\ \text{теорему
\ref{0<->infty}}\end{pmatrix}}= 0
$$
Это верно для всякой последовательности $x_n \underset{n\to
\infty}{\longrightarrow}\infty$, поэтому $\lim\limits_{x\to
\infty}\frac{1}{x^3}=0$.
\end{ex}

\begin{ex}
Покажем что функция $f(x)=\cos x$ не имеет предела при $x\to +\infty$:
$$
\nexists \lim_{x\to +\infty}\cos x
$$
Возьмем сначала последовательность аргументов
$$
x_n =2\pi n \underset{n\to \infty}{\longrightarrow} +\infty
$$
Тогда мы получим
$$
f(x_n)=\cos x_n= \cos 2\pi n=1 \underset{n\to \infty}{\longrightarrow} 1
$$
С другой стороны, если взять
$$
x_n =\pi+2\pi n \underset{n\to \infty}{\longrightarrow} +\infty
$$
то мы получим
$$
f(x_n)=\cos x_n= \cos (\pi+2\pi n)=-1 \underset{n\to \infty}{\longrightarrow}
-1
$$
Таким образом, мы получаем, что если взять одну последовательность аргументов
$x_n \underset{n\to \infty}{\longrightarrow} +\infty$, то соотвествующая
последовательность $f(x_n)$ будет стремиться к одному числу ($A=1$), а если
взять другую последовательность $x_n \underset{n\to \infty}{\longrightarrow}
+\infty$, то $f(x_n)$ будет стремиться к другому числу ($A=-1$).

Это означает, что функция $f(x)=\cos x$ не имеет предела при $x\to +\infty$,
потому что иначе $f(x_n)$ должно было бы стремиться к этому конкретному пределу
$A$, независимо от того, какую берешь последовательность $x_n \underset{n\to
\infty}{\longrightarrow} +\infty$.
\end{ex}

\begin{ex}
Покажем что
$$
\lim_{x\to +\infty} |x|x=+\infty \qquad \lim_{x\to -\infty} |x|x=-\infty
$$
Действительно, если взять последовательность аргументов $\{ x_n \}$,
стремящуюся к $+\infty$,
$$
x_n \underset{n\to \infty}{\longrightarrow}  +\infty
$$
то мы получим
 \begin{multline*}
\lim_{n\to \infty} f(x_n)=\lim_{n\to \infty} |x_n| x_n=\\=
 {\smsize
 \begin{pmatrix} |x_n|=x_n, \\ \text{потому что}\\ \text{почти
все}\,\, x_n>0
 \end{pmatrix}}=\lim_{n\to \infty} x_n^2 =+\infty
 \end{multline*}
Это верно для всякой последовательности $x_n \underset{n\to
\infty}{\longrightarrow} +\infty$, поэтому $\lim\limits_{x\to +\infty}
|x|x=+\infty$.

Возьмем теперь последовательность аргументов $\{ x_n \}$, стремящуюся к
$-\infty$,
$$
x_n \underset{n\to \infty}{\longrightarrow}  -\infty
$$
Тогда
 \begin{multline*}
\lim_{n\to \infty} f(x_n)=\lim_{n\to \infty} |x_n| x_n=\\=
 {\smsize \begin{pmatrix} |x_n|=-x_n, \text{потому что}\\
\text{почти все}\,\, x_n<0
\end{pmatrix}}=-\lim_{n\to \infty} x_n^2 =-\infty
 \end{multline*}
Это верно для всякой последовательности $x_n \underset{n\to
\infty}{\longrightarrow} -\infty$, поэтому $\lim\limits_{x\to -\infty}
|x|x=-\infty$.
\end{ex}

\end{multicols}\noindent\rule[10pt]{160mm}{0.1pt}

\begin{tm}[\bf о связи между нулевым и бесконечным пределами]\label{TH:0<->infty-x}
Предел функции $f$ при $x\to a$ равен нулю в том и только в том случае, если
предел обратной функции $\frac{1}{f}$ при $x\to a$ равен бесконечности:
 \beq\label{0<->infty-x}
\lim_{x\to a}f(x)=0\qquad\Longleftrightarrow\qquad \lim_{x\to
a}\frac{1}{f(x)}=\infty.
 \eeq
\end{tm}
\bpr Это следует из теоремы \ref{0<->infty} о связи между бесконечно большими и
бесконечно малыми последовательностями: если $x_n\to a$, то соотношение
$y_n=f(x_n)\to 0$ эквивалентно соотношению $\frac{1}{y_n}=\frac{1}{f(x_n)}\to
\infty$.
 \epr

\noindent\rule{160mm}{0.1pt}\begin{multicols}{2}

 \bex\label{EX:x^k->0(x->0)}
 \beq\label{x^k->0(x->0)}
\boxed{\quad\lim_{x\to 0} x^k=
 \begin{cases}0, & \text{\rm если $k>0$}\\ \infty, & \text{\rm если
 $k<0$}\end{cases}\quad}
 \eeq
\eex
 \bpr
1. Пусть $k>0$. Тогда
$$
x_n\underset{x\to 0}{\longrightarrow}0
$$
$$
\phantom{\text{\scriptsize (свойство $3^0$ на с.\pageref{lim-ariphm}})}\quad
\Downarrow\quad\text{\scriptsize (свойство $3^0$ на с.\pageref{lim-ariphm})}
$$
$$
(x_n)^k\underset{x\to0}{\longrightarrow}0
$$
Поскольку это верно для любой последовательности $x_n\to 0$, мы получаем, что
$$
x^k\underset{x\to 0}{\longrightarrow}0 \qquad (k>0)
$$

2. Пусть $k<0$. Тогда $m=-k>0$, поэтому:
$$
\phantom{\text{\scriptsize (уже доказано)}}\quad x^m\underset{x\to
0}{\longrightarrow}0\quad\text{\scriptsize (уже доказано)}
$$
$$
\phantom{\text{\scriptsize \eqref{0<->infty-x}}}\quad
\Downarrow\quad\text{\scriptsize \eqref{0<->infty-x}}
$$
$$
x^k=\frac{1}{x^m}\underset{x\to0}{\longrightarrow}\infty
$$
 \epr

\bex\label{EX:x^k->8(x->8)}
 \beq\label{x^k->8(x->8)}
\boxed{\quad\lim_{x\to\infty} x^k=
 \begin{cases}\infty, & \text{\rm если $k>0$}\\ 0, & \text{\rm если
 $k<0$}\end{cases}\quad}
 \eeq
\eex
 \bpr
1. Пусть $k>0$ и
$$
x_n\underset{x\to 0}{\longrightarrow}\infty
$$
Обозначим $y_n=\frac{1}{x_n}$. Тогда
$$
\phantom{\text{\scriptsize (теорема \ref{0<->infty})}}\quad y_n\underset{x\to
0}{\longrightarrow}0\quad\text{\scriptsize (теорема \ref{0<->infty})}
$$
$$
\phantom{\text{\scriptsize (свойство $3^0$ на с.\pageref{lim-ariphm}})}\quad
\Downarrow\quad\text{\scriptsize (свойство $3^0$ на с.\pageref{lim-ariphm})}
$$
$$
(y_n)^k\underset{x\to0}{\longrightarrow}0
$$
$$
\phantom{\text{\scriptsize (теорема \ref{0<->infty})}}\quad
\Downarrow\quad\text{\scriptsize (теорема \ref{0<->infty})}
$$
$$
(x_n)^k=\frac{1}{(y_n)^k}\underset{x\to 0}{\longrightarrow}\infty
$$
Поскольку это верно для любой последовательности $x_n\to\infty$, мы получаем,
что
$$
x^k\underset{x\to\infty}{\longrightarrow}\infty \qquad (k>0)
$$

2. Пусть $k<0$. Тогда $m=-k>0$, поэтому:
$$
\phantom{\text{\scriptsize (уже доказано)}}\quad x^m\underset{x\to
\infty}{\longrightarrow}\infty\quad\text{\scriptsize (уже доказано)}
$$
$$
\phantom{\text{\scriptsize \eqref{0<->infty-x}}}\quad
\Downarrow\quad\text{\scriptsize \eqref{0<->infty-x}}
$$
$$
x^k=\frac{1}{x^m}\underset{x\to\infty}{\longrightarrow}0
$$
 \epr

\bex Для любого $n\in\N$
 \begin{align}
&\lim_{x\to-\infty} x^{2n}=+\infty \label{x^(2n)->8(x->-8)} \\
&\lim_{x\to+\infty} x^{2n}=+\infty \label{x^(2n)->8(x->+8)}
 \end{align}
 \eex

\bex Для любого $n\in\N$
 \begin{align}
&\lim_{x\to0} x^{-2n}=+\infty \label{x^(-2n)->+8(x->0)} \\
&\lim_{x\to\infty} x^{-2n}=0 \label{x^(-2n)->0(x->8)}
 \end{align}
 \eex

\bex Для любого $n\in\N$
 \begin{align}
&\lim_{x\to-\infty} x^{2n-1}=-\infty \label{x^(2n-1)->-8(x->-8)}\\
&\lim_{x\to+\infty} x^{2n-1}=+\infty \label{x^(2n-1)->+8(x->+8)}
 \end{align}
 \eex

\bex Для любого $n\in\N$
 \begin{align}
&\lim_{x\to-0} x^{-(2n-1)}=-\infty \label{x^(2n-1)->-8(x->-8)}\\
&\lim_{x\to+0} x^{-(2n-1)}=+\infty \label{x^(2n-1)->+8(x->+8)} \\
&\lim_{x\to \infty} x^{-(2n-1)}=0 \label{x^(2n-1)->-8(x->-8)}
 \end{align}
 \eex

\begin{ers} Покажите что
 \biter{
\item[1)] $\lim\limits_{x\to +\infty}\frac{|x|}{x}=1$

\item[2)] $\lim\limits_{x\to -\infty}\frac{|x|}{x}=-1$
 }\eiter
\end{ers}

\end{multicols}\noindent\rule[10pt]{160mm}{0.1pt}

\subsection{Связь между односторонними и двусторонними
пределами}

\begin{tm}
Величина $A$ будет пределом функции $f$ в точке $a$
$$
  A=\lim_{x\to a} f(x)
$$
тогда и только тогда, когда $A$ является пределом слева и пределом справа
функции $f$ в точке $a$:
$$
  \lim_{x\to a-0} f(x)=A=\lim_{x\to a+0} f(x)
$$
\end{tm}\begin{proof} 1. Если $A=\lim\limits_{x\to a} f(x)$, то для
любой последовательности $x_n\underset{n\to \infty}{\longrightarrow} a, \,
x_n\ne a$ должно выполняться $f(x_n)\underset{n\to \infty}{\longrightarrow} A$.
В частности, если $x_n\underset{n\to \infty}{\longrightarrow} a, \, x_n<a$ то
$f(x_n)\underset{n\to \infty}{\longrightarrow} A$, и поэтому
$A=\lim\limits_{x\to a-0} f(x)$. Аналогично получается $A=\lim\limits_{x\to
a+0} f(x)$.

2. Наоборот, если $\lim\limits_{x\to a-0} f(x)=A=\lim\limits_{x\to a+0} f(x)$,
то это означает, что
$$
\forall x_n\underset{n\to \infty}{\longrightarrow} a, \, x_n< a \quad
\Rightarrow f(x_n)\underset{n\to \infty}{\longrightarrow} A
$$
и
$$
\forall x_n\underset{n\to \infty}{\longrightarrow} a, \, x_n>a \quad
\Rightarrow f(x_n)\underset{n\to \infty}{\longrightarrow} A
$$
Поэтому если взять последовательность $x_n\underset{n\to
\infty}{\longrightarrow} a, \, x_n\ne a$ то, разложив ее на две
подпоследовательности
$$
\Big( x_{n_k}\underset{k\to \infty}{\longrightarrow} a, \, x_{n_k}<a \Big),
\quad \Big( x_{m_k}\underset{k\to \infty}{\longrightarrow} a, \, x_{m_k}>a
\Big)
$$
мы получим
$$
f(x_{n_k})\underset{k\to \infty}{\longrightarrow} A, \quad
f(x_{m_k})\underset{k\to \infty}{\longrightarrow} A
$$
Это означает, что
$$
f(x_n)\underset{n\to \infty}{\longrightarrow} A,
$$
Поскольку это верно для любой последовательности $x_n\underset{n\to
\infty}{\longrightarrow} a, \, x_n\ne a$, мы получаем $A=\lim\limits_{x\to a}
f(x)$. \end{proof}

\begin{tm}
Величина $A$ будет пределом функции $f$ на бесконечности
$$
  A=\lim_{x\to \infty} f(x)
$$
тогда и только тогда, когда $A$ является пределом функции $f$ при $x$
стремящимся к $-\infty$ и к $+\infty$:
$$
  \lim_{x\to -\infty} f(x)=A=\lim_{x\to +\infty} f(x)
$$
\end{tm}\begin{proof} 1. Если $A=\lim\limits_{x\to \infty} f(x)$, то для
любой последовательности $x_n\underset{n\to \infty}{\longrightarrow}\infty$
должно выполняться $f(x_n)\underset{n\to \infty}{\longrightarrow} A$. В
частности, если $x_n\underset{n\to \infty}{\longrightarrow} -\infty$ то
$f(x_n)\underset{n\to \infty}{\longrightarrow} A$, и поэтому
$A=\lim\limits_{x\to -\infty} f(x)$. Аналогично получается $A=\lim\limits_{x\to
+\infty} f(x)$.

2. Наоборот, если $\lim\limits_{x\to -\infty} f(x)=A=\lim\limits_{x\to +\infty}
f(x)$, то это означает, что
$$
\forall x_n\underset{n\to \infty}{\longrightarrow} -\infty \quad \Rightarrow
f(x_n)\underset{n\to \infty}{\longrightarrow} A
$$
и
$$
\forall x_n\underset{n\to \infty}{\longrightarrow} +\infty \quad \Rightarrow
f(x_n)\underset{n\to \infty}{\longrightarrow} A
$$
Поэтому если взять последовательность $x_n\underset{n\to
\infty}{\longrightarrow}\infty$ то, разложив ее на две подпоследовательности
$$
x_{n_k}\underset{k\to \infty}{\longrightarrow} -\infty, \quad
x_{m_k}\underset{k\to \infty}{\longrightarrow} +\infty,
$$
мы получим
$$
f(x_{n_k})\underset{k\to \infty}{\longrightarrow} A, \quad
f(x_{m_k})\underset{k\to \infty}{\longrightarrow} A
$$
Это означает, что
$$
f(x_n)\underset{n\to \infty}{\longrightarrow} A,
$$
Поскольку это верно для любой последовательности $x_n\underset{n\to
\infty}{\longrightarrow}\infty$, мы получаем $A=\lim\limits_{x\to \infty}
f(x)$. \end{proof}

\noindent\rule{160mm}{0.1pt}\begin{multicols}{2}

\begin{ex}
Покажем что
$$
\lim_{x\to 0-0}\frac{|x|}{x}=-1 \qquad \lim_{x\to 0+0}\frac{|x|}{x}=1
$$
Действительно, если взять последовательность аргументов $\{ x_n \}$,
стремящуюся к 0 слева,
$$
x_n \underset{n\to \infty}{\longrightarrow}  0 \qquad x_n < 0
$$
то мы получим
 \begin{multline*}\lim_{n\to \infty}\frac{|x_n|}{x_n}=
 {\smsize {\smsize\begin{pmatrix}
|x_n|=-x_n, \\ \text{поскольку}\\ x_n<0 \end{pmatrix}}}=
 \lim_{n\to\infty}\frac{-x_n}{x_n}=\\= \lim_{n\to \infty} (-1)=-1
 \end{multline*}
Это верно для всякой последовательности $x_n \underset{n\to
\infty}{\longrightarrow} 0\quad x_n < 0 $, поэтому $\lim\limits_{x\to
0-0}\frac{|x|}{x}=-1$.

Аналогично, если
$$
x_n \underset{n\to \infty}{\longrightarrow} 0 \qquad x_n > 0
$$
то мы получим
 $$
 \lim_{n\to \infty}\frac{|x_n|}{x_n}=
 {\smsize
 {\smsize\begin{pmatrix} |x_n|=x_n, \\ \text{поскольку}\\ x_n>0
 \end{pmatrix}}}=
 \lim_{n\to \infty}\frac{x_n}{x_n}= \lim_{n\to \infty} 1=1
 $$
Это верно для всякой последовательности $x_n \underset{n\to
\infty}{\longrightarrow} 0\quad x_n > 0 $, поэтому $\lim\limits_{x\to
0-0}\frac{|x|}{x}=1$.
\end{ex}

\begin{ex}
Покажем что
$$
\lim_{x\to 1-0}\frac{x^2-1}{|x-1|}=-2 \qquad \lim_{x\to
1+0}\frac{x^2-1}{|x-1|}=2
$$
Действительно, если взять последовательность аргументов $\{ x_n \}$,
стремящуюся к 1 слева,
$$
x_n \underset{n\to \infty}{\longrightarrow} 1 \qquad x_n < 1
$$
то мы получим
 \begin{multline*}\lim_{n\to \infty}\frac{x_n^2-1}{|x_n-1|}= \lim_{n\to \infty}\frac{(x_n+1)(x_n-1)}{|x_n-1|}=\\=
 {\smsize
 {\smsize\begin{pmatrix} |x_n-1|=-(x_n-1),
 \\ \text{поскольку $x_n-1<0$}\end{pmatrix}}}=
 \lim_{n\to \infty}\frac{(x_n+1)(x_n-1)}{-(x_n-1)}=\\= -\lim_{n\to
\infty} (
 \boxed{x_n}\put(-5,-8){\vector(1,-2){5}\put(2,-16){1}}
 +1)= -2
  \end{multline*}
 Это верно для всякой
последовательности $x_n \underset{n\to \infty}{\longrightarrow} 1\quad x_n < 1
$, поэтому $\lim\limits_{x\to 1-0}\frac{x^2-1}{|x-1|}=-2$.

Аналогично, если
$$
x_n \underset{n\to \infty}{\longrightarrow} 1 \qquad x_n > 1
$$
то мы получим
 \begin{multline*}\lim_{n\to \infty}\frac{x_n^2-1}{|x_n-1|}= \lim_{n\to \infty}\frac{(x_n+1)(x_n-1)}{|x_n-1|}=\\=
 {\smsize {\smsize\begin{pmatrix} |x_n-1|=x_n-1,
\\ \text{поскольку $x_n-1>0$}
 \end{pmatrix}}}=\\=
 \lim_{n\to \infty}\frac{(x_n+1)(x_n-1)}{x_n-1}=
 \lim_{n\to \infty}\big(
 {\smsize
 \boxed{x_n}\put(-5,-8){\vector(1,-2){5}\put(2,-16){1}}
 }
 +1\big)=2
  \end{multline*}
 Это верно для всякой последовательности $x_n \underset{n\to
\infty}{\longrightarrow} 1$, $x_n > 1 $, поэтому $\lim\limits_{x\to
1-0}\frac{x^2-1}{|x-1|}=2$.
\end{ex}

\begin{ex}
Покажем что
$$
\lim_{x\to 0}\frac{x^2+x^3}{|x|}=0
$$
Возьмем сначала последовательность аргументов $\{ x_n \}$, стремящуюся к 0
слева,
$$
x_n \underset{n\to \infty}{\longrightarrow} 0 \qquad x_n < 0
$$
тогда
 \begin{multline*}
\lim_{n\to \infty}\frac{x_n^2+x_n^3}{|x_n|}=
 {\smsize
 {\smsize\begin{pmatrix} |x_n|=-x_n, \\ \text{поскольку $x_n<0$}
 \end{pmatrix}}}=\\=
 \lim_{n\to \infty}\frac{x_n^2+x_n^3}{-x_n}=
\lim_{n\to \infty} (-x_n-x_n^2)=0
 \end{multline*}
Это верно для всякой последовательности $x_n \underset{n\to
\infty}{\longrightarrow} 0\quad x_n < 0 $, поэтому $\lim\limits_{x\to
0-0}\frac{x^2+x^3}{|x|}=0$.

Аналогично, если
$$
x_n \underset{n\to \infty}{\longrightarrow} 0 \qquad x_n > 0
$$
то мы получим
 \begin{multline*}\lim_{n\to \infty}\frac{x_n^2+x_n^3}{|x_n|}=
 {\smsize
 {\smsize\begin{pmatrix}
 |x_n|=x_n, \\ \text{поскольку $x_n>0$}
 \end{pmatrix}}}=\\=
 \lim_{n\to \infty}\frac{x_n^2+x_n^3}{x_n}=
\lim_{n\to \infty} (x_n+x_n^2)=0
 \end{multline*}
Это верно для всякой последовательности $x_n \underset{n\to
\infty}{\longrightarrow} 0\quad x_n > 0 $, поэтому $\lim\limits_{x\to
0+0}\frac{x^2+x^3}{|x|}=1$.

Мы получаем, что для функции $f(x)=\frac{x^2+x^3}{|x|}$ левый и правый предел в
точке 0 совпадают, поэтому
$$
\lim_{x\to 0}\frac{x^2+x^3}{|x|}= \lim_{x\to 0-0}\frac{x^2+x^3}{|x|}=\lim_{x\to
0+0}\frac{x^2+x^3}{|x|}=0
$$
\end{ex}

\begin{ex}
Покажем что
$$
\lim_{x\to 0}\frac{1}{|x|}=+\infty
$$
Действительно, если взять последовательность аргументов $\{ x_n \}$,
стремящуюся к 0 слева,
$$
x_n \underset{n\to \infty}{\longrightarrow}  0 \qquad x_n < 0
$$
то мы получим
 \begin{multline*}\lim_{n\to \infty}\frac{1}{|x_n|}=
 {\smsize
 {\smsize\begin{pmatrix}
 |x_n|=-x_n, \\ \text{потому что}\,\, x_n<0
 \end{pmatrix}}}=
 -\lim_{n\to
\infty}\frac{1}{x_n}=\\=
 {\smsize {\smsize\begin{pmatrix}\text{применяем}\\
 \text{теорему \ref{0<->infty}}\end{pmatrix}}}=
 +\infty
 \end{multline*}
Это верно для всякой последовательности $x_n \underset{n\to
\infty}{\longrightarrow} 0\quad x_n < 0 $, поэтому $\lim\limits_{x\to
0-0}\frac{1}{x}=+\infty$.

Аналогично, если
$$
x_n \underset{n\to \infty}{\longrightarrow} 0 \qquad x_n > 0
$$
то мы получим
 \begin{multline*}\lim_{n\to \infty}\frac{1}{|x_n|}
  \quad = {\smsize {\smsize\begin{pmatrix}
 |x_n|=x_n, \\ \text{потому что}\,\, x_n>0
 \end{pmatrix}}}=
 \lim_{n\to \infty}\frac{1}{x_n}=\\=
 {\smsize {\smsize\begin{pmatrix}\text{применяем}\\
 \text{теорему \ref{0<->infty}}\end{pmatrix}}}=
 +\infty
 \end{multline*}
Это верно для всякой последовательности $x_n \underset{n\to
\infty}{\longrightarrow} 0$ $x_n > 0 $, поэтому $\lim\limits_{x\to
0+0}\frac{1}{x}=+\infty$.

Мы получили
$$
\lim_{x\to 0-0}\frac{1}{|x|}=+\infty=\lim_{x\to 0+0}\frac{1}{|x|}
$$
и это означает что
$$
\lim_{x\to 0}\frac{1}{|x|}=+\infty
$$
\end{ex}

\end{multicols}\noindent\rule[10pt]{160mm}{0.1pt}

\subsection{Предел функции и монотонные последовательности}

Определение предела по Коши упрощает доказательство следующего факта.

\begin{tm}\label{lim-monot}\footnote{Эта теорема используется в главе \ref{ch-lopital} при
доказательстве правила Лопиталя для раскрытия неопределенностей типа
$\frac{\infty}{\infty}$} Пусть функция $f$ определена на некотором интервале
$(a,b)$ и $C$ -- произвольное число. Тогда следующие условия эквивалентны:
 \bit{
\item[(i)] $\lim\limits_{x\to b-0} f(x)=C$;

\item[(ii)] для любой монотонной последовательности $x_n$
\begin{equation}
x_1<x_2<...<x_n<...<b, \quad x_n\underset{n\to\infty}{\longrightarrow} b
\label{6.2.1}\end{equation} выполняется
\begin{equation}
f(x_n)\underset{n\to\infty}{\longrightarrow} C \label{6.2.2}\end{equation}
 }\eit
 \end{tm}

\begin{proof} Ясно, что из $(i)$ следует $(ii)$, потому что
если $\lim\limits_{x\to b-0} f(x)=C$, то
$f(x_n)\underset{n\to\infty}{\longrightarrow} C$ для любой последовательности
$x_n\underset{n\to\infty}{\longrightarrow} b, \, x_n<b$ (необязательно
монотонной). Докажем, что наоборот если $(i)$ не выполняется, то не выполняется
и $(ii)$. Пусть $\lim\limits_{x\to b-0} f(x)\ne C$, то есть существует точка
последовательность $x_n\in M$ такая что
$$
x_n\in (a,b)\qquad \&\qquad x_n\underset{n\to\infty}{\longrightarrow} a\qquad
\&\qquad f(x_n)\underset{n\to \infty}{\notarrow} C
$$
План наших действий состоит в том, чтобы, переходя от $x_n$ к ее
подпоследовательностям, построить строго монотонную последовательность с теми
же свойствами.

1. Прежде всего, по свойству подпоследовательностей $2^\circ$ на
с.\pageref{podposledovatelnosti}, второе условие означает, что существует
$\e>0$ и последовательность натуральных чисел
$n_k\underset{k\to\infty}{\longrightarrow}\infty$, такие что
$$
\forall k\in\N\qquad f(x_{n_k})\notin (C-\e;C+\e)
$$
Обозначив $y_k=x_{n_k}$, мы получим последовательность с такими свойствами:
$$
y_k\in (a,b)\qquad \&\qquad y_k\underset{k\to\infty}{\longrightarrow} a\qquad
\&\qquad \forall k\in\N\qquad f(y_k)\notin (C-\e;C+\e)
$$

2. После этого мы индуктивно определим две последовательности $k_m$ и $t_m$.
Сначала полагаем
$$
k_1=1,\qquad t_1=y_1
$$
Затем, если для $m=1,...,l$ числа $k_m$ и $t_m$ определены, мы замечаем вот
что. Поскольку $a<y_k$, $y_k\underset{k\to\infty}{\longrightarrow} a$ и
$a<t_l$, найдется такое $k_{l+1}$, что $y_{k_{l+1}}$ лежит в интервале
$(a;t_l)$:
$$
y_{k_{l+1}}<t_l
$$
Выбираем такое $k_{l+1}$ и полагаем
$$
t_{l+1}=y_{k_{l+1}}
$$
Полученная последовательность $t_m=z_{m+1}$ обладает нужными нам свойствами.
Во-первых, (как и для $x_n$ и $y_k$) все числа $t_m$ лежат в интервале $(a,b)$:
$$
t_m\in (a,b)
$$
Во-вторых, последовательность $t_m$, как подпоследовательность $y_k$, должна
стремиться к $a$:
$$
t_m=k_{k_m}\underset{m\to\infty}{\longrightarrow} a
$$
В-третьих, последовательность значений $f$ на ней не заходит в интервал
$(C-\e;C+\e)$:
$$
\forall m\in\N\qquad  f(t_m)=f(z_{k_m})\notin (C-\e;C+\e)
$$
И, в-четвертых, она строго монотонна:
$$
t_m<y_{k_{m+1}}=t_{m+1}
$$
\end{proof}

Заменой переменной доказывается симметричное утверждение:

\begin{tm}\label{lim-monot>} Пусть функция $f$ определена на некотором интервале
$(a,b)$ и $C$ -- произвольное число. Тогда следующие условия эквивалентны:
 \bit{
\item[(i)] $\lim\limits_{x\to a+0} f(x)=C$;

\item[(ii)] для любой монотонной последовательности $x_n$
 $$
a<...<x_n<...<x_2<x_1<b, \quad x_n\underset{n\to\infty}{\longrightarrow} a
 $$
выполняется
 $$
f(x_n)\underset{n\to\infty}{\longrightarrow} C
 $$
 }\eit
\end{tm}

\subsection{Предел монотонной функции}

 \begin{tm}\label{lim-monot-func}
Пусть функция $f$ определена на некотором интервале $(a,b)$ и монотонна на нем.
Тогда
 \bit{
\item[---] если $f$ неубывает, то
 \begin{align}\label{lim(x-to-a+0)f(x)=inf(x-in-(a;b))f(x)}
&\lim_{x\to a+0} f(x)=\inf_{x\in(a;b)} f(x) && \lim_{x\to b-0}
f(x)=\sup_{x\in(a;b)} f(x)
 \end{align}

\item[---] если $f$ невозрастает, то
 \begin{align}\label{lim(x-to-a+0)f(x)=sup(x-in-(a;b))f(x)}
&\lim_{x\to a+0} f(x)=\sup_{x\in(a;b)} f(x) && \lim_{x\to b-0}
f(x)=\inf_{x\in(a;b)} f(x)
 \end{align}
 }\eit
 \end{tm}

\subsection{Классические теоремы о пределах}

\paragraph{Связь между понятием предела и непрерывностью.}

\begin{tm}[\bf критерий непрерывности]\label{cont-crit}
Пусть функция $f$ определена на некотором интервале $(\alpha;\beta)$,
содержащем точку $x=a$. Тогда следующие условия эквивалентны:
 \bit{
\item[(i)] функция $f$ непрерывна в точке $a$; \item[(ii)] функция $f$ имеет
конечный предел в точке $a$, равный $f(a)$:
$$
\exists \, \lim\limits_{x\to a} f(x)=f(a)
$$
 }\eit
\end{tm}\begin{proof} 1. $(i)\Longrightarrow (ii)$. Пусть
функция $f$ непрерывна в точке $a$, то есть для всякой последовательности
аргументов $\{ x_n \}$, стремящейся к точке $a$ соответствующая
последовательность значений $\{ f(x_n) \}$ стремится к значению $f(a)$. Тогда,
в частности, для всякой последовательности
$$
x_n\underset{n\to \infty}{\longrightarrow} a, \qquad x_n\ne a
$$
выполняется
$$
f(x_n)\underset{n\to \infty}{\longrightarrow} f(a)
$$
По определению предела это означает, что предел функции $f$ в точке $x=a$
равен $f(a)$:
$$
\lim_{x\to a} f(x)=f(a)
$$

2. $(ii)\Longrightarrow (i)$. Пусть наоборот
$$
\lim_{x\to a} f(x)=f(a)
$$
Это означает, что для любой последовательности
$$
x_n\underset{n\to \infty}{\longrightarrow} a, \qquad x_n\ne a
$$
выполняется
$$
f(x_n)\underset{n\to \infty}{\longrightarrow} f(a)
$$
Нам нужно проверить то же самое для последовательностей $\{ x_n \}$, у которых
некоторые элементы $x_n$ совпадают с $a$. Пусть $\{ x_n \}$ -- как раз такая
последовательность, то есть $x_n\underset{n\to \infty}{\longrightarrow} a$,
причем $x_n=a$ для некоторых $n$. Тогда $\{ x_n \}$ можно разложить на две
подпоследовательности $\{ x_{n_k}\}$ и $\{ x_{m_k}\}$ такие что
$$
x_{n_k}\underset{k\to \infty}{\longrightarrow} a \quad x_{n_k}\ne a, \qquad
x_{m_k}=a
$$
Для них мы получаем
$$
f(x_{n_k})\underset{k\to \infty}{\longrightarrow} f(a), \qquad f(x_{m_k})=f(a)
$$
откуда
$$
f(x_n)\underset{n\to \infty}{\longrightarrow} f(a),
$$
Мы показали, что для всякой последовательности аргументов $\{ x_n \}$,
стремящейся к точке $a$ соответствующая последовательность значений $\{ f(x_n)
\}$ стремится к значению $f(a)$. Это означает, что функция $f$ непрерывна в
точке $a$. \end{proof}

\begin{tm}[\bf о перестановочности предела с непрерывной функцией]
\label{lim_f(g)=f(lim g)} Если существует конечный предел
$$
\lim_{x\to a} f(x)=A
$$
и функция $g(y)$ определена в некоторой окрестности точки $A$ и непрерывна в
$A$, то
\begin{equation}\exists \quad \lim_{x\to a} g\Big( f(x) \Big)= g\Big(\lim_{x\to a}
f(x) \Big) \label{5.4.1}\end{equation}\end{tm}\begin{proof} Если
$x_n\underset{n\to \infty}{\longrightarrow} a, \quad x_n\ne a$, то, поскольку
$\lim_{x\to a} f(x)=A$, получаем $f(x_n)\underset{n\to \infty}{\longrightarrow}
A$, и поскольку $g(y)$ непрерывна в точке $A$, $g\Big(
f(x_n)\Big)\underset{n\to \infty}{\longrightarrow} g(A)$. Поскольку это верно
для произвольной последовательности $x_n\underset{n\to \infty}{\longrightarrow}
a, \quad x_n\ne a$, мы получаем $\lim_{x\to a} g\Big( f(x)
\Big)=g(A)=g\Big(\lim_{x\to a} f(x) \Big)$
\end{proof}

\paragraph{Теорема о замене переменной под знаком предела.}

\begin{tm}[\bf о замене переменной под знаком предела]\label{TH:zamena-perem-v-lim}
Пусть
$$
\lim_{x\to a} f(x)=b \qquad \& \qquad \lim_{y\to b} g(y)=c
$$
причем
$$
\forall x\ne a \quad f(x)\ne b
$$
Тогда
$$
\lim_{x\to a} g(f(x))=c
$$
\end{tm}\begin{proof} Возьмем какую-нибудь
последовательность
$$
x_n\underset{n\to \infty}{\longrightarrow} a, \qquad x_n\ne a
$$
Тогда мы получим
$$
f(x_n)\underset{n\to \infty}{\longrightarrow} b, \qquad f(x_n)\ne b
$$
Поэтому
$$
g(f(x_n))\underset{n\to \infty}{\longrightarrow} c
$$
Поскольку это верно для всякой последовательности $x_n\underset{n\to
\infty}{\longrightarrow} a\quad (x_n\ne a)$, мы получаем, что
\newline
$\lim_{x\to a} g(f(x))=c$. \end{proof}

Эта теорема бывает полезна при вычислении предела от сложной функции. Если нам
нужно найти предел
$$
\lim_{x\to a} g(f(x))
$$
и известно, что $f(x)\underset{x\to a}{\longrightarrow} b$ причем $f(x)\ne b$
при $x\ne a$, то можно сделать замену переменных:
$$
\lim_{x\to a} g(f(x))=
\left| \begin{array}{c} f(x)=y \\
y\underset{x\to a}{\longrightarrow} b \\
y\ne b \, \text{при}\, x\ne a \end{array}\right| = \lim_{y\to b} g(y)
$$

\noindent\rule{160mm}{0.1pt}\begin{multicols}{2}

\begin{ex}
Найти предел
$$
\lim_{x\to \infty}\frac{2x^2-x}{x^2+10}
$$
Решение:
 \begin{multline*}
\lim_{x\to \infty}\frac{2x^2-x}{x^2+10}=
\left| \begin{array}{c} x=\frac{1}{y}, \, y=\frac{1}{x}\\
y\underset{x\to \infty}{\longrightarrow} 0
\end{array}\right| =\\=
\lim_{y\to 0}\frac{\frac{2}{y^2}-\frac{1}{y}}{\frac{1}{y^2}+10}= \lim_{y\to
0}\frac{2-y}{1+10 y^2}=\\= \left(\text{подстановка}\right)= \frac{2-0}{1+10
\cdot 0}=2
 \end{multline*}\end{ex}

\begin{ers} Найти пределы
 \biter{

\item[1)] $\lim\limits_{x\to \infty}\frac{2x^2-x}{x^2+10}$,

\item[2)] $\lim\limits_{x\to \infty}\frac{1-3x}{2x+3}$,

\item[3)] $\lim\limits_{x\to \infty}\frac{(x-1)^3}{2x^3+3x+1}$,

\item[4)] $\lim\limits_{x\to \infty}\frac{10x+5}{0.01 x^2-6x}$,

\item[5)] $\lim\limits_{x\to \infty}\frac{x^2-5x+1}{3 x+7}$,

\item[6)] $\lim\limits_{x\to \infty}\frac{x^3-5x}{x+10^{10}}$
 }\eiter
\eers

\end{multicols}\noindent\rule[10pt]{160mm}{0.1pt}

\paragraph{Арифметические операции над пределами.}

\begin{tm}
Если существуют конечные пределы $\lim\limits_{x\to a} f(x)$ и
$\lim\limits_{x\to a} g(x)$, то справедливы формулы:
 \begin{align}
& \lim\limits_{x\to a} (f(x)+g(x))=\lim\limits_{x\to a} f(x) +
\lim\limits_{x\to a} g(x) && \label{5.7.1}
 \\
& \lim\limits_{x\to a} (f(x)-g(x)) =\lim\limits_{x\to a} f(x) -
\lim\limits_{x\to a} g(x) && \label{5.7.2}
 \\
&\lim\limits_{x\to a} (C\cdot f(x))=C\cdot \lim\limits_{x\to a} f(x) &&
(\text{если}\, \, 0\ne C\ne \infty) \label{5.7.3}
 \\
& \lim\limits_{x\to a} f(x)\cdot g(x)=\lim\limits_{x\to a} f(x) \cdot
\lim\limits_{x\to a} g(x) && \label{5.7.4}
 \\
&\lim\limits_{x\to a}\frac{g(x)}{f(x)} = \frac{\lim\limits_{x\to a}
g(x)}{\lim\limits_{x\to a} f(x)} && (\text{если}\, \, \lim\limits_{x\to a}
f(x)\ne 0) \label{5.7.5}
 \\
&\lim\limits_{x\to a} f(x)^{g(x)}= \lim\limits_{x\to a} f(x)^{\lim\limits_{x\to
a} g(x)} && (\text{если}\, \, \lim\limits_{x\to a} f(x)>0) \label{5.7.6}
 \end{align}
 \end{tm}\begin{proof}

Формулы \eqref{5.7.1} -- \eqref{5.7.5} следуют из арифметических свойств
последовательностей (то есть из формул \eqref{5.7.1} -- \eqref{5.7.5} главы
\ref{ch-x_n}). Докажем например \eqref{5.7.4}. Если $\lim\limits_{x\to a}
f(x)=A$ и $\lim\limits_{x\to a} g(x)=B$, то для всякой последовательности
$x_n\underset{n\to \infty}{\longrightarrow} a, \quad x_n\ne a$ возникает
логическая цепочка:
$$
x_n\underset{n\to \infty}{\longrightarrow} a, \quad x_n\ne a
$$
$$
\Downarrow
$$
$$
f(x_n)\underset{n\to \infty}{\longrightarrow} A, \qquad g(x_n)\underset{n\to
\infty}{\longrightarrow} B
$$
$$
\Downarrow\put(20,0){\smsize \text{$\begin{pmatrix}\text{используем свойство}\\
\text{$3^0\, \S 4$ главы \ref{ch-x_n}}\end{pmatrix}$}}
$$
$$
f(x_n)\cdot g(x_n) \underset{n\to \infty}{\longrightarrow} A\cdot B
$$
Поскольку это верно для любой последовательности $x_n\underset{n\to
\infty}{\longrightarrow} a, \quad x_n\ne a$, мы получаем $\lim\limits_{x\to a}
f(x)\cdot g(x)= \lim\limits_{x\to a} f(x)\cdot \lim\limits_{x\to a} g(x)$.

Отдельно нужно доказать формулу \eqref{5.7.6}.

\begin{multline*}\lim\limits_{x\to a} f(x)^{g(x)}= \lim\limits_{x\to a} 2^{\log_2
f(x)^{g(x)}}= \lim\limits_{x\to a} 2^{g(x)\cdot \log_2 f(x)}=
{\smsize{\smsize\begin{pmatrix}\text{применяем формулу}\\
\text{\eqref{5.4.1} этой главы}\end{pmatrix}}} =\\= 2^{\lim\limits_{x\to
a}\Big( g(x)\cdot \log_2 f(x) \Big)}=
{\smsize{\smsize\begin{pmatrix}\text{применяем уже доказанное}\\
\text{свойство $3^0$ этого параграфа}\end{pmatrix}}} = 2^{\lim\limits_{x\to a}
g(x)\cdot \lim\limits_{x\to a}\log_2
f(x)}=\\= {\smsize{\smsize\begin{pmatrix}\text{снова применяем}\\
\text{формулу \eqref{5.4.1} этой главы}\end{pmatrix}}} = 2^{\lim\limits_{x\to
a} g(x)\cdot \log_2 \l \lim\limits_{x\to a} f(x)\r}= 2^{\log_2 \lll \l
\lim\limits_{x\to a} f(x)\r^{\lim\limits_{x\to a} g(x)}\rrr }=\\=
\lim\limits_{x\to a} f(x)^{\lim\limits_{x\to a} g(x)}
\end{multline*}

\end{proof}

\paragraph{Критерий Коши существования предела функции.}

В следующих двух теоремах, называемых критериями Коши, объясняется, что
конечный предел функции $f$ при $x\to a$ существует тогда и только тогда,
когда разность значений этой функции в соседних точках $f(s)-f(t)$ стремится к
нулю при $s,t\to a$.

\begin{tm}[\bf критерий Коши существования двустороннего предела функции]
\label{Cauchy-crit-x->a}\index{критерий!Коши!существования предела функции}
Пусть функция $f$ определена в некоторой проколотой окрестности $U$ величины
$a$ (где $a$ -- произвольное число или символ бесконечности $\infty$). Тогда
следующие условия эквивалентны:
 \bit{
\item[(i)] существует конечный предел
$$
  \lim_{x\to a} f(x)
$$
\item[(ii)] для любых двух последовательностей $s_n, t_n\in U$ стремящихся к
$a$
$$
s_n\underset{n\to \infty}{\longrightarrow} a, \quad t_n\underset{n\to
\infty}{\longrightarrow} a,
$$
разность между соответствующими значениями функции $f$ стремится к нулю:
$$
f(t_n)-f(s_n)\underset{n\to \infty}{\longrightarrow} 0
$$
 }\eit
\end{tm}

\begin{tm}[\bf критерий Коши существования одностороннего предела функции]
\label{Cauchy-crit-x->a-0}\footnote{Этот результат понадобится при
доказательстве критерия Коши сходимости несобственного интеграла (теорема
\ref{tm-17.6.1}).} Пусть функция $f$ определена на интервале $(\alpha;a)$,
где $a$ -- произвольное число или символ бесконечности $\infty$. Тогда
следующие условия эквивалентны:
 \bit{
\item[(i)] существует конечный предел
$$
  \lim_{x\to a-0} f(x)
$$
\item[(ii)] для любых двух последовательностей $s_n, t_n\in (\alpha;a)$
стремящихся к $a$
$$
s_n\underset{n\to \infty}{\longrightarrow} a, \quad t_n\underset{n\to
\infty}{\longrightarrow} a,
$$
разность между соответствующими значениями функции $f$ стремится к нулю:
$$
f(t_n)-f(s_n)\underset{n\to \infty}{\longrightarrow} 0
$$
 }\eit
\end{tm}\begin{proof} Эти две теоремы доказываются
одинаково, поэтому мы докажем лишь вторую.

1. Легко проверяется, что из условия $(i)$ следует условие $(ii)$.
Действительно, если существует конечный предел $\lim\limits_{x\to a-0} f(x)=C$,
то есть для всякой последовательности $x_n\in (\alpha;a), \,\,
x_n\underset{n\to \infty}{\longrightarrow} a$ выполняется соотношение
$$
f(x_n)\underset{n\to \infty}{\longrightarrow} C
$$
то взяв произвольные последовательности $s_n, t_n\in (\alpha;a), \,\,
s_n\underset{n\to \infty}{\longrightarrow} a, \quad t_n\underset{n\to
\infty}{\longrightarrow} a$, мы получим
$$
f(t_n)-f(s_n)\underset{n\to \infty}{\longrightarrow} C-C=0
$$

2. Теперь докажем, что из условия $(ii)$ следует условие $(i)$. Пусть
выполняется $(ii)$.

a) Покажем сначала, что тогда для любой последовательности $x_n\in (\alpha;a),
\,\, x_n\underset{n\to \infty}{\longrightarrow} a$ существует конечный предел
$$
\lim_{n\to \infty} f(x_n)
$$
Действительно, обозначим $y_n=f(x_n)$. Тогда если взять любые две
последовательности индексов
$$
p_i\underset{i\to \infty}{\longrightarrow}\infty, \quad q_i\underset{i\to
\infty}{\longrightarrow}\infty
$$
то положив $t_i=x_{p_i}, \, s_i=x_{q_i}$, мы получим, что, в силу условия
$(ii)$
$$
y_{p_i}-y_{q_i}=f(x_{p_i})-f(x_{q_i})= f(t_i)-f(s_i)\underset{i\to
\infty}{\longrightarrow} 0
$$
Это означает, что последовательность $y_n=f(x_n)$ удовлетворяет критерию Коши
(то есть обладает свойством $(i)$ теоремы \ref{Cauchy-crit-seq}). Значит,
$y_n=f(x_n)$ имеет конечный предел.

b) Мы показали, что для любой последовательности $x_n\in (\alpha;a), \,\,
x_n\underset{n\to \infty}{\longrightarrow} a$ существует конечный предел
$\lim\limits_{n\to \infty} f(x_n)$. Теперь из условия $(ii)$ следует, что все
такие пределы совпадают, потому что для любых двух последовательностей $s_n,
t_n\in (\alpha;a)$ стремящихся к $a$
$$
s_n\underset{n\to \infty}{\longrightarrow} a, \quad t_n\underset{n\to
\infty}{\longrightarrow} a
$$
выполняется равенство
$$
\lim_{n\to \infty}\l f(t_n)-f(s_n) \r=0
$$
и следовательно,
$$
\lim_{n\to \infty} f(t_n)=\lim_{n\to \infty} f(s_n)
$$
Это в свою очередь означает, что существует число $C$ такое, что
$$
\lim_{n\to \infty} f(x_n)=C
$$
для любой последовательности $x_n\in (\alpha;a), \,\, x_n\underset{n\to
\infty}{\longrightarrow} a$, то есть
$$
\lim_{x\to a-0} f(x)=C
$$
Мы убедились, что из условия $(ii)$ следует условие $(i)$. Теорема 1.1
доказана. \end{proof}

\section{Язык $\varepsilon$-$\delta$ Коши}\label{e-d-cauchy}

\subsection{Определение предела функции по Коши}

Данное нами на с.\pageref{Opred-predela-funktsii} определение предела функции
по Гейне -- не единственно возможное. В некоторых случаях бывает полезно знать
еще одно, эквивалентное определение, предложенное известным французским
математиком Огюстеном Коши (1789-1857). Чтобы его сформулировать, нам, как и в
случае с определением Гейне, понадобится вспомнить несколько терминов.

Пусть, как и в \ref{sec-pred-funct}\ref{subsec-opred-pred-funct}, $a$ и $A$
обозначают числа, или символы бесконечности, возможно с определенным знаком
$$
a, A \in \R \quad \text{или}\quad a,A=\infty, \, -\infty, \, +\infty
$$
и пусть для всякого такого $a$
 \bit{
\item[---] {\it проколотой окрестностью} величины $a$ называется всякое
множество вида
$$
U=\left\{ \begin{array}{c} (a-\delta,a)\cup (a,a+\delta), \quad \text{где}\,
\delta>0 \quad
\Big(\text{если $a$ -- число}\Big) \\
(-\infty,-E)\cup (E,+\infty), \quad \text{где}\, E>0 \quad
\Big(\text{если}\, a=\infty \Big);\\
(-\infty,-E), \quad \text{где}\, E>0 \quad
\Big(\text{если}\, a=-\infty \Big); \\
(E,+\infty), \quad \text{где}\, E>0 \quad \Big(\text{если}\, a=+\infty \Big);
\end{array}\right\}
$$
\item[---] {\it левой полуокрестностью} $a$ называется всякий интервал вида
$$
U=(a-\delta,a), \quad \text{где}\, \delta>0 \quad \Big(\text{если $a$ --
число}\Big)
$$
\item[---] {\it правой полуокрестностью} $a$ называется всякий интервал вида
$$
U=(a,a+\delta), \quad \text{где}\, \delta>0 \quad \Big(\text{если $a$ --
число}\Big)
$$
\item[---] {\it окрестностью} величины $a$ называется
 \bit{
\item[---] всякий интервал вида $U=(a-\varepsilon,a+\varepsilon)$ (где
$\varepsilon>0$), если $a$ -- точка; \item[---] всякая проколотая окрестность
величины $a$, если $a$ -- символ бесконечности.
 }\eit
 }\eit

Теперь даем само определение по Коши:

\bit{\item[$\bullet$] Величина $A$ называется {\it пределом
функции}\index{предел!функции!по Коши} $f(x)$ при $x$ стремящимся к
$a$\label{DEF:Cauchy-lim(x->a)}
$$
A=\lim_{x\to a} f(x) \quad \l f(x)\underset{x\to a}{\longrightarrow} A \r
$$
если
 \bit{
\item[1)] $f(x)$ определена в некоторой выколотой окрестности $U$ величины $a$
и \item[2)] для любой окрестности $V$ величины $A$ найдется выколотая
окрестность $W$ величины $a$ такая, что
$$
  \forall x\in W\quad f(x)\in V
$$
 }\eit
}\eit

Отдельно определяются односторонние пределы:

\bit{\item[$\bullet$] Величина $A$ называется {\it пределом слева (справа)
функции} $f(x)$ в точке $a$\label{DEF:Cauchy-lim(x->a-0)}
$$
  A=\lim_{x\to a-0} f(x) \quad \l A=\lim_{x\to a+0} f(x) \r
$$
если
 \bit{
\item[1)] $f(x)$ определена в некоторой левой (правой) полуокрестности $U$
точки $a$ и \item[2)] для любой окрестности $V$ величины $A$ найдется левая
(правая) полуокрестность $W$ точки $a$ такая, что
$$
  \forall x\in W\quad f(x)\in V
$$
 }\eit
}\eit

\noindent\rule{160mm}{0.1pt}\begin{multicols}{2}

Покажем, как эта общая схема работает в конкретных ситуациях.

\biter{\item[$\bullet$]  {\bf Конечный предел функ\-ции в точке по Коши.} Число
$A$ называется {\it пределом функции $f$ в точке $a$}
$$
\lim_{x\to a} f(x)=A
$$
если
 \biter{
\item[1)] функция $f$ определена на некотором множестве вида $(a-\eta,a)\cup
(a,a+\eta)$ (где $\eta>0$), и \item[2)] для всякого числа $\varepsilon>0$
существует число $\delta>0$ такое что для любого $x\in (a-\delta;a)\cup
(a,a+\delta)$ выполняется включение $f(x)\in (A-\varepsilon; A+\varepsilon)$.
 }\eiter

\item[$\bullet$] {\bf Бесконечный предел функции в точке по Коши.} Говорят, что
функция $f$ {\it имеет бесконечный предел в точке $a$}
$$
\lim_{x\to a} f(x)=\infty
$$
если
 \biter{
\item[1)] функция $f$ определена на некотором множестве вида $(a-\eta,a)\cup
(a,a+\eta)$ (где $\eta>0$), и \item[2)] для всякого числа $E>0$ существует
число $\delta>0$ такое что для любого $x\in (a-\delta;a)\cup (a,a+\delta)$
выполняется включение $f(x)\in (-\infty,-E)\cup (E,+\infty)$.
 }\eiter

\item[$\bullet$] {\bf Конечный предел функции слева в точке по Коши.} Число $A$
называется {\it пределом слева функции $f$ в точке $a$}
$$
\lim_{x\to a-0} f(x)=A
$$
если
 \biter{
\item[1)] функция $f$ определена на некотором интервале вида $(a-\eta; a)$ (где
$\eta>0$, и \item[2)] для всякого числа $\varepsilon>0$ существует число
$\delta>0$ такое что для любого $x\in (a-\delta;a)$ выполняется включение
$f(x)\in (A-\varepsilon; A+\varepsilon)$.
 }\eiter

\item[$\bullet$] {\bf Конечный предел функции при $x\to +\infty$ по Коши.}
Число $A$ называется {\it пределом функции $f$ при $x\to +\infty$}
$$
\lim_{x\to +\infty} f(x)=A
$$
если
 \biter{
\item[1)] функция $f$ определена на некотором интервале вида $(\alpha;
+\infty)$, и \item[2)] для всякого числа $\varepsilon>0$ существует число
$\beta>\alpha$ такое что для любого $x\in (\beta;+\infty)$ выполняется
включение $f(x)\in (A-\varepsilon; A+\varepsilon)$.
 }\eiter

}\eiter

\end{multicols}\noindent\rule[10pt]{160mm}{0.1pt}

\begin{tm}[\bf об эквивалентности определений предела по Коши и по Гейне]
\label{Heine-Cauchy}\footnote{ Эквивалентность определений предела по Коши и по
Гейне используется ниже в главе 17 при доказательстве признаков сравнения
несобственных интегралов (теоремы 5.7 и 7.5).} Для произвольной функции $f$
и любых величин $A$ и $a$ следующие условия эквивалентны:
 \bit{
\item[(i)] $A$ является пределом функции $f$ при $x\to a$ по Коши;
\item[(ii)] $A$ является пределом функции $f$ при $x\to a$ по Гейне.
 }\eit
(Аналогичное утверждение справедливо для односторонних пределов.)
\end{tm}\begin{proof} Мы докажем эту теорему для конечных
пределов слева, поскольку в каждом случае доказательство отличается лишь
несущественными деталями.

1. $(i)\Rightarrow (ii)$. Пусть $A$ является пределом по Коши функции $f$
при $x\to a-0$, то есть для всякого $\varepsilon>0$ существует $\delta>0$ такое
что для любого $x\in [a-\delta;a)$ выполняется включение $f(x)\in
(A-\varepsilon; A+\varepsilon)$. Покажем, что тогда $A$ является пределом по
Гейне $f(x)$ при $x\to a-0$, то есть что для любой последовательности
аргументов $x_n\underset{n\to \infty}{\longrightarrow} a-0$ выполняется
$f(x_n)\underset{n\to \infty}{\longrightarrow} A$.

Возьмем произвольную такую последовательноть $x_n\underset{n\to
\infty}{\longrightarrow} a-0$. и зафиксируем какое-нибудь $\varepsilon>0$. По
предположению, для него можно выбрать $\delta>0$ так чтобы
\begin{equation}\forall x\in [a-\delta;a) \quad f(x)\in (A-\varepsilon;
A+\varepsilon) \label{6.1.1}\end{equation} Зафиксируем это число $\delta$.
Получаем следующую логическую цепочку:
$$
x_n\underset{n\to \infty}{\longrightarrow} a-0
$$
$$
\Downarrow
$$
$$
x_n\in [a-\delta;a) \quad \text{для почти всех}\,\, n\in \mathbb{N}
$$
$$
\Downarrow
$$
$$
(\text{применяем \eqref{6.1.1}})
$$
$$
\Downarrow
$$
$$
f(x_n)\in (A-\varepsilon;A+\varepsilon) \quad \text{для почти всех}\,\, n\in
\mathbb{N}
$$
Это верно для произвольного $\varepsilon>0$. Значит, $f(x_n)\underset{n\to
\infty}{\longrightarrow} A$, а именно это нам и нужно было доказать.

2. $(i)\kern-9pt\Big\backslash \Rightarrow (ii)\kern-11pt\Big\backslash$. Пусть
наоборот, $A$ не является пределом по Коши функции $f$ при $x\to a-0$, то
есть существует такое $\varepsilon>0$, что для любого $\delta>0$ найдется такой
$x\in [a-\delta;a)$, для которого выполняется $f(x)\notin (A-\varepsilon;
A+\varepsilon)$. Тогда, в частности, для любого числа вида $\delta=\frac{1}{n}$
мы получим, что
\begin{equation}\exists x_n\in [a-\delta;a)=[a-\frac{1}{n};a) \quad \text{такие
что}\quad f(x_n)\notin (A-\varepsilon; A+\varepsilon)
\label{6.1.2}\end{equation}

Зафиксируем такие числа $x_n$. Получается:
$$
a-\frac{1}{n}<x_n<a
$$
$$
\Downarrow
$$
$$
x_n\underset{n\to \infty}{\longrightarrow} a-0
$$
Но, с другой стороны, в силу \eqref{6.1.2},
$$
f(x_n)\underset{n\to \infty}{\notarrow} A
$$
Значит, $A$ не является пределом по Гейне функции $f$ при $x\to a-0$.

Мы получили, что если не выполняется $(i)$ то автоматически не выполняется
$(ii)$. \end{proof}

 \subsection{Равномерная непрерывность функции по Коши}
 \label{SEC-opredelenie-Cauchy-ravnomernoj-nepreryvnosti}

В \ref{SEC-th-o-nepr-func}\ref{SEC-teorema-Kantora} главы \ref{ch-cont-f(x)} мы
дали определение равномерно непрерывной функции на множестве. Это определение
легко переводится на язык $\varepsilon-\delta$ Коши. В дальнейшем при
доказательстве результатов такая переформулировка нам не понадобится, но мы
полагаем, что она сама по себе будет интересна читателю.

\bit{\item[$\bullet$] Функция $f(x)$ называется {\it равномерно непрерывной} на
множестве $E$, если $f(x)$ определена на $E$, и для всякого числа
$\varepsilon>0$ найдется такое число $\delta>0$, что для любых $x$ и $y$ из
$E$, расстояние между которыми меньше $\delta$, соответствующее расстояние
между значениями функции $f$ и $f(y)$ будет меньше $\varepsilon$:
 \begin{equation}\forall
\varepsilon>0 \quad \exists \delta>0 \quad \forall x,y\in E \qquad |x-y|<\delta
\quad \Longrightarrow \quad |f(x)-f(y)|<\varepsilon
 \label{ravnomernaya nepreryvnost po Cauchy}\end{equation}}\eit

\begin{proof}[Докажем равносильность этих определений.]

1. Пусть выполняется \eqref{ravnomernaya nepreryvnost po Cauchy}. Возьмем
произвольные последовательности $x_n,y_n\in E$ такие, что
$$
x_n-y_n\underset{n\to\infty}{\longrightarrow} 0
$$
Покажем, что тогда
$$
f(x_n)-f(y_n)\underset{n\to\infty}{\longrightarrow} 0
$$
Для этого возьмем $\varepsilon>0$. В силу \eqref{ravnomernaya nepreryvnost po
Cauchy}, должен существовать такой $\delta>0$, что
$$
\forall x,y\in E \qquad |x-y|<\delta \quad \Longrightarrow \quad
|f(x)-f(y)|<\varepsilon
$$
В частности,
$$
\underbrace{|x_n-y_n|<\delta}_{\text{выполняется для почти всех $n$}}\quad
\Longrightarrow \quad \underbrace{|f(x_n)-f(y_n)|<\varepsilon}_{\text{значит, и
это выполняется для почти всех $n$}}
$$
Последнее верно при любом $\varepsilon>0$, поэтому
$f(x_n)-f(y_n)\underset{n\to\infty}{\longrightarrow} 0$.

2. Наоборот, предположим, что \eqref{ravnomernaya nepreryvnost po Cauchy} не
выполняется, то есть существует некоторое $\varepsilon>0$ такое что
$$
\forall \delta>0 \quad \exists x,y\in E \qquad |x-y|<\delta \quad \& \quad
|f(x)-f(y)|\ge \varepsilon
$$
В частности, это должно выполняться для чисел $\delta=\frac{1}{n}$:
$$
\forall n\in\N \quad \exists x_n,y_n\in E \qquad |x_n-y_n|<\frac{1}{n}\quad \&
\quad |f(x_n)-f(y_n)|\ge \varepsilon
$$
Мы получили последовательности $x_n$ и $y_n$ такие, что $0\le
|x_n-y_n|<\frac{1}{n}\underset{n\to\infty}{\longrightarrow} 0$, и поэтому
$$
x_n-y_n\underset{n\to\infty}{\longrightarrow} 0
$$
но при этом $|f(x_n)-f(y_n)|\ge \varepsilon$, и поэтому
$$
f(x_n)-f(y_n)\underset{n\to\infty}{\notarrow} 0
$$
Значит, \eqref{ravnomernaya nepreryvnost po Heine} не может выполняться.
\end{proof}

\btm[критерий равномерной непрерывности]\label{TH-krit-ravn-nepr} Функция
$f:E\to\R$ тогда и только тогда равномерно непрерывная на $E$, когда
выполняется соотношение
 \beq
 \sup_{x,y\in E: \; |x-y|\le \delta} |f(x)-f(y)|
 \underset{\delta\to  0}{\longrightarrow} 0
 \label{kriterij-ravnomernoj-nepreryvnosti}
 \eeq
\etm

\chapter{СТАНДАРТНЫЕ ФУНКЦИИ}
\label{ch-ELEM-FUNCTIONS}

Среди функций, рассматриваемых в математическом анализе, имеется несколько,
используемых особенно часто. Эти функции называются {\it элементарными}, и их
свойства, а также свойства функций, получаемых из них с помощью алгебраических
операций и композиции (такие функции мы будем называть {\it стандартными})
служат предметом изучения в разделе анализа, называемом {\it Исчислением}. Что
это такое мы начнем обсуждать в следующей главе
(с.\pageref{SEC-proizv-kak-oper-nad-simv}), а здесь мы определим элементарные
функции и обсудим их свойства, которые понадобятся нам в дальнейшем. Список
элементарных функций выглядит так:\footnote{Смысл встречающихся в этой таблице
обозначений типа $\frac{\Z}{2\N-1}$, $\frac{\N}{2\N-1}$ объяснялся в конце
главы \ref{ch-R&N} на с.\pageref{Z/2N-1}.}

{\smsize

\begin{center}\label{spisok-elem-functsij}
\begin{tabular}{|l|l|l|} \hline
  &  &  \\
 название & $\begin{array}{l} \text{обозначение:} \\ x \\ \text{\rotatebox{-90}{$\mapsto$}} \end{array}$ & область определения \\
  &  &  \\ \hline
  &  &  \\
 степенная функция & $x^b$ &
$\begin{cases}x\in\R, & b\in\frac{\Z_+}{2\N-1} \\ x\ne 0, &
b\in-\frac{\N}{2\N-1} \\ x\ge 0, & b\in(0,+\infty)\setminus\frac{\Z}{2\N-1} \\
x>0, & b\in(-\infty,0)\setminus\frac{\Z}{2\N-1}
\end{cases}$
 \\
   &  &  \\
 \hline
  &  &  \\
 показательная функция & $a^x$ &
$\begin{cases} x\in\R, & a>0 \\
x\ge 0, & a=0 \\
x\in\frac{\Z}{2\N-1}, & a<0
\end{cases}$
 \\
   &  &  \\
 \hline
  &  &  \\
 логарифм & $\log_a x$ & $x>0$, $\ a\in(0;1)\cup(1;+\infty)$ \\
   &  &  \\
 \hline
   &  &  \\
синус & $\sin x$ & $x\in\R$ \\
  &  &  \\
\hline
  &  &  \\
 косинус & $\cos x$ & $x\in\R$ \\
   &  &  \\
\hline
  &  &  \\
 тангенс & $\tg x$ & $x\notin \frac{\pi}{2}+\pi\Z$  \\
   &  &  \\
\hline
  &  &  \\
 котангенс & $\ctg x$ & $x\notin \pi\Z$ \\
   &  &  \\
\hline
  &  &  \\
 арксинус & $\arcsin x$ & $x\in[-1;1]$ \\
   &  &  \\
\hline
  &  &  \\
 арккосинус & $\arccos x$ & $x\in[-1;1]$ \\
   &  &  \\
\hline
  &  &  \\
 арктангенс & $\arctg x$ & $x\in\R$ \\
   &  &  \\
\hline
  &  &  \\
 арккотангенс & $\arcctg x$ & $x\in\R$ \\
   &  &  \\
\hline
\end{tabular}
\end{center}

} \noindent Как может заметить читатель, только одна из этих функций -- $x^b$
-- была определена нами к настоящему моменту, и то только для целых значений
параметра $b$ (см. определения на с.\pageref{def:a^n} и \pageref{DF:a^n-n<0}).
Остальные функции в этом списке мы определим необычным приемом -- с помощью
двух так называемых {\it избыточных аксиом}\index{аксиомы!избыточные}.

Это вот что такое. {\it Избыточной аксиомой} называется такая аксиома теории,
про которую впоследствии становится известно, что ее можно вывести из остальных
аксиом этой теории. Как следствие, избыточные аксиомы можно (и с точки зрения
логической стройности правильнее) оформлять не как аксиомы, а как теоремы этой
теории. Однако в математике встречаются время от времени ситуации, когда
знакомство с теоремой бывает желательно задолго до того, как появляется
возможность эту теорему строго доказать. В таких случаях полезно
объявить эту теорему аксиомой, сообщив слушателям, что в действительности эта
аксиома избыточна.

Элементарные функции -- как раз такой пример. Чисто логически ничто не мешает
молчать о них в курсе анализа до тех пор, пока не появится возможность
аккуратно их определить и доказать их характеристические свойства. Однако тогда
возможность решать задачи с элементарными функциями, отличными от рациональных,
появится у преподавателя и у студентов только в конце второго семестра. Поскольку
выбрасывание из упражнений всех иррациональных элементарных функций
эквивалентно просто отказу от Исчисления в классическом его понимании, разумно
поступить по-другому. В этой главе мы сформулируем две избыточные аксиомы --
Аксиому степеней (аксиома B1 на с.\pageref{TH-o-step-otobr}) и Аксиому
тригонометрии (аксиома B2 на с.\pageref{osn-teor-trigonometrii}).
Доказательство их избыточности мы сможем привести очень нескоро -- лишь в
\ref{SEC:prilozh-step-ryadov} главы \ref{CH-step-ryady}. Поэтому, в
соответствии с общим принципом, мы переносим эти утверждения вместе с их
следствиями, в иллюстративный материал.

\noindent\rule{160mm}{0.1pt}\begin{multicols}{2}

\section{Элементарные функции}\label{SUBSEC-elem-funktsii}

\subsection{Степени с нецелым показателем и логарифмы}

\paragraph{Степенное отображение.}

На страницах \pageref{def:a^n} и \pageref{DF:a^n-n<0} мы уже определили степень
числа с целым показателем:
$$
a^n \qquad (n\in\Z)
$$
Для нецелых показателей это понятие вводится с помощью следующей избыточной
аксиомы:

 \biter{\it
\item[\rm B1.] {\bf Аксиома
степеней.}\label{TH-o-step-otobr}\footnote{Избыточность Аксиомы степеней
доказывается в \ref{SEC:prilozh-step-ryadov}\ref{PR:teor-o-step-otobr} главы
\ref{CH-step-ryady}.} Существует единственное отображение
$$
(a,b)\mapsto a^b,
$$
со следующими свойствами:
 \biter{
\item[$P_0$:] оно определено в следующих трех ситуациях:
\label{usloviya-sushestv-x^a}
 \biter{
\item[---] при $a>0$ и произвольном $b\in\R$,

\item[---] при $a=0$ и $b\ge 0$,

\item[---] при $a<0$ и $b\in\frac{\Z}{2\N-1}$,
 }\eiter\noindent
или, что то же самое, в следующих четырех:
 \biter{
\item[---] при $b\in\frac{\Z_+}{2\N-1}$ и тогда $a\in\R$ может быть
произвольным,

\item[---] при $b\in-\frac{\N}{2\N-1}$ и тогда $a\ne 0$,

\item[---] при $b\in(0,+\infty)\setminus\frac{\Z}{2\N-1}$, и тогда $a\ge 0$,

\item[---] при $b\in(-\infty,0)\setminus\frac{\Z}{2\N-1}$, и тогда $a>0$,
 }\eiter\noindent

\item[$P_1$:] следующие две группы тождеств выполняются всякий раз, когда  обе
части тождества определены:
 \biter{
\item[---] {\bf показательные законы}\index{закон!показательный}:
 \begin{align}
&a^0=1,\label{a^0} \\
&a^{-x}=\frac{1}{a^x}, \label{a^(-x)} \\
&a^{x+y}=a^x\cdot a^y;\label{a^(x+y)}
 \end{align}

\item[---] {\bf степенные законы}\index{закон!степенной}:
 \begin{align}
&1^b=1,\label{1^b} \\
&\left(\frac{1}{x}\right)^b=\frac{1}{x^b},\label{(1/x)^b} \\
&(x\cdot y)^b=x^b\cdot y^b \label{(xy)^b}
 \end{align}
 }\eiter\noindent

\item[$P_2$:] тождество, называемое {\bf накопительным
законом}\index{закон!накопительный},
 \begin{align}\label{nakop-zakon}
&(a^x)^y=a^{x\cdot y}
 \end{align}
выполняется для следующих значений переменных:
 \biter{
\item[---] при $a>0$ и $x,y\in\R$,

\item[---] при $a=0$ и $x\ge 0$, $y\ge 0$,

\item[---] при $a<0$ и $x,y\in\frac{\Z}{2\N-1}$,
 }\eiter

\item[$P_3$:] для нулевого основания степень, в случаях, когда она определена,
описывается формулой
 \beq\label{0^b}
0^b=\begin{cases}1,& b=0 \\ 0,& b>0\end{cases}
 \eeq
а для положительных выполняется следующее {\bf условие сохранения знака}:
 \beq\label{x^b>0}
\kern-10pt a>0\quad\Rightarrow\quad a^b>0\quad (b\in\R)
 \eeq

\item[$P_4$:] для положительных оснований выполняются следующие {\bf условия
монотонности}:
 \biter{
\item[---] если $b>0$, то возведение в степень $b$ не меняет знак неравенства:
 \beq\label{monot-x^b-b>0}
 0<x<y \quad\Longrightarrow\quad x^b<y^b
 \eeq
\item[---] если $b<0$, то возведение в степень $b$ меняет знак неравенства:
 \beq\label{monot-x^b-b<0}
 0<x<y \quad\Longrightarrow\quad x^b>y^b
 \eeq

\item[---] если $a>1$, то потенцирование с основанием $a$ не меняет знак
неравенства:
 \beq\label{monot-a^x-a>1}
x<y \quad\Longrightarrow\quad a^x<a^y
 \eeq

\item[---] если $0<a<1$, то потенцирование с основанием $a$ меняет знак
неравенства:
 \beq\label{monot-a^x-0<a<1}
x<y \quad\Longrightarrow\quad a^x>a^y
 \eeq
 }\eiter
 }\eiter

 }\eiter

\brem Для натуральных значений $b$ тождество \eqref{a^(x+y)} влечет за собой
индуктивное тождество \eqref{def:a^n}:
$$
a^0=1,\quad a^{n+1}=a^n\cdot a,
$$
поэтому $a^n$ в этом случае совпадает со степенью числа $a$, определенной на
странице \pageref{opr-stepeni}. Как следствие, выполняется тождество
 \beq\label{x^1=x}
x^1=x
 \eeq
 \erem

\brem Накопительный закон \eqref{nakop-zakon} не обязан выполняться в случаях,
не предусмотренных условиями Аксиомы степеней B1. Например, равенство
$$
\underbrace{\left((-1)^2\right)^\frac{1}{2}}_{\scriptsize\begin{matrix} \| \\
1^\frac{1}{2}\\ \| \\ 1
\end{matrix}}=\underbrace{(-1)^{2\cdot\frac{1}{2}}}_{\scriptsize\begin{matrix} \| \\
(-1)^1\\ \| \\ -1
\end{matrix}}
$$
конечно, не будет верным.
 \erem

\brem\label{REM:0^b} Формула \eqref{0^b} вытекает из остальных утверждений
теоремы: если $b>0$, то
$$
0\cdot 0=0
$$
$$
\Downarrow
$$
$$
0^b\cdot 0^b=\eqref{(xy)^b}=(0\cdot 0)^b=0^b
$$
$$
\Downarrow
$$
$$
0^b\cdot 0^b-0^b=0
$$
$$
\Downarrow
$$
$$
\left[\begin{matrix}0^b=0 \\ 0^b=1\end{matrix}\right]
$$
Здесь второе равенство -- $0^b=1$ -- невозможно, потому, что, если бы оно было
верно, мы получили бы противоречие:
 \begin{multline*}
0=0^1=0^{b\cdot\frac{1}{b}}=\eqref{nakop-zakon}=\\=
(\kern-22pt\underbrace{0^b}_{\scriptsize\begin{matrix}\| \\ \phantom{,}1,\\
\text{по предположению}\end{matrix}}\kern-22pt)^{\frac{1}{b}}
=1^{\frac{1}{b}}=\eqref{1^b}=1
 \end{multline*}
Значит, должно быть верно первое равенство: $0^b=0$.  \erem

\ber Можно заметить, что условие сохранения знака \eqref{x^b>0} также следует
из остальных утверждений Аксиомы степеней B1. Нам этот факт в дальнейшем не
понадобится, однако читатель в качестве упражнения может доказать его
самостоятельно. В качестве подсказки мы сообщим только, что здесь удобно
использовать теорему о степенном отображении \ref{TH-o-step-otobr-malaya} и
предложение \ref{LM:x^(1/2n-1)-nepr} на с.\pageref{LM:x^(1/2n-1)-nepr}.
 \eer

\bcor\label{COR:a-ge-0=>a^b-ge-0} Справедлива импликация:
 \beq\label{a-ge-0=>a^b-ge-0}
\Big(a\ge 0\quad \&\quad b>0\Big)\quad\Longrightarrow\quad a^b\ge 0
 \eeq
\ecor
 \bpr
Если $a=0$, то
$$
a^b=0^b\overset{\eqref{0^b}}{=}0
$$
Если же $a>0$, то
$$
a^b\overset{\eqref{x^b>0}}{>}0
$$
В обоих случаях получается $a^b\ge 0$.
 \epr

\paragraph{Корни.}\label{PAR:korni}

 \biter{

\item[$\bullet$] {\it Корнем степени $n\in\N$} называется функция
 \beq
\sqrt[n]{x}:=x^{\frac{1}{n}}
 \eeq
Из аксиомы степеней B1 следует, что эта функция определена
 \biter{
\item[--] на всей прямой $\R$ для $n\in 2\N-1$ (то есть когда $n$ нечетно),

\item[--] на полуинтервале $[0;+\infty)$ для $n\in 2\N$ (то есть когда $n$
четно).
 }\eiter

\item[$\bullet$] Корень степени $2$ называется {\it квадратным корнем}, и его
обозначение упрощается до $\sqrt{\phantom{x}}$:
 \beq
\sqrt{x}:=\sqrt[2]{x}=x^{\frac{1}{2}}
 \eeq

\item[$\bullet$] Корень степени $3$ называется {\it кубическим корнем}:
 $$
\sqrt[3]{x}=x^{\frac{1}{3}}
 $$
 }\eiter

\bprop\label{LM:x^(1/2n-1)-nepr} Корень нечетной степени $x\mapsto
\sqrt[2n-1]{x}$
 \biter{
\item[---] определен всюду на $\R$,

\item[---] возрастает на $\R$,

\item[---] непрерывен на $\R$,

\item[---] равен нулю в нуле:
 \beq\label{sqrt[2n-1](0)=0}
\sqrt[2n-1]{0}=0
 \eeq

\item[---] имеет бесконечные пределы на бесконечностях:
 \begin{align}
&\lim_{x\to-\infty} \sqrt[2n-1]{x}=-\infty \label{x^(1/(2n-1))->8(x->-8)} \\
&\lim_{x\to+\infty} \sqrt[2n-1]{x}=+\infty \label{x^(1/(2n-1))->8(x->+8)}
 \end{align}

\item[---] удовлетворяет следующему правилу, которое может считаться его
определением:
 \beq\label{DF:sqrt[2n-1]}
y=\sqrt[2n-1]{x}\quad\Longleftrightarrow\quad y^{2n-1}=x
 \eeq
 }\eiter
\eprop
 \bpr
1. Как мы уже отмечали, тот факт, что функция $x\mapsto
\sqrt[2n-1]{x}=x^{\frac{1}{2n-1}}$ определена на $\R$ постулируется в Аксиоме
степеней B1 на с.\pageref{TH-o-step-otobr}, поэтому его доказывать не нужно.

2. Докажем правило \eqref{DF:sqrt[2n-1]}. Оно означает, что функции
$$
f(x)=x^\frac{1}{2n-1},\quad g(y)=y^{2n-1}
$$
(определенные на всей прямой $\R$, по условиям Аксиомы степеней B1 на с.
\pageref{TH-o-step-otobr}) должны быть обратны друг другу. Это доказывается с
помощью накопительного закона \eqref{nakop-zakon}: заметим, что цепочка
равенств
 \beq\label{(a^(1/(2n-1)))^(2n-1)=a}
(a^\frac{1}{2n-1})^{2n-1}=a^{\frac{1}{2n-1}\cdot
(2n-1)}=a^1\overset{\eqref{x^1=x}}{=}a
 \eeq
выполняется при любых значениях $a\in\R$:
 \biter{
\item[---] для $a>0$ первое равенство в \eqref{(a^(1/(2n-1)))^(2n-1)=a}
выполняется, независимо от свойств степеней $\frac{1}{2n-1}$ и $2n-1$,

\item[---] для $a=0$ оно выполняется, потому что $\frac{1}{2n-1}\ge 0$ и
$2n-1\ge 0$,

\item[---] для $a<0$ оно выполняется, потому что
$\frac{1}{2n-1}\in\frac{\Z}{2\N-1}$ и $2n-1\in\frac{\Z}{2\N-1}$.
 }\eiter
Из \eqref{(a^(1/(2n-1)))^(2n-1)=a} следует тождество:
 $$
g(f(x))=(x^\frac{1}{2n-1})^{2n-1}=x,\qquad x\in\R
 $$
По тем же причинам цепочка равенств
 $$
(a^{2n-1})^\frac{1}{2n-1}=a^{(2n-1)\cdot\frac{1}{2n-1}}=a^1=a
 $$
выполняется при любых значениях $a\in\R$, и из нее следует тождество:
 $$
f(g(y))=(y^{2n-1})^\frac{1}{2n-1}=y,\qquad y\in\R
 $$
Вместе это и означает, что функции $f$ и $g$ обратны друг другу.

3. В примере \ref{EX:x^(2n-1)-vozrastaet} мы показали, что степенная функция
$$
g(y)=y^{2n-1}
$$
возрастает на всей прямой $\R$ (а для интервала $(0;+\infty)$ это утверждается
в первом условии монотонности \eqref{monot-x^b-b>0}). По теореме о монотонной
функции \ref{Teor-ob-obratnoi-funktsii-na-int-v-R^1} отсюда следует, что ее
обратная функция
$$
f(x)=x^\frac{1}{2n-1}
$$
должна возрастать и быть непрерывной на своей области определения, то есть на
множестве
 \begin{multline*}
\Big(\inf_{y\in\R}g(y);\sup_{y\in\R}g(y)\Big)=\\=
\Big(\lim_{y\to-\infty}g(y);\lim_{y\to+\infty}g(y)\Big)=\\=
\eqref{x^(2n-1)->-8(x->-8)},\eqref{x^(2n-1)->+8(x->+8)}= (-\infty;+\infty)=\R
 \end{multline*}

4. Равенство \eqref{sqrt[2n-1](0)=0} следует из тождества \eqref{0^b}.

5. Формулы \eqref{x^(1/(2n-1))->8(x->-8)}-\eqref{x^(1/(2n-1))->8(x->+8)}
следуют из формул \eqref{inf(y-in-J)f^(-1)(y)=inf-I}:
 \begin{multline*}
\lim_{x\to-\infty}
x^{\frac{1}{2n-1}}=\eqref{lim(x-to-a+0)f(x)=inf(x-in-(a;b))f(x)}=\inf_{x\in\R}
x^{\frac{1}{2n-1}}=\\=\eqref{inf(y-in-J)f^(-1)(y)=inf-I}=\inf\D(x\mapsto
x^{2n-1})=\inf\R=-\infty
 \end{multline*}
и
 \begin{multline*}
\lim_{x\to+\infty}
x^{\frac{1}{2n-1}}=\eqref{lim(x-to-a+0)f(x)=inf(x-in-(a;b))f(x)}=\sup_{x\in\R}
x^{\frac{1}{2n-1}}=\\=\eqref{inf(y-in-J)f^(-1)(y)=inf-I}=\sup\D(x\mapsto
x^{2n-1})=\sup\R=+\infty
 \end{multline*}
 \epr

\bex\label{EX:graph-of-x^(1/2n-1)} Предложение \ref{LM:x^(1/2n-1)-nepr} дает
достаточно информации для построения графика корня нечетной степени, за
исключением одной детали, а именно свойства графика иметь выпуклость вверх или
вниз на разных участках. Смысл этих слов интуитивно очевиден, однако аккуратное
объяснение им удобно давать после того, как будет определена производная,
поэтому мы отложим разговор о выпуклости до
\ref{SEC:grafik}\ref{SUBSEC:vypuklost} главы \ref{ch-graph-f(x)}. Если же не
требовать здесь объяснений относительно выпуклости, то в остальном из
предложения \ref{LM:x^(1/2n-1)-nepr} следует, что по виду графика (то есть по
тому, где у графика имеются интервалы возрастания и убывания, интервалы, где
сохраняется выпуклость, и какие функция имеет пределы на концах этих
интервалов) все корни нечетной степени ведут себя одинаково:
 \biter{

\item[---] график функции $x\mapsto\sqrt[2n-1]{x}$ выглядит как график любого
представителя из этого семейства функций, например, как график {\it кубического
корня}
$$
\sqrt[3]{x}=x^\frac{1}{3} \qquad (x\in\R)
$$

%\picture{120pt}{0pt}{x^2.pcx}
\vglue40pt
 }\eiter
 \eex

\bprop\label{LM:x^(1/2n)-nepr} Корень четной степени $x\mapsto \sqrt[2n]{x}$
 \biter{
\item[---] определен всюду на полуинтервале $[0;+\infty)$,

\item[---] возрастает на $[0;+\infty)$,

\item[---] непрерывен на $[0;+\infty)$,

\item[---] равен нулю в нуле:
 \beq\label{sqrt[2n](0)=0}
\sqrt[2n]{0}=0
 \eeq
\item[---] имеет бесконечный предел на бесконечности:
 \begin{align}
&\lim_{x\to+\infty} \sqrt[2n]{x}=+\infty \label{x^(1/(2n))->8(x->+8)}
 \end{align}

\item[---] удовлетворяет следующему правилу, которое может считаться его
определением:
 \beq\label{DF:sqrt[2n]}
y=\sqrt[2n]{x}\quad\Longleftrightarrow\quad y^{2n}=x\quad\&\quad y\ge 0
 \eeq
 }\eiter

\eprop
 \bpr

1. Тот факт, что функция $x\mapsto \sqrt[2n]{x}=x^{\frac{1}{2n}}$ определена на
$[0;+\infty)$ постулируется в Аксиоме степеней B1 на
с.\pageref{TH-o-step-otobr}, поэтому его доказывать не нужно.

2. Докажем соотношение \eqref{x^(1/(2n))->8(x->+8)}. При $x>1$ мы получим:
$$
\frac{1}{2n}>\frac{1}{2n+1}
$$
$$
\Downarrow
$$
$$
x^\frac{1}{2n}>\eqref{monot-a^x-0<a<1}>x^\frac{1}{2n+1}\overset{\eqref{x^(1/(2n-1))->8(x->+8)}}{\underset{x\to+\infty}{\longrightarrow}}+\infty
$$
$$
\Downarrow
$$
$$
x^\frac{1}{2n}\underset{x\to+\infty}{\longrightarrow}+\infty
$$

3. Равенство \eqref{sqrt[2n](0)=0} следует из тождества \eqref{0^b}.

4. Проверим, что функция $x\mapsto\sqrt[2n]{x}=x^\frac{1}{2n}$ возрастает на
$[0;+\infty)$. Для подмножества $(0;+\infty)$ этот факт постулируется в условии
монотонности \eqref{monot-x^b-b>0}. Поэтому нам нужно лишь проверить, что он
будет верен для случая пары $x,y\in[0;+\infty)$, $x<y$, в которой $x=0$. Здесь
применяется следствие \ref{COR:a>0=>a^b>0}:
$$
0<y\quad\Longrightarrow\quad 0<\eqref{x^b>0}<y^\frac{1}{2n}
$$

5. Докажем правило \eqref{DF:sqrt[2n]}. В прямую сторону мы получаем: если
$$
y=\sqrt[2n]{x}=x^{\frac{1}{2n}}
$$
то, во-первых,
$$
y^{2n}=\left(x^{\frac{1}{2n}}\right)^{2n}=\eqref{nakop-zakon}=x^{\frac{1}{2n}\cdot
2n}=x^1=x
$$
и, во-вторых, по Аксиоме степеней B1 на с.\pageref{TH-o-step-otobr}, из того,
что $x^{\frac{1}{2n}}$ существует, следует, что $x\ge 0$. Поэтому
$$
y=x^{\frac{1}{2n}}\overset{\eqref{a-ge-0=>a^b-ge-0}}{\ge}0
$$
Мы доказали импликацию
$$
y=\sqrt[2n]{x}\quad\Longrightarrow\quad y^{2n}=x\quad\&\quad y\ge 0
$$
Теперь докажем обратную импликацию:
$$
y=\sqrt[2n]{x}\quad\Longleftarrow\quad y^{2n}=x\quad\&\quad y\ge 0
$$
Пусть
$$
y^{2n}=x\quad\&\quad y\ge 0
$$
Тогда
$$
x=y^{2n}\overset{\eqref{a-ge-0=>a^b-ge-0}}{\ge}0
$$
и поэтому, по Аксиоме степеней B1 нас.\pageref{TH-o-step-otobr},
$$
\exists\
x^{\frac{1}{2n}}=(y^{2n})^{\frac{1}{2n}}=\eqref{nakop-zakon}=y^{2n\cdot\frac{1}{2n}}=y^1=y
$$
То есть
$$
y=\sqrt[2n]{x}
$$

6. Остается доказать непрерывность функции $x\mapsto \sqrt[2n]{x}$. Заметим,
что правило \eqref{DF:sqrt[2n]} можно интерпретировать так: функции
$$
f(x)=\sqrt[2n]{x}=x^\frac{1}{2n},\qquad g(y)=y^{2n},
$$
если их считать определенными на полуинтервале $[0;+\infty)$ (эта оговорка
нужна, потому что функция $g(y)=y^{2n}$ формально определена на всей прямой
$\R$, что больше, чем $[0;+\infty)$) обратны друг другу. Поэтому можно
воспользоваться замечанием
\ref{REM:Teor-ob-obratnoi-funktsii-na-poluint-v-R^1}: функции $f$ и $g$
определены на полуинтервале $[0;+\infty)$, взаимно обратны, причем функция $f$
возрастает на нем (мы это уже доказали на первом шаге), поэтому они должны быть
непрерывны.

Однако утверждение, упомянутое нами в замечании
\ref{REM:Teor-ob-obratnoi-funktsii-na-poluint-v-R^1} мы не доказали (и даже
аккуратно не сформулировали). Поэтому для строгости можно поступить иначе.
Сначала заметим, что функции $f$ и $g$, если их считать определенными на
интервале $(0;+\infty)$, взаимно обратны, и строго монотонны (в силу условия
монотонности \eqref{monot-x^b-b>0}). Значит, по теореме
\ref{Teor-ob-obratnoi-funktsii-na-int-v-R^1} они должны быть непрерывны на
$(0;+\infty)$. Поэтому для непрерывности на $[0;+\infty)$ нам остается
доказать, что $f$ непрерывна в нуле. При $0<x<1$ мы получим:
$$
\frac{1}{2n}>\frac{1}{2n+1}
$$
$$
\phantom{{\scriptsize \eqref{monot-a^x-0<a<1}}}\quad\Downarrow\quad{\scriptsize
\eqref{monot-a^x-0<a<1}}
$$
$$
0\overset{\eqref{x^b>0}}{<}x^\frac{1}{2n}<\underbrace{x^\frac{1}{2n+1}\underset{x\to
0}{\longrightarrow}0^\frac{1}{2n+1}}_{\scriptsize\begin{matrix}\text{функция
$x\mapsto x^\frac{1}{2n+1}$}\\ \text{непрерывна}\\ \text{по предложению
\ref{LM:x^(1/2n-1)-nepr}}\end{matrix}}\overset{\scriptsize\eqref{sqrt[2n](0)=0}}{=}0
$$
$$
\Downarrow
$$
$$
x^\frac{1}{2n}\underset{x\to 0}{\longrightarrow}0
$$
 \epr

\bex\label{EX:graph-of-x^(1/2n)} Предложение \ref{LM:x^(1/2n)-nepr} позволяет
строить графики корней четной степени. Опять если не интересоваться участками
выпуклости (о которых мы упоминали в примере \ref{EX:graph-of-x^(1/2n-1)}), мы
получаем, что по виду графика все корни четной степени ведут себя одинаково:
 \biter{

\item[---] график функции $x\mapsto\sqrt[2n]{x}$ выглядит как график любого
представителя из этого семейства функций, например, как график {\it квадратного
корня}
$$
\sqrt{x}=x^\frac{1}{2} \qquad (x\in\R)
$$

%\picture{120pt}{0pt}{x^2.pcx}
\vglue40pt
 }\eiter
 \eex

\bprop Справедливо тождество:
 \beq\label{|x|=sqrt-x^2}
|x|=\sqrt{x^2}
 \eeq
 \eprop
\bpr Здесь нужно рассмотреть три случая: во-первых, если $x=0$, то
$\sqrt{x^2}=0=|x|$ -- очевидно. Во-вторых, если $x>0$, то получаем
 $$
y=\sqrt{x^2}
 $$
 $$
 \Updownarrow
 $$
 $$
y^2=x^2\quad\&\quad y\ge 0
 $$
 $$
 \Updownarrow
 $$
 $$
y^2-x^2=0\quad\&\quad y\ge 0
 $$
 $$
 \Updownarrow
 $$
 $$
(y-x)(y+x)=0\quad\&\quad y\ge 0
 $$
 $$
 \Updownarrow
 $$
 $$
\left[\begin{matrix}y-x=0 \\ y+x=0\end{matrix}\right] \quad\&\quad y\ge 0
 $$
 $$
 \Updownarrow
 $$
 $$
\left[\begin{matrix}y=x \\
\boxed{y=-x}\end{matrix}\right] \put(-20,-30){\vector(-1,3){5}}\put(-50,-47){
 \boxed{\scriptsize\begin{matrix} \text{невозможно,}
\\ \text{потому что $x>0$} \end{matrix}}}
 \quad\&\quad y\ge 0
 $$
 $$
 \Updownarrow
 $$
 $$
y=x=|x|
 $$

В-третьих, если $x<0$, то
 $$
y=\sqrt{x^2}
 $$
 $$
 \Updownarrow
 $$
 $$
y^2=x\quad\&\quad y\ge 0
 $$
 $$
 \Updownarrow
 $$
 $$
y^2-x^2=0\quad\&\quad y\ge 0
 $$
 $$
 \Updownarrow
 $$
 $$
(y-x)(y+x)=0\quad\&\quad y\ge 0
 $$
 $$
 \Updownarrow
 $$
 $$
\left[\begin{matrix} y+x=0 \\ y-x=0\end{matrix}\right] \quad\&\quad y\ge 0
 $$
 $$
 \Updownarrow
 $$
 $$
\left[\begin{matrix} y=-x \\ \boxed{y=x}
\end{matrix}\right]\put(-20,-30){\vector(-1,3){5}}\put(-50,-47){
 \boxed{\scriptsize\begin{matrix} \text{невозможно,}
\\ \text{потому что $x<0$} \end{matrix}}}
 \quad\&\quad y\ge 0
 $$
 $$
 \Updownarrow
 $$
 $$
y=-x=|x|
 $$
\epr

\paragraph{Степенная функция.} Корни, о которых мы говорили в предыдущем разделе,
являются частными случаями функций, называемых степенными.

 \biter{

\item[$\bullet$] При каждом фиксированном значении $b$ функция
$$
x\mapsto x^b
$$
называется {\it степенной} (с показателем $b$). По Аксиоме степеней B1 на с.
\pageref{TH-o-step-otobr}, эту функцию можно считать определенной
 \biter{
\item[---] на множестве $\R$, если $b\in\frac{\Z_+}{2\N-1}$;

\item[---] на множестве $\R\setminus\{0\}$, если $b\in-\frac{\N}{2\N-1}$;

\item[---] на множестве $[0,+\infty)$, если
$b\in(0,+\infty)\setminus\frac{\Z}{2\N-1}$.

\item[---] на множестве $(0,+\infty)$, если и
$b\in(-\infty,0)\setminus\frac{\Z}{2\N-1}$.
 }\eiter
 }\eiter

\bprop\label{PROP:x^b-b-in-2N-1/2N-1} При $b\in\frac{2\N-1}{2\N-1}$ степенная
функция $x\mapsto x^b$
 \biter{

\item[--] определена всюду на $\R$,

\item[--] непрерывна на $\R$,

\item[--] возрастает на $\R$,

\item[--] равна нулю в нуле,
 $$
0^b=0,
 $$
\item[--] имеет бесконечные пределы на бесконечностях:
 \begin{align}
&\lim_{x\to-\infty} x^b=-\infty \label{x^b->8(x->-8)-b-in-2N-1/2N-1} \\
&\lim_{x\to+\infty} x^b=+\infty \label{x^b->8(x->+8)-b-in-2N-1/2N-1}
 \end{align}
 }\eiter
 \eprop
\bpr

1. Поскольку $b>0$, по Аксиоме степеней B1 на с.\pageref{TH-o-step-otobr},
функция $x\mapsto x^b$ действительно определена на $\R$.

2. Если $b=\frac{2m-1}{2n-1}$, $m,n\in\N$, то степенная функция $h(x)=x^b$
является композицией степенных функций $g(y)=y^{2m-1}$ и
$f(x)=x^{\frac{1}{2n-1}}$:
 \begin{multline*}
h(x)=g(f(x))=y^{2m-1}\Big|_{y=x^{\frac{1}{2n-1}}}=\\=(x^{\frac{1}{2n-1}})^{2m-1}=\eqref{nakop-zakon}=x^{\frac{2m-1}{2n-1}}
 \end{multline*}

3. Поскольку функция $f(x)=x^{\frac{1}{2n-1}}$ непрерывна на $\R$ по
предложению \ref{LM:x^(1/2n-1)-nepr}, а функция $g(y)=y^{2m-1}$ непрерывна на
$\R$ в силу примера \ref{nepr-x^n}, их композиция $h=g\circ f$ должна быть
непрерывна по теореме \ref{cont-composition}.

4. Функция $f(x)=x^{\frac{1}{2n-1}}$ возрастает всюду на $\R$ в силу
предложения \ref{LM:x^(1/2n-1)-nepr}, а функция $g(y)=y^{2m-1}$ возрастает
всюду на $\R$ в силу примера \ref{EX:x^(2n-1)-vozrastaet}. Значит, по свойству
$5^\circ$ на с.\pageref{vozrastanie=>neubyvanie}, композиция $h=g\circ f$ также
возрастает.

5. Равенство нулю в нуле опять следует из тождества \eqref{0^b}.

6. А пределы
\eqref{x^b->8(x->-8)-b-in-2N-1/2N-1}-\eqref{x^b->8(x->+8)-b-in-2N-1/2N-1}
следуют из \eqref{x^(1/(2n-1))->8(x->-8)}-\eqref{x^(1/(2n-1))->8(x->+8)} и
\eqref{x^(2n-1)->-8(x->-8)}-\eqref{x^(2n-1)->+8(x->+8)}. \epr

\bex Как и в случае с предложением \ref{PROP:x^b-b-in-2N/2N-1}, предложение
\ref{PROP:x^b-b-in-2N-1/2N-1} позволяет строить график степенной функции,
теперь уже с показателем $b\in\frac{2\N-1}{2\N-1}$. Опять вид графика зависит
от того, больше или меньше единицы степень $b$ (здесь она может быть равной 1,
поскольку $1\in\frac{2\N-1}{2\N-1}$), и, забегая вперед, мы употребляем термин
выпуклость, рассчитывая на интуицию (или эрудицию) читателя:
 \biter{

\item[---] в случае $b>1$, $b\in\frac{2\N-1}{2\N-1}$, график функции $x\mapsto
x^b$ является выпуклым вверх на полуинтервале $(-\infty;0]$ и выпуклым вниз на
полуинтервале $[0;+\infty)$, и выглядит как график любого представителя из
этого семейства функций, например, как {\it кубическая парабола}, то есть как
график функции:
$$
f(x)=x^3 \qquad (x\in\R)
$$

%\picture{120pt}{0pt}{x^2.pcx}
\vglue40pt

\item[---] в случае $b=1$ график функции $x\mapsto x^b=x^1=x$, понятное дело,
представляет собой {\it прямую}:
$$
f(x)=x \qquad (x\in\R)
$$

%\picture{120pt}{0pt}{x^2.pcx}
\vglue40pt

\item[---] если же $b<1$, $b\in\frac{2\N-1}{2\N-1}$, то график функции
$x\mapsto x^b$ является выпуклым вниз на полуинтервале $(-\infty;0]$ и выпуклым
вверх на полуинтервале $[0;+\infty)$, и выглядит как график любого
представителя из этого семейства, например, как график {\it кубического корня}:
$$
f(x)=x^{\frac{1}{3}}=\sqrt[3]{x} \qquad (x\in\R)
$$

%\picture{120pt}{0pt}{x^2.pcx}
 \vglue40pt
 }\eiter
\eex

\bprop\label{PROP:x^b-b-in-2N/2N-1} При $b\in\frac{2\N}{2\N-1}$ степенная
функция $x\mapsto x^b$
 \biter{

\item[--] определена всюду на $\R$,

\item[--] непрерывна на $\R$,

\item[--] убывает на полуинтервале $(-\infty;0]$,

\item[--] возрастает  на полуинтервале $[0;+\infty)$,

\item[--] равна нулю в нуле,
 $$
0^b=0,
 $$
\item[--] имеет бесконечные пределы на бесконечностях:
 \begin{align}
&\lim_{x\to-\infty} x^b=+\infty \label{x^b->+8(x->-8)-b-in-2N/2N-1} \\
&\lim_{x\to+\infty} x^b=+\infty \label{x^b->+8(x->+8)-b-in-2N/2N-1}
 \end{align}
 }\eiter
 \eprop
\bpr Здесь с очевидными изменениями мы применяем те же рассуждения, что и при
доказательстве предложения \ref{PROP:x^b-b-in-2N-1/2N-1}.

1. Поскольку $b>0$, по Аксиоме степеней B1 на с.\pageref{TH-o-step-otobr},
функция $x\mapsto x^b$ действительно определена на $\R$.

2. Заметим, что для $b=\frac{2m}{2n-1}$, $m,n\in\N$, степенная функция
$h(x)=x^b$ является композицией степенных функций $g(y)=y^{2m}$ и
$f(x)=x^{\frac{1}{2n-1}}$:
 \begin{multline*}
h(x)=g(f(x))=y^{2m}\Big|_{y=x^{\frac{1}{2n-1}}}=\\=(x^{\frac{1}{2n-1}})^{2m}=\eqref{nakop-zakon}=x^{\frac{2m}{2n-1}}
 \end{multline*}

3. Теперь проверим непрерывность. Поскольку функция $g(y)=y^{2m}$ непрерывна на
$\R$ в силу примера \ref{nepr-x^n}, а функция $f(x)=x^{\frac{1}{2n-1}}$
непрерывна на $\R$ по предложению \ref{LM:x^(1/2n-1)-nepr}, их композиция
$h=g\circ f$ должна быть непрерывна по теореме \ref{cont-composition}.

4. Далее разберемся с монотонностью. Функция $f(x)=x^{\frac{1}{2n-1}}$
возрастает на полуинтервале $(-\infty;0]$ в силу предложения
\ref{LM:x^(1/2n-1)-nepr}, а функция $g(y)=y^{2m}$ убывает на полуинтервале
$(-\infty;0]\supseteq f((-\infty;0])$ в силу примера
\ref{EX:x^(2n)-vozrastaet-x>0}. Значит, по свойству $5^\circ$ на
с.\pageref{vozrastanie=>neubyvanie}, композиция $h=g\circ f$ убывает на
полуинтервале $(-\infty;0]$.

Наоборот, функция $f(x)=x^{\frac{1}{2n-1}}$ возрастает на $[0;+\infty)$ в силу
предложения \ref{LM:x^(1/2n-1)-nepr}, а функция $g(y)=y^{2m}$ возрастает на
$[0;+\infty)\supseteq f([0;+\infty))$ в силу примера
\ref{EX:x^(2n)-vozrastaet-x>0}. Значит, по свойству $5^\circ$ на
с.\pageref{vozrastanie=>neubyvanie}, композиция $h=g\circ f$ также возрастает
на $[0;+\infty)$.

5. Равенство нулю в нуле следует из тождества \eqref{0^b}.

6. Пределы
\eqref{x^b->+8(x->-8)-b-in-2N/2N-1}-\eqref{x^b->+8(x->+8)-b-in-2N/2N-1} следуют
из \eqref{x^(1/(2n-1))->8(x->-8)}-\eqref{x^(1/(2n-1))->8(x->+8)} и
\eqref{x^(2n)->8(x->-8)}-\eqref{x^(2n)->8(x->+8)}. Например, первый из них:
$$
x\to-\infty
$$
$$
\phantom{{\scriptsize\eqref{x^b->+8(x->-8)-b-in-2N/2N-1}}}
\quad\Downarrow\quad{\scriptsize\eqref{x^b->+8(x->-8)-b-in-2N/2N-1}}
$$
$$
f(x)=x^{\frac{1}{2n-1}}\to-\infty
$$
$$
\phantom{{\scriptsize\text{\eqref{x^(2n)->8(x->-8)}, теорема
\ref{TH:zamena-perem-v-lim}}}}
\quad\Downarrow\quad{\scriptsize\text{\eqref{x^(2n)->8(x->-8)}, теорема
\ref{TH:zamena-perem-v-lim}}}
$$
$$
g(f(x))=x^{\frac{2n}{2n-1}}\to+\infty
$$
 \epr

\bex\label{EX:graph-of-x^(2m/2n-1)} Предложение \ref{PROP:x^b-b-in-2N/2N-1}
позволяет строить графики степенных функций с показателями
$b\in\frac{2\N}{2\N-1}$. Эти графики имеют уловимое взглядом различие, в
зависимости от того, больше или меньше единицы степень $b$ (равной 1 она быть
не может, из-за того, что $1\notin\frac{2\N}{2\N-1}$). Различие же состоит в
том, как устроены участки выпуклости, которые мы еще не успели определить, и о
которых речь пойдет в \ref{SEC:grafik}\ref{SUBSEC:vypuklost} главы
\ref{ch-graph-f(x)}. Поскольку интуитивный смысл выпуклости очевиден, мы
думаем, что читателю будет понятно, если, забегая вперед, сообщить ему, что при
$b>1$ график степенной функции является выпуклым вниз на прямой $\R$, а при
$0<b<1$ прямая разбивается на два полуинтервала $(-\infty,0]$ и $[0;+\infty)$,
на которых график является выпуклым вверх:
 \biter{

\item[---] в случае $b>1$, $b\in\frac{2\N}{2\N-1}$, график функции $x\mapsto
x^b$, выглядит как график любого представителя из этого семейства функций,
например, как {\it квадратичная парабола}, то есть как график функции
$$
f(x)=x^2 \qquad (x\in\R)
$$

%\picture{120pt}{0pt}{x^2.pcx}
\vglue40pt

\item[---] точно так же, если $b<1$, $b\in\frac{2\N}{2\N-1}$, график функции
$x\mapsto x^b$, выглядит как график любого представителя из этого семейства,
например, как {\it парабола со степенью $\frac{2}{3}$}, то есть как график
функции
$$
f(x)=x^{\frac{2}{3}} \qquad (x\in\R)
$$

%\picture{120pt}{0pt}{x^2.pcx}
 \vglue40pt
 }\eiter
 \eex

\bprop\label{PROP:x^b-b-in-(-2N/2N-1)} При $b\in-\frac{2\N}{2\N-1}$ степенная
функция $x\mapsto x^b$
 \biter{

\item[--] определена на множестве $\R\setminus\{0\}$,

\item[--] непрерывна на $\R\setminus\{0\}$,

\item[--] возрастает на интервале $(-\infty;0)$,

\item[--] убывает на интервале $(0;+\infty)$,

\item[--] имеет следующие пределы в граничных точках области определения и на
бесконечности:
 \begin{align}
&\lim_{x\to-0} x^b=+\infty \label{x^b->+8(x->-0)-b-in-(-2N/2N-1)} \\
&\lim_{x\to+0} x^b=+\infty \label{x^b->+8(x->+0)-b-in-(-2N/2N-1)} \\
&\lim_{x\to\infty} x^b=0 \label{x^b->0(x->8)-b-in-(-2N/2N-1)}
 \end{align}
 }\eiter
 \eprop
\bpr Если обозначить $c=-b$, то мы получим $c\in\frac{2\N}{2\N-1}$, и наша
функция будет композицией функции $f(x)=x^c$, рассмотренной в предложении
\ref{PROP:x^b-b-in-2N/2N-1} и функции $g(y)=y^{-1}$, рассмотренной в примерах
\ref{EX:monotonnost-1/x} и \ref{nepr-x^n}. Все нужные нам свойства функции
$x\mapsto x^b$ следуют из этих утверждений. \epr

\bex Предложение \ref{PROP:x^b-b-in-(-2N/2N-1)} позволяет строить график
степенной функции с показателем $b\in-\frac{2\N}{2\N-1}$. В отличие от
предыдущих случаев, здесь вид графика (в тех деталях, которые нас интересуют, а
именно, монотонность, выпуклость и пределы в крайних токах интервалов
монотонности и выпуклости) не зависит от значения $b$:
 \biter{

\item[---] если $b\in-\frac{2\N}{2\N-1}$, то график функции $x\mapsto x^b$,
выглядит как график любого представителя из этого семейства функций, например,
как {\it квадратичная гипербола}, то есть как график функции
$$
f(x)=x^{-2} \qquad (x\in\R)
$$

%\picture{120pt}{0pt}{x^2.pcx}
\vglue40pt
 }\eiter
\eex

\bprop\label{PROP:x^b-b-in-(-2N-1/2N-1)} При $b\in-\frac{2\N-1}{2\N-1}$
степенная функция $x\mapsto x^b$
 \biter{

\item[--] определена на множестве $\R\setminus\{0\}$,

\item[--] непрерывна на $\R\setminus\{0\}$,

\item[--] убывает на интервале $(-\infty;0)$,

\item[--] убывает на интервале $(0;+\infty)$,

\item[--] имеет следующие пределы в граничных точках области определения и на
бесконечности:
 \begin{align}
&\lim_{x\to-0} x^b=-\infty \label{x^b->+8(x->-0)-b-in-(-2N-1/2N-1)} \\
&\lim_{x\to+0} x^b=+\infty \label{x^b->+8(x->+0)-b-in-(-2N-1/2N-1)} \\
&\lim_{x\to\infty} x^b=0 \label{x^b->0(x->8)-b-in-(-2N-1/2N-1)}
 \end{align}
 }\eiter
 \eprop
\bpr Если обозначить $c=-b$, то мы получим $c\in\frac{2\N}{2\N-1}$, и наша
функция будет композицией функции $f(x)=x^c$, рассмотренной в предложении
\ref{PROP:x^b-b-in-2N-1/2N-1} и функции $g(y)=y^{-1}$, рассмотренной в примерах
\ref{EX:monotonnost-1/x} и \ref{nepr-x^n}. Все нужные нам свойства функции
$x\mapsto x^b$ следуют из этих утверждений. \epr

\bex Предложение \ref{PROP:x^b-b-in-(-2N-1/2N-1)} позволяет строить график
степенной функции с показателем $b\in-\frac{2\N-1}{2\N-1}$. Как и в предыдущем
случае, здесь вид графика не зависит от значения $b$:
 \biter{

\item[---] если $b\in-\frac{2\N-1}{2\N-1}$, то график функции $x\mapsto x^b$,
выглядит как график любого представителя из этого семейства функций, например,
как {\it кубическая гипербола}, то есть как график функции
$$
f(x)=x^{-3} \qquad (x\in\R)
$$

%\picture{120pt}{0pt}{x^2.pcx}
\vglue40pt
 }\eiter
\eex

\bprop\label{PROP:x^b-b-in-(0,8)-setminus-Z/2N-1} При
$b\in(0;+\infty)\setminus\frac{\Z}{2\N-1}$ степенная функция $x\mapsto x^b$
 \biter{

\item[--] определена на полуинтервале $[0;+\infty)$,

\item[--] непрерывна на $[0;+\infty)$,

\item[--] возрастает на $[0;+\infty)$,

\item[--] равна нулю в нуле,
 \beq\label{0^b=0,b-in-(0,8)}
0^b=0,
 \eeq
\item[--] имеет бесконечный предел в бесконечности:
 \beq\label{x->8=>x^b->8}
\lim_{x\to+\infty} x^b=+\infty
 \eeq
 }\eiter
 \eprop
\bpr Здесь мы с небольшими изменениями повторяем рассуждения, применявшиеся при
доказательстве предложения \ref{LM:x^(1/2n)-nepr}.

1. Поскольку $b>0$, по Аксиоме степеней B1 на с.\pageref{TH-o-step-otobr}
функция $x\mapsto x^b$ действительно определена на $[0;+\infty)$.

2. Докажем соотношение \eqref{x->8=>x^b->8}. Для этого подберем $n\in\N$ так,
чтобы $b>\frac{1}{2n-1}$. Тогда при $x>1$ мы получим:
$$
x^b>\eqref{monot-a^x-0<a<1}>x^\frac{1}{2n-1}\underset{x\to+\infty}{\longrightarrow}+\infty
$$
$$
\Downarrow
$$
$$
x^b\underset{x\to+\infty}{\longrightarrow}+\infty
$$

3. Равенство \eqref{0^b=0,b-in-(0,8)} следует из тождества \eqref{0^b}.

4. Проверим, что наша функция возрастает на $[0;+\infty)$. Для подмножества
$(0;+\infty)$ этот факт постулируется в условии монотонности
\eqref{monot-x^b-b>0}. Поэтому нам нужно лишь проверить, что он будет верен для
случая пары $x,y\in[0;+\infty)$, $x<y$, в которой $x=0$. Здесь применяется
следствие \ref{COR:a>0=>a^b>0}:
$$
0<y\quad\Longrightarrow\quad 0<\eqref{x^b>0}<y^b
$$

5. Заметим далее, что функции
$$
f(x)=x^b,\qquad g(y)=y^\frac{1}{b}
$$
если их считать определенными на полуинтервале $[0;+\infty)$ (эта оговорка
нужна, потому что может случиться, что вторую из них можно определить на
множестве, большем, чем $[0;+\infty)$; например, при
$b=\frac{1}{2}\in(0;+\infty)\setminus\frac{\Z}{2\N-1}$ мы получаем
$\frac{1}{b}=2\in\frac{\Z}{2\N-1}$, и функция $y\mapsto y^\frac{1}{b}=y^2$
оказывается определенной на всей прямой $\R$) обратны друг другу.
Действительно, с одной стороны, для любого $x\in[0;+\infty)$ мы получаем
 \begin{multline*}
g(f(x))=(x^b)^\frac{1}{b}=\eqref{nakop-zakon}=\\=x^{b\cdot\frac{1}{b}}=x^1=\eqref{x^1=x}=x
 \end{multline*}
А с другой -- для любого $y\in[0;+\infty)$
 \begin{multline*}
f(g(y))=(y^\frac{1}{b})^b=\eqref{nakop-zakon}=\\=y^{\frac{1}{b}\cdot
b}=y^1=\eqref{x^1=x}=y
 \end{multline*}

6. После этого можно воспользоваться замечанием
\ref{REM:Teor-ob-obratnoi-funktsii-na-poluint-v-R^1}: функции $f$ и $g$
определены на полуинтервале $[0;+\infty)$, взаимно обратны, причем функция $f$
возрастает на нем (мы это уже доказали на первом шаге), поэтому они должны быть
непрерывны.

Однако поскольку утверждение, упомянутое нами в замечании
\ref{REM:Teor-ob-obratnoi-funktsii-na-poluint-v-R^1} мы не доказали, мы для
строгости поступим иначе. Сначала заметим, что функции $f$ и $g$, если их
считать определенными на интервале $(0;+\infty)$, взаимно обратны, и строго
монотонны (в силу условия монотонности \eqref{monot-x^b-b>0}). Значит, по
теореме \ref{Teor-ob-obratnoi-funktsii-na-int-v-R^1} они должны быть непрерывны
на $(0;+\infty)$. Поэтому для непрерывности на $[0;+\infty)$ нам остается
доказать, что $f$ непрерывна в нуле. Для этого снова возьмем $n\in\N$ так,
чтобы $\frac{1}{2n-1}<b$. Тогда при $0<x<1$ мы получим:
$$
0<\eqref{x^b>0}<x^b<\eqref{monot-a^x-0<a<1}<x^\frac{1}{2n-1}\underset{x\to
0}{\longrightarrow}0
$$
$$
\Downarrow
$$
$$
x^b\underset{x\to 0}{\longrightarrow}0
$$
 \epr

\bex Предложение \ref{PROP:x^b-b-in-(0,8)-setminus-Z/2N-1} позволяет строить
график степенной функций с показателем
$b\in(0;+\infty)\setminus\frac{\Z}{2\N-1}$. Здесь как и в примере
\ref{EX:graph-of-x^(2m/2n-1)}, вид графика зависит от того, больше или меньше
единицы степень $b$ (равной 1 она опять быть не может, из-за того, что
$1\in\frac{\Z}{2\N-1}$):
 \biter{

\item[---] в случае $b>1$, $b\in(0;+\infty)\setminus\frac{\Z}{2\N-1}$, график
функции $x\mapsto x^b$, выглядит как график любого представителя из этого
семейства функций, например, как график функции
$$
f(x)=x^{\frac{3}{2}} \qquad (x\ge 0)
$$

%\picture{120pt}{0pt}{x^2.pcx}
\vglue40pt

\item[---] точно так же, если $b<1$,
$b\in(0;+\infty)\setminus\frac{\Z}{2\N-1}$, график функции $x\mapsto x^b$,
выглядит как график любого представителя из этого семейства, например, как
график квадратного корня
$$
f(x)=x^{\frac{1}{2}}=\sqrt{x} \qquad (x\ge 0)
$$

%\picture{120pt}{0pt}{x^2.pcx}
 \vglue40pt
 }\eiter
 \eex

\bprop\label{PROP:x^b-b-in-(-8,0)-setminus-Z/2N-1} При
$b\in(-\infty;0)\setminus\frac{\Z}{2\N-1}$ степенная функция $x\mapsto x^b$
 \biter{
\item[--] определена на интервале $(0;+\infty)$.

\item[--] непрерывна на $(0;+\infty)$.

\item[--] убывает на $(0;+\infty)$

\item[--] имеет следующие пределы на концах этого интервала:
 \begin{align}
&\lim_{x\to+0} x^b=+\infty \label{x^b->8(x->+0)-b<0} \\
&\lim_{x\to+\infty} x^b=0 \label{x^b->0(x->+8)-b<0}
 \end{align}
 }\eiter
 \eprop
\bpr Если обозначить $c=-b$, то мы получим $c\in\frac{2\N}{2\N-1}$, и наша
функция будет композицией функции $f(x)=x^c$, рассмотренной в предложении
\ref{PROP:x^b-b-in-(0,8)-setminus-Z/2N-1} и функции $g(y)=y^{-1}$,
рассмотренной в примерах \ref{EX:monotonnost-1/x} и \ref{nepr-x^n}. Все нужные
нам свойства функции $x\mapsto x^b$ следуют из этих утверждений. \epr

\bex Предложение \ref{PROP:x^b-b-in-(-8,0)-setminus-Z/2N-1} позволяет строить
график степенной функций с показателем
$b\in(-\infty;0)\setminus\frac{\Z}{2\N-1}$. Здесь вид графика (в тех деталях,
которые мы здесь уже перечисляли) не зависит от значения $b$:
 \biter{

\item[---] при $b\in(-\infty;0)\setminus\frac{\Z}{2\N-1}$, график функции
$x\mapsto x^b$, выглядит как график любого представителя из этого семейства
функций, например, как график функции
$$
f(x)=x^{-\frac{1}{2}} \qquad (x\in\R)
$$

%\picture{120pt}{0pt}{x^2.pcx}
 \vglue40pt
 }\eiter
 \eex

\paragraph{Показательная функция.}
При каждом фиксированном значении $a$ отображение
$$
x\mapsto a^x
$$
называется {\it показательной функцией} (с основанием $a$). По Аксиоме степеней
B1 на с.\pageref{TH-o-step-otobr}, эту функцию можно считать определенной
 \bit{
\item[---] на множестве $\R$, если $a>0$;

\item[---] на множестве $[0;+\infty)$, если $a=0$;

\item[---] на множестве $\frac{\Z}{2\N-1}$, если $a<0$.
 }\eit

\blm\label{LM:a^(1/n)->0} При любом $a>0$ последовательности $a^{\frac{1}{n}}$
и $a^{-\frac{1}{n}}$ стремятся к единице:
 \begin{align}\label{a^(1/n)->1}
& a^{\frac{1}{n}}\underset{n\to\infty}{\longrightarrow} 1, &&
a^{-\frac{1}{n}}\underset{n\to\infty}{\longrightarrow} 1
 \end{align}
\elm
 \bpr
1. Пусть сначала $a>1$. Обозначим
$$
x_n=a^{\frac{1}{n}}-1
$$
и заметим, что, во-первых,
$$
a>1
$$
$$
\phantom{\text{\scriptsize\eqref{monot-x^b-b>0}}}\quad\Downarrow\quad\text{\scriptsize\eqref{monot-x^b-b>0}}
$$
$$
a^{\frac{1}{n}}>1^{\frac{1}{n}}=\eqref{1^b}=1
$$
$$
\Downarrow
$$
$$
x_n=a^{\frac{1}{n}}-1>0
$$
И, во-вторых, в силу неравенства Бернулли \eqref{nerav-Bernoulli},
$$
a=(1+x_n)^n\ge \eqref{nerav-Bernoulli}\ge 1+n\cdot x_n
$$
$$
\Downarrow
$$
$$
\frac{a-1}{n}\ge x_n
$$
Записав эти два неравенства вместе, мы получим оценку для $x_n$, из которой
получается первый из двух пределов \eqref{a^(1/n)->1}:
$$
0<x_n\le\frac{a-1}{n}
$$
$$
\Downarrow
$$
$$
a^{\frac{1}{n}}-1=x_n\underset{n\to\infty}{\longrightarrow} 0
$$
$$
\Downarrow
$$
$$
a^{\frac{1}{n}}\underset{n\to\infty}{\longrightarrow} 1
$$
После этого второй предел получается применением тождества \eqref{a^(-x)} и
свойства пределов сохранять дробь:
$$
a^{-\frac{1}{n}}=\eqref{a^(-x)}=\frac{1}{a^{\frac{1}{n}}}\overset{\begin{matrix}\text{\scriptsize
$4^0$, с.\pageref{lim-ariphm}}\\
\end{matrix}}{\underset{n\to\infty}{\longrightarrow}} \frac{1}{1}=1
$$

2. Если $a=1$, то пределы \eqref{a^(1/n)->1} становятся тривиальными
следствиями тождества \eqref{1^b}:
$$
a^{\frac{1}{n}}=1^{\frac{1}{n}}=\eqref{1^b}=1\underset{n\to\infty}{\longrightarrow}
1
$$
$$
a^{-\frac{1}{n}}=1^{-\frac{1}{n}}=\eqref{1^b}=1\underset{n\to\infty}{\longrightarrow}
1
$$

3. Если $0<a<1$, то положив $A=\frac{1}{a}$, мы получим $A>1$, и поэтому, как
мы уже доказали,
$$
A^{\frac{1}{n}}\underset{n\to\infty}{\longrightarrow} 1
$$
Как следствие,
$$
a^{\frac{1}{n}}=\left(\frac{1}{A}\right)^{\frac{1}{n}}=\eqref{(1/x)^b}=
\frac{1}{A^{\frac{1}{n}}} \underset{n\to\infty}{\longrightarrow} \frac{1}{1}=1
$$
а отсюда уже
$$
a^{-\frac{1}{n}}=\eqref{a^(-x)}=\frac{1}{a^{\frac{1}{n}}}\underset{n\to\infty}{\longrightarrow}
\frac{1}{1}=1
$$
 \epr

\blm\label{LM:a^(1/n)->0} При любом $a>0$ справедливо соотношение:
 \beq\label{a^x->1}
 a^x\underset{x\to 0}{\longrightarrow} 1
 \eeq
\elm
 \bpr
Нам нужно показать, что для произвольной бесконечно малой последовательности
 \beq\label{x_k->0}
 x_k\underset{k\to\infty}{\longrightarrow} 0
 \eeq
выполняется
 \beq\label{a^(x_k)->1}
a^{x_k}\underset{k\to\infty}{\longrightarrow} 1
 \eeq
Рассмотрим несколько случаев.

1. Если $a=1$, то \eqref{a^(x_k)->1} выполняется тривиально, потому что в этом
случае $a^{x_k}=1$.

2. Пусть $a>1$. Из соотношений \eqref{a^(1/n)->1} следует, что для любого
$\e>0$ найдется номер $N$ такой, что
$$
1-\e<a^{-\frac{1}{N}}<a^{\frac{1}{N}}<1+\e
$$
Из \eqref{x_k->0} следует, что для данного $N$ можно подобрать номер $K$ такой,
что
$$
\forall k>K\qquad |x_k|<\frac{1}{N}
$$
Поэтому:
$$
\forall k>K\qquad -\frac{1}{N}<x_k<\frac{1}{N}
$$
$$
\phantom{\text{\scriptsize\eqref{monot-a^x-a>1}}}\quad\Downarrow\quad\text{\scriptsize\eqref{monot-a^x-a>1}}
$$
$$
\forall k>K\qquad 1-\e<a^{-\frac{1}{N}}<a^{x_k}<a^{\frac{1}{N}}<1+\e
$$
$$
\Downarrow
$$
$$
\forall k>K\qquad a^{x_k}\in(1-\e,1+\e)
$$
Мы получили, что какой ни возьми $\e>0$ почти все элементы последовательности
$a^{x_k}$ лежат в $\e$-окрестности точки $1$. Это и означает, что справедливо
\eqref{a^(x_k)->1}.

3. Пусть $0<a<1$. Тогда из \eqref{a^(1/n)->1} следует, что для любого $\e>0$
найдется номер $N$ такой, что
$$
1-\e<a^{\frac{1}{N}}<a^{-\frac{1}{N}}<1+\e
$$
Опять из \eqref{x_k->0} получаем, что для данного $N$ можно подобрать номер $K$
такой, что
$$
\forall k>K\qquad |x_k|<\frac{1}{N}
$$
и поэтому
$$
\forall k>K\qquad -\frac{1}{N}<x_k<\frac{1}{N}
$$
$$
\phantom{\text{\scriptsize\eqref{monot-a^x-0<a<1}}}\quad\Downarrow\quad\text{\scriptsize\eqref{monot-a^x-0<a<1}}
$$
$$
\forall k>K\qquad 1-\e<a^{\frac{1}{N}}<a^{x_k}<a^{-\frac{1}{N}}<1+\e
$$
$$
\Downarrow
$$
$$
\forall k>K\qquad a^{x_k}\in(1-\e,1+\e)
$$
То есть какой ни возьми $\e>0$, почти все элементы последовательности $a^{x_k}$
лежат в $\e$-окрестности точки $1$. Опять это означает, что справедливо
\eqref{a^(x_k)->1}.
 \epr

\bprop\label{PROP:nepr-a^x} При любом $a>0$ функция $x\mapsto a^x$ непрерывна.
\eprop
 \bpr
Если $x_k\underset{k\to\infty}{\longrightarrow} x$, то
$x_k-x\underset{k\to\infty}{\longrightarrow} 0$ и поэтому
$$
a^{x_k}-a^x=a^{x_k-x}\cdot
a^x-a^x=(\kern-13pt\underbrace{a^{x_k-x}}_{\scriptsize\begin{matrix}\downarrow\\
\phantom{,}1, \\ \text{в силу \eqref{a^x->1}}\end{matrix}}\kern-13pt-1)\cdot
a^x \underset{k\to\infty}{\longrightarrow} 0
$$
 \epr

\bprop\label{PROP:predely-a^x} Справедливы следующие соотношения:
 \biter{
\item[---] если $a>1$, то
 \begin{align}\label{a>1=>a^x->8}
& \lim_{x\to-\infty} a^x= 0, && \lim_{x\to+\infty} a^x=+\infty
 \end{align}

\item[---] если $0<a<1$, то
 \begin{align}\label{0<a<1=>a^x->0}
& \lim_{x\to-\infty} a^x=+\infty, && \lim_{x\to+\infty} a^x= 0
 \end{align}
 }\eiter
 \eprop
 \bpr
1. Пусть $a>1$. Тогда в силу \eqref{a>1=>a^n->infty},
$$
a^n\underset{n\to\infty}{\longrightarrow}+\infty,
$$
поэтому для любого $E>0$ найдется $N\in\N$ такое, что
$$
a^N>E
$$
Для любого $x>N$ мы теперь получаем:
$$
a^x>\eqref{monot-a^x-a>1}>a^N>E
$$
Мы получили, что для любого $E>0$ найдется $N$ такое, что при $x>N$ выполняется
неравенство $a^x>E$. Это означает, что справедливо второе соотношение в
\eqref{a>1=>a^x->8}:
$$
\lim_{x\to+\infty} a^x=+\infty
$$
После того, как оно доказано, первое выводится как его следствие:
 \begin{multline*}
\lim_{x\to-\infty} a^x=\begin{pmatrix}y=-x \\ y\to+\infty \end{pmatrix}=
\lim_{y\to+\infty} a^{-y}=\\=\lim_{y\to+\infty}\frac{1}{{\smsize
 \boxed{a^y}\put(2,-10){\vector(1,-1){10}\put(2,-15){$+\infty$}}
 }
}=\eqref{0<->infty-x}=0
 \end{multline*}

2. Пусть $0<a<1$. Тогда, положив $A=\frac{1}{a}$, мы получим $A>1$, поэтому
 \begin{multline*}
\lim_{x\to-\infty} a^x=\begin{pmatrix}y=-x \\ y\to+\infty \\ a^x=A^{-x}=A^y
\end{pmatrix}=\\=\lim_{y\to+\infty} A^y=\eqref{a>1=>a^x->8}=+\infty
 \end{multline*}
и
 \begin{multline*}
\lim_{x\to+\infty} a^x=\begin{pmatrix}y=-x \\ y\to-\infty \\ a^x=A^{-x}=A^y
\end{pmatrix}=\\=\lim_{y\to-\infty} A^y=\eqref{a>1=>a^x->8}=0
 \end{multline*}
 \epr

\bprop\label{PROP:a^x>0} При любых $a>0$ и $x\in\R$ справедливо неравенство
 \beq\label{a^x>0}
 a^x>0
 \eeq
 \eprop
\bpr Если $a=1$, то это неравенство будет очевидно,  в силу первого степенного
закона \eqref{1^b}: $1^x=1$. Поэтому важно рассмотреть случаи $a>1$ и $0<a<1$.

1. Пусть $a>1$. Тогда рассмотрев целую часть $[x]$ числа $x$ (мы ее определили
выше формулой \eqref{tselaya-chast-chisla}), мы получим:
 \begin{align*}
[x]&\le x && \text{\scriptsize \eqref{opr-tsel-chasti}}\\
&\Downarrow && \text{\scriptsize \eqref{monot-a^x-a>1}} \\
0\underset{\eqref{x^n>0}}{<}a^{[x]} & \le a^x && \\
&\Downarrow && \\
  0 &<a^x &&
 \end{align*}

2. Наоборот, если $0<a<1$, то
 \begin{align*}
x&<[x]+1 && \text{\scriptsize \eqref{opr-tsel-chasti}}\\
&\Downarrow && \text{\scriptsize \eqref{monot-a^x-0<a<1}} \\
a^x & > a^{[x]+1}\underset{\eqref{x^n>0}}{>}0 && \\
&\Downarrow && \\
 a^x &>0 &&
 \end{align*}
 \epr

\bex Предложений \ref{PROP:nepr-a^x}, \ref{PROP:predely-a^x} и \ref{PROP:a^x>0}
достаточно, чтобы построить график функции $x\mapsto a^x$ в тех случаях, когда
это возможно, то есть когда $a\ge 0$:

\begin{itemize}
\item[--] если $a>1$, график выглядит так:

%\picture{120pt}{0pt}{1_5^x.pcx}
\vglue40pt

\item[--] если $a=1$, график становится прямой:

%\picture{120pt}{0pt}{1^x.pcx}
\vglue40pt

\item[--] если $0<a<1$, график выглядит так:

%\picture{120pt}{0pt}{1_5^(-x).pcx}
\vglue40pt

\item[--] если $a=0$, график выглядит так:

%\picture{120pt}{0pt}{1^x.pcx}
\vglue40pt

\end{itemize}
\eex

\bex Отдельно нужно сказать, что при $a<0$, график изобразить невозможно,
оттого, что он перестает быть непрерывной линией; можно пытаться рисовать
отдельные точки на нем, например

%\picture{120pt}{0pt}{1^x.pcx}
\vglue40pt

увеличивая число точек, мы будем добиваться усложнения картинки

%\picture{120pt}{0pt}{1^x.pcx}
\vglue40pt

но, в отличие от предыдущих случаев, здесь не получится, чтобы все точки
принадлежали какой-то одной непрерывной линии.

Частным случаем будет график функции
$$
f(x)=(-1)^x
$$
В зависимости от степени подробности он может изображаться так

%\picture{120pt}{0pt}{1^x.pcx}
\vglue40pt

или так

%\picture{120pt}{0pt}{1^x.pcx}
\vglue40pt

или так

%\picture{120pt}{0pt}{1^x.pcx}
\vglue40pt

\eex

\paragraph{Логарифм.}\label{usloviya-sushestvovaniya-log}

 \biter{

\item[$\bullet$] Для всякого числа $a$, удовлетворяющего двойному неравенству
$$
0<a\ne 1,
$$
правило
 \beq\label{DF:log_a}
y=\log_a x\quad\Longleftrightarrow\quad a^y=x
 \eeq
определяет функцию
$$
x\mapsto \log_a x\qquad (x>0)
$$
называемую {\it логарифмом} по основанию $a$. Понятно, что функции $x\mapsto
a^x$ и $\log_a$ связаны тождествами
\begin{align}
& a^{\log_a x}=x && (x>0) \label{a^(log_a-x)=x}\\
& \log_a (a^x)=x && (x\in\R) \label{log_a(a^x)=x}
\end{align}
Второе из них влечет за собой тождество
 \begin{align}\label{log_a(a)=1}
& \log_a a=1 && (0<a\ne 1)
\end{align}
 }\eiter
 \bpr
Мы докажем это утверждение для случая $a>1$ (случай $0<a<1$ рассматривается
аналогично).

1. Неравенство \eqref{a^x>0} означает, что отображение $x\mapsto a^x$
действительно переводит $\R$ в $(0;+\infty)$.

2. Далее из условия возрастания \eqref{monot-a^x-a>1} следует, что отображение
$x\mapsto a^x$ инъективно.

3. Покажем, что $x\mapsto a^x$ сюръективно отображает $\R$ на $(0;+\infty)$.
Пусть $C\in(0;+\infty)$. Поскольку $\lim_{x\to-\infty} a^x=0$ (первое равенство
в \eqref{a>1=>a^x->8}), найдется $x\in\R$ такое, что $a^x<C$. С другой стороны,
$\lim_{x\to+\infty}a^x=+\infty$ (второе равенство в \eqref{a>1=>a^x->8}),
поэтому найдется $y\in\R$ такое, что $a^y>C$. Мы получаем, что
$$
a^x<C<a^y
$$
причем $x\mapsto a^x$ -- непрерывная функция в силу предложения
\ref{PROP:nepr-a^x}. Значит, по теореме Коши о среднем значении \ref{Cauchy-I},
найдется точка $c\in\R$ такая, что
$$
a^c=C
$$
Мы получили, что какое ни возьми $C\in(0;+\infty)$, для него найдется точка
$c\in\R$, в которой функция $x\mapsto a^x$ принимает значение $C$. Это и
означает сюръективность $x\mapsto a^x$ как отображения из $\R$ в $(0;+\infty)$.

4. Инъективность и сюръективность вместе означают биективность $x\mapsto
a^x:\R\to (0;+\infty)$.
 \epr

 \bprop\label{PROP:log-monot} Функция $\log_a$, определенная правилом \eqref{DF:log_a},
 \biter{
\item[---] возрастает при $a>1$:
 \beq\label{log_a-monot-a>1}
0<x<y\quad\Longrightarrow\quad \log_a x<\log_a y
 \eeq

\item[---] убывает при $0<a<1$:
 \beq\label{log_a-monot-0<a<1}
0<x<y\quad\Longrightarrow\quad \log_a x>\log_a y
 \eeq
 }\eiter
 \eprop
 \bpr
Пусть $a>1$ и $x<y$. Если бы оказалось, что $\log_a x\ge \log_a y$, то, в силу
возрастания $x\mapsto a^x$, мы получили бы
$$
\ln x\ge \ln y
$$
$$
\Downarrow
$$
$$
\underbrace{a^{\log_a x}}_{x}\ge\underbrace{a^{\log_a y}}_{y}
$$
$$
\Downarrow
$$
$$
x\ge y
$$
Случай $0<a<1$ рассматривается аналогично.
 \epr

 \bprop\label{PROP:log(xy)} Функция $\log_a$, определенная правилом \eqref{DF:log_a},
удовлетворяет следующему равенству и двум тождествам:
 \begin{align}\label{log(xy)}
& \log_a 1=0, \\
& \log_a\frac{1}{x}=-\log_a x \\ & \log_a(x\cdot y)=\log_a x+
\log_a y,
 \end{align}
 \eprop
\bpr Первое равенство следует непосредственно из определения \eqref{DF:log_a}:
$$
a^0=\eqref{a^0}=1 \qquad\Longleftrightarrow\qquad 0=\log_a 1
$$

Чтобы доказать второе, достаточно вычислить значения функции $x\mapsto a^x$ в
обеих его частях:
 $$
a^{\log_a\frac{1}{x}}=\frac{1}{x}=\frac{1}{a^{\log_a
x}}=\eqref{a^(-x)}=a^{-\log_a x}
 $$
Поскольку $x\mapsto a^x$ -- инъективное отображение, равенство
$a^{\log_a\frac{1}{x}}=a^{-\log_a x}$ означает, что аргументы должны быть
равны: $\log_a\frac{1}{x}=-\log_a x$.

Точно так же доказывается последнее равенство:
 \begin{multline*}
a^{\log_a(x\cdot y)}=x\cdot y=a^{\log_a x}\cdot a^{\log_a
y}=\\=\eqref{a^(x+y)}=a^{\log_a x+\log_a y}
 \end{multline*}
Опять, поскольку $x\mapsto a^x$ -- инъективное отображение, равенство
$a^{\log_a(x\cdot y)}=a^{\log_a x+\log_a y}$ возможно только если аргументы
равны: $\log_a(x\cdot y)=\log_a x+\log_a y$. \epr

\bprop\label{PROP:predely-log_a(x)} Справедливы следующие соотношения:
 \biter{
\item[---] если $a>1$, то
 \begin{align}\label{a>1=>log_a(x)->8}
& \lim_{x\to+0} \log_a x=-\infty, && \lim_{x\to+\infty}\log_a x=+\infty
 \end{align}

\item[---] если $0<a<1$, то
 \begin{align}\label{0<a<1=>log_a(x)->8}
& \lim_{x\to+0} \log_a x=+\infty, && \lim_{x\to+\infty}\log_a x=-\infty
 \end{align}
 }\eiter
 \eprop
\bpr Пусть $a>0$. Если $E>0$, то для любого $x\in(0;a^{-E})$ получаем в силу
\eqref{log_a-monot-a>1}:
$$
\log_a x<\log_a a^{-E}=-E
$$
Это означает, что справедливо первое соотношение в \eqref{a>1=>log_a(x)->8}:
$$
\lim_{x\to+0} \log_a x=-\infty
$$
Наоборот, для любого $x>a^E$ получаем в силу \eqref{log_a-monot-a>1}:
$$
\log_a x>\log_a a^E=E
$$
Это означает, что справедливо второе соотношение в \eqref{a>1=>log_a(x)->8}:
$$
\lim_{x\to+\infty} \log_a x=+\infty
$$
 \epr

\bex Предложений \ref{PROP:log-monot}, \ref{PROP:log(xy)} и
\ref{PROP:predely-log_a(x)} достаточно, чтобы построить график функции
$\log_a$:
\begin{itemize}
\item[---] если $a>1$, то картинка получается такой:

%\picture{120pt}{0pt}{ln_x.pcx}
\vglue40pt

\item[---] если $0<a<1$, то такой:

%\picture{120pt}{0pt}{-ln_x.pcx}
\vglue40pt

\end{itemize}
\eex

\paragraph{Формулы, связывающие показательные функции и логарифмы с разными основаниями.}

\bprop Справедливы тождества:
 \begin{align}
& a^x=b^{x\cdot \log_b a}   \label{a^x=b^(x-log_b-a)} \\
& \log_b (a^x)=x\cdot \log_b a  \label{log_b(a^x)=x-log_b-a} \\
&\log_a x=\frac{\log_b x}{\log_b a} \label{log_a-x=log_b-x/log_b-a}
 \end{align}
(в первых двух тождествах $0<x$, $0<a$, $0<b\ne 1$, а в третьем $0<x$, $0<a\ne
1$, $0<b\ne 1$) \eprop
 \bpr
Первое тождество доказывается цепочкой:
$$
b^{x\cdot \log_b a}=\eqref{nakop-zakon}=(\underbrace{b^{\log_b
a}}_{\scriptsize\begin{matrix}\| \\ \eqref{a^(log_a-x)=x}
\\ \| \\ a
\end{matrix}})^x=a^x
$$
Из него следует второе:
$$
\log_b (a^x)\overset{\eqref{a^x=b^(x-log_b-a)}}{=}\log_b (b^{x\cdot \log_b
a})\overset{\eqref{log_a(a^x)=x}}{=}x\cdot \log_b a
$$
А из второго следует третье:
$$
\log_b a\cdot \log_a x\overset{\eqref{log_b(a^x)=x-log_b-a}}{=}\log_b
(a^{\log_a x})=\log_b x
$$
$$
\Downarrow
$$
$$
\log_a x=\frac{\log_b x}{\log_b a}
$$
 \epr

\paragraph{Экспонента $e^x$ и натуральный логарифм $\ln x$.}

Напомним, что число $e$ мы определили в \ref{SEC:e} главы \ref{ch-x_n}, и там
же мы показали, что оно иррационально и лежит в интервале
 $$
2,708<e\le 2,720
 $$
Среди всех показательных функций  $x\mapsto a^x$ и всех логарифмов
$x\mapsto\log_a x$ функции с основанием $e$ занимают особое место, потому что,
с одной стороны, формулы для вычисления их значений, их производных и
интегралов, оказываются особенно простыми, а с другой -- все остальные
показательные функции и логарифмы выражаются через них с помощью формул
\eqref{a^x=b^(x-log_b-a)}-\eqref{log_a-x=log_b-x/log_b-a}.
 \biter{
\item[$\bullet$] Показательная функция $x\mapsto e^x$ с основанием $e$ имеет
специальное название -- {\it экспонента}. А логарифм $x\mapsto\log_e x$ с
основанием $e$ называется {\it натуральным логарифмом} и у него есть
специальное обозначение -- $\ln$:
 $$
\ln x=\log_e x
 $$
 }\eiter
Из \eqref{a^x=b^(x-log_b-a)}-\eqref{log_a-x=log_b-x/log_b-a} следует, что любую
показательную функцию и любой логарифм можно выразить через экспоненту и
натуральный логарифм так:
 \begin{align}
& a^x=e^{x\cdot \ln a} && (0<a)  \label{5.10.1} \\
&\log_a x=\frac{\ln x}{\ln a} && (0<a\ne 1)
 \end{align}
Графики эти функций выглядят как обычная показательная функция и обычный
логарифм:

%\picture{0pt}{0pt}{99.pcx}

\vglue140pt \noindent

\subsection{Тригонометрические функции}

\paragraph{Аксиома тригонометрии.}
Напомним, что синус $\sin t$ и косинус $\cos t$ определяются в школе как
ордината и абсцисса точки на единичной окружности, получаемой поворотом крайней
правой точки $(1;0)$ на угол $t$ против часовой стрелки (при этом угол обычно
измеряется в радианах, то есть в качестве угла берется длина дуги единичной
окружности):

\vglue150pt

 \noindent
Поскольку здесь не объясняется, что такое длина дуги, такое <<определение>>
синуса и косинуса, конечно, нельзя считать строгим. Однозначно эти функции
описываются следующей избыточной аксиомой:

\break

 \biter{\it

\item[\rm B2.] {\bf Аксиома тригонометрии.}
\label{osn-teor-trigonometrii}\footnote{Избыточность Аксиомы тригонометрии
доказывается в \ref{SEC:prilozh-step-ryadov}\ref{PR:osn-teor-trig} главы
\ref{CH-step-ryady}.} Существуют и единственны числовые функции синус $\sin$ и
косинус $\cos$, удовлетворяющие следующим условиям:
 \biter{
\item[$T_0$:] функции $\sin$ и $\cos$ определены всюду на $\R$ и являются
периодическими,

\item[$T_1$:] справедливы следующие три тождества $(x,y\in\R)$:
 \begin{align}
&\sin^2x+\cos^2x=1 \label{sin^t+cos^2t=1}\\
&\sin(x+y)=\sin x\cdot\cos y+\cos x\cdot\sin y \label{sin(x+y)} \\
&\cos(x+y)=\cos x\cdot\cos y-\sin x\cdot\sin y \label{cos(x+y)}
 \end{align}

\item[$T_2$:] при $x\in(0;1)$ выполняется следующее тройное неравенство:
 \beq\label{0<sin-x<x}
0<x\cos x<\sin x<x.
 \eeq
 }\eiter
 }\eiter

\paragraph{Синус и косинус.}
Это может быть неожиданным, но все остальные свойства синуса и косинуса
выводятся из этого списка. Мы покажем здесь, как это делается.

\bprop\label{PROP:sin-2x} Справедливы тождества:
\medskip
 \begin{align}
& \sin 2x=2\sin x\cdot\cos x \label{sin-2x} \\
& \cos 2x=2\cos^2x-1=1-2\sin^2 x \label{cos-2x}
 \end{align}\eprop
 \bpr Во-первых,
 \begin{multline*}
\sin 2x=\sin(x+x)=\eqref{sin(x+y)}=\\=\sin x\cdot\cos x+\sin x\cdot \cos x=
2\cdot\sin x\cdot\cos x
 \end{multline*}
Во-вторых,
 \begin{multline*}
\cos 2x=\cos(x+x)=\eqref{cos(x+y)}=\\=\cos x\cdot\cos x-\sin x\cdot\sin x=\\=
\cos^2 x-\underbrace{\sin^2 x}_{\scriptsize\begin{matrix}\| \\
\eqref{sin^t+cos^2t=1}\\ \| \\ 1-\cos^2 x \end{matrix}}=2\cos^2 x-1
 \end{multline*}
И отсюда уже следует:
 $$
\cos 2x=2\underbrace{\cos^2 x}_{\scriptsize\begin{matrix}\| \\
\eqref{sin^t+cos^2t=1}\\ \| \\ 1-\sin^2 x \end{matrix}}-1= 1-2\sin^2 x
 $$
 \epr

 \bprop
Функция $\sin$ обладает наименьшим периодом (а значит и наименьшим
полупериодом).
 \eprop
 \biter{
\item[$\bullet$] Наименьший полупериод функции $\sin$ обозначается буквой
$\pi$:
 \begin{multline}\label{DEF:pi}
\pi=\min\Big\{h>0:\ \forall x\in\R\\ \sin(x+2h)=\sin x\Big\}
 \end{multline}
Как следствие,
 \beq\label{2-pi-period-dlya-sin}
\sin(x+2\pi)=\sin x.
 \eeq
 }\eiter

\brem Формула для вычисления $\pi$ будет выведена нами только на
с.\pageref{ryad-dlya-pi}. Сейчас мы лишь отметим, что число $\pi$ иррационально
и удовлетворяет оценке
$$
3,14<\pi<3,15
$$
\erem

 \bpr
Рассмотрим число
$$
C=\sin\frac{1}{2}
$$
Из \eqref{0<sin-x<x} следует, что оно лежит в интервале
$\left(0;\frac{1}{2}\right)$:
$$
0<\sin\frac{1}{2}=C<\frac{1}{2}
$$
Положим $\e=\frac{C}{2}$ и рассмотрим интервал $(0;\e)$. Если предположить, что
у функции $\sin$ нет наименьшего периода, то по теореме
\eqref{TH:period-f-bez-naim-perioda}, найдется точка $t\in(0,\e)$ такая, что
$$
\sin t=C
$$
Но в то же время из $t\in(0,\e)\subseteq(0;1)$ в силу \eqref{0<sin-x<x}
следует:
$$
\sin t<t<\e=\frac{C}{2}<C
$$
Полученное противоречие означает, что наше предположение (что у $\sin$ нет
наименьшего периода) неверно. \epr

\bprop В точках $0$, $\pi$, $2\pi$ синус и косинус принимают следующие
значения:

\bigskip
 \vbox{\tabskip=0pt\offinterlineskip \halign to \hsize{ \vrule#\tabskip=2pt
 plus3pt minus1pt & \strut\hfil\;#\hfil & \vrule#&\hfil#\hfil &
 \vrule#&\hfil#\hfil &
 \vrule#&\hfil#\hfil & \vrule#\tabskip=0pt\cr \noalign{\hrule}
 height2pt&\omit&&\omit&&\omit&&\omit&\cr

 & $x$ && $0$ && $\pi$ && $2\pi$   &\cr \noalign{\hrule}
 height2pt&\omit&&\omit&&\omit&&\omit&\cr

 & $\sin x$ && $0$ && $0$ && $0$ &\cr \noalign{\hrule}
 height2pt&\omit&&\omit&&\omit&&\omit&\cr
 & $\cos x$ && $1$ && $-1$ && $1$ &\cr \noalign{\hrule}
 } }
\label{tablitsa-sin-cos-0-pi-2pi} \eprop
 \bpr
1. Сначала рассмотрим случай $x=0$. Во-первых,
$$
\sin 0=\sin(2\cdot 0)=\eqref{sin-2x}=2\cdot\sin 0\cdot\cos 0
$$
$$
\Downarrow
$$
$$
\sin 0-2\cdot\sin 0\cdot\cos 0=0
$$
$$
\Downarrow
$$
$$
\sin 0\cdot(1-2\cdot\cos 0)=0
$$
$$
\Downarrow
$$
 \beq\label{sin-0=0-ili-cos-0=1/2}
\left[\begin{matrix}\sin 0=0 \\ \cos 0=\frac{1}{2}\end{matrix}\right]
 \eeq
Но с другой стороны,
 \beq\label{cos-0=2-cos^2-0-1}
\cos 0=\cos(2\cdot 0)=\eqref{cos-2x}=2\cos^2 0-1
 \eeq
и отсюда следует, что $\cos 0$ не может быть равен $\frac{1}{2}$, потому что
при подстановке в \eqref{cos-0=2-cos^2-0-1} мы получаем неверное равенство:
$$
\frac{1}{2}=2\left(\frac{1}{2}\right)^2-1
$$
Таким образом, в совокупности \eqref{sin-0=0-ili-cos-0=1/2} второе равенство
невозможно, и остается первое:
$$
\sin 0=0.
$$

Отсюда уже получаем значение для $\cos 0$:
$$
\overbrace{\sin^2}^{0} 0+\cos^2 0=\eqref{sin^t+cos^2t=1}=1
$$
$$
\Downarrow
$$
$$
\cos^2 0=1
$$
$$
\Downarrow
$$
 \beq\label{cos-0=1-ili-cos-0=-1}
\left[\begin{matrix}\cos 0=1 \\ \cos 0=-1\end{matrix}\right]
 \eeq
Здесь снова если попробовать подставить в \eqref{cos-0=2-cos^2-0-1} вместо
$\cos 0$ значение $-1$, получится неверное равенство:
$$
-1=2\cdot(-1)^2-1,
$$
поэтому второе равенство в совокупности \eqref{cos-0=1-ili-cos-0=-1} неверно, и
остается первое:
$$
\cos 0=1
$$

2. После того, как значения в $x=0$ вычислены, из тождества периодичности для
синуса \eqref{2-pi-period-dlya-sin} сразу получается значение $\sin$ в точке
$x=2\pi$:
 $$
 \sin 2\pi=\sin(0+2\pi)=\sin 0=0
 $$
Отсюда следует тождество
 \begin{multline*}
\sin x=\sin(x+2\pi)=\sin x\cdot\cos 2\pi+\\+\cos x\cdot \underbrace{\sin
2\pi}_{0}=\sin x\cdot\cos 2\pi
 \end{multline*}
В частности, при $x=\frac{1}{2}$ получаем:
$$
\sin \frac{1}{2}=\sin \frac{1}{2}\cdot\cos 2\pi
$$
Из \eqref{0<sin-x<x} следует, что $\sin\frac{1}{2}>0$, потому на него можно
сократить:
$$
1=\cos 2\pi
$$
Таким образом, мы заполнили третий столбец в таблице.

3. Только после этого можно получить значения в точке $x=\pi$. Во-первых,
$$
0\kern-23pt\overset{\scriptsize\begin{matrix}\text{(уже вычислено)}\\
\downarrow
\end{matrix}}{=}\kern-23pt\sin2\pi=\eqref{sin-2x}=2\sin\pi\cdot\cos\pi
$$
$$
\Downarrow
$$
 $$
\left[\begin{matrix}\sin\pi=0 \\ \cos\pi=0\end{matrix}\right]
 $$
Второе из этих равенств -- $\cos\pi=0$ -- неверно, потому что если это значение
подставить в формулу \eqref{cos-2x} при $x=\pi$,
$$
\underbrace{\cos2\pi}_{\scriptsize\begin{matrix}\| \\ 1 \\ \text{(уже
вычислено)}
\end{matrix}}\kern-15pt=2\cos^2\pi-1
$$
получается неверное равенство:
$$
1=2\cdot 0^2-1
$$
Значит остается первое:
$$
\sin\pi=0
$$

Из этого равенства теперь получаем:
$$
\underbrace{\sin^2\pi}_{0}+\cos^2\pi=1
$$
$$
\Downarrow
$$
$$
\cos^2\pi=1
$$
$$
\Downarrow
$$
 \beq\label{cos-pi=1-ili-cos-pi=-1}
\left[\begin{matrix}\cos\pi=1 \\ \cos\pi=-1\end{matrix}\right]
 \eeq
Предположим, что верно первое: $\cos\pi=1$. Тогда в качестве следствия мы
получим, во-первых,
 \begin{multline*}
\sin(x+\pi)=\eqref{sin(x+y)}=\\=\sin x\cdot\kern-15pt\underbrace{\cos\pi}_{\scriptsize\begin{matrix}\|\\ 1\\
\text{(предположение)}\end{matrix}}\kern-15pt+\cos
x\cdot\kern-15pt\underbrace{\sin\pi}_{\scriptsize\begin{matrix}\|\\ 0\\
\text{(уже доказано)}\end{matrix}}\kern-15pt=\sin x
 \end{multline*}
То есть число $\pi>0$ должно быть периодом функции $\sin$. Но это невозможно,
потому что наименьшим периодом для $\sin$ является число $2\pi$. Полученное
противоречие означает, что наше предположение о равенстве $\cos\pi=1$ было
неверным.

Значит, в совокупности \eqref{cos-pi=1-ili-cos-pi=-1} должно быть верно второе
равенство:
$$
\cos\pi=-1
$$
 \epr

\bprop Число $\pi$ является полупериодом и для функции $\cos$:
 \beq\label{2-pi-period-dlya-cos}
\cos(x+2\pi)=\cos x.
 \eeq
\eprop
 \bpr Здесь используются уже найденные значения $\sin$ и $\cos$ в точке $2\pi$:
 \begin{multline*}
\cos(x+2\pi)=\eqref{cos(x+y)}=\\=\cos x\cdot\underbrace{\cos2\pi}_{1}-\sin
x\cdot\underbrace{\sin2\pi}_{0}=\cos x
 \end{multline*}
 \epr

\bprop Справедливы тождества:
\medskip
 \begin{align}\label{sin(-x)}
& \sin(-x)=-\sin x, && \cos(-x)=\cos x
 \end{align}\eprop
 \bpr
Зафиксируем точку $x\in\R$ и обозначим
$$
A=\sin(-x),\qquad B=\cos(-x)
$$
Тогда
 \begin{multline*}
 0=(\text{\scriptsize таблица на с.\pageref{tablitsa-sin-cos-0-pi-2pi}})=\sin 0=\\
 =\sin(x-x)=\sin x\cdot\underbrace{\cos(-x)}_{B}+\cos
 x\cdot\underbrace{\sin(-x)}_{A}=\\=
 B\cdot \sin x+A\cdot \cos x
 \end{multline*}
И, аналогично,
 \begin{multline*}
 1=(\text{\scriptsize таблица на с.\pageref{tablitsa-sin-cos-0-pi-2pi}})=\cos 0=\\
 =\cos(x-x)=\cos x\cdot\underbrace{\cos(-x)}_{B}-\sin x\cdot\underbrace{\sin(-x)}_{A}=\\=
 B\cdot \cos x-A\cdot \sin x
 \end{multline*}
Вместе эти равенства дают систему линейных уравнений для $A$ и $B$:
 $$
\begin{cases}B\cdot \sin x+A\cdot \cos x=0 \\
 B\cdot \cos x-A\cdot \sin x=1 \end{cases}
 $$
Ее можно по-разному решать, например, так:
 $$
\begin{cases}B\cdot \sin x+A\cdot \cos x=0 & \text{\scriptsize $\leftarrow$ умножаем на $\sin x$}\\
 B\cdot \cos x-A\cdot \sin x=1 & \text{\scriptsize $\leftarrow$ умножаем на $\cos x$} \end{cases}
 $$
 $$
 \Downarrow
 $$
 $$
 \begin{cases}
 B\cdot \sin^2 x+A\cdot\cos x\cdot\sin x=0 \\
 B\cdot \cos^2 x-A\cdot\sin x\cdot\cos x=\cos x
 \end{cases}
 $$
 $$
\phantom{\text{\scriptsize (складываем равенства)}}\quad \Downarrow\quad
\text{\scriptsize (складываем равенства)}
 $$
 $$
 \underbrace{B\cdot \sin^2 x+ B\cdot \cos^2 x}_{\scriptsize
 \begin{matrix}\| \\ B\cdot(\sin^2 x+ \cos^2 x)\\ \| \\ B \end{matrix}}=\cos x
 $$
 $$
 \Downarrow
 $$
 $$
 B=\cos x
 $$
Это нам и нужно было узнать про $B$. Точно так же и с $A$:
 $$
\begin{cases}B\cdot \sin x+A\cdot \cos x=0 & \text{\scriptsize $\leftarrow$ умножаем на $\cos x$}\\
 B\cdot \cos x-A\cdot \sin x=1 & \text{\scriptsize $\leftarrow$ умножаем на $\sin x$} \end{cases}
 $$
 $$
 \Downarrow
 $$
 $$
 \begin{cases}
 B\cdot \sin x\cdot\cos x+A\cdot\cos^2 x=0 \\
 B\cdot \cos x\cdot\sin x-A\cdot\sin^2 x=\sin x
 \end{cases}
 $$
 $$
\phantom{\text{\scriptsize (вычитаем второе из перво)}}\ \Downarrow\
\text{\scriptsize (вычитаем второе из первого)}
 $$
 $$
 \underbrace{A\cdot \cos^2 x+ A\cdot\sin^2 x}_{\scriptsize
 \begin{matrix}\| \\ A\cdot(\cos^2 x+\sin^2 x)\\ \| \\ A \end{matrix}}=-\sin x
 $$
 $$
 \Downarrow
 $$
 $$
 A=-\sin x
 $$
 \epr

Очевидным следствием \eqref{sin(-x)} будет

\bprop Справедливы тождества:
 \begin{align}
&\sin(x-y)=\sin x\cdot\cos y-\cos x\cdot\sin y \label{sin(x-y)} \\
&\cos(x-y)=\cos x\cdot\cos y+\sin x\cdot\sin y \label{cos(x-y)}
 \end{align}
 \eprop

Из этого в свою очередь следует еще несколько тождеств. Во-первых, тождества
для произведений:

\bprop Справедливы тождества:
 \begin{align}
& \sin x\cdot \sin y=\frac{1}{2}\Big(\cos(x-y)-\cos(x+y)\Big)\label{sin(x)sin(y)} \\
& \cos x\cdot \cos y=\frac{1}{2}\Big(\cos(x-y)+\cos(x+y)\Big)\label{cos(x)cos(y)} \\
& \sin x\cdot\cos
y=\frac{1}{2}\Big(\sin(x-y)+\sin(x+y)\Big)\label{sin(x)cos(y)}
 \end{align}
 \eprop
\bpr Мы докажем первое из этих тождеств, потому что остальные два доказываются
по аналогии. Здесь нужно просто вычислить выражение в скобках в правой части:
 \begin{multline*}
\cos(x-y)-\cos(x+y)=\eqref{cos(x-y)},\eqref{cos(x+y)}=\\= \underbrace{\cos
x\cdot\cos y}
\put(-19.3,-14){\put(-6.9,-6){\line(1,0){115}}\put(-9,-4){$\uparrow$}\put(106,-4){$\uparrow$}}
+\sin x\cdot\sin y\underbrace{-\cos x\cdot\cos y}+\\+\sin x\cdot\sin y=2\sin
x\cdot\sin y
 \end{multline*}
Поделив на 2, мы получим \eqref{sin(x)sin(y)}.
  \epr

И, во-вторых, тождества для сумм:

\bprop Справедливы тождества:
 \begin{align}
& \sin x+\sin y=2\cdot\sin\frac{x+y}{2}\cdot\cos\frac{x-y}{2}\label{sin-x+sin-y} \\
& \cos x+\cos y =2\cdot\cos\frac{x+y}{2}\cdot\cos\frac{x-y}{2}\label{cos-x+cos-y}\\
& \sin x-\sin y=2\cdot\cos\frac{x+y}{2}\cdot\sin\frac{x-y}{2}\label{sin-x-sin-y} \\
& \cos x-\cos y
=-2\cdot\sin\frac{x+y}{2}\cdot\sin\frac{x-y}{2}\label{cos-x-cos-y}
 \end{align}
 \eprop
 \bpr
Ниже нам понадобятся только два последних тождества, поэтому мы докажем их. Два
первых доказываются по аналогии. Обозначим
$$
\alpha=\frac{x+y}{2},\qquad \beta=\frac{x-y}{2}
$$
Тогда
$$
\alpha+\beta=x,\qquad \alpha-\beta=y
$$
Поэтому
 \begin{multline*}
\sin x-\sin y=\sin(\alpha+\beta)-\sin(\alpha-\beta)=\\=
\eqref{sin(x+y)},\eqref{sin(x-y)}=
\sin\alpha\cdot\cos\beta+\cos\alpha\cdot\sin\beta-\\-
\sin\alpha\cdot\cos\beta+\cos\alpha\cdot\sin\beta=
2\cos\alpha\cdot\sin\beta=\\=2\cos\frac{x+y}{2}\cdot\sin\frac{x-y}{2}
 \end{multline*}
И точно так же,
 \begin{multline*}
\cos x-\cos y=\cos(\alpha+\beta)-\cos(\alpha-\beta)=\\=
\eqref{cos(x+y)},\eqref{cos(x-y)}=
\cos\alpha\cdot\cos\beta-\sin\alpha\cdot\sin\beta-\\-
\cos\alpha\cdot\cos\beta-\sin\alpha\cdot\sin\beta=
-2\sin\alpha\cdot\sin\beta=\\=-2\sin\frac{x+y}{2}\cdot\sin\frac{x-y}{2}
 \end{multline*}
 \epr

\bprop Справедливы следующие неравенства $(x\in\R)$:
 \begin{align}
&|\sin x|\le 1,   \label{|sin-x|<1}\\
&|\cos x|\le 1, \label{|cos-x|<1}\\
&|\sin x|\le |x|. \label{|sin-x|<|x|}
 \end{align}
 \eprop
 \bpr
Из \eqref{sin^t+cos^2t=1} получаем
$$
\sin^2 x+\cos^2 x=1
$$
$$
\Downarrow
$$
$$
\sin^2 x\le 1
$$
$$
\Downarrow
$$
$$
|\sin x|=\sqrt{\sin^2 x}\le 1
$$
и точно так же доказывается \eqref{|cos-x|<1}.

Последнее неравенство доказывается отдельно для случаев $x=0$, $x\in (0;1)$,
$x\in(-1,0)$ и $x\notin(-1,1)$:
 \biter{

\item[--] при $x=0$ получаем
$$
|\sin 0|=0\le |0|
$$

\item[--] при $0<x<1$:
$$
0\le \sin x<x=|x|
$$
 $$
\phantom{\text{\scriptsize \eqref{0<sin-x<x}}}\Downarrow\ \text{\scriptsize
\eqref{0<sin-x<x}}
 $$
$$
|\sin x| \le |x|
$$

\item[--] при $-1<x<0$ получаем:
 $$
x\in (-1,0)
$$
 $$
 \Downarrow
 $$
 $$
-x\in (0,1)
$$
 $$
\phantom{\text{\scriptsize \eqref{0<sin-x<x}}}\Downarrow\ \text{\scriptsize
\eqref{0<sin-x<x}}
 $$
$$
0< \underbrace{\sin(-x)}_{\scriptsize\begin{matrix}\| \\ -\sin x \\ \| \\
|\sin x|
\end{matrix}}<\underbrace{-x}_{\scriptsize\begin{matrix}\| \\ |x|
\end{matrix}}
$$
 $$
 \Downarrow
 $$
 $$
 |\sin x| \le |x|
 $$

\item[--] и, наконец, если $x\notin(-1,+1)$, то $1\le |x|$ и поэтому
$$
|\sin x|\le 1\le |x|.
$$
 }\eiter
 \epr

 \bprop Справедливы неравенства
 \begin{align}
&|\sin x-\sin a|\le |x-a| \label{3.3.3} \\
&|\cos x-\cos a|\le |x-a| \label{3.3.4}
 \end{align}
 \eprop
 \bpr
 \begin{multline*}
|\sin x-\sin a|=\eqref{sin(x-y)}=\\= \left|\ 2 \cdot \cos \frac{x+a}{2}\cdot
\sin \frac{x-a}{2}\right|=\\= 2 \cdot \left|\ \cos \frac{x+a}{2}\right|\cdot
\left|\ \sin \frac{x-a}{2}\right|\le \\ \le
\eqref{|cos-x|<1},\eqref{|sin-x|<|x|} \le 2 \cdot 1\cdot
\left|\frac{x-a}{2}\right|= |x-a|
 \end{multline*}

\begin{multline*}
|\cos x-\cos a|= \eqref{cos(x-y)}=\\= \left|-2 \cdot \sin \frac{x+a}{2}\cdot
\sin \frac{x-a}{2}\right|=\\= 2 \cdot \left|\ \sin \frac{x+a}{2}\right|\cdot
\left|\ \sin \frac{x-a}{2}\right|\le\\ \le
\eqref{|sin-x|<1},\eqref{|sin-x|<|x|}\le 2 \cdot 1\cdot
\left|\frac{x-a}{2}\right|= |x-a|
\end{multline*}

 \epr

 \bprop\label{PROP:sin-cos-nepreryvny}
Функции $\sin$ и $\cos$ непрерывны на $\R$.
 \eprop
 \bpr
Для любого числа $a\in \R$ и любой последовательности $x_n\underset{n\to
\infty}{\longrightarrow} a$ получаем:
 \begin{gather*}\quad x_n\underset{n\to \infty}{\longrightarrow} a \\
\Downarrow \put(20,0){\text{\smsize$\begin{pmatrix}
\text{теорема \ref{x->a<=>x-a->0}}\end{pmatrix}$}}\\
\quad x_n-a\underset{n\to \infty}{\longrightarrow} 0 \\
\Downarrow\put(20,0){\text{\smsize
$\begin{pmatrix}\text{вспоминаем}\\ \text{определение}\\
\text{бесконечно малой}\\
\text{пос-ти на с.\pageref{->0-&-->infty}}\end{pmatrix}$}}\\
|x_n-a|\underset{n\to \infty}{\longrightarrow} 0\\
\Downarrow\put(20,0){\text{\smsize $\begin{pmatrix}\text{применяем}\\
\text{неравенство \eqref{3.3.3}}\end{pmatrix}$}}\\
0\le |\sin x_n-\sin a|\le |x_n-a| \underset{n\to
\infty}{\longrightarrow} 0\\
\Downarrow\put(20,0){\text{\smsize $\begin{pmatrix}
\text{теорема о двух}\\ \text{милиционерах \ref{milit}}\end{pmatrix}$}}\\
|\sin x_n-\sin a| \underset{n\to \infty}{\longrightarrow} 0\\
\Downarrow\put(20,0){\text{\smsize $\begin{pmatrix}
\text{определение}\\ \text{бесконечно малой}\\
\text{пос-ти на с.
\pageref{->0-&-->infty}}\end{pmatrix}$}}\\
\sin x_n-\sin a \underset{n\to \infty}{\longrightarrow} 0\\
\Downarrow\put(20,0){\text{\smsize $\begin{pmatrix} \text{теорема \ref{x->a<=>x-a->0}}\end{pmatrix}$}}\\
\sin x_n \underset{n\to \infty}{\longrightarrow}\sin a
 \end{gather*}
Мы получили, что если $x_n\underset{n\to \infty}{\longrightarrow} a$, то $\sin
x_n\underset{n\to \infty}{\longrightarrow}\sin a$. Это означает, что функция
$f(x)=\sin x$ непрерывна в точке $a\in \R$. Поскольку $a$ -- произвольная точка
на $\R$, функция $\sin$ непрерывна всюду на $\R$.

Для функции $\cos$ этот факт доказывается аналогично, только нужно вместо
неравенства \eqref{3.3.3} применять \eqref{3.3.4}. \epr

\blm\label{LM:poluperiod-sin-cos} Если в некоторой точке $h>0$ синус равен
нулю,
$$
\sin h=0,
$$
то $h$ -- полупериод для функций $\sin$ и $\cos$. \elm
 \bpr
$$
\bigg\{\begin{matrix}
\sin 2h=2\cdot\overbrace{\sin h}^{0}\cdot\cos h=0 \\
\cos 2h=1-2\cdot\underbrace{\sin^2 h}_{0}=1
\end{matrix}\bigg\}
$$
$$
\Downarrow
$$
$$
\bigg\{\begin{matrix}
\sin(x+2h)=\sin x\cdot\overbrace{\cos 2h}^{1}+\cos x\cdot\overbrace{\sin 2h}^{0}=\sin x \\
\cos(x+2h)=\cos x\cdot\underbrace{\cos 2h}_{1}-\sin x\cdot\underbrace{\sin
2h}_{0}=\cos x
\end{matrix}\bigg\}
$$
 \epr

\bprop\label{PROP:sin>0} На интервале $(0,\pi)$ синус положителен, а на
интервале $(-\pi,0)$ -- отрицателен:
 \begin{align}
& x\in(0,\pi) && \Longrightarrow &&\sin x>0 \label{sin-x>0}\\
& x\in(-\pi,0)&& \Longrightarrow &&\sin x<0 \label{sin-x<0}
 \end{align}
 \eprop
\bpr Из-за тождества нечетности \eqref{sin(-x)} здесь достаточно доказать
первое утверждение. Предположим, что оно неверно, то есть что для некоторого
$x\in(0,\pi)$ выполняется $\sin x\le 0$. Рассмотрим отдельно два случая: когда
$\sin x=0$, и когда $\sin x<0$.

1) Если $\sin x=0$, то по лемме \ref{LM:poluperiod-sin-cos} число $x$ должно
быть полупериодом для функций $\sin$ и $\cos$. Это невозможно, потому что
$x<\pi$, то есть $x$ меньше наименьшего полупериода $\pi$.

2) Если же $\sin x<0$, то выберем какое-нибудь число $\e>0$ так, чтобы
выполнялись два условия:
$$
\e<1,\qquad \e<x
$$
Из первого неравенства, в силу \eqref{0<sin-x<x}, мы получим:
$$
\sin\e>0
$$
То есть на отрезке $[\e,x]$ функция $\sin$ меняет знак:
$$
\sin\e>0>\sin x
$$
Значит, по теореме Коши о промежуточном значении \ref{Cauchy-I}, в некоторой
точке $h\in(\e,x)$ она обращается в нуль:
$$
\sin h=0
$$
Но тогда снова по лемме \ref{LM:poluperiod-sin-cos} получается, что число $h$
должно быть полупериодом для функции $\sin$, что невозможно, потому что
$h<x<\pi$, то есть $h$ меньше наименьшего полупериода $\pi$.
 \epr

\bprop На отрезке $[0,\pi]$ косинус строго убывает, а на отрезке $[-\pi,0]$
монотонно возрастает:
 \begin{align}
& x,y\in [0,\pi], && x<y && \Longrightarrow && \cos x>\cos y \\
& x,y\in [-\pi,0], && x<y && \Longrightarrow && \cos x<\cos y
 \end{align}
\eprop
 \bpr
Заметим сначала, что из-за тождества четности \eqref{sin(-x)} здесь достаточно
доказать только утверждение для отрезка $[0,\pi]$. Кроме того, поскольку
косинус непрерывен, его пределы в граничных точках этого отрезка равны его
значениям в этих точках:
$$
\lim_{x\to 0}\cos x=\cos 0,\qquad \lim_{x\to\pi}\cos x=\cos\pi
$$
Отсюда следует, что можно в доказательстве заменить отрезок $[0,\pi]$ на
интервал $(0,\pi)$: если косинус строго убывает на интервале $(0,\pi)$, то он
строго убывает и на отрезке $[0,\pi]$. Для интервала $(0,\pi)$ получаем такую
цепочку:
$$
x,y\in (0,\pi),\qquad x<y
$$
$$
\Downarrow
$$
$$
\frac{x+y}{2}\in(0,\pi),\qquad \frac{x-y}{2}\in(-\pi,0)
$$
$$
\Downarrow
$$
 \begin{multline*}
\cos x-\cos y =\eqref{cos-x-cos-y}=\\=-2\cdot
\underbrace{\sin\frac{x+y}{2}}_{\scriptsize\begin{matrix}\text{\rotatebox{90}{$<$}}
\\ 0\end{matrix}}\cdot\underbrace{\sin\frac{x-y}{2}}_{\scriptsize\begin{matrix}\text{\text{\rotatebox{90}{$>$}}}
\\ 0\end{matrix}}>0
 \end{multline*}
$$
\Downarrow
$$
$$
\cos x>\cos y
$$
 \epr

\bprop\label{tablitsa-sin-cos} Справедлива следующая таблица значений синуса и
косинуса в некоторых характерных точках на отрезке $[0,\frac{\pi}{2}]$:

 \bigskip
 \vbox{\tabskip=0pt\offinterlineskip \halign to \hsize{ \vrule#\tabskip=2pt
 plus3pt minus1pt & \strut\hfil\;#\hfil & \vrule#&\hfil#\hfil &
 \vrule#&\hfil#\hfil & \vrule#&\hfil#\hfil & \vrule#&\hfil#\hfil &
 \vrule#&\hfil#\hfil & \vrule#\tabskip=0pt\cr \noalign{\hrule}
 height2pt&\omit&&\omit&&\omit&&\omit&&\omit&&\omit&\cr

 & $t$ && $0$ && $\frac{\pi}{6}$ && $\frac{\pi}{4}$ & & $\frac{\pi}{3}$ && $\frac{\pi}{2}$  &\cr \noalign{\hrule}
 height2pt&\omit&&\omit&&\omit&&\omit&&\omit&&\omit&\cr

 & $\cos t$ && $1$ && $\frac{\sqrt{3}}{2}$ && $\frac{\sqrt{2}}{2}$ && $\frac{1}{2}$ & & $0$ &\cr \noalign{\hrule}
 height2pt&\omit&&\omit&&\omit&&\omit&&\omit&&\omit&\cr

 & $\sin t$ && $0$ && $\frac{1}{2}$ && $\frac{\sqrt{2}}{2}$ && $\frac{\sqrt{3}}{2}$ && $1$ &\cr \noalign{\hrule} } }

 \bigskip
 \eprop
\bpr Мы покажем, как выводятся значения этих функций в двух точках, оставив
вывод для остальных точек читателю.

1. Вычислим значения в точке $\frac{\pi}{2}$: во-первых,
$$
0=\sin\pi=2\cdot\sin\frac{\pi}{2}\cdot\cos\frac{\pi}{2}
$$
$$
\Downarrow
$$
$$
\left[\begin{matrix}\sin\frac{\pi}{2}=0 & \text{\scriptsize $\leftarrow$ невозможно, в силу \eqref{sin-x>0}}\\
\cos\frac{\pi}{2}=0 &
\end{matrix}\right]
$$
$$
\Downarrow
$$
$$
\cos\frac{\pi}{2}=0
$$
И, во-вторых,
$$
\sin^2\frac{\pi}{2}+\underbrace{\cos^2\frac{\pi}{2}}_{0}=1
$$
$$
\Downarrow
$$
$$
\sin^2\frac{\pi}{2}=1
$$
$$
\Downarrow
$$
$$
\left[\begin{matrix}\sin\frac{\pi}{2}=1 & \\
\sin\frac{\pi}{2}=-1 & \text{\scriptsize $\leftarrow$ невозможно, в силу
\eqref{sin-x>0}}
\end{matrix}\right]
$$
$$
\Downarrow
$$
$$
\sin\frac{\pi}{2}=1
$$

2. Заметим, что $\sin\frac{\pi}{4}=\cos\frac{\pi}{4}$:
 \begin{multline*}
\sin\frac{\pi}{4}=\sin\left(\frac{\pi}{2}-\frac{\pi}{4}\right)=\underbrace{\sin\frac{\pi}{2}}_{1}\cdot\cos\frac{\pi}{4}
-\\- \underbrace{\cos\frac{\pi}{2}}_{0}\cdot\sin\frac{\pi}{4}=\cos\frac{\pi}{4}
 \end{multline*}
Отсюда получаем:
$$
1=\sin\frac{\pi}{2}=\sin\left(\frac{\pi}{4}+\frac{\pi}{4}\right)=2\cdot\sin\frac{\pi}{4}\cdot\cos\frac{\pi}{4}=2\cdot
\sin^2\frac{\pi}{4}
$$
$$
\Downarrow
$$
$$
\sin^2\frac{\pi}{4}=\frac{1}{2}
$$
$$
\Downarrow
$$
$$
\left|\sin\frac{\pi}{4}\right|=\sqrt{\frac{1}{2}}=\frac{\sqrt{2}}{2}
$$
$$
\Downarrow
$$
$$
\left[\begin{matrix} \sin\frac{\pi}{4}=\frac{\sqrt{2}}{2}& \\
\sin\frac{\pi}{4}=-\frac{\sqrt{2}}{2}& {\scriptsize \gets \text{невозможно, в
силу \eqref{sin-x>0}}}\end{matrix}\right]
$$
$$
\Downarrow
$$
$$
\sin\frac{\pi}{4}=\cos\frac{\pi}{4}=\frac{\sqrt{2}}{2}
$$
 \epr

 \bprop
Справедливы тождества:
 \begin{align}
& \sin(x+\pi)=-\sin x \label{sin(x+pi)} \\
& \cos(x+\pi)=-\cos x \label{cos(x+pi)} \\
& \sin\left(x+\frac{\pi}{2}\right)=\cos x \label{sin(x+pi/2)} \\
& \cos\left(x+\frac{\pi}{2}\right)=-\sin x \label{cos(x+pi/2)} \\
& \sin\left(x-\frac{\pi}{2}\right)=-\cos x \label{sin(x-pi/2)} \\
& \cos\left(x-\frac{\pi}{2}\right)=\sin x \label{cos(x-pi/2)} \\
& \sin\left(\frac{\pi}{2}-x\right)=\cos x \label{sin(pi/2-x)} \\
& \cos\left(\frac{\pi}{2}-x\right)=\sin x \label{cos(x-pi/2)}
 \end{align}
 \eprop
 \bpr
Мы докажем только первое тождество, остальные доказываются по аналогии:
$$
\sin(x+\pi)=\sin x\cdot\underbrace{\cos\pi}_{-1}+\cos
x\cdot\underbrace{\sin\pi}_{0}=\sin x
$$
 \epr

Индукцией из формул \eqref{sin(x+pi)} и \eqref{cos(x+pi)} выводится

\bcor Справедливы тождества:
 \begin{align}
& \cos (\pi n)=\cos (-\pi n)=(-1)^n, \label{cos(pi-n)} \\
& \sin (\pi n)=\sin (-\pi n)=0. \label{sin(pi-n)}
 \end{align}
\ecor

\bprop\label{PROP:cos>0} На интервале
$\left(-\frac{\pi}{2},\frac{\pi}{2}\right)$ косинус положителен, а на интервале
$\left(\frac{\pi}{2},\frac{3\pi}{2}\right)$ -- отрицателен:
 \begin{align}
& x\in\left(-\frac{\pi}{2},\frac{\pi}{2}\right) && \Longrightarrow &&\cos x>0 \label{cos-x>0}\\
& x\in\left(\frac{\pi}{2},\frac{3\pi}{2}\right)&& \Longrightarrow &&\cos x<0
\label{cos-x<0}
 \end{align}
 \eprop
 \bpr
Здесь все следует из предложения \ref{PROP:sin>0}:
$$
x\in\left(-\frac{\pi}{2},\frac{\pi}{2}\right)
$$
$$
\Downarrow
$$
$$
x+\frac{\pi}{2}\in (0,\pi)
$$
$$
\Downarrow
$$
$$
\cos x=\eqref{sin(x+pi/2)}=\sin\left(x+\frac{\pi}{2}\right)>0
$$
И точно так же
$$
x\in\left(\frac{\pi}{2},\frac{3\pi}{2}\right)
$$
$$
\Downarrow
$$
$$
x+\frac{\pi}{2}\in (\pi,2\pi)
$$
$$
\Downarrow
$$
$$
\cos x=\eqref{sin(x+pi/2)}=\sin\left(x+\frac{\pi}{2}\right)<0
$$
 \epr

\bprop\label{PROP:int-monot-sin} На отрезке
$\left[-\frac{\pi}{2},\frac{\pi}{2}\right]$ синус монотонно возрастает, а на
отрезке $\left[\frac{\pi}{2},\frac{3\pi}{2}\right]$ монотонно убывает:
 \begin{align}
& x,y\in \left[-\frac{\pi}{2},\frac{\pi}{2}\right], && x<y && \Longrightarrow && \sin x>\sin y \\
& x,y\in \left[\frac{\pi}{2},\frac{3\pi}{2}\right], && x<y && \Longrightarrow
&& \sin x<\sin y
 \end{align}
\eprop
 \bpr
Для отрезка $\left[-\frac{\pi}{2},\frac{\pi}{2}\right]$ получаем:
$$
x,y\in\left[-\frac{\pi}{2},\frac{\pi}{2}\right],\qquad x<y
$$
$$
\Downarrow
$$
$$
x+\frac{\pi}{2},y+\frac{\pi}{2}\in[0,\pi],\qquad
x+\frac{\pi}{2}<y+\frac{\pi}{2}
$$
$$
\Downarrow
$$
$$
\cos\left(x+\frac{\pi}{2}\right)>\cos\left(y+\frac{\pi}{2}\right)
$$
$$
\Downarrow
$$
 \begin{multline*}
\sin
x=\eqref{cos(x+pi/2)}=-\cos\left(x+\frac{\pi}{2}\right)<\\
<-\cos\left(y+\frac{\pi}{2}\right)=\eqref{cos(x+pi/2)}=\sin
y
 \end{multline*}
А для $\left[\frac{\pi}{2},\frac{3\pi}{2}\right]$:
$$
x,y\in\left[\frac{\pi}{2},\frac{3\pi}{2}\right],\qquad x<y
$$
$$
\Downarrow
$$
$$
x+\frac{\pi}{2},y+\frac{\pi}{2}\in(\pi,2\pi),\qquad
x+\frac{\pi}{2}<y+\frac{\pi}{2}
$$
$$
\Downarrow
$$
$$
\cos\left(x+\frac{\pi}{2}\right)<\cos\left(y+\frac{\pi}{2}\right)
$$
$$
\Downarrow
$$
 \begin{multline*}
\sin
x=\eqref{cos(x+pi/2)}=-\cos\left(x+\frac{\pi}{2}\right)>\\
>-\cos\left(y+\frac{\pi}{2}\right)=\eqref{cos(x+pi/2)}=\sin y
 \end{multline*}
 \epr

\bex Полученной теперь информации достаточно для построения графиков синуса и
косинуса. Эти картинки общеизвестны, однако мы предлагаем читателю проверить
самостоятельно, что, если ставить себе целью указать на графике интервалы
монотонности и знакоопределенности, то из того, что мы уже успели доказать про
синус и косинус никаких других картинок получиться не может. Для синуса график
должен выглядеть так:

%\picture{120pt}{0pt}{sin_x.pcx}
\vglue120pt

А для косинуса так:

%\picture{120pt}{0pt}{cos_x.pcx}
\vglue120pt

\eex

\paragraph{Тангенс и котангенс.}

 \bit{
\item[$\bullet$]Тангенс $\tg$ и котангенс $\ctg$ определяются формулами
 \begin{align}
\tg x&=\frac{\sin x}{\cos x} && \left(x\notin\frac{\pi}{2}+\pi\Z\right)
\\
\ctg x&=\frac{\cos x}{\sin x} && (x\notin\pi\Z)
 \end{align}
 }\eit

\bprop Справедливы следующие тождества:
 \begin{align}
&\tg(x+\pi)=\tg x, && \ctg(x+\pi)=\ctg x \label{tg(x+pi)=tg(x)}\\
&\tg\left(x-\frac{\pi}{2}\right)=-\ctg x, && \ctg\left(x-\frac{\pi}{2}\right)=-\tg x \label{tg(x+pi/2)=tg(x)}\\
& \tg(-x)=-\tg x, &&  \ctg(-x)=-\ctg x \label{tg(-x)=-tg(x)}
 \end{align}
 \eprop
 \bpr
Это сразу вытекает из свойств синуса и косинуса, например, первая формула
доказывается так:
 \begin{multline*}
\tg(x+\pi)=\frac{\sin(x+\pi)}{\cos(x+\pi)}=\eqref{sin(x+pi)},\eqref{cos(x+pi)}=\\=
\frac{-\sin x}{-\cos x}= \frac{\sin x}{\cos x}=\tg x
 \end{multline*}
 \epr

 \bprop На интервале $\left(-\frac{\pi}{2},\frac{\pi}{2}\right)$ функция
 $\tg$ возрастает:
 \begin{align}
x,y\in \left(-\frac{\pi}{2},\frac{\pi}{2}\right), && x<y && \Longrightarrow &&
\tg x<\tg y \label{tg-vozrast}
 \end{align}
 \eprop
 \bpr
В силу нечетности тангенса \eqref{tg(-x)=-tg(x)}, достаточно рассмотреть
интервал $\left(0,\frac{\pi}{2}\right)$. Тогда все следует из монотонности
синуса и косинуса:
$$
x,y\in \left(0,\frac{\pi}{2}\right), \qquad x<y
$$
$$
\Downarrow
$$
$$
0<\sin x<\sin y,\qquad \cos x>\cos y>0
$$
$$
\Downarrow
$$
$$
0<\sin x<\sin y,\qquad 0<\frac{1}{\cos x}<\frac{1}{\cos y}
$$
$$
\Downarrow
$$
$$
\tg x=\frac{\sin x}{\cos x}<\frac{\sin y}{\cos y}=\tg y
$$
 \epr

 \bprop На интервале $(0,\pi)$ функция $\ctg$ убывает:
 \begin{align}
x,y\in (0,\pi), && x<y && \Longrightarrow && \ctg x>\ctg y \label{ctg-ubyv}
 \end{align}
 \eprop
 \bpr
$$
x,y\in (0,\pi), \qquad x<y
$$
$$
\Downarrow
$$
$$
x-\frac{\pi}{2},y-\frac{\pi}{2}\in \left(-\frac{\pi}{2},\frac{\pi}{2}\right),
\qquad x-\frac{\pi}{2}<y-\frac{\pi}{2}
$$
$$
\Downarrow
$$
$$
\tg\left(x-\frac{\pi}{2}\right)<\tg\left(y-\frac{\pi}{2}\right)
$$
$$
\Downarrow
$$
$$
-\ctg x=\tg\left(x-\frac{\pi}{2}\right)<\tg\left(y-\frac{\pi}{2}\right)=-\ctg y
$$
$$
\Downarrow
$$
$$
\ctg x>\ctg y
$$
 \epr

\bprop На границах интервалов, где $\tg$ и $\ctg$ определены, эти функции имеют
бесконечные пределы:
 \begin{align}
&\lim_{x\to\frac{\pi}{2}+\pi n}\tg x=\infty, && (n\in\Z) \\
& \lim_{x\to\pi n}\ctg x=\infty, && (n\in\Z)
 \end{align}
\eprop
 \bpr
Из-за $\pi$-периодичности и формул \eqref{tg(x+pi/2)=tg(x)}, связывающих $\tg$
и $\ctg$, здесь достаточно рассмотреть только одну точку. Например, при
$x\to\frac{\pi}{2}$ мы, в силу непрерывности $\sin$ и $\cos$, получаем:
$$
\sin x\underset{x\to\frac{\pi}{2}}{\longrightarrow} \sin\frac{\pi}{2}= 1,\qquad
\cos x\underset{x\to\frac{\pi}{2}}{\longrightarrow} \cos\frac{\pi}{2}=0
$$
$$
\Downarrow
$$
$$
\tg x=\frac{\sin x}{\cos x}\underset{x\to\frac{\pi}{2}}{\longrightarrow} \infty
$$
 \epr

\bex Для построения графиков нам теперь не хватает нескольких значений в
характерных точках. Из таблицы для синуса и косинуса на
с.\pageref{tablitsa-sin-cos} мы прямым вычислением получаем таблицу для
тангенса и котангенса:

\bprop Справедлива следующая таблица значений тангенса и котангенса:

 \bigskip
 \vbox{\tabskip=0pt\offinterlineskip \halign to \hsize{ \vrule#\tabskip=2pt
 plus3pt minus1pt & \strut\hfil\;#\hfil & \vrule#&\hfil#\hfil &
 \vrule#&\hfil#\hfil & \vrule#&\hfil#\hfil & \vrule#&\hfil#\hfil &
 \vrule#&\hfil#\hfil & \vrule#\tabskip=0pt\cr \noalign{\hrule}
 height2pt&\omit&&\omit&&\omit&&\omit&&\omit&&\omit&\cr

 & $x$ && $0$ && $\frac{\pi}{6}$ && $\frac{\pi}{4}$ & & $\frac{\pi}{3}$ && $\frac{\pi}{2}$  &\cr \noalign{\hrule}
 height2pt&\omit&&\omit&&\omit&&\omit&&\omit&&\omit&\cr

 & $\tg x$ && $0$ && $\frac{1}{\sqrt{3}}$ && $1$ && $\sqrt{3}$ & & $\not\exists$ &\cr \noalign{\hrule}
 height2pt&\omit&&\omit&&\omit&&\omit&&\omit&&\omit&\cr

 & $\ctg x$ && $\not\exists$ && $\sqrt{3}$ && $1$ && $\frac{1}{\sqrt{3}}$ && $0$ &\cr \noalign{\hrule} } }
 \label{tablitsa-tg-ctg}
 \bigskip
 \eprop

Теперь можно строить графики. Для тангенса график выглядит так:

%\picture{120pt}{0pt}{tg_x.pcx}
\vglue120pt

А для котангенса так:

%\picture{120pt}{0pt}{ctg_x.pcx}
\vglue120pt
 \eex

\paragraph{Обратные тригонометрические функции.}

 \biter{
\item[$\bullet$] Следующие правила определяют функции, называемые
соответственно арксинусом, арккосинусом, арктангенсом и арккотангенсом:
 }\eiter
 \begin{align} t &=\arcsin
x\;\Longleftrightarrow\; x=\sin t\;\&\; t\in
\left[-\frac{\pi}{2}; \frac{\pi}{2}\right] \label{DEF:arcsin}\\
t &=\arccos x\;\Longleftrightarrow\; x=\cos t\;\&\; t\in
\left[0; \pi\right] \label{DEF:arccos}\\
t &=\arctg x\;\Longleftrightarrow\; x=\tg t\;\&\; t\in
\left(-\frac{\pi}{2}; \frac{\pi}{2}\right) \label{DEF:arctg} \\
t &=\arcctg x\;\Longleftrightarrow\; x=\ctg t\;\&\; t\in \left(0; \pi\right)
\label{DEF:arcctg}
 \end{align}
\bpr Здесь используется теоремы о монотонной функции на отрезке
\ref{Teor-ob-obratnoi-funktsii-na-R^1} и на интервале
\ref{Teor-ob-obratnoi-funktsii-na-int-v-R^1}. Например, в случае с синусом мы
получаем вот что: по предложению \ref{PROP:int-monot-sin}, синус строго
возрастает на отрезке $\left[-\frac{\pi}{2}; \frac{\pi}{2}\right]$. Вдобавок
$$
\sin\left(-\frac{\pi}{2}\right)=-1,\qquad \sin\left(\frac{\pi}{2}\right)=1
$$
поэтому по теореме \ref{Teor-ob-obratnoi-funktsii-na-R^1} правило
\eqref{DEF:arcsin} определяет некую функцию на отрезке
$[-1;1]=\sin\left(\left[-\frac{\pi}{2}; \frac{\pi}{2}\right]\right)$
 \epr

Из теорем \ref{Teor-ob-obratnoi-funktsii-na-R^1} и
\ref{Teor-ob-obratnoi-funktsii-na-int-v-R^1} мы автоматически получаем

\bprop\label{PROP:nepr-arcsin-arccos} Функции $\arcsin$, $\arccos$, $\arctg$ и
$\arcctg$ непрерывны и монотонны на своей области определения.
 \eprop

А из предложений \ref{PROP:sin>0} и \ref{PROP:cos>0} получаем

 \bprop Функции $\arcsin$, $\arccos$, $\arctg$ и
$\arcctg$ имеют следующие знаки на характерных интервалах:
 \begin{align}
& \arcsin x>0, && x\in(0,1] \label{arcsin(x)>0}\\
& \arcsin x<0, && x\in[-1,0)\label{arcsin(x)<0} \\
& \arccos x>0, && x\in[-1,1) \label{arcsoc(x)>0} \\
& \arctg x>0, && x>0 \label{arctg(x)>0} \\
& \arctg x<0, && x<0 \label{arctg(x)<0} \\
& \arcctg x>0, && x\in\R \label{arcctg(x)>0}
  \end{align}
 \eprop

 Остается привести графики этих функций.

\bex Арксинус:

%\picture{120pt}{0pt}{arcsin_x.pcx}
\vglue120pt

\eex

\bex
 Арккосинус:

%\picture{120pt}{0pt}{arccos_x.pcx}
\vglue120pt \eex

\bex

Арктангенс:

%\picture{120pt}{0pt}{arctg_x.pcx}
\vglue120pt \eex

\bex Арккотангенс:

%\picture{120pt}{0pt}{arcctg_x.pcx}
\vglue120pt \eex

\section{Числовые выражения и стандартные функции}
\label{Chislovye-termy-i-stand-func}

\subsection{Собственные обозначения для функций}

Как, наверное, уже заметил читатель, подобно числам (с.
\pageref{SUBSEC-sobstv-obozn-chisla}), некоторые функции также имеют
собственные обозначения. Скажем, символ $\sin$ используется только для
обозначения тригонометрической функции <<синус>>, и не бывает такого, чтобы
математик обозначал им какую-нибудь другую функцию на $\R$ (как не бывает,
чтобы символом 1 обозначалось число, отличное от единицы).

В этом пункте мы перечислим некоторые собственные обозначения для функций.

\bexs Понятно, что символы тригонометрических функций будут собственными
обозначениями:
 \begin{align*}
f(x)&=\sin x \\
f(x)&=\cos x \\
f(x)&=\tg x \\
f(x)&=\ctg x
 \end{align*}
Точно так же, обратные тригонометрические функции имеют собственные обозначения:
 \begin{align*}
f(x)&=\arcsin x \\
f(x)&=\arccos x \\
f(x)&=\arctg x \\
f(x)&=\arcctg x
 \end{align*}
Еще один пример собственного обозначения для функции -- логарифм, например, с
основанием 2:
$$
f(x)=\log_2 x
$$

Очень важно при этом, что собственные обозначения для функций позволяют
определять функции без участия переменных. Например, запись
 \beq\label{korotkaya-zapis'-sin}
f=\sin
 \eeq
означает, что под $f$ понимается функция <<синус>>, которая в точке $x\in\R$
будет равна
 \beq\label{dlinnaya-zapis'-sin}
\forall x\in\R\quad f(x)=\sin x
 \eeq
Таким образом, формулу \eqref{korotkaya-zapis'-sin} можно считать просто
короткой записью формулы \eqref{dlinnaya-zapis'-sin}:
$$
f=\sin\;\Longleftrightarrow\; \Big(\forall x\in\R\quad f(x)=\sin x\Big)
$$
Точно так же, запись
 $$
f=\cos
 $$
означает, что под $f$ понимается функция <<косинус>>, и расшифровывается она
так:
 $$
\forall x\in\R\quad f(x)=\cos x
 $$

Из уже имеющихся собственных обозначений можно с помощью алгебраических
операций строить новые. Например, функция
$$
\sin+\cos
$$
действует по правилу
$$
(\sin+\cos)(x)=\sin x+\cos x\qquad (x\in\R)
$$
Точно так же функции
$$
\sin^2,\quad \cos\cdot\log_2
$$
действуют по правилам
$$
(\sin^2)(x)=(\sin x)^2\qquad (x\in\R)
$$
$$
(\cos\cdot\log_2)(x)=\cos x\cdot\log_2 x\qquad (x\in\R)
$$

\eexs

\paragraph{Числовые выражения.}

Числовые выражения служат главным предметом изучения в школе, поэтому читателю
они должны быть хорошо знакомы. Грубо говоря, это конечные последовательности
символов (букв, цифр, или символов операций), которые можно записывать как
левые или правые части уравнений, или которые определяют функции (возможно, от
нескольких переменных) с помощью элементарных функций и алгебраических
операций. Например, выражениями будут
$$
25,\quad t+2,\quad \sin x,\quad \frac{a+b}{a-b},\quad \sqrt{x^2+y^2}.
$$
Строгое определение этому понятию выглядит так:

\biter{

\item[$\bullet$] {\it Числовыми
выражениями}\index{выражение!числовое}\label{opredelenie-chislovogo-terma}
называются конечные последовательности символов, определенные следующими
индуктивными правилами:

 \biter{

\item[1)] всякое собственное обозначение для числа является числовым
выражением, например
$$
1,\; 2,\; 15,\; 7/11
$$
в частности, символы $e$ и $\pi$, употреблявшиеся выше для обозначения чисел, определенных формулами \eqref{5.9.1} и \eqref{DEF:pi}, также считаются выражениями;

\item[2)] всякая буква латинского или греческого алфавита, строчная или прописная, считается числовым выражением:
$$
a,b,c,...,x,y,z,\qquad A,B,C,...,X,Y,Z
$$
$$
\alpha,\beta,\gamma,...,\chi, \psi,\omega,\qquad
A,B,\varGamma,...,X,\varPsi,\varOmega
$$
 (в этом списке не должно быть только букв $e$ и $\pi$, зарезервированных как собственные обозначения двух конкретных чисел, и поэтому уже упомянутых в предыдущем пункте);

\item[3)] если $\mathcal P$ и $\mathcal Q$ --- числовые выражения, то числовыми
выражениями будут также
\begin{multline*}
(\mathcal P)+(\mathcal Q),\quad (\mathcal P)-(\mathcal Q),\quad (\mathcal
P)\cdot(\mathcal Q),\\ ({\mathcal P})/({\mathcal Q}),\quad (\mathcal
P)^{(\mathcal Q)},\quad \log_{(\mathcal P)}({\mathcal Q})
\end{multline*}

\item[4)] если $\mathcal P$ --- числовое выражение, то числовыми выражениями
будут также выражения, получаемые подстановкой $\mathcal P$ как аргумента под
знак элементарных функций, в частности,
\begin{multline*}
\sin (\mathcal P), \quad \cos (\mathcal P), \\
 \tg (\mathcal P), \quad \ctg
(\mathcal P), \\
\arcsin (\mathcal P), \quad \arccos (\mathcal P), \\
 \arctg (\mathcal P),
\quad \arcctg (\mathcal P)
\end{multline*}

\item[5)] Если $\mathcal Q$ --- числовое выражение, $y$ --- входящая в него
буква, и $\mathcal P$ --- какое-то другое числовое выражение, то, подставив
всюду в $\mathcal Q$ вместо буквы $y$ выражение $\mathcal P$, мы получим новое
числовое выражение, обозначаемое
 \beq\label{podstanovka-v-vyrazhenie}
{\mathcal Q}\Big|_{y={\mathcal P}},
 \eeq
и называемое {\it результатом подстановки} в $\mathcal Q$ вместо буквы $y$
выражения $\mathcal P$.
 }\eiter
 }\eiter

 \biter{
\item[$\bullet$] Если выражение $\mathcal R$ получено из какого-то выражения
$\mathcal Q$ подстановкой выражения $\mathcal P$,
$$
{\mathcal Q}\Big|_{y={\mathcal P}}={\mathcal R}
$$
то говорят, что выражение $\mathcal P$ {\it подчинено} выражению $\mathcal R$
(или {\it является частью} $\mathcal R$).

\item[$\bullet$] {\it Порядок выражения} $\mathcal R$ --- это сколько раз в нем
встречаются алгебраические операции (сложение, вычитание, умножение, деление,
возведение в степень) и символы элементарных функций.

}\eiter

\paragraph{Переменные и параметры.}

Буквы, входящие в числовое выражение, делятся на {\it переменные} и {\it
параметры}. Это значит, что, задавая выражение, Вы обязаны указать, какие из
входящих в него букв следует считать переменными, а какие параметрами.
Например, в выражении
$$
x^\alpha
$$
можно объявить букву $x$ переменной, а $\alpha$ --- параметром, и формально это
будет не то же самое выражение, как если бы мы заявили, что наоборот, $x$ ---
параметр, а $\alpha$ --- переменная.

Разница между переменными и параметрами станет понятной только на с.
\pageref{SEC-proizv-kak-oper-nad-simv}, где будет введено понятие производной
выражения
--- там мы объясним, что, в отличие от параметра, переменная обладает тем
свойством, что по ней можно дифференцировать, причем производная любого
параметра по любой переменной по определению считается равной нулю: если
$\alpha$ --- параметр, а $x$ --- переменная, то
$$
\frac{\d \alpha}{\d x}=0
$$
Позже, на с.\pageref{diff-vyrazh-step-1} это формальное различие приобретет
вполне простую и аксиоматически удобную форму с введением понятия
дифференциального выражения
--- там оно выразится в постулировании правил вычисления
дифференциала:
 \begin{align*}
&\d x\ne 0, \qquad \text{если $x$ --- переменная} \\
&\d x=0, \qquad \text{если $x$ --- параметр}
 \end{align*}

До той поры мы будем считать, что разница между переменной и параметром чисто
номинальна: если, например, дан выражение
$$
x^\alpha,
$$
и указано, что буква $x$ в нем является переменной, а $\alpha$ --- параметром,
то мы просто должны это помнить. Интуитивно же такую запись полезно
представлять себе не как один выражение, а как целое <<семейство выражений,
зависящее от параметра $\alpha$>>: при разных значениях $\alpha$ мы будем
получать разные выражения:
$$
x^1,\qquad x^{-2},\qquad x^\frac{3}{5},\qquad x^{\sqrt{2}},\qquad ...
$$
Если параметров два, например,
$$
k\cdot x+b,\qquad (\text{$x$ --- переменная, $k$, $b$ --- параметры})
$$
то запись можно считать <<двухпараметрическим семейством выражений>>: при
разных значениях $k$ и $b$ получаем выражения:
$$
2\cdot x+1,\qquad -3\cdot x+0,\qquad 0\cdot x+1,\qquad x+(-{\sqrt{2}}),\qquad
...
$$
И так далее.

\bit{

\item[$\bullet$] Если число переменных в числовом выражении $\mathcal P$ не
превышает $n$, то (независимо от того, сколько в нем параметров) такой
выражение называется {\it $n$-местным}. В частности, если $\mathcal P$ содержит
только одну переменную, или не содержит никакой переменной, то $\mathcal P$
называется {\it одноместным выражением}.

\item[$\bullet$] Одноместное числовое выражение $\mathcal P$ называется {\it
выражением от переменной} $x$, если он содержит переменную $x$ (либо не
содержит никакой).
 }\eit

\subsection{Область допустимых значений переменных. Стандартная
функция. Решение уравнений и неравенств.}

С числовыми выражениями связаны три важных понятия, о которых приходится
говорить одновременно из-за невозможности их определить отдельно друг от друга
-- это область допустимых значений переменной, стандартная функция,
определяемая выражением и множество решений уравнений (а также неравенств, или,
в общем случае, условий), определяемых выражением. Точное определение эти
понятий (``трехсложное определение'' на с.\pageref{3-opredelenie}) выглядит
довольно громоздко, поэтому мы предпошлем ему интуитивное объяснение, чтобы к
моменту, когда дается формальная конструкция читатель психологически был к ней
готов.

\paragraph{Область допустимых значений.}
Прежде всего, всякое числовое выражение $\mathcal P$ от одной переменной,
скажем $x$, и без параметров, определяет на прямой $\R$ так называемую область
допустимых значений переменной $x$, неформально понимаемую как множество
значений $x$, для которых <<выражение $\mathcal P$ имеет смысл>>.

\bex Для выражения
$$
\frac{1}{x}
$$
областью допустимых значений переменной будет множество
$(-\infty,0)\cup(0,+\infty)$, то есть множество $\{x\in\R:\; x\ne0\}$, потому
что при подстановке $x=0$ это выражение <<утрачивает смысл>>
$$
\frac{1}{0},
$$
(делить на ноль нельзя), а при любых других подстановках получающемуся
выражению можно приписать смысл. Например, при $x=2$
$$
\frac{1}{2}
$$
мы получаем выражение, смысл которого --- имя числа. \eex

\bex Для выражения
$$
\frac{1}{\log_2 x}
$$
область допустимых значений переменной определяется условиями
$$
x>0
$$
(для того, чтобы существовал $\log_2 x$) и
$$
\log_2 x\ne 0
$$
(для того, чтобы можно было делить на $\log_2 x$). Поскольку эти условия должны
быть выполнены оба сразу, их нужно объединить в систему, решая которую получим
ответ:
 \begin{multline*}
 \left\{\begin{matrix} x>0 \\ \log_2 x\ne 0\end{matrix}\right\}
\quad\Longleftrightarrow\quad
 \left\{\begin{matrix} x>0 \\ 0< x\ne 1 \end{matrix}\right\}
\quad\Longleftrightarrow\\ \Longleftrightarrow\quad x\in (0,1)\cup(1,+\infty)
 \end{multline*}
Множество $(0,1)\cup(1,+\infty)$ будет областью допустимых значений переменной
в этом выражении.
 \eex

В более общем случае, когда выражение $\mathcal P$ содержит параметры (и
по-прежнему, зависит от одной переменной $x$), его область допустимых значений
представляет собой не одно множество, а семейство множеств на прямой $\R$,
зависящее от параметров.

\bex В соответствии со сказанным на с.\pageref{usloviya-sushestv-x^a}, областью
допустимых значений переменной $x$ в выражении
$$
x^\alpha
$$
будет семейство множеств
$$
D_\alpha=\begin{cases}
\R,& \alpha\in\frac{\N}{2\N-1};\\
\R\setminus\{0\},& \alpha\in -\frac{\Z_+}{2\N-1};\\
[0,+\infty),& \alpha\in(0,+\infty)\setminus\frac{\N}{2\N-1}\\
\varnothing,& \alpha\in (-\infty,0)\setminus\frac{-\N}{2\N-1}\end{cases}
$$
\eex

\bex Аналогично, областью допустимых значений переменной $x$ в выражении
$$
a^x
$$
будет семейство множеств
$$
D_a=\begin{cases}
\R,& a>0\\
(0;+\infty),& a=0 \\
\frac{\Z}{2\N-1},& a<0
\end{cases}
$$
\eex

\bex Для выражения
$$
\log_a x
$$
область допустимых значений переменной $x$ имеет вид
$$
D_a=\begin{cases}
(0;+\infty),& a\in(0;1)\cup(1;+\infty);\\
\varnothing,& \alpha\in (-\infty,0]\cup\{1\}\end{cases}
$$

\eex

\paragraph{Стандартная функция, определяемая выражением.}\label{standartnye-funktsii}

Второе, на что мы хотим обратить внимание читателя, это что всякое числовое
выражение $\mathcal P$ без параметров от переменной $x$ определяет некоторую
функцию формулой
 \beq\label{functsiya-opred-termom}
f(x)={\mathcal P}
 \eeq
Такие функции (то есть, определенные выражением) мы будем называть {\it
стандартными}.

\bex Выражение
$$
\frac{1}{\log_2 x}
$$
определяет функцию одного переменного
$$
f(x)=\frac{1}{\log_2 x}
$$
а выражение
$$
\frac{1}{x+y}
$$
--- функцию двух переменных
$$
f(x,y)=\frac{1}{x+y}
$$
Областью определения такой функции считается область допустимых значений
переменных выражения. В первом случае это будет множество на прямой
$$
\{x:\; 0<x\ne 1\}=(0,1)\cup(1,+\infty)
$$
а во втором --- множество на плоскости
$$
\{(x,y):\; x+y\ne 0\}.
$$
\eex

\bex С другой стороны, выражения
$$
\frac{1}{0},\qquad  \sqrt{x-1}+\sqrt{-1-x}
$$
определяют пустую функцию, потому что область допустимых значений переменных у
них пуста. \eex

Если же выражение $\mathcal P$ содержит параметры (и, по-прежнему, только одну
переменную $x$), то формула \eqref{functsiya-opred-termom} определяет некое
семейство функций, зависящее от параметров.

\bex Формула
$$
f(x)=x^\alpha
$$
определяет семейство функций с областями определения, зависящими от параметра
$\alpha$:
$$
D_\alpha=\begin{cases}
\R,& \alpha\in\frac{\N}{2\N-1};\\
\R\setminus\{0\},& \alpha\in -\frac{\Z_+}{2\N-1};\\
[0,+\infty),& \alpha\in(0,+\infty)\setminus\frac{\N}{2\N-1}\\
\varnothing,& \alpha\in (-\infty,0)\setminus\frac{-\N}{2\N-1}\end{cases}
$$
\eex

\bex Формула
$$
f(x)=a^x
$$
определяет семейство функций с областями определения
$$
D_a=\begin{cases}
\R,& a>0\\
(0;+\infty),& a=0 \\
\frac{\Z}{2\N-1},& a<0
\end{cases}
$$
\eex

\bex Формула
$$
f(x)=\log_a x
$$
определяет семейство функций на множествах
$$
D_a=\begin{cases}
(0;+\infty),& a\in(0;1)\cup(1;+\infty);\\
\varnothing,& \alpha\in (-\infty,0]\cup\{1\}\end{cases}
$$
\eex

\bex Функция {\it модуль} является стандартной функцией, потому что представима
в виде
$$
|x|=\sqrt{x^2}
$$
\eex

\bex Функция {\it сигнум}, определяемая правилом
$$
\sgn x=\begin{cases} 1,& x>0 \\ -1,& x<0\end{cases}
$$
является стандартной, потому что совпадает со стандартной функцией
$$
\sgn x=\frac{x}{|x|}
$$
 \eex

\bex Покажем, что для всякого открытого интервала $(\alpha,\beta)$ (возможно,
бесконечного), функция
 $$
f(x)=\begin{cases} 1,& x\in (\alpha,\beta) \\ 0,& x\notin
[\alpha,\beta]\end{cases}
 $$
будет стандартной. Для бесконечных интервалов это верно в силу формул
 \beq
\left\{\begin{matrix} 0,& x<\alpha \\ 1,&
x>\alpha\end{matrix}\right\}=\frac{1+\sgn (x-\alpha)}{2}
 \eeq
и
 \beq
\left\{\begin{matrix} 1,& x<\alpha \\ 0,&
x>\alpha\end{matrix}\right\}=\frac{1-\sgn (x-\alpha)}{2}
 \eeq
В частности, для случая положительной и отрицательной полупрямой формулы
выглядят так:
 \beq
\left\{\begin{matrix} 0,& x<0 \\ 1,& x>0\end{matrix}\right\}=\frac{1+\sgn x}{2}
 \eeq
 \beq
\left\{\begin{matrix} 1,& x<0 \\ 0,& x>0\end{matrix}\right\}=\frac{1-\sgn x}{2}
 \eeq
А для конечных интервалов -- в силу формулы
 \begin{multline}
\left\{\begin{matrix} 1,& x\in(\alpha,\beta) \\ 0,&
x\notin[\alpha,\beta]\end{matrix}\right\}=\\=
\frac{1+\sgn(x-\alpha)}{2}-\frac{1+\sgn(x-\beta)}{2}
 \end{multline}

 \eex

\bex Всякая функция вида
 $$
f(x)=1,\qquad  x\in (\alpha,\beta),
 $$
(то есть определенная на заданном интервале $(\alpha,\beta)$ и равная на нем
единице) является стандартной, потому что получается из какой-то функции
предыдущего примера возведением в степень $-1$. Для бесконечных интервалов
формулы выглядят так:
 \beq
\left\{\begin{matrix} 1,& x>\alpha\end{matrix}\right\}=\frac{2}{1+\sgn
(x-\alpha)},
 \eeq
 \beq
\left\{\begin{matrix} 1,& x<\alpha \end{matrix}\right\}=\frac{2}{1-\sgn
(x-\alpha)}
 \eeq
в частности, для случая положительной и отрицательной полупрямой,
 \beq
\left\{\begin{matrix} 1,& x>0\end{matrix}\right\}=\frac{2}{1+\sgn x}
 \eeq
 \beq
\left\{\begin{matrix} 1,& x<0\end{matrix}\right\}=\frac{2}{1-\sgn x}
 \eeq
А, для конечных интервалов --- так:
 \begin{multline}\label{stand-func-1-(a,b)}
\left\{\begin{matrix} 1,& x\in(\alpha,\beta) \end{matrix}\right\}=\\=
\left(\frac{1+\sgn(x-\alpha)}{2}-\frac{1+\sgn(x-\beta)}{2}\right)^{-1}
 \end{multline}
\eex

\bex Из предыдущего примера следует, что если нам дана некая стандартная
функция $f$, то умножая ее на стандартную функцию вида
\eqref{stand-func-1-(a,b)}, где интервал $(\alpha,\beta)$ лежит в области
определения $f$, мы получаем стандартную функцию, определенную строго на
$(\alpha,\beta)$, и совпадающую с $f$ на этом интервале:
$$
\left\{\begin{matrix} f(x),& x\in(\alpha,\beta) \end{matrix}\right\}=f(x)\cdot
\left\{\begin{matrix} 1,& x\in(\alpha,\beta) \end{matrix}\right\}
$$
Таким образом, {\it ограничение стандартной функции на любой интервал
(возможно, бесконечный) также будет стандартной функцией}.

В частности, функция
$$
\left\{\begin{matrix} \frac{1}{x},& x>0\end{matrix}\right\}=\frac{1}{x}\cdot
\left\{\begin{matrix} 1,& x>0
\end{matrix}\right\}
$$
тоже является стандартной.
 \eex

\ber Проверьте, что для любых стандартных функций $f_1$,...,$f_k$ и любого
набора непересекающихся интервалов
$(\alpha_1,\beta_1)$,...,$(\alpha_k,\beta_k)$, функция, определенная правилом
$$
f(x)=\left\{\begin{matrix} f_1(x),& x\in(\alpha_1,\beta_1)\\...\\ f_k(x),&
x\in(\alpha_k,\beta_k)\end{matrix}\right\}
$$
также будет стандартной.
 \eer

\paragraph{Уравнения, неравенства и условия, определяемые выражением.}
Третье наблюдение состоит в том, что, имея дело с выражениями $\mathcal P$ от
переменной $x$, постоянно приходится решать всякие уравнения вида
$$
{\mathcal P}=C,
$$
неравенства
$$
{\mathcal P}<C,\qquad {\mathcal P}\le C,\qquad {\mathcal P}> C,\qquad {\mathcal
P}\ge C
$$
и часто даже более общие задачи, в которых дано какое-то множество $E$ на
прямой $\R$ и требуется найти такие значения переменной $x$, для которых
выполняется условие
$$
{\mathcal P}\in E
$$

В школе учат решать уравнения и неравенства разных типов, и для этого
оказывается полезным употребление понятий системы и совокупности. Напомним, что
 \bit{
\item[$\bullet$] {\it системой условий} $\mathcal A$ и $\mathcal B$ на
переменную $x$ называется условие, записываемое как $\mathcal A$ и $\mathcal B$
в фигурных скобках
$$
\left\{\begin{matrix} {\mathcal A} \\ {\mathcal B}
\end{matrix}\right\},
$$
и представляющее собой логическое пересечение $\mathcal A\&\mathcal B$ условий
$\mathcal A$ и $\mathcal B$;

\item[$\bullet$] {\it совокупностью условий} $\mathcal A$ и $\mathcal B$ на
переменную $x$ называется условие, записываемое как $\mathcal A$ и $\mathcal B$
в квадратных скобках
$$
\left[\begin{matrix} {\mathcal A} \\ {\mathcal B}
\end{matrix}\right]
$$
и представляющее собой логическое объединение $\mathcal A\vee\mathcal B$
условий $\mathcal A$ и $\mathcal B$.
 }\eit
Как пользоваться этими понятиями мы напомним следующих примерах.

\bex {\bf Уравнения} --- формулы вида
$$
{\mathcal P}={\mathcal Q}
$$
где $\mathcal P$ и $\mathcal Q$ --- выражения. Частным случаем будет уравнение
вида
$$
{\mathcal P}=C
$$
где $\mathcal P$ --- выражение, а $C$ --- параметр не содержащийся в $\mathcal
P$, или собственное обозначение для числа.

Решение уравнения представляет собой последовательное сведение его к
равносильным системам и совокупностям (обычно многократно вложенным друг в
друга) простейших условий на переменную и параметры. Вот некоторые примеры
общих правил решения уравнений:
$$
ax+b=0\quad\Longleftrightarrow\quad \left[\begin{matrix} \left\{\begin{matrix}
a\ne 0 \\ x=-\frac{b}{a}
\end{matrix}\right\}
\\
\left\{\begin{matrix} a=0 \\ \left[\begin{matrix} \left\{\begin{matrix} b\ne 0
\\ x\in\varnothing\end{matrix}\right\}
\\
\left\{\begin{matrix} b=0 \\ x\in \R \end{matrix}\right\}
\end{matrix}\right]
\end{matrix}\right\}
\end{matrix}\right]
$$
 \begin{multline*}
ax^2+bx+c=0
 \quad\Longleftrightarrow \\
 \Longleftrightarrow\quad
 \left[\begin{matrix}
  \left\{\begin{matrix} a=0 \\ bx+c=0
  \end{matrix}\right\}
\\
 \left\{\begin{matrix} a\ne 0 \\
  \left[\begin{matrix}
   \left\{\begin{matrix} b^2-4ac<0 \\ x\in\varnothing
   \end{matrix}\right\}
   \\
   \left\{\begin{matrix} b^2-4ac=0 \\ x=-\frac{b}{2a}
   \end{matrix}\right\}
   \\
   \left\{\begin{matrix} b^2-4ac>0 \\ x=\frac{-b\pm\sqrt{b^2-4ac}}{2a}
   \end{matrix}\right\}
  \end{matrix}\right]
 \end{matrix}\right\}
\end{matrix}\right]
 \end{multline*}
$$
\log_a x=C\quad\Longleftrightarrow\quad \left\{\begin{matrix} 0<a\ne 1 \\ x=a^C
\end{matrix}\right\}
$$
 \begin{multline*}
\sin x=C
 \quad\Longleftrightarrow \\
 \Longleftrightarrow\quad
 \left[\begin{matrix}
  \left\{\begin{matrix}
   C\in[-1;1]  \\
   \left[\begin{matrix}
   x\in\arcsin C+2\pi\Z \\
   x\in\pi-\arcsin C+2\pi\Z
   \end{matrix}\right]
  \end{matrix}\right\}
 \\
  \left\{\begin{matrix}
  C\notin[-1;1] \\ x\in\varnothing
  \end{matrix}\right\}
 \end{matrix}\right]
\end{multline*}

 \begin{multline*}
\cos x=C
 \quad\Longleftrightarrow \\
 \Longleftrightarrow\quad
 \left[\begin{matrix}
  \left\{\begin{matrix}
   C\in[-1;1]  \\
   x\in\pm\arccos C+2\pi\Z
  \end{matrix}\right\}
 \\
  \left\{\begin{matrix}
  C\notin[-1;1] \\ x\in\varnothing
  \end{matrix}\right\}
 \end{matrix}\right]
\end{multline*}
\eex

$$
\arcsin x=C
 \quad\Longleftrightarrow\quad
 \left[\begin{matrix}
  \left\{\begin{matrix}
   C\in\left[-\frac{\pi}{2};\frac{\pi}{2}\right]  \\
   x\in\sin C
  \end{matrix}\right\}
 \\
  \left\{\begin{matrix}
  C\notin\left[-\frac{\pi}{2};\frac{\pi}{2}\right]
  \\ x\in\varnothing
  \end{matrix}\right\}
 \end{matrix}\right]
$$

$$
\arccos x=C
 \quad\Longleftrightarrow\quad
 \left[\begin{matrix}
  \left\{\begin{matrix}
   C\in[\,0;\pi]  \\
   x\in\cos C
  \end{matrix}\right\}
 \\
  \left\{\begin{matrix}
  C\notin[\,0;\pi]
  \\ x\in\varnothing
  \end{matrix}\right\}
 \end{matrix}\right]
$$

\ber Найдите общие правила решения уравнений:
$$
x^\alpha=C,\quad a^x=C,
$$
$$
\tg x=C,\quad \ctg x=C,
$$
$$
\arctg x=C,\quad \arcctg x=C.
$$
\eer

\bex {\bf Неравенства} --- формулы вида
$$
{\mathcal P}<{\mathcal Q},\quad {\mathcal P}\le {\mathcal Q},\quad {\mathcal
P}>{\mathcal Q},\quad {\mathcal P}\ge {\mathcal Q},
$$
где $\mathcal P$ и $\mathcal Q$ --- выражения. Частными случаями будут
неравенства вида
$$
{\mathcal P}<C,\quad {\mathcal P}\le C,\quad {\mathcal P}>C,\quad {\mathcal
P}\ge C,
$$
где $\mathcal P$ --- выражение, а $C$ --- параметр не содержащийся в $\mathcal
P$, или собственное обозначение для числа.

В школе неравенства обычно учат решать двумя способами: либо тоже сведением их
к равносильным системам и совокупностям простых условий, либо методом
интервалов. Второй способ мы разбирать не будем, полагаясь на память читателя,
а первый проиллюстрируем лишь несколькими правилами:
$$
\log_a x<C \quad\Longleftrightarrow\quad
 \left[\begin{matrix}
  \left\{\begin{matrix}  0<a<1 \\ x>a^C
  \end{matrix}\right\}
  \\
  \left\{\begin{matrix}  a>1 \\ x<a^C
  \end{matrix}\right\}
 \end{matrix}\right]
$$

 \begin{multline*}
\sin x<C
 \quad\Longleftrightarrow \\
 \Longleftrightarrow\quad
 \left[\begin{matrix}
  \left\{\begin{matrix}
   C<-1  \\ x\in\varnothing
  \end{matrix}\right\}
 \\
  \left\{\begin{matrix}
   C\in[-1;1]  \\
   x\in[-\arcsin С-\pi;\arcsin C]+2\pi\Z
  \end{matrix}\right\}
 \\
  \left\{\begin{matrix}
  C>1 \\ x\in\R
  \end{matrix}\right\}
 \end{matrix}\right]
\end{multline*}

$$
\arcsin x<C
 \quad\Longleftrightarrow\quad
 \left[\begin{matrix}
  \left\{\begin{matrix}
   C\le -\frac{\pi}{2}
    \\
   x\in\varnothing
  \end{matrix}\right\}
 \\
  \left\{\begin{matrix}
   C\in\left(-\frac{\pi}{2};\frac{\pi}{2}\right]  \\
   x\in[-1,\sin C)
  \end{matrix}\right\}
 \\
  \left\{\begin{matrix}
  C>\frac{\pi}{2}
  \\ x\in[-1;1]
  \end{matrix}\right\}
 \end{matrix}\right]
$$
 \eex

\ber Найдите общие правила решения неравенств:
$$
ax+b<0,
$$
$$
ax^2+bx+c<0,\quad ax^2+bx+c\le 0,
$$
$$
ax^2+bx+c>0,\quad ax^2+bx+c\ge 0,
$$
$$
x^\alpha<C,\quad a^x<C,
$$
$$
\cos x<C,\quad \arccos x<C,
$$
$$
\tg x<C,\quad \arctg x<C,
$$
$$
\ctg x<C,\quad \arcctg x<C.
$$
\eer

\paragraph{``Трехсложное определение''.}\label{3-opredelenie}

Как мы уже говорили, строгое определение области допустимых значений переменных
для выражения $\mathcal U$ невозможно дать иначе как по индукции, причем увязав
его одновременно с определением множества решений условия
 \beq\label{U-in-E}
{\mathcal U}\in E
 \eeq
(где $E$ --- произвольное множество на прямой $\R$), а также с замечанием, что
если область допустимых значений $D$ переменной $x$ в $\mathcal U$ определена,
то формула
$$
f(x)={\mathcal U}
$$
определяет функцию $f:D\to\R$ (с областью определения $D$).

Точное определение выглядит так:
 \bit{
\item[0)] Если $\mathcal U$ --- выражение порядка 0 от переменной $x$, то есть
$$
{\mathcal U}=x
$$
или
$$
{\mathcal U}=a,
$$
где $a$ --- какой-нибудь параметр, то в обоих случаях областью допустимых
значений переменной считается вся числовая прямая $\R$, причем условие
\eqref{U-in-E} в первом случае приобретает вид
$$
x\in E
$$
и, понятное дело, решением тогда будет само множество $E$; а во втором случае
---
 \beq\label{a-in-E}
a\in E
 \eeq
и решение (то есть множеством таких $x$, для которых это верно) будет зависеть
от $a$:
 \biter{
\item[---] если $a$ лежит в $E$, то это будет верно при любом $x\in\R$, и
поэтому решением \eqref{a-in-E} будет множество $\R$;

\item[---] если $a$ не лежит в $E$, то наоборот, при любом $x\in\R$ это будет
неверно, и поэтому решением \eqref{a-in-E} будет множество $\varnothing$.
 }\eiter

\item[n)] Если область допустимых значений переменной и множество решений
условия \eqref{U-in-E} определены для всех выражений порядка, не превышающего
некоторого $n\in\N$, то для выражения $\mathcal U$ порядка $n+1$ областью
допустимых значений переменных в числовом выражении $\mathcal U$ считается
множество решений связанной с этим выражением системы условий (уравнений,
неравенств или включений) $\varPhi$, определяемой следующими правилами:
 \biter{
\item[---] если $\mathcal U$ содержит в качестве подчиненного выражения
какое-то выражение вида $\frac{\mathcal P}{\mathcal Q}$, то в систему $\varPhi$
добавляется условие
$$
{\mathcal Q}\ne 0
$$

\item[---] если $\mathcal U$ содержит в качестве подчиненного выражения
какое-то выражение вида ${(\mathcal P)}^{(\mathcal Q)}$, то в систему $\varPhi$
добавляется следующая совокупность условий:\footnote{См. объяснение на с.
\pageref{usloviya-sushestv-x^a} по поводу условий существования выражения
$a^b$.}
$$
\left[\begin{matrix}
{\mathcal P}>0                    \\
\left\{\begin{matrix}{\mathcal P}=0 \\ {\mathcal Q}>0 \end{matrix}
\right\}    \\
\left\{\begin{matrix} {\mathcal P}<0 \\ {\mathcal Q}\in\frac{\Z}{2\N-1}
\end{matrix}\right\}
\end{matrix}\right]
$$

\item[---] если $\mathcal U$ содержит в качестве подчиненного выражения
какое-то выражение вида $\log_{(\mathcal P)} ({\mathcal Q})$, то в систему
$\varPhi$ добавляется следующая система условий:\footnote{См. аналогичное
объяснение на с. \pageref{usloviya-sushestvovaniya-log} по поводу $\log_a x$.}
$$
\left\{\begin{matrix}
{\mathcal P}>0                    \\
{\mathcal P}\ne 1 \\
{\mathcal Q}>0
\end{matrix}\right\}
$$

\item[---] если $\mathcal U$ содержит в качестве подчиненного выражения
какое-то выражение вида $\tg {\mathcal P}$, то в систему $\varPhi$ добавляется
условие
$$
{\mathcal P}\notin \frac{\pi}{2}+\pi\Z
$$

\item[---] если $\mathcal U$ содержит в качестве подчиненного выражения
какое-то выражение вида $\ctg {\mathcal P}$, то в систему $\varPhi$ добавляется
условие
$$
{\mathcal P}\notin \pi\Z
$$

\item[---] если $\mathcal U$ содержит в качестве подчиненного выражения
какое-то выражение вида $\arcsin {\mathcal P}$ или $\arccos {\mathcal P}$, то в
систему $\varPhi$ добавляется условие
$$
{\mathcal P}\in [-1,1]
$$
 }\eiter
После того, как определена область $D$ допустимых значений переменной в
$\mathcal U$, становится осмысленным равенство
$$
f(x)={\mathcal U},
$$
определяющее функцию $f:D\to\R$. Тогда уже множество решений условия
\eqref{U-in-E} для этого выражения определяется формулой
$$
\{x\in D:\; f(x)\in E\}
$$

}\eit

 \biter{
\item[$\bullet$] Функция $f$ (одного или нескольких переменных) называется {\it
стандартной}\label{def-stand-func}, если она определяется некоторым числовым
выражением $\mathcal P$:
$$
f(x)={\mathcal P}
$$
(или $f(a)={\mathcal P}$, или $f(x,y)={\mathcal P}$,...)

\item[$\bullet$] В частном случае, когда выражение $\mathcal P$ определяет
элементарную функцию
\begin{multline*}
x,\quad x^\alpha,\quad a^x,\quad \log_a x,\quad \sin x,\quad \cos x,\\
\tg
x,\quad \ctg x,\quad \arcsin x,\quad \arccos x,\\
 \arctg x,\quad \arcctg x
\end{multline*}
оно называется {\it элементарным}.

\item[$\bullet$] Область допустимых значений переменной в одноместном числовом
выражении $\mathcal P$ обозначается $\D(\mathcal P)$ (по аналогии с
обозначениями с. \pageref{def-domain-range}).

 }\eiter

\bprop\label{PROP:D(P+Q),...} Справедливы равенства:
\begin{align}
& \D\big({\mathcal P}+{\mathcal Q}\big)=\D({\mathcal P})\cap\D({\mathcal Q}) \label{D(P+Q)}\\
& \D\big({\mathcal P}-{\mathcal Q}\big)=\D({\mathcal P})\cap\D({\mathcal Q}) \label{D(P-Q)} \\
& \D\big({\mathcal P}\cdot{\mathcal Q}\big)=\D({\mathcal P})\cap\D({\mathcal Q}) \label{D(PQ)}\\
& \D\l \frac{\mathcal P}{\mathcal Q}\r=\D({\mathcal P})\cap\D({\mathcal Q})\setminus\{x: {\mathcal Q}=0\} \label{D(P/Q)}\\
& \D\big({\mathcal Q}\Big|_{y={\mathcal P}}\big)=\Big\{x:\ {\mathcal P}\in\D({\mathcal Q})\Big\}\cap\D({\mathcal P}) \label{D(Q-0-P)}
\end{align}
\eprop

\subsection{Равенство выражений}

Различаясь орфографически (то есть как последовательности символов), выражения могут быть {\it равны по смыслу}. Например, из формулы бинома Ньютона, которую мы доказали на с.\pageref{binom-Newtona} следует равенство
 \beq\label{(x+1)^2=x^2+2x+1}
(x+1)^2=x^2+2x+1.
 \eeq
Здесь степень в левой части фиксирована и невелика, поэтому эту формулу нетрудно доказать напрямую, без ссылок на бином Ньютона:
\begin{multline*}
(x+1)^2=(x+1)\cdot (x+1)=\eqref{(a+b)-cdot-c=a-cdot-c+b-cdot-c}=\\=x\cdot (x+1)+1\cdot (x+1)=\eqref{a-cdot-b=b-cdot-a},\eqref{(a+b)-cdot-c=a-cdot-c+b-cdot-c}=\\=
x^2+x\cdot 1+1\cdot x+1\cdot 1=\eqref{a-cdot-b=b-cdot-a},\eqref{a-cdot-1=a}=\\=
x^2+1\cdot x+1\cdot x+1=x^2+(1+1)\cdot x+1=\eqref{DEF:2}=\\=x^2+2x+1.
\end{multline*}

Глядя на эту цепочку, можно заметить, что здесь совершенно неважно, какой смысл вкладывается в символы в формуле \eqref{(x+1)^2=x^2+2x+1}: если справедливы тождества \eqref{(a+b)-cdot-c=a-cdot-c+b-cdot-c}, \eqref{a-cdot-b=b-cdot-a}, \eqref{a-cdot-1=a} и \eqref{DEF:2}, то  \eqref{(x+1)^2=x^2+2x+1} будет верно независимо от того, что мы понимаем под 1, 2, $x$ и +.

И действительно, искушенный читатель заметит, что равенство \eqref{(x+1)^2=x^2+2x+1} будет справедливо не только для вещественных чисел. Например, под 1 и $x$ могут пониматься комплексные числа, или матрицы, или вообще элементы произвольной ассоциативной алгебры, и все равно формула \eqref{(x+1)^2=x^2+2x+1} будет верна.

В соответствии с общим духом математики, такой уровень абстракции в вопросе о том, когда справедливо то или иное равенство, естественно было бы применить и к числовым выражениям. То есть вопрос, верно или неверно равенство числовых выражений
$$
{\mathcal P}={\mathcal Q}
$$
-- логично было бы решать, выясняя, можно ли это равенство вывести из заранее выписанного списка элементарных равенств (тождеств).
Такой формальный подход позволил бы сразу отделить Исчисление от Анализа, и, как всегда бывает в случаях, когда изложение удается разбить на независимые части, которые затем, складываясь вместе, обнаруживают неожиданные связи друг с другом, давая возможность поглядеть на предмет с разных точек зрения, это сильно упростило бы идеологию изучаемой дисциплины и повысило бы ее эстетическую ценность.

Однако (и это не делает чести специалистам по компьютерной алгебре и авторам бесчисленных руководств по дифференциальному и интегральному исчислению), к настоящему времени никаких исследований на этот счет не проводилось. Каким должен быть список элементарных тождеств, и можно ли при таком подходе добиться, чтобы равенство функций, определяемых выражениями, автоматически влекло за собой формальное равенство выражений (то есть выводимость этого равенства из элементарных тождеств для выражений)  --  все это предмет исследования для математиков будущего (и повод для упреков математикам настоящего).

Как следствие, нам, в соответствии с подходом наиболее близкой к этой области математической дисциплины, компьютерной алгебры, остается только определить равенство выражений просто как равенство определяемых этими выражениями функций.

 \biter{

\item[$\bullet$] Два числовых выражения ${\mathcal P}$ и ${\mathcal Q}$ считаются {\it
равными на множестве $x\in X$}, и изображается это формулой
$$
{\mathcal P}\underset{x\in X}{=}{\mathcal Q},
$$
или, в случае, если $X=\D(\mathcal P)\cap\D(\mathcal Q)$, то формулой
$$
{\mathcal P}={\mathcal Q},
$$
и, наконец, в случае, если $X=\D(\mathcal P)=\D(\mathcal Q)$, то формулой
$$
{\mathcal P}\equiv{\mathcal Q},
$$
если они порождают тождественно равные функции на множестве $X$:
$$
p(x)\underset{x\in X}{=}q(x),
$$
где $p$ и $q$ определяются равенствами
$$
p(x):={\mathcal P},\quad q(x):={\mathcal Q}.
$$

 }\eiter

\section{Стандартные функции в пройденных темах}\label{SEC:stand-func-v-analize}

Добавление стандартных функций к общему списку объектов, рассматриваемых в математическом анализе, стремительно расширяет возможности иллюстрировать идеи этой дисциплины. В этом параграфе мы приведем серию таких иллюстраций к трем основным пройденным на этот момент темам: последовательности, непрерывность и пределы. Одновременно здесь будут доказаны некоторые утверждения, которые понадобятся затем в иллюстративном материале последующих глав, в частности формулы первого и второго замечательных пределов.

\subsection{Числовые последовательности, определяемые стандартными функциями}

Примеры последовательностей, рассматривавшиеся в главах \ref{ch-R&N} и \ref{ch-x_n}, полезно дополнить следующими, определяемыми с помощью стандартных функций.

\paragraph{Метод математической индукции.}

\bers Докажите методом математической индукции (используя бином Ньютона):
 \biter{
\item[1)] $\frac{1}{2}\cdot \frac{3}{4}\cdot \frac{5}{6}\cdot ...
\frac{2n-1}{2n}\le \frac{1}{\sqrt{3n+1}}$

\item[2)] $\sqrt{n} < \sum_{k=1}^n \frac{1}{\sqrt{k}} < 2\sqrt{n} \quad (n\ge
2)$.
 }\eiter
 \eers

\paragraph{Предел последовательности.}

\bers Какими должны быть величины $\alpha, \beta$, чтобы почти все
 элементы последовательности $x_n$ лежали в интервале  $(\alpha;\beta)$?
 \biter{
 \item[1)] $x_n=\cos\frac{2\pi n}{3}$ \\
{\smsize (Ответ: $\alpha<-\frac{1}{2}, \, \beta>1$)}

\item[2)] $x_n=\sin\frac{2\pi n}{3}$ \\
{\smsize (Ответ: $\alpha<-\frac{\sqrt{3}}{2}, \, \beta>\frac{\sqrt{3}}{2}$)}

\item[3)] $x_n=\cos\frac{\pi n}{3}$ \\
{\smsize (Ответ: $\alpha<-1,\,\beta>1$)}

\item[4)] $x_n=\sin\frac{\pi n}{3}$ \\
{\smsize (Ответ: $\alpha<-\frac{\sqrt{3}}{2}, \, \beta>\frac{\sqrt{3}}{2}$)}
 }\eiter
\eers

\begin{ers}
Доказать что
 \biter{
\item[1)] $0\ne \lim\limits_{n\to \infty} \sin \frac{\pi n}{3}$
 }\eiter
\end{ers}

\begin{ers}
Доказать что
 \biter{
\item[1)] $\lim\limits_{n\to \infty} \sqrt{n}=+\infty$ (здесь расстояние между
соседними точками последова\-тельности убывает!)

\item[2)] $\lim\limits_{n\to \infty} \log_2 n=+\infty$ (здесь расстояние между
соседними точками последова\-тельности убывает!)

\item[3)] $\lim\limits_{n\to \infty} \log_3 \frac{1}{n}=-\infty$

\item[4)] $\lim\limits_{n\to \infty} (-n^\frac{2}{3})=-\infty$

\item[5)] $\lim\limits_{n\to \infty} (-1)^n \cdot n=\infty$
 }\eiter
\end{ers}

\begin{er}
Какие из последовательностей являются бесконечно малыми, а какие -- бесконечно
большими?
 \biter{
 \item[1)] $x_n=\frac{1}{\sqrt{n}}$
 \item[2)] $x_n=\sqrt{n}$
 }\eiter
\end{er}

\begin{ex}
Проверьте, что следующие последовательности ограничены:
$$
x_n= \sin n, \, x_n= \cos n, \, x_n= \arctg n.
$$
Какие из них сходятся (то есть имеют конечный предел)?
\end{ex}

\begin{ers}
Найти пределы
$$
 \lim\limits_{n\to \infty} \frac{\sin n}{n},
\quad \lim\limits_{n\to \infty} \frac{n+5}{n+\sin n}, \quad \lim\limits_{n\to
\infty} \frac{n+\arctg n}{n-\cos n},
$$
\end{ers}

\paragraph{Последовательности, заданные рекуррентно.}

\begin{ex}
Покажем, что последовательность, заданная рекуррентно
$$
x_1=0, \quad x_{n+1}=\sqrt{6+x_n}
$$
сходится, и найдем ее предел.

1. Вычислим сначала несколько первых элементов последовательности $\{ x_n \}$,
и нарисуем картинку:

%\picture{0pt}{0pt}{82.pcx}

\vglue100pt

2. Заметим, что все числа $\{ x_n \}$ неотри\-ца\-тель\-ны (это можно доказать
отдельно по индукции):
$$
x_n \ge 0
$$

3. Из картинки видно, что на первых нескольких элементах $\{ x_n \}$ наша
последовательность неубывает. Попробуем доказать, что она неубывает всегда:
\begin{equation}
x_n\le x_{n+1}\label{2.2.2}\end{equation} Поймем сначала, что означает это
неравенство:
 \begin{multline*}
x_n\le x_{n+1}\quad \Leftrightarrow \quad x_n\le \sqrt{6+x_n}\quad
 \kern-35pt\overset{\smsize \begin{matrix}
 \text{используем, что}\, x_n\ge 0 \\ \downarrow \\ \phantom{\cdot}
 \end{matrix}}{\Leftrightarrow}\\
 \Leftrightarrow \quad x_n^2\le 6+x_n \quad
\Leftrightarrow  \quad x_n^2-x_n-6\le 0 \quad \Leftrightarrow \\
\Leftrightarrow \quad
 \underbrace{(x_n+2)(x_n-3)\le 0}_{\smsize
 \begin{matrix}\uparrow \\ \text{делим на $x_n+2$;}\\
 \text{и, поскольку $x_n+2>0$}\\ \text{знак неравенства не меняется}\end{matrix}}\quad
 \Leftrightarrow \\ \Leftrightarrow \quad x_n-3\le 0 \quad \Leftrightarrow \quad
x_n\le 3 \end{multline*} Теперь докажем, последнее неравенство:
\begin{equation}
x_n\le 3 \label{2.2.3}\end{equation} Это делается математической индукцией.

a) Сначала проверяем наше неравенство при $n=1$:
$$
x_1=0\le 3
$$

b) Предполагаем, что оно верно при $n=k$:
\begin{equation}
x_k\le 3 \label{2.2.4}\end{equation}

c) Доказываем, что тогда оно будет верно при $n=k+1$.
\begin{equation}
x_{k+1}\le 3 \label{2.2.5}\end{equation} Действительно,
 \begin{gather*}
x_{k+1}\le 3 \quad \Leftrightarrow \quad \sqrt{6+x_k}\le 3 \quad
\Leftrightarrow \\ \Leftrightarrow \quad 6+x_k\le 9 \quad \Leftrightarrow \quad
x_k\le 3
 \end{gather*}
а последнее неравенство выполняется в силу \eqref{2.2.4}. Таким образом,
получилось что из \eqref{2.2.4} следует \eqref{2.2.5}, а это нам и нужно было.

Итак, мы доказали \eqref{2.2.3}, а значит и равносильное ему условие
\eqref{2.2.2}. Эти неравенства означают, что $\{ x_n \}$ неубывает и ограничена
сверху. Значит, по теореме Вейерштрасса \ref{Wei-I} она имеет предел. Обозначим
его буквой $c$:
$$
\lim_{n\to \infty} x_n =c
$$
и заметим, что $c\ge 0$, поскольку $x_n\ge 0$ (по теореме \ref{x_n<_y_n} о
переходе к пределу в неравенстве). Тогда из формулы $x_{n+1}=\sqrt{6+x_n}$ мы
получаем
$$
(x_{n+1})^2=6+x_n
$$
откуда
$$
c^2=\lim_{n\to \infty} (x_{n+1})^2=\lim_{n\to \infty} (6+x_n)=6+c
$$
То есть $c$ удовлетворяет квадратному уравнению
$$
c^2=6+c
$$
решая которое (с учетом, что $c\ge 0$), мы получаем $c=3$.

Ответ: $\lim\limits_{n\to \infty} x_n =3$
\end{ex}

\begin{ers}
Докажите, что последовательность, заданная рекуррентно сходится, и найдите ее
предел:
 \biter{
\item[1)] $x_1=1, \, x_{n+1}=\sqrt{x_n(6-x_n)}$;

\item[2)] $x_1=3, \, x_{n+1}=\frac{4x_n}{2+\sqrt{2x_n^2-4}}$;

\item[3)] $x_1=9, \, x_{n+1}=2\sqrt{x_n}$;

\item[4)] $x_1=1, \, x_{n+1}=3\sqrt{x_n}$;

\item[5)] $x_{n+1}=\frac{10x_n}{5+\sqrt{2x_n^2-25}}\; ,\; x_1=6$

\item[6)] $x_{n+1}=\frac{x_n 3\sqrt{2}}{\sqrt{9+x_n^4}}\; ,\; x_1=1$.
 }\eiter
\end{ers}

\paragraph{Подпоследовательности.}

\begin{ers}
Является ли данная последовательность $\{ y_k \}$ подпоследовательностью
последовательности $\{ x_n \}$? (Если да, то указать соответствующую
последовательность индексов $\{ n_k \}$.)
 \biter{
\item[1)] $x_n=\sqrt{n}, \, y_k=k$

 }\eiter
\end{ers}

\begin{ers}
Для данной последовательности $\{ x_n \}$ укажите какую-нибудь сходящуюся
подпоследовательность $\{ x_{n_k}\}$, если она существует.
 \biter{
\item[1)] $x_n=\cos\frac{2\pi n}{3}$
 }\eiter
\end{ers}

\begin{ers}
Проверьте с помощью критерия Коши сходимость последовательностей:
 \biter{
\item[1)] $x_n=\sin \frac{\pi n}{4}$;

\item[2)] $x_n=\cos \frac{\pi n}{7}$;

\item[3)] $x_n=\arctg (-1)^n$.
 }\eiter
\end{ers}

\subsection{Непрерывность стандартных функций}

Можно заметить, что в списке элементарных функций, который мы приводили на с.
\pageref{spisok-elem-functsij}, не все функции будут непрерывны. Покажем это на
примерах.

\bex Функция
$$
f(x)=0^x
$$
разрывна в точке $x=0$. Действительно, для последовательности
$$
x_n=\frac{1}{n}\underset{n\to\infty}{\longrightarrow} 0
$$
мы получим
$$
f(x_n)=0^{\frac{1}{n}}=0\underset{n\to\infty}{\longrightarrow}0\ne 1=0^0=f(0)
$$
\eex

\bex Функция
$$
f(x)=(-1)^x
$$
разрывна в любой точке своей области определения. Например, для точки $a=0$
можно взять последовательность
$$
x_n=\frac{1}{3n}\underset{n\to\infty}{\longrightarrow} 0
$$
для которой мы получим
$$
f(x_n)=(-1)^{\frac{1}{3n}}=-1\underset{n\to\infty}{\longrightarrow}-1\ne
1=(-1)^0=f(0)
$$
\eex

Доказательство следующего факта мы предоставляем читателю:

\btm\label{TH-nepr-elem-func} Следующие элементарные функции (и только они)
непрерывны в своей области определения:\etm
 {\smsize
 \begin{center}
\begin{tabular}{|l|l|l|} \hline
  &  &  \\
 название & $\begin{array}{l}  x \\ \text{\rotatebox{-90}{$\mapsto$}} \end{array}$ & область определения \\
  &  &  \\ \hline
  &  &  \\
 $\begin{matrix}\text{степенная}\\ \text{функция}\end{matrix}$ & $x^b$ &
$\begin{cases}x\in\R, & b\in\frac{\Z_+}{2\N-1} \\ x\ne 0, &
b\in-\frac{\N}{2\N-1} \\ x\ge 0, & b\in(0,+\infty)\setminus\frac{\Z}{2\N-1} \\
x>0, & b\in(-\infty,0)\setminus\frac{\Z}{2\N-1}
\end{cases}$
 \\
   &  &  \\
 \hline
  &  &  \\
 $\begin{matrix}\text{показательная}\\ \text{функция}\end{matrix}$ & $a^x$ &
$x\in\R$, $\ a>0$
 \\
   &  &  \\
 \hline
  &  &  \\
 логарифм & $\log_a x$ & $x>0$, $\ a\in(0;1)\cup(1;+\infty)$ \\
   &  &  \\
 \hline
   &  &  \\
синус & $\sin x$ & $x\in\R$ \\
  &  &  \\
\hline
  &  &  \\
 косинус & $\cos x$ & $x\in\R$ \\
   &  &  \\
\hline
  &  &  \\
 тангенс & $\tg x$ & $x\notin \frac{\pi}{2}+\pi\Z$  \\
   &  &  \\
\hline
  &  &  \\
 котангенс & $\ctg x$ & $x\notin \pi\Z$ \\
   &  &  \\
\hline
  &  &  \\
 арксинус & $\arcsin x$ & $x\in[-1;1]$ \\
   &  &  \\
\hline
  &  &  \\
 арккосинус & $\arccos x$ & $x\in[-1;1]$ \\
   &  &  \\
\hline
  &  &  \\
 арктангенс & $\arctg x$ & $x\in\R$ \\
   &  &  \\
\hline
  &  &  \\
 арккотангенс & $\arcctg x$ & $x\in\R$ \\
   &  &  \\
\hline
\end{tabular}
\end{center}

}

 \bit{
\item[$\bullet$] Функции из этого списка называются {\it непрерывными
элементарными функциями}.\label{DEF:nepr-elem-funktsii}
 }\eit

\begin{tm}\label{cont-elem}
Стандартная функция, составленная из непрерывных элементарных функций (на своих
областях определения), также непрерывна (на своей области определения).
\end{tm}

 \bit{
\item[$\bullet$] Стандартные функции, описываемые в этой теореме, называются
{\it непрерывными стандартными функциями}.
 }\eit

\begin{proof} По теоремам
\ref{cont-alg} и \ref{cont-composition}, если из непрерывных элементарных
функций конструировать новые функции с помощью алгебраических операций и
операции композиции, то будут получаться также непрерывные функции.
\end{proof}

\subsection{Вычисление пределов}\label{SUBSEC:vychisl-predelov}

Рассмотрим теперь вопрос, как вычисляются пределы стандартных функций:
$$
\lim_{x\to a} G(x)=?
$$

\paragraph{Существование предела.}
Прежде всего, конечно, нужно понимать, что предел существует не всегда.

\begin{ex}
Покажем что функция $f(x)=\sin \frac {1}{x}$ не имеет предела в точке $x=0$:
$$
\nexists \lim_{x\to 0}\sin \frac {1}{x}
$$
Возьмем сначала последовательность аргументов
$$
x_n =\frac{1}{\frac{\pi}{2}+2\pi n}\underset{n\to \infty}{\longrightarrow} 0
$$
тогда мы получим
$$
f(x_n)=\sin \frac {1}{x_n}= \sin \l \frac{\pi}{2}+2\pi n\r=1 \underset{n\to
\infty}{\longrightarrow} 1
$$
С другой стороны, если взять
$$
x_n =\frac{1}{-\frac{\pi}{2}+2\pi n}\underset{n\to \infty}{\longrightarrow} 0
$$
то мы получим
$$
f(x_n)=\sin \frac {1}{x_n}= \sin \l -\frac{\pi}{2}+2\pi n\r=-1 \underset{n\to
\infty}{\longrightarrow} -1
$$
Таким образом, мы получаем, что если взять одну последовательность аргументов
$x_n \underset{n\to \infty}{\longrightarrow} 0$, то соотвествующая
последовательность $f(x_n)$ будет стремиться к одному числу ($A=1$), а если
взять другую последовательность $x_n \underset{n\to \infty}{\longrightarrow}
0$, то $f(x_n)$ будет стремиться к другому числу ($A=-1$).

Это означает, что функция $f(x)=\sin \frac {1}{x}$ не имеет предела в точке
$x=0$, потому что иначе $f(x_n)$ должно было бы стремиться к этому конкретному
пределу $A$, независимо от того, какую берешь последовательность $x_n
\underset{n\to \infty}{\longrightarrow} 0$.
\end{ex}

\paragraph{Принцип подстановки.}

Во-вторых, нужно помнить, что если функция непрерывна, то по критерию непрерывности
\ref{cont-crit} ее предел равен ее значению в точке:

\begin{tm}[\bf принцип подстановки]
Пусть $G$ -- непрерывная стандартная функция, определенная в некоторой
окрестности точки $a$. Тогда предел функции $G$ в точке $x=a$ существует и
равен ее значению в этой точке:
$$
\exists \, \lim_{x\to a} G(x)=G(a).
$$
\end{tm}

Рассмотрим примеры.

\begin{ex}
Вычислить предел
$$
  \lim_{x\to 4}\frac {\sqrt{x}+x^2}{1+x}
$$
Решение:
 \begin{multline*}
\lim_{x\to 4}\frac {\sqrt{x}+x^2}{1+x}= \left(\text{подстановка}\right)=\\=
\frac {\sqrt{4}+4^2}{1+4}=\frac {2+16}{5}=\frac {18}{5}
 \end{multline*}
Число $\frac {18}{5}$ и будет ответом.
\end{ex}

\begin{ex}
Найти предел
$$
  \lim_{x\to 1} (x^2-2x+3)
$$
Решение:
 \begin{multline*} \lim_{x\to 1} (x^2-2x+3)= \left(\text{подстановка}\right)=\\= 1^2-2\cdot 1+3=1-2+3=2
 \end{multline*}\end{ex}

\paragraph{Вычисление упрощением.}

Обсудим теперь ситуацию, когда, по-прежнему, требуется найти предел
$$
\lim_{x\to a} G(x)=?
$$
но функция $G(x)$ не определена в точке $x=a$ (хотя, например, определена слева
и справа от $a$). Тогда говорят, что имеется так называемая {\it
неопределенность}, и алгоритм действий в таких случаях состоит в том, чтобы
подобрать новую функцию $F(x)$, которая совпадала бы с $G(x)$ вблизи точки $a$,
и при этом была бы непрерывной в точке $a$.

Зрительно изобразить это руководство можно так:
\begin{multline*}
\lim_{x\to a} G(x)=...\\...= {\smsize\begin{pmatrix}\text{подбираем функцию}\, F(x) \\
\text{совпадающую с}\, G(x) \\
\text{вблизи точки}\, a \\
\text{и непрерывную в точке}\, a
\end{pmatrix}}=...\\...= \lim_{x\to a} F(x)=F(a)
\end{multline*}
Проиллюстрируем это на примерах.

\begin{ex}[\bf упрощение + сокращение]
Найти предел
$$
\lim_{x\to 1}\left(\frac{1}{1-x} - \frac{2}{1-x^2}\right)
$$
При подстановке возникает неопределенность:
$$
\lim_{x\to 1}\left(\frac{1}{1-x} - \frac{2}{1-x^2}\right)= \frac{1}{0} -
\frac{2}{0}=?
$$
Значит, нужны преобразования. Попробуем упростить наше выражение:
 \begin{multline*}\lim_{x\to 1}\left(\frac{1}{1-x} - \frac{2}{1-x^2}\right)=\\=
\lim_{x\to 1}\left(\frac{1}{1-x} - \frac{2}{(1-x)(1+x)}\right)=\\= \lim_{x\to
1}\frac{1+x-2}{(1-x)(1+x)}=\lim_{x\to 1}\frac{x-1}{(1-x)(1+x)} =\\= \lim_{x\to
1}\frac{-1}{1+x} = \lim_{x\to 1}\frac{-1}{1+x} =\\=
(\text{подстановка})=\frac{-1}{1+1}=-\frac{1}{2}\end{multline*}\end{ex}

\begin{ex}[\bf разложение на множители]
Найти предел
$$
  \lim_{x\to -1}\frac{x^2-1}{x^2-x-2}
$$
Здесь возникает неопределенность:
$$
\lim_{x\to -1}\frac{x^2-1}{x^2-x-2}=\frac{1-1}{1+1-2}=\frac{0}{0}=?
$$
Это означает, что необходимы преобразования. В данном случае они состоят в том,
что числитель и знаменатель нужно разложить на множители, после чего сократить
дробь:
 \begin{multline*}
\lim_{x\to -1}\frac{x^2-1}{x^2-x-2}= \lim_{x\to
-1}\frac{(x-1)(x+1)}{(x-2)(x+1)}=\\= \lim_{x\to -1}\frac{x-1}{x-2}=
(\text{\scriptsize снова подстановка})=\frac{-1-1}{-1-2}=\frac{2}{3}
 \end{multline*}\end{ex}

\begin{ex}[\bf использование сопряженного радикала]
Вычислить предел
$$
\lim_{x\to 0}\frac{\sqrt{1+x}-\sqrt{1-x}}{x}
$$
При подстановке получается неопределенность:
$$
\lim_{x\to 0}\frac{\sqrt{1+x}-\sqrt{1-x}}{x}=\frac{0}{0}
$$
Поэтому нужны какие-то преобразования. В примерах с корнями, как этот, в
качестве такого преобразования можно использовать {\it домножение на
сопряженный радикал}.

Под этим понимается вот что. В числителе имеется иррациональное выражение
$\sqrt{1+x}-\sqrt{1-x}$, называемое {\it радикалом}. Подберем такое выражение,
которое при умножении на него ``уничтожает корни''. Оно называется {\it
сопряженным радикалом}, и в нашем случае имеет вид $\sqrt{1+x}+\sqrt{1-x}$.
Помножим числитель и знаменатель дроби на сопряженный радикал:
 \begin{multline*}\lim_{x\to 0}\frac{\sqrt{1+x}-\sqrt{1-x}}{x}=\\= \lim_{x\to
0}\frac{(\sqrt{1+x}-\sqrt{1-x})\cdot (\sqrt{1+x}+\sqrt{1-x})}{x\cdot
(\sqrt{1+x}+\sqrt{1-x})}=\\= \lim_{x\to
0}\frac{\left(\sqrt{1+x}\right)^2-\left(\sqrt{1-x}\right)^2} {x\cdot
(\sqrt{1+x}+\sqrt{1-x})}=\\= \lim_{x\to 0}\frac{(1+x)-(1-x)}{x\cdot
(\sqrt{1+x}+\sqrt{1-x})}=\\= \lim_{x\to 0}\frac{2x}{x\cdot
(\sqrt{1+x}+\sqrt{1-x})}=\\= \lim_{x\to 0}\frac{2}{\sqrt{1+x}+\sqrt{1-x}}=
\left(\text{подстановка}\right)=\\= \frac{2}{\sqrt{1+0}+\sqrt{1-0}}=
\frac{2}{1+1}=1
\end{multline*}\end{ex}

\begin{ex}[\bf использование сопряженного радикала]
Вычислить предел
$$
\lim_{x\to 0}\frac{\sqrt[3]{1+x^2}-1}{x^2}
$$
При подстановке имеем неопределенность:
$$
\lim_{x\to 0}\frac{\sqrt[3]{1+x^2}-1}{x^2}= \frac{\sqrt[3]
{1+0}-1}{0}=\frac{0}{0}
$$
Значит нужны преобразования, которые здесь опять состоят в умножении на
сопряженный радикал. В нашем примере радикалом является выражение $\sqrt[3]
{1+x^2}-1$. Сопряженный же радикал -- то есть выражение, которое при умножении
на наш радикал ``сокращает корни'' -- имеет вид $\left(\sqrt[3]
{1+x^2}\right)^2+\sqrt[3] {1+x^2}+1$. Умножим на него числитель и знаменатель:
 \begin{multline*}\lim_{x\to 0}\frac{\sqrt[3] {1+x^2}-1}{x^2}=\\ {\smsize =\lim_{x\to
0}\frac{\left(\sqrt[3] {1+x^2}-1 \right) \left(\left(\sqrt[3]
{1+x^2}\right)^2+\sqrt[3] {1+x^2}+1 \right)} {x^2 \left(\left(\sqrt[3]
{1+x^2}\right)^2+\sqrt[3] {1+x^2}+1 \right)}=}
\\=
\lim_{x\to 0}\frac{1+x^2-1}{x^2 \left(\left(\sqrt[3] {1+x^2}\right)^2+ \sqrt[3]
{1+x^2}+1 \right)}=\\= \lim_{x\to 0}\frac{x^2}{x^2 \left(\left(\sqrt[3]
{1+x^2}\right)^2+ \sqrt[3] {1+x^2}+1 \right)}=\\= \lim_{x\to
0}\frac{1}{\left(\sqrt[3] {1+x^2}\right)^2+ \sqrt[3] {1+x^2}+1}=\\=
\left(\text{подстановка}\right)= \frac{1}{\left(\sqrt[3]
{1+0}\right)^2+\sqrt[3] {1+0}+1}=\\=
\frac{1}{1+1+1}=\frac{1}{3}\end{multline*}\end{ex}

\begin{ers} Найти пределы
 \biter{
\item[1)] $\lim\limits_{x\to \pi /2}\sin x$

\item[2)] $\lim\limits_{x\to 0}\arctg{x}-1+x^3$

\item[3)] $\lim\limits_{x\to 0}\frac{\sqrt{x}+1}{\arcsin{x}-4}$

\item[4)] $\lim\limits_{x\to \sqrt{3}}\frac{x^2-3}{x^4+x^2+1}$,

\item[5)] $\lim\limits_{x\to 2}\frac{x^2-5x+6}{x^2-2x}$,

\item[6)] $\lim\limits_{x\to 0}\frac{\sqrt{1+x^2}-1}{x^2}$

\item[7)] $\lim\limits_{x\to 0}\frac{\sqrt{1+x^2}-1}{\sqrt{16+x^2}-4}$

\item[8)] $\lim\limits_{x\to 0}\frac{\sqrt[3] {1+x}-\sqrt[3] {1-x}}{x}$

\item[9)] $\lim\limits_{x\to 1}\frac{\sqrt{x}-1}{\sqrt[3] {x}-1}$

\item[10)] $\lim\limits_{x\to -8}\frac{\sqrt{1-x}-3}{2+\sqrt[3] {x}}$
 }\eiter
\end{ers}

\paragraph{Замена переменной под знаком предела.}
В следующих примерах используется теорема \ref{TH:zamena-perem-v-lim} о замене
переменной под знаком предела.

\begin{ex}
Найти предел
$$
\lim_{x\to 16}\frac{\sqrt{x}-4}{\sqrt{x}-\sqrt[4]{x}-2}
$$
Решение:
 \begin{multline*}\lim_{x\to 16}\frac{\sqrt{x}-4}{\sqrt{x}-\sqrt[4]{x}-2}=
\left| \begin{array}{c}\sqrt[4]{x}=y, \, x=y^4 \\
y\underset{x\to 16}{\longrightarrow} 4 \\
y\ne 4 \, \text{при}\, x\ne 16 \end{array}\right| =\\= \lim_{y\to
4}\frac{y^2-4}{y^2-y-2}= \left(\frac{0}{0}\right)=\lim_{y\to
4}\frac{(y-2)(y+2)}{(y-2)(y+1)}=\\= \lim_{y\to 4}\frac{y+2}{y+1}=
\frac{4+2}{4+1}= \frac{6}{5}\end{multline*}\end{ex}

\begin{ex}
Найти предел
$$
\lim_{x\to 1}\frac{\sqrt{x}-1}{\sqrt[3]{x}-1}
$$
Решение:
 \begin{multline*}\lim_{x\to 1}\frac{\sqrt{x}-1}{\sqrt[3]{x}-1}=
\left| \begin{array}{c}\sqrt[6]{x}=y, \, x=y^6 \\
y\underset{x\to 1}{\longrightarrow} 1 \\
y\ne 1 \, \text{при}\, x\ne 1 \end{array}\right| =\\= \lim_{y\to
1}\frac{y^3-1}{y^2-1}= \left(\frac{0}{0}\right)=\\= \lim_{y\to
1}\frac{(y-1)(y^2+y+1)}{(y-1)(y+1)}= \lim_{y\to 1}\frac{y^2+y+1}{y+1}=\\=
\frac{1+1+1}{1+1}=\frac{3}{2}\end{multline*}\end{ex}

\begin{ex}
Найти предел
$$
\lim_{x\to 0}\frac{\cos^2 x -1}{\cos^2 x+2 \cos x -3}
$$
Решение:
 \begin{multline*}\lim_{x\to 0}\frac{\cos^2 x -1}{\cos^2 x+2 \cos x -3}=
\left| \begin{array}{c}\cos x=y, \\
y\underset{x\to 0}{\longrightarrow} 1
\end{array}\right| =\\=
\lim_{y\to 1}\frac{y^2 -1}{y^2 +2 y -3}= \lim_{y\to 1}\frac{(y
-1)(y+1)}{(y-1)(y+3)}=\\= \lim_{y\to 1}\frac{y+1}{y+3}=\frac{2}{4}=\frac{1}{2}
 \end{multline*}\end{ex}

\begin{ex}
Вычислить предел
$$
\lim_{x\to +\infty}\Big(\sqrt{1+x^2}-x\Big)
$$
Решение:
 \begin{multline*}\lim_{x\to +\infty}\Big(\sqrt{1+x^2}-x\Big)=
\left| \begin{array}{c} x=\frac{1}{y}, \, y=\frac{1}{x}\\
y\underset{x\to +\infty}{\longrightarrow} 0+0
\end{array}\right| =\\=
\lim_{y\to 0+0}\left(\sqrt{1+\frac{1}{y^2}}-\frac{1}{y}\right)=\\= \lim_{y\to
0+0}\left(\frac{1}{\sqrt{y^2}}\sqrt{y^2+1}-\frac{1}{y}\right)=
\\=
{\smsize\begin{pmatrix}\text{используем тождество}\\
\sqrt{y^2}=|y| \end{pmatrix}}=\\= \lim_{y\to 0+0}\left(\frac{1}{|y|}\sqrt{y^2+1}-\frac{1}{y}\right)=\\= {\smsize\begin{pmatrix} y\to 0+0, \\ \text{значит}\, y>0, \\
\text{поэтому}\, |y|=y \end{pmatrix}}=\\= \lim_{y\to
0+0}\left(\frac{1}{y}\sqrt{y^2+1}-\frac{1}{y}\right)= \lim_{y\to
0+0}\frac{\sqrt{y^2+1}-1}{y}=\\= \lim_{y\to
0+0}\frac{(\sqrt{y^2+1}-1)(\sqrt{y^2+1}+1)}{y (\sqrt{y^2+1}+1)}=\\= \lim_{y\to
0+0}\frac{y^2+1-1}{y (\sqrt{y^2+1}+1)}= \lim_{y\to 0+0}\frac{y^2}{y
(\sqrt{y^2+1}+1)}=\\= \lim_{y\to 0+0}\frac{y}{\sqrt{y^2+1}+1}=\\= \lim_{y\to
0+0}\frac{y}{\sqrt{y^2+1}+1}= \left(\text{подстановка}\right)=\\=
\frac{0}{\sqrt{0+1}+1}=0 \end{multline*}\end{ex}

\begin{ex}
Вычислить предел
$$
\lim_{x\to -\infty}\frac{x}{\sqrt{x^2+5}}
$$
Решение:
 \begin{multline*}\lim_{x\to -\infty}\frac{x}{\sqrt{x^2+5}}=
\left| \begin{array}{c} x=\frac{1}{y}, \, y=\frac{1}{x}\\
y\underset{x\to -\infty}{\longrightarrow} 0-0
\end{array}\right| =\\=
\lim_{y\to 0-0}\frac{\frac{1}{y}}{\sqrt{\frac{1}{y^2}+5}}= \lim_{y\to
0-0}\frac{\frac{1}{y}}{\sqrt{\frac{1}{y^2}+5}}=\\= \lim_{y\to
0-0}\frac{\frac{1}{y}}{\frac{1}{\sqrt{y^2}}\sqrt{1+5
y^2}}= {\smsize {\smsize\begin{pmatrix}\text{используем}\\
\text{тождество}\\
\sqrt{y^2}=|y| \end{pmatrix}}}=\\= \lim_{y\to
0-0}\frac{\frac{1}{y}}{\frac{1}{|y|}\sqrt{1+5 y^2}}= \left({\smsize
\begin{array}{c}
y\to 0-0, \\ \text{значит}\, y<0,  \\ \text{поэтому}\, |y|=-y
\end{array}}\right)=\\= \lim_{y\to 0-0}\frac{\frac{1}{y}}{\frac{1}{-y}\sqrt{1+5
y^2}}= \lim_{y\to 0-0}\frac{1}{-\sqrt{1+5 y^2}}=\\= \frac{1}{-\sqrt{1+5\cdot
0}}=-1 \end{multline*}\end{ex}

\begin{ers} Найти пределы
 \biter{
\item[1)] $\lim\limits_{x\to -\frac{\pi}{2}}\frac{\sin^2 x-\sin x-2}{\sin^3
x+1}$

\item[2)] $\lim\limits_{x\to +\infty}\frac{2^{3x+2}-3\cdot 2^x}{1-2^{3x+1}}$

\item[3)] $\lim\limits_{x\to -\infty}\frac{2^{3x+2}-3\cdot 2^x}{1-2^{3x+1}}$

\item[4)] $\lim\limits_{x\to 0}\frac{\sqrt{1+\sin x}-1}{\sin x}$

\item[5)] $\lim\limits_{x\to +\infty}\frac{\sqrt{1+9^x}}{1+3^x}$

\item[6)] $\lim\limits_{x\to -\infty}\frac{\sqrt{1+9^x}}{1+3^x}$

\item[7)] $\lim\limits_{x\to 0}\frac{16^x-1}{8^x-1}$

\item[8)] $\lim\limits_{x\to \infty}\left(\sqrt{x^2+1}-\sqrt{x^2-1}\right)$

\item[9)] $\lim\limits_{x\to \infty} x (\sqrt{x^2+1}-x)$

\item[10)] $\lim\limits_{x\to \infty}\frac{2x^2-3x-4}{\sqrt{x^4+1}}$

\item[11)] $\lim\limits_{x\to
+\infty}\frac{\sqrt{x}}{\sqrt{x+\sqrt{x+\sqrt{x}}}}$
 }\eiter
\end{ers}

\paragraph{Первый замечательный предел.} При вычислении пределов с
тригонометрическими функциями оказывается полезным следующий факт.

\begin{tm}\label{sin_x/x}
 \begin{equation}
 \lim_{x\to 0}\frac{\sin x}{x}=1
 \label{5.8.1}
 \end{equation}\end{tm}

 \bit{\item[$\bullet$] Это равенство называется
{\it первым замечательным пределом}\index{предел!первый замечательный}. }\eit

\begin{proof}
Пусть сначала
$$
x_n\underset{n\to \infty}{\longrightarrow} 0, \qquad x_n>0
$$
Почти все $x_n$ лежат в интервале $(0;\frac{\pi}{2})$, поэтому для них из
тройного неравенства \eqref{0<sin-t<t} мы получаем цепочку:
$$
0<\sin x_n< x_n< \tg x_n
$$
$$
\phantom{(\text{\scriptsize делим на $\sin x_n$})} \Downarrow
(\text{\scriptsize делим на $\sin x_n$})
$$
$$
1<\frac{x_n}{\sin x_n} <\frac{1}{\cos x_n}
$$
$$
\phantom{(\text{\scriptsize возводим в степень $-1$})} \Downarrow
(\text{\scriptsize возводим в степень $-1$})
$$
$$
1>\frac{\sin x_n}{x_n} >\kern-105pt \underbrace{\cos
x_n}_{\scriptsize\begin{matrix}\phantom{\text{($\cos$ -- непрерывная функция)}}
\quad  \downarrow \quad \text{($\cos$ -- непрерывная функция)} \\ \cos 0 \\ \|
\\ 1
\end{matrix}}
$$
$$
\phantom{(\text{\scriptsize теорема
\ref{milit}})} \Downarrow (\text{\scriptsize теорема \ref{milit}})
$$
$$
\frac{\sin x_n}{x_n}\underset{n\to \infty}{\longrightarrow} 1.
$$
Это верно для всякой последовательности $x_n\underset{n\to
\infty}{\longrightarrow} 0$, $x_n>0$, значит
$$
\lim_{x\to +0}\frac{\sin x}{x}=1
$$
Отсюда уже следует, что
 \begin{multline*}
\lim_{x\to -0}\frac{\sin x}{x}=\left|\begin{matrix}t=-x \\
t\to+0\end{matrix}\right|=\lim_{t\to +0}\frac{\sin(-t)}{-t}=\\=\lim_{t\to
+0}\frac{-\sin t}{-t}=\lim_{t\to +0}\frac{\sin t}{t}=1
 \end{multline*}
Вместе эти равенства дают \eqref{5.8.1}. \end{proof}

Обсудим теперь применения первого замечательного предела.

\begin{ex}
Найти предел
$$
\lim_{x\to 0}\frac{\sin 2x}{x}
$$
Здесь можно сделать замену переменных, после чего наш предел сведется к первому
замечательному:
 \begin{multline*}
 \lim_{x\to 0}\frac{\sin 2x}{x}={\smsize
\left| \begin{array}{c} 2x=y, \, x=\frac{y}{2}\\
y\underset{x\to 0}{\longrightarrow} 0 \\
y\ne 0 \, \text{при}\, x\ne 0 \end{array}\right|} =\\= \lim_{y\to 0}\frac{2\sin
y}{y}= 2 \lim_{y\to 0}\frac{\sin y}{y}=2\cdot 1=2
 \end{multline*}\end{ex}

\begin{ex}
 \begin{multline*}
\lim_{x\to \infty} x\cdot \sin \frac{1}{x}={\smsize
\left| \begin{array}{c}\frac{1}{x}=y, \, x=\frac{1}{y}\\
y\underset{x\to \infty}{\longrightarrow} 0 \\
y\ne 0 \, \text{при любом}\, x \end{array}\right|} =\\= \lim_{y\to 0}\frac{\sin
y}{y}=1
 \end{multline*}\end{ex}

\begin{ex}
 \begin{multline*}
 \lim_{x\to 0}\frac{\tg x}{x}= \lim_{x\to 0}\frac{\sin x}{x \cos x}=
\lim_{x\to 0}\frac{\sin x}{x}\cdot \lim_{x\to 0}\frac{1}{\cos x}=\\= 1\cdot 1=1
 \end{multline*}\end{ex}

\begin{ex}

\begin{multline*}\lim_{x\to 0}\frac{1-\cos x}{x^2}= \lim_{x\to 0}\frac{2\sin^2
\frac{x}{2}}{x^2}= 2 \lim_{x\to 0}\left(\frac{\sin \frac{x}{2}}{x}\right)^2=\\=
2 \left(\lim_{x\to 0}\frac{\sin \frac{x}{2}}{x}\right)^2={\smsize
\left| \begin{array}{c}\frac{x}{2}=y, \, x=2y \\
y\underset{x\to 0}{\longrightarrow} 0 \\
y\ne 0 \, \text{при}\, x\ne 0 \end{array}\right|} =\\= 2 \left(\lim_{y\to
0}\frac{\sin y}{2y}\right)^2= \frac {1}{2}\left(\lim_{y\to 0}\frac{\sin
y}{y}\right)^2= \frac {1}{2}\cdot 1= \frac {1}{2}\end{multline*}\end{ex}

\begin{ers}
Найти пределы
 \biter{
\item[1)] $\lim\limits_{x\to 0}\frac{\sin 5x}{3x}$

\item[2)] $\lim\limits_{x\to 0}\frac{\sin \alpha x}{\sin \beta x}$

\item[3)] $\lim\limits_{x\to \infty} (x-3) \sin \frac{2}{3x}$

\item[4)] $\lim\limits_{x\to 0} x \ctg x$

\item[5)] $\lim\limits_{x\to \infty} x^2 \sin \frac{1}{x}$

\item[6)] $\lim\limits_{x\to 0}\frac{\tg x -\sin x}{\sin^3 x}$
 }\eiter
\end{ers}

\paragraph{Второй замечательный предел.} Следующее утверждение оказывается также весьма полезным при вычислении пределов:

\begin{tm}\label{II-lim}
 \begin{equation}\lim_{x\to
0}\left(1+x \right)^\frac{1}{x}=e= \lim_{x\to
\infty}\left(1+\frac{1}{x}\right)^x \label{5.9.2}\end{equation}
\end{tm}

\bit{\item[$\bullet$] Эти (эквивалентные) равенства называются {\it вторым
замечательным пределом}\index{предел!второй замечательный}. }\eit

\begin{proof} Равенство между этими двумя пределами получается заменой переменной, поэтому нам
достаточно доказать только второе равенство в этой цепочке. Для начала будем
считать, что $x$ положительно, то есть докажем равенство
 \beq\label{5.9.2+}
\lim_{x\to+\infty}\left(1+\frac{1}{x}\right)^x=e
 \eeq
Зафиксируем последовательность $\{ x_k \}$, стремящуюся к $+\infty$,
$$
x_k\underset{k\to \infty}{\longrightarrow} +\infty
$$
и покажем, что
$$
\lim_{k\to \infty}\left( 1+\frac{1}{x_k}\right)^{x_k}=e
$$
Действительно, всякое число $x_k$ лежит между какими-нибудь двумя соседними
натуральными числами, поэтому можно найти число $n_k\in \mathbb{N}$ такое, что
\begin{equation}
n_k\le x_k <n_k+1 \label{5.9.3}\end{equation} Возникающая таким образом
последовательность натуральных чисел $\{ n_k \}$ будет бесконечно большой
$$
n_k \underset{k\to \infty}{\longrightarrow} +\infty
$$
потому что $n_k>x_k-1 \underset{k\to \infty}{\longrightarrow} +\infty$. Значит,
можно рассмотреть подпоследовательность $\{ r_{n_k}\}$, последовательности $\{
r_n \}$, рассматривавшейся в лемме \ref{e}, и поскольку $\{ r_n \}$ стремится к
$e$, по свойству $1^0$ пункта \ref{SEC-BOL-WEI} мы получим, что $\{ r_{n_k}\}$
тоже должна стремится к $e$:
\begin{equation}\left( 1+\frac{1}{n_k}\right)^{n_k}\underset{k\to
\infty}{\longrightarrow} e \label{5.9.4}\end{equation} Аналогично доказывается
что
\begin{equation}\left( 1+\frac{1}{n_k+1}\right)^{n_k+1}\underset{k\to
\infty}{\longrightarrow} e \label{5.9.5}\end{equation}

Теперь рассмотрим двойное неравенство \eqref{5.9.3}. Из него следует
$$
\frac{1}{n_k}\ge \frac{1}{x_k} >\frac{1}{n_k+1}
$$
а отсюда, в свою очередь, получается
$$
1+\frac{1}{n_k}\ge 1+\frac{1}{x_k} >1+\frac{1}{n_k+1}
$$
и
$$
\left(1+\frac{1}{n_k}\right)^{x_k}\ge \left(1+\frac{1}{x_k}\right)^{x_k}
> \left(1+\frac{1}{n_k+1}\right)^{x_k}
$$
Расширим это двойное неравенство влево и вправо:

\begin{multline*}\left(1+\frac{1}{n_k}\right)^{n_k+1}
> {\smsize\begin{pmatrix}\text{используем неравенство}\\
n_k+1>x_k \end{pmatrix}}>\\ > \left(1+\frac{1}{n_k}\right)^{x_k}\ge
\left(1+\frac{1}{x_k}\right)^{x_k} >\\> \left(1+\frac{1}{n_k+1}\right)^{x_k}\ge
{\smsize\begin{pmatrix}\text{используем неравенство}\\ x_k \ge n_k
\end{pmatrix}} \ge\\ \ge
\left(1+\frac{1}{n_k+1}\right)^{n_k}\end{multline*} Теперь выбросим ненужные
звенья в этой цепочке:
$$
\underbrace{\left(1+\frac{1}{n_k}\right)^{n_k+1}}_{\scriptsize\begin{matrix}\phantom{\eqref{5.9.4}}\ \downarrow \ \eqref{5.9.4}
\\ e \end{matrix}}
\kern-7pt>
\left(1+\frac{1}{x_k}\right)^{x_k}
\kern-3pt
>\kern-4pt
\underbrace{\left(1+\frac{1}{n_k+1}\right)^{n_k}}_{\scriptsize\begin{matrix}\phantom{\eqref{5.9.5}}\ \downarrow \ \eqref{5.9.5}
\\ e \end{matrix}}
$$
$$
\phantom{(\text{\scriptsize теорема
\ref{milit}})} \Downarrow (\text{\scriptsize теорема \ref{milit}})
$$
$$
\left( 1+\frac{1}{x_k}\right)^{x_k}\underset{k\to \infty}{\longrightarrow} e
$$
Мы доказали равенство \eqref{5.9.2+}. Из него теперь получаем
 \begin{multline*}
\lim_{x\to-\infty}\left(1+\frac{1}{x}\right)^x=\left|\begin{matrix}y=-x
\\ y\to+\infty\end{matrix}\right|=\\=
\lim_{y\to+\infty}\left(1-\frac{1}{y}\right)^{-y}=
\lim_{y\to+\infty}\left(\frac{y-1}{y}\right)^{-y}=\\=
\lim_{y\to+\infty}\left(\frac{y}{y-1}\right)^y=
\lim_{y\to+\infty}\left(1+\frac{1}{y-1}\right)^y=\\=\left|\begin{matrix}z=y-1
\\ z\to+\infty\end{matrix}\right|=\lim_{z\to+\infty}\left(1+\frac{1}{z}\right)^{z+1}=\\=
\underbrace{\lim_{z\to+\infty}\left(1+\frac{1}{z}\right)^z}_{\scriptsize\begin{matrix}\text{\rotatebox{90}{$=$}}\\
\phantom{,} e,\\ \text{в силу \eqref{5.9.2+}}\end{matrix}}\cdot
\underbrace{\lim_{z\to+\infty}\left(1+\frac{1}{z}\right)}_{\scriptsize\begin{matrix}\text{\rotatebox{90}{$=$}}\\
1\end{matrix}}=e
 \end{multline*}
В результате мы доказали
 \beq\label{5.9.2-}
\lim_{x\to-\infty}\left(1+\frac{1}{x}\right)^x=e
 \eeq
и вместе с \eqref{5.9.2+} это дает второе равенство в \eqref{5.9.2}.
\end{proof}

Теорема \ref{II-lim} позволяет вычислять сложные пределы, в которых возникает
неопределенность типа $1^\infty$.

\begin{ex}
Найти предел
$$
\lim_{x\to 0}\left( 1+\tg x \right)^{\ctg x}
$$
Решение:
 \begin{multline*}
\lim_{x\to 0}\left( 1+\tg x \right)^{\ctg x}={\smsize
\left| \begin{array}{c}\tg x=y, \, \ctg x=\frac{1}{y}\\
y\underset{x\to 0}{\longrightarrow}\infty
\end{array}\right|} =\\=
\lim_{y\to 0}\left( 1+y \right)^{\frac{1}{y}}= e
 \end{multline*}\end{ex}

\begin{ex}
Найти предел
$$
\lim_{x\to \infty}\left(\frac{x+1}{x-1}\right)^x
$$
Решение:

\begin{multline*}\lim_{x\to \infty}\left(\frac{x+1}{x-1}\right)^x= \lim_{x\to
\infty}\left(1+ \frac{x+1}{x-1} -1 \right)^x=\\= \lim_{x\to \infty}\left(1+
\frac{2}{x-1}\right)^x={\smsize
\left| \begin{array}{c}\frac{2}{x-1}=y, \, x=1+\frac{2}{y}\\
y\underset{x\to \infty}{\longrightarrow} 0
\end{array}\right|} = \\=
\lim_{y\to 0}\left( 1+y \right)^{1+\frac{2}{y}}= \lim_{y\to 0}\left( 1+y
\right) \left[\left( 1+y \right)^\frac{1}{y}\right]^2=\\= \lim_{y\to 0}\left(
1+y \right) \left[\lim_{y\to 0}\left( 1+y \right)^\frac{1}{y}\right]^2= e^2
\end{multline*}\end{ex}

\begin{ex}
Найти предел
$$
\lim_{x\to 0}\left(\cos x \right)^\frac{1}{x^2}
$$
Решение:

\begin{multline*}\lim_{x\to 0}\left(\cos x \right)^\frac{1}{x^2}= \lim_{x\to 0}\left( 1-2\sin^2 \frac{1}{x}\right)^\frac{1}{x^2}=\\= \lim_{x\to 0}\left[ \left( 1-2\sin^2 x \right)^\frac{1}{-2\sin^2 x}\right]^\frac{-2\sin^2 x}{x^2}=\\= \left[ \lim_{x\to 0}\left(
1-2\sin^2 x \right)^\frac{1}{-2\sin^2 x}\right]^{\lim_{x\to 0}\frac{-2\sin^2
x}{x^2}}=\\={\smsize
\left| \begin{array}{c} -2\sin^2 x=y, \\
y\underset{x\to 0}{\longrightarrow} 0
\end{array}\right|} =\\=
\left[ \lim_{y\to 0}\left( 1+y \right)^\frac{1}{y}\right]^{-2 \lim\limits_{x\to
0}\left(\frac{\sin x}{x}\right)^2}= e^{-2}\end{multline*}\end{ex}

\begin{ers} Найти пределы

1) $\lim\limits_{x\to 0} (1+3\tg^2 x)^{\ctg^2 x}$

2) $\lim\limits_{x\to \infty}\left(\frac{2x}{2x-3}\right)^{3x}$

3) $\lim\limits_{x\to +\infty}\left(\frac{5x^3+2}{5x^3}\right)^{\sqrt{x}}$

4) $\lim\limits_{x\to 0}\left( 1+2x^2 \right)^\frac{1}{x^2}$

5) $\lim\limits_{x\to \infty}\left(\frac{2x+3}{2x+1}\right)^{x+1}$

6) $\lim\limits_{x\to 0}\left( 1+\tg x \right)^{\ctg x}$

7) $\lim\limits_{x\to 0}\left( 1+3\tg^2 x \right)^{\ctg^2 x}$

8) $\lim\limits_{x\to 0}\left(\cos x \right)^{\frac{1}{\sin^2 x}}$

9) $\lim\limits_{x\to 0}\left(\cos^2 x \right)^{\frac{1}{\sin^2 x}}$

10) $\lim\limits_{x\to 0}\left(\cos x \right)^{\ctg x}$

11) $\lim\limits_{x\to \infty}\left(\frac{2x+3}{2x-1}\right)^{x+1}$

12) $\lim\limits_{x\to \infty}\left(\frac{x^2+3x+1}{x^2-x+2}\right)^{x^2}$

13) $\lim\limits_{x\to \infty}\left(\frac{5x^2+1}{5x^2+x+1}\right)^{x}$

14) $\lim\limits_{x\to \infty}\left(\frac{x+1}{x+3}\right)^{x^2}$
\end{ers}

\bcor Справедливы соотношения:
 \begin{equation}\lim_{x\to 0}\frac{\ln (1+x)}{x}=1
 \label{5.10.2}
 \end{equation}
и
 \begin{equation}\lim_{x\to 0}\frac{e^x-1}{x}=1
 \label{5.10.3}
 \end{equation}
 \ecor
 \bpr Сразу получается \eqref{5.10.2}:
 \begin{multline*}
\lim_{x\to 0}\frac{\ln (1+x)}{x}= \lim_{x\to 0}\ln (1+x)^{\frac{1}{x}}=\\= \ln
\lim_{x\to 0} (1+x)^{\frac{1}{x}}=\ln e=1
 \end{multline*}
А отсюда уже следует \eqref{5.10.3}:
 \begin{multline*}
\lim_{x\to 0}\frac{e^x-1}{x}={\smsize \left| \begin{array}{c} e^x-1=t \\
x=\ln (1+t)\end{array}\right|}= \lim_{t\to 0}\frac{t}{\ln(1+t)}=\\=
\left(\lim_{t\to 0}\frac{\ln(1+t)}{t}\right)^{-1}=\eqref{5.10.2}= 1^{-1}=1
 \end{multline*}
 \epr

\begin{ers} Найдите пределы

1. $\lim\limits_{x\to 0}\frac{e^{2x}-1}{x}$;

2. $\lim\limits_{x\to 0}\frac{\ln (1-3x)}{x}$;

3. $\lim\limits_{x\to \infty} x\l e^{\frac{5}{x}}-1\r$;

4. $\lim\limits_{x\to \infty}\frac{1}{x}\ln \frac{x-1}{x+1}$;

5. $\lim\limits_{x\to -\infty}\frac{\ln (1-7e^x)}{e^x}$;

6. $\lim\limits_{x\to -\infty}\frac{1-e^{\sin x}}{1-\cos^2 x}$.
\end{ers}

\end{multicols}\noindent\rule[10pt]{160mm}{0.1pt}

\chapter{ПРОИЗВОДНАЯ}\label{ch-f'(x)}

Понятие производной появилось в анализе как инструмент для решения задач на
экстремум (то есть на максимум и минимум), поэтому обсуждение этой темы полезно
начать с какого-нибудь простого примера, иллюстрирующего эту связь.

\noindent\rule{160mm}{0.1pt}\begin{multicols}{2}

\bex\label{EX:f(x)=x^2-x^4} Пусть нам дана функция
\beq\label{f(x)=x^2-x^4}
  f(x)=x^2-x^4
\eeq
и требуется найти ее наибольшее значение на прямой $\R$.

Нулями этой функции (то есть решениями уравнения $f(x)=0$) являются числа $\{
-1,0,1 \}$, поэтому легко сообразить, что график $f$ будет выглядеть примерно
так:

%\picture{50pt}{-30pt}{99-1.pcx}

\vglue120pt \noindent Отсюда можно заключить, что $f$ действительно будет иметь
максимум в каких-то точках, лежащих на интервалах $(-1,0)$ и $(0,1)$, и нам
нужно придумать, как найти эти точки.

Чтобы это понять, нарисуем касательную к графику функции $f$ в какой-нибудь
точке $a\in\R$ (что такое касательная, мы надеемся, что читатель знает по
школе):
%\picture{50pt}{-30pt}{99-2.pcx}

\vglue120pt \noindent Как любая линейная функция (мы определили линейные функции формулой \eqref{DEF:lineynaya-functsiya}), касательная описывается
уравнением вида
$$
y=k\cdot x+b
$$
в котором $k$ -- угловой коэффициент, связанный с углом наклона $\ph$
касательной формулой
$$
k=\tg\ph
$$
Если точку $a$ сдвинуть в какое-нибудь другое место, то угловой коэффициент
$k$, вообще говоря, изменится (как и вся касательная). Таким образом, мы
получаем некую зависимость, или, точнее сказать, функцию
$$
a\mapsto k(a)
$$
которая каждой точке $a\in\R$ ставит в соответствие угловой коэффициент
$k=k(a)$ касательной к графику функции $f$ в точке $a$.

Заметим далее, что в точках максимума (а также, минимума) касательная к графику
функции $f$ горизонтальна, то есть угловой коэффициент касательной равен нулю
 \beq\label{k(a)=0}
  k(a)=0.
 \eeq
Поэтому если б мы знали функцию $a\mapsto k(a)$, то решив уравнение
\eqref{k(a)=0}, мы как раз смогли бы получить точки максимума.

Если теперь вообразить себя очень умным, настолько, что по формуле, задающей
функцию $x\mapsto f(x)$, ты всегда можешь определить, как устроена функция
$a\mapsto k(a)$, то останется совсем небольшой шаг, чтобы понять, что у
конкретной функции $f(x)=x^2-x^4$ соответствующая функция $k$ вычисляется по
формуле
 \beq\label{k(a)=2a-4a^3}
  k(a)=2a-4a^3.
 \eeq
Как можно достичь такого уровня интеллекта станет понятно после \ref{SEC-proizv-i-diff-functsii}\ref{SEC-proizv-kak-oper-nad-simv} настоящей главы (а конкретно для функции \eqref{f(x)=x^2-x^4} эта задача будет решена в примере \ref{EX:proizv-f(x)=x^2-x^4}), для нас же сейчас важно,
что, зная \eqref{k(a)=2a-4a^3}, мы легко можем решить уравнение \eqref{k(a)=0}:
 \begin{multline*}
k(a)=0 \quad \Leftrightarrow \quad 2a-4a^3=0 \quad \Leftrightarrow \\
\Leftrightarrow \quad a\in \left\{ -\frac{1}{\sqrt{2}}, 0,
\frac{1}{\sqrt{2}}\right\},
 \end{multline*}
После этого становится понятно, что максимум нашей функции достигается в точках
$\pm \frac{1}{\sqrt{2}}$, и равен
$$
f\l \pm
\frac{1}{\sqrt{2}}\r=\frac{1}{2}-\frac{1}{4}=\frac{1}{4}
$$

%\picture{50pt}{-30pt}{99-3.pcx}

\vglue120pt

\eex

\end{multicols}\noindent\rule[10pt]{160mm}{0.1pt}

Функция $a\mapsto k(a)$, о которой мы говорили в этом примере, на
математическом языке называется производной функции $f$ и обозначается
$$
k(a)=f'(a)
$$
Это одно из главных понятий в математическом анализе, и, помимо нахождения
экстремумов (максимумов и минимумов), оно используется в разных других задачах,
например,
 \bit{
\item[--] при построении графика функции (это описывается в \ref{SEC:grafik}
этой главы),

\item[--] при вычислении сложных пределов с помощью правила Лопиталя (см.
\ref{SEC:Lopital} этой главы),

\item[--]  при определении еще одного важного объекта в математическом анализе
-- {\it неопределенного интеграла} (об этом речь пойдет в главе
\ref{CH-indef-integral}).
 }\eit\noindent
В этой главе мы дадим аккуратное определение этому понятию, научимся вычислять
производные стандартных функций и обсудим приложения производной.

\section{Производная и ее вычисление}
\label{SEC-proizv-i-diff-functsii}

\subsection{Производная и дифференцируемость} \label{CH-opredelenie-proizvodnoi}

 \bit{
\item[$\bullet$] Пусть функция $f$ определена в некоторой окрестности точки $a$. Если существует конечный предел
\beq\label{DEF:f'(a)}
f'(a)=\lim_{t\to 0}\frac{f(a+t)-f(a)}{t}=\lim_{x\to a}\frac{f(x)-f(a)}{x-a}\in\R,
\eeq
то он называется {\it производной} функции $f$ в точке $a$\index{производная}, а про функцию $f$ говорят, что она {\it дифференцируема в точке}
\index{функция!дифференцируемая!в
точке} $a$.

\item[$\bullet$] Функция $f$ называется {\it
дифференцируемой}\index{функция!дифференцируемая!на
множестве} {\it на множестве} $E$, если она (определена в окрестности каждой точки этого множества и) дифференцируема
в каждой точке этого множества:
$$
\forall a\in E \qquad \exists f'(a)= \lim_{t\to 0}\frac{f(a+t)-f(a)}{t}\in \R
$$
 }\eit

\noindent\rule{160mm}{0.1pt}\begin{multicols}{2}

\bex Для функции
$$
f(x)=x^2
$$
по определению, получаем:
 \begin{multline*}
f'(a)=\lim_{t\to 0}\frac{(a+t)^2-a^2}{t}=\\= \lim_{t\to
0}\frac{a^2+2at+t^2-a^2}{t}= \lim_{t\to 0} (2a+t)=2 a
 \end{multline*}
То есть, в любой точке $a$ производная функции
$f(x)=x^2$ равна $2 a$:
$$
f'(a)=2 a, \qquad a\in \R
$$
Ясно, что смысл нашего утверждения не изменится, если
заменить $a$ на любую другую переменную, например, на $x$:
$$
f'(x)=2 x, \qquad x\in \R
$$
-- эта формула будет означать, что производная функции $f$
в любой точке $x$ равна $2x$.
\eex

\bex\label{EX:proizv-f(x)=x^2-x^4} Рассмотрим функцию из примера \ref{EX:f(x)=x^2-x^4}:
$$
f(x)=x^2-x^4
$$
опять по определению, получаем:
 \begin{multline*}
f'(a)=\lim_{t\to 0}\frac{\Big\{(a+t)^2-(a+t)^4\Big\}-\Big\{a^2-a^4\Big\}}{t}=\\=\eqref{binom-Newtona}=\\=
\lim_{t\to 0}\frac{1}{t}\cdot\Big\{\big(a^2+2at+t^2\big)-\\-\big(a^4+4a^3t+6a^2t^2+4at^3+t^4\big)-a^2+a^4\Big\}= \\=
\lim_{t\to 0}\frac{1}{t}\cdot\Big\{\big(2at+t^2\big)-\big(4a^3t+6a^2t^2+4at^3+t^4\big)\Big\}=\\=
\lim_{t\to 0}\Big\{\big(2a+t\big)-\big(4a^3+6a^2t+4at^2+t^3\big)\Big\}=\\=2a-4a^3
 \end{multline*}
Это совпадает с тем, что мы заявляли в примере \ref{EX:f(x)=x^2-x^4}:
$$
f'(a)=k(a)=2 a-4a^3, \qquad a\in \R
$$
Если поменять переменную на $x$, то
$$
f'(x)=2 x-4x^3, \qquad x\in \R.
$$
\eex

\bex Наоборот, функция
$$
f(x)=|x|
$$
не дифференцирема в точке $a=0$,  потому что предел
$$
f'(0)=\lim_{t\to 0}\frac{|t|}{t}
$$
не существует.
\eex

\end{multicols}\noindent\rule[10pt]{160mm}{0.1pt}

\begin{tm}[\bf о связи между непрерывностью и дифференцируемостью]
\label{cont-diff} Если функция $f$ дифференцируема в точке $a$, то она
непрерывна в этой точке.
\end{tm}\begin{proof} Если функция $f$ определена в
некоторой окрестности точки $a$ и дифференцируема в $a$, то это значит, что
существует конечный предел
$$
\lim_{t\to 0}\frac{f(a+t)-f(a)}{t}=f'(a)
$$
Отсюда следует, что существует предел
$$
\lim_{t\to 0}\left(f(a+t)-f(t)\right)= \lim_{t\to
0}\left(\frac{f(a+t)-f(a)}{t}\cdot t \right)= \lim_{t\to
0}\frac{f(a+t)-f(a)}{t}\cdot \lim_{t\to 0} t= f'(a)\cdot
0=0
$$
То есть,
$$
\lim_{x\to a} f(x)=\lim_{t\to 0} f(a+t)=f(a)
$$
а это как раз означает, что функция $f$ непрерывна в
точке $a$.
\end{proof}

\noindent\rule{160mm}{0.1pt}\begin{multicols}{2}

\paragraph*{Геометрический смысл производной.}

У производной имеется естественная геометрическая интерпретация: число $f'(a)$, определенное формулой \eqref{DEF:f'(a)} как раз и есть {\it угловой коэффициент касательной} к графику функции $f$ в точке $a$, о котором мы говорили в примере \ref{EX:f(x)=x^2-x^4}. Чтобы понять, почему это так, зафиксируем сначала число $t$ и проведем секущую к графику функции $f$ через точки $(a; f(a))$ и $(a+t; f(a+t))$:

%\picture{0pt}{0pt}{106.pcx}

\vglue110pt \noindent Из треугольника ABC видно, что
угловой коэффициент (то есть тангенс угла наклона) секущей
будет равен
$$
k_{t}=\frac{f(a+t)-f(a)}{t}
$$
Теперь ``расфиксируем'' число $t$ и устремим его к нулю.
Тогда наша секущая тоже станет двигаться, ``приближаясь к
касательной'', а ``в пределе'' -- просто превратится в
касательную.

%\picture{0pt}{0pt}{107.pcx}

\vglue110pt \noindent Угловой коэффициент (тангенс угла
наклона) этой касательной будет как раз равен производной:
$$
k=\lim_{t\to 0}\frac{f(a+t)-f(a)}{t}=f'(a)
$$

\paragraph*{Наглядный смысл дифферецируемости.} На
интуитивном уровне дифференцируемость функции $f$  в точке $a$ означает, что в
точке $a$ к графику функции $f$ можно провести (однозначно определенную)
касательную (и это не будет вертикальная прямая).

Для понимания этих терминов следует уяснить два важных
момента:
 \biter{
\item[--] во-первых, производная существует не всегда, и глядя на график обычно
бывает нетрудно сообразить, имеет ли функция $f$ в данной точке производную (то
есть, дифференцируема ли она в этой точке), или нет; например, функция
$f(x)=|\arctg x|$

%\picture{0pt}{0pt}{108.pcx}

\vglue120pt

\item[] дифференцируема в любой точке $x_0\ne 0$ (потому
что в любой точке $x_0\ne 0$ к графику можно провести
касательную)

%\picture{0pt}{0pt}{109.pcx}

\vglue120pt

\item[] но в точке $x_0=0$ она не дифференцируема (потому
что в $x_0=0$ однозначно определенную касательную провести
невозможно)

%\picture{0pt}{0pt}{110.pcx}

\vglue120pt

\item[--] во-вторых, по графику функции можно сообразить, в
каких точках производная больше, а в каких -- меньше;
например, у той же самой функции $f(x)=|\arctg x|$
производная в точке $x_0=1$ больше, чем в точке $x_1=2$
(потому что угол наклона касательной в точке $x_0=1$
больше, чем в точке $x_1=2$):

%\picture{0pt}{0pt}{111.pcx}

\vglue120pt

\item[] в частности, если производная в какой-то точке
$x_0$ равна нулю
$$
f'(x_0)=0
$$
то это означает, что касательная в этой точке $x_0$
горизонтальна

%\picture{0pt}{0pt}{112.pcx}

\vglue120pt

\item[] если же производная в $x_0$ положительна
$$
f'(x_0)>0
$$
то это означает, что касательная в этой точке $x_0$
возрастает

%\picture{0pt}{0pt}{113.pcx}

\vglue140pt

\item[] а если производная в $x_0$ отрицательна
$$
f'(x_0)<0
$$
то касательная в $x_0$ убывает

%\picture{0pt}{0pt}{114.pcx}

\vglue120pt
 }\eiter

В качестве упражнения полезно порешать задачи следующего
типа.

\begin{er}
Нарисуйте график какой-нибудь функции $f$ со следующими
свойствами:
 \bit{
\item[1)] $f(x)$ дифференцируема в любой точке, кроме
$x=-2$ и $x=1$; \item[2)] $f'(0)=0$, $f'(2)=-1$,
$f'(-1)=1$; \item[3)] $f(x)$ убывает на множестве $(0;
+\infty)$; \item[4)] $f(x)$ возрастает на множестве
$(-\infty; -3)$.
 }\eit
\end{er}

\begin{er}
Нарисуйте график какой-нибудь функции $f$ со следующими
свойствами:
 \bit{
\item[1)] $f(x)$ дифференцируема в любой точке, кроме
$x=0$; \item[2)] $f'(2)=0$, $f'(3)=1$, $f'(1)=-1$;
\item[3)] $f(x)$ ограничена на множестве $(0; +\infty)$;
\item[4)] $f(x)$ не ограничена на множестве $(-\infty; 0)$.
 }\eit
\end{er}

\begin{er}
По графику функции

%\picture{0pt}{0pt}{115.pcx}

\vglue160pt \noindent определите, в каких точках она
дифференцируема, и будет ли она дифференцируемой

1) на интервале $(-\infty; 1)$?

2) на интервале $(-\infty;0)$?

3) на интервале $(1;+\infty)$?

4) на интервале $(0;1)$?
\end{er}

\end{multicols}\noindent\rule[10pt]{160mm}{0.1pt}

\subsection{Производные элементарных функций.}\label{SUBSEC-proizv-elem-func}

Естественно теперь задаться вопросом, каковы производные других элементарных
функций. Они, конечно, все посчитаны (в тех случаях, когда это возможно), и
результаты можно привести в следующей серии примеров.

\noindent\rule{160mm}{0.1pt}\begin{multicols}{2}

\bex {\bf Производная линейной функции:}
 \beq\label{7.1.1}
f(x)=kx+b\quad\Longrightarrow\quad f'(x)=k
 \eeq
В частности, производная от константы:
 \beq\label{7.1.14}
f(x)=C\quad\Longrightarrow\quad f'(x)=0
 \eeq
Действительно,
 \begin{multline*}
f'(x)=\lim_{t\to 0}\frac{f(x+t)-f(x)}{t}=\\=\lim_{t\to
0}\frac{(k\cdot (x+t) + b)-(k\cdot x + b)}{t}
=\\=\lim_{t\to 0}\frac{k\cdot t}{t}=\lim_{t\to 0} k=k
 \end{multline*}
\eex

\bex{\bf Производная степенной функции:}
 \beq\label{7.1.2}
f(x)=x^\alpha\quad\Longrightarrow\quad f'(x)=\alpha\cdot
x^{\alpha-1}
 \eeq
В частности,
 \beq\label{7.1.15}
f(x)=\frac{1}{x}\quad\Longrightarrow\quad
f'(x)=-\frac{1}{x^2}
 \eeq
Действительно,
 \begin{multline*}
f'(x)=\lim_{t\to 0}\frac{(x+t)^\alpha-x^\alpha}{t}=\\=
\left|
\begin{array}{c} x+t=x e^s\\ t=x(e^s-1) \\{s=\ln (1+\frac{t}{x})}\end{array}\right|=\\=
\lim_{s\to 0}\frac{x^\alpha \cdot e^{\alpha
s}-x^\alpha}{x(e^s-1)}= \frac{x^\alpha}{x}\cdot \lim_{s\to
0}\frac{e^{\alpha s}-1}{e^s-1}=\\= x^{\alpha-1}\cdot
\lim_{s\to 0}\frac{e^{\alpha s}-1}{\alpha s}\frac{\alpha
s}{e^s-1}=\\= x^{\alpha-1}\cdot \lim_{s\to
0}\frac{e^{\alpha s}-1}{\alpha s}\lim_{s\to 0}\frac{\alpha
s}{e^s-1}=\\=
{\smsize {\smsize\begin{pmatrix}\text{в первом пределе}\\
\text{делаем замену}\; y=\alpha s \end{pmatrix}}}=\\=
x^{\alpha-1}\cdot \lim_{y\to 0}\frac{e^y-1}{y}\lim_{s\to
0}\frac{\alpha s}{e^s-1}=\\={\smsize
{\smsize\begin{pmatrix}\text{выносим $\alpha$ из второго предела}\\
\text{и преобразуем этот предел}\end{pmatrix}}}=\\= \alpha
x^{\alpha-1}\cdot \lim_{y\to
0}\frac{e^y-1}{y}\left(\lim_{s\to
0}\frac{e^s-1}{s}\right)^{-1}=\\= {\smsize
{\smsize\begin{pmatrix}\text{дважды применяем}\\
\text{формулу \eqref{5.10.3}}\end{pmatrix}}}=\\= \alpha
x^{\alpha-1}\cdot 1\cdot 1^{-1}=\alpha x^{\alpha-1}
 \end{multline*}
\eex

\bex\label{EX:(e^x)'} {\bf Производная показательной функции:}
 \beq\label{7.1.3}
f(x)=a^x\quad\Longrightarrow\quad f'(x)=a^x\cdot \ln a
 \eeq
В частности,
 \beq\label{7.1.16}
f(x)=e^x\quad\Longrightarrow\quad f'(x)=e^x
 \eeq
Действительно,
 \begin{multline*}
f'(x)=\lim_{t\to 0}\frac{a^{x+t}-a^x}{t}= a^x \cdot
\lim_{t\to 0}\frac{a^t-1}{t}=\\=
{\smsize {\smsize\begin{pmatrix}\text{применяем}\\
\text{формулу \eqref{5.10.1}}\end{pmatrix}}}= a^x \cdot
\lim_{t\to 0}\frac{e^{t \ln a}-1}{t}=\\={\smsize \left|
\begin{array}{c} s=t \ln a  \\ t=\frac{s}{\ln a}\end{array}\right|}=
a^x \ln a \cdot \lim_{s\to 0}\frac{e^s-1}{s}=\\={\smsize
{\smsize\begin{pmatrix}\text{применяем}\\
\text{формулу \eqref{5.10.3}}\end{pmatrix}}}=a^x \ln a
\cdot 1= a^x \ln a
\end{multline*}
\eex

\bex {\bf Производная логарифма:}
 \beq\label{7.1.4}
f(x)=\log_a x\quad\Longrightarrow\quad
f'(x)=\frac{1}{x\cdot \ln a}
 \eeq
В частности,
 \beq\label{7.1.17}
f(x)=\ln x\quad\Longrightarrow\quad f'(x)=\frac{1}{x}
 \eeq
Действительно,
 \begin{multline*}
f'(x)=\lim_{t\to 0}\frac{\log_a (x+t)-\log_a x}{t}=\\=
{\smsize
{\smsize\begin{pmatrix}\text{применяем}\\
\text{формулу \eqref{5.10.1}}\end{pmatrix}}}= \lim_{t\to
0}\frac{\ln (x+t)-\ln x}{t \ln a}=\\= \lim_{t\to
0}\frac{\ln \frac{x+t}{x}}{t \ln a}= \lim_{t\to 0}\frac{\ln
(1+\frac{t}{x})}{t \ln a}= {\smsize \left|
\begin{matrix}
s=\frac{t}{x}\\
t=s\cdot x
\end{matrix}\right|}=\\=
\lim_{s\to 0}\frac{\ln (1+s)}{s\cdot x \cdot \ln a}=
\frac{1}{x\cdot \ln a}\cdot \lim_{s\to 0}\frac{\ln
(1+s)}{s}=\\= {\smsize
{\smsize\begin{pmatrix}\text{применяем}\\
\text{формулу
\eqref{5.10.2}}\end{pmatrix}}}=\frac{1}{x\cdot \ln a}\cdot
1=\frac{1}{x\cdot \ln a}
 \end{multline*}
\eex

\bex {\bf Производная синуса:}
 \beq\label{7.1.5}
f(x)=\sin x\quad\Longrightarrow\quad f'(x)=\cos x
 \eeq
Для доказательства нам понадобится следующая
вспомогательная формула:
 \begin{equation}\label{7.1.13}
 \lim_{t\to 0}\frac{\cos t -1}{t}=0, \end{equation}
которую можно вывести из первого замечательного предела:
 \begin{multline*}
\lim_{t\to 0}\frac{\cos t -1}{t}= \lim_{t\to
0}\frac{-2\sin^2 \frac{t}{2}}{t}= \left| \begin{array}{c}
s=\frac{t}{2}\\ t=2 s
\end{array}\right|=\\= \lim_{s\to 0}\frac{-2\sin^2 s}{2 s}=
-\lim_{s\to 0}\frac{\sin s}{s}\cdot \lim_{s\to 0}\sin s = -
1\cdot 0 =0
 \end{multline*}
Мы можем теперь доказать \eqref{7.1.5}:
 \begin{multline*}
f'(x)=\lim_{t\to 0}\frac{\sin (x+t)-\sin x}{t}=\\=
\lim_{t\to 0}\frac{\sin x\cdot \cos t+\cos x\cdot \sin
t-\sin x}{t}=\\= \lim_{t\to 0}\frac{\sin x\cdot (\cos
t-1)+\cos x\cdot \sin t}{t}=\\= \lim_{t\to 0}\frac{\sin
x\cdot (\cos t-1)}{t}+ \lim_{t\to 0}\frac{\cos x\cdot \sin
t}{t}=\\= \sin x\cdot \lim_{t\to 0}\frac{\cos t-1}{t}+ \cos
x\cdot \lim_{t\to 0}\frac{\sin t}{t}=\\={\smsize
{\smsize\begin{pmatrix}\text{применяем формулу \eqref{7.1.13}}\\
\text{и первый замечательный предел}\end{pmatrix}}}=\\=
\sin x\cdot 0+ \cos x\cdot 1=\cos x
 \end{multline*}
\eex

\bex {\bf Производная косинуса:}
 \beq\label{7.1.6}
f(x)=\cos x\quad\Longrightarrow\quad f'(x)=-\sin x
 \eeq
Действительно,
 \begin{multline*}
f'(x)=\lim_{t\to 0}\frac{\cos (x+t)-\cos x}{t}=\\=
\lim_{t\to 0}\frac{\cos x\cdot \cos t-\sin x\cdot \sin
t-\cos x}{t}=\\= \lim_{t\to 0}\frac{\cos x\cdot (\cos
t-1)-\sin x\cdot \sin t}{t}=\\= \cos x\cdot \lim_{t\to
0}\frac{\cos t-1}{t}- \sin x\cdot \lim_{t\to
0}\frac{\sin t}{t}=\\={\smsize {\smsize\begin{pmatrix}\text{применяем формулу \eqref{7.1.13}}\\
\text{и первый замечательный предел}\end{pmatrix}}}=\\=
\cos x\cdot 0- \sin x\cdot 1=-\sin x
 \end{multline*}
 \eex

\bex {\bf Производная тангенса:}
 \beq\label{7.1.7}
f(x)=\tg x\quad\Longrightarrow\quad f'(x)=\frac{1}{\cos^2
x}
 \eeq
Действительно,
 \begin{multline*}
(\tg \boldsymbol{x})'(x)=\lim_{t\to 0}\frac{\tg (x+t)-\tg x}{t}=\\= \lim_{t\to
0}\frac{\frac{\sin (x+t)}{\cos (x+t)}-\frac{\sin x}{\cos x}}{t}=\\= \lim_{t\to
0}\frac{\sin (x+t)\cdot \cos x- \sin x\cdot \cos (x+t)} {t\cdot \cos (x+t)\cos
x}=\\= \lim_{t\to 0}\frac{\sin \left((x+t)-x\right)} {t\cdot \cos (x+t)\cos x}=
\lim_{t\to 0}\frac{\sin t} {t\cdot \cos (x+t)\cos x}=\\= \lim_{t\to
0}\frac{\sin t}{t}\cdot \lim_{t\to 0}\frac{1}{\cos (x+t)\cos x}= 1 \cdot
\frac{1}{\cos^2 x}= \frac{1}{\cos^2 x}\end{multline*} \eex

\bex {Производная котангенса:}
 \beq\label{7.1.8}
f(x)=\ctg x\quad\Longrightarrow\quad f'(x)=-\frac{1}{\sin^2
x}
 \eeq
Действительно,
 \begin{multline*}
f'(x)=\lim_{t\to 0}\frac{\ctg (x+t)-\ctg x}{t}=\\=
\lim_{t\to 0}\frac{\frac{\cos (x+t)}{\sin (x+t)}-\frac{\cos
x}{\sin x}}{t}=\\= \lim_{t\to 0}\frac{\cos (x+t)\cdot \sin
x- \cos x\cdot \sin (x+t)} {t\cdot \sin (x+t)\sin x}=\\=
\lim_{t\to 0}\frac{\sin \left(x-(x+t) \right)} {t\cdot \sin
(x+t)\sin x}= \lim_{t\to 0}\frac{-\sin t} {t\cdot \sin
(x+t)\sin x}=\\= - \lim_{t\to 0}\frac{\sin t}{t}\cdot
\lim_{t\to 0}\frac{1}{\sin (x+t)\sin x}=\\= - 1 \cdot
\frac{1}{\sin^2 x}= - \frac{1}{\sin^2 x}\end{multline*}
 \eex

\bex {\bf Производная арксинуса:}
 \beq\label{7.1.9}
f(x)=\arcsin x\quad\Longrightarrow\quad
f'(x)=\frac{1}{\sqrt{1-x^2}}
 \eeq
Действительно,
 \begin{multline*}
(\arcsin \boldsymbol{x})'(x)=\lim_{t\to 0}\frac{\arcsin
(x+t)-\arcsin x}{t}=\\= {\smsize \left|
\begin{array}{c}
s=\arcsin (x+t)-\arcsin x
\\
t=\sin s\cdot \cos(\arcsin x)+x\cdot (\cos s -1)
\end{array}\right|}=\\=
\lim_{s\to 0}\frac{s}{\sin s\cdot \cos(\arcsin x)+x\cdot
(\cos s
-1)}=\\={\smsize {\smsize\begin{pmatrix}\text{используем формулу}\\
\cos(\arcsin x)=\sqrt{1-x^2}\end{pmatrix}}}=\\= \lim_{s\to
0}\frac{s}{\sin s\cdot \sqrt{1-x^2}+x\cdot (\cos s -1)}=\\=
\lim_{s\to 0}\frac{1}{\frac{\sin s}{s}\cdot \sqrt{1-x^2}+
x\cdot \frac{\cos s -1}{s}}=\\= \frac{1}{\lim_{s\to
0}\frac{\sin s}{s}\cdot \sqrt{1-x^2}+x\cdot \lim_{s\to
0}\frac{\cos s -1}{s}}=\\= {\smsize
{\smsize\begin{pmatrix}\text{используем первый}\\
\text{замечательный предел}\\
\text{и формулу \eqref{7.1.13}}\end{pmatrix}}}=\\=
\frac{1}{1\cdot \sqrt{1-x^2}+x\cdot 0}=
\frac{1}{\sqrt{1-x^2}}\end{multline*}
 \eex

 \bex {\bf Производная арккосинуса:}
 \beq\label{7.1.10}
f(x)=\arccos x\quad\Longrightarrow\quad
f'(x)=-\frac{1}{\sqrt{1-x^2}}
 \eeq
Действительно,
 \begin{multline*}
f'(x)=\lim_{t\to 0}\frac{\arccos (x+t)-\arccos x}{t}=\\=
{\smsize \left|
\begin{array}{c}
s=\arccos (x+t)-\arccos x
\\
t=(\cos s -1)\cdot x-\sin s\cdot \sin(\arccos x)
\end{array}\right|}=\\=
\lim_{s\to 0}\frac{s}{(\cos s -1)\cdot x-\sin s\cdot
\sin(\arccos
x)}=\\= {\smsize {\smsize\begin{pmatrix}\text{используем формулу}\\
\sin(\arccos x)=\sqrt{1-x^2}\end{pmatrix}}}=\\= \lim_{s\to
0}\frac{s}{(\cos s -1)\cdot x-\sin s\cdot \sqrt{1-x^2}}=\\=
\lim_{s\to 0}\frac{1}{\frac{\cos s -1}{s}\cdot x-
\frac{\sin s}{s}\cdot \sqrt{1-x^2}}=\\= \frac{1}{\lim_{s\to
0}\frac{\cos s -1}{s}\cdot x- \lim_{s\to 0}\frac{\sin
s}{s}\cdot
\sqrt{1-x^2}}=\\={\smsize {\smsize\begin{pmatrix}\text{используем первый}\\
\text{замечательный предел}\\
\text{и формулу \eqref{7.1.13}}\end{pmatrix}}}=\\=
\frac{1}{0\cdot x- 1\cdot \sqrt{1-x^2}}=
-\frac{1}{\sqrt{1-x^2}}\end{multline*}
 \eex

\bex {\bf Производная арктангенса:}
 \beq\label{7.1.11}
f(x)=\arctg x\quad\Longrightarrow\quad
f'(x)=\frac{1}{1+x^2}
 \eeq
Действительно,
 \begin{multline*}
f'(x)=\lim_{t\to 0}\frac{\arctg (x+t)-\arctg x}{t}=\\=
{\smsize \left|
\begin{array}{c}
s=\arctg (x+t)-\arctg x
\\
t=\frac{(1+x^2)\cdot \tg s}{1-x \tg s}\end{array}\right|}=
\lim_{s\to 0}\frac{s (1-x \tg s)}{\tg s(1+x^2)}=\\=
\lim_{s\to 0}\frac{s \cos s (1-x \tg s)}{\sin s(1+x^2)}=
\lim_{s\to 0}\frac{\cos s (1-x \tg s)}{\frac{\sin
s}{s}\cdot (1+x^2)}=\\= \frac{\lim_{s\to 0}\cos s \cdot
(1-x \cdot \lim_{s\to 0}\tg s)} {\lim_{s\to 0}\frac{\sin
s}{s}\cdot (1+x^2)}=\\= \frac{1\cdot (1-x \cdot 0)} {1
\cdot (1+x^2)}= \frac{1}{1+x^2}\end{multline*}
 \eex

 \bex {\bf Производная арккотангенса:}
 \beq\label{7.1.12}
f(x)=\arcctg x\quad\Longrightarrow\quad
f'(x)=-\frac{1}{1+x^2}
 \eeq
Действительно,
 \begin{multline*}
f'(x)=\lim_{t\to 0}\frac{\arcctg (x+t)-\arcctg x}{t}=\\=
{\smsize \left|
\begin{array}{c}
s=\arcctg (x+t)-\arcctg x
\\
t=-\frac{1+x^2}{x+ \ctg s}\end{array}\right|}=\\=
\lim_{s\to 0}\frac{-s\cdot (x+\ctg s)}{1+x^2}=\\=
-\lim_{s\to 0}\frac{s\cdot (x\cdot \sin s+\cos
s)}{(1+x^2)\cdot \sin s}=\\= -\lim_{s\to 0}\frac{x\cdot
\sin s+\cos s}{(1+x^2)\cdot \frac{\sin s}{s}}=\\= -
\frac{x\cdot \lim_{s\to 0}\sin s+\lim_{s\to 0}\cos s}
{(1+x^2)\cdot \lim_{s\to 0}\frac{\sin s}{s}}=\\= -
\frac{x\cdot 0+1} {(1+x^2)\cdot 1}= -
\frac{1}{1+x^2}\end{multline*}
 \eex

Из этих примеров видно, что (как и в случае с непрерывностью, о чем мы говорили в теореме \ref{TH-nepr-elem-func}), не все
элементарные функции дифференцируемы на своей области определения:

\btm\label{TH-diff-elem-func} Среди элементарных функций все
непрерывные\footnote{Непрерывные элементарные функции были определены нами на
с.\pageref{DEF:nepr-elem-funktsii}.} (и только они) дифференцируемы на всяком
интервале в своей области определения.\etm

\end{multicols}\noindent\rule[10pt]{160mm}{0.1pt}

\subsection{Правила вычисления производных}\label{SEC-pravila-proizvodnoi}

\paragraph{Арифметические действия с производными.}

\begin{tm}\label{diff-alg}
Если функции $f$ и $g$ дифференцируемы в точке $x$ то их сумма, разность и
произведение тоже дифференцируемы в точке $x$, причем
 \begin{align}
& \Big( C\cdot f\Big)'(x)=C\cdot f'(x),\qquad C\in\R
\label{7.3.1}
 \\
& \Big( f+g\Big)'(x)=f'(x)+g'(x) \label{7.3.2}
 \\
& \Big( f-g\Big)'(x)=f'(x)-g'(x) \label{7.3.3}
\\
& \Big( f\cdot g\Big)'(x)=f'(x)\cdot g(x)+ f(x)\cdot g'(x)
\label{7.3.4}
 \end{align}
а если, кроме того $g(x)\ne 0$ то функция $\frac{f}{g}$
тоже дифференцируема в точке $x$, причем
 \begin{equation}
\l\frac{f}{g}\r'(x)= \frac{f'(x)\cdot g(x)-f(x)\cdot
g'(x)}{g^2(x)}\label{7.3.5}
 \end{equation}\end{tm}

\begin{proof}
Докажем \eqref{7.3.2}:

\begin{multline*}\Big(f(x)+g(x)\Big)'= \lim_{t\to 0}\frac{\Big(f(x+t)+g(x+t)\Big)-\Big(f(x)+g(x)\Big)}{t}=\\=
\lim_{t\to
0}\frac{\Big(f(x+t)-f(x)\Big)+\Big(g(x+t)-g(x)\Big)}{t}=
\lim_{t\to 0}\frac{f(x+t)-f(x)}{t}+\lim_{t\to
0}\frac{g(x+t)-g(x)}{t}=\\= f'(x)+g'(x) \end{multline*}
Аналогично доказывается \eqref{7.3.3}:

\begin{multline*}\Big(f(x)-g(x)\Big)'= \lim_{t\to 0}\frac{\Big(f(x+t)-g(x+t)\Big)-\Big(f(x)-g(x)\Big)}{t}=\\=
\lim_{t\to
0}\frac{\Big(f(x+t)-f(x)\Big)-\Big(g(x+t)-g(x)\Big)}{t}=
\lim_{t\to 0}\frac{f(x+t)-f(x)}{t}-\lim_{t\to
0}\frac{g(x+t)-g(x)}{t}=\\= f'(x)-g'(x) \end{multline*} Для
\eqref{7.3.4} получаем:

\begin{multline*}\Big(f(x)\cdot g(x)\Big)'= \lim_{t\to 0}\frac{f(x+t)\cdot
g(x+t)-f(x)\cdot g(x)}{t}= {\smsize\begin{pmatrix}\text{вычтем и прибавим}\\
\text{слагаемое}\,\, f(x)\cdot g(x+t)
\end{pmatrix}}=\\= \lim_{t\to 0}\frac{f(x+t)\cdot g(x+t) -f(x)\cdot
g(x+t)+f(x)\cdot g(x+t)-f(x)\cdot g(x)}{t}=\\= \lim_{t\to
0}\Big(\frac{f(x+t)\cdot g(x+t) -f(x)\cdot
g(x+t)}{t}+\frac{f(x)\cdot g(x+t)-f(x)\cdot
g(x)}{t}\Big)=\\= \lim_{t\to 0}\frac{\Big( f(x+t) - f(x)
\Big)\cdot g(x+t)}{t}+ \lim_{t\to 0}\frac{f(x)\cdot \Big(
g(x+t)-g(x)\Big)}{t}=\\= \lim_{t\to 0}\frac{\Big( f(x+t) -
f(x) \Big)}{t}\cdot \lim_{t\to 0} g(x+t)+ f(x)\cdot
\lim_{t\to 0}\frac{\Big(
g(x+t)-g(x)\Big)}{t}=\\= {\smsize\begin{pmatrix}\text{по теореме \ref{cont-diff}, функция}\\
\text{$g(x)$ непрерывна, значит}\\
\lim_{t\to 0} g(x+t)=g(x)
\end{pmatrix}}= f'(x)\cdot g(x)+ f(x)\cdot g'(x) \end{multline*}
И, наконец, для \eqref{7.3.5}:

\begin{multline*}\Big(\frac{f(x)}{g(x)}\Big)'= \lim_{t\to 0}\frac{\frac{f(x+t)}{g(x+t)}-\frac{f(x)}{g(x)}}{t}= \lim_{t\to 0}\frac{f(x+t)\cdot g(x)-f(x)\cdot g(x+t)}{t\cdot g(x+t)\cdot g(x)}=
{\smsize\begin{pmatrix}\text{вычтем и прибавим}\\
\text{слагаемое}\,\, f(x)\cdot g(x)
\end{pmatrix}}=\\= \lim_{t\to 0}\frac{f(x+t)\cdot g(x)-f(x)\cdot g(x)
+f(x)\cdot g(x)-f(x)\cdot g(x+t)}{t\cdot g(x+t)\cdot
g(x)}=\\= \lim_{t\to 0}\Big(\frac{f(x+t)\cdot
g(x)-f(x)\cdot g(x)}{t\cdot g(x+t)\cdot g(x)}
+\frac{f(x)\cdot g(x)-f(x)\cdot g(x+t)}{t\cdot g(x+t)\cdot
g(x)}\Big)=\\= \lim_{t\to
0}\left(\Big(\frac{f(x+t)-f(x)}{t}\cdot g(x) - f(x)\cdot
\frac{g(x+t)-g(x)}{t}\Big) \cdot \frac{1}{g(x+t)\cdot
g(x)}\right) =\\= \Big(\lim_{t\to
0}\frac{f(x+t)-f(x)}{t}\cdot g(x) - f(x)\cdot
\lim_{t\to 0}\frac{g(x+t)-g(x)}{t}\Big) \cdot \lim_{t\to 0}\frac{1}{g(x+t)\cdot g(x)} =\\= {\smsize\begin{pmatrix}\text{по теореме \ref{cont-diff}, функция}\\
\text{$g(x)$ непрерывна, значит}\\
\lim_{t\to 0} g(x+t)=g(x)
\end{pmatrix}}= \Big( f'(x)\cdot g(x) - f(x)\cdot g'(x) \Big) \cdot
\frac{1}{g(x)^2} = \frac{f'(x)\cdot g(x)- f(x)\cdot
g'(x)}{g(x)^2}\end{multline*}\end{proof}

\noindent\rule{160mm}{0.1pt}\begin{multicols}{2}

\bex\label{EX:(ln-x-sin-x)'} Вычислим производную функции
$$
h(x)=\ln x\cdot\sin x
$$
Для этого введем обозначения
$$
f(x)=\ln x,\qquad g(x)=\sin x
$$
Тогда $h=f\cdot g$, и по формуле \eqref{7.3.4} получаем:
$$
h'(x)=f'(x)\cdot g(x)+f(x)\cdot g(x)=\frac{1}{x}\cdot\sin x
+\ln x\cdot\cos x
$$
\eex

\bex Чтобы найти производную функции
$$
h(x)=\frac{e^x}{\arctg x}
$$
введем обозначения
$$
f(x)=e^x,\qquad g(x)=\arctg x,
$$
Тогда $h=\frac{f}{g}$, и по формуле \eqref{7.3.4} получаем:
 \begin{multline*}
h'(x)=\frac{f'(x)\cdot g(x)-f(x)\cdot g(x)}{g(x)^2}=\\=
\frac{e^x\cdot\arctg x-e^x\cdot\frac{1}{1+x^2}}{\arctg^2 x}
 \end{multline*}
\eex

\end{multicols}\noindent\rule[10pt]{160mm}{0.1pt}

\paragraph{Производная композиции.}

Напомним, что {\it композиция $g\circ f$ функций} $g$ и $f$ (или, что то же
самое, {\it сложная функция}, составленная из $f$ и $g$) была определена нами
на с.\pageref{DEF:kompozitsiya-funktsij} формулой
$$
(g\circ f)(x)=g(f(x)).
$$

\begin{tm}\label{proizv-slozhn-functsii}
Пусть функция $f$ дифференцирема в точке $x=a$, а функция $g$ дифференцируема в
точке $f(a)$. Тогда сложная функция
$$
h(x)=g(f(x))
$$
дифференцируема в точке $x=a$, причем
 \begin{equation}\label{7.4.1}
h'(a)=g'(f(a))\cdot f'(a)
 \end{equation}
 \end{tm}

\brem Формулу \eqref{7.4.1} удобно записывать в виде
 \beq\label{proizvodnaya-kompozitsii}
(g\circ f)'=(g'\circ f)\cdot f'
 \eeq

\erem

 \begin{proof} Перепишем равенство
\eqref{7.4.1} по-другому:
$$
\lim_{x\to a}\frac{g(f(x))-g(f(a))}{x-a}=g'(f(a))\cdot
f'(a)
$$
Чтобы это доказать, нужно взять произвольную последовательность
\begin{equation}
x_n\underset{n\to \infty}{\longrightarrow} a, \quad x_n\ne
a \label{7.4.2}\end{equation} и убедиться, что
\begin{equation}\lim_{n\to \infty}\frac{g(f(x_n))-g(f(a))}{x_n-a}=g'(f(a))\cdot
f'(a) \label{7.4.3}\end{equation} Для этого нужно
рассмотреть два случая:
 \bit{
\item[--] возможна ситуация, когда для почти всех $n\in \mathbb{N}$
выполняется $f(x_n)\ne f(a)$;

\item[--] если же это не так, то для некоторой подпоследовательности аргументов
$x_{n_k}$ соответствующая последовательность значений $f(x_{n_k})$ будет
обладать свойством $f(x_{n_k})=f(a)$.
 }\eit

1. Рассмотрим сначала первый случай, то есть когда для
почти всех $n\in \mathbb{N}$ выполняется $f(x_n)\ne f(a)$.
Тогда можно делить и умножать на $f(x_n)-f(a)\ne 0$, и мы
получим

\begin{multline*}\lim_{n\to \infty}\frac{g(f(x_n))-g(f(a))}{x_n-a}= \lim_{n\to
\infty}\frac{g(f(x_n))-g(f(a))}{f(x_n)-f(a)}\cdot
\frac{f(x_n)-f(a)}{x_n-a}=\\= \lim_{n\to
\infty}\frac{g(f(x_n))-g(f(a))}{f(x_n)-f(a)}\cdot
\lim_{n\to \infty}\frac{f(x_n)-f(a)}{x_n-a}=
{\smsize\begin{pmatrix}\text{в левом пределе делаем}\\
\text{замену}\,\, y_n=f(x_n)
\end{pmatrix}}=\\= \lim_{n\to \infty}\frac{g(y_n)-g(f(a))}{y_n-f(a)}\cdot \lim_{n\to \infty}\frac{f(x_n)-f(a)}{x_n-a}= g'(f(a))\cdot f'(a) \end{multline*}

2. После этого рассмотрим случай, когда для некоторой подпоследовательности
аргументов $x_{n_k}$ соответствующая последовательность значений $f(x_{n_k})$
обладает свойством $f(x_{n_k})=f(a)$. Тогда последовательность $x_n$ можно
разбить на две подпоследовательности $x_{n_k}$ и $x_{n_i}$, такие что
$$
f(x_{n_k})=f(a), \qquad f(x_{n_i})\ne f(a)
$$
Для подпоследовательности $x_{n_i}$ мы получаем то же самое, что и в предыдущем
случае:
$$
\lim_{i\to
\infty}\frac{g(f(x_{n_i}))-g(f(a))}{x_{n_i}-a}=g'(f(a))\cdot
f'(a)
$$
А для $x_{n_k}$ получаем:
$$
\lim_{k\to \infty}\frac{g(f(x_{n_k}))-g(f(a))}{x_{n_k}-a}=
{\smsize\begin{pmatrix}\text{вспоминаем, что}\\
f(x_{n_k})=f(a)
\end{pmatrix}}= \lim_{k\to \infty}\frac{g(f(a))-g(f(a))}{x_{n_k}-a}=
\lim_{k\to \infty}\frac{0}{x_{n_k}-a}=0
$$
А с другой стороны,
$$
f'(a)=\lim_{k\to \infty}\frac{f(x_{n_k})-f(a)}{x_{n_k}-a}=
{\smsize\begin{pmatrix}\text{вспоминаем, что}\\
f(x_{n_k})=f(a)
\end{pmatrix}}= \lim_{k\to \infty}\frac{0}{x_{n_k}-a}=0
$$
И, таким образом,
$$
\lim_{k\to \infty}\frac{g(f(x_{n_k}))-g(f(a))}{x_{n_k}-a}
=0=g'(f(a))\cdot f'(a)
$$

Мы получили, что для любой последовательности \eqref{7.4.2} выполняется
равенство \eqref{7.4.3}, а это нам и нужно было доказать.
\end{proof}

\noindent\rule{160mm}{0.1pt}\begin{multicols}{2}

\begin{ex}\label{primery-vych-proizv-slozhnoi-f}
Вычислим производную функции $h(x)=\sin (x^2)$. Для этого
представим ее как композицию двух элементарных функций:
$$
h(x)=\sin (x^2)=g(f(x)),\ g(y)=\sin y,\ f(x)=x^2
$$
Тогда
$$
g'(y)=\cos y,\quad f'(x)=2x
$$
Применяя формулу \eqref{7.4.4}, получаем:
 $$
h'(x)=g'(f( x))\cdot f'( x)= \cos (x^2)\cdot 2x
 $$
 \end{ex}

\begin{ex}
Чтобы посчитать производную функции
$h(x)=\sin^2 x$, нужно представить ее в виде композиции
элементарных функций
$$
h(x)=\sin^2 x=g(f(x)),\ g(y)=y^2,\ f(x)=\sin x
$$
Тогда
$$
g'(y)=2y,\quad f'(x)=\cos x
$$
и, применяя \eqref{7.4.4}, получаем:
 $$
h'(x)=g'(f( x))\cdot f'( x)= 2\sin
 x\cdot \cos  x
 $$
 \end{ex}

\bex
Вычислим производную функции $h(x)=e^{\arctg  x}$. Для этого
представим ее как композицию двух элементарных функций:
$$
h(x)=e^{\arctg x}=g(f(x)),\ g(y)=e^y,\ f(x)=\arctg x
$$
Тогда
$$
g'(y)=e^y,\quad f'(x)=\frac{1}{1+x^2}
$$
Применяя формулу \eqref{7.4.4}, получаем:
 $$
h'(x)=g'(f( x))\cdot f'( x)= e^{\arctg  x}\cdot \frac{1}{1+
x^2}
 $$
\eex

Одним из следствий теоремы о производной сложной
функции является возможность вычислять производные функций вида $f^{g}$. Это делается с помощью следующей очевидной формулы:
 \beq\label{EQ:f^g}
f^{g}=e^{\ln f^{g}}=e^{g \cdot \ln f}
 \eeq
(справедливой, разумеется, только при $f>0$).

 \begin{ex}\label{EX-(x^x)'} Вычислим производную функции
 $$
h(x)=x^x, \quad x>0
 $$
Для этого представим нашу функцию, как композицию двух функций с помощью формулы \eqref{EQ:f^g}:
$$
h(x)=x^x=e^{\ln x^x}=e^{x\cdot \ln x}
$$
Тогда
$$
h(x)=g(f(x)),\qquad g(y)=e^y,\quad f(x)=x\cdot\ln x
$$
Производная функции $f$ вычисляется по формуле \eqref{proizvodnaya-kompozitsii} и равна
$$
f'(x)=\ln x+1
$$
а для функции $g$ она посчитана в примере \ref{EX:(e^x)'}:
$$
g'(y)=e^y
$$
Поэтому
 \begin{multline*}
h'(x)=g'(f(x))\cdot f'(x)=\\=e^{x\cdot \ln x}\cdot\Big(\ln x+1\Big)= x ^{ x }\cdot \Big(\ln x +1\Big).
 \end{multline*}
\end{ex}
\end{multicols}\noindent\rule[10pt]{160mm}{0.1pt}

\subsection{Дифференциальное исчисление: производная, как формальная операция}
\label{SEC-proizv-kak-oper-nad-simv}

В примерах на с.\pageref{EX:(ln-x-sin-x)'} и с.\pageref{primery-vych-proizv-slozhnoi-f}, как наверняка заметил читатель, для
вычисления производной стандартной функции $h$ нам приходилось вводить новые
обозначения для <<более простых>> функций, через которые $h$ выражается
алгебраически или в виде композиции. Например, чтобы вычислить производную функции
$h(x)=\sin (x^2)$ в примере \ref{primery-vych-proizv-slozhnoi-f} нам пришлось
рассмотреть функции $g(y)=\sin y$ и $f(x)=x^2$, выписать их производные и
подставить их в формулу \eqref{7.4.1}.

Этот способ вычислений выглядит довольно громоздким в сравнении с приемами, употреблявшимися для таких задач в школе, где, если помнит читатель, вычисление производной оформляется в виде цепочки равенств (иногда длинной цепочки, но все же гораздо короче, чем запись, которую приходится употреблять, следуя алгоритму действий, описанному на страницах \pageref{EX:(ln-x-sin-x)'} и \pageref{primery-vych-proizv-slozhnoi-f}). Например, для функции $h(x)=\sin (x^2)$ можно (используя обозначение \eqref{podstanovka-v-vyrazhenie}) записать эту цепочку так:
 \beq\label{sin(x^2)'}
\Big(\sin (x^2)\Big)'=\Big(\sin y\Big|_{y=x^2}\Big)'=(\sin y)'\Big|_{y=x^2}\cdot (x^2)'=\cos x^2\cdot 2x
 \eeq
Естественно поэтому задуматься, отчего школьные приемы приводят к правильным результатам, и как их следует аккуратно сформулировать, чтобы грамотно использовать для вычислений?

Раздел математического анализа, занимающийся этими вопросами, называется {\it исчислением}, и мы упоминали уже о нем в начале \ref{SUBSEC-elem-funktsii} главы \ref{ch-ELEM-FUNCTIONS}. Его идея состоит в том, чтобы поглядеть на производную, как на формальную
операцию над числовыми выражениями (о которых мы вели речь в \ref{Chislovye-termy-i-stand-func} той же главы \ref{ch-ELEM-FUNCTIONS}).
Такой отстраненный взгляд позволяет формализовать приемы вычисления производной таким образом, что цепочки вида \eqref{sin(x^2)'} становятся результатом применения одного общего алгоритма, обоснованность которого отдельно доказывается, и поэтому не вызывает сомнений.

В этом пункте мы опишем этот алгоритм (и докажем правомерность его применения), чтобы иметь формальные основания пользоваться им в примерах.

\noindent\rule{160mm}{0.1pt}\begin{multicols}{2}

\paragraph{Производная числового выражения.}
\label{subsec-diff-chisl-terma}

Всякому числовому выражению $\mathcal P$ от переменной $x$, можно поставить в соответствие, причем единственным с точностью до равенства выражений способом, некое новое числовое выражение, обозначаемое $\frac{\d}{\d x}({\mathcal P})$ или $\frac{\d ({\mathcal P})}{\d x}$ (если это не приводит к недоразумениям, скобки вокруг ${\mathcal P}$ могут не писаться) и называемое {\it производным выражением} таким образом, чтобы выполнялись следующие
условия:
 \biter{
\item[0)] если $a$  --- параметр, то по определению считается, что
 \beq\label{da/dx}
\frac{\d a}{\d x}=0,
 \eeq
и это равенство выполняется всякий раз, когда в числителе под дифференциалом
(то есть после символа $\d$) стоит параметр, а в знаменателе --- переменная,

\item[1)] если $a$  --- собственное обозначение числа, то опять же по
определению считается, что
 \beq\label{d const/dx}
\frac{\d a}{\d x}=0,
 \eeq

\item[2)] на элементарных выражениях операция дифференцирования $\frac{\d}{\d
x}$ определяется формулами, представляющими из себя сокращенную запись
утверждений \eqref{7.1.1}-\eqref{7.1.12} (здесь всюду $x$ --- единственная
переменная):

 \bigskip
 \centerline{\bf Таблица производных:}
 \label{tablitsa-proizvodnyh}

\begin{align*}
 &\frac{\d}{\d x}\Big( C\Big)=k
 \\
 &\frac{\d}{\d x}\Big(x^\alpha \Big)=\alpha \cdot x^{\alpha-1}
 \\
 &\frac{\d}{\d x}\Big(a^{x} \Big)=\ln a \cdot a^{x}\quad(a>0)
 \\
 &\frac{\d}{\d x}\Big(e^{x} \Big)=e^{x}
 \\
 &\frac{\d}{\d x}\Big(\log_a x \Big)=\frac{1}{x\cdot \ln a}\quad (x>0)
 \\
  &\frac{\d}{\d x}\Big(\ln x \Big)=\frac{1}{x}\quad (x>0)
 \\
 &\frac{\d}{\d x}\Big(\sin x \Big)= \cos x
 \\
 &\frac{\d}{\d x}\Big(\cos x \Big)= -\sin x
 \\
 &\frac{\d}{\d x}\Big(\tg x \Big)=\frac{1}{\cos^2 x}
 \\
 &\frac{\d}{\d x}\Big(\ctg x \Big)=-\frac{1}{\sin^2 x}
 \\
 &\frac{\d}{\d x}\Big(\arcsin x \Big)=\frac{1}{\sqrt{1-x^2}}
 \\
 &\frac{\d}{\d x}\Big(\arccos x \Big)=-\frac{1}{\sqrt{1-x^2}}
 \\
 &\frac{\d}{\d x}\Big(\arctg x \Big)=\frac{1}{1+x^2}
 \\
 &\frac{\d}{\d x}\Big(\arcctg x \Big)=-\frac{1}{1+x^2}
 \end{align*}\medskip

\item[3)] для произвольных выражений $\mathcal P$ и $\mathcal Q$
 \begin{align}
&\kern-25pt\frac{\d}{\d x}\Big({\mathcal P}+{\mathcal Q}\Big)
=\frac{\d}{\d x}{\mathcal P}+\frac{\d}{\d x}{\mathcal Q} \label{DEF:d(P+Q)/dx}\\
&\kern-25pt\frac{\d}{\d x}\Big({\mathcal P}-{\mathcal Q}\Big)
=\frac{\d}{\d x}{\mathcal P}-\frac{\d}{\d x}{\mathcal Q} \label{DEF:d(P-Q)/dx} \\
&\kern-25pt\frac{\d}{\d x}\Big({\mathcal P}\cdot{\mathcal Q}\Big) =\left(\frac{\d}{\d
x}{\mathcal P}\right)\cdot{\mathcal Q}+{\mathcal P}\cdot\left(\frac{\d}{\d
x}{\mathcal
Q}\right) \label{DEF:d(P-cdot-Q)/dx} \\
&\kern-25pt\frac{\d}{\d x}\l\frac{{\mathcal P}}{{\mathcal Q}}\r =
\frac{\Big(\frac{\d{\mathcal P}}{\d x}\Big)\cdot {\mathcal Q} -{\mathcal
P}\cdot \Big(\frac{\d{\mathcal Q}}{\d x}\Big)}{{\mathcal Q}^2} \label{DEF:d(P/Q)/dx}
 \end{align}
в частности, если $\mathcal P$ не содержит переменной $x$, то
 \begin{align*}
&\frac{\d}{\d x}\Big({\mathcal P}\cdot{\mathcal Q}\Big) ={\mathcal
P}\cdot\left(\frac{\d}{\d x}{\mathcal Q}\right)
 \end{align*}

\item[4)] если ${\mathcal P}$ --- одноместное выражение от переменной $x$, а
${\mathcal Q}$ --- одноместное выражение от переменной $y$, то для выражения,
полученного подстановкой, выполняется равенство
 \begin{align}\label{7.4.4}
&\frac{\d}{\d x}\left({\mathcal Q}\Big|_{y={\mathcal P}}\right)=
\left(\frac{\d}{\d y}{\mathcal Q}\right)\Big|_{y={\mathcal P}}\cdot
\frac{\d}{\d x}{\mathcal P}
 \end{align}

 }\eiter

\brem К этим правилам следует добавить, что операции дифференцирования
выражений по другим переменным определяются по аналогии. Поэтому, если,
например, ${\mathcal P}$ -- одноместное выражение от переменной $t$, то
производную $\frac{\d}{\d t}{\mathcal P}$ можно вычислить, заменив сначала $t$
на $x$, затем вычислив $\frac{\d}{\d t}{\mathcal P}$ по уже сформулированному
алгоритму, а после этого снова заменив $x$ на $t$:
$$
\frac{\d}{\d t}{\mathcal P}=\left[\frac{\d}{\d x}\left({\mathcal
P}\Big|_{t=x}\right)\right]\Bigg|_{x=t}
$$
\erem

\brem Из формул \eqref{DEF:d(P+Q)/dx}-\eqref{7.4.4} и предложения \ref{PROP:D(P+Q),...} следует, что область допустимых значений переменной для производного выражения подчиняется равенствам:
\begin{align}
& \D\l\frac{\d({\mathcal P}+{\mathcal Q})}{\d x}\r=\D\l\frac{\d{\mathcal P}}{\d x}\r\cap\D\l\frac{\d{\mathcal Q}}{\d x}\r \label{D(d(P+Q)/dx)}\\
& \D\l\frac{\d({\mathcal P}-{\mathcal Q})}{\d x}\r=\D\l\frac{\d{\mathcal P}}{\d x}\r\cap\D\l\frac{\d{\mathcal Q}}{\d x}\r \label{D(d(P-Q)/dx)}\\
& \D\l\frac{\d({\mathcal P}\cdot{\mathcal Q})}{\d x}\r=\D\l\frac{\d{\mathcal P}}{\d x}\r\cap\D\l\frac{\d{\mathcal Q}}{\d x}\r \label{D(d(P-cdot-Q)/dx)}\\
& \begin{array}{l}\kern-3pt \D\l\frac{\d(\frac{\mathcal P}{\mathcal Q})}{\d x}\r=\\
\kern20pt=\D\l\frac{\d{\mathcal P}}{\d x}\r\cap\D\l\frac{\d{\mathcal Q}}{\d x}\r\setminus\{x: {\mathcal Q}=0\}\end{array} \label{D(d(P/Q)/dx)}\\
& \begin{array}{l}\kern-3pt \D\l\frac{\d\l {\mathcal Q}\big|_{y={\mathcal P}}\r}{\d x}\r=\\ \kern20pt=
\left\{x:\ {\mathcal P}\in\D\l\frac{\d{\mathcal Q}}{\d y}\r\right\}\cap\D\l\frac{\d{\mathcal P}}{\d x}\r\end{array} \label{D(d(Q-0-P)/dx)}
\end{align}
\erem

\bprop Для всякого одноместного числового выражения $\mathcal P$ от переменной
$x$ область допустимых значений переменной $x$ в производном выражении $\frac{\d
{\mathcal P}}{\d x}$ является открытым множеством на прямой $\R$, содержащимся
в области допустимых значений $x$ в $\mathcal P$:
$$
\Int\left(\D\left(\frac{\d {\mathcal P}}{\d x}\right)\right)=\D\left(\frac{\d
{\mathcal P}}{\d x}\right)\subseteq \D({\mathcal P})
$$
\eprop
 \bpr
Заметим, что это очевидно для выражений порядка 0 и 1. После этого предположим,
что это верно для всех выражений до порядка $n$ включительно, и покажем, что
тогда это утверждение верно и для выражений порядка $n+1$. Пусть $\mathcal U$
-- такое выражение. По определению, $\mathcal U$ будет либо суммой, либо
разностью, либо произведением, либо композицией выражений порядка не больше
$n$. Рассмотрим каждый из этих случаев.

Если ${\mathcal U}={\mathcal P}+{\mathcal Q}$, где ${\mathcal P}$ и ${\mathcal
Q}$ -- выражения порядка не больше $n$, то
 \begin{multline*}
\D\left(\frac{\d {\mathcal U}}{\d x}\right)=\D\left(\frac{\d ({\mathcal
P}+{\mathcal Q})}{\d x}\right)=\\=\underbrace{\D\left(\frac{\d {\mathcal P}}{\d
x}\right)\cap \D\left(\frac{\d {\mathcal Q}}{\d
x}\right)}_{\scriptsize\begin{matrix}\text{пересечение}\\ \text{открытых
множеств}\end{matrix}}\subseteq \D({\mathcal P})\cap \D({\mathcal
Q})=\\=\D({\mathcal P}+{\mathcal Q})=\D({\mathcal U})
 \end{multline*}
То же самое верно для случая разности и произведения выражений. Для отношения
получаем: если ${\mathcal U}=\frac{\mathcal P}{\mathcal Q}$, где ${\mathcal P}$
и ${\mathcal Q}$ -- выражения порядка не больше $n$, то
 \begin{multline*}
\D\left(\frac{\d {\mathcal U}}{\d x}\right)=\D\left(\frac{\d}{\d
x}\left(\frac{\mathcal P}{\mathcal Q}\right)\right)=\\= \D\left( \frac{
\frac{\d{\mathcal P}}{\d x}\cdot {\mathcal Q}-{\mathcal
P}\cdot\frac{\d{\mathcal Q}}{\d x}}{{\mathcal
Q}^2}\right)=\\=\left\{\D\left(\frac{\d{\mathcal P}}{\d x}\right)\cap
\D({\mathcal Q}) \cap \D\left(\frac{\d{\mathcal Q}}{\d x}\right)\cap
\D({\mathcal P})\right\}\setminus\\ \setminus \{x:\ {\mathcal Q}=0\}=\\=
\bigg\{\underbrace{\D\left(\frac{\d{\mathcal P}}{\d x}\right)\cap
\D\left(\frac{\d{\mathcal Q}}{\d
x}\right)}_{\scriptsize\begin{matrix}\text{пересечение}\\ \text{открытых
множеств}\end{matrix}}\bigg\}\setminus \underbrace{\{x:\ {\mathcal
Q}=0\}}_{\scriptsize\begin{matrix}\text{замкнутое}\\
\text{множество в $\D({\mathcal Q})$}\end{matrix}}\subseteq\\
\subseteq \big\{\D({\mathcal P})\cap \D({\mathcal Q})\big\} \setminus \{x:\
{\mathcal Q}=0\}=\D\left(\frac{\mathcal P}{\mathcal Q}\right)=\D({\mathcal U})
 \end{multline*}
Остается рассмотреть случай композиции: если ${\mathcal U}={\mathcal
Q}\Big|_{y={\mathcal P}}$, где ${\mathcal P}$ и ${\mathcal Q}$ -- выражения
порядка не больше $n$, то
 \begin{multline*}
\D\left(\frac{\d{\mathcal U}}{\d x}\right)= \D\left(\frac{\d}{\d
x}\left({\mathcal Q}\Big|_{y={\mathcal P}}\right)\right)=\\= \D\left(
\frac{\d{\mathcal Q}}{\d y}\bigg|_{y={\mathcal P}}\cdot \frac{\d{\mathcal
P}}{\d x}\right)=\\=\D\left( \frac{\d{\mathcal Q}}{\d y}\bigg|_{y={\mathcal
P}}\right)\cap \D\left(\frac{\d{\mathcal P}}{\d x}\right)=\\=
\underbrace{\underbrace{{\mathcal P}^{-1}\Big(\overbrace{\D\left(
\frac{\d{\mathcal Q}}{\d
y}\right)}^{\scriptsize\begin{matrix}\text{открытое}\\
\text{множество} \\ \text{в
$\R$}\end{matrix}}\Big)}_{\scriptsize\begin{matrix}\text{открытое}\\
\text{множество в $\D({\mathcal P})$}\end{matrix}}\cap
\underbrace{\D\left(\frac{\d{\mathcal P}}{\d x}\right)}_{\scriptsize\begin{matrix}\text{открытое}\\
\text{множество} \\ \text{в $\R$}\\ \text{и подмножество}\\
\text{в $\D({\mathcal P})$}\end{matrix}}}_{\text{открытое множество в
$\R$}}\subseteq\\ \subseteq {\mathcal P}^{-1}\big(\D({\mathcal Q})\big)\cap\D({\mathcal
P})=\D({\mathcal U})
 \end{multline*}

 \epr

\btm\label{TM:proizv-stand-func} Стандартная функция $f$, определяемая
выражением ${\mathcal U}$,
$$
f(x)={\mathcal U}
$$
дифференцируема всюду на области допустимых значений $x$ в производном
выражении $\frac{\d{\mathcal U}}{\d x}$, и производная $f$ задается на этом
множестве выражением $\frac{\d{\mathcal U}}{\d x}$:
$$
f'(x)=\frac{\d{\mathcal U}}{\d x},\qquad x\in \D\left(\frac{\d{\mathcal U}}{\d
x}\right)
$$
\etm
 \bpr
Для элементарных выражений это верно. Предположим, что это верно для выражений
до порядка $n$ включительно. Пусть ${\mathcal U}$ -- выражение порядка $n+1$.
Тогда по определению, $\mathcal U$ будет либо суммой, либо разностью, либо
произведением, либо отношением, либо композицией выражений порядка не больше $n$. Рассмотрим
каждый из этих случаев.

Если ${\mathcal U}={\mathcal P}+{\mathcal Q}$, где ${\mathcal P}$ и ${\mathcal
Q}$ -- выражения порядка не больше $n$, то можно рассмотреть функции
 \beq\label{p=P,q=Q}
p(x)={\mathcal P},\qquad q(x)={\mathcal Q}
 \eeq
и по предположению индукции, они будут дифференцируемы на множествах
$\D\left(\frac{\d{\mathcal P}}{\d x}\right)$ и $\D\left(\frac{\d{\mathcal
Q}}{\d x}\right)$ соответственно, причем
 \beq\label{p'=P',q'=Q'}
 \begin{split}
& p'(x)=\frac{\d{\mathcal P}}{\d x},\qquad x\in \D\left(\frac{\d{\mathcal P}}{\d
x}\right),\\
& q'(x)=\frac{\d{\mathcal Q}}{\d x},\qquad x\in
\D\left(\frac{\d{\mathcal Q}}{\d x}\right).
\end{split}
 \eeq
С другой стороны, функция $f$ будет их суммой на множестве
$\D\left(\frac{\d{\mathcal U}}{\d x}\right)=\D\left(\frac{\d{\mathcal P}}{\d
x}\right)\cap \D\left(\frac{\d{\mathcal Q}}{\d x}\right)$:
$$
f(x)=p(x)+q(x),\qquad x\in \D\left(\frac{\d{\mathcal U}}{\d x}\right)
$$
И значит, ее производная на этом множестве будет равна сумме производных $p$ и
$q$:
$$
f'(x)=p'(x)+q'(x),\qquad x\in \D\left(\frac{\d{\mathcal U}}{\d x}\right)
$$
Следовательно, функция $f'$ должна задаваться формулой
$$
f'(x)=p'(x)+q'(x)=\frac{\d{\mathcal P}}{\d x}+\frac{\d{\mathcal Q}}{\d
x}=\frac{\d{\mathcal U}}{\d x}
$$

То же самое верно для случая разности и произведения выражений. Для отношения
получаем: если ${\mathcal U}=\frac{\mathcal P}{\mathcal Q}$, где ${\mathcal P}$
и ${\mathcal Q}$ -- выражения порядка не больше $n$, то снова рассмотрев
функции \eqref{p=P,q=Q}, мы по предположению индукции получим тождества
\eqref{p'=P',q'=Q'}, а с другой стороны из тождества
$$
f(x)=\frac{p(x)}{q(x)}
$$
будет следовать, что функция $f$ дифференцируема на множестве
\begin{multline*}
\D\left(\frac{\d{\mathcal P}}{\d x}\right)\cap \D\left(\frac{\d{\mathcal Q}}{\d
x}\right)\setminus \{x:\ q(x)={\mathcal Q}=0\}=\\=\D\left(\frac{\d}{\d
x}\left(\frac{\mathcal P}{\mathcal Q}\right)\right)= \D\left(\frac{\d{\mathcal
U}}{\d x}\right)
\end{multline*}
и ее производная на этом множестве равна
\begin{multline*}
f'(x)=\frac{p'(x)\cdot q(x)-p(x)\cdot q'(x)}{q(x)^2}=\\=
\frac{\Big(\frac{\d{\mathcal P}}{\d x}\Big)\cdot {\mathcal Q} -{\mathcal
P}\cdot \Big(\frac{\d{\mathcal Q}}{\d x}\Big)}{{\mathcal Q}^2}=\\=\frac{\d}{\d
x}\left(\frac{\mathcal P}{\mathcal Q}\right)=\frac{\d{\mathcal U}}{\d x}
\end{multline*}

Остается рассмотреть случай композиции. Пусть ${\mathcal U}={\mathcal
Q}\Big|_{y={\mathcal P}}$, где ${\mathcal P}$ и ${\mathcal Q}$ -- выражения
порядка не больше $n$. Снова рассмотрев функции \eqref{p=P,q=Q}, мы по
предположению индукции получим тождества \eqref{p'=P',q'=Q'}, а с другой
стороны из тождества
$$
f(x)=q\big(p(x)\big)
$$
будет следовать, что функция $f$ дифференцируема на множестве
\begin{multline*}
\left\{x:\ {\mathcal P}\in\D\l\frac{\d{\mathcal Q}}{\d y}\r\right\}\cap\D\l\frac{\d{\mathcal P}}{\d x}\r=\\=\D\left(\frac{\d}{\d
x}\left({\mathcal Q}\Big|_{y={\mathcal P}}\right)\right)=
\D\left(\frac{\d{\mathcal U}}{\d x}\right)
\end{multline*}
и ее производная на этом множестве равна
$$
f'(x)=q'(p(x))\cdot p'(x)=\frac{\d{\mathcal Q}}{\d y}\Big|_{y={\mathcal
P}}\cdot \frac{\d{\mathcal P}}{\d x}=\frac{\d{\mathcal U}}{\d x}
$$

 \epr

\bex Покажем, что {\it при выборе выражения, представляющего функцию, область допустимых значений переменной в его производном выражении может меняться.} Рассмотрим функцию
$$
f(x)=x.
$$
Выражение в правой части этой формулы представляет собой переменную, поэтому
его производная будет единицей, и значит, производная функции $f$ также равна
единице:
$$
f'(x)=\frac{\d x}{\d x}=1,\quad x\in\R
$$
Но функцию $f$ можно представить другим выражением:
$$
f(x)=\sqrt[3]{x^3}
$$
Производная этого нового выражения будет определена только на множестве
$(-\infty,0)\cup (0,+\infty)$, поэтому производная функции $f$, если ее
вычислять таким образом, также должна быть определена только на этом множестве:
$$
f'(x)=\frac{\d\sqrt[3]{x^3}}{\d x}=\frac{1}{3 (\sqrt[3]{x^3})^2}\cdot 3\cdot
x^2=1, \quad x\ne 0
$$
Таким образом, при изменении выражения, которым представляется данная функция $f$, множество точек, в которых действует формула, выражающая производную функции $f$, может меняться. Замена данного представления другим может добавлять или выбрасывать какие-то точки, в которых формула для производной $f$ будет верна.
\eex

\bex {\it Композиция выражений, дифференцируемых на любом интервале своей области допустимых значений переменной, может быть недифференцируема на каких-то интервалах своей области допустимых значений переменной.}
Таким свойством обладает выражение
$$
\arcsin\sin x
$$
График определяемой им функции
$$
h(x)=\arcsin\sin x
$$
представляет собой ``пилу'':

\vglue60pt

 \noindent
-- и видно, что в точках $x=\frac{\pi}{2}+\pi n$, $n\in\Z$, эта функция не дифференцируема. Поэтому ее нельзя назвать дифференцируемой на прямой $\R$.
Но при этом составляющие ее функции
$$
f(x)=\sin x, \qquad g(y)=\arcsin y
$$
дифференцируемы на каждом (открытом) интервале своей области определения.
\eex

\paragraph{Вычисление производной.}

Всюду, как может заметить читатель, вычисление производной
теперь сжимается до единственной цепочки равенств. Вот как это выглядит в
примерах, аналогичных рассмотренным на с.\pageref{EX:(ln-x-sin-x)'}.

\bex
 \begin{multline*}
\frac{\d}{\d x}\Big(3 x ^2-  x +5\Big)=\\= 3\frac{\d}{\d x} \Big( x ^2\Big)-
\frac{\d}{\d x}\Big( x \Big)+\frac{\d}{\d x}\Big(5\Big)=\\= 3\cdot 2 x -1+0= 6
x -1
 \end{multline*}
\eex

\bex
 \begin{multline*}
\frac{\d}{\d x}\Big( x \ln x \Big)=\\= \frac{\d}{\d x}\Big( x \Big)\cdot \ln x
+  x \cdot\frac{\d}{\d x}\Big(\ln x \Big)=\\= 1\cdot \ln x +
 x \cdot\frac{1}{ x }=
\ln x + 1
 \end{multline*}
\eex

\bex
 \begin{multline*}
\frac{\d}{\d x}\Big(\sqrt{ x }\cdot \sin x \Big)=\\= \frac{\d}{\d x}\Big(\sqrt{
x }\Big)\cdot \sin x + \sqrt{ x }\cdot \frac{\d}{\d x}\Big(\sin x \Big)=\\=
\frac{1}{2\sqrt{ x }}\cdot\sin x + \sqrt{ x }\cdot\cos x
 \end{multline*}
\eex

\bex
 \begin{multline*}
\frac{\d}{\d x}\left(\frac{\arctg x }{ x ^2+1}\right)=\\= \frac{\frac{\d}{\d
x}(\arctg x )\cdot ( x ^2+1) - \frac{\d}{\d x}(x ^2+1)\cdot\arctg x }{( x
^2+1)^2}=\\= \frac{\frac{1}{ x ^2+1}\cdot ( x ^2+1) - 2 x \cdot\arctg x }{( x
^2+1)^2}=\\= \frac{1-2 x \cdot\arctg x }{( x ^2+1)^2}
\end{multline*} \eex

\bers Продифференцируйте выражения по переменной $x$:
 \bit{

\item[1)] $\arctg x + x +\arcctg x $;

\item[2)] $ x \arcsin x $;

\item[3)] $\frac{\arccos x }{\arcsin x }$;

\item[4)] $\big( x ^2-7 x +8\Big)\cdot e^{ x }$;

\item[5)] $\ln x ^3- \frac{9}{ x }-\frac{27}{2 x ^2}$;

\item[6)] $\frac{\sin x +\cos x }{\sin x - \cos x }$;

\item[7)] $(\sqrt{2})^{ x }+(\sqrt{5})^{- x }$;

\item[8)] $2^{ x }\cdot\ln| x |$;

\item[9)] $e^{ x }\cdot\log_2 x $;

\item[10)] $\log_2 x \cdot\ln x \cdot\log_3 x $;

 }\eit
\eers

Вычисления типа тех, что проводились на с.\pageref{primery-vych-proizv-slozhnoi-f}, трансформируются
в цепочку равенств, в которой первым действием служит представление сложной
функции в виде подходящей подстановки, например
$$
e^{\arctg x}=e^y|_{y=\arctg x},
$$
после чего применяется формула \eqref{7.4.4}.

\begin{ex}
 \begin{multline*}
\frac{\d }{\d x}\Big( e^{\arctg x}\Big)= \frac{\d }{\d x}\Big(
e^y\Big|_{y=\arctg x}\Big)=\eqref{7.4.4}=\\= \frac{\d e^y}{\d y}\Big|_{y=\arctg
x}\cdot \frac{\d \arctg x}{\d x}=\\= e^y\Big|_{y=\arctg x}\cdot
\frac{1}{1+x^2}= e^{\arctg x}\cdot \frac{1}{1+x^2}
 \end{multline*}\end{ex}

\begin{ex}
 \begin{multline*}
\frac{\d }{\d x}\Big(\sin^2 x\Big)= \frac{\d }{\d x}\Big( y^2\Big|_{y=\sin
x}\Big)=\eqref{7.4.4}=\\= \frac{\d y^2}{\d y}\Big|_{y=\sin x}\cdot \frac{\d
\sin x}{\d x}=\\= 2y\Big|_{y=\sin x}\cdot \cos x= 2\sin x\cdot \cos x
 \end{multline*}\end{ex}

\begin{ers}
Найдите производную по переменной $x$:

\begin{multicols}{2}

1) $\arctg e^{ x}$;

2) $\sin  x^2$;

3) $\tg \ln  x$;

4) $\ln \tg  x$.

5) $\ln \sin e^{ x }$;

6) $\sin e^{\arcsin  x }$;

7) $e^{\arcsin \ln  x }$;

8) $\arcsin \ln \sin  x $;

9) $\ln^2 (1+ x ^3)$;

10) $\cos \sqrt{ x + x ^3}$;

11) $\sqrt[3] {1+\cos  x }$.

\end{multicols}

\end{ers}

Для задач, подобных рассматривавшимся на с.\pageref{EX-(x^x)'}, полезно заметить
формулу:
 \beq\label{EQ:P^Q}
{\mathcal P}^{\mathcal Q}=e^{\ln {\mathcal P}^{\mathcal Q}}=e^{{\mathcal Q} \cdot \ln {\mathcal P}}
 \eeq
(справедливоую при ${\mathcal P}>0$).

\begin{ex} Действия в примере
\ref{EX-(x^x)'} выглядят теперь так:
 \begin{multline*}
\frac{\d}{\d x}\left( x ^{ x }\right)=\eqref{EQ:P^Q}=\frac{\d}{\d x}\left( e^{x\cdot\ln x }\right)=\\=
\frac{\d}{\d x}\left( e^y\Big|_{y=x\cdot\ln x}\right)=
\frac{\d e^y}{\d y}\Big|_{y=x\cdot\ln x}\cdot \frac{\d(x\cdot\ln x )}{\d x}=\\=
e^{x\cdot\ln x }\cdot \Big(1\cdot\ln x +x\cdot\frac{1}{x}\Big)
= x ^{ x }\cdot \Big(\ln x +1\Big).
 \end{multline*}
\end{ex}

 \begin{ex}
 \begin{multline*}
\frac{\d}{\d x}\left(\Big(\cos  x \Big)^{\sin  x }\right)=\eqref{EQ:P^Q}=\\=
\frac{\d}{\d x}\left(e^{\sin x\cdot \ln \cos x}\right)=\frac{\d}{\d x}\left( e^y\Big|_{y=\sin x\cdot \ln \cos x}\right)=\\=
\frac{\d e^y}{\d y}\Big|_{y=\sin x\cdot \ln \cos x}\cdot \frac{\d(\sin x\cdot \ln \cos x)}{\d x}=\\=
 e^y\Big|_{y=\sin x\cdot \ln \cos x} \cdot \Big(\frac{\d \sin x}{\d x}\cdot \ln \cos x+\\+\sin x\cdot\frac{\d\ln\cos x}{\d x}\Big)=
\Big(\cos x \Big)^{\sin  x }\cdot\\ \cdot \bigg(\cos x\cdot
\ln \cos x + \sin x \cdot \frac{\frac{\d}{\d x}\big(\cos x \big)}{\cos x
}\bigg) =\\= \Big(\cos x \Big)^{\sin  x }\cdot \Big(\cos x \ln \cos
 x -\frac{\sin^2
 x }{\cos x }\Big) \end{multline*}
 \end{ex}

\begin{ers} Вычислите производную:

\begin{multicols}{2}

1) $\ln^{ x }  x $;

2) $ x ^{\sqrt{ x }}$;

3) $\Big(\sin
 x \Big)^{\frac{1}{ x }}$;

4) $ x ^{e^{ x }}$.

\end{multicols}
\end{ers}

\end{multicols}\noindent\rule[10pt]{160mm}{0.1pt}

\section{Классические теоремы о дифференцируемых функциях}

\subsection{Теорема Ферма}

\begin{tm}[\bf Ферма]\label{Fermat}\index{теорема!Ферма}\footnote{Эта теорема используется в
следующем параграфе при доказательстве теоремы Ролля, а также в $\S 1$ главы
\ref{ch-graph-f(x)} при доказательстве необходимого условия локального
экстремума} Пусть функция $f$ определена на интервале $(a,b)$ и в некоторой
точке $\xi\in (a,b)$ достигает экстремума на этом интервале, то есть
$$
f(\xi)=\max_{x\in (a,b)} f(x),
$$
или
$$
f(\xi)=\min_{x\in (a,b)} f(x).
$$
Тогда, если в точке $\xi$ функция $f$ дифференцируема, то
$$
f'(\xi)=0
$$
\end{tm}

%\picture{120pt}{-40pt}{fermat.pcx}

\vglue120pt \begin{proof} Пусть для определенности
$$
f(\xi)=\max_{x\in (a,b)} f(x)
$$
то есть
\begin{equation}\forall x\in (a,b) \quad f(x)\le f(\xi) \label{8.1.1}\end{equation}
Тогда, с одной стороны,

$$
f'(\xi)=\lim_{x\to \xi}\frac{f(x)-f(\xi)}{x-\xi}={\smsize
{\smsize\begin{pmatrix}\text{если предел существует,}\\
\text{то он равен пределу слева}\end{pmatrix}}}= \lim_{x\to
\xi-0}\frac{\overbrace{f(x)-f(\xi)}^{\scriptsize\begin{matrix} 0
\\
\phantom{\tiny \eqref{8.1.1}} \ \text{\rotatebox{90}{$\ge$}} \ {\tiny
\eqref{8.1.1}}\end{matrix}}}{\underbrace{x-\xi}_{\scriptsize\begin{matrix}
\text{\rotatebox{90}{$<$}} \\ 0 \\ \Uparrow \\ x<\xi \end{matrix}}}\ge 0
$$
А, с другой стороны,
$$
f'(\xi)=\lim_{x\to \xi}\frac{f(x)-f(\xi)}{x-\xi}={\smsize
{\smsize\begin{pmatrix}\text{если предел существует,}\\
\text{то он равен пределу справа}\end{pmatrix}}}= \lim_{x\to
\xi+0}\frac{\overbrace{f(x)-f(\xi)}^{\scriptsize\begin{matrix} 0
\\
\phantom{\tiny \eqref{8.1.1}} \ \text{\rotatebox{90}{$\ge$}} \ {\tiny
\eqref{8.1.1}}\end{matrix}}}{\underbrace{x-\xi}_{\scriptsize\begin{matrix}
\text{\rotatebox{90}{$>$}} \\ 0 \\ \Uparrow \\ x>\xi \end{matrix}}}\le 0
$$
Итак, мы получили, что $f'(\xi)\ge 0$ и, одновременно, $f'(\xi)\le 0$. Значит,
$f'(\xi)=0$.
\end{proof}

\subsection{Теорема Ролля}

\begin{tm}[\bf Ролля]\label{Roll}\index{теорема!Ролля}\footnote{Эта теорема используется в
следующих трех параграфах при доказательстве теорем Лагранжа, Коши и Тейлора}
Пусть функция $f$ определена на отрезке $[a,b]$, причем
 \bit{
\item[(i)] $f$ непрерывна на отрезке $[a,b]$,

\item[(ii)] $f$ дифференцируема на интервале $(a,b)$,

\item[(iii)] $f(a)=f(b)$.
 }\eit
Тогда существует точка $\xi\in (a,b)$, такая что
$$
f'(\xi)=0
$$
\end{tm}

%\picture{120pt}{-40pt}{roll.pcx}

\vglue120pt \begin{proof} Поскольку функция $f$ непрерывна на отрезке
$[a,b]$, по теореме Вейерштрасса об экстремумах (теорема \ref{Wei-III}), $f(x)$
имеет максимум и минимум на $[a,b]$:
$$
m=f(x_m)=\min_{x\in [a,b]} f(x) \qquad M=f(x_M)=\max_{x\in [a,b]} f(x)
$$
Ясно, что $m\le M$, поэтому возможны два случая:
 \bit{
\item[1)] $m=M$, или \item[2)] $m<M$.
 }\eit

В первом случае мы получаем, что функция $f$ должна быть постоянной
$$
 \forall x\in [a,b] \quad m=f(x)=M,
$$
поэтому ее производная $f'(x)$ будет равна нулю в любой точке $x\in [a,b]$.

Во втором же случае мы получаем, что с одной стороны $f(a)=f(b)$, а с другой,
$m=f(x_m)<f(x_M)=M$. Это означает, что $x_m$ или $x_M$ не не может быть концом
отрезка $[a,b]$:
 \bit{
\item[--] либо $a<x_m<b$, \item[--] либо $a<x_M<b$.
 }\eit
Значит, функция $f$ имеет экстремум на интервале $(a,b)$ (если $a<x_m<b$ --
то минимум, а если $a<x_M<b$ -- то максимум). Следовательно, по теореме Ферма
\ref{Fermat}, найдется точка $\xi\in (a,b)$ такая, что $f'(\xi)=0$.
\end{proof}

\subsection{Теорема Лагранжа}

\begin{tm}[\bf Лагранжа]\label{Lagrange}\index{теорема!Лагранжа}\footnote{Эта теорема используется в
главе \ref{ch-graph-f(x)} при доказательстве теоремы \ref{suff-conv}
(достаточное условие строгой выпуклости), а также в главе 12 при доказательстве
леммы 1.1 (о постоянности функции с нулевой производной).} Пусть функция $f$
определена на отрезке $[a,b]$, причем
 \bit{
\item[(i)] $f$ непрерывна на отрезке $[a,b]$,

\item[(ii)] $f$ дифференцируема на интервале $(a,b)$.
 }\eit
Тогда существует точка $\xi\in (a,b)$, такая что
 \beq\label{Lagrange-xi}
f'(\xi)=\frac{f(b)-f(a)}{b-a}
 \eeq
\end{tm}

%\picture{120pt}{-40pt}{lagrange.pcx}

\vglue120pt \begin{proof} Рассмотрим вспомогательную функцию
$$
F(x)=f(x)-f(a)-\frac{f(b)-f(a)}{b-a} (x-a)
$$
Эта функция удовлетворяет всем трем условиям теоремы Ролля \ref{Roll}:
 \bit{
\item[(1)] $F$ непрерывна на отрезке $[a,b]$ (по теореме \ref{cont-alg}, потому
что $F(x)$ составлена из непрерывных функций с помощью операции вычитания);

\item[(2)] $F$ дифференцируема на интервале $(a,b)$ (по теореме \ref{diff-alg},
потому что $F$ составлена из дифференцируемых функций с помощью операции
вычитания);

\item[(3)] $F(a)=F(b)$ (потому что $F(a)=f(a)-f(a)-\frac{f(b)-f(a)}{b-a}
(a-a)=0-0=0$ и $F(b)=f(b)-f(a)-\frac{f(b)-f(a)}{b-a} (b-a)=
f(b)-f(a)-\frac{f(b)-f(a)}{1} =0$).
 }\eit
Значит, по теореме Ролля \ref{Roll}, найдется такая точка $\xi\in (a,b)$, что
$F'(\xi)=0$. Иными словами,
$$
f'(\xi)-\frac{f(b)-f(a)}{b-a}=0
$$
или
$$
f'(\xi)=\frac{f(b)-f(a)}{b-a}\qquad $$ \end{proof}

Из теоремы Лагранжа следует утверждение, которое понадобится нам в главе \ref{CH-indef-integral} в разговоре о неопределенных интегралах:

\bcor\label{LM:o-konstantah} Пусть функция $f$ дифференцируема на интервале $I$ и всюду на нем имеет нулевую производную:
 \beq\label{f'(x)=0}
f'(x)=0,\quad x\in I.
 \eeq
Тогда $f$ постоянна на $I$:
$$
f(a)=f(b),\quad a,b\in I
$$
\ecor
 \bpr Возьмем произвольные точки $a,b\in I$, причем пусть $a<b$. Тогда функция $f$
 \bit{
\item[--] дифференцируема на интервале $(a;b)$ (потому что она дифференцируема
на $I$), и

\item[--] непрерывна на отрезке $[a;b]$ (по предложению
\ref{PROP:diff=>cont-[a,b]})
 }\eit\noindent
поэтому, мы можем применить к функции $f$ на отрезке $[a;b]$ теорему Лагранжа
\ref{Lagrange}, и тогда получим, что существует такая точка $\xi \in [a;b]$,
что
$$
\frac{f(b)-f(a)}{b-a}=f'(\xi)= \eqref{f'(x)=0}=0
$$
отсюда
$$
f(b)-f(a)=0\quad\Longrightarrow\quad f(b)=f(a)
$$
 \epr

\subsection{Теорема Коши об отношении приращений}

\begin{tm}[\bf Коши]\label{Cauchy-II}\index{теорема!Коши!об отношении приращений}\footnote{Эта теорема
используется в главе \ref{ch-lopital} при доказательстве теоремы
\ref{Lopital_x->a_0/0} (правило Лопиталя для предела в точке с
неопределенностью $\frac{0}{0}$)} Пусть функции $f$ и $g$ определены на отрезке
$[a,b]$, причем
 \bit{
\item[(i)] $f$ и $g$ непрерывны на отрезке $[a,b]$,

\item[(ii)] $f$ и $g$ дифференцируемы на интервале $(a,b)$,

\item[(iii)] $g'(x)\ne 0$ при $x\in (a,b)$.
 }\eit
Тогда существует точка $\xi\in (a,b)$, такая что
\begin{equation}\frac{f'(\xi)}{g'(\xi)}=\frac{f(b)-f(a)}{g(b)-g(a)}\label{8.5.1}\end{equation}\end{tm}\begin{proof} Сначала нужно убедиться, что
$g(a)\ne g(b)$, то есть что формула \eqref{8.5.1} имеет смысл. Действительно,
если бы оказалось, что $g(a)=g(b)$, то это означало бы, что $g(x)$
удовлетворяет всем трем условиям теоремы Ролля \ref{Roll}, поэтому нашлась бы
точка $\xi \in (a,b)$ такая, что $g'(\xi)=0$. Это невозможно по условию $(iii)$
нашей теоремы.

Теперь перейдем к доказательству формулы \eqref{8.5.1}. Рассмотрим
вспомогательную функцию
$$
F(x)=f(x)-f(a)-\frac{f(b)-f(a)}{g(b)-g(a)}\cdot (g(x)-g(a))
$$
Эта функция удовлетворяет всем трем условиям теоремы Ролля \ref{Roll}:
 \bit{
\item[(1)] $F$ непрерывна на отрезке $[a,b]$ (по теореме \ref{cont-alg}, потому
что $F$ составлена из непрерывных функций с помощью операции вычитания);

\item[(2)] $F$ дифференцируема на интервале $(a,b)$ (по теореме \ref{diff-alg},
потому что $F$ составлена из дифференцируемых функций с помощью операции
вычитания);

\item[(3)] $F(a)=F(b)$ (потому что
$F(a)=f(a)-f(a)-\frac{f(b)-f(a)}{g(b)-g(a)}\cdot (g(a)-g(a))=0-0=0$ и
$F(b)=f(b)-f(a)-\frac{f(b)-f(a)}{g(b)-g(a)} (g(b)-g(a))=
f(b)-f(a)-\frac{f(b)-f(a)}{1} =0$).
 }\eit
Значит, по теореме Ролля \ref{Roll}, найдется такая точка $\xi\in (a,b)$, что
$$
F'(\xi)=0
$$
Иными словами,
$$
f'(\xi)-\frac{f(b)-f(a)}{g(b)-g(a)}\cdot g'(\xi)=0
$$
или, учитывая что $g'(\xi)\ne 0$,
$$
\frac{f'(\xi)}{g'(\xi)}=\frac{f(b)-f(a)}{g(b)-g(a)}\qquad
$$ \end{proof}

\section{Построение графика}\label{SEC:grafik}

В главе \ref{ch-f'(x)} мы уже говорили, что производная бывает нужна при
построении графика функции. Мы можем теперь, наконец, объяснить, как она
используется, и доказать необходимые для этого математические утверждения.

\subsection{Монотонность и экстремум}

\paragraph{Монотонность.}

\btm[\bf о нестрогой монотонности]\label{TH:o-netsr-monot} Пусть $f$ --
дифференцируемая функция на интервале $(a;b)$. Тогда
 \bit{
\item[(a)] $f$ неубывает на $(a;b)$ в том и только в том случае, если ее
производная неотрицательна всюду на $(a;b)$:
$$
f'(x)\ge 0,\qquad x\in(a;b)
$$

\item[(b)] $f$ невозрастает на $(a;b)$ в том и только в том случае, если ее
производная неположительна всюду на $(a;b)$:
$$
f'(x)\le 0,\qquad x\in(a;b)
$$
 }\eit
\etm
 \bpr
Мы докажем утверждение (a), оставив читателю выполнить то же самое для (b) по
аналогии. Если $f$ неубывает,
$$
\forall x,y\in (a;b)\qquad \Big(x<y\qquad\Longrightarrow\qquad f(x)\le
f(y)\Big)
$$
то, переходя при вычислении производной в точке $x\in(a;b)$ к одностороннему
пределу, мы получим:
$$
f'(x)=\lim_{y\to x}\frac{f(y)-f(x)}{y-x}=\lim_{y\to x+0}
\frac{\overbrace{f(y)-f(x)}^{\scriptsize\begin{matrix}0\\\text{\rotatebox{90}{$\ge$}}\end{matrix}}}
{\underbrace{y-x}_{\scriptsize\begin{matrix}\text{\rotatebox{90}{$<$}}\\0\end{matrix}}}\ge
0
$$
Наоборот, если в любой точке $x\in(a;b)$ производная неотрицательна, то для
любых $x,y\in (a;b)$ таких, что $x<y$, подобрав по теореме Лагранжа
\ref{Lagrange} точку $\xi\in(a;b)$ со свойством \eqref{Lagrange-xi}, мы
получим:
$$
\frac{f(y)-f(x)}{\underbrace{y-x}_{\scriptsize\begin{matrix}\text{\rotatebox{90}{$<$}}\\0\end{matrix}}}=f'(\xi)\ge
0 \qquad\Longrightarrow\qquad f(y)-f(x)\ge 0\qquad\Longrightarrow\qquad f(x)\le
f(y)
$$
 \epr

\begin{tm}[\bf о строгой монотонности]
\label{Roll-cons}\index{теорема!о строгой монотонности}\footnote{Теорема
\ref{Roll-cons} используется в главе \ref{ch-lopital} при доказательстве
правила Лопиталя для раскрытия неопределенностей типа $\frac{\infty}{\infty}$}
Пусть $f$ -- дифференцируемая функция на интервале $(a,b)$, причем ее
производная нигде не обращается в нуль:
 \beq\label{f'(x)-ne-0}
  \forall x\in (a,b) \quad f'(x)\ne 0
 \eeq
Тогда $f$ строго монотонна на $(a,b)$, а именно:
 \bit{
\item[---] либо, $f'(x)>0$ всюду на $(a,b)$, и в этом случае $f$ возрастает на
$(a,b)$:
$$
\forall x,y\in (a;b)\qquad \Big(x<y\quad\Longrightarrow\quad f(x)<f(y)\Big)
$$

\item[---] либо, $f'(x)<0$ всюду на $(a,b)$, и в этом случае $f$ убывает на
$(a,b)$:
$$
\forall x,y\in (a;b)\qquad \Big(x<y\quad\Longrightarrow\quad f(x)>f(y)\Big).
$$

 }\eit
\end{tm}

Доказательству этой теоремы мы предпошлем две леммы.

\blm\label{LM:Roll-cons-1} В условиях теоремы \ref{Roll-cons} для любых трех
точек $x,y,z\in(a;b)$, связанных неравенством
$$
x<y<z,
$$
выполняется следующее:
 \bit{
\item[(a)] условие $f(x)<f(z)$ влечет за собой условие $f(x)<f(y)<f(z)$:
$$
f(x)<f(z)\qquad\Longrightarrow\qquad f(x)<f(y)<f(z)
$$

\item[(b)] условие $f(x)>f(z)$ влечет за собой условие $f(x)>f(y)>f(z)$:
$$
f(x)>f(z)\qquad\Longrightarrow\qquad f(x)>f(y)>f(z).
$$
 }\eit
\elm
 \bpr Мы докажем утверждение (a) (а (b) читатель сможет доказать по аналогии).
Пусть $f(x)<f(z)$. Нам нужно показать, что число $f(y)$ лежит в интервале
$(f(x);f(z))$:
$$
f(y)\in \big(f(x);f(z)\big)
$$
Если это не так, то либо $f(y)$ лежит ниже этого интервала, то есть $f(y)\le
f(x)$, либо, наоборот, выше, то есть $f(y)\ge f(z)$. Убедимся, что ни то, ни
другое невозможно.

1. Предположим, что $f(y)\le f(x)$. Если $f(y)=f(x)$, то по теореме Ролля
\ref{Roll} найдется точка $\xi\in(x;y)$ такая, что
$$
f'(\xi)=0,
$$
что невозможно. Значит, должно быть $f(y)<f(x)$. Тогда, если обозначить
$C=f(x)$, то мы получим
$$
f(y)<\underbrace{f(x)}_{\scriptsize\begin{matrix}\|\\ C\end{matrix}}<f(z),
$$
и, по теореме Коши о промежуточном значении \ref{Cauchy-I}, на интервале
$(y;z)$ должна найтись такая точка $c\in(y;z)$, что
$$
f(c)=C=f(x)
$$
Теперь опять по теореме Ролля \ref{Roll} должна существовать точка
$\xi\in(x;c)$ такая, что
$$
f'(\xi)=0,
$$
и это снова противоречит условию \eqref{f'(x)-ne-0}. Мы получаем, что
неравенство $f(y)\le f(x)$ невозможно.

2. Проверим, что будет, если $f(y)\ge f(z)$. Если $f(y)=f(z)$, то по теореме
Ролля \ref{Roll} найдется точка $\xi\in(y;z)$ такая, что
$$
f'(\xi)=0,
$$
что невозможно. Значит, должно быть $f(y)>f(z)$. Тогда, если обозначить
$C=f(z)$, то мы получим
$$
f(x)<\underbrace{f(z)}_{\scriptsize\begin{matrix}\|\\ C\end{matrix}}<f(y),
$$
и, по теореме Коши о промежуточном значении \ref{Cauchy-I}, на интервале
$(x;y)$ должна найтись такая точка $c\in(x;y)$, что
$$
f(c)=C=f(z)
$$
Теперь опять по теореме Ролля \ref{Roll} должна существовать точка
$\xi\in(c;z)$ такая, что
$$
f'(\xi)=0,
$$
и это снова противоречит условию \eqref{f'(x)-ne-0}. Мы получаем, что
неравенство $f(y)\ge f(x)$ тоже невозможно.
 \epr

\blm\label{LM:Roll-cons-2} В условиях теоремы \ref{Roll-cons} пусть
$\alpha,\beta\in(a,b)$ -- произвольные точки, связанные неравенством
$$
\alpha<\beta,
$$
Тогда:
 \bit{
\item[(a)] если $f(\alpha)<f(\beta)$, то функция $f$ возрастает на интервале
$(\alpha;\beta)$:
$$
f(\alpha)<f(\beta)\quad\Longrightarrow\quad\bigg\{\forall x,y\in
(\alpha;\beta)\qquad \Big(x<y\quad\Longrightarrow\quad f(x)<f(y)\Big)\bigg\}
$$

\item[(b)] если $f(\alpha)>f(\beta)$, то функция $f$ убывает на интервале
$(\alpha;\beta)$:
$$
f(\alpha)>f(\beta)\quad\Longrightarrow\quad\bigg\{\forall x,y\in
(\alpha;\beta)\qquad \Big(x<y\quad\Longrightarrow\quad f(x)>f(y)\Big)\bigg\}
$$

\item[(c)] равенство $f(\alpha)=f(\beta)$ невозможно.
 }\eit
 \elm
\bpr 1. Если $f(\alpha)<f(\beta)$, то мы получаем логическую цепочку:
$$
\underbrace{f(\alpha)<f(\beta),\qquad \alpha<x<\beta}_{ \text{лемма
\ref{LM:Roll-cons-1}(a)}\quad{\text{\large$\Downarrow$}}\quad\phantom{\text{лемма
\ref{LM:Roll-cons-1}(a)}}}
$$
$$
\phantom{x<y<\beta}\qquad f(\alpha)<\underbrace{f(x)<f(\beta),\qquad
x<y<\beta}_{ \text{лемма
\ref{LM:Roll-cons-1}(a)}\quad{\text{\large$\Downarrow$}}\quad\phantom{\text{лемма
\ref{LM:Roll-cons-1}(a)}}}
$$
$$
\kern100pt \begin{matrix} \underbrace{f(x)<f(y)<f(\beta)}\\ \text{\large$\Downarrow$} \\
f(x)<f(y)\end{matrix}
$$

2. Наоборот, если $f(\alpha)>f(\beta)$, то мы получаем цепочку:
$$
\underbrace{f(\alpha)>f(\beta),\qquad \alpha<x<\beta}_{ \text{лемма
\ref{LM:Roll-cons-1}(b)}\quad{\text{\large$\Downarrow$}}\quad\phantom{\text{лемма
\ref{LM:Roll-cons-1}(b)}}}
$$
$$
\phantom{x<y<\beta}\qquad f(\alpha)>\underbrace{f(x)>f(\beta),\qquad
x<y<\beta}_{ \text{лемма
\ref{LM:Roll-cons-1}(b)}\quad{\text{\large$\Downarrow$}}\quad\phantom{\text{лемма
\ref{LM:Roll-cons-1}(b)}}}
$$
$$
\kern100pt \begin{matrix} \underbrace{f(x)>f(y)>f(\beta)}\\ \text{\large$\Downarrow$} \\
f(x)>f(y)\end{matrix}
$$

3. Остается убедиться, что равенство $f(\alpha)=f(\beta)$ невозможно.
Действительно, если бы оно выполнялось, то мы получили бы по теореме Ролля
\ref{Roll}, что  $f'(\xi)=0$ для некоторого $\xi\in (\alpha,\beta)$, а это
противоречит условию \eqref{f'(x)-ne-0}.
 \epr

\begin{proof}[Доказательство теоремы \ref{Roll-cons}]

Зафиксируем произвольные точки $\alpha,\beta\in (a,b)$, связанные неравенством
$$
\alpha<\beta
$$
По лемме \ref{LM:Roll-cons-2}(c), $f(\alpha)\ne f(\beta)$, значит должно быть
либо $f(\alpha)<f(\beta)$, либо $f(\alpha)>f(\beta)$. Рассмотрим каждый случай
отдельно.

1. Пусть $f(\alpha)<f(\beta)$. Тогда по лемме \ref{LM:Roll-cons-2}(a), функция
$f$ возрастает на интервале $(\alpha;\beta)$. Если теперь выбрать какие-нибудь
точки $x,y$ так, чтобы
 \beq\label{a<x<alpha<beta<y<b}
a<x<\alpha<\beta<y<b
 \eeq
то снова по лемме \ref{LM:Roll-cons-2} функция $f$ должна быть строго
монотонной на интервале $(x;y)$. При этом, как мы уже поняли, $f$ возрастает на
меньшем интервале $(\alpha;\beta)$, значит, она возрастает и на большем
интервале $(x;y)$. Поскольку это верно для любых $x,y$, удовлетворяющих условию
\eqref{a<x<alpha<beta<y<b}, мы получаем, что $f$ возрастает на всем интервале
$(a;b)$.

Далее получается цепочка:
 $$
 \text{$f$ возрастает на интервале $(a;b)$}
 $$
 $$
\phantom{\text{\scriptsize теорема \ref{TH:o-netsr-monot}(a)}}
\qquad\Downarrow\qquad\text{\scriptsize теорема \ref{TH:o-netsr-monot}(a)}
 $$
 $$
 f'(x)\ge 0,\qquad x\in(a;b)
 $$
 $$
\phantom{\text{\scriptsize \eqref{f'(x)-ne-0}}}
\qquad\Downarrow\qquad\text{\scriptsize \eqref{f'(x)-ne-0}}
 $$
 $$
 f'(x)>0,\qquad x\in(a;b)
 $$

2. Пусть наоборот, $f(\alpha)>f(\beta)$. Тогда по лемме
\ref{LM:Roll-cons-2}(b), функция $f$ убывает на интервале $(\alpha;\beta)$.
Если теперь выбрать какие-нибудь точки $x,y$ так, чтобы выполнялось
\eqref{a<x<alpha<beta<y<b}, то снова по лемме \ref{LM:Roll-cons-2} функция $f$
должна быть строго монотонной на интервале $(x;y)$. При этом, как мы уже
поняли, $f$ убывает на меньшем интервале $(\alpha;\beta)$, значит, она убывает
и на большем интервале $(x;y)$. Поскольку это верно для любых $x,y$,
удовлетворяющих условию \eqref{a<x<alpha<beta<y<b}, мы получаем, что $f$
убывает на всем интервале $(a;b)$.

Далее получается цепочка:
 $$
 \text{$f$ убывает на интервале $(a;b)$}
 $$
 $$
\phantom{\text{\scriptsize теорема \ref{TH:o-netsr-monot}(b)}}
\qquad\Downarrow\qquad\text{\scriptsize теорема \ref{TH:o-netsr-monot}(b)}
 $$
 $$
 f'(x)\le 0,\qquad x\in(a;b)
 $$
 $$
\phantom{\text{\scriptsize \eqref{f'(x)-ne-0}}}
\qquad\Downarrow\qquad\text{\scriptsize \eqref{f'(x)-ne-0}}
 $$
 $$
 f'(x)<0,\qquad x\in(a;b)
 $$
\end{proof}

\paragraph{Экстремум.}

\bit{ \item[$\bullet$] Точка $x_0$ называется
 \bit{
\item[--] точкой {\it строгого локального максимума} функции $f$, если в
некоторой окрестности $(x_0-\delta;x_0+\delta)$ точки $x_0$ выполняется
неравенство
$$
  f(x)<f(x_0)
$$
\item[--] точкой {\it строгого локального минимума} функции $f$, если в
некоторой окрестности $(x_0-\delta;x_0+\delta)$ точки $x_0$ выполняется
неравенство
$$
  f(x)>f(x_0)
$$
\item[--] точкой {\it  локального экстремума} функции $f$, если она является
точкой строгого локального минимума, или строгого локального максимума.
 }\eit
 }\eit

\begin{tm}[\bf необходимое условие локального экстремума]
\label{nes-loc-extr} Если функция $f$ имеет в точке $x_0$ локальный
экстремум и дифференцируема в этой точке, то
$$
f'(x_0)=0
$$
\end{tm}\begin{proof} В точке $x_0$ функция $f$
достигает максимума или минимума на некотором интервале
$(x_0-\delta;x_0+\delta)$, значит, по теореме Ферма \ref{Fermat}, $f'(x_0)=0$.
\end{proof}

\begin{tm}[\bf достаточное условие локального экстремума]
Пусть функция $f$ диф\-ферен\-цируема в некоторой окрестности
$(x_0-\delta;x_0+\delta)$ точки $x_0$. Тогда
 \bit{
\item[--] если $f'(x)>0$ при $x\in (x_0-\delta;x_0)$, и $f'(x)<0$ при $x\in
(x_0;x_0+\delta)$ (то есть при переходе через точку $x_0$ производная меняет
знак с $+$ на $-$), то $x_0$ -- точка строгого локального максимума; \item[--]
если $f'(x)<0$ при $x\in (x_0-\delta;x_0)$, и $f'(x)>0$ при $x\in
(x_0;x_0+\delta)$ (то есть при переходе через точку $x_0$ производная меняет
знак с $-$ на $+$), то $x_0$ -- точка строгого локального минимума.
 }\eit
\end{tm}\begin{proof} Пусть при переходе через точку
$x_0$ производная меняет знак с $-$ на $+$. Тогда по теореме \ref{Roll-cons},
на интервале $(x_0-\delta;x_0)$ функция $f$ строго убывает, а на интервале
$(x_0;x_0+\delta)$ функция $f$ строго возрастает, значит $x_0$ -- точка
локального минимума.

Аналогично рассматривается случай, когда при переходе через точку $x_0$
производная меняет знак с $-$ на $+$. \end{proof}

\subsection{Выпуклость и перегиб}\label{SUBSEC:vypuklost}

\paragraph{Выпуклость.}

Пусть функция $f$ определена на интервале $(a;b)$ и пусть
$a<\alpha<\beta<b$. Проведем хорду через точки $(\alpha; f(\alpha))$ и $(\beta;
f(\beta))$ на графике функции $f$:

%\picture{0pt}{0pt}{119.pcx}

\vglue120pt

Уравнение этой хорды имеет вид
$$
  y=\frac{f(\beta)\cdot (x-\alpha)+f(\alpha)\cdot (\beta-x)}{\beta-\alpha}
$$

 \bit{
\item[$\bullet$] Функция $f(x)$ называется
 \bit{
\item[--] {\it выпуклой вверх} на интервале $(a;b)$, если для любых двух точек
$\alpha,\beta : \,\, a<\alpha<\beta<b$ график функции $f$ на интервале
$(\alpha,\beta)$ лежит выше хорды, то есть
 \begin{equation}\frac{f(\beta)\cdot (x-\alpha)+f(\alpha)\cdot
(\beta-x)}{\beta-\alpha}< f(x), \quad a<\alpha<x<\beta<b
 \label{10.2.1}\end{equation}
соответствующая картинка выглядит так:

%\picture{0pt}{0pt}{120.pcx}

\vglue120pt

\item[--] {\it выпуклой вниз} на интервале $(a;b)$, если для любых двух точек
$\alpha,\beta : \,\, a<\alpha<\beta<b$ график функции $f$ на интервале
$(\alpha,\beta)$ лежит ниже хорды, то есть
\begin{equation}
f(x)<\frac{f(\beta)\cdot (x-\alpha)+f(\alpha)\cdot (\beta-x)}{\beta-\alpha},
\quad a<\alpha<x<\beta<b \label{10.2.2}\end{equation} соответствующая картинка
выглядит так:

%\picture{0pt}{0pt}{pp2.pcx}

\vglue120pt
 }\eit
 }\eit

\begin{tm}[\bf достаточное условие выпуклости]
\label{suff-conv} Пусть функция $f$ дважды дифференцируема на интервале
$(a;b)$. Тогда
 \bit{
\item[--] если $f''(x)<0$ при $x\in (a;b)$, то функция $f$ выпукла вверх на
интервале $(a;b)$; \item[--] если $f''(x)>0$ при $x\in (a;b)$, то функция $f$ выпукла вниз на интервале $(a;b)$.
 }\eit
\end{tm}\begin{proof} Рассмотрим разность между хордой и
функцией:
 \begin{multline*}\frac{f(\beta)\cdot (x-\alpha)+f(\alpha)\cdot
(\beta-x)}{\beta-\alpha}-f(x)= \frac{f(\beta)\cdot (x-\alpha)+f(\alpha)\cdot
(\beta-x)- f(x)\cdot (\beta-\alpha)}{\beta-\alpha}=\\= \frac{f(\beta)\cdot
(x-\alpha)+f(\alpha)\cdot (\beta-x)- f(x)\cdot (\beta-\alpha)+f(x)\cdot
(x-x)}{\beta-\alpha}=\\= \frac{[f(\beta)-f(x)]\cdot
(x-\alpha)-[f(x)-f(\alpha)]\cdot (\beta-x)} {\beta-\alpha}=\\=
{\smsize\begin{pmatrix}\text{применяем
теорему Лагранжа \ref{Lagrange}}\\
\text{для $f(x)$ на отрезках $[x;\beta]$ и $[\alpha;x]$:}\\
\text{найдутся $\theta \in [x;\beta]$ и $\eta \in [\alpha;x]$}\\
\text{такие что}\,\, f(\beta)-f(x)=f'(\theta)\cdot (\beta-x)
\\
\text{и}\,\, f(x)-f(\alpha)=f'(\eta)\cdot (x-\alpha)
\end{pmatrix}}=\\= \frac{f'(\theta)\cdot (\beta-x)\cdot (x-\alpha)-
f'(\eta)\cdot (x-\alpha)\cdot (\beta-x)} {\beta-\alpha}=
\frac{[f'(\theta)-f'(\eta)]\cdot (\beta-x)\cdot (x-\alpha)} {\beta-\alpha}=\\=
{\smsize\begin{pmatrix}\text{снова применяем
теорему Лагранжа \ref{Lagrange}}\\
\text{но уже для $f'(x)$ на отрезке $[\eta;\theta]$:}\\
\text{найдется $\xi \in [\eta;\theta]$}\\
\text{такое что}\,\, f'(\theta)-f'(\eta)=f''(\xi)\cdot (\theta-\eta)
\end{pmatrix}}= f''(\xi)\cdot \frac{(\theta-\eta)\cdot (\beta-x)\cdot
(x-\alpha)} {\beta-\alpha}
 \end{multline*}
 Заметим теперь, что
$(\theta-\eta)>0, \quad (\beta-x)>0, \quad (x-\alpha)>0, \quad
(\beta-\alpha)>0$. Значит, если $f''(x)<0$ при $x\in (a;b)$, то, в частности,
$f''(\xi)<0$, и поэтому
$$
\frac{f(\beta)\cdot (x-\alpha)+f(\alpha)\cdot (\beta-x)}{\beta-\alpha}-f(x)<0
$$
(то есть функция $f$ выпукла вверх на  интервале $(a;b)$). Если же
$f''(x)>0$ при $x\in (a;b)$, то, в частности, $f''(\xi)>0$, и поэтому
$$
\frac{f(\beta)\cdot (x-\alpha)+f(\alpha)\cdot (\beta-x)}{\beta-\alpha}-f(x)<0
$$
(то есть  функция $f$ выпукла вниз на интервале $(a;b)$).
\end{proof}

\paragraph{Перегиб.}

\bit{ \item[$\bullet$] Точка $x_0$ называется {\it точкой перегиба} функции $f$, если при переходе через $x_0$ выпуклость вверх меняется на выпуклость
вниз (или наоборот). Иными словами, имеется некоторая окрестность
$(x_0-\delta;x_0+\delta)$ такая, что на интервале $(x_0-\delta;x_0)$ функция $f$ выпукла вверх, а на интервале $(x_0;x_0+\delta)$ она выпукла вниз (или
наоборот, на $(x_0-\delta;x_0)$ функция $f$ выпукла вниз, а на интервале
$(x_0;x_0+\delta)$ она выпукла вверх).
 }\eit

\begin{tm}[\bf о точках перегиба]\label{tochki-peregiba}
Пусть функция $f$ дважды дифференцируема на интервале $(a;b)$, причем ее
вторая производная $f''$ непрерывна на $(a;b)$. Тогда если точка $x_0\in
(a;b)$ -- точка перегиба функции $f$, то $f''(x_0)=0$.
\end{tm}\begin{proof} Предположим, что $f''(x_0)\ne 0$,
тогда либо $f''(x_0)<0$, либо $f''(x_0)>0$. Если $f''(x_0)<0$, то, по теореме о
сохранении знака \ref{sign-pres}, $f''(x)<0$ для всех $x$ из некоторой
окрестности $(x_0-\delta;x_0+\delta)$ точки $x_0$. Поэтому по теореме
\ref{suff-conv} этой главы, $f(x)$ будет выпукла вверх всюду в окрестности
$(x_0-\delta;x_0+\delta)$, и значит $x_0$ не может быть точкой перегиба.
Аналогично, если $f''(x_0)>0$, то по теореме о сохранении знака
\ref{sign-pres}, $f''(x)>0$ для всех $x$ из некоторой окрестности
$(x_0-\delta;x_0+\delta)$ точки $x_0$. Поэтому по теореме \ref{suff-conv} этой
главы, $f(x)$ будет выпукла вверх всюду в окрестности
$(x_0-\delta;x_0+\delta)$, и опять $x_0$ не может быть точкой перегиба.
\end{proof}

\subsection{Асимптоты}

Асимптотой называется прямая, к которой график функции ``бесконечно
приближается''. Более точное определение этому понятию выглядит так:

\bit{ \item[$\bullet$] Пусть функция $f$ обладает следующим свойством:
$$
\lim_{x\to x_0-0} f(x)=\infty \quad \text{или}\quad \lim_{x\to x_0+0}
f(x)=\infty
$$
Тогда прямая
$$
x=x_0
$$
называется {\it вертикальной асимптотой}\index{асимптота!вертикальная} функции $f$.
 }\eit

\bit{ \item[$\bullet$] Пусть функция $f$ обладает следующим свойством:
$$
\lim_{x\to +\infty}\left( f(x)-(k\cdot x+b) \right)=0 \quad \text{или}\quad
\lim_{x\to -\infty}\left( f(x)-(k\cdot x+b) \right)=0
$$
Тогда прямая
$$
y=k\cdot x+b
$$
называется {\it наклонной асимптотой}\index{асимптота!наклонная} функции $f$
(при $x\to +\infty$ или $x\to -\infty$).
 }\eit

\noindent\rule{160mm}{0.1pt}\begin{multicols}{2}

\begin{ex}
Функция $f(x)=\frac{1}{x}$ имеет вертикальную асимпоту $x=0$.

%\picture{0pt}{0pt}{121.pcx}

\vglue80pt
\end{ex}

\begin{ex} Функция $f(x)=\frac{1}{x}+x$
имеет наклонную асимпоту $y=x$.

%\picture{0pt}{0pt}{122.pcx}

\vglue100pt
\end{ex}

\end{multicols}\noindent\rule[10pt]{160mm}{0.1pt}

\begin{tm}[\bf об асимптотах]
Если прямая $y=k\cdot x+b$ является наклонной асимптотой функции $f$ при
$x\to +\infty$ (или $x\to -\infty$), то
\begin{equation}
 \boxed{
 k=\lim_{x\to +\infty}\frac{f(x)}{x}
 }
 \quad
 \left(\text{или}\quad
  \boxed{
  k=\lim_{x\to -\infty}\frac{f(x)}{x}
  }\quad
 \right)
 \label{10.3.1}\end{equation}
и
\begin{equation}
 \boxed{
 b=\lim_{x\to +\infty} (f(x)-k\cdot x)
 }
 \quad \left(\text{или}\quad
  \boxed{
 b=\lim_{x\to -\infty} (f(x)-k\cdot x)
 }\quad \right)
 \label{10.3.2}\end{equation}\end{tm}\begin{proof} Если
$$
\lim_{x\to +\infty}\left( f(x)-(k\cdot x+b) \right)=0
$$
то
$$
f(x)=(k\cdot x+b)+\alpha(x) \qquad \left(\text{где}\quad
\alpha(x)\underset{x\to +\infty}{\longrightarrow} 0 \right)
$$
поэтому
$$
\lim_{x\to +\infty}\frac{f(x)}{x}= \lim_{x\to +\infty}\left(
k+\frac{b}{x}+\frac{\alpha(x)}{x}\right)=k
$$
и, кроме того,
$$
\lim_{x\to +\infty} (f(x)-k\cdot x)= \lim_{x\to +\infty} ( b+\alpha(x) )=b
$$
Аналогично рассматривается случай $x\to -\infty$.
\end{proof}

\noindent\rule{160mm}{0.1pt}\begin{multicols}{2}

\subsection{Общая схема построения графика функции}

Построение графика функции удобно проводить по следующей схеме.

 \biter{
\item[1.] Находится область определения; обычно она состоит из интервалов (или
отрезков), концы которых удобно называть специальным термином, например, {\it
особыми точками}.\footnote {Более точное определение: {\it особой точкой}
функции $f$ называется любая граничная точка области определения функции $f$, то есть такая точка $a$, в любой окрестности которой
$(a-\varepsilon,a+\varepsilon)$ содержатся точки, в которых $f(x)$ определена и
точки, в которых  $f(x)$ не определена.}

\item[2.] Находятся пределы функции в особых точках и вертикальные асимптоты.

\item[3.] Находятся пределы на бесконечности и наклонные асимптоты.

\item[4.] Вычисляется первая производная, с помощью которой затем находятся
интервалы монотонности и точки экстремума.

\item[5.] Вычисляется вторая производная, с помощью которой затем находятся
интервалы выпуклости вверх и вниз и точки перегиба.

\item[6.] С учетом полученных данных строится график.
 }\eiter

\noindent Объясним эту схему на примерах.

\begin{ex} Построим график следующей функции:
$$
f(x)=\frac{x^2+1}{x-1}
$$

1. Сначала находим область определения. Она определяется неравенством $x-1\ne
0$, и поэтому имеет вид $(-\infty;1)\cup (1;+\infty)$. Точка $x=1$ является
особой точкой функции $f$.

2. Затем находим пределы функции в особых точках:
 \begin{multline*}\lim_{x\to 1-0} f(x)=\lim_{x\to 1-0}\frac{x^2+1}{x-1}=
\frac{1^2+1}{1-0-1}=\\=\left(\frac{2}{-0}\right)=-\infty
 \end{multline*}
 \begin{multline*}\lim_{x\to 1+0} f(x)=\lim_{x\to 1+0}\frac{x^2+1}{x-1}=
\frac{1^2+1}{1+0-1}=\\=\left(\frac{2}{+0}\right)=+\infty
 \end{multline*}
Из этих равенств следует, что прямая $x=1$ является вертикальной асимптотой
нашей функции.

3. Вычисляем пределы в бесконечности:
 \begin{multline*}\lim_{x\to +\infty} f(x)=\lim_{x\to +\infty}\frac{x^2+1}{x-1}=
\lim_{x\to +\infty}\frac{x+\frac{1}{x}}{1-\frac{1}{x}}=\\=
\left(\frac{+\infty+\frac{1}{+\infty}}{1-\frac{1}{+\infty}}\right)= +\infty
 \end{multline*}
 \begin{multline*}\lim_{x\to -\infty} f(x)=\lim_{x\to -\infty}\frac{x^2+1}{x-1}=
\lim_{x\to -\infty}\frac{x+\frac{1}{x}}{1-\frac{1}{x}}=\\=
\left(\frac{-\infty+\frac{1}{-\infty}}{1-\frac{1}{-\infty}}\right)= -\infty
 \end{multline*}
Ищем наклонные асимптоты по формулам \eqref{10.3.1} и \eqref{10.3.2}:
 \begin{multline*}
k_+=\lim_{x\to +\infty}\frac{f(x)}{x}= \lim_{x\to
+\infty}\frac{x^2+1}{x^2-x}=\\= \lim_{x\to
+\infty}\frac{1+\frac{1}{x^2}}{1-\frac{1}{x}}=1
 \end{multline*}
 \begin{multline*}
b_+=\lim_{x\to +\infty} [f(x)-k_+\cdot x]=\\= \lim_{x\to
+\infty}\left[\frac{x^2+1}{x-1}-1\cdot x \right]=\\= \lim_{x\to
+\infty}\frac{x^2+1-x(x-1)}{x-1}=\\= \lim_{x\to
+\infty}\frac{x^2+1-x^2+x}{x-1}= \lim_{x\to +\infty}\frac{1+x}{x-1}=\\=
\lim_{x\to +\infty}\frac{\frac{1}{x}+1}{1-\frac{1}{x}}=1
\end{multline*}
Значит, прямая
$$
y=x+1
$$
является асимптотой при $x\to +\infty$. Аналогично получаем, что $k_-=1$ и
$b_-=1$, и поэтому та же прямая $y=x+1$ является асимптотой при $x\to -\infty$.

4. Находим производную:
 \begin{multline*}
f'(x)=\left(\frac{x^2+1}{x-1}\right)'=\\= \frac{(x^2+1)'\cdot (x-1)-(x-1)'\cdot
(x^2+1)}{(x-1)^2} =\\= \frac{2x\cdot (x-1)-1\cdot (x^2+1)}{(x-1)^2} =
\frac{x^2-2x-1}{(x-1)^2}
 \end{multline*}
Решаем уравнение $f'(x)=0$:
 \begin{multline*}
f'(x)=0 \quad \Leftrightarrow \quad \frac{x^2-2x-1}{(x-1)^2}=0 \quad
\Leftrightarrow \\ \Leftrightarrow \quad
\begin{cases} {x^2-2x-1=0}\\ {(x-1)^2\ne 0}\end{cases}\quad \Leftrightarrow \\ \Leftrightarrow \quad x=1-\sqrt{2}\,\,
\text{или}\,\, x=1+\sqrt{2}
 \end{multline*}
Рисуем картинку, изображающую интервалы знакопостоянства производной, а значит,
интервалы монотонности функции $f$:

%\picture{0pt}{0pt}{123.pcx}

\vglue80pt

5. Находим вторую производную:
 \begin{multline*}
f''(x)=\left(\frac{x^2-2x-1}{(x-1)^2}\right)'=\\=\text{\smsize
$\frac{(x^2-2x-1)'\cdot (x-1)^2- ((x-1)^2)'\cdot (x^2-2x-1)}{(x-1)^4}$}=\\=
\frac{(2x-2)\cdot (x-1)^2- 2(x-1)\cdot (x^2-2x-1)}{(x-1)^4}=\\=
\frac{4}{(x-1)^3}\end{multline*} Решаем уравнение $f''(x)=0$:
$$
f''(x)=0 \quad \Leftrightarrow \quad \frac{4}{(x-1)^3}=0 \quad \Leftrightarrow
\quad x\in \varnothing
$$
Рисуем картинку, изображающую интервалы знакопостоянства второй производной, а
значит, интервалы выпуклости вверх и вниз функции $f$:

%\picture{0pt}{0pt}{124.pcx}

\vglue80pt

6. Рисуем график:

%\picture{0pt}{0pt}{125.pcx}

\vglue160pt
\end{ex}

\begin{ex} Построим график функции
$$
f(x)=\sqrt[3] {x^3-3x}
$$

1. Очевидно, функция определена всюду, поэтому область определения есть вся
числовая прямая $\R=(-\infty;+\infty)$.

2. Функция непрерывна везде на $\R=(-\infty;+\infty)$, поэтому вертикальных
асимптот нет.

3. Вычисляем пределы в бесконечности:
 \begin{multline*}\lim_{x\to +\infty} f(x)=\lim_{x\to +\infty}\sqrt[3] {x^3-3x}=\\=
\lim_{x\to +\infty} x\cdot \sqrt[3] {1-\frac{3}{x^2}}= (+\infty \cdot
1)=+\infty
 \end{multline*}
 \begin{multline*}\lim_{x\to -\infty} f(x)=\lim_{x\to -\infty}\sqrt[3] {x^3-3x}=\\=
\lim_{x\to -\infty} x\cdot \sqrt[3] {1-\frac{3}{x^2}}= (-\infty \cdot
1)=-\infty
 \end{multline*}
Ищем наклонные асимптоты по формулам \eqref{10.3.1} и \eqref{10.3.2}:
 \begin{multline*}
k_+=\lim_{x\to +\infty}\frac{f(x)}{x}= \lim_{x\to +\infty}\frac{\sqrt[3]
{x^3-3x}}{x}=\\= \lim_{x\to +\infty}\sqrt[3] {1-\frac{3}{x^2}}=1
 \end{multline*}
 \begin{multline*}
b_+=\lim_{x\to +\infty} [f(x)-k_+\cdot x]=\\= \lim_{x\to +\infty}\left[
\sqrt[3] {x^3-3x}-1\cdot x \right]=\\=\text{\tiny $\lim_{x\to
+\infty}\frac{\left(\sqrt[3] {x^3-3x}-x \right) \cdot \left(\sqrt[3]
{(x^3-3x)^2}+x\cdot \sqrt[3] {x^3-3x}+x^2 \right)} {\sqrt[3]
{(x^3-3x)^2}+x\cdot \sqrt[3] {x^3-3x}+x^2}$}=\\= \lim_{x\to +\infty}\frac{
\sqrt[3] {(x^3-3x)^3}-x^3 } {\sqrt[3] {(x^3-3x)^2}+x\cdot \sqrt[3]
{x^3-3x}+x^2}=\\= \lim_{x\to +\infty}\frac{ -3x } {\sqrt[3] {(x^3-3x)^2}+x\cdot
\sqrt[3] {x^3-3x}+x^2}=\\= \lim_{x\to +\infty}\frac{ -3} {\sqrt[3]
{(1-\frac{3}{x^2})\cdot (x^3-3x)}+ \sqrt[3] {x^3-3x}+x}=\\= \lim_{x\to
+\infty}\frac{ -3} {x\cdot \left(\sqrt[3] {(1-\frac{3}{x^2})^2}+ \sqrt[3]
{1-\frac{3}{x^2}}+1 \right)}=0
\end{multline*}
Значит, прямая $y=x$ является асимптотой при $x\to +\infty$. Аналогично
получаем, что $k_-=1$ и $b_-=0$, и поэтому та же прямая $y=x$ является
асимптотой при $x\to -\infty$.

4. Находим производную:
 \begin{multline*}
f'(x)=\left(\sqrt[3] {x^3-3x}\right)'= \frac{(x^3-3x)'}{3 \left(\sqrt[3]
{x^3-3x}\right)^2}=\\= \frac{3x^2-3}{3 \left(\sqrt[3] {x^3-3x}\right)^2}=
\frac{x^2-1}{\left(x^3-3x \right)^\frac{2}{3}}
 \end{multline*}
Решаем уравнение $f'(x)=0$:
 \begin{multline*}
f'(x)=0 \quad \Leftrightarrow \quad \frac{x^2-1}{\left(x^3-3x
\right)^\frac{2}{3}}=0 \quad \Leftrightarrow \\ \Leftrightarrow \quad
\begin{cases} {x^2-1=0}\\ {x^3-3x \ne 0}\end{cases}\quad \Leftrightarrow \quad x\in \{ -1;1 \}
 \end{multline*}
Рисуем картинку, изображающую интервалы знакопостоянства производной, а значит,
интервалы монотонности функции $f$:

%\picture{0pt}{0pt}{p9.pcx}

\vglue80pt

5. Находим вторую производную:
$$
f''(x)=\left(\frac{x^2-1}{\left(x^3-3x \right)^\frac{2}{3}}\right)'=...= -2
\frac{2x^4+x^2+1}{\left(x^3-3x \right)^\frac{5}{3}}
$$
Решаем уравнение $f''(x)=0$:
$$
f''(x)=0 \quad \Leftrightarrow \quad -2 \frac{2x^4+x^2+1}{\left(x^3-3x
\right)^\frac{5}{3}}=0 \; \Leftrightarrow \; x\in \varnothing
$$
Рисуем картинку, изображающую интервалы знакопостоянства второй производной, а
значит, интервалы выпуклости вверх и вниз функции $f$:

%\picture{0pt}{0pt}{p10.pcx}

\vglue80pt

6. Рисуем график:

%\picture{0pt}{-10pt}{p11.pcx}

\vglue170pt
\end{ex}

\begin{ers}
Постройте графики следующих функций:
 \biter{
 \item[1.] $y=\sqrt[3]{x^2-x^3}\quad$ ({\smsize здесь $y'=\frac{2-3x}{3
\sqrt[3]{x(1-x)^2}}$ и $y''=-\frac{2}{9 x^{4/3} (1-x)^{5/3}}$}).

 \item[2.] $y=\frac{x^2}{x-1}$.

 \item[3.] $y=\sqrt[5] {\frac{x^6}{x-a}}$.

 \item[4.] $y=\frac{x^3}{4(x-2)^2}$.

 \item[5.] $y=x\sqrt[3]{(x+1)^2}$.

 \item[6.] $y=\sqrt{\frac{x^3-2x^2}{x-3}}$.

 \item[7.] $y=\frac{x^2-4}{x}\cdot e^{-\frac{5}{3x}}$.

 \item[8.] $y=x^3+3x^2-9x-3$.

 \item[9.] $y=x+\sqrt{x^2-1}$.
 }\eiter
\end{ers}

\subsection{Графики с симметриями}

\paragraph{Четные и нечетные функции.}

Напомним, что четные и нечетные функции мы определили в
\ref{SEC-numb-function}\ref{func-s-simmetr} главы \ref{ch-R&N} как
удовлетворяющие тождествам
$$
f(-x)=f(x)\quad \text{(четность)},
$$
$$
f(-x)=-f(x)\quad \text{(нечетность)}
$$
Там же отмечалось, что график четной функции должен быть симметричен
относительно оси ординат, а график нечетной -- центрально симметричен
относительно начала координат.

Отсюда следует, что если нам нужно построить график четной функции, то его
достаточно построить на правой полуплоскости, а затем, чтобы получить
окончательный результат, отобразить полученную картинку налево симметрично
относительно оси ординат. Точно так же график нечетной функции достаточно
строить на правой полуплоскости, а потом отобразить центрально симметрично
относительно начала координат. Здесь мы рассмотрим примеры.

\begin{ex} Постройте график функции:
$$
f(x)=2x-\arcsin x
$$

1. Область определения: $[-1;1]$. Точки $x=-1$  и $x=1$ являются особыми.
Функция является нечетной, поэтому нам достаточно построить график только на
отрезке $x\in [0;1]$.

2. Предел в особой точке конечный
$$
\lim_{x\to 1-0} (2x-\arcsin x)=2-\arcsin 1=2-\frac{\pi}{2},
$$
значит вертикальных асимптот нет.

3. Наклонных асимптот тоже нет, поскольку $x$ не может стремиться к
бесконечности.

4. Производная:
$$
f'(x)=(2x-\arcsin x)'=2-\frac{1}{\sqrt{1-x^2}}
$$
При $x\ge 0$ получаем $f'(x)=0 \, \Leftrightarrow \, 2-\frac{1}{\sqrt{1-x^2}}=0
\, \Leftrightarrow \, 2=\frac{1}{\sqrt{1-x^2}}\, \Leftrightarrow \,
\sqrt{1-x^2}=\frac{1}{2}\, \Leftrightarrow \, 1-x^2=\frac{1}{4}\,
\Leftrightarrow \, x^2=\frac{3}{4}\, \Leftrightarrow \, x=\frac{\sqrt{3}}{2} $.
Интервалы знакопостоянства производной:

%\picture{0pt}{0pt}{p18.pcx}

\vglue90pt

5. Вторая производная:
 \begin{multline*}
f''(x)=\Big( 2-\frac{1}{\sqrt{1-x^2}}\Big)'= \Big(
2-(1-x^2)^{-\frac{1}{2}}\Big)'=\\= \frac{1}{2}(1-x^2)^{-\frac{3}{2}}\cdot
(-2x)= -x\cdot (1-x^2)^{-\frac{3}{2}}
 \end{multline*}
При $x\ge 0$ получаем $f''(x)=0 \, \Leftrightarrow \, -x\cdot
(1-x^2)^{-\frac{3}{2}}=0
 \, \Leftrightarrow \, x=0$.
Интервалы знакопостоянства второй производной:

%\picture{0pt}{0pt}{p19.pcx}

\vglue90pt

6. График функции для $x\ge 0$:

%\picture{0pt}{0pt}{p20.pcx}

\vglue110pt \noindent Отобразив эту картинку симметрично относительно начала
координат, мы получим полный график нашей функции:

%\picture{0pt}{0pt}{p21.pcx}

\vglue110pt
\end{ex}

\paragraph{Периодические функции.}

Напомним, что периодические функции мы определили в
\ref{SEC-numb-function}\ref{func-s-simmetr} главы \ref{ch-R&N} как
удовлетворяющие тождеству
$$
f(x+T)=f(x)
$$
при некотором $T>0$. Там же отмечалось, что график периодической функции $f$
должен накладываться на себя при сдвиге вправо (и влево) на число $T$ (период).

Отсюда следует, что если нам нужно построить график периодической функции, то
его достаточно построить на отрезке вида $[a;a+T]$, а затем, чтобы получить
окончательный результат, сдвинуть полученную картинку влево и вправо нужное
число раз (чтобы зрителю стало понятно, что сдвигать можно сколько угодно раз с
тем же результатом). Здесь мы рассмотрим примеры.

\begin{ex} Постройте график функции:
$$
f(x)=\sin^3 x +\cos^3 x
$$

1. Область определения: $(-\infty;+\infty)$. Функция является периодической с
периодом $T=2\pi$, поэтому нам достаточно построить график только на отрезке
$x\in [0;2\pi]$.

2. Особых точек нет, значит и вертикальных асимптот нет.

3. Наклонные асимптоты мы не ищем, поскольку нас интересует только отрезок
$x\in [0;2\pi]$.

4. Производная:
 \begin{multline*}
f'(x)=3\sin^2 x \cos x -3 \cos^2 x \sin x=\\=3\sin x \cos x (\sin x - \cos x)
 \end{multline*}
При $x\in [0;2\pi]$ получаем
 \begin{multline*}
f'(x)=0 \, \Leftrightarrow \, 3\sin x \cos x (\sin x - \cos x)=0 \,
\Leftrightarrow \\
\text{\smsize $\Leftrightarrow  \, \left[\aligned & \sin x=0
\\
& \cos x=0
\\
& \sin x = \cos x
\endaligned\right]
\, \Leftrightarrow \, \left[\aligned & x\in \{ 0; \pi \}\\
& x\in \left\{\frac{\pi}{2}; \frac{3\pi}{2}\right\}\\
& x\in \left\{\frac{\pi}{4}; \frac{5\pi}{4}\right\}\endaligned\right] \,
\Leftrightarrow$}\\ \Leftrightarrow  \, x\in \left\{0; \frac{\pi}{4};
\frac{\pi}{2}; \pi; \frac{5\pi}{4}; \frac{3\pi}{2}\right\}
 \end{multline*}
Интервалы знакопостоянства производной:

%\picture{0pt}{0pt}{p23.pcx}

\vglue90pt

5. Вторая производная:
 \begin{multline*}
f''(x)= (3\sin x \cos x (\sin x - \cos x))'=\\= 3\cos^2 x (\sin x - \cos
x)-3\sin^2 x (\sin x - \cos x)+\\+3\sin x \cos x (\cos x + \sin x) =\\=
3(\cos^2 x -\sin^2 x) (\sin x - \cos x)+\\+ 3\sin x \cos x (\cos x + \sin
x)=\\= -3(\cos x +\sin x) (\cos x - \sin x)^2+\\+ 3\sin x \cos x (\cos x + \sin
x)=\\= 3(\cos x +\sin x) [-(\cos x - \sin x)^2+\sin x \cos x]=\\= 3(\cos x
+\sin x) [-\cos^2 x - \sin^2 x+3\sin x \cos x]=\\= 3(\cos x +\sin x) [3\sin x
\cos x-1]
 \end{multline*}
При $x\in [0;2\pi]$ получаем
 \begin{multline*}
f''(x)=0 \, \Leftrightarrow\\ \Leftrightarrow \, 3(\cos x +\sin x) [3\sin x
\cos x-1]=0
 \, \Leftrightarrow \\
\text{\smsize $\Leftrightarrow\, \left[\aligned & \cos x +\sin x=0
\\
& 3\sin x \cos x-1=0
\endaligned\right]
 \, \Leftrightarrow \,
\left[\aligned & \cos x =-\sin x
\\
& 3\sin x \cos x=1
\endaligned\right]
\Leftrightarrow$}\\
\text{\smsize $\Leftrightarrow\, \left[\aligned & \cos x =-\sin x
\\
& \frac{3}{2}\sin 2x =1
\endaligned\right]\,
\Leftrightarrow\, \left[\aligned & \cos x =-\sin x
\\
& \sin 2x =\frac{2}{3}\endaligned\right] \Leftrightarrow$}\\
\text{\smsize $\Leftrightarrow\,
 \left[\aligned & x\in \left\{ \frac{3\pi}{4};
\frac{7\pi}{4}\right\}\\
& 2x=\arcsin \frac{2}{3}+2\pi n
\\
& 2x=\pi-\arcsin \frac{2}{3}+2\pi n
\endaligned\right]
\Leftrightarrow$}\\
\text{\smsize $\Leftrightarrow\,\left[\aligned & x\in \left\{
\frac{3\pi}{4}; \frac{7\pi}{4}\right\}\\
& x=\frac{1}{2}\arcsin \frac{2}{3}+\pi n
\\
& x=\frac{\pi}{2}-\frac{1}{2}\arcsin \frac{2}{3}+\pi n
\endaligned\right]\, \Leftrightarrow$}\\
 \text{\smsize $\Leftrightarrow\,
 {\smsize\begin{pmatrix}
 \text{оставляем только}\\
 \text{те значения $x$, которые}\\
 \text{лежат в отрезке $[0; 2\pi]$}\end{pmatrix}}
 \,\Leftrightarrow$}\\
\Leftrightarrow x\in \Big\{ \frac{1}{2}\arcsin \frac{2}{3};
\frac{\pi}{2}-\frac{1}{2}\arcsin \frac{2}{3}; \frac{3\pi}{4};\\
\pi+\frac{1}{2}\arcsin \frac{2}{3}; \frac{3\pi}{2}-\frac{1}{2}\arcsin
\frac{2}{3}; \frac{7\pi}{4}\Big\}
 \end{multline*}

Интервалы знакопостоянства второй производной:

%\picture{0pt}{0pt}{p24.pcx}

\vglue90pt

6. График функции на отрезке $x\in [0;2\pi]$:

%\picture{0pt}{0pt}{p25.pcx}

\vglue120pt \noindent Сдвинув эту картинку несколько раз влево и вправо на
период $T=2\pi$, мы получим полный график нашей функции:

%\picture{0pt}{0pt}{p26.pcx}

\vglue120pt
\end{ex}

\begin{ers}
Постройте графики следующих функций:

1. $y=\frac{x^2+2}{2x^2+1}$.

2. $y=\sqrt[3] {(x+1)^2} - \sqrt[3] {(x-1)^2}$.

3. $y=e^{\sin x+\cos x}$.

4. $y=\sqrt[3] {(x+1)^2} + \sqrt[3] {(x-1)^2}$.

5. $y=\frac{-8x}{x^2+4}$.

6. $y=\arctg (\sin x)$.

7. $y=\frac{x^2}{\sqrt{x^2+2}}$.

8. $y=x \ln |x|$.

9. $y=\ln (\sin x+\cos x)$.

10. $y=\frac{x^2+1}{\sqrt{4x^2-3}}$.

11. $y=2x-\arctg x$.

12. $y=\frac{1}{\sin x-\cos x}$.
\end{ers}

\end{multicols}\noindent\rule[10pt]{160mm}{0.1pt}

\section{Правило Лопиталя}\label{SEC:Lopital}

Теоремы о дифференцируемых функциях, обсуждавшиеся в этой главе,
позволяют доказать важное {\it правило Лопиталя}, с помощью которого удается
вычислять ``сложные'' пределы (которые иначе сосчитать бывает невозможно).

\subsection{Раскрытие неопределенностей типа $\frac{0}{0}$}

\begin{tm}
[правило Лопиталя для предела в точке с неопределенностью $\frac{0}{0}$]
\label{Lopital_x->a_0/0}\footnote{Эта теорема используется в главе
\ref{ch-o(f(x))} при доказательстве асимптотической формулы Тейлора-Пеано
\ref{Taylor-Peano}}

Пусть функции $f$ и $g$ обладают следующими свойствами:
 \bit{
\item[(i)] $f$ и $g$ определены и дифференцируемы в некоторой выколотой
окрестности $(a-\varepsilon,a)\cup (a,a+\varepsilon)$ точки $a$,

\item[(ii)] $\lim\limits_{x\to a} f(x)=0$ и $\lim\limits_{x\to a} g(x)=0$,

\item[(iii)] $g'(x)\ne 0$ при $x\ne a$,

\item[(iv)] существует предел $\lim\limits_{x\to a}\frac{f'(x)}{g'(x)}$
(конечный или бесконечный).
 }\eit
Тогда существует предел $\lim\limits_{x\to a}\frac{f(x)}{g(x)}$, причем
\begin{equation}\label{9.1.1}
\boxed{ \lim_{x\to a}\frac{f(x)}{g(x)}=\lim_{x\to a}\frac{f'(x)}{g'(x)} }
\end{equation}
 \end{tm}
 \begin{proof}
Чтобы доказать \eqref{9.1.1}, нужно
взять произвольную последовательность $x_n\underset{n\to
\infty}{\longrightarrow} a, \quad x_n\ne a$, и убедиться, что
$$
\lim_{n\to \infty}\frac{f(x_n)}{g(x_n)}=\lim_{x\to a}\frac{f'(x)}{g'(x)}
$$
Доопределим функции $f$ и $g(x)$ в точке $a$ формулами
\begin{equation}\tilde f(x)=\begin{cases} f(x), & \text{если}\, x\in (a,b)
\\ 0, &  \text{если}\, x=a   \end{cases}\qquad \tilde g(x)=\begin{cases} g(x), & \text{если}\, x\in (a,b)
\\ 0, & \text{если}\, x=a   \end{cases}\label{9.1.2}\end{equation}
Тогда на каждом отрезке $[a,x_n]$ функции $\tilde f(x)$ и $\tilde g(x)$
удовлетворяют условиям теоремы Коши \ref{Cauchy-II}:
 \bit{
\item[1)] $\tilde f(x)$ и $\tilde g(x)$ непрерывны на отрезке $[a,x_n]$,
\item[2)] $\tilde f(x)$ и $\tilde g(x)$ дифференцируемы на интервале $(a,x_n)$,
\item[3)] $g'(x)\ne 0$ при $x\in (a,x_n)$.
 }\eit
Значит, по теореме Коши \ref{Cauchy-II}, существует точка $\xi_n\in (a,x_n)$,
такая что
$$
\frac{\tilde f(x_n)-\tilde f(a)}{\tilde g(x_n)-\tilde
g(a)}=\frac{f'(\xi_n)}{g'(\xi_n)}\qquad \xi_n\in (a,x_n)
$$
Но, по формуле \eqref{9.1.2}, $\tilde f(a)=\tilde g(a)=0$, поэтому
$$
\frac{\tilde f(x_n)}{\tilde g(x_n)}=\frac{f'(\xi_n)}{g'(\xi_n)}, \qquad
\xi_n\in (a,x_n)
$$
Теперь устремим $n$ к бесконечности:

\begin{multline*}\lim_{n\to \infty}\frac{f(x_n)}{g(x_n)}= \lim_{n\to \infty}\frac{\tilde f(x_n)}{\tilde g(x_n)}= \lim_{n\to \infty}\frac{f'(\xi_n)}{g'(\xi_n)}= {\smsize\begin{pmatrix}\xi_n\in (a,x_n),
\\
\text{поэтому}\\
\xi_n\underset{n\to \infty}{\longrightarrow} a
\end{pmatrix}}= \lim_{\xi\to a}\frac{f'(\xi)}{g'(\xi)}=\\=
{\smsize\begin{pmatrix}\text{неважно, какая}\\
\text{буква стоит под}\\
\text{знаком предела}\end{pmatrix}}= \lim_{x\to
a}\frac{f'(x)}{g'(x)}\end{multline*} Как раз это нам и нужно было проверить.
\end{proof}

\begin{tm}[\bf правило Лопиталя для предела в бесконечности
с неопределенностью $\frac{0}{0}$]\label{Lopital_x->infty_0/0}

Пусть функции $f$ и $g(x)$ обладают следующими свойствами:
 \bit{
\item[(i)] $f(x)$ и $g(x)$ определены и дифференцируемы на некотором множестве
$(-\infty,-E)\cup (E;+\infty)$, \item[(ii)] $\lim\limits_{x\to \infty} f(x)=0$
и $\lim\limits_{x\to \infty} g(x)=0$, \item[(iii)] $g'(x)\ne 0$ при $x\in
(-\infty,-E)\cup (E;+\infty)$, \item[(ii)] существует предел $\lim\limits_{x\to
\infty}\frac{f'(x)}{g'(x)}$ (конечный или бесконечный).
 }\eit
Тогда существует предел $\lim\limits_{x\to \infty}\frac{f(x)}{g(x)}$, причем
\begin{equation}\boxed{ \lim_{x\to \infty}\frac{f(x)}{g(x)}=\lim_{x\to \infty}\frac{f'(x)}{g'(x)} }\label{9.1.3}\end{equation}\end{tm}\begin{proof} Вычислим этот предел с помощью
теоремы \ref{Lopital_x->infty_0/0}:
\begin{multline*}\lim_{x\to \infty}\frac{f(x)}{g(x)}= \left|
\begin{array}{c}
x=\frac{1}{t}, \, t=\frac{1}{x}\\
t\underset{x\to \infty}{\longrightarrow} 0
\end{array}\right|= \lim_{t\to 0}\frac{f(\frac{1}{t})}{g(\frac{1}{t})}=
\left(\text{применяем теорему \ref{Lopital_x->infty_0/0}}\right)= \lim_{t\to
0}\frac{\left(f(\frac{1}{t}) \right)'}{\left(g(\frac{1}{t})\right)'}=\\=
\lim_{t\to 0}\frac{f'(\frac{1}{t}) \cdot \left( -\frac{1}{t^2}\right)}
{g'(\frac{1}{t}) \cdot \left( -\frac{1}{t^2}\right)}= \lim_{t\to
0}\frac{f'(\frac{1}{t})}{g'(\frac{1}{t})}= \left|
\begin{array}{c}
t=\frac{1}{x}, \, x=\frac{1}{t},
\\
x\underset{t\to 0}{\longrightarrow}\infty
\end{array}\right|= \lim_{x\to \infty}\frac{f'(x)}{g'(x)}  \end{multline*}\end{proof}

\noindent\rule{160mm}{0.1pt}\begin{multicols}{2}

\begin{ex} Вычислим предел
$$
\lim_{x\to 0}\frac{\sin 2x}{\arctg 3x}
$$
Для этого можно воспользоваться теоремой \ref{Lopital_x->a_0/0}:
 \begin{multline*}\lim_{x\to 0}\frac{\sin 2x}{\arctg 3x}= \left(\frac{0}{0}\right)={\smsize {\smsize\begin{pmatrix}\text{применяем}\\
\text{теорему \ref{Lopital_x->a_0/0}}\end{pmatrix}}}=\\= \lim_{x\to
0}\frac{(\sin 2x)'}{(\arctg 3x)'}= \lim_{x\to 0}\frac{\cos 2x \cdot
2}{\frac{3}{1+9x^2}}= \frac{\cos 0 \cdot 2}{\frac{3}{1+0}}=\\= \frac{1 \cdot
2}{3}=\frac{2}{3}
 \end{multline*}
 \end{ex}

\begin{ex} Докажем еще раз формулу \eqref{5.10.3}:
$$
\lim_{x\to 0}\frac{e^x-1}{x}=1
$$
Действительно,
 \begin{multline*}\lim_{x\to 0}\frac{e^x-1}{x}= \left(\frac{0}{0}\right)= {\smsize
{\smsize\begin{pmatrix}\text{применяем}\\
\text{теорему \ref{Lopital_x->a_0/0}}\end{pmatrix}}}=\\= \lim_{x\to
0}\frac{(e^x-1)'}{(x)'}= \lim_{x\to 0}\frac{e^x}{1}=1
 \end{multline*}\end{ex}

\begin{ex} Правило Лопиталя можно применять несколько раз подряд.
Покажем это на следующем примере:
 \begin{multline*}\lim_{x\to 0}\frac{1-\cos x}{x^2}= \left(\frac{0}{0}\right)=
{\smsize {\smsize\begin{pmatrix}\text{применяем}\\
\text{теорему \ref{Lopital_x->a_0/0}}\end{pmatrix}}}=\\= \lim_{x\to
0}\frac{(1-\cos x)'}{(x^2)'}= \lim_{x\to 0}\frac{\sin x}{2x}=
\left(\frac{0}{0}\right)=\\={\smsize {\smsize\begin{pmatrix}\text{снова}\\ \text{применяем}\\
\text{теорему \ref{Lopital_x->a_0/0}}\end{pmatrix}}}= \lim_{x\to 0}\frac{(\sin
x)'}{(2x)'}= \lim_{x\to 0}\frac{\cos x}{2}= \frac{1}{2}\end{multline*}\end{ex}

\begin{ex} Еще один пример, где правило Лопиталя
применяется несколько раз:
 \begin{multline*}\lim_{x\to 0}\frac{\ln \cos x}{x^2}= \left(\frac{0}{0}\right)={\smsize {\smsize\begin{pmatrix}\text{применяем}\\
\text{теорему \ref{Lopital_x->a_0/0}}\end{pmatrix}}}=\\= \lim_{x\to
0}\frac{(\ln \cos x)'}{(x^2)'}= \lim_{x\to 0}\frac{\frac{1}{\cos x}\cdot (-\sin
x)}{2x}=\\= \lim_{x\to 0}\frac{-\sin x}{2x\cos x}=
\left(\frac{0}{0}\right)={\smsize
{\smsize\begin{pmatrix}\text{снова}\\ \text{применяем}\\
\text{теорему \ref{Lopital_x->a_0/0}}\end{pmatrix}}}=\\= \lim_{x\to
0}\frac{(-\sin x)'}{(2x \cos x)'}= \lim_{x\to 0}\frac{-\cos x}{2\cos x-2x\sin
x}=\\= \frac{-1}{2-0}= -\frac{1}{2}\end{multline*}\end{ex}

\begin{ex} Теперь пример, где применяется теорема
\ref{Lopital_x->infty_0/0}:
 \begin{multline*}\lim_{x\to +\infty}\frac{\frac{\pi}{2}-\arctg x} {\ln \left(
1+\frac{1}{x^2}\right)}= \left(\frac{0}{0}\right)={\smsize
{\smsize\begin{pmatrix}\text{применяем}\\
\text{теорему \ref{Lopital_x->infty_0/0}}\end{pmatrix}}}=\\= \lim_{x\to
+\infty}\frac{(\frac{\pi}{2}-\arctg x)'} {\left(\ln \left(
1+\frac{1}{x^2}\right) \right)'}= \lim_{x\to +\infty}\frac{-\frac{1}{1+x^2}}
{\frac{1}{1+\frac{1}{x^2}}\cdot \left( -\frac{2}{x^3}\right) }=\\= \lim_{x\to
+\infty}\frac{1}{1+x^2}\cdot \frac{x^2+1}{x^2}\cdot \frac{x^3}{2} = \lim_{x\to
+\infty}\frac{x}{2} =+\infty
\end{multline*}\end{ex}

\begin{ers}
Вычислите пределы:

1. $\lim\limits_{x\to 1}\frac{\sin 3\pi x}{\sin 2\pi x}.$

2. $\lim\limits_{x\to 0}\frac{\tg x - x}{\sin x - x^2}.$

3. $\lim\limits_{x\to 0}\frac{\tg x - x}{\sin x - x}.$
\end{ers}

\end{multicols}\noindent\rule[10pt]{160mm}{0.1pt}

\subsection{Раскрытие неопределенностей типа $\frac{\infty}{\infty}$}

\begin{tm}[\bf правило Лопиталя для предела в точке с неопределенностью
$\frac{\infty}{\infty}$]\label{Lopital_x->a_infty/infty}

Пусть функции $f$ и $g(x)$ обладают следующими свойствами:
 \bit{
\item[(i)] $f(x)$ и $g(x)$ определены и дифференцируемы в некоторой окрестности
точки $a$,

\item[(ii)] $\lim\limits_{x\to a} f(x)=\infty$ и  $\lim\limits_{x\to a}
g(x)=\infty$,

\item[(iii)] $g'(x)\ne 0$ в некоторой выколотой окрестности точки $a$,

\item[(iv)] существует предел $\lim\limits_{x\to a}\frac{f'(x)}{g'(x)}$
(конечный или бесконечный).
 }\eit
Тогда существует предел $\lim\limits_{x\to a}\frac{f(x)}{g(x)}$, причем
 \begin{equation}\label{9.2.1}
\boxed{ \lim_{x\to a}\frac{f(x)}{g(x)}=\lim_{x\to a}\frac{f'(x)}{g'(x)} }
 \end{equation}
 \end{tm}
\begin{proof} 1. Рассмотрим сначала случай,
когда предел $\lim\limits_{x\to a}\frac{f'(x)}{g'(x)}$ конечен:
 \begin{equation}
\lim_{x\to a}\frac{f'(x)}{g'(x)}=A\in \R \label{9.2.2}
 \end{equation}
Зафиксируем интервал $(a-\delta,a)$, на котором $g'(x)\ne 0$. В силу теоремы
\ref{Roll-cons}, функция $g(x)$ должна быть строго монотонна на $(a-\delta,a)$.
Поэтому, если взять возрастающую последовательность
$$
x_1<x_2<...<x_n<...<a, \quad x_n\underset{n\to\infty}{\longrightarrow} a
$$
то мы получим, что $g(x_n)$ -- строго монотонная бесконечно большая
последовательность:
$$
g(x_n)\underset{n\to\infty}{\longrightarrow} +\infty, \qquad
g(x_1)<g(x_2)<...<g(x_n)<... \quad \Big(\text{или}\quad
g(x_1)>g(x_2)>...>g(x_n)>... \Big)
$$
Положив
$$
y_n=g(x_n), \quad z_n=f(x_n)
$$
мы теперь получим, что последовательности $y_n$ и $z_n$ удовлетворяют всем
условиям теоремы Штольца \ref{Stoltz}: во первых,
$$
y_n\underset{n\to\infty}{\longrightarrow}\infty, \qquad y_1<y_2<...<y_n<...
\quad \Big(\text{или}\quad y_1>y_2>...>y_n>...\Big)
$$
и во вторых,

\begin{multline*}\lim_{n\to\infty}\frac{z_n-z_{n-1}}{y_n-y_{n-1}}=
\lim_{n\to\infty}\frac{f(x_n)-f(x_{n-1})}{g(x_n)-g(x_{n-1})}=\\= {\smsize\begin{pmatrix}\text{по теореме Коши \ref{Cauchy-II},}\\
\exists \xi_n\in (x_{n-1},x_n) \quad
\frac{f(x_n)-f(x_{n-1})}{g(x_n)-g(x_{n-1})}=\frac{f'(\xi_n)}{g'(\xi_n)}\end{pmatrix}}
=\lim_{n\to\infty}\frac{f'(\xi_n)}{g'(\xi_n)}=(2.2)=A
\end{multline*} Значит, по теореме Штольца \ref{Stoltz},
$$
\lim_{n\to\infty}\frac{f(x_n)}{g(x_n)}= \lim_{n\to\infty}\frac{z_n}{y_n}=
\lim_{n\to\infty}\frac{z_n-z_{n-1}}{y_n-y_{n-1}}=A
$$
Это верно для всякой возрастающей последовательности
$x_n\underset{n\to\infty}{\longrightarrow} a$, значит, по теореме
\ref{lim-monot},
$$
\lim_{x\to a-0}\frac{f(x)}{g(x)}=A
$$
Аналогично получаем
$$
\lim_{x\to a+0}\frac{f(x)}{g(x)}=A
$$

2. Рассмотрим теперь случай, когда
$$
\lim\limits_{x\to a}\frac{f'(x)}{g'(x)}=\infty
$$
Тогда $\frac{f'(x)}{g'(x)}\ne 0$ в некоторой выколотой окрестности точки $a$.
(Это следует, например, из определения предела по Коши.) Значит, $f'(x)\ne 0$ в
некоторой выколотой окрестности точки $a$. Таким образом, получается:
 \bit{
\item[1)] $f'(x)\ne 0$ в некоторой выколотой окрестности точки $a$, и \item[2)]
$\lim\limits_{x\to a}\frac{g'(x)}{f'(x)}=0$
 }\eit
\noindent То есть выполняется посылка теоремы, но только с переставленными
символами $f$ и $g$, и с условием, что предел отношения производных конечен.
Поскольку для случая конечного предела теорема уже доказана, имеем
$$
\lim_{x\to a}\frac{g(x)}{f(x)}=0=\lim_{x\to a}\frac{g'(x)}{f'(x)}
$$
откуда
$$
\lim_{x\to a}\frac{f(x)}{g(x)}=\infty=\lim_{x\to a}\frac{f'(x)}{g'(x)}\qquad $$
\end{proof}

\begin{tm}[\bf правило Лопиталя для предела в бесконечности
с неопределенностью $\frac{\infty}{\infty}$]
\label{Lopital_x->infty_infty/infty}

Пусть функции $f$ и $g(x)$ обладают следующими свойствами:
 \bit{
\item[(i)] $f(x)$ и $g(x)$ определены и дифференцируемы на некотором множестве
$(-\infty,-E)\cup (E,+\infty)$, \item[(ii)] $\lim\limits_{x\to \infty}
f(x)=\infty$ и $\lim\limits_{x\to \infty} g(x)=\infty$, \item[(iii)] $g'(x)\ne
0$ при $x\in (-\infty,-E)\cup (E,+\infty)$, \item[(ii)] существует предел
$\lim\limits_{x\to \infty}\frac{f'(x)}{g'(x)}$ (конечный или бесконечный).
 }\eit
Тогда существует предел $\lim\limits_{x\to \infty}\frac{f(x)}{g(x)}$, причем
\begin{equation}\boxed{ \lim_{x\to \infty}\frac{f(x)}{g(x)}=\lim_{x\to \infty}\frac{f'(x)}{g'(x)} }\label{9.2.3}\end{equation}\end{tm}\begin{proof} Здесь используется тот же прием,
что и в доказательстве теоремы \ref{Lopital_x->infty_0/0}:
\begin{multline*}\lim_{x\to \infty}\frac{f(x)}{g(x)}= \left|
\begin{array}{c}
x=\frac{1}{t}, \, t=\frac{1}{x}\\
t\underset{x\to \infty}{\longrightarrow} 0
\end{array}\right|= \lim_{t\to 0}\frac{f(\frac{1}{t})}{g(\frac{1}{t})}=
\left(\text{применяем теорему \ref{Lopital_x->a_infty/infty}}\right)=
\lim_{t\to 0}\frac{\left(f(\frac{1}{t})
\right)'}{\left(g(\frac{1}{t})\right)'}=\\= \lim_{t\to 0}\frac{f'(\frac{1}{t})
\cdot \left( -\frac{1}{t^2}\right)} {g'(\frac{1}{t}) \cdot \left(
-\frac{1}{t^2}\right)}= \lim_{t\to 0}\frac{f'(\frac{1}{t})}{g'(\frac{1}{t})}=
\left|
\begin{array}{c}
t=\frac{1}{x}, \, x=\frac{1}{t},
\\
x\underset{t\to 0}{\longrightarrow}\infty
\end{array}\right|= \lim_{x\to \infty}\frac{f'(x)}{g'(x)}
 \end{multline*}\end{proof}

\noindent\rule{160mm}{0.1pt}\begin{multicols}{2}

\begin{ex} Вычислим предел
$$
\lim_{x\to +0}\frac{\ctg x}{\ln x}
$$
Для этого воспользуемся теоремами \ref{Lopital_x->a_infty/infty} и
\ref{Lopital_x->a_0/0}:
 \begin{multline*}\lim_{x\to +0}\frac{\ctg x}{\ln x}= \left(\frac{\infty}{\infty}\right)={\smsize {\smsize\begin{pmatrix}\text{применяем}\\
\text{теорему \ref{Lopital_x->a_infty/infty}}\end{pmatrix}}}=\\= \lim_{x\to
+0}\frac{(\ctg x)'}{(\ln x)'}= \lim_{x\to +0}\frac{-\frac{1}{\sin^2
x}}{\frac{1}{x}}=\\= -\lim_{x\to +0}\frac{x}{\sin^2 x}=
\left(\frac{0}{0}\right)={\smsize
{\smsize\begin{pmatrix}\text{применяем}\\
\text{теорему \ref{Lopital_x->a_0/0}}\end{pmatrix}}}=\\= -\lim_{x\to
+0}\frac{(x)'}{(\sin^2 x)'}= -\lim_{x\to +0}\frac{1}{2 \sin x \cos x}=\\=
\left( -\frac{1}{+0}\right)=-\infty \end{multline*}\end{ex}

\begin{ex} Вычислим предел
$$
\lim_{x\to +\infty}\frac{\ln x}{x}
$$
Для этого воспользуемся теоремой \ref{Lopital_x->infty_infty/infty}:
 \begin{multline*}\lim_{x\to +\infty}\frac{\ln x}{x}= \left(\frac{\infty}{\infty}\right)={\smsize {\smsize\begin{pmatrix}\text{применяем}\\
\text{теорему \ref{Lopital_x->infty_infty/infty}}\end{pmatrix}}}=\\= \lim_{x\to
+\infty}\frac{(\ln x)'}{(x)'}= \lim_{x\to +\infty}\frac{\frac{1}{x}}{1}=+\infty
 \end{multline*}\end{ex}

\begin{ex}
Вычислим предел
 \begin{multline*}\lim_{x\to +\infty}\frac{x^3}{e^x}= \left(\frac{\infty}{\infty}\right)={\smsize {\smsize\begin{pmatrix}\text{применяем}\\
\text{теорему \ref{Lopital_x->infty_infty/infty}}\end{pmatrix}}}=\\= \lim_{x\to
+\infty}\frac{(x^3)'}{(e^x)'}= \lim_{x\to +\infty}\frac{3x^2}{e^x}=
\left(\frac{\infty}{\infty}\right)=\\={\smsize {\smsize\begin{pmatrix}\text{снова}\\ \text{применяем}\\
\text{теорему \ref{Lopital_x->infty_infty/infty}}\end{pmatrix}}}= \lim_{x\to
+\infty}\frac{(3x^2)'}{(e^x)'}=\\= \lim_{x\to +\infty}\frac{6x}{e^x}=
\left(\frac{\infty}{\infty}\right)= {\smsize
{\smsize\begin{pmatrix}\text{еще раз}\\ \text{применяем}\\
\text{теорему \ref{Lopital_x->infty_infty/infty}}\end{pmatrix}}}=\\= \lim_{x\to
+\infty}\frac{(6x)'}{(e^x)'}= \lim_{x\to +\infty}\frac{6}{e^x}=
\left(\frac{6}{\infty}\right)=0
\end{multline*}\end{ex}

\begin{ers}
Вычислите пределы:

1. $\lim\limits_{x\to \infty}\frac{\ln (1+x^2)}{\ln \left(\frac{\pi}{2}-\arctg
x \right)}$.

2. $\lim\limits_{x\to a}\frac{\ln (x-a)}{\ln (e^x-e^a)}$.

3. $\lim\limits_{x\to +\infty}\frac{e^x+x^2+ x}{x^3-e^x}$.
\end{ers}

\end{multicols}\noindent\rule[10pt]{160mm}{0.1pt}

\subsection{Раскрытие неопределенностей $0\cdot \infty$,
$\infty-\infty$, $0^0$, $1^\infty$, $\infty^0$}

Покажем на примерах, как с помощью правил Лопиталя можно вычислять пределы с
неопределенностями $0\cdot \infty$, $\infty-\infty$, $0^0$, $1^\infty$,
$\infty^0$.

\noindent\rule{160mm}{0.1pt}\begin{multicols}{2}

\begin{ex}[\bf неопределенность вида $0\cdot \infty$] Вычислим предел
 \begin{multline*}\lim_{x\to +0} x\cdot \ln x= \left( 0\cdot \infty \right)= \lim_{x\to
+0}\frac{\ln x}{\frac{1}{x}}=
\left(\frac{\infty}{\infty}\right)=\\={\smsize {\smsize\begin{pmatrix}\text{применяем}\\
\text{теорему \ref{Lopital_x->a_infty/infty}}\end{pmatrix}}}= \lim_{x\to
+0}\frac{(\ln x)'}{(\frac{1}{x})'}= \lim_{x\to
+0}\frac{\frac{1}{x}}{-\frac{1}{x^2}}=\\= -\lim_{x\to +0} x=0
\end{multline*}\end{ex}

\begin{ex}[\bf неопределенность вида $\infty-\infty$] Вычислим предел
 \begin{multline*}\lim_{x\to \frac{\pi}{2}}\left(\frac{1}{\cos x}-\tg x \right)=
\left(\infty-\infty \right)=\\= \lim_{x\to \frac{\pi}{2}}\frac{1-\sin
x}{\cos x}= \left(\frac{0}{0}\right)={\smsize {\smsize\begin{pmatrix}\text{применяем}\\
\text{теорему \ref{Lopital_x->a_0/0}}\end{pmatrix}}}=\\= \lim_{x\to
\frac{\pi}{2}}\frac{(1-\sin x)'}{(\cos x)'}= \lim_{x\to
\frac{\pi}{2}}\frac{-\cos x}{-\sin x}=0
\end{multline*}\end{ex}

\begin{ex}[\bf неопределенность вида $0^0$] Вычислим предел
 \begin{multline*}\lim_{x\to +0} x^x= \left( 0^0 \right)= \lim_{x\to +0} e^{\ln (x^x)}=
\lim_{x\to +0} e^{x\ln x}=\\= e^{\lim_{x\to +0} x\ln x}={\smsize
{\smsize\begin{pmatrix}\text{применяем}\\
\text{пример 3.1}\end{pmatrix}}}= e^0=1
 \end{multline*}\end{ex}

\begin{ex}[\bf неопределенность вида $1^\infty$] Вычислим предел
 \begin{multline*}\lim_{x\to 0} {\cos x}^{\ctg^2 x}= \left( 1^\infty \right)=
\lim_{x\to 0} e^{\ln \left( {\cos x}^{\ctg^2 x}\right)}=\\= \lim_{x\to 0}
e^{\ctg^2 x \ln \cos x}= e^{\lim_{x\to 0}\ctg^2 x \ln \cos x}=\\= e^{\lim_{x\to
0}\cos^2 x \cdot \frac{\ln \cos x}{\sin^2 x}}= e^{\lim_{x\to 0}\cos^2 x \cdot
\lim_{x\to 0}\frac{\ln \cos x}{\sin^2 x}}=\\= e^{1\cdot \lim_{x\to 0}\frac{\ln
\cos x}{\sin^2
x}}= \left( e^\frac{0}{0}\right)={\smsize {\smsize\begin{pmatrix}\text{применяем}\\
\text{теорему \ref{Lopital_x->a_0/0}}\end{pmatrix}}}=\\= e^{\lim_{x\to
0}\frac{(\ln \cos x)'}{(\sin^2 x)'}}= e^{\lim_{x\to 0}\frac{\frac{-\sin x}{\cos
x}}{2\sin x \cos x}}=\\= e^{\lim_{x\to 0}\frac{-1}{2\cos^2 x}}=
e^{\frac{-1}{2}}=\frac{1}{\sqrt{e}}\end{multline*}\end{ex}

\begin{ex}[\bf неопределенность вида $\infty^0$] Вычислим предел
 \begin{multline*}\lim_{x\to +\infty} (\ln x)^{\frac{1}{x}}= \left(\infty^0 \right)=
\lim_{x\to +\infty} e^{\ln (\ln x)^{\frac{1}{x}}}=\\= \lim_{x\to +\infty}
e^{\frac{\ln (\ln x)}{x}}= e^{\lim_{x\to +\infty}\frac{\ln (\ln x)}{x}}= \left(
e^\frac{\infty}{\infty}\right)=\\={\smsize
{\smsize\begin{pmatrix}\text{применяем}\\
\text{теорему \ref{Lopital_x->infty_infty/infty}}\end{pmatrix}}}= e^{\lim_{x\to
+\infty}\frac{(\ln (\ln x))'}{(x)'}}=\\= e^{\lim_{x\to
+\infty}\frac{\frac{1}{x\ln x}}{1}}= e^{\lim_{x\to +\infty}\frac{1}{x\ln x}}=
e^0=1
\end{multline*}\end{ex}

\begin{ers}
Вычислите пределы:

1. $\lim\limits_{x\to +\infty} x^2 \cdot e^{-2x}$.

2. $\lim\limits_{x\to 1+0}\ln x \cdot \ln(x-1)$.

3. $\lim\limits_{x\to 0}\left(\frac{1}{x}-\frac{1}{e^x-1}\right)$.

4. $\lim\limits_{x\to +\infty}\left(\sqrt{x}-\ln x \right)$.

5. $\lim\limits_{x\to +0}\left(\ln \frac{1}{x}\right)^x$.

6. $\lim\limits_{x\to +\infty}\left(\ln 2x \right)^{\frac{1}{\ln x}}$.

7. $\lim\limits_{x\to 0}\left( 1+\sin^2 x \right)^{\frac{1}{\tg^2 x}}$.

8. $\lim\limits_{x\to 1-0} (1-x)^{\ln x}$.

9. $\lim\limits_{x\to +0}\left(\ln(x+e) \right)^{\frac{1}{x}}$.
\end{ers}

\end{multicols}\noindent\rule[10pt]{160mm}{0.1pt}

\subsection{Правило Лопиталя для построения графиков}

Правило Лопиталя используется везде, где нужно вычислить пределы, которые
элементарными способами найти трудно. Такое, в частности, случается при
построении графиков, и вот типичный пример.

\noindent\rule{160mm}{0.1pt}\begin{multicols}{2}

\begin{ex} Постройте график функции:
$$
f(x)=x^2\cdot \ln|x|
$$

1. Область определения: $(-\infty;0)\cup (0;+\infty)$. Точка $x=0$ является
особой. Функция является четной, поэтому нам достаточно построить график только
на интервале $x\in (0;+\infty)$.

2. Поскольку нас интересует только область $x\in (0;+\infty)$, в особой точке
$x=0$ мы вычисляем односторонний предел:
 \begin{multline*}\lim_{x\to +0} x^2\cdot \ln|x|= \lim_{x\to +0} x^2\cdot \ln x=
(0\cdot \infty)=\\= \lim_{x\to +0}\frac{\ln x}{x^{-2}}= \left(\frac{\infty}{\infty}\right)=\\={\smsize {\smsize\begin{pmatrix}\text{применяем}\\
\text{правило}\\ \text{Лопиталя}\end{pmatrix}}}= \lim_{x\to +0}\frac{x^{-1}}{-2
x^{-3}}= \lim_{x\to +0}\frac{x^2}{-2}=0
 \end{multline*}
Значит, вертикальных асимптот нет.

3. Наклонная асимптота нас тоже интересует только при $x\to +\infty$:
 \begin{multline*}
k_+=\lim_{x\to +\infty}\frac{f(x)}{x}=\lim_{x\to +\infty} x\cdot \ln|x|=\\=
\lim_{x\to +\infty} x\cdot \ln x=\infty
 \end{multline*}
Это означает, что и наклонных асимптот нет.

4. Производная:
 \begin{multline*}
f'(x)=(x^2\cdot \ln|x|)'=2x\cdot \ln|x|+x^2\cdot \frac{1}{x}=\\= x\cdot
(2\ln|x|+1)
 \end{multline*}
При $x>0$ получаем $f'(x)=0 \, \Leftrightarrow \, 2\ln x+1=0
 \, \Leftrightarrow \, \ln x=-\frac{1}{2}
 \, \Leftrightarrow \, x=e^{-\frac{1}{2}}$.
Интервалы знакопостоянства производной:

%\picture{0pt}{0pt}{p13.pcx}

\vglue90pt

5. Вторая производная:
 \begin{multline*}
f''(x)=\Big( x\cdot (2\ln|x|+1)\Big)'=(2\ln|x|+1)+x\cdot \frac{2}{x}=\\=
2\ln|x|+2=2(\ln|x|+1)
 \end{multline*}
При $x>0$ получаем $f''(x)=0 \, \Leftrightarrow \, \ln x+1=0
 \, \Leftrightarrow \, \ln x=-1
 \, \Leftrightarrow \, x=e^{-1}$.
Интервалы знакопостоянства второй производной:

%\picture{0pt}{0pt}{p14.pcx}

\vglue90pt

6. График функции для $x>0$:

%\picture{0pt}{-10pt}{p15.pcx}

\vglue110pt \noindent Отобразив эту картинку симметрично относительно оси
ординат, мы получим полный график нашей функции:

%\picture{0pt}{-10pt}{p16.pcx}

\vglue110pt
\end{ex}

\end{multicols}\noindent\rule[10pt]{160mm}{0.1pt}

\chapter{НЕОПРЕДЕЛЕННЫЙ ИНТЕГРАЛ}\label{CH-indef-integral}

Операция дифференцирования каждой функции $F$ ставит в соответствие ее
производную $F'$. В этой главе мы поговорим об обратной операции, которая по
производной $F'$ позволяет определить исходную функцию $F$.

\noindent\rule{160mm}{0.1pt}\begin{multicols}{2}

\section{Неопределенный интеграл}

При вычислении неопределенных интегралов от стандартных функций оказывается
удобным представлять эти интегралы как операции не над функциями, а над
однородными дифференциальными выражениями степени 1 (см. определение на
с.\pageref{odnor-diff-term}). Например, чтобы вычислить
интеграл\footnote{Интеграл $\int \ln x\d x$ будет подсчитан нами ниже в примере
\ref{ex-int-ln-x}.}
$$
\int \ln x\d x
$$
удобно считать, что это такая операция, которая однородному дифференциальному
выражению степени 1
$$
\ln x\d x
$$
(а не функции $f(x)=\ln x$) ставит в соответствие некое числовое выражение
$\int \ln x\d x$. Хотя разница кажется чисто формальной, тем не менее
удивительный факт состоит в том, что при таком подходе вычисления резко
упрощаются. В этом параграфе мы объясним правила нахождения интегралов от
дифференциальных термов (а с ними и от функций) и приведем примеры таких
вычислений.

\subsection{Однородные дифференциальные выражения степени 1}

\paragraph{Определение.}

\biter{

\item[$\bullet$] {\it Однородными дифференциальными выражениями степени 1}\label{diff-vyrazh-step-1} (в этой главе мы для краткости будем иногда называть их просто {\it дифференциальными выражениями}) называются конечные последовательности символов, составленные их числовых выражений (определенных нами на с.\pageref{opredelenie-chislovogo-terma}) и специального символа $\d$, называемого {\it формальным дифференциалом}, по следующим индуктивным правилам:
 \biter{
\item[1)] для любого одноместного числового выражения ${\mathcal P}$ от переменной $x$ запись
$$
\d ({\mathcal P})
$$
считается дифференциальным выражением степени 1; {\it областью допустимых значений} такого дифференциального выражения называется область допустимых значений его производной ${\mathcal P}$:
 \beq\label{D(dP)}
\D\Big(\d ({\mathcal P})\Big):=\D\l\frac{\d {\mathcal P}}{\d x}\r
 \eeq

\item[2)] если ${\mathcal P}$ -- числовое выражение, а ${\mathcal E}$ --- дифференциальное выражение степени 1, то записи
$$
({\mathcal P})\cdot({\mathcal E}),\quad ({\mathcal E})\cdot({\mathcal P}),\quad \frac{{\mathcal E}}{{\mathcal P}}
$$
также считаются дифференциальными выражениями степени 1; {\it область допустимых значений переменной} в первых двух из этих дифференциальных выражений определяется как пересечение областей допустимых значений для ${\mathcal P}$ и для ${\mathcal D}$:
\begin{align}
&\D\Big(({\mathcal P})\cdot({\mathcal E})\Big):=\D\Big({\mathcal P}\Big)\cap\D\Big({\mathcal E}\Big) \label{D(PE)}\\
&\D\Big(({\mathcal E})\cdot({\mathcal P})\Big):=\D\Big({\mathcal E}\Big)\cap\D\Big({\mathcal P}\Big) \label{D(EP)}
\end{align}
а в третьем -- формулой
\beq\label{D(E/P)}
 \D\l \frac{\mathcal E}{\mathcal P}\r=\D({\mathcal E})\cap\D({\mathcal P})\setminus\{x: {\mathcal P}=0\}
\eeq

\item[3)] если ${\mathcal D}$ и ${\mathcal E}$ --- два дифференциальных выражения степени 1, то записи
$$
({\mathcal D})+({\mathcal E}),\qquad ({\mathcal D})-({\mathcal E})
$$
также считаются дифференциальными выражениями степени 1;  {\it область допустимых значений} для таких дифференциальных выражений определяется как пересечение областей допустимых значений для ${\mathcal D}$ и для ${\mathcal E}$:
\begin{align}
&\D\Big(({\mathcal D})+({\mathcal E})\Big):=\D\Big({\mathcal D}\Big)\cap\D\Big({\mathcal E}\Big) \label{D(P+E)}\\
&\D\Big(({\mathcal D})-({\mathcal P})\Big):=\D\Big({\mathcal D}\Big)\cap\D\Big({\mathcal P}\Big) \label{D(P-E)}
\end{align}

}\eiter

\item[] Как и в случае с числовым выражением, каждое дифференциальное выражение
должно быть снабжено указанием, какие из входящих в него букв являются
переменными, а какие параметрами.

}\eiter

\paragraph{Равенство однородных дифференциальных выражений степени 1.}

\biter{

\item[$\bullet$] Два дифференциальных выражения ${\mathcal E}$ и ${\mathcal F}$ степени 1 считаются {\it
равными на множестве $x\in X$}, и изображается это формулой
$$
{\mathcal E}\underset{x\in X}{=}{\mathcal F},
$$
или, в случае, если $X=\D(\mathcal E)\cap\D(\mathcal F)$, то формулой
$$
{\mathcal E}={\mathcal F},
$$
и, наконец, в случае, если $X=\D(\mathcal E)=\D(\mathcal F)$, то формулой
$$
{\mathcal E}\equiv{\mathcal F},
$$
если их можно отождествить с помощью следующих индуктивных правил:

 \biter{
\item[(a)] умножение на числовое выражение подчинено правилам (здесь ${\mathcal P}$ и ${\mathcal Q}$ -- числовые выражения, а ${\mathcal E}$ -- дифференциальное степени 1):
\kern-20pt \begin{align}
& {\mathcal P}\cdot {\mathcal E}\equiv{\mathcal E}\cdot {\mathcal P} \label{kommut-diff-vyrazh} \\
& \frac{1}{\mathcal P}\cdot {\mathcal E}\equiv\frac{\mathcal E}{\mathcal P}\equiv{\mathcal E}\cdot \frac{1}{\mathcal P} \label{drob-diff-vyrazh} \\
& {\mathcal P}\cdot \Big({\mathcal Q}\cdot {\mathcal E}\Big)\equiv\Big({\mathcal P}\cdot {\mathcal Q}\Big)\cdot {\mathcal E}
\label{assoc-diff-vyrazh} \\
& \Big({\mathcal P}+{\mathcal Q}\Big)\cdot {\mathcal E}\equiv{\mathcal P}\cdot {\mathcal E}+ {\mathcal Q}\cdot {\mathcal E} \label{distribut-diff-vyrazh} \\
& 1\cdot {\mathcal E}\equiv{\mathcal E} \label{edin-v-diff-vyrazh}
 \end{align}

\item[(b)] если $x$ -- переменная, а ${\mathcal P}$ и ${\mathcal Q}$ -- числовые выражения, то
 \beq\label{ravenstvo-diff-vyrazh-step-1}
{\mathcal P}\d x\underset{x\in X}{=}{\mathcal Q}\d x \quad\Longleftrightarrow
\quad {\mathcal P}\underset{x\in X}{=}{\mathcal Q}
 \eeq

\item[(c)] для элементарных числовых выражений дифференциалы задаются следующей
таблицей (здесь всюду $x$ --- единственная переменная):

\bigskip
 \centerline{\bf Таблица дифференциалов:}\label{tablitsa-differentsialov}

\begin{align}
 & \d \Big( C \Big) \equiv 0\cdot\d x \label{dC}
 \\
 & \d \Big(x^\alpha \Big)\equiv \alpha \cdot x^{\alpha-1}\cdot \d x \label{d(x^a)}
 \\
 & \d \Big(a^{x} \Big)\equiv \ln a \cdot a^{x}\cdot \d x \label{d(a^x)}
 \\
 & \d\Big(e^{x} \Big)\equiv e^{x}\cdot \d x \label{d(e^x)}
 \\
 & \d\Big(\log_a x \Big)\equiv \frac{1}{x\cdot \ln a}\cdot \d x \label{d(log_a-x)}
 \\
 & \d\Big(\ln x \Big)\equiv \frac{1}{x}\cdot \d x \label{d(ln-x)}
 \\
 & \d\Big(\sin x \Big)\equiv  \cos x\cdot \d x \label{d(sin-x)}
 \\
 & \d\Big(\cos x \Big)\equiv  -\sin x\cdot \d x \label{d(cos-x)}
 \\
 & \d\Big(\tg x \Big)\equiv \frac{1}{\cos^2 x}\cdot \d x \label{d(tg-x)}
 \\
 & \d\Big(\ctg x \Big)\equiv -\frac{1}{\sin^2 x}\cdot \d x \label{d(ctg-x)}
 \\
 & \d\Big(\arcsin x \Big)\equiv \frac{1}{\sqrt{1-x^2}}\cdot \d x \label{d(arcsin-x)}
 \\
 & \d\Big(\arccos x \Big)\equiv -\frac{1}{\sqrt{1-x^2}}\cdot \d x \label{d(arccos-x)}
 \\
 & \d\Big(\arctg x \Big)\equiv \frac{1}{1+x^2}\cdot \d x \label{d(arctg-x)}
 \\
 & \d\Big(\arcctg x \Big)\equiv -\frac{1}{1+x^2}\cdot \d x \label{d(arcctg-x)}
 \end{align}\medskip

\item[(d)] дифференциал от суммы, разности, произведения и частного
дифференциальных выражений подчиняется правилам:
\begin{align}
& \kern-20pt \d \Big({\mathcal P}\pm {\mathcal Q}\Big)\equiv
\d {\mathcal P}\pm \d {\mathcal Q}  \label{differentsial-ot-summy} \\
& \kern-20pt \d \Big({\mathcal P}\cdot {\mathcal Q}\Big)\equiv  \Big(\d {\mathcal P}\Big)\cdot
{\mathcal Q}+ {\mathcal P}\cdot \Big(\d {\mathcal Q} \Big)
\label{differentsial-ot-proizvedeniya} \\
& \kern-20pt \d \l \frac{{\mathcal P}}{{\mathcal Q}}\r \equiv  \frac{\Big(\d {\mathcal
P}\Big)\cdot {\mathcal Q} -{\mathcal P}\cdot \Big(\d {\mathcal
Q}\Big)}{{\mathcal Q}^2} \label{differentsial-ot-drobi}
\end{align}

\item[(e)] если ${\mathcal P}$ --- одноместное числовое выражение от переменной
$x$, а ${\mathcal Q}$ --- одноместное числовое выражение от переменной $y$, то
для выражения, полученного подстановкой, выполняется равенство
 \begin{align}\label{differentsial-slozhnogo-terma}
& \kern-20pt \d \left({\mathcal Q}\Big|_{y={\mathcal P}}\right)\equiv  \left(\frac{\d}{\d
y}{\mathcal Q}\right)\Big|_{y={\mathcal P}}\cdot \d {\mathcal P}
 \end{align}

 }\eiter
 }\eiter

\btm\label{TH:dU=U'dx} Для всякого числового выражения ${\mathcal U}$ от переменной
$x$ справедливо равенство
 \beq\label{d(P)=d(P)/d(x)-cdot-d(x)}
\d{\mathcal U}\equiv \frac{\d {\mathcal U}}{\d x}\cdot \d x
 \eeq
\etm
\bpr
Это  доказывается индукцией. Прежде всего, из \eqref{dC}-\eqref{d(arcctg-x)} следует, что эта формула верна для элементарных выражений.
Предположим, что она верна для выражений до порядка $n$ включительно. Пусть ${\mathcal U}$ -- выражение порядка $n+1$.
Тогда по определению, $\mathcal U$ будет либо суммой, либо разностью, либо произведением, либо отношением, либо композицией выражений порядка не больше $n$. Рассмотрим каждый из этих случаев.

Если ${\mathcal U}\equiv {\mathcal P}+{\mathcal Q}$, где ${\mathcal P}$ и ${\mathcal Q}$ -- выражения порядка не больше $n$, то
\begin{multline*}
\d{\mathcal U}\equiv \d\Big({\mathcal P}+{\mathcal Q}\Big)\equiv \eqref{differentsial-ot-summy}\equiv \d{\mathcal P}+\d{\mathcal Q}\equiv \\ \equiv
{\scriptsize\begin{pmatrix}\text{предположение}\\ \text{индукции}\end{pmatrix}}\equiv \frac{\d {\mathcal P}}{\d x}\cdot \d x+\frac{\d {\mathcal Q}}{\d x}\cdot \d x\equiv \eqref{distribut-diff-vyrazh}\equiv \\ \equiv \l\frac{\d {\mathcal P}}{\d x}+\frac{\d {\mathcal Q}}{\d x}\r\cdot \d x\equiv \eqref{DEF:d(P+Q)/dx}\equiv \\ \equiv
\frac{\d \Big({\mathcal P}+{\mathcal Q}\Big)}{\d x}\cdot \d x\equiv \frac{\d {\mathcal U}}{\d x}\cdot \d x
\end{multline*}
И точно так же рассматривается случай ${\mathcal U}\equiv {\mathcal P}-{\mathcal Q}$.

Пусть далее
${\mathcal U}\equiv {\mathcal P}\cdot{\mathcal Q}$, где ${\mathcal P}$ и ${\mathcal Q}$ -- выражения порядка не больше $n$.
Тогда
\begin{multline*}
\d{\mathcal U}\equiv \d\Big({\mathcal P}\cdot{\mathcal Q}\Big) \equiv \eqref{differentsial-ot-proizvedeniya} \equiv \d{\mathcal P}\cdot{\mathcal Q}+{\mathcal P}\cdot\d{\mathcal Q} \equiv \\ \equiv
{\scriptsize\begin{pmatrix}\text{предположение}\\ \text{индукции}\end{pmatrix}} \equiv \\ \equiv
\l\frac{\d {\mathcal P}}{\d x}\cdot \d x\r\cdot{\mathcal Q}+{\mathcal P}\cdot\l\frac{\d {\mathcal Q}}{\d x}\cdot \d x\r \equiv
\eqref{kommut-diff-vyrazh} \equiv \\ \equiv
{\mathcal Q}\cdot\l\frac{\d {\mathcal P}}{\d x}\cdot \d x\r+{\mathcal P}\cdot\l\frac{\d {\mathcal Q}}{\d x}\cdot \d x\r \equiv
\eqref{assoc-diff-vyrazh} \equiv \\ \equiv
\l{\mathcal Q}\cdot\frac{\d {\mathcal P}}{\d x}\r\cdot \d x+\l{\mathcal P}\cdot\frac{\d {\mathcal Q}}{\d x}\r\cdot \d x \equiv
\eqref{distribut-diff-vyrazh} \equiv \\ \equiv
\l{\mathcal Q}\cdot\frac{\d {\mathcal P}}{\d x}+{\mathcal P}\cdot\frac{\d {\mathcal Q}}{\d x}\r\cdot \d x \equiv
\eqref{DEF:d(P-cdot-Q)/dx} \equiv \\ \equiv
\frac{\d \Big({\mathcal P}\cdot{\mathcal Q}\Big)}{\d x}\cdot \d x \equiv \frac{\d {\mathcal U}}{\d x}\cdot \d x
\end{multline*}

Теперь пусть ${\mathcal U} \equiv \frac{\mathcal P}{\mathcal Q}$, где опять ${\mathcal P}$ и ${\mathcal Q}$ -- выражения порядка не больше $n$. Тогда
\begin{multline*}
\d{\mathcal U} \equiv \d\l\frac{\mathcal P}{\mathcal Q}\r \equiv \eqref{differentsial-ot-drobi} \equiv \\ \equiv
\frac{\Big(\d {\mathcal
P}\Big)\cdot {\mathcal Q} -{\mathcal P}\cdot \Big(\d {\mathcal
Q}\Big)}{{\mathcal Q}^2} \equiv \\ \equiv
{\scriptsize\begin{pmatrix}\text{предположение}\\ \text{индукции}\end{pmatrix}} \equiv \\ \equiv
\frac{\l\frac{\d {\mathcal P}}{\d x}\cdot \d x\r\cdot{\mathcal Q}-{\mathcal P}\cdot\l\frac{\d {\mathcal Q}}{\d x}\cdot \d x\r}{{\mathcal Q}^2} \equiv
\eqref{kommut-diff-vyrazh} \equiv \\ \equiv
\frac{{\mathcal Q}\cdot\l\frac{\d {\mathcal P}}{\d x}\cdot \d x\r-{\mathcal P}\cdot\l\frac{\d {\mathcal Q}}{\d x}\cdot \d x\r}{{\mathcal Q}^2} \equiv \eqref{assoc-diff-vyrazh} \equiv \\ \equiv
\frac{\l{\mathcal Q}\cdot\frac{\d {\mathcal P}}{\d x}\r\cdot \d x-\l{\mathcal P}\cdot\frac{\d {\mathcal Q}}{\d x}\r\cdot \d x}{{\mathcal Q}^2} \equiv \eqref{distribut-diff-vyrazh} \equiv \\ \equiv
\frac{\l{\mathcal Q}\cdot\frac{\d {\mathcal P}}{\d x}-{\mathcal P}\cdot\frac{\d {\mathcal Q}}{\d x}\r\cdot \d x}{{\mathcal Q}^2} \equiv \eqref{drob-diff-vyrazh} \equiv \\ \equiv
\frac{{\mathcal Q}\cdot\frac{\d {\mathcal P}}{\d x}-{\mathcal P}\cdot\frac{\d {\mathcal Q}}{\d x}}{{\mathcal Q}^2}\cdot \d x \equiv  \eqref{DEF:d(P/Q)/dx} \equiv \\ \equiv
\frac{\d \frac{\mathcal P}{\mathcal Q}}{\d x}\cdot \d x \equiv \frac{\d {\mathcal U}}{\d x}\cdot \d x
\end{multline*}

И остается последний случай, когда
${\mathcal U} \equiv {\mathcal Q}\Big|_{y \equiv {\mathcal P}}$, и снова ${\mathcal P}$ и ${\mathcal Q}$ -- выражения порядка не больше $n$. Тогда
\begin{multline*}
\d{\mathcal U} \equiv \d\l{\mathcal Q}\Big|_{y={\mathcal P}}\r \equiv \eqref{differentsial-slozhnogo-terma} \equiv \frac{\d{\mathcal Q}}{\d
y}\Big|_{y={\mathcal P}}\cdot \d {\mathcal P} \equiv \\ \equiv
{\scriptsize\begin{pmatrix}\text{предположение}\\ \text{индукции}\end{pmatrix}} \equiv \frac{\d{\mathcal Q}}{\d
y}\Big|_{y={\mathcal P}}\cdot \l\frac{\d{\mathcal P}}{\d x}\d x\r \equiv \\ \equiv \eqref{assoc-diff-vyrazh} \equiv \l\frac{\d{\mathcal Q}}{\d
y}\Big|_{y={\mathcal P}}\cdot \frac{\d{\mathcal P}}{\d x}\r\cdot\d x \equiv \eqref{7.4.4} \equiv \\ \equiv \frac{\d}{\d x}\l{\mathcal Q}\Big|_{y={\mathcal P}}\r\cdot\d x \equiv \frac{\d {\mathcal U}}{\d x}\cdot \d x.
\end{multline*}
\epr

\bcor\label{COR:E=Sdx}
Любое однородное дифференциальное выражение ${\mathcal E}$ степени 1 от переменной $x$ имеет вид
 \beq\label{odnor-diff-term-1}
 {\mathcal E} \equiv {\mathcal S}\d x
 \eeq
где ${\mathcal S}$ --- некоторое числовое выражение (с той же областью допустимых значений переменной).
\ecor
\bpr
Из теоремы \ref{TH:dU=U'dx} и формулы \eqref{D(dP)} следует, что это верно для дифференциальных выражений степени 1 вида ${\mathcal E}=\d {\mathcal U}$. Все остальные дифференциальные выражения степени 1 получаются из таких либо домножением (слева или справа), либо делением на числовое выражение, либо прибавлением каких-то других дифференциальных выражений степени 1, полученных таким же способом. Но, во-первых, из формул \eqref{kommut-diff-vyrazh}-\eqref{distribut-diff-vyrazh} следует, что при всех этих операциях множитель $\d x$ можно вынести за скобки, а внутри скобок останется какое-то числовое выражение. Это и будет ${\mathcal S}$.

А, во-вторых, из формул \eqref{D(dP)}-\eqref{D(P-E)} следует, что при каждой из этих операций (умножение или деление на числовое выражение, или прибавление и вычитание дифференциального выражения) получающиееся в скобках при вынесении $\d x$ числовое выражение ${\mathcal S}$ имеет ту же область допустимых значений переменной, что дифференциальное выражение ${\mathcal E}$ (в котором $\d x$ еще не вынесено за скобки).
\epr

\paragraph{Нулевое однородное дифференциальное выражение степени 1.}

\biter{

\item[$\bullet$] {\it Нулевым дифференциальным выражением степени 1} от переменной $x$ называется выражение $0\cdot \d x$. Оно для удобства обозначается символом 0:
    $$
    0:=0\cdot \d x.
    $$
    (это означает, что если где-то встречается запись ${\mathcal E}=0$ или ${\mathcal E}+0$ или ${\mathcal E}-0$, где ${\mathcal E}$ -- какое-то дифференциальное выражение степени 1, то под символом $0$ здесь понимается $0\cdot \d x$).
}\eiter

\bigskip

\centerline{\bf Свойства нулевого}
\centerline{\bf дифференциального  выражения:}

 \biter{\it

\item[$1^\circ.$] Умножение на нулевое дифференциальное выражение всегда дает нулевое дифференциальное выражение:
\beq\label{umnozh-na-nul-diff-vyrazh}
{\mathcal P}\cdot 0=0=0\cdot{\mathcal P}
\eeq

\item[$2^\circ.$] Умножение на нулевое числовое выражение всегда дает нулевое дифференциальное выражение:
\beq\label{umnozh-na-nul-chisl-vyrazh}
0\cdot {\mathcal E}=0={\mathcal E}\cdot 0
\eeq

\item[$3^\circ.$] Прибавление нулевого дифференциального выражения не меняет дифференциальное выражение:
\beq\label{slozh-s-nul-diff-vyrazh}
{\mathcal E}+ 0\equiv{\mathcal E}\equiv 0+{\mathcal E}
\eeq

\item[$4^\circ.$] Дифференциал от произвольного параметра равен нулю:
\beq\label{differentsial-ot-parametra}
\d C\equiv 0
\eeq

\item[$5^\circ.$] Дифференциал от переменной наоборот, не может быть равен нулю:
\beq\label{differentsial-ot-peremennoi}
\d x\ne 0
\eeq

}\eiter
\bpr 1. Умножение на нулевое дифференциальное выражение:
 $$
{\mathcal P}\cdot 0={\mathcal P}\cdot \Big(0\cdot \d x\Big)=\eqref{assoc-diff-vyrazh}=\Big({\mathcal P}\cdot 0\Big)\cdot \d x=0\cdot x=0
 $$

2. Умножение на нулевое числовое выражение:
 \begin{multline*}
 0\cdot {\mathcal E}=\eqref{odnor-diff-term-1}= 0\cdot \Big({\mathcal S}\cdot \d x\Big)=\\=\eqref{assoc-diff-vyrazh}=\Big(0\cdot {\mathcal S}\Big)\cdot \d x=0\cdot x=0
 \end{multline*}

3. Прибавление нулевого дифференциального выражения:
 \begin{multline*}
{\mathcal E}+ 0=\eqref{odnor-diff-term-1}={\mathcal S}\cdot\d x+0\cdot \d x=\eqref{distribut-diff-vyrazh}=\\=\Big({\mathcal S}+0\Big)\cdot \d x={\mathcal S}\cdot\d x={\mathcal E}
 \end{multline*}

4. Дифференциал от параметра:
$$
\d C=\eqref{dC}=0\cdot\d x=0
$$

5. Если бы $\d x=0=0\cdot\d x$, то мы получили бы
$$
1\cdot\d x=\d x=0=0\cdot\d x
$$
и в силу \eqref{ravenstvo-diff-vyrazh-step-1} это означало бы, что $\d x$ можно отбрость:
$$
1=0
$$
Поскольку это неверно, наше исходное предположение, что $\d x=0$, тоже неверно.
\epr

\bcor
Всякий параметр можно выносить за знак
дифференциала:
 \beq\label{vynesenie-konstaty-za-differentsial}
\d\big( C\cdot {\mathcal P}\big)=C\cdot \d {\mathcal P}
 \eeq
\ecor
\bpr
 \begin{multline*}
\d\big( C\cdot {\mathcal P}\big)=\eqref{differentsial-ot-proizvedeniya}= \kern-18pt\overbrace{\d C}^{\scriptsize\begin{matrix}0\\
\phantom{\tiny\eqref{differentsial-ot-parametra}} \ \text{\rotatebox{90}{$=$}}
\ {\tiny\eqref{differentsial-ot-parametra}}\end{matrix}}\kern-18pt\cdot {\mathcal P}+ C\cdot\d {\mathcal P}=\\=
\kern-18pt\underbrace{0\cdot {\mathcal P}}_{\scriptsize\begin{matrix}
\phantom{\tiny\eqref{umnozh-na-nul-diff-vyrazh}} \ \text{\rotatebox{90}{$=$}}
\ {\tiny\eqref{umnozh-na-nul-diff-vyrazh}}\\ 0\end{matrix}}\kern-18pt + C\cdot\d {\mathcal P}=0+ C\cdot\d {\mathcal P}=\eqref{slozh-s-nul-diff-vyrazh}=C\cdot\d {\mathcal P}
 \end{multline*}

\epr

\paragraph{Константы.}

 \biter{

\item[$\bullet$] Числовое выражение ${\mathcal P}$ от переменной $x$ называется {\it константой} на интервале $I$, если определяемая им функция
$$
f(x)={\mathcal P}
$$
постоянна на этом интервале:
$$
\forall x,y\in I\quad f(x)=f(y).
$$
 }\eiter

\bex Всякие нульместное (то есть не содержащее переменной) числовое выражение будет, очевидно, константой, например:
$$
\frac{3}{2},\quad \sqrt{2},\quad e^2
$$
\eex

\bex Нетрудно придумать константы, содержащие переменную:
$$
(x+1)^2-x^2-2x,\quad \frac{e^{x+2}}{e^x},\quad \sin^2x+\cos^2x.
$$
\eex

\btm\label{TH:osn-sv-konstant} Для всякой константы ${\mathcal P}$ на произвольном интервале $I$ найдется нульместное числовое выражение ${\mathcal C}$, совпадающее с ${\mathcal P}$ на $I$:
$$
{\mathcal P}\underset{x\in I}{=}{\mathcal C}.
$$
\etm
\bpr
Обозначим
$$
f(x)={\mathcal P},\qquad x\in I.
$$
Выберем в интервале $I$ какую-нибудь точку $a$, определяемую каким-нибудь нульместным числовым выражением ${\mathcal A}$ (такая точка найдется, например, можно взять какую-нибудь рациональную точку в $I$ -- именно для этого требуется, чтобы $I$ был интервалом):
$$
{\mathcal A}=a\in I
$$
После этого полагаем
$$
{\mathcal C}={\mathcal P}\Big|_{x={\mathcal A}}
$$
\epr

\btm\label{TH:konst<=>d/dx=0} Для числового выражения ${\mathcal R}$ от переменной $x$ следующие условия эквивалентны:
 \biter{

 \item[(i)] ${\mathcal R}$ -- константа на интервале $I$,

 \item[(ii)] производная ${\mathcal R}$ по переменной $x$ равна нулю на интервале $I$:
 \beq
 \frac{\d {\mathcal R}}{\d x}\underset{x\in I}{=}0,
 \eeq

 \item[(iii)] дифференциал ${\mathcal R}$ равен нулю на интервале $I$:
 \beq
 \d {\mathcal R}\underset{x\in I}{=}0.
 \eeq

 }\eiter

\etm
\begin{proof} Эквивалентность условий (ii) и (iii) очевидна, импликация $(i)\Rightarrow (ii)$ тоже, поэтому нужно только доказать импликацию $(ii)\Rightarrow (i)$. Пусть производная ${\mathcal R}$ равна нулю на интервале $I$. Обозначим через $f$ функцию, определяемую выражением ${\mathcal R}$:
$$
f(x)={\mathcal R}
$$
И, как в теореме \ref{TH:osn-sv-konstant}, выберем в интервале $I$ какую-нибудь точку $a$, определяемую каким-нибудь нульместным числовым выражением ${\mathcal A}$ (например, можно взять какую-нибудь рациональную точку в $I$):
$$
{\mathcal A}=a\in I
$$
Положим
$$
{\mathcal C}={\mathcal R}\Big|_{x={\mathcal A}}
$$
По лемме \ref{LM:o-konstantah}, функция $f$ всюду на интервале $I$ будет равна своему значению в точке $a$,
$$
f(x)=f(a),\quad x\in I
$$
Применительно к выражениям $\mathcal R$ и $\mathcal C$ это означает,  что они совпадают на $I$:
$$
{\mathcal R}\underset{x\in I}{=}{\mathcal C}.
$$
 \end{proof}

\paragraph{Эквивалентность числовых выражений.}

\biter{

\item[$\bullet$] Два числовых выражения $\mathcal P$ и $\mathcal Q$ мы будем называть {\it эквивалентными}, и обозначать это записью
 $$
 {\mathcal P}\approx {\mathcal Q}
 $$
если их дифференциалы имеют одинаковую область допустимых значений переменной
$$
\D(\d\mathcal P)=\D(\d\mathcal Q)
$$
и на каждом интервале $I\subseteq \D(\d\mathcal P)=\D(\d\mathcal Q)$ выражения $\mathcal P$ и $\mathcal Q$ отличаются на константу:
$$
{\mathcal P}\underset{x\in I}{=} {\mathcal Q}+{\mathcal C},\quad \d {\mathcal C}\underset{x\in I}{=} 0.
$$
}\eiter

\btm\label{COR:dQ=dP=>Q=P+C} Числовые выражения ${\mathcal P}$ и ${\mathcal Q}$ имеют один и тот же дифференциал тогда и только тогда, когда они эквивалентны:
$$
\d {\mathcal Q}\equiv \d {\mathcal P}\quad\Longleftrightarrow\quad {\mathcal Q}\approx{\mathcal P}.
$$
\etm
 \bpr Обозначив
 $$
 {\mathcal C}\equiv {\mathcal Q}-{\mathcal P},
 $$
 мы получим:
 $$
\d {\mathcal C}\equiv \d{\mathcal Q}-\d{\mathcal P}\underset{x\in I}{=}0,
 $$
поэтому по теореме \ref{TH:konst<=>d/dx=0}, ${\mathcal C}$ -- константа на $I$.
 \epr

\subsection{Первообразная и неопределен\-ный интеграл}

 \biter{

\item[$\bullet$] Числовое выражение ${\mathcal P}$ называется {\it первообразной} для дифференциального выражения ${\mathcal E}$ степени 1, если дифференциал от ${\mathcal P}$ тождественно равен ${\mathcal E}$:
$$
\d {\mathcal P}\equiv {\mathcal E}
$$
(напомним, что это предполагает равенство $\D(\d {\mathcal P})=\D({\mathcal E})$).

\item[$\bullet$] Числовое выражение ${\mathcal P}$ называется {\it неопределенным интегралом}, или просто {\it интегралом}, от дифференциального выражения ${\mathcal E}$ степени 1 с {\it параметром интегрирования} $C$, и обозначается это записью
    \beq\label{DEF:neopr-int}
    {\mathcal P}=\int {\mathcal E},
    \eeq
    если

    \biter{
\item[(a)] ${\mathcal P}$ является первообразной для ${\mathcal E}$
 \beq\label{DEF:neopr-int-1}
\d {\mathcal P}\equiv {\mathcal E}
 \eeq

\item[(b)] $C$ является параметром в ${\mathcal P}$, причем любая другая первообразная ${\mathcal Q}$  для ${\mathcal E}$
$$
\d {\mathcal Q}\equiv {\mathcal E}
$$
на любом интервале $I\subseteq \D({\mathcal E})$ получается из ${\mathcal P}$ заменой параметра $C$ на некоторое нульместное числовое выражение ${\mathcal C}$:
 \beq\label{DEF:neopr-int-2}
{\mathcal Q}\underset{x\in I}{=}{\mathcal P}\Big|_{C={\mathcal C}}.
 \eeq
    }\eiter
 }\eiter

\btm\label{DF:int-terma} Пусть числовое выражение $\mathcal P$ -- какая-нибудь первообразная для $\mathcal E$,
$$
{\mathcal E}\equiv \d {\mathcal P}
$$
причем множество параметров у этих выражений одно и то же, и буква $C$ в него не входит. Тогда выражение ${\mathcal P}+C$ является неопределенным интегралом для $\mathcal E$ с константой интегрирования $C$:
$$
\int {\mathcal E}={\mathcal P}+C
$$
\etm
\bpr
Пусть $\mathcal Q$ -- какая-то другая первообразная для $\mathcal E$:
$$
\d {\mathcal Q}\equiv {\mathcal E}\equiv \d {\mathcal P}
$$
Тогда по теореме \ref{COR:dQ=dP=>Q=P+C}, на всяком интервале $I\subseteq \D\big(\d{\mathcal Q}\big)=\D\big(\d{\mathcal P}\big)$
выражения ${\mathcal Q}$ и ${\mathcal P}$ отличаются друг от друга на некоторую константу ${\mathcal C}$:
$$
{\mathcal Q}\underset{x\in I}{=}{\mathcal P}+{\mathcal C}
$$
По теореме \ref{TH:osn-sv-konstant}, константу ${\mathcal C}$ можно на интервале $I$ считать нульместной. Поэтому
$$
{\mathcal P}+C\Big|_{C={\mathcal C}}={\mathcal P}+{\mathcal C}\underset{x\in I}{=}{\mathcal Q}.
$$
\epr

\bigskip
 \centerline{\bf Таблица интегралов:}\label{TAB:int-dlya-termov}
 \label{tablitsa-integralov}

\begin{align}
&\int x^\alpha \, \d x=\frac{x^{\alpha+1}}{\alpha+1}+C \qquad (\alpha\ne -1) \label{int-x^a}
\\
&\int \frac{d x}{x} =\ln|x|+C \label{int-ln-x} \\
&\int a^x \, \d x=\frac{a^x}{\ln a}+C, \qquad (0<a\ne 1) \label{int-a^x}
\\
&\int e^x \, \d x= e^x+C \label{int-e^x} \\
&\int \sin x \, \d x=-\cos x+C \label{int-sin-x}
\\
&\int \cos x \, \d x=\sin x+C  \label{int-cos-x} \\
&\int \frac{d x}{\cos^2 x} =\tg x+C \label{int-1/cos^2-x}
\\
&\int \frac{d x}{\sin^2 x} =-\ctg x+C \label{int-1/sin^2-x} \\
&\int \frac{d x}{a^2+x^2} =\frac{1}{a}\arctg\frac{x}{a}+C \label{int-1/a^2+x^2}
\\
&\int \frac{d x}{\sqrt{a^2-x^2}}=\arcsin \frac{x}{a}+C  \label{int-1/sqrt(a^2-x^2)} \\
&\int \frac{d x}{a^2-x^2} =\frac{1}{2a}\ln \left| \frac{a+x}{a-x}\right|+C \label{int-1/a^2-x^2}
\\
& \int \frac{d x}{\sqrt{x^2+A}}=\ln \left| x+\sqrt{x^2+A}\right|+C \label{int-1/x^2+A}
 \end{align}

\bpr Все формулы здесь выводятся с помощью теоремы \ref{DF:int-terma} взятием
дифференциала. Например, из
$$
\d\left(\frac{x^{\alpha+1}}{\alpha+1}\right)= \frac{(\alpha+1)\cdot
x^{\alpha}}{\alpha+1}\,\d x=x^{\alpha}\,\d x
$$
следует, формула \eqref{int-x^a}:
$$
\int x^\alpha \, \d x= \frac{x^{\alpha+1}}{\alpha+1}+C
$$

Формула \eqref{int-ln-x} содержит модуль, поэтому для ее доказательства нужно рассмотреть два случая:
при $x>0$ получаем:
$$
\frac{\d}{\d x}(\ln |x| )=\frac{\d}{\d x}(\ln x )=\frac{1}{x},
$$
И при $x<0$ получается то же самое:
 \begin{multline*}
\frac{\d}{\d x}(\ln |x| )=\frac{\d}{\d x}\Big(\ln (-x) \Big)=
\frac{\d}{\d x}\Big(\ln y\Big|_{y=-x} \Big)=\\=\frac{\d}{\d y}\Big(\ln y\Big)\cdot \frac{\d (-x)}{\d x} =\frac{1}{-x}\cdot (-1)=\frac{1}{x}.
 \end{multline*}
Вместе это дает \eqref{int-ln-x}.
\epr

\brem Важно понимать, что во всех формулах правую часть можно записывать по-разному (потому что необязательно в записи следовать предписаниям теоремы \ref{DF:int-terma}). Например, заменив параметр $C$ на любой другой новый параметр, мы получим равносильные формулы, скажем, такую:
$$
\int x^\alpha \, \d x=\frac{x^{\alpha+1}}{\alpha+1}+D
$$
Точно так же, можно прибавлять к правой части разные нульместные константы, и опять будут получаться верные формулы, например:
$$
\int x^\alpha \, \d x=\frac{x^{\alpha+1}}{\alpha+1}+C+10
$$
Более того, в случае, если область допустимых значений переменной не совпадает со всей прямой, правая часть допускает совсем экзотические варианты записи. Например, формулу \eqref{int-ln-x} вполне можно заменить на формулу
$$
\int \frac{d x}{x} =\begin{cases}\ln x+C, & x>0 \\ \ln(-x)+D, & x<0
\end{cases},
$$
или на формулу
$$
\int \frac{d x}{x} =\ln|x|+C\cdot \sgn x.
$$
Мы советуем читателю самостоятельно проверить справедливость этих формул.
\erem

\bex Рассмотрим дифференциальное выражение 1 степени
$$
{\mathcal U}=x^2\d (x^3)
$$
Если вынести $x^3$ из-под знака дифференциала
$$
\d (x^3)=3x^2\d x,
$$
то мы получим каноническое представление
$$
{\mathcal U}=x^2\d (x^3)=x^2\cdot 3x^2\cdot\d x=3x^4\d x
$$
Отсюда следует, что интеграл от ${\mathcal U}$ можно вычислить с помощью
таблицы на с. \pageref{TAB:int-dlya-termov}:
 \begin{multline*}
\int{\mathcal U}=\int x^2\d (x^3)=\int 3x^4\d x=\\=3\int x^4\d
x=\frac{3}{5}\cdot x^5+C
 \end{multline*}
 \eex

\bex Можно рассмотреть похожее выражение, в котором выражения перед и после
дифференциала переставлены:
$$
{\mathcal U}=x^3\d (x^2)
$$
Для него каноническое представление будет выглядеть так:
$$
{\mathcal U}=x^3\d (x^2)=x^3\cdot 2x\cdot\d x=2x^4\d x
$$
И интеграл получится другим:
 \begin{multline*}
\int{\mathcal U}=\int x^3\d (x^2)=\int 2x^4\d x=\\=2\int x^4\d
x=\frac{2}{5}\cdot x^5+C
 \end{multline*}
 \eex

\bex Следующий пример показывает, что интеграл можно вычислять от выражения с
параметром (в данном случае $x$ будет переменной, а $y$ --- параметром):
$$
\int x^2\cdot y \, \d x= y\cdot \int x^2\, \d x= \frac{1}{3}x^3\cdot y+C,
$$
Если же, наоборот, $x$ --- параметр, а $y$ --- переменная, то ответ меняется:
$$
\int x^2\cdot y \, \d y=  x^2\cdot \int y \, \d y= \frac{1}{2}x^2\cdot y^2+C
$$
\eex

\begin{samepage}

\bigskip

\centerline{\bf Свойства неопределенного интеграла}

 \biter{\it
\item[$1^\circ.$] Если дифференциальное выражение $\mathcal E$ обладает интегралом, то дифференциал от этого интеграла тождественно равен $\mathcal E$:
 \beq\label{d(int)}
\d \int {\mathcal E}\equiv {\mathcal E}
 \eeq

\item[$2^\circ.$] Дифференциал $\d{\mathcal R}$ произвольного
числового выражения ${\mathcal R}$ всегда обладает интегралом, который эквивалентен ${\mathcal R}$:
 \beq\label{int(d)}
\int\d{\mathcal R}\approx {\mathcal R}
 \eeq

\item[$3^\circ$] {\bf Аддитивность:} интеграл от суммы эквивалентен сумме
интегралов
 \beq\label{additivnost-int-ot-terma}
\int \big( {\mathcal P}+{\mathcal Q}\big) \d x\approx \int {\mathcal P}\d x
+\int {\mathcal Q}\d x
 \eeq

\item[$4^\circ$] {\bf Однородность:} константу можно выносить за знак интеграла
 \beq\label{vynesenie-konst-v-int-ot-terma}
\d {\mathcal C}\equiv 0\quad\Longrightarrow\quad
 \int \big( {\mathcal C}\cdot {\mathcal P}\big) \d x\approx
 {\mathcal C}\cdot \int {\mathcal P}\d x
 \eeq

\item[$5^\circ$] {\bf Замена переменной:} если ${\mathcal P}$ --- числовое
выражение от переменной $x$, а ${\mathcal Q}$ --- числовое выражение от
переменной $y$, то для выражения, полученного подстановкой, выполняется
эквивалентность
 \beq\label{zamena-perem-v-int-ot-terma}
 \int \left( {\mathcal Q}\Big|_{y={\mathcal P}}\right) \d {\mathcal P}
 \approx
 \left( \int {\mathcal Q}\d y\right)\Big|_{y={\mathcal P}}
 \eeq

\item[$6^\circ$] {\bf Интегрирование по частям:}
 \beq\label{int-po-chastyam-int-ot-terma}
 \int {\mathcal P}\d {\mathcal Q}
 \approx
 {\mathcal P}\cdot {\mathcal Q}-\int {\mathcal Q}\d {\mathcal P}
 \eeq

}\eiter

\end{samepage}

\bpr

1. Формула \eqref{d(int)} -- просто следствие формулы \eqref{DEF:neopr-int}: если $\mathcal P=\int {\mathcal E}$, то
$\d{\mathcal P}\equiv {\mathcal E}$.

2. Если $\mathcal P=\int \d {\mathcal R}$, то по формуле \eqref{DEF:neopr-int},
$$
\d{\mathcal P}\equiv \d {\mathcal R}
$$
$$
\phantom{\scriptsize \text{теорема \ref{COR:dQ=dP=>Q=P+C}}}\ \Downarrow\ {\scriptsize \text{теорема \ref{COR:dQ=dP=>Q=P+C}}}
$$
$$
{\mathcal P}\approx {\mathcal R}
$$

3. Для доказательства \eqref{additivnost-int-ot-terma} обозначим
$$
{\mathcal U}=\int {\mathcal P}\d x,
$$
$$
{\mathcal V}= \int
{\mathcal Q}\d x,
$$
$$
{\mathcal W}=\int \Big({\mathcal P}+{\mathcal
Q}\Big) \d x.
$$
Тогда
$$
\d {\mathcal W}\equiv{\mathcal P}\d x+{\mathcal Q}\d x\equiv\d{\mathcal U}+\d {\mathcal
V}\equiv\d\Big({\mathcal U}+{\mathcal V}\Big)
$$
$$
\phantom{\scriptsize \text{теорема \ref{COR:dQ=dP=>Q=P+C}}}\ \Downarrow\ {\scriptsize \text{теорема \ref{COR:dQ=dP=>Q=P+C}}}
$$
$$
{\mathcal W}\approx {\mathcal U}+{\mathcal V}.
$$

4. Для доказательства \eqref{vynesenie-konst-v-int-ot-terma} положим
$$
\d{\mathcal C}\equiv 0,\qquad {\mathcal U}=\int {\mathcal P}\d x,\qquad
{\mathcal V}=\int \Big({\mathcal C}\cdot{\mathcal P}\Big)\d x.
$$
Тогда
\begin{multline*}
\d {\mathcal V}\equiv\Big({\mathcal C}\cdot{\mathcal P}\Big)\d x\equiv {\mathcal
C}\cdot\Big({\mathcal P}\d x\Big)\equiv\\ \equiv{\mathcal C}\cdot \Big(\d {\mathcal
U}\Big)\equiv\eqref{vynesenie-konstaty-za-differentsial}\equiv\d \Big({\mathcal
C}\cdot{\mathcal U}\Big)
\end{multline*}
$$
\phantom{\scriptsize \text{теорема \ref{COR:dQ=dP=>Q=P+C}}}\ \Downarrow\ {\scriptsize \text{теорема \ref{COR:dQ=dP=>Q=P+C}}}
$$
$$
{\mathcal V}\approx {\mathcal C}\cdot  {\mathcal U}.
$$

5. Для \eqref{zamena-perem-v-int-ot-terma} положим
$$
{\mathcal U}=\int \left( {\mathcal Q}\Big|_{y={\mathcal P}}\right) \d
{\mathcal P},\qquad {\mathcal V}=\int {\mathcal Q}\d y
$$
Тогда
$$
\frac{\d {\mathcal V}}{\d y}\equiv{\mathcal Q}
$$
и поэтому
\begin{multline*}
\d {\mathcal U}\equiv \left( {\mathcal Q}\Big|_{y={\mathcal P}}\right) \d {\mathcal
P}\equiv\left( \frac{\d {\mathcal V}}{\d y}\Bigg|_{y={\mathcal P}}\right) \d
{\mathcal P}\equiv \\ \equiv\eqref{differentsial-slozhnogo-terma}\equiv \d \left({\mathcal
V}\Big|_{y={\mathcal P}}\right)
\end{multline*}
$$
\phantom{\scriptsize \text{теорема \ref{COR:dQ=dP=>Q=P+C}}}\ \Downarrow\ {\scriptsize \text{теорема \ref{COR:dQ=dP=>Q=P+C}}}
$$
$$
{\mathcal U}\approx
{\mathcal V}\Big|_{y={\mathcal P}}.
$$

6. Остается \eqref{int-po-chastyam-int-ot-terma}. Обозначим
$$
{\mathcal U}= \int {\mathcal P}\d {\mathcal Q},\qquad {\mathcal
V}= \int {\mathcal Q}\d {\mathcal P}
$$
Тогда
 \begin{multline*}
\d\Big({\mathcal U}+{\mathcal V}\Big)\equiv\d {\mathcal U}+\d {\mathcal V}\equiv{\mathcal
P}\d {\mathcal Q}+{\mathcal Q}\d {\mathcal P}\equiv\\ \equiv\eqref{differentsial-ot-summy}\equiv\d
\Big({\mathcal P}\cdot{\mathcal Q}\Big)
 \end{multline*}
$$
\phantom{\scriptsize \text{теорема \ref{COR:dQ=dP=>Q=P+C}}}\ \Downarrow\ {\scriptsize \text{теорема \ref{COR:dQ=dP=>Q=P+C}}}
$$
$$
{\mathcal U}+{\mathcal V} \approx {\mathcal
P}\cdot{\mathcal Q}
$$
$$
\Downarrow
$$
$$
{\mathcal U}\approx {\mathcal
P}\cdot{\mathcal Q}-{\mathcal V}.
$$

 \epr

\subsection{Вычисление интеграла}

Покажем теперь на примерах, как используются формулы
\eqref{additivnost-int-ot-terma}-\eqref{int-po-chastyam-int-ot-terma}.

\paragraph{Применение формул аддитивности и однородности.}
Начнем с \eqref{additivnost-int-ot-terma} и
\eqref{vynesenie-konst-v-int-ot-terma}.

\begin{ex}\label{ex-12.3.2}
 \begin{multline*}
\int \left[ 2+x^3-3\sin x-\frac{5}{1+x^2}\right] \, \d x \approx  \\
 \approx \int 2 \, \d x+\int x^3 \, \d x-3\int \sin x \, \d x-\\
-5\int \frac{1}{1+x^2}\, \d x = 2 x+ \frac{x^4}{4} +3\cos x -5\arctg x+C
 \end{multline*}
\end{ex}

\begin{ex}\label{ex-12.3.3}
 \begin{multline*}
\int \left(\sqrt{x}+\frac{2}{\cos^2 x}\right) \, \d x \approx  \int
x^{\frac{1}{2}}\, \d x+\\+2\int \frac{1}{\cos^2 x}\, \d x =
\frac{2}{3}x^{\frac{3}{2}} +2\tg x+C
 \end{multline*}
\end{ex}

\paragraph{Применение формулы замены переменной.}
Иллюстрации на применение \eqref{zamena-perem-v-int-ot-terma}:

\begin{ex}\label{ex-12.4.1}
Как мы уже видели,
 \begin{multline*}
\int x^4 \, \d x^2 = \int t^2 \, \d t \, \Big|_{t=x^2} = \\ =
\frac{1}{3}t^3+C \, \Big|_{t=x^2} =  \frac{1}{3}x^6+C
 \end{multline*}
Этот интеграл можно вычислить, вынеся $x^2$ из-под знака дифференциала, и от
этого результат не изменится:
 \begin{multline*}
\int x^4 \, \d x^2 = \int x^4 2x \, \d x =  2\int x^5 \, \d x =
\\ =  \frac{2}{6}x^6+C =  \frac{1}{3}x^6+C.
 \end{multline*}
\end{ex}

\begin{ex}\label{ex-12.4.2}
Точно также можно двумя способами сосчитать интеграл
$$
\int x^2 \, \d \sqrt{x}
$$
и результат будет одним и тем же: с одной стороны,
 \begin{multline*}
\int x^2 \, \d \sqrt{x} = \int t^4 \, \d t \, \Big|_{t=\sqrt{x}} =
\\ =  \frac{1}{5}t^5+C \, \Big|_{t=\sqrt{x}} =  \frac{1}{5}x^\frac{5}{2}+C
 \end{multline*}
а с другой,
 \begin{multline*}
\int x^2 \, \d \sqrt{x} = \int x^2 \frac{1}{2\sqrt{x}}\, \d x = \\
= \frac{1}{2}\int x^\frac{3}{2}\, \d x =  \frac{1}{2}\cdot
\frac{2}{5}\cdot x^\frac{5}{2} +C =  \frac{1}{5}\cdot x^\frac{5}{2} +C
 \end{multline*}
\end{ex}

Вынесение из-под знака дифференциала оказывается полезным средством, если при
этом выражения упрощаются:

\begin{exs} Иногда вынесение
из-под знака дифференциала упрощает интеграл и позволяет сразу его вычислить:
$$
\int x \, \d \ln x = \int x\cdot \frac{1}{x} \, \d x = \int 1 \,\d x
= x+C,
$$
Но так случается не всегда. Например, в следующем интеграле вынесение из-под
знака дифференциала ничего не дает:
$$
\int e^x \, \d\arctg x  = \int \frac{e^x}{1+x^2} \, \d x
$$
\end{exs}

\begin{ex}\label{ex-12.5.5}
 \begin{multline*}
\int (2x-5)^{10}\, \d x   =  \left| \begin{array}{c} y=2x-5 \\
x=\frac{y+5}{2}\end{array}\right| = \\ =  \int y^{10}\,  d
\frac{y+5}{2}\, \Big|_{y=2x-5} =
{\smsize\begin{pmatrix}\text{выносим}\, \, \frac{y+5}{2}\\
\text{из-под знака}\\ \text{дифференциала}\end{pmatrix}} = \\ =
\frac{1}{2}\int y^{10}\,  d y \, \Big|_{y=2x-5} =  \frac{1}{2\cdot 11}
y^{11}+C \, \Big|_{y=2x-5} = \\ =  \frac{1}{2\cdot 11} (2x-5)^{11}+C
= \frac{1}{22} (2x-5)^{11}+C
\end{multline*}\end{ex}

\begin{ex}\label{ex-12.5.6}
 \begin{multline*}
\int \frac{d x}{1+\sqrt{x}}  =  \left| \begin{array}{c} y=1+\sqrt{x}\\
x=(y-1)^2 \end{array}\right| = \\ =  \int \frac{d (y-1)^2 }{y}  \,
\Big|_{y=1+\sqrt{x}} = \\ =  \int \frac{2(y-1) d y }{y}  \,
\Big|_{y=1+\sqrt{x}} = \\ =  2\int \left( 1-\frac{1}{y}\right) \, \d
y \,\Big|_{y=1+\sqrt{x}} = \\ =  2 y-2\ln y +C \,\Big|_{y=1+\sqrt{x}}
= \\ =  2 (1+\sqrt{x})-2\ln (1+\sqrt{x}) +C
\end{multline*}\end{ex}

\begin{ers}\label{ers-12-5.7} Найдите неопределенные интегралы:

1. $\int \frac{\cos \ln x}{x}\, \d x $

2. $\int \frac{\cos x}{\sqrt{\sin^2 x+3}}\, \d x $

3. $\int \frac{2 x}{1+x^4}\, \d x $

4. $\int \frac{e^{\tg x}}{\cos^2 x}\, \d x $

5. $\int \frac{2x-1}{\sqrt{1-x+x^2}}\, \d x $

6. $\int \frac{3x^2}{x^3+1}\ln (x^3+1) \, \d x $

7. $\int \frac{e^{\arcsin x}}{\sqrt{1-x^2}}\, \d x $

8. $\int \frac{e^x}{\sqrt{4-e^{2x}}}\, \d x $

9. $\int \frac{\sin x}{\sqrt{\cos x}}\, \d x $

10. $\int x^3 \sqrt{2x^4-3}\, \d x $

11. $\int \frac{x^5}{\sqrt{7-x^6}}\, \d x $

12. $\int \frac{\cos x}{\sqrt{3+\sin x}}\, \d x $
 \end{ers}

Внесение под знак дифференциала с последующей заменой переменных позволяет
иногда упростить интеграл.

\begin{ex}\label{ex-12.5.2}
 \begin{multline*}
\int e^{\sin x}\cdot \cos x \, \d x = {\smsize\begin{pmatrix}\text{вносим
$\cos x$ под}\\ \text{знак дифференциала}\end{pmatrix}} = \\ =  \int
e^{\sin x}\, \d \sin x =  \int e^y \, \d y \, \Big|_{y=\sin x} =
\\ = e^y+C\, \Big|_{y=\sin x} = e^{\sin x}+C
\end{multline*}\end{ex}

\begin{ex}\label{ex-12.5.3}
 \begin{multline*}
\int \cos 3x  \, \d x =
{\smsize\eqref{vynesenie-konst-v-int-ot-terma}} = \\ = \frac{1}{3}\int
\cos 3x  \cdot 3 \, \d x =  {\smsize\begin{pmatrix}\text{вносим $3$ под}\\
\text{знак дифференциала}\end{pmatrix}} = \\ =  \frac{1}{3}\int \cos
3x \, \d 3x =  \frac{1}{3}\int \cos y  \, \d y \, \Big|_{y=3x} = \\
= \frac{1}{3}\sin y +C \, \Big|_{y=3x} =  \frac{1}{3}\sin 3x+C
\end{multline*}\end{ex}

\begin{ex}\label{ex-12.5.4}
 \begin{multline*}\int \tg x  \, \d x =
\int \frac{\sin x}{\cos x}  \, \d x = \\ =
{\smsize\eqref{vynesenie-konst-v-int-ot-terma}} =  -\int \frac{-\sin
x}{\cos x}  \, \d x = \\ =  {\smsize\begin{pmatrix}\text{вносим $-\sin x$ под}\\
\text{знак дифференциала}\end{pmatrix}} =  -\int \frac{1}{\cos x} \, \d
\cos x = \\ =  -\int \frac{1}{y}  \, \d y \, \Big|_{y=\cos x} =
-\ln |y|+C \, \Big|_{y=\cos x} = \\ =  -\ln |\cos x|+C
\end{multline*}\end{ex}

\begin{exs} Внесение под знак диф\-ферен\-циала:
$$
\int e^x \cdot \sin x \, \d x = -\int e^x \, \d \cos x,
$$
$$\int x
\cdot \ln x \, \d x = \frac{1}{2}\int \ln x \, \d x^2
$$
\end{exs}

\begin{er}\label{ex-12-4.3} Внесите под знак диф\-ферен\-циала второй множитель:
$$
\int e^x \cdot \sin x \, \d x, \quad \int x^4 \cdot \cos x \, \d x, \quad \int
x\cdot \frac{1}{1+x^2}\, \d x
$$
\end{er}

\begin{er}\label{ex-12-4.4} В упражнении \ref{ex-12-4.3} внесите под знак
диф\-ферен\-циала первый множитель.
\end{er}

\paragraph{Интегрирование по частям.}
Покажем, как применяется формула \eqref{int-po-chastyam-int-ot-terma}.

\begin{ex}\label{ex-int-ln-x}
 \begin{multline*}
\int \ln x \, \d x  = {\smsize\begin{pmatrix}\text{применяем}\\
\text{интегрирование}\\ \text{по частям}\end{pmatrix}} = \\ =  \ln
x\cdot
x- \int x\d \ln x = {\smsize\begin{pmatrix}\text{выносим $\ln x$}\\
\text{из-под знака}\\
\text{дифференциала}\end{pmatrix}} = \\
 = \ln x\cdot x-\int x\cdot\frac{1}{x} \d x =  \ln x\cdot x-\int 1 \d
x = \\ = \ln x\cdot x-x+C
\end{multline*}\end{ex}

\begin{ex}\label{ex-12.6.2}\begin{multline*}\int x e^x \, \d x  =  {\smsize\begin{pmatrix}\text{вносим $e^x$}\\ \text{под знак}\\ \text{дифференциала}\end{pmatrix}} =
\int x  \, \d e^x  = \\ =  {\smsize\begin{pmatrix}\text{применяем}\\
\text{интегрирование}\\ \text{по частям}\end{pmatrix}} =  x e^x-\int e^x
\, \d x  = x e^x- e^x +C
\end{multline*}\end{ex}

\begin{ex}\label{ex-12.6.3}
 \begin{multline*}\int x \sin x \, \d x  =
{\smsize\begin{pmatrix}\text{вносим $\sin x$}\\ \text{под знак}\\
\text{дифференциала}\end{pmatrix}} = \\ =
-\int x  \, \d \cos x  =  {\smsize\begin{pmatrix}\text{применяем}\\
\text{интегрирование}\\ \text{по частям}\end{pmatrix}} = \\ =
-\left(x \cos x+\int \cos x  \, \d x \right) = \\ =  -\left(x \cos x+
\sin x +C_1\right) =  -x \cos x- \sin x +C
\end{multline*}\end{ex}

\begin{ex}\label{ex-12.6.4}
 \begin{multline*}
\int x^2 \ln x \, \d x  =  {\smsize\begin{pmatrix}\text{вносим $x^2$}\\
\text{под знак}\\ \text{дифференциала}\end{pmatrix}} = \\ =
\frac{1}{3}\int \ln x \, \d x^3   =
{\smsize\begin{pmatrix}\text{применяем}\\ \text{интегрирование}\\
\text{по частям}\end{pmatrix}} = \\ =  \frac{1}{3}\left(\ln x \cdot
x^3 -\int x^3 \, \d \ln x \right) = \\ =
{\smsize\begin{pmatrix}\text{выносим $\ln x$}\\ \text{из-под знака}\\
\text{дифференциала}\end{pmatrix}} = \\ =  \frac{1}{3}\left(\ln x
\cdot x^3 -\int x^3 \frac{1}{x}\, \d x \right) = \\ =
\frac{1}{3}\left(\ln x \cdot x^3 -\int x^2 \, \d x \right) = \\ =
\frac{1}{3}\left(\ln x \cdot x^3 -\frac{1}{3} x^3 +C \right) =
\frac{1}{3}\ln x \cdot x^3 -\frac{1}{9} x^3 +C
\end{multline*}\end{ex}

\begin{ers}\label{ers-12.6.5} Найдите интегралы:
 \begin{multicols}{2}
1) $\int \frac{\ln x}{x^2}\, \d x$

2) $\int x e^{2x}\, \d x$

3) $\int \arctg x \, \d x$

4) $\int x\cos x \, \d x$

5) $\int x\arctg x \, \d x$

6) $\int \frac{x}{\sin^2 x}\, \d x$

7) $\int \sqrt{x}\ln x \, \d x$

8) $\int x^n \ln x \, \d x$

9) $\int \arcsin x \, \d x$
 \end{multicols}\end{ers}

Иногда бывает необходимо применять формулу интегрирования по частям несколько
раз.

\begin{ex}\label{ex-12.6.6}
 \begin{multline*}
\int x^2 e^x \, \d x  =  {\smsize\begin{pmatrix}\text{вносим $e^x$}\\
\text{под знак}\\ \text{дифференциала}\end{pmatrix}} =
\int x^2  \, \d e^x  = \\ = {\smsize\begin{pmatrix}\text{применяем}\\
\text{интегрирование}\\ \text{по частям}\end{pmatrix}} =  x^2 e^x-\int e^x
\, \d x^2  = \\ =  {\smsize\begin{pmatrix}\text{выносим
$x^2$}\\ \text{из-под знака}\\
\text{дифференциала}\end{pmatrix}} =  x^2 e^x-2\int e^x  x\, \d x  =
\\ =
{\smsize\begin{pmatrix}\text{вносим $e^x$}\\ \text{под знак}\\
\text{дифференциала}\end{pmatrix}} =  x^2 e^x-2\int  x\, \d e^x
 = \\ =  {\smsize\begin{pmatrix}\text{применяем}\\ \text{интегрирование}\\
\text{по частям}\end{pmatrix}} =  x^2 e^x-2\left( x e^x -\int e^x\, \d x
\right)   = \\ =  x^2 e^x-2\left( x e^x - e^x -C_1 \right)   =
\\ =  x^2 e^x-2x e^x - 2e^x +C
\end{multline*}\end{ex}

\begin{ex}\label{ex-12.6.7}
 \begin{multline*}
\int x^2 \cos x \, \d x  =  {\smsize\begin{pmatrix}\text{вносим $\cos
x$}\\ \text{под знак}\\ \text{дифференциала}\end{pmatrix}} = \\ =
\int
x^2  \, \d \sin x  =  {\smsize\begin{pmatrix}\text{применяем}\\
\text{интегрирование}\\ \text{по частям}\end{pmatrix}} = \\ =  x^2
\sin x-\int \sin x  \, \d x^2  =  {\smsize\begin{pmatrix}\text{выносим
$x^2$}\\ \text{из-под знака}\\
\text{дифференциала}\end{pmatrix}} = \\ =  x^2 \sin x-2\int x \sin x
\,
\d x  =  {\smsize\begin{pmatrix}\text{вносим $\sin x$}\\ \text{под знак}\\
\text{дифференциала}\end{pmatrix}} = \\ =  x^2 \sin x+2\int
x\, \d \cos x   =  {\smsize\begin{pmatrix}\text{применяем}\\
\text{интегрирование}\\ \text{по частям}\end{pmatrix}} = \\ =  x^2
\sin x+2\left( x \cos x -\int  \cos x\, \d x \right)   = \\ =  x^2
\sin x+2\left( x \cos x - \sin x -C_1 \right)   = \\ =  x^2 \sin x-2x
\cos x - 2\sin x +C
\end{multline*}\end{ex}

\begin{ers}\label{ers-12.6.8} Найдите интегралы:

1) $\int x^2\sin x \, \d x$

2) $\int x^2  e^{-x}\, \d x$

3) $\int x  \ln^2 x \, \d x$

4) $\int x^3 \, \d \cos x$,

5) $\int \cos x \, \d x^3$,

6) $\int x \, \d e^x$.
\end{ers}

Случается, что интегрирование по частям приводит к исходным интегралам. В таких
случаях нужно поступать следующим образом.

\begin{ex}\label{ex-12.6.9}
 \begin{multline*}
\int e^x \sin x \, \d x  =  {\smsize\begin{pmatrix}\text{вносим $e^x$}\\
\text{под знак}\\ \text{дифференциала}\end{pmatrix}} = \\ =
\int \sin x  \, \d e^x  =  {\smsize\begin{pmatrix}\text{применяем}\\
\text{интегрирование}\\ \text{по частям}\end{pmatrix}} = \\ =  e^x
\sin x-\int e^x  \, \d \sin x  =  {\smsize\begin{pmatrix}\text{выносим
$\sin x$}\\ \text{из-под знака}\\
\text{дифференциала}\end{pmatrix}} = \\ =  e^x \sin x-\int e^x \cos x
\,
\d x  =  {\smsize\begin{pmatrix}\text{вносим $e^x$}\\ \text{под знак}\\
\text{дифференциала}\end{pmatrix}} = \\ =  e^x \sin x-\int \cos x \,
\d e^x   =  {\smsize\begin{pmatrix}\text{применяем}\\
\text{интегрирование}\\ \text{по частям}\end{pmatrix}} = \\ =  e^x
\sin x-\left(e^x \cos x -\int  e^x \, \d \cos x\right)   = \\ =  e^x
\sin x-e^x \cos x +\int  e^x \, \d \cos x = \\ =
{\smsize\begin{pmatrix}\text{выносим $\cos x$}\\ \text{из-под знака}\\
\text{дифференциала}\end{pmatrix}} = \\ =  e^x \sin x-e^x \cos x
-\int e^x \sin x\, \d x
\end{multline*}
Мы вернулись к исходному интегралу. Если теперь обозначить
$$
I=\int  e^x \sin x\, \d x,
$$
то мы получим уравнение
$$
I = e^x \sin x-e^x \cos x -I,
$$
из которого легко находится $I$:
$$
2I = e^x \sin x-e^x \cos x
$$
$$
\Downarrow
$$
$$
I = \frac{1}{2}\left( e^x \sin x-e^x \cos x \right)
$$
Теперь нужно вспомнить, что неопределенный интеграл вычисляется с точностью до
константы, поэтому окончательный ответ будет таким:
$$
\int  e^x \sin x\, \d x=\frac{1}{2}\left( e^x \sin x-e^x \cos x \right)+С
$$
\end{ex}

\begin{ex}\label{ex-12.6.10}
 \begin{multline*}
\int \sqrt{1-x^2}\, \d x  =  {\smsize\begin{pmatrix}\text{применяем}\\
\text{интегрирование}\\ \text{по частям}\end{pmatrix}} = \\ =  x
\sqrt{1-x^2} - \int x \d\sqrt{1-x^2}  =
{\smsize\begin{pmatrix}\text{выносим $\sqrt{1-x^2}$}\\ \text{из-под знака}\\
\text{дифференциала}\end{pmatrix}} = \\ =  x \sqrt{1-x^2} + \int
\frac{x^2}{\sqrt{1-x^2}}\, \d x  = \\ =  x \sqrt{1-x^2} - \int
\frac{1-x^2}{\sqrt{1-x^2}}\, \d x +\int \frac{1}{\sqrt{1-x^2}}\, \d x  =
\\ =  x \sqrt{1-x^2} - \int \sqrt{1-x^2}\, \d x +\arcsin x+C
\end{multline*}
Мы вернулись к исходному интегралу. Если теперь обозначить
$$
I=\int \sqrt{1-x^2}\, \d x ,
$$
то получается алгебраическое уравнение, из которого находится $I$:
$$
I = x \sqrt{1-x^2} - I +\arcsin x+C
$$
$$
\Downarrow
$$
$$
2I = x \sqrt{1-x^2} +\arcsin x+C
$$
$$
\Downarrow
$$
$$
I = \frac{1}{2}\left( x \sqrt{1-x^2} +\arcsin x +C\right)
$$
Выбрав новую константу, получаем окончательный ответ:
$$
\int \sqrt{1-x^2}\, \d x =\frac{1}{2}\left( x \sqrt{1-x^2} +\arcsin x \right)+A
$$
\end{ex}

\begin{ers}\label{ERS:int-po-chastyam} Найдите интегралы:

1) $\int \sqrt{1+x^2}\, \d x$

2) $\int e^x  \cos x \, \d x$

3) $\int e^{2x}  \sin^2 x \, \d x$

4) $\int \cos (\ln x) \, \d x$

5) $\int e^{\arcsin x}\, \d x$
 \end{ers}

\section{Специальные методы вычисления неопределенных интегралов}\label{CH-comp-indef-integr}

В этом параграфе мы опишем некоторые приемы вычисления неопределенных
интегралов.

\subsection{Рекуррентные формулы}

Некоторые интегралы удается найти только с помощью так называемых рекуррентных
интегральных формул. Что это такое станет ясно из следующего примера.

\begin{tm}\label{tm-13.1.1} Для интегралов
$$
  I_n=\int \frac{d x}{(x^2+a^2)^n}\label{13.1.1}
$$
первый из которых, по формуле \eqref{12.2.9}, равен
$$
I_1=\frac{1}{a}\arctg \frac{x}{a}+C \label{13.1.2}
$$
справедлива следующая рекуррентная формула
$$
I_{n+1}=\frac{1}{2n a^2}\cdot \frac{x}{(x^2+a^2)^n}+\frac{2n-1}{2n a^2}\cdot
I_n \qquad (n\in \mathbb{N}) \label{13.1.3}
$$
\end{tm}
\begin{proof}
 \begin{multline*}
I_{n+1}=\int \frac{d x}{(x^2+a^2)^{n+1}}=\\= \frac{1}{a^2}\int
\frac{a^2}{(x^2+a^2)^{n+1}}\, \d x=\\= \frac{1}{a^2}\int
\frac{x^2+a^2-x^2}{(x^2+a^2)^{n+1}}\, \d x=\\= \frac{1}{a^2}\int
\frac{1}{(x^2+a^2)^n}d x- \frac{1}{a^2}\int \frac{x^2}{(x^2+a^2)^{n+1}}\, \d x=\\=
\frac{1}{a^2} I_n-\frac{1}{a^2}\int x^2\cdot (x^2+a^2)^{-(n+1)}\, \d x=\\=
{\smsize\begin{pmatrix}\text{вносим $x$}\\
\text{под знак}\\ \text{дифференциала}\end{pmatrix}}=\\=
\frac{1}{a^2}I_n-\frac{1}{2a^2}\int x\cdot (x^2+a^2)^{-(n+1)}\, \d
(x^2+a^2)=\\=
{\smsize\begin{pmatrix}\text{вносим $(x^2+a^2)^{-(n+1)}$}\\
\text{под знак дифференциала}\end{pmatrix}}=\\= \frac{1}{a^2} I_n+\frac{1}{2n
a^2}\int x \, \d (x^2+a^2)^{-n}=
{\smsize\begin{pmatrix}\text{интегрируем}\\
\text{по частям}\end{pmatrix}}=\\= \frac{1}{a^2} I_n+\frac{1}{2n a^2} x \cdot
(x^2+a^2)^{-n} -\\-\frac{1}{2n a^2}\int (x^2+a^2)^{-n}\, \d x=\\= \frac{1}{a^2}
I_n+\frac{1}{2n a^2}\frac{x}{(x^2+a^2)^n} -\frac{1}{2n a^2}\int  \frac{d
x}{(x^2+a^2)^n}=\\= \frac{1}{a^2} I_n+\frac{1}{2n a^2}\frac{x}{(x^2+a^2)^n}
-\frac{1}{2n a^2} I_n=\\= \frac{1}{2n a^2}\cdot \frac{x}{(x^2+a^2)^n}
+\frac{2n-1}{2n a^2}\cdot I_n
\end{multline*}\end{proof}

Теперь покажем, как применяется рекуррентная формула \eqref{13.1.3}.

\begin{ex}\label{ex-13.1.2} Вычислим, например, интеграл $I_2$.
Для этого в формулу \eqref{13.1.3} подставим $n=1$:
 \begin{multline*}
I_2=\frac{1}{2 a^2}\cdot \frac{x}{(x^2+a^2)}+\frac{1}{2 a^2}\cdot I_1={\smsize
\eqref{13.1.2}}=\\= \frac{x}{2 a^2(x^2+a^2)}+\frac{1}{2
a^3}\cdot \arctg \frac{x}{a}+C
 \end{multline*}
То есть,
 \begin{multline}
I_2=\int \frac{d x}{(x^2+a^2)^2} =\\= \frac{x}{2 a^2(x^2+a^2)}+\frac{1}{2
a^3}\cdot \arctg \frac{x}{a}+C \label{13.1.4}
 \end{multline}
\end{ex}

\begin{ex}\label{ex-13.1.3} Вычислим, интеграл $I_3$,
для чего в формулу \eqref{13.1.3} подставим $n=2$:
 \begin{multline*}
I_3=\frac{1}{4 a^2}\cdot \frac{x}{(x^2+a^2)^2}+\frac{3}{4 a^2}\cdot
I_2={\smsize \eqref{13.1.4}}=\\= \frac{x}{4
a^2(x^2+a^2)^2}+\\+\frac{3}{4 a^2}\cdot \left(\frac{x}{2
a^2(x^2+a^2)}+\frac{1}{2 a^3}\cdot \arctg \frac{x}{a}+C \right)=\\= \frac{x}{4
a^2(x^2+a^2)^2}+ \frac{3x}{8 a^4(x^2+a^2)}+\frac{3}{8 a^5}\cdot \arctg
\frac{x}{a}+C
 \end{multline*}
То есть,
 \begin{multline}
I_3=\int \frac{d x}{(x^2+a^2)^3} =\\= \frac{x}{4 a^2(x^2+a^2)^2}+ \frac{3x}{8
a^4(x^2+a^2)}+\frac{3}{8 a^5}\cdot \arctg \frac{x}{a}+C \label{13.1.5}
 \end{multline}
\end{ex}

\begin{ers} Выведите рекуррентные формулы для следующих
интегралов, и вычислите их при $n=2,3$:

1) $J_n=\int x^k \ln^n x \, \d x \quad (n\in \mathbb{N})$

2) $K_n=\int \sin^n x \, \d x \quad (n\in \mathbb{N})$
 \end{ers}

\subsection{Интегрирование рациональных функций}

\biter{

\item[$\bullet$] Функции вида
$$
f(x)=\sum_{k=0}^n a_k\cdot x^k\qquad (x\in\R),
$$
где $n\in\N$ и $a_k\in\R$, называются {\it алгебраическими многочленами}\label{EX:algebr-mnogochl} или
просто {\it многочленами} (от одной переменной). Если $a_n\ne 0$, то число $n$
называется {\it степенью многочлена} $f$.

\item[$\bullet$] {\it Рациональной функцией}\index{функция!рациональная} называется всякая
функция $R(x)$, которую можно представить в виде отношения двух многочленов:
$$
R(x)=\frac{P(x)}{Q(x)}
$$
}\eiter

\bex Например, функции
$$
\frac{x^2+x+1}{x-2}, \quad \frac{2x+1}{x^3-x+5}, \quad x^5-2x+1
$$
являются рациональными, а функции
$$
\sin x, \quad \frac{x}{\sqrt{x}-3}, 2^x
$$
не являются рациональными.
\eex

В этом параграфе мы покажем, как вычисляются интегралы от рациональных функций.

\paragraph{Интегрирование простейших дробей}

Самыми простыми рациональными функциями являются так называемые {\it простейшие
дроби}\index{дробь!простейшая}, то есть дроби следующих четырех типов.

\bigskip
\centerline{\bf Типы простейших дробей}
 \biter{
\item[---] простейшие дроби 1 типа:
$$\frac{A}{x-a}$$
\item[---] простейшие дроби 2 типа:
$$\frac{A}{(x-a)^n}\quad (n\in \mathbb{N}, \, n>1)$$
\item[---] простейшие дроби 3 типа:
$$\frac{Bx+C}{x^2+px+q}\quad (p^2-4q<0)$$
\item[---] простейшие дроби 4 типа:
$$
\frac{Bx+C}{(x^2+px+q)^n}\quad (p^2-4q<0, \, n\in \mathbb{N}, \, n>1)
$$
 }\eiter

Покажем, как они интегрируются.

\begin{ex}[\bf простейшая дробь 1 типа]\label{ex-13.2.1}
Для таких дробей можно вывести общую формулу:
 \begin{multline*}
\int \frac{A}{x-a}\, \d x= \left|
\begin{array}{c}
y=x-a \\ x=y+a
\end{array}\right|=\\= \int \frac{A}{y}\, \d (y+a) \, \Big|_{y=x-a}= A \int
\frac{d y}{y}\, \Big|_{y=x-a}=\\= A \ln |y| +C \, \Big|_{y=x-a}= A \ln |x-a| +C
 \end{multline*}\end{ex}

\begin{ex}[\bf простейшая дробь 2 типа]\label{ex-13.2.2}
Здесь также можно найти общую формулу
\begin{multline*}\int \frac{A}{(x-a)^n}\, \d x= \left|
\begin{array}{c}
y=x-a \\ x=y+a
\end{array}\right|=\\=
\int \frac{A}{y^n}\, \d (y+a) \, \Big|_{y=x-a}= A \int \frac{d y}{y^n}\,
\Big|_{y=x-a}=\\= \frac{A}{(n-1)y^{n-1}} +C \, \Big|_{y=x-a}=
\frac{A}{(n-1)(x-a)^{n-1}} +C
\end{multline*}\end{ex}

\begin{ex}[\bf простейшая дробь 3 типа]\label{ex-13.2.3}
Здесь, хотя и можно вывести общую формулу, но запомнить ее уже будет вряд ли
возможно. Поэтому мы приведем лишь общий план действий. Чтобы вычислить
интеграл
$$
\int \frac{Bx+C}{x^2+px+q}\, \d x \qquad (p^2-4q<0)
$$
нужно сначала {\it выделить полный квадрат}\index{полного квадрата выделение} в
знаменателе, то есть переписать знаменатель так, чтобы переменная $x$
встречалась только один раз:
 \begin{multline*}
x^2+px+q=x^2+2\cdot \frac{p}{2}\cdot x+
\left(\frac{p}{2}\right)^2-\left(\frac{p}{2}\right)^2+q=\\= \left(
x+\frac{p}{2}\right)^2+\left(q-\frac{p^2}{4}\right)
 \end{multline*}
а затем выражение в скобках заменить на новую переменную:
$$
x+\frac{p}{2}=y
$$
Оказывается, что при этом мы получим табличные интегралы.

Приведем иллюстрацию. Пусть требуется вычислить интеграл
$$
\int \frac{3x-1}{x^2-x+1}\, \d x
$$
Поскольку дискриминант знаменателя отрицателен,
$$
p^2-4q=1-4=-3<0
$$
выражение под интегралом действительно является простейшей дробью третьего
типа. Выделим полный квадрат в знаменателе:
 \begin{multline*}
x^2-x+1=x^2-2\cdot \frac{1}{2}\cdot x+
\left(\frac{1}{2}\right)^2-\left(\frac{1}{2}\right)^2+1=\\= \left(
x-\frac{1}{2}\right)^2+\left(1-\frac{1}{4}\right)= \left(
x-\frac{1}{2}\right)^2+\frac{3}{4}
 \end{multline*}
Теперь делаем замену переменной в интеграле:
 \begin{multline*}
\int \frac{3x-1}{x^2-x+1}\, \d x= \int \frac{3x-1}{\left(
x-\frac{1}{2}\right)^2+\frac{3}{4}}\, \d x=\\= \left|
\begin{array}{c}
y=x-\frac{1}{2}\\
x=y+\frac{1}{2}\end{array}\right|= \int
\frac{3\left(y+\frac{1}{2}\right)-1}{y^2+\frac{3}{4}}\, \d
\left(y+\frac{1}{2}\right)=\\= \int \frac{3y+\frac{1}{2}}{y^2+\frac{3}{4}}\, \d
y=\\= 3\int \frac{y}{y^2+\frac{3}{4}}\, \d y+ \frac{1}{2}\int
\frac{1}{y^2+\frac{3}{4}}\, \d y=\\= {\smsize\begin{pmatrix}\text{в первом интеграле}\\
\text{вносим $y$ под знак}\\
\text{дифференциала}\end{pmatrix}}=\\= \frac{3}{2}\int \frac{d
y^2}{y^2+\frac{3}{4}} + \frac{1}{2}\int \frac{1}{y^2+\frac{3}{4}}\,
\d y=\\= {\smsize\begin{pmatrix}\text{в первом интеграле}\\
\text{делаем замену}\\
y^2=z-\frac{3}{4}\\
z=y^2+\frac{3}{4}\end{pmatrix}}=\\= \frac{3}{2}\int \frac{d
\left(z-\frac{3}{4}\right)}{z}\, \Big|_{y=z-\frac{3}{4}z=y+\frac{3}{4}} +
\frac{1}{2}\int \frac{d y}{y^2+\frac{3}{4}} =\\= \frac{3}{2}\int \frac{d
z}{z}\, \Big|_{z=y^2+\frac{3}{4}} + \frac{1}{2}\int \frac{d
y}{y^2+\left(\frac{\sqrt{3}}{2}\right)^2} =\\
= {\smsize\begin{pmatrix}\text{применяем}\\
\text{формулы}\\ \text{\eqref{12.2.2} и \eqref{12.2.9}}\end{pmatrix}}=\\=
\frac{3}{2}\ln |z|  \, \Big|_{z=y^2+\frac{3}{4}} + \frac{1}{2}\cdot
\frac{2}{\sqrt{3}}\cdot \arctg \frac{2y}{\sqrt{3}}+C=\\= \frac{3}{2}\ln
\left|y^2+\frac{3}{4}\right| + \frac{1}{\sqrt{3}}\cdot \arctg
\frac{2y}{\sqrt{3}}+C=\\= {\smsize\begin{pmatrix}\text{возвращаемся}\\
\text{к исходной}\\
\text{переменной $x$}\end{pmatrix}}=\\= \frac{3}{2}\ln
\left|\left(x-\frac{1}{2}\right)^2+\frac{3}{4}\right| + \frac{1}{\sqrt{3}}\cdot
\arctg \frac{2\left( x-\frac{1}{2}\right)}{\sqrt{3}}+C=\\= \frac{3}{2}\ln
\left|x^2-x+1 \right| + \frac{1}{\sqrt{3}}\cdot \arctg
\frac{2x-1}{\sqrt{3}}+C=\\=
{\smsize\begin{pmatrix}\text{всегда}\\
x^2-x+1>0
\end{pmatrix}}=\\= \frac{3}{2}\ln \left( x^2-x+1 \right) +
\frac{1}{\sqrt{3}}\cdot \arctg \frac{2x-1}{\sqrt{3}}+C
\end{multline*}\end{ex}

\begin{ex}[\bf простейшая дробь 4 типа]\label{ex-13.2.4}
 Для таких дробей
вывести общую формулу и вовсе невозможно. План действий же состоит в следующем.
Чтобы вычислить интеграл
$$
\int \frac{Bx+C}{(x^2+px+q)^n}\, \d x \qquad (n\in \mathbb{N}, \, n>1, \,
p^2-4q<0)
$$
нужно, как и в предыдущем случае, сначала выделить полный квадрат у квадратного
трехчлена в знаменателе,
 \begin{multline*}
x^2+px+q=x^2+2\cdot \frac{p}{2}\cdot x+
\left(\frac{p}{2}\right)^2-\left(\frac{p}{2}\right)^2+q=\\= \left(
x+\frac{p}{2}\right)^2+\left(q-\frac{p^2}{4}\right)
 \end{multline*}
затем выражение в скобках заменить на новую переменную:
$$
x+\frac{p}{2}=y
$$
При этом получится несколько интегралов, из которых одни окажутся табличными, а
другие можно будет вычислить с помощью рекуррентных формул $\S 7$.

Приведем иллюстрацию. Пусть требуется вычислить интеграл
$$
\int \frac{x-1}{x^2+2x+3}\, \d x
$$
Поскольку дискриминант знаменателя отрицателен,
$$
p^2-4q=4-12=-8<0
$$
выражение под интегралом действительно является простейшей дробью четвертого
типа. Выделим полный квадрат в знаменателе:
$$
x^2+2x+3=x^2+2\cdot 1\cdot x+1^2-1^2+3=(x+1)^2+2
$$
Теперь делаем замену переменной в интеграле:
 \begin{multline*}
\int \frac{x-1}{(x^2+2x+3)^2}\, \d x= \int
\frac{x-1}{\left((x+1)^2+2\right)^2}\, \d x=\\= \left|
\begin{array}{c}
y=x+1 \\
x=y-1
\end{array}\right|= \int \frac{y-2}{(y^2+2)^2}\, \d (y-2)=\\= \int
\frac{y-2}{(y^2+2)^2}\, \d y= \int \frac{y \,  d y}{(y^2+2)^2} -2\int \frac{d
y}{(y^2+2)^2}=\\= \frac{1}{2}\int \frac{d y^2}{(y^2+2)^2}-2\int \frac{d
y}{(y^2+2)^2}=\\= \frac{1}{2}\int \frac{d (z-2)}{z^2}\, \Big|_{z=y^2+2}-2\int
\frac{d y}{(y^2+2)^2}=\\= \frac{1}{2}\int \frac{d z}{z^2}\, \Big|_{z=y^2+2}
-2\int \frac{d y}{\left(y^2+\left(\sqrt{2}\right)^2 \right)^2}=\\=
{\smsize\begin{pmatrix}\text{применяем}\\
\text{формулы}\\ \text{\eqref{12.2.1} и
\eqref{13.1.4}}\end{pmatrix}}=-\frac{1}{2z}\, \Big|_{z=y^2+2} -\\-2
\left(\frac{y}{4(y^2+2)}+\frac{1}{2 (\sqrt{2})^3}\cdot \arctg
\frac{y}{\sqrt{2}}+C \right)=\\= -\frac{1}{2(y^2+2)} -\\-2
\left(\frac{y}{4(y^2+2)}+\frac{1}{2 (\sqrt{2})^3}\cdot \arctg
\frac{y}{\sqrt{2}}+C \right)=\\= -\frac{1}{2(y^2+2)} -
\frac{y}{2(y^2+2)}-\frac{1}{2 \sqrt{2}}\cdot \arctg \frac{y}{\sqrt{2}}+C =\\=
-\frac{y+1}{2(y^2+2)}-\frac{\sqrt{2}}{4}\cdot \arctg
\frac{y}{\sqrt{2}}+C =\\= {\smsize\begin{pmatrix}\text{возвращаемся}\\
\text{к исходной}\\
\text{переменной $x$}\end{pmatrix}}=\\= -\frac{x+2}{2(x^2+2x+3)}-
\frac{\sqrt{2}}{4}\cdot \arctg \frac{x+1}{\sqrt{2}}+C
\end{multline*}\end{ex}

\begin{ers} Вычислите интегралы:
 \begin{multicols}{2}
1) $\int \frac{d x}{x-7}$;

2) $\int \frac{d x}{3-x}$;

3) $\int \frac{d x}{1+5x}$;

4) $\int \frac{d x}{(3x+2)^2}$;

5) $\int \frac{d x}{x^2+2x+1}$;

6) $\int \frac{d x}{4x^2-4x+1}$;

7) $\int \frac{x \, \d x}{x^2+7x+13}$;

8) $\int \frac{4x+8}{3x^2+2x+5}\, \d x$;

9) $\int \frac{x+5}{2x^2+2x+3}\, \d x$;

10) $\int \frac{5x-7}{8x^2+x+1}\, \d x$;

11) $\int \frac{3x+5}{(x^2+2x+2)^2}\, \d x$;

12) $\int \frac{3x+2}{(x^2-3x+3)^2}\, \d x$.
\end{multicols}\end{ers}

\paragraph{Интегрирование правильных дробей}

\biter{

\item[$\bullet$]
Рациональная функция
$$
R(x)=\frac{P(x)}{Q(x)}
$$
называется {\it правильной дробью}\index{дробь!правильная}, если степень
многочлена в числителе меньше степени многочлена в знаменателе:
$$
  \deg P<\deg Q
$$
}\eiter

\bex

Например, дроби
$$
  \frac{x+1}{x^2+x+1}, \quad
  \frac{x^3}{(x^2-1)(x^2+4)}, \quad
  \frac{x^2-2x+1}{x^3}, \quad
$$
будут правильными, а дроби
$$
  \frac{x^2+1}{x^2+x+1}, \quad
  \frac{x^3}{(x-1)(x+4)}, \quad
  \frac{x^2-2x+1}{x}, \quad
$$
-- неправильными.

\eex

Для интегрирования правильных дробей нужно научиться пользоваться следующей
классической теоремой.

\begin{tm}\label{tm-13.2.6} Всякая правильная дробь однозначным образом
раскладывается в сумму простейших дробей.
\end{tm}

Мы не будем доказывать этот результат, а лишь объясним как он работает. Для
этого полезно запомнить следующие

\bigskip

\centerline{\bf Правила разложения на простейшие
дроби:}\label{pravila-razlozh-na-prost-drobi}
 \biter{
\item[---] если в разложении знаменателя $Q(x)$ на множители имеется множитель
$(x-a)^n$
$$
  Q(x)=... \cdot (x-a)^n \cdot ...
$$
то в разложении правильной дроби $\frac{P(x)}{Q(x)}$ в сумму простейших дробей
должны присутствовать слагаемые вида $\frac{A_1}{x-a}, \frac{A_2}{(x-a)^2},
..., \frac{A_n}{(x-a)^n}$:
$$
\frac{P(x)}{Q(x)}=
...+\frac{A_1}{x-a}+\frac{A_2}{(x-a)^2}+...+\frac{A_n}{(x-a)^n}+...
$$
\item[---] если в разложении знаменателя $Q(x)$ на множители имеется множитель
$(x^2+px+q)^n \, (p^2-4q<0)$
$$
  Q(x)=... \cdot (x^2+px+q)^n \cdot ...
$$
то в разложении правильной дроби $\frac{P(x)}{Q(x)}$ в сумму простейших дробей
должны присутствовать слагаемые вида $\frac{B_1x+C_1}{x^2+px+q}, \frac{B_2
x+C_2}{(x^2+px+q)^2}, ..., \frac{B_n x+C_n}{(x^2+px+q)^n}$:
\begin{multline*}
\frac{P(x)}{Q(x)}= ...+\frac{B_1x+C_1}{x^2+px+q}+\frac{B_2
x+C_2}{(x^2+px+q)^2}+...\\
...+ \frac{B_n x+C_n}{(x^2+px+q)^n}+...
\end{multline*}
 }\eiter

\bigskip

Покажем теперь на примерах, что эти правила означают.

\begin{ex}[\bf случай неповторяющихся мно\-жи\-телей первой степени]
\label{ex-13.2.7} Рассмотрим интеграл
$$
\int \frac{9x^2-2x-8}{x^3-4x}\, \d x
$$
Разложим знаменатель дроби на множители:
$$
x^3-4x=x(x^2-4)=x(x-2)(x+2)
$$
В соответствии с первым правилом, это означает, что в разложении нашей
(правильной) дроби в сумму простейших дробей должны присутствовать следующие
слагаемые:
$$
\frac{A}{x}, \, \frac{B}{x-2}, \, \frac{C}{x+2}
$$
То есть,
$$
\frac{9x^2-2x-8}{x^3-4x}=\frac{A}{x}+\frac{B}{x-2}+\frac{C}{x+2}
$$
Теперь необходимо найти коэффициенты $A,B,C$. Для этого приведем правую часть к
общему знаменателю
 \begin{multline*}
\frac{9x^2-2x-8}{x^3-4x}=\\=\frac{A(x-2)(x+2)+Bx(x+2)+Cx(x-2)}{x(x-2)(x+2)}
 \end{multline*}
Эти дроби действительно будут совпадать, если у них совпадают числители:
$$
9x^2-2x-8=A(x-2)(x+2)+Bx(x+2)+Cx(x-2) \label{13.2.1}
$$

Далее коэффициенты $A,B,C$ можно искать двумя разными методами.

Первый из них называется {\it методом неопределенных
коэффициентов}\index{метод!неопределенных коэффициентов}. Он состоит в
следующем. Сначала приведем подобные слагаемые справа:
$$
9x^2-2x-8=(A+B+C)x^2+(2B-2C)x-4A
$$
Слева и справа записаны многочлены от переменной $x$. Они будут совпадать
только если у них совпадают коэффициенты:
 \begin{align*}
x^2: & \quad 9=A+B+C  \\
x^1: & \quad -2=2B-2C \\
x^0: & \quad -8=-4A
\end{align*}\noindent
Эта запись означает, что число 9 является коэффициентом перед $x^2$ слева, и
оно совпадает с числом $A+B+C$, которое является коэффициентом перед $x^2$
справа, и т.д. Мы получаем линейную систему, решая которую находим числа
$A,B,C$:
 \begin{multline*}
\begin{cases}{A+B+C=9}\\{2B-2C=-2}\\{-4A=-8}\end{cases}
\; \Leftrightarrow \;
\begin{cases}{A+B+C=9}\\{B-C=-1}\\{A=2}\end{cases}
\; \Leftrightarrow
\\ \Leftrightarrow \;
\begin{cases}{2+B+C=9}\\{B-C=-1}\\{A=2}\end{cases}
\; \Leftrightarrow \;
\begin{cases}{B+C=7}\\{B-C=-1}\\{A=2}\end{cases}
\; \Leftrightarrow
\\
\Leftrightarrow \;
\begin{cases}{B+C=7}\\{2B=6}\\{A=2}\end{cases}\; \Leftrightarrow \;
\begin{cases}{B+C=7}\\{B=3}\\{A=2}\end{cases}
\; \Leftrightarrow
\\
\Leftrightarrow \;
\begin{cases}{3+C=7}\\{B=3}\\{A=2}\end{cases}\; \Leftrightarrow \;
\begin{cases}{C=4}\\{B=3}\\{A=2}\end{cases}
 \end{multline*}

Мы нашли коэффициенты $A,B,C$ методом неопределенных коэффициентов. Теперь
покажем, как их можно было найти по-другому, а именно {\it методом частных
значений}\index{метод!частных значений}.

Еще раз рассмотрим равенство \eqref{13.2.1}:
$$
9x^2-2x-8=A(x-2)(x+2)+Bx(x+2)+Cx(x-2) \label{13.2.2}
$$
Оно должно выполняться при любом $x$. В частности, при  $x=0$, $x=2$ и $x=-2$.
Посмотрим, во что превращаются левая и правая части при этих значениях
переменной:
 \begin{align*}
x=0: & \quad -8=-4A \qquad \Rightarrow \qquad A=2  \\
x=2: & \quad 24=8B  \qquad \Rightarrow \qquad B=3 \\
x=-2: & \quad 32=8C \qquad \Rightarrow \qquad C=4
\end{align*}\noindent Эта запись означает, что при $x=0$ равенство
\eqref{13.2.2} превращается в равенство $-8=-4A$, откуда следует что $A=2$, и
т.д.

В общем, каким способом ни считай, получается
$$
\begin{cases}{A=2}\\{B=3}\\{C=4}\end{cases}
$$
Это означает, что
$$
\frac{9x^2-2x-8}{x^3-4x}=\frac{2}{x}+\frac{3}{x-2}+\frac{4}{x+2}
$$
и теперь, наконец, можно вычислить наш интеграл:
 \begin{multline*}
\int \frac{9x^2-2x-8}{x^3-4x}\, \d x =\\= \int \frac{2}{x}\, \d x + \int
\frac{3}{x-2}\, \d x + \int \frac{4}{x+2}\, \d x=\\= 2\ln |x| + 3\ln |x-2| +
4\ln |x+2|+D
 \end{multline*}
($D$ -- постоянная интегрирования).
\end{ex}

\begin{ex}[\bf случай повторяющихся мно\-жи\-телей первой степени]\label{ex-13.2.8}
Рассмотрим интеграл
$$
\int \frac{3x+2}{x(x+1)^3}\, \d x
$$
Знаменатель дроби уже разложен на множители:
$$
x(x+1)^3
$$
В соответствии с первым правилом, в разложении нашей (правильной) дроби в сумму
простейших дробей должны присутствовать следующие слагаемые:
$$
\frac{A}{x}, \quad \frac{B}{x+1}, \quad \frac{C}{(x+1)^2}, \quad
\frac{D}{(x+1)^3}
$$
Значит,
$$
\frac{3x+2}{x(x+1)^3}=
\frac{A}{x}+\frac{B}{x+1}+\frac{C}{(x+1)^2}+\frac{D}{(x+1)^3}
$$
Теперь необходимо найти коэффициенты $A,B,C,D$. Для этого приведем правую часть
к общему знаменателю
 \begin{multline*}
\frac{3x+2}{x(x+1)^3}=\\= \frac{A(x+1)^3+Bx(x+1)^2+Cx(x+1)+Dx}{x(x+1)^3}
 \end{multline*}
Эти дроби действительно будут совпадать, если у них совпадают числители:
$$
3x+2=A(x+1)^3+Bx(x+1)^2+Cx(x+1)+Dx
$$

Теперь коэффициенты $A,B,C,D$ лучше всего искать, комбинируя метод частных
значений с методом неопределенных коэффициентов.
 \begin{align*}
x=0: & \quad 2=A \quad\Rightarrow\quad A=2  \\
x=-1: & \quad -1=-D  \quad\Rightarrow\quad D=1\\
x^3: & \quad 0=A+B=2+B \quad\Rightarrow\quad B=-2 \\
x^2: & \quad 0=3A+2B+C=2+C \quad\Rightarrow \\
 & \kern100pt \Rightarrow\quad C=-2
 \end{align*}
Таким образом,
$$
\begin{cases}{A=2}\\{B=-2}\\{C=-2}\\{D=1}\end{cases}
$$
Значит,
$$
\frac{3x+2}{x(x+1)^3}=
\frac{2}{x}-\frac{2}{x+1}-\frac{2}{(x+1)^2}+\frac{1}{(x+1)^3}
$$
и теперь можно вычислять наш интеграл:
 \begin{multline*}
\int \frac{3x+2}{x(x+1)^3}\, \d x = \int \frac{2}{x}\, \d x - \int
\frac{2}{x+1}\, \d x -\\- \int \frac{2}{(x+1)^2}\, \d x + \int
\frac{1}{(x+1)^3}\, \d x=\\= 2\ln |x| - 2\ln |x+1| + \frac{2}{x+1} -
\frac{1}{2(x+1)^2}+E \end{multline*} ($E$ -- постоянная интегрирования).
\end{ex}

\begin{ex}[\bf случай неповторяющихся множителей второй степени]\label{ex-13.2.9}
Рассмотрим интеграл
$$
\int \frac{x^2-2x+2}{(x-1)^2 (x^2+1)}\, \d x
$$
Знаменатель дроби уже разложен на множители:
$$
(x-1)^2 (x^2+1)
$$
В соответствии с первым и вторым правилами, в разложении нашей (правильной)
дроби в сумму простейших дробей должны присутствовать следующие слагаемые:
$$
\frac{A}{x-1}, \quad \frac{B}{(x-1)^2}, \quad \frac{Cx+D}{x^2+1}
$$
Значит,
$$
\frac{x^2-2x+2}{(x-1)^2 (x^2+1)}=
\frac{A}{x-1}+\frac{B}{(x-1)^2}+\frac{Cx+D}{x^2+1}
$$
Чтобы найти коэффициенты $A,B,C,D$, приведем правую часть к общему знаменателю
 \begin{multline*}
\frac{x^2-2x+2}{(x-1)^2 (x^2+1)}=\\={\smsize\text{$
\frac{A(x-1)(x^2+1)+B(x^2+1)+(Cx+D)(x-1)^2}{(x-1)^2 (x^2+1)}$}}
 \end{multline*}
Эти дроби совпадают, если у них совпадают числители:
 \begin{multline*}
x^2-2x+2=\\= A(x-1)(x^2+1)+B(x^2+1)+(Cx+D)(x-1)^2
 \end{multline*}

Теперь коэффициенты $A,B,C,D$ ищем, комбинируя метод частных значений с методом
неопределенных коэффициентов.
 \begin{align*}
x=1: & \quad 1=2B  \\
x^3: & \quad 1=A+C \\
x^2: & \quad 0=-A+B-2C+D  \\
x^0: & \quad 2=-A+B+D
 \end{align*}
Получается система
$$
\begin{cases}{2B=1}\\{A+C=1}\\{-A+B-2C+D=0}\\{-A+B+D=2}\end{cases}
\;\Leftrightarrow\; ... \; \Leftrightarrow \;
\begin{cases}{A=0}\\{B=\frac{1}{2}}\\{C=1}\\{D=\frac{3}{2}}\end{cases}
$$
Значит,
$$
\frac{x^2-2x+2}{(x-1)^2 (x^2+1)}= \frac{1}{2(x-1)^2}+\frac{2x+3}{2(x^2+1)}
$$
и теперь можно вычислять наш интеграл:
\begin{multline*}\int \frac{x^2-2x+2}{(x-1)^2 (x^2+1)}\, \d x=\\= \int
\frac{1}{2(x-1)^2}\, \d x + \int \frac{2x+3}{2(x^2+1)}\, \d x=\\=
-\frac{1}{2(x-1)} + \int \frac{2x}{2(x^2+1)}\, \d x+ \int \frac{3}{2(x^2+1)}\,
\d x=\\= -\frac{1}{2(x-1)} + \int \frac{d (x^2+1)}{2(x^2+1)}+ \frac{3}{2}\int
\frac{d x}{x^2+1}=\\= -\frac{1}{2(x-1)} + \frac{1}{2}\ln (x^2+1)+
\frac{3}{2}\arctg x+E
\end{multline*} ($E$ -- постоянная интегрирования).
\end{ex}

\begin{ex}[\bf еще один случай не\-пов\-то\-ряю\-щихся множителей второй
степени]\label{ex-13.2.10} Рассмотрим интеграл
$$
\int \frac{2x^2-3x+1}{x^3+1}\, \d x
$$
Разложим знаменатель на множители:
$$
x^3+1=(x+1)(x^2-x+1)
$$
В соответствии с первым и вторым правилами, в разложении нашей (правильной)
дроби в сумму простейших дробей должны присутствовать следующие слагаемые:
$$
\frac{A}{x+1}, \quad \frac{Bx+C}{x^2-x+1}
$$
Значит,
$$
\frac{2x^2-3x+1}{x^3+1}=\frac{A}{x+1}+\frac{Bx+C}{x^2-x+1}
$$
Чтобы найти коэффициенты $A,B,C$, приведем правую часть к общему знаменателю
$$
\frac{2x^2-3x+1}{x^3+1}=\frac{A(x^2-x+1)+(Bx+C)(x+1)}{(x+1)(x^2-x+1)}
$$
Эти дроби совпадают, если у них совпадают числители:
$$
2x^2-3x+1=A(x^2-x+1)+(Bx+C)(x+1)
$$

Теперь коэффициенты $A,B,C$ ищем, комбинируя метод частных значений с методом
неопределенных коэффициентов.
 \begin{align*}
x=-1: & \quad 6=3A  \\
x^2: & \quad 2=A+B \\
x^0: & \quad 1=A+C
 \end{align*}
Получается система
$$
\begin{cases}{3A=6}\\{A+B=2}\\{A+C=1}\end{cases}\quad \Leftrightarrow \quad ... \quad \Leftrightarrow \quad
\begin{cases}{A=2}\\{B=0}\\{C=-1}\end{cases}
$$
Значит,
$$
\frac{2x^2-3x+1}{x^3+1}=\frac{2}{x+1}-\frac{1}{x^2-x+1}
$$
и теперь можно вычислять наш интеграл:
 \begin{multline*}
\int\frac{2x^2-3x+1}{x^3+1}\, \d x =\\= \int \frac{2}{x+1}\, \d x - \int
\frac{1}{x^2-x+1}\, \d x=\\= 2\ln |x+1| - \int
\frac{1}{\left(x+\frac{1}{2}\right)^2+ \left(\frac{\sqrt{3}}{2}\right)^2}\, \d
x=\\= \left|
\begin{array}{c} y=x+\frac{1}{2}\\ x=y-\frac{1}{2}\end{array}\right|=
2\ln |x+1| -\\
-\int\frac{1}{y^2+\left(\frac{\sqrt{3}}{2}\right)^2}\d \left(
y-\frac{1}{2}\right) \, \Big|_{y=x+\frac{1}{2}}=\\= 2\ln |x+1| - \int
\frac{1}{y^2+\left(\frac{\sqrt{3}}{2}\right)^2}\, \d y \,
\Big|_{y=x+\frac{1}{2}}=\\
={\smsize\begin{pmatrix}\text{применяем}\\
\text{формулу \eqref{12.2.9}}\end{pmatrix}}=\\
=2\ln |x+1| - \frac{2}{\sqrt{3}}\arctg \frac{2y}{\sqrt{3}}\,
\Big|_{y=x+\frac{1}{2}}=\\= 2\ln |x+1| - \frac{2}{\sqrt{3}}\arctg
\frac{2 \left(x+\frac{1}{2}\right)}{\sqrt{3}}+D=\\
=2\ln |x+1| - \frac{2}{\sqrt{3}}\arctg \frac{2 x+1}{\sqrt{3}}+D
\end{multline*} ($D$ -- постоянная интегрирования).
\end{ex}

\begin{ex}[\bf случай повторяющихся множителей второй степени]\label{ex-13.2.11}
Рассмотрим интеграл
$$
\int \frac{x^3+3}{(x+1) (x^2+1)^2}\, \d x
$$
Знаменатель дроби уже разложен на множители:
$$
(x+1) (x^2+1)^2
$$
В соответствии с первым и вторым правилами, в разложении нашей (правильной)
дроби в сумму простейших дробей должны присутствовать следующие слагаемые:
$$
\frac{A}{x+1}, \quad \frac{Bx+C}{x^2+1}, \quad \frac{Dx+E}{(x^2+1)^2}
$$
Значит,
$$
\frac{x^3+3}{(x+1) (x^2+1)^2}= \frac{A}{x+1} + \frac{Bx+C}{x^2+1} +
\frac{Dx+E}{(x^2+1)^2}
$$
Чтобы найти коэффициенты $A,B,C,D,E$, приведем правую часть к общему
знаменателю
 \begin{multline*}
\frac{x^3+3}{(x+1) (x^2+1)^2}=\\={\smsize\text{
$\frac{A(x^2+1)^2+(Bx+C)(x+1)(x^2+1)+(Cx+D)(x+1)}{(x+1)(x^2+1)^2}$}}
 \end{multline*}
Эти дроби совпадают, если у них совпадают числители:
 \begin{multline*}
x^3+3=A(x^2+1)^2+\\+(Bx+C)(x+1)(x^2+1)+(Cx+D)(x+1)
 \end{multline*}

Теперь коэффициенты $A,B,C,D,E$ ищем, комбинируя метод частных значений с
методом неопределенных коэффициентов.
 \begin{align*}
x=-1: & \quad 2=4A  \\
x^4: & \quad 0=A+B \\
x^3: & \quad 1=B+C  \\
x^2: & \quad 0=2A+B+C+D \\
x^0: & \quad 3=A+C+E
 \end{align*}
Получается система
$$
\begin{cases}{4A=2}\\{A+B=0}\\{B+C=1}\\{2A+B+C+D=0}\\{A+C+E=3}\end{cases}
\;\Leftrightarrow\; ... \;\Leftrightarrow\;
\begin{cases}{A=\frac{1}{2}}\\{B=-\frac{1}{2}}\\{C=\frac{3}{2}}\\{D=-2}\\{E=1}\end{cases}
$$
Значит,
$$
\frac{x^3+3}{(x+1) (x^2+1)^2}= \frac{\frac{1}{2}}{x+1} +
\frac{-\frac{1}{2}x+\frac{3}{2}}{x^2+1} + \frac{-2x+1}{(x^2+1)^2}
$$
и теперь можно вычислять наш интеграл:
\begin{multline*}\int \frac{x^3+3}{(x+1) (x^2+1)^2}\, \d x= \int
\frac{\frac{1}{2}}{x+1}\, \d x +\\+ \int
\frac{-\frac{1}{2}x+\frac{3}{2}}{x^2+1}\, \d x + \int \frac{-2x+1}{(x^2+1)^2}\,
\d x=\\= \frac{1}{2}\int \frac{d x}{x+1} -\frac{1}{2}\int \frac{x \, \d
x}{x^2+1}  + \frac{3}{2}\int \frac{d x}{x^2+1} -\\-\int \frac{2x\, \d
x}{(x^2+1)^2} + \int \frac{d x}{(x^2+1)^2}=\\= \frac{1}{2}\ln |x+1|
-\frac{1}{4}\int \frac{d (x^2+1)}{x^2+1}+ \frac{3}{2}\arctg x -\\-\int \frac{d
(x^2+1)}{(x^2+1)^2} + \int \frac{d x}{(x^2+1)^2}=
{\smsize\begin{pmatrix}\text{в последнем}\\
\text{интеграле}\\
\text{применяем}\\
\text{формулу \eqref{13.1.4}}\end{pmatrix}}=\\= \frac{1}{2}\ln |x+1|
-\frac{1}{4}\ln (x^2+1)+ \frac{3}{2}\arctg x + \frac{1}{x^2+1} +\\+ \frac{x}{2
(x^2+1)}+\frac{1}{2}\cdot \arctg x+H=\\= \frac{1}{2}\ln |x+1| -\frac{1}{4}\ln
(x^2+1)+ 2 \arctg x +\\+ \frac{1}{x^2+1} + \frac{x}{2 (x^2+1)}+H
\end{multline*} ($H$ -- постоянная интегрирования).
\end{ex}

\begin{ers} Вычислите интегралы:
 \begin{multicols}{2}
 \biter{
\item[1)] $\int \frac{x^2+x-1}{x^3+x^2-6x}\, \d x$;

\item[2)] $\int \frac{d x}{(x+1)(x-2)}$;

\item[3)] $\int \frac{x d x}{3x^2-3x-2}$;

\item[4)] $\int \frac{2x+11}{(x^2+6x+13)^2}\, \d x$;

\item[5)] $\int \frac{x^2 \, \d x}{(x+2)^2(x+4)^2}$;

\item[6)] $\int \frac{d x}{x^3(x-1)}$;

\item[7)] $\int \frac{d x}{(x^2+2)(x-1)^2}$;

\item[8)] $\int \frac{x^3+3x^2-3x+1}{(x+1)^2(x^2+1)}\, \d x$;

\item[9)] $\int \frac{3x^3-x^2-4x+13}{x^2(x^2-4x+13)}\, \d x$;

\item[10)] $\int \frac{4x^2-8x}{(x-1)^2(x^2+1)^2}\, \d x$;

\item[11)] $\int \frac{x^2}{(x^2+2x+2)^2}\, \d x$;

\item[12)] $\int \frac{x^3+2x^2+3x+4}{x^4+x^3+2x^2}\, \d x$.
 }\eiter
 \end{multicols}
 \end{ers}

\paragraph{Интегрирование неправильных дробей.}

Напомним (см. предыдущий подпункт (b)) , что рациональная функция
$$
R(x)=\frac{P(x)}{Q(x)}
$$
называется {\it неправильной дробью}\index{дробь!неправильная}, если степень
многочлена в числителе больше или равна степени многочлена в знаменателе:
$$
  \deg P\ge \deg Q
$$
Оказывается, что справедлива

\begin{tm}\label{tm-13.2.13} Всякая неправильная дробь представима в виде
суммы многочлена и правильной дроби.
\end{tm}

Она позволяет вычислять интегралы от любых рациональных функций (потому что
любую рациональную функцию можно представить в виде суммы многочлена и
правильной дроби, а их мы уже умеем интегрировать).

\begin{ex}\label{ex-13.2.14}
Рассмотрим интеграл
$$
\int \frac{5x^3+9x^2-22x-8}{x^3-4x}\, \d x
$$
Разделив числитель на знаменатель с остатком, получим
$$
\frac{5x^3+9x^2-22x-8}{x^3-4x}= 5+\frac{9x^2-2x-8}{x^3-4x}
$$
Поэтому
 \begin{multline*}
\int \frac{5x^3+9x^2-22x-8}{x^3-4x}\, \d x =\\= \int 5 \, \d x+\int
\frac{9x^2-2x-8}{x^3-4x}\, \d x=\\= {\smsize\begin{pmatrix}\text{второй интеграл}\\
\text{мы уже нашли}\\
\text{в примере \ref{ex-13.2.7}}\end{pmatrix}} =\\= 5x+ \left ( 2\ln |x| + 3\ln
|x-2| + 4\ln |x+2|+D \right)
\end{multline*}\end{ex}

\begin{ex}\label{ex-13.2.15}
Рассмотрим интеграл
$$
\int \frac{2x^4+5x^2-2}{2x^3-x-1}\, \d x
$$
Разделив числитель на знаменатель с остатком, получим
$$
\frac{2x^4+5x^2-2}{2x^3-x-1} = x+ \frac{6x^2+x-2}{2x^3-x-1}
$$
Теперь надо подумать о том, как вычислить интеграл от второго слагаемого (то
есть от полученной правильной дроби). Для этого разложим знаменатель на
множители:
$$
2x^3-x-1=(x-1)(2x^2+2x+1)
$$
В соответствии с первым и вторым правилами, в разложении нашей (правильной)
дроби в сумму простейших дробей должны присутствовать следующие слагаемые:
$$
\frac{A}{x-1}, \quad \frac{Bx+C}{2x^2+2x+1}
$$
Значит,
$$
\frac{6x^2+x-2}{2x^3-x-1}=\frac{A}{x-1}+\frac{Bx+C}{2x^2+2x+1}
$$
Чтобы найти коэффициенты $A,B,C$, приведем правую часть к общему знаменателю
$$
\frac{6x^2+x-2}{2x^3-x-1}=\frac{A(2x^2+2x+1)+(Bx+C)(x-1)}{(x-1)(2x^2+2x+1)}
$$
Эти дроби совпадают, если у них совпадают числители:
$$
6x^2+x-2=A(2x^2+2x+1)+(Bx+C)(x-1)
$$

Теперь коэффициенты $A,B,C$ ищем, комбинируя метод частных значений с методом
неопределенных коэффициентов.
 \begin{align*}
x=1: & \quad 5=5A \\
x^2: & \quad 6=2A+B  \\
x^0: & \quad -2=A-C
 \end{align*} Получается система
$$
\begin{cases}{A=1}\\{2A+B=6}\\{A-C=-2}\end{cases}\quad \Leftrightarrow \quad ... \quad \Leftrightarrow \quad
\begin{cases}{A=1}\\{B=4}\\{C=3}\end{cases}
$$
Значит,
$$
\frac{6x^2+x-2}{2x^3-x-1}=\frac{1}{x-1}+\frac{4x+3}{2x^2+2x+1}
$$
Отсюда
 \begin{multline*}
\frac{2x^4+5x^2-2}{2x^3-x-1}=\\
=x+\frac{6x^2+x-2}{2x^3-x-1}=x+\frac{1}{x-1}+\frac{4x+3}{2x^2+2x+1}
 \end{multline*}
и теперь можно вычислять исходный интеграл:
 \begin{multline*}
\int \frac{2x^4+5x^2-2}{2x^3-x-1}\, \d x = \int x \, \d x+\\+ \int
\frac{1}{x-1}\, \d x+ \int \frac{4x+3}{2x^2+2x+1}\, \d x=\\= \frac{x^2}{2}+\ln
|x-1|+ \int \frac{2x+\frac{3}{2}}{x^2+x+\frac{1}{2}}\, \d x=\\=
\frac{x^2}{2}+\ln |x-1|+\\+ \int \frac{2x+\frac{3}{2}}{x^2+2\cdot
\frac{1}{2}x+\left(\frac{1}{2}\right)^2
-\left(\frac{1}{2}\right)^2+\frac{1}{2}}\, \d x=\\= \frac{x^2}{2}+\ln |x-1|+
\int \frac{2x+\frac{3}{2}}{\left(x+\frac{1}{2}\right)^2 +\frac{1}{4}}\, \d
x=\\= \left|
\begin{array}{c}
y=x+\frac{1}{2}\\ x=y-\frac{1}{2}\end{array}\right|= \frac{x^2}{2}+\ln
|x-1|+\\+ \int \frac{2\left(y-\frac{1}{2}\right)+\frac{3}{2}}
{y^2+\frac{1}{4}}\, \d \left(y-\frac{1}{2}\right) \, \Big|_{y=x+\frac{1}{2}}=
\frac{x^2}{2}+\\+\ln |x-1|+ \int \frac{2y+\frac{1}{2}} {y^2+\frac{1}{4}}\, \d y
\, \Big|_{y=x+\frac{1}{2}}=\\= \frac{x^2}{2}+\ln |x-1|+ \int \frac{2y\, \d y}
{y^2+\frac{1}{4}}\, \Big|_{y=x+\frac{1}{2}}+\\+ \frac{1}{2}\int \frac{d y}
{y^2+\frac{1}{4}}\, \Big|_{y=x+\frac{1}{2}}= \frac{x^2}{2}+\ln |x-1|+\\+ \int
\frac{d \left(y^2+\frac{1}{4}\right)} {y^2+\frac{1}{4}}\,
\Big|_{y=x+\frac{1}{2}}+ \frac{1}{2}\int \frac{d y}
{y^2+\left(\frac{1}{2}\right)^2}\, \Big|_{y=x+\frac{1}{2}}=\\=
{\smsize\begin{pmatrix}\text{применяем}\\
\text{формулы}\\ \text{\eqref{12.2.2} и \eqref{12.2.9}}\end{pmatrix}}=
\frac{x^2}{2}+\ln |x-1|+\\+ \ln \left(y^2+\frac{1}{4}\right) \,
\Big|_{y=x+\frac{1}{2}}+ \arctg (2y) \, \Big|_{y=x+\frac{1}{2}}+D=\\=
\frac{x^2}{2}+\ln |x-1|+ \ln
\left(\left(x+\frac{1}{2}\right)^2+\frac{1}{4}\right) +\\+\arctg
2\left(x+\frac{1}{2}\right)+D=\\= \frac{x^2}{2}+\ln |x-1|+ \ln
\left(x^2+x+\frac{1}{2}\right) +\arctg (2x+1)+D
\end{multline*} ($D$ -- постоянная интегрирования).
\end{ex}

\begin{ers} Вычислите интегралы:
 \begin{multicols}{2}
1) $\int \frac{x^5+x^4-8}{x^3-4x}\, \d x$;

2) $\int \frac{3x^3+5x^2-25x-1}{(x+2)(x-1)^2}$;

3) $\int \frac{x^5-2x^2+3}{x^2-4x+4}$;

4) $\int \frac{x^3+x^2+x+3}{(x+3)(x^2+x+1)}$;

5) $\int \frac{x^9}{(x^4-1)^2}$; \end{multicols}\end{ers}

\subsection{Интегрирование рациональных тригонометрических функций}

{\it Рациональной функцией от двух переменных}\index{функция!рациональная!от
двух переменных} $u$ и $v$ называется любая функция $R(u,v)$, составленная из
$u$ и $v$ с помощью алгебраических операций $+,-,\cdot,/$. Например, функция
$$
R(u,v)=\frac{3u^2v-2uv^2+7v-1}{v^2+\sqrt{5}}
$$
является рациональной.

{\it Рациональной тригонометрической
функцией}\index{функция!рациональная!тригонометрическая} называется любая
рациональная функция от переменных $\sin x$ и $\cos x$, например,
\begin{multline*}
f(x)=R(\sin x, \cos x)=\\=\frac{3\sin^2 x \cos x-2\sin x \cos^2 x+7\cos x-1}
{\cos^2 x+\sqrt{5}}
\end{multline*}
В этом параграфе мы рассмотрим способы интегрирования рациональных
тригонометрических функций.

\paragraph{Универсальная тригонометрическая подстановка}

Существует один универсальный способ превратить интеграл от рациональной
тригонометрической функции в интеграл от обычной рациональной функции (который
мы уже умеем вычислять). Он называется {\it универсальной тригонометрической
подстановкой}\index{подстановка!тригонометрическая!универсальная} и состоит в
следующем. Если необходимо вычислить интеграл
$$
  \int R(\sin x, \cos x) \, \d x
$$
то сделать это можно заменой
$$
x=2 \arctg t \qquad \left( t=\tg \frac{x}{2}\right) \label{13.3.1}
$$
тогда получится
$$
\sin x=\frac{2t}{1+t^2}, \quad \cos x=\frac{1-t^2}{1+t^2}, \quad
dx=\frac{2dt}{1+t^2}\label{13.3.2}
$$
Рассмотрим примеры.

\begin{ex}\label{ex-13.3.1}
 \begin{multline*}\int \frac{d x}{\sin x}= \left|
\begin{matrix}
t=\tg \frac{x}{2}\\
\sin x=\frac{2t}{1+t^2}\\
dx=\frac{2dt}{1+t^2}\end{matrix}\right|= \int
\frac{\frac{2dt}{1+t^2}}{\frac{2t}{1+t^2}}=\\= \int \frac{dt}{t}=\ln
|t|+C= {\smsize\begin{pmatrix}\text{возвращаемся}\\
\text{к переменной $x$}\end{pmatrix}}=\\=\ln \left|\tg \frac{x}{2}\right| +C
\end{multline*}\end{ex}

\begin{ex}\label{ex-13.3.2}
 \begin{multline*}\int \frac{d x}{\cos x}= \left|
\begin{array}{c}
t=\tg \frac{x}{2}\\
\cos x=\frac{1-t^2}{1+t^2}\\
dx=\frac{2dt}{1+t^2}\end{array}\right|= \int
\frac{\frac{2dt}{1+t^2}}{\frac{1-t^2}{1+t^2}}=\\= \int \frac{2dt}{1-t^2}= -\int
\frac{d(1-t^2)}{1-t^2}= -\ln |1-t^2|+C=\\=
{\smsize\begin{pmatrix}\text{возвращаемся}\\
\text{к переменной $x$}\end{pmatrix}}=-\ln \left|1-\tg^2 \frac{x}{2}\right| +C
\end{multline*}\end{ex}

\begin{ex}\label{ex-13.3.3}
 \begin{multline*}\int \frac{\d x}{8-4\sin x+7\cos x}= \left|
\begin{array}{c}
t=\tg \frac{x}{2}\\
\sin x=\frac{2t}{1+t^2}\\
\cos x=\frac{1-t^2}{1+t^2}\\
\d x=\frac{2\d t}{1+t^2}\end{array}\right|=\\= \int \frac{\frac{2\d
t}{1+t^2}}{8-4\frac{2t}{1+t^2}+7\frac{1-t^2}{1+t^2}}=\\= \int \frac{2\d
t}{8(1+t^2)-8t+7(1-t^2)}=\\=
 \int\frac{2\d t}{t^2-8t+15}=
 \int \frac{2}{(t-3)(t-5)}\, \d t=\\=
 {\smsize\begin{pmatrix}
 \text{здесь мы ''в уме''}
 \\
 \text{раскладываем дробь}
 \\
 \text{под интегралом}
 \\
 \text{на простейшие}
 \end{pmatrix}}=\\=
 \int \left(\frac{1}{t-5}-\frac{1}{t-3}\right)\, \d t=\\=
 \ln |t-5|-\ln |t-3|+C=
 \left(
 \begin{array}{c}
 \text{возвращаемся}
 \\
 \text{к переменной $x$}
 \end{array}
 \right)=\\=
 \ln \left|\tg \frac{x}{2}-5\right|-\ln \left|\tg
 \frac{x}{2}-3\right|+C
\end{multline*}\end{ex}

\begin{ers} Найдите интегралы
 \begin{multicols}{2}
1) $\int \frac{d x}{5+4\sin x}$

2) $\int \frac{d x}{3\sin x-4\cos x}$

3) $\int \frac{d x}{5+\sin x+3\cos x}$

4) $\int \frac{\sin x \cdot \cos x}{(3+\cos x)^2}\, \d x$

5) $\int \frac{\sin x}{1+\tg x}\, \d x$ \end{multicols}\end{ers}

\paragraph{Частные тригонометрические подстановки.}

Часто бывает, что универсальная тригонометрическая подстановка приводит к
сложным выкладкам, которых можно избежать, если пользоваться другими
подстановками. Такие подстановки называются {\it
частными}\index{подстановка!тригонометрическая!частная}, и их имеется три:

 \biter{
\item[---] если функция $R(\sin x, \cos x)$ нечетная относительно синуса, то
есть
$$
R(-\sin x, \cos x)=-R(\sin x, \cos x)
$$
то применима подстановка
 \beq\label{13.3.3}
 \begin{cases}
\cos x=t \\
\sin x=\sqrt{1-t^2}, \\
x=\arccos t, \\
d x=-\frac{\d t}{\sqrt{1-t^2}}
\end{cases}
 \eeq
\item[---] если функция $R(\sin x, \cos x)$ нечетная относительно косинуса, то
есть
$$
R(\sin x, -\cos x)=-R(\sin x, \cos x)
$$
то применима подстановка
 \beq
  \begin{cases}
\sin x=t \\
\cos x=\sqrt{1-t^2}, \\
x=\arcsin t, \\
d x=\frac{\d t}{\sqrt{1-t^2}}
  \end{cases}
 \label{13.3.4}
 \eeq
\item[---] если функция $R(\sin x, \cos x)$ четная относительно синуса и
косинуса, то есть
$$
R(-\sin x, -\cos x)=R(\sin x, \cos x)
$$
то применима подстановка
 \beq
  \begin{cases}
\tg x=t, \\
\sin x=\frac{t}{\sqrt{1+t^2}}, \\
\cos x=\frac{1}{\sqrt{1+t^2}}, \\
x=\arctg t, \\
d x=\frac{d t}{1+t^2}
  \end{cases}
\label{13.3.5}
 \eeq
 }\eiter

\begin{ex}\label{ex-13.3.5}
 \begin{multline*}\int \frac{\sin^3 x}{\cos x-3}\, \d x=
{\smsize\begin{pmatrix}\text{подынтегральная}\\
\text{функция нечетна}\\
\text{относительно $\sin x$}\\
\text{поэтому используем}\\
\text{подстановку}\, \cos x=t
\end{pmatrix}}=\\
=-\int \frac{\left(\sqrt{1-t^2}\right)^3}{t-3}\, \frac{d t}{\sqrt{1-t^2}}= \int
\frac{t^2-1}{t-3}\, \d t=\\=
{\smsize\begin{pmatrix}\text{делим многочлены}\\
\text{с остатком}\end{pmatrix}}= \int \left(t+3+\frac{8}{t-3}\right) \, \d
t=\\= \frac{t^2}{2}+3t+8\ln |t-3|+C
= {\smsize\begin{pmatrix}\text{возвращаемся}\\
\text{к переменной $x$}\end{pmatrix}}=\\= \frac{\cos^2 x}{2}+3\cos x+8\ln |\cos
x-3|+C \end{multline*}\end{ex}

\begin{ex}\label{ex-13.3.6}
 \begin{multline*}\int \frac{\cos^3 x}{\sin^4 x}\, \d x=
{\smsize\begin{pmatrix}\text{подынтегральная}\\
\text{функция нечетна}\\
\text{относительно $\cos x$}\\
\text{поэтому используем}\\
\text{подстановку}\, \sin x=t
\end{pmatrix}}=\\
=\int \frac{\left(\sqrt{1-t^2}\right)^3}{t^4}\, \frac{d t}{\sqrt{1-t^2}}= \int
\frac{1-t^2}{t^4}\, \d t=\\= \int \left(\frac{1}{t^4}-\frac{1}{t^2}\right) \,
\d t=
-\frac{1}{3t^3}+\frac{1}{t}+C=\\
={\smsize\begin{pmatrix}\text{возвращаемся}\\
\text{к переменной $x$}\end{pmatrix}}= -\frac{1}{3\sin^3 x}+\frac{1}{\sin x}+C
\end{multline*}\end{ex}

\begin{ex}\label{ex-13.3.7}
 \begin{multline*}\int \frac{d x}{\sin^2 x-4\sin x\cos x+5\cos^2
 x}=\\=
{\smsize\begin{pmatrix}\text{подынтегральная}\\
\text{функция четна}\\
\text{относительно $\sin x$ и $\cos x$}\\
\text{поэтому используем}\\
\text{подстановку}\, \tg x=t
\end{pmatrix}}=\\= \int \frac{\frac{dt}{1+t^2}}
{\left(\frac{t}{\sqrt{1+t^2}}\right)^2- 4\cdot \frac{t}{\sqrt{1+t^2}}\cdot
\frac{1}{\sqrt{1+t^2}} +5\left(\frac{1}{\sqrt{1+t^2}}\right)^2}=\\= \int
\frac{dt}{t^2-4t+5}= \int \frac{dt}{(t-2)^2+1}=\\= \left|
\begin{array}{c}
t-2=y
\\
t=y+2
\end{array}\right|= \int \frac{dy}{y^2+1}\, \Big|_{y=t-2}=\\
=\arctg y+C \, \Big|_{y=t-2}=\arctg (t-2)+C=\\=
{\smsize\begin{pmatrix}\text{возвращаемся}\\
\text{к переменной $x$}\end{pmatrix}}= \arctg (\tg x-2)+C
\end{multline*}\end{ex}

\begin{ers} Найдите интегралы

1) $\int \cos^4 x\sin^3 x \, \d x$

2) $\int \cos^5 x\sin^2 x \, \d x$

3) $\int \frac{\sin^3 x}{1+\cos^2 x}\, \d x$

4) $\int \frac{\cos^3 x}{4\sin^2 x-1}\, \d x$

5) $\int \frac{d x}{3\cos^2 x+4\sin^2 x}$

6) $\int \frac{d x}{19\sin^2 x-8\sin x \cos x-3}$
 \end{ers}

\subsection{Интегрирование тригонометрических произведений}

Для интегрирования некоторых специальных тригонометрических функций полезно
знать "дополнительные хитрости", которые мы обсудим в этом параграфе.

\paragraph{Интегрирование функций $\sin ax \cdot \cos bx,
\sin ax \cdot \sin bx, \cos ax \cdot \cos bx$}

Для интегрирования таких произведений применяются тригонометрические тождества
\eqref{sin(x)sin(y)}-\eqref{sin(x)cos(y)}. Рассмотрим примеры.

\begin{ex}\label{ex-13.4.1}
 \begin{multline*}\int \sin 4x \sin 6x \, \d x=
{\smsize \eqref{sin(x)sin(y)}}=\\= \int
\frac{1}{2}\left(\cos (4x-6x)-\cos (4x+6x) \right) \, \d x=\\= \frac{1}{2}\int
\cos 2x \, \d x - \frac{1}{2}\int \cos 10x \, \d x=\\= \frac{1}{4}\int \cos 2x
\, \d (2x) - \frac{1}{20}\int \cos 10x \, \d (10x)=\\= \frac{1}{4}  \sin 2x -
\frac{1}{20}\sin 10x +C
\end{multline*}\end{ex}

\begin{ers} Найдите интегралы

1) $\int \cos \frac{x}{3}\sin  \frac{x}{4}  \, \d x$

2) $\int \cos 4x \cos 7x \, \d x$

3) $\int \sin  \frac{x}{3}\sin  \frac{x}{2}  \, \d x$

4) $\int \sin 3x \cos 10x \, \d x$

5) $\int \sin 3x\sin^2 x \, \d x$
 \end{ers}

\paragraph{Интегрирование функций $\sin^\alpha x \cdot \cos^\beta x$}

Для интегрирования таких произведений нужно запомнить следующие правила:
 \biter{
\item[---] если $\alpha$ -- нечетное положительное целое число (то есть
$\alpha\in 2\mathbb{N}-1$), то делается подстановка
 \beq\label{13.4.4}
 \begin{cases}
\cos x=t, \\
\sin x=\sqrt{1-t^2}, \\
x=\arccos t, \\
\d x=-\frac{\d t}{\sqrt{1-t^2}}
 \end{cases}
 \eeq
\item[---] если $\beta$ -- нечетное положительное целое число (то есть  $\beta
\in 2\mathbb{N}-1$), то делается подстановка
 \beq\label{13.4.5}
  \begin{cases}
\sin x=t, \\
\cos x=\sqrt{1-t^2}, \\
x=\arcsin t, \\
\d x=\frac{\d t}{\sqrt{1-t^2}}
 \end{cases}
 \eeq
\item[---] если $\alpha+\beta$ -- четное отрицательное целое число
($\alpha+\beta \in -2\mathbb{N}$), то делается подстановка
 \beq\label{13.4.6}
   \begin{cases}
\tg x=t, \\
\sin x=\frac{t}{\sqrt{1+t^2}}, \\
\cos x=\frac{1}{\sqrt{1+t^2}}, \\
 x=\arctg t, \\
 d x=\frac{d t}{1+t^2}
  \end{cases}
 \eeq
\item[---] если $\alpha$ и $\beta$ -- четные неотрицательные целые числа
($\alpha,\beta \in 2(\mathbb{N}-1)$), то применяются формулы понижения степени:
 \beq
  \sin^2 x=\frac{1-\cos 2x}{2}, \qquad \cos^2 x=\frac{1+\cos 2x}{2}\label{13.4.7}
 \eeq
 }\eiter

\begin{ex}\label{ex-13.4.3}
 \vglue-20pt
 \begin{multline*}\int \sqrt[3] {\sin^2 x}\cos^3 x \, \d x=
 {\smsize \left|
\begin{array}{c}\text{поскольку}\\
\beta=3\in 2\mathbb{N}-1,
\\
\text{подстановка}\\
\sin x=t
\end{array}\right|}=\\= \int \sqrt[3]{t^2}\left(\sqrt{1-t^2}\right)^3 \, \frac{d
t}{\sqrt{1-t^2}}=\\= \int t^{\frac{2}{3}} (1-t^2) \, \d t= \int
\left(t^{\frac{2}{3}} -t^{\frac{8}{3}}\right) \, \d t=\\=
\frac{3}{5}t^{\frac{5}{3}} -\frac{3}{11}t^{\frac{11}{3}}+C=
{\smsize\begin{pmatrix}\text{возвращаемся}\\
\text{к переменной $x$}\end{pmatrix}}=\\= \frac{3}{5}(\sin x)^{\frac{5}{3}}
-\frac{3}{11}(\sin x)^{\frac{11}{3}}+C=\\= \frac{3}{5}\sqrt[3]{\sin^5 x}
-\frac{3}{11}\sqrt[3]{\sin^{11} x}+C
\end{multline*}\end{ex}

\begin{ex}\label{ex-13.4.4}
 \begin{multline*}\int \frac{dx}{\sqrt {\cos^7 x \sin x}} =
 {\smsize \left|
\begin{array}{c}\text{поскольку}\\
\alpha+\beta=-\frac{7}{2}-\frac{1}{2}=\\= -4\in -2\mathbb{N},
\\
\text{подстановка}\\
\tg x=t
\end{array}\right|}=\\= \int \frac{\frac{d t}{1+t^2}} {\sqrt
{\left(\frac{1}{\sqrt{1+t^2}}\right)^7 \frac{t}{\sqrt{1+t^2}}}} = \int
\frac{\frac{d t}{1+t^2}} {\left(\frac{1}{1+t^2}\right)^2 \sqrt{t}} =\\= \int
\frac{1+t^2}{\sqrt{t}}\, \d t= \int \left(
t^{-\frac{1}{2}}+t^{\frac{3}{2}}\right) \, \d t=\\=
2t^{\frac{1}{2}}+\frac{2}{5} t^{\frac{5}{2}}+C=
{\smsize\begin{pmatrix}\text{возвращаемся}\\
\text{к переменной $x$}\end{pmatrix}}=\\= 2\sqrt{\tg x}+\frac{2}{5}\sqrt{\tg^5
x}+C \end{multline*}\end{ex}

\begin{ex}\label{ex-13.4.5}
 \begin{multline*}
\int \sin^2 x \cos^4 x \, \d x =
{\smsize\left|\begin{matrix}\text{поскольку}\\
\alpha=2\in 2(\mathbb{N}-1),
\\
\beta=4\in 2(\mathbb{N}-1),
\\
\text{применяем формулы} \\ \eqref{13.4.7}\end{matrix}\right|}=\\
=\int \frac{1-\cos 2x}{2}\cdot  \left(\frac{1+\cos 2x}{2}\right)^2 \, \d x =\\=
\int \frac{1-\cos 2x}{2}\cdot \frac{1+2\cos 2x+\cos^2 2x}{4}  \, \d x =\\=
\frac{1}{8}\int (1+\cos 2x-\cos^2 2x-\cos^3 2x) \, \d x
=\\={\smsize\text{$\frac{1}{8}\int \left(1+\cos 2x-\frac{1+\cos 4x}{2}-
\frac{1+\cos 4x}{2}\cos 2x \right)  \, \d x$}}=\\ = \frac{1}{16}\int (1+\cos
2x-\cos
4x-\cos 4x\cos 2x)  \, \d x =\\= {\smsize\begin{pmatrix}\text{применяем}\\
\text{формулу \eqref{cos(x)cos(y)}}\end{pmatrix}}=\\
={\smsize\text{$\frac{1}{16}\int \left(1+\cos 2x-\cos 4x- \frac{1}{2}(\cos
2x+\cos 6x) \right)  \, \d x$}} =\\= \frac{1}{32}\int (2+\cos 2x-2\cos 4x-\cos
6x )  \, \d x =\\= \frac{1}{32}\left(2x+ \frac{1}{2}\sin 2x-\frac{1}{2}\sin
4x-\frac{1}{6}\sin 6x+C\right) =\\= \frac{x}{16}+ \frac{\sin 2x}{64}-\frac{\sin
4x}{64}-\frac{\sin 6x}{192}+C
\end{multline*}\end{ex}

\begin{ers} Найдите интегралы
 \begin{multicols}{2}
1) $\int \sin^5 x \cos^4 x  \, \d x$

2) $\int \sin^{\frac{3}{5}} x \cos^5 x \, \d x$

3) $\int  \frac{\sin^3 x}{\sqrt[3]{\cos^4 x}}  \, \d x$

4) $\int \frac{dx}{\sin^3 x \cos x}$

5) $\int \frac{dx}{\sqrt{\sin x \cos^3 x}}$

6) $\int \frac{dx}{\sqrt[4]{\sin^3 x \cos^5 x}}$

7) $\int \cos^2 x \, \d x$

8) $\int \cos^4 x \, \d x$

9) $\int \sin^4 x \, \d x$

10) $\int \sin^4 x \cos^2 x \, \d x$

11) $\int \cos^6 x \, \d x$ \end{multicols}\end{ers}

\subsection{Интегрирование некоторых алгебраических иррациональностей}

\paragraph{Интегрирование функций $R(x, \sqrt[n]{x}), \, n\in \mathbb{N}$}

Если $R(x, y)$ -- рациональная функция от переменных $x$ и $y$, то функция
$R(x, \sqrt[n]{x})$ интегрируется заменой
$$
t=\sqrt[n]{x}\qquad (x=t^n,\, \d x=nt^{n-1} dt )
$$

\begin{ex}\label{ex-13.5.1}\begin{multline*}\int \frac{\sqrt{x}}{x-\sqrt[3] {x^2}}\, \d x= \int
\frac{\left(\sqrt[6]{x}\right)^3}{x-\left(\sqrt[6]{x}\right)^4}\, \d
x=\\={\smsize \left|
\begin{array}{c}
t=\sqrt[6]{x}\\
x=t^6
\\
dx=6t^5 dt
\end{array}\right|}=
\int \frac{t^3}{t^6-t^4}\, 6t^5 dt=\\
=6\int \frac{t^4}{t^2-1}\, \d t={\smsize\begin{pmatrix}\text{делим}\\
\text{многочлены}\\
\text{с остатком}\end{pmatrix}}=\\
= 6\int \left(t^2+1+\frac{1}{t^2-1}\right)\, \d t=\\
=6\int \left(t^2+1+\frac{1}{(t-1)(t+1)}\right)\, \d t=\\=
{\smsize\begin{pmatrix}\text{раскладываем}\\
\text{на простейшие}\\
\text{дроби}\end{pmatrix}}=\\
=6\int \left(t^2+1+\frac{1}{2}\frac{1}{t-1}-\frac{1}{2}\frac{1}{t+1}\right)\,
\d t=\\
=2t^3+6t+3\ln |t-1|-3\ln |t+1| +C =\\=
{\smsize\begin{pmatrix}\text{возвращаемся}\\
\text{к переменной $x$}\end{pmatrix}}=
2\left(\sqrt[6]{x}\right)^3+6\sqrt[6]{x}+\\+ 3\ln
\left|\sqrt[6]{x}-1\right|-3\ln \left|\sqrt[6]{x}+1\right| +C
\end{multline*}\end{ex}

\begin{ers} Найдите интегралы
 \begin{multicols}{2}
1) $\int \frac{\sqrt{x}}{1+\sqrt{x}}  \, \d x$

2) $\int \frac{d x}{\left(1+\sqrt[3]{x}\right)\cdot \sqrt{x}}$

3) $\int \frac{\sqrt{x}}{1+\sqrt[4]{x}}  \, \d x$
\end{multicols}\end{ers}

\paragraph{Интегрирование функций
$R\left(x, \sqrt[n]{\frac{ax+b}{cx+d}}\right), \, n\in \mathbb{N}$.}

Интегралы такого вида вычисляются заменой
$$
t=\sqrt[n]{\frac{ax+b}{cx+d}}
$$

\begin{ex}\label{ex-13.5.3}
 \begin{multline*}
\int \frac{1}{(1-x)^2}\cdot \sqrt{\frac{1-x}{1+x}}\, \d x={\smsize
\left|\begin{array}{c} t=\sqrt{\frac{1-x}{1+x}}\\
x=\frac{1-t^2}{1+t^2}\\
dx=-\frac{4t dt}{(1+t^2)^2}\end{array}\right|}=\\
=-\int \frac{1}{\left(1-\frac{1-t^2}{1+t^2}\right)^2}\cdot t \, \frac{4t
dt}{(1+t^2)^2}=\\= -\int \frac{(1+t^2)^2}{4t^4}\, \frac{4t^2 dt}{(1+t^2)^2}=
-\int \frac{dt}{t^2}= \frac{1}{t}+C=\\=
{\smsize\begin{pmatrix}\text{возвращаемся}\\
\text{к переменной $x$}\end{pmatrix}}= \sqrt{\frac{1+x}{1-x}}+C
\end{multline*}\end{ex}

\begin{ers} Найдите интегралы

1) $\int \frac{\sqrt{1+x}}{x}  \, \d x$

2) $\int \frac{\sqrt[3] {3x+4}}{1+\sqrt[3] {3x+4}}  \, \d x$

3) $\int \frac{1}{(1-x)(1+x)^2}\cdot \sqrt{\frac{1+x}{1-x}}  \, \d x$

4) $\int  \sqrt[3]{\frac{1-x}{1+x}}\cdot   \frac{dx}{(1+x)^2}$
 \end{ers}

\paragraph{Интегрирование функций $\frac{1}{\sqrt{ax^2+bx+c}}$}

После выделения полного квадрата под знаком радикала интегралы этого вида
сводятся к следующим:
 \bit{
\item[---] если $a>0$, то получаем табличный интеграл
$$
\int \frac{1}{\sqrt{x^2+A}}\, \d x
$$
\item[---] если $a<0$, то получаем табличный интеграл
$$
\int \frac{1}{\sqrt{B^2-x^2}}\, \d x
$$
 }\eit

\begin{ex}\label{ex-13.5.5}
 \begin{multline*}
\int \frac{dx}{2x^2-x+3} = \frac{1}{\sqrt{2}}\int \frac{dx}
{\sqrt{x^2-\frac{1}{2}x+\frac{3}{2}}} =\\
=\frac{1}{\sqrt{2}}\int \frac{d x}
{\sqrt{\left(x-\frac{1}{4}\right)^2+\frac{23}{16}}}={\smsize \left|
\begin{matrix} t=x-\frac{1}{4}\\ x=t+\frac{1}{4}\\ \d x=\d t
\end{matrix}\right|}=\\= \frac{1}{\sqrt{2}}\int \frac{d
t}{\sqrt{t^2+\frac{23}{16}}}= \frac{1}{\sqrt{2}}\ln \left|
t+\sqrt{t^2+\frac{23}{16}}\right|+C=\\=
{\smsize\begin{pmatrix}\text{возвращаемся}\\
\text{к переменной $x$}\end{pmatrix}}=\\= \frac{1}{\sqrt{2}}\ln \left|
x-\frac{1}{4}+\sqrt{ \left( x-\frac{1}{4}\right)^2+\frac{23}{16}}\right|+C=\\=
\frac{1}{\sqrt{2}}\ln \left| x-\frac{1}{4}+
\sqrt{x^2-\frac{1}{2}x+\frac{3}{2}}\right|+C
\end{multline*}\end{ex}

\begin{ex}\label{ex-13.5.6}
 \begin{multline*}
\int \frac{ dx}{\sqrt{5-2x-3x^2}}= \frac{1}{\sqrt{3}}\int \frac{dx}
{\sqrt{\frac{5}{3}-\frac{2}{3}x-x^2}}=\\
=\frac{1}{\sqrt{3}}\int \frac{dx} {\sqrt{\frac{16}{9}-\left(\frac{1}{3}+x
\right)^2}}= {\smsize \left|
\begin{array}{c}
t=\frac{1}{3}+x
\\
x=t-\frac{1}{3}\\
dx=dt
\end{array}\right|}=\\= \frac{1}{\sqrt{3}}\int \frac{dt}
{\sqrt{\left(\frac{4}{3}\right)^2-t^2}}= \frac{1}{\sqrt{3}}\arcsin
\frac{3t}{4}+C= \\
={\smsize\begin{pmatrix}\text{возвращаемся}\\
\text{к переменной $x$}\end{pmatrix}}= \frac{1}{\sqrt{3}}\arcsin
\frac{1+3x}{4}+C \end{multline*}\end{ex}

\begin{ers} Найдите интегралы
 \begin{multicols}{2}
1) $\int \frac{dx}{\sqrt{5x^2+3x+2}}$

2) $\int \frac{dx}{\sqrt{2+3x-2x^2}}$ \end{multicols}\end{ers}

\paragraph{Интегрирование функций $R\left(x, \sqrt{ax^2+bx+c}\right)$}

Выделением полного квадрата и последующей заменой функция $R\left(x,
\sqrt{ax^2+bx+c}\right)$ приводится к виду
$$
R\left(x, \sqrt{a^2-x^2}\right),
$$
$$
R\left(x, \sqrt{x^2+a^2}\right),
$$
$$
R\left(x, \sqrt{x^2-a^2}\right)
$$
а затем делается тригонометрическая подстановка:
 \bit{
\item[---] для
$$
\int R\left(x, \sqrt{a^2-x^2}\right) \, \d x
$$
используется подстановка
$$
x=a\sin t \qquad \left(\text{или}\,\, x=a\cos t\right)
$$
\item[---] для
$$
\int R\left(x, \sqrt{x^2+a^2}\right)      \, \d x
$$
используется подстановка
$$
x=a\tg t  \qquad \left(\text{или}\,\, x=a\ctg t\right)
$$
\item[---] для
$$
\int R\left(x, \sqrt{x^2-a^2}\right)      \, \d x
$$
используется подстановка
$$
x=\frac{a}{\cos t}  \qquad \left(\text{или}\,\, x=\frac{a}{\sin t}\right)
$$
 }\eit

\begin{ex}\label{ex-13.5.8}
\begin{multline*}\int \frac{ dx}{(x^2+16)\sqrt{9-x^2}}={\smsize \left|
\begin{array}{c}
x=3\sin t
\\
t=\arcsin \frac{x}{3}\\
dx=3\cos t dt
\end{array}\right|}=\\
=\int \frac{3\cos t dt}{(9\sin^2 t+16)\sqrt{9-9\sin^2 t}}=\\
=\int \frac{3\cos t dt}{(9\sin^2 t+16) 3\cos t}=\\= \int \frac{dt}{9\sin^2
t+16}={\smsize \left|
\begin{array}{c}\text{применяем частную}\\
\text{подстановку \eqref{13.3.5}}\\
\tg t=s
\\
\sin t=\frac{s}{\sqrt{1+s^2}}\\
\d t=\frac{\d s}{\sqrt{1+s^2}}\end{array}\right|}=\\
=\int \frac{\frac{\d s}{\sqrt{1+s^2}}}
{9\left(\frac{s}{\sqrt{1+s^2}}\right)^2+16}= \int \frac{ds}{9s^2+16(1+s^2)}=\\=
\int \frac{ds}{25s^2+16}=\\= \frac{1}{25}\int
\frac{ds}{s^2+\left(\frac{4}{5}\right)^2}= \frac{1}{20}\arctg
\frac{5s}{4}+C=\\=
{\smsize\begin{pmatrix}\text{возвращаемся}\\
\text{к переменной $t$}\end{pmatrix}}= \frac{1}{20}\arctg \frac{5
\tg t}{4}+C=\\= {\smsize\begin{pmatrix}\text{возвращаемся}\\
\text{к переменной $x$}\end{pmatrix}}= \frac{1}{20}\arctg \frac{5 \tg \arcsin
\frac{x}{3}}{4}+C=\\= \frac{1}{20}\arctg \frac{5 \frac{\frac{x}{3}}
{\sqrt{1-\left(\frac{x}{3}\right)^2}}}{4}+C= \frac{1}{20}\arctg \frac{5 x}
{4\sqrt{9-x^2}}+C
\end{multline*}\end{ex}

\begin{ex}\label{ex-13.5.9}
\begin{multline*}\int \frac{x^2-x+1}{(x^2+1)\sqrt{x^2+1}}\, \d x=
{\smsize \left|\begin{array}{c} x=\ctg t
\\
t=\arcctg x
\\
dx=-\frac{dt}{\sin^2 t}\end{array}\right|}=\\
=-\int \frac{\ctg^2 t-\ctg t+1}{(\ctg^2 t+1)\sqrt{\ctg^2 t+1}}\cdot
\frac{dt}{\sin^2 t}=\\= -\int \frac{\ctg^2 t-\ctg t+1}{\frac{1}{\sin^2 x}\cdot
\frac{1}{\sin x}}\cdot \frac{dt}{\sin^2
t}=\\
=-\int (\ctg^2 t-\ctg t+1)\sin x \, \d
t=\\
={\smsize\begin{pmatrix}\text{заметим, что}\\
\ctg^2 t+1=\frac{1}{\sin^2 t}\end{pmatrix}}=\\= -\int
\left(\frac{1}{\sin^2 t}-\ctg t \right)\sin x \, \d t=\\
=-\int \left(\frac{1}{\sin t}-\cos t \right) \, \d t=\\
=-\int \frac{dt}{\sin t}+\int \cos t \, \d t=
{\smsize\begin{pmatrix}\text{первый
интеграл}\\
\text{мы уже вычисляли}\\
\text{в примере \ref{ex-13.3.1}}\end{pmatrix}}=\\
=-\ln \left|\tg \frac{t}{2}\right|+\sin t +C=
{\smsize\begin{pmatrix}\text{возвращаемся}\\
\text{к переменной $x$}\end{pmatrix}}=\\= -\ln \left|\tg
\frac{\arcctg x}{2}\right|+\sin \arcctg x +C=\\
={\smsize\begin{pmatrix}\text{применяем}\\
\text{тригонометрические}\\
\text{тождества}\end{pmatrix}}=\\
= \frac{1}{\sqrt{1+x^2}} +\ln \left| x+\sqrt{1+x^2}\right|+C
\end{multline*}\end{ex}

\begin{ex}\label{ex-13.5.10}
 \begin{multline*}\int \frac{dx}{(x^2+2)\sqrt{x^2-1}}=
{\smsize \left|\begin{array}{c}
x=\frac{1}{\cos t}\\
t=\arccos \frac{1}{x}\\
dx=\frac{\sin t}{\cos^2 t}\, \d t
\end{array}\right|}=\\
=\int \frac{\frac{\sin t}{\cos^2 t}\, \d t} {\left(\frac{1}{\cos^2
t}+2\right)\sqrt{\frac{1}{\cos^2 t}-1}}= \int \frac{\cos t \, \d t}{1+2\cos^2
t}=\\= \int \frac{\cos t \, \d t}{3-2\sin^2 t}=
{\smsize\left|\begin{array}{c}\sin t=s
\\
ds=\cos t dt\end{array}\right|}= \int \frac{ds}{3-2s^2}=\\= \frac{1}{2}\int
\frac{1}{\left(\sqrt{\frac{3}{2}}\right)^2-s^2}\, \d s=
{\smsize\begin{pmatrix}\text{применяем}\\
\text{формулу \eqref{12.2.11}}\end{pmatrix}}=\\= \frac{1}{2\sqrt{6}}\ln \left|
\frac{\sqrt{\frac{3}{2}}+s}{\sqrt{\frac{3}{2}}-s}\right|+C=
{\smsize\begin{pmatrix}\text{возвращаемся}\\
\text{к переменной $t$}\end{pmatrix}}=\\= \frac{1}{2\sqrt{6}}\ln \left|
\frac{\sqrt{\frac{3}{2}}+\sin t}{\sqrt{\frac{3}{2}}-\sin
t}\right|+C= {\smsize\begin{pmatrix}\text{возвращаемся}\\
\text{к переменной $x$}\end{pmatrix}}=\\= \frac{1}{2\sqrt{6}}\ln \left|
\frac{\sqrt{\frac{3}{2}}+\sin \arccos \frac{1}{x}} {\sqrt{\frac{3}{2}}-\sin
\arccos \frac{1}{x}}\right|+C=\\= \frac{1}{2\sqrt{6}}\ln \left|
\frac{\sqrt{\frac{3}{2}}+\sqrt{1-\frac{1}{x^2}}}
{\sqrt{\frac{3}{2}}-\sqrt{1-\frac{1}{x^2}}}\right|+C=\\= \frac{1}{2\sqrt{6}}\ln
\left| \frac{x\sqrt{3}+\sqrt{2}\sqrt{x^2-1}}
{x\sqrt{3}-\sqrt{2}\sqrt{x^2-1}}\right|+C
\end{multline*}\end{ex}

\begin{ers} Найдите интегралы
 \begin{multicols}{2}
 \biter{
\item[1)] $\int \frac{dx}{(9+x^2)^\frac{3}{2}}$

\item[2)] $\int \frac{dx}{\sqrt{(5-x^2)^3}}$

\item[3)] $\int \frac{\sqrt{x^2-4}}{x^3}\, \d x$

\item[4)] $\int \frac{x^2}{(x^2+5)^\frac{3}{2}}\, \d x$

\item[5)] $\int \frac{dx}{\sqrt{(1+x+x^2)^3}}$

\item[6)] $\int \frac{\sqrt{1+x^2}}{2+x^2}\, \d x$

\item[7)] $\int \frac{dx}{(x^2+4)\sqrt{4x^2+1}}\, \d x$
 }\eiter
\end{multicols}\end{ers}

\paragraph{Интегрирование функций $x^m \cdot (a+bx^n)^p, \, m,n,p\in
\mathbb{Q}$.}

Интегралы от таких функций выражаются через элементарные функции только если
одно из чисел
$$
  p, \, \frac{m+1}{n}, \, \frac{m+1}{n}+p
$$
является целым. В соответствии с этим следует употреблять три различных
подстановки:
 \bit{
\item[---] если $p$  целое отрицательное ($p\in -\mathbb{N}$) и $m=\frac{q}{s},
\, n=\frac{r}{s}$, то используется подстановка
$$
  x=t^s
$$
\item[---] если $\frac{m+1}{n}$  целое ($\frac{m+1}{n}\in \mathbb{Z}$), то
используется подстановка
$$
a+bx^n=t
$$
\item[---] если $\frac{m+1}{n}+p$  целое ($\frac{m+1}{n}+p\in \mathbb{Z}$), то
используется подстановка
$$
ax^{-n}+b=t
$$
 }\eit

\begin{ex}\label{ex-13.5.12}
\begin{multline*}\int x^3(1-x^2)^{-\frac{3}{2}}\, \d x=\\
={\smsize \left|\begin{array}{c}
\text{здесь}\,\, m=3, n=2, p=-\frac{3}{2}\\
\text{и}\,\, \frac{m+1}{n}=2\in Z\\
\text{поэтому применяем}\\
\text{подстановку}\\
t=1-x^2, \, xdx=-\frac{1}{2} dt
\end{array}\right|}=\\= -\frac{1}{2}\int (1-t)t^{-\frac{3}{2}}\, \d t=
-\frac{1}{2}\int (t^{-\frac{3}{2}}-t^{-\frac{1}{2}})\, \d t=\\=
t^{-\frac{1}{2}}+t^{\frac{1}{2}}+C=
{\smsize\begin{pmatrix}\text{возвращаемся}\\
\text{к переменной}\, \, x
\end{pmatrix}}=\\= \frac{1}{\sqrt{1-x^2}}+\sqrt{1-x^2}+C
\end{multline*}\end{ex}

\begin{ex}\label{ex-13.5.13}
\begin{multline*}\int x\sqrt{1+x^4}\, \d x
={\smsize \left|
\begin{array}{c}\text{здесь}\,\, m=1, n=4, p=\frac{1}{2}\\
\text{и}\,\, \frac{m+1}{n}+p=1\in Z
\\
\text{поэтому применяем}\\
\text{подстановку}\\
t=x^{-4}+1, \,
x=(t-1)^{-\frac{1}{4}}\\
dx=-\frac{1}{4}(t-1)^{-\frac{5}{4}} dt
\end{array}\right|}=\\
=-\frac{1}{4}\int (t-1)^{-\frac{1}{4}}\sqrt{1+(t-1)^{-\frac{1}{4}\cdot 4}}\,
(t-1)^{-\frac{5}{4}} dt=\\= -\frac{1}{4}\int
(t-1)^{-\frac{3}{2}}\sqrt{1+(t-1)^{-1}}\, \d t=\\
=-\frac{1}{4}\int (t-1)^{-\frac{3}{2}}\frac{\sqrt{t}}{\sqrt{t-1}}\,
\d t=\\
=-\frac{1}{4}\int (t-1)^{-\frac{3}{2}}\frac{\sqrt{t}}{\sqrt{t-1}}\, \d t=\\=
-\frac{1}{4}\int \frac{\sqrt{t}}{(t-1)^2}\, \d t= {\smsize \left|
\begin{array}{c}
y=\sqrt{t}\\
t=y^2 \\
dt=2y dy
\end{array}\right|}=\\
=-\frac{1}{4}\int \frac{y}{(y^2-1)^2}\, 2y dy= -\frac{1}{2}\int
\frac{y^2}{(y^2-1)^2}\, \d y=\\=
{\smsize\begin{pmatrix}\text{раскладываем дробь}\\
\text{на простейшие}\end{pmatrix}}=\\
{\smsize\text{$=-\frac{1}{2}\int \left(\frac{\frac{1}{4}}{y-1} +
\frac{\frac{1}{4}}{(y-1)^2}- \frac{\frac{1}{4}}{y+1}+
\frac{\frac{1}{4}}{(y+1)^2}\right)\, \d y=$}}\\
{\smsize\text{$=\frac{1}{8}\ln |y-1| -\frac{1}{8(y-1)}+
\frac{1}{8}\ln |y+1|+\frac{1}{8(y+1)}+C=$}}\\
=\frac{1}{8}\ln \left|\frac{y+1}{y-1}\right|+\frac{y}{4(y^2-1)}+C=\\=
{\smsize\begin{pmatrix}\text{возвращаемся}\\
\text{к переменной}\, \, t
\end{pmatrix}}=
\frac{1}{8}\ln \left|
\frac{\sqrt{t}+1}{\sqrt{t}-1}\right|+\frac{\sqrt{t}}{4(t-1)}+C=
\\={\smsize\begin{pmatrix}\text{возвращаемся}\\
\text{к переменной}\, \, x
\end{pmatrix}}=\\= \frac{1}{8}\ln \left|
\frac{\sqrt{x^{-4}+1}+1}{\sqrt{x^{-4}+1}-1}\right|+
\frac{\sqrt{x^{-4}+1}}{4x^{-4}}+C=\\
=\frac{1}{8}\ln \left|\frac{\sqrt{1+x^4}+x^2}{\sqrt{1+x^4}-x^2}\right|+
\frac{x^2}{4}\sqrt{1+x^4}+C
\end{multline*}\end{ex}

\begin{ers} Найдите интегралы
 \begin{multicols}{2}
1) $\int \frac{\sqrt{x}}{\left( 1+\sqrt[3]{x}\right)^2}\, \d x$

2) $\int x^5 (1+x^2)^\frac{2}{3}\, \d x$

3) $\int \frac{\sqrt{4+\sqrt[3]{x}}}{\sqrt[3]{x^2}}\, \d x$

4) $\int \sqrt[3]{x}\cdot \sqrt[3]{1+3\sqrt[3]{x^2}}\, \d x$

5) $\int \sqrt[4]{\left(1+x^\frac{1}{2}\right)^3}\, \d x$

6) $\int \frac{\sqrt{2-\sqrt[3]{x}}}{\sqrt[3]{x}}\, \d x$

7) $\int \frac{dx}{x^2\sqrt[3]{(1+x^3)^2}}\, \d x$
\end{multicols}\end{ers}

\end{multicols}\noindent\rule[10pt]{160mm}{0.1pt}

\chapter{ОПРЕДЕЛЕННЫЙ ИНТЕГРАЛ}\label{CH-definite-integral}

Определенный интеграл функции $f$ на отрезке $[a;b]$ есть величина,
геометрический смысл которой -- площадь криволинейной трапеции, ограничиваемой
графиком функции $f$ (при условии, что $f$ неотрицательна на $[a;b]$):

%\pucture{0pt}{0pt}{11-1.pcx}

\vglue100pt

Точный смысл этих слов вряд ли может быть сейчас понятен читателю, поскольку в
школьном курсе геометрии не объясняется, что такое площадь криволинейной
трапеции (мы же даем соответствующие строгие определения лишь в главе
\ref{CH-mera-Jordana}). Но всякий раз, когда вводится новое определение,
полезно, чтобы слушатель держал в голове какую-то геометрическую картинку,
упрощающую понимание соответствующих формальных определений и математических
результатов. В случае с определенным интегралом такой картинкой, несомненно,
будет площадь криволинейной трапеции (понимаемая пока на интуитивном уровне).

В этой главе мы даем формальное определение этому новому понятию, приводим
способ его вычисления и описываем некоторые приложения.

\section{Определенный интеграл}

\subsection{Определение определенного интеграла}

\paragraph{Разбиения отрезка.}

Перед тем как давать определение определенному интегралу нужно
проделать некую предварительную работу -- сказать несколько слов
про разбиения отрезка.

{\it Разбиением}\index{разбиение отрезка} отрезка $[a;b]$
называется произвольная конечная система точек $\{
x_0;x_1;...;x_k\}$ такая что
$$
 a=x_0<x_1<x_2<...<x_{k-1}<x_k=b
$$
Если нам дано какое-то разбиение $\{ x_0;x_1;...;x_k\}$, то его
удобно обозначать одной буквой, например, $\tau$:
$$
\tau=\{ x_0;x_1;...;x_k\}
$$
Наибольшее расстояние между соседними точками разбиения $\tau$
обозначается
$$
\diam\ \tau=\max_{i=1,...,k} (x_i-x_{i-1})
$$
и называется {\it диаметром разбиения}\index{диаметр разбиения}
$\tau=\{ x_0;x_1;...;x_k\}$.

\noindent\rule{160mm}{0.1pt}\begin{multicols}{2}

\begin{ex}\label{ex-14.1.1} Система точек $\tau=\{ 0;\frac{1}{2};1;\sqrt{2};2
\}$ является разбиением отрезка  $[0;2]$, потому что эти числа образуют строго
возрастающую конечную последовательность, у которой концы совпадают с концами
отрезка $[0;2]$:
$$
0<\frac{1}{2}<1<\sqrt{2}<2
$$
Диаметр этого разбиения равен
$$
\diam\ \tau=\max \left\{ \frac{1}{2};\frac{1}{2};\sqrt{2}-1;2-\sqrt{2}\right\}=
2-\sqrt{2}
$$
\end{ex}

\begin{er}\label{er-14.1.2} Проверьте, что следующие системы точек
являются разбиениями данных отрезков и найдите диаметр этих
разбиений:

1) $\tau=\{ 0;\frac{2}{3};\frac{5}{4};\frac{11}{5};3 \}, \quad
[a;b]=[0;3]$

2) $\tau=\{
-1;-\frac{1}{2};-\frac{1}{4};0;\frac{1}{4};\frac{1}{2};1 \}, \quad
[a;b]=[-1;1]$

3) $\tau=\{ 1;\sqrt{2};\sqrt{3};2 \}, \quad [a;b]=[1;2]$
\end{er}

\end{multicols}\noindent\rule[10pt]{160mm}{0.1pt}

Понятно, что у любого отрезка имеется много всяких разбиений. Поэтому можно
рассматривать {\it последовательности разбиений}\index{последовательность
разбиений}.

\noindent\rule{160mm}{0.1pt}\begin{multicols}{2}

\begin{ex}\label{ex-14.1.3} Рассмотрим, например, отрезок $[0;1]$ и его
разбиения:
 \begin{align*}
& \tau^{(1)}=\{ 0; 1 \}\\
& \tau^{(2)}=\left\{ 0; \frac{1}{2} ; 1\right\}\\
& \tau^{(3)}=\left\{ 0; \frac{1}{3} ; \frac{2}{3} ; 1\right\}\\
& ... \\
& \tau^{(n)}=\left\{ 0; \frac{1}{n} ; \frac{2}{n};...; \frac{n-1}{n} ; 1\right\} \\
& ...
 \end{align*}
Ясно, что для любого $n$ точки $\tau^{(n)}=\{ 0; \frac{1}{n} ; \frac{2}{n};...;
\frac{n-1}{n} ; 1\}$ действительно являются разбиением отрезка $[0;1]$. Диаметр
такого разбиения равен
$$
\diam\ \tau^{(n)}=\frac{1}{n}
$$
Можно заметить, что при $n\to \infty$ эта величина стремится к
нулю:
$$
\diam\ \tau^{(n)}=\frac{1}{n}\underset{n\to \infty}{\longrightarrow} 0
$$
\end{ex}

\end{multicols}\noindent\rule[10pt]{160mm}{0.1pt}

\bit{ \item[$\bullet$] Последовательность разбиений $\tau^{(n)}$, у которой
диаметр стремится к нулю
$$
\diam\ \tau^{(n)}\underset{n\to \infty}{\longrightarrow} 0
$$
называется {\it измельчающейся}\index{последовательность
разбиений!измельчающаяся}.
 }\eit

\noindent\rule{160mm}{0.1pt}\begin{multicols}{2}

\begin{er}\label{er-14.1.4} Проверьте,
будут ли сле\-дую\-щие последовательности разбиений отрезка $[0;1]$
измельчающимися:
 \biter{
\item[1)] $\tau^{(n)}=\left\{ 0;\frac{1}{n};\frac{n-1}{n};1 \right\}, \quad
n>2$ (здесь каждое разбиение состоит из четырех точек);

\item[2)] $\tau^{(n)}=\left\{
0;\frac{1}{2^n};\frac{2}{2^n};...;\frac{i}{2^n};...; \frac{2^n-1}{2^n};1
\right\}$

\item[3)] $\tau^{(n)}=\left\{ 0; \frac{1}{2n} ;\frac{1}{2n-1};
\frac{1}{2n-2};...; \frac{1}{n}; 1\right\}$

\item[4)] $\tau^{(n)}=\left\{ 0; \sqrt{\frac{1}{n}} ; \sqrt{\frac{2}{n}};...;
\sqrt{\frac{n-1}{n}} ; 1\right\}$
 }\eiter
\end{er}

\end{multicols}\noindent\rule[10pt]{160mm}{0.1pt}

\paragraph{Определение определенного интеграла.}

Пусть функция $f$ определена на отрезке $[a;b]$ и пусть дано разбиение этого
отрезка $\tau=\{ x_0;x_1;...;x_k\}$. В каждом маленьком отрезке $[x_{i-1};x_i]$
зафиксируем произвольную точку $\xi_i\in [x_{i-1};x_i]$. Такая система точек
$\xi_i$ называется {\it системой выделенных точек} в данном разбиении $\tau$.
Рассмотрим сумму
$$
\sum_{i=1}^k f(\xi_i)\cdot \Delta x_i \qquad (\text{где}\, \Delta
x_i=x_i-x_{i-1})
$$
Она называется {\it интегральной суммой}\index{сумма!интегральная} функции $f$,
соответствующей разбиению $\tau=\{ x_0;x_1;...;x_k\}$ с выделенными точками
$\xi_i$. Легко видеть, что эта величина равна площади ступенчатой фигуры,
построенной на точках $\xi_i\in [x_{i-1};x_i]$:

%\pucture{0pt}{0pt}{11-2.pcx}

\vglue100pt

Представим теперь, что разбиения $\tau=\{ x_0;x_1;...;x_k\}$
отрезка $[a;b]$ меняются, становясь все более мелкими. Тогда
площадь ступенчатой фигуры должна приближаться к числу, которое
естественно считать "площадью криволинейной трапеции", о которой
мы говорили в начале этой главы.

%\pucture{0pt}{0pt}{11-3.pcx}

\vglue100pt

\noindent Это число называется определенным интегралом функции $f$ на отрезке
$[a;b]$. Более точное определение выглядит так:

\bit{ \item[$\bullet$] Число $I$ называется {\it определенным
интегралом}\index{интеграл!определенный} функции $f$ на отрезке $[a;b]$, и
обозначается
$$
  I=\int_{[a,b]} f(x) \, \d x
$$
если для всякой измельчающейся последовательности разбиений отрезка $[a;b]$
$$
\tau^{(n)}: \quad a=x_0^{(n)}<x_1^{(n)}<...
<x_{k^{(n)}-1}^{(n)}<x_{k^{(n)}}^{(n)}=b, \qquad \diam\
\tau^{(n)}\underset{n\to \infty}{\longrightarrow} 0
$$
и любой системы выделенных точек
$$
\xi_i^{(n)}\in [x_{i-1}^{(n)};x_i^{(n)}]
$$
соответствующие интегральные суммы стремятся к числу $I$:
$$
\sum_{i=1}^{k^{(n)}} f(\xi_i^{(n)})\cdot \Delta x_i \underset{n\to
\infty}{\longrightarrow} I
$$
Если такое число $I$ существует, то функция $f$ называется {\it
интегрируемой}\index{функция!интегрируемая} (по Риману) на отрезке $[a;b]$.
 }\eit\noindent
Это определение можно коротко записать в виде формулы:
 \beq \label{14.2.1}
\int_{[a,b]} f(x) \, \d x=\lim_{\diam\ \{x_0,...,x_k\}\to 0}\sum_{i=1}^k
f(\kern-10pt\underset{\scriptsize\begin{matrix}\text{\rotatebox{-90}{$\in$}}\\
[x_{i-1}; x_i]\end{matrix}}{\xi_i}\kern-10pt)\cdot \Delta x_i
 \eeq
Ее удобно применять в случаях, когда не нужно следить за
строгостью рассуждений.

\noindent\rule{160mm}{0.1pt}\begin{multicols}{2}

Проверим, что в простейших случаях наше определение приводит к правильным
результатам, то есть что, вычисляя интеграл, мы действительно получаем площадь
соответствующей фигуры.

\begin{ex} {\bf Интеграл от константы.}\label{ex-14.2.1}
Рассмотрим постоянную функцию $f(x)=C$ и вычислим определенный
интеграл от нее на произвольном отрезке $[a;b]$. Соответствующая
картинка выглядит так:

%\pucture{0pt}{0pt}{11-4.pcx}

\vglue100pt

Поскольку криволинейная трапеция в данном случае представляет
собой прямоугольник, ответ нам известен заранее -- интеграл должен
быть равен $C\cdot (b-a)$. Проверим это:
 \begin{multline*}
\int_{[a,b]} f(x) \, \d x= \lim_{\diam\ \tau\to 0}\sum_{i=1}^k
f(\xi_i)\cdot \Delta x_i=\\
=\lim_{\diam\ \tau\to 0}\sum_{i=1}^k C\cdot \Delta x_i= \lim_{\diam\ \tau\to 0}
C\cdot \underbrace{\sum_{i=1}^k \Delta x_i}_{b-a}=\\= \lim_{\diam\ \tau\to 0}
C\cdot (b-a)= C\cdot (b-a)
\end{multline*}
Запишем вывод:
 \beq\label{14.2.2}
  \int_{[a,b]} C \, \d x=C\cdot (b-a)
 \eeq
\end{ex}

\begin{ex} {\bf Интеграл от линейной функции.}\label{ex-14.2.2}
Рассмотрим теперь функцию $f(x)=x$ и вычислим определенный
интеграл от нее на отрезке $[0;1]$. Соответствующая картинка
выглядит так:

%\pucture{0pt}{0pt}{11-5.pcx}

\vglue100pt

Криволинейная трапеция в данном случае представляет собой
треугольник, поэтому ответом должна быть площадь треугольника
$S=\frac{1}{2}\cdot 1\cdot 1=\frac{1}{2}$ (половина произведения
основания на высоту).

Чтобы проверить это, возьмем какое-нибудь разбиение
$\tau=\{x_0;x_1;...;x_k \}$ отрезка $[0;1]$ и какие-нибудь точки
$$
\xi_i\in [x_{i-1};x_i]
$$
Заметим, что тогда площадь нашего треугольника $S=\frac{1}{2}$
можно представить, как сумму площадей маленьких трапеций:

%\pucture{0pt}{0pt}{11-6.pcx}

\vglue100pt

 \begin{multline}\label{14.2.3}
S=\frac{1}{2}=\sum_{i=1}^{k}\frac{f(x_{i-1})+f(x_i)}{2}\cdot \Delta
x_i=\\
=\sum_{i=1}^{k}\frac{x_{i-1}+x_i}{2}\cdot \Delta x_i
 \end{multline}
Теперь мы получим
 \begin{multline*}
\bigg| \sum_{i=1}^{k} \overbrace{f(\xi_i)}^{\scriptsize\begin{matrix}\xi_i\\
\text{\rotatebox{90}{$=$}}\end{matrix}}\cdot \Delta
x_i-\kern-25pt\overbrace{\frac{1}{2}}^{\scriptsize\begin{matrix}
\sum_{i=1}^{k}\frac{x_{i-1}+x_i}{2}\cdot \Delta x_i\\ \phantom{\eqref{14.2.3}}\
\text{\rotatebox{90}{$=$}}
\ \eqref{14.2.3} \end{matrix}}\kern-25pt\bigg| =\\
=\left| \sum_{i=1}^{k}\xi_i\cdot \Delta x_i-
\sum_{i=1}^{k}\frac{x_{i-1}+x_i}{2}\cdot \Delta x_i \right|=\\= \left|
\sum_{i=1}^{k}\left(\xi_i- \frac{x_{i-1}+x_i}{2}\right) \cdot
\Delta x_i \right|=\\
=\left| \sum_{i=1}^{k}\frac{2\xi_i-x_{i-1}-x_i}{2}\cdot \Delta x_i \right|=\\=
\left| \sum_{i=1}^{k}\frac{(\xi_i-x_{i-1})+(\xi_i-x_i)}{2}\cdot \Delta x_i
\right|\le \\ \le\eqref{ind-nerav-s-modulem}\le\\
 \le \sum_{i=1}^{k}\Bigg|
\frac{(\xi_i-x_{i-1})+(\xi_i-x_i)}{2}\cdot \underbrace{\Delta
x_i}_{\scriptsize\begin{matrix}
\text{\rotatebox{90}{$\le$}}\\ 0\end{matrix}} \Bigg|=\\
= \sum_{i=1}^{k}\frac{1}{2}\cdot\Big| (\xi_i-x_{i-1})+(\xi_i-x_i)\Big| \cdot
\Delta x_i \le\\ \le \eqref{module-2^0}\le\\ \le \sum_{i=1}^{k}\frac{1}{2}\cdot
\Big(\underbrace{\left|\xi_i-x_{i-1}\right|}_{\scriptsize\begin{matrix}
\text{\rotatebox{90}{$\ge$}}\\ \phantom{,}\Delta x_i,\\
\text{поскольку} \\ \xi_i\in [x_{i-1};x_i]\end{matrix}}+
\underbrace{\left|\xi_i-x_i\right|}_{\scriptsize\begin{matrix}
\text{\rotatebox{90}{$\ge$}}\\ \phantom{,}\Delta x_i,\\
\text{поскольку} \\ \xi_i\in [x_{i-1};x_i]\end{matrix}}\Big)\cdot \Delta x_i
\le \\ \le \sum_{i=1}^{k}\frac{\Delta x_i+\Delta x_i }{2}\cdot \Delta x_i=
\sum_{i=1}^{k}\underbrace{\Delta
x_i}_{\scriptsize\begin{matrix}\text{\rotatebox{90}{$\ge$}}\\
\diam\ \tau\end{matrix}}\cdot \Delta x_i
\le\\
 \le \sum_{i=1}^{k}\diam\
\tau\cdot \Delta x_i=\\= \diam\ \tau\cdot \underbrace{\sum_{i=1}^{k}\Delta
x_i}_{\scriptsize\begin{matrix}\text{\rotatebox{90}{$=$}}\\
1\end{matrix}}=\diam\ \tau
\end{multline*}
Мы получили неравенство, справедливое для любого разбиения
$\tau=\{x_0;x_1;...;x_k \}$ отрезка $[0;1]$ и любых точек
$\xi_i\in [x_{i-1};x_i]$:
$$
\left| \sum_{i=1}^{k} f(\xi_i)\cdot \Delta x_i-\frac{1}{2}\right| \le \diam\
\tau
$$
Если теперь $\tau^{(n)}=\{x_0^{(n)};x_1^{(n)};...;x_{k_n}^{(n)}\}$ --
какая-нибудь измельчающая последовательность разбиений с выделенными точками
$\xi_i^{(n)}\in [x_{i-1}^{(n)};x_i^{(n)}]$, то мы получаем
$$
\left| \sum_{i=1}^{k^{(n)}} f(\xi_i^{(n)})\cdot \Delta
x_i^{(n)}-\frac{1}{2}\right| \le \diam\ \tau^{(n)} \underset{n\to
\infty}{\longrightarrow} 0
$$
Отсюда
$$
\sum_{i=1}^{k^{(n)}} f(\xi_i^{(n)})\cdot \Delta
x_i^{(n)}-\frac{1}{2}\underset{n\to \infty}{\longrightarrow} 0
$$
то есть
$$
\sum_{i=1}^{k^{(n)}} f(\xi_i^{(n)})\cdot \Delta
x_i^{(n)}\underset{n\to \infty}{\longrightarrow}\frac{1}{2}
$$
Это верно для любой измельчающейся последовательности
$\tau^{(n)}=\{x_0^{(n)};x_1^{(n)};...;x_k^{(n)}\}$ разбиений отрезка $[0;1]$ и
любой системы выделенных точек $\xi_i^{(n)}\in [x_{i-1}^{(n)};x_i^{(n)}]$,
значит число $\frac{1}{2}$ действительно является определенным интегралом
$\int_{[a,b]} x \, \d x$.

Запишем вывод:
$$
  \int_0^1 x \, \d x=\frac{1}{2}
$$
\end{ex}

\begin{ex} {\bf Интеграл от функции Дирихле.}\label{ex-14.2.3}
Покажем, что определенный интеграл не всегда существует. Классическим примером
здесь является функция Дирихле, которую мы определили выше формулой
\eqref{func-Dirichle}:
$$
D(x)=\begin{cases} 1, \,\, \text{если $x$ -- рациональное число}\\
0, \,\, \text{если $x$ -- иррациональное число}\end{cases}
$$
Попробуем понять, что будет интегралом этой функции на
произвольном отрезке $[a;b]$. Выберем какое-нибудь разбиение
$\tau=\{x_0;x_1;...;x_k\}$ отрезка $[a;b]$.

Точки $\xi_i\in [x_{i-1};x_i]$ можно выбирать по-разному, но нас
будут интересовать два способа:
 \biter{
\item[--] первый способ заключается в том, чтобы в каждом отрезке
$[x_{i-1};x_i]$ выбрать какую-нибудь {\it рациональную} точку
$\xi_i$. Тогда получится $D(\xi_i)=1$, и интегральные суммы
окажутся равны $b-a$:
$$
\sum_{i=1}^{k} D(\xi_i)\cdot \Delta x_i= \sum_{i=1}^{k} 1\cdot
\Delta x_i=b-a
$$
\item[--] второй способ заключается в том, чтобы в каждом отрезке
$[x_{i-1};x_i]$ выбрать какую-нибудь {\it иррациональную} точку
$\xi_i$. Тогда получится $D(\xi_i)=0$, и интегральные суммы
окажутся равны $0$:
$$
\sum_{i=1}^{k} D(\xi_i)\cdot \Delta x_i= \sum_{i=1}^{k} 0\cdot
\Delta x_i=0
$$
 }\eiter
 Если теперь взять измельчающуюся последовательность
разбиений $\tau^{(n)}=\{x_0^{(n)};x_1^{(n)};...;x_k^{(n)}\}$ отрезка
$[a;b]$,
$$
\diam\ \tau^{(n)}\underset{n\to \infty}{\longrightarrow} 0
$$
то выбрав точки $\xi_i^{(n)}$ первым способом, мы получим, что
интегральные суммы стремятся к числу $b-a$:
$$
\sum_{i=1}^{k^{(n)}} D(\xi_i^{(n)})\cdot \Delta x_i^{(n)}=b-a
\underset{n\to \infty}{\longrightarrow} b-a
$$
а для второго способа получим, что они стремятся к нулю:
$$
\sum_{i=1}^{k^{(n)}} D(\xi_i^{(n)})\cdot \Delta x_i^{(n)}=0
\underset{n\to \infty}{\longrightarrow} 0
$$
Спросим себя теперь: существует ли единое число $I$, к которому бы стремились
все интегральные суммы, независимо от того, какая выбрана измельчающаяся
последовательность разбиений $\tau^{(n)}=\{x_0;x_1;...;x_k^{(n)}\}$ и какие
выбраны точки $\xi_i^{(n)}\in [x_{i-1}^{(n)};x_i^{(n)}]$? Понятно, что такого
числа $I$ (то есть интеграла $\int_{[a,b]} D(x) \, \d x$) не существует.
Значит, мы можем сделать

Вывод: {\it функция Дирихле не интегрируема ни на каком отрезке $[a;b]$}.
\end{ex}

\end{multicols}\noindent\rule[10pt]{160mm}{0.1pt}

\subsection{Когда существует определенный интеграл?}

Из примера \ref{ex-14.2.3} видно, что не всякую функцию можно интегрировать. В
этом пункте мы обсудим вопрос о том, когда определенный интеграл существует.

\paragraph{Суммы Дарбу.}

Пусть функция $f$ определена и ограничена на отрезке $[a;b]$ и пусть $\tau=\{
x_0;x_1;...;x_k \}$ -- какое-нибудь его разбиение. На каждом отрезке
$[x_{i-1};x_i]$ найдем нижнюю и верхнюю грань функции $f$
$$
m_i=\inf_{\xi\in [x_{i-1};x_i]} f(\xi), \quad M_i=\sup_{\xi\in
[x_{i-1};x_i]} f(\xi)
$$
и обозначим
$$
s_\tau=\sum_{i=1}^k m_i\cdot \Delta x_i, \quad S_\tau=\sum_{i=1}^k
M_i\cdot \Delta x_i
$$
Поскольку $f$ ограничена на $[a;b]$, $m_i$ и $M_i$ являются обычными числами, и
поэтому $s_\tau$ и $S_\tau$ тоже являются числами, причем
 \begin{equation}
 s_\tau \le S_\tau
 \label{14.4.1}
 \end{equation}

\bit{ \item[$\bullet$] Величина $s_\tau$ называется {\it
нижней}\index{сумма!Дарбу!нижняя}, а величина $S_\tau$ -- {\it верхней суммой
Дарбу}\index{сумма!Дарбу!верхняя} (функции $f$ на отрезке $[a;b]$ при разбиении
$\tau$).
 }\eit

Перечислим главные

\bigskip

\centerline{\bf Свойства сумм Дарбу}\index{свойства!сумм!Дарбу}

{\it
 \bit{
\item[$1^0.$]
Если разбиение $\sigma$ получено из разбиения $\tau$ добавлением
новых точек
$$
\tau \subseteq \sigma
$$
то
 \begin{equation} \label{14.4.2}
 s_{\tau}\le s_{\sigma}\le S_{\sigma}\le S_{\tau}
 \end{equation}
(то есть, при добавлении к данному разбиению $\tau$ новых точек
нижняя сумма Дарбу увеличивается, а верхняя уменьшается).
\item[$2^0.$]
Для любых двух разбиений $\sigma$ и $\tau$ данного отрезка
выполняется неравенство
 \begin{equation}
s_{\sigma}\le S_{\tau}\label{14.4.3}
 \end{equation}
(то есть, любая нижняя сумма Дарбу меньше любой верхней).
\item[$3^0.$]
Любая интегральная сумма лежит между нижней и верхней суммами
Дарбу
 \begin{equation}\label{14.4.4}
s_{\tau}\le \sum_{i=1}^k f(\xi_i)\cdot \Delta x_i \le S_{\tau}
 \end{equation}
причем нижняя сумма Дарбу есть нижняя грань всех интегральных сумм
при данном разбиении, а верхняя сумма Дарбу -- верхняя грань:
 \begin{equation}
s_{\tau}=\inf_{\xi_i\in [x_{i-1};x_i]}\sum_{i=1}^k f(\xi_i)\cdot
\Delta x_i, \qquad S_{\tau}=\sup_{\xi_i\in [x_{i-1};x_i]}\sum_{i=1}^k
f(\xi_i)\cdot \Delta x_i \label{14.4.5}
 \end{equation}

\item[$4^0.$] Если функция $f$ интегрируема на $[a,b]$, то для произвольного
разбиения $\tau$ отрезка $[a,b]$ интеграл от $f$ лежит между нижней и верхней
суммами Дарбу:
 \begin{equation}\label{s_tau<int<S_tau}
s_{\tau}\le \int_{[a,b]} f(x)\ \d x \le S_{\tau}
 \end{equation}

 }\eit
}\begin{proof}

1. Для доказательства $1^0$ достаточно рассмотреть случай, когда к
данному разбиению $\tau=\{ x_0;x_1;...;x_k \}$ добавляется всего
одна точка $x^*$:
$$
\sigma=\tau \cup \{ x^* \}=\{ x_0;x_1;...;x_k \}\cup \{ x^* \}
$$
Тогда $x^*$ лежит в каком-то отрезке $[x_{j-1};x_j]$:
 \begin{equation}
 x^*\in [x_{j-1};x_j]
\label{14.4.6}
 \end{equation}
и поэтому
\begin{multline*}
s_\tau=\sum_{i=1}^k m_i\cdot \Delta x_i= m_j\cdot \Delta
x_j+\sum_{i\ne j} m_i\cdot \Delta x_i= m_j\cdot
(x_j-x_{j-1})+\sum_{i\ne j} m_i\cdot \Delta x_i=
 \eqref{14.4.6}
 =\\=
 \underbrace{\inf_{\xi\in [x_{j-1};x_j]} f(\xi)}_{\scriptsize\begin{matrix}
 \text{\rotatebox{90}{$\ge$}}\\ \inf\limits_{\xi\in [x^*,x_j]} f(\xi),\\ \text{потому что} \\
[x_{j-1};x_j]\supseteq [x^*,x_j] \end{matrix}}\cdot (x_j-x^*)+
 \underbrace{\inf_{\xi\in [x_{j-1};x_j]} f(\xi)}_{\scriptsize\begin{matrix}
 \text{\rotatebox{90}{$\ge$}}\\ \inf\limits_{\xi\in [x_{j-1},x^*]} f(\xi),\\ \text{потому что} \\
[x_{j-1};x_j]\supseteq [x_{j-1},x^*] \end{matrix}}\cdot (x^*-x_{j-1})+
\sum_{i\ne j} m_i\cdot \Delta x_i\le \\ \le
 \inf_{\xi\in [x^*,x_j]} f(\xi)\cdot (x_j-x^*)+
\inf_{\xi\in [x_{j-1},x^*]} f(\xi)\cdot (x^*-x_{j-1})+ \sum_{i\ne j} m_i\cdot
\Delta x_i=s_\sigma
 \end{multline*}
То есть, $s_\tau \le s_\sigma$. Аналогично доказывается. что
$S_\tau \ge S_\sigma$:
\begin{multline*}
S_\tau=\sum_{i=1}^k M_i\cdot \Delta x_i= M_j\cdot \Delta
x_j+\sum_{i\ne j} M_i\cdot \Delta x_i= M_j\cdot
(x_j-x_{j-1})+\sum_{i\ne j} M_i\cdot \Delta x_i
 =\eqref{14.4.6}=\\=
 \underbrace{\sup_{\xi\in [x_{j-1};x_j]} f(\xi)}_{\scriptsize\begin{matrix}
 \text{\rotatebox{90}{$\le$}}\\ \sup\limits_{\xi\in [x^*,x_j]} f(\xi),\\ \text{потому что} \\
[x_{j-1};x_j]\supseteq [x^*,x_j] \end{matrix}}\cdot (x_j-x^*)+
  \underbrace{\sup_{\xi\in [x_{j-1};x_j]} f(\xi)}_{\scriptsize\begin{matrix}
 \text{\rotatebox{90}{$\le$}}\\ \sup\limits_{\xi\in [x_{j-1},x^*]} f(\xi),\\ \text{потому что} \\
[x_{j-1};x_j]\supseteq [x_{j-1},x^*] \end{matrix}}\cdot (x^*-x_{j-1})+
\sum_{i\ne j} M_i\cdot \Delta x_i\ge \\ \ge
 \sup_{\xi\in [x_{j-1};x^*]} f(\xi)\cdot (x_j-x^*)+
\sup_{\xi\in [x^*;x_j]} f(\xi)\cdot (x^*-x_{j-1})+ \sum_{i\ne j}
M_i\cdot \Delta x_i=S_\sigma
 \end{multline*}
Оставшееся неравенство в \eqref{14.4.2}, $s_\sigma \le S_\sigma$, есть просто
неравенство \eqref{14.4.1}, в котором $\tau$ заменено на $\sigma$.

2. Если $\sigma$ и $\tau$ -- два разбиения отрезка $[a;b]$, то, добавлением
новых точек, можно сделать из них третье разбиение $\rho$, которое бы содержало
и $\sigma$, и $\tau$:
$$
\sigma \subseteq \rho, \qquad \tau \subseteq \rho
$$
Для него мы получим:
$$
\underbrace{s_{\sigma}\overset{\eqref{14.4.2}}{\le} s_{\rho}}_{\sigma \subseteq
\rho}\overset{\eqref{14.4.1}}{\le}
\underbrace{S_{\rho}\overset{\eqref{14.4.2}}{\le} S_{\tau}}_{\rho\supseteq
\tau}
$$

3. Неравенства \eqref{14.4.4} очевидны,
$$
s_{\tau}= \sum_{i=1}^k \inf_{\zeta_i \in [x_{i-1};x_i]} f(\zeta_i)\cdot \Delta
x_i \le \sum_{i=1}^k f(\xi_i)\cdot \Delta x_i \le \sum_{i=1}^k \sup_{\zeta_i
\in [x_{i-1};x_i]} f(\zeta_i)\cdot \Delta x_i \le S_{\tau},
$$
и из них сразу следуют неравенства
$$
s_{\tau}\le\inf_{\xi_i\in [x_{i-1};x_i]}\sum_{i=1}^k f(\xi_i)\cdot \Delta x_i,
\qquad \sup_{\xi_i\in [x_{i-1};x_i]}\sum_{i=1}^k f(\xi_i)\cdot \Delta x_i\le
S_{\tau}.
$$
Поэтому, чтобы, например, доказать первое равенство в \eqref{14.4.5},
достаточно доказать обратное неравенство:
 \beq\label{s_tau-ge-inf-int-sum}
s_{\tau}\ge\inf_{\xi_i\in [x_{i-1};x_i]}\sum_{i=1}^k f(\xi_i)\cdot \Delta x_i
 \eeq
Для этого (зафиксировав разбиение $\tau=\{x_0,...,x_k\}$) для каждого
$i=1,...,k$ выберем последовательность точек $\xi_i^{(n)}\in [x_{i-1};x_i]$
такую, что
$$
f(\xi_i^{(n)}) \underset{n\to \infty}{\longrightarrow}\inf_{\xi_i\in
[x_{i-1};x_i]} f(\xi_i)
$$
Тогда \eqref{s_tau-ge-inf-int-sum} получается так:
 \begin{multline*}
\inf_{\xi_i\in [x_{i-1};x_i]}\sum_{i=1}^k f(\xi_i)\cdot \Delta x_i \le
\inf_{n\to \infty}\sum_{i=1}^k f\l \xi_i^{(n)}\r\cdot \Delta x_i= \lim_{n\to
\infty}\sum_{i=1}^k f\l\xi_i^{(n)}\r\cdot \Delta x_i=\\= \sum_{i=1}^k\lim_{n\to
\infty} f\l\xi_i^{(n)}\r\cdot \Delta x_i= \sum_{i=1}^k \inf_{\xi_i\in
[x_{i-1};x_i]} f(\xi_i)\cdot \Delta x_i=s_\tau
 \end{multline*}
Второе равенство в \eqref{14.4.5} доказывается аналогично.

4. Для доказательства $4^0$, зафиксируем разбиение $\tau$. Для всякого более
мелкого разбиения $\sigma\supseteq\tau$ и любой системы выделенных точек
$\xi_i$ разбиения $\sigma$ мы получим:
$$
 s_{\tau}\overset{\eqref{14.4.2}}{\le} s_{\sigma}\overset{\eqref{14.4.4}}{\le}
\sum_{i=1}^k f(\xi_i)\cdot \Delta x_i
 \overset{\eqref{14.4.4}}{\le} S_{\sigma}\overset{\eqref{14.4.2}}{\le} S_{\tau}
$$
Если теперь заморозить $\tau$, а $\sigma$ выбирать все более измельчающимся, то
мы получим:
$$
 s_{\tau}\le
 \underbrace{\sum_{i=1}^k f(\xi_i)\cdot \Delta x_i}_{\scriptsize\begin{matrix}
 \phantom{\tiny\begin{matrix}\diam\sigma\\ \downarrow \\ 0\end{matrix}} \ \downarrow\ {\tiny\begin{matrix}\diam\sigma\\ \downarrow \\ 0\end{matrix}}
 \\
 \int_{[a,b]} f(x) \ \d x
 \end{matrix}}
 \le S_{\tau}
$$
Отсюда и следует \eqref{s_tau<int<S_tau}.
\end{proof}

\bit{ \item[$\bullet$] Верхняя грань всех нижних сумм Дарбу для функции $f$ на
отрезке $[a;b]$
$$
I_*=\sup_\tau s_\tau
$$
называется {\it нижним интегралом
Дарбу}\index{инетграл!Дарбу!нижний}, а нижняя грань всех верхних
сумм Дарбу
$$
I^*=\inf_\tau S_\tau
$$
называется {\it верхним интегралом
Дарбу}\index{инетграл!Дарбу!верхний}. Из неравенства \eqref{14.4.3}
следует, что для любых двух разбиений $\sigma$ и $\tau$ данного
отрезка выполняются неравенства
 \beq\label{14.4.7}
s_{\sigma}\le I_*\le I^*\le S_{\tau}
 \eeq
Если вдобавок функция $f$ интегрируема на $[a,b]$, то эту цепочку можно
дополнить так:
 \beq\label{14.4.7-*}
s_{\sigma}\le I_*\le\int_{[a,b]} f(x) \ \d x\le I^*\le S_{\tau}
 \eeq
 }\eit

\noindent\rule{160mm}{0.1pt}\begin{multicols}{2}

\begin{er}\label{er-14.4.1} Найдите верхний и нижний интегралы Дарбу
на отрезке $[-1;1]$ для функции сигнум $\sgn x$ (определенной выше формулой
\eqref{DEF:sgn}) и функции Дирихле (определенной формулой
\eqref{func-Dirichle}).
\end{er}

\end{multicols}\noindent\rule[10pt]{160mm}{0.1pt}

\paragraph{Критерий интегрируемости.}\index{критерий!интегрируемости}

\begin{tm}[\bf критерий интегрируемости]\footnote{Эта теорема используется ниже при
доказательстве теорем \ref{tm-14.3.1} и \ref{tm-14.3.2}.}\label{tm-14.5.1}
Ограниченная на отрезке $[a;b]$ функция $f$ тогда и только тогда интегрируема
на $[a;b]$, когда для всякой измельчающейся последовательности $\tau^{(n)}$
разбиений отрезка $[a;b]$
$$
\diam\ \tau^{(n)}\underset{n\to \infty}{\longrightarrow} 0
$$
соответствующие верхняя и нижняя суммы Дарбу стремятся друг к
другу:
 \beq
S_{\tau^{(n)}}-s_{\tau^{(n)}}\underset{n\to \infty}{\longrightarrow}
0 \label{14.5.1}
 \eeq
В этом случае (для всякой измельчающейся последовательности $\tau^{(n)}$
разбиений отрезка $[a;b]$) суммы Дарбу необходимо стремятся к интегралу от $f$
на $[a,b]$:
 \beq\label{S_tau_n->int<-s_tau_n}
S_{\tau^{(n)}}\underset{n\to\infty}{\longrightarrow}\int_{[a,b]}f(x)\ \d x
\underset{\infty\gets n}{\longleftarrow} s_{\tau^{(n)}}
 \eeq
\end{tm}\begin{proof} 1. Необходимость. Пусть
$f$ интегрируема на $[a;b]$, покажем что тогда выполняется \eqref{14.5.1}.
Возьмем какую-нибудь измельчающуюся последовательность разбиений $\tau^{(n)}$
отрезка $[a;b]$. Из левой формулы \eqref{14.4.5}
$$
s_{\tau}=\inf_{\xi_i\in [x_{i-1};x_i]}\sum_{i=1}^k f(\xi_i)\cdot
\Delta x_i,
$$
следует, что для всякого $n$ можно выбрать точки $\xi_i^{(n)}\in
[x_{i-1};x_i]$ так, чтобы
$$
\left| s_{\tau}- \sum_{i=1}^k f\left(\xi_i^{(n)}\right)\cdot
\Delta x_i \right| \le \frac{1}{n}
$$
Тогда мы получим
$$
s_{\tau^{(n)}}- \sum_{i=1}^k f\left(\xi_i^{(n)}\right)\cdot \Delta
x_i \underset{n\to \infty}{\longrightarrow} 0
$$
Теперь вспомним, что, по предположению, функция $f$ интегрируема на $[a;b]$, то
есть существует такое число $I$, что
 \beq\label{14.5.2}
\sum_{i=1}^k f\left(\xi_i^{(n)}\right)\cdot \Delta x_i \underset{n\to
\infty}{\longrightarrow} I
 \eeq
для всякой измельчающейся последовательности $\tau^{(n)}$ и любой системы точек
$\xi_i^{(n)}\in [x_{i-1};x_i]$. В частности, это должно быть верно для нашей
последовательности $\tau^{(n)}$ и выбранной нами системы точек $\xi_i^{(n)}\in
[x_{i-1};x_i]$. Отсюда
$$
s_{\tau^{(n)}}= \left(s_{\tau^{(n)}}- \sum_{i=1}^k
f\left(\xi_i^{(n)}\right)\cdot \Delta x_i\right) + \sum_{i=1}^k
f\left(\xi_i^{(n)}\right)\cdot \Delta x_i \underset{n\to
\infty}{\longrightarrow} 0+I=I
$$

Точно также из правой формулы \eqref{14.4.5}
$$
S_{\tau}=\sup_{\xi_i\in [x_{i-1};x_i]}\sum_{i=1}^k f(\xi_i)\cdot
\Delta x_i,
$$
следует, что для всякого $n$ можно выбрать (новые) точки
$\xi_i^{(n)}\in [x_{i-1};x_i]$ так, чтобы
$$
\left| S_{\tau}- \sum_{i=1}^k f\left(\xi_i^{(n)}\right)\cdot
\Delta x_i \right| \le \frac{1}{n}
$$
Тогда мы получим
$$
S_{\tau^{(n)}}- \sum_{i=1}^k f\left(\xi_i^{(n)}\right)\cdot \Delta
x_i \underset{n\to \infty}{\longrightarrow} 0
$$
Поскольку функция $f$ интегрируема, выполняется \eqref{14.5.2}, значит
$$
S_{\tau^{(n)}}= \left(S_{\tau^{(n)}}- \sum_{i=1}^k
f\left(\xi_i^{(n)}\right)\cdot \Delta x_i\right) + \sum_{i=1}^k
f\left(\xi_i^{(n)}\right)\cdot \Delta x_i \underset{n\to
\infty}{\longrightarrow} 0+I=I
$$

Таким образом, мы получили, что
$$
s_{\tau^{(n)}}\underset{n\to \infty}{\longrightarrow} I, \qquad
S_{\tau^{(n)}}\underset{n\to \infty}{\longrightarrow} I
$$
откуда и следует \eqref{14.5.1}:
$$
S_{\tau^{(n)}}-s_{\tau^{(n)}}\underset{n\to \infty}{\longrightarrow}
0
$$

2. Докажем достаточность. Пусть для всякой измельчающейся последовательности
$\tau^{(n)}$ разбиений отрезка $[a;b]$
$$
\diam\ \tau^{(n)}\underset{n\to \infty}{\longrightarrow} 0
$$
соответствующие верхняя и нижняя суммы Дарбу стремятся друг к
другу:
$$
S_{\tau^{(n)}}-s_{\tau^{(n)}}\underset{n\to \infty}{\longrightarrow}
0
$$
Заметим, что тогда верхний и нижний интегралы Дарбу должны
совпадать
$$
  I_*=I^*,
$$
потому что из неравенств \eqref{14.4.7} следует
$$
0\le I^*-I_*\le
S_{\tau^{(n)}}-s_{\tau^{(n)}}\underset{n\to\infty}{\longrightarrow}0.
$$
Теперь положим
$$
I=I_*=I^*
$$
и покажем, что $I$ является интегралом функции $f$ на отрезке $[a;b]$.

Для этого опять возьмем измельчающуюся последовательность $\tau^{(n)}$
разбиений отрезка $[a;b]$. Из неравенств \eqref{14.4.7} получаем две цепочки:
 \begin{align*}
&
\begin{matrix}
0\le I-s_{\tau^{(n)}}\le
\underbrace{S_{\tau^{(n)}}-s_{\tau^{(n)}}}_{\scriptsize\begin{matrix}
\phantom{\tiny
\begin{matrix}n\\ \downarrow\\ \infty\end{matrix}} \ \downarrow \
{\tiny
\begin{matrix}n\\ \downarrow\\ \infty\end{matrix}}\\ 0
\end{matrix}}
\\ \Downarrow \\
I-s_{\tau^{(n)}}\underset{n\to\infty}{\longrightarrow}0
\\ \Downarrow \\
s_{\tau^{(n)}}\underset{n\to\infty}{\longrightarrow} I
\end{matrix}
&&
\begin{matrix}
0\le S_{\tau^{(n)}}-I\le
\underbrace{S_{\tau^{(n)}}-s_{\tau^{(n)}}}_{\scriptsize\begin{matrix}
\phantom{\tiny
\begin{matrix}n\\ \downarrow\\ \infty\end{matrix}} \ \downarrow \
{\tiny
\begin{matrix}n\\ \downarrow\\ \infty\end{matrix}}\\ 0
\end{matrix}}
\\ \Downarrow \\
S_{\tau^{(n)}}-I\underset{n\to\infty}{\longrightarrow}0
\\ \Downarrow \\
S_{\tau^{(n)}}\underset{n\to\infty}{\longrightarrow} I
\end{matrix}
 \end{align*}
Теперь из формулы \eqref{14.4.4} при произвольном выборе точек $\xi_i^{(n)}\in
[x_{i-1};x_i]$ мы получаем
$$
\underbrace{s_{\tau^{(n)}}}_{\scriptsize\begin{matrix} {\tiny
\begin{matrix}n\\ \downarrow\\ \infty\end{matrix}} \ \downarrow \
\phantom{\tiny
\begin{matrix}n\\ \downarrow\\ \infty\end{matrix}}\\ I
\end{matrix}}\le
\sum_{i=1}^k f(\xi_i^{(n)})\cdot \Delta x_i^{(n)}\le
\underbrace{S_{\tau^{(n)}}}_{\scriptsize\begin{matrix} \phantom{\tiny
\begin{matrix}n\\ \downarrow\\ \infty\end{matrix}} \ \downarrow \
{\tiny
\begin{matrix}n\\ \downarrow\\ \infty\end{matrix}}\\ I
\end{matrix}}
$$
Последовательности по бокам этого двойного неравенства стремятся к $I$, значит,
интегральная сумма тоже стремится к $I$:
$$
\sum_{i=1}^k f(\xi_i^{(n)})\cdot \Delta x_i^{(n)}\underset{n\to
\infty}{\longrightarrow} I
$$
Это верно для всякой измельчающейся последовательности $\tau^{(n)}$ разбиений
отрезка $[a;b]$ и любой системы выделенных точек $\xi_i^{(n)}\in
[x_{i-1}^{(n)};x_i^{(n)}]$, значит $I$ действительно является интегралом для
$f$ на $[a;b]$. Это нам и нужно было доказать. \end{proof}

\paragraph{Ограниченность интегрируемой функции.}

\begin{tm}\label{tm-14.3.1} Если функция $f$ интегрируема на отрезке $[a;b]$, то она ограничена на нем.
\end{tm}

\begin{proof} Докажем это утверждение в эквивалентной формулировке: {\it если функция $f$ (определена и) НЕ
ограничена на отрезке $[a;b]$, то она НЕ интегрируема на нем.} Пусть функция
$f$ не ограничена на отрезке $[a;b]$. Это означает, что она или неограничена
сверху, или неограничена снизу:
$$
\sup_{x\in [a;b]} f(x)=+\infty \quad \text{или}\quad \inf_{x\in [a;b]}
f(x)=-\infty
$$
Предположим, для определенности, что она неограничена сверху:
$$
\sup_{x\in [a;b]} f(x)=+\infty
$$
Зафиксируем какое-нибудь число $M$ и возьмем произвольное разбиение
$\tau=\{x_0;x_1;...;x_k\}$ отрезка $[a;b]$.

Поскольку $f$ неограничена на отрезке $[a;b]$, она должна быть неограничена на
каком-нибудь отрезке $[x_{i-1};x_i]$. То есть, существует индекс $j$ такой что
 \begin{equation}
  \sup_{\xi\in [x_{j-1};x_j]} f(\xi)=+\infty
  \label{14.3.3}
 \end{equation}
Зафиксируем этот индекс $j$ и выберем произвольным образом точки $\xi_i\in
[x_{i-1};x_i], \quad i\ne j$. Тогда в интегральной сумме
$$
\sum_{i=1}^{k} f(\xi_i)\cdot \Delta x_i= f(\xi_j)\cdot \Delta x_j+\sum_{i\ne j}
f(\xi_i)\cdot \Delta x_i
$$
величина $\sum_{i\ne j} f(\xi_i)\cdot \Delta x_i$ фиксирована, а слагаемое
$f(\xi_j)\cdot \Delta x_j$ -- переменная, поскольку точка $\xi_j\in
[x_{j-1};x_j]$ еще не выбрана. Из \eqref{14.3.3} следует, что недостающую точку
$\xi_j\in [x_{j-1};x_j]$ можно выбрать так, чтобы вся сумма была больше числа
$M$:
$$
\sum_{i=1}^{k} f(\xi_i)\cdot \Delta x_i= f(\xi_j)\cdot \Delta x_j+\sum_{i\ne j}
f(\xi_i)\cdot \Delta x_i \ge M
$$

Мы получили, что для всякого разбиения  $\tau=\{x_0;x_1;...;x_k\}$ отрезка
$[a;b]$ можно выбрать числа $\xi_i\in [x_{i-1};x_i]$ так, чтобы интегральная
сумма оказалась больше числа $M$.

Значит, если взять измельчающуюся последовательность разбиений $\tau^{(n)}$, то
для нее тоже можно выбрать числа $\xi_i^{(n)}\in [x_{i-1}^{(n)};x_i^{(n)}]$
так, чтобы интегральная сумма оказалась больше числа $M$:
$$
\sum_{i=1}^{k^{(n)}} f(\xi_i)\cdot \Delta x_i\ge M
$$
Отсюда следует, что интеграл (то есть предел таких интегральных сумм), если он
существует, не может быть меньше $M$:
$$
I=\lim_{n\to\infty}\sum_{i=1}^{k^{(n)}} f(\xi_i)\cdot \Delta x_i \ge M
$$
Но число $M$ мы выбирали произвольным: получается, что число $I$ должно быть
больше какого хочешь другого числа $M$, а это абсурд.
\end{proof}

\paragraph{Интегрируемость монотонной функции.}

\begin{tm}\label{tm-14.3.2} Если функция $f$ (определена и)
монотонна на отрезке $[a;b]$, то она интегрируема на нем.
\end{tm}
\begin{proof} Пусть
функция $f$ (определена и) монотонна на отрезке $[a;b]$, например, неубывает на
нем:
$$
s\le t \quad \Longrightarrow \quad f(s)\le f(t)
$$
покажем, что тогда $f$ интегрируема на $[a;b]$. Для этого возьмем какое-нибудь
разбиение $\tau$ отрезка $[a;b]$ и покажем, что суммы Дарбу удовлетворяют
неравенству:
 \begin{equation} S_{\tau}-s_{\tau}\le \diam\ \tau\cdot \left\{f(b)- f(a)\right\}
 \label{S_tau-s_tau le Delta tau left f(b)- f(a)}
 \end{equation}
Действительно,
\begin{multline*}
 S_{\tau}-s_{\tau}=
\sum_{i=1}^{k}
 \underbrace{\sup_{\xi\in [x_{i-1};x_i]} f(\xi)}_{\tiny
 \begin{array}{c}
 \| \\ f\left(x_i\right), \\ \text{поскольку $f$}\\
 \text{монотонно неубывает}
 \end{array}
 }
 \cdot \Delta x_i-
\sum_{i=1}^{k}
  \underbrace{
  \inf_{\xi\in [x_{i-1};x_i]} f(\xi)}_{\tiny
   \begin{array}{c}
   \|
   \\
   f\left(x_{i-1}\right),
   \\ \text{поскольку $f$}
   \\
   \text{монотонно неубывает}
   \end{array}
  }
  \cdot \Delta x_i
  =
  \sum_{i=1}^{k} f\left(x_i\right)\cdot \Delta x_i-
  \sum_{i=1}^{k} f\left(x_{i-1}\right)\cdot \Delta x_i
  =\\=
  \sum_{i=1}^{k}\left(f\left(x_i\right)-f\left(x_{i-1}\right)\right)\cdot
  \underbrace{\Delta x_i}_{\scriptsize\begin{matrix}
  \text{\rotatebox{90}{$\ge$}} \\
  \diam\ \tau\end{matrix}}
  \le \sum_{i=1}^{k}\left(f\left(x_i\right)-
 f\left(x_{i-1}\right)\right)\cdot \kern-25pt
 \underbrace{\diam\ \tau}_{\tiny
 \begin{array}{c}
 \text{не зависит}\\
 \text{от индекса}\,\, i,  \\
 \text{поэтому можно}\\
 \text{вынести за знак суммы}
 \end{array}
 }\kern-25pt
=\diam\ \tau\cdot \sum_{i=1}^{k}\left(f\left(x_i\right)- f\left(x_{i-1}\right)\right)=\\
=\diam\ \tau \cdot
  \Big\{
  \underbrace{
  \begin{array}[b]{c} f(b) \\ \text{\rotatebox{90}{$=$}} \\ f(x_k)
  \end{array}
  -
 f(x_{k-1})
 \put(-18,20){
 \text{\tiny сокращаются}
 \put(-41,-2){\line(1,0){40}}
 \put(-43,-8){$\downarrow\kern35pt\downarrow$}}
 }_{i=k}
  +
  \underbrace{
  f(x_{k-1})
  -
  f(x_{k-2})
   \put(-21,20){
 \text{\tiny сокращаются}
 \put(-38,-2){\line(1,0){33}}
 \put(-40,-8){$\downarrow\kern28pt\downarrow$}}
    }_{i=k-1}
  +
  \phantom{f(x)}
  ...
  \phantom{f(x)}
   \put(-16,20){
 \text{\tiny сокращаются}
 \put(-36,-2){\line(1,0){30}}
 \put(-38,-8){$\downarrow\kern25pt\downarrow$}}
  +
  \underbrace{
  f(x_2)
   -
  f(x_1)
   \put(-18,20){
 \text{\tiny сокращаются}
 \put(-36,-2){\line(1,0){30}}
 \put(-38,-8){$\downarrow\kern25pt\downarrow$}}
  }_{i=2}
  +
  \underbrace{
  f(x_1)
  -
  \begin{array}[b]{c} f(a) \\ \text{\rotatebox{90}{$=$}} \\ f(x_0)
  \end{array}
  }_{i=1}
  \Big\}
  =\\=
  \diam\ \tau\cdot \left\{f(b)-
  f(a)\right\}
 \end{multline*}
 Если теперь взять измельчающуюся последовательность разбиений $\tau^{(n)}$  отрезка
$[a;b]$, то мы получим
$$
 0\le S_{\tau^{(n)}}-s_{\tau^{(n)}}\le
 \eqref{S_tau-s_tau le Delta tau left f(b)- f(a)} \le
 \diam\ \tau^{(n)}\cdot
 \left\{f(b)- f(a)\right\}\underset{n\to \infty}{\longrightarrow}
 0
$$
значит,
$$
S_{\tau^{(n)}}-s_{\tau^{(n)}}\underset{n\to \infty}{\longrightarrow} 0
$$
Это верно для любой измельчающейся последовательности разбиений $\tau^{(n)}$,
значит, по теореме \ref{tm-14.5.1}, $f$ интегрируема на $[a;b]$.\end{proof}

\paragraph{Интегрируемость непрерывной функции.}

\begin{tm}\label{tm-14.3.3} Если функция $f$ (определена и)
непрерывна на отрезке $[a;b]$, то она интегрируема на нем.
\end{tm}
Для доказательства теоремы \ref{tm-14.3.3} нам понадобится следующая

\begin{lm}\label{LM-o-razbieniyah-nepreryvnij-funktsii}
Если $f$ -- непрерывная функция на отрезке $[a;b]$, то для любого разбиения
$\tau$ отрезка $[a;b]$ существуют числа $\alpha, \beta \in [a;b]$, такие что
расстояние между ними не превышает диаметра разбиения $\tau$
\begin{equation}
 |\alpha-\beta|\le \diam\ \tau,
 \label{alpha-beta-le-Delta-tau}\end{equation}
а разность между верхней и нижней интегральными суммами Дарбу этого разбиения
оценивается неравенством:
\begin{equation}
S_{\tau}-s_{\tau}\le \Big|f(\alpha)-f(\beta)\Big|\cdot (b-a)
 \label{S_tau-s_tau le f(alpha)-f(beta)}\end{equation}\end{lm}
 \begin{proof} По теореме Вейерштрасса об экстремумах \ref{Wei-III},
 функция $f$ достигает максимума и минимума на каждом отрезке
 $[x_{i-1};x_i]$, то есть существуют такие точки $\eta_i,\theta_i \in [x_{i-1};x_i]$, что
 \beq\label{14.7.2}
f(\eta_i) =\sup_{x\in [x_{i-1};x_i]} f(x) \qquad \& \qquad f(\theta_i)
=\inf_{x\in [x_{i-1};x_i]} f(x)
 \eeq
Рассмотрим конечный набор чисел
$$
f\left(\eta_1\right)-f\left(\theta_1\right), \,\,
f\left(\eta_2\right)-f\left(\theta_2\right), \,\, ..., \,\,
f\left(\eta_i\right)-f\left(\theta_i\right), \,\, ..., \,\,
f\left(\eta_{k}\right)-f\left(\theta_{k}\right)
$$
и выберем среди них максимальное: пусть соответствующий индекс обозначается
$j$:
 \beq\label{14.7.3}
f\left(\eta_{j}\right)-f\left(\theta_{j}\right)= \max_{i=1,2,...,k}\left(
f\left(\eta_i\right)-f\left(\theta_i \right)\right)
 \eeq
Обозначив теперь
$$
\alpha=\eta_{j}, \qquad \beta=\theta_{j}
$$
мы получим нужные числа: с одной стороны, поскольку $\alpha, \beta \in
[x_{j-1};x_{j}]$, получается
$$
|\alpha-\beta|\le \Delta x_j \le \Delta \Big(\tau \Big)
$$
А, с другой стороны,
\begin{equation}
 f(\alpha)-f(\beta)=\max_{1\le i\le k}\Big(
 f(\eta_i)-f(\theta_i)\Big)=\max_{1\le i\le k}\l  \sup_{x\in
 [x_{i-1};x_i]} f(x) - \inf_{x\in [x_{i-1};x_i]} f(x) \r
 \label{f(alpha)-f(beta)}\end{equation}
Поэтому
\begin{multline*}
 S_{\tau}-s_{\tau}=
\sum_{i=1}^{k}\sup_{\xi\in [x_{i-1};x_i]} f(\xi)\cdot \Delta x_i-
\sum_{i=1}^{k}\inf_{\xi\in [x_{i-1};x_i]} f(\xi)\cdot \Delta x_i =
 \sum_{i=1}^{k}
 \underbrace{\left[ \sup_{\xi\in [x_{i-1};x_i]} f(\xi) -
 \inf_{\xi\in [x_{i-1};x_i]} f(\xi) \right]}_{
 \begin{array}[t]{c}
 \phantom{\ \eqref{f(alpha)-f(beta)}}
 \IA
 \ \eqref{f(alpha)-f(beta)}
 \\
 f(\alpha) - f(\beta)
 \end{array}
 }
 \cdot \Delta x_i
 \le \\ \le
 \sum_{i=1}^{k}
 \underbrace{\Big[ f(\alpha) - f(\beta) \Big]}_{\scriptsize
 \begin{array}[t]{c}
 \text{не зависит}
 \\
 \text{от индекса $i$}
 \end{array}
 }
 \cdot \Delta x_i =
 \Big[ f(\alpha)-f(\beta) \Big] \cdot
 \underbrace{\sum_{i=1}^{k}\Delta x_i}_{\scriptsize
 \text{длина отрезка [a,b]}
 }
 = \Big[ f(\alpha)-f(\beta) \Big]\cdot (b-a)
 \end{multline*}\end{proof}

\begin{proof}[Доказательство теоремы \ref{tm-14.3.3}] Пусть
$\tau^{(n)}$ -- измельчающаяся последовательность разбиений отрезка $[a;b]$:
$$
\diam\ \tau^{(n)}\underset{n\to \infty}{\longrightarrow} 0
$$
По лемме \ref{LM-o-razbieniyah-nepreryvnij-funktsii}, найдутся числа
$\alpha^{(n)},\beta^{(n)}\in [a,b]$ такие, что
$$
 |\alpha^{(n)}-\beta^{(n)}|\le \diam\ \tau^{(n)}
$$
и
$$
S_{\tau^{(n)}}-s_{\tau^{(n)}}\le \Big[f(\alpha^{(n)})-f(\beta^{(n)})\Big]\cdot
(b-a)
$$
Из первого условия следует, что последовательности $\alpha^{(n)}$ и
$\beta^{(n)}$ стремятся друг к другу:
$$
\alpha^{(n)}-\beta^{(n)}\underset{n\to \infty}{\longrightarrow} 0
$$
Поэтому по теореме Кантора \ref{Kantor} (функция $f$ непрерывна на $[a;b]$,
значит она должна быть равномерно непрерывна на нем), получаем:
$$
f(\alpha^{(n)})-f(\beta^{(n)}) \underset{n\to \infty}{\longrightarrow} 0
$$
Отсюда уже следует
$$
0\le S_{\tau^{(n)}}-s_{\tau^{(n)}}\le
\Big[f(\alpha^{(n)})-f(\beta^{(n)})\Big]\cdot (b-a)
 \underset{n\to \infty}{\longrightarrow} 0
 \quad \Longrightarrow \quad
 S_{\tau^{(n)}}-s_{\tau^{(n)}}\underset{n\to \infty}{\longrightarrow} 0
$$
Это верно для любой измельчающейся последовательности разбиений $\tau^{(n)}$,
значит, по теореме \ref{tm-14.5.1}, $f$ интегрируема на $[a;b]$.  \end{proof}

\subsection{Свойства определенного
интеграла}\label{SEC-svojstva-opredelyonnogo-integrala}

\bigskip

\centerline{\bf Свойства определенного интеграла}

 \bit{\it

\item[$1^0.$] {\bf Линейность:} если функции $f$ и $g$ интегрируемы на отрезке
$[a;b]$ то для любых чисел $\alpha,\beta \in \R$ функция $\alpha\cdot
f+\beta\cdot g$ тоже интегрируема на отрезке $[a;b]$, причем
 \begin{equation}\label{14.6.1}
\int_{[a,b]} \Big( \alpha\cdot f(x)+\beta\cdot g(x) \Big)\, \d x= \alpha\cdot
\int_{[a,b]} f(x)\, \d x+\beta\cdot \int_{[a,b]} g(x)\, \d x
 \end{equation}

\item[$2^0.$]\label{additivnost-integrala} {\bf Аддитивность:} пусть $a<b<c$ и
функция $f$ интегрируема на отрезках $[a;b]$ и $[b;c]$, тогда она интегрируема
на отрезке $[a;c]$, и при этом
 \begin{equation}\label{14.6.2}
\int_{[a,b]} f(x)\, \d x+\int_{[b,c]} f(x)\, \d x= \int_{[a,c]}  f(x) \, \d x
 \end{equation}

\item[$3^0.$] {\bf Монотонность:} если функции $f$ и $g$ интегрируемы на
отрезке $[a;b]$ и при этом
 \begin{equation}\label{14.6.3}
  f(x)\le g(x),  \qquad x\in [a;b]
 \end{equation}
то
 \begin{equation}\label{14.6.4}
\int_{[a,b]}  f(x) \, \d x\le \int_{[a,b]} g(x)\, \d x
 \end{equation}

\item[$4^0.$] {\bf Выпуклость:} если функция $f$ интегрируема на отрезке
$[a;b]$, то функция $|f|$ тоже интегрируема на отрезке $[a;b]$, и при этом
 \begin{equation}\label{14.6.5}
\left|\int_{[a,b]}  f(x) \, \d x\right| \le \int_{[a,b]} \left|f(x)\right| \,
\d x
 \end{equation}

\item[$5^0.$] {\bf Оценка сверху:} если функция $f$ интегрируема на отрезке
$[a;b]$, то
 \begin{equation}\label{otsenka-integrala}
\left|\int_{[a,b]}  f(x) \, \d x\right|\le (b-a)\cdot\sup_{x\in[a,b]}|f(x)|
 \end{equation}

\item[$6^0.$] {\bf Непрерывность:} если функция $f$ интегрируема на отрезке
$[a;b]$, то для всякой точки $c\in[a,b]$
 \begin{align}\label{nepreryvnost-integrala-1}
&\int_{[a,y]}  f(x) \, \d x\underset{y\to c}{\longrightarrow} \int_{[a,c]}
f(x) \, \d
x \\
\label{nepreryvnost-integrala-2} &\int_{[y,b]}  f(x) \, \d x\underset{y\to
c}{\longrightarrow} \int_c^b  f(x) \, \d x
\end{align}

\item[$7^0.$] {\bf Теорема о среднем:} если функция $f$ непрерывна на отрезке
$[a;b]$, то найдется такая точка $\xi\in [a;b]$, что
 \begin{equation}\label{14.6.6}
\int_{[a,b]}  f(x) \, \d x=f(\xi)\cdot (b-a)
 \end{equation}

\item[$8^0.$]\label{integrir-proizvedeniya} {\bf Интегрируемость произведения:}
если функции $f$ и $g$ интегрируемы на отрезке $[a;b]$, то их произведение
$f\cdot g$ -- тоже интегрируемая функция на отрезке $[a;b]$.

\item[$9^0.$]\label{tm-15.4.1-1} {\bf Замена переменной:} пусть
 \bit{
\item[1)] функция $f$ интегрируема на отрезке $[a; b]$, и

\item[2)] функция $\ph$ определена на отрезке $[\alpha; \beta]$ и обладает
свойствами:
 \bit{
\item[a)] $\ph$ -- гладкая на отрезке $[\alpha; \beta]$;

\item[b)] $\ph$ -- строго монотонная на отрезке $[\alpha; \beta]$;

\item[c)] $\ph$ сюръективно отображает отрезок $[\alpha; \beta]$ на отрезок
$[a;b]$:
$$
\ph\Big([\alpha; \beta]\Big)=[a,b],
$$
 }\eit
 }\eit
тогда функция $(f\circ\ph)\cdot \ph'$ интегрируема на отрезке $[\alpha,\beta]$,
и
 \beq\label{15.4.1-1}
\int_{[a,b]} f (x) \, \d x= \int_{[\alpha,\beta]} f (\ph(t))\cdot |\ph'(t)| \,
\d t
 \eeq

 }\eit

\bigskip

\begin{proof}

1. Линейность. Пусть функции $f$ и $g$ интегрируемы на отрезке $[a;b]$, и
$$
  F=\int_{[a,b]} f(x) \, \d x, \qquad G=\int_{[a,b]} g(x) \, \d x
$$
Для всякого разбиения $\tau=\{x_0,...,x_k\}$ отрезка $[a;b]$ и любой системы
выделенных точек $\xi_i\in [x_{i-1};x_i]$ мы получим:
$$
\sum_{i=1}^{k}\Big( \alpha \cdot f\left(\xi_i\right) +\beta \cdot
g\left(\xi_i\right) \Big)\cdot \Delta x_i= \alpha \cdot
\underbrace{\sum_{i=1}^{k} f\left(\xi_i\right) \cdot \Delta
x_i}_{\scriptsize\begin{matrix} \phantom{\tiny\begin{matrix}\diam\tau
\\ \downarrow\\ 0\end{matrix}} \downarrow
{\tiny\begin{matrix}\diam\tau
\\ \downarrow\\ 0\end{matrix}}
\\ F\end{matrix}} +\beta \cdot
\underbrace{\sum_{i=1}^{k} g\left(\xi_i\right) \cdot \Delta
x_i}_{\scriptsize\begin{matrix} \phantom{\tiny\begin{matrix}\diam\tau
\\ \downarrow\\ 0\end{matrix}} \downarrow
{\tiny\begin{matrix}\diam\tau
\\ \downarrow\\ 0\end{matrix}}
\\ F\end{matrix}}\underset{\diam\tau\to 0}{\longrightarrow}\alpha\cdot F+\beta\cdot
G
$$
Это нужно понимать так: если взять последовательность разбиений $\tau^{(n)}$,
диаметры которых стремятся к нулю, то какую ни выбирай систему выделенных точек
$\xi_i^{(n)}$, интегральная сумма функции $\alpha \cdot f +\beta \cdot g$ будет
стремится к числу $\alpha\cdot F+\beta\cdot G$. То есть, функция $\alpha \cdot
f +\beta \cdot g$ получается интегрируемой, и
$$
\int_{[a,b]} \Big(\alpha\cdot f(x)+\beta\cdot g(x) \Big)\, \d x= \alpha\cdot
F+\beta\cdot G= \alpha\cdot \int_{[a,b]} f(x)\, \d x+\beta\cdot \int_{[a,b]}
g(x)\, \d x
$$

2. Аддитивность. Пусть $f$ интегрируема на отрезках $[a;b]$ и $[b;c]$. Тогда по
теореме \ref{tm-14.3.1}, $f$ ограничена на отрезках $[a;b]$ и $[b;c]$, а значит
и на отрезке $[a;c]$:
 \beq\label{14.6.7}
  \exists M\in \R\quad \forall x\in [a;c] \quad |f(x)|\le M
 \eeq
Обозначим через $A$ и $B$ интегралы по отрезкам $[a;b]$ и $[b;c]$
$$
  A=\int_{[a,b]} f(x) \, \d x, \quad B=\int_{[b,c]} f(x) \, \d x
$$
Нам нужно показать, что $f$ интегрируема на  $[a;c]$, причем
$$
  \int_{[a,c]} f(x) \, \d x=A+B
$$
Возьмем разбиение $\tau=\{x_0,...,x_k\}$ отрезка $[a;c]$. Точка $b$ попадет в
какой-то полуинтервал $(x_{j-1},x_{j}]$:
$$
\tau=\{a=x_0<x_1<...<x_{j-1}<b\le x_{j}< x_{j+1}<...<x_{k}=c \}
$$
Запомним этот индекс $j$ и выберем произвольные точки $\xi_i\in [x_{i-1};
x_i]$. Тогда:
 \begin{multline*}
\overbrace{\sum_{i=1}^{k} f\left(\xi_i\right)\cdot \Delta
x_i}^{\scriptsize\begin{matrix}\text{интегральная сумма}\\ \text{на отрезке
$[a;c]$}\end{matrix}}= \sum_{i=1}^{j-1} f\left(\xi_i\right)\cdot \Delta x_i+
f\left(\xi_{j}\right)\cdot \Delta x_{j}+ \sum_{i=j+1}^{k}
f\left(\xi_i\right)\cdot \Delta x_i=\\= \sum_{i=1}^{j-1}
f\left(\xi_i\right)\cdot \Delta x_i+ \Big( f(b)+
f\left(\xi_{j}\right)-f(b)\Big) \cdot \Delta x_{j}+ \sum_{i=j+1}^{k}
f\left(\xi_i\right)\cdot \Delta x_i=\\= \sum_{i=1}^{j-1}
f\left(\xi_i\right)\cdot \Delta x_i+ f(b) \cdot \Delta x_{j}+ \sum_{i=j+1}^{k}
f\left(\xi_i\right)\cdot \Delta x_i+ \Big( f\left(\xi_{j}\right)-f(b)\Big)
\cdot \Delta x_{j} =\\= \sum_{i=1}^{j-1} f\left(\xi_i\right)\cdot \Delta x_i+
f(b) \cdot \l x_{j}-x_{j-1}\r+ \sum_{i=j+1}^{k} f\left(\xi_i\right)\cdot \Delta
x_i+ \left[ f\left(\xi_{j}\right)-f(b)\right] \cdot \Delta x_{j}
 =\\=
 \underbrace{\sum_{i=1}^{j-1} f\left(\xi_i
 \right)\cdot \Delta x_i+ f(b) \cdot \l
 b-x_{j-1}\r}_{\text{интегральная сумма на отрезке $[a;b]$}}
 +
 \underbrace{f(b)\cdot
 \l x_{j}-b \r+ \sum_{i=j+1}^{k} f\left(
 \xi_i\right)\cdot \Delta x_i}_{\text{интегральная сумма на
 отрезке $[b;c]$}}
  +
 \left[ f\left(\xi_{j}\right)-f(b)\right] \cdot
 \Delta x_{j}
 \end{multline*}
$$
\Downarrow
$$
 \begin{multline*}
\Bigg|\overbrace{\sum_{i=1}^{k} f\left(\xi_i\right)\cdot \Delta
x_i}^{\scriptsize\begin{matrix}\text{интегральная сумма}\\ \text{на отрезке
$[a;c]$}\end{matrix}}- \underbrace{\overbrace{\l\sum_{i=1}^{j-1} f\left(\xi_i
 \right)\cdot \Delta x_i+ f(b) \cdot (
 b-x_{j-1})\r}^{\text{интегральная сумма на отрезке $[a;b]$}}}_{\scriptsize\begin{matrix}
 \phantom{\tiny\begin{matrix}\diam\tau\\ \downarrow\\ 0\end{matrix}}\
 \downarrow \ {\tiny\begin{matrix}\diam\tau\\ \downarrow\\ 0\end{matrix}}
 \\ A\end{matrix}}
 -
 \underbrace{\overbrace{\l f(b)\cdot
 \l x_{j}-b \r+ \sum_{i=j+1}^{k} f\left(
 \xi_i\right)\cdot \Delta x_i\r}^{\text{интегральная сумма на
 отрезке $[b;c]$}}}_{\scriptsize\begin{matrix}
 \phantom{\tiny\begin{matrix}\diam\tau\\ \downarrow\\ 0\end{matrix}}\
 \downarrow \ {\tiny\begin{matrix}\diam\tau\\ \downarrow\\ 0\end{matrix}}
 \\ B\end{matrix}}\Bigg|
 =\\=
 \Big| f\left(\xi_{j}\right)-f(b)\Big| \cdot
 \Delta x_j\le
 \Big(\underbrace{\big| f(\xi_j)\big|}_{\scriptsize\begin{matrix}\text{\rotatebox{90}{$\ge$}}
 \\ M\end{matrix}}+\underbrace{\big| f(b)\big|}_{\scriptsize\begin{matrix}\text{\rotatebox{90}{$\ge$}}
 \\ M\end{matrix}}\Big) \cdot
 \underbrace{\Delta x_j}_{\scriptsize\begin{matrix}\text{\rotatebox{90}{$\ge$}}
 \\ \diam\tau\end{matrix}}\le 2M\cdot\diam\
 \tau\underset{\diam\tau\to0}{\longrightarrow}0
 \end{multline*}
$$
\Downarrow
$$
$$
\sum_{i=1}^{k} f\left(\xi_i\right)\cdot \Delta
x_i\underset{\diam\tau\to0}{\longrightarrow} A+B
$$
Опять же, это нужно понимать так, что какую ни возьми последовательность
разбиений $\tau^{(n)}$ отрезка $[a,c]$, с диаметрами стремящимися к нулю, то
при любом выборе выделенных точек $\xi_i^{(n)}$, интегральные суммы функции $f$
на отрезке $[a,c]$ будут стремиться к числу $A+B$. Это нам и нужно было
доказать.

3. Монотонность.  Пусть функции $f$ и $g$ интегрируемы на отрезке $[a;b]$,
причем
$$
  f(x)\le g(x),  \qquad x\in [a;b]
$$
Тогда для любых разбиений $\tau=\{x_0,...,x_k\}$ отрезка $[a;c]$ и любых точек
$\xi_i\in [x_{i-1};x_i]$ мы получим
$$
\forall i \qquad f\left(\xi_i\right) \le g\left(\xi_i\right)
$$
$$
\Downarrow
$$
$$
\forall i \qquad f\left(\xi_i\right)\cdot \Delta x_i \le
g\left(\xi_i\right)\cdot \Delta x_i
$$
$$
\Downarrow
$$
$$
\int_{[a,b]} f(x) \, \d
x\underset{0\gets\diam\tau}{\longleftarrow}\sum_{i=1}^{k}
f\left(\xi_i\right)\cdot \Delta x_i \le \sum_{i=1}^{k} g\left(\xi_i\right)\cdot
\Delta x_i\underset{\diam\tau\to 0}{\longrightarrow}\int_{[a,b]} g(x) \, \d x
$$
$$
\Downarrow
$$
$$
\int_{[a,b]} f(x) \, \d x\le\int_{[a,b]} g(x) \, \d x
$$

4. Выпуклость. Здесь при доказательстве используются свойства точных граней
функции, о которых мы говорили на с.\pageref{SEC-tochnaya-gran-funktsii}.
Сначала нужно доказать следующее неравенство, связывающее верхние и нижние
суммы Дарбу для $f$ и $|f|$:
 \begin{equation}
 S_{\tau}(|f|)-s_{\tau}(|f|)\le S_{\tau}(f)-s_{\tau}(f)
 \label{S_tau(|f|)}
 \end{equation}
Действительно:
\begin{multline*}
 S_{\tau}(|f|) -s_{\tau}(|f|)=
 \sum_{i=1}^{k}
 \left(
 \sup_{\eta\in
 [x_{i-1};x_i]}|f(\eta )|- \inf_{\theta\in [x_{i-1};x_i]}|f(\theta)|
 \right) \cdot \Delta x_i
 =\eqref{antiodnorodnost-inf-i-sup}=\\=
 \sum_{i=1}^{k}
 \left(
 \sup_{\eta\in
 [x_{i-1};x_i]}|f(\eta )|+\sup_{\theta\in
 [x_{i-1};x_i]}\Big(-|f(\theta)|\Big)
 \right) \cdot \Delta x_i
 =
 \sum_{i=1}^{k}
 \sup_{\eta,\theta\in [x_{i-1};x_i]}
 \Big( |f(\eta)|-|f(\theta)| \Big) \cdot \Delta x_i
 \le \\ \le \eqref{module-3^0}
 \le
 \sum_{i=1}^{k}
 \sup_{\eta,\theta\in [x_{i-1};x_i]}
 \left|f(\eta )-f(\theta ) \right|
 \cdot \Delta x_i=
 \eqref{sup_x,y in X |f(x)-f(y)| sup_x,y in X (f(x)-f(y))}
 =\sum_{i=1}^{k}
 \sup_{\eta,\theta\in [x_{i-1};x_i]}\left(f(\eta )-f(\theta ) \right)
 \cdot \Delta x_i
 =\\=
  \sum_{i=1}^{k}
 \left(
 \sup_{\eta\in
 [x_{i-1};x_i]}f(\eta )+\sup_{\theta\in
 [x_{i-1};x_i]}\Big(-f(\theta)\Big)
 \right) \cdot \Delta x_i
 =\eqref{antiodnorodnost-inf-i-sup}=\\=
 \sum_{i=1}^{k}\left(\sup_{\eta\in
 [x_{i-1};x_i]} f(\eta )- \inf_{\theta\in [x_{i-1};x_i]} f(\theta
 ) \right) \cdot \Delta x_i = S_{\tau}(f)-s_{\tau}(f)
  \end{multline*}
Если теперь функция $f$ интегрируема на отрезке $[a;b]$, то по критерию
интегрируемости (теорема \ref{tm-14.5.1} этой главы), для любой измельчающейся
последовательности $\tau^{(n)}$ разбиений отрезка $[a;b]$ разность между
верхними и нижними суммами Дарбу для функции $f$ должна стремиться к нулю,
поэтому из \eqref{S_tau(|f|)} получаем:
$$
0\le S_{\tau^{(n)}}(|f|) -s_{\tau^{(n)}}(|f|)\le
S_{\tau^{(n)}}(f)-s_{\tau^{(n)}}(f) \underset{n\to \infty}{\longrightarrow} 0
$$
откуда
$$
S_{\tau^{(n)}}(|f|) -s_{\tau^{(n)}}(|f|) \underset{n\to
\infty}{\longrightarrow} 0
$$
То есть функция $|f|$ тоже должна быть интегрируема на отрезке $[a;b]$. Нам
остается доказать неравенство \eqref{14.6.5}:
$$
\left|\int_{[a,b]} f(x) \, \d x \right|
\underset{0\gets\diam\tau}{\longleftarrow}\left| \sum_{i=1}^{k}
f\left(\xi_i\right)\cdot \Delta x_i \right| \le\eqref{ind-nerav-s-modulem}\le
\sum_{i=1}^{k}\left|f\left(\xi_i\right) \right|\cdot \Delta
x_i\underset{\diam\tau\to 0}{\longrightarrow} \int_{[a,b]} |f(x)| \, \d x
$$
$$
\Downarrow
$$
$$
\left|\int_{[a,b]} f(x) \, \d x \right|\le\int_{[a,b]} |f(x)| \, \d x
$$

5. Оценка сверху. Если функция $f$ интегрируема на отрезке $[a,b]$, то по
теореме \ref{tm-14.3.1}, она ограничена. Обозначим $C=\sup_{x\in[a,b]}|f(x)|$.
Тогда
$$
\left|\int_{[a,b]}  f(x) \, \d x\right| \le\eqref{14.6.5}\le
\underbrace{\int_{[a,b]} |f(x)| \, \d x\le \int_{[a,b]} C \, \d
x}_{\scriptsize\begin{matrix}\phantom{\eqref{14.6.4}}\ \Uparrow\ \eqref{14.6.4}
\\|f(x)|\le C\end{matrix}}=\eqref{14.2.2}=(b-a)\cdot C
$$

6. Непрерывность. Из двух формул \eqref{nepreryvnost-integrala-1} и
\eqref{nepreryvnost-integrala-2} мы докажем первую, имея в виду, что вторая
доказывается по аналогии:
$$
\int_{[a,y]}  f(x) \, \d x\underset{y\to c}{\longrightarrow} \int_{[a,c]}  f(x)
\, \d x
$$
Ее, в свою очередь можно разбить на два односторонних предела,
$$
\int_{[a,y]}  f(x) \, \d x\underset{y\to c-0}{\longrightarrow} \int_{[a,c]}
f(x) \, \d x,\qquad \int_{[a,y]}  f(x) \, \d x\underset{y\to
c+0}{\longrightarrow} \int_{[a,c]} f(x) \, \d x,
$$
и доказать, например, второй, заявив, что первый доказывается аналогично:
$$
\bigg|\kern-14pt\underbrace{\int_{[a,y]}  f(x) \, \d
x}_{\scriptsize\begin{matrix}\text{\phantom{\eqref{14.6.2}}\ \rotatebox{90}{$=$}}\ \eqref{14.6.2}\\
\int_{[a,c]} f(x) \, \d x+\int_c^y  f(x) \, \d x
\end{matrix}}\kern-14pt-\int_{[a,c]} f(x) \, \d x\bigg|= \left|\int_c^y f(x) \, \d
x\right|\le\eqref{otsenka-integrala}\le
\underbrace{(y-c)}_{\scriptsize\begin{matrix}\downarrow\\
0\end{matrix}}\cdot\kern-25pt\underbrace{\sup_{x\in[c,y]}|f(x)|}_{\scriptsize\begin{matrix}
\text{\rotatebox{90}{$\ge$}}\\
\sup\limits_{x\in[c,y]}|f(x)| \\
\phantom{\text{теорема \ref{tm-14.3.1}}}\ \text{\rotatebox{90}{$>$}}\ \text{теорема \ref{tm-14.3.1}}\\
\infty\end{matrix}}\kern-25pt\underset{y\to c+0}{\longrightarrow} 0
$$

7. Теорема о среднем. Пусть $f$ непрерывна на отрезке $[a;b]$. Тогда по теореме
Вейерштрасса об экстремумах (теорема \ref{Wei-III}), найдутся такие точки
$\alpha, \beta \in [a;b]$, что
$$
\forall x\in [a;b] \qquad f(\alpha)\le f(x)\le f(\beta)
$$
Теперь по уже доказанному свойству монотонности $3^0$,
$$
\int_{[a,b]} f(\alpha) \, \d x \le \int_{[a,b]} f(x) \, \d x\le \int_{[a,b]}
f(\beta) \, \d x
$$
Вспомнив теперь \eqref{14.2.2}, получаем
$$
f(\alpha) \cdot (b-a) \le \int_{[a,b]} f(x) \, \d x \le f(\beta) \cdot (b-a)
$$
откуда
$$
f(\alpha) \le \frac{1}{b-a}\cdot \int_{[a,b]} f(x) \, \d x \le f(\beta)
$$
Таким образом, число $M=\frac{1}{b-a}\cdot \int_{[a,b]} f(x)$ лежит между
значениями $f(\alpha)$ и $f(\beta)$ непрерывной функции $f$ на отрезке
$[\alpha;\beta]$. Значит, по теореме Коши о промежуточном значении (теорема
\ref{Cauchy-I}), найдется точка $\xi \in [\alpha;\beta]\subseteq [a;b]$ такая
что
$$
f(\xi)=\frac{1}{b-a}\cdot \int_{[a,b]} f(x) \, \d x
$$
Отсюда и получается \eqref{14.6.6}:
$$
\int_{[a,b]}  f(x) \, \d x=f(\xi)\cdot (b-a)
$$

8. Произведение интегрируемых функций. Пусть функции $f$ и $g$ интегрируемы на
отрезке $[a;b]$. Тогда по теореме \ref{tm-14.3.1} они должны быть ограничены на
$[a;b]$:
$$
  \exists A, B\in \R\quad \forall x\in [a;b] \quad
|f(x)|\le A,\quad |g(x)|\le B \label{14.6.10}
$$
Заметим, что для любых $\eta,\theta\in [a,b]$ выполняется неравенство
 \begin{equation}\label{14.6.11}
\left| f(\eta)\cdot g(\eta)-f(\theta)\cdot g(\theta) \right|\le B\cdot
|f(\eta)-f(\theta)|+A\cdot |g(\eta)-g(\theta)|
 \end{equation}
Действительно,
 \begin{multline*}
| f(\eta)\cdot g(\eta)-f(\theta)\cdot g(\theta)|= | f(\eta)\cdot
g(\eta)-f(\theta)\cdot g(\eta)+ f(\theta)\cdot g(\eta)-f(\theta)\cdot
g(\theta)|\le \\ \le | f(\eta)\cdot g(\eta)-f(\theta)\cdot g(\eta)|+ |
f(\theta)\cdot g(\eta)-f(\theta)\cdot g(\theta)|= | f(\eta)-f(\theta)|\cdot
|g(\eta)| +| f(\theta)|\cdot |g(\eta)-g(\theta)| \le \\ \le
|f(\eta)-f(\theta)|\cdot B+A\cdot |g(\eta)-g(\theta)|
 \end{multline*}
Теперь, чтобы доказать, что $f\cdot g$ интегрируема на $[a;b]$, возьмем
разбиение $\tau$ отрезка $[a;b]$ и оценим разность между верхними и нижними
суммами Дарбу для функции $f(x)\cdot g(x)$:
\begin{multline*}
S_{\tau}(f\cdot g) -s_{\tau}(f\cdot g)= \sum_{i=1}^{k}\left(\sup_{\eta\in
[x_{i-1};x_i]}f(\eta )\cdot g(\eta )- \inf_{\theta\in [x_{i-1};x_i]}f(\theta
)\cdot g(\theta ) \right) \cdot \Delta x_i
 =\\=
 \sum_{i=1}^{k}
 \sup_{\eta,\theta\in [x_{i-1};x_i]}
 \Big| f(\eta )\cdot g(\eta )-f(\theta)\cdot g(\theta ) \Big|
 \cdot \Delta x_i
 \le \eqref{14.6.11} \le\\ \le
 \sum_{i=1}^{k}
 \sup_{\eta,\theta\in [x_{i-1};x_i]}
 \Big( B\cdot \Big|f(\eta)-f(\theta)\Big|+A\cdot \Big|g(\eta)-g(\theta)\Big| \Big)
 \cdot \Delta x_i =\\=
 B\cdot \sum_{i=1}^{k}
 \sup_{\eta,\theta\in [x_{i-1};x_i]}
 \Big|f(\eta)-f(\theta)\Big| \cdot \Delta x_i +
 A\cdot \sum_{i=1}^{k}
 \sup_{\eta,\theta\in [x_{i-1};x_i]}
 \Big|g(\eta)-g(\theta)\Big| \cdot \Delta x_i
 =\\=
B\cdot \underbrace{\Big( S_{\tau}(f)
-s_{\tau}(f)\Big)}_{\scriptsize\begin{matrix} \phantom{\tiny
\begin{matrix}\diam \tau\\ \downarrow \\  0\end{matrix}}\
\text{\rotatebox{90}{$\longleftarrow$}} \ {\tiny \begin{matrix}\diam \tau\\ \downarrow \\  0\end{matrix}} \\
 \phantom{,}0, \\ \text{поскольку $f$}\\ \text{интегрируема}\end{matrix}}
 +
 A\cdot \underbrace{\Big( S_{\tau}(g) -s_{\tau}(g) \Big)}_{\scriptsize\begin{matrix} \phantom{\tiny
\begin{matrix}\diam \tau\\ \downarrow \\  0\end{matrix}}\
\text{\rotatebox{90}{$\longleftarrow$}} \ {\tiny \begin{matrix}\diam \tau\\ \downarrow \\  0\end{matrix}} \\
 \phantom{,}0, \\ \text{поскольку $g$}\\ \text{интегрируема}\end{matrix}}
 \underset{\diam \tau\to 0}{\longrightarrow} 0 \end{multline*}

9. Замена переменной. Пусть функции $f:[a,b]\to\R$ и
$\ph:[\alpha,\beta]\to[a,b]$ обладают свойствами, описанными на
с.\pageref{tm-15.4.1-1}. Поскольку $f$ интегрируема, по теореме \ref{tm-14.3.1}
должна быть конечна величина
$$
M=\sup_{x\in[a,b]}|f(x)|
$$
Далее, функция $\ph$ строго монотонна, значит она либо возрастает, либо
убывает. Будем считать для начала, что $\ph$ возрастает. Тогда
$$
\ph(\alpha)=a,\qquad \ph(\beta)=b,\qquad \ph'(x)\ge 0,\quad x\in[\alpha,\beta]
$$
Рассмотрим произвольное разбиение $\tau=\{t_0,...,t_i,...,t_k\}$ отрезка
$[\alpha,\beta]$ с выделенной системой точек $\zeta_i\in[t_{i-1},t_i]$. Положив
$x_i=\ph(t_i)$ и $\xi_i=\ph(\zeta_i)$, мы получим разбиение отрезка $[a,b]$,
причем по теореме Лагранжа \ref{Lagrange}, для некоторых точек
$\eta_i\in[t_{i-1},t_i]$ будет выполняться равенство
$$
\Delta
x_i=x_i-x_{i-1}=\ph(t_i)-\ph(t_{i-1})=\eqref{Lagrange-xi}=\ph'(\eta_i)\cdot
(t_i-t_{i-1})=\ph'(\eta_i)\cdot\Delta t_i
$$
Если теперь взять $\e>0$, то по теореме Кантора \ref{Kantor} найдется
$\delta>0$ такое, что для любых $s,t\in[\alpha,\beta]$ выполняется импликация:
$$
|s-t|<\delta\quad\Longrightarrow\quad
|\ph'(s)-\ph'(t)|<\frac{\e}{M\cdot(\beta-\alpha)}
$$
Поэтому если диаметр разбиения $\tau$ меньше $\delta$, мы получим
 \begin{multline*}
\Big|\sum_{i=1}^k
\underbrace{f\big(\ph(\zeta_i)\big)}_{\scriptsize\begin{matrix}\text{\rotatebox{90}{$=$}}\\
f(\xi_i)\end{matrix}}\cdot\underbrace{|\ph'(\zeta_i)|}_{\scriptsize\begin{matrix}\text{\rotatebox{90}{$=$}}\\
\ph'(\zeta_i)\end{matrix}}\cdot\Delta t_i- \sum_{i=1}^k
f(\xi_i)\cdot\kern-10pt\underbrace{\Delta x_i}_{\scriptsize\begin{matrix}\text{\rotatebox{90}{$=$}}\\
\ph'(\eta_i)\cdot\Delta t_i\end{matrix}}\kern-10pt\Big|=\Big|\sum_{i=1}^k
f(\xi_i)\cdot\big( \ph'(\zeta_i)-\ph'(\eta_i)\big)\cdot\Delta t_i\Big|\le \\
\le \sum_{i=1}^k \underbrace{|f(\xi_i)|}_{\scriptsize\begin{matrix}\text{\rotatebox{90}{$\ge$}}\\
\sup\limits_{x\in[a,b]}|f(x)| \\ \text{\rotatebox{90}{$=$}} \\ M
\end{matrix}}\cdot\underbrace{\big|
\ph'(\zeta_i)-\ph'(\eta_i)\big|}_{\scriptsize\begin{matrix}\text{\rotatebox{90}{$>$}}\\
\frac{\e}{M\cdot(\beta-\alpha)} \end{matrix}}\cdot\Delta t_i<M\cdot
\frac{\e}{M\cdot(\beta-\alpha)} \cdot\underbrace{\sum_{i=1}^k \Delta
t_i}_{\scriptsize\begin{matrix} \text{\rotatebox{90}{$=$}} \\ \beta-\alpha
\end{matrix}}=\e
 \end{multline*}
Это означает, что при стремлении к нулю диаметра разбиения $\tau$ отрезка
$[\alpha,\beta]$ интегральные суммы стремятся друг к другу:
$$
\sum_{i=1}^k f\big(\ph(\zeta_i)\big)\cdot\ph'(\zeta_i)\cdot\Delta t_i-
\sum_{i=1}^k f(\xi_i)\cdot \Delta x_i\underset{\diam\tau\to
0}{\longrightarrow}0
$$
С другой стороны, диаметр соответствующего разбиения
$\ph(\tau)=\{x_0,...,x_k\}$ отрезка $[a,b]$ тоже при этом будет стремиться к
нулю, потому что
$$
\diam\ph(\tau)=\max_{i}\Delta x_i=\max_{i}\ph'(\eta_i)\cdot\Delta t_i\le
\max_{t\in[\alpha,\beta]}\ph'(t)\cdot\max_{i}\Delta
t_i=\max_{t\in[\alpha,\beta]}\ph'(t)\cdot\diam \tau\underset{\diam\tau\to
0}{\longrightarrow}0
$$
Значит интегральные суммы $\sum_{i=1}^k f(\xi_i)\cdot \Delta x_i$ должны
стремиться к интегралу $\int_{[a,b]} f(x) \ \d x$, и мы получаем:
$$
\sum_{i=1}^k f\big(\ph(\zeta_i)\big)\cdot\ph'(\zeta_i)\cdot\Delta t_i-
\underbrace{\sum_{i=1}^k f(\xi_i)\cdot \Delta x_i}_{\scriptsize\begin{matrix}
\phantom{\tiny \begin{matrix} \diam \ph(\tau) \\ \downarrow \\ 0
\end{matrix}} \quad \downarrow \quad {\tiny \begin{matrix} \diam \ph(\tau) \\ \downarrow \\
0
\end{matrix}}
\\ \int\limits_{[a,b]} f(x) \ \d x
\end{matrix}}\underset{\diam\tau\to
0}{\longrightarrow}0
$$
То есть
$$
\sum_{i=1}^k f\big(\ph(\zeta_i)\big)\cdot\ph'(\zeta_i)\cdot\Delta
t_i\underset{\diam\tau\to 0}{\longrightarrow} \int_{[a,b]} f(x) \ \d x
$$

Остается рассмотреть случай, когда $\ph$ убывает. Тогда возрастающая
последовательность
$$
\alpha=t_0<...<t_i<...<t_k=\beta
$$
превращается в убывающую
$$
b=x_0>...>x_i>...>x_k=a
$$
Поэтому в интегральной сумме величину $\Delta x_i$ нужно определять как
разность $x_{i-1}-x_i$, чтобы она получилась положительной:
$$
\Delta
x_i=-\Big(x_i-x_{i-1}\Big)=-\Big(\ph(t_i)-\ph(t_{i-1})\Big)=\eqref{Lagrange-xi}=-\ph'(\eta_i)\cdot
(t_i-t_{i-1})=-\ph'(\eta_i)\cdot\Delta t_i
$$
Вдобавок, поскольку $\ph$ убывает, ее производная должна быть отрицательной
$$
\ph'(t)\le 0\qquad \Longrightarrow\quad |\ph'(\zeta_i)|=-\ph'(\zeta_i)
$$
и в вычислениях мы получим
 $$
\Big|\sum_{i=1}^k
\underbrace{f\big(\ph(\zeta_i)\big)}_{\scriptsize\begin{matrix}\text{\rotatebox{90}{$=$}}\\
f(\xi_i)\end{matrix}}\cdot\underbrace{|\ph'(\zeta_i)|}_{\scriptsize\begin{matrix}\text{\rotatebox{90}{$=$}}\\
-\ph'(\zeta_i)\end{matrix}}\cdot\Delta t_i- \sum_{i=1}^k
f(\xi_i)\cdot\kern-10pt\underbrace{\Delta x_i}_{\scriptsize\begin{matrix}\text{\rotatebox{90}{$=$}}\\
-\ph'(\eta_i)\cdot\Delta t_i\end{matrix}}\kern-10pt\Big|=\Big|-\sum_{i=1}^k
f(\xi_i)\cdot\big( \ph'(\zeta_i)-\ph'(\eta_i)\big)\cdot\Delta t_i\Big|\le ...<
\e
 $$
Таким образом, после взятия модуля ничто не меняется, и все рассуждения можно
оставить прежними.
 \end{proof}

\section{Формула Ньютона-Лейбница}\label{CH-Newton-Leibnitz}

Между двумя независимыми на первый взгляд понятиями -- первообразной и
определенным интегралом, -- как оказывается, существует глубокая связь. Эта
связь устанавливается знаменитой формулой Ньютона-Лейбница, которая считается
главной в интегральном исчислении. Здесь мы поговорим об этой формуле.

\subsection{Основные результаты}

\paragraph{Дифференцируемые функции и первообразная на отрезке.}

 \bit{
\item[$\bullet$]
 Пусть функция $h$ определена на отрезке $[a;b]$. Тогда
 \bit{
 \item[---] ее {\it правой производной}\index{производная!правая}
 в точке $a$ называется предел
  \beq\label{h_+'(a)}
 h_+'(a):=\lim_{x\to a+0}\frac{h(x)-h(a)}{x-a}
  \eeq
 \item[---] ее {\it левой производной}\index{производная!левая}
 в точке $b$ называется предел
 \beq\label{h_-'(a)}
 h_-'(b):=\lim_{x\to b-0}\frac{h(x)-h(b)}{x-b}
 \eeq
 }\eit

\item[$\bullet$] Если эти величины конечны, и кроме того функция $h$ имеет
 (обычную) конечную производную в любой точке интервала $(a;b)$, то
 говорят, что $h$ {\it дифференцируема на отрезке}\label{func-diff-na-otrezke}
 \index{функция!дифференцируемая!на отрезке} $[a;b]$, а ее {\it производной на
 отрезке} $[a;b]$ называется функция
 \beq\label{DEF:proizvodnaya-na-otrezke}
 h'(x)=
 \begin{cases}
 h_+'(a), \quad x=a \\
 h'(x), \quad x\in (a;b) \\
 h_-'(a), \quad x=b
 \end{cases}
 \eeq
 }\eit

\bprop\label{PROP:diff=>cont-[a,b]} Всякая дифференцируемая функция на отрезке
$[a,b]$ непрерывна на нем. \eprop

\bpr Если функция $h$ дифференцируема на отрезке $[a;b]$, то по теореме
\ref{cont-diff}, она будет непрерывна в каждой внутренней точке этого отрезка.
Остается убедиться, что на концах она тоже непрерывна. Это делается тем же
приемом, что и в теореме \ref{cont-diff}: если пределы \eqref{h_+'(a)} и
\eqref{h_-'(a)} существуют и конечны, то функция $h$ непрерывна в точках $a$ и
$b$.
 \epr

\bit{ \item[$\bullet$] Функция $F$ называется {\it
первообразной}\label{DEF:pervoobraznaya}\index{первообразная!на отрезке} для
функции $f$ {\it на отрезке} $[a;b]$, если $F$ определена на $[a;b]$,
дифференцируема на нем и ее производная на нем равна $f$:
 \beq\label{15.1.4}
F'(x)=f(x), \qquad x\in [a;b]
 \eeq
При этом, в соответствии с определением производной на отрезке формулой
\eqref{DEF:proizvodnaya-na-otrezke}, на концах отрезка формула \eqref{15.1.4}
интерпретируется, как условие на односторонние производные:
 \beq\label{15.1.4+}
F'_+(a)=f(a),\qquad F'_-(b)=f(b)
 \eeq
 }\eit

\begin{lm}\label{lm-12.1.3} Пусть функция $H$ определена и дифференцируема на
отрезке $[a,b]$, причем ее производная равна нулю на $[a,b]$:
 \begin{equation}
H'(x)=0, \quad x\in [a,b] \label{12.1.1}
 \end{equation}
Тогда функция $H(x)$ постоянна на $[a,b]$:
$$
H(x)=C, \quad x\in [a,b] \label{12.1.2}
$$
для некоторого $C\in \R$.
\end{lm}\begin{proof} Положим $C=H(a)$. Для любого $x\in(a,b]$ мы получим, что
функция $H$
 \bit{
\item[--] дифференцируема на интервале $(a;x)$ (потому что она дифференцируема
на $(a,b)$), и

\item[--] непрерывна на отрезке $[a;x]$ (по предложению
\ref{PROP:diff=>cont-[a,b]})
 }\eit
поэтому, мы можем применить к функции $H$ на отрезке $[a;x]$ теорему Лагранжа
\ref{Lagrange}, и тогда получим, что существует такая точка $\xi \in [a;x]$,
что
$$
\frac{H(x)-H(a)}{x-a}=H'(\xi)= \eqref{12.1.1}=0
$$
отсюда
$$
H(x)-H(a)=0\quad\Longrightarrow\quad H(x)=H(a)=C
$$
 \end{proof}

\medskip

\centerline{\bf Свойства первообразных на
отрезке}\index{свойства!первообразных}

{\it

 \bit{

\item[$1^0.$]\label{svoistvo-pervoobraznyh-1} Если $F$ -- первообразная для
некоторой функции $f$ на отрезке $[a;b]$, то $F$ является непрерывной функцией
на $[a;b]$.

\item[$2^0.$]\label{svoistvo-pervoobraznyh-2} Если $F$ -- какая-нибудь
первообразная для функции $f$ на отрезке $[a;b]$, то все первообразные для
функции $f$ на отрезке $[a;b]$ имеют вид
$$
F(x)+C
$$
где $C$ -- произвольные константы.
 }\eit
}

\begin{proof}

1. Если $F$ -- первообразная для $f$ на $[a;b]$, то $F$ дифференцируема на
$[a;b]$ (и ее производной будет $f$), поэтому в силу предложения
\ref{PROP:diff=>cont-[a,b]}, $F$ непрерывна на $[a;b]$.

2. Пусть $F$ и $G$ -- две первообразные для $f$ на отрезке $[a;b]$:
$$
F'(x)=G'(x)=f(x), \quad x\in [a,b]
$$
Тогда их разность $H=G-F$ будет иметь нулевую производную на $[a;b]$
$$
H'(x)=G'(x)-F'(x)=f(x)-f(x)=0, \quad x\in [a;b]
$$
Значит, по лемме \ref{lm-12.1.3}, функция $H$ постоянна:
$$
H(x)=G(x)-F(x)=C, \quad x\in [a;b]
$$
(где $C$ -- некоторая константа). Следовательно,
$$
G(x)=F(x)+C, \quad x\in [a;b]
$$
У нас получилось, что если $F$ -- какая-нибудь первообразная, то любая другая
первообразная $G$ имеет вид $G(x)=F(x)+C$. \end{proof}

\noindent\rule{160mm}{0.1pt}\begin{multicols}{2}

\brem Из свойства $2^0$ можно сделать вывод, что если функция $f$ определяется
(одноместным) числовым выражением от переменной $x$ без параметров (или с $m$
параметрами)
$$
f(x)={\mathcal P},
$$
то ее неопределенный интеграл можно описать одноместным числовым выражением от
$x$ с одним параметром (соответственно, с $m+1$ параметрами). В качестве такого
нового параметра обычно берется $C$, но, разумеется, вместо $C$ можно брать
любую другую букву, кроме $x$. Этот выражение обозначается символом
$$
\int {\mathcal P} \d x
$$
и определяется правилом
 \begin{multline}\label{opred-int-f(x)-dx}
\int {\mathcal P} \d x={\mathcal Q}\qquad\Longleftrightarrow\qquad \d {\mathcal
Q}={\mathcal P}\cdot\d x \qquad\Longleftrightarrow
\\ \Longleftrightarrow\qquad \frac{\d}{\d
x}{\mathcal Q}={\mathcal P}
 \end{multline}
--- здесь важно помнить, что, в соответствии с правилами
дифференцирования выражений на с.\pageref{subsec-diff-chisl-terma}, дифференциал и
производная от параметра $C$ по $x$ всегда считается равной нулю (формула
\eqref{d/dx-bez-x}):
 \beq\label{dC/dx=0}
\d C=0,\qquad \frac{\d C}{\d x}=0
 \eeq
К этому определению следует еще добавить, что если функция $f$ определяется
каким-то выражением $\mathcal P$ с параметром (или параметрами), например,
$$
f(x)=x^\alpha,
$$
то, разумеется, при интегрировании в качестве нового параметра нужно выбирать
какую-нибудь букву из тех, что не содержатся в $\mathcal P$ (то есть в данном
случае не $x$ и не $\alpha$). \erem

\end{multicols}\noindent\rule[10pt]{160mm}{0.1pt}

\paragraph{Гладкие функции и интеграл с переменным верхним пределом на отрезке.}

 \bit{
\item[$\bullet$] Функция $h$ называется {\it гладкой}\index{функция!гладкая!на
отрезке} {\it на отрезке} $[a;b]$, если она дифференцируема на $[a;b]$ (в
смысле определения на с.\pageref{func-diff-na-otrezke}), и ее производная $h'$
на этом отрезке (определенная формулой \eqref{DEF:proizvodnaya-na-otrezke})
непрерывна на $[a;b]$. Множество всех гладких функций на $[a;b]$ обозначается
символом ${\mathcal C}^1[a,b]$. Из предложения \ref{PROP:diff=>cont-[a,b]}
следует, что любая гладкая функция автоматически непрерывна, то есть
выполняется включение:
 \beq\label{C^1[a,b]-subset-C[a,b]}
{\mathcal C}^1[a,b]\subset {\mathcal C}[a,b]
 \eeq
 }\eit

\begin{tm}[\bf об интеграле с переменным верхним пределом]\label{tm-15.1.1}
Пусть функция $f$ непрерывна на отрезке $[a;b]$, тогда функция
 \beq\label{15.1.6}
  F(x)=\int_{[a,x]} f(t) \, \d t
 \eeq
является гладкой на $[a;b]$ и будет первообразной для $f$ на $[a;b]$:
 \beq\label{15.1.6-1}
F'(c)=f(c),\qquad c\in[a,b].
 \eeq
 \end{tm}
\brem В соответствии с определением производной на отрезке формулой
\eqref{DEF:proizvodnaya-na-otrezke}, соотношение \eqref{15.1.6-1} следует
понимать так:
 \bit{
\item[---] для точек внутри отрезка $[a;b]$ эта формула расшифровывается как
равенство
 \beq\label{15.1.7}
  F'(c)=\lim_{x\to c}\frac{1}{x-c}\l\int_{[a,x]} f(t) \, \d t-\int_{[a,c]} f(t) \, \d
  t\r=f(c), \quad c\in (a;b)
 \eeq
\item[---] а на концах отрезка $[a;b]$ имеются в виду односторонние
производные:
 \beq\label{15.1.7+}
  F'_+(a)=\lim_{x\to a+0}\frac{1}{x-a}\int_{[a,x]} f(t)\ \d t=f(a),
 \eeq
 \beq\label{15.1.7-}
  F'_-(b)=\lim_{x\to b-0}\frac{1}{x-b}\l \int_{[a,x]} f(t)\ \d t-\int_{[a,b]} f(t)\ \d t\r=
  \lim_{x\to b-0}\frac{1}{b-x}\int_x^b f(t)\ \d t=f(b).
 \eeq
 }\eit
 \erem
\begin{proof} Заметим сразу, что здесь достаточно доказать тождество \eqref{15.1.6-1}, потому что остальное
-- то есть утверждение, что функция $F$ является гладкой -- будет следовать из
\eqref{15.1.6-1} и того факта, что $f$ непрерывна. В свою очередь
\eqref{15.1.6-1} разбивается на условия \eqref{15.1.7}-\eqref{15.1.7-}, и
каждое нужно проверить. Зафиксируем для этого какую-нибудь точку $c\in [a;b]$ и
найдем в ней производную функции $F$. Нам придется рассмотреть несколько
случаев.

1. Пусть для начала точка $c$ лежит внутри отрезка $[a;b]$, то есть $a<c<b$.
Чтобы найти производную функции $F(x)$ в точке $c$,
$$
F'(c)=\lim_{x\to c}\frac{F(x)-F(c)}{x-c},
$$
вычислим по отдельности правый и левый пределы в этой точке.

 \bit{
\item[a)] Пусть $x>c$, то есть $a<c<x<b$. Тогда
\begin{multline*}\frac{F(x)-F(c)}{x-c}= \frac{1}{x-c}\left\{ \int_{[a,x]} f(t) \, \d
t- \int_{[a,c]} f(t) \, \d t \right\}=\\= \frac{1}{x-c}\left\{ \int_{[a,c]}
f(t) \, \d t+ \int_{[c,x]} f(t) \, \d t- \int_{[a,c]} f(t) \, \d t \right\}=\\=
\frac{1}{x-c}\cdot \underbrace{\int_{[c,x]} f(t) \, \d
t}_{\scriptsize\begin{matrix}\phantom{\tiny \eqref{15.4.1-1}}\quad
\text{\rotatebox{90}{$=$}}\quad {\tiny \eqref{15.4.1-1}} \\
f(\xi)\cdot (x-c),
\\ \xi\in[c,x] \end{matrix}} =  \frac{1}{x-c}\cdot
f(\xi)\cdot (x-c)  =\underbrace{f(\xi) \underset{x\to c+0}{\longrightarrow}
f(c)}_{\scriptsize\begin{matrix}\Uparrow \\ \xi\underset{x\to c+0}{\longrightarrow} c\\ \Uparrow \\
\xi\in[c,x]\end{matrix}}
\end{multline*}
Таким образом,
 \beq\label{15.1.8}
\lim_{x\to c+0}\frac{F(x)-F(c)}{x-c}=f(c)
 \eeq

\item[b)] Аналогично, если $x<c$, то есть $a<x<c<b$, то
 \begin{multline*}
\frac{F(x)-F(c)}{x-c}= \frac{1}{x-c}\left\{ \int_{[a,x]} f(t) \, \d t-
\int_{[a,c]} f(t) \, \d t \right\}=\\= \frac{1}{x-c}\left\{ \int_{[a,x]} f(t)
\, \d t- \left(\int_{[a,x]} f(t) \, \d t+\int_{[x,c]} f(t) \, \d t \right)
\right\}=\\= -\frac{1}{x-c}\cdot \underbrace{\int_{[x,c]} f(t) \, \d
t}_{\scriptsize\begin{matrix}\phantom{\tiny \eqref{15.4.1-1}}\quad
\text{\rotatebox{90}{$=$}}\quad {\tiny \eqref{15.4.1-1}} \\
f(\xi)\cdot (c-x),
\\ \xi\in[x,c] \end{matrix}} =  -\frac{1}{x-c}\cdot
f(\xi)\cdot (c-x)  = f(\xi)\underbrace{f(\xi) \underset{x\to
c-0}{\longrightarrow}
f(c)}_{\scriptsize\begin{matrix}\Uparrow \\ \xi\underset{x\to c-0}{\longrightarrow} c\\ \Uparrow \\
\xi\in[x,c]\end{matrix}}
\end{multline*}
Таким образом,
 \beq\label{15.1.9}
\lim_{x\to c-0}\frac{F(x)-F(c)}{x-c}=f(c)
 \eeq
 }\eit

Формулы \eqref{15.1.8} и \eqref{15.1.9} вместе дают
$$
F'(c)=\lim_{x\to c}\frac{F(x)-F(c)}{x-c}=f(c), \qquad c\in (a;b)
$$

2. Теперь рассмотрим случай, когда $c=a$. Тогда аналогично пункту a) получаем
равенство
$$
\lim_{x\to a+0}\frac{F(x)-F(a)}{x-a}=f(a)
$$

3. Остается случай, когда $c=b$, который нужно рассмотреть по аналогии с
пунктом b), и тогда получится формула
$$
\lim_{x\to b-0}\frac{F(x)-F(b)}{x-b}=f(b)
$$
\end{proof}

\paragraph{Формула Ньютона-Лейбница.}

\begin{tm}[\bf Ньютона-Лейбница]\label{tm-15.3.1}
Пусть функция $f$ непрерывна на отрезке $[a;b]$ тогда для любой ее
первообразной $\varPhi$ на этом отрезке выполняется равенство
 \beq\label{15.3.1}
 \int_{[a,b]} f(x) \, \d x=\varPhi (b)-\varPhi (a)
 \eeq
\end{tm}
 \bit{
\item[$\bullet$] Формула \eqref{15.3.1} называется {\it формулой
Ньютона-Лейбница}\index{формула!Ньютона-Лейбница}.
 }\eit\begin{proof} Рассмотрим функцию
$$
 F(x)=\int_{[a,x]} f(t) \, \d t
$$
По теореме \ref{tm-15.1.1}, она является первообразной для функции $f$ на
отрезке $[a;b]$. Значит, $F$ и $\varPhi$ -- две первообразные для одной и той
же функции $f$ на отрезке $[a;b]$. Поэтому, по свойству $2^0$ на
с.\pageref{svoistvo-pervoobraznyh-2}, эти функции отличаются друг от друга на
какую-то константу:
$$
F(x)=\varPhi(x)+C, \qquad x\in [a;b]
$$
То есть,
 \beq\label{15.3.2}
\int_{[a,x]} f(t) \, \d t=\varPhi(x)+C, \qquad x\in [a;b]
 \eeq
 Взяв $x=a$, мы получим
$$
0=\int_{[a,a]} f(t) \, \d t=\varPhi(a)+C,
$$
откуда
$$
C=-\varPhi(a),
$$
поэтому, подставив это в \eqref{15.3.2}, получаем
$$
\int_{[a,x]} f(t) \, \d t=\varPhi(x)-\varPhi(a), \qquad x\in [a;b]
$$
и, если сюда подставить $x=b$, то получится как раз \eqref{15.3.1}.
\end{proof}

\subsection{Интеграл по ориентированному отрезку и вычисления}

\paragraph{Ориентированные отрезки в $\R$.}

 \bit{
\item[$\bullet$] {\it Ориентированным отрезком} на прямой $\R$ называется
отрезок $[x,y]$, относительно крайних точек которого $x$ и $y$ имеется соглашение, какая из них
считается {\it началом}, а какая {\it концом} (при этом, необязательно $x$ должна считаться началом, а $y$ концом, возможно и наоборот).
Частный случай, когда $x=y$, также считается ориентированным отрезком, в
котором конец и начало совпадают, и такой ориентированный отрезок называется
{\it вырожденным}. Запись $\overrightarrow{ab}$, или $\overrightarrow{a;b}$,
означает ориентированный отрезок с началом $a$ и концом $b$. {\it Носителем}
$[\overrightarrow{ab}]$ ориентированного отрезка $\overrightarrow{ab}$
называется тот же самый отрезок, но уже обычный, неориентированный (то есть такой, у которого мы ``забыли'', какая из его крайних точек была началом, а какая концом).
Понятно, что
$$
[\overrightarrow{ab}]=\begin{cases}[a,b], & a<b \\
\{a\}, & a=b \\ [b,a], & a>b
\end{cases}
$$

\item[$\bullet$] {\it Суммой ориентированных отрезков} $\overrightarrow{ab}$ и
$\overrightarrow{bc}$ называется ориентированный отрезок $\overrightarrow{ac}$:
 \beq\label{summa-orient-otrezkov}
\overrightarrow{ab}+\overrightarrow{bc}:=\overrightarrow{ac}
 \eeq
(подразумевается, что сумма определена только если конец предыдущего отрезка
совпадает с началом следующего).

\item[$\bullet$] {\it Противоположным отрезком} к отрезку $\overrightarrow{ab}$
называется отрезок $\overrightarrow{ba}$, и обозначается он символом
$-\overrightarrow{ab}$:
 \beq
-\overrightarrow{ab}:=\overrightarrow{ba}
 \eeq
Сумму вида
$\overrightarrow{ac}+(-\overrightarrow{bc})=\overrightarrow{ac}+\overrightarrow{cb}=\overrightarrow{ab}$
принято записывать как $\overrightarrow{ac}-\overrightarrow{bc}$ и называть
{\it разностью ориентированных отрезков} $\overrightarrow{ac}$ и
$\overrightarrow{bc}$:
 \beq
\overrightarrow{ac}-\overrightarrow{bc}:=\overrightarrow{ac}+(-\overrightarrow{bc})=\overrightarrow{ac}+\overrightarrow{cb}=\overrightarrow{ab}
 \eeq

\item[$\bullet$] Говорят, что функция $\ph$ {\it является преобразованием
ориентированного отрезка $\overrightarrow{\alpha\beta}$ в ориентированный
отрезок $\overrightarrow{ab}$}, и обозначают это записью
 \beq\label{DEF:preobr-orient-otrezkov}
\ph:\overrightarrow{\alpha\beta}\to\overrightarrow{ab},
 \eeq
если
 \bit{
\item[(i)] $\ph$ гладко отображает носитель $\overrightarrow{\alpha\beta}$ в
носитель $\overrightarrow{ab}$:
 \beq\label{DEF:preobr-orient-otrezkov-1}
\ph:[\overrightarrow{\alpha\beta}]\to [\overrightarrow{ab}]
 \eeq
(это, в частности, означает, что для всякого $t\in
[\overrightarrow{\alpha\beta}]$ значение $\ph(t)$ лежит в
$[\overrightarrow{ab}]$);

\item[(ii)] $\ph$ сохраняет ориентацию в следующем смысле:
 \beq\label{ph(alpha)=a}
\ph(\alpha)=a, \qquad \ph(\beta)=b
 \eeq
 }\eit

 }\eit

\noindent\rule{160mm}{0.1pt}\begin{multicols}{2}

\bex Запись $\overrightarrow{0;1}$ означает отрезок $[0;1]$, в котором 0
считается началом, а 1 концом. Наоборот, $\overrightarrow{1;0}$ означает
отрезок $[0,1]$, но в котором 0 считается концом, а 1 началом. Носители у этих
отрезков одинаковы
$$
[\overrightarrow{0;1}]=[0;1]=[\overrightarrow{1;0}]
$$
но ориентация у них противоположная, и поэтому
$$
\overrightarrow{1;0}=-\overrightarrow{0;1}
$$
\eex

\bex Очевидно, справедливы равенства:
$$
\overrightarrow{0;1}+\overrightarrow{1;2}=\overrightarrow{0;2}
$$
$$
\overrightarrow{0;2}-\overrightarrow{1;2}=\overrightarrow{0;1}
$$
$$
\overrightarrow{0;1}-\overrightarrow{2;1}=\overrightarrow{0;2}
$$
\eex

\bex Функция
$$
\ph(t)=\cos t
$$
является преобразованием ориентированного отрезка $\overrightarrow{0;\pi}$ в
ориентированный отрезок $\overrightarrow{1;-1}$:
$$
\cos:\overrightarrow{0;\pi}\to\overrightarrow{1;-1}
$$
потому что она гладко отображает носитель первого отрезка в носитель второго
$$
\cos:[\overrightarrow{0;\pi}]=[0,\pi]\to[-1;1]=[\overrightarrow{1;-1}]
$$
и сохраняет ориентацию:
$$
\cos 0=1,\qquad \cos\pi=-1.
$$
\eex

\bex Но та же функция
$$
\ph(t)=\cos t
$$
не будет преобразованием отрезка $\overrightarrow{0;\frac{3\pi}{2}}$ в отрезок
$\overrightarrow{1;0}$, хотя она по-прежнему гладкая, и сохраняет ориентацию
этих отрезков:
$$
\cos 0=1,\qquad \cos\frac{3\pi}{2}=0
$$
Дело в том, что эта функция не будет отображением носителя первого отрезка в
носитель второго:
$$
\cos:
\left[\overrightarrow{0;\frac{3\pi}{2}}\right]=\left[0,\frac{3\pi}{2}\right]\not\kern-1pt\to
[0;1]=[\overrightarrow{1;0}]
$$
потому что, например, в точке $t=\pi\in \left[0,\frac{3\pi}{2}\right]$ ее
значения вылезают из отрезка $[0;1]$:
$$
\cos\pi=-1\notin[0;1].
$$
  \eex

\end{multicols}\noindent\rule[10pt]{160mm}{0.1pt}

\paragraph{Интеграл по ориентированному отрезку.}

 \bit{
\item[$\bullet$] Говорят, что числовая функция $f$ определена (интегрируема,
непрерывна, гладка) на ориентированном отрезке $\overrightarrow{ab}$, если она
определена (интегрируема, непрерывна, гладка) на его носителе
$[\overrightarrow{ab}]$. При этом

 \bit{
\item[---] {\it перемещением функции $f$ на ориентированном отрезке}
$\overrightarrow{ab}$ называется число
 \beq\label{15.3.3}
f(x) \Big|_{x=a}^{x=b}:=f(b)-f(a)
 \eeq

\item[---] {\it интегралом функции $f$ по ориентированному отрезку}
$\overrightarrow{ab}$, или {\it анизотропным интегралом} по
$\overrightarrow{ab}$, называется число
 \beq\label{int-s-razdelennymi-predelami}
\int_a^b f(x)\ \d x:=\begin{cases}\int_{[a,b]}f(x)\ \d x, & a<b \\ 0, & a=b \\
-\int_{[b,a]}f(x)\ \d x, & b<a
\end{cases}
 \eeq
 }\eit
 }\eit

\bigskip

\centerline{\bf Свойства интеграла по ориентированному отрезку:}

 \bit{\it

\item[$1^0.$] {\bf Линейность:} интеграл от линейной комбинации функций равен
линейной комбинации интегралов
 \begin{equation}\label{linein-aniz-int-v-R^1}
\int_a^b \Big( \alpha\cdot f(x)+\beta\cdot g(x) \Big)\, \d x= \alpha\cdot
\int_a^b f(x)\, \d x+\beta\cdot \int_a^b g(x)\, \d x
 \end{equation}

\item[$2^0.$] {\bf Аддитивность:} интеграл по сумме ориентированных отрезков
равен сумме интегралов по этим отрезкам:
 \begin{equation}\label{addit-aniz-int-v-R^1}
\int_a^b f(x)\, \d x+\int_b^c f(x)\, \d x= \int_a^c  f(x) \, \d x
 \end{equation}
 а интеграл по разности равен разности интегралов:
 \begin{equation}\label{raznost-aniz-int-v-R^1}
\int_a^c  f(x) \, \d x-\int_a^b f(x)\, \d x=\int_b^c f(x)\, \d x
 \end{equation}

\item[$3^0.$] {\bf Теорема о среднем:} если функция $f$ непрерывна на отрезке
$\overrightarrow{ab}$, то найдется такая точка $\xi\in\overrightarrow{ab}$, что
 \begin{equation}\label{teor-o-srednem-v-aniz-int}
\int_a^b  f(x) \, \d x=f(\xi)\cdot (b-a)
 \end{equation}

\item[$4^\circ$] {\bf Формула Ньютона-Лейбница:} интеграл от непрерывной
функции по ориентированному отрезку равен перемещению какой-нибудь ее
первообразной на этом отрезке:
 \beq\label{Newton-Leibnitz-v-aniz-integr-R^1}
 \int_a^b f(x) \, \d x=\varPhi(x)\Big|_{x=a}^{x=b}
 \eeq

\item[$5^\circ$] {\bf Формула замены переменной:} \label{tm-15.4.1} пусть даны
 \bit{
\item[(i)] преобразование $\ph$ ориентированного отрезка
$\overrightarrow{\alpha\beta}$ в ориентированный отрезок $\overrightarrow{ab}$
(удовлетворяющее условиям \eqref{DEF:preobr-orient-otrezkov-1} и
\eqref{ph(alpha)=a}), и

\item[(ii)] функция $f$, непрерывная на ориентированном отрезке
$\overrightarrow{ab}$;
 }\eit
тогда
 \beq\label{15.4.1}
\int_a^b f (x) \, \d x= \int_\alpha^\beta f (\ph(t)) \, \d \ph(t)=
\int_\alpha^\beta f (\ph(t))\cdot \ph'(t) \, \d t
 \eeq

 }\eit
\bpr

1. Первое свойство следует из \eqref{14.6.1}, и для его доказательства нужно
просто рассмотреть два случая взаимного расположения точек $a$ и $b$ на прямой.
Если $a<b$, то
 \begin{multline*}
\int_a^b \Big( \alpha\cdot f(x)+\beta\cdot g(x)\Big) \, \d
x=\eqref{int-s-razdelennymi-predelami}=\int_{[a,b]} \Big( \alpha\cdot
f(x)+\beta\cdot g(x)\Big) \, \d x=\eqref{14.6.1}=\\=\alpha\cdot\int_{[a,b]}
f(x)\, \d x+\beta\cdot\int_{[a,b]} g(x)\, \d
x=\eqref{int-s-razdelennymi-predelami}= \alpha\cdot\int_a^b f(x)\, \d
x+\beta\cdot\int_a^b g(x)\, \d x
 \end{multline*}
А если $a>b$, то
 \begin{multline*}
\int_a^b \Big( \alpha\cdot f(x)+\beta\cdot g(x)\Big) \, \d
x=\eqref{int-s-razdelennymi-predelami}=-\int_{[b,a]} \Big( \alpha\cdot
f(x)+\beta\cdot g(x)\Big) \, \d x=\eqref{14.6.1}=\\=-\alpha\cdot\int_{[b,a]}
f(x)\, \d x-\beta\cdot\int_{[b,a]} g(x)\, \d
x=\eqref{int-s-razdelennymi-predelami}= \alpha\cdot\int_a^b f(x)\, \d
x+\beta\cdot\int_a^b g(x)\, \d x
 \end{multline*}

2. Равенство \eqref{addit-aniz-int-v-R^1} эквивалентно
\eqref{raznost-aniz-int-v-R^1} и следует из \eqref{14.6.2}. Оно доказывается
тем же приемом: нужно рассмотреть несколько вариантов расположения точек $a$,
$b$ и $c$. В простейшем случае, если $a<b<c$, мы получаем:
 \begin{multline*}
\int_a^b f(x)\, \d x+\int_b^c f(x)\, \d x=\eqref{int-s-razdelennymi-predelami}=
\int_{[a,b]} f(x)\, \d x+\int_{[b,c]} f(x)\, \d x=\eqref{14.6.2}=\\=
\int_{[a,c]} f(x) \, \d x=\eqref{int-s-razdelennymi-predelami}=\int_a^c f(x) \,
\d x
 \end{multline*}
В случае, если, например, $a<c<b$, мы получаем:
 \begin{multline*}
\int_a^b f(x)\, \d x+\int_b^c f(x)\, \d x=\eqref{int-s-razdelennymi-predelami}=
\int_{[a,b]} f(x)\, \d x-\int_{[c,b]} f(x)\, \d x=\eqref{14.6.2}=\\=
\int_{[a,c]} f(x) \, \d x=\eqref{int-s-razdelennymi-predelami}=\int_a^c f(x) \,
\d x
 \end{multline*}
И в таком духе. Остальные варианты мы предлагаем читателю проверить
самостоятельно.

3. Свойство $3^\circ$ следует из \eqref{14.6.6}: если $a<b$, то
$$
\int_a^b  f(x) \, \d x=\int_{[a,b]}  f(x) \, \d x=\eqref{14.6.6}= f(\xi)\cdot
(b-a)
$$
Если же $a>b$, то
$$
\int_a^b  f(x) \, \d x=-\int_{[b,a]}  f(x) \, \d x=\eqref{14.6.6}=-f(\xi)\cdot
(a-b)=f(\xi)\cdot (b-a)
$$

4. Свойство $4^\circ$ следует из \eqref{15.3.1}: если $a<b$, то
$$
\int_a^b f(x) \, \d x=\eqref{int-s-razdelennymi-predelami}=\int_{[a,b]} f(x) \,
\d x=\eqref{15.3.1}=\varPhi(b)-\varPhi(a)=\eqref{15.3.3}=\varPhi (x)
\Big|_{x=a}^{x=b}
$$
Если же $a>b$, то
 \begin{multline*}
\int_a^b f(x) \, \d x=\eqref{int-s-razdelennymi-predelami}=-\int_{[b,a]} f(x)
\, \d
x=\eqref{15.3.1}=\\=-\Big(\varPhi(a)-\varPhi(b)\Big)=\varPhi(b)-\varPhi(a)=\eqref{15.3.3}=\varPhi
(x) \Big|_{x=a}^{x=b}
 \end{multline*}

5. Докажем $5^\circ$. Пусть $F$ -- первообразная функции $f$ на отрезке
$[\overrightarrow{ab}]$:
 \beq\label{15.4.2}
F'(x)=f(x), \qquad  x\in [\overrightarrow{ab}]
 \eeq
Тогда композиция $F\circ\ph$ будет первообразной для функции $(f\circ \ph)\cdot
\ph'$, на отрезке $[\overrightarrow{\alpha\beta}]$, по формуле
\eqref{proizvodnaya-kompozitsii}:
 \beq\label{15.4.3}
(F\circ\ph)'(t)=F'(\ph(t))\cdot \ph'(t), \qquad t \in
[\overrightarrow{\alpha\beta}]
 \eeq
Поэтому
 \begin{multline*}
\int_\alpha^\beta f (\ph(t))\cdot \ph'(t) \, \d
t=\eqref{Newton-Leibnitz-v-aniz-integr-R^1}=(F\circ\ph)(t)\Big|_{t=\alpha}^{t=\beta}=
F\big(\kern-20pt\underbrace{\ph(\beta)}_{\scriptsize\begin{matrix}
\phantom{\tiny \eqref{ph(alpha)=a}}\quad\text{\rotatebox{90}{$=$}}\quad{\tiny
\eqref{ph(alpha)=a}}\\ b
\end{matrix}}\kern-20pt\big)-F\big(\kern-20pt\underbrace{\ph(\alpha)}_{\scriptsize\begin{matrix} \phantom{\tiny
\eqref{ph(alpha)=a}}\quad\text{\rotatebox{90}{$=$}}\quad{\tiny
\eqref{ph(alpha)=a}}\\ a
\end{matrix}}\kern-20pt\big)=\\=
F(b)-F(a)=\eqref{Newton-Leibnitz-v-aniz-integr-R^1}= \int_a^b f (x) \, \d x
\end{multline*}
 \epr

\noindent\rule{160mm}{0.1pt}\begin{multicols}{2}

Формула Ньютона-Лейбница считается такой важной в анализе потому, что позволяет
вычислять определенные интегралы от стандартных функций, сводя их к
неопределенным интегралам (вычислять которые мы уже научились). Вот несколько
примеров.

\begin{ex}\label{ex-15.3.2} В начале этой главы мы говорили, что геометрический смысл определенного
интеграла -- площадь фигуры под графиком функции. Найдем в качестве примера
площадь фигуры, ограниченной аркой синусоиды (то есть, интеграл \break
$\int_0^\pi \sin x \, \d x$):

%\pucture{0pt}{0pt}{11-7.pcx}

\vglue100pt

 \begin{multline*}
\int_0^\pi \sin x \, \d x= {\smsize\begin{pmatrix}\text{вспоминаем, что}\\
\text{$-\cos x$ есть}\\
\text{первообразная для $\sin x$}\\
\text{и применяем}\\
\text{формулу \eqref{Newton-Leibnitz-v-aniz-integr-R^1}}\end{pmatrix}}=\\=
-\cos x \Big|_{x=0}^{x=\pi}=-\cos \pi-(-\cos 0)=1+1=2
 \end{multline*}
\end{ex}

\begin{ex}\label{ex-15.3.3} Поучительно было бы найти площадь какой-нибудь
известной фигуры, например, круга, с помощью теоремы Ньютона-Лейбница, и
убедиться, что при этом получается правильный ответ. В качестве примера найдем
интеграл
$$
\int_{-1}^1 \sqrt{1-x^2}\, \d x
$$
(то есть площадь полукруга радиуса 1).

%\pucture{0pt}{0pt}{11-8.pcx}

\vglue100pt

 \begin{multline*}
\int_{-1}^1 \sqrt{1-x^2}\, \d x=
{\smsize\begin{pmatrix}\text{этот неопределенный}\\
\text{интеграл был вычислен}\\
\text{нами в примере \ref{ex-12.6.10}}\end{pmatrix}}=\\= \frac{1}{2}\left( x
\sqrt{1-x^2} +\arcsin x \right)\, \Big|_{x=-1}^{x=1}=\\=
\frac{1}{2}\left(\arcsin 1-\arcsin(-1)
\right)=\frac{1}{2}\left(\frac{\pi}{2}+\frac{\pi}{2}\right)=\frac{\pi}{2}
 \end{multline*}
Видно, что ответ правильный: площадь полукруга радиуса 1 равна $\frac{\pi}{2}$.
\end{ex}

\begin{ers} Найдите определенные интегралы:
 \begin{multicols}{2}
1) $\int_1^4 \frac{1+\sqrt{x}}{x^2}\, \d x$

2) $\int_0^1 \sqrt{1+x}\, \d x$

3) $\int_0^\frac{\pi}{2}\sin^3 x \, \d x$

4) $\int_{-1}^1 \frac{x\, \d x}{(x^2+1)^2}$

5) $\int_1^e \frac{1+\ln x}{x}\, \d x$

6) $\int_{-\frac{\pi}{2}}^{-\frac{\pi}{4}}\frac{\cos^3 x} {\sqrt[3]{\sin x}}\,
\d x$

7) $\int_0^1 \frac{d x}{x^2+4x+5}$

8) $\int_0^1 \frac{d x}{\sqrt{3+2x-x^2}}$ \end{multicols}\end{ers}

\end{multicols}\noindent\rule[10pt]{160mm}{0.1pt}

\paragraph{Интегрирование по частям.}

 \bit{
\item[$\bullet$] Если $f$ -- интегрируемая функция, а $g$ -- гладкая функция на
ориентированном отрезке $\overrightarrow{ab}$, то {\it интегралом от функции
$f$ вдоль функции $g$ по ориентированному отрезку $\overrightarrow{ab}$}
называется число
 \beq\label{DEF:int_a^b-u(x)-d-v(x)}
\int_a^b f(x) \, \d g(x):= \int_a^b f(x) \cdot g'(x) \, \d x
 \eeq
В частном случае, когда $f(x)=1$, это число обозначается также
 \beq\label{DEF:int_a^b-d-v(x)}
\int_a^b \d g(x):= \int_a^b g'(x) \, \d x
 \eeq
 }\eit

\btm[{\bf о перемещении гладкой
функции}]\label{TH:peremeshenie-gladkoi-funktsii} Перемещение гладкой функции
$g$ на ориентированном отрезке $\overrightarrow{ab}$ равно интегралу от единицы
вдоль $g$ по этому отрезку:
 \beq\label{peremeshenie-gladkoi-funktsii}
 \int_a^b \d g(x) = \int_a^b g'(x) \, \d x=g(x)\Big|_{x=a}^{x=b}
 \eeq
\etm \bpr Это следует из формулы Ньютона-Лейбница
\eqref{Newton-Leibnitz-v-aniz-integr-R^1}, если положить $f=g'$, $\varPhi=g$.
\epr

\bcor\label{TH:stroenie-glad-func} Функция $g:[a,b]\to\R$ является гладкой
тогда и только тогда, когда она представима в виде интеграла с переменным
верхним пределом от некоторой непрерывной функции $f:[a,b]\to\R$:
$$
g(x)=g(a)+\int_a^x f(t)\ \d t
$$
\ecor
 \bpr
В качестве $f$ нужно взять производную $g'$ функции $g$, доопределенную
произвольным образом в точках недифференцируемости $g$.
 \epr

\btm[{\bf об интегрировании по частям}]\label{TH:integrir-po-chastyam} Если $f$
и $g$ -- гладкие функции на ориентированном отрезке $\overrightarrow{ab}$, то
 \beq\label{integrir-po-chastyam}
\int_a^b f(x) \, \d g(x)= f(x)\cdot g(x) \, \Big|_{x=a}^{x=b}- \int_a^b g(x) \,
\d f(x)
 \eeq
 \etm
 \bpr
 \begin{multline*}
\int_a^b f(x)\ \d g(x)+\int_a^b g(x)\ \d
f(x)=\eqref{DEF:int_a^b-u(x)-d-v(x)}=\int_a^b f(x)\cdot g'(x)\ \d x+\int_a^b
g(x)\cdot f'(x)\ \d
x=\eqref{linein-aniz-int-v-R^1}=\\=\int_a^b\Big(\underbrace{f(x)\cdot
g'(x)+f'(x)\cdot g(x)}_{\scriptsize\begin{matrix}\text{\rotatebox{90}{$=$}}
\\ (f\cdot g)'(x)\end{matrix}}\Big)\ \d x=\int_a^b (f\cdot g)'(x)\ \d
x=\eqref{peremeshenie-gladkoi-funktsii}=(f\cdot
g)(x)\Big|_{x=a}^{x=b}=f(x)\cdot g(x)\Big|_{x=a}^{x=b}
 \end{multline*}
$$
 \Downarrow\put(20,0){\smsize
 (\text{переносим $\int_a^b  g(x)\, \d f(x)$ направо})}
$$
$$
\int_a^b f(x)\, \d g(x)= f(x)\cdot g(x) \, \Big|_{x=a}^{x=b}-\int_a^b g(x)\, \d
f(x)
$$
\epr

\noindent\rule{160mm}{0.1pt}\begin{multicols}{2}

\paragraph*{Интегрирование вдоль гладкой функции.}
Чтобы привыкнуть к формуле \eqref{DEF:int_a^b-u(x)-d-v(x)}, полезно вычислить
несколько раз определяемую этой формулой величину.

\begin{ex}\label{ex-15.4.1}
 \begin{multline*}
\int_0^1 x^3 \, \d x^5= \int_0^1 x^3 5 x^4 \, \d x=\\
=5\int_0^1 x^7 \, \d x= \frac{5}{8} x^8 \Big|_{x=0}^{x=1}=\frac{5}{8}
 \end{multline*}
\end{ex}

\begin{ex}
 \begin{multline*}
\int_0^\pi \cos x \, \d \sin x= \int_0^\pi \cos^2 x \, \d x=
\\
=\int_0^\pi \frac{1+\cos 2 x}{2} \, \d x= \frac{1}{2}\cdot x+
\frac{1}{4}\cdot\sin 2x \Big|_0^\pi=\frac{\pi}{2}
 \end{multline*}
\end{ex}

\paragraph*{Применение формулы замены переменной.} Покажем теперь, как
применяется формула \eqref{15.4.1}.

\begin{ex}\label{ex-15.4.3}
 \begin{multline*}
\int_0^2 x\cdot e^{x^2}\, \d x= \frac{1}{2}\int_0^2 e^{x^2}\, \d
x^2=\\={\smsize \left|
\begin{array}{c}
x=\sqrt{t}, \quad t=x^2
\\
x\in \overrightarrow{0;2} \Leftrightarrow t\in \overrightarrow{0;4}
\end{array}\right|}=\\= \frac{1}{2}\int_0^4 e^t \, \d t= \frac{1}{2} e^t
\Big|_0^4= \frac{e^4-1}{2}
 \end{multline*}
\end{ex}

\begin{ex}\label{ex-15.4.4}
 \begin{multline*}\int_1^4 \frac{x \, \d x}{1+\sqrt{x}} ={\smsize \left|
\begin{array}{c}
1+\sqrt{x}=t, \quad x=t^2-1
\\
x\in \overrightarrow{1;4} \Leftrightarrow t\in \overrightarrow{2;3}
\end{array}\right|}=\\= \int_2^3 \frac{(t^2-1) \d(t^2-1)}{t} = \int_2^3
\frac{(t^2-1)\cdot 2t \, \d t}{t} =\\= 2\int_2^3 (t^2-1)\, \d t = 2
\left\{\frac{t^3}{3}-t\right\}\, \Big|_{t=2}^{t=3} =\\= 2 \left\{
(9-3)-\left(\frac{8}{3}-2\right) \right\}= \frac{32}{3}\end{multline*}\end{ex}

\begin{ers} Найдите определенные интегралы:

1) $\int_0^1 \sqrt{4-x^2}\, \d x$

2) $\int_0^{\ln 5}\frac{e^x \sqrt{e^x-1}}{e^x+3}\, \d x$

3) $\int_0^\frac{a}{2}\frac{d x}{\sqrt{a^2-x^2}}$

4) $\int_0^1 x^2\sqrt{1-x^2}\, \d x$

5) $\int_0^1 \frac{d x}{e^x+e^{-x}}$

6) $\int_0^{\frac{\pi}{2}}\frac{d x}{1+\sin x+\cos x}$

7) $\int_0^{\frac{\pi}{4}}\frac{d x}{1+2\sin^2 x}$

\end{ers}

\paragraph*{Применение формулы интегрирования по частям.} Покажем, как
применяется формула \eqref{integrir-po-chastyam}.

\begin{ex}\label{ex-15.5.2}
 \begin{multline*}
\int_1^e \ln x \, \d x= \eqref{integrir-po-chastyam}= x\cdot \ln x \,
\Big|_{x=1}^{x=e} -\\- \int_1^e x \, \d \ln x= x\cdot \ln x \,
\Big|_{x=1}^{x=e} - \int_1^e x \cdot \frac{1}{x}\, \d x=\\= x\cdot \ln x \,
\Big|_{x=1}^{x=e} - \int_1^e 1\, \d x= x\cdot \ln x \, \Big|_{x=1}^{x=e} - x \,
\Big|_{x=1}^{x=e} =\\= e\cdot \ln e - 1\cdot \ln 1 -(e-1) =1
\end{multline*}\end{ex}

\begin{ex}\label{ex-15.5.3}
 \begin{multline*}\int_0^2 x e^x \, \d x= \int_0^2 x  \, \d e^x=
 \eqref{integrir-po-chastyam}=\\= x\cdot e^x \, \Big|_{x=0}^{x=2}
- \int_0^2 e^x \, \d x= x\cdot e^x \, \Big|_{x=0}^{x=2} -  e^x \,
\Big|_{x=0}^{x=2}=\\= 2\cdot e^2 -0\cdot e^0  -  (e^2-e^0)=e^2+1
\end{multline*}\end{ex}

\paragraph*{Использование рекуррентных соотношений.}

\bex Некоторые определенные интегралы вычисляются с помощью рекуррентных
соотношений. Например, следующие равенства\footnote{Равенства
\eqref{int_0^(pi/2)-sin^n-cos^n} понадобятся нам на с.\pageref{Wallis} при
доказательстве формулы Валлиса.}:
 \begin{multline}\label{int_0^(pi/2)-sin^n-cos^n}
\int_{0}^{\frac{\pi}{2}}\sin^nx\ \d x= \int_{0}^{\frac{\pi}{2}}\cos^nx\ \d
x=\\=\begin{cases}\frac{(n-1)!!}{n!!}\cdot\frac{\pi}{2}, & n\in 2\N-1\\
\frac{(n-1)!!}{n!!}, & n\in 2\N
\end{cases}
 \end{multline}
(напомним, что двойной факториал $n!!$ был определен формулами
\eqref{DEF:(2m)!!} и \eqref{DEF:(2m-1)!!}). Сначала проверим, что эти интегралы
совпадают:
 \begin{multline*}
 \int_{0}^{\frac{\pi}{2}}\sin^nx\ \d x=
 \left|\begin{matrix}x=\frac{\pi}{2}-y
 \\
 x=0\quad\Leftrightarrow\quad y=\frac{\pi}{2}
 \\
 x=\frac{\pi}{2} \quad\Leftrightarrow\quad y=0
\end{matrix}\right|
=\\=
 \int_{\frac{\pi}{2}}^0\sin^n\l \frac{\pi}{2}-y\r\ \d\l \frac{\pi}{2}-y \r=\\=
-\int_{\frac{\pi}{2}}^0\cos^ny\ \d y=\int_0^{\frac{\pi}{2}}\cos^ny\ \d y
 \end{multline*}
После этого выведем рекуррентную формулу:
 \begin{multline*}
I_n=\int_{0}^{\frac{\pi}{2}}\sin^nx\ \d x=\int_{0}^{\frac{\pi}{2}}\sin^{n-1}x\
\d(-\cos x)=\\=\underbrace{-\sin^{n-1}x\cdot\cos
x\Big|_0^{\frac{\pi}{2}}}_{\scriptsize\begin{matrix}\text{\rotatebox{90}{$=$}}\\
0\end{matrix}}+\\+(n-1)\cdot\int_{0}^{\frac{\pi}{2}}\sin^{n-2}x\cdot\cos^2 x\
\d x=\\=(n-1)\cdot\int_{0}^{\frac{\pi}{2}}\sin^{n-2}x\cdot(1-\sin^2 x)\ \d
x=\\=(n-1)\cdot I_{n-2}-(n-1)\cdot I_n
 \end{multline*}
$$
\Downarrow
$$
$$
I_n=\frac{n-1}{n}\cdot I_{n-2}
$$
Теперь для $n=0$ получаем:
$$
I_0=\int_{0}^{\frac{\pi}{2}}\ \d x=\frac{\pi}{2}
$$
и поэтому при $n=2k$,
 \begin{multline*}
I_{2k}=\frac{2k-1}{2k}\cdot
I_{2k-2}=\frac{(2k-1)\cdot(2k-3)}{2k\cdot(2k-2)}\cdot I_{2k-4}=\\=...=
\frac{(2k-1)\cdot(2k-3)\cdot...\cdot 1}{2k\cdot(2k-2)\cdot...\cdot 2}\cdot I_0=
\frac{(2k-1)!!}{(2k)!!}\cdot\frac{\pi}{2}
 \end{multline*}
А для $n=1$ получаем:
$$
I_1=\int_{0}^{\frac{\pi}{2}}\sin x\ \d x=-\cos x\Big|_0^{\frac{\pi}{1}}=1,
$$
поэтому при $n=2k+1$,
 \begin{multline*}
I_{2k+1}=\frac{2k}{2k+1}\cdot
I_{2k-1}=\frac{2k\cdot(2k-2)}{(2k+1)\cdot(2k-1)}\cdot I_{2k-3}=\\=...=
\frac{(2k)\cdot(2k-2)\cdot...\cdot 2}{(2k+1)\cdot(2k-1)\cdot...\cdot 1}\cdot
I_1= \frac{(2k)!!}{(2k+1)!!}
 \end{multline*}
\eex

\begin{ers} Найдите определенные интегралы:
 \begin{multicols}{2}
1) $\int_0^1 \arctg x \, \d x$

2) $\int_0^1 x e^{-x}\, \d x$

3) $\int_0^\frac{1}{2}\arcsin x \, \d x$

4) $\int_1^2 x \ln x\, \d x$

5) $\int_0^\pi x^3 \sin x \, \d x$

6) $\int_0^\pi e^x \sin x \, \d x$

7) $\int_0^\frac{\pi}{2} e^{2x}\cos x \, \d x$
\end{multicols}\end{ers}

\end{multicols}\noindent\rule[10pt]{160mm}{0.1pt}

\subsection{Вычисление геометрических величин} \label{CH-applic-def-integral}

В примерах \ref{ex-15.3.2} и \ref{ex-15.3.3} выше мы видели, что с помощью
определенного интеграла можно вычислять площади плоских фигур. Помимо этого,
как оказывается, определенный интеграл имеет много других полезных приложений,
например, с его помощью вычисляются длина кривой, объем и площадь фигуры
вращения, и разные другие геометрические величины. Здесь мы перечислим
некоторые из этих приложений в геометрии, не объясняя, однако, ни происхождение
применяемых нами формул, ни даже точный смысл самих терминов ``площадь'',
``объем'' и т.д. Наша цель -- дать возможность читателю поупражняться в
применении определенного интеграла, не обременяя его пока строгими
определениями и формальными доказательствами. Мы предполагаем, что интуитивный
смысл соот\-вест\-вующих терминов ему понятен уже сейчас, а точные объяснения
по этому поводу мы дадим позже в главах \ref{CH-mera-Jordana},
\ref{CH-dlina-krivoj} и \ref{CH-ploshad-poverhnosti}.

\noindent\rule{160mm}{0.1pt}\begin{multicols}{2}

\paragraph{Площадь правильной области на
плоскости.}\label{SEC-ploshad-pravil'noi-oblasti}

Пусть функции $f$ и $g$ непрерывны на отрезке $[a,b]$ и удовлетворяют
неравенству
$$
  f(x)\le g(x), \quad x\in [a;b]
$$
Тогда система неравенств
$$
\begin{cases}{a\le x\le b}\\{f(x)\le y\le g(x)}\end{cases}
$$
задает на координатной плоскости некоторую область $D$

%\pucture{0pt}{0pt}{12-1.pcx}

\vglue100pt

Всякая такая область называется {\it правильной}\index{область!правильная}
(относительно оси OX).

\begin{tm}\label{tm-16.1.1} Площадь $S_D$ правильной области
$$
D: \quad \begin{cases}{a\le x\le b}\\{f(x)\le y\le g(x)}\end{cases}
$$
равна
 \begin{equation} \label{16.1.1}
 S_D=\int_a^b \Big[ g(x)-f(x) \Big] \, \d x
 \end{equation}\end{tm}

\begin{ex}\label{ex-16.1.2} Найдем площадь фигуры, ограниченной кривыми
$$
  y=x, \quad y=2-x^2
$$
Найдем сначала точки пересечения этих линий:
 \begin{multline*}\left\{ \aligned
& y=x \\
& y=2-x^2
\endaligned
\right\}\; \Leftrightarrow \; \left\{ \aligned &
y=x \\
& x=2-x^2
\endaligned
\right\}\; \Leftrightarrow \\ \Leftrightarrow \; \left\{ \aligned
& y=x \\
& x^2+x-2=0
\endaligned
\right\}\; \Leftrightarrow \; \left\{ \aligned & y=x
\\
& \left[\aligned &
x=-2\\
& x=1\endaligned\right]
\endaligned
\right\}\; \Leftrightarrow \\ \Leftrightarrow\; \left[\aligned & \left\{
\aligned & x=-2 \\ & y=-2
\endaligned
\right\}\\
& \left\{ \aligned & x=1\\ & y=1
\endaligned
\right\}\endaligned\right]
\end{multline*}
После этого рисуем область:

%\pucture{0pt}{0pt}{12-2.pcx}

\vglue100pt

\noindent Теперь интересующая нас область в правильном виде записывается
следующим образом:
$$
D: \quad \begin{cases}{-2\le x\le 1}\\{x\le y\le 2-x^2}\end{cases}
$$
и остается применить формулу \eqref{16.1.1}:
 \begin{multline*}
S_D= \int_{-2}^1 \{ 2-x^2-x \}\, \d x =\\= \left\{
2x-\frac{x^3}{3}-\frac{x^2}{2}\right\}\, \Big|_{-2}^1 =\\= \left\{
2-\frac{1}{3}-\frac{1}{2}\right\}- \left\{ -4+\frac{8}{3}-2 \right\}
=\frac{9}{2}
 \end{multline*}
\end{ex}

\begin{ers}\label{ers-16.1.3} Найдите площади фигур, ограниченных линиями:
  \biter{
\item[1)] $y=e^x, \, y=e^{-x}, \, x=1$,

\item[2)] $y=2px^2, \, x=2py^2$,

\item[3)] $y=-x, \, y=2x-x^2$.
 }\eiter
 \end{ers}

\paragraph{Площадь области, ограниченной параметризованной кривой.}\label{SEC-ploshad-oblasti-ogranich-param-krivoi}

Для любых непрерывных функций $\ph,\chi:[\alpha,\beta]\to\R$ система уравнений
вида
$$
  \begin{cases}{x=\ph(t)}\\{y=\chi(t)}\\{t\in [\alpha;\beta]}\end{cases}
$$
задает на плоскости некоторое множество точек, называемых {\it (компактной)
параметризованной кривой}\index{кривая!параметризованная}.

\begin{ex}[\bf эллипс]\label{ex-16.2.1} Можно проверить, что уравнения
$$
\begin{cases}{x=a\cos t}\\{y=b\sin t}\end{cases}, \quad t\in \R
$$
определяют эллипс:

%\pucture{0pt}{0pt}{12-3.pcx}

\vglue100pt \noindent Действительно,
$$
\frac{x^2}{a}+\frac{y^2}{b}=\cos^2 t+\sin^2 t=1
$$
\end{ex}

\begin{ex}[\bf астроида]\label{ex-16.2.2} Уравнения
$$
  \begin{cases}{x=a\cos^3 t}\\{y=a\sin^3 t}\end{cases}, \quad t\in \R
$$
задают на плоскости кривую, называемую {\it астроидой}\index{астроида}. Для ее
построения можно составить таблицу

\bigskip

\vbox{\tabskip=0pt\offinterlineskip \halign to \hsize{ \vrule#\tabskip=2pt
plus3pt minus1pt & \strut\hfil\;#\hfil & \vrule#&\hfil#\hfil &
\vrule#&\hfil#\hfil & \vrule#&\hfil#\hfil & \vrule#&\hfil#\hfil &
\vrule#&\hfil#\hfil & \vrule#\tabskip=0pt\cr \noalign{\hrule}
height2pt&\omit&&\omit&&\omit&&\omit&&\omit&&\omit&\cr & $t$ && 0 &&
$\pm\frac{\pi}{6}$ && $\pm\frac{\pi}{4}$ & & $\pm\frac{\pi}{3}$ &&
$\pm\frac{\pi}{2}$  &\cr \noalign{\hrule}
height2pt&\omit&&\omit&&\omit&&\omit&&\omit&&\omit&\cr & $x$ && $a$ &&
$\frac{3\sqrt{3}}{8}a$ && $\frac{\sqrt{2}}{4}a$ & & $\frac{1}{8}a$ & & $0$ &\cr
\noalign{\hrule} height2pt&\omit&&\omit&&\omit&&\omit&&\omit&&\omit&\cr & $y$
&& $0$ && $\pm\frac{1}{8}a$ && $\pm\frac{\sqrt{2}}{4}a$ & &
$\pm\frac{3\sqrt{3}}{8}a$ && $\pm a$ &\cr \noalign{\hrule} } }

\bigskip

\bigskip

\vbox{\tabskip=0pt\offinterlineskip \halign to \hsize{ \vrule#\tabskip=2pt
plus3pt minus1pt & \strut\hfil\;#\hfil & \vrule#&\hfil#\hfil &
\vrule#&\hfil#\hfil & \vrule#&\hfil#\hfil & \vrule#&\hfil#\hfil &
\vrule#\tabskip=0pt\cr \noalign{\hrule}
height2pt&\omit&&\omit&&\omit&&\omit&&\omit&\cr & $t$ && $\pm\frac{2\pi}{3}$ &
& $\pm\frac{3\pi}{4}$ && $\pm\frac{5\pi}{6}$ && $\pm\pi$ &\cr \noalign{\hrule}
height2pt&\omit&&\omit& &\omit&&\omit&&\omit&\cr & $x$ && $-\frac{1}{8}a$ & &
$-\frac{\sqrt{2}}{4}a$ && $-\frac{3\sqrt{3}}{8}a$ & & $-a$ &\cr
\noalign{\hrule} height2pt&\omit& &\omit&&\omit&&\omit&&\omit&\cr & $y$ &&
$\pm\frac{3\sqrt{3}}{8}a$ & & $\pm\frac{\sqrt{2}}{4}a$ && $\pm\frac{1}{8}a$ &&
$0$ &\cr \noalign{\hrule} } }

\bigskip

\noindent затем отметить все эти точки на рисунке:

%\pucture{0pt}{0pt}{12-4.pcx}

\vglue100pt \noindent и соединить их ``плавной линией''. Получится картинка

%\pucture{0pt}{0pt}{12-5.pcx}

\vglue100pt
\end{ex}

\begin{ex}[\bf циклоида]\label{ex-16.2.3} Уравнения
$$
\begin{cases}{x=a(t-\sin t)}\\{y=a(1-\cos t)}\end{cases}, \quad
t\in \R
$$
задают на плоскости кривую, называемую {\it циклоидой}\index{циклоида}. Ее
также можно построить по точкам, и соответсвующая картинка будет выглядеть так:

%\pucture{0pt}{0pt}{12-6.pcx}

\vglue100pt
\end{ex}

Область на плоскости можно определять параметризованными кривыми. Например,
система
$$
D: \quad \begin{cases}{x=\ph(t)}\\{0\le y\le \psi(t)}\end{cases}, \quad  t\in
[\alpha;\beta]
$$
задает на плоскости примерно такую область

%\pucture{0pt}{0pt}{12-7.pcx}

\vglue100pt \noindent если $\ph)$ возрастает, и

%\pucture{0pt}{0pt}{12-8.pcx}

\vglue100pt \noindent если $\ph$ убывает.

\begin{tm}\label{tm-16.2.4} Пусть

 \bit{
\item[1)] $\ph$ -- гладкая монотонная функция на отрезке $t\in [\alpha;\beta]$;

\item[2)] $\chi$ -- непрерывная неотрицательная функция на отрезке $t\in
[\alpha;\beta]$: $\chi(t)\ge 0$.
 }\eit

Тогда площадь $S_D$ области, ограниченной параметризованной кривой
$$
D: \quad \begin{cases}{x=\ph(t)}\\{0\le y\le \chi(x)}\end{cases}, \quad  t\in
[\alpha;\beta]
$$
равна
 \beq
S_D=\left|\int_\alpha^\beta \chi(t) \, \d \ph(t)\right| \label{16.2.1}
 \eeq
\end{tm}

\begin{ex}[\bf площадь эллипса]\label{ex-16.2.5} Найдем площадь эллипса из примера
\ref{ex-16.2.1}.

%\pucture{0pt}{0pt}{12-9.pcx}

\vglue100pt \noindent Для этого выделим четвертинку этой области, определяемую
системой
$$
D: \quad \begin{cases}{x=a\cos t}\\{0\le y \le b\sin t}\end{cases}, \quad t\in
[0;\frac{\pi}{2}]
$$
%\pucture{0pt}{0pt}{12-10.pcx}

\vglue100pt \noindent и найдем ее площадь:
 \begin{multline*}
 S_D=\left|\int_0^{\frac{\pi}{2}} a\cos t \, \d a\sin t\right|=
 ab \left|\int_0^{\frac{\pi}{2}}\cos^2 t \, \d t\right|=\\=
 ab \left|\int_0^{\frac{\pi}{2}}\frac{1+\cos 2 t}{2}\, \d t\right|=
 \frac{ab}{2}\left|\l t+\frac{\sin 2 t}{2}\r \Big|_0^{\frac{\pi}{2}}\right|
 =\\=
 \frac{\pi ab}{4}
 \end{multline*}
Умножив это число на 4, получаем площадь эллипса:
$$
S=\pi ab
$$
\end{ex}

\begin{ex}\label{ex-16.2.6} Найдем площадь фигуры, ограниченной астроидой
из примера \ref{ex-16.2.2}.

%\pucture{0pt}{0pt}{12-11.pcx}

\vglue100pt \noindent Для этого выделим четвертинку, определяемую системой
$$
D: \quad \begin{cases}{x=a\cos^3 t}\\{0\le y \le a\sin^3 t}\end{cases}, \quad
t\in \left[0;\frac{\pi}{2}\right]
$$
%\pucture{0pt}{0pt}{12-12.pcx}

\vglue100pt \noindent и найдем ее площадь:
 \begin{multline*}
 S_D= \left|\int_0^{\frac{\pi}{2}} a\sin^3 t  \, \d \l a\cos^3 t\r\right|=
 \\=
 3a^2\left|\int_0^{\frac{\pi}{2}}\sin^4 t  \cos^2 t \, \d
 t\right|=\\=
 3a^2\left|\int_0^{\frac{\pi}{2}}\l \frac{1-\cos 2t}{2}\r^2\frac{1+\cos 2t}{2}\, \d t\right|
 =\\=
 \frac{3a^2}{8}\left|\int_0^{\frac{\pi}{2}}\l 1-\cos 2t \r  \l 1-\cos^2 2t \r \, \d
 t\right|=\\=
 \frac{3a^2}{8}\left|\int_0^{\frac{\pi}{2}}\l 1-\cos 2t\r \l 1-\frac{1+\cos 4t}{2}\r \, \d t
 \right|=\\=
\frac{3a^2}{16}\left|\int_0^{\frac{\pi}{2}}\l 1-\cos 2t \r \l 1-\cos 4t \r \,
\d t\right|=\\=
 \frac{3a^2}{16}\left|\int_0^{\frac{\pi}{2}}\l 1-\cos 2t-\cos 4t+\cos 2t \cos 4t \r \, \d t
 \right|
 =\\=
\frac{3a^2}{16}\Bigg|\int_0^{\frac{\pi}{2}}\Big( 1-\cos 2t -\cos
4t+\\+\frac{1}{2}\l \cos 6t +\cos 2t \r \Big) \, \d t \Bigg|=\\=
 \frac{3a^2}{32}\left|\int_0^{\frac{\pi}{2}}\l 2-\cos 2t -2\cos 4t+\cos
 6t \r \, \d t\right|=\\=
 \frac{3a^2}{32}\left|\l 2t-\frac{1}{6}\sin  2t-\frac{1}{2}\sin 4t+ \frac{1}{6}\sin 6t \r
 \Big|_0^{\frac{\pi}{2}}\right|=\\= \frac{3a^2 \pi}{32}
 \end{multline*}
Умножив это число на 4, получим площадь всей фигуры:
$$
  S=\frac{3a^2 \pi}{8}
$$
\end{ex}

\begin{ex}\label{ex-16.2.7}
Найдем площадь фигуры, ограниченной аркой циклоиды
$$
D: \quad \begin{cases}{x=a(t-\sin t)}\\{0\le y \le a(1-\cos t)}\end{cases},
\quad t\in [0;2\pi]
$$
%\pucture{0pt}{0pt}{12-13.pcx}

\vglue100pt \noindent Вычисляем площадь:
 \begin{multline*}
 S_D= \left|\int_0^{2\pi} a(1-\cos t) \, \d a(t-\sin t)\right|=\\=
 a^2\left|\int_0^{2\pi} (1-\cos t)^2 \, \d t\right|=
 \\=a^2 \left|\int_0^{2\pi}\l 1-2\cos t+\cos^2 t \r \, \d t\right|
 =\\= a^2
 \left|\int_0^{2\pi}\l 1-2\cos t+\frac{1+\cos 2 t}{2}\r \, \d
 t\right|
 =\\=
 a^2  \left|\l \frac{3}{2} t-2\sin t+\frac{\sin 2 t}{4}\r \Big|_0^{2\pi}\right|= 3\pi a^2
 \end{multline*}\end{ex}

\begin{ers}\label{ers-16.2.7} Найдите площади фигур, ограниченных
параметризованными кривыми:

 \begin{multicols}{2}
1) $x=a(t^2+1), \, y=b(t^3-3t)$ (площадь петли)

2) $x=\frac{1}{3}(3-t^2), \, y=t^2$ (площадь петли)

3) $x=2t-t^2, \, y=2t^2-t^3$ (площадь петли)
\end{multicols}\end{ers}

\paragraph{Площадь области в полярных координатах.}\label{SEC-ploshad-v-polyarnyh-koordinatah}

Положение точки A на плоскости можно задавать не только декартовыми
координатами (то есть обычными значениями переменных $x$ и $y$), но также
полярными координатами -- значением угла $\ph$ между лучом OX и лучом OA,
соединяющим начало координат O с точкой A, и расстояния $\rho$ от A до начала
координат O:

%\pucture{0pt}{0pt}{12-14.pcx}

\vglue100pt

При этом, полярные координаты связаны с декартовыми формулами
$$
\begin{cases}{x=\rho\cdot \cos \ph}\\{y=\rho\cdot \sin \ph}\end{cases},
\qquad
\begin{cases}{\rho=\sqrt{x^2+y^2}}\\{\tg \ph=\frac{y}{x}}\end{cases},
$$

Как и в случае с декартовыми координатами, всякое уравнение, связывающее
полярные координаты,
$$
  \rho=R(\ph)
$$
задает некоторую кривую на плоскости. Рассмотрим примеры.

\begin{ex}[\bf кардиоида]\label{ex-16.3.1}
Построим кривую, заданную уравнением в полярных координатах
$$
  \rho=a\cdot(1+\cos \ph)
$$
Это можно сделать по точкам. Сначала составляется таблица

\bigskip

\vbox{\tabskip=0pt\offinterlineskip \halign to \hsize{ \vrule#\tabskip=2pt
plus3pt minus1pt & \strut\hfil\;#\hfil & \vrule#&\hfil#\hfil &
\vrule#&\hfil#\hfil & \vrule#&\hfil#\hfil & \vrule#&\hfil#\hfil &
\vrule#&\hfil#\hfil & \vrule#\tabskip=0pt\cr \noalign{\hrule}
height2pt&\omit&&\omit&&\omit& &\omit&&\omit&&\omit&\cr & $\ph$ && 0 &&
$\pm\frac{\pi}{6}$ && $\pm\frac{\pi}{4}$ & & $\pm\frac{\pi}{3}$ & &
$\pm\frac{\pi}{2}$ &\cr \noalign{\hrule} height2pt&\omit&&\omit&
&\omit&&\omit&&\omit&&\omit&\cr & $\rho$ && $2a$ &&
$a\left(1+\frac{\sqrt{3}}{2}\right)$ & & $a\left(1+\frac{\sqrt{2}}{2}\right)$ &
& $\frac{3a}{2}$ && $a$ &\cr \noalign{\hrule} } }

\bigskip

\vbox{\tabskip=0pt\offinterlineskip \halign to \hsize{ \vrule#\tabskip=2pt
plus3pt minus1pt & \strut\hfil\;#\hfil & \vrule#&\hfil#\hfil &
\vrule#&\hfil#\hfil & \vrule#&\hfil#\hfil & \vrule#&\hfil#\hfil &
\vrule#\tabskip=0pt\cr \noalign{\hrule} height2pt&\omit&&\omit&
&\omit&&\omit&&\omit&\cr & $\ph$ && $\pm\frac{2\pi}{3}$ & & $\pm\frac{3\pi}{4}$
&& $\pm\frac{5\pi}{6}$ & & $\pm\pi$ &\cr \noalign{\hrule}
height2pt&\omit&&\omit&&\omit&&\omit&&\omit&\cr & $\rho$ && $\frac{a}{2}$ & &
$a\left(1-\frac{\sqrt{2}}{2}\right)$ & & $a\left(1-\frac{\sqrt{3}}{2}\right)$
&& $0$ &\cr \noalign{\hrule} } }

\bigskip

\noindent Потом соответствующие точки строятся на плоскости:

%\pucture{0pt}{0pt}{12-15.pcx}

\vglue100pt \noindent Эти точки соединяются плавной линией, и получается
кривая:

%\pucture{0pt}{0pt}{12-16.pcx}

\vglue100pt
\end{ex}

\begin{ex}[\bf лемниската Бернулли]\label{ex-16.3.2}
По\-стро\-им кривую, заданную уравнением в полярных координатах
$$
  \rho=a\sqrt{\cos 2\ph}
$$
Это также делается по точкам. Сначала составляется таблица (в которой знак
$\nexists$ говорит, что для данного значения $\ph$ соответсвующее значение
параметра $\rho$ не определено):

\bigskip

\vbox{\tabskip=0pt\offinterlineskip \halign to \hsize{ \vrule#\tabskip=2pt
plus3pt minus1pt & \strut\hfil\;#\hfil & \vrule#&\hfil#\hfil &
\vrule#&\hfil#\hfil & \vrule#&\hfil#\hfil & \vrule#&\hfil#\hfil &
\vrule#&\hfil#\hfil & \vrule#\tabskip=0pt\cr \noalign{\hrule}
height2pt&\omit&&\omit&&\omit&&\omit&&\omit&&\omit&\cr & $\ph$ && 0 &&
$\pm\frac{\pi}{6}$ && $\pm\frac{\pi}{4}$ & & $\pm\frac{\pi}{3}$ & &
$\pm\frac{\pi}{2}$ &\cr \noalign{\hrule} height2pt&\omit&&\omit&
&\omit&&\omit&&\omit&&\omit&\cr & $\rho$ && $a$ && $\frac{a}{\sqrt{2}}$ && 0 &
& $\nexists$ && $\nexists$ &\cr \noalign{\hrule} } }

\bigskip

\vbox{\tabskip=0pt\offinterlineskip \halign to \hsize{ \vrule#\tabskip=2pt
plus3pt minus1pt & \strut\hfil\;#\hfil & \vrule#&\hfil#\hfil &
\vrule#&\hfil#\hfil & \vrule#&\hfil#\hfil & \vrule#&\hfil#\hfil &
\vrule#\tabskip=0pt\cr \noalign{\hrule} height2pt&\omit&&\omit&
&\omit&&\omit&&\omit&\cr & $\ph$ && $\pm\frac{2\pi}{3}$ & & $\pm\frac{3\pi}{4}$
&& $\pm\frac{5\pi}{6}$ & & $\pm\pi$ &\cr \noalign{\hrule} height2pt&\omit&
&\omit&&\omit&&\omit&&\omit&\cr & $\rho$ && $\nexists$ & & 0 &&
$\frac{a}{\sqrt{2}}$ && $a$ &\cr \noalign{\hrule} } }

\bigskip

\noindent Потом соответствующие точки строятся на плоскости:

%\pucture{0pt}{0pt}{12-17.pcx}

\vglue100pt \noindent И в результате получается кривая:

%\pucture{0pt}{0pt}{12-18.pcx}

\vglue100pt
\end{ex}

\begin{ex}[\bf окружность]\label{ex-16.3.3}
Построим кривую, заданную уравнением в полярных координатах
$$
  \rho=a\sin \ph
$$
Составляем таблицу

\bigskip

\vbox{\tabskip=0pt\offinterlineskip \halign to \hsize{ \vrule#\tabskip=2pt
plus3pt minus1pt & \strut\hfil\;#\hfil & \vrule#&\hfil#\hfil &
\vrule#&\hfil#\hfil & \vrule#&\hfil#\hfil & \vrule#&\hfil#\hfil &
\vrule#&\hfil#\hfil & \vrule#\tabskip=0pt\cr \noalign{\hrule}
height2pt&\omit&&\omit&&\omit& &\omit&&\omit&&\omit&\cr & $\ph$ && $-\pi$ &&
$-\frac{3\pi}{4}$ && $-\frac{\pi}{2}$ & & $-\frac{\pi}{4}$ & & $0$ &\cr
\noalign{\hrule} height2pt&\omit&&\omit& &\omit&&\omit&&\omit&&\omit&\cr &
$\rho$ && $0$ && $\nexists$ && $\nexists$ & & $\nexists$ && $0$ &\cr
\noalign{\hrule} } }

\bigskip

\vbox{\tabskip=0pt\offinterlineskip \halign to \hsize{ \vrule#\tabskip=2pt
plus3pt minus1pt & \strut\hfil\;#\hfil & \vrule#&\hfil#\hfil &
\vrule#&\hfil#\hfil & \vrule#&\hfil#\hfil & \vrule#&\hfil#\hfil &
\vrule#\tabskip=0pt\cr \noalign{\hrule} height2pt&\omit&&\omit&
&\omit&&\omit&&\omit&\cr & $\ph$ && $\frac{\pi}{4}$ & & $\frac{\pi}{2}$ &&
$\frac{3\pi}{4}$ & & $\pi$ &\cr \noalign{\hrule} height2pt&\omit&
&\omit&&\omit&&\omit&&\omit&\cr & $\rho$ && $\frac{a}{\sqrt{2}}$ & & $a$ &&
$\frac{a}{\sqrt{2}}$ && $0$ &\cr \noalign{\hrule} } }

\bigskip

\noindent (здесь знак $\nexists$ мы пишем в случае, когда значение $\rho$
получается отрицательным; это должно означать, что при данном $\ph$
соответствующая точка не существует, поскольку $\rho$ есть расстояние до нуля,
и значит не может быть отрицательным).

Потом на плоскости строятся точки:

%\pucture{0pt}{0pt}{12-19.pcx}

\vglue100pt \noindent В результате получается кривая:

%\pucture{0pt}{0pt}{12-20.pcx}

\vglue100pt

\noindent Нетрудно показать, что эта кривая -- окружность. Для этого нужно
перейти к декартовым координатам:
 \begin{multline*}
\rho=a\sin \ph \;\Leftrightarrow\; \rho^2=a \rho \sin \ph
\;\Leftrightarrow\; x^2+y^2=a y \;\Leftrightarrow\\
\Leftrightarrow\; x^2+y^2-a y+\frac{a^2}{4}=\frac{a^2}{4} \;\Leftrightarrow\;
x^2+\left(y-\frac{a}{2}\right)^2=\left(\frac{a}{2}\right)^2
 \end{multline*}
\end{ex}

\begin{tm}\label{tm-16.3.4}
Площадь области на плоскости, заданной неравенствами
$$
D:  \begin{cases}{\alpha\le \ph \le \beta}\\{\rho\le R(\ph)}\end{cases}
$$
%\pucture{0pt}{0pt}{12-21.pcx}

\vglue100pt

\noindent вычисляется по формуле
 \beq
 S_D=\frac{1}{2}\int_\alpha^\beta R(\ph)^2 \, \d \ph
 \label{16.3.1}
 \eeq
\end{tm}

\begin{ex}\label{ex-16.3.5} Вычислим площадь области, ограниченной кардиоидой
(линией, построенной нами в примере \ref{ex-16.3.1}):
$$
  \rho=a\cdot(1+\cos \ph)
$$

%\pucture{0pt}{0pt}{12-22.pcx}

\vglue100pt \noindent По теореме \ref{tm-16.3.4} получаем:
 \begin{multline*}
S_D=\frac{1}{2}\int_0^{2\pi} a^2\cdot(1+\cos \ph)^2 \, \d \ph=\\=
\frac{a^2}{2}\cdot \int_0^{2\pi} (1+2\cos \ph+\cos^2 \ph) \, \d \ph=\\=
\frac{a^2}{2}\cdot \int_0^{2\pi}\l 1+2\cos \ph+\frac{1+\cos 2 \ph}{2}\r \, \d
\ph=\\= \frac{a^2}{2}\cdot \int_0^{2\pi}\l \frac{3}{2}+2\cos \ph+\frac{\cos 2
\ph}{2}\r \, \d \ph=\\= \frac{a^2}{2}\cdot  \l \frac{3}{2}\ph +2\sin
\ph+\frac{\sin 2 \ph}{4}\r \Big|_{\ph=0}^{\ph=2\pi}= \frac{3\pi
a^2}{2}\end{multline*}\end{ex}

\begin{ex}\label{ex-16.3.6} Вычислим площадь области, ограниченной лемнискатой
Бернулли (линией, построенной нами в примере \ref{ex-16.3.2}):
$$
  \rho=a\sqrt{\cos 2\ph}
$$

%\pucture{0pt}{0pt}{12-23.pcx}

\vglue100pt \noindent Для этого достаточно найти площадь одного ``лепестка'':
 \begin{multline*}
S_A= \frac{1}{2}\int_{-\frac{\pi}{4}}^{\frac{\pi}{4}} a^2\cdot \cos 2 \ph \, \d
\ph= \frac{a^2}{2}\cdot \frac{\sin 2
\ph}{2}\Big|_{\ph=-\frac{\pi}{4}}^{\ph=\frac{\pi}{4}}=\\= \frac{a^2}{4}\cdot \l
1-(-1) \r= \frac{a^2}{2}
 \end{multline*}
Общая площадь теперь получается удвоением площади одного лепестка:
$$
  S_D=2 S_A=a^2
$$
\end{ex}

\begin{ex}\label{ex-16.3.7} Убедимся, что формула \eqref{16.3.1}, будучи
применена к вычислению площадь круга, построенного в примере \ref{ex-16.3.3},
приводит к правильному результату : $\pi R^2=\frac{\pi a^2}{4}$.
$$
  \rho=a \sin \ph
$$

%\pucture{0pt}{0pt}{12-24.pcx}

\vglue100pt

 \begin{multline*}
S_D= \frac{1}{2}\int_0^{\pi} a^2\cdot \sin^2 \ph \, \d \ph=\\=
\frac{a^2}{2}\cdot \int_0^{\pi}\l \frac{1-\cos 2 \ph}{2}\r \, \d \ph=\\=
\frac{a^2}{2}\cdot \l \frac{\ph}{2} -\frac{\sin 2 \ph}{2}\r
\Big|_{\ph=0}^{\ph=\pi} = \frac{\pi a^2}{4}
 \end{multline*}
\end{ex}

\begin{ers} Найдите площади фигур, ограниченных линиями:
 \begin{multicols}{2}
1) $\rho=4\sin^2 \ph$

2) $\rho=a(2+\sin \ph)$

3) $\rho=a \cos \ph$

4) $\rho=a \sin 3 \ph$

5) $\rho=a \cos 4 \ph$

6) $\rho=a \sin^2 3\ph$

7) $\rho=a \cos^2 3\ph$ \end{multicols}\end{ers}

\paragraph{Длина кривой.}\label{SEC-dlina-krivoi}

\begin{tm}\label{tm-16.4.1}
Если кривая задана уравнением
$$
  y=f(x), \quad x\in [a;b]
$$
(в котором $f$ -- гладкая функция на отрезке $[a;b]$), то ее длина равна
 \beq\label{16.4.1}
  l=\int_a^b \sqrt{1+\l f'(x) \r^2}\, \d x
 \eeq
\end{tm}

\begin{ex}\label{ex-16.4.2} Найдем длину куска полукубической параболы:
$$
  y=x^\frac{3}{2}, \quad x\in [0;1]
$$
По формуле \eqref{16.4.1} имеем:
\begin{multline*}
l= \int_0^1 \sqrt{1+\l \frac{3}{2} x^\frac{1}{2}\r^2}\, \d x= \int_0^1 \sqrt{1+
\frac{9}{4} x}\, \d x=\\={\smsize \left|
\begin{array}{c}
t=\sqrt{1+ \frac{9}{4} x}\\ x=\frac{4}{9}\l t^2-1 \r \\
x\in [0;1] \Leftrightarrow t\in [1;\frac{\sqrt{13}}{2}]
\end{array}\right|}=\\
=\int_1^\frac{\sqrt{13}}{2} t \, \d \lll \frac{4}{9}\l t^2-1 \r \rrr =
\frac{8}{9}\int_1^\frac{\sqrt{13}}{2} t^2 \, \d t =\\=
\frac{8}{9}\frac{t^3}{3}\, \Big|_1^\frac{\sqrt{13}}{2} = \frac{8}{27}\lll \l
\frac{\sqrt{13}}{2}\r^3 - 1 \rrr  = \frac{13\sqrt{13} -
8}{27}\end{multline*}\end{ex}

\begin{tm}\label{tm-16.4.3}
Если кривая задана параметрическими уравнениями
$$
  \begin{cases}{x=\ph(t)}\\{y=\chi(t)}\end{cases}, \quad t\in [\alpha;\beta]
$$
 причем $\ph(t), \, \chi(t)$ -- гладкие функции на $[\alpha;\beta]$, то ее длина равна
 \beq\label{16.4.2}
  l=\int_\alpha^\beta \sqrt{\l \ph'(t) \r^2+\l \chi'(t) \r^2}\, \d t
 \eeq
\end{tm}

\begin{ex}[\bf длина астроиды]\label{ex-16.4.4} Найдем длину астроиды:
$$
  \begin{cases}{x=a\cos^3 t}\\{y=a\sin^3 t}\end{cases}
$$
По формуле \eqref{16.4.2}, длина четвертинки равна
 \begin{multline*}
\frac{l}{4}=\\= \int_0^\frac{\pi}{2}\sqrt{\l -3a\cos^2 t\sin t \r^2+ \l
3a\sin^2 t\cos t  \r^2}\, \d t=\\= 3a\int_0^\frac{\pi}{2}\sqrt{\cos^4 t\sin^2
t+ \sin^4 t\cos^2 t}\, \d t=\\= 3a\int_0^\frac{\pi}{2}\left| \cos t\sin t
\right| \sqrt{\cos^2 t+\sin^2 t}\, \d t=\\= 3a\int_0^\frac{\pi}{2}  \cos t\sin
t \sqrt{1}\, \d t= 3a\int_0^\frac{\pi}{2}  \sin t \, \d \sin t= \\=3a
\frac{\sin^2 t}{2}\, \Big|_0^\frac{\pi}{2} = \frac{3a }{2}\end{multline*}
Умножив на 4, получаем ответ:
$$
  l=6a
$$
\end{ex}

\begin{tm}\label{tm-16.4.5}
Если кривая задана уравнением в полярных координатах,
$$
\rho=R(\ph), \quad \ph\in [\alpha;\beta]
$$
(где $R(\ph)$ -- гладкая функция на $[\alpha;\beta]$), то ее длина равна
 \beq\label{16.4.3}
  l=\int_\alpha^\beta \sqrt{R(\ph)^2+\l R'(\ph) \r^2}\, \d \ph
 \eeq
\end{tm}

\begin{ex}[\bf длина кардиоиды]\label{ex-16.4.5} Найдем длину кардиоиды:
$$
  \rho=a\cdot(1+\cos \ph)
$$
По формуле \eqref{16.4.3},
 \begin{multline*}
l= \int_0^{2\pi}\sqrt{\lll a\cdot (1+\cos \ph) \rrr ^2+\l a\sin \ph \r^2}\, \d
\ph=\\= a\int_0^{2\pi}\sqrt{1+2\cos \ph+\cos^2 \ph+\sin^2 \ph }\, \d \ph=\\=
a\int_0^{2\pi}\sqrt{2+2\cos \ph}\, \d \ph=\\= a\int_0^{2\pi}\sqrt{4\cos^2
\frac{\ph}{2}}\, \d \ph=
2a\int_0^{2\pi}\left| \cos \frac{\ph}{2}\right| \, \d \ph=\\= {\smsize\begin{pmatrix}\text{разбиваем отрезок $[0;2\pi]$}\\
\text{на два отрезка, где}\,\, \cos \frac{\ph}{2}\\
\text{имеет постоянный знак}\end{pmatrix}}=\\= 2a\int_0^\pi \left| \cos
\frac{\ph}{2}\right| \, \d \ph+ 2a\int_\pi^{2\pi}\left| \cos
\frac{\ph}{2}\right| \, \d \ph=\\=
{\smsize\begin{pmatrix}\text{при $\ph\in [0;\pi]$ имеем}  \\
\left| \cos \frac{\ph}{2}\right|= \cos \frac{\ph}{2}\\
\text{а при $\ph\in [\pi;2\pi]$ имеем}  \\
\left| \cos \frac{\ph}{2}\right|=-\cos \frac{\ph}{2}\end{pmatrix}}=
2a\int_0^\pi \cos \frac{\ph}{2}\, \d \ph -\\- 2a\int_\pi^{2\pi}\cos
\frac{\ph}{2}\, \d \ph= 4a\sin \frac{\ph}{2}\, \Big|_0^\pi - 4a\sin
\frac{\ph}{2}\, \Big|_\pi^{2\pi} =\\= 4a \l \sin \frac{\pi}{2}-0 \r - 4a \l
\sin \frac{2\pi}{2} -\sin \frac{\pi}{2}\r =\\= 4a + 4a =8a
\end{multline*}\end{ex}

\begin{ers} Найдите длину кривой:
 \biter{
\item[1)] $y=\frac{1}{4}x^2-\frac{1}{2}\ln x, \, x\in [1;e]$

\item[2)] $y=\frac{1}{3}(3-x)\sqrt{x}, \, x\ge 0$

\item[3)] $x=a(t-\sin t), y=a(1-\cos t)\, t\in [0;2\pi]$

\item[4)] $x=e^t\cos t, y=e^t\sin t)\, t\in [0;\pi]$

\item[5)] $\rho=a\sqrt{\cos 2\ph}$

\item[6)] $\rho=a\sin 2\ph$

\item[7)] $\rho=a\sin^2 \ph$
 }\eiter
 \end{ers}

\paragraph{Объем тела вращения.}

Пусть на отрезке $[a;b]$ задана непрерывная неотрицательная функция $f$

%\pucture{0pt}{0pt}{12-25.pcx}

\vglue100pt \noindent Если ее график заставить вращаться вокруг оси OX,
получится поверхность, описываемая уравнением
$$
y^2+z^2=f(x)^2, \quad x\in [a;b]
$$
и называемая {\it поверхностью вращения}\index{поверхность!вращения}.

%\pucture{0pt}{0pt}{12-26.pcx}

\vglue100pt \noindent Эта поверхность ограничивает вместе с плоскостями $x=a$ и
$x=b$ область в терхмерном пространстве, которая называется {\it телом
вращения}\index{тело!вращения} и описывается системой неравенств
 \beq
\begin{cases}{y^2+z^2\le f(x)^2}\\{a\le x\le b}\end{cases}\label{16.5.1}
 \eeq

\begin{tm}\label{tm-16.5.1}
Объем тела вращения \eqref{16.5.1} вычисляется по формуле
 \beq \label{16.5.2}
 V=\pi \int_a^b f(x)^2 \, \d x
 \eeq
\end{tm}

\begin{ex}\label{ex-16.5.2} Найдем объем тела, образо\-ван\-но\-го вращением
параболы $y=\sqrt{x}$ и плоскостью $x=b$:

%\pucture{0pt}{0pt}{12-27.pcx}

\vglue100pt \noindent По формуле \eqref{16.5.2} получаем
$$
V=\pi \int_0^b \l \sqrt{x}\r ^2 \, \d x= \pi \int_0^b x \, \d x=\frac{\pi}{2}
x^2 \, \Big|_0^b =\frac{\pi}{2} b^2
$$
\end{ex}

\begin{ex}\label{ex-16.5.3} Найдем объем тела, образован\-но\-го вращением
полуволны синусоиды $y=\sin x, \, x\in [0;\pi]$:

%\pucture{0pt}{0pt}{12-28.pcx}

\vglue100pt \noindent По формуле \eqref{16.5.2} получаем
 \begin{multline*}
V=\pi \int_0^\pi  \sin^2 x \, \d x=\pi \int_0^\pi  \frac{1-\cos 2 x}{2}\, \d
x=\\= \frac{\pi}{2}\l x -\frac{\sin 2 x}{2}\r \Big|_0^\pi=\frac{\pi^2}{2}
 \end{multline*}
\end{ex}

Кривая, вращаемая вокруг оси, может быть задана параметрически, и тогда
соответствующее тело вращения будет описываться системой
 \beq
D: \quad \begin{cases}{x=\ph(t)}\\{y^2+z^2\le \chi(t)^2}\end{cases}, \quad t\in
[\alpha;\beta] \label{16.5.3}
 \eeq

\begin{tm}\label{tm-16.5.4}
Объем тела вращения \eqref{16.5.3} вычисляется по формуле
 \beq
 V=\pi \left|\int_\alpha^\beta \chi(t)^2 \, \d\ph(t)\right|
 \label{16.5.4}
 \eeq
\end{tm}

\begin{ex}\label{ex-16.5.5} Найдем объем эллипсоида, то есть тела,
образованного вращением эллипса
$$
\begin{cases}{x=a\cos t}\\{y=b\sin t}\end{cases}
$$
вокруг оси OX:

%\pucture{0pt}{0pt}{12-29.pcx}

\vglue100pt \noindent По формуле \eqref{16.5.4} получаем
\begin{multline*}
 V=\pi\left|\int_0^\pi \l b\sin t \r ^2 \, \d \l a\cos t
 \r\right|=\\=
 a b^2 \pi \left|\int_0^\pi \sin^2 t \, \d \cos t \right|=\\=
 a b^2 \pi \left|\int_0^\pi \l 1- \cos^2 t \r \, \d \cos t\right| =\\=
 a b^2 \pi  \left|\l \cos t- \frac{\cos^3 t}{3}\r  \Big|_0^\pi
 \right|=\\=
 a b^2 \pi  \l 2-\frac{2}{3}\r = \frac{4\pi}{3} a b^2
 \end{multline*}\end{ex}

\begin{ers} Найдите объем тел вращения:
 \biter{
\item[1)] $y=x^3, \, x\in [0;b]$

\item[2)] $y=e^x, \, x\in [a;b]$

\item[3)] $x=a(t-\sin t), \, y=a(1-\cos t),$ $t\in [0;2\pi]$

\item[4)] $x=a\cos^3 t, \, y=a\sin^3 t$
 }\eiter
\end{ers}

\paragraph{Площадь поверхности вращения.}

\begin{tm}\label{tm-16.6.1}
Площадь поверхности вращения
$$
y^2+z^2=f(x)^2, \quad x\in [a;b] \label{16.6.1}
$$

%\pucture{0pt}{0pt}{12-26.pcx}

\vglue100pt \noindent вычисляется по формуле
$$
S= 2\pi \int_a^b f(x)\cdot \sqrt{1+ \l f'(x) \r^2}\, \d x \label{16.6.2}
$$
\end{tm}

\begin{ex}\label{ex-16.6.2} Найдем площадь поверхности, образованной вращением
параболы $y=\sqrt{x}, \quad x\in [0;b]$:

%\pucture{0pt}{0pt}{12-27.pcx}

\vglue100pt \noindent По формуле \eqref{16.6.2} получаем
\begin{multline*}
S=2\pi \int_0^b \sqrt{x}\cdot \sqrt{1+ \l \frac{1}{2\sqrt{x}}\r^2}\, \d x=\\=
2\pi \int_0^b \sqrt{x}\cdot \sqrt{1+ \frac{1}{4x}}\, \d x= 2\pi \int_0^b
\sqrt{x+ \frac{1}{4}}\, \d x=\\= 2\pi \int_0^b \sqrt{x+ \frac{1}{4}}\, \d \l x+
\frac{1}{4}\r= 2\pi \frac{2}{3}\l x+ \frac{1}{4}\r^\frac{3}{2}\, \Big|_0^b=\\=
\frac{4\pi}{3}\lll \l b+ \frac{1}{4}\r^\frac{3}{2}-
\l\frac{1}{4}\r^\frac{3}{2}\rrr =\\= \frac{4\pi}{3}\lll \l b+
\frac{1}{4}\r^\frac{3}{2}- \frac{1}{8}\rrr
\end{multline*}\end{ex}

\begin{ex}\label{ex-16.6.3} Найдем площадь поверхности, образованной вращением
полуволны синусоиды $y=\sin x, \, x\in [0;\pi]$:

%\pucture{0pt}{0pt}{12-28.pcx}

\vglue100pt \noindent По формуле \eqref{16.5.2} получаем
 \begin{multline*}
S= 2\pi \int_0^\pi \sin x\cdot \sqrt{1+\cos^2 x}\, \d x=\\= -2\pi \int_0^\pi
\sqrt{1+\cos^2 x}\, \d \cos x=\\={\smsize \ml
\begin{array}{c}\cos x=t
\\
x\in [0;\pi] \Leftrightarrow t\in [1;-1]
\end{array}\mr}= -2\pi \int_1^{-1}\sqrt{1+t^2}\, \d t=\\=
{\smsize\begin{pmatrix}\text{меняем местами}\\
\text{пределы интегрирования}\end{pmatrix}}= 2\pi \int_{-1}^1
\sqrt{1+t^2}\, \d t=\\= {\smsize\begin{pmatrix}\text{интегрируем по частям}\\
\text{с возвращением к исходному}\\
\text{интегралу (как в примере \ref{ex-12.6.10})}\end{pmatrix}}=\\=
\frac{1}{2}\l t\sqrt{1+t^2}+\ln \ml t+\sqrt{1+t^2}\mr \r \, \Big|_{-1}^1=\\=
\frac{1}{2}\l 2\sqrt{2}+\ln \ml 1+\sqrt{2}\mr -\ln \ml -1+\sqrt{2}\mr \r=\\=
\sqrt{2}+\frac{1}{2}\ln \ml \frac{\sqrt{2}+1}{\sqrt{2}-1}\mr=
\sqrt{2}+\frac{1}{2}\ln \l \sqrt{2}+1 \r^2=\\= \sqrt{2}+\ln \l \sqrt{2}+1 \r
\end{multline*}\end{ex}

\begin{tm}\label{tm-16.6.4}
Площадь поверхности вращения, образованной параметризованной кривой
$$
D: \quad \begin{cases}{x=\ph(t)}\\{y^2+z^2=\chi(t)^2}\end{cases}, \quad  t\in
[\alpha;\beta] \label{16.6.3}
$$
вычисляется по формуле
$$
S= 2\pi \int_\alpha^\beta \chi(t)\cdot \sqrt{\l \ph'(t) \r^2+ \l \chi'(t)
\r^2}\, \d t \label{16.6.4}
$$
\end{tm}

\begin{ex}\label{ex-16.6.5} Найдем площадь поверхности,
образованной вращением астроиды
$$
\begin{cases}{x=a\cos^3 t}\\{y=a\sin^3 t}\end{cases}
$$
вокруг оси OX:

%\pucture{0pt}{0pt}{12-30.pcx}

\vglue100pt \noindent По формуле \eqref{16.6.4} получаем площадь половинки
равна
 \begin{multline*}
\frac{S}{2}=\\={\smsize\text{$2\pi \int_0^\frac{\pi}{2} a\sin^3 t\cdot \sqrt{\l
-3a\cos^2 t \sin t\r^2+ \l 3a\sin^2 t\cos t \r^2}\,
\d t$}}=\\
={\smsize\text{$2\pi \int_0^\frac{\pi}{2} a\sin^3 t\cdot
\sqrt{9a^2\cos^4 t \sin^2 t+ 9a^2\sin^4 t \cos^2 t}\, \d t$}}=\\
=6a^2\pi \int_0^\frac{\pi}{2}\sin^3 t\cdot |\cos t\sin t| \sqrt{\cos^2 t +
\sin^2 t}\, \d t=\\= 6a^2\pi \int_0^\frac{\pi}{2}\sin^3 t\cdot \cos t\sin t
\cdot 1 \cdot \, \d t=\\= 6a^2\pi \int_0^\frac{\pi}{2}\sin^4 t\cdot \cos t \,
\d t=\\= 6a^2\pi \int_0^\frac{\pi}{2}\sin^4 t \, \d \sin t= 6a^2\pi
\frac{\sin^5 t}{5}\, \Big|_0^\frac{\pi}{2} =\\= \frac{6}{5}a^2\pi
 \end{multline*}
 Отсюда полная площадь
$$
  S=\frac{12}{5}a^2\pi
$$
\end{ex}

\begin{ers} Найдите площади поверх\-ностей, образованных
вращением следующих кривых вокруг оси $OX$:
 \biter{
\item[1)] $x=a\cos t, \, y=b\sin t$

\item[2)] $x=a(t-\sin t), \, y=b(1-\cos t), \, t\in [0;2\pi]$

\item[3)] $y=\ch x=\frac{e^x+e^{-x}}{2}, x\in [a;b]$
 }\eiter\end{ers}

\end{multicols}\noindent\rule[10pt]{160mm}{0.1pt}

\subsection{Интегрирование кусочно-непрерывных и кусочно-гладких функций}
\label{SEC:kusochno-neprer-i-kusochno-differents-func}

Интегралы, с которыми нам придется иметь дело в дальнейшем, будут, как правило,
интегралами от непрерывных функций, возможно, вдоль гладких функций. Однако в
главе \ref{CH-trig-series}, где речь пойдет о рядах Фурье, а также в главе
\ref{ch-o(f(x))}, где будут изучаться асимптотики, нам придется интегрировать
функции из несколько более широкого класса, а именно, кусочно-непрерывные
функции, возможно, вдоль непрерывных кусочно-гладких функций. В этом пункте мы
поговорим о таких интегралах.

\paragraph{Кусочно-непрерывные функции.}

\bit{ \item[$\bullet$]\label{DEF:kusochno-nepr-func} Функция $f$ на отрезке
$[a,b]$ называется {\it
кусочно-непрерывной}\index{функция!кусочно-непрерывная}, если
 \bit{

\item[1)] на $[a,b]$ функция $f$ непрерывна во всех точках, кроме, возможно,
конечного набора;

\item[2)] в каждой точке разрыва $c\in \R$ функция $f$ имеет конечные левый и
правый предел
 \beq\label{14.3.1}
f(c-0)=\lim_{x\to c-0} f(x), \qquad f(c+0)=\lim_{x\to c+0} f(x)
 \eeq
 }\eit
 }\eit

\btm\label{TH:osn-svoistvo-kusoch-neprer-func} Функция $f$ кусочно-непрерывна
на отрезке $[a,b]$ тогда и только тогда, когда существует разбиение
$\tau=\{c_0,...,c_k\}$ этого отрезка такое, что на каждом отрезке
$[c_{i-1},c_i]$ можно определить непрерывную функцию $\ph_i$, совпадающую с $f$
на интервале $(c_{i-1},c_i)$:
 $$
f(x)=\ph_i(x),\qquad x\in(c_{i-1},c_i)
 $$
\etm \bpr Занумеруем по возрастанию точки разрыва функции $f$ и добавим к ним,
если этого не произошло сразу, концы отрезка $[a,b]$. У нас получится некоторое
разбиение:
$$
  a=c_0<c_1<...<c_k=b
$$
На всяком отрезке $[c_{i-1},c_i]$ формула
$$
\ph_i(x)=\begin{cases} f(x), \quad x\in (c_{i-1},c_i)
\\
f(c_{i-1}+0), \quad x=c_{i-1}\\
f(c_i-0), \quad x=c_i
\end{cases}
$$
определяет непрерывную функцию $\ph_i$, совпадающую с $f$ на интервале
$(c_{i-1},c_i)$. \epr

\noindent\rule{160mm}{0.1pt}\begin{multicols}{2}

\begin{er}\label{ex-14.3.4} Проверьте, что функция
$$
  f(x)=\sgn \sin x
$$
имеющая график

%\pucture{0pt}{0pt}{ii-23.pcx}

\vglue100pt \noindent является кусочно-непрерывной на любом отрезке $[a,b]$.
\end{er}

\ber Наоборот, функция
$$
  f(x)=\begin{cases}\sgn \sin \frac{1}{x}, \quad x\ne 0 \\ 0,\quad x=0
\end{cases}
$$
с графиком

%\pucture{0pt}{0pt}{ii-23-1.pcx}

\vglue100pt \noindent не является кусочно-непрерывной, например, на отрезке
$[0,1]$, потому что имеет на нем бесконечное число точек разрыва
($x=\frac{1}{\pi n}$). \eer

\ber Функция
$$
  f(x)=\begin{cases}\frac{1}{x}, \quad x\ne 0 \\ 0,\quad x=0 \end{cases}
$$
с графиком

%\pucture{0pt}{0pt}{ii-23-2.pcx}

\vglue100pt \noindent тоже не будет кусочно-непрерывной на отрезке $[0,1]$,
потому что имеет бесконечный предел в точке 0.
\end{er}

\end{multicols}\noindent\rule[10pt]{160mm}{0.1pt}

\begin{tm}\label{tm-14.3.5} Всякая кусочно-непрерывная функция
$f$ на отрезке $[a,b]$ интегрируема на нем.
\end{tm}

Для доказательства теоремы \ref{tm-14.3.5} нам понадобится следующая

\begin{lm}\label{lm-14.7.1} Пусть функции
$f$ и $\ph$ определены на отрезке $[a,b]$ и совпадают на интервале $(a,b)$:
$$
  f(x)=\ph (x), \qquad x\in (a,b)
$$
Тогда если $\ph$ интегрируема на отрезке $[a,b]$, то и $f$ интегрируема на
$[a,b]$, причем
$$
  \int_a^b f(x) \, \d x= \int_a^b \ph (x) \, \d x
\label{14.7.7}
$$
\end{lm}\begin{proof} Заметим сразу, что обе функции $f$ и $\ph$
должны быть ограничены на $[a,b]$: функция $\ph$ ограничена по теореме
\ref{tm-14.3.1}, как интегрируемая на $[a,b]$, а функция $f$ отличается от нее
не более чем в двух точках. Мы можем даже считать, что они ограничены одной
константой:
$$
\sup_{x\in[a,b]}|f(x)|\le M, \qquad \sup_{x\in[a,b]}|\ph(x)|\le M
$$
Теперь для произвольных разбиений $\tau$ отрезка $[a;b]$ и выделенных точек
$\xi_i$ мы получим:
 \begin{multline*}
\Bigg|\overbrace{\sum_{i=1}^{k} f (\xi_i)\cdot\Delta
x_i}^{\scriptsize\begin{matrix}\text{интегральная}\\ \text{сумма для
$f$}\end{matrix}}-\underbrace{\overbrace{\sum_{i=1}^{k}\ph (\xi_i) \cdot \Delta
x_i}^{\scriptsize\begin{matrix}\text{интегральная}\\ \text{сумма для
$\ph$}\end{matrix}}}_{\scriptsize\begin{matrix}
\phantom{\tiny\begin{matrix}\diam\tau\\ \downarrow \\ 0\end{matrix}}
\ \downarrow \ {\tiny\begin{matrix}\diam\tau\\ \downarrow \\ 0\end{matrix}}\\
 \int_a^b \ph (x) \ \d x\end{matrix}}\Bigg|=\\ =
 \Bigg|\bigg( f(\xi_1)\cdot\Delta x_1+
\underbrace{\sum_{i=2}^{k-1} f\l \xi_i\r \cdot \Delta x_i}
 \put(-18,-30){
 \text{\scriptsize сокращаются, поскольку $f(\xi)=\ph(\xi)$ при $\xi\in(a,b)$}}
  \put(-33,-30){
 \put(0,-2){\line(1,0){219}}
 \put(-2,0){$\uparrow\kern214pt\uparrow$}}
+ f\l \xi_{k}\r \cdot \Delta x_{k}\bigg)-\bigg( \ph(\xi_1)\cdot \Delta x_1+
\underbrace{\sum_{i=2}^{k-1}\ph \l \xi_i\r \cdot \Delta x_i}+ \ph\l \xi_{k}\r
\cdot \Delta
x_{k}\bigg)\Bigg|=\\
=\Bigg|\Big(f(\xi_1)-\ph(\xi_1)\Big)\cdot\Delta
x_1+\Big(f(\xi_k)-\ph(\xi_k)\Big)\cdot\Delta x_k\Bigg|\le\\ \le
\underbrace{\Big|f(\xi_1)-\ph(\xi_1)\Big|}_{\scriptsize\begin{matrix}
\text{\rotatebox{-90}{$\le$}}\\ |f(\xi_1)|-|\ph(\xi_1)| \\
\text{\rotatebox{-90}{$\le$}}\\
2M\end{matrix}}\cdot\underbrace{\Delta x_1}_{\scriptsize\begin{matrix}
\text{\rotatebox{-90}{$\le$}}\\ \diam\tau
\end{matrix}}+\underbrace{\Big|f(\xi_k)-\ph(\xi_k)\Big|}_{\scriptsize\begin{matrix}
\text{\rotatebox{-90}{$\le$}}\\ |f(\xi_k)|-|\ph(\xi_k)| \\
\text{\rotatebox{-90}{$\le$}}\\
2M\end{matrix}}\cdot\underbrace{\Delta x_k}_{\scriptsize\begin{matrix}
\text{\rotatebox{-90}{$\le$}}\\ \diam\tau
\end{matrix}}\le 4M\cdot\diam\ \tau\underset{\diam\tau\to 0}{\longrightarrow}0
\end{multline*}
Отсюда получаем:
$$
\sum_{i=1}^{k} f (\xi_i)\cdot\Delta x_i\underset{\diam\tau\to
0}{\longrightarrow}\int_a^b \ph (x) \, \d x
$$
Это значит, что функция $f$ интегрируема на $[a,b]$ и ее интеграл равен
$\int\limits_a^b \ph (x) \, \d x$.
\end{proof}

\begin{proof}[Доказательство теоремы \ref{tm-14.3.5}]
Подберем с помощью теоремы \ref{TH:osn-svoistvo-kusoch-neprer-func} разбиение
отрезка $[a,b]$
$$
  a=c_0<c_1<...<c_k=b
$$
такое, что на всяком интервале $(c_{i-1},c_i)$ функция $f$ совпадает с
некоторой функцией $\ph_i$, непрерывной на отрезке $[c_{i-1},c_i]$. Тогда по
лемме \ref{lm-14.7.1}, функция $f$ будет интегрируемой на всяком отрезке
$[c_{i-1},c_i]$. Значит, по свойству аддитивности интеграла ($2^0$ на
с.\pageref{additivnost-integrala}), $f$ интегрируема на отрезке
$[a,b]=[a,c_1]\cup [c_1,c_2]\cup ... \cup [c_n,b]$.
\end{proof}

\paragraph{Кусочно-гладкие функции.}

\bit{ \item[$\bullet$]\label{DEF:kusochno-gladkaya-func} Функция $f$ на отрезке
$[a,b]$ называется {\it кусочно-гладкой}\index{функция!кусочно-гладкая}, если
 \bit{
\item[1)] на отрезке $[a,b]$ она дифференцируема во всех точках, кроме,
возможно, конечного набора;

\item[2)] в каждой точке $c\in[a,b]$, где $f$ дифференцируема, ее производная
непрерывна,
 \beq\label{22.1.2-0}
f'(c)=\lim_{x\to c}f'(x)
 \eeq

\item[3)] в каждой точке $c\in[a,b]$, где $f$ не дифференцируема, она обладает
следующими свойствами:
 \bit{
\item[(i)] $f$ имеет конечные левый и правый пределы (либо один из них, если
$c$ совпадает с каким-то из концов отрезка $[a,b]$)
 \beq\label{22.1.2-1}
f(c-0)=\lim_{x\to c-0} f(x), \quad f(c+0)=\lim_{x\to c+0} f(x)
 \eeq

\item[(ii)] $f$ имеет конечные левую и правую производные (либо одну из них,
если $c$ совпадает с каким-то из концов отрезка $[a,b]$)
 \beq\label{22.1.3}
f'_{-}(c)=\lim_{x\to c-0}\frac{f(x)-f(c-0)}{x-c}, \quad f'_{+}(c)=\lim_{x\to
c+0}\frac{f(x)-f(c+0)}{x-c}
 \eeq

\item[(iii)] при приближении аргумента $x$ к $c$ справа или слева, значение
$f'(x)$ стремится соответственно к $f'_{-}(c)$ или $f'_{+}(c)$:
 \beq\label{22.1.3-*}
f'_{-}(c)=\lim_{x\to c-0} f'(x), \quad f'_{+}(c)=\lim_{x\to c+0}f'(x)
 \eeq

 }\eit
 }\eit
 }\eit

Следующее утверждение доказывается в точности, как теорема
\ref{TH:osn-svoistvo-kusoch-neprer-func}:

\btm\label{TH:osn-svoistvo-kusoch-gladkih-func} Функция $f$ является
кусочно-гладкой на отрезке $[a,b]$ тогда и только тогда, когда существует
разбиение $\tau=\{c_0,...,c_k\}$ этого отрезка такое, что на каждом отрезке
$[c_{i-1},c_i]$ можно определить гладкую функцию $\ph_i$, совпадающую с $f$ на
интервале $(c_{i-1},c_i)$:
 $$
f(x)=\ph_i(x),\qquad x\in(c_{i-1},c_i)
 $$
\etm

Следующее утверждение мы считаем очевидным:

\btm\label{TH:interg-ot-v'-dlya-kusoch-glad-v} Пусть $g$ -- непрерывная
кусочно-гладкая функция на отрезке $[a,b]$. Тогда
 \bit{
\item[(i)] произвольным образом доопределяя производную $g'$ функции $g$ в
точках недифференцируемости $g$ на $[a,b]$, мы получим кусочно-непрерывную
функцию на $[a,b]$;

\item[(ii)] для любой интегрируемой функции $f$ на $[a,b]$ интеграл
$$
\int_a^b f(x)\cdot g'(x)\ \d x
$$
не зависит от того, как доопределена функция $g'$ в точках недифференцируемости
$g$.
 }\eit
\etm

 \bit{
\item[$\bullet$] Как следствие, для любой интегрируемой функции $f$ и любой
непрерывной\footnote{Объяснение, почему функцию $g$ в определении интеграла
$\int_a^b f(x) \ \d g(x)$ нужно брать непрерывной кусочно-гладкой (а не просто
кусочно-гладкой) заключается в том, что только тогда этот интеграл можно
определить формулой \eqref{DEF:int_a^b-u(x)-d-v(x)-kusoch-gladkie}. В общем
случае эта конструкция называется {\it интегралом Римана-Стилтьеса} и
определяется она по аналогии с интегралом Римана но с заменой $\Delta
x_i=x_i-x_{i-1}$ на $\Delta g_i=g(x_i)-g(x_{i-1})$. Для сходимости частичных
сумм необходимо требовать, чтобы $g$ была {\it функцией ограниченной вариации}.
Подробности можно найти, например, в учебнике: А.Н.Колмогоров, С.В.Фомин,
Элементы теории функций и функционального анализа, М.: Наука, 1981.}
кусочно-гладкой функции $g$ на отрезке $[a,b]$ формула
 \beq\label{DEF:int_a^b-u(x)-d-v(x)-kusoch-gladkie}
\int_a^b f(x) \, \d g(x):= \int_a^b f(x) \cdot g'(x) \, \d x
 \eeq
однозначно определяет число, называемое {\it определенным интегралом от функции
$f$ вдоль функции $g$ на отрезке $[a,b]$}. В частном случае, когда $f=1$ эта
величина обозначается
 \beq\label{peremeshenie-kusoch-glad-funktsii}
 \int_a^b \d g(x):= \int_a^b g'(x) \, \d x
 \eeq
 }\eit

Целью наших рассмотрений здесь является обобщение теорем
\ref{TH:peremeshenie-gladkoi-funktsii} и \ref{TH:integrir-po-chastyam} об
интегрировании вдоль функции на случай, когда эта функция является непрерывной
кусочно-гладкой (а не просто гладкой).

\btm[{\bf о перемещении непрерывной кусочно-гладкой
функции}]\label{TH:peremeshenie-gladkoi-funktsii} Перемещение непрерывной
кусочно-гладкой функции $g$ на ориентированном отрезке $\overrightarrow{ab}$
равно интегралу от единицы вдоль $g$ по этому отрезку:
 \beq\label{peremeshenie-nepreryvnoi-kusochno-gladkoi-funktsii}
 \int_a^b \d g(x) = \int_a^b g'(x) \, \d x=g(x)\Big|_{x=a}^{x=b}
 \eeq
\etm
\begin{proof}
Понятно, что здесь достаточно считать, что $a<b$, потому что противоположный
случай получается умножением равенства
\eqref{peremeshenie-nepreryvnoi-kusochno-gladkoi-funktsii} на $-1$.

1. Пусть сначала функция $g$ дифференцируема везде внутри отрезка $[a,b]$ (а
недифференцируема может быть только на концах этого отрезка). Тогда для любых
точек $\alpha$ и $\beta$ таких, что
$$
a<\alpha<\beta<b,
$$
функция $v$ будет гладкой на отрезке $[\alpha,\beta]$, и значит к ним применима
теорема \ref{TH:peremeshenie-gladkoi-funktsii}:
$$
\int_\alpha^\beta \d g(x)= g(\beta)-g(\alpha)
$$
Переходя к пределам при $\alpha\to a$ и $\beta\to b$, мы получаем как раз
\eqref{peremeshenie-nepreryvnoi-kusochno-gladkoi-funktsii}:
$$
\underbrace{\int_\alpha^\beta \d g(x)}_{\scriptsize\begin{matrix}
\phantom{\eqref{DEF:int_a^b-u(x)-d-v(x)-kusoch-gladkie}}\
\text{\rotatebox{90}{$=$}}\ \eqref{DEF:int_a^b-u(x)-d-v(x)-kusoch-gladkie}\\
\int_\alpha^\beta  g'(x)\ \d x \\
\downarrow \\
\int_a^b g'(x)\ \d x \\
\phantom{\eqref{DEF:int_a^b-u(x)-d-v(x)-kusoch-gladkie}}\
\text{\rotatebox{90}{$=$}}\ \eqref{DEF:int_a^b-u(x)-d-v(x)-kusoch-gladkie}\\
\int_a^b \d g(x)
\end{matrix}}= \underbrace{g(\beta)-g(\alpha)}_{\scriptsize\begin{matrix}
\downarrow \\
g(b)-g(a) \\ \begin{pmatrix}\text{здесь используется} \\
\text{непрерывность $g$} \end{pmatrix}
\end{matrix}}
$$

2. В общем случае мы выбираем по теореме
\ref{TH:osn-svoistvo-kusoch-gladkih-func} разбиение отрезка $[a,b]$
$$
a=c_0<c_1<...<c_k=b
$$
так, чтобы внутри каждого отрезка $[c_{i-1},c_i]$ функция $g$ была
дифференцируема, и, по уже доказанному, мы получим:
$$
\int_a^b \d g(x)=\sum_{i=1}^k \int_{c_{i-1}}^{c_i} \d g(x)= \sum_{i=1}^k  g(x)
\, \Big|_{x=c_{i-1}}^{x=c_i}= g(x) \ \Big|_{x=b}^{x=a}
$$
 \end{proof}

\bcor\label{TH:stroenie-nepr-kusoch-glad-func} Функция $g:[a,b]\to\R$ является
непрерывной кусочно-гладкой тогда и только тогда, когда она представима в виде
интеграла с переменным верхним пределом от некоторой кусочно-непрерывной
функции $f:[a,b]\to\R$:
$$
g(x)=g(a)+\int_a^x f(t)\ \d t
$$
\ecor
 \bpr
В качестве $f$ нужно взять производную $g'$ функции $g$, доопределенную
произвольным образом в точках недифференцируемости $g$.
 \epr

\btm[{\bf об интегрировании по частям для непрерывных кусочно-гладких
функций}]\label{TH:integrir-po-chastyam-dlya-nepr-kus-glad-funk} Если $f$ и $g$
-- непрерывные кусочно-гладкие функции на ориентированном отрезке
$\overrightarrow{ab}$, то
 \beq\label{integrir-po-chastyam-dlya-nepr-kus-glad-funk}
\int_a^b f(x) \, \d g(x)= f(x)\cdot g(x) \, \Big|_{x=a}^{x=b}- \int_a^b g(x) \,
\d f(x)
 \eeq
 \etm
\begin{proof} Как и в доказательстве теоремы \ref{TH:peremeshenie-gladkoi-funktsii}, здесь достаточно считать, что $a<b$.

1. Пусть для начала функции $f$ и $g$ дифференцируемы везде внутри отрезка
$[a,b]$ (а недифференцируемы могут быть только на концах этого отрезка). Тогда
для любых точек $\alpha$ и $\beta$ таких, что
$$
a<\alpha<\beta<b,
$$
функции $f$ и $g$ будут гладкими на отрезке $[\alpha,\beta]$, и значит к ним
применима теорема \ref{tm-15.5.1}:
$$
\int_\alpha^\beta f(x) \, \d g(x)= f(x)\cdot g(x) \,
\Big|_{x=\alpha}^{x=\beta}- \int_\alpha^\beta g(x) \, \d f(x)
$$
Переходя к пределам при $\alpha\to a$ и $\beta\to b$ по формулам
\eqref{nepreryvnost-integrala-1} и \eqref{nepreryvnost-integrala-2}, мы
получаем как раз \eqref{integrir-po-chastyam-dlya-nepr-kus-glad-funk}:
$$
\underbrace{\int_\alpha^\beta f(x) \ \d g(x)}_{\scriptsize\begin{matrix}
\phantom{\eqref{DEF:int_a^b-u(x)-d-v(x)-kusoch-gladkie}}\
\text{\rotatebox{90}{$=$}}\ \eqref{DEF:int_a^b-u(x)-d-v(x)-kusoch-gladkie}\\
\int_\alpha^\beta f(x)\cdot g'(x)\ \d x \\
\downarrow \\
\int_a^b f(x)\cdot g'(x)\ \d x \\
\phantom{\eqref{DEF:int_a^b-u(x)-d-v(x)-kusoch-gladkie}}\
\text{\rotatebox{90}{$=$}}\ \eqref{DEF:int_a^b-u(x)-d-v(x)-kusoch-gladkie}\\
\int_a^b f(x) \ \d g(x)
\end{matrix}}= \underbrace{f(x)\cdot g(x) \,
\Big|_{x=\alpha}^{x=\beta}}_{\scriptsize\begin{matrix}
\downarrow \\
f(x)\cdot g(x) \, \Big|_{x=a}^{x=b} \\ \begin{pmatrix}\text{здесь используется} \\
\text{непрерывность $f$ и $g$} \end{pmatrix}
\end{matrix}}-\underbrace{\int_\alpha^\beta g(x) \, \d f(x)}_{\scriptsize\begin{matrix}
\phantom{\eqref{DEF:int_a^b-u(x)-d-v(x)-kusoch-gladkie}}\
\text{\rotatebox{90}{$=$}}\ \eqref{DEF:int_a^b-u(x)-d-v(x)-kusoch-gladkie}\\
\int_\alpha^\beta g(x)\cdot f'(x)\ \d x \\
\downarrow \\
\int_a^b g(x)\cdot f'(x)\ \d x \\
\phantom{\eqref{DEF:int_a^b-u(x)-d-v(x)-kusoch-gladkie}}\
\text{\rotatebox{90}{$=$}}\ \eqref{DEF:int_a^b-u(x)-d-v(x)-kusoch-gladkie}\\
\int_a^b g(x) \ \d f(x)
\end{matrix}}
$$

2. В общем случае мы выбираем по теореме
\ref{TH:osn-svoistvo-kusoch-gladkih-func} разбиение отрезка $[a,b]$
$$
a=c_0<c_1<...<c_k=b
$$
так, чтобы внутри каждого отрезка $[c_{i-1},c_i]$ функции $f$ и $g$ были
дифференцируемы, и, по уже доказанному, мы получим:
 \begin{multline*}
\int_a^b f(x) \, \d g(x)=\sum_{i=1}^k \int_{c_{i-1}}^{c_i} f(x) \, \d g(x)=
\sum_{i=1}^k f(x)\cdot g(x) \, \Big|_{x=c_{i-1}}^{x=c_i}- \sum_{i=1}^k
\int_{c_{i-1}}^{c_i} g(x) \, \d f(x)=\\= f(x)\cdot g(x) \, \Big|_{x=a}^{x=b}-
\int_a^b g(x) \, \d f(x)
 \end{multline*}
 \end{proof}

\section{Некоторые следствия формулы Ньютона-Лейбница}

\subsection{Первообразная на интервале}

Напомним, что понятие первообразной для функции, определенной на отрезке было
введено нами на с.\pageref{DEF:pervoobraznaya}. Для случая интервала
определение ничем не отличается:

  \bit{
\item[$\bullet$] Функция $F$ называется {\it первообразной для функции $f$ на
интервале $(\alpha,\beta)$}, если производная $F$ на этом интервале равна $f$:
 \beq\label{pervoobr-nepr-f-na-intervale}
F'(x)=f(x),\qquad x\in(\alpha,\beta).
 \eeq
  }\eit

\bprop\label{PROP:pervoobr-nepr-f-na-intervale} Для любой непрерывной функции
$f$ на интервале $(\alpha,\beta)$ и любой точки $c\in (\alpha;\beta)$ формула
$$
F(x)= \int_c^x f(t) \, \d t, \qquad x\in(\alpha,\beta)
$$
определяет первообразную функции $f$ на интервале $(\alpha,\beta)$.
  \eprop
\begin{proof} Здесь используется теорема \ref{tm-15.1.1} об интеграле с переменным верхним
пределом.

1. Зафиксируем $x_0\in(c,\beta)$. Взяв какое-нибудь $b\in (x_0,\beta)$, мы
можем считать, что переменная $x$ бегает по отрезку $[c,b]$, на котором $x_0$
будет внутренней точкой, поэтому производную функции $F$ в $x_0$ можно
вычислить по формуле \eqref{15.1.7}:
 \begin{multline}\label{15.2.2}
F'(x_0)=\lim_{x\to x_0}\frac{1}{x-x_0}\l\int_c^x f(t) \, \d t-\int_c^{x_0} f(t)
\, \d  t\r=\\= \lim_{x\to x_0}\frac{1}{x-x_0}\l\int_{[c,x]} f(t) \, \d
t-\int_{[c,x_0]} f(t) \, \d  t\r=\eqref{15.1.7}=f(x_0).
 \end{multline}
Более того, представив $F$ тем же интегралом с $x$ бегающим по отрезку $[c,b]$,
мы можем найти правую производную $F$ в точке $c$ по формуле \eqref{15.1.7+}:
 \begin{multline}\label{15.2.4}
F_+'(c)=\lim_{x\to c+0}\frac{F(x)-F(c)}{x-c}=\lim_{x\to
c+0}\frac{1}{x-c}\int_c^x f(t) \ \d t=\\=\lim_{x\to
c+0}\frac{1}{x-c}\int_{[c,x]} f(t) \ \d t=\eqref{15.1.7+}=f(c).
 \end{multline}

2. Если $\alpha<x_0<c$, то взяв какое-нибудь $a: \, \alpha<a<x_0$,  мы получим
 \begin{multline}\label{15.2.3}
F'(x_0)=\lim_{x\to x_0}\frac{1}{x-x_0}\l\int_c^x f(t) \, \d t-\int_c^{x_0} f(t)
\, \d  t\r=\eqref{raznost-aniz-int-v-R^1}=\lim_{x\to
x_0}\frac{1}{x-x_0}\l\int_{x_0}^x f(t) \, \d t\r=
\\=\eqref{raznost-aniz-int-v-R^1}=\lim_{x\to x_0}\frac{1}{x-x_0}\l\int_a^x f(t) \, \d t-\int_a^{x_0}f(t) \, \d
  t\r=\\=\lim_{x\to x_0}\frac{1}{x-x_0}\l\int_{[a,x]} f(t) \, \d t-\int_{[a,x_0]}f(t) \, \d
  t\r=\eqref{15.1.7}=f(x_0)
 \end{multline}
И представив $F$ тем же интегралом с $x$ бегающим по отрезку $[a,c]$, мы можем
найти левую производную $F$ в точке $c$ по формуле \eqref{15.1.7-}:
\begin{multline}\label{15.2.5}
F_-'(c)=\lim_{x\to c-0}\frac{F(x)-F(c)}{x-c}=\lim_{x\to c-0}
\frac{1}{x-c}\int_c^x f(t) \ \d t=\\=\lim_{x\to c-0}\frac{1}{c-x}\int_x^c f(t)
\ \d t= \eqref{15.1.7-}=f(c).
 \end{multline}

3. Теперь мы получаем, что формулы \eqref{15.2.4} и \eqref{15.2.5} вместе дают
\eqref{pervoobr-nepr-f-na-intervale} для случая $x_0=c$, а для всех остальных
точек $x_0\ne c$ то же самое утверждается в формулах \eqref{15.2.2} и
\eqref{15.2.3}.
\end{proof}

\subsection{Лемма Адамара}

\blm[Адамар] Для любой гладкой функции $f$ на отрезке $[a,b]$ и любой точки
$c\in[a,b]$ найдется непрерывная функция $g$ на $[a,b]$ такая, что
 \beq\label{Hadamard}
f(x)-f(c)=(x-c)\cdot g(x),\qquad x\in[a,b]
 \eeq
\elm
 \bpr
Положим
$$
g(x)=\int_0^1 f'\Big((x-c)\cdot s+c\Big)\ \d s,\qquad x\in[a,b]
$$
Из непрерывности функции $f'$ следует, что функция $g$ тоже непрерывна на
$[a,b]$. Действительно, пусть $\e>0$. По теореме Кантора
\ref{Kantor}, найдется $\delta>0$ такое, что для любых
$x,y\in[a,b]$
$$
|x-y|<\delta\qquad \Longrightarrow\qquad |f'(x)-f'(y)|<\e
$$
Поэтому справедлива цепочка:
$$
|x-y|<\delta
$$
$$
\Downarrow
$$
$$
\forall s\in[0,1]\qquad \Big|\big((x-c)\cdot s+c\big)-\big((y-c)\cdot
s+c\big)\Big|=\Big|(x-y)\cdot s\Big|=|x-y|\cdot |s|\le |x-y| <\delta
$$
$$
\Downarrow
$$
$$
\forall s\in[0,1]\qquad \Big|f'\big((x-c)\cdot s+c\big)-f'\big((y-c)\cdot
s+c\big)\Big|<\e
$$
$$
\Downarrow
$$
 \begin{multline*}
|g(x)-g(y)|=\left|\int_0^1f'\big((x-c)\cdot s+c\big)\ \d s
-\int_0^1f'\big((y-c)\cdot s+c\big)\ \d
s\right|=\\=\left|\int_{[0;1]}f'\big((x-c)\cdot s+c\big)\ \d s
-\int_{[0;1]}f'\big((y-c)\cdot s+c\big)\ \d s\right|=\\=
\left|\int_{[0;1]}\Big( f'\big((x-c)\cdot s+c\big)-f'\big((y-c)\cdot
s+c\big)\Big)\ \d s\right|\le\\ \le \int_{[0;1]}\underbrace{\Big|
f'\big((x-c)\cdot s+c\big)-f'\big((y-c)\cdot
s+c\big)\Big|}_{\scriptsize\begin{matrix}\text{\rotatebox{90}{$>$}}\\ \e
\end{matrix}}\ \d s< \int_{[0;1]}\e \ \d s=\e
 \end{multline*}
Остается проверить тождество \eqref{Hadamard}. Для любого $x>c$ мы получаем:
 \begin{multline*}
f(x)-f(c)=\int_c^x f'(t)\ \d t={\scriptsize \left|\begin{matrix} t=(x-c)\cdot
s+c \\ s\in[0,1] \\ \d t=(x-c)\ \d s \\ t\in\overrightarrow{cx}\
\Leftrightarrow \ s\in\overrightarrow{0;1}
\end{matrix}\right|}=\int_0^1 f'\Big((x-c)\cdot s+c\Big)\cdot(x-c)\ \d s=\\=
(x-c)\cdot\int_0^1 f'\Big((x-c)\cdot s+c\Big)\ \d s=(x-c)\cdot g(x)
 \end{multline*}
 \epr

\chapter{НЕСОБСТВЕННЫЕ ИНТЕГРАЛЫ}\label{CH-improper-integral}

Несобственный интеграл формализует в математическом анализе идею площади
неограниченной фигуры.

\noindent\rule{160mm}{0.1pt}\begin{multicols}{2}

\bex\label{EX:ploshad-1/sqrt(x)} В предыдущей главе мы изучали определенный
интеграл, геометрический смысл которого -- площадь криволинейной трапеции.
Рассмотрим теперь какую-нибудь неограниченную криволинейную трапецию (то есть
такую, которая не лежит ни в каком прямоугольнике)

%\pucture{0pt}{0pt}{ii-1.pcx}

\vglue120pt \noindent Ее площадь будет непонятно чему равна, потому
что определенный интеграл от неограниченной функции не существует (по
теореме \ref{tm-14.3.1}).

Тем более неожиданным должно быть наблюдение, что у некоторых
неограниченных криволинейных трапеций все-таки имеется конечная
площадь. Точнее, понятие площади можно обобщить таким образом, что
некоторые неограниченные криволинейные трапеции будут иметь
конечную площадь. Например, если Вы задумаетесь, чему должна быть
равна площадь криволинейной трапеции, ограниченной кривой
$y=\frac{1}{\sqrt{x}}$

%\pucture{0pt}{0pt}{ii-2.pcx}

\vglue120pt
 \noindent
 то Вас наверняка удивит наблюдение, что эта площадь конечна и равна 2 (хотя сама фигура
 неограничена).

Дело в том, что эту величину можно посчитать как предел площадей
ограниченных криволинейных трапеций:

%\pucture{0pt}{0pt}{ii-3.pcx}

\vglue120pt \noindent
 \begin{multline*}
S=\lim_{t\to +0}\int_t^1 \frac{1}{\sqrt{x}}\, \d x= \lim_{t\to +0} 2\sqrt{x}\,
\Big|_t^1=\\= \lim_{t\to +0}  \l 2-2\sqrt{t}\r=2
 \end{multline*}
\eex

\bex\label{EX:ploshad-1/x} С другой стороны, если вместо кривой
$y=\frac{1}{\sqrt{x}}$ взять ничем особенным, как будто, не отличающуюся от нее
кривую $y=\frac{1}{x}$,

%\pucture{0pt}{0pt}{ii-4.pcx}

\vglue120pt \noindent то окажется, что соответствующая площадь,
посчитанная тем же способом бесконечна:
 \begin{multline*}
S=\lim_{t\to +0}\int_t^1 \frac{1}{x}\, \d x= \lim_{t\to +0}  \ln x \,
\Big|_t^1=\\= \lim_{t\to +0}  \l \ln 1- \ln t \r=+\infty
 \end{multline*}
\eex

\begin{ex}\label{ex-17.2.1} Еще один пример -- площадь криволинейной трапеции,
ограниченной кривой $y=\frac{1}{1+x^2}$ на полуинтервале $[0;+\infty)$:

%\pucture{0pt}{0pt}{ii-5.pcx}

 \vglue120pt \noindent
Ее можно вычислить как предел площадей получающихся, когда $x$ бегает по
отрезкам $[0,t]$, при $t\to\infty$,

 \vglue120pt \noindent
и оказывается, что она равна $\frac{\pi}{2}$:
 \begin{multline*}
\int_0^{+\infty}\frac{\, \d x}{1+x^2}= \lim_{t\to +\infty}\int_0^t \frac{\, \d
x}{1+x^2}=\\=\lim_{t\to +\infty}  \arctg x \, \Big|_0^t=\lim_{t\to +\infty}
\arctg t=\frac{\pi}{2}
 \end{multline*}
\end{ex}

\end{multicols}\noindent\rule[10pt]{160mm}{0.1pt}

После этих примеров естественно спросить, какой должна быть неограниченная
криволинейная трапеция, чтобы ее площадь была конечна? И чему равны площади
конкретных неограниченных криволинейных трапеций? Об этом мы поговорим в
настоящей главе.

\section{Определения несобственного интеграла и его свойства}

Начнем со следующего определения.

 \bit{
\item[$\bullet$]\label{DEF:loc-int-func} Функция $f$ называется {\it локально
интегрируемой на множестве $E$}, если она определена на $E$ и интегрируема на
любом отрезке $[a,b]$, содержащемся в $E$.
 }\eit

\noindent\rule{160mm}{0.1pt}\begin{multicols}{2}

\bex Функция
$$
f(x)=\begin{cases}\frac{1}{x}, & x\ne 0\\ 0,& x=0
\end{cases}
$$
конечно, не интегрируема на отрезке $[0,1]$ (потому что не ограничена на нем).
Но она локально интегрируема на полуинтервале $(0,1]$ (потому что непрерывна на
$(0,1]$). \eex

\bex Функция Дирихле, которую мы определяли формулой \eqref{func-Dirichle},
$$
D(x)=\begin{cases}1,& x\in\Q \\ 0,& x\notin\Q \end{cases}
$$
не будет локально интегрируемой, например на $\R$, потому что она не
интегрируема ни на каком отрезке. \eex

\end{multicols}\noindent\rule[10pt]{160mm}{0.1pt}

\subsection{Виды несобственных интегралов}

\paragraph{Несобственный интеграл по конечному промежутку.}

\bit{

\item[$\bullet$] Пусть функция $f$ локально интегрируема на полуинтервале
$(a;b]$, где $a,b$ -- произвольные числа. Тогда предел
$$
\lim_{t\to a+0}\int_t^b f(x) \, \d x
$$
называется {\it несобственным интегралом от $f$ по конечному промежутку
$(a;b]$}\index{интеграл!несобственный!по конечному промежутку} и обозначается
$$
\int_a^b f(x) \, \d x
$$
Если этот предел существует и конечен, то говорят, что {\it
несобственный интеграл сходится}, а если не существует или
бесконечен, то говорят, что {\it несобственный интеграл
расходится}.

\item[$\bullet$] Аналогично, если $f$  локально интегрируема на полуинтервале
$[a;b)$, то предел
$$
\lim_{t\to b-0}\int_a^t f(x) \, \d x
$$
называется {\it несобственным интегралом от $f$ по конечному промежутку
$[a;b)$}\index{интеграл!несобственный!по конечному промежутку} и обозначается
$$
\int_a^b f(x) \, \d x
$$
и опять если этот предел существует и конечен, то говорят, что
{\it несобственный интеграл сходится}, а если не существует или
бесконечен, то говорят, что {\it несобственный интеграл
расходится}.
 }\eit

Рассмотрим примеры.

\noindent\rule{160mm}{0.1pt}\begin{multicols}{2}

\begin{ex}\label{ex-17.1.1}
 \begin{multline*}
\int_0^1 \frac{\, \d x}{x^2}= \lim_{t\to +0}\int_t^1 \frac{\, \d
x}{x^2}= \lim_{t\to +0}\l -\frac{1}{x}\r \, \Big|_{x=t}^{x=1}=\\=
\lim_{t\to +0}\l -1+\frac{1}{t}\r =+\infty
 \end{multline*}
Вывод: интеграл расходится.
\end{ex}

\begin{ex}\label{ex-17.1.2}
  \begin{multline*}
\int_0^1 \frac{\, \d x}{\sqrt{1-x^2}}= \lim_{t\to 1-0}\int_0^t
\frac{\, \d x}{\sqrt{1-x^2}}= \\=\lim_{t\to 1-0}  \arcsin x \;
\Big|_{x=0}^{x=t}= \lim_{t\to 1-0}  \l \arcsin t - 0\r =\\=
 \arcsin 1 =\frac{\pi}{2}
 \end{multline*}
Вывод: интеграл сходится и равен $\frac{\pi}{2}$.
\end{ex}

\begin{ex}\label{ex-17.1.3}
 \begin{multline*}
\int_0^1 \frac{\sin \frac{1}{x}\, \d x}{x^2}= \lim_{t\to +0}\int_t^1
\frac{\sin \frac{1}{x}\, \d x}{x^2}=\\= -\lim_{t\to +0}\int_t^1 \sin
\frac{1}{x}\, \d \frac{1}{x}= \lim_{t\to +0}\cos \frac{1}{x}\,
\Big|_t^1 =\\= \lim_{t\to +0}\lll \cos 1-\cos \frac{1}{t}\rrr = \ml
\begin{array}{c}\frac{1}{t}=y
\\
t=\frac{1}{y}\\
t \to +0
\\
y\to +\infty
\end{array}\mr=\\= \cos 1-\kern-15pt\underbrace{\lim_{y\to +\infty}\cos y}_{
\text{предел не существует}}
 \end{multline*}
Вывод: интеграл расходится.
\end{ex}

\begin{ex}\label{ex-17.1.4}
  \begin{multline*}
\int_0^\frac{1}{2}\frac{\, \d x}{x\ln^2 x}= \lim_{t\to
+0}\int_t^\frac{1}{2}\frac{\, \d x}{x\ln^2 x}=\\= \lim_{t\to
+0}\int_t^\frac{1}{2}\frac{\, \d (\ln x)}{\ln^2 x}= -\lim_{t\to
+0}\frac{1}{\ln x}\, \Big|_t^\frac{1}{2}=\\= \lim_{t\to +0}\l
\frac{1}{\ln t}-\frac{1}{\ln \frac{1}{2}}\r= \frac{1}{\ln 2}
 \end{multline*}
Вывод: интеграл сходится и равен $\frac{1}{\ln 2}$.
\end{ex}

\begin{ers} Вычислите интегралы или установите их расходимость:
 \biter{
\item[1)] $\int_0^\frac{\pi}{4}\ctg x \, \d x$

\item[2)] $\int_0^1 \frac{x \, \d x}{\sqrt{1-x^2}}$

\item[3)] $\int_0^4 \frac{\, \d x}{x+\sqrt{x}}$

\item[4)] $\int_0^\frac{1}{2}\frac{\, \d x}{x\ln x}$

\item[5)] $\int_{-1}^0 e^\frac{1}{x}\frac{\, \d x}{x^2}$

\item[6)] $\int_0^1 e^\frac{1}{x}\frac{\, \d x}{x^2}$

\item[7)] $\int_0^\frac{\pi}{4}\frac{\sin x+\cos x}{\sqrt[3]{\sin x-\cos x}}$
 }\eiter
\end{ers}
\end{multicols}\noindent\rule[10pt]{160mm}{0.1pt}

\paragraph{Несобственный интеграл по бесконечному промежутку.}

Выше мы определили несобственный интеграл {\it по конечному промежутку},
формализующий понятие площади криволинейной трапеции специального вида:
основанием такой трапеции является конечный полуинтервал, а ``по вертикали''
эта фигура может быть неограничена (такие фигуры мы описывали в примерах
\ref{EX:ploshad-1/sqrt(x)} и \ref{EX:ploshad-1/x}).

Но помимо таких криволинейных трапеций можно рассматривать другие, основаниями
которых являются бесконечные полуинтервалы, и которые поэтому ``неограничены по
горизонтали'' (эту ситуацию мы описывали в примере \ref{ex-17.2.1}).
Соответствующий математический объект, формализующий это понятие также
называется неопределенным интегралом, но уже {\it по бесконечному промежутку}.

\bit{

\item[$\bullet$] Пусть функция $f$  локально интегрируема на полуинтервале
$[a;+\infty)$, где $a$ -- произвольное число. Тогда предел
$$
\lim_{t\to +\infty}\int_a^t f(x) \, \d x
$$
называется {\it несобственным интегралом от $f$ по бесконечному промежутку
$[a;+\infty)$}\index{интеграл!несобственный!по бесконечному промежутку} и
обозначается
$$
\int_a^{+\infty} f(x) \, \d x
$$
Если этот предел существует и конечен, то говорят, что {\it
несобственный интеграл сходится}, а если не существует или
бесконечен, то говорят, что {\it несобственный интеграл
расходится}.

\item[$\bullet$] Аналогично, если $f$  локально интегрируема на полуинтервале
$(-\infty;a]$, то предел
$$
\lim_{t\to -\infty}\int_t^b f(x) \, \d x
$$
называется {\it несобственным интегралом от $f$ по бесконечному промежутку
$(-\infty;a]$}\index{интеграл!несобственный!по бесконечному промежутку} и
обозначается
$$
\int_{-\infty}^b f(x) \, \d x
$$
и опять если этот предел существует и конечен, то говорят, что
{\it несобственный интеграл сходится}, а если не существует или
бесконечен, то говорят, что {\it несобственный интеграл
расходится}.
 }\eit

Рассмотрим примеры.

\noindent\rule{160mm}{0.1pt}\begin{multicols}{2}

\begin{ex}\label{ex-17.2.2}
 \begin{multline*}
\int_1^{+\infty}\frac{\, \d x}{x}= \lim_{t\to +\infty}\int_1^t
\frac{\, \d x}{x}=\\= -\lim_{t\to +\infty}\ln x\, \Big|_1^t =
-\lim_{t\to +\infty}  \ln t  =\infty
 \end{multline*}
Вывод: интеграл расходится.
\end{ex}

\begin{ex}\label{ex-17.2.3}
 \begin{multline*}
\int_{-\infty}^0 e^x \, \d x= \lim_{t\to -\infty}\int_t^0 e^x \, \d
x= \lim_{t\to -\infty}  e^x \, \Big|_t^0=\\= \lim_{t\to -\infty}\l
1-e^t \r=1
 \end{multline*}
Вывод: интеграл сходится и равен 1.
\end{ex}

\begin{ex}\label{ex-17.2.4}
 \begin{multline*}
\int_0^{+\infty}\cos x \, \d x= \lim_{t\to +\infty}\int_0^t \cos x \,
\d x=\\= -\lim_{t\to +\infty}  \sin x \, \Big|_0^t=
 -\kern-15pt\underbrace{\lim_{t\to+\infty}\sin t}_{\text{предел не существует}}
 \end{multline*}
Вывод: интеграл расходится.
\end{ex}

\begin{ers} Вычислите интегралы или установите их расходимость:
 \begin{multicols}{2}
1) $\int_2^{+\infty}\frac{\arctg x \, \d x}{x^2+1}$

2) $\int_0^{+\infty} e^x \, \d x$

3) $\int_{-\infty}^0 e^x \, \d x$

4) $\int_0^{+\infty} xe^x \, \d x$

5) $\int_e^{+\infty}\frac{\, \d x}{x\ln x}$

6) $\int_3^{+\infty}\frac{x\, \d x}{x^2-1}$

7) $\int_{-\infty}^0 \frac{\, \d x}{x^2+2x+2}$
\end{multicols}\end{ers}

\end{multicols}\noindent\rule[10pt]{160mm}{0.1pt}

\subsection{Несобственные интегралы от степенной и показательной функций}

Мы видели, что несобственный интеграл может сходиться, а может расходиться.
Здесь мы поймем, в каких случаях сходится интеграл от конкретных двух функций
-- $\frac{1}{x^\alpha}$ и $a^x$.

\noindent\rule{160mm}{0.1pt}\begin{multicols}{2}

\begin{tm}
[о несобственном интеграле по конечному промежутку от степенной
функции]\label{tm-17.3.1}
$$
\int_0^b \frac{\, \d x}{x^\alpha}\quad\longleftarrow\;
\begin{cases}{\text{сходится, если}\quad \alpha<1}\\
{\text{расходится, если}\quad \alpha\ge 1}\end{cases}
$$
\end{tm}\begin{proof} Рассмотрим три случая.

1) Если $\alpha<1$, то интеграл сходится:
 \begin{multline*}
\int_0^b \frac{\, \d x}{x^\alpha}= \int_0^b x^{-\alpha}\, \d x=
\lim_{t\to +0}\int_t^b x^{-\alpha}\, \d x=\\= \lim_{t\to
+0}\frac{x^{1-\alpha}}{(1-\alpha)}\, \Big|_{x=t}^{x=b}=\\= \lim_{t\to
+0}\lll \frac{b^{1-\alpha}}{(1-\alpha)} -
\frac{t^{1-\alpha}}{(1-\alpha)}\rrr =\\= {\smsize\begin{pmatrix}
1-\alpha<0
\\
\text{поэтому}\\
t^{1-\alpha}\to 0
\end{pmatrix}}= \frac{b^{1-\alpha}}{(1-\alpha)}
  \end{multline*}

2) Если $\alpha=1$, то интеграл расходится:
 \begin{multline*}
\int_0^b \frac{\, \d x}{x^\alpha}= \int_0^b \frac{\, \d x}{x}=
\lim_{t\to +0}\int_t^b \frac{\, \d x}{x}=\\= \lim_{t\to +0}\ln x \,
\Big|_{x=t}^{x=b}= \lim_{t\to +0}\lll \ln b-\ln t \rrr = \infty
  \end{multline*}

3) Если $\alpha>1$, то интеграл расходится:
 \begin{multline*}
\int_0^b \frac{\, \d x}{x^\alpha}= \int_0^b x^{-\alpha}\, \d x=
\lim_{t\to +0}\int_t^b x^{-\alpha}\, \d x=\\= \lim_{t\to
+0}\frac{x^{1-\alpha}\, \d x}{(1-\alpha)}\, \Big|_{x=t}^{x=b}=\\=
\lim_{t\to +0}\lll \frac{b^{1-\alpha}\, \d x}{(1-\alpha)} -
\frac{t^{1-\alpha}\, \d x}{(1-\alpha)}\rrr = {\smsize\begin{pmatrix}
1-\alpha<0
\\
\text{поэтому}\\
t^{1-\alpha}\to \infty
\end{pmatrix}}=\\= \infty
 \end{multline*}
 \end{proof}

\begin{tm}[\bf о несобственном интеграле по
бесконечному промежутку от степенной функции]\label{tm-17.3.2}
$$
\int_a^{+\infty}\frac{\, \d x}{x^\alpha}\quad\longleftarrow\;
\begin{cases}{\text{сходится, если}\quad \alpha>1}\\
{\text{расходится, если}\quad \alpha\le 1}\end{cases}
$$
\end{tm}\begin{proof} Рассмотрим три случая.

1) Если $\alpha>1$, то интеграл сходится:
 \begin{multline*}
\int_a^{+\infty}\frac{\, \d x}{x^\alpha}= \int_a^{+\infty}
x^{-\alpha}\, \d x= \lim_{t\to +\infty}\int_a^t x^{-\alpha}\, \d
x=\\= \lim_{t\to +\infty}\frac{x^{1-\alpha}\, \d x}{(1-\alpha)}\,
\Big|_{x=a}^{x=t}=\\= \lim_{t\to +\infty}\lll
\frac{t^{1-\alpha}}{(1-\alpha)} - \frac{a^{1-\alpha}}{(1-\alpha)}\rrr
= {\smsize\begin{pmatrix} 1-\alpha<0
\\
\text{поэтому}\\
t^{1-\alpha}\to 0
\end{pmatrix}}=\\= \frac{t^{1-\alpha}}{(1-\alpha)}
 \end{multline*}

2) Если $\alpha=1$, то интеграл расходится:
 \begin{multline*}
\int_a^{+\infty}\frac{\, \d x}{x^\alpha}= \int_a^{+\infty}\frac{\, \d
x}{x}=\\= \lim_{t\to +\infty}\int_a^t \frac{\, \d x}{x}= \lim_{t\to
+\infty}\ln x \, \Big|_{x=a}^{x=t}=\\= \lim_{t\to +\infty}\lll \ln
t-\ln a \rrr = \infty
 \end{multline*}

3) Если $\alpha<1$, то интеграл расходится:
 \begin{multline*}
\int_a^{+\infty}\frac{\, \d x}{x^\alpha}= \int_a^{+\infty}
x^{-\alpha}\, \d x=\\= \lim_{t\to +\infty}\int_a^t x^{-\alpha}\, \d
x= \lim_{t\to +\infty}\frac{x^{1-\alpha}\, \d x}{(1-\alpha)}\,
\Big|_{x=a}^{x=t}=\\= \lim_{t\to +\infty}\lll
\frac{t^{1-\alpha}}{(1-\alpha)} - \frac{a^{1-\alpha}}{(1-\alpha)}\rrr
= {\smsize\begin{pmatrix} 1-\alpha>0
\\
\text{поэтому}\\
t^{1-\alpha}\to \infty
\end{pmatrix}}= \infty
 \end{multline*}
 \end{proof}

\begin{tm}
[о несобственном интеграле от показательной
функции]\label{tm-17.3.3}
$$
\int_a^{+\infty} A^x \, \d x \quad\longleftarrow\;
\begin{cases}{\text{сходится, если}\quad 0<A<1}\\
{\text{расходится, если}\quad A\ge 1}\end{cases}
$$
\end{tm}\begin{proof} Рассмотрим три случая.

1) Если $0<A<1$, то интеграл сходится:
 \begin{multline*}
\int_a^{+\infty} A^x \, \d x= \lim_{t\to +\infty}\int_a^t A^x \, \d
x=\\= \lim_{t\to +\infty}\frac{A^x}{\ln A}\, \Big|_a^t = \frac{1}{\ln
A}\lim_{t\to +\infty} (A^t-A^a)=\\= {\smsize\begin{pmatrix} 0<A<1
\\
\text{поэтому}\\
A^t \to 0
\end{pmatrix}}= -\frac{A^a}{\ln A}
 \end{multline*}

2) Если $A=1$, то интеграл расходится:
 \begin{multline*}
\int_a^{+\infty} A^x \, \d x= \int_a^{+\infty} 1 \, \d x= \lim_{t\to
+\infty}\int_a^t 1 \, \d x=\\= \lim_{t\to +\infty}\frac{x}{\ln A}\,
\Big|_a^t = \lim_{t\to +\infty}\frac{t-a}{\ln A}=+\infty
 \end{multline*}

3) Если $A>1$, то интеграл расходится:
 \begin{multline*}
\int_a^{+\infty} A^x \, \d x= \lim_{t\to +\infty}\int_a^t A^x \, \d
x=\\= \lim_{t\to +\infty}\frac{A^x}{\ln A}\, \Big|_a^t = \frac{1}{\ln
A}\lim_{t\to +\infty} (A^t-A^a)=\\= {\smsize\begin{pmatrix} A>1
\\
\text{поэтому}\\
A^t \to \infty
\end{pmatrix}}= \infty
 \end{multline*}
 \end{proof}

\end{multicols}\noindent\rule[10pt]{160mm}{0.1pt}

\subsection{Замена переменной в несобственном интеграле}

Если функции $f$ и $g$ определены на полуинтервале $[a;b)$, то символом
$$
\int_a^b f(x) \, \d g(x)
$$
обозначается несобственный интеграл от функции $f \cdot g'$ по $[a;b)$:
$$
\int_a^b f(x) \, \d g(x):= \int_a^b f(x) \cdot g'(x) \, \d x=
\lim_{t\to b-0}\int_a^t f(x) \cdot g'(x) \, \d x
$$

\noindent\rule{160mm}{0.1pt}\begin{multicols}{2}

\begin{ex}\label{ex-17.4.1}
 \begin{multline*}
\int_1^{+\infty}\frac{1}{\sqrt{x}}\, \d \ln x=
\int_1^{+\infty}\frac{1}{x^\frac{3}{2}}\, \d x= \lim_{t\to +\infty}\int_1^t
\frac{1}{x^\frac{3}{2}}\, \d x=\\= \lim_{t\to +\infty}\l
-\frac{2}{x^\frac{1}{2}}\r \Big|_1^t= \lim_{t\to +\infty}\l
2-\frac{2}{t^\frac{1}{2}}\r=2
 \end{multline*}
\end{ex}

\end{multicols}\noindent\rule[10pt]{160mm}{0.1pt}

\begin{tm}[\bf о замене переменной в несобственном интеграле]\label{tm-17.4.1}
Пусть даны:
 \bit{
\item[1)] функция $f$, непрерывная на полуинтервале $[a; b)$

\item[2)] функция $\ph$, определенная на полуинтервале $[\alpha; \beta)$, со
следующими свойствами:
 \bit{
\item[a)] $\ph$ дифференцируема на полуинтервале $[\alpha; \beta)$ (то есть
$\ph$ дифференцируема на интервале $(\alpha; \beta)$, и в точке $\alpha$ имеет
левую производную);

\item[b)] $\forall t\in [\alpha; \beta) \qquad \ph (t)\in [a; b)$;

\item[c)] $\ph(\alpha)=a, \quad \ph(t)\underset{t\to
\beta-0}{\longrightarrow}b$.
 }\eit
 }\eit
Тогда
 \beq\label{17.4.1}
\int_a^b f (x) \, \d x= \int_\alpha^\beta f (\ph(t)) \, \d \ph(t)=
\int_\alpha^\beta f (\ph(t))\cdot \ph'(t) \, \d t
 \eeq
\end{tm}\noindent
 \bpr
 $$
\int_\alpha^\beta f (\ph(t))\, \d \ph(t)= \lim_{t\to \beta-0}\int_\alpha^t f
(\ph(t))\, \d \ph(t)=\eqref{15.4.1}= \lim_{t\to
\beta-0}\int_{\ph(\alpha)}^{\ph(t)}f (x) \, \d x= \int_a^b f (x) \, \d x
 $$
 \epr

\noindent\rule{160mm}{0.1pt}\begin{multicols}{2}

\begin{ex}\label{ex-17.4.3}
 \begin{multline*}
\int_2^{+\infty}\frac{\, \d x}{x\sqrt{x^2-1}}={\smsize \left|
\begin{array}{c}
x=\frac{1}{t},\qquad t=\frac{1}{x}\\
x=2 \Leftrightarrow t=\frac{1}{2}\\
x\to +\infty \Leftrightarrow t\to +0
\end{array}\right|}=\\= \int_\frac{1}{2}^0 \frac{-\, \d t}{\sqrt{1-t^2}}=
\int_0^\frac{1}{2}\frac{\, \d t}{\sqrt{1-t^2}}= \arcsin t
\Big|_0^\frac{1}{2}=\frac{\pi}{6}
 \end{multline*}
\end{ex}

\begin{ex}\label{ex-17.4.4}
 \begin{multline*}
\int_0^1 \frac{\, \d x}{(2-x)\sqrt{1-x}}={\smsize \left|
\begin{array}{c}
t=\sqrt{1-x},\qquad x=1-t^2\\
x=0 \Leftrightarrow t=1\\
x\to 1-0 \Leftrightarrow t\to +0
\end{array}\right|}=\\= \int_1^0 \frac{-t\, \d t}{t(t^2+1)}= \int_0^1 \frac{\, \d
t}{t^2+1}= \arctg t \Big|_0^1=\frac{\pi}{2}
 \end{multline*}
\end{ex}

\end{multicols}\noindent\rule[10pt]{160mm}{0.1pt}

\section{Признаки сходимости несобственных интегралов}

Не всякий несобственный интеграл можно явно вычислить, даже если
известно, что он сходится. Например, попробуйте найти
несобственный интеграл
 \beq
  \int_0^1 \frac{\sin^2 \frac{1}{x}}{\sqrt{x}}\, \d x
\label{17.5.1}
 \eeq
и Вы увидите, что это нелегко. Поэтому часто бывает важно просто понять,
сходится или нет данный несобственный интеграл, не вычисляя его. Например,
интеграл \eqref{17.5.1}, как будет показано в примере \ref{ex-17.5.2} сходится
(хотя и непонятно, чему равен).

В этом параграфе мы приведем некоторые признаки сходимости несобственных
интегралов. Впоследствии в \ref{SEC:asymp-integr} главы \ref{ch-o(f(x))} мы
дополним этот список еще двумя асимптотическими признаками.

\subsection{Признаки сходимости знакопостоянных интегралов}

 \bit{

\item[$\bullet$] Несобственный интеграл $\int_a^b f(x) \d x$ называется

 \bit{

\item[---] {\it знакопостоянным}\index{интеграл!несобственный!знакопостоянный},
если его подынтегральная функция $f$ не меняет знак на промежутке
интегрирования;

\item[---] {\it
знакоположительным}\index{интеграл!несобственный!знакоположительный}, если его
подынтегральная функция $f$ неотрицательная на промежутке интегрирования
$$
f(x)\ge 0,\qquad x\in(a,b)
$$
 }\eit
 }\eit

Понятно, что исследование на сходимость знакопостоянных интегралов сводится
исследованию знакоположительных: если $f(x)\le 0$, то можно взять функцию
$g(x)=-f(x)\ge 0$, и окажется, что
$$
\int_a^b f(x) \, \d x= \lim_{t\to b-0}\int_a^t f(x) \, \d x= -\lim_{t\to
b-0}\int_a^t g(x) \, \d x= -\int_a^b g(x) \, \d x
$$
откуда и будет следовать, что интеграл
$$
\int_a^b f(x) \, \d x
$$
сходится тогда и только тогда, когда сходится интеграл
$$
\int_a^b g(x) \, \d x
$$

\paragraph{Критерий сходимости знакоположительного интеграла.}

\begin{tm}\label{TH:int<infty}
Пусть функция $f$ локально интегрируема и неотрицательна на полуинтервале
$[a;b)$
$$
f(x)\ge 0,\qquad x\in[a,b)
$$
Тогда сходимость интеграла $\int_a^b f(x) \, \d x$ эквивалентна ограниченности
интегралов $\int_a^t f(x) \, \d x$, $t\in[a,b)$:
$$
\int_a^b f(x) \, \d x \quad \text{сходится}\quad \Longleftrightarrow \quad
\sup_{t\in[a,b)}\int_a^t f(x) \, \d x<\infty
$$
\end{tm}

\brem Условие справа (то есть утверждение, что интегралы $\int_a^b f(x) \, \d
x$ ограничены) принято записывать неравенством
$$
\int_a^b f(x) \, \d x<\infty
$$
и, в силу теоремы \ref{TH:int<infty}, такая запись считается эквивалентной
утверждению, что интеграл $\int_a^b f(x) \, \d x$ сходится (при неотрицательной
$f$).
 \erem

\bpr Обозначим
$$
F(t)=\int_a^t f(x) \, \d x
$$
и заметим, что это будет монотонно неубывающая функция, поскольку под
интегралом стоит неотрицательная функция. Тогда мы получим цепочку:
$$
\text{интеграл}\quad \int_a^b f(x) \, \d x \quad \text{сходится}
$$
$$
\Updownarrow
$$
$$
\text{$F$ имеет конечный предел $\lim_{t\to b-0} F(t)=\int_a^b f(x) \, \d x$}
$$
$$
\Updownarrow
$$
$$
\text{$F$ -- монотонная и ограниченная функция}
$$
$$
\Updownarrow
$$
$$
\sup_{t\in[a,b)}\int_a^t f(x) \, \d x=\sup_{t\in[a,b)} F(t)<\infty
$$
\epr

\paragraph{Признак сравнения интегралов.}

\begin{tm}[\bf признак сравнения несобственных интегралов]\label{tm-17.5.1}
Пусть $f$ и $g$ -- локально интегрируемые функции на полуинтервале $[a;c)$ (где
$c$ -- число или символ бесконечности $+\infty$), причем
 \beq\label{17.5.1-1}
0\le f(x)\le g(x),\quad x\in [a;c)
 \eeq
Соответствующая зависимость между несобственными интегралами
коротко записывается следующим образом:
$$
\int_a^c f(x) \, \d x \le \int_a^c g(x) \, \d x
$$
Тогда
 \bit{
\item[1)] из сходимости большего интеграла $\int_a^c g(x) \, \d x$
следует сходимость меньшего интеграла $\int_a^c f(x) \, \d x$:
$$
\int_a^c f(x) \, \d x \quad \text{сходится}\quad \Longleftarrow \quad
\int_a^c g(x) \, \d x \quad \text{сходится}  \quad
$$
\item[2)] из расходимости меньшего интеграла $\int_a^c f(x) \, \d x$
следует расходимость большего интеграла $\int_a^c g(x) \, \d x$:
$$
\int_a^c f(x) \, \d x \quad \text{расходится}\quad \Longrightarrow
\quad \int_a^c g(x) \, \d x \quad \text{расходится}  \quad
$$
 }\eit
\end{tm}\begin{proof} Обозначим
$$
  F(t)=\int_a^t f(x) \, \d x, \quad  G(t)=\int_a^t g(x) \, \d x
$$
тогда из \eqref{17.5.1-1} следует, что, во-первых,
$$
0\le F(t)\le G(t),\quad t\in [a;c) \label{17.5.2}
$$
и, во-вторых,
$$
F(t) \, \, \text{и}\, \, G(t) - \text{неубывающие функции на
промежутке}\,\, t\in [a;c) \label{17.5.3}
$$
потому что если $\alpha\le \beta$, то
$$
F(\alpha)=\int_a^\alpha f(x) \, \d x \le \int_a^\alpha f(x) \, \d
x+\int_\alpha^\beta f(x) \, \d x= \int_a^\beta f(x) \, \d
x=F(\beta)
$$
и аналогично,
$$
  G(\alpha)\le G(\beta)
$$

1. Предположим теперь, что интеграл
$$
\int_a^c g(x) \, \d x =\lim_{t\to c} G(t)
$$
сходится, то есть существует конечный предел
$$
\lim_{t\to c} G(t)=C
$$
Тогда, поскольку $G(t)$ -- неубывающая функция на промежутке $t\in
[a;c)$, она должна быть ограничена на этом промежутке:
$$
\sup_{t\in [a;c)} G(t)=C<+\infty
$$
Значит, функция $F(t)$ тоже должна быть ограничена на промежутке
$t\in [a;c)$:
$$
\sup_{t\in [a;c)} F(t)=B\le C<+\infty \label{17.5.4}
$$
Покажем, что
$$
\lim_{t\to c} F(t)=B \label{17.5.5}
$$
Действительно, возьмем произвольную последовательность $t_n\underset{n\to
\infty}{\longrightarrow} c$. Из \eqref{17.5.4} следует, что для всякого числа
$\varepsilon>0$ найдется такое $b\in [a;c)$, что
$$
B-\varepsilon<F(b)\le B
$$
При этом, поскольку $b<c$ и $t_n\underset{n\to
\infty}{\longrightarrow} c$, почти все числа $t_n$ должны быть
больше $b$:
$$
b<t_n \,\, \text{для почти всех}\,\, n\in \mathbb{N}
$$
$$
\phantom{\scriptsize(\text{вспоминаем, что $F$
неубывает})}\quad\Downarrow\quad{\scriptsize(\text{вспоминаем, что $F$
неубывает})}
$$
$$
F(b)<F(t_n) \,\, \text{для почти всех}\,\, n\in \mathbb{N}
$$
$$
\Downarrow
$$
$$
B-\varepsilon<F(b)<F(t_n)\le \sup_{t\in [a;c)} F(t)=B \,\, \text{для
почти всех}\,\, n\in \mathbb{N}
$$
$$
\Downarrow
$$
$$
F(t_n)\in (B-\varepsilon;B] \,\, \text{для почти всех}\,\, n\in
\mathbb{N}
$$
$$
\Downarrow
$$
$$
F(t_n)\in (B-\varepsilon;B+\varepsilon) \,\, \text{для почти
всех}\,\, n\in \mathbb{N}
$$
Последнее верно для любого числа $\varepsilon>0$, значит
$$
\lim_{n\to \infty} F(t_n)=B
$$
Это верно для всякой последовательности $t_n\underset{n\to
\infty}{\longrightarrow} c$, значит выполняется \eqref{17.5.5}.

Формула \eqref{17.5.5} означает, что существует конечный предел
$$
\int_a^c f(x) \, \d x =\lim_{t\to c} F(t)
$$
то есть несобственный интеграл $\int_a^c f(x) \, \d x$ сходится.

2. Мы доказали утверждение 2 (A). Утверждение 2 (B) является его
следствием: если интеграл $\int_a^c f(x) \, \d x$ расходится, то
интеграл $\int_a^c g(x) \, \d x$ не может сходиться (потому что
иначе мы получили бы что $\int_a^c f(x) \, \d x$ расходится, а
$\int_a^c g(x) \, \d x$ сходится, а это невозможно в силу уже
доказанного утверждения 2 (A)). \end{proof}

\noindent\rule{160mm}{0.1pt}\begin{multicols}{2}

\begin{ex}\label{ex-17.5.2} Чтобы понять, сходится ли интеграл
$$
  \int_0^1 \frac{\sin^2 \frac{1}{x}}{\sqrt{x}}\, \d x
$$
заметим, что подынтегральная функция положительна и меньше функции
$\frac{1}{\sqrt{x}}$:
$$
  0\le \frac{\sin^2 \frac{1}{x}}{\sqrt{x}}\le \frac{1}{\sqrt{x}}
$$
Отсюда получаем
$$
0\le \int_0^1 \frac{\sin^2 \frac{1}{x}}{\sqrt{x}}\, \d x \le
\int_0^1\frac{1}{\sqrt{x}}\, \d x
$$
причем больший интеграл сходится по теореме \ref{tm-17.3.1}.
Значит, по теореме \ref{tm-17.5.1}, наш исходный интеграл тоже
сходится.

Вывод: интеграл $\int_0^1 \frac{\sin^2 \frac{1}{x}}{\sqrt{x}}\, \d x$ сходится.
\end{ex}

\begin{ex}\label{ex-17.5.3} Для следующего интеграла
$$
  \int_0^1 \frac{1}{\sin x}\, \d x
$$
мы аналогичные рассуждения запишем более коротко:
$$
\int_0^1 \frac{1}{\sin x}\, \d x \ge \underbrace{\int_0^1 \frac{1}{x}\, \d
x}_{\scriptsize\begin{matrix}\text{расходится по}\\ \text{теореме
\ref{tm-17.3.1}}\end{matrix}}
$$

Вывод: интеграл $\int_0^1 \frac{1}{\sin x}\, \d x$ расходится.
\end{ex}

\begin{ex}\label{ex-17.5.4} Теперь рассмотрим интеграл по бесконечному
промежутку:
 \begin{multline*}
\int_\pi^{+\infty}\frac{3+\cos x}{x\sqrt{x}}\, \d x \le
\int_\pi^{+\infty}\frac{4}{x\sqrt{x}}\, \d x=\\=
4\underbrace{\int_\pi^{+\infty}\frac{1}{x^\frac{3}{2}}\, \d
x}_{\scriptsize\begin{matrix}\text{сходится по}\\ \text{теореме
\ref{tm-17.3.2}}\end{matrix}}
 \end{multline*}

Вывод: интеграл $\int_\pi^{+\infty}\frac{3+\cos x}{x\sqrt{x}}\, \d x$ сходится.
\end{ex}

\begin{ex}\label{ex-17.5.5} Снова рассмотрим интеграл по бесконечному
промежутку:
 \begin{multline*}
\int_\pi^{+\infty}\frac{3+\cos x}{\sqrt{x}}\, \d x \ge
\int_\pi^{+\infty}\frac{2}{\sqrt{x}}\, \d x=\\=
2\underbrace{\int_\pi^{+\infty}\frac{1}{x^\frac{1}{2}}\, \d
x}_{\scriptsize\begin{matrix}\text{расходится по}\\ \text{теореме
\ref{tm-17.3.2}}\end{matrix}}
 \end{multline*}

Вывод: интеграл $\int_\pi^{+\infty}\frac{3+\cos x}{\sqrt{x}}\, \d x$
расходится.
\end{ex}

\begin{ers} Исследуйте на сходимость интегралы:
 \begin{multicols}{2}
1) $\int_\pi^{+\infty}\frac{\sin^2 x}{x^2}\, \d x$;

2) $\int_e^{+\infty}\frac{\ln x}{\sqrt[3]{x}}\, \d x$;

3) $\int_2^{+\infty}\frac{\arctg x}{x\sqrt{x}}\, \d x$;

4) $\int_0^1 \frac{\cos \frac{1}{x}}{\sqrt[3]{x}}\, \d x$;

5) $\int_0^1 \frac{e^\frac{1}{x}}{x^3}\, \d x$;

6) $\int_{-1}^0 \frac{e^\frac{1}{x}}{x^3}\, \d x$
\end{multicols}\end{ers}

\end{multicols}\noindent\rule[10pt]{160mm}{0.1pt}

\subsection{Критерий Коши сходимости несобственного интеграла}

Теперь от знакопостоянных интегралов мы возвращаемся к произвольным.

\begin{tm}
[\bf критерий Коши сходимости несобственного
интеграла]\label{tm-17.6.1}\footnote{Этот результат понадобится ниже при
доказательстве признака абсолютной сходимости несобственного интеграла (теорема
\ref{tm-17.7.1})} Пусть функция $f$ локально интегрируема на полуинтервале
$[a;b)$. Тогда следующие условия эквивалентны:
 \bit{
\item[(i)]
несобственный интеграл
$$
\int_a^b f(x) \, \d x
$$
сходится;
\item[(ii)]
для любых двух числовых последовательностей $s_n, t_n\in [a;b)$,
стремящихся к $b$
$$
s_n\underset{n\to \infty}{\longrightarrow} b, \quad
t_n\underset{n\to \infty}{\longrightarrow} b,
$$
интеграл от функции $f$ по отрезку $[s_n;t_n]$ стремится к нулю при $n\to
\infty$:
$$
\int_{s_n}^{t_n} f(x) \, \d x \underset{n\to
\infty}{\longrightarrow} b
$$
 }\eit
\end{tm}\begin{proof} Рассмотрим функцию
$$
F(t)=\int_a^t f(x) \, \d x
$$
Тогда мы получим следующую логическую цепочку:
$$
\text{несобственный интеграл}\quad \int_a^b f(x) \, \d x \quad
\text{сходится}
$$
$$
\Updownarrow
$$
$$
\text{существует конечный предел}\quad \lim_{t\to b-0} F(t)
$$
$$
\Updownarrow\put(20,0){\smsize (\text{теорема \ref{Cauchy-crit-x->a-0}})}
$$
$$
\forall s_n\underset{n\to \infty}{\longrightarrow} b, \quad \forall
t_n\underset{n\to \infty}{\longrightarrow} b \qquad
\underbrace{F(t_n)-F(s_n)}_{\scriptsize\begin{matrix}
\text{\rotatebox{90}{$=$}} \\ \int_a^{t_n} f(x) \, \d x-\int_a^{s_n} f(x) \, \d
x \\ \text{\rotatebox{90}{$=$}} \\ \int_{s_n}^{t_n} f(x) \, \d x
\end{matrix}}\underset{n\to
\infty}{\longrightarrow} 0
$$
$$
\Updownarrow
$$
$$
\forall  s_n\underset{n\to \infty}{\longrightarrow} b, \quad \forall
t_n\underset{n\to \infty}{\longrightarrow} b \qquad \int_{s_n}^{t_n} f(x) \, \d
x \underset{n\to \infty}{\longrightarrow} 0
$$
\end{proof}

\noindent\rule{160mm}{0.1pt}\begin{multicols}{2}
\begin{ex}\label{ex-17.6.2} Критерий Коши обычно бывает удобен для
доказательства
расходимости несобственного интеграла (если он действительно
расходится). Рассмотрим, например, интеграл
$$
  \int_0^{+\infty}\sqrt{x}\sin x \, \d x
$$
Выберем следующие две числовые последовательности:
$$
  s_n=\frac{\pi}{6}+2\pi n, \quad t_n=\frac{5\pi}{6}+2\pi n
$$
Тогда
$$
  \forall x\in [s_n;t_n] \quad \sin x\ge \frac{1}{2}
$$
поэтому
 \begin{multline*}
  \int_{s_n}^{t_n}\sqrt{x}\sin x \, \d x\ge
\int_{s_n}^{t_n}\sqrt{x}\frac{1}{2}\, \d x=\\= \frac{1}{3}
x^\frac{3}{2}\Big|_{\frac{\pi}{6}+2\pi n}^{\frac{5\pi}{6}+2\pi n}=\\=
\frac{1}{3}\lll \l \frac{5\pi}{6}+2\pi n\r^\frac{3}{2} - \l
\frac{\pi}{6}+2\pi n\r^\frac{3}{2}\rrr
 \end{multline*}
Вычислим предел получившейся последовательности:
 \begin{multline*}
\lim_{n\to \infty}\lll \l \frac{5\pi}{6}+2\pi n\r^\frac{3}{2} - \l
\frac{\pi}{6}+2\pi n\r^\frac{3}{2}\rrr =\\= {\smsize \left|
\begin{array}{c}
t=\frac{1}{2\pi n}\\
n\to \infty \Leftrightarrow t\to +0
\end{array}\right|}=\\= \lim_{t\to +0}\lll \l
\frac{5\pi}{6}+\frac{1}{t}\r^\frac{3}{2} - \l
\frac{\pi}{6}+\frac{1}{t}\r^\frac{3}{2}\rrr =\\= \lim_{t\to
+0}\frac{\l \frac{5\pi t}{6}+1\r^\frac{3}{2} - \l \frac{\pi
t}{6}+1\r^\frac{3}{2} }{t^\frac{3}{2}}=\\= {\smsize \text{$\lim_{t\to
+0}\frac{\l 1+\frac{3}{2}\frac{5\pi t}{6}+ \underset{t\to
+0}{\bold{o}} (t)\r - \l 1+\frac{3}{2}\frac{\pi t}{6}+ \underset{t\to
+0}{\bold{o}} (t)\r }{t^\frac{3}{2}}$}}=\\= \lim_{t\to
+0}\frac{\frac{\pi t}{6}+\underset{t\to +0}{\bold{o}} (t)}
{t^\frac{3}{2}}=\\= \lim_{t\to +0}\frac{\frac{\pi}{6}+\underset{t\to
+0}{\bold{o}} (1)} {t^\frac{1}{2}}=+\infty
 \end{multline*}
Таким образом, мы получаем, что
 \begin{multline*}
\int_{s_n}^{t_n}\sqrt{x}\sin x \, \d x\ge\\ \ge \frac{1}{3}\lll \l
\frac{5\pi}{6}+2\pi n\r^\frac{3}{2} - \l \frac{\pi}{6}+2\pi
n\r^\frac{3}{2}\rrr \underset{n\to \infty}{\longrightarrow} +\infty
 \end{multline*}
Следовательно,
$$
\int_{s_n}^{t_n}\sqrt{x}\sin x \, \d x \underset{n\to \infty}{
\longrightarrow\kern-15pt{\Big/} } 0
$$
и по признаку Коши это означает, что наш интеграл расходится.

Вывод: интеграл $\int_0^{+\infty}\sqrt{x}\sin x \, \d x$ расходится.
\end{ex}

\begin{er}\label{er-17.6.3} Проверьте по критерию Коши, что интеграл
$$
  \int_0^{+\infty}\sqrt[3] {x}\cos x \, \d x
$$
расходится.
\end{er}

\end{multicols}\noindent\rule[10pt]{160mm}{0.1pt}

\subsection{Признак абсолютной сходимости}

\begin{tm}[\bf признак абсолютной сходимости]\label{tm-17.7.1}
Пусть функция $f$ локально интегрируема на полуинтервале $[a;b)$. Тогда
справедливы импликации:
$$
\int_a^b f(x) \, \d x \quad \text{сходится}\quad \Longleftarrow \quad \int_a^b
|f(x)| \, \d x \quad \text{сходится}\quad \Longleftarrow \quad
\sup_{t\in[a,b)}\int_a^t |f(x)| \, \d x<\infty
$$
\end{tm}\begin{proof}
$$
\text{интеграл}\quad \int_a^b |f(x)| \, \d x \quad \text{сходится}
$$
$$
\Downarrow\put(20,0){\smsize (\text{теорема \ref{tm-17.6.1}})}
$$
$$
\forall  s_n\underset{n\to \infty}{\longrightarrow} b, \quad \forall
t_n\underset{n\to \infty}{\longrightarrow} b \qquad
\underbrace{\int_{s_n}^{t_n} |f(x)| \, \d
x}_{\scriptsize\begin{matrix}\text{\rotatebox{90}{$\le$}}\\
\left|\int_{s_n}^{t_n} f(x) \, \d x \right|\\ \text{\rotatebox{90}{$\le$}}\\ 0
\end{matrix}}
\underset{n\to \infty}{\longrightarrow} 0
$$
$$
\Downarrow
$$
$$
\forall  s_n\underset{n\to \infty}{\longrightarrow} b, \quad \forall
t_n\underset{n\to \infty}{\longrightarrow} b \qquad \left|\int_{s_n}^{t_n} f(x)
\, \d x \right| \underset{n\to \infty}{\longrightarrow} 0
$$
$$
\Downarrow
$$
$$
\forall  s_n\underset{n\to \infty}{\longrightarrow} b, \quad \forall
t_n\underset{n\to \infty}{\longrightarrow} b \qquad
\int_{s_n}^{t_n} f(x) \, \d
x \underset{n\to \infty}{\longrightarrow} 0
$$
$$
\Downarrow\put(20,0){\smsize (\text{теорема \ref{tm-17.6.1}})}
$$
$$
\text{несобственный интеграл}\quad \int_a^b f(x) \, \d x \quad
\text{сходится}\qquad $$ \end{proof}

\noindent\rule{160mm}{0.1pt}\begin{multicols}{2}

\begin{ex}\label{ex-17.7.2} Чтобы исследовать на сходимость несобственный
интеграл
$$
  \int_0^1 \sin \frac{1}{x}\, \d x
$$
рассмотрим соответствующий интеграл от модуля
$$
0\le \int_0^1 \left|\sin \frac{1}{x}\right| \, \d x \le \int_0^1 1
\, \d x
$$
и применим признак сравнения для знакопостоянных интегралов
(теорема \ref{tm-17.5.1}): поскольку больший интеграл -- $\int_0^1
1\, \d x=1$ -- сходится, меньший -- $\int_0^1 \left|\sin
\frac{1}{x}\right| \, \d x$ -- тоже должен сходиться. Отсюда по
признаку абсолютной сходимости (теорема \ref{tm-17.7.1}) получаем

Вывод: Интеграл $\int_0^1 \sin \frac{1}{x}\, \d x$ сходится.
\end{ex}

\begin{ex}\label{ex-17.7.3} Интеграл
$$
  \int_1^{+\infty}\frac{\sin x}{x^2}\, \d x
$$
исследуется тем же способом:
$$
0\le \int_1^{+\infty}\left| \frac{\sin x}{x^2}\right| \, \d x \le
\int_1^{+\infty}\frac{1}{x^2}\, \d x
$$
По признаку сравнения для знакопостоянных интегралов (теорема
\ref{tm-17.5.1}), из сходимости большего интеграла --
$\int_1^{+\infty}\frac{1}{x^2}\, \d x$ -- следует сходимость меньшего
-- $\int_1^{+\infty}\left| \frac{sin x}{x^2}\right| \, \d x$. Отсюда
по признаку абсолютной сходимости (теорема \ref{tm-17.7.1})

Вывод: интеграл $\int_1^{+\infty}\frac{\sin x}{x^2}\, \d x$ сходится.
\end{ex}

\begin{ers} Исследуйте на сходимость интегралы:
 \biter{
\item[1)] $\int_0^{+\infty} e^{-x}\cos x \, \d x$;

\item[2)] $\int_2^{+\infty}\frac{\sqrt{x}\cos x}{x^2+x+\ln x}\, \d
x$;

\item[3)] $\int_2^{+\infty}\frac{(1+\sqrt{x}+\ln x)\sin
x}{x^2+x^3+x^4}\, \d x$;

\item[4)] $\int_0^1 \frac{\cos{\frac{1}{x}}}{\sqrt{x}+x+x^2}\, \d x$;

\item[5)] $\int_0^1 \frac{\cos{\frac{1}{x}}}{e^{\sqrt{x}}-1}\, \d x$;

\item[6)] $\int_0^1\frac{\sqrt{x}\sin{\frac{1}{x^2}}}{\ln(1+x+x^2)}\,
\d x$.
 }\eiter\end{ers}

\end{multicols}\noindent\rule[10pt]{160mm}{0.1pt}

\subsection{Формула Бонне и признаки Дирихле и Абеля для интегралов.}

\blm Если $f$ -- интегрируемая, а $g$ -- монотонная функции на ориентированном
отрезке $\overrightarrow{ab}$, то найдется точка $\xi\in\overrightarrow{ab}$
такая, что
 \beq\label{Bonnet}
\int_a^b f(x)\cdot g(x)\ \d x=g(a)\cdot \int_a^\xi f(x)\ \d x+g(b)\cdot
\int_\xi^b f(x)\ \d x
 \eeq
\elm
 \bit{
\item[$\bullet$] Равенство \eqref{Bonnet} называется {\it формулой Бонне}.
 }\eit
\bpr Пусть для начала $a<b$. Рассмотрим функцию
$$
F(t)=g(a)\cdot \int_a^t f(x)\ \d x+g(b)\cdot \int_t^b f(x)\ \d x
$$
По свойству непрерывности интеграла $6^0$ на
с.\pageref{nepreryvnost-integrala-1}, это будет непрерывная функция на отрезке
$[a,b]$. На концах этого отрезка она принимает такие значения:
$$
F(a)=g(b)\cdot \int_a^b f(x)\ \d x,\qquad F(b)=g(a)\cdot \int_a^b f(x)\ \d x
$$

Предположим теперь, что $g$ -- неубывающая функция на $[a,b]$. Тогда
$$
g(a)\le g(x)\le g(b),\qquad x\in[a,b]
$$
$$
\Downarrow
$$
$$
g(a)\cdot f(x)\le g(x)\cdot f(x)\le g(b)\cdot f(x),\qquad x\in[a,b]
$$
$$
\Downarrow
$$
$$
\underbrace{g(a)\cdot \int_a^b f(x)\ \d
x}_{\scriptsize\begin{matrix}\text{\rotatebox{90}{$=$}} \\ F(b)\end{matrix}}\le
\int_a^b g(x)\cdot f(x)\ \d x\le \underbrace{g(b)\cdot \int_a^b f(x)\ \d
x}_{\scriptsize\begin{matrix}\text{\rotatebox{90}{$=$}} \\ F(a)\end{matrix}}
$$
Таким образом, число $\int_a^b g(x)\cdot f(x)\ \d x$ лежит между двумя
значениями $F(a)$ и $F(b)$ функции $F$, непрерывной на отрезке $[a,b]$. По
теореме Коши о промежуточном значении \ref{Cauchy-I}, найдется точка
$\xi\in[a,b]$ такая, что
$$
F(\xi)=\int_a^b g(x)\cdot f(x)\ \d x,
$$
а это и есть формула \eqref{Bonnet}.

Если же функция $g$ невозрастает на $[a,b]$, то мы получим ту же цепочку
рассуждений, только в ней неравенства будут глядеть в противоположную сторону.

Остается рассмотреть случай, когда $a>b$. По уже доказанному, мы получим, что
для некоторой точки $\xi\in[b,a]$ должно выполняться равенство
$$
\int_b^a f(x)\cdot g(x)\ \d x=g(b)\cdot \int_b^\xi f(x)\ \d x+g(a)\cdot
\int_\xi^a f(x)\ \d x
$$
Умножив это на $-1$, получим:
$$
\underbrace{-\int_b^a f(x)\cdot g(x)\ \d x}_{\int_a^b f(x)\cdot g(x)\ \d
x}=\underbrace{-g(b)\cdot \int_b^\xi f(x)\ \d x}_{g(b)\cdot \int_\xi^b f(x)\ \d
x}\underbrace{-g(a)\cdot \int_\xi^a f(x)\ \d x}_{g(a)\cdot \int_a^\xi f(x)\ \d
x}
$$
-- и это будет формула \eqref{Bonnet} для случая $a>b$.
 \epr

 \begin{tm}[\bf признак Дирихле]\label{tm-17.7.11} Пусть даны:
 \bit{
\item[(i)] непрерывная функция $f:[a;+\infty)\to\R$ с ограниченной
первообразной на $[a;+\infty)$:
 \beq
  \sup_{t\ge a}\left| \int_a^t f(x) \, \d x \right|=M<\infty
\label{17.7.2}
 \eeq

\item[(ii)] гладкая функция $g:[a;+\infty)\to\R$, монотонно стремящаяся к нулю;
 \beq\label{17.7.3}
g(x)\overset{\text{\rm монотонно}}{\underset{x\to +\infty}{\longrightarrow}} 0
 \eeq
 }\eit
\noindent Тогда интеграл
$$
\int_a^{+\infty} f(x) \cdot g(x) \, \d x
$$
сходится.
 \end{tm}

\bpr\label{PROOF:Dirichlet-int} Обозначим
 $$
  C=\sup_{t\ge a}\left| \int_a^t f(x) \, \d x \right|<\infty
 $$
и заметим, что
 \beq\label{PROOF:Dirichlet-int-1}
\forall s,t\in[a,\infty)\qquad \left| \int_s^t f(x) \, \d x \right|\le 2C
 \eeq
Действительно,
$$
\bigg|\int_s^t f(x) \, \d x\bigg|=\bigg|\int_a^t f(x) \, \d x-\int_a^s f(x) \,
\d x\bigg|\le\underbrace{\bigg|\int_a^t f(x) \, \d x\bigg|}_{\scriptsize\begin{matrix}\text{\rotatebox{90}{$\ge$}}\\
C\end{matrix}}+\underbrace{\bigg|\int_a^s f(x) \, \d x\bigg|}_{\scriptsize\begin{matrix}\text{\rotatebox{90}{$\ge$}}\\
C\end{matrix}}\le 2C
$$
Чтобы убедиться, что интеграл $\int_a^\infty f(x)\cdot g(x) \ \d x$ сходится,
воспользуемся критерием Коши (теорема \ref{tm-17.6.1}): выберем произвольные
последовательности
$$
s_n\underset{n\to \infty}{\longrightarrow}+\infty, \quad t_n\underset{n\to
\infty}{\longrightarrow} +\infty.
$$
По формуле Бонне \eqref{Bonnet}, для любого $n$ найдется точка
$\xi_n\in\overrightarrow{s_n,t_n}$ такая, что
$$
\int_{s_n}^{t_n} f(x)\cdot g(x)\ \d x=g(s_n)\cdot \int_{s_n}^{\xi_n} f(x)\ \d
x+g(t_n)\cdot \int_{\xi_n}^{t_n} f(x)\ \d x
$$
Поскольку при этом $\xi_n\underset{n\to \infty}{\longrightarrow}+\infty$ (это
следует из условия $\xi_n\in\overrightarrow{s_n,t_n}$), мы получим:
 $$
\Big|\int_{s_n}^{t_n} f(x)\cdot g(x)\ \d x\Big|\le
\underbrace{|g(s_n)|}_{\scriptsize\begin{matrix}\downarrow
\\ 0\end{matrix}}\cdot
\underbrace{\Big|\int_{s_n}^{\xi_n} f(x)\ \d
x\Big|}_{\scriptsize\begin{matrix}\phantom{\tiny\eqref{PROOF:Dirichlet-int-1}}
\ \text{\rotatebox{90}{$\ge$}}\ {\tiny\eqref{PROOF:Dirichlet-int-1}}\\
2C\end{matrix}}+\underbrace{|g(t_n)|}_{\scriptsize\begin{matrix}\downarrow
\\ 0\end{matrix}}\cdot\underbrace{\Big|\int_{\xi_n}^{t_n}
f(x)\ \d
x\Big|}_{\scriptsize\begin{matrix}\phantom{\tiny\eqref{PROOF:Dirichlet-int-1}}
\ \text{\rotatebox{90}{$\ge$}}\ {\tiny\eqref{PROOF:Dirichlet-int-1}}\\
2C\end{matrix}}\underset{n\to \infty}{\longrightarrow} 0
 $$
Это верно для любых последовательностей $s_n$, $t_n$, стремящихся к $+\infty$.
Значит, интеграл $\int_a^\infty f(x)\cdot g(x) \ \d x$ сходится. \epr

\noindent\rule{160mm}{0.1pt}\begin{multicols}{2}

\begin{ex}\label{ex-17.7.12} Для исследования интеграла
$$
  \int_1^{+\infty}\frac{\sin x}{x}\, \d x
$$
можно воспользоваться признаком Дирихле, положив $f(x)=\sin x$ и
$g(x)=\frac{1}{x}$, и мы получим

Вывод: интеграл $\int_1^{+\infty}\frac{\sin x}{x}\, \d x$ сходится.
\end{ex}

\begin{ex}\label{ex-17.7.13} Интеграл
 \begin{multline*}
\int_1^{+\infty} x^2\cdot \sin x^2 \, \d x= {\smsize \left|
\begin{array}{c}
x^2=t, \quad x=\sqrt{t}\\
x=1 \Leftrightarrow t=1 \\
x\to +\infty \Leftrightarrow t\to +\infty
\end{array}\right|}=\\= \int_1^{+\infty} t\cdot \sin t \, \d \sqrt{t}=
\int_1^{+\infty}\frac{\sin t}{2\sqrt{t}}\, \d t
 \end{multline*}
тоже сходится по признаку Дирихле (с $f(t)=\sin t, \,
g(t)=\frac{1}{2\sqrt{t}}$).
\end{ex}

\begin{ex}\label{ex-17.7.14} Интеграл
 \begin{multline*}
\int_0^1 \frac{\cos \frac{1}{x}}{x\sqrt[3] {x}}\, \d x={\smsize
\left|
\begin{array}{c}\frac{1}{x}=t, \quad x=\frac{1}{t}\\
x=1 \Leftrightarrow t=1 \\
x\to +0 \Leftrightarrow t\to +\infty
\end{array}\right|}=\\= \int_{+\infty}^1 t \sqrt[3] {t}\cdot \cos t  \, \d \l
\frac{1}{t}\r= -\int_{+\infty}^1 \frac{\cos t}{t^\frac{2}{3}}  \, \d
t=\\= \int_1^{+\infty}\frac{\cos t}{t^\frac{2}{3}}  \, \d t
 \end{multline*}
тоже сходится по признаку Дирихле.
\end{ex}

\end{multicols}\noindent\rule[10pt]{160mm}{0.1pt}

\begin{tm}[\bf признак Абеля]\label{TH:priznak-Abelya-dlya-nesobst-int}
Пусть даны:
 \bit{
\item[$(i)$] сходящийся интеграл $\int_a^\infty f(x) \ \d x$ и

\item[$(ii)$] монотонная и ограниченная функция $g:[a,\infty)\to\R$.
 }\eit
Тогда интеграл $\int_a^\infty f(x)\cdot g(x) \ \d x$ сходится.
\end{tm}

\bpr\label{PROOF:Abel-int} Воспользуемся критерием Коши (теорема
\ref{tm-17.6.1}): выберем произвольные последовательности
$$
s_n\underset{n\to \infty}{\longrightarrow}+\infty, \quad t_n\underset{n\to
\infty}{\longrightarrow} +\infty.
$$
По формуле Бонне \eqref{Bonnet}, для любого $n$ найдется точка
$\xi_n\in\overrightarrow{s_n,t_n}$ такая, что
$$
\int_{s_n}^{t_n} f(x)\cdot g(x)\ \d x=g(s_n)\cdot \int_{s_n}^{\xi_n} f(x)\ \d
x+g(t_n)\cdot \int_{\xi_n}^{t_n} f(x)\ \d x
$$
Поскольку при этом $\xi_n\underset{n\to \infty}{\longrightarrow}+\infty$ (это
следует из условия $\xi_n\in\overrightarrow{s_n,t_n}$), мы получим:
$$
\int_{s_n}^{t_n} f(x)\cdot g(x)\ \d
x=\overbrace{g(s_n)}^{\text{ограничена}}\cdot\kern-15pt
\underbrace{\int_{s_n}^{\xi_n} f(x)\ \d
x}_{\scriptsize\begin{matrix}\phantom{\tiny\text{теорема \ref{tm-17.6.1}}}\
\downarrow \ {\tiny\text{теорема \ref{tm-17.6.1}}} \\
\phantom{,}0,\\ \text{поскольку} \\ \int_a^\infty f(x)\ \d x \\ \text{сходится}
\end{matrix}}\kern-15pt +\overbrace{g(t_n)}^{\text{ограничена}}\cdot\kern-15pt
\underbrace{\int_{\xi_n}^{t_n} f(x)\ \d
x}_{\scriptsize\begin{matrix}\phantom{\tiny\text{теорема \ref{tm-17.6.1}}}\
\downarrow \ {\tiny\text{теорема \ref{tm-17.6.1}}} \\
\phantom{,}0,\\ \text{поскольку} \\ \int_a^\infty f(x)\ \d x \\ \text{сходится}
\end{matrix}}\kern-15pt\underset{n\to \infty}{\longrightarrow} 0
$$
Это верно для любых последовательностей $s_n$, $t_n$, стремящихся к $+\infty$.
Значит, интеграл $\int_a^\infty f(x)\cdot g(x) \ \d x$ сходится. \epr

\noindent\rule{160mm}{0.1pt}\begin{multicols}{2}

\begin{ex} Для исследования интеграла
$$
  \int_1^{+\infty}\frac{\sin x}{x}\cdot\arctg x\, \d x
$$
вспомним, что в примере \eqref{ex-17.7.12} мы убедились, что интеграл
$$
  \int_1^{+\infty}\frac{\sin x}{x}\, \d x
$$
сходится. Поскольку функция $x\mapsto \arctg x$ монотонна, по признаку Абеля мы
получаем, что интеграл $\int_1^{+\infty}\frac{\sin x}{x}\cdot\arctg x\, \d x$
тоже сходится.
\end{ex}

\begin{ex} Для исследования интеграла
$$
\int_1^{+\infty} x^2\cdot \sin x^2\cdot\cos\frac{1}{x} \, \d x
$$
вспомним, что в примере \ref{ex-17.7.13} мы доказали сходимость интеграла
$$
\int_1^{+\infty} x^2\cdot \sin x^2 \, \d x
$$
Поскольку функция $x\mapsto\cos\frac{1}{x}$ монотонна, по признаку Абеля
интеграл $ \int_1^{+\infty} x^2\cdot \sin x^2\cdot\cos\frac{1}{x} \, \d x $
тоже сходится.
\end{ex}

\bex Следующий интеграл заменой переменных превращается в интеграл, который
можно исследовать по признаку Абеля:
 \begin{multline*}
\int_0^1\frac{\arccos x}{x}\cdot\cos\frac{1}{x}\ \d x=\left|\begin{matrix}
t=\frac{1}{x}\\ x=\frac{1}{t}\\ \d x=-\frac{\d t}{t^2}\end{matrix}\right|=\\=
-\int_{+\infty}^1 t\cdot\arccos\frac{1}{t}\cdot\cos t\cdot\frac{\d
t}{t^2}=\\=\int_1^\infty\frac{\cos t}{t}\cdot\arccos\frac{1}{t}\ \d t
 \end{multline*}
По признаку Дирихле интеграл $\int_1^\infty\frac{\cos t}{t}\ \d t$ сходится. А
функция $x\mapsto \arccos\frac{1}{t}$ монотонна и ограничена на $[0,+\infty)$.
Значит, по признаку Абеля, сходится интеграл $\int_1^\infty\frac{\cos
t}{t}\cdot\arccos\frac{1}{t}\ \d t$, а вместе с ним и интеграл
$\int_0^1\frac{\arccos x}{x}\cdot\cos\frac{1}{x}\ \d x$. \eex

\end{multicols}\noindent\rule[10pt]{160mm}{0.1pt}

\chapter{ЧИСЛОВЫЕ РЯДЫ}\label{CH-number-series}

Представьте себе, что у Вас есть числовая последовательность $\{ a_n \}$ и Вам
захотелось сосчитать (бесконечную) сумму ее элементов:
$$
\sum_{n=1}^{\infty} a_n=a_1+a_2+a_3+...
$$
В таких случаях говорят, что нужно вычислить сумму числового ряда. Интересное
наблюдение, с которого начинается вся теория рядов, состоит в том, что иногда
такая бесконечная сумма оказывается конечной величиной.

\noindent\rule{160mm}{0.1pt}\begin{multicols}{2}

\bex Например, нетрудно убедиться (это следует из приводимой ниже формулы
\eqref{18.1.1-1}), что
$$
\sum_{n=1}^{\infty}\frac{1}{2^n}=1
$$
Труднее доказать, что
 \beq\label{sum-1/n^2=pi^2/6}
\sum_{n=1}^{\infty}\frac{1}{n^2}=\frac{\pi^2}{6}
 \eeq
(мы это сделаем в замечании \ref{REM:sum-1/n^2} на с.\pageref{REM:sum-1/n^2}).
С другой стороны, сумма может оказаться и бесконечной, например,
$$
\sum_{n=1}^{\infty}\frac{1}{n}=\infty
$$
(это будет доказано в примере \ref{ex-18.3.2}).
 \eex

\end{multicols}\noindent\rule[10pt]{160mm}{0.1pt}

 \noindent
Заинтриговав читателя этими заявлениями,
мы можем перейти теперь к точным формулировкам.

\section{Определение числового ряда и его свойства}

\subsection{Определение числового ряда}

 \bit{

\item[$\bullet$] Пусть задана числовая последовательность $\{ a_n \}$ и из нее
составлена новая последовательность $\{ S_N \}$ по формуле
$$
S_N=\sum_{n=1}^N a_n=a_1+a_2+a_3+...+a_N
$$
Тогда такая пара последовательностей $\{ a_n \}$ и $\{ S_N \}$
называется {\it числовым рядом} и обозначается
$$
\sum_{n=1}^{\infty} a_n=a_1+a_2+a_3+...
$$
Числа $ a_n $ называются {\it слагаемыми} ряда
$\sum\limits_{n=1}^{\infty} a_n$,
а числа $ S_N $ -- {\it частичными суммами} этого ряда.

\item[$\bullet$] Если частичные суммы стремятся к некоторому конечному пределу
$$
S_N \underset{N\to \infty}{\longrightarrow} S, \qquad (S\in \R)
$$
то этот предел $S$ называется {\it суммой ряда}
$\sum\limits_{n=1}^{\infty} a_n$, а сам ряд называется
{\it сходящимся} (или про него говорят, что он {\it сходится}).
Коротко это записывают формулой
$$
\sum_{n=1}^{\infty} a_n=S
$$

\item[$\bullet$] Если же последовательность частичных сумм не имеет конечного
предела
$$
\not\exists \lim_{N\to \infty} S_N \in \R,
$$
то ряд $\sum\limits_{n=1}^{\infty} a_n$ называется {\it расходящимся}
(или про него говорят, что он {\it расходится}).

 }\eit

\noindent\rule{160mm}{0.1pt}\begin{multicols}{2}

\begin{ex}\label{ex-18.1.1} Рассмотрим ряд
$$
\sum_{n=1}^{\infty}\frac{1}{n(n+1)}
$$
Чтобы понять, сходится он или нет, заметим, что слагаемые можно разложить на
простейшие дроби:
$$
\frac{1}{n(n+1)}=\frac{1}{n}-\frac{1}{n+1}
$$
Теперь найдем частичные суммы:
\begin{multline*}
S_N=\sum_{n=1}^N \frac{1}{n(n+1)}=\sum_{n=1}^N \l
\frac{1}{n}-\frac{1}{n+1}\r=\\=
 \underbrace{
 \frac{1}{1}-
 \frac{1}{2}
    \put(-17,28){
 \text{\tiny сокращаются}
 \put(-31,-3){\line(1,0){18}}
 \put(-33,-10){$\downarrow\kern13pt\downarrow$}}
 }_{n=1}
 +
 \underbrace{
 \frac{1}{2}
 -
 \frac{1}{3}
    \put(27,28){
 \put(-31,-3){\line(1,0){18}}
 \put(-33,-10){$\downarrow\kern13pt\downarrow$}}
 }_{n=2}
 +
 \underbrace{
  \frac{1}{3}
  -
  \frac{1}{4}
    \put(27,28){
 \put(-31,-3){\line(1,0){15}}
 \put(-33,-10){$\downarrow\kern10pt\downarrow$}}
 }_{n=3}
 +
  ...
 \put(30,28){
 \put(-31,-3){\line(1,0){18}}
 \put(-33,-10){$\downarrow\kern13pt\downarrow$}}
 +
 \underbrace{
   \frac{1}{N}
 -
 \frac{1}{N+1}
 }_{n=N}
 =\\=1-\frac{1}{N+1}\underset{N\to \infty}{\longrightarrow} 1=S
 \end{multline*}

Вывод: ряд $\sum_{n=1}^{\infty}\frac{1}{n(n+1)}$ сходится, и его сумма равна 1.
\end{ex}

\begin{ex}\label{ex-18.1.2} Чтобы исследовать на сходимость ряд
$$
\sum_{n=1}^{\infty}\frac{1}{\sqrt{n}+\sqrt{n+1}}
$$
преобразуем его слагаемые:
\begin{multline*}\frac{1}{\sqrt{n}+\sqrt{n+1}}=
{\smsize\begin{pmatrix}\text{умножаем на}\\
\text{сопряженный радикал}\end{pmatrix}}=\\=
\frac{\sqrt{n}-\sqrt{n+1}}{(\sqrt{n}+\sqrt{n+1})(\sqrt{n}-\sqrt{n+1})}=\\=
\frac{\sqrt{n}-\sqrt{n+1}}{n-(n+1)}=
\frac{\sqrt{n}-\sqrt{n+1}}{-1}=\\=-\sqrt{n}+\sqrt{n+1}
  \end{multline*}
Теперь найдем частичные суммы:
\begin{multline*}
S_N=\sum_{n=1}^N \frac{1}{\sqrt{n}+\sqrt{n+1}} =\\
=\sum_{n=1}^N \l-\sqrt{n}+\sqrt{n+1}\r=\\=
 \underbrace{
  -\sqrt{1}
  +
  \sqrt{2}
 \put(-6,21){\put(2,0){\line(1,0){23}}
 \put(0,-7){$\downarrow$}\put(23,-7){$\downarrow$}}
  }_{n=1}
  \underbrace{
  -
  \sqrt{2}
  +
  \sqrt{3}
 \put(-6,21){\put(2,0){\line(1,0){23}}
 \put(0,-7){$\downarrow$}\put(23,-7){$\downarrow$}}
  }_{n=2}
  \underbrace{
  -
  \sqrt{3}
  +
  \sqrt{4}
 \put(-6,21){\put(2,0){\line(1,0){16}}
 \put(0,-7){$\downarrow$}\put(16,-7){$\downarrow$}}
  }_{n=3}
 +
 ... \\ ...
 \put(-3,21){\put(2,0){\line(1,0){21}}
 \put(0,-7){$\downarrow$}\put(21,-7){$\downarrow$}}
  \underbrace{
  -
  \sqrt{N}
  +
  \sqrt{N+1}
  }_{n=N}
 =
 -1+\sqrt{N+1}\underset{N\to \infty}{\longrightarrow} +\infty
 \end{multline*}
Вывод: ряд $\sum_{n=1}^{\infty}\frac{1}{\sqrt{n}+\sqrt{n+1}}$ расходится.
\end{ex}

\begin{ex}\label{ex-18.1.3} Рассмотрим ряд
$$
\sum_{n=1}^{\infty} (-1)^n
$$
Найдем первые несколько частичных сумм:
$$
  S_1=\sum_{n=1}^1 (-1)^n=(-1)^1=-1
$$
$$
  S_2=\sum_{n=1}^2 (-1)^n=(-1)^1+(-1)^2=-1+1=0
$$
 \begin{multline*}
  S_3=\sum_{n=1}^3 (-1)^n=(-1)^1+(-1)^2+(-1)^3=\\=-1+1-1=-1
 \end{multline*}
Теперь видна закономерность:
$$
S_N=\begin{cases}-1, & \text{если $N$ нечетное}\\
1, & \text{если $N$ четное}\end{cases}
$$
Ясно, что такая последовательность не имеет предела:
$$
  \not\exists \lim_{N\to \infty} S_N
$$

Вывод: ряд $\sum_{n=1}^{\infty} (-1)^n$ расходится.
\end{ex}

Следующий пример особенно важен.

\begin{tm}[\bf геометрическая прогрессия]\label{ex-18.1.4}
$$
\sum_{n=0}^\infty q^n \quad\leftarrow\;
\begin{cases}\text{сходится,} & \text{если $|q|<1$}\\
\text{расходится,} & \text{если $|q|\ge 1$}\end{cases}
$$
причем
 \begin{align}
& \sum_{n=0}^\infty q^n =\frac{1}{1-q}, \quad |q|<1 \label{18.1.1} \\
& \sum_{n=1}^\infty q^n =\frac{q}{1-q}, \quad |q|<1 \label{18.1.1-1}
 \end{align}
\end{tm}\begin{proof}
Здесь исполь\-зуется формула \eqref{Geom-progr-N} суммы конечного числа членов
геометрической прогрессии, которую мы доказали в главе \ref{ch-R&N}:
 \beq
S_N=\sum_{n=0}^N q^n =\frac{1-q^{N+1}}{1-q}, \quad q\ne 1
\label{18.1.2}
 \eeq
Рассмотрим три случая:

а) если $|q|<1$, то
$$
S_N=\frac{1-
\boxed{q^{N+1}}\put(0,16){\vector(1,2){6}\put(0,14){$0$}}}{1-q}\underset{N\to
\infty}{\longrightarrow}\frac{1}{1-q}
$$
и, значит в этом случае ряд сходится, причем мы сразу получаем формулу
\eqref{18.1.1}, из которой в свою очередь тут же следует формула
\eqref{18.1.1-1}:
$$
\sum_{n=1}^\infty q^n =q\cdot\sum_{n=1}^\infty q^{n-1}=
q\cdot\underbrace{\sum_{k=0}^\infty q^k}_{\scriptsize\begin{matrix}\text{\rotatebox{90}{$=$}}
\\
\frac{1}{1-q}\end{matrix}}=\frac{q}{1-q}
$$

б) если $|q|>1$, то
$$
S_N=\frac{1-\boxed{q^{N+1}}\put(0,16){\vector(1,2){6}\put(-2,15){$\infty$}}}{1-q}\underset{N\to
\infty}{\longrightarrow}\infty
$$
и, значит в этом случае ряд расходится.

в) если $|q|=1$, то это означает, что $q=1$ или $q=-1$,
и тогда
 \biter{
\item[---] при $q=1$ получается
$$
S_N=\sum_{n=0}^N 1=\boxed{N}\put(0,16){\vector(1,2){6}\put(-2,15){$\infty$}}+1
\underset{N\to \infty}{\longrightarrow}\infty;
$$

\item[---] а если $q=1$, то этот случай мы уже рассмотрели в примере
\ref{ex-18.1.3}, и показали, что получающийся ряд $\sum_{n=0}^\infty (-1)^n$
расходится.
 }\eiter
\end{proof}

\begin{ers} Исследуйте на сходимость ряды:
 \biter{
\item[1)] $\sum\limits_{n=1}^\infty \frac{1}{4n^2-1}$;

\item[2)] $\sum\limits_{n=1}^\infty \frac{2n+1}{n^2(n+1)^2}$;

\item[3)] $\sum\limits_{n=1}^\infty
(\sqrt{n}-2\sqrt{n+1}+\sqrt{n+2})$.
 }\eiter
 \end{ers}

\end{multicols}\noindent\rule[10pt]{160mm}{0.1pt}

\subsection{Арифметические свойства числовых рядов}

 \bit{\it
\item[$1^0$.] Если ряд $\sum\limits_{n=1}^\infty a_n$ сходится, то для всякого
числа $\lambda \in \R$ ряд $\sum\limits_{n=1}^\infty \lambda \cdot  a_n$ тоже
сходится, причем
 \beq
\sum\limits_{n=1}^\infty \lambda \cdot  a_n=
\lambda \cdot \sum\limits_{n=1}^\infty a_n
\label{18.2.1}
 \eeq
\item[$2^0$.]
Если ряды $\sum\limits_{n=1}^\infty a_n$
и $\sum\limits_{n=1}^\infty b_n$
сходятся,
то ряд $\sum\limits_{n=1}^\infty (a_n + b_n)$ тоже сходится, причем
 \beq
\sum\limits_{n=1}^\infty (a_n + b_n)=
\sum\limits_{n=1}^\infty a_n+\sum\limits_{n=1}^\infty b_n
\label{18.2.2}
 \eeq
 }\eit

\begin{proof}
1.
Если ряд $\sum\limits_{n=1}^\infty a_n$ сходится, то это означает, что
его частичные суммы $S_N=\sum\limits_{n=1}^N a_n$
стремятся к какому-то числу $S$:
$$
\sum\limits_{n=1}^N a_n=
S_N\underset{N\to \infty}{\longrightarrow} S
$$
Поэтому для всякого $\lambda \in \R$
$$
\sum\limits_{n=1}^N \lambda  a_n= \lambda \sum\limits_{n=1}^N  a_n=
\lambda S_N\underset{N\to \infty}{\longrightarrow}\lambda S
$$
Это означает, что ряд $\sum\limits_{n=1}^\infty \lambda a_n$
сходится, и его сумма равна $\lambda S=\lambda
\sum\limits_{n=1}^\infty a_n$, то есть справедлива формула
\eqref{18.2.1}.

2.
Пусть ряды $\sum\limits_{n=1}^\infty a_n$ и $\sum\limits_{n=1}^\infty b_n$
сходятся. Это означает, что если обозначить через $A_N$ и $B_N$
их частичные суммы, то они стремятся к каким-то числам $A$ и $B$:
$$
\sum\limits_{n=1}^N a_n=
A_N\underset{N\to \infty}{\longrightarrow} A, \quad
\sum\limits_{n=1}^N b_n=
B_N\underset{N\to \infty}{\longrightarrow} B
$$
Поэтому
$$
\sum\limits_{n=1}^N (a_n+b_n)=
\sum\limits_{n=1}^N a_n+\sum\limits_{n=1}^N b_n=
A_N+B_N\underset{N\to \infty}{\longrightarrow} A+B
$$
Это означает, что ряд $\sum\limits_{n=1}^\infty (a_n+b_n)$ сходится,
и его сумма равна $A+B=\sum\limits_{n=1}^\infty
a_n+\sum\limits_{n=1}^\infty b_n$, то есть справедлива формула
\eqref{18.2.2}.
\end{proof}

\noindent\rule{160mm}{0.1pt}\begin{multicols}{2}

\begin{ex}\label{ex-18.2.1}
 \begin{multline*}
\sum\limits_{n=1}^\infty \frac{5}{2^n}=
\frac{5}{2}\sum\limits_{n=1}^\infty \frac{1}{2^{n-1}}=
{\smsize\begin{pmatrix}\text{делаем замену}\\ n-1=k \\
\text{тогда}\,\, k=0,1,2,...
\end{pmatrix}}=\\=
\frac{5}{2}\cdot \sum\limits_{k=0}^\infty \frac{1}{2^k}=
\frac{5}{2}\cdot \frac{1}{1-\frac{1}{2}}=5
 \end{multline*}
\end{ex}

\begin{ex}\label{ex-18.2.2}
 \begin{multline*}\sum\limits_{n=1}^\infty \frac{2^n+3^n}{5^n}=
\sum\limits_{n=1}^\infty \lll \l \frac{2}{5}\r^n+\l \frac{3}{5}\r^n \rrr =\\=
\sum\limits_{n=1}^\infty \l \frac{2}{5}\r^n+ \sum\limits_{n=1}^\infty \l
\frac{3}{5}\r^n= \eqref{18.1.1-1}=\\=
\frac{1}{1-\frac{2}{5}}+\frac{1}{1-\frac{3}{5}}= \frac{5}{3}+\frac{5}{2}=
\frac{25}{6}\end{multline*}\end{ex}

\end{multicols}\noindent\rule[10pt]{160mm}{0.1pt}

\section{Признаки сходимости рядов}

Как правило, сумму ряда точно вычислить невозможно, даже если известно, что ряд
сходится. Например, формула, упомянутая нами в начале этой главы
$$
\sum_{n=1}^{\infty}\frac{1}{n^2}=\frac{\pi^2}{6},
$$
является результатом случайного наблюдения, а не какого-то найденного
математиками способа вычисления сумм, с помощью которого ее и другие подобные
формулы можно было бы выводить. Из-за этого в теории рядов становится важно
просто понять, сходится данный ряд или нет, не вычисляя его суммы.

Правила, объясняющие, когда данный ряд сходится, а когда нет, называются {\it
признаками сходимости}. В этом параграфе мы приведем некоторые признаки
сходимости рядов. В \ref{SEC:asymp-summ} главы \ref{ch-o(f(x))} мы дополним
этот список еще двумя асимптотическими признаками. Мы начнем со знакопостоянных
рядов.

\subsection{Сходимость знакопостоянных рядов}

 \bit{

\item[$\bullet$] Числовой ряд $\sum\limits_{n=1}^\infty a_n$ называется

 \bit{

\item[---] {\it знакопостоянным}\index{ряд!знакопостоянный}, если его члены
$a_n$ не меняют знак:
$$
\Big( \forall n\in\N a_n\ge0 \Big)\quad \text{или}\quad \Big( \forall n\in\N
a_n\le 0 \Big)
$$

\item[---] {\it знакоположительным}\index{ряд!знакоположительный}, если его
члены $a_n$ неотрицательны:
$$
\forall n\in\N a_n\ge0
$$
 }\eit
 }\eit

Как и в случае с несобственными интегралами, исследование на сходимость
знакопостоянных рядов сводится исследованию знакоположительных: если $a_n\le
0$, то можно взять последовательность $b_n=-a_n\ge 0$, и окажется, что
$$
\sum_{n=1}^\infty a_n= \lim_{N\to +\infty}\sum_{n=1}^N a_n= -\lim_{N\to
+\infty}\sum_{n=1}^N b_n= -\sum_{n=1}^\infty a_n
$$
откуда и будет следовать, что ряд
$$
\sum_{n=1}^\infty a_n
$$
сходится тогда и только тогда, когда сходится ряд
$$
\sum_{n=1}^\infty b_n
$$

\paragraph{Критерий сходимости знакоположительного ряда.}

\begin{tm}\label{TH:sum<infty}
Пусть последовательность $a_n$ неотрицательна
$$
a_n\ge 0,\qquad n\in\N
$$
Тогда сходимость ряда $\sum\limits_{n=1}^\infty a_n$ эквивалентна
ограниченности его частичных сумм $\sum\limits_{n=1}^N a_n$, $N\in\N$:
$$
\sum\limits_{n=1}^\infty a_n \quad \text{сходится}\quad \Longleftrightarrow
\quad \sup_{N\in\N}\sum\limits_{n=1}^N a_n <\infty
$$
\end{tm}

\brem Условие справа (то есть утверждение, что частичные суммы
$\sum\limits_{n=1}^N a_n$, $N\in\N$ ограничены) принято записывать неравенством
$$
\sum\limits_{n=1}^\infty a_n<\infty
$$
и, в силу теоремы \ref{TH:sum<infty}, такая запись считается эквивалентной
утверждению, что ряд $\sum\limits_{n=1}^\infty a_n$ сходится (при
неотрицательных $a_n$).
 \erem

\bpr Поскольку $a_n\ge 0$, частичные суммы
$$
S_N=\sum_{n=1}^N a_n
$$
образуют монотонно неубывающую последовательность. Поэтому справедлива цепочка
$$
\text{ряд}\quad \sum_{n=1}^\infty a_n \quad \text{сходится}
$$
$$
\Updownarrow
$$
$$
\text{$S_N$ имеет конечный предел $\lim_{N\to\infty} S_N=\sum_{n=1}^\infty
a_n$}
$$
$$
\Updownarrow
$$
$$
\text{$S_N$ -- монотонная и ограниченная последовательность}
$$
$$
\Updownarrow
$$
$$
\sup_{N\in\N} \sum_{n=1}^N a_n=\sup_{N\in\N} S_N<\infty
$$
\epr

\paragraph{Интегральный признак Коши и постоянная Эйлера.}

\begin{tm}[\bf интегральный признак Коши]\label{tm-18.3.1}
Пусть $f$ -- неотрицательная и невозрастающая функция на полуинтервале
$[1;+\infty)$. Тогда существует такое число $C\in\R$, что
 \beq\label{int-priznak-Cauchy}
\sum_{n=1}^Nf(n)-\int_1^N f(x)\ \d x\underset{N\to\infty}{\longrightarrow} C
 \eeq
 \bit{\rm
\item[$\bullet$] Число $C$ в этой формуле называется {\it постоянной
Эйлера}\label{post-Euler} ряда $\sum\limits_{n=1}^\infty f(n)$.
 }\eit\noindent
Как следствие,
$$
\text{ряд}\quad \sum_{n=1}^\infty f(n)\quad \text{сходится}\quad
\Longleftrightarrow \quad \text{интеграл}\quad \int_1^\infty f(x) \, \d x \quad
\text{сходится}
$$
\end{tm}
 \bit{
\item[$\bullet$] Такую связь между рядом и несобственным интегралом (когда
сходимость одного эквивалентна сходимости другого) коротко записывают формулой
$$
\sum_{n=1}^\infty f(n)\sim\int_1^\infty f(x) \, \d x.
$$
и говорят при этом, что ряд $\sum\limits_{n=1}^\infty f(n)$ {\it эквивалентен}
интегралу $\int\limits_1^\infty f(x) \, \d x$.
 }\eit

\begin{proof}
1. Рассмотрим сначала вспомогательную последовательность
$$
F_N=\sum_{n=1}^Nf(n)-\int_1^{N+1} f(x)\ \d x
$$
и покажем, что она сходится. Для этого нужно просто заметить, что она монотонна
и ограничена. Действительно, во-первых,
$$
\overbrace{f(n)\ge f(x)\ge f(n+1)}^{\forall x\in[n,n+1]}
$$
$$
\Downarrow
$$
 \beq\label{f(n)-ge-int_n^(n+1)f(x)-dx-ge-f(n+1)}
f(n)\ge \int_n^{n+1} f(x)\ \d x\ge f(n+1)
 \eeq
$$
\Downarrow
$$
$$
-f(n)\le -\int_n^{n+1} f(x)\ \d x\le -f(n+1)
$$
$$
\Downarrow
$$
 \beq\label{0-le-f(n)-int_n^(n+1)f(x)-dx-le-f(n)-f(n+1)}
0\le f(n)-\int_n^{n+1} f(x)\ \d x\le f(n)-f(n+1)
 \eeq
$$
\Downarrow
$$
 \begin{multline*}
0\le\overbrace{\sum_{n=1}^N\Big(f(n)-\int_n^{n+1} f(x)\ \d
x\Big)}^{\scriptsize\begin{matrix} F_N \\ \text{\rotatebox{90}{$=$}}
\\ \sum_{n=1}^Nf(n)-\int_1^{N+1} f(x)\ \d x\\
\text{\rotatebox{90}{$=$}}
\\
\sum_{n=1}^Nf(n)-\sum_{n=1}^N\int_n^{n+1} f(x)\ \d x\\
\text{\rotatebox{90}{$=$}}\end{matrix}}\le \sum_{n=1}^N\Big(
f(n)-f(n+1)\Big)=\\=\underbrace{f(1)-\overbrace{f(2)}}_{n=1}
\put(19.5,28){\put(-30.5,-3.5){\line(1,0){28}}\put(-33,-10){$\downarrow\kern23pt\downarrow$}}
+\underbrace{\overbrace{f(2)}-\overbrace{f(3)}}_{n=2}
\put(19.5,28){\put(-30.5,-3.5){\line(1,0){19}}\put(-33,-10){$\downarrow\kern14pt\downarrow$}}
+...
\put(30,28){\put(-30.5,-3.5){\line(1,0){20}}\put(-33,-10){$\downarrow\kern15pt\downarrow$}}
+\underbrace{\overbrace{f(N)}-f(N+1)}_{n=N}=f(1)-\underbrace{f(N+1)}_{\scriptsize\begin{matrix}
\text{\rotatebox{90}{$\le$}}
\\ 0\end{matrix}}\le f(1)
 \end{multline*}
$$
\Downarrow
$$
$$
0\le F_N\le f(1)
$$
И, во-вторых,
 \begin{multline*}
F_{N+1}-F_N=\l\sum_{n=1}^{N+1}f(n)-\int_1^{N+2} f(x)\ \d x\r
-\l\sum_{n=1}^Nf(n)-\int_1^{N+1} f(x)\ \d
x\r=\\=\l\sum_{n=1}^{N+1}f(n)-\sum_{n=1}^N f(n)\r-\l\int_1^{N+2} f(x)\ \d
x-\int_1^{N+1} f(x)\ \d x \r=f(N+1)-\int_N^{N+1} f(x)\ \d
x\overset{\eqref{0-le-f(n)-int_n^(n+1)f(x)-dx-le-f(n)-f(n+1)}}{\ge}0
 \end{multline*}
$$
\Downarrow
$$
$$
F_{N+1}\ge F_N
$$
Итак, $F_N$ монотонна и ограничена, и значит, имеет предел:
$$
F_N\underset{N\to\infty}{\longrightarrow} A
$$

2. Заметим далее, что последовательность $f(n)$ тоже монотонна и ограничена, и
значит тоже имеет предел:
$$
f(n)\underset{n\to\infty}{\longrightarrow} B
$$
Отсюда, в силу \eqref{f(n)-ge-int_n^(n+1)f(x)-dx-ge-f(n+1)}, следует, что
$$
\int_n^{n+1} f(x)\ \d x\underset{n\to\infty}{\longrightarrow} B
$$
(потому что $\int_n^{n+1} f(x)\ \d x$ заключена между двумя милиционерами,
стремящимися к $B$). Теперь получаем:
$$
\sum_{n=1}^N f(n)-\int_1^N f(x)\ \d x= \underbrace{\sum_{n=1}^N
f(n)-\int_1^{N+1} f(x)\ \d x}_{\scriptsize\begin{matrix}
\text{\rotatebox{90}{$=$}}
\\ F_N \\ \downarrow \\ A \end{matrix}}+\underbrace{\int_N^{N+1} f(x)\ \d x}_{\scriptsize\begin{matrix}
\downarrow \\ B \end{matrix}}\underset{N\to\infty}{\longrightarrow} A+B=C
$$
Это доказывает \eqref{int-priznak-Cauchy}.

3. Из \eqref{int-priznak-Cauchy} следует, что если ряд
$\sum\limits_{n=1}^\infty f(n)$ сходится, то есть
$\sum\limits_{n=1}^Nf(n)\underset{N\to\infty}{\longrightarrow}S\in\R$, то
$$
\int_1^N f(x)\ \d x= \underbrace{\sum_{n=1}^Nf(n)}_{\scriptsize\begin{matrix}
\downarrow \\ S \end{matrix}}-\underbrace{\l \sum_{n=1}^Nf(n)-\int_1^N f(x)\ \d
x\r}_{\scriptsize\begin{matrix} \downarrow \\ C
\end{matrix}}\underset{N\to\infty}{\longrightarrow} S-C
$$
Отсюда следует, что
$$
\int_1^y f(x)\ \d x\underset{y\to\infty}{\longrightarrow} S-C
$$
поскольку функция $F(y)=\int_1^y f(x)\ \d x$ монотонна (из-за неотрицательности
$f$). То есть, интеграл $\int_1^\infty f(x)\ \d x$ сходится.

4. Наоборот, если интеграл $\int_1^\infty f(x)\ \d x$ сходится, то есть
$\int_1^y f(x)\ \d x\underset{y\to\infty}{\longrightarrow} I\in\R$, то
$$
\int_1^N f(x)\ \d x\underset{N\to\infty}{\longrightarrow} I
$$
и поэтому
$$
\sum_{n=1}^Nf(n)=\underbrace{\l \sum_{n=1}^Nf(n)-\int_1^N f(x)\ \d
x\r}_{\scriptsize\begin{matrix} \downarrow \\ C
\end{matrix}}+\underbrace{\int_1^N f(x)\ \d
x}_{\scriptsize\begin{matrix} \downarrow \\ I
\end{matrix}}\underset{N\to\infty}{\longrightarrow} C-I
$$
То есть, ряд $\sum\limits_{n=1}^Nf(n)$ сходится.
 \end{proof}

\noindent\rule{160mm}{0.1pt}\begin{multicols}{2}

\begin{ex} {\bf Гармонический ряд.}\label{ex-18.3.2}
Числовой ряд
$$
\sum_{n=1}^\infty \frac{1}{n}
$$
называется {\it гармоническим}. Он расходится, потому что эквивалентен
расходящемуся интегралу:
$$
\sum_{n=1}^\infty \frac{1}{n}\sim \int_1^\infty \frac{1}{x}\, \d x= \ln x
\Big|_1^\infty=\infty
$$
Несмотря на то, что суммы у этого ряда нет, у него есть постоянная Эйлера:
 \beq\label{harmonic-Euler}
C=\lim_{N\to\infty}\bigg(\sum_{n=1}^N \frac{1}{n}-\underbrace{\ln
N}_{\scriptsize\begin{matrix}\text{\rotatebox{90}{$=$}}\\
\int\limits_1^N \frac{1}{x}\, \d x
\end{matrix}}\bigg)
 \eeq
Это число привлекает внимание специалистов по теории чисел тем, что за
прошедшие с его открытия два с половиной века до сих пор так и не удалось
понять, будет ли оно рациональным.
\end{ex}

\bex {\bf Ряд Дирихле.}\label{ex-18.3.3} Покажем, что
 \beq\label{ryad-Dirichlet}
\sum_{n=1}^\infty \frac{1}{n^\alpha}\quad\leftarrow\;
\begin{cases}\text{сходится}, &\text{если $\alpha>1$}\\
\text{расходится}, & \text{если $\alpha\le 1$}\end{cases}
 \eeq
Этот ряд называется {\it рядом Дирихле}. При $\alpha>0$ к нему применима
теорема \ref{tm-18.3.1}, и в этом случае ему соответствуют постоянные Эйлера
 \beq\label{Dirichlet-Euler}
C_{\alpha}=\lim_{N\to\infty}\bigg(\sum_{n=1}^N \frac{1}{n^{\alpha}}-
 \underbrace{\frac{N^{1-\alpha}-1}{1-\alpha}}_{\scriptsize
 \begin{matrix} \text{\rotatebox{90}{$=$}}\\
 \int\limits_1^N \frac{1}{x^{\alpha}}\, \d x \end{matrix}}
\bigg)
 \eeq
 \eex
 \bpr
Если $\alpha\le 0$, то общий член этого ряда не меньше единицы
$$
\frac{1}{n^\alpha}\ge 1,
$$
и тогда
$$
\sum_{n=1}^N \frac{1}{n^\alpha}\ge \sum_{n=1}^N
1=N\underset{N\to\infty}{\longrightarrow}\infty,
$$
и значит, ряд расходится. Если же $\alpha>0$, то можно применить теорему
\ref{tm-18.3.1}, и мы получим:
$$
\sum_{n=1}^\infty \frac{1}{n^\alpha}\sim \underbrace{\int_1^\infty
\frac{1}{x^\alpha}\, \d x}_{\scriptsize\begin{matrix}\text{сходится}\\
\text{только при  $\alpha>1$,}\\ \text{по теореме
\ref{tm-17.3.2}}\end{matrix}}.
$$
\epr

\begin{ex}\label{ex-18.3.4}
 \begin{multline*}
\sum_{n=2}^\infty \frac{1}{n \ln n}\sim \int_2^\infty \frac{1}{x\ln
x}\, \d x=\\= \int_2^\infty \frac{1}{\ln x}\, \d \ln x= \ln \ln x
\Big|_2^\infty=\infty
 \end{multline*}
Вывод: ряд $\sum_{n=2}^\infty \frac{1}{n \ln n}$ расходится.
\end{ex}

\begin{ex}\label{ex-18.3.5}
 \begin{multline*}
\sum_{n=2}^\infty \frac{1}{n \ln^2 n}\sim \int_2^\infty
\frac{1}{x\ln^2 x}\, \d x=\\= \int_2^\infty \frac{1}{\ln^2 x}\, \d
\ln x= -\frac{1}{\ln x}\Big|_2^\infty=\frac{1}{\ln 2}
 \end{multline*}
Вывод: ряд $\sum_{n=2}^\infty \frac{1}{n \ln^2 n}$ сходится.
\end{ex}

\begin{ers} Исследуйте на сходимость ряды:
 \begin{multicols}{2}
1) $\sum_{n=2}^\infty \frac{1}{n \ln^\alpha n}$

2) $\sum_{n=2}^\infty \frac{1}{n \ln n \ln^\alpha \ln n}$

3) $\sum_{n=1}^\infty \frac{e^{-\sqrt{n}}}{\sqrt{n}}$
\end{multicols}\end{ers}

\end{multicols}\noindent\rule[10pt]{160mm}{0.1pt}

\paragraph{Признак сравнения рядов.}

\begin{tm}[\bf признак сравнения рядов]\label{tm-18.3.7}
Пусть $0\le a_n\le b_n$.
Соответствующая зависимость между рядами коротко
записывается следующим образом:
 \beq\label{DEF:sravnenie-ryadov}
\sum_{n=1}^\infty a_n \le \sum_{n=1}^\infty b_n
 \eeq
Тогда
 \bit{
\item[1)] из сходимости большего ряда $\sum\limits_{n=1}^\infty b_n$
следует сходимость меньшего ряда $\sum\limits_{n=1}^\infty a_n$:
$$
\sum\limits_{n=1}^\infty a_n\quad \text{сходится}\quad \Longleftarrow
\quad \sum\limits_{n=1}^\infty b_n \quad \text{сходится}  \quad
$$
\item[2)] из расходимости меньшего ряда $\sum\limits_{n=1}^\infty a_n$
следует расходимость большего ряда $\sum\limits_{n=1}^\infty b_n$:
$$
\sum\limits_{n=1}^\infty a_n \quad \text{расходится}\quad
\Longrightarrow \quad \sum\limits_{n=1}^\infty b_n \quad
\text{расходится}  \quad
$$
 }\eit
\end{tm}\begin{proof}
Обозначим через $A_N$ и $B_N$ частичные суммы рядов
$\sum\limits_{n=1}^\infty a_n$ и $\sum\limits_{n=1}^\infty b_n$:
$$
  A_N=\sum\limits_{n=1}^N a_n, \qquad B_N=\sum\limits_{n=1}^N b_n
$$
Тогда
$$
\begin{array}{ccccccccccc}
B_1 & \le &  B_2 & \le & B_3 & \le & ... & \le & B_N & \le & ...\\
\VI &   & \VI &     & \VI &     &     &     & \VI &     &    \\
A_1 & \le & A_2 & \le  &  A_3& \le & ... & \le  & A_N & \le & ...
 \end{array}
$$
1. Теперь возникает логическая цепочка:
$$
\text{если ряд $\sum\limits_{n=1}^\infty b_n$ сходится, то}
$$
$$
\Downarrow
$$
$$
\text{существует конечный предел $\lim\limits_{N\to \infty} B_N=B$}
$$
$$
\Downarrow
$$
$$
\text{последовательность $B_N$ ограничена}: B_1\le B_2\le B_3\le ... \le B_N
\le ...\le B
$$
$$
\Downarrow
$$
$$
\text{последовательность $A_N$ (монотонна и) ограничена}: A_1\le A_2\le A_3\le
... \le A_N \le ...\le B
$$
$$
\Downarrow
$$
$$
\text{по теореме Вейерштрасса,
существует конечный предел $\lim\limits_{N\to \infty} A_N=A$}
$$
$$
\Downarrow
$$
$$
\text{ряд $\sum\limits_{n=1}^\infty a_n$ сходится}
$$

2. Мы доказали, что если $\sum\limits_{n=1}^\infty b_n$ сходится, то
$\sum\limits_{n=1}^\infty a_n$ тоже сходится.
То есть,

$$
\text{\it невозможна ситуация, когда}\, \sum\limits_{n=1}^\infty a_n
\, \text{\it расходится. а}\, \sum\limits_{n=1}^\infty b_n \,
\text{\it сходится}
$$
Отсюда следует, что если $\sum\limits_{n=1}^\infty a_n$ расходится, то
$\sum\limits_{n=1}^\infty b_n$ тоже должен расходиться.
\end{proof}

\noindent\rule{160mm}{0.1pt}\begin{multicols}{2}

\begin{ex}\label{ex-18.3.8}
$$
\sum_{n=1}^\infty \frac{\cos^2 n}{2^n}\le \kern-50pt
 \underbrace{\sum_{n=1}^\infty \frac{1}{2^n}}_{\smsize
 \begin{matrix}
 \text{сходится,}\\ \text{как сумма геометрической прогрессии} \\
 \text{со знаменателем $|q|<1$}
 \end{matrix}}
$$
Больший ряд сходится, значит и меньший ряд сходится.

Вывод: ряд $\sum\limits_{n=1}^\infty \frac{\cos^2 n}{2^n}$ сходится.
\end{ex}

\begin{ex}\label{ex-18.3.9}
$$
\sum_{n=1}^\infty \frac{2+\cos n}{n}\ge \kern-20pt\underbrace{\sum_{n=1}^\infty
\frac{1}{n}}_{\smsize \begin{matrix}\text{расходится,}\\
\text{в силу примера \ref{ex-18.3.2}}\end{matrix}}
$$
Меньший ряд расходится, значит и больший ряд расходится.

Вывод: ряд $\sum\limits_{n=1}^\infty \frac{2+\cos n}{n}$ расходится.
\end{ex}

\begin{ex}\label{ex-18.3.10}
$$
\sum_{n=1}^\infty \frac{2+\cos n}{n^2}\le \sum_{n=1}^\infty
\frac{3}{n^2}= 3\cdot \kern-20pt\underbrace{\sum_{n=1}^\infty
\frac{1}{n^2}}_{\smsize \begin{matrix} \text{сходится,}\\ \text{в
силу примера \ref{ex-18.3.3}}\end{matrix}}
$$

Вывод: ряд $\sum\limits_{n=1}^\infty \frac{2+\cos n}{n^2}$ сходится.
\end{ex}

\begin{ex}\label{ex-18.3.11}
 \begin{multline*}
\sum_{n=1}^\infty \frac{1}{\sqrt{n}\cdot \arctg (3+\sin n)}\ge
\sum_{n=1}^\infty \frac{1}{\sqrt{n}\cdot \arctg 4}=\\=
\frac{1}{\arctg 4}\cdot \kern-10pt\underbrace{\sum_{n=1}^\infty
\frac{1}{\sqrt{n}}}_{\smsize\begin{matrix} \text{расходится,}\\
\text{в силу примера \ref{ex-18.3.3}}\end{matrix}}
 \end{multline*}

Вывод: ряд $\sum\limits_{n=1}^\infty \frac{1}{\sqrt{n}\cdot \arctg (3+\sin n)}$
расходится.
\end{ex}

\begin{ers} Исследуйте на сходимость ряды:
 \begin{multicols}{2}
1) $\sum\limits_{n=1}^\infty e^{-n^2}$;

2) $\sum\limits_{n=1}^\infty \frac{\ln n}{\sqrt[3] {n^2}}$;

3) $\sum\limits_{n=1}^\infty \frac{2^n}{n\cdot 3^n}$;

4) $\sum\limits_{n=1}^\infty \frac{4+\cos n}{n^\alpha}$;

5) $\sum\limits_{n=1}^\infty (2+\sin n)\cdot n^\alpha$.
\end{multicols}\end{ers}

\end{multicols}\noindent\rule[10pt]{160mm}{0.1pt}

\paragraph{Признак Даламбера.}

\begin{tm}[\bf признак Даламбера]\label{tm-18.3.20}
Пусть $a_n>0$, и
 \beq
D=\lim_{n\to \infty}\frac{a_{n+1}}{a_n}\label{18.3.3}
 \eeq
Тогда
 \bit{
\item[1)] если $D<1$, то ряд $\sum\limits_{n=1}^\infty a_n$ сходится;
\item[2)] если $D>1$, то ряд $\sum\limits_{n=1}^\infty a_n$ расходится.
 }\eit
\end{tm}
\bit{ \item[$\bullet$] Число $D$, определяемое формулой \eqref{18.3.3},
называется {\it числом Даламбера} ряда \break $\sum\limits_{n=1}^\infty a_n$.
 }\eit
\begin{proof}
1. Предположим, что
 \beq
D=\lim_{n\to \infty}\frac{a_{n+1}}{a_n}<1 \label{18.3.4}
 \eeq
Возьмем какое-нибудь число $\varepsilon>0$ так чтобы
$D+\varepsilon<1$. Из \eqref{18.3.4} следует, что
$\varepsilon$-окрестность числа $D$ содержит почти все числа
$\frac{a_{n+1}}{a_n}$, то есть
$$
\exists M \quad \forall n\ge M \quad
D-\varepsilon<\frac{a_{n+1}}{a_n}<D+\varepsilon
$$
Если обозначить $q=D+\varepsilon$, то мы получим что $q<1$ и
$$
\forall n\ge M \quad
\frac{a_{n+1}}{a_n}<q
$$
то есть
$$
\forall n\ge M \quad
a_{n+1}<q\cdot a_n
$$
В частности,
$$
a_{M+1}<q\cdot a_M
$$
$$
a_{M+2}<q\cdot a_{M+1}<q\cdot (q\cdot a_M)=q^2\cdot a_M
$$
$$
a_{M+3}<q\cdot a_{M+2}<q\cdot (q^2\cdot a_M)=q^3\cdot a_M
$$
$$
...
$$
$$
a_{M+k}<q\cdot a_{M+k-1}<...<q^k\cdot a_M
$$
Отсюда
$$
\sum_{n=M}^\infty a_n=
\sum_{k=0}^\infty a_{M+k}<
\sum_{k=0}^\infty \l q^k\cdot a_M \r=
a_M\cdot \sum_{k=0}^\infty q^k
$$
Последний ряд сходится как геометрическая прогрессия со знаменателем
$q<1$. Значит, по признаку сравнения, меньший ряд $\sum_{n=M}^\infty
a_n$ тоже сходится. Отсюда по лемме об остатке \ref{lm-18.3.13}, ряд
$\sum_{n=1}^\infty a_n$ сходится.

2. Предположим, что наоборот,
 \beq
D=\lim_{n\to \infty}\frac{a_{n+1}}{a_n}>1 \label{18.3.5}
 \eeq
Возьмем какое-нибудь число $\varepsilon>0$ так чтобы
$D-\varepsilon>1$. Из \eqref{18.3.5} следует, что
$\varepsilon$-окрестность числа $D$ содержит почти все числа
$\frac{a_{n+1}}{a_n}$, то есть
$$
\exists M \quad \forall n\ge M \quad
D-\varepsilon<\frac{a_{n+1}}{a_n}<D+\varepsilon
$$
Если обозначить $q=D-\varepsilon$, то мы получим что $q>1$ и
$$
\forall n\ge M \quad
q<\frac{a_{n+1}}{a_n}
$$
то есть
$$
\forall n\ge M \quad
a_{n+1}>q\cdot a_n
$$
В частности,
$$
a_{M+1}>q\cdot a_M
$$
$$
a_{M+2}>q\cdot a_{M+1}>q\cdot (q\cdot a_M)=q^2\cdot a_M
$$
$$
a_{M+3}>q\cdot a_{M+2}>q\cdot (q^2\cdot a_M)=q^3\cdot a_M
$$
$$
...
$$
$$
a_{M+k}>q\cdot a_{M+k-1}>...>q^k\cdot a_M
$$
Отсюда
$$
\sum_{n=M}^\infty a_n=
\sum_{k=0}^\infty a_{M+k}>
\sum_{k=0}^\infty \l q^k\cdot a_M \r=
a_M\cdot \sum_{k=0}^\infty q^k
$$
Последний ряд расходится как геометрическая прогрессия со
знаменателем $q>1$. Значит, по признаку сравнения, больший ряд
$\sum_{n=M}^\infty a_n$ тоже расходится. Отсюда по лемме об остатке
\ref{lm-18.3.13}, ряд $\sum_{n=1}^\infty a_n$ расходится.
\end{proof}

\noindent\rule{160mm}{0.1pt}\begin{multicols}{2}

\begin{ex}\label{ex-18.13.21} Чтобы исследовать на сходи\-мость ряд
$$
\sum_{n=1}^\infty \frac{n\cdot 2^n}{3^n}
$$
найдем его число Даламбера:
 \begin{multline*}
D=\lim_{n\to \infty}\frac{a_{n+1}}{a_n}= \lim_{n\to
\infty}\frac{\frac{(n+1)\cdot 2^{n+1}}{3^{n+1}}} {\frac{n\cdot
2^n}{3^n}}=\\= \lim_{n\to \infty}\frac{(n+1)\cdot 2^{n+1}\cdot 3^n}
{n\cdot 2^n\cdot 3^{n+1}}=\\= \lim_{n\to \infty}\frac{(n+1)\cdot
2}{n\cdot 3}=\frac{2}{3}<1
 \end{multline*}
Вывод: ряд $\sum\limits_{n=1}^\infty \frac{n\cdot 2^n}{3^n}$ сходится.
\end{ex}

\begin{ex}\label{ex-18.13.22} Рассмотрим ряд
$$
\sum_{n=1}^\infty \frac{n^n}{n!}
$$
Его число Даламбера:
\begin{multline*}
D=\lim_{n\to \infty}\frac{a_{n+1}}{a_n}= \lim_{n\to
\infty}\frac{\frac{(n+1)^{n+1}}{(n+1)!}}{\frac{n^n}{n!}}=\\=
\lim_{n\to \infty}\frac{(n+1)^{n+1}\cdot n!}{n^n\cdot (n+1)!}=
\lim_{n\to \infty}\frac{(n+1)^{n+1}\cdot n!}{n^n\cdot (n+1)\cdot
n!}=\\= \lim_{n\to \infty}\frac{(n+1)^n}{n^n}= \lim_{n\to \infty}\l
\frac{n+1}{n}\r^n=\\= \lim_{n\to \infty}\l 1+\frac{1}{n}\r^n=e>1
\end{multline*}
Вывод: ряд $\sum\limits_{n=1}^\infty \frac{n^n}{n!}$ расходится.
\end{ex}

\begin{ers} Исследуйте на сходимость ряды:
 \begin{multicols}{2}
1) $\sum_{n=1}^\infty \frac{1}{n!}$;

2) $\sum_{n=1}^\infty \frac{e^n}{n!}$;

3) $\sum_{n=1}^\infty \frac{(n!)^2}{(2n)!}$;

4) $\sum_{n=1}^\infty \frac{2^n\cdot n!}{n^n}$;

5) $\sum_{n=1}^\infty \frac{3^n\cdot n!}{n^n}$;

6) $\sum_{n=1}^\infty \frac{(n!)^2}{2^{n^2}}$;

7) $\sum_{n=1}^\infty \frac{n^5}{2^n+3^n}$.
\end{multicols}\end{ers}

\end{multicols}\noindent\rule[10pt]{160mm}{0.1pt}

\paragraph{Радикальный признак Коши.}

\begin{tm}[\bf радикальный признак Коши]\label{tm-18.3.24}
Пусть $a_n\ge 0$, и
 \beq
C=\lim_{n\to \infty}\sqrt[n]{a_n}\label{18.3.6}
 \eeq
Тогда
 \bit{
\item[1)] если $C<1$, то ряд $\sum\limits_{n=1}^\infty a_n$ сходится;
\item[2)] если $C>1$, то ряд $\sum\limits_{n=1}^\infty a_n$ расходится.
 }\eit
\end{tm}

\bit{ \item[$\bullet$] Число $C$, определяемое формулой \eqref{18.3.6},
называется {\it числом Коши} ряда $\sum\limits_{n=1}^\infty a_n$.
 }\eit

\begin{proof}
1. Предположим, что
 \beq
C=\lim_{n\to \infty}\sqrt[n]{a_n}<1 \label{18.3.7}
 \eeq
Возьмем какое-нибудь число $\varepsilon>0$ так чтобы
$C+\varepsilon<1$. Из \eqref{18.3.7} следует, что
$\varepsilon$-окрестность числа $C$ содержит почти все числа
$\sqrt[n]{a_n}$, то есть
$$
\exists M \quad \forall n\ge M \quad
C-\varepsilon<\sqrt[n]{a_n}<C+\varepsilon
$$
Если обозначить $q=C+\varepsilon$, то мы получим что $q<1$ и
$$
\forall n\ge M \quad
\sqrt[n]{a_n}<q
$$
то есть
$$
\forall n\ge M \quad
a_n<q^n
$$
Отсюда
$$
\sum_{n=M}^\infty a_n<
\sum_{n=M}^\infty q^n
$$
Последний ряд сходится как геометрическая прогрессия со знаменателем
$q<1$. Значит, по признаку сравнения, меньший ряд $\sum_{n=M}^\infty
a_n$ тоже сходится. Отсюда по лемме об остатке \ref{lm-18.3.13}, ряд
$\sum_{n=1}^\infty a_n$ сходится.

2. Предположим, что наоборот
 \beq
C=\lim_{n\to \infty}\sqrt[n]{a_n}>1 \label{18.3.8}
 \eeq
Возьмем какое-нибудь число $\varepsilon>0$ так чтобы
$C-\varepsilon>1$. Из \eqref{18.3.8} следует, что
$\varepsilon$-окрестность числа $C$ содержит почти все числа
$\sqrt[n]{a_n}$, то есть
$$
\exists M \quad \forall n\ge M \quad
C-\varepsilon<\sqrt[n]{a_n}<C+\varepsilon
$$
Если обозначить $q=C-\varepsilon$, то мы получим что $q>1$ и
$$
\forall n\ge M \quad
\sqrt[n]{a_n}>q
$$
то есть
$$
\forall n\ge M \quad
a_n>q^n
$$
Отсюда
$$
\sum_{n=M}^\infty a_n>
\sum_{n=M}^\infty q^n
$$
Последний ряд расходится как геометрическая прогрессия со
знаменателем $q>1$. Значит, по признаку сравнения, меньший ряд
$\sum_{n=M}^\infty a_n$ тоже расходится. Отсюда по лемме об остатке
\ref{lm-18.3.13}, ряд $\sum_{n=1}^\infty a_n$ расходится.
\end{proof}

\noindent\rule{160mm}{0.1pt}\begin{multicols}{2}

\begin{ex}\label{ex-18.13.25} Чтобы исследовать на сходи\-мость ряд
$$
\sum_{n=1}^\infty \l \frac{n-1}{2n+1}\r^n
$$
найдем его число Коши:
$$
C=\lim_{n\to \infty}\sqrt[n]{a_n}= \lim_{n\to
\infty}\frac{n-1}{2n+1}=\frac{1}{2}<1
$$
Вывод: ряд $\sum\limits_{n=1}^\infty \l \frac{n-1}{2n+1}\r^n$ сходится.
\end{ex}

\begin{ex}\label{ex-18.13.26} Рассмотрим ряд
$$
\sum_{n=1}^\infty \l 1+\frac{1}{n}\r^{n^2}
$$
Его число Коши:
$$
C=\lim_{n\to \infty}\sqrt[n]{a_n}= \lim_{n\to \infty}\l
1+\frac{1}{n}\r^n=e>1
$$
Вывод: ряд $\sum\limits_{n=1}^\infty \l 1+\frac{1}{n}\r^{n^2}$ расходится.
\end{ex}

\begin{ers} Исследуйте на сходимость ряды:
 \biter{
\item[1)] $\sum\limits_{n=1}^\infty \frac{1}{2^n}\cdot \l
1+\frac{1}{n}\r^{n^2}$;

\item[2)] $\sum\limits_{n=1}^\infty \frac{1}{3^n}\cdot \l
1+\frac{1}{n}\r^{n^2}$;

\item[3)] $\sum\limits_{n=1}^\infty n\cdot \l \arcsin
\frac{1}{n}\r^n$;

\item[4)] $\sum\limits_{n=1}^\infty
\frac{n^{n+\frac{1}{n}}}{(n+1)^n}$ (здесь признак Коши не дает
результатов);

\item[5)] $\sum\limits_{n=1}^\infty 2^{-n+\sin n}$;

\item[6)] $\sum\limits_{n=1}^\infty \frac{2+(-1)^n}{2^n}$ (этот ряд
исследуется по признаку Коши, но не по признаку Даламбера).
 }\eiter \end{ers}

\end{multicols}\noindent\rule[10pt]{160mm}{0.1pt}

\subsection{Критерий Коши сходимости ряда}

Теперь от знакопостоянных рядов мы возвращаемся к произвольным.

\begin{tm}[\bf критерий Коши сходимости ряда]\label{tm-18.4.1}
Пусть $\{ a_n \}$ -- произвольная числовая последовательность. Тогда следующие
условия эквивалентны:
 \bit{
\item[(i)]
ряд
$$
\sum_{n=1}^\infty a_n
$$
сходится;

\item[(ii)] для любой последовательности $l_i\in \mathbb{N}$, и любой
бесконечно большой последовательности $k_i\in \mathbb{N}$
$$
k_i\underset{i\to \infty}{\longrightarrow}\infty,
$$
сумма $\sum\limits_{n=k_i+1}^{k_i+l_i} a_n$
стремится к нулю при $i\to \infty$:
$$
\sum\limits_{n=k_i+1}^{k_i+l_i} a_n
\underset{i\to \infty}{\longrightarrow} 0
$$
 }\eit
\end{tm}\begin{proof}
Рассмотрим частичные суммы ряда $\sum\limits_{n=1}^\infty a_n$:
$$
S_N=\sum_{n=1}^N a_n
$$
Мы получаем следующую логическую цепочку:
$$
\text{ряд $\sum\limits_{n=1}^\infty a_n$ сходится}
$$
$$
\Updownarrow
$$
$$
\text{последовательность $S_N$ сходится}
$$
$$
\Updownarrow\put(20,0){\smsize \text{$\begin{pmatrix}
\text{вспоминаем критерий Коши сходимости}\\
\text{последовательности -- теорему \ref{Cauchy-crit-seq}}
\end{pmatrix}$}}
$$
$$
\begin{array}{c}\text{для любой последовательности $l_i\in \mathbb{N}$,}\\
\text{и любой бесконечно большой последовательности $k_i\in \mathbb{N}$}\quad
k_i\underset{i\to \infty}{\longrightarrow}\infty,\\
\text{выполняется соотношение:}\quad
S_{k_i+l_i}-S_{k_i}\underset{i\to \infty}{\longrightarrow} 0
\end{array}
$$
$$
\Updownarrow\put(20,0){\smsize \text{$\begin{pmatrix} \text{замечаем,
что}\\ S_{k_i+l_i}-S_{k_i}=\sum\limits_{n=k_i+1}^{k_i+l_i} a_n
\end{pmatrix}$}}
$$
$$
\begin{array}{c}\text{для любой последовательности $l_i\in \mathbb{N}$,}\\
\text{и любой бесконечно большой последовательности $k_i\in \mathbb{N}$}\quad
k_i\underset{i\to \infty}{\longrightarrow}\infty,\\
\text{выполняется соотношение:}\quad \sum\limits_{n=k_i+1}^{k_i+l_i}
a_n \underset{i\to \infty}{\longrightarrow} 0
\end{array}
$$
\end{proof}

\subsection{Общие признаки сходимости рядов}

\paragraph{Необходимое условие сходимости.}

\begin{tm}\label{tm-18.5.1}
Если ряд $\sum\limits_{n=1}^\infty a_n$ сходится, то
$a_n\underset{n\to \infty}{\longrightarrow} 0$.
\end{tm}\begin{tm}[\bf эквивалентная формулировка]
Если $a_n\underset{n\to \infty}{\not\kern-5pt\longrightarrow} 0$, то ряд
$\sum\limits_{n=1}^\infty a_n$ расходится.
\end{tm}\begin{proof}
Если ряд $\sum\limits_{n=1}^\infty a_n$ сходится, то
для последовательностей $k_i=i$ и $l_i=1$,
по критерию Коши (теорема \ref{tm-18.4.1}), мы получим
$$
\sum\limits_{n=k_i+1}^{k_i+1} a_n=a_{k_i+1}=a_i
\underset{i\to \infty}{\longrightarrow} 0
\qquad
$$ \end{proof}

\noindent\rule{160mm}{0.1pt}\begin{multicols}{2}

\begin{ex}\label{ex-18.5.2} Чтобы исследовать на сходи\-мость ряд
$$
\sum\limits_{n=1}^\infty \frac{1}{\sqrt[n]{n}}
$$
достаточно вычислить предел
 \begin{multline*}
\lim_{n\to \infty}\frac{1}{\sqrt[n]{n}}= \lim_{n\to \infty}
n^{-\frac{1}{n}}= \lim_{n\to \infty} e^{\ln \l
n^{-\frac{1}{n}}\r}=\\= \lim_{n\to \infty} e^{-\frac{\ln
n}{n}}=e^0=1\ne 0
 \end{multline*}
Вывод: ряд $\sum\limits_{n=1}^\infty \frac{1}{\sqrt[n]{n}}$ расходится.
\end{ex}

\begin{ex}\label{ex-18.5.3} Чтобы исследовать ряд
$$
\sum\limits_{n=1}^\infty \frac{n^{n+\frac{1}{n}}}{(n+1)^n}
$$
вычисляем предел
 \begin{multline*}
\lim_{n\to \infty}\frac{n^{n+\frac{1}{n}}}{(n+1)^n}= {\smsize\begin{pmatrix}\text{делим числитель}\\
\text{и знаменатель на $n^n$}\end{pmatrix}}=\\= \lim_{n\to
\infty}\frac{n^{\frac{1}{n}}}{\l 1+\frac{1}{n}\r^n}= \frac{1}{e}\ne 0
 \end{multline*}
Вывод: ряд $\sum\limits_{n=1}^\infty \frac{n^{n+\frac{1}{n}}}{(n+1)^n}$
расходится.
\end{ex}

\begin{ex}\label{ex-18.5.4} Чтобы исследовать на сходи\-мость ряд
$$
\sum\limits_{n=1}^\infty (-1)^n\cdot n\cdot \sin \frac{1}{n}
$$
достаточно вычислить предел
$$
\lim_{n\to \infty} |a_n|=
\lim_{n\to \infty} n\cdot \sin \frac{1}{n}=
\lim_{n\to \infty} n\cdot \frac{1}{n}=1\ne 0
$$
То есть  $|a_n|\underset{n\to \infty}{
\longrightarrow\kern-15pt{\Big/}
} 0$,
откуда $a_n\underset{n\to \infty}{
\longrightarrow\kern-15pt{\Big/}
} 0$.

Вывод: ряд $\sum\limits_{n=1}^\infty (-1)^n\cdot n\cdot \sin \frac{1}{n}$
расходится.
\end{ex}

\begin{ers} Исследуйте на сходимость ряды:
 \biter{
\item[1)] $\sum\limits_{n=1}^\infty \frac{(-1)^n}{\sqrt[n]{n}}$;

\item[2)] $\sum\limits_{n=1}^\infty \frac{1}{\sqrt[n]{\ln n}}$;

\item[3)] $\sum\limits_{n=1}^\infty (-1)^n \cdot n^2\cdot \l 1-\cos
\frac{1}{n}\r$.
 }\eiter \end{ers}

\end{multicols}\noindent\rule[10pt]{160mm}{0.1pt}

\paragraph{Признак абсолютной сходимости.}

\begin{tm}\label{tm-18.5.6}
Если
ряд $\sum\limits_{n=1}^\infty |a_n|$ сходится, то
ряд $\sum\limits_{n=1}^\infty a_n$ тоже сходится.
\end{tm}
\begin{proof}
$$
\text{ряд $\sum\limits_{n=1}^\infty |a_n|$ сходится}
$$
$$
\Downarrow\put(20,0){\smsize \text{$\begin{pmatrix} \text{применяем
критерий Коши}\\ \text{сходимости ряда -- теорему
\ref{tm-18.4.1}}\end{pmatrix}$}}
$$
$$
\forall \, l_i, k_i\in \mathbb{N}\quad \l k_i\underset{n\to
\infty}{\longrightarrow}\infty \r \quad
\sum\limits_{n=k_i+1}^{k_i+l_i} |a_n| \underset{i\to
\infty}{\longrightarrow} 0
$$
$$
\Downarrow
$$
$$
\l
0\le \left| \sum\limits_{n=k_i+1}^{k_i+l_i} a_n \right| \le
\sum\limits_{n=k_i+1}^{k_i+l_i} |a_n|
\underset{i\to \infty}{\longrightarrow} 0
\quad
\Rightarrow
\quad
\left| \sum\limits_{n=k_i+1}^{k_i+l_i} a_n \right|
\underset{i\to \infty}{\longrightarrow} 0
\quad
\Rightarrow
\quad
\sum\limits_{n=k_i+1}^{k_i+l_i} a_n
\underset{i\to \infty}{\longrightarrow} 0
\r
$$
$$
\Downarrow
$$
$$
\forall \, l_i, k_i\in \mathbb{N}\quad \l k_i\underset{n\to
\infty}{\longrightarrow}\infty \r \quad
\sum\limits_{n=k_i+1}^{k_i+l_i} a_n \underset{i\to
\infty}{\longrightarrow} 0
$$
$$
\Downarrow\put(20,0){\smsize \text{$\begin{pmatrix} \text{снова
применяем критерий Коши}\\ \text{сходимости ряда}\end{pmatrix}$}}
$$
$$
\text{ряд $\sum\limits_{n=1}^\infty a_n$ сходится}\quad
$$ \end{proof}

\noindent\rule{160mm}{0.1pt}\begin{multicols}{2}

\begin{ex}\label{ex-18.5.7} Чтобы исследовать ряд
$$
\sum\limits_{n=1}^\infty \frac{(-1)^n}{n\sqrt{n}}
$$
рассмотрим соответствующий ряд из модулей:
 $$
\sum\limits_{n=1}^\infty \ml \frac{(-1)^n}{n\sqrt{n}}\mr=
\sum\limits_{n=1}^\infty \frac{1}{n\sqrt{n}}=\kern-35pt
\underbrace{\sum\limits_{n=1}^\infty
\frac{1}{n^\frac{3}{2}}}_{\smsize
\begin{matrix} \text{сходится как ряд Дирихле}\\
\text{с показателем степени $\alpha>1$}\end{matrix}}
 $$
Вывод: ряд $\sum\limits_{n=1}^\infty \frac{(-1)^n}{n\sqrt{n}}$ сходится.
\end{ex}

\begin{ex}\label{ex-18.5.8} Чтобы исследовать ряд
$$
\sum\limits_{n=1}^\infty \frac{\sin n}{2^n}
$$
рассмотрим соответствующий ряд из модулей:
$$
\sum_{n=1}^\infty \ml \frac{\sin n}{2^n}\mr \le\kern-50pt
 \underbrace{\sum_{n=1}^\infty \frac{1}{2^n}}_{\smsize
 \begin{matrix}
 \text{сходится,}\\ \text{как сумма геометрической прогрессии} \\
 \text{со знаменателем $|q|<1$}
 \end{matrix}}
$$
Вывод: ряд $\sum\limits_{n=1}^\infty \frac{\sin n}{2^n}$ сходится.
\end{ex}

\begin{ers} Исследуйте на сходимость ряды:

1) $\sum\limits_{n=1}^\infty \frac{\cos n}{n^2+n+\ln n}$;

2) $\sum\limits_{n=1}^\infty (-1)^\frac{n(n-1)}{2}\frac{n^5}{3^n}$.
 \end{ers}

\end{multicols}\noindent\rule[10pt]{160mm}{0.1pt}

\subsection{Специальные признаки сходимости
рядов}\label{SEC-spetsialnye-priznaki-shodimosti}

\paragraph{Признак Лейбница.}

\begin{tm}\label{tm-18.6.1}
Пусть последовательность $\{b_n\}$ обладает свойствами:
 \bit{
\item[(i)] она неотрицательна и невозрастает,
 \beq\label{18.6.1}
b_0\ge b_1\ge b_2\ge ... \ge 0
 \eeq

\item[(ii)] и стремится к нулю:
 \beq
b_n\underset{n\to \infty}{\longrightarrow} 0
 \eeq
 }\eit
Тогда
 \bit{
\item[(a)] ряд $\sum\limits_{n=0}^\infty (-1)^n b_n$ сходится,

\item[(b)] его остаток оценивается сверху своим первым членом:
 \beq\label{otsenka-ostatka-Leibnitz}
\left|\sum\limits_{n=N+1}^\infty (-1)^n b_n\right|\le b_{N+1}
 \eeq
 }\eit
 \end{tm}

\brem То же справедливо и для ряда, у которого нижний предел суммирования равен
1:
$$
\sum_{n=1}^\infty (-1)^n b_n
$$
\erem

\begin{proof}
Рассмотрим частичные суммы нашего ряда:
$$
  S_N=\sum\limits_{n=1}^N (-1)^n b_n
$$
Мы покажем, что $S_N$ образуют последовательность, которую можно изоборазить
картинкой

%\pucture{0pt}{0pt}{ii-8.pcx}

\vglue120pt \noindent (нечетные суммы образуют возрастающую последовательность,
а четные -- убывающую, причем расстояние между четными и нечетными стремится к
нулю).

Действительно, рассмотрим подпоследовательности нечетных и четных сумм:
$$
A_k=S_{2k-1}, \quad B_k=S_{2k}
$$

1. Покажем сначала, что $A_k$ монотонна и ограничена:
 \beq \label{18.6.2}
A_1\le A_2\le A_3\le ... \le b_0
 \eeq
Действительно, с одной стороны,
 $$
A_{k+1}=S_{2k+1}=\sum_{n=0}^{2k+1} (-1)^n b_n=\underbrace{ \sum_{n=0}^{2k-1}
(-1)^n b_n}_{\tiny\begin{matrix}\|
\\ S_{2k-1}\end{matrix}} +\underbrace{
b_{2k}-b_{2k+1}}_{\tiny\begin{matrix}\VI
\\ 0\end{matrix}} \ge S_{2k-1}=A_k
$$
А, с другой, --
$$
A_{k+1}=S_{2k+1}=b_0 -\underbrace{(b_1-b_2)}_{\tiny\begin{matrix}\VI
\\ 0\end{matrix}}-
\underbrace{(b_3-b_4)}_{\tiny\begin{matrix}\VI
\\ 0\end{matrix}}-...-\underbrace{(b_{2k-1}-b_{2k})}_{\tiny\begin{matrix}\VI
\\ 0\end{matrix}}-\underbrace{(b_{2k+1})}_{\tiny\begin{matrix}\VI
\\ 0\end{matrix}}\le b_0
$$

2. Докажем монотонность и ограниченность последовательности $\{ B_k \}$:
 \beq \label{18.6.2-1}
0\le...\le B_2\le B_1\le B_0
 \eeq
С одной стороны, получаем:
 $$
B_{k+1}=S_{2k+2}=\sum_{n=0}^{2k+2} (-1)^n b_n=\underbrace{ \sum_{n=0}^{2k}
(-1)^n b_n}_{\tiny\begin{matrix}\|
\\ S_{2k}\end{matrix}} -\underbrace{
(b_{2k+1}-b_{2k+2})}_{\tiny\begin{matrix}\VI
\\ 0\end{matrix}} \le S_{2k}=B_k
$$
А, с другой, --
$$
B_k=S_{2k}=\underbrace{(b_0-b_1)}_{\tiny\begin{matrix}\VI
\\ 0\end{matrix}}+
\underbrace{(b_2-b_3)}_{\tiny\begin{matrix}\VI
\\ 0\end{matrix}}+...+\underbrace{(b_{2k-1}-b_{2k})}_{\tiny\begin{matrix}\VI
\\ 0\end{matrix}}\ge 0
$$

3. Из \eqref{18.6.2} и \eqref{18.6.2-1} следует, что последовательности $\{ A_k
\}$ и $\{ B_k \}$ сходятся:
 \beq \label{18.6.3}
A_k\underset{n\to \infty}{\longrightarrow} A=\sup_{k\in\Z_+} A_k,\qquad
B_k\underset{n\to \infty}{\longrightarrow} B=\inf_{k\in\Z_+} B_k
 \eeq
При этом,
 $$
B-A=\lim_{k\to\infty} B_k-\lim_{k\to\infty} A_k=\lim_{k\to\infty}
(B_k-A_k)=\lim_{k\to\infty} (S_{2k}-S_{2k-1})=\lim_{k\to\infty} b_{2k}=0,
 $$
то есть
$$
A=B,
$$
и поэтому соотношения \eqref{18.6.3} удобно переписать так:
 \beq \label{18.6.3-1}
S_{2k-1}\underset{k\to \infty}{\longrightarrow} \sup_{k\in\Z_+}
S_{2k-1}=S=\inf_{k\in\Z_+} S_{2k}\underset{\infty\gets k}{\longleftarrow}
S_{2k}
 \eeq
где $S=A=B$. Отсюда сразу следует, что последовательность $S_N$ сходится к $S$,
$$
S_N \underset{N\to \infty}{\longrightarrow} S,
$$
то есть ряд $\sum\limits_{n=0}^\infty (-1)^n b_n$ сходится. С другой стороны,
из \eqref{18.6.3-1} следует двойное неравенство
 \beq\label{S_(2k-1)-le-S-le-S_(2l)}
\forall k,l\in\N\qquad S_{2k-1}\le S\le S_{2l}
 \eeq
из которого мы получаем, во-первых,
 \beq \label{18.6.3-2}
 |S-S_{2k-1}|=S-S_{2k-1}\le S_{2k}-S_{2k-1}=(-1)^{2k}\cdot b_{2k}=b_{2k}
 \eeq
и, во-вторых,
 \beq \label{18.6.3-3}
|S-S_{2k}|=S_{2k}-S\underset{\scriptsize\begin{matrix}\Uparrow\\ -S\le
-S_{2k+1}\\ \Uparrow\\ S_{2k+1}\le S \\ \Uparrow\\
\eqref{S_(2k-1)-le-S-le-S_(2l)}
\end{matrix}}{\le} S_{2k}-S_{2k+1}=-(-1)^{2k+1}\cdot b_{2k+1}=b_{2k+1}
 \eeq
Цепочки \eqref{18.6.3-2} и \eqref{18.6.3-3} вместе дают оценку
\eqref{otsenka-ostatka-Leibnitz}.
\end{proof}

\noindent\rule{160mm}{0.1pt}\begin{multicols}{2}

\begin{ex}\label{ex-18.6.2} Ряд
$$
\sum\limits_{n=1}^\infty \frac{(-1)^n}{n}
$$
можно представить в виде $\sum\limits_{n=1}^\infty (-1)^n b_n$, где
$$
b_n=\frac{1}{n}\overset{\text{монотонно}}{\underset{n\to
\infty}{\longrightarrow}} 0
$$
поэтому по теореме Лейбница получаем

Вывод: ряд $\sum\limits_{n=1}^\infty \frac{(-1)^n}{n}$ сходится.
\end{ex}

\begin{ex}\label{ex-18.6.3} Ряд
$$
\sum\limits_{n=1}^\infty \frac{(-1)^n \ln n}{n}
$$
можно представить в виде $\sum\limits_{n=1}^\infty (-1)^n b_n$, где
$$
b_n=\frac{\ln n}{n}
$$
Эта последовательность стремится к нулю по теореме о шкале бесконечностей, но
непонятно, будет ли это стремление монотоным. Для того, чтобы это проверить,
можно рассмотреть вспомогательную функцию
$$
  f(x)=\frac{\ln x}{x}
$$
Ее производная будет отрицательной на промежутке $x\ge 3$
$$
f'(x)=\frac{1-\ln x}{x^2}<0, \qquad x\ge 3
$$
Поэтому $f(x)$ монотонно убывает при $x\ge 3$, значит, наша последовательность
$b_n=f(n)$ монотонно убывает начиная с номера $n=3$. Отсюда по теореме Лейбница
получаем, что должен сходиться ряд $\sum\limits_{n=3}^\infty \frac{(-1)^n \ln
n}{n}$. Он является остатком ряда $\sum\limits_{n=1}^\infty \frac{(-1)^n \ln
n}{n}$, значит этот ряд $\sum\limits_{n=1}^\infty \frac{(-1)^n \ln n}{n}$ тоже
сходится (по лемме об остатке \ref{lm-18.3.13}).

Вывод: ряд $\sum\limits_{n=1}^\infty \frac{(-1)^n \ln n}{n}$ сходится.
\end{ex}

\begin{ers} Исследуйте на сходимость ряды:
 \begin{multicols}{2}
1) $\sum\limits_{n=1}^\infty \frac{(-1)^n}{\ln n}$;

2) $\sum\limits_{n=1}^\infty \frac{(-1)^n}{n^\alpha}$;

3) $\sum\limits_{n=1}^\infty \frac{(-1)^n}{n\cdot (\ln n)^\alpha}$.
\end{multicols}\end{ers}

\end{multicols}\noindent\rule[10pt]{160mm}{0.1pt}

\paragraph{Преобразование Абеля и признаки Дирихле и Абеля для рядов.}\label{subsec-preobr-Abelya}

\blm\label{lm-18.6.13} Для любых числовых последовательностей $\{ a_n \}$ и $\{
b_n \}$ справедливо равенство
 \beq\label{18.6.6}
\sum_{n=p}^q a_n\cdot b_n= \sum_{n=p}^{q-1} A_n\cdot (b_n-b_{n+1})+A_q\cdot
b_q,\qquad p\le q
 \eeq
где
$$
A_N=\sum\limits_{n=p}^N a_n
$$
Если вдобавок $\{ b_n \}$ -- монотонная последовательность, то справедливо
неравенство
 \beq \label{18.6.7}
\ml \sum\limits_{n=p}^q a_n \cdot b_n \mr \le 3\cdot\max_{ p\le N\le q}\ml
\sum\limits_{n=p}^N a_n \mr\cdot\max\Big\{|b_p|,|b_q|\Big\}
 \eeq
\elm
 \bit{
\item[$\bullet$] Равенство \eqref{18.6.6} называется {\it преобразованием
Абеля}\footnote{В \eqref{18.6.6} можно ввести третий параметр, число $M<p$, и
тогда в обозначении $A_N=\sum_{n=M}^N a_n$ формула примет более привычный вид:
$$
\sum_{n=p}^q a_n\cdot b_n= \sum_{n=p}^{q-1} A_n\cdot (b_n-b_{n+1})+A_q\cdot
b_q-A_{p-1}\cdot  b_p,
$$
}, а неравенство \eqref{18.6.7} -- {\it неравенством Абеля}.
 }\eit
 \begin{proof} 1. Сначала докажем формулу \eqref{18.6.6}.
 \begin{multline*}
\sum_{n=p}^q a_n\cdot b_n= \overbrace{a_p}^{\scriptsize\begin{matrix}A_p\\
\text{\rotatebox{90}{$=$}}\end{matrix}}\cdot b_p+\kern-10pt\overbrace{a_{p+1}}^{\scriptsize\begin{matrix}-A_p+A_{p+1}\\
\text{\rotatebox{90}{$=$}}\end{matrix}}\kern-10pt\cdot b_{p+1}+...+\kern-10pt\overbrace{a_{q-1}}^{\scriptsize
\begin{matrix}-A_{q-2}+A_{q-1}\\
\text{\rotatebox{90}{$=$}}\end{matrix}}\kern-10pt\cdot
b_{q-1}+\kern-10pt\overbrace{a_q}^{\scriptsize
\begin{matrix}-A_{q-1}+A_q\\
\text{\rotatebox{90}{$=$}}\end{matrix}}\kern-10pt\cdot b_q=\\= A_p
b_p+(-A_p+A_{p+1}) b_{p+1}+...+(-A_{q-2}+A_{q-1}) b_{q-1} +(-A_{q-1}+A_q)
b_q=\\= \underbrace{A_p b_p-A_p
b_{p+1}}_{\scriptsize\begin{matrix}\text{\rotatebox{90}{$=$}}
\\ A_p
(b_p-b_{p+1})\end{matrix}}+\underbrace{A_{p+1} b_{p+1}+}...\underbrace{
-A_{q-2}b_{q-1}}+\underbrace{A_{q-1} b_{q-1}-A_{q-1}
b_q}_{\scriptsize\begin{matrix}\text{\rotatebox{90}{$=$}}
\\ A_{q-1}
(b_{q-1}-b_q)\end{matrix}}+A_q b_q=\\= \sum_{n=p}^{q-1} A_n\cdot
(b_n-b_{n+1})+A_q\cdot b_q
\end{multline*}

2. Теперь предположим, что последовательность $\{ b_n \}$ монотонна и докажем
неравенство \eqref{18.6.7}. Здесь можно считать, что $\{ b_n \}$ невозрастает,
потому что случай, когда она неубывает, сводится к этому умножением $\{ b_n \}$
на $-1$:
 \beq\label{18.6.9}
b_1\ge b_2\ge b_3\ge ...
 \eeq
Обозначим
 \beq
A=\max_{ p\le N\le q}|A_N|=\max_{ p\le N\le q}\ml \sum\limits_{n=p}^N a_n \mr
\label{18.6.8}
 \eeq
Тогда:
 \begin{multline*}
\ml \sum\limits_{n=p}^q a_n \cdot b_n \mr =\eqref{18.6.6}= \ml
\sum\limits_{n=p}^{q-1} A_n\cdot (b_n-b_{n+1})+A_q\cdot b_q \mr \le
\sum_{n=p}^{q-1} \kern-17pt\underbrace{|A_n|}_{\scriptsize\begin{matrix}
\phantom{\tiny\eqref{18.6.8}} \ \text{\rotatebox{90}{$\ge$}}\
{\tiny\eqref{18.6.8}} \\ A \end{matrix}}\kern-17pt\cdot |b_n-b_{n+1}|+
\underbrace{|A_q|}_{\scriptsize\begin{matrix} \text{\rotatebox{90}{$\ge$}} \\ A
\end{matrix}}\cdot |b_q|
 \le \\ \le
A\cdot \bigg(\sum\limits_{n=p}^{q-1}
\underbrace{|\kern-7pt\overbrace{b_n-b_{n+1}}^{\scriptsize\begin{matrix}0\\
\phantom {\tiny\eqref{18.6.9}} \ \text{\rotatebox{90}{$\ge$}} \
{\tiny\eqref{18.6.9}}\end{matrix}}\kern-7pt|}_{\scriptsize\begin{matrix}
 \text{\rotatebox{90}{$=$}} \\ b_n-b_{n+1} \end{matrix}}+ |b_q|
 \bigg)=
 A\cdot \bigg(\sum_{n=p}^{q-1} (b_n-b_{n+1})+ |b_q|\bigg)=\\= A\cdot \bigg(
\underbrace{b_p-b_{p+1}\put(20,22){\put(-30.5,-3.5){\line(1,0){30}}\put(-33,-10){$\downarrow\kern25pt\downarrow$}}
}_{n=p}+\underbrace{b_{p+1}-b_{p+2}\put(20,22){\put(-30.5,-3.5){\line(1,0){20}}\put(-33,-10){$\downarrow\kern15pt\downarrow$}}
}_{n=p+1}+...\put(30,22){\put(-30.5,-3.5){\line(1,0){20}}\put(-33,-10){$\downarrow\kern15pt\downarrow$}}
+\underbrace{b_{q-1}-b_q}_{n=q-1}+ |b_q| \bigg)=\\=A\cdot(b_p-b_q+|b_q|)\le
A\cdot(|b_p|+|b_q|+|b_q|)\le 3\cdot A\cdot\max\{|b_p|,|b_q|\}
 \end{multline*}
 \end{proof}

\begin{tm}[\bf признак Дирихле]\label{tm-18.6.5}
Пусть последовательности $\{ a_n \}$ и $\{ b_n \}$ обладают следующими
свойствами:
 \bit{
\item[$(i)$] частичные суммы ряда $\sum\limits_{n=1}^\infty a_n$
ограничены:
$$
\sup_{N\in \mathbb{N}}\ml \sum\limits_{n=1}^N a_n \mr <\infty
$$
\item[$(ii)$] последовательность $\{ b_n \}$ монотонно стремится к нулю:
$$
b_n\overset{\text{монотонно}}{\underset{n\to \infty}{\longrightarrow}} 0
$$
 }\eit
Тогда ряд $\sum\limits_{n=1}^\infty a_n\cdot b_n$ сходится.
\end{tm}

\begin{proof}
Обозначим
$$
\sup_{N\in \mathbb{N}}\ml \sum_{n=1}^N a_n \mr=C <\infty
$$
и заметим, что
 \beq\label{sup_p,q|sum_pq|<infty}
\forall p,q\in \mathbb{N}\qquad \ml \sum_{n=p}^q a_n \mr\le 2C
 \eeq
Действительно,
$$
\bigg|\sum_{n=p}^q a_n\bigg|=\bigg|\sum_{n=1}^q a_n-\sum_{n=1}^{p-1}
a_n\bigg|\le\underbrace{\bigg|\sum_{n=1}^q
a_n\bigg|}_{\scriptsize\begin{matrix}\text{\rotatebox{90}{$\ge$}}\\
C\end{matrix}}+\underbrace{\bigg|\sum_{n=1}^{p-1} a_n\bigg|}_{\scriptsize\begin{matrix}\text{\rotatebox{90}{$\ge$}}\\
C\end{matrix}}\le 2C
$$
Чтобы убедиться, что ряд $\sum_{n=1}^\infty a_n\cdot b_n$ сходится,
воспользуемся критерием Коши (теорема \ref{tm-18.4.1}): выберем произвольные
последовательности  $l_i\in \mathbb{N}$ и $k_i\in \mathbb{N}$
($k_i\underset{n\to \infty}{\longrightarrow}\infty$). По неравенству Абеля
\eqref{18.6.7} получаем:
$$
\ml \sum\limits_{n=k_i+1}^{k_i+l_i} a_n \cdot b_n \mr \le 3\cdot
\underbrace{\max_{ k_i\le N\le k_i+l_i}\ml \sum\limits_{n=k_i+1}^N a_n
\mr}_{\scriptsize\begin{matrix}\phantom{\tiny\eqref{sup_p,q|sum_pq|<infty}}\
\text{\rotatebox{90}{$\ge$}}\ {\tiny\eqref{sup_p,q|sum_pq|<infty}}\\
2C\end{matrix}}\cdot\max\Big\{\underbrace{|b_{k_i+1}|}_{\scriptsize
\begin{matrix}\downarrow \\ 0\end{matrix}},\underbrace{|b_{k_i+l_i}|}_{\scriptsize
\begin{matrix}\downarrow \\ 0\end{matrix}}\Big\}\underset{i\to\infty}{\longrightarrow}0
$$
Это верно для любых последовательностей $l_i\in \mathbb{N}$ и $k_i\in
\mathbb{N}$ ($k_i\underset{n\to \infty}{\longrightarrow}\infty$), значит ряд
$\sum\limits_{n=1}^\infty a_n\cdot b_n$ сходится.
\end{proof}

Признак Дирихле часто используется для исследования рядов вида
$$
\sum_{n=1}^\infty b_n\cdot \sin n x, \quad
\sum_{n=1}^\infty b_n\cdot \cos n x
$$
В этих случаях бывают полезны следующие тригонометрические формулы:
 \beq\label{18.6.5}
\sum_{n=1}^N \sin n x= \frac{\sin \frac{N}{2} x \cdot \sin
\frac{N+1}{2} x}{\sin \frac{x}{2}}\qquad \sum_{n=1}^N \cos n x =
\frac{\sin \frac{N}{2} x \cdot \cos \frac{N+1}{2} x}{\sin
\frac{x}{2}}\qquad (x\ne 2\pi k)
 \eeq
Покажем как они применяются.

\noindent\rule{160mm}{0.1pt}\begin{multicols}{2}

\begin{ex}\label{ex-18.6.6} Исследуем на сходимость ряд
$$
\sum_{n=1}^\infty \frac{\sin n}{n}
$$
Если положить
$$
a_n=\sin n, \quad b_n=\frac{1}{n}
$$
то мы получим $b_n\overset{\text{монотонно}}{\underset{n\to
\infty}{\longrightarrow}} 0$, и
 \begin{multline*}
\sup_{N\in \mathbb{N}}\ml \sum\limits_{n=1}^N a_n \mr = \sup_{N\in
\mathbb{N}}\ml \sum\limits_{n=1}^N \sin n \mr =\\= {\smsize \eqref{18.6.5}}=\\=
\sup_{N\in \mathbb{N}}\ml \frac{\sin \frac{N}{2}\cdot \sin \frac{N+1}{2}}{\sin
\frac{1}{2}}\mr \le \frac{1}{\sin \frac{1}{2}}<+\infty
 \end{multline*}
По признаку Дирихле, получаем

Вывод: ряд $\sum\limits_{n=1}^\infty \frac{\sin n}{n}$ сходится.
\end{ex}

\begin{ex}\label{ex-18.6.7} Исследуем на сходимость ряд
$$
\sum_{n=1}^\infty \frac{\sin n x}{n}
$$
Если положить
$$
a_n=\sin n x, \quad b_n=\frac{1}{n}
$$
то мы получим $b_n\overset{\text{монотонно}}{\underset{n\to
\infty}{\longrightarrow}} 0$, и
 \begin{multline*}
\sup_{N\in \mathbb{N}}\ml \sum\limits_{n=1}^N a_n \mr = \sup_{N\in
\mathbb{N}}\ml \sum\limits_{n=1}^N \sin n x \mr =\\= {\smsize
\eqref{18.6.5}}=\\= \sup_{N\in \mathbb{N}}\ml \frac{\sin \frac{N}{2} x \cdot
\sin \frac{N+1}{2} x}{\sin \frac{x}{2}}\mr \le \frac{1}{\ml \sin
\frac{x}{2}\mr}<+\infty
 \end{multline*}
По признаку Дирихле, получаем, что наш ряд сходится при $x\ne 2\pi k$.

Остается проверить, будет ли он сходиться при $x=2\pi k$.
Подставим это значение в наш ряд:
$$
\sum_{n=1}^\infty \frac{\sin 2\pi kn}{n}=
\sum_{n=1}^\infty 0
$$
Ясно, что этот ряд сходится.

Вывод: ряд $\sum\limits_{n=1}^\infty \frac{\sin n x}{n}$ сходится при любом
$x\in \R$.
\end{ex}

\begin{ers} Исследуйте на сходимость ряды:
 \biter{
\item[1)] $\sum\limits_{n=1}^\infty \frac{\cos n}{n}$;

\item[2)] $\sum\limits_{n=1}^\infty \frac{\cos n x}{n}$;

\item[3)] $\sum\limits_{n=1}^\infty \frac{\sin \frac{\pi n}{12}}{\ln
n}$;

\item[4)] $\sum\limits_{n=1}^\infty \frac{\ln n}{n}\cdot \cos
\frac{\pi n}{12}$;

\item[5)] $\sum\limits_{n=1}^\infty (-1)^\frac{n(n+1)}{2}\cdot
\frac{1}{\sqrt{n}}$.
 }\eiter\end{ers}

\end{multicols}\noindent\rule[10pt]{160mm}{0.1pt}

\begin{tm}[\bf признак Абеля]\label{tm-18.6.9}
Пусть последовательности  $\{ a_n \}$ и $\{ b_n \}$ обладают следующими
свойствами:
 \bit{
\item[$(i)$] ряд $\sum\limits_{n=1}^\infty a_n$ сходится;

\item[$(ii)$] последовательность $\{ b_n \}$ монотонна и ограничена.
 }\eit
Тогда ряд $\sum\limits_{n=1}^\infty a_n\cdot b_n$ сходится.
\end{tm}

\begin{proof}
Заметим, что из сходимости ряда $\sum\limits_{n=1}^\infty a_n$ следует, что для
произвольных последовательностей $l_i\in \mathbb{N}$ и $k_i\in \mathbb{N}$
($k_i\underset{n\to \infty}{\longrightarrow}\infty$) выполняется соотношение
 \beq\label{PROOF:Abel-sum-1}
\max_{ k_i\le N\le k_i+l_i}\ml \sum\limits_{n=k_i+1}^N a_n
\mr\underset{i\to\infty}{\longrightarrow}0
 \eeq
Действительно, если бы это было не так,
$$
\max_{ k_i\le N\le k_i+l_i}\ml \sum\limits_{n=k_i+1}^N a_n
\mr\underset{i\to\infty}{\not\kern-5pt\longrightarrow}0
$$
то, переходя к подпоследовательности, мы получили бы, что для некоторого $\e>0$
выполняется соотношение
$$
\max_{ k_i\le N\le k_i+l_i}\ml \sum\limits_{n=k_i+1}^N a_n \mr>\e
$$
$$
\Downarrow
$$
$$
\exists N_i\in[k_i,k_i+l_i]\qquad \underbrace{\ml \sum_{n=k_i+1}^{N_i}
a_n \mr}_{\scriptsize\begin{matrix}\text{здесь $N_i>k_i$,}\\
\text{потому что}\\ \text{иначе было бы}\\ \sum\limits_{n=k_i+1}^{N_i} a_n=0
\end{matrix}}>\e
$$
$$
\Downarrow
$$
$$
\exists m_i=N_i-k_i\qquad \underbrace{\ml \sum\limits_{n=k_i+1}^{k_i+m_i} a_n
\mr>\e}_{\scriptsize\begin{matrix}\uparrow\\ \text{это невозможно,}\\
\text{потому что}\\ \text{по критерию Коши}\\
\text{(теорема \ref{tm-18.4.1})}\\
\ml \sum\limits_{n=k_i+1}^{k_i+m_i} a_n
\mr\underset{i\to\infty}{\longrightarrow}0
\end{matrix}}
$$

Теперь обозначим
 $$
B=\sup_{n\in\N}|b_n|<\infty
 $$
Чтобы доказать сходимость ряда $\sum\limits_{n=1}^\infty a_n\cdot b_n$, снова
воспользуемся критерием Коши (теоремой \ref{tm-18.4.1}): для произвольных
последовательностей  $l_i\in \mathbb{N}$ и $k_i\in \mathbb{N}$
($k_i\underset{n\to \infty}{\longrightarrow}\infty$) по неравенству Абеля
\eqref{18.6.7} мы получим
$$
\ml \sum\limits_{n=k_i+1}^{k_i+l_i} a_n \cdot b_n \mr \le 3\cdot
\underbrace{\max_{ k_i\le N\le k_i+l_i}\ml \sum\limits_{n=k_i+1}^N a_n
\mr}_{\scriptsize\begin{matrix}
\phantom{\tiny\eqref{PROOF:Abel-sum-1}}\ \downarrow\ {\tiny\eqref{PROOF:Abel-sum-1}} \\
0\end{matrix}}\cdot\max\Big\{\underbrace{|b_{k_i+1}|}_{\scriptsize
\begin{matrix}\text{\rotatebox{90}{$\ge$}} \\ B\end{matrix}},\underbrace{|b_{k_i+l_i}|}_{\scriptsize
\begin{matrix}\text{\rotatebox{90}{$\ge$}} \\ B\end{matrix}}\Big\}\underset{i\to\infty}{\longrightarrow}0
$$
Это верно для любых последовательностей $l_i\in \mathbb{N}$ и $k_i\in
\mathbb{N}$ ($k_i\underset{n\to \infty}{\longrightarrow}\infty$), значит  ряд
$\sum\limits_{n=1}^\infty a_n\cdot b_n$ сходится.
\end{proof}

\noindent\rule{160mm}{0.1pt}\begin{multicols}{2}

\begin{ex}\label{ex-18.6.10} Рассмотрим ряд
$$
\sum_{n=1}^\infty \frac{\sin n}{n}\cos \frac{\pi}{n}
$$
Если положить
$$
a_n=\frac{\sin n}{n}, \quad b_n=\cos \frac{\pi}{n}
$$
то мы получим, что $\sum\limits_{n=1}^\infty a_n$ -- сходящийся ряд (мы уже
доказали это в примере \ref{ex-18.6.6}), а $b_n$ -- монотонная ограниченная
последовательность. Значит, по признаку Абеля, ряд $\sum\limits_{n=1}^\infty \frac{\sin n}{n}\cos \frac{\pi}{n}$
сходится.
\end{ex}

\begin{ex}\label{ex-18.6.11} Следующий ряд
$$
\sum_{n=2}^\infty \frac{\sin n x\cdot \cos n x}{\ln n}\arctg (nx)
$$
исследуется таким же образом. Если положить
$$
a_n=\frac{\sin n x\cdot \cos n x}{\ln n}, \quad b_n= \arctg (nx)
$$
то мы получим, что
$$
\sum\limits_{n=2}^\infty a_n=
\sum\limits_{n=2}^\infty \frac{\sin n x\cdot \cos n x}{\ln n}=
\sum\limits_{n=2}^\infty \frac{\sin 2n x}{2\ln n}
$$
-- сходящийся ряд при любом $x\in \R$ (это доказывается так же, как в примере
\ref{ex-18.6.7}), а $b_n$ -- монотонная ограниченная последовательность.
Значит, по признаку Абеля, ряд $\sum\limits_{n=2}^\infty \frac{\sin n x\cdot \cos n x}{\ln n}\arctg
(nx)$ сходится при любом $x\in \R$.
\end{ex}

\begin{ers} Исследуйте на сходимость ряды:
 \biter{
\item[1)] $\sum\limits_{n=1}^\infty \frac{\cos n}{n}\cdot \arcsin
\frac{n}{n+1}$;

\item[2)] $\sum\limits_{n=1}^\infty \frac{\cos n x}{n}\cdot \arctg
(nx)$;

\item[3)] $\sum\limits_{n=1}^\infty \frac{\sin \frac{\pi n}{12}}{\ln
n}\cdot \arcctg (-n)$;

\item[4)] $\sum\limits_{n=1}^\infty \frac{\ln n}{n}\cdot \cos
\frac{\pi n}{12}\cdot \arccos \l -\frac{n}{n+1}\r$;

\item[5)] $\sum\limits_{n=1}^\infty (-1)^\frac{n(n+1)}{2}\cdot
\frac{\cos \frac{1}{n}}{\sqrt{n}}$.
 }\eiter\end{ers}

\end{multicols}\noindent\rule[10pt]{160mm}{0.1pt}

\chapter{ФУНКЦИОНАЛЬНЫЕ ПОСЛЕДОВАТЕЛЬНОСТИ И РЯДЫ}\label{CH-functional-sequen}

Объяснить, зачем нужны функциональные последовательности и функциональные ряды,
о которых пойдет речь в этой главе, можно, начав, например, со следующего
вопроса: {\it может ли непрерывная функция на отрезке быть нигде не
дифференцируемой?}

Это нужно вот для чего. Начиная с главы \ref{ch-f'(x)}, мы доказывали
утверждения о дифференцируемых функциях, постоянно оговаривая в формулировках
условия дифференцирумости или гладкости, но может быть природа устроена проще,
и эти условия в подходящей интерпретации выполняются автоматически? Вот,
например, в теореме \ref{TH:integrir-po-chastyam-dlya-nepr-kus-glad-funk} об
интегрировании по частям речь идет о непрерывных кусочно гладких функциях, но,
может быть условие кусочной гладкости автоматически следует из непрерывности, и
поэтому отдельно его требовать не нужно? Мы ведь до сих пор не задумывались,
может ли непрерывная функция быть существенно негладкой (точнее, не быть
кусочно гладкой). Вдруг такое невозможно? Тогда можно было бы упростить теорему
\ref{TH:integrir-po-chastyam-dlya-nepr-kus-glad-funk} и некоторые подобные
утверждения, выбросив из формулировок упоминание о кусочной гладкости функций,
и это легче было бы запомнить.

Точно так же, если бы обнаружилось, например, что любая непрерывная функция на
отрезке дифференцируема во всех точках, кроме, может быть, конечного набора (а
других непрерывных функций мы пока не встречали), то это означало бы, что
вообще половину утверждений, представленных нами здесь начиная с главы
\ref{ch-f'(x)}, необходимо срочно переписывать, потому что все оказывается
намного проще.

Вопрос, {\it не может ли это быть проще}, является по-видимому, главным в
науке, поскольку сам научный метод, если пытаться объяснить его суть просто,
представляет собой не что иное, как способ упорядочивать поступающие знания о
мире таким образом, чтобы, не теряя практической ценности, они выстраивались в
структуру, достаточно простую, чтобы ее можно было легко запомнить.
Естественно, в математике этот вопрос также важен\footnote{Обычно по этому
поводу цитируется упоминавшийся нами в главе \ref{ch-R&N} Давид Гильберт.
Констанс Рид в своей знаменитой биографии Гильберта приводит такие его слова из
разговора с Гаральдом Бором: «То, что мне удалось что-то сделать в математике,
объясняется, на самом деле, тем, что я всегда находил всё очень сложным. Когда
я читаю или когда мне что-то рассказывают, мне почти всегда это кажется очень
трудным и практически невозможным понять. Тогда я не могу не задать себе
вопрос, а не может ли это быть проще. И в некоторых случаях оказывалось, что
это действительно намного проще».}, причем применительно к математическим
утверждениям он имеет две важные частные формулировки:
 \bit{\it
\item[1)] Насколько существенны условия на объекты, оговариваемые в теоремах?

\item[2)] Не будут ли какое-то аксиомы, в частности, аксиомы, постулирующие
существование объектов с заданными свойствами, данной теории ее теоремами?
 }\eit
Здесь в обоих случаях речь идет об объектах с нужными свойствами, и понятно,
что для того, чтобы судить о существовании таких объектов, желательно иметь
операции, позволяющие строить новые объекты из уже имеющихся. В анализе такой
операцией является предельный переход, причем если под искомым объектом
понимается число, то таким предельным переходом будет предел числовой
последовательности, о котором речь шла в главе \ref{ch-x_n}, а если под искомым
объектом понимается функция, то операцией, позволяющей ее построить, будет
предел функциональной последовательности, и об этой конструкции мы как раз
поговорим в настоящей главе. В качестве приложения мы, в частности, дадим ответ
на сформулированный здесь в самом начале вопрос о существовании нигде не
гладкой непрерывной функции (см.
\ref{SEC:ravnom-shodimost}\ref{nepr-no-nedifferent-func}). Одновременно
описываемое в этой главе понятие функционального ряда позволит затем в
\ref{SEC:prilozh-step-ryadov} главы \ref{CH-step-ryady} доказать избыточность
аксиомы степеней и аксиомы тригонометрии, о которых речь шла в главе
\ref{ch-ELEM-FUNCTIONS}.

\section{Поточечная сходимость}

\subsection{Функциональная последовательность и ее область сходимости}

 \bit{
\item[$\bullet$] {\it Функциональной последовательностью}
\index{последовательность!функциональная} называется последовательность,
элементами которой являются функции:
$$
  f_1, f_2, f_3, ...
$$
{\it Пределом функциональной последовательности}\index{предел!функциональной
последовательности} считается функция,
$$
  f(x)=\lim_{n\to \infty} f_n(x)
$$
определенная для тех значений переменной $x$, для которых предел
$\lim\limits_{n\to \infty} f_n(x)$ существует и конечен. Это множество
$$
D=\{ x\in \R:\quad \text{предел $\lim\limits_{n\to \infty} f_n(x)$ существует и
конечен}\}
$$
называется {\it областью сходимости}\index{область сходимости!функциональной
последовательности} функциональной последовательности $\{ f_n \}$.

\item[$\bullet$] Говорят, что функциональная последовательность $\{ f_n \}$
стремится к функции $f$ {\it поточечно}\index{сходимость!функциональной
последовательности!поточечная} на множестве $E$, если для всякого $x\in E$
числовая последовательность $\{ f_n(x) \}$ стремится к числу $f(x)$:
 \beq\label{19.3.1}
\forall x\in E \quad f_n(x)\underset{n\to \infty}{\longrightarrow} f(x)
 \eeq
коротко это записывается так:
$$
f_n(x)\overset{x\in E}{\underset{n\to \infty}{\longrightarrow}} f(x)
$$

 }\eit

\noindent\rule{160mm}{0.1pt}\begin{multicols}{2}

\begin{ex}\label{er-3.4.1}
Рассмотрим функцию, заданную формулой
$$
f(x)=\lim\limits_{n\to \infty}\frac{1}{1+x^n}
$$
Покажем, что, несмотря на такой необычный способ задания функции, можно
построить ее график (и заодно найти точки разрыва). Здесь важно заметить, что
$$
x^n\underset{n\to \infty}{\longrightarrow}\begin{cases}\infty, & \text{если}\;
|x|>1\\0, & \text{если}\; |x|<1\end{cases}
$$
После этого надо рассмотреть несколько случаев:

1) если $|x|>1$, то $f(x)=\lim\limits_{n\to \infty}\frac{1}{1+x^n}=
\left(\frac{1}{1+\infty}\right)=0$;

2) если $|x|<1$, то $f(x)=\lim\limits_{n\to \infty}\frac{1}{1+x^n}=
\frac{1}{1+0}=0$;

3) если $x=1$, то получаем $f(1)=\lim\limits_{n\to \infty}\frac{1}{1+x^n}=
\lim\limits_{n\to \infty}\frac{1}{1+1^n}=\frac{1}{1+1}=\frac{1}{2}$;

4) если $x=-1$, то $\lim\limits_{n\to \infty}\frac{1}{1+x^n}$ не существует, и
поэтому $f(-1)$ не определено.

В результате получаем
$$
f(x)=
\begin{cases}0, & \text{если}\; |x|>1
\\ 1, & \text{если}\; |x|<1
\\ \frac{1}{2}, & \text{если}\; x=1
\end{cases}
$$
и график этой функции выглядит следующим образом:

%\picture{0pt}{0pt}{p7.pcx}

\vglue100pt \noindent Видно, что $x=1$ здесь является точкой разрыва, а все
остальные точки, в которых функция определена -- то есть точки $x\in
(-\infty;-1)\cup (-1;1)\cup (1;+\infty)$ -- являются точками непрерывности
$f(x)$. (При этом, в точке $x=-1$ функция $f$ не определена, и поэтому не имеет
смысла говорить о непрерывности $f$ в этой точке.)
\end{ex}

\begin{ex}\label{er-3.4.2}
Решим ту же задачу для функции
$$
f(x)=\lim\limits_{n\to \infty}\frac{n^x-n^{-x}}{n^x+n^{-x}}
$$
Заметим, что
$$
n^\alpha\underset{n\to \infty}{\longrightarrow}\begin{cases}\infty, &
\text{если}\; \alpha>0\\0, & \text{если}\; \alpha<0\end{cases}
$$
и рассмотрим несколько случаев:

1) если $x>0$, то $f(x)=\lim\limits_{n\to \infty}\frac{n^x-n^{-x}}{n^x+n^{-x}}=
\lim\limits_{n\to \infty}\frac{1-n^{-2x}}{1+n^{-2x}}= \frac{1-0}{1+0}=1$;

2) если $x<0$, то $f(x)=\lim\limits_{n\to \infty}\frac{n^x-n^{-x}}{n^x+n^{-x}}=
\lim\limits_{n\to \infty}\frac{n^{2x}-1}{n^{2x}+1}= \frac{0-1}{0+1}=-1$;

2) если $x=0$, то $f(x)=\lim\limits_{n\to \infty}\frac{n^0-n^0}{n^0+n^0}=
\frac{1-1}{1+1}=0$.

В результате получаем
$$
f(x)=
\begin{cases}\; 1, & \text{если}\; x>0
\\ \; 0, & \text{если}\; x=0
\\ -1, & \text{если}\; x<0
\end{cases}
$$
то есть $f(x)$ совпадает с функцией сигнум:

%\picture{0pt}{0pt}{96.pcx}

\vglue100pt \noindent Видно, что $x=0$ здесь является точкой разрыва, а все
остальные точки, -- то есть точки $x\in (-\infty;0)\cup (0;+\infty)$ --
являются точками непрерывности $f(x)$.
\end{ex}

\begin{ex}\label{ex-19.2.3}
Рассмотрим функциональную последовательность
$$
f_n(x)=\sqrt[n]{1+x^{2n}}
$$
Ее областью сходимости будет вся числовая прямая
$$
D=\R
$$
а пределом будет функция
$$
f(x)=
\begin{cases}
 1, & \text{если}\, |x|<1 \\
 1, & \text{если}\, |x|=1 \\
 x^2, & \text{если}\, |x|>1
\end{cases}
$$
\end{ex}

\begin{ex}\label{ex-19.2.4}
Функциональная последовательность
$$
f_n(x)=x^n
$$
имеет область сходимости
$$
D=(-1;1]
$$
и ее пределом будет функция
$$
f(x)=
\begin{cases}
 0, & \text{если}\,\, x\in (-1;1)
 \\
 1, & \text{если}\,\, x=1
 \end{cases}
$$
\end{ex}

\begin{ex}\label{ex-19.2.5}
Функциональная последовательность
$$
f_n(x)=n^x
$$
имеет область сходимости
$$
D=(-\infty;0]
$$
и ее пределом будет функция
$$
f(x)=
\begin{cases}
 0, & \text{если}\,\, x<0
 \\
 1, & \text{если}\,\, x=0
 \end{cases}
$$
\end{ex}

\begin{ers}
Постройте графики и найдите точки разрыва следующих функций:
 \biter{
\item[1)] $f(x)=\lim\limits_{n\to \infty}\cos^{2n} x$;

\item[2)] $f(x)=\lim\limits_{n\to \infty}\frac{x}{1+ (2\sin x)^{2n}}$;

\item[3)] $f(x)=\lim\limits_{n\to \infty} x\cdot \arctg(n\cdot \ctg x)$;

\item[4)] $f(x)=\lim\limits_{n\to \infty}\frac{x+x^2 2^{n x}}{1+2^{n x}}$;

\item[5)] $f(x)=\lim\limits_{n\to \infty}\frac{\log_2(1+2^{n
x})}{\log_2(1+2^{n})}$;

\item[6)] $f(x)=\lim\limits_{n\to \infty}\frac{x^{-n}-1}{x^{-n}+1}$;

\item[7)] $f(x)=\lim\limits_{n\to \infty} (\sin^{2n} x- \cos^{2n} x)$;

\item[8)] $f(x)=\lim\limits_{n\to \infty}\frac{2^{-n x}-2^{n x}}{2^{-n x}+2^{n
x}}$;

\item[9)] $f(x)=\lim\limits_{n\to \infty} x^2\cdot \arctg(x^{2n})$;

\item[10)] $f(x)=\lim\limits_{n\to \infty}\arcsin \frac{n^x-1}{n^x+1}$.

\item[11)] $f_n(x)=\arctg (nx)$;

\item[12)] $f_n(x)=\arctg \frac{x}{n}$;

\item[13)] $f_n(x)=x^n-x^{n+1}$;

\item[14)] $f_n(x)=\sqrt[n]{x^2+\frac{1}{n^2}}$.
 }\eiter\end{ers}

\bex Докажем следующую формулу, выражающую функцию Дирихле из
\eqref{func-Dirichle} как двойной поточечный предел стандартных функций:
 \beq\label{vyrazhenie-Dirichlet-cherez-cos}
D(x)=\left\{\begin{matrix}1,& x\in\Q \\ 0,& x\notin\Q \end{matrix}\right\}=
\lim_{n\to\infty}\lim_{m\to\infty}\cos^m2\pi\cdot n!\cdot x
 \eeq
(эту формулу надо понимать так: чтобы вычислить значение в точке $x$ нужно
сначала взять предел при $m\to\infty$, а затем предел при $n\to\infty$).
 \eex
\bpr Пусть $x\in\Q$, то есть $x=\frac{p}{q}$, $p\in\Z$, $q\in\N$, тогда при
$n\ge q$ мы получим:
$$
n!\div q
$$
$$
\Downarrow
$$
$$
n!\cdot x=n!\cdot\frac{p}{q}\in \Z
$$
$$
\Downarrow
$$
$$
\cos 2\pi\cdot n!\cdot x=1
$$
$$
\Downarrow
$$
$$
\cos^m 2\pi\cdot n!\cdot x=1
$$
$$
\Downarrow
$$
$$
\lim_{m\to\infty}\cos^m 2\pi\cdot n!\cdot x=1
$$
Это верно для $n\ge q$, то есть для почти всех $n$. Поэтому
$$
\lim_{n\to\infty}\lim_{m\to\infty}\cos^m 2\pi\cdot n!\cdot x=1=D(x)
$$

Наоборот, если $x\notin\Q$, то при любом $n\in\N$ получим:
$$
n!\cdot x\notin \Z
$$
$$
\Downarrow
$$
$$
|\cos 2\pi\cdot n!\cdot x|<1
$$
$$
\Downarrow
$$
$$
\lim_{m\to\infty}\cos^m 2\pi\cdot n!\cdot x=0
$$
Это верно для всех $n$, значит
$$
\lim_{n\to\infty}\lim_{m\to\infty}\cos^m 2\pi\cdot n!\cdot x=0=D(x)
$$
\epr

\end{multicols}\noindent\rule[10pt]{160mm}{0.1pt}

\subsection{Функциональный ряд и его область сходимости}

 \bit{
\item[$\bullet$] Пусть задана функциональная последовательность $\{ a_n \}$ и
из нее составлена новая функциональная последовательность $\{ S_N \}$ по
формуле
$$
S_N(x)=\sum_{n=1}^N a_n(x)=a_1(x)+a_2(x)+a_3(x)+...+a_N(x)
$$
Тогда такая пара последовательностей $\{ a_n\}$ и $\{S_N\}$ называется {\it
функциональным рядом}\index{ряд!функциональный} и обозначается
$$
\sum_{n=1}^{\infty} a_n(x)=a_1(x)+a_2(x)+a_3(x)+...
$$
Функции $a_n$ называются {\it слагаемыми}\index{слагаемое!функционального ряда}
ряда $\sum_{n=1}^{\infty} a_n$, а функции $S_N$ -- {\it частичными
суммами}\index{сумма!частичная!функционального ряда} этого ряда.

При каждом фиксированном значении переменной $x=x_0$ функциональная
последовательность $\{ S_N(x) \}$ превращается в числовую последовательность
$\{ S_N(x_0) \}$. Если эта числовая последовательность сходится
$$
\exists \lim_{N\to \infty} S_N(x_0)=S(x_0)=\sum_{n=1}^\infty a_n(x_0) \, \in \R
$$
$\l \text{то есть, если числовой ряд $\sum\limits_{n=1}^\infty a_n(x_0)$
сходится}\r$, то говорят, что {\it функциональный ряд $\sum\limits_{n=1}^\infty
a_n(x)$ сходится в точке}\index{сходимость!функционального ряда!поточечная}
$x=x_0$.

Множество всех точек $x$, в которых ряд $\sum\limits_{n=1}^\infty a_n(x)$
сходится (то есть, область сходимости функциональной последовательности $\{
S_N(x) \}$) называется {\it областью сходимости функционального
ряда}\index{область сходимости!функционального ряда} $\sum\limits_{n=1}^\infty
a_n(x)$.

\item[$\bullet$] Говорят, что функциональный ряд $\sum_{n=1}^{\infty} a_n$ {\it
сходится к функции $S$ поточечно} на множестве $E\subseteq\R$, , если для
всякого $x\in E$ числовой ряд $\sum_{n=1}^{\infty} a_n(x)$ сходится к числу
$S(x)$:
 \beq\label{19.3.1*}
\forall x\in E \quad \sum_{n=1}^N a_n(x)\underset{N\to \infty}{\longrightarrow}
S(x).
 \eeq
 }\eit

\noindent\rule{160mm}{0.1pt}\begin{multicols}{2}

\begin{ex}\label{ex-20.1.1} Областью сходимости ряда
$$
  \sum\limits_{n=1}^\infty x^n
$$
будет, как мы уже отмечали в теореме \ref{ex-18.1.4}, множество $(-1;1)$. Для
наглядности в таких случаях рисуют картинку:

%\pucture{0pt}{0pt}{ii-11.pcx}

\vglue80pt \noindent
\end{ex}

\begin{ex}\label{ex-20.1.2} Областью сходимости ряда Дирихле
$$
  \sum\limits_{n=1}^\infty \frac{1}{n^x}
$$
будет, как уже говорилось в примере \ref{ex-19.3.3}, множество $(1;+\infty)$.

%\pucture{0pt}{0pt}{ii-12.pcx}

\vglue80pt \noindent
\end{ex}

\begin{ex}\label{ex-20.1.3} Найдем область сходимости ряда
$$
  \sum\limits_{n=1}^\infty \frac{(-1)^n}{n^x}
$$
Для этого заметим, что
$$
\ml \frac{(-1)^n}{n^x}\mr=\frac{1}{n^x}\underset{n\to
\infty}{\longrightarrow}\begin{cases}
 0, & \text{если}\,\, x>0 \\
 1, & \text{если}\,\, x=0 \\
 \infty, & \text{если}\,\, x<0
\end{cases}
$$
Отсюда получаем:

1) если $x>0$, то $\frac{1}{n^x}\underset{n\to \infty}{\longrightarrow} 0$, и
ряд $\sum\limits_{n=1}^\infty \frac{(-1)^n}{n^x}$ будет сходиться по признаку
Лейбница;

2) если $x\le 0$, то $\ml \frac{(-1)^n}{n^x}\mr \underset{n\to
\infty}{\longrightarrow\kern-15pt{\Big/}} 0$, и ряд $\sum\limits_{n=1}^\infty
\frac{(-1)^n}{n^x}$ будет расходиться в силу необходимого условия сходимости.

%\pucture{0pt}{0pt}{ii-13.pcx}

\vglue80pt \noindent Вывод: областью сходимости ряда $\sum\limits_{n=1}^\infty
\frac{(-1)^n}{n^x}$ является множество $(0;+\infty)$.
\end{ex}

\begin{ex}\label{ex-20.1.4} Найдем область сходимости ряда
$$
  \sum\limits_{n=1}^\infty \frac{(-1)^n}{n\cdot 3^n\cdot \sqrt{(x+2)^n}}
$$
Для этого сразу заметим, что область определения нашего ряда (то есть общая
область определения слагаемых ряда) есть множество $x>-2$.

1. Исследуем теперь ряд на абсолютную сходимость (то есть пробуем применить к
нему теорему \ref{tm-18.5.6}):
 \begin{multline*}
\sum\limits_{n=1}^\infty \ml \frac{(-1)^n}{n\cdot 3^n\cdot
\sqrt{(x+2)^n}}\mr=\\= \sum\limits_{n=1}^\infty \frac{1}{n\cdot 3^n\cdot
\sqrt{(x+2)^n}}=\\= \sum\limits_{n=1}^\infty \frac{1}{n\cdot \l 3\cdot
\sqrt{x+2}\r^n}
 \end{multline*}
Находим число Коши:
$$
C=\lim_{n\to \infty}\sqrt[n]{ \frac{1}{n\cdot \l 3\cdot \sqrt{x+2}\r^n}}=
\frac{1}{3\cdot \sqrt{x+2}}
$$
Отсюда
 \begin{multline*}
C=\frac{1}{3\cdot \sqrt{x+2}}<1 \quad \Leftrightarrow \quad
\sqrt{x+2}>\frac{1}{3}\quad \Leftrightarrow \\ \Leftrightarrow \quad
x+2>\frac{1}{9}\quad \Leftrightarrow \quad x>\frac{1}{9}-2=-\frac{17}{9}
 \end{multline*}
Теперь можно сделать вывод, что ряд сходится при $x>-\frac{17}{9}$:

%\pucture{0pt}{0pt}{ii-14.pcx}

\vglue80pt

2. Посмотрим, что получается при $x<-\frac{17}{9}$:
\begin{multline*}
x<-\frac{17}{9}\; \Rightarrow \; \sqrt{x+2}<\frac{1}{3}\;
\Rightarrow \; \frac{1}{3\cdot \sqrt{x+2}}>1 \; \Rightarrow \\
\Rightarrow \; \ml \frac{(-1)^n}{3\cdot \sqrt{x+2}}\mr=\frac{1}{3\cdot
\sqrt{x+2}}\underset{n\to \infty}{\longrightarrow\kern-15pt{\Big/}} 0 \;
\Rightarrow
\\
\Rightarrow \; {\smsize \text{$\begin{pmatrix}\text{вспоминаем}
\\ \text{необходимое условие}\\ \text{сходимости}\end{pmatrix}$}}
\; \Rightarrow\\ \Rightarrow \; \text{ряд}\,\, \sum\limits_{n=1}^\infty
\frac{(-1)^n}{n\cdot 3^n\cdot \sqrt{(x+2)^n}}\,\, \text{расходится}
 \end{multline*} Отмечаем это на картинке:

%\pucture{0pt}{0pt}{ii-15.pcx}

\vglue80pt

3. Нам остается посмотреть, что будет при $x=-\frac{17}{9}$:
 \begin{multline*}
\sum\limits_{n=1}^\infty \frac{(-1)^n}{n\cdot 3^n\cdot \sqrt{(x+2)^n}}=
\sum\limits_{n=1}^\infty \frac{(-1)^n} {n\cdot 3^n\cdot
\sqrt{\l\frac{1}{9}\r^n}}=\\= \sum\limits_{n=1}^\infty \frac{(-1)^n} {n\cdot
3^n\cdot \frac{1}{3^n}}= \sum\limits_{n=1}^\infty \frac{(-1)^n}{n}
 \end{multline*}
Этот ряд сходится по признаку Лейбница:

%\pucture{0pt}{0pt}{ii-16.pcx}

\vglue80pt \noindent

Вывод: областью сходимости ряда $\sum\limits_{n=1}^\infty \frac{(-1)^n}{n\cdot
3^n\cdot \sqrt{(x+2)^n}}$ является множество $\left[ -\frac{17}{9};+\infty \r$.
\end{ex}

\begin{ers} Найдите область сходимости функциональных рядов:
 \biter{
 \item[1)] $\sum\limits_{n=1}^\infty \frac{\cos nx}{n}$ (признак Дирихле +
необходимое условие сходимости);

 \item[2)] $\sum\limits_{n=1}^\infty \frac{\cos nx}{n\sqrt{n}}$ (теорема об
абсолютной сходимости);

 \item[3)] $\sum\limits_{n=1}^\infty \frac{n!}{x^n}$ (необходимое условие
сходимости);

 \item[4)] $\sum\limits_{n=1}^\infty \frac{x^n}{1-x^n}$ (теорема об
абсолютной сходимости + признак Даламбера + необходимое условие сходимости);

 \item[5)] $\sum\limits_{n=1}^\infty \frac{2^n\cdot \sin^n x}{n^2}$
(теорема об абсолютной сходимости + признак Даламбера + необходимое условие
сходимости);

 \item[6)] $\sum\limits_{n=1}^\infty \l x^n+\frac{1}{n\cdot 2^n\cdot x^n}\r$ (теорема об абсолютной сходимости + признак Даламбера +
необходимое условие сходимости);

 \item[7)] $\sum\limits_{n=1}^\infty \frac{n}{n+1}\cdot \l \frac{x}{2x+1}\r^n$ (теорема об абсолютной сходимости + признак Даламбера +
необходимое условие сходимости);

 \item[8)] $\sum\limits_{n=1}^\infty \frac{n\cdot 3^{2n}}{2^n}\cdot
x^n\cdot (1-x)^n$ (теорема об абсолютной сходимости + признак Даламбера +
необходимое условие сходимости).
 }\eiter
 \end{ers}

\end{multicols}\noindent\rule[10pt]{160mm}{0.1pt}

\section{Равномерная сходимость}\label{SEC:ravnom-shodimost}

\subsection{Равномерная норма функции}

\bit{

\item[$\bullet$] {\it Равномерной нормой функции $f$ на множестве
$E$}\index{норма!равномерная}, называется величина
 \beq\label{ravnomernaya-norma}
\norm{f}_E=\norm{f(x)}_{x\in E}:=\sup_{x\in E} |f(x)|
 \eeq
(здесь предполагается, что функция $f$ определена на множестве $E$).
 }\eit

\noindent\rule{160mm}{0.1pt}\begin{multicols}{2}

\begin{exs}
 \begin{align*}
& \parallel x^2 \parallel_{x\in [0;2]}=\sup_{x\in [0;2]} |x^2|=4 \\
& \parallel x^2 \parallel_{x\in \R}=\sup_{x\in \R} |x^2|=\infty \\
& \parallel\sin x \parallel_{x\in \R}=\sup_{x\in \R} |\sin x|=1 \\
& \norm{\arctg x}_{x\in \R}=\sup_{x\in \R} |\arctg x|=\frac{\pi}{2}
 \end{align*}
\end{exs}
\end{multicols}\noindent\rule[10pt]{160mm}{0.1pt}

\bigskip

\centerline{\bf Свойства равномерной нормы}\label{svoistva-ravnom-normy}
 \bit{\it

\item[$1^\circ.$] Неотрицательность:
 \beq\label{neotr-normy}
\norm{f}_E\ge 0
 \eeq
причем
 \beq\label{norma=0}
\norm{f}_E=0\quad\Longleftrightarrow\quad  \forall x\in E\quad f(x)=0
 \eeq

\item[$2^\circ.$] Однородность:
  \beq\label{odnorod-normy}
\norm{\lambda\cdot f}_E = |\lambda|\cdot\norm{f}_E, \qquad \lambda\in\R
 \eeq

\item[$3^\circ.$] Полуаддитивность:
  \beq\label{poluaddit-normy}
\norm{f+g}_E\le \norm{f}_E+\norm{g}_E
 \eeq

\item[$4^\circ.$] Полумультпликативность:
  \beq\label{polumultipl-normy}
 \norm{f\cdot g}_E\le \norm{f}_E\cdot \norm{g}_E
 \eeq

\item[$5^\circ.$] Монотонность по параметру $E$:
 \beq\label{monot-normy-po-E}
 D\subseteq E\qquad\Longrightarrow\qquad \norm{f}_D\le \norm{f}_E
 \eeq

}\eit

\begin{proof} Все эти свойства напрямую следуют из свойств модуля. Например, полуаддитивность доказывается так:
 \begin{multline*}
\norm{f+g}_E = \sup_{x\in E} |f(x)+g(x)|\le\eqref{module-2^0}\le \sup_{x\in
E}\Big(|f(x)|+|g(x)|\Big)=\eqref{poluaddit-sup}=\\= \sup_{x\in E}
|f(x)|+\sup_{x\in E} |g(x)|= \norm{f}_E+\norm{g}_E
 \end{multline*}
 \end{proof}

\subsection{Равномерная сходимость функциональных последовательностей}

\bit{

\item[$\bullet$] Говорят, что {\it последовательность функций $\{ f_n \}$
стремится к функции $f$ равномерно\index{сходимость!функциональной
последовательности!равномерная} на множестве $E$}, если равномерная норма
разности $f_n-f$ на множестве $E$ стремится к нулю при $n\to \infty$:
 \beq
\norm{f_n-f}_E\underset{n\to \infty}{\longrightarrow} 0 \label{19.3.2}
 \eeq
Коротко это записывается так:
$$
f_n(x)\overset{x\in E}{\underset{n\to \infty}{\rightrightarrows}}
f(x)
$$
 }\eit

В следующих примерах мы увидим, что поточечная и равномерная
сходимость функциональных последовательностей -- не одно и то же.

\noindent\rule{160mm}{0.1pt}\begin{multicols}{2}

\begin{ex}\label{ex-19.3.1}
Очевидно,
$$
x^n\underset{n\to \infty}{\longrightarrow}\begin{cases}
 0, & \text{если}\,\, x\in (-1;1)
 \\
 1, & \text{если}\,\, x=1
\end{cases}.
$$
Поэтому если взять $E=\left[ -\frac{1}{2};\frac{1}{2}\right]$, то
последовательность $f_n(x)=x^n$  будет стремиться к функции $f(x)=0$ поточечно
на множестве $E$:
$$
x^n\overset{x\in \left[
-\frac{1}{2};\frac{1}{2}\right]}{\underset{n\to
\infty}{\longrightarrow}} 0
$$
Проверим, будет ли $f_n(x)=x^n$ стремиться к функции $f(x)=0$
равномерно на $E=\left[ -\frac{1}{2};\frac{1}{2}\right]$. Для этого
надо вычислить норму их разности:
 \begin{multline*}
||f_n(x)-f(x)||_{x\in E}= ||x^n-0||_{x\in \left[
-\frac{1}{2};\frac{1}{2}\right]}=\\= \sup_{x\in \left[
-\frac{1}{2};\frac{1}{2}\right]}|x^n-0|= \sup_{x\in \left[
-\frac{1}{2};\frac{1}{2}\right]}|x^n|= \l\frac{1}{2}\r^n
\underset{n\to \infty}{\longrightarrow} 0
 \end{multline*}
Вывод: $f_n(x)=x^n$ стремится к $f(x)=0$ поточечно и равномерно на множестве
$E=\left[ -\frac{1}{2};\frac{1}{2}\right]$:
$$
x^n\overset{x\in \left[ -\frac{1}{2};\frac{1}{2}\right]}{\underset{n\to
\infty}{\longrightarrow}} 0 \quad \& \quad x^n\overset{x\in \left[
-\frac{1}{2};\frac{1}{2}\right]}{\underset{n\to \infty}{\rightrightarrows}} 0
$$
\end{ex}

\begin{ex}\label{ex-19.3.2}
Та же самая последовательность $f_n(x)=x^n$ стремится к функции $f(x)=0$
поточечно и на более широком множестве $E=(-1;1)$:
$$
x^n\overset{x\in (-1;1)}{\underset{n\to \infty}{\longrightarrow}} 0
$$
Проверим, будет ли $f_n(x)=x^n$ стремиться к $f(x)=0$ равномерно
на $E=(-1;1)$. Для этого снова вычисляем норму их разности на этом
множестве:
 \begin{multline*}
||f_n(x)-f(x)||_{x\in E}= ||x^n-0||_{x\in (-1;1)}=\\= \sup_{x\in
(-1;1)}|x^n-0|= \sup_{x\in (-1;1)}|x^n|=1 \underset{n\to \infty}{
\longrightarrow\kern-15pt{\Big/} } 0
 \end{multline*}
Вывод: последовательность $f_n(x)=x^n$ стремится к функции $f(x)=0$ поточечно,
но не равномерно на множестве $E=(-1;1)$:
$$
x^n\overset{x\in (-1;1)}{\underset{n\to \infty}{\longrightarrow}} 0 \quad \&
\quad x^n\overset{x\in (-1;1)}{\underset{n\to \infty}{
\rightrightarrows\kern-12pt{\Big/} }} 0
$$
\end{ex}

\begin{ex}\label{ex-19.3.3}
Нетрудно увидеть, что последовательность $f_n(x)=\frac{\sin x}{n}$ стремится к
функции $f(x)=0$ поточечно на множестве $E=\R$:
$$
\frac{\sin x}{n}\overset{x\in \R}{\underset{n\to \infty}{ \longrightarrow }} 0
$$
Проверим, будет ли $f_n(x)=\frac{\sin x}{n}$ стремиться к $f(x)=0$ равномерно
на $E=\R$:
 \begin{multline*}
||f_n(x)-f(x)||_{x\in E}=
\parallel \frac{\sin x}{n}-0\parallel_{x\in \R}=\\=
\sup_{x\in \R}\ml \frac{\sin x}{n}-0\mr= \frac{1}{n}\cdot \sup_{x\in \R}|\sin
x|= \frac{1}{n}\underset{n\to \infty}{\longrightarrow} 0
 \end{multline*}
Вывод: $f_n(x)=\frac{\sin x}{n}$ стремится к $f(x)=0$ поточечно и равномерно на
множестве $E=\R$:
$$
\frac{\sin x}{n}\overset{x\in \R}{\underset{n\to \infty}{ \longrightarrow }} 0
\quad \& \quad \frac{\sin x}{n}\overset{x\in \R}{\underset{n\to \infty}{
\rightrightarrows }} 0
$$
\end{ex}

\begin{ex}\label{ex-19.3.4}
Похожая последовательность $f_n(x)=\sin \frac{x}{n}$ тоже стремится к функции
$f(x)=0$ поточечно на множестве $E=\R$:
$$
\sin \frac{x}{n}\overset{x\in \R}{\underset{n\to \infty}{ \longrightarrow }} 0
$$
Проверим, будет ли $f_n(x)=\sin \frac{x}{n}$ стремиться к $f(x)=0$ равномерно
на $E=\R$:
 \begin{multline*}
||f_n(x)-f(x)||_{x\in E}=
\parallel \sin \frac{x}{n}-0\parallel_{x\in \R}=\\=
\sup_{x\in \R}\ml \sin \frac{x}{n}-0\mr= \sup_{x\in \R}\ml\sin
\frac{x}{n}\mr=1 \underset{n\to \infty}{
\longrightarrow\kern-15pt{\Big/} } 0
 \end{multline*}
Вывод: последовательность $f_n(x)=\sin \frac{x}{n}$ стремится к функции
$f(x)=0$ поточечно, но не равномерно на множестве $E=\R$:
$$
\sin \frac{x}{n}\overset{x\in \R}{\underset{n\to \infty}{ \longrightarrow }} 0
\quad \& \quad \sin \frac{x}{n}\overset{x\in \R}{\underset{n\to \infty}{
\rightrightarrows\kern-12pt{\Big/} }} 0
$$
\end{ex}

\end{multicols}\noindent\rule[10pt]{160mm}{0.1pt}

Примеры 3.2 и 3.4 показывают, что из поточечной сходимости не
всегда следует равномерная:
$$
f_n(x)\overset{x\in E}{\underset{n\to \infty}{\longrightarrow}} f(x)
\quad \Longrightarrow\kern-15pt{\Big/}\quad f_n(x)\overset{x\in
E}{\underset{n\to \infty}{\rightrightarrows}} f(x)
$$
Но обратное всегда верно:

\begin{tm}
[о связи между поточечной и равномерной сходимостью]\label{tm-19.3.5} Если
последовательность функций $f_n$ стремится к функции $f$ равномерно на
множестве $E$, то $f_n$ стремится к $f$ поточечно на $E$:
$$
f_n(x)\overset{x\in E}{\underset{n\to \infty}{\longrightarrow}}
f(x) \quad \Longleftarrow \quad f_n(x)\overset{x\in
E}{\underset{n\to \infty}{\rightrightarrows}} f(x)
$$
\end{tm}\begin{proof} Доказательство представляет собой простую
логическую цепочку:
$$
f_n(x)\overset{x\in E}{\underset{n\to \infty}{\rightrightarrows}}
f(x)
$$
$$
\Downarrow
$$
$$
\parallel f_n(x)-f(x)\parallel_{x\in E} =
\sup_{x\in E}\ml f_n(x)-f(x) \mr \underset{n\to
\infty}{\longrightarrow}   0
$$
$$
\Downarrow
$$
$$
\forall x\in E \quad 0\le |f_n(x)-f(x)|\le \sup_{x\in E}\ml
f_n(x)-f(x) \mr \underset{n\to \infty}{\longrightarrow}  0
$$
$$
\Downarrow
$$
$$
\forall x\in E \quad |f_n(x)-f(x)| \underset{n\to
\infty}{\longrightarrow}  0
$$
$$
\Downarrow
$$
$$
\forall x\in E \quad f_n(x)-f(x) \underset{n\to
\infty}{\longrightarrow} 0
$$
$$
\Downarrow
$$
$$
\forall x\in E \quad f_n(x) \underset{n\to
\infty}{\longrightarrow} f(x)
$$
$$
\Downarrow
$$
$$
f_n(x)\overset{x\in E}{\underset{n\to \infty}{\longrightarrow}}
f(x) \qquad $$ \end{proof}

Из этой теоремы становится понятно, что в примерах \ref{ex-19.3.2}
и \ref{ex-19.3.4} необязательно было писать в выводах, что
$f_n(x)$ стремится к $f(x)$ поточечно на $E$, потому что это
автоматически следует из равномерного стремления $f_n(x)$ к $f(x)$
на $E$.

Рассмотрим теперь еще несколько примеров.

\noindent\rule{160mm}{0.1pt}\begin{multicols}{2}

\begin{ex}\label{ex-19.3.6} Рассмотрим функциональную последовательность
$f_n(x)=x^n-x^{n+1}$ на множестве $E=[0;1]$. Ясно, что она
поточечно стремится к нулю:
$$
x^n-x^{n+1}\overset{x\in [0;1]}{\underset{n\to \infty}{\longrightarrow}} 0
$$
Чтобы проверить стремится ли она равномерно, нам нужно найти норму
 \begin{multline*}
\norm{f_n(x)-f(x)}_{x\in [0;1]}= \norm{x^n-x^{n+1}}_{x\in [0;1]}=\\= \sup_{x\in
[0;1]}\abs{x^n-x^{n+1}}
 \end{multline*}
Это можно сделать, честно исследовав функцию $f_n(x)=x^n-x^{n+1}$
на возрастание и убываение с помощью производной:
 \begin{multline*}
f_n'(x)=nx^{n-1}-(n+1)x^n=\\=x^{n-1}\lll n-(n+1)x\rrr=0 \quad
\Leftrightarrow \quad x\in \lll 0;\frac{n}{n+1}\rrr
 \end{multline*}
Интервалы возрастания и убывания на отрезке $[0;1]$:

%\pucture{0pt}{0pt}{ii-9.pcx}

\vglue70pt \noindent Значение функции в точке максимума:
 \begin{multline*}
f_n \l
\frac{n}{n+1}\r=\frac{n^n}{(n+1)^n}-\frac{n^{n+1}}{(n+1)^{n+1}}=\\=
\frac{n^n}{(n+1)^n}\l 1-\frac{n}{n+1}\r= \frac{n^n}{(n+1)^n}\cdot
\frac{1}{n+1}=\\= \frac{1}{\l 1+\frac{1}{n}\r^n}\cdot \frac{1}{n+1}
 \end{multline*}
График:

%\pucture{0pt}{0pt}{ii-10.pcx}

\vglue120pt \noindent Норма на отрезке $[0;1]$:
 \begin{multline*}
\norm{f_n(x)-f(x)}_{x\in [0;1]}= \sup_{x\in [0;1]}\abs{x^n-x^{n+1}}=\\=f_n \l
\frac{n}{n+1}\r= \frac{1}{\l 1+\frac{1}{n}\r^n}\cdot
\frac{1}{n+1}\underset{n\to \infty}{\longrightarrow}\frac{1}{e}\cdot 0=0
 \end{multline*}
Вывод: последовательность $f_n(x)=x^n-x^{n+1}$ стремится к $f(x)=0$ (поточечно
и) равномерно на множестве $E=[0;1]$:
$$
x^n-x^{n+1}\overset{x\in [0;1]}{\underset{n\to \infty}{ \rightrightarrows }} 0
$$
\end{ex}

\begin{ex}\label{ex-19.3.7} Рассмотрим функциональную последовательность
$f_n(x)=\sqrt{x^2+\frac{1}{n^2}}$ на множестве $E=\R$. Ясно, что она поточечно
стремится к функции $f(x)=|x|$:
$$
\sqrt{x^2+\frac{1}{n^2}}\overset{x\in \R}{\underset{n\to
\infty}{\longrightarrow}}\sqrt{x^2+0}=|x|
$$
Проверим, будет ли это стремление равномерным:
 \begin{multline*}
\parallel f_n(x)-f(x) \parallel_{x\in \R}=
\parallel \sqrt{x^2+\frac{1}{n^2}}-|x| \parallel_{x\in \R}=\\=
\sup_{x\in \R}\ml \sqrt{x^2+\frac{1}{n^2}}-|x|  \mr=
{\smsize\begin{pmatrix}\text{домножаем на}\\
\text{сопряженный}\\ \text{радикал}\end{pmatrix}}=\\= \sup_{x\in
\R}\ml \frac{\l\sqrt{x^2+\frac{1}{n^2}}\r^2-|x|^2}
{\sqrt{x^2+\frac{1}{n^2}}+|x|}  \mr=\\= \sup_{x\in \R}\ml
\frac{x^2+\frac{1}{n^2}-x^2} {\sqrt{x^2+\frac{1}{n^2}}+|x|}  \mr=
\sup_{x\in \R}\ml \frac{\frac{1}{n^2}}
{\sqrt{x^2+\frac{1}{n^2}}+|x|}\mr=\\= \frac{1}{n^2}\cdot
\frac{1}{\inf_{x\in \R}\l \sqrt{x^2+\frac{1}{n^2}}+|x|\r} =\\=
\frac{1}{n^2}\cdot \frac{1}{\sqrt{0+\frac{1}{n^2}}+0}=
\frac{\frac{1}{n^2}}{\frac{1}{n}}=\frac{1}{n}\underset{n\to
\infty}{\longrightarrow} 0
\end{multline*}
Вывод: последовательность $f_n(x)=\sqrt{x^2+\frac{1}{n^2}}$ стремится к
$f(x)=|x|$ (поточечно и) равномерно на множестве $E=\R$:
$$
\sqrt{x^2+\frac{1}{n^2}}\overset{x\in \R}{\underset{n\to \infty}{
\rightrightarrows }} |x|
$$
\end{ex}

\begin{ex}\label{ex-19.3.8} Рассмотрим функциональную последовательность
$f_n(x)=\arctg \frac{x}{n}$ на множестве $E=[-a;a]$. Ясно, что она
поточечно стремится к функции $f(x)=0$:
$$
\arctg \frac{x}{n}\qquad \overset{x\in [-a;a]}{\underset{n\to
\infty}{\longrightarrow}}\qquad \arctg 0=0
$$
Проверим, будет ли это стремление равномерным:
\begin{multline*}
\parallel f_n(x)-f(x) \parallel_{x\in [-a;a]}=
\parallel \arctg \frac{x}{n}\parallel_{x\in [-a;a]}=\\=
\sup_{x\in [-a;a]}\ml \arctg \frac{x}{n}  \mr=\\= \sup_{x\in
[-a;a]}\ml \int_{-a}^x \l \arctg \frac{t}{n}\r'_t \, \d t  \mr=\\=
\sup_{x\in [-a;a]}\ml \int_{-a}^x
\frac{\frac{1}{n}}{1+\frac{t^2}{n^2}}\, \d t \mr=\\= \sup_{x\in
[-a;a]}\int_{-a}^x \frac{\frac{1}{n}}{1+\frac{t^2}{n^2}}\, \d t \le\\
\le \sup_{x\in [-a;a]}\int_{-a}^x  \frac{\frac{1}{n}}{1+0}\, \d t =
\frac{1}{n}\cdot (a-a) \underset{n\to \infty}{\longrightarrow} 0
\end{multline*}
Вывод: последовательность $f_n(x)=\arctg \frac{x}{n} $ стремится к $f(x)=0$
(поточечно и) равномерно на множестве $E=[-a;a]$:
$$
\arctg \frac{x}{n}\overset{x\in [-a;a]}{\underset{n\to \infty}{
\rightrightarrows }} 0
$$
\end{ex}

\begin{ers} Проверьте соотношения:
 \biter{
\item[1)] $ \arctg \frac{x}{n}\overset{x\in \R}{\underset{n\to \infty}{
\rightrightarrows }} 0 $;

\item[2)] $ \sin \frac{x}{n}\overset{x\in [-a;a]}{\underset{n\to
\infty}{ \rightrightarrows }} 0 $ {\smsize (воспользоваться тем же
приемом, что и в примере \ref{ex-19.3.8})};

\item[3)] $ n^x \overset{x\in (-\infty;0)}{\underset{n\to \infty}{
\rightrightarrows }} 0 $;

\item[4)] $ \frac{1}{1+x^n}\overset{|x|<\frac{3}{4}}{\underset{n\to
\infty}{ \rightrightarrows }} 1 $;

\item[5)] $ \frac{1}{1+x^n}\overset{|x|<1}{\underset{n\to \infty}{
\rightrightarrows }} 1 $;

\item[6)] $ \frac{n^x-n^{-x}}{n^x+n^{-x}}\overset{x>0}{\underset{n\to
\infty}{ \rightrightarrows }} 1 $;

\item[7)] $ \sin^{2n} x\overset{|x|<1}{\underset{n\to \infty}{
\rightrightarrows }} 0 $;

\item[8)] $ \sin^{2n} x\overset{|x|<\frac{\pi}{2}}{\underset{n\to
\infty}{ \rightrightarrows }} 0 $;
 }\eiter
\end{ers}

\end{multicols}\noindent\rule[10pt]{160mm}{0.1pt}

\paragraph{Свойства равномерно сходящихся функциональных последовательностей.}

Перечислим некоторые свойства равномерно сходящихся
последовательностей.

{\it
 \bit{
\item[$1^0$.]\label{nepreryvnost-ravnom-predela} Если функции $f_n$ непрерывны
на множестве $E$ и равномерно сходятся на нем к функции $f$
$$
f_n(x)\overset{x\in E}{\underset{n\to \infty}{\rightrightarrows}}
f(x),
$$
то функция $f$ непрерывна на множестве $E$.

\item[$2^0$.] Если функции $f_n$ непрерывны и равномерно сходятся на интервале
$(a,b)$, то для всякой точки $c\in (a,b)$ справедливо равенство
 \beq
\lim_{x\to c}\lim_{n\to \infty} f_n(x) =\lim_{n\to \infty}\lim_{x\to
c}f_n (x) \label{19.4.1}
 \eeq

\item[$3^0$.] Если функции $f_n$ непрерывны и равномерно сходятся на отрезке
$[a,b]$, то
 \beq
\int_a^b \lim_{n\to \infty} f_n(x) \, \d x =\lim_{n\to
\infty}\int_a^b  f_n (x) \, \d x \label{19.4.2}
 \eeq

\item[$4^0$.] Пусть $f_n$ -- гладкие функции на $[a,b]$, их производные $f_n'$
равномерно сходятся на отрезке $[a,b]$, и хотя бы для одной точки $c\in [a,b]$
существует конечный предел
 \beq
  \lim_{n\to \infty} f_n(c)
\label{19.4.3}
 \eeq
Тогда функции $f_n$ равномерно сходятся на отрезке $[a,b]$, причем
 \beq
\lll \lim_{n\to \infty} f_n(x) \rrr' =\lim_{n\to \infty} f_n' (x)
\label{19.4.4}
 \eeq
 }\eit
}

\begin{proof} 1. Пусть функции $f_n(x)$ непрерывны на множестве
$E$ и равномерно сходятся к функции $f$ на $E$. Покажем, что $f(x)$
непрерывна на $E$, то есть (см. определение $\S 1$ главы 3) что для любой
последовательности $x_i\in E$, сходящейся к точке $c\in E$
$$
  x_i\underset{n\to\infty}{\longrightarrow} c
$$
обязательно выполняется соотношение
$$
 f(x_i)\underset{i\to\infty}{\longrightarrow} f(c)
$$
Для этого зафиксируем какой-нибудь $\varepsilon>0$ и покажем, что
почти все числа $f(x_i)$ содержатся в интервале
$(f(c)-\varepsilon,f(c)+\varepsilon)$, то есть, что для почти всех
номеров $i\in \mathbb{N}$ выполняется неравенство
 \beq
 |f(x_i)-f(c)|<\varepsilon
\label{19.4.5}
 \eeq
Для этого заметим, что из равномерной сходимости $f_n(x)$ к $f(x)$
то есть из соотношения
$$
  \parallel f_n(x)-f(x) \parallel_{x\in E}\underset{n\to\infty}{\longrightarrow} 0
$$
следует что найдется номер $N\in \mathbb{N}$, для которого
$$
\parallel f_N(x)-f(x) \parallel_{x\in E}<\frac{\varepsilon}{3}
$$
Поскольку функция $f_N(x)$ непрерывна на множестве $E$, должно
выполняться соотношение
$$
 f_N(x_i)\underset{i\to\infty}{\longrightarrow} f_N(c)
$$
Значит, найдется такой номер $I\in \mathbb{N}$, что
$$
\forall i\ge I \quad |f_N(x_i)-f_N(c)|<\frac{\varepsilon}{3}
$$
Теперь получаем $\forall i\ge I$
\begin{multline*}
|f(x_i)-f(c)|= |f(x_i)-f_N(x_i)+f_N(x_i)-f_N(c)+f_N(c)-f(c)|\le \\
\le |f(x_i)-f_N(x_i)|+|f_N(x_i)-f_N(c)|+|f_N(c)-f(c)| \le
\\
\le
\parallel f(x)-f_N(x)\parallel_{x\in E}+|f_N(x_i)-f_N(c)|+
\parallel f_N(c)-f(c) \parallel_{x\in E}
<\frac{\varepsilon}{3}+\frac{\varepsilon}{3}+\frac{\varepsilon}{3}=
\varepsilon
\end{multline*}
Мы получили то, что хотели: формула \eqref{19.4.5} выполняется для
почти всех $i\in \mathbb{N}$.

2. Пусть функции $f_n(x)$ непрерывны и равномерно сходятся на
интервале $(a,b)$ к какой-то функции $f$:
$$
f_n(x)\overset{x\in (a,b)}{\underset{n\to
\infty}{\rightrightarrows}} f(x)
$$
Тогда функция  $f(x)$ будет непрерывна в силу уже доказанного
свойства $1^0$, поэтому для всякой точки $c\in (a,b)$ мы получим
$$
\lim_{x\to c}\lim_{n\to \infty} f_n(x)= \lim_{x\to c}
f(x)=f(c)=\lim_{n\to \infty} f_n (c) =\lim_{n\to \infty}\lim_{x\to c}
f_n (x)
$$

3. Пусть функции $f_n(x)$ непрерывны и равномерно сходятся на
отрезке $[a,b]$ к какой-нибудь функции $f$. Тогда
\begin{multline*}\left| \int_a^b f_n(x) \, \d x - \int_a^b f(x)\, \d x \right|=
\left| \int_a^b \lll f_n(x)-f(x) \rrr \, \d x \right|\le \int_a^b
\left|  f_n(x)-f(x)  \right| \, \d x\le \\ \le \int_a^b \sup_{x\in
[a,b]}\left|  f_n(x)-f(x)  \right| \, \d x= \sup_{x\in [a,b]}\left|
f_n(x)-f(x)  \right|\cdot \int_a^b 1 \, \d x=
\parallel  f_n(x)-f(x) \parallel_{x\in [a,b]}\cdot (b-a)
\underset{n\to\infty}{\longrightarrow} 0
\end{multline*}
То есть
$$
\int_a^b f_n(x) \, \d x - \int_a^b f(x)\, \d x
\underset{n\to\infty}{\longrightarrow} 0
$$
или
$$
\int_a^b f_n(x) \, \d x
\underset{n\to\infty}{\longrightarrow}\int_a^b f(x)\, \d x
$$
а это уже означает, что выполняется \eqref{19.4.2}:
$$
\lim_{n\to \infty}\int_a^b  f_n (x) \, \d x= \int_a^b \lim_{n\to
\infty} f_n(x) \, \d x
$$

4. Пусть $f_n$ -- гладкие функции на $[a,b]$, их производные $f_n'$ равномерно
сходятся на отрезке $[a,b]$ к какой-то функции $g(x)$,
 \beq
f_n'(x) \overset{x\in E}{\underset{n\to
\infty}{\rightrightarrows}} g(x), \label{19.4.6}
 \eeq
и хотя бы для одной точки $c\in [a,b]$ существует конечный предел
\eqref{19.4.3}:
$$
 H=\lim_{n\to \infty} f_n(c).
$$
Тогда, в силу доказанного уже свойства $1^0$, функция $g(x)$
непрерывна на отрезке $[a,b]$. Поэтому можно рассмотреть функцию
 \beq
f(x)=H+\int_c^x g(t) \, \d t, \quad x\in [a,b] \label{19.4.7}
 \eeq
Покажем, что функции $f_n(x)$ равномерно сходятся к функции $f$
на отрезке $[a,b]$. Заметим, что
\begin{multline*}
|f_n(x)-f(x)|= \left| \lll f_n(c)+\int_c^x f_n'(t) \, \d x \rrr -
\lll H+\int_c^x g(t) \, \d t \rrr \right|=\\= \left| \lll
f_n(c)-H\rrr+\lll \int_c^x f_n'(t) \, \d x - \int_c^x g(t) \, \d t
\rrr \right|\le \left| f_n(c)-H \right| + \left| \int_c^x \lll
f_n'(t) - g(t) \rrr \, \d t \right|\le \\ \le \left| f_n(c)-H
\right| + \int_c^x \left| f_n'(t) -  g(t) \right|\, \d t\le \left|
f_n(c)-H \right| + \int_c^x \sup_{t\in [a,b]}\left| f_n'(t) -
g(t) \right|\, \d t=\\= \left| f_n(c)-H \right| + \sup_{t\in
[a,b]}\left| f_n'(t) -  g(t) \right| \cdot |x-c|= \left| f_n(c)-H
\right| +
\parallel f_n'(t) -  g(t) \parallel_{t\in [a,b]}\cdot |x-c|
\end{multline*}
Отсюда
\begin{multline*}\parallel f_n(x)-f(x)\parallel_{x\in [a,b]} =
\sup_{x\in [a,b]} |f_n(x)-f(x)|= |f_n(c)-H| +
\parallel f_n'(t) -  g(t) \parallel_{t\in [a,b]}\cdot
\sup_{x\in [a,b]} |x-c|\le \\ \le |f_n(c)-H|+
\parallel f_n'(t) -  g(t) \parallel_{t\in [a,b]}\cdot |b-a|
\underset{n\to \infty}{\longrightarrow} 0+0=0
\end{multline*}
То есть,
 \beq
f_n(x) \overset{x\in E}{\underset{n\to \infty}{\rightrightarrows}}
f(x). \label{19.4.8}
 \eeq
Теперь формула \eqref{19.4.4} становится очевидной:
$$
\lll \lim_{n\to \infty} f_n(x) \rrr'=(4.8)= \lll f(x)
\rrr'=(4.7)=g(x)=(4.6)=\lim_{n\to \infty} f_n' (x) \quad $$
\end{proof}

\paragraph{Критерий Коши равномерной сходимости функциональной последовательности.}

\begin{tm}
[критерий Коши равномерной сходимости функциональной последовательности]
\label{tm-19.5.1} Функциональная последовательность $\{ f_n(x) \}$ тогда и
только тогда равномерно сходится на множестве $E$, когда она удовлетворяет
следующим двум эквивалентным условиям:
 \bit{
\item[(i)]
для любых двух бесконечно больших последовательностей индексов $\{
p_i \},\, \{ q_i \}\subseteq \mathbb{N}$
 \beq
p_i \underset{i\to \infty}{\longrightarrow}\infty, \quad q_i
\underset{i\to \infty}{\longrightarrow}\infty \label{19.5.1}
 \eeq
соответствующие подпоследовательности $\{ f_{p_i}(x) \}$ и $\{ f_{q_i}(x) \}$
последовательности $\{ f_n(x) \}$ равномерно стремятся друг к другу на
множестве $E$:
 \beq
\parallel f_{p_i}(x)-f_{q_i}(x)\parallel_{x\in E}\underset{i\to \infty}{\longrightarrow} 0 \label{19.5.2}
 \eeq
\item[(ii)] для любой последовательности $\{ l_i \}\subseteq \mathbb{N}$ и
любой бесконечно большой последовательности $\{ k_i \}\subseteq \mathbb{N}$
 \beq
   k_i \underset{i\to \infty}{\longrightarrow}\infty
\label{19.5.3}
 \eeq
выполняется соотношение
 \beq
\parallel f_{k_i+l_i}(x)-f_{k_i}(x)\parallel_{x\in E}\underset{i\to \infty}{\longrightarrow} 0 \label{19.5.4}
 \eeq
 }\eit
\end{tm}

\begin{proof} Убедимся сначала, что условия $(i)$ и $(ii)$ действительно
эквивалентны.

1. Импликация $(i)\Rightarrow (ii)$ очевидна: если выполняется $(i)$ то есть
любые две подпоследовательности $\{ f_{p_i}(x) \}$ и $\{ f_{q_i}(x) \}$
последовательности $\{ f_n(x) \}$ стремятся друг к другу
$$
\parallel f_{p_i}(x)-f_{q_i}(x)\parallel_{x\in E}\underset{i\to \infty}{\longrightarrow} 0,
$$
то автоматически выполняется и $(ii)$, потому что для любой последовательности
$\{ l_i \}\subseteq \mathbb{N}$ и любой бесконечно большой последовательности
$\{ k_i \}\subseteq \mathbb{N}$ мы можем положить $p_i=k_i+l_i$ и $q_i=k_i$, и
тогда будет выполняться соотношение
$$
\parallel f_{k_i+l_i}(x)-f_{k_i}(x)\parallel_{x\in E}=
\parallel f_{p_i}(x)-f_{q_i}(x)\parallel_{x\in E}\underset{i\to \infty}{\longrightarrow} 0
$$

2. Докажем импликацию $(ii)\Rightarrow (i)$. Пусть выполняется $(ii)$, то есть
для любой последовательности $\{ l_i \}\subseteq \mathbb{N}$ и любой бесконечно
большой последовательности $\{ k_i \}\subseteq \mathbb{N}$ выполняется
соотношение \eqref{19.5.4}. Возьмем какие-нибудь две бесконечно большие
последовательности индексов
$$
p_i \underset{i\to \infty}{\longrightarrow}\infty, \quad q_i
\underset{i\to \infty}{\longrightarrow}\infty
$$
и положим
$$
k_i=\min\{ p_i; q_i \}, \quad l_i=\max\{ p_i; q_i \}-k_i
$$
Тогда
$$
k_i+l_i=\max\{ p_i; q_i \}
$$
и поэтому
 \begin{multline*}
f_{p_i}(x)-f_{q_i}(x)= \lll
\begin{array}{c}
{ f_{\max\{ p_i; q_i \}}(x)-f_{\min\{ p_i; q_i \}}(x), \quad
\text{если}\,\, p_i>q_i }\\{ f_{\min\{ p_i; q_i \}}(x)-f_{\max\{
p_i; q_i \}}(x), \quad \text{если}\,\, p_i<q_i
}\\
{ 0, \quad \text{если}\,\, p_i=q_i }\end{array}\rrr=\\= \lll
\begin{array}{c}{
f_{ k_i+l_i}(x)-f_{k_i}(x), \quad \text{если}\,\, p_i>q_i }\\{
f_{k_i}(x)-f_{k_i+l_i}(x), \quad \text{если}\,\, p_i<q_i
}\\
{ 0, \quad \text{если}\,\, p_i=q_i }\end{array}\rrr
 \end{multline*}
В любом случае получается
$$
\parallel f_{p_i}(x)-f_{q_i}(x)\parallel_{x\in E}
=
\parallel f_{k_i+l_i}(x)-f_{k_i}(x)\parallel_{x\in E}
$$
поэтому
$$
\parallel f_{p_i}(x)-f_{q_i}(x)\parallel_{x\in E} =
\parallel f_{k_i+l_i}(x)-f_{k_i}(x)\parallel_{x\in E}\underset{i\to \infty}{\longrightarrow} 0
$$
и значит,
$$
\parallel f_{p_i}(x)-f_{q_i}(x)\parallel_{x\in E}\underset{i\to \infty}{\longrightarrow} 0
$$
Последовательности $\{ p_i \}$ и $\{ q_i \}$ здесь выбирались с
самого начала произвольными. Это значит, что выполняется $(i)$.
Таким образом, мы доказали, что из $(ii)$ следует $(i)$.

3. Докажем теперь, что если последовательность $\{ f_n(x) \}$ сходится
равномерно к некоторой функции $f$, то есть
$$
\parallel  f_n(x)-f(x)\parallel_{x\in E}\underset{i\to \infty}{\longrightarrow} 0
$$
то $\{ f_n(x) \}$ обязательно должна обладать свойством $(i)$. Возьмем любые
две бесконечно большие последовательности индексов $\{ p_i \},\, \{ q_i
\}\subseteq \mathbb{N}$
$$
p_i \underset{i\to \infty}{\longrightarrow}\infty, \quad q_i
\underset{i\to \infty}{\longrightarrow}\infty
$$
Соответствующие подпоследовательности $\{ f_{p_i}(x) \}$ и $\{ f_{q_i}(x) \}$
тоже сходится равномерно к $f(x)$, потому что
$$
\parallel  f_{p_i}(x)-f(x)\parallel_{x\in E}\underset{i\to \infty}{\longrightarrow} 0, \quad
\parallel  f_{q_i}(x)-f(x)\parallel_{x\in E}\underset{i\to \infty}{\longrightarrow} 0
$$
поэтому
\begin{multline*}
0\le \parallel  f_{p_i}(x)-f_{q_i}(x) \parallel_{x\in E} =
\parallel  f_{p_i}(x)-f(x)+f(x)-f_{q_i}(x) \parallel_{x\in
E}\le\eqref{poluaddit-normy}\le
\\
\le
\parallel  f_{p_i}(x)-f(x) \parallel_{x\in E}
+
\parallel f(x)-f_{q_i}(x) \parallel_{x\in E}\underset{i\to \infty}{\longrightarrow} 0+0=0
\end{multline*}
То есть,
$$
\parallel  f_{p_i}(x)-f_{q_i}(x) \parallel_{x\in E}\underset{i\to \infty}{\longrightarrow} 0
$$
и это значит, что $\{ f_n(x) \}$ действительно должна обладать
свойством $(i)$.

4. Пусть выполнено $(i)$, то есть для любых последовательностей
индексов
$$
p_i \underset{i\to \infty}{\longrightarrow}\infty, \quad q_i
\underset{i\to \infty}{\longrightarrow}\infty
$$
выполняется соотношение
 \beq
\parallel f_{p_i}(x)-f_{q_i}(x)\parallel_{x\in E}\underset{i\to \infty}{\longrightarrow} 0 \label{19.5.5}
 \eeq
Тогда для всякой точки $a\in E$ мы получим
$$
0\le |f_{p_i}(a)-f_{q_i}(a)|\le \sup_{x\in E}
|f_{p_i}(x)-f_{q_i}(x)|=
\parallel f_{p_i}(x)-f_{q_i}(x)\parallel_{x\in E}\underset{i\to \infty}{\longrightarrow} 0
$$
откуда по теореме о двух милиционерах,
$$
f_{p_i}(a)-f_{q_i}(a) \underset{i\to \infty}{\longrightarrow} 0
$$
Это значит, что для всякой точки $a\in E$ числовая последовательность $\{
f_n(a)\}$ сходится по критерию Коши \ref{Cauchy-crit-seq}. Обозначим через
$f(a)$ ее предел:
 \beq
  f(a)=\lim_{n\to \infty} f_n(a), \quad a\in E
\label{19.5.6}
 \eeq
Мы получили функцию $f$, определенную на множестве $x\in E$.

Покажем, что функциональная последовательность $f_n(x)$ равномерно сходится к
функции $f$ на множестве $E$.
$$
\parallel f_n(x)-f(x) \parallel_{x\in E}\underset{n\to \infty}{\longrightarrow} 0
$$
Предположим, что это не так, то есть что числовая последовательность
$c_n=\parallel f_n(x)-f(x)
\parallel_{x\in E}$ не стремится к нулю:
 \beq
c_n=\parallel f_n(x)-f(x) \parallel_{x\in E}\underset{n\to
\infty}{\not\longrightarrow} 0 \label{19.5.7}
 \eeq
Тогда, по свойству $2^0 \S 3$ главы 2, существуют подпоследовательность
$p_i\underset{i\to \infty}{\longrightarrow}\infty$ и число $\varepsilon>0$
такие, что
$$
c_{p_i}>\varepsilon
$$
То есть,
$$
\parallel f_{p_i}(x)-f(x) \parallel_{x\in E}=
\sup_{x\in E}\ml f_{p_i}(x)-f(x) \mr > \varepsilon
$$
Отсюда следует, что для всякого $i$ существует точка $a_i\in E$
такая что
$$
\ml f_{p_i}(a_i)-f(a_i) \mr > \varepsilon
$$
то есть
$$
f_{p_i}(a_i)-\varepsilon< f(a_i) <f_{p_i}(a_i)-\varepsilon
$$
Вспомним теперь формулу \eqref{19.5.6}: число $f(a_i)$ является пределом
последовательности $f_n(a_i)$ и лежит в интервале $\l f_{p_i}(a_i)-\varepsilon,
f_{p_i}(a_i)-\varepsilon \r$ Значит, почти все числа $f_n(a_i)$ лежат в этом
интервале:
$$
f_{p_i}(a_i)-\varepsilon< f_n(a_i) <f_{p_i}(a_i)-\varepsilon
$$
Обозначим через $q_i$ какое-нибудь значение $n$ с таким свойством:
$$
f_{p_i}(a_i)-\varepsilon< f_{q_i} (a_i) <f_{p_i}(a_i)-\varepsilon
$$
Из этого двойного неравенства получаем
$$
|f_{p_i}(a_i)-f_{q_i} (a_i)|>\varepsilon
$$
а отсюда
$$
\parallel f_{p_i}(x)-f_{q_i} (x)\parallel_{x\in E}\ge
|f_{p_i}(a_i)-f_{q_i} (a_i)|>\varepsilon
$$
Это противоречит формуле \eqref{19.5.5}. Таким образом, наше
предположение (5.7) неверно. \end{proof}

\subsection{Равномерная сходимость функционального ряда}

Пусть нам дан функциональный ряд
$$
\sum\limits_{n=1}^\infty a_n(x)
$$
и пусть $S_N$ -- последовательность его частичных сумм:
$$
S_N(x)=\sum\limits_{n=1}^N a_n(x)
$$
Говорят, что
 \bit{
\item[---] ряд $\sum\limits_{n=1}^\infty a_n$ сходится {\it
поточечно}\index{сходимость!функционального ряда!поточечная на множестве} на
множестве $E$ к сумме $S$, если $S_N$ стремится к $S$ поточечно на $E$:
 \beq
S_N(x)\overset{x\in E}{\underset{N\to \infty}{\longrightarrow}} S(x)
\label{20.2.1}
 \eeq
(то есть $\forall x\in E \quad S_N(x)\underset{N\to \infty}{\longrightarrow}
S(x)$);

\item[---] ряд $\sum\limits_{n=1}^\infty a_n$ сходится {\it
равномерно}\index{сходимость!функционального ряда!равномерная} на множестве $E$
к сумме $S$, если $S_N$ стремится к $S$ равномерно на $E$:
 \beq
S_N(x)\overset{x\in E}{\underset{N\to \infty}{\rightrightarrows}} S(x)
\label{20.2.2}
 \eeq
(то есть $||S_N(x)-S(x)||_{x\in E}\underset{N\to \infty}{\longrightarrow} 0$).
 }\eit

\noindent\rule{160mm}{0.1pt}\begin{multicols}{2}

\begin{ex}\label{ex-20.2.1} Выясним, будет ли ряд
$$
  \sum_{n=1}^\infty (1-x)\cdot x^n
$$
сходиться равномерно на множестве $E=(0;1)$. Для этого найдем предел его
частичных сумм:
\begin{multline*}
S_N(x)=\sum_{n=1}^N (1-x)\cdot x^n=\sum_{n=1}^N (x^n-x^{n+1})=\\=
(x-x^2)+(x^2-x^3)+(x^3-x^4)+...+(x^N-x^{N+1})=\\= x-x^2
 \put(-9,20){\put(2,0){\line(1,0){22}}
 \put(0,-7){$\downarrow$}\put(22,-7){$\downarrow$}}
+x^2 -x^3
 \put(-9,20){\put(2,0){\line(1,0){22}}
 \put(0,-7){$\downarrow$}\put(22,-7){$\downarrow$}}
+x^3 -x^4
 \put(-9,20){\put(2,0){\line(1,0){18}}
 \put(0,-7){$\downarrow$}\put(18,-7){$\downarrow$}}
 +...  \put(-1,20){\put(2,0){\line(1,0){14}}
 \put(0,-7){$\downarrow$}\put(14,-7){$\downarrow$}}
+x^N -x^{N+1}=\\=x-x^{N+1}\overset{x\in (0;1)}{\underset{N\to
\infty}{\longrightarrow}} x=S(x)
\end{multline*}
Теперь найдем норму остатка
 \begin{multline*}
||S_N(x)-S(x)||_{x\in (0;1)}=\\=|| \l x-x^{N+1}\r-x ||_{x\in (0;1)}=\\=
\sup_{x\in (0;1)}\left|x^{N+1}\right|=1 \underset{N\to
\infty}{\longrightarrow\kern-15pt{\Big/}} 0
 \end{multline*}
Вывод: ряд $\sum\limits_{n=1}^\infty (1-x)\cdot x^n$ сходится на множестве
$(0;1)$ поточечно, но не равномерно к сумме $S(x)=x$.
\end{ex}

\begin{ex}\label{ex-20.2.2} Выясним, будет ли ряд
$$
  \sum_{n=0}^\infty \frac{x}{(nx-x+1)(nx+1)}
$$
сходиться равномерно на множестве $E=(1;+\infty)$. Для этого сразу заметим, что
слагаемые раскладываются на простейшие дроби:
$$
\frac{x}{(nx-x+1)(nx+1)}=\frac{1}{nx-x+1}-\frac{1}{nx+1}
$$
Вычислим частичные суммы:
\begin{multline*}
S_N(x)=\sum_{n=0}^N \l \frac{1}{nx-x+1}-\frac{1}{nx+1}\r=\\=
\underbrace{\frac{1}{-x+1}-\frac{1}{1}
 \put(-6,26){\put(2,0){\line(1,0){18}}
 \put(0,-7){$\downarrow$}\put(18,-7){$\downarrow$}}
}_{n=0}+ \underbrace{\frac{1}{1}-\frac{1}{x+1}
 \put(-15,26){\put(2,0){\line(1,0){36}}
 \put(0,-7){$\downarrow$}\put(36,-7){$\downarrow$}}
}_{n=1}+ \underbrace{\frac{1}{x+1}-\frac{1}{2x+1}
 \put(-18,26){\put(2,0){\line(1,0){25}}
 \put(0,-7){$\downarrow$}}
}_{n=2}+\\+
 \put(-20,26){\put(2,0){\line(1,0){18}}\put(18,-7){$\downarrow$}}
 ... \put(-1,26){\put(2,0){\line(1,0){39}}
 \put(0,-7){$\downarrow$}\put(39,-7){$\downarrow$}}
 + \underbrace{\frac{1}{(N-1)x+1}-\frac{1}{Nx+1}}_{n=N} =\\=
\frac{1}{1-x}-\frac{1}{Nx+1}\overset{x\in (0;1)}{\underset{N\to
\infty}{\longrightarrow}}\frac{1}{1-x}=S(x)
\end{multline*}
Теперь найдем норму
 \begin{multline*}
||S_N(x)-S(x)||_{x\in (1;+\infty)}= \sup_{x\in (1;+\infty)}\left|
\frac{1}{Nx+1}\right|=\\=\frac{1}{N+1}\underset{N\to \infty}{\longrightarrow} 0
 \end{multline*}
Вывод: ряд $\sum\limits_{n=1}^\infty \frac{x}{(nx-x+1)(nx+1)}$ сходится на
множестве $(1;+\infty)$ равномерно к сумме $S(x)=\frac{1}{1-x}$.
\end{ex}

\begin{ers} Выясните, будет ли данный ряд сходиться равномерно
на множестве $E$:
 \bit{
 \item[1)] $\sum\limits_{n=1}^\infty (1-x)\cdot x^n$, $E=\l
0;\frac{1}{2}\r$;

 \item[2)] $\sum\limits_{n=1}^\infty \frac{x}{(nx-x+1)(nx+1)}$, $E=(0;1)$;

 \item[3)] $\sum\limits_{n=1}^\infty (1-x^n)\cdot x^n$, $E=(0;1)$,
 $E=\l 0;\frac{1}{2}\r$ (представить как сумму двух геометрических
прогрессий);

 \item[4)] $\sum\limits_{n=1}^\infty \frac{x}{\sqrt{nx+x}+\sqrt{nx}}=
\sum\limits_{n=1}^\infty \l \sqrt{nx+x}-\sqrt{nx}\r$, $E=(0;1)$;

 \item[5)] $\sum\limits_{n=1}^\infty \l
\sqrt{nx+2x}-2\sqrt{nx+x}+\sqrt{nx}\r= \sum\limits_{n=1}^\infty \l
\frac{x}{\sqrt{nx+2x}+\sqrt{nx+x}}-\frac{x}{\sqrt{nx+x}+\sqrt{nx}}\r$,
$E=(0;+\infty),$ $E=(1;+\infty)$;

 \item[6)] $\sum\limits_{n=1}^\infty \frac{x^n+(x+1)^n}{x^n\cdot (x+1)^n}=
\sum\limits_{n=1}^\infty \l \frac{1}{(x+1)^n}+\frac{1}{x^n}\r$, $E=(0;1)$,
$E=(1;+\infty)$.
 }\eit
 \end{ers}

\end{multicols}\noindent\rule[10pt]{160mm}{0.1pt}

\paragraph{Свойства равномерно сходящихся функциональных рядов.}

Перечислим некоторые свойства равномерно сходящихся рядов.

\bigskip

{\it
 \bit{
\item[$1^0$.] Если функции $a_n$ непрерывны на множестве $E$ и ряд
$$
  \sum_{n=1}^\infty a_n(x)
$$
равномерно сходится на множестве $E$ то его сумма
$$
S(x)=\sum_{n=1}^\infty a_n(x)
$$
непрерывна на множестве $E$.

\item[$2^0$.] Если функции $a_n$ непрерывны на интервале $(a,b)$ и ряд
$$
  \sum_{n=1}^\infty a_n(x)
$$
равномерно сходится на $(a,b)$, то для всякой точки $c\in (a,b)$ справедливо
равенство
 \beq
\lim_{x\to c}\lll \sum_{n=1}^\infty a_n(x) \rrr =\sum_{n=1}^\infty \lll
\lim_{x\to c} a_n (x) \rrr \label{20.3.1}
 \eeq

\item[$3^0$.]\label{PROP:int-func-ryada} Если функции $a_n$ непрерывны на
отрезке $[a,b]$ и ряд
$$
  \sum_{n=1}^\infty a_n(x)
$$
равномерно сходится на $[a,b]$, то
 \beq
\int_a^b \lll \sum_{n=1}^\infty a_n(x) \rrr \, \d x =\sum_{n=1}^\infty \lll
\int_a^b a_n (x) \, \d x \rrr \label{20.3.2}
 \eeq

\item[$4^0$.]\label{PROP:diff-func-ryada} Пусть $a_n$ -- гладкие функции на
отрезке $[a,b]$, ряд из производных
$$
  \sum_{n=1}^\infty a_n'(x)
$$
равномерно сходится на отрезке $[a,b]$, и хотя бы для одной точки $c\in [a,b]$
сходится числовой ряд
 \beq\label{20.3.3}
  \sum_{n=1}^\infty a_n(c)
 \eeq
Тогда функциональный ряд
$$
  \sum_{n=1}^\infty a_n(x)
$$
равномерно сходится на отрезке $[a,b]$, причем его сумма будет гладкой на
$[a,b]$, и ее производную можно вычислить почленным дифференцированием:
 \beq\label{20.3.4}
\lll \sum_{n=1}^\infty a_n(x) \rrr' =\sum_{n=1}^\infty a_n'(x)
 \eeq
 }\eit
}\bigskip

\begin{proof}

1. Пусть дан ряд $\sum\limits_{n=1}^\infty a_n(x)$, коэффициенты которого
$a_n(x)$ -- непрерывные функции на множестве $E$. Тогда частичные суммы
$$
S_N(x)=\sum\limits_{n=1}^N a_n(x)
$$
будут непрерывны на множестве $E$ (по теореме \ref{cont-alg}), как конечные
суммы непрерывных функций). Если ряд сходится равномерно на $E$, то есть
частичные суммы стремятся равномерно на $E$ к некоторой функции $S(x)$
$$
S_N(x)\overset{x\in E}{\underset{N\to \infty}{\rightrightarrows}} S(x),
$$
то по свойству $1^0 \, \S 4$ главы 19, эта функция
$$
S(x)=\sum_{n=1}^\infty a_n(x)
$$
должна быть непрерывна на множестве $E$.

2. Пусть функции $a_n(x)$ непрерывны и ряд
$$
  \sum_{n=1}^\infty a_n(x)
$$
равномерно сходятся на интервале $(a,b)$. Это значит, что его частичные суммы
$$
S_N(x)=\sum\limits_{n=1}^N a_n(x)
$$
непрерывны на $(a,b)$ и стремятся равномерно на $(a,b)$ к некоторой функции
$S(x)$
$$
S_N(x)\overset{x\in (a,b)}{\underset{N\to \infty}{\rightrightarrows}} S(x),
$$
По свойству $2^0 \, \S 4$ главы 19, отсюда следует, что для всякой точки $c\in
(a,b)$ справедливо равенство
$$
\lim_{x\to c}\lim_{N\to \infty} S_N(x) =\lim_{N\to \infty}\lim_{x\to c} S_N (x)
$$
то есть равенство
$$
\lim_{x\to c}\lll \sum_{n=1}^\infty a_n(x) \rrr =\sum_{n=1}^\infty \lll
\lim_{x\to c} a_n (x) \rrr
$$

3. Пусть функции $a_n(x)$ непрерывны и ряд
$$
  \sum_{n=1}^\infty a_n(x)
$$
равномерно сходятся на отрезке $[a,b]$.  Это означает, что частичные суммы ряда
$S_N(x)=\sum\limits_{n=1}^N a_n(x)$ непрерывны на $[a,b]$ и стремятся
равномерно на $[a,b]$ к некоторой функции $S(x)$
$$
S_N(x)\overset{x\in (a,b)}{\underset{N\to \infty}{\rightrightarrows}} S(x),
$$
По свойству $3^0 \, \S 4$ главы 19, отсюда следует, что
$$
\int_a^b \lim_{N\to \infty} S_N(x) \, \d x =\lim_{N\to \infty}\int_a^b  S_N (x)
\, \d x
$$
то есть,
$$
\int_a^b \lll \sum_{n=1}^\infty a_n(x) \rrr \, \d x =\sum_{n=1}^\infty \lll
\int_a^b a_n (x) \, \d x \rrr
$$

4. Пусть $a_n(x)$ -- гладкие функции на $[a,b]$, ряд из производных
$$
  \sum_{n=1}^\infty a_n'(x)
$$
равномерно сходятся на отрезке $[a,b]$, и для какой-нибудь точки $c\in [a,b]$
сходится ряд
$$
\sum_{n=1}^\infty a_n(c)
$$
Тогда частичные суммы $S_N(x)=\sum\limits_{n=1}^N a_n(x)$ -- гладкие функции
на $[a,b]$, последовательность их производных равномерно сходится на $[a,b]$ к
некоторой функции $g(x)$
$$
S_N'(x)\overset{x\in [a,b]}{\underset{N\to \infty}{\rightrightarrows}} g(x),
$$
и для точки $c\in [a,b]$ существует конечный предел
$$
\lim_{N\to\infty} S_N(c)
$$
По свойству $4^0 \, \S 4$ главы 19, отсюда следует, что функции $S_N(x)$
равномерно на $[a,b]$ сходится (к некоторой гладкой функции $S(x)$),
$$
S_N(x)\overset{x\in [a,b]}{\underset{N\to \infty}{\rightrightarrows}} S(x),
$$
причем
$$
\lll \lim_{N\to \infty} S_N(x) \rrr' =\lim_{N\to \infty} S_N' (x)
$$
Это означает, что ряд
$$
  \sum_{n=1}^\infty a_n(x)
$$
равномерно сходятся на отрезке $[a,b]$, и
$$
\lll \sum_{n=1}^\infty a_n(x) \rrr' =\sum_{n=1}^\infty a_n'(x) \quad $$
\end{proof}

\paragraph{Критерий Коши равномерной сходимости функционального ряда.}

\begin{tm}[\bf критерий Коши равномерной сходимости ряда]\label{tm-20.4.1}
Пусть $\{ a_n \}$ -- произвольная последовательность функций, определенных на
множестве $x\in E$. Тогда следующие условия эквивалентны:
 \bit{
\item[(i)] функциональный ряд
$$
\sum_{n=1}^\infty a_n
$$
сходится равномерно на множестве $E$; \item[(ii)] для любой последовательности
$l_i\in \mathbb{N}$, и любой бесконечно большой последовательности $k_i\in
\mathbb{N}$
$$
k_i\underset{i\to \infty}{\longrightarrow}\infty,
$$
сумма $\sum\limits_{n=k_i+1}^{k_i+l_i} a_n(x)$ стремится к нулю при $i\to
\infty$ равномерно на множестве $E$:
$$
\norm{\sum\limits_{n=k_i+1}^{k_i+l_i} a_n}_E \underset{i\to
\infty}{\longrightarrow} 0
$$
 }\eit
\end{tm}\begin{proof} Для частичных суммы
$S_N(x)=\sum_{n=1}^N a_n(x)$ мы получаем следующую логическую цепочку:
$$
\text{ряд $\sum\limits_{n=1}^\infty a_n(x)$ сходится равномерно на $E$}
$$
$$
\Updownarrow
$$
$$
\text{последовательность $S_N(x)$ сходится равномерно на $E$}
$$
$$
\Updownarrow\put(20,0){\smsize\text{$\begin{pmatrix}\text{вспоминаем
критерий Коши}\\ \text{сходимости функциональной}\\
\text{последовательности -- теорему \ref{tm-19.5.1}}\end{pmatrix}$}}
$$
$$
\begin{array}{c}\text{для любой последовательности $l_i\in \mathbb{N}$,}\\
\text{и любой бесконечно большой последовательности $k_i\in \mathbb{N}$}\quad
k_i\underset{i\to \infty}{\longrightarrow}\infty,\\
\text{выполняется соотношение:}\quad
\parallel S_{k_i+l_i}(x)-S_{k_i}(x)\parallel_{x\in E}\underset{i\to \infty}{\longrightarrow} 0
\end{array}
$$
$$
\Updownarrow\put(20,0){\smsize\text{$\begin{pmatrix}
 \text{замечаем, что}\\
S_{k_i+l_i}(x)-S_{k_i}(x)=\sum\limits_{n=k_i+1}^{k_i+l_i} a_n(x)
\end{pmatrix}$}}
$$
$$
\begin{array}{c}\text{для любой последовательности $l_i\in \mathbb{N}$,}\\
\text{и любой бесконечно большой последовательности $k_i\in \mathbb{N}$}\quad
k_i\underset{i\to \infty}{\longrightarrow}\infty,\\
\text{выполняется соотношение:}\quad \left|\left|
\sum\limits_{n=k_i+1}^{k_i+l_i} a_n(x) \right|\right|_{x\in E}\underset{i\to
\infty}{\longrightarrow} 0
\end{array}
$$
\end{proof}

\paragraph{Необходимое условие равномерной сходимости ряда.}

Часто бывает, что частичные суммы функционального ряда не поддаются вычислению,
как, например, в случае с рядом
 \beq\label{20.5.1}
  \sum_{n=1}^\infty \sin \frac{x}{n^2}
 \eeq
Чтобы понять, сходится этот ряд равномерно, или нет, используются теоремы,
называемые признаками равномерной сходимости. Первый из них звучит так:

\begin{tm}\label{tm-20.5.1}
Если функциональный ряд
$$
  \sum_{n=1}^\infty a_n(x)
$$
равномерно сходится на множестве $E$, то его общий член равномерно стремится к
нулю на $E$:
$$
a_n(x)\overset{x\in E}{\underset{n\to \infty}{\rightrightarrows}} 0
$$
\end{tm}\begin{proof} По критерию Коши (теорема \ref{tm-20.4.1}), если
этот ряд сходится равномерно на $E$, то для последовательностей $l_i=1$, и
$k_i=i$ должно выполняться условие
$$
\left|\left| \sum\limits_{n=k_i+1}^{k_i+l_i} a_n(x) \right|\right|_{x\in E}=
\left|\left| \sum\limits_{n=i+1}^{i+1} a_n(x) \right|\right|_{x\in E}=
\left|\left| a_{i+1}(x) \right|\right|_{x\in E}\underset{i\to
\infty}{\longrightarrow} 0 \quad $$ \end{proof}

\noindent\rule{160mm}{0.1pt}\begin{multicols}{2}

\begin{ex}\label{ex-20.5.2} Рассмотрим ряд \eqref{20.5.1}:
$$
  \sum_{n=1}^\infty \sin \frac{x}{n^2}
$$
По признаку абсолютной сходимости (теорема \ref{tm-18.5.6}), он сходится при
любом $x\in\R$:
$$
\sum_{n=1}^\infty \left|\sin \frac{x}{n^2}\right|\le\eqref{|sin-x|<|x|}\le
\sum_{n=1}^\infty \left|\frac{x}{n^2}\right|= |x|\cdot \sum_{n=1}^\infty
\frac{1}{n^2}<\infty
$$
Поймем, будет ли он сходиться равномерно на множестве $E=\R$. Для этого
вычислим норму общего члена:
$$
\norm{\sin\frac{x}{n^2}}_{x\in \R}= \sup_{x\in \R}\ml \sin
\frac{x}{n^2}\mr=1\underset{n\to \infty}{\not\longrightarrow} 0
$$
Вывод: ряд сходится на $\R$, но не равномерно.
\end{ex}

\ber Проверьте, будет ли ряд сходится равномерно на множестве $\R$:
$$
\sum\limits_{n=1}^\infty \arctg\frac{x}{n^2}
$$
Здесь предварительно нужно доказать формулу, аналогичную \eqref{|sin-x|<|x|}:
 \beq\label{|arctg(x)|<|x|}
|\arctg x|\le |x|,\qquad x\in\R
 \eeq
Это можно сделать, например, с помощью формулы Ньютона-Лейбница: при $x>0$ мы
получаем
 \begin{multline*}
|\kern-20pt\overbrace{\arctg x}^{\scriptsize\begin{matrix}0\\
\phantom{\eqref{arctg(x)>0}}\ \text{\rotatebox{90}{$>$}}\ \eqref{arctg(x)>0}
\end{matrix}}\kern-20pt|=\arctg x=\int_0^x\frac{\d t}{1+t^2}\le \int_0^x 1\ \d t=\\=x=|x|
 \end{multline*}
и отсюда уже для $x<0$ получаем
$$
\kern-20pt\overbrace{\arctg x}^{\scriptsize\begin{matrix}0\\
\phantom{\eqref{arctg(x)<0}}\ \text{\rotatebox{90}{$<$}}\ \eqref{arctg(x)<0}
\end{matrix}}\kern-20pt=-\arctg x=\arctg(\overbrace{-x}^{\scriptsize\begin{matrix}0\\
\text{\rotatebox{90}{$>$}}
\end{matrix}})\le -x=|x|
$$
 \eer

\begin{ers} Проверьте, будут ли данные ряды сходится равномерно
на множестве $E$:
 \biter{
 \item[1.] $\sum\limits_{n=1}^\infty x^n$, $E=(-1;1)$;

 \item[2.] $\sum\limits_{n=1}^\infty x^n\cdot (1-x^n)$, $E=[0;1]$
(здесь необходимо исследование функции $x^n\cdot (1-x^n)$ на экстремум на
множестве $[0;1]$).
 }\eiter
 \end{ers}

\end{multicols}\noindent\rule[10pt]{160mm}{0.1pt}

\paragraph{Признак Вейерштрасса равномерной сходимости.}

\begin{tm}\label{tm-20.5.4}
Если числовой ряд, состоящий из норм функций $a_n(x)$ на множестве $E$
 \beq
  \sum_{n=1}^\infty \left|\left| a_n(x)\right|\right|_{x\in E}\label{20.5.2}
 \eeq
сходится, то функциональный ряд
$$
  \sum_{n=1}^\infty a_n(x)
$$
равномерно сходится на множестве $E$.
\end{tm}\begin{proof} Если ряд \eqref{20.5.2} сходится, то, по критерию
Коши сходимости числового ряда (теорема \ref{tm-18.4.1}), для любой
последовательности $l_i\in \mathbb{N}$, и любой бесконечно большой
последовательности $k_i\in \mathbb{N}$ $k_i\underset{i\to
\infty}{\longrightarrow}\infty$,
$$
\sum\limits_{n=k_i+1}^{k_i+l_i}\left|\left| a_n(x)\right|\right|_{x\in
E}\underset{i\to \infty}{\longrightarrow} 0
$$
Отсюда
$$
0\le \left|\left| \sum\limits_{n=k_i+1}^{k_i+l_i} a_n(x) \right|\right|_{x\in
E}\le \sum\limits_{n=k_i+1}^{k_i+l_i}\left|\left| a_n(x)\right|\right|_{x\in
E}\underset{i\to \infty}{\longrightarrow} 0
$$
$$
\Downarrow
$$
$$
\left|\left| \sum\limits_{n=k_i+1}^{k_i+l_i} a_n(x) \right|\right|_{x\in
E}\underset{i\to \infty}{\longrightarrow} 0
$$
Поскольку это верно для любых последовательностей $l_i\in \mathbb{N}$ и $k_i\in
\mathbb{N}$ ($k_i\underset{i\to \infty}{\longrightarrow}\infty$), по критерию
Коши равномерной сходимости функционального ряда (теорема \ref{tm-20.4.1}), ряд
$\sum\limits_{n=1}^\infty a_n(x)$ равномерно сходится на множестве $E$.
\end{proof}

\noindent\rule{160mm}{0.1pt}\begin{multicols}{2}

\begin{ex}\label{ex-20.5.5} Проверим, будет ли ряд \eqref{20.5.1} сходиться равномерно на
множестве $[0;1]$:
$$
  \sum_{n=1}^\infty \sin \frac{x}{n^2}
$$
Для этого вычислим норму общего члена:
$$
\norm{\sin\frac{x}{n^2}}_{x\in [0;1]}= \sup_{x\in [0;1]}\ml
\sin\frac{x}{n^2}\mr=\sin\frac{1}{n^2}
$$
Теперь получаем
 \begin{multline*}
\sum_{n=1}^\infty \norm{\sin\frac{x}{n^2}}_{x\in [0;1]}= \sum_{n=1}^\infty
\sin\frac{1}{n^2}\le\\ \le \eqref{|sin-x|<|x|}\le \sum_{n=1}^\infty
\frac{1}{n^2}<\infty
 \end{multline*}
Вывод: ряд сходится равномерно на отрезке $[0;1]$.
\end{ex}

\begin{ers} Проверьте, будут ли данные ряды сходится равномерно
на множестве $E$:
 \biter{
 \item[1.] $\sum\limits_{n=1}^\infty \frac{1}{x^2+n^2}, \quad
E=\R$;

 \item[2.] $\sum\limits_{n=1}^\infty x^n, \quad
E=\left[-\frac{1}{2};\frac{1}{2}\right]$;

 \item[3.] $\sum\limits_{n=1}^\infty \frac{x}{1+n^4\cdot x^2},$
$E=\R$ (здесь необходимо исследование функции $\frac{x}{1+n^4\cdot x^2}$ на
экстремум);

 \item[4.] $\sum\limits_{n=1}^\infty \frac{n\cdot x}{1+n^5\cdot x}, \quad
E=\R$;

 \item[5.] $\sum\limits_{n=1}^\infty \l
\frac{x^n}{n}-\frac{x^{n+1}}{n+1}\r, \quad E=\R$;

 \item[6.] $\sum\limits_{n=1}^\infty \frac{n^2}{\sqrt{n!}}\l x^n+x^{-n}\r, \quad E=\left[ \frac{1}{2}, 2 \right]$.
 }\eiter
\end{ers}

\end{multicols}\noindent\rule[10pt]{160mm}{0.1pt}

\paragraph{Признак Лейбница равномерной сходимости.}

\begin{tm}\label{TH:priznak-Leibnitza-ravnom-shod}
Пусть функциональная последовательность $\{b_n(x);\ x\in E\}$ обладает
следующими свойствами:
 \bit{
\item[(i)] она неотрицательна и невозрастает,
 \beq\label{18.6.1}
b_0(x)\ge b_1(x)\ge b_2(x)\ge ... \ge 0
 \eeq

\item[(ii)] она стремится к нулю равномерно на $E$:
 \beq
||b_n||_E\underset{n\to \infty}{\longrightarrow} 0
 \eeq
 }\eit\noindent
Тогда
 \bit{
\item[(a)] функциональный ряд $\sum\limits_{n=0}^\infty (-1)^n b_n(x)$ сходится
равномерно на множестве $E$,

\item[(b)] его остаток оценивается сверху нормой своего первого члена:
 \beq\label{otsenka-ostatka-ravnom-Leibnitz}
\left|\left|\sum\limits_{n=N+1}^\infty (-1)^n b_n(x)\right|\right|_{x\in E}\le
||b_{N+1}||_E
 \eeq
 }\eit
 \end{tm}

\bpr Из признака Лейбница сходимости числового ряда (теорема \ref{tm-18.6.1})
следует, что при каждом фиксированном $x\in E$ числовой ряд
$\sum\limits_{n=0}^\infty (-1)^n b_n(x)$ сходится. Обозначим через $S(x)$ его
сумму, а через $S_N(x)$ его частичные суммы:
$$
S(x)=\sum\limits_{n=0}^\infty (-1)^n b_n(x),\qquad S_N(x)=\sum\limits_{n=0}^N
(-1)^n b_n(x)
$$
В силу \eqref{otsenka-ostatka-Leibnitz}, мы получаем:
$$
\forall x\in E\qquad |S(x)-S_N(x)|=\left|\sum\limits_{n=N+1}^\infty (-1)^n
b_n(x)\right|\le b_{N+1}(x)\le ||b_{N+1}||_E
$$
$$
\Downarrow
$$
$$
||S(x)-S_N(x)||_{x\in E}=\sup_{x\in E}|S(x)-S_N(x)|= \sup_{x\in
E}\left|\sum\limits_{n=N+1}^\infty (-1)^n b_n(x)\right|\le ||b_{N+1}||_E
$$
То есть справедлива оценка \eqref{otsenka-ostatka-ravnom-Leibnitz}. Из нее
сразу следует равномерная сходимость ряда:
$$
||S(x)-S_N(x)||_{x\in E}\le
||b_{N+1}||_E\underset{N\to\infty}{\longrightarrow}0
$$
$$
\Downarrow
$$
$$
||S(x)-S_N(x)||_{x\in E}\underset{N\to\infty}{\longrightarrow}0
$$
$$
\Downarrow
$$
$$
S_N(x)\overset{x\in E}{\underset{N\to\infty}{\rightrightarrows}} S(x)
$$

\epr

\paragraph{Признаки Дирихле и Абеля равномерной сходимости.}

\begin{tm}[\bf признак Дирихле]\label{tm-20.5.7}\footnote{Теорема \ref{tm-20.5.7} часто называется также признаком Харди.}
Пусть функциональные последовательности $\{ a_n \}$ и $\{ b_n \}$ обладают
следующими свойствами:
 \bit{
\item[$(i)$] частичные суммы ряда $\sum\limits_{n=1}^\infty a_n$ равномерно
ограничены на множестве $E$:
$$
\sup_{N\in \mathbb{N}}\norm{\sum\limits_{n=1}^N a_n(x)}_{x\in E} <\infty
$$
\item[$(ii)$] последовательность $\{ b_n (x)\}$ монотонна при любом $x\in E$;

\item[$(iii)$] последовательность $\{ b_n \}$ равномерно стремится к нулю на
множестве $E$:
 \beq\label{tm-20.5.7-0}
\norm{b_n(x)}_{x\in E}\underset{n\to \infty}{\longrightarrow} 0
 \eeq
 }\eit
Тогда ряд $\sum\limits_{n=1}^\infty a_n(x)\cdot b_n(x)$ равномерно сходится на
множестве $E$.
 \end{tm}
 \begin{proof} Здесь с очевидными изменениями копируется доказательство
теоремы \ref{tm-18.6.5}. Обозначим
$$
C=\sup_{N\in \mathbb{N}}\norm{\sum_{n=1}^N a_n(x)}_{x\in E}<\infty
$$
и заметим, что
 \beq\label{tm-20.5.7-1}
\sup_{p,q\in \mathbb{N}}\norm{\sum_{n=p}^q a_n(x)}_{x\in E}\le 2C <\infty
 \eeq
Действительно,
$$
\norm{\sum_{n=p}^q a_n(x)}_{x\in E}=\norm{\sum_{n=1}^q a_n(x)-\sum_{n=1}^{p-1}
a_n(x)}_{x\in E}\le\underbrace{\norm{\sum_{n=1}^q
a_n(x)}_{x\in E}}_{\scriptsize\begin{matrix}\text{\rotatebox{90}{$\ge$}}\\
C\end{matrix}}+\underbrace{\norm{\sum_{n=1}^{p-1} a_n(x)}_{x\in E}}_{\scriptsize\begin{matrix}\text{\rotatebox{90}{$\ge$}}\\
C\end{matrix}}\le 2C
$$
Воспользуемся теперь критерием Коши (теорема \ref{tm-20.4.1}): выберем
произвольные последовательности  $l_i\in \mathbb{N}$ и $k_i\in \mathbb{N}$
($k_i\underset{n\to \infty}{\longrightarrow}\infty$). По неравенству Абеля
\eqref{18.6.7} получаем:
$$
\ml \sum\limits_{n=k_i+1}^{k_i+l_i} a_n(x) \cdot b_n(x) \mr \le 3\cdot \max_{
k_i\le N\le k_i+l_i}\underbrace{\ml \sum\limits_{n=k_i+1}^N a_n(x)
\mr}_{\scriptsize\begin{matrix}
\text{\rotatebox{90}{$\ge$}}\\
\norm{\sum\limits_{n=k_i+1}^N a_n(x)}_{x\in E}\\
\phantom{\tiny\eqref{tm-20.5.7-1}}\
\text{\rotatebox{90}{$\ge$}}\ {\tiny\eqref{tm-20.5.7-1}}\\
2C
\end{matrix}}\cdot\max\Big\{\underbrace{|b_{k_i+1}(x)|}_{\scriptsize\begin{matrix}
\text{\rotatebox{90}{$\ge$}}\\
\norm{b_{k_i+1}(x)}_{x\in E}
\end{matrix}},\underbrace{|b_{k_i+l_i}(x)|}_{\scriptsize\begin{matrix}
\text{\rotatebox{90}{$\ge$}}\\
\norm{b_{k_i+l_i}(x)}_{x\in E}
\end{matrix}}\Big\}
$$
$$
\Downarrow
$$
$$
\norm{\sum\limits_{n=k_i+1}^{k_i+l_i} a_n(x) \cdot b_n(x)}_{x\in E} \le 6\cdot
C\cdot\max\Big\{\underbrace{\norm{b_{k_i+1}(x)}_{x\in E}}_{\scriptsize
\begin{matrix}\phantom{\tiny\eqref{tm-20.5.7-0}}\ \downarrow\ {\tiny\eqref{tm-20.5.7-0}} \\ 0\end{matrix}},
\underbrace{\norm{b_{k_i+l_i}(x)}_{x\in E}}_{\scriptsize
\begin{matrix}\phantom{\tiny\eqref{tm-20.5.7-0}}\ \downarrow\ {\tiny\eqref{tm-20.5.7-0}} \\ 0\end{matrix}}
\Big\}\underset{i\to\infty}{\longrightarrow}0
$$
Это верно для любых последовательностей $l_i\in \mathbb{N}$ и $k_i\in
\mathbb{N}$ ($k_i\underset{n\to \infty}{\longrightarrow}\infty$), значит ряд
$\sum\limits_{n=1}^\infty a_n\cdot b_n$ сходится.

\end{proof}

\noindent\rule{160mm}{0.1pt}\begin{multicols}{2}

\begin{ex}\label{ex-20.5.8} Проверим, будет ли ряд
$$
  \sum_{n=1}^\infty \frac{x^n}{\sqrt{n}}
$$
сходиться равномерно на множестве $(-1;0)$. Обозначим
$$
  a_n(x)=x^n, \quad b_n(x)=\frac{1}{\sqrt{n}}
$$
Тогда
 \begin{multline*}
\norm{ \sum_{n=1}^N a_n(x) }_{x\in (-1;0)}= \sup_{x\in (-1;0)}\ml \sum_{n=1}^N
a_n(x) \mr=\\= \sup_{x\in (-1;0)}\ml \sum_{n=1}^N x^n  \mr= \sup_{x\in
(-1;0)}\ml x\cdot \frac{1-x^N}{1-x}
\mr\le\\ \le \sup_{x\in (-1;0)}\frac{|x|\cdot (1+|x^N|)}{|1-x|}\le\\
\le \frac{\sup_{x\in (-1;0)} |x|\cdot (1+|x^N|)}{\inf_{x\in (-1;0)} |1-x|}\le
\frac{2}{1}=2
\end{multline*}
откуда
$$
\sup_{N\in \mathbb{N}}\norm{\sum_{n=1}^N a_n(x)}_{x\in (-1;0)}\le 2<\infty
$$
А, с другой стороны,
$$
b_n(x)=\frac{1}{\sqrt{n}}\overset{\smsize \begin{matrix}\text{монотонно}\\
\text{и равномерно}\\ \text{по $x\in (-1,0)$}\end{matrix}}
{\underset{n\to\infty}{\rightrightarrows}} 0
$$
Вывод: ряд сходится равномерно на множестве $(-1;0)$.
\end{ex}

\begin{ers} Проверьте, будут ли данные ряды сходится равномерно
на множестве $E$:
 \biter{
 \item[1.] $\sum\limits_{n=1}^\infty \frac{(-1)^\frac{n(n-1)}{2}}
{\sqrt[3]{n^2+e^x}}$, $E=\R$;

 \item[2.] $\sum\limits_{n=1}^\infty \frac{\cos \frac{2\pi
n}{3}}{\sqrt{n^2+x^2}}$, $E=\R$ (здесь нужно применить формулы \eqref{18.6.5});

 \item[3.] $\sum\limits_{n=1}^\infty \frac{\sin x\cdot \sin
nx}{\sqrt{n+1}}$, $E=\R$ (тоже нужны формулы \eqref{18.6.5}).
 }\eiter
 \end{ers}

\end{multicols}\noindent\rule[10pt]{160mm}{0.1pt}

\begin{tm}[\bf признак Абеля]\label{TH:priznak-Abelya-ravnom-shod-ryadov}
Пусть функциональные последовательности  $\{ a_n \}$ и $\{ b_n \}$ на множестве
$E$ обладают следующими свойствами:
 \bit{
\item[$(i)$] ряд $\sum\limits_{n=1}^\infty a_n$ сходится равномерно на $E$;

\item[$(ii)$] последовательность $\{ b_n \}$ монотонна и равномерно ограничена
на $E$:
 $$
\sup_{n\in\N}\norm{b_n(x)}_{x\in E}<\infty
 $$
 }\eit
Тогда ряд $\sum\limits_{n=1}^\infty a_n\cdot b_n$ сходится равномерно на $E$.
\end{tm}

\begin{proof} Здесь повторяются рассуждения, применявшиеся при доказательстве
теоремы \ref{tm-18.6.9}. Сначала нужно заметить, что из равномерной сходимости
ряда $\sum\limits_{n=1}^\infty a_n$ следует, что для произвольных
последовательностей $l_i\in \mathbb{N}$ и $k_i\in \mathbb{N}$
($k_i\underset{n\to \infty}{\longrightarrow}\infty$) выполняется соотношение
 \beq\label{PROOF:Abel-sum-ravn-1}
\max_{ k_i\le N\le k_i+l_i}\norm{\sum\limits_{n=k_i+1}^N a_n
}_{E}\underset{i\to\infty}{\longrightarrow}0
 \eeq
После этого, обозначив
 $$
B=\sup_{n\in\N}\norm{b_n}_E=\sup_{n\in\N}\sup_{x\in E}|b_n(x)|<\infty
 $$
мы пользуемся критерием Коши (теоремой \ref{tm-20.4.1}): для произвольных
последовательностей  $l_i\in \mathbb{N}$ и $k_i\in \mathbb{N}$
($k_i\underset{n\to \infty}{\longrightarrow}\infty$) по неравенству Абеля
\eqref{18.6.7} мы получим
$$
\ml \sum\limits_{n=k_i+1}^{k_i+l_i} a_n(x) \cdot b_n(x) \mr \le 3\cdot \max_{
k_i\le N\le k_i+l_i}\underbrace{\ml \sum\limits_{n=k_i+1}^N a_n(x)
\mr}_{\scriptsize\begin{matrix} \text{\rotatebox{90}{$\ge$}} \\
\norm{\sum\limits_{n=k_i+1}^N a_n
}_E\end{matrix}}\cdot\max\Big\{\underbrace{|b_{k_i+1}(x)|}_{\scriptsize
\begin{matrix}\text{\rotatebox{90}{$\ge$}} \\ B\end{matrix}},\underbrace{|b_{k_i+l_i}(x)|}_{\scriptsize
\begin{matrix}\text{\rotatebox{90}{$\ge$}} \\ B\end{matrix}}\Big\}
$$
$$
\Downarrow
$$
$$
\ml \sum\limits_{n=k_i+1}^{k_i+l_i} a_n \cdot b_n \mr \le 3 B\cdot
\underbrace{\max_{ k_i\le N\le k_i+l_i}\ml \sum\limits_{n=k_i+1}^N a_n
\mr}_{\scriptsize\begin{matrix}
\phantom{\tiny\eqref{PROOF:Abel-sum-ravn-1}}\ \downarrow\ {\tiny\eqref{PROOF:Abel-sum-ravn-1}} \\
0\end{matrix}}\underset{i\to\infty}{\longrightarrow}0
$$
Это верно для любых последовательностей $l_i\in \mathbb{N}$ и $k_i\in
\mathbb{N}$ ($k_i\underset{n\to \infty}{\longrightarrow}\infty$), значит  ряд
$\sum\limits_{n=1}^\infty a_n\cdot b_n$ сходится равномерно.
\end{proof}

\subsection{Всюду непрерывная, но нигде не дифференцируемая функция}\label{nepr-no-nedifferent-func}

\noindent\rule{160mm}{0.1pt}\begin{multicols}{2}

\bex {\it Пусть числа $a$ и $b$ удовлетворяют условиям:
$$
0<b<1,\quad a\in 2\N-1, \quad ab>1+\frac{3}{2}\pi
$$
Тогда функция
 \beq\label{primer-Weierstrassa}
f(x)=\sum_{n=0}^\infty b^n\cos(a^n\pi x)
 \eeq
определена и непрерывна всюду на $\R$, но не дифференцируема ни в одной точке.}
\eex

 \bpr
Равномерные нормы слагаемых этого ряда образуют сходящийся ряд
$$
\sum_{n=0}^\infty b^n\norm{\cos(a^n\pi x)}_{x\in\R}=\sum_{n=0}^\infty
b^n<\infty
$$
поэтому по признаку Вейерштрасса (теорема \ref{tm-20.5.4}), функциональный ряд
в \eqref{primer-Weierstrassa} сходится равномерно, и значит, определяет
непрерывную функцию. Зафиксируем точку $x\in\R$ и убедимся, что $f$ не
дифференцируема в ней.

Заметим прежде всего, что
$$
ab>1+\frac{3}{2}\pi
$$
$$
\Downarrow
$$
 \beq\label{3/2>pi/(ab-1)}
\frac{3}{2}>\frac{\pi}{ab-1}
 \eeq

Далее рассмотрим три последовательности:
 \begin{align*}
& \xi_k=\left\{a^kx+\frac{1}{2}\right\}-\frac{1}{2}\in\left[-\frac{1}{2},\frac{1}{2}\right) \\
& \alpha_k=\left[a^kx+\frac{1}{2}\right]=a^kx-\xi_k\in\Z \\
&
h_k=\frac{1-\xi_k}{a^k}=\frac{\frac{3}{2}-\left\{a^kx+\frac{1}{2}\right\}}{a^k}
 \end{align*}
(здесь $[\cdot]$ и $\{\cdot\}$ -- целая и дробная части числа, определенные
формулами \eqref{tselaya-chast-chisla} и \eqref{drobnaya-chast-chisla}).
Очевидно, что
 \beq\label{0<h_k<3/2a^k}
0<h_k\le \frac{3}{2a^k}
 \eeq
$$
\Downarrow
$$
$$
h_k \underset{k\to\infty}{\longrightarrow}0,
$$
поэтому наша задача будет выполнена, если мы докажем, что
 $$
\frac{f(x+h_k)-f(x)}{h_k}\underset{k\to\infty}{\longrightarrow}\infty
 $$

Обозначим
 \begin{align*}
& S=\sum_{n=0}^{k-1} b^n\cdot\frac{\cos\big(a^n\pi (x+h_k)\big)-\cos(a^n\pi
x)}{h_k}\\
& R=\sum_{n=k}^\infty b^n\cdot\frac{\cos\big(a^n\pi (x+h_k)\big)-\cos(a^n\pi
x)}{h_k}.
 \end{align*}

Тогда
 \begin{multline*}
\frac{f(x+h_k)-f(x)}{h_k}=\\= \sum_{n=0}^\infty b^n\cdot\frac{\cos\big(a^n\pi
(x+h_k)\big)-\cos(a^n\pi x)}{h_k}=\\=\sum_{n=0}^{k-1}+\sum_{n=k}^\infty=S+R
 \end{multline*}

Сначала оценим сверху $S$:
 \begin{multline*}
\Big|\cos\big(a^n\pi (x+h_k)\big)-\cos(a^n\pi x)\Big|=\\=\Big|\cos\big(a^n\pi
x+a^n\pi h_k\big)-\cos(a^n\pi x)\Big|=\eqref{cos-x-cos-y}=\\=
2\underbrace{\left|\sin\frac{2a^n\pi x+a^n\pi
h_k}{2}\right|}_{\scriptsize\begin{matrix} \IA\\
1\end{matrix}}\cdot\underbrace{\left|\sin\frac{a^n\pi
h_k}{2}\right|}_{\scriptsize\begin{matrix} \phantom{\eqref{|sin-x|<|x|}}\ \IA \ \eqref{|sin-x|<|x|} \\
\left|\frac{a^n\pi h_k}{2}\right|\\ \text{\rotatebox{90}{$=$}} \\ \frac{a^n\pi
h_k}{2} \end{matrix}}\le a^n\pi h_k
 \end{multline*}
 $$
 \Downarrow
 $$
 \begin{multline*}
|S|=\left|\sum_{n=0}^{k-1} b^n\cdot\frac{\cos\big(a^n\pi
(x+h_k)\big)-\cos(a^n\pi x)}{h_k}\right|\le\\ \le \sum_{n=0}^{k-1}
b^n\cdot\left|\frac{\cos\big(a^n\pi (x+h_k)\big)-\cos(a^n\pi x)}{h_k}\right|\le
\\ \le \sum_{n=0}^{k-1} b^n\cdot\frac{a^n\pi h_k}{h_k}=\pi\cdot\sum_{n=0}^{k-1}
(ab)^n=\\=\eqref{Geom-progr-N}=\pi\cdot\frac{(ab)^k-1}{ab-1}\le \frac{\pi a^k
b^k}{ab-1}
 \end{multline*}

После этого оценим $R$. Пусть $n\ge k$. Тогда, во-первых,
 \begin{multline*}
a^k\cdot(x+h_k)=a^k\cdot x+a^k\cdot h_k=\\= a^k\cdot
x+\frac{3}{2}-\left\{a^kx+\frac{1}{2}\right\}=\\= a^k\cdot
x+\frac{1}{2}-\left\{a^kx+\frac{1}{2}\right\}+1=\\=
\left[a^kx+\frac{1}{2}\right]+1=\alpha_k+1
 \end{multline*}
 $$
 \Downarrow
 $$
 \begin{multline*}
a^n\cdot\pi\cdot(x+h_k)=a^{n-k}\cdot a^k (x+h_k)\cdot\pi=\\=
a^{n-k}\cdot(\alpha_k+1)\cdot\pi
 \end{multline*}
 $$
 \Downarrow
 $$
 \begin{multline*}
\cos\Big(a^n\cdot\pi\cdot(x+h_k)\Big)=
\cos\Big(\underbrace{a^{n-k}\cdot(\alpha_k+1)}_{\scriptsize\begin{matrix}
\text{\rotatebox{90}{$\owns$}}\\ \Z \end{matrix}}\cdot\pi\Big)=\\=
(-1)^{a^{n-k}\cdot(\alpha_k+1)}=
\Big((-1)^{\overbrace{a^{n-k}}^{\scriptsize\begin{matrix} 2\N-1 \\
\text{\rotatebox{90}{$\in$}}\end{matrix}}}\Big)^{\alpha_k+1}=(-1)^{\alpha_k+1}
 \end{multline*}
Во-вторых,
 $$
a^k\cdot x=\alpha_k+\xi_k
 $$
 $$
 \Downarrow
 $$
 $$
a^n\cdot x= a^{n-k}\cdot(\alpha_k+\xi_k)
 $$
 $$
 \Downarrow
 $$
 \begin{multline*}
\cos\Big(a^n\pi x\Big)=\cos\Big(a^{n-k}\cdot\pi\cdot(\alpha_k+\xi_k)\Big)=\\=
\cos\Big(a^{n-k}\pi\alpha_k+a^{n-k}\pi\xi_k\Big)=\eqref{cos(x+y)}=\\=
\underbrace{\cos\Big(a^{n-k}\alpha_k\pi\Big)}_{\scriptsize\begin{matrix}
\text{\rotatebox{90}{$=$}}\\ (-1)^{a^{n-k}\cdot\alpha_k}
\end{matrix}}\cdot\cos\Big(a^{n-k}\pi\xi_k\Big)-\\-
\underbrace{\sin\Big(a^{n-k}\alpha_k\pi\Big)}_{\scriptsize\begin{matrix}
\text{\rotatebox{90}{$=$}}\\ 0
\end{matrix}}\cdot\sin\Big(a^{n-k}\pi\xi_k\Big)=\\=
(-1)^{a^{n-k}\cdot\alpha_k}\cdot\cos\Big(a^{n-k}\pi\xi_k\Big)=\\=
\Big((-1)^{\overbrace{a^{n-k}}^{\scriptsize\begin{matrix} 2\N-1 \\
\text{\rotatebox{90}{$\in$}}\end{matrix}}}\Big)^{\alpha_k}\cdot\cos\Big(a^{n-k}\pi\xi_k\Big)=\\=
(-1)^{\alpha_k}\cdot\cos\big(a^{n-k}\pi\xi_k\big)
 \end{multline*}

Вместе мы получаем:
 \begin{multline*}
\cos\big(a^n\pi (x+h_k)\big)-\cos(a^n\pi x)=\\=
(-1)^{\alpha_k+1}-(-1)^{\alpha_k}\cdot\cos\big(a^{n-k}\pi\xi_k\big)=\\=
(-1)^{\alpha_k+1}\Big(1+\cos\big(a^{n-k}\pi\xi_k\big)\Big)
 \end{multline*}
 $$
 \Downarrow
 $$
 \begin{multline*}
|R|=\left|(-1)^{\alpha_k}\cdot\sum_{n=k}^\infty
b^n\cdot\frac{1+\cos\big(a^{n-k}\pi\xi_k\big)}{h_k}\right|=\\=
\bigg|\sum_{n=k}^\infty
\underbrace{b^n\cdot\frac{1+\cos\big(a^{n-k}\pi\xi_k\big)}{h_k}}_{\scriptsize\begin{matrix}
\VI\\ 0
\end{matrix}}\bigg|= \\ =
\sum_{n=k}^\infty b^n\cdot\frac{1+\cos\big(a^{n-k}\pi\xi_k\big)}{h_k}\ge \\ \ge
\underbrace{b^k\cdot\frac{1+\cos\big(a^{\overbrace{k-k}^{\scriptsize\begin{matrix} 0 \\
\text{\rotatebox{90}{$=$}}\end{matrix}}}\pi\xi_k\big)}{h_k}}_{n=k}=\\=
b^k\cdot\frac{1+\cos\big(\kern-5pt\overbrace{\pi\xi_k}^{\scriptsize\begin{matrix} \left[-\frac{\pi}{2},\frac{\pi}{2}\right) \\
\text{\rotatebox{90}{$\in$}}\end{matrix}}\kern-5pt\big)}{h_k}\ge
\frac{b^k}{h_k}\eqref{0<h_k<3/2a^k}\ge\frac{2a^kb^k}{3}
 \end{multline*}

Теперь получаем:
 \begin{multline*}
\left|\frac{f(x+h_k)-f(x)}{h_k}\right|\ge |R|-|S|\ge \\ \ge
\frac{2a^kb^k}{3}-\frac{\pi a^k b^k}{ab-1}=\\=
\underbrace{\l\frac{2}{3}-\frac{\pi}{ab-1}\r}_{\scriptsize\begin{matrix}
\phantom{\eqref{3/2>pi/(ab-1)}}\ \text{\rotatebox{90}{$<$}}\
\eqref{3/2>pi/(ab-1)}\\ 0\end{matrix}}
 (ab)^k\underset{k\to\infty}{\longrightarrow}\infty
 \end{multline*}

 \epr

\end{multicols}\noindent\rule[10pt]{160mm}{0.1pt}

\section{Равномерная по производным сходимость}

\subsection{Равномерная по производным норма функции}

Пусть $E$ -- интервал или отрезок в $\R$.

 \bit{
\item[$\bullet$]  Функция $f$ называется {\it гладкой на $E$ порядка $m$}, если
она имеет $m$ непрерывных производных на $E$, то есть если существует конечная
последовательность непрерывных функций $f_0,f_1,...,f_m$ на $E$ таких, что
$$
f_0=f
$$
и для всякого $k=0,...,m-1$ функция $f_k$ дифференцируема на отрезке $E$ (если
$E$ -- отрезок, то дифференцируемость на нем понимается в смысле определения на
с.\pageref{func-diff-na-otrezke}) и ее производная равна $f_{k+1}$:
$$
f'_k=f_{k+1},\qquad k<m
$$
Множество всех функций, гладких на $E$ порядка $m$, обозначается ${\mathcal
C}^m(E)$.

\item[$\bullet$]  Функция $f$ называется {\it бесконечно гладкой на $E$}, если
для всякого $m\in\N$ она является гладкой на $E$ порядка $m$:
$$
\forall m\in\N\qquad f\in{\mathcal C}^m(E).
$$
Множество всех бесконечно гладких функций на $E$ обозначается ${\mathcal
C}^\infty(E)$. Таким образом,
$$
{\mathcal C}^\infty(E)=\bigcap_{m=0}^\infty {\mathcal C}^m(E)
$$
Можно заметить, что бесконечная гладкость $f$ эквивалентна просто существованию
производной любого порядка, то есть существованию бесконечной
последовательности функций $\{f_k\}$ на $E$ такой, что
$$
f_0=f, \qquad \forall k\in\Z_+\quad \exists f'_k=f_{k+1}.
$$

\item[$\bullet$] {\it Нормой функции $f\in {\mathcal C}^m(E)$, равномерной на
множестве $E$ по производным до порядка $m\in\Z_+$}\index{норма!равномерная!по
производным}, называется величина
 \beq\label{ravnomernaya-norma-po-proizvodnym}
\norm{f}_E^{(m)}=\norm{f(x)}_{x\in E}^{(m)}:=\sum_{k=0}^m
\frac{1}{k!}\cdot\norm{f^{(k)}}_E=\sum_{k=0}^m \frac{1}{k!}\cdot\sup_{x\in I}
|f^{(k)}(x)|
 \eeq
Если для какого-нибудь $k$ функция $f^{(k)}$ не ограничена на $E$, то эта
величина считается бесконечной. Однако такое может случиться только если $E$ --
интервал, потому что если $E$ -- отрезок, то по теореме Вейерштрасса
\ref{Wei-II}, все функции $f^{(k)}$ будут ограничены на $E$.

Отметим, что равномерная по производным норма мажорирует обычную равномерную
норму, которую мы определяли формулой \eqref{ravnomernaya-norma}:
$$
\norm{f}_E\le\norm{f}_E^{(m)}
$$

 }\eit
\noindent\rule{160mm}{0.1pt}\begin{multicols}{2}

\begin{exs}
 \begin{align*}
& \norm{x^2}_{x\in [0;2]}^{(0)}=\sup_{x\in [0;2]} |x^2|=4 \\
& \norm{x^2}_{x\in [0;2]}^{(1)}=\sup_{x\in [0;2]} |x^2|+\sup_{x\in [0;2]}
|2x|=4+4=8 \\
& \norm{x^2}_{x\in [0;2]}^{(2)}=\sup_{x\in [0;2]} |x^2|+\sup_{x\in [0;2]} |2x|+
\\
&\qquad +\frac{1}{2}\cdot \sup_{x\in [0;2]} |2|=4+4+1=9
\\
& \norm{\sin x}_{x\in \R}^{(1)}=\sup_{x\in \R} |\sin x|+\sup_{x\in \R} |\cos
x|=1+1=2 \\
&\norm{\sin x}_{x\in\R}^{(2)}=\sup_{x\in\R} |\sin x|+\sup_{x\in\R} |\cos x|+
\\
&\qquad +\frac{1}{2}\cdot \sup_{x\in\R} |-\sin x|=1+1+\frac{1}{2}=2,5
 \end{align*}
\end{exs}
\end{multicols}\noindent\rule[10pt]{160mm}{0.1pt}

\bigskip

\centerline{\bf Свойства равномерной по производным нормы}
 \bit{\it

\item[$1^\circ.$] Неотрицательность:
 \beq\label{neotr-normy-po-proizv}
\norm{f}_E^{(m)}\ge 0
 \eeq
причем
 \beq\label{norma-po-proizv=0}
\norm{f}_E^{(m)}=0\quad\Longleftrightarrow\quad  \forall k=0,...,m \quad
\forall x\in E\quad f^{(k)}(x)=0
 \eeq

\item[$2^\circ.$] Однородность:
  \beq\label{odnorod-normy-po-proizv}
\norm{\lambda\cdot f}_E^{(m)} = |\lambda|\cdot\norm{\lambda\cdot f}_E^{(m)},
\qquad \lambda\in\R
 \eeq

\item[$3^\circ.$] Полуаддитивность:
  \beq\label{poluaddit-normy-po-proizv}
\norm{f+g}_E^{(m)}\le \norm{f}_E^{(m)}+\norm{g}_E^{(m)}
 \eeq

\item[$4^\circ.$] Полумультпликативность:
  \beq\label{polumultipl-normy-po-proizv}
 \norm{f\cdot g}_E^{(m)}\le \norm{f}_E^{(m)}\cdot \norm{g}_E^{(m)}
 \eeq

\item[$5^\circ.$] Монотонность по параметру $E$:
 \beq\label{monot-normy-po-proizv-po-E}
 D\subseteq E\qquad\Longrightarrow\qquad \norm{f}_D^{(m)}\le \norm{f}_E^{(m)}
 \eeq

\item[$6^\circ.$] Монотонность по параметру $m$:
  \beq\label{monot-normy-po-proizv-po-m}
 m\le n\qquad\Longrightarrow\qquad \norm{f}_E^{(m)}\le \norm{f}_E^{(n)}
 \eeq

\item[$7^\circ.$] Норма функции мажорирует норму ее производной: для любого
$k\le m$
  \beq\label{monot-normy-po-proizv-po-m}
\norm{f^{(k)}}_E^{(m)}\le \norm{f}_E^{(k+m)}
 \eeq

 }\eit

\bpr Как и в случае со свойствами равномерной нормы на
с.\pageref{svoistva-ravnom-normy}, все эти утверждения следуют напрямую из
свойств модуля. Мы предоставляем читателю самостоятельно проверить их
справедливость. \epr

\subsection{Равномерная по производным сходимость функциональных последовательностей}

\bit{

\item[$\bullet$] Говорят, что {\it последовательность функций $\{ f_n \}$
стремится к функции $f$ равномерно\index{сходимость!функциональной
последовательности!равномерная!по производным} по производным до порядка $m$ на
множестве $E$}, если
 \beq\label{ravnom-shod-po-proizv}
\norm{f_n-f}_E^{(m)}\underset{n\to \infty}{\longrightarrow} 0
 \eeq

 }\eit

Следующее утверждение аналогично теореме \ref{tm-19.5.1} и доказывается так же:

\begin{tm}
[\bf критерий Коши равномерной по производным сходимости функциональной
последовательности] \label{TH:cauchy-dlya-ravnom-po-proizv-shodim}
Функциональная последовательность $f_n\in{\mathcal C}^m(E)$ тогда и только
тогда сходится к некоторой функции $f\in{\mathcal C}^m(E)$ равномерно по
производным до порядка $m$ на множестве $E$, когда она удовлетворяет следующим
двум эквивалентным условиям:
 \bit{
\item[(i)] для любых двух бесконечно больших последовательностей индексов $\{
p_i \},\, \{ q_i \}\subseteq \mathbb{N}$
 \beq
p_i \underset{i\to \infty}{\longrightarrow}\infty, \quad q_i \underset{i\to
\infty}{\longrightarrow}\infty \label{19.5.1-*}
 \eeq
соответствующие подпоследовательности $\{ f_{p_i} \}$ и $\{ f_{q_i} \}$
равномерно по производным до порядка $m$ стремятся друг к другу на множестве
$E$:
 \beq
\norm{f_{p_i}-f_{q_i}}^{(m)}_E\underset{i\to \infty}{\longrightarrow} 0
\label{19.5.2-*}
 \eeq

\item[(ii)] для любой последовательности $\{ l_i \}\subseteq \mathbb{N}$ и
любой бесконечно большой последовательности $\{ k_i \}\subseteq \mathbb{N}$
 \beq
   k_i \underset{i\to \infty}{\longrightarrow}\infty
\label{19.5.3-*}
 \eeq
выполняется соотношение
 \beq
\norm{f_{k_i+l_i}-f_{k_i}}^{(m)}_E\underset{i\to \infty}{\longrightarrow} 0
\label{19.5.4-*}
 \eeq
 }\eit
\end{tm}

\begin{tm}[\bf о почленном дифференцировании последовательности, сходящейся равномерно по производным]
\label{TH:pochlennoe-differentsir-ravnom-po-proizv-shodysh-posledov} Если
последовательность функций $f_n$, гладких порядка $m$ на множестве $E$,
сходится на этом множестве равномерно по производным до порядка $m$, то
 \bit{

\item[(i)] предел этой последовательности
$$
f(x)=\lim_{n\to \infty}f_n(x),\qquad x\in E
$$
является функцией, гладкой порядка $m$ на $E$,

\item[(ii)] для всякого $k\le m$ производная порядка $k$ от предела совпадает с
пределом последовательности из производных порядка $k$ от ее членов на
множестве $E$:
$$
  f^{(k)}(x)=\lim_{n\to\infty} f_n^{(k)}(x),\qquad x\in E.
$$
 }\eit
\end{tm}

\subsection{Равномерная по производным сходимость функциональных рядов}

\bit{

\item[$\bullet$] Говорят, что {\it функциональный ряд $\sum_{n=0}^\infty a_n$
сходится на множестве $E$ равномерно по производным до порядка $m$ к функции
$S$}\index{сходимость!функционального ряда!равномерная!по производным}, если
частичные суммы этого ряда сходятся к $S$ на $E$ равномерно по производным до
порядка $m$:
 \beq\label{ravnom-shod-ryada-po-proizv}
\norm{\sum_{k=0}^n a_k-S}_E^{(m)}\underset{n\to \infty}{\longrightarrow} 0
 \eeq

 }\eit

Следующие две теоремы аналогичны теоремам \ref{tm-20.4.1} и \ref{tm-20.5.4}, и
доказывается так же:

\begin{tm}[\bf критерий Коши равномерной по производным сходимости ряда]\label{TH:cauchy-dlya-ravnom-po-proizv-shodim-ryada}
Для последовательности функций $a_n\in{\mathcal C}^m(E)$ следующие условия
эквивалентны:
 \bit{
\item[(i)] функциональный ряд
$$
\sum_{n=1}^\infty a_n(x)
$$
сходится равномерно по производным до порядка $m$ на множестве $E$ к некоторой
сумме $S\in{\mathcal C}^m(E)$;

\item[(ii)] для любой последовательности $l_i\in \mathbb{N}$, и любой
бесконечно большой последовательности $k_i\in \mathbb{N}$
$$
k_i\underset{i\to \infty}{\longrightarrow}\infty,
$$
сумма $\sum\limits_{n=k_i+1}^{k_i+l_i} a_n$ стремится к нулю при $i\to \infty$
равномерно по производным до порядка $m$ на множестве $E$:
$$
\norm{\sum\limits_{n=k_i+1}^{k_i+l_i} a_n}^{(m)}_E \underset{i\to
\infty}{\longrightarrow} 0
$$
 }\eit
\end{tm}

\begin{tm}[\bf признак Вейерштрасса для равномерной по производным сходимости]\label{TH:Weierstrass-dlya-ravnom-po-proizv-shodim-ryada}
Пусть $a_n\in{\mathcal C}^m(E)$, и числовой ряд
 $$
  \sum_{n=1}^\infty \norm{a_n(x)}_{x\in E}^{(m)}
 $$
сходится. Тогда функциональный ряд
$$
  \sum_{n=1}^\infty a_n(x)
$$
сходится равномерно по производным до порядка $m$ на множестве $E$.
\end{tm}

\begin{tm}[\bf о почленном дифференцировании ряда, сходящегося равномерно по производным]
\label{TH:pochlennoe-differentsir-ravnom-po-proizv-shodysh-ryada} Если ряд из
функций $a_n$, гладких порядка $m$ на множестве $E$
 $$
  \sum_{n=1}^\infty a_n(x)
 $$
сходится равномерно по производным до порядка $m$ на множестве $E$, то
 \bit{

\item[(i)] сумма этого ряда
$$
  S(x)=\sum_{n=1}^\infty a_n(x),\qquad x\in E,
$$
является функцией, гладких порядка $m$ на множестве $E$;

\item[(ii)] для всякого $k\le m$ производная порядка $k$ от суммы совпадает с
суммой ряда из производных порядка $k$ на множестве $E$:
$$
  S^{(k)}(x)=\sum_{n=1}^\infty a_n^{(k)}(x),\qquad x\in E.
$$
 }\eit
\end{tm}

\noindent\rule{160mm}{0.1pt}\begin{multicols}{2}

\subsection{Контрпримеры в классе гладких функций}

\bex\label{EX:gladkaya-func-0<f(x)<1} {\it Существует бесконечно гладкая
функция $f\in{\mathcal C}^\infty(\R)$ со следующими свойствами:}
$$
\begin{cases}
f(x)=0, & x\le 0 \\
f(x)>0,& x>0
\end{cases}
$$
\eex
 \bpr Рассмотрим функцию
$$
f(x)=\begin{cases}e^{-\frac{1}{x}}, & x>0\\ 0, & x\le 0
\end{cases}
$$
Все объявленные ее свойства очевидны, кроме бесконечной гладкости. Очевидно,
$f$ будет бесконечно гладкой в некоторой окрестности каждой точки $x\ne 0$,
поэтому нужно только доказать, что она является бесконечно гладкой в
окрестности точки $x=0$. Чтобы это понять, вычислим сначала односторонние
производные в этой точке. Ясно, что левая производная равна нулю:
$$
f'_{-}(0)=0
$$
Найдем правую:
$$
f'_+(0)=\lim_{x\to +0}\frac{e^{-\frac{1}{x}}}{x}=\left|\begin{matrix}
\frac{1}{x}=t \\ t\to+\infty\end{matrix}\right|= \lim_{t\to +\infty} t e^{-t}=0
$$

Мы получили, что первая производная существует в каждой точке. Существование
остальных проверяется по индукции. Предположим, что мы доказали, что существует
производная порядка $n$. Очевидно, она имеет вид
$$
f^{(n)}(x)=\begin{cases}\frac{P_n(x)}{Q_n(x)}e^{-\frac{1}{x}}, & x>0\\ 0, &
x\le 0
\end{cases}
$$
где $P_n(x)$ и $Q_n(x)$ некоторые многочлены. Отсюда получаем
 \begin{multline*}
f^{(n+1)}_+(0)=\lim_{x\to +0}\frac{P_n(x)e^{-\frac{1}{x}}}{xQ_n(x)}=\\=
\left|\begin{matrix} \frac{1}{x}=t \\
t\to+\infty\end{matrix}\right|= \lim_{t\to +\infty} \frac{R_n(t)}{S_n(t)}
e^{-t}=0
 \end{multline*}
где $R_n(t)$ и $S_n(t)$ -- некоторые новые многочлены.
 \epr

\bex\label{EX:gladkaya-func-a-b} {\it Для любого интервала $(a,b)$ на $\R$
существует бесконечно гладкая функция $g\in{\mathcal C}^\infty(\R)$ со
следующими свойствами:}
$$
\begin{cases}
g_{a,b}(x)=0, & x\notin(a,b)\\
g_{a,b}(x)>0,& x\in(a,b)
\end{cases}
$$
\eex
 \bpr
Можно положить
$$
g_{a,b}(x)=f(b-x)\cdot f(x-a),
$$
где $f$ -- функция из примера \ref{EX:gladkaya-func-0<f(x)<1}.
 \epr

\bex\label{EX:gladkaya-func-a->b} {\it Для любого интервала $(a,b)$ на $\R$
существует бесконечно гладкая функция $h_{a,b}\in{\mathcal C}^\infty(\R)$ со
следующими свойствами:}
$$
\begin{cases}
h_{a,b}(x)=0, & x\le a\\
0<h_{a,b}(x)<1,& x\in(a,b)\\
h_{a,b}(x)=1, & x\ge b
\end{cases}
$$
\eex
 \bpr
Этими свойствами будет обладать функция
$$
h_{a,b}(x)=\frac{\int_a^x g_{a,b}(t)\d t}{\int_a^b g_{a,b}(t)\d t},
$$
где $g_{a,b}$ -- функция из примера \ref{EX:gladkaya-func-a-b}.
 \epr

\bex\label{EX:gladkaya-func-a-alpha-beta-b} {\it Для любой последовательности
чисел
$$
a<\alpha<\beta<b
$$
существует бесконечно гладкая функция $\eta\in{\mathcal C}^\infty(\R)$ со
следующими свойствами:}
$$
\begin{cases}
\eta(x)=0, & x\notin(a,b)\\
\eta(x)=1, & x\in(\alpha,\beta)\\
0<\eta(x)<1,& x\in(a,\alpha)\cup(\beta,b)
\end{cases}
$$
\eex
 \bpr
Можно взять
$$
\eta(x)=h_{a,\alpha}(x)+1-h_{\beta,b}(x),
$$
где $h_{a,\alpha}$ и $h_{\beta,b}$ -- функции из примера
\ref{EX:gladkaya-func-a->b}.
 \epr

\blm\label{LM:func-dlya-Borelya} {\it Для любых $\e>0$ и $n\in\Z_+$ существует
гладкая функция $\ph\in{\mathcal C}^\infty[-1;1]$ со следующими свойствами:}
 \biter{
\item[1)] {\it в точке $x=0$ все ее производные равны нулю, кроме производной
порядка $n$, которая равна единице,}
$$
\ph^{(k)}(0)=\begin{cases}0, & k\ne n\\ 1,& k=n\end{cases}
$$

\item[2)] {\it ее равномерная норма на отрезке $[-1,1]$ по производным до
порядка $n-1$ меньше $\e$,}
$$
\norm{\ph}_{[-1,1]}^{(n-1)}<\e.
$$
 }\eiter
\elm
 \bpr
Положим $\delta=\frac{\e}{e}$ (где $e$ -- число Непера) и рассмотрим функцию
$\eta$ из примера \ref{EX:gladkaya-func-a-alpha-beta-b} со следующими
свойствами:
$$
\begin{cases}
\eta(x)=0, & x\notin(-\delta,\delta)\\
\eta(x)=1, & x\in\l-\frac{\delta}{2},\frac{\delta}{2}\r\\
0<\eta(x)<1,& \text{-- в остальных случаях}
\end{cases}
$$
Рассмотрим последовательность функций $f_n$, определенную рекуррентными
соотношениями:
 \begin{align*}
& f_0=\eta, \\
& f_m(x)=\begin{cases}
\int_0^x f_{m-1}(t)\d t, & x\notin[0,1]\\
-\int_x^0 f_{m-1}(t)\d t, & x\notin[-1,0)
\end{cases},\quad m\ge 1
 \end{align*}
По предложению \ref{PROP:pervoobr-nepr-f-na-intervale}, все функции $f_m$
являются гладкими на отрезке $[-1,1]$, причем
$$
f_m'=f_{m-1}, \qquad m\ge 1
$$
Отсюда следует формула:
 \beq\label{f_m^(k)=f_(m-k)}
f_m^{(k)}=\begin{cases}f_{m-k},& 0\le k<m \\
\eta,& k=m \\
\eta^{(k-m)},& k>m
\end{cases}
 \eeq
Из нее, в свою очередь, следует три важных для нас соотношения:
 \begin{align}
& f_m^{(k)}(0)=\begin{cases}0,& k\ne m \\ 1,& k=m \end{cases}, \label{f_m^(k)(0)=0-1} \\
& \norm{f_m}_{[-1,1]}=\sup_{x\in[-1,1]}|f_m(x)|\le\delta,\quad m\ge 1 \label{norm(f_m)_[-1,1]-le-delta} \\
& \norm{f_m}_{[-1,1]}^{(m-1)}<\e \label{norm(f_m)_[-1,1]^(m-1)<e}
 \end{align}

Здесь равенство \eqref{f_m^(k)(0)=0-1} доказывается перечислением случаев: если
$0\le k<m$, то
$$
f_m^{(k)}(0)=f_{m-k}(0)=\int_0^0 f_{m-1}(t)\d t=0,
$$
если $k=m$, то
$$
f_m^{(k)}(0)=\eta(0)=1,
$$
и если же $k>m$, то
$$
f_m^{(k)}(0)=\kern-20pt\underset{\scriptsize\begin{matrix}\uparrow\\
\text{постоянная} \\ \text{на интервале} \\
(-\frac{\delta}{2},\frac{\delta}{2})\end{matrix}}{\eta}\kern-20pt^{(k-m)}(0)=0.
$$

Неравенство \eqref{norm(f_m)_[-1,1]-le-delta} доказывается индукцией. При $m=1$
получаем, с одной стороны,
 \begin{multline*}
\norm{f_1}_{[0,1]}=\sup_{x\in[0,1]}|f_1(x)|=
\sup_{x\in[0,1]}\left|\int_0^x\eta(t)\ \d t\right|=\\=
\sup_{x\in[0,1]}\int_0^x\eta(t)\ \d t= \int_0^{\delta}\eta(t)\ \d t  \le\\ \le
\delta\cdot\sup_{x\in[0,\delta]}\eta(x)=\delta
 \end{multline*}
С другой стороны,
 \begin{multline*}
\norm{f_1}_{[-1,0]}=\sup_{x\in[-1,0]}|f_1(x)|=\\=
\sup_{x\in[-1,0]}\left|-\int_x^0\eta(t)\ \d t\right|=
\sup_{x\in[-1,0]}\int_x^0\eta(t)\ \d t=\\= \int_{-\delta}^0\eta(t)\ \d t  \le
\delta\cdot\sup_{x\in[-\delta,0]}\eta(x)=\delta
 \end{multline*}
Вместе это дает
$$
\norm{f_1}_{[-1,1]}=\max\left\{\norm{f_1}_{[-1,0]},\norm{f_1}_{[0,1]}\right\}\le
\delta
$$
Далее, если для $m-1$ неравенство \eqref{norm(f_m)_[-1,1]-le-delta} уже
доказано, то для $m$ получаем: с одной стороны,
 \begin{multline*}
\norm{f_m}_{[0,1]}=\sup_{x\in[0,1]}|f_m(x)|=\\=
\sup_{x\in[0,1]}\left|\int_0^xf_{m-1}(t)\ \d t\right|\le \\ \le
\sup_{x\in[0,1]}(x-0)\cdot \sup_{x\in[0,1]}|f_{m-1}(x)|\le\\ \le
1\cdot\underbrace{\norm{f_{m-1}}_{[-1,1]}}_{\scriptsize\begin{matrix}
\text{\rotatebox{90}{$\ge$}}
\\
\phantom{,}\delta, \\
\text{по предположению} \\ \text{индукции}\end{matrix}}\le \delta
 \end{multline*}
А с другой стороны,
 \begin{multline*}
\norm{f_m}_{[-1,0]}=\sup_{x\in[-1,0]}|f_m(x)|=\\=
\sup_{x\in[-1,0]}\left|\int_x^0 f_{m-1}(t)\ \d t\right|\le \\ \le
\sup_{x\in[-1,0]}(0-x)\cdot \sup_{x\in[-1,0]}|f_{m-1}(x)| \le\\ \le
1\cdot\underbrace{\norm{f_{m-1}}_{[-1,1]}}_{\scriptsize\begin{matrix}
\text{\rotatebox{90}{$\ge$}}
\\
\phantom{,}\delta, \\
\text{по предположению} \\ \text{индукции}\end{matrix}}\le \delta
 \end{multline*}
И вместе получается \eqref{norm(f_m)_[-1,1]-le-delta} для $m$:
$$
\norm{f_m}_{[-1,1]}=\max\left\{\norm{f_m}_{[-1,0]},\norm{f_m}_{[0,1]}\right\}\le
\delta
$$

Неравенство \eqref{norm(f_m)_[-1,1]^(m-1)<e} доказывается вычислением:
 \begin{multline*}
\norm{f_m}_{[-1,1]}^{(m-1)}=\sum_{k=0}^{m-1}\frac{1}{k!}\cdot\norm{f_m^{(k)}}_{[-1,1]}=\eqref{f_m^(k)=f_(m-k)}=\\=
\sum_{k=0}^{m-1}\frac{1}{k!}\cdot\norm{f_{m-k}}_{[-1,1]}\le\eqref{norm(f_m)_[-1,1]-le-delta}\le\\
\le
\sum_{k=0}^{m-1}\frac{1}{k!}\cdot\delta<\l\sum_{k=0}^\infty\frac{1}{k!}\r\cdot\delta=e\cdot\frac{\e}{e}=\e
 \end{multline*}

После того, как соотношения
\eqref{f_m^(k)(0)=0-1}--\eqref{norm(f_m)_[-1,1]^(m-1)<e} доказаны, нам остается
положить
$$
\ph=f_n
$$
и мы получим
$$
\ph^{(k)}=\eqref{f_m^(k)(0)=0-1}=\begin{cases}0,& k\ne n \\ 1,& k=n \end{cases}
$$
и
$$
\norm{\ph}_{[-1,1]}^{(n-1)}<\eqref{norm(f_m)_[-1,1]^(m-1)<e}<\e
$$
 \epr

\btm[Борель]\label{TH:Borel} Для любой последовательности чисел $a=\{a_n;\
n\in\Z_+\}$ существует гладкая функция $f$ на $\R$ такая, что
 \beq
\forall n\in\Z_+\qquad f^{(n)}(0)=a_n
 \eeq
\etm

\bpr Выберем последовательность функций $\ph_n$ со следующими свойствами:
 \biter{
\item[1)] если $a_n=0$, то $\ph_n=0$,

\item[2)] если $a_n\ne 0$, то $\ph_n$ выбирается по лемме
\ref{LM:func-dlya-Borelya} так, чтобы выполнялись условия:
 \begin{align}
& \ph_n^{(k)}=\begin{cases}0,& k\ne n \\ 1,& k=n \end{cases}, \label{Borel-1}\\
& \norm{\ph_n}_{[-1,1]}^{(n-1)}<\frac{1}{2^{n-1}\cdot |a_n|} \label{Borel-2}
 \end{align}
 }\eiter
Заметим, что функциональный ряд
$$
\sum_{n=0}^\infty a_n\cdot\ph_n
$$
сходится на $[-1,1]$ равномерно по производным до произвольного порядка $k$.
Действительно, при любом фиксированном $k$ мы получаем:
 \begin{multline*}
\sum_{n\ge k+1}
\underbrace{\norm{a_n\cdot\ph_n}_{[-1,1]}^{(k)}}_{\scriptsize\begin{matrix}\text{\rotatebox{90}{$\ge$}}
\\ \norm{a_n\cdot\ph_n}_{[-1,1]}^{(n-1)}, \\ \Uparrow \\
\eqref{monot-normy-po-proizv-po-m}\\ \Uparrow \\ k\le n-1
\end{matrix}}\le \sum_{n\ge k+1}
\norm{a_n\cdot\ph_n}_{[-1,1]}^{(n-1)}=\\=  \sum_{n\ge k+1}
|a_n|\cdot\norm{\ph_n}_{[-1,1]}^{(n-1)}\le\eqref{Borel-2}\le\\ \le \sum_{n\ge
k+1} |a_n|\cdot\frac{1}{2^{n-1}\cdot |a_n|}=\sum_{n\ge k+1}
\frac{1}{2^{n-1}}<\infty
 \end{multline*}
Из сходимости этого ряда по признаку Вейерштрасса
\ref{TH:Weierstrass-dlya-ravnom-po-proizv-shodim-ryada} следует, что его сумма
$$
f=\sum_{n=0}^\infty a_n\cdot\ph_n
$$
является гладкой функцией порядка $k$. Это справедливо для любого $k$, значит
$f$ является гладкой. При этом, по теореме о почленном дифференцировании
\ref{TH:pochlennoe-differentsir-ravnom-po-proizv-shodysh-ryada}
$$
f^{(k)}=\sum_{n=0}^\infty a_n\cdot\ph_n^{(k)}
$$
и отсюда
$$
f^{(k)}(0)=\sum_{n=0}^\infty
a_n\cdot\kern-15pt\underbrace{\ph_n^{(k)}(0)}_{\scriptsize\begin{matrix}\uparrow\\
\text{среди этих чисел}\\ \text{отлично от нуля}\\ \text{только с индексом}
\\ n=k
\end{matrix}}\kern-15pt =\eqref{Borel-1} =a_k
$$
\epr

\end{multicols}\noindent\rule[10pt]{160mm}{0.1pt}

\section{Аппроксимация}

\subsection{Приближение интегрируемых функций}

\paragraph{Интегральная норма функции.}

\bit{

\item[$\bullet$] {\it Интегральной нормой функции $f$ на отрезке
$[a,b]$}\index{норма!интегральная}, называется величина
 \beq\label{integralnaya-norma}
\inorm{f}_{[a,b]}=\inorm{f(x)}_{x\in[a,b]}:=\int_{[a,b]} |f(x)|\ \d x
 \eeq
(здесь предполагается, что функция $f$ определена и интегрируема на отрезке
$[a,b]$).
 }\eit

\noindent\rule{160mm}{0.1pt}\begin{multicols}{2}

\begin{exs}
 $$
 \inorm{x^2}_{x\in [0;2]}=\int_{[0;2]} |x^2|\ \d
x=\frac{x^3}{3}\Big|_{x=0}^{x=2}=\frac{8}{3}
 $$
 \begin{multline*}
\inorm{\sin x}_{x\in [0;\pi]}=\int_{[0;\pi]} |\sin x|\ \d x=\\= \int_0^{\pi}
\sin x\ \d x=-\cos x\Big|_{x=0}^{x=\pi}=2
 \end{multline*}
 \begin{multline*}
\inorm{\sin x}_{x\in [0,2\pi]}=\int_{[0,2\pi]} |\sin x|\ \d x=\\= \int_0^{\pi}
\sin x\ \d x-\int_{\pi}^{2\pi} \sin x\ \d x =\\=-\cos x\Big|_{x=0}^{x=\pi}
+\cos x\Big|_{x=\pi}^{x=2\pi}=2+2=4
 \end{multline*}
\end{exs}
\end{multicols}\noindent\rule[10pt]{160mm}{0.1pt}

\bigskip

\centerline{\bf Свойства интегральной нормы}\label{svoistva-integr-normy}
 \bit{\it

\item[$1^\circ.$] Неотрицательность:
 \beq\label{neotr-int-normy}
\inorm{f}_{[a,b]}\ge 0
 \eeq

\item[$2^\circ.$] Однородность:
  \beq\label{odnorod-int-normy}
\inorm{\lambda\cdot f}_{[a,b]} = |\lambda|\cdot\inorm{\lambda\cdot f}_{[a,b]},
\qquad \lambda\in\R
 \eeq

\item[$3^\circ.$] Полуаддитивность:
  \beq\label{poluaddit-int-normy}
\inorm{f+g}_{[a,b]}\le \inorm{f}_{[a,b]}+\inorm{g}_{[a,b]}
 \eeq

\item[$4^\circ.$] Монотонность по отрезку $[a,b]$:
 \beq\label{monot-int-normy-po-E}
 [a,b]\subseteq [c,d]\qquad\Longrightarrow\qquad \inorm{f}_{[a,b]}\le \inorm{f}_{[c,d]}
 \eeq

}\eit

\begin{proof} Как и в случаях с равномерной нормой и с равномерной по
производным нормой эти свойства следуют из свойств модуля (и в данном случае
еще из свойств интеграла).
 \end{proof}

\paragraph{Приближение интегрируемой функции кусочно постоянными.}

 \bit{
\item[$\bullet$] Функция $g$ на отрезке $[a,b]$ называется {\it
кусочно-постоянной}, если существует разбиение $\tau=\{x_0,...,x_k\}$ отрезка
$[a,b]$ такое, что на каждом интервале $(x_{i-1},x_i)$ функция $g$ постоянна:
$$
g(s)=g(t),\qquad s,t\in(x_{i-1},x_i)
$$
}\eit

Понятно, что всякая кусочно-постоянная функция $g$ будет кусочно-гладкой, и
значит, интегрируемой на $[a,b]$. Если через $C_i$ обозначить значения $g$ на
интервалах постоянства $(x_{i-1},x_i)$,
$$
g(x)=C_i,\qquad x\in(x_{i-1},x_i),
$$
то по лемме \ref{lm-14.7.1}, интеграл по каждому отрезку $[x_{i-1},x_i]$ от $g$
будет равен интегралу от $C_i$ по этому отрезку
$$
\int_{x_{i-1}}^{x_i} g(x)\ \d x=C_i\cdot\big(x_i-x_{i-1}\big)
$$
Как следствие, интеграл от $g$ по всему $[a,b]$ будет равен
 \beq\label{int-ot-kusochno-post}
\int_a^b g(x)\ \d x=\sum_{i=1}^k C_i\cdot\big(x_i-x_{i-1}\big)
 \eeq

\btm\label{TH:o-priblizh-integrir-func-kusoch-postoyann} Для любой
интегрируемой функции $f$ на отрезке $[a,b]$ и любого $\e>0$ можно найти
кусочно-постоянную функцию $g$ на $[a,b]$ такую, что
 \beq\label{priblizh-integrir-func-kusoch-postoyann}
\inorm{f-g}_{[a,b]}=\int_a^b |f(x)-g(x)|\ \d x<\e
 \eeq
Функцию $g$ можно выбрать так, чтобы выполнялось условие:
 \beq\label{norm(g)-le-norm(f)}
\sup_{x\in[a,b]}|g(x)|\le \sup_{x\in[a,b]}|f(x)|
 \eeq
\etm \bpr Выберем разбиение $\tau=\{x_0,...,x_k\}$ отрезка $[a,b]$ так, чтобы
нижняя сумма Дарбу $s_\tau$ отличалась от интеграла меньше чем на $\e$:
$$
0\le \int_a^b f(x)\ \d x-s_\tau<\e
$$
(этого всегда можно добиться, в силу соотношений
\eqref{S_tau_n->int<-s_tau_n}). Положим
$$
g(t)=\inf_{x\in[x_{i-1},x_i]}f(x),\qquad t\in [x_{i-1},x_i)
$$
(а в последней точке $x_k=b$ функцию $g$ можно определить, например, положив
$g(b)=f(b)$). Тогда мы получим:
$$
g(x)\le f(x),\qquad x\in[a,b]
$$
и поэтому
 \begin{multline*}
\int_a^b |\underbrace{f(x)-g(x)}_{\scriptsize\begin{matrix} \VI \\ 0
\end{matrix}}|\
\d x=\int_a^b \Big(f(x)-g(x)\Big)\ \d x= \int_a^b f(x)\ \d x-\int_a^b g(x)\ \d
x=\eqref{int-ot-kusochno-post}=\\=\int_a^b f(x)\ \d x- \underbrace{\sum_{i=1}^k
\inf_{x\in[x_{i-1},x_i]}f(x)\cdot\Delta x_i}_{\scriptsize\begin{matrix} \text{\rotatebox{90}{$=$}} \\
s_\tau \end{matrix}}=\int_a^b f(x)\ \d x-s_\tau<\e
 \end{multline*}
Условие \eqref{norm(g)-le-norm(f)} выполняется по построению $g$.
 \epr

\paragraph{Приближение интегрируемой функции непрерывными.}

\btm\label{TH:o-priblizh-integrir-func-nepreryvnyni} Для любой интегрируемой
функции $f$ на отрезке $[a,b]$ и любого $\e>0$ можно найти непрерывную функцию
$g$ на $[a,b]$ такую, что
 \beq\label{priblizh-integrir-func-nepreryvnyni}
\inorm{f-g}_{[a,b]}=\int_a^b |f(x)-g(x)|\ \d x<\e
 \eeq
Функцию $g$ можно выбрать так, чтобы выполнялось условие:
 \beq\label{norm(g)-le-norm(f)-2}
\sup_{x\in[a,b]}|g(x)|\le \sup_{x\in[a,b]}|f(x)|
 \eeq
 а также условие
 \beq\label{g(a)=g(b)}
g(a)=0=g(b)
 \eeq
 \etm
 \bpr 1. Рассмотрим сначала случай, когда функция
$f$ представляет собой ступеньку, то есть имеет вид
$$
f(x)=\begin{cases}C,& x\in[\alpha,\beta]\\ 0,& x\notin[\alpha,\beta]\end{cases}
$$
для некоторого отрезка $[\alpha,\beta]\subseteq[a,b]$ и некоторого числа
$C\in\R$. Тогда для любого $\e>0$ можно поправить функцию $f$ линейно в
окрестностях точек разрыва $\alpha$ и $\beta$ (или в окрестности одной из них,
если другая является концом отрезка $[a,b]$) так, чтобы получилась непрерывная
функция $g$, отличающаяся от $f$ только в этих окрестностях, и интеграл от
модуля разности будет маленьким.

На языке формул это можно выразить так: если $a\ne\alpha$ и $\beta\ne b$, то
положим
$$
\delta=\min\left\{\frac{\e}{3|C|},\alpha-a,b-\beta\right\},\qquad
g(x)=\begin{cases}0,&x<\alpha-\delta\\ \frac{C}{\delta}\cdot
x+C-\frac{C}{\delta}\cdot\alpha, & x\in[\alpha-\delta,\alpha]\\
C,& x\in[\alpha,\beta]\\
-\frac{C}{\delta}\cdot
x+C+\frac{C}{\delta}\cdot\beta, & x\in[\beta,\beta+\delta]\\
0,& x>\beta+\delta
\end{cases}
$$
и тогда $g$ непрерывна, отличается от $f$ только на интервалах
$(\alpha-\delta,\alpha)$ и $(\beta,\beta+\delta)$, и разность $f-g$ нигде не
превосходит $|C|$, поэтому
$$
\int_a^b |f(x)-g(x)|\ \d x=\int_{\alpha-\delta}^{\alpha} |f(x)-g(x)|\ \d
x+\int_{\beta}^{\beta+\delta} |f(x)-g(x)|\ \d x\le 2|C|\delta\le
2|C|\cdot\frac{\e}{3|C|}=\frac{2\e}{3}<\e
$$
Если же $a=\alpha$ или $\beta=b$, то в определении $g$ интервал
$(\alpha-\delta,\alpha)$ или интервал $(\beta,\beta+\delta)$ соответственно
нужно исключить, и результат будет тем же. Условие \eqref{norm(g)-le-norm(f)-2}
при этом выполняется автоматически.

2. После этого предположим, что функция $f$ кусочно-постоянна. Тогда, если не
считать ее значений в точках разрыва, она будет совпадать с суммой некоторого
конечного набора ступенек $f_1,...,f_n$:
$$
f=\sum_{i=1}^n f_i
$$
Для каждой ступеньки $f_i$ подберем непрерывную функцию $g_i$ так, чтобы
$$
\int_a^b|f_i(x)-g_i(x)|\ \d x<\frac{\e}{n}
$$
Тогда, положив $g=\sum_{i=1}^n g_i$, мы получим:
 \begin{multline*}
\int_a^b|f(x)-g(x)|\ \d x=\int_a^b\left|\sum_{i=1}^n f_i(x)-\sum_{i=1}^n
g_i(x)\right|\ \d x=\int_a^b\left|\sum_{i=1}^n \Big(f_i(x)-g_i(x)\Big)\right|\
\d x\le\\ \le \int_a^b\sum_{i=1}^n \Big|f_i(x)-g_i(x)\Big|\ \d x=
\sum_{i=1}^n\int_a^b\Big|f_i(x)-g_i(x)\Big|\ \d x<\sum_{i=1}^n\frac{\e}{n}=\e
 \end{multline*}
Если вдобавок $g_i$ выбирались так, чтобы выполнялось условие
\eqref{norm(g)-le-norm(f)-2},
$$
\sup_{x\in[a,b]}|g_i(x)|\le \sup_{x\in[a,b]}|f_i(x)|
$$
то это условие будет выполняться и для $g$:
$$
\sup_{x\in[a,b]}|g(x)|=\max_{1\le i\le n}\sup_{x\in[a,b]}|g_i(x)|\le \max_{1\le
i\le n}\sup_{x\in[a,b]}|f_i(x)|=\sup_{x\in[a,b]}|f(x)|
$$

3. Наконец, если $f$ -- произвольная интегрируемая функция, и $\e>0$, то,
подберем для нее по теореме \ref{TH:o-priblizh-integrir-func-kusoch-postoyann}
кусочно-постоянную функцию, обозначим ее $h$, так, чтобы
$$
\int_a^b|f(x)-h(x)|\ \d x<\frac{\e}{2}
$$
Затем, по уже доказанному, мы можем подобрать для функции $h$ непрерывную
функцию $g$ так, чтобы
$$
\int_a^b|h(x)-g(x)|\ \d x<\frac{\e}{2}
$$
И тогда мы получим:
 \begin{multline*}
\int_a^b|f(x)-g(x)|\ \d x=\int_a^b|f(x)-h(x)+h(x)-g(x)|\ \d x\le
\int_a^b\Big(|f(x)-h(x)|+|h(x)-g(x)|\Big)\ \d x=\\= \int_a^b|f(x)-h(x)|\ \d
x+\int_a^b|h(x)-g(x)|\ \d x<\frac{\e}{2}+\frac{\e}{2}=\e
 \end{multline*}
Если кроме того, $h$ и $g$ выбирались так, чтобы удовлетворять условиям
\eqref{norm(g)-le-norm(f)-2} и \eqref{norm(g)-le-norm(f)}, то мы получим
$$
\sup_{x\in[a,b]}|g(x)|\overset{\eqref{norm(g)-le-norm(f)-2}}{\le}
\sup_{x\in[a,b]}|h(x)|\overset{\eqref{norm(g)-le-norm(f)}}{\le}
\sup_{x\in[a,b]}|f(x)|
$$

4. Последнее условие \eqref{g(a)=g(b)} можно обеспечить, поправив еще немного
функцию $g$ в окрестности точки $a$ или $b$ так, как мы это делали в пункте 1
(то есть линейно).
 \epr

\subsection{Свертка}

 \bit{
\item[$\bullet$] {\it Сверткой} функций $f:\R\to\R$ и $g:\R\to\R$ называется
функция $f*g:\R\to\R$, определенная формулой
 \beq\label{DEF:svertka}
(f*g)(x)=\int_{-\infty}^{+\infty} f(y)\cdot g(x-y)\ \d
y=\lim_{\scriptsize\begin{matrix}A\to-\infty \\
B\to+\infty\end{matrix}}\int_A^B f(y)\cdot g(x-y)\ \d y, \qquad y\in\R
 \eeq
 }\eit
Здесь в правой части стоит несобственный интеграл, который, конечно, не для
любых $f$, $g$ и $x$ существует, и поэтому свертка $f*g$ не всегда определена.
Среди достаточных условий для существования свертки нас будет интересовать
только одно, а именно, то, в котором функции $f$ и $g$ локально интегрируемы,
причем одна из них финитна. Понятие локально интегрируемой функции было введено
нами на с.\pageref{DEF:loc-int-func}, а что такое финитная функция объясняется
в следующем определении:
 \bit{
\item[$\bullet$] Функция $f:\R\to\R$ называется {\it финитной}, если существует
отрезок $[a,b]$, вне которого функция $f$ обращается в нуль:
 $$
\forall x\notin[a,b]\qquad f(x)=0
 $$
Наименьший из таких отрезков мы будем называть {\it выпуклым носителем функции}
$f$ и обозначать $\convsupp f$:
$$
\convsupp f=\bigcap \Big\{[a,b]: \quad \forall x\notin[a,b]\qquad f(x)=0 \Big\}
$$
(если $f$ финитная и не везде нулевая, то такой отрезок всегда существует).
 }\eit

 \btm\label{TH:korrekt-opr-svertki}
Если функции  $f$ и $g$ локально интегрируемы, причем одна из них финитна, то
свертка $f*g:\R\to\R$ корректно определена и непрерывна.
 \etm
 \bpr Если финитна функция $f$, то интеграл в правой части \eqref{DEF:svertka}
будет равен интегралу по ее выпуклому носителю $[a,b]=\convsupp f$:
 \begin{multline*}
(f*g)(x)=\int_{-\infty}^{+\infty} f(y)\cdot g(x-y)\ \d y=
\lim_{\scriptsize\begin{matrix}A\to-\infty \\ B\to+\infty\end{matrix}}\int_A^B
f(y)\cdot g(x-y)\ \d y=\\= \lim_{\scriptsize\begin{matrix}A\to-\infty \\
B\to+\infty\end{matrix}}\Bigg\{\underbrace{\int_A^a f(y)\cdot g(x-y)\ \d
y}_{\scriptsize\begin{matrix}\text{\rotatebox{90}{$=$}}\\ 0\end{matrix}}
+\int_a^b f(y)\cdot g(x-y)\ \d y+\underbrace{\int_b^B f(y)\cdot g(x-y)\ \d
y}_{\scriptsize\begin{matrix}\text{\rotatebox{90}{$=$}}\\ 0\end{matrix}}\Bigg\}
=\\=\int_a^b f(y)\cdot g(x-y)\ \d y
 \end{multline*}
Если же финитна функция $g$, и $[a,b]=\convsupp g$ -- ее выпуклый носитель, то
для всякого фиксированного $x\in\R$ функция $y\mapsto g(x-y)$ будет тоже
финитной, и ее выпуклым носителем будет отрезок $[x-b,x-a]$:
$$
x-y\in[a,b]\quad\Leftrightarrow\quad a\le x-y\le b \quad\Leftrightarrow\quad
-a\ge y-x\ge -b\quad\Leftrightarrow\quad x-a\ge y\ge
x-b\quad\Leftrightarrow\quad y\in[x-b,x-a]
$$
Поэтому интеграл в \eqref{DEF:svertka} будет интегралом по отрезку $[x-b,x-a]$:
 \begin{multline*}
(f*g)(x)=\int_{-\infty}^{+\infty} f(y)\cdot g(x-y)\ \d y=
\lim_{\scriptsize\begin{matrix}A\to-\infty \\ B\to+\infty\end{matrix}}\int_A^B
f(y)\cdot g(x-y)\ \d y=\\= \lim_{\scriptsize\begin{matrix}A\to-\infty \\
B\to+\infty\end{matrix}}\Bigg\{\underbrace{\int_A^{x-b} f(y)\cdot g(x-y)\ \d
y}_{\scriptsize\begin{matrix}\text{\rotatebox{90}{$=$}}\\ 0\end{matrix}}
+\int_{x-b}^{x-a} f(y)\cdot g(x-y)\ \d y+\underbrace{\int_{x-a}^B f(y)\cdot
g(x-y)\
\d y}_{\scriptsize\begin{matrix}\text{\rotatebox{90}{$=$}}\\
0\end{matrix}}\Bigg\} =\\=\int_{x-b}^{x-a} f(y)\cdot g(x-y)\ \d y
 \end{multline*}
Остается проверить, что функция $f*g$ непрерывна -- мы это сделаем позже, на
с.\pageref{okonch-dok-TH:korrekt-opr-svertki}.
 \epr

 \bit{
\item[$\bullet$] {\it Сдвигом функции} $f:\R\to\R$ на величину $a\in\R$
называется функция $T_af:\R\to\R$, определенная правилом
 \beq\label{DEF:sdvig-functsii}
T_af(x):=f(x+a),\qquad x\in\R
 \eeq
 }\eit

\bigskip

\centerline{\bf Свойства свертки}
 \bit{\it

\item[$1^\circ.$] Дистрибутивность:
 \beq\label{distrib-svertki}
(f+g)*h=f*h+g*h
 \eeq

\item[$2^\circ.$] Коммутативность:
 \beq\label{kommut-svertki}
f*g=g*f
 \eeq

\item[$3^\circ.$] Перестановочность со сдвигом:
 \beq\label{perest-svertki-so-sdvigom}
(T_af)*g=T_a(f*g)=f*(T_ag)
 \eeq

\item[$4^\circ.$] Оценка равномерной нормы: если функция $g$ (локально
интегрируема и) сосредоточена на отрезке $[a,b]$,
$$
\forall x\notin [a,b]\quad g(x)=0,
$$
то для любой (локально интегрируемой) функции $f$ и любого отрезка $[A,B]$
выполняются неравенства:
 \begin{align}
& \norm{f*g}_{[A,B]}\le \norm{f}_{[A-b,B-a]}\cdot \inorm{g}_{[a,b]} \label{norma-svertki-1} \\
& \norm{f*g}_{[A,B]}\le \inorm{f}_{[A-b,B-a]}\cdot \norm{g}_{[a,b]}
\label{norma-svertki-2}
 \end{align}

 }\eit

\bpr 1. Первое свойство доказывается просто вычислением:
 \begin{multline*}
\Big((f+g)*h\Big)(x)=\int_{-\infty}^{+\infty} \Big(f(y)+g(y)\Big)\cdot h(x-y)\
\d y=\\= \int_{-\infty}^{+\infty} f(y)\cdot h(x-y)\ \d
y+\int_{-\infty}^{+\infty} g(y)\cdot h(x-y)\ \d y=(f*h)(x)+(g*h)(x)
 \end{multline*}

2. Следующие два свойства доказываются заменой переменной. Коммутативность:
 \begin{multline*}
f*g(x)=\int_{-\infty}^{+\infty} f(y)\cdot g(x-y)\ \d y=\lim_{\scriptsize\begin{matrix}A\to-\infty\\
B\to+\infty\end{matrix}}\int_A^B f(y)\cdot g(x-y)\ \d y={\scriptsize\left|\begin{matrix}x-y=t \\
y=x-t
\\ \d y=-\d
t\end{matrix}\right|}=\eqref{15.4.1-1}=\\=-\lim_{\scriptsize\begin{matrix}A\to-\infty\\
B\to+\infty\end{matrix}}\int_{x-A}^{x-B} f(x-t)\cdot g(t)\ \d
t=\lim_{\scriptsize\begin{matrix}A\to-\infty\\
B\to+\infty\end{matrix}}\int^{x-A}_{x-B} f(x-t)\cdot g(t)\ \d
t={\scriptsize\left|\begin{matrix}
\widetilde{B}=x-A\underset{A\to-\infty}{\longrightarrow}+\infty \\
\widetilde{A}=x-B\underset{B\to+\infty}{\longrightarrow}-\infty
\end{matrix}\right|}=\\=\lim_{\scriptsize\begin{matrix}\widetilde{A}\to-\infty\\
\widetilde{B}\to+\infty\end{matrix}}\int_{\widetilde{A}}^{\widetilde{B}}
f(x-t)\cdot g(t)\ \d t=\int_{-\infty}^{+\infty} f(x-t)\cdot g(t)\ \d t=(g*f)(x)
 \end{multline*}

3. Перестановочность со сдвигами:
 \begin{multline*}
(T_af)*g(x)=\int_{-\infty}^{+\infty} T_af(y)\cdot g(x-y)\ \d y=\lim_{\scriptsize\begin{matrix}A\to-\infty\\
B\to+\infty\end{matrix}}\int_A^B f(y+a)\cdot g(x-y)\ \d y={\scriptsize\left|\begin{matrix}y+a=t \\
y=t-a
\\ \d y=\d
t\end{matrix}\right|}=\eqref{15.4.1-1}=\\=\lim_{\scriptsize\begin{matrix}A\to-\infty\\
B\to+\infty\end{matrix}}\int_{A+a}^{B+a} f(t)\cdot g(x+a-t)\ \d
t={\scriptsize\left|\begin{matrix}
\alpha=A+a\underset{A\to-\infty}{\longrightarrow}-\infty \\
\beta=B+a\underset{B\to+\infty}{\longrightarrow}+\infty
\end{matrix}\right|}=\\=\lim_{\scriptsize\begin{matrix}\alpha\to-\infty\\
\beta\to+\infty\end{matrix}}\int_{\alpha}^{\beta} f(t)\cdot g(x+a-t)\ \d
t=\int_{-\infty}^{+\infty} f(t)\cdot g(x+a-t)\ \d
t=(g*f)(x+a)=\big(T_a(g*f)\big)(x)
 \end{multline*}
И точно так же доказывается вторая формула в \eqref{perest-svertki-so-sdvigom}.

4. Докажем \eqref{norma-svertki-1}:
 \begin{multline*}
\forall x\in[A,B]\qquad \Big|(f*g)(x)\Big|=\bigg|\int_{-\infty}^{\infty}
f(x-y)\cdot\kern-25pt\underbrace{g(y)}_{\scriptsize\begin{matrix}\uparrow\\
\text{отлично от нуля} \\
\text{только при $y\in[a,b]$}\end{matrix}}\kern-25pt\ \d y\bigg|=\left|\int_a^b
f(x-y)\cdot g(y)\ \d y\right|\le\\ \le \int_a^b
\underbrace{|f(\kern-10pt\overbrace{x-y}^{\scriptsize\begin{matrix}[A-b,B-a]
\\ \text{\rotatebox{90}{$\subseteq$}} \\ [x-b,x-a]\\
\text{\rotatebox{90}{$\in$}}\end{matrix}}\kern-10pt)|}_{\scriptsize\begin{matrix}\text{\rotatebox{90}{$\ge$}}\\
\norm{f}_{[A-b,B-a]}
\end{matrix}}\cdot\Big|g(y)\Big|\
\d y \le \norm{f}_{[A-b,B-a]}\cdot \int_a^b \Big|g(y)\Big|\ \d y=
\norm{f}_{[A-b,B-a]}\cdot\inorm{g}_{[a,b]}
 \end{multline*}
$$
\Downarrow
$$
$$
\norm{f*g}_{[A,B]}=\sup_{x\in[A,B]}\Big|(f*g)(x)\Big|\le
\norm{f}_{[A-b,B-a]}\cdot\inorm{g}_{[a,b]}
$$
Теперь \eqref{norma-svertki-2}:
 \begin{multline*}
\forall x\in[A,B]\qquad \Big|(f*g)(x)\Big|=\bigg|\int_{-\infty}^{\infty}
f(x-y)\cdot\kern-28pt\underbrace{g(y)}_{\scriptsize\begin{matrix}\uparrow\\
\text{отлично от нуля} \\
\text{только при $y\in[a,b]$}\end{matrix}}\kern-28pt \ \d
y\bigg|=\\=\left|\int_{a}^{b} f(x-y)\cdot g(y)\ \d y\right|\le\int_{a}^{b}
\Big|f(x-y)\Big|\cdot \underbrace{|g(y)|}_{\scriptsize\begin{matrix}
\text{\rotatebox{90}{$\ge$}}\\
\norm{g}_{[a,b]}\end{matrix}}\ \d y\le\eqref{14.6.3}\le \\
\le  \int_{a}^{b}
\Big|f(x-y)\Big|\ \d y\cdot \norm{g}_{[a,b]}={\scriptsize\left|\begin{matrix}x-y=t\\
y=x-t\\ \d y=-\d t \end{matrix}\right|}=-\norm{g}_{[a,b]}\cdot \int_{x-a}^{x-b}
\Big|f_n(t)-f(t)\Big|\ \d t=\\= \int_{x-b}^{x-a} \Big|f(t)\Big|\ \d t\cdot
\norm{g}_{[a,b]}\kern-40pt\underset{\scriptsize\begin{matrix}\uparrow\\
[x-b,x-a]\subseteq [A-b,B-a]\end{matrix}}{\le}\kern-40pt \int_{A-b}^{B-a}
\Big|f(t)\Big|\ \d t\cdot \norm{g}_{[a,b]}=\inorm{f}_{[A-b,B-a]}\cdot
\norm{g}_{[a,b]}
 \end{multline*}
$$
\Downarrow
$$
$$
\norm{f*g}_{[A,B]}=\sup_{x\in[A,B]}\Big|(f*g)(x)\Big|\le
\inorm{f}_{[A-b,B-a]}\cdot \norm{g}_{[a,b]}
$$
\epr

\blm\label{LM:T_s_nf->f} Пусть $f$ -- непрерывная функция на $\R$. Тогда для
любого отрезка $[\alpha,\beta]$ справедливо соотношение:
 \beq\label{T_s_nf->f-0}
\norm{T_{s_n}f-f}_{[\alpha,\beta]}=\norm{f(x+s)-f(x)}_{x\in[\alpha,\beta]}\underset{s\to
0}{\longrightarrow} 0
 \eeq \elm
 \bpr
Функция $f$ непрерывна на отрезке $[\alpha-1,\beta+1]$, значит по теореме
Кантора \ref{Kantor} она равномерно непрерывна на нем, то есть для любого
$\e>0$ существует $\delta>0$ такое что для любых $x$ и $y$ будет справедлива
импликация:
 $$
\left\{\begin{matrix} y\in[\alpha-1,\beta+1]\\ x\in[\alpha-1,\beta+1] \\
|y-x|<\delta \end{matrix}\right\} \quad\Longrightarrow\quad |f(y)-f(x)|<\e
 $$
Здесь при необходимости число $\delta$ можно уменьшить так, чтобы выполнялось
дополнительное неравенство
$$
\delta\le 1
$$
Тогда для любых $s\in(-\delta,\delta)$ и $x\in[\alpha,\beta]$ мы получим:
 $$
|f(x+s)-f(x)|\kern-40pt \overset{\scriptsize\begin{matrix}
\left\{\begin{matrix} x+s\in[\alpha-1,\beta+1]\\ x\in[\alpha-1,\beta+1] \\
|(x+s)-x|=|s|<\delta \end{matrix}\right\}
\\ \Downarrow\end{matrix}}{<}\kern-40pt\e
 $$
Поэтому
 \beq\label{T_s_nf->f-1}
\norm{T_sf-f}_{[\alpha,\beta]}=\sup_{x\in[\alpha,\beta]}|f(x+s)-f(x)|\le\e
 \eeq
То есть по произвольному данному $\e>0$ мы нашли такое $\delta>0$, что для всех
$s\in(-\delta,\delta)$ выполняется \eqref{T_s_nf->f-1}. Это и нужно было
доказать.
 \epr

\blm\label{LM:1/s_n(T_s_nf-f)->f'} Пусть $f$ -- гладкая функция на $\R$. Тогда
для любого отрезка $[\alpha,\beta]$ справедливо соотношение:
 \beq
\norm{\frac{T_{s}f-f}{s}-f'}_{[\alpha,\beta]}=\norm{\frac{f(x+s)-f(x)}{s}-f'(x)}_{x\in[\alpha,\beta]}\underset{s\to
0}{\longrightarrow}0
 \eeq
 \elm
\bpr Здесь тот же прием примененяется не к функции $f$, а к ее производной
$f'$. Функция $f'$ непрерывна на отрезке $[\alpha-1,\beta+1]$, значит по
теореме Кантора \ref{Kantor} она равномерно непрерывна на нем, то есть для
любого $\e>0$ существует $\delta>0$ такое что для любых $x$ и $y$ будет
справедлива импликация:
 $$
\left\{\begin{matrix} y\in[\alpha-1,\beta+1]\\ x\in[\alpha-1,\beta+1] \\
|y-x|<\delta \end{matrix}\right\} \quad\Longrightarrow\quad |f'(y)-f'(x)|<\e
 $$
Опять замечаем, что при необходимости число $\delta$ можно уменьшить так, чтобы
выполнялось дополнительное неравенство
$$
\delta\le 1
$$
Тогда для любых $s\in(-\delta,\delta)$ и $x\in[\alpha,\beta]$ мы получим:
 $$
\Big|\frac{1}{s}\cdot\underbrace{
\Big(f(x+s)-f(x)\Big)}_{\scriptsize\begin{matrix}
\phantom{\tiny\quad\eqref{Newton-Leibnitz-v-aniz-integr-R^1}}\text{\rotatebox{90}{$=$}}{\tiny\quad\eqref{Newton-Leibnitz-v-aniz-integr-R^1}}
\\
\int\limits_x^{x+s} f'(y)\ \d y
\\
\text{\rotatebox{90}{$=$}}
\\
\int\limits_0^s f'(x+t)\ \d t \\
\phantom{\tiny\quad\eqref{teor-o-srednem-v-aniz-int}}
\text{\rotatebox{90}{$=$}}{\tiny\quad\eqref{teor-o-srednem-v-aniz-int}}\\
f'(x+\tau)\cdot s,\\ \tau\in\overrightarrow{0s}
\end{matrix}}- f'(x)\Big|=|f'(x+\underset{\scriptsize\begin{matrix}\text{\rotatebox{90}{$\owns$}}\\
\overrightarrow{0s}\end{matrix}}{\tau})-f'(x)|\kern-50pt
\overset{\scriptsize\begin{matrix}
\left\{\begin{matrix} x+\tau\in[\alpha-1,\beta+1]\\ x\in[\alpha-1,\beta+1] \\
|(x+\tau)-x|=|\tau|\le|s|<\delta \end{matrix}\right\}
\\ \Downarrow\end{matrix}}{<}\kern-50pt\e
 $$
Поэтому
 \beq\label{1/s_n(T_s_nf-f)->f'-1}
\norm{\frac{T_{s}f-f}{s}-f'}_{[\alpha,\beta]}=\sup_{x\in[\alpha,\beta]}\left|\frac{f(x+s)-f(x)}{s}-f'(x)\right|\le\e
 \eeq
То есть по произвольному данному $\e>0$ мы нашли такое $\delta>0$, что для всех
$s\in(-\delta,\delta)$ выполняется \eqref{1/s_n(T_s_nf-f)->f'-1}. Это и нужно
было доказать.
 \epr

\bpr[Окончание доказательства теоремы
\ref{TH:korrekt-opr-svertki}.]\label{okonch-dok-TH:korrekt-opr-svertki}
Вспомним, что в доказательстве теоремы \ref{TH:korrekt-opr-svertki} мы отложили
на будущее проверку того, что свертка $f*g$ действительно является непрерывной
функцией. Теперь можем в этом убедиться.

1. Предположим сначала, что функция $g$ непрерывна и финитна. Пусть $[a,b]$ --
ее выпуклый носитель:
$$
\convsupp g=[a,b]
$$
Тогда для любой последовательности $s_n\to 0$, по модулю не превосходящей 1,
$$
|s_n|\le 1
$$
сдвиги $T_{s_n}g$ функции $g$ будут сосредоточены на отрезке $[a-1,b+1]$:
 \beq
\convsupp T_{s_n}g\subseteq [a-1,b+1]
 \eeq
Из леммы \ref{LM:T_s_nf->f} мы получим цепочку:
$$
\phantom{\text{\scriptsize (лемма \ref{LM:T_s_nf->f})}}\qquad
\norm{T_{s_n}g-g}_{[a-1,b+1]}\underset{n\to\infty}{\longrightarrow} 0\qquad
\text{\scriptsize (лемма \ref{LM:T_s_nf->f})}
$$
$$
\Downarrow
$$
 \begin{multline*}
\forall x\in\R\qquad
|(f*g)(x+s_n)-(f*g)(x)|=|T_{s_n}(f*g)(x)-(f*g)(x)|=\norm{T_{s_n}(f*g)-f*g}_{[x,x]}=\\=\eqref{perest-svertki-so-sdvigom}=
\norm{f*(T_{s_n}g)-f*g}_{[x,x]}=\eqref{distrib-svertki}=\norm{f*\Big(T_{s_n}g-g\Big)}_{[x,x]}\le\\
\le\eqref{norma-svertki-2}\le
\inorm{f}_{[x-b-1,x-a+1]}\cdot\underbrace{\norm{T_{s_n}g-g}_{[a-1,b+1]}}_{\scriptsize\begin{matrix}
\downarrow\\ 0\end{matrix}}\underset{n\to\infty}{\longrightarrow}0
 \end{multline*}
$$
\Downarrow
$$
$$
\forall x\in\R\qquad (f*g)(x+s_n)\underset{n\to\infty}{\longrightarrow}(f*g)(x)
$$
Это верно для любой последовательности $s_n\to 0$ со свойством $|s_n|\le 1$,
значит,
$$
\forall x\in\R\qquad (f*g)(x+s)\underset{s\to 0}{\longrightarrow}(f*g)(x)
$$
то есть функция $f*g$ непрерывна.

2. Теперь пусть $g$ -- произвольная локально интегрируемая финитная функция и
$[a,b]$ -- ее выпуклый носитель. По теореме
\ref{TH:o-priblizh-integrir-func-nepreryvnyni} можно подобрать
последовательность непрерывных функций $g_n$ так, чтобы выполнялись условия
$$
\inorm{g_n-g}_{[a,b]}\underset{n\to\infty}{\longrightarrow}0,\qquad g(a)=0=g(b)
$$
Второе из этих условий означает, что функции $g_n$ можно продолжить нулем вне
отрезка $[a,b]$, и они станут непрерывными функциями, определенными на всей
прямой $\R$, и, как и $g$, сосредоточенными на отрезке $[a,b]$, поэтому
$$
\forall n\qquad \convsupp (g_n-g)\subseteq[a,b]
$$
После этого для всякого отрезка $[A,B]\subset\R$ мы получим:
$$
\norm{f*g_n-f*g}_{[A,B]}=\eqref{distrib-svertki}=\norm{f*(g_n-g)}_{[A,B]}\le\eqref{norma-svertki-1}\le
\norm{f}_{[A-b,B-a]}\cdot\inorm{g_n-g}_{[a,b]}\underset{n\to\infty}{\longrightarrow}0
$$
$$
\Downarrow
$$
$$
(f*g_n)(x)\overset{x\in[A,B]}{\underset{n\to\infty}{\rightrightarrows}}
(f*g)(x)
$$
При этом, поскольку функции $g_n$ непрерывны, по уже доказанному, свертки
$f*g_n$ -- тоже непрерывные функции. Значит, по свойству $1^\circ$ на
с.\pageref{nepreryvnost-ravnom-predela}, равномерный на отрезке $[A,B]$ предел
$f*g$ функций $f*g_n$ должен быть непрерывной функцией на $[A,B]$. Это верно
для любого отрезка $[A,B]$, значит функция $f*g$ должна быть непрерывна всюду
на $\R$.
 \epr

 \btm\label{TH:gladkost-f*g} Если одна из функций $f$ или $g$ является гладкой, то свертка $f*g$ тоже
является гладкой, причем если гладкой является функция $f$, то
 \beq\label{diff-svertki-1}
(f*g)'=f'*g
 \eeq
если же гладкой является функция $g$, то
 \beq\label{diff-svertki-2}
(f*g)'=f*g'
 \eeq
 \etm
 \bpr Пусть функция $f$ -- гладкая, а $g$ -- финитная (и локально интегрируемая)
с носителем $[a,b]$. Тогда для произвольной последовательности
$s_n\underset{n\to\infty}{\longrightarrow}0$, $s_n\ne 0$, мы получим:
$$
\phantom{\text{\scriptsize (лемма \ref{LM:1/s_n(T_s_nf-f)->f'})}}\qquad \forall
[\alpha,\beta]\subset\R\qquad
\norm{\frac{1}{s_n}\Big(T_{s_n}f-f\Big)-f'}_{[\alpha,\beta]}\underset{n\to\infty}{\longrightarrow}0\qquad
\text{\scriptsize (лемма \ref{LM:1/s_n(T_s_nf-f)->f'})}
$$
$$
\Downarrow
$$
 \begin{multline*}
\forall x\in\R\qquad
\left|\frac{(f*g)(x+s_n)-(f*g)(x)}{s_n}-(f'*g)(x)\right|=\\=
\left|\frac{1}{s_n}\Big((T_{s_n}f*g)(x)-(f*g)(x)\Big)-(f'*g)(x)\right|=
\norm{\frac{1}{s_n}\Big(T_{s_n}f*g-f*g\Big)-f'*g}_{[x,x]}=\\=
\norm{\l\frac{1}{s_n}\Big(T_{s_n}f-f\Big)-f'\r*g}_{[x,x]}
\le\eqref{norma-svertki-1}\le
\underbrace{\norm{\frac{1}{s_n}\Big(T_{s_n}f-f\Big)-f'}_{[x-b,x-a]}}_{\scriptsize\begin{matrix}
\phantom{\tiny \begin{matrix}n\\ \downarrow\\ \infty\end{matrix}}
\quad \downarrow\quad {\tiny \begin{matrix}n\\ \downarrow\\
\infty\end{matrix}}\\ 0
\end{matrix}}\cdot\inorm{g}_{[a,b]}
\underset{n\to\infty}{\longrightarrow}0
 \end{multline*}
$$
\Downarrow
$$
$$
\forall x\in\R\qquad
\frac{(f*g)(x+s_n)-(f*g)(x)}{s_n}\underset{n\to\infty}{\longrightarrow}(f'*g)(x)
$$
Это верно для любой последовательности $s_n\to 0$. Поэтому
$$
\forall x\in\R\qquad \frac{(f*g)(x+s)-(f*g)(x)}{s}\underset{s\to
0}{\longrightarrow}(f'*g)(x)
$$
Мы получаем, что функция $f*g$ дифференцируема и ее производная будет равна
$f'*g$. Остается вспомнить, что в силу только что законченного доказательства
теоремы \ref{TH:korrekt-opr-svertki}, функция $f'*g$ будет непрерывна.

По аналогии рассматривается случай, когда функция $f$ -- произвольная (локально
интегрируемая), а $g$ -- гладкая и финитная. \epr

\bcor\label{COR:besk-gladkost-f*g} Если одна из функций $f$ или $g$ является
бесконечно гладкой, то их свертка тоже будет бесконечно гладкой. \ecor
 \bpr Пусть $f$ бесконечно гладкая. Тогда по теореме \ref{TH:gladkost-f*g}
 свертка $f*g$ будет гладкой, причем
 $$
(f*g)'=f'*g
 $$
Поскольку $f'$ тоже гладкая, опять по теореме \ref{TH:gladkost-f*g} получаем,
что свертка $f'*g$, то есть функция $(f*g)'$, будет гладкой, причем
$$
(f*g)''=(f'*g)'=f''*g
$$
И так далее.
 \epr

\subsection{Приближение непрерывных функций}

\paragraph{Аппроксимативная единица.}

 \bit{
\item[$\bullet$] Последовательность локально интегрируемых функций
$\varDelta_n:\R\to\R$ называется {\it аппроксимативной единицей}, если
выполняются следующие условия:
 \biter{
\item[(i)] функции $\varDelta_n$ имеют общий компактный носитель, то есть
существует отрезок $[\alpha,\beta]$, вне которого все они равны нулю:
 $$
\forall x\notin[\alpha,\beta]\qquad \forall n\in\N \qquad \varDelta_n(x)=0
 $$
если это верно, то в качестве $[\alpha,\beta]$ всегда можно выбрать отрезок
вида $[-D,D]$, где $D>0$, и тогда будет выполняться условие
 \beq\label{komp-nosit-approx-edinitsy}
\forall x\in\R \qquad\Big(\quad |x|>D \quad\Longrightarrow\quad \forall n\in\N
\qquad \varDelta_n(x)=0\quad\Big)
 \eeq

\item[(ii)] все функции $\varDelta_n$ неотрицательны:
 \beq\label{neotric-approx-edinitsy}
 \varDelta_n\ge 0
 \eeq

\item[(iii)] интеграл от всякой функции $\varDelta_n$ равен единице:
 \beq\label{integral-ot-approx-edinitsy}
\int_{-\infty}^{+\infty}\varDelta_n(x)\ \d x=1
 \eeq

\item[(iv)] для всякого $\delta>0$ выполняется соотношение:
 \beq\label{asymp-approx-edinitsy-1}
\int_{-\delta}^{+\delta}\varDelta_n(x)\ \d
x\underset{n\to\infty}{\longrightarrow} 1
 \eeq
 или, что равносильно, соотношение
 \beq\label{asymp-approx-edinitsy-0}
\int_{|x|\ge\delta}\varDelta_n(x)\ \d x\underset{n\to\infty}{\longrightarrow} 0
 \eeq
 }\eiter
 }\eit

\noindent\rule{160mm}{0.1pt}\begin{multicols}{2}

\bex Покажем, что функции
 \beq\label{algebr-approx-1}
\varDelta_n(x)=\begin{cases}\frac{1}{\int\limits_{-1}^1(1-x^2)^n\, \d x}\cdot
(1-x^2)^n,& |x|\le 1 \\ 0,& |x|>1
\end{cases}
 \eeq
образуют аппроксимативную единицу. Условия (i)-(iii) здесь очевидны. Докажем
(iv). Для всякого $\delta>0$ мы получим:
 \begin{multline*}
\int_{-1}^1(1-x^2)^n\ \d x=2\int_0^1(1-x^2)^n\ \d x\ge\\ \ge 2\int_0^1(1-x)^n\
\d x=-2\cdot\frac{(1-x)^{n+1}}{n+1}\Big|_{x=0}^{x=1}=\frac{2}{n+1}
 \end{multline*}
 $$
 \Downarrow
 $$
 $$
\frac{2}{\int\limits_{-1}^1(1-x^2)^n\, \d x}\le n+1
 $$
 $$
 \Downarrow
 $$
 \begin{multline*}
\int_{|x|\ge\delta}\varDelta_n(x)\ \d x=2\cdot\int_{\delta}^1\varDelta_n(x)\ \d
x=\\= \underbrace{\frac{2}{\int\limits_{-1}^1(1-x^2)^n\, \d
x}}_{\scriptsize\begin{matrix}\IA
\\ n+1\end{matrix}}\cdot\underbrace{\int_{\delta}^1(1-x^2)^n\, \d x}_{\scriptsize\begin{matrix}\IA
\\ \int_{\delta}^1(1-\delta^2)^n\, \d x\\ \text{\rotatebox{90}{$=$}} \\
(1-\delta^2)^n\cdot (1-\delta)
\end{matrix}}\le\\ \le (n+1)\cdot(1-\delta^2)^n\cdot
(1-\delta)\overset{\eqref{n^k/a^n->?}}{\underset{n\to\infty}{\longrightarrow}}0
 \end{multline*}
 \eex

\bex Покажем, что функции
 \beq\label{trigonom-approx-1}
\varDelta_n(x)=\begin{cases}\frac{1}{\int\limits_{-\pi}^{\pi}\cos^{2n}\frac{x}{2}\,
\d x}\cdot \cos^{2n}\frac{x}{2},& |x|\le \pi \\ 0,& |x|>\pi
\end{cases}
 \eeq
образуют аппроксимативную единицу. Как и в предыдущем примере, здесь нужно
проверить только условие (iv). Для всякого $\delta>0$ мы получим:
 \begin{multline*}
\int_{-\pi}^{\pi}\cos^{2n}\frac{x}{2}\ \d x=\left|\begin{matrix}\frac{x}{2}=y\\
x=2y\end{matrix}\right|=\\=2\int_{-\frac{\pi}{2}}^{\frac{\pi}{2}}\cos^{2n}y\ \d
y=4\int_0^{\frac{\pi}{2}}
\kern-18pt\underbrace{\cos^{2n}y}_{\scriptsize\begin{matrix}\VI
\\ 1-\frac{2y}{\pi}, \\ \text{потому что} \\
\text{$\cos$ -- выпуклая} \\
\text{функция на $\left[0,\frac{\pi}{2}\right]$}
\end{matrix}}\kern-18pt\ \d y \ge\\
\ge 4\int_0^{\frac{\pi}{2}}\l 1-\frac{2y}{\pi}\r^{2n}\ \d
y=\\=-\frac{2\pi}{2n+1} \l
1-\frac{2y}{\pi}\r^{2n+1}\bigg|_{y=0}^{y=\frac{\pi}{2}}=\frac{2\pi}{2n+1}
 \end{multline*}
 $$
 \Downarrow
 $$
 $$
\frac{1}{\int_{-\pi}^{\pi}\cos^{2n}\frac{x}{2}\ \d x}\le\frac{2n+1}{2\pi}
 $$
 $$
 \Downarrow
 $$
 \begin{multline*}
\int_{|x|\ge\delta}\varDelta_n(x)\ \d x=2\int_{\delta}^{\pi}\varDelta_n(x)\ \d
x=\\= \underbrace{\frac{1}{\int_{-\pi}^{\pi}\cos^{2n}\frac{x}{2}\ \d
x}}_{\scriptsize\begin{matrix}\IA
\\ \frac{2n+1}{2\pi}\end{matrix}}\cdot 2\int_{\delta}^{\pi}
\Big(\cos\frac{x}{2}\Big)^{2n}\ \d x\le \\ \le \frac{2n+1}{2\pi}\cdot
2\int_{\delta}^{\pi} \Big(\cos\frac{x}{2}\Big)^{2n}\ \d x=\left|\begin{matrix}\frac{x}{2}=y\\
x=2y\end{matrix}\right|=\\= \frac{2n+1}{2\pi}\cdot
4\int_{\delta}^{\frac{\pi}{2}} \Big(\kern-10pt\underbrace{\cos
y}_{\scriptsize\begin{matrix}\IA
\\ \cos\delta, \\ \text{потому что} \\
\text{$\cos$ убывает} \\
\text{на $\left[0,\frac{\pi}{2}\right]$}
\end{matrix}}\kern-10pt\Big)^{2n}\ \d y \le\\
\le \frac{4n+2}{\pi}\cdot\int_{\delta}^{\frac{\pi}{2}}\cos^{2n}\delta\ \d x=\\=
\frac{4n+2}{\pi}\cdot\Big(\underbrace{\cos^2\delta}_{\scriptsize\begin{matrix}
\text{\rotatebox{90}{$>$}}
\\ 1\end{matrix}}\Big)^n\cdot
\l\frac{\pi}{2}-\delta\r\overset{\eqref{n^k/a^n->?}}{\underset{n\to\infty}{\longrightarrow}}0
 \end{multline*}
 \eex

\bex\label{EX:gladkaya-approx-edin} {\bf Бесконечно гладкая аппроксимативная
единица.} Следуя примеру \ref{EX:gladkaya-func-a-b} построим гладкую функцию
$g\in{\mathcal C}^\infty(\R)$ со следующими свойствами:
$$
\begin{cases}g(x)=0,& x\notin(-1,1)\\
g(x)>0, & x\in(-1,1) \end{cases}
$$
Положив
$$
h(x)=\frac{g(x)}{\int_{-1}^1 g(t) \ \d t}
$$
мы получим функцию со следующими свойствами:
$$
 \begin{cases}
 h(x)=0,& x\notin(-1,1)\\
 h(x)>0, & x\in(-1,1) \\
 \int_{-1}^1 h(x)\ \d x=1
 \end{cases}
$$
Теперь положив
$$
\varDelta_n(x)=n\cdot h(nx)
$$
мы получим последовательность неотрицательных функций $\varDelta_n\in{\mathcal
C}^\infty(\R)$, имеющих в качестве общего носителя отрезок $[-1,1]$:
 \begin{multline*}
x\notin(-1,1)\quad\Longrightarrow\quad nx\notin(-1,1)\quad\Longrightarrow\\
\Longrightarrow\quad \varDelta_n(x)=h(nx)=0
 \end{multline*}
Интеграл от каждой такой функции будет равен единице:
 \begin{multline*}
\int_{-\infty}^{\infty}\varDelta_n(x)\ \d x= \int_{-\infty}^{\infty} n\cdot
h(nx)\ \d x=\\= \int_{-\infty}^{\infty} h(nx)\ \d(nx)= \int_{-\infty}^{\infty}
h(y)\ \d y=1
 \end{multline*}
А для каждого $\delta>0$ и любого $n>\frac{1}{\delta}$ мы получим:
 \begin{multline*}
\int_{-\delta}^{\delta}\varDelta_n(x)\ \d x= \int_{-\delta}^{\delta} n\cdot
h(nx)\ \d x=\\= \int_{-\delta}^{\delta} h(nx)\ \d(nx)=
\underbrace{\int_{-n\delta}^{n\delta} h(y)\ \d y}_{n\delta>1}=\\=\int_{-1}^{1}
h(y)\ \d y=1
 \end{multline*}
Таким образом, выполняются все четыре условия (i)-(iv) на
с.\pageref{komp-nosit-approx-edinitsy}, и наша последовательность $\varDelta_n$
является аппроксимативной единицей.
 \eex

\end{multicols}\noindent\rule[10pt]{160mm}{0.1pt}

\btm\label{TH:approx-1} Пусть $f:\R\to\R$ -- непрерывная функция на $\R$. Тогда
для всякой аппроксимативной единицы $\varDelta_n$ последовательность сверток
$f*\varDelta_n$ стремится к $f$ равномерно на каждом отрезке $[a,b]$:
 \beq
f*\varDelta_n(x)\overset{x\in[a,b]}{\underset{n\to\infty}{\rightrightarrows}}
f(x)
 \eeq
 \etm
\bpr Зафиксируем отрезок $[a,b]$ и число $D$ со свойством
\eqref{komp-nosit-approx-edinitsy}. Поскольку функция $f$ непрерывна на отрезке
$[a-D,b+D]$, по теореме Вейерштрасса \ref{Wei-III} должна быть конечна величина
 \beq\label{approx-1-0}
 M=\norm{f}_{[a-D,b+D]}=\sup_{s\in[a-D,b+D]}|f(s)|
 \eeq
Пусть далее $\e>0$. Поскольку функция $f$ непрерывна на $\R$, она непрерывна и
на отрезке $[a-D,b+D]$. По теореме Кантора \ref{Kantor} это означает, что $f$
равномерно непрерывна на $[a-D,b+D]$, поэтому для числа $\frac{\e}{2}>0$ должно
существовать $\delta>0$, которое можно подчинить дополнительному условию
$\delta<D$, такое, что
$$
\forall s,t\in [a-D,b+D]\qquad |s-t|<\delta\quad \Longrightarrow\quad
|f(s)-f(t)|<\frac{\e}{2}
$$
В частности, при $t=x\in[a,b]$ и $s=x-y$, где $y\in(-\delta,\delta)$, мы
получим такое соотношение:
 \beq\label{approx-1-1}
\forall x\in [a,b]\quad \forall y\in(-\delta,\delta)\qquad
|f(x-y)-f(x)|<\frac{\e}{2}
 \eeq
Теперь для любого $x\in[a,b]$ мы получаем:
 \begin{multline*}
\big|f*\varDelta_n(x)-f(x)\big|=\bigg|\overbrace{\int_{-\infty}^{\infty}
f(x-y)\cdot\varDelta_n(y)\d y}^{\scriptsize\begin{matrix} f*\varDelta_n(x) \\
\phantom{\tiny\eqref{DEF:svertka}}\text{\rotatebox{90}{$=$}}{\tiny\eqref{DEF:svertka}}\end{matrix}}-f(x)\cdot
\overbrace{\int_{-\infty}^{\infty} \varDelta_n(y)\d
y}^{\scriptsize\begin{matrix} 1 \\
\phantom{\tiny\eqref{integral-ot-approx-edinitsy}}\text{\rotatebox{90}{$=$}}{\tiny\eqref{integral-ot-approx-edinitsy}}
\end{matrix}}
\bigg|=\\=\bigg|\int_{-\infty}^{\infty} f(x-y)\cdot\varDelta_n(y)\d y-
\int_{-\infty}^{\infty} f(x-y)\cdot\varDelta_n(y)\d y \bigg|=
\bigg|\int_{-\infty}^{\infty}\Big(
f(x-y)-f(x)\Big)\cdot\kern-5pt\underbrace{\varDelta_n(y)}_{\scriptsize\begin{matrix}\text{\rotatebox{90}{$=$}}\\
\phantom{,}0,\\ \text{при $|y|>D$}\end{matrix}}\kern-5pt\d y\bigg|=\\=
\bigg|\int_{-D}^D\Big( f(x-y)-f(x)\Big)\cdot\varDelta_n(y)\d y\bigg|\le
\int_{-D}^D\Big| f(x-y)-f(x)\Big|\cdot\varDelta_n(y)\d y\le \\ \le
\int_{-\delta}^{\delta}\underbrace{\Big|
f(x-y)-f(x)\Big|}_{\scriptsize\begin{matrix}\phantom{\tiny\eqref{approx-1-1}}
\text{\rotatebox{90}{$\ge$}}{\tiny\eqref{approx-1-1}}\\ \frac{\e}{2}
\end{matrix}}\cdot\varDelta_n(y)\d y+ \int_{|y|\ge \delta}\underbrace{\Big|
f(x-y)-f(x)\Big|}_{\scriptsize\begin{matrix} \text{\rotatebox{90}{$\ge$}}\\
|f(x-y)|+|f(x)|\\ \phantom{\tiny\eqref{approx-1-0}} \text{\rotatebox{90}{$\ge$}}{\tiny\eqref{approx-1-0}}\\ M+M \\
\text{\rotatebox{90}{$=$}}\\ 2M
\end{matrix}}\cdot\varDelta_n(y)\d y\le \\ \le
\frac{\e}{2}\cdot \int_{-\delta}^{\delta}\varDelta_n(y)\d y+2M\cdot\int_{|y|\ge
\delta}\varDelta_n(y)\d y
 \end{multline*}
Это верно для любого $x\in[a,b]$. Значит,
$$
\norm{f*\varDelta_n-f}_{[a,b]}=\sup_{x\in[a,b]}\big|f*\varDelta_n(x)-f(x)\big|\le
\frac{\e}{2}\cdot \underbrace{\int_{-\delta}^{\delta}\varDelta_n(y)\d
y}_{\scriptsize\begin{matrix}\phantom{\tiny
\begin{matrix}n\\ \downarrow\\ \infty \end{matrix}}{\tiny \eqref{asymp-approx-edinitsy-1}} \ \downarrow \ {\tiny
\begin{matrix}n\\ \downarrow\\ \infty \end{matrix}}\phantom{\tiny \eqref{asymp-approx-edinitsy-1}}
\\ 1\end{matrix}}+2M\cdot\underbrace{\int_{|y|\ge
\delta}\varDelta_n(y)\d y}_{\scriptsize\begin{matrix}\phantom{\tiny
\begin{matrix}n\\ \downarrow\\ \infty \end{matrix}}{\tiny \eqref{asymp-approx-edinitsy-0}} \ \downarrow \ {\tiny
\begin{matrix}n\\ \downarrow\\ \infty \end{matrix}}\phantom{\tiny \eqref{asymp-approx-edinitsy-0}}
\\ 0\end{matrix}}\underset{n\to\infty}{\longrightarrow}\frac{\e}{2}<\e
$$
И это верно для любого $\e>0$. Мы получаем, что, какое ни возьми $\e>0$, для
почти всех $n\in\N$ будет верно неравенство
$$
\norm{f*\varDelta_n-f}_{[a,b]}<\e
$$
То есть,
$$
\norm{f*\varDelta_n-f}_{[a,b]}\underset{n\to\infty}{\longrightarrow}0,
$$
а это нам и нужно было доказать.
 \epr

\paragraph{Аппроксимация гладкими функциями.}

\btm\label{TH:approx-nepr-func-gladkimi} Для любой непрерывной функции
$f:\R\to\R$ можно подобрать последовательность бесконечно гладких функций
$\ph_n:\R\to\R$, равномерно сходящуюся к $f$ на каждом отрезке
$[a,b]\subset\R$:
 \beq
\forall [a,b]\subset\R\qquad
\ph_n(x)\overset{x\in[a,b]}{\underset{n\to\infty}{\rightrightarrows}} f(x)
 \eeq
\etm \bpr Пусть $\varDelta_n$ -- бесконечно гладкая аппроксимативная единица,
построенная в примере \ref{EX:gladkaya-approx-edin}. Тогда по следствию
\ref{COR:besk-gladkost-f*g}, функции
$$
\ph_n=f*\varDelta_n
$$
будут бесконечно гладкими, а по теореме \ref{TH:approx-1}, они будут стремиться
к $f$ равномерно на каждом отрезке в $\R$. \epr

\bcor Для любой непрерывной функции $f:[a,b]\to\R$ можно подобрать
последовательность бесконечно гладких функций $\ph_n:[a,b]\to\R$, равномерно
сходящуюся к $f$ на $[a,b]$:
 \beq
\ph_n(x)\overset{x\in[a,b]}{\underset{n\to\infty}{\rightrightarrows}} f(x)
 \eeq
\ecor \bpr Нужно сначала продолжить $f$ до непрерывной функции на всю прямую
$\R$. Это можно сделать, например, положив
$$
f(x)=\begin{cases}f(a),& x< a \\ f(b),& x>b\end{cases}
$$
После этого по теореме \ref{TH:approx-nepr-func-gladkimi}, найдется
последовательность гладких функций $\ph_n$, приближающая $f$ равномерно на
каждом отрезке, в частности, на отрезке $[a,b]$. \epr

\paragraph{Аппроксимация алгебраическими многочленами.} Напомним определение, данное на с.\pageref{EX:algebr-mnogochl}:

\bit{

\item[$\bullet$] Функции вида
$$
f(x)=\sum_{k=0}^n a_k\cdot x^k\qquad (x\in\R),
$$
где $n\in\N$ и $a_k\in\R$, называются {\it алгебраическими многочленами}\label{EX:algebr-mnogochl*} или
просто {\it многочленами} (от одной переменной). Если $a_n\ne 0$, то число $n$
называется {\it степенью многочлена} $f$.
}\eit

\btm[\bf Вейерштрасса об аппроксимации алгебраическими многочленами] Для любой
непрерывной функции $f:[a,b]\to\R$ можно подобрать последовательность
многочленов $\ph_n:[a,b]\to\R$, равномерно сходящуюся к $f$ на $[a,b]$:
 \beq
f_n(x)\overset{x\in[a,b]}{\underset{n\to\infty}{\rightrightarrows}} f(x)
 \eeq
\etm \bpr 1. Прежде всего заметим, что линейной заменой переменной утверждение
сводится к случаю $[a,b]\subset(0,1)$. Например, функция
$$
\ph(t)=3(b-a)\cdot t+2a-b
$$
превращает отрезок $\left[\frac{1}{3},\frac{2}{3}\right]$ в отрезок $[a,b]$:
$$
\ph\l\left[\frac{1}{3},\frac{2}{3}\right]\r=[a,b]
$$
А функция $f$ в композиции с $\ph$ превращается в функцию $f\circ\ph$ на
отрезке $\left[\frac{1}{3},\frac{2}{3}\right]$. Если мы сможем построить
последовательность многочленов $g_n$, приближающую $f\circ\ph$ равномерно на
отрезке $\left[\frac{1}{3},\frac{2}{3}\right]$ (содержащимся в интервале
$(0,1)$)
$$
g_n(t)\overset{t\in\left[\frac{1}{3},\frac{2}{3}\right]}{\underset{n\to\infty}{\rightrightarrows}}
f\Big(\ph(t)\Big),
$$
то функции $f_n=g_n\circ\ph^{-1}$ будут многочленами, приближающими $f$
равномерно на $[a,b]$:
$$
f_n(x)=g_n\Big(\ph^{-1}(x)\Big)\overset{x\in[a,b]}{\underset{n\to\infty}{\rightrightarrows}}
f\Big(\ph\big(\ph^{-1}(x)\big)\Big)=f(x)
$$

2. Итак, можно считать, что отрезок $[a,b]$ лежит в интервале $(0;1)$:
$$
[a,b]\subset(0,1)
$$
Тогда $f$ можно доопределить непрерывно на всю прямую $\R$ так, чтобы ее
выпуклым носителем был отрезок $[0,1]$. Например, можно на оставшихся кусках
интервала $(0,1)$ доопределить функцию $f$ линейно, то есть положить
$$
f(x)=\begin{cases}\frac{f(a)}{a}\cdot x,& x\in(0,a)\\
\frac{f(b)}{b-1}\cdot(x-b)+f(b),& x\in(b,1)\\ 0, &
x\in(-\infty,0]\cup[1,+\infty)\end{cases}
$$
Запомним это, и рассмотрим аппроксимативную единицу \eqref{algebr-approx-1}:
$$
\varDelta_n(x)=\begin{cases}C_n\cdot (1-x^2)^n,& |x|\le 1 \\ 0,& |x|>1
\end{cases},\qquad C_n=\frac{1}{\int\limits_{-1}^1(1-x^2)^n\, \d x}
$$
По теореме \ref{TH:approx-1}, выполняется соотношение
$$
f*\varDelta_n(x)\overset{x\in[a,b]}{\underset{n\to\infty}{\rightrightarrows}}
f(x)
$$
Нам остается только заметить, что на отрезке $[a,b]$ функции $f*\varDelta_n$
представляют собой многочлены:
 \begin{multline*}
f*\varDelta_n(x)=\int_{-\infty}^{+\infty}
f(y)\cdot\kern-15pt\underbrace{\varDelta_n(x-y)}_{\scriptsize\begin{matrix}\text{\rotatebox{90}{$=$}}\\
\phantom{,}0, \\
\text{при $y\notin[x-1,x+1]$,}\\ \text{поскольку}\\
\convsupp\varDelta_n=[-1,1]\end{matrix}}\kern-15pt\ \d y=\int_{x-1}^{x+1}
f(y)\cdot\varDelta_n(x-y)\ \d y=\\= C_n\cdot\int_{x-1}^{x+1}
\kern-25pt\underbrace{f(y)}_{\scriptsize\begin{matrix}\text{\rotatebox{90}{$=$}}\\
\phantom{,}0, \\
\text{при $y\notin[0,+1]$,}\\ \text{причем}\\
x-1\le 0<1\le x+1\end{matrix}}\kern-25pt\cdot \Big(1-(x-y)^2\Big)^n\ \d y=
C_n\cdot\int_0^1 f(y)\cdot\underbrace{\Big(1-(x-y)^2\Big)^n}_{\text{многочлен
по $y$}}\ \d y=\\=C_n\cdot\int_0^1 f(y)\cdot\sum_{k=0}^{2n}p_k(x)\cdot y^n\ \d
y= \sum_{k=0}^{2n} \overbrace{\l C_n\cdot\int_0^1 f(y)\cdot y^n\ \d
y\r}^{\text{число}}\cdot\kern-5pt
\underbrace{p_k(x)}_{\scriptsize\begin{matrix}\text{многочлен}\\ \text{по
$x$}\end{matrix}}
 \end{multline*}

\epr

\paragraph{Аппроксимация тригонометрическими многочленами.}

\bit{

\item[$\bullet$] Функции вида
$$
f(x)=c+\sum_{k=1}^n \Big\{a_k\cdot \cos kx+b_k\cdot\sin kx\Big\}\qquad
(x\in\R),
$$
где $n\in\N$, $c,a_k,b_k\in\R$, называются {\it тригонометрическими
многочленами}\label{EX:trigonom-mnogochl}. Если $a_n\ne 0$ или $b_n\ne 0$, то число $n$ называется {\it
степенью тригонометрического многочлена} $f$.
 }\eit

\btm[\bf Вейерштрасса об аппроксимации тригонометрическими
многочленами]\label{TH:Weierstrass-approx-trig-mnogochl} Для любой непрерывной
$2\pi$-периодической функции $f:\R\to\R$ можно подобрать последовательность
тригонометрических многочленов $f_n:\R\to\R$, равномерно сходящуюся к $f$ на
$\R$:
 \beq
f_n(x)\overset{x\in\R}{\underset{n\to\infty}{\rightrightarrows}} f(x)
 \eeq
\etm \bpr Здесь нужно рассмотреть аппроксимативную единицу
\eqref{trigonom-approx-1}. По теореме \ref{TH:approx-1}, выполняется
соотношение
$$
f*\varDelta_n(x)\overset{x\in[-\pi,\pi]}{\underset{n\to\infty}{\rightrightarrows}}
f(x)
$$
и наша задача -- убедиться, что функции $f*\varDelta_n$ являются
тригонометрическими многочленами. Для этого нужно заметить, что на отрезке
$[-\pi,\pi]$ функции $\varDelta_n$ сами будут тригонометрическими многочленами:
$$
\varDelta_n(x)=c+\sum_{k=1}^N\Big\{ a_k\cos kx+b_k\sin kx\Big\}
$$
Отсюда следует цепочка:
 \begin{multline*}
f*\varDelta_n(x)=\int_{-\infty}^{+\infty}
f(y)\cdot\kern-15pt\underbrace{\varDelta_n(x-y)}_{\scriptsize\begin{matrix}\text{\rotatebox{90}{$=$}}\\
\phantom{,}0, \\
\text{при $y\notin[x-\pi,x+\pi]$,}\\ \text{поскольку}\\
\convsupp\varDelta_n=[-\pi,\pi]\end{matrix}}\kern-15pt\ \d
y=\int_{x-\pi}^{x+\pi}
f(y)\cdot\varDelta_n(\underbrace{x-y}_{\scriptsize\begin{matrix}\text{\rotatebox{90}{$\owns$}}\\
[-\pi,\pi]
\end{matrix}})\ \d y=\\= \int_{x-\pi}^{x+\pi} \underbrace{f(y)\cdot\l
c+\sum_{k=1}^N\Big\{a_k\cdot\cos k(x-y)+b_k(x)\cdot\sin
k(x-y)\Big\}\r}_{\text{$2\pi$-периодическая функция от $y$}}\ \d
y=\left|{\scriptsize\begin{matrix}\text{интеграл по сдвинутому периоду}\\
\text{равен интегралу по периоду}\end{matrix}}\right|=\\= \int_{-\pi}^{+\pi}
f(y)\cdot\l c+\sum_{k=1}^N\Big\{a_k\cdot\cos k(x-y)+b_k(x)\cdot\sin
k(x-y)\Big\}\r\ \d y=\eqref{cos(x-y)},\eqref{sin(x-y)}=\\= \int_{-\pi}^{+\pi}
f(y)\cdot\bigg(
C+\sum_{k=1}^N\Big\{\kern-5pt\underbrace{A_k(x)}_{\scriptsize\begin{matrix}\text{тригоно-}\\ \text{метрический}\\
\text{многочлен}\\
\text{от $x$}\end{matrix}}\kern-5pt\cdot\cos ky+\kern-5pt\underbrace{B_k(x)}_{\scriptsize\begin{matrix}\text{тригоно-}\\ \text{метрический}\\
\text{многочлен}\\
\text{от $x$}\end{matrix}}\kern-5pt\cdot\sin ky\Big\}\bigg)\ \d y=\\=
C\cdot\overbrace{\int_{-\pi}^{+\pi} f(y)\ \d
y}^{\text{число}}+\sum_{k=1}^N\Big\{\kern-5pt\underbrace{A_k(x)}_{\scriptsize\begin{matrix}\text{тригоно-}\\ \text{метрический}\\
\text{многочлен}\\
\text{от $x$}\end{matrix}}\kern-5pt\cdot\overbrace{\int_{-\pi}^{+\pi}
f(y)\cdot\cos ky\ \d y}^{\text{число}}+\kern-5pt\underbrace{B_k(x)}_{\scriptsize\begin{matrix}\text{тригоно-}\\ \text{метрический}\\
\text{многочлен}\\
\text{от $x$}\end{matrix}}\kern-5pt\cdot\overbrace{\int_{-\pi}^{+\pi}
f(y)\cdot\sin ky\ \d y}^{\text{число}}\Big\}
 \end{multline*}
\epr

\chapter{СТЕПЕННЫЕ РЯДЫ И АНАЛИТИЧЕСКИЕ ФУНКЦИИ}\label{CH-step-ryady}

 \bit{
\item[$\bullet$] {\it Степенным рядом}\index{ряд!степенной} называется
функциональный ряд следующего специального вида:
 \beq\label{DEF:step-ryad}
  \sum_{n=0}^\infty c_n\cdot (x-x_0)^n
 \eeq
Число $x_0$ называется {\it центром}\index{центр степенного ряда} этого
степенного ряда, а числа $c_n$ -- его {\it коэффициентами}\index{коэффициенты
степенного ряда}.

\item[$\bullet$] Частным случаем степенного ряда будет ряд, в котором
коэффициенты $c_n$ вычисляются по формулам
 \beq\label{coefficienty-Teilora}
c_0=\frac{f^{(n)}(x_0)}{n!}
 \eeq
где $f$ -- некоторая функция, бесконечно гладкая в окрестности точки $x_0$.
Такой степенной ряд называется {\it рядом Тейлора} функции $f$ в точке $x_0$.
Его коэффициенты \eqref{coefficienty-Teilora} называются {\it коэффициентами
Тейлора} функции $f$ в точке $x_0$, а частичные суммы
 \beq\label{mnogochlen-Taylor}
T_N(x)=\sum_{n=0}^\infty c_n\cdot (x-x_0)^n
 \eeq
-- {\it многочленами Тейлора} функции $f$ в точке $x_0$.
 }\eit

Среди всевозможных функций в математическом анализе важный класс образуют
функции, являющиеся суммами своего ряда Тейлора:
 $$
f(x)=\sum_{n=0}^\infty c_n\cdot (x-x_0)^n.
 $$
Такие функции называются аналитическими (точное определение см. на
с.\pageref{DEF:analitich-func}), и в этой главе мы поговорим о них.

\section{Степенные ряды}\label{power-series}

\subsection{Область сходимости степенного ряда}

\begin{tm}[\bf об области сходимости степенного ряда]\label{tm-21.1.1}
Область сходимости $D$ всякого степенного ряда
 \beq
  \sum_{n=0}^\infty c_n\cdot (x-x_0)^n
\label{21.1.1}
 \eeq
имеет следующий вид:
 \bit{
\item[---] либо $D$ состоит из одной точки: $D=\{ x_0 \}$; \item[---] либо $D$
совпадает со всей числовой прямой: $D=(-\infty,+\infty)$; \item[---] либо
существует такое число $R>0$, что $D$ совпадает с отрезком $[x_0-R,x_0+R]$,
исключая, может быть, точки на границе: $(x_0-R,x_0+R)\subseteq D\subseteq
[x_0-R,x_0+R]$.
 }\eit
\end{tm}

\bit{ \item[$\bullet$] Число $R$ при этом называется {\it радиусом
сходимости}\index{радиус сходимости степенного ряда}, а интервал
$(x_0-R,x_0+R)$ -- {\it интервалом сходимости}\index{интервал сходимости
степенного ряда} степенного ряда \eqref{21.1.1}.
 }\eit

Доказательство этой теоремы использует следующую лемму.

\begin{lm}[\bf Абель]\label{lm-Abel}
Пусть нам дан степенной ряд
 \beq\label{21.1.2}
  \sum_{n=0}^\infty c_n\cdot z^n.
 \eeq
Тогда:
 \bit{
\item[---] если ряд \eqref{21.1.2} сходится в какой-то точке $z=\zeta$, то он
сходится также в любой точке $z=\eta$ по модулю меньшей $\zeta$:
$$
|\eta|<|\zeta|
$$

\item[---] если ряд \eqref{21.1.2} расходится в какой-то точке $z=\eta$, то он
расходится также в любой точке $z=\zeta$ по модулю большей $\eta$:
$$
|\eta|<|\zeta|.
$$
 }\eit

\end{lm}\begin{proof} 1. Докажем сначала первую часть. Здесь используется некий стандартный прием,
постоянно применяемый при изучении степенных рядов, и называемый {\it
логической цепочкой Абеля}\index{Абеля!логическая
цепочка}\label{Abel-log-chain}. Он состоит в следующем.
$$
\text{Ряд \eqref{21.1.2} сходится в точке $z=\zeta$}
$$
$$
  \Downarrow
$$
$$
\text{числовой ряд}\quad \sum_{n=0}^\infty c_n\cdot \zeta^n \quad
\text{сходится}
$$
$$
  \Downarrow
$$
$$
\text{общий член стремится к нулю:}\quad c_n\cdot \zeta^n\underset{n\to
\infty}{\longrightarrow} 0
$$
$$
  \Downarrow
$$
 \beq
\text{последовательность $c_n\cdot \zeta^n$ ограничена:}\quad \exists M>0 \quad
\forall n\in \mathbb{N}\quad |c_n\cdot \zeta^n|\le M \label{21.1.3}
 \eeq
$$
  \Downarrow
$$
$$
\text{для всякого}\quad \eta: \, |\eta| < | \zeta | \quad \text{получаем:}
$$
$$
|c_n\cdot \eta^n|=\left|c_n\cdot \zeta^n\cdot \frac{\eta^n}{\zeta^n}\right|=
|c_n\cdot \zeta^n|\cdot \left|\frac{\eta}{\zeta}\right|^n \le \l
\text{применяем \eqref{21.1.3}}\r \le M\cdot \left|\frac{\eta}{\zeta}\right|^n,
\quad \text{где}\quad \left|\frac{\eta}{\zeta}\right|<1
$$
$$
  \Downarrow
$$
$$
\sum_{n=0}^\infty |c_n\cdot \eta^n|\le \sum_{n=0}^\infty M\cdot
\left|\frac{\eta}{\zeta}\right|^n= M\cdot \sum_{n=0}^\infty
\left|\frac{\eta}{\zeta}\right|^n \quad - \quad \text{сходится, поскольку}\quad
\left|\frac{\eta}{\zeta}\right|<1
$$
$$
  \Downarrow
$$
$$
\quad \sum\limits_{n=0}^\infty |c_n\cdot \eta^n| \quad - \quad \text{сходится,
по признаку сравнения (теорема \ref{tm-18.3.7})}
$$
$$
  \Downarrow
$$
$$
\quad \sum\limits_{n=0}^\infty c_n\cdot \eta^n \quad - \quad \text{сходится, по
признаку абсолютной сходимости (теорема \ref{tm-18.5.6})} $$

2. Для доказательства второй части достаточно переформулировать доказанное
утверждение так: {\it не бывает, чтобы $|\eta|<|\zeta|$, и при этом ряд
$\sum_{n=0}^\infty c_n\cdot \eta^n$ расходился, а ряд $\sum_{n=0}^\infty
c_n\cdot \zeta^n$ сходился.} Отсюда сразу получается импликация
$$
|\eta|<|\zeta|\quad\&\quad \text{ряд $\sum_{n=0}^\infty c_n\cdot \eta^n$
расходится}\quad\Longrightarrow\quad \text{ряд $\sum_{n=0}^\infty c_n\cdot
\zeta^n$ расходится}
$$

\end{proof}

\begin{proof}[Доказательство теоремы \ref{tm-21.1.1}]
 Заметим сразу, что нам достаточно рассмотреть ряд
\eqref{21.1.2}, потому что ряд \eqref{21.1.1} превращается в него заменой
переменной
$$
  x-x_0=z
$$
Обозначим через $D$ его область сходимости. Нам нужно доказать, что либо $D=\{
0\}$, либо $D=\R$, либо существует такое число $R>0$, что $(-R,R)\subseteq
D\subseteq [-R,R]$.

Ясно, что $D\ne\varnothing$, поскольку $0\in D$. Кроме того, из леммы Абеля
\ref{lm-Abel} следует утверждение:
 \beq
\text{если}\quad z\in \R\quad \text{и}\quad \exists \zeta\in D: \quad
|z|<|\zeta| \quad \text{то}\quad z\in D \label{21.1.4}
 \eeq
А из следствия 1.3 -- утверждение
 \beq
\text{если}\quad  z\in \R\quad \text{и}\quad \exists \eta\notin D: \quad
|\eta|<|z| \quad \text{то}\quad z\notin D \label{21.1.5}
 \eeq
Положим
$$
  R=\sup \{ |\zeta|; \zeta\in D \}
$$
и рассмотрим несколько возможных случаев.

1. Если $R=0$, то это означает, что область сходимости состоит из одной точки:
$D=\{ 0 \}$.

2. Если $R=\infty$, то это означает, что для всякой точки $z\in \R$ найдется
$\zeta\in D$ такое, что $|z|<|\zeta|$. Поэтому, в силу \eqref{21.1.4}, $z$
принадлежит $D$. Это верно для любой $z\in \R$, поэтому $D=\R$.

3. Пусть $0<R<\infty$, тогда мы получим
 \bit{
\item[--] если $z\in (-R,R)$, то есть $|z|<R$, то это означает, что найдется
$\zeta\in D$ такое, что $|z|<|\zeta|$, поэтому, в силу утверждения
\eqref{21.1.4}, $z\in D$; \item[--] если же $z\notin [-R,R]$, то есть $|z|>R$,
то это означает, что найдется $\eta\notin D$ такое, что $|z|>|\eta|$, поэтому,
в силу утверждения \eqref{21.1.5}, $z\notin D$;
 }\eit
Таким образом, мы получаем, что $D$ содержит интервал $(-R,R)$ и не содержит
точек, не лежащих в отрезке $[-R,R]$. То есть,
$$
(-R,R)\subseteq D\subseteq [-R,R]
$$
Это нам и нужно было доказать. \end{proof}

\begin{tm}[\bf о радиусе сходимости степенного ряда]\label{tm-21.1.4}
Радиус сходимости $R$ произвольного степенного ряда
$$
  \sum_{n=0}^\infty c_n\cdot (x-x_0)^n
$$
можно вычислить по формулам
$$
R=\lim_{n\to\infty}\left|\frac{c_n}{c_{n+1}}\right|=
\lim_{n\to\infty}\frac{1}{\sqrt[n] {|c_n|}}
$$
(если такие пределы существуют). При этом, равенство
$$
  R=0
$$
означает, что область сходимости $D$ этого ряда имеет вид $D=\{ x_0 \}$, а
равенство
$$
  R=\infty
$$
-- что $D=\R$.
\end{tm}\begin{proof} Достаточно рассмотреть ряды с нулевым центром, то есть рядов вида \eqref{21.1.2}.
Мы докажем лишь первую формулу, сказав про вторую лишь, что она доказывается
аналогично. Пусть
$$
R=\lim_{n\to\infty}\left|\frac{c_n}{c_{n+1}}\right|,
$$
причем предел справа существует и конечен. Тогда
 \bit{
\item[---] если $|z|<R$, то ряд $\sum\limits_{n=0}^\infty |c_n\cdot z^n|$ будет
сходиться по признаку Даламбера (теорема \ref{tm-18.3.20}), потому что число
Даламбера будет меньше единицы:
$$
D= \lim_{n\to\infty}\frac{\ml c_{n+1}\cdot z^{n+1}\mr}{\ml c_n\cdot z^n\mr}=
\lim_{n\to\infty}\frac{\ml c_{n+1}\mr}{\ml c_n\mr}\cdot |z|=
\frac{|z|}{\lim_{n\to\infty}\frac{\ml c_n\mr}{\ml c_{n+1}\mr}}=
\frac{|z|}{R}<1;
$$
значит, ряд $\sum\limits_{n=0}^\infty c_n\cdot z^n$ тоже должен сходиться;

\item[---] если
$|z|>R=\lim\limits_{n\to\infty}\left|\frac{c_n}{c_{n+1}}\right|$, то есть
$\lim\limits_{n\to\infty}\left|\frac{c_{n+1}\cdot z}{c_n}\right|>1$, то это
означает, что существует такое число $C>1$, что, для некоторого номера $N$
выполняется следующее:
$$
\forall n\ge N \quad \left|\frac{c_{n+1} z}{c_n}\right|\ge C
$$
$$
\Downarrow
$$
$$
\forall n\ge N \quad |c_{n+1} z|\ge C |c_n|
$$
$$
\Downarrow
$$
$$
\forall n\ge N \quad \ml c_{n+1} z^{n+1}\mr\ge C \ml c_n z^n\mr
$$
то есть, последовательность $\ml c_{n+1} z^{n+1}\mr$ должна быть неубывающей;
отсюда следует, что
$$
c_n z^n \underset{n\to\infty}{\not\longrightarrow} 0,
$$
значит, ряд $\sum\limits_{n=0}^\infty c_{n+1} z^{n+1}$ не может сходиться.

Мы доказали, что
$$
(-R,R)\subseteq D\subseteq [-R,R]
$$
\item[] и поэтому $R$ действительно должен быть радиусом сходимости ряда
\eqref{21.1.2}.
 }\eit

Нам нужно только еще рассмотреть случай $R=\lim\limits_{n\to\infty}\frac{\ml
c_n\mr}{\ml c_{n+1}\mr}=\infty$. Тогда мы получим для любого $z\in \R$, то ряд
$\sum\limits_{n=0}^\infty |c_n\cdot z^n|$ снова будет сходиться по признаку
Даламбера:
$$
D= \lim_{n\to\infty}\frac{\ml c_{n+1}\cdot z^{n+1}\mr}{\ml c_n\cdot z^n\mr}=
\lim_{n\to\infty}\frac{\ml c_{n+1}\mr}{\ml c_n\mr}\cdot |z|=
\frac{|z|}{\lim\limits_{n\to\infty}\frac{\ml c_n\mr}{\ml c_{n+1}\mr}}= \l
\frac{|z|}{\infty}\r=0
$$
поэтому ряд $\sum_{n=0}^\infty c_n\cdot z^n$ тоже сходится, и таким образом,
$D=\R$.  \end{proof}

Перейдем, наконец, к примерам.

\noindent\rule{160mm}{0.1pt}\begin{multicols}{2}

\begin{ex}\label{ex-21.1.5} Рассмотрим ряд
$$
  \sum_{n=1}^\infty n!\cdot x^n
$$
Радиус сходимости в этом случае оказывается равным нулю
$$
R=\lim_{n\to\infty}\left|\frac{c_n}{c_{n+1}}\right|=
\lim_{n\to\infty}\frac{n!}{(n+1)!}= \lim_{n\to\infty}\frac{1}{n+1}=0
$$
поэтому область сходимости будет состоять из одной точки -- центра нашего ряда
$x_0=0$.

Вывод: Область сходимости $D=\{ 0\}$.
\end{ex}

\begin{ex}\label{ex-21.1.6} Рассмотрим ряд
$$
  \sum_{n=1}^\infty \frac{x^n}{n!}
$$
Радиус сходимости в этом случае оказывается равным бесконечности
$$
R=\lim_{n\to\infty}\left|\frac{c_n}{c_{n+1}}\right|=
\lim_{n\to\infty}\frac{(n+1)!}{n!}= \lim_{n\to\infty} (n+1)=\infty
$$
поэтому область сходимости будет совпадать со всей прямой $\R$.

Вывод: Область сходимости $D=\R$.
\end{ex}

\begin{ex}\label{ex-21.1.7} Рассмотрим ряд
$$
  \sum_{n=1}^\infty \frac{x^n}{n}
$$
Вычисляем радиус сходимости:
$$
R=\lim_{n\to\infty}\left|\frac{c_n}{c_{n+1}}\right|=
\lim_{n\to\infty}\frac{n+1}{n}=1
$$
Это значит, что область сходимости $D$ является отрезком $[-1,1]$, кроме, может
быть, точек на границе. Наглядно это удобно показать следующей картинкой:

%\pucture{0pt}{0pt}{ii-17.pcx}

\vglue80pt \noindent Нам остается лишь проверить точки $x=-1$ и $x=1$.

1) При $x=-1$ получаем ряд
$$
  \sum_{n=1}^\infty \frac{(-1)^n}{n}
$$
который сходится по признаку Лейбница (теорема \ref{tm-18.6.1}).

2) При $x=1$ получаем гармонический ряд
$$
  \sum_{n=1}^\infty \frac{1^n}{n}=\sum_{n=1}^\infty \frac{1}{n}
$$
который расходится, в силу примера \ref{ex-18.3.2}.

Таким образом, нашу картинку можно поправить следующим образом:

%\pucture{0pt}{0pt}{ii-18.pcx}

\vglue80pt

Вывод: Область сходимости $D=[-1,1)$.
\end{ex}

\begin{ex}\label{ex-21.1.8} Рассмотрим ряд с ненулевым центром:
$$
  \sum_{n=1}^\infty \frac{(x+1)^n}{2^n}
$$
Радиус сходимости:
$$
R=\lim_{n\to\infty}\left|\frac{c_n}{c_{n+1}}\right|=
\lim_{n\to\infty}\frac{2^{n+1}}{2^n}=2
$$
Картинка:

%\pucture{0pt}{0pt}{ii-19.pcx}

\vglue80pt \noindent Проверяем точки $x=-3$ и $x=1$.

1) При $x=-3$ получаем ряд
$$
  \sum_{n=1}^\infty \frac{(-2)^n}{2^n}=\sum_{n=1}^\infty (-1)^n
$$
расходящийся по необходимому условию сходимости (теорема \ref{tm-18.5.1}).

2) При $x=1$ получаем ряд
$$
  \sum_{n=1}^\infty \frac{2^n}{2^n}=\sum_{n=1}^\infty 1
$$
также расходящийся по необходимому условию сходимости.

Поправленная картинка:

%\pucture{0pt}{0pt}{ii-20.pcx}

\vglue80pt

Вывод: Область сходимости $D=(-3,1)$.
\end{ex}

\begin{ers} Найдите область сходимости:
 \begin{multicols}{2}
1. $\sum\limits_{n=1}^\infty \frac{(x-2)^n}{n^2}$;

2. $\sum\limits_{n=1}^\infty n^2\cdot (x-2)^n$;

3. $\sum\limits_{n=1}^\infty 3^n\cdot (x+2)^n$;

4. $\sum\limits_{n=1}^\infty n^n\cdot x^n$;

5. $\sum\limits_{n=1}^\infty \frac{x^n}{n^n}$.
\end{multicols}\end{ers}

\end{multicols}\noindent\rule[10pt]{160mm}{0.1pt}

\subsection{Равномерная сходимость степенного ряда и ее следствия}

\paragraph{Равномерная сходимость степенного ряда.}

\begin{tm}\label{tm-21.2.1} Всякий степенной ряд
$$
  \sum_{n=0}^\infty c_n\cdot (x-x_0)^n
$$
сходится равномерно на каждом отрезке $[a,b]$ внутри интервала сходимости
$(x_0-R, x_0+R)$.
\end{tm}

%\pucture{0pt}{0pt}{ii-21.pcx}

\vglue80pt \begin{proof} Поскольку $[a,b]\subseteq (x_0-R, x_0+R)$, можно найти
число $C>0$ такое, что
$$
x_0-R<x_0-C<a<b<x_0+C<x_0+R
$$

%\pucture{0pt}{0pt}{ii-22.pcx}

\vglue80pt \noindent Зафиксируем его и заметим, что
 \beq
\sup_{x\in [a,b]}\frac{|x-x_0|}{C}= \max \lll \frac{|a-x_0|}{C};
\frac{|b-x_0|}{C}\rrr=\lambda<1 \label{21.2.1}
 \eeq
Наш ряд сходится в точке $x_0+C\in (x_0-R, x_0+R)$, и это приводит к логической
цепочке Абеля (описывавшейся при доказательстве леммы \ref{lm-Abel}):
$$
\sum_{n=0}^\infty c_n\cdot C^n \quad - \quad \text{сходится}
$$
$$
  \Downarrow
$$
$$
c_n\cdot C^n \underset{n\to \infty}{\longrightarrow} 0
$$
$$
  \Downarrow
$$
$$
\text{последовательность}\quad  c_n\cdot C^n  \quad \text{ограничена}:
$$
 \beq
\exists M>0 \quad \forall n \quad |c_n\cdot C^n|\le M \label{21.2.2}
 \eeq
$$
  \Downarrow
$$
\begin{multline*}\ml\ml c_n\cdot (x-x_0)^n \mr\mr_{x\in [a,b]}= \sup_{x\in
[a,b]}\ml c_n\cdot (x-x_0)^n \mr= \sup_{x\in [a,b]}\ml c_n\cdot C^n\cdot
\frac{(x-x_0)^n}{C^n}\mr=\\= \ml c_n\cdot C^n\mr\cdot \sup_{x\in
[a,b]}\ml\frac{(x-x_0)^n}{C^n}\mr\le (\text{применяем \eqref{21.2.2}})\le \\
\le M\cdot \l \sup_{x\in [a,b]}\frac{|x-x_0|}{C}\r^n \le (\text{применяем
\eqref{21.2.1}})\le M\cdot \lambda^n \quad (\text{где}\quad \lambda<1)
\end{multline*}
$$
  \Downarrow
$$
$$
\sum_{n=0}^\infty \ml\ml c_n\cdot (x-x_0)^n \mr\mr_{x\in [a,b]}\le
\sum_{n=0}^\infty M\cdot \lambda^n = M\cdot \sum_{n=0}^\infty \lambda^n \quad
(\text{где}\quad \lambda<1)
$$
$$
  \Downarrow
$$
$$
\sum_{n=0}^\infty \ml\ml c_n\cdot (x-x_0)^n \mr\mr_{x\in [a,b]}\quad - \quad
\text{сходится}
$$
$$
  \Downarrow \put(20,0){\smsize\text{$\begin{pmatrix}
\text{вспоминаем признак Вейерштрасса}\\ \text{равномерной сходимости (теорему
\ref{tm-20.5.4})}\end{pmatrix}$}}
$$
$$
\sum_{n=0}^\infty c_n\cdot (x-x_0)^n  \quad - \quad \text{сходится равномерно
на отрезке}\quad [a,b] \quad $$ \end{proof}

\paragraph{Непрерывность суммы степенного ряда.}

\begin{tm}\label{tm-21.2.2} Сумма всякого степенного ряда
$$
S(x)=\sum_{n=0}^\infty c_n\cdot (x-x_0)^n
$$
непрерывна на интервале сходимости $(x_0-R, x_0+R)$.
\end{tm}\begin{proof} По теореме \ref{tm-21.2.1}, этот ряд сходится
равномерно на любом отрезке $[a,b]\subseteq (x_0-R, x_0+R)$. С другой стороны,
он состоит их непрерывных функций $c_n\cdot (x-x_0)^n$, поэтому, по свойству
$1^0 \S 3$ главы 20, сумма $S(x)$ должна быть непрерывна на $[a,b]$. Это верно
для любого отрезка $[a,b]\subseteq (x_0-R, x_0+R)$, поэтому сумма $S(x)$ должна
быть непрерывна на всем интервале $(x_0-R, x_0+R)$. \end{proof}

\paragraph{Дифференцирование и интегрирование степенных рядов.}

Зафиксируем произвольный степенной ряд
 \beq
\sum_{n=0}^\infty c_n\cdot (x-x_0)^n \label{21.2.3}
 \eeq
и выпишем ряд, получающийся его почленным дифференцированием
 \beq
\sum_{n=0}^\infty \lll c_n\cdot (x-x_0)^n \rrr'= \sum_{n=1}^\infty c_n\cdot
n\cdot (x-x_0)^{n-1}\label{21.2.4}
 \eeq
и ряд, получающийся  почленным интегрированием:
 \beq
\sum_{n=0}^\infty \int_{x_0}^x c_n\cdot (t-x_0)^n \, \d t= \sum_{n=0}^\infty
\frac{c_n}{n+1}\cdot (x-x_0)^{n+1}\label{21.2.5}
 \eeq

\begin{tm}\label{TH:diff+int-step-ryada} Радиусы
(и интервалы) сходимости рядов \eqref{21.2.3}, \eqref{21.2.4} и \eqref{21.2.5}
совпадают, и
 \bit{
\item[(i)] сумма ряда \eqref{21.2.3} является гладкой функцией на интервале
сходимости, и ее производную на нем можно вычислить почленным
дифференцированием:
 \beq\label{diff-step-ryada}
\frac{\d}{\d x}\l\sum_{n=0}^\infty c_n\cdot (x-x_0)^n\r=\sum_{n=0}^\infty
\frac{\d}{\d x}\Big(c_n\cdot (x-x_0)^n \Big)= \sum_{n=1}^\infty c_n\cdot n\cdot
(x-x_0)^{n-1}
 \eeq

\item[(ii)] интеграл от суммы ряда \eqref{21.2.3} можно вычислить почленным
интегрированием:
 \beq\label{int-step-ryada}
\int_{x_0}^x \l\sum_{n=0}^\infty c_n\cdot (t-x_0)^n\r\d t=\sum_{n=0}^\infty
\int_{x_0}^x c_n\cdot (t-x_0)^n \, \d t= \sum_{n=0}^\infty \frac{c_n}{n+1}\cdot
(x-x_0)^{n+1}
 \eeq
 }\eit
 \end{tm}
\begin{proof}

1. Убедимся сначала, что радиусы сходимости рядов \eqref{21.2.3},
\eqref{21.2.4} и \eqref{21.2.5} совпадают. Заметим прежде всего, что для этого
нам достаточно проверить совпадение радиусов сходимости для рядов
\eqref{21.2.3} и \eqref{21.2.4}. Это будет означать, что при дифференцировании
степенного ряда его радиус сходимости не меняется, значит, поскольку ряд
\eqref{21.2.3} тоже получается из ряда \eqref{21.2.5} почленным
дифференцированием, их радиусы сходимости тоже будут совпадать.

Сделаем после этого замену переменной
$$
  x-x_0=z
$$
Тогда ряды \eqref{21.2.3} и \eqref{21.2.4} перепишутся следующим образом:
 \beq
\sum_{n=0}^\infty c_n\cdot z^n \label{21.2.6}
 \eeq
 \beq
\sum_{n=1}^\infty c_n\cdot n\cdot z^{n-1}\label{21.2.7}
 \eeq
Обозначим через $R_1$ радиус сходимости ряда \eqref{21.2.6}, а через $R_2$ --
радиус сходимости ряда \eqref{21.2.7}.

A) Докажем сначала, что $R_1\le R_2$. Это равносильно тому, что если ряд
$\sum\limits_{n=0}^\infty c_n\cdot \zeta^n$ сходится при каком-нибудь $\zeta$,
то при любом $z: \, 0<|z|<|\zeta|$ сходится ряд $\sum\limits_{n=1}^\infty
c_n\cdot n\cdot z^{n-1}$. Применяем логическую цепочку Абеля:
$$
\sum\limits_{n=0}^\infty c_n\cdot \zeta^n \quad - \quad \text{сходится}
$$
$$
\Downarrow
$$
$$
c_n\cdot \zeta^n \underset{n\to \infty}{\longrightarrow} 0
$$
$$
\Downarrow
$$
$$
\exists M>0 \quad \forall n \quad |c_n\cdot \zeta^n|\le M
$$
$$
\Downarrow
$$
$$
\text{для всякого}\quad z: \, 0<|z|<|\zeta|
$$
$$
|c_n\cdot n\cdot z^{n-1}|= \left|c_n\cdot \zeta^n \cdot \frac{n}{z}\cdot
\frac{z^n}{\zeta^n}\right|= |c_n\cdot \zeta^n| \cdot \frac{n}{|z|}\cdot
\left|\frac{z}{\zeta}\right|^n\le M \cdot \frac{n}{|z|}\cdot
\left|\frac{z}{\zeta}\right|^n
$$
$$
\Downarrow
$$
$$
\sum\limits_{n=1}^\infty |c_n\cdot n\cdot z^{n-1}|\le \sum\limits_{n=1}^\infty
M \cdot \frac{n}{|z|}\cdot \left|\frac{z}{\zeta}\right|^n= M \cdot
\sum\limits_{n=1}^\infty \frac{n}{|z|}\cdot \left|\frac{z}{\zeta}\right|^n
$$
$$
- \quad \text{сходится по признаку Даламбера, поскольку}\quad
\left|\frac{z}{\zeta}\right|<1
$$
$$
\Downarrow
$$
$$
\sum\limits_{n=1}^\infty |c_n\cdot n\cdot z^{n-1}| \quad - \quad \text{сходится
по признаку сравнения}
$$
$$
\Downarrow
$$
$$
\sum\limits_{n=1}^\infty c_n\cdot n\cdot z^{n-1}\quad - \quad \text{сходится по
признаку абсолютной сходимости}
$$

B) Теперь докажем, что $R_2\le R_1$. Это равносильно тому, что если ряд
$\sum\limits_{n=1}^\infty c_n\cdot n\cdot \zeta^{n-1}$ сходится при
каком-нибудь $\zeta$, то при любом $z: \, 0<|z|<|\zeta|$ сходится ряд
$\sum\limits_{n=0}^\infty c_n\cdot z^n$. Применяем логическую цепочку Абеля:
$$
\sum\limits_{n=1}^\infty c_n\cdot n\cdot \zeta^{n-1}\quad - \quad
\text{сходится}
$$
$$
\Downarrow
$$
$$
c_n\cdot n\cdot \zeta^{n-1}\underset{n\to \infty}{\longrightarrow} 0
$$
$$
\Downarrow
$$
$$
\exists M>0 \quad \forall n \quad |c_n\cdot n\cdot \zeta^{n-1}|\le M
$$
$$
\Downarrow
$$
$$
\text{для всякого}\quad z: \, 0<|z|<|\zeta|
$$
$$
|c_n\cdot z^n|= \ml c_n\cdot n\cdot \zeta^{n-1}\cdot \frac{\zeta}{n}\cdot
\frac{z^n}{\zeta^n}\mr= |c_n\cdot n\cdot \zeta^{n-1}| \cdot
\frac{|\zeta|}{n}\cdot \ml \frac{z}{\zeta}\mr^n\le M \cdot
\frac{|\zeta|}{n}\cdot \ml \frac{z}{\zeta}\mr^n
$$
$$
\Downarrow
$$
$$
\sum\limits_{n=0}^\infty |c_n\cdot z^n|\le \sum\limits_{n=0}^\infty M \cdot
\frac{|\zeta|}{n}\cdot \ml \frac{z}{\zeta}\mr^n= M \cdot
\sum\limits_{n=0}^\infty \frac{|\zeta|}{n}\cdot \ml \frac{z}{\zeta}\mr^n
$$
$$
- \quad \text{сходится по признаку Даламбера, поскольку}\quad
\left|\frac{z}{\zeta}\right|<1
$$
$$
\Downarrow
$$
$$
\sum\limits_{n=0}^\infty |c_n\cdot z^n| \quad - \quad \text{сходится по
признаку сравнения}
$$
$$
\Downarrow
$$
$$
\sum\limits_{n=0}^\infty c_n\cdot z^n \quad - \quad \text{сходится по признаку
абсолютной сходимости}
$$

Итак, мы доказали, что
$$
  R_1\le R_2 \le R_1
$$
то есть
$$
  R_1=R_2
$$

2. Докажем далее тождество \eqref{diff-step-ryada}. Возьмем какие-нибудь точки
$a,b$ так чтобы
$$
  x_0-R<a<x_0<b<x_0+R
$$
Тогда, по теореме \ref{tm-21.2.1}, мы получим, что ряд \eqref{21.2.4}
$$
\sum_{n=1}^\infty c_n\cdot n\cdot (x-x_0)^{n-1}= \sum_{n=0}^\infty \l c_n\cdot
n\cdot (x-x_0)^n \r'
$$
сходится равномерно на отрезке $[a,b]$, а, с другой стороны, ряд
$$
\sum_{n=0}^\infty c_n\cdot (x-x_0)^n
$$
сходится в точке $x_0$, потому что все его члены, кроме нулевого равны нулю в
этой точке:
$$
\sum_{n=0}^\infty c_n\cdot n\cdot (x-x_0)^n \Big|_{x=x_0}= c_0\cdot 1+c_1\cdot
0+c_2\cdot 0+...
$$
По свойству $4^0$ на с.\pageref{PROP:diff-func-ryada}, это означает, что сумма
ряда \eqref{21.2.3} будет дифференцируема на $[a,b]$, и ее производную в каждой
точке $x\in [a,b]$ можно вычислить почленным дифференцированием:
$$
\frac{\d}{\d x}\l\sum_{n=0}^\infty c_n\cdot (x-x_0)^n\r=\sum_{n=0}^\infty
\frac{\d}{\d x}\Big(c_n\cdot (x-x_0)^n \Big)= \sum_{n=1}^\infty c_n\cdot n\cdot
(x-x_0)^{n-1}, \quad x\in [a,b]
$$
То есть тождество \eqref{diff-step-ryada} выполняется для $x\in [a,b]$.
Поскольку числа $a,b$ здесь выбирались произвольными, удовлетворяющими условию
$x_0-R<a<x_0<b<x_0+R$, мы получаем, что это тождество выполняется при любом
$x\in (x_0-R,x_0+R)$.

3. Теперь докажем \eqref{int-step-ryada}. Зафиксируем $x$ из интервала
сходимости (общего для рядов \eqref{21.2.3} и \eqref{21.2.5}). По теореме
\ref{tm-21.2.1}, мы получим, что ряд
$$
\sum_{n=0}^\infty c_n\cdot n\cdot (t-x_0)^n
$$
сходится равномерно на отрезке $[x_0,x]$. Поэтому, по свойству $3^0$ на
с.\pageref{PROP:int-func-ryada}, должно выполняться равенство
$$
\int_{x_0}^x \sum_{n=0}^\infty c_n\cdot (t-x_0)^n \, \d t= \sum_{n=0}^\infty
\int_{x_0}^x c_n\cdot (t-x_0)^n\, \d t= \sum_{n=0}^\infty \frac{c_n}{n+1}\cdot
(x-x_0)^{n+1}
$$
То есть справедливо \eqref{int-step-ryada}.
\end{proof}

\begin{cor}\label{cor-21.2.5}
Сумма всякого степенного ряда является бесконечно гладкой функцией на интервале
сходимости $(x_0-R, x_0+R)$ этого ряда, и производная порядка $k$ этой функции
вычисляется по формуле
 \beq
\frac{\d^k}{\d x^k}\l\sum_{n=0}^\infty c_n\cdot (x-x_0)^n\r=\sum_{n=k}^\infty
\frac{n!}{(n-k)!}\cdot c_n \cdot (x-x_0)^{n-k}, \qquad x\in (x_0-R, x_0+R)
\label{21.2.12}
 \eeq
\end{cor}\begin{proof} Обозначим для удобства сумму нашего степенного ряда буквой
$P$:
$$
P(x)=\sum_{n=0}^\infty c_n\cdot (x-x_0)^n, \qquad x\in (x_0-R, x_0+R)
$$
По теореме \ref{tm-21.2.4}, $P$ дифференцируема на $(x_0-R, x_0+R)$, причем ее
производная является суммой степенного ряда с тем же интервалом сходимости:
$$
P'(x)=\sum_{n=1}^\infty c_n\cdot n\cdot (x-x_0)^{n-1}, \qquad x\in (x_0-R,
x_0+R)
$$
Значит, опять по теореме \ref{tm-21.2.4}, $P'$ дифференцируема на $(x_0-R,
x_0+R)$, причем ее производная является суммой степенного ряда с тем же
интервалом сходимости:
$$
P''(x)=\sum_{n=2}^\infty c_n\cdot n\cdot (n-1)\cdot (x-x_0)^{n-2}, \qquad x\in
(x_0-R, x_0+R)
$$
И так далее. Индукцией по $k$ получаем формулу \eqref{21.2.12}.
\end{proof}

\noindent\rule{160mm}{0.1pt}\begin{multicols}{2}

\subsection{Вычисление суммы степенного ряда}

Теорема \ref{TH:diff+int-step-ryada} позволяет в некоторых случаях явно
вычислить сумму степенного ряда. Здесь мы покажем на примерах, как это
делается.

Все начинается с формулы для суммы членов бесконечной геометрической
прогрессии, которую мы выписывали в примере \ref{ex-18.1.4}:
 \beq\label{21.3.11}
\boxed{\quad\sum_{n=0}^\infty x^n=\frac{1}{1-x},\quad |x|<1}
 \eeq
Из нее выводится целый ряд других полезных формул для сумм степенных рядов.
Прежде всего, заменой $x$ на $-x$ мы получаем тождество
 \beq\label{21.3.10}
\boxed{\quad\sum_{n=0}^\infty (-1)^n \cdot x^n=\frac{1}{1+x},\quad  |x|<1}
 \eeq
Интегрируя его при $|x|<1$ мы получаем
 \begin{multline*}
\sum_{n=1}^\infty (-1)^{n-1}\cdot \frac{x^n}{n}=\sum_{n=0}^\infty (-1)^n\cdot
\frac{x^{n+1}}{n+1}=\\=\sum_{n=0}^\infty \int_0^x (-1)^n \cdot t^n \d t=
\ln(1+x)=\int_0^x \frac{\d t}{1+t}
 \end{multline*}
То есть,
 \beq\label{21.3.7}
\boxed{\quad\sum_{n=1}^\infty (-1)^{n-1}\cdot \frac{x^n}{n}=\ln (1+x),\quad
|x|<1}
 \eeq
Далее, заменив в \eqref{21.3.10} $x$ на $x^2$, мы получим
 \beq\label{21.3.10-1}
\boxed{\quad\sum_{n=0}^\infty (-1)^n \cdot x^{2n}=\frac{1}{1+x^2},\quad |x|<1}
 \eeq
Интегрируя \eqref{21.3.10-1} по теореме \ref{TH:diff+int-step-ryada} (при
$|x|<1$) мы получаем:
 \begin{multline*}
\sum_{n=0}^\infty (-1)^n\cdot \frac{x^{2n+1}}{2n+1}=\sum_{n=0}^\infty \int_0^x
(-1)^n \cdot t^{2n} \d t=\\=\int_0^x \frac{\d t}{1+t^2} =\arctg x
 \end{multline*}
То есть
 \beq\label{21.3.12}
\boxed{\quad\sum_{n=0}^\infty (-1)^n\cdot \frac{x^{2n+1}}{2n+1}=\arctg x,\quad
|x|<1}
 \eeq

Полученных формул уже достаточно, чтобы перейти к решению задач.

\begin{ex}\label{ex-21.4.1} Предположим, нам нужно найти сумму степенного ряда
 \beq
  \sum_{n=1}^\infty \frac{x^n}{n}\label{21.4.1}
 \eeq
Рецепт решения заключается в том, чтобы сообразить, как этот ряд можно получить
из формул \eqref{21.3.11}-\eqref{21.3.12} с помощью операций дифференцирования,
интегрирования, замены переменной и умножения на $x^k$. Для данного ряда
последовательность действий выглядит следующим образом. Сначала мы выписываем
формулу \eqref{21.3.11}:
$$
\sum_{n=0}^\infty x^n=\frac{1}{1-x}
$$
Затем меняем в ней индекс $n$ на $k-1$:
$$
\sum_{k=1}^\infty x^{k-1}=\frac{1}{1-x}
$$
После этого интегрируем:
$$
\sum_{k=1}^\infty \frac{x^k}{k}=\int_0^x \sum_{k=1}^\infty t^{k-1}\, \d
t=\int_0^x \frac{1}{1-t}\, \d t=-\ln(1-x)
$$
Остается увидеть, что полученный ряд отличается от \eqref{21.4.1} только
индексом. Если теперь поменять $k$ на $n$, то получится

Ответ:
$$
\sum_{n=1}^\infty \frac{x^n}{n}=-\ln(1-x)
$$
\end{ex}

\begin{ex}\label{ex-21.4.2}  Найти сумму степенного ряда
$$
  \sum_{n=1}^\infty n \cdot x^n
$$
Снова выписываем формулу \eqref{21.3.11}:
$$
\sum_{n=0}^\infty x^n=\frac{1}{1-x}
$$
Дифференцируем ее:
$$
\lll \sum_{n=0}^\infty x^n\rrr'=\lll\frac{1}{1-x}\rrr'
$$
$$
\Downarrow
$$
$$
\sum_{n=1}^\infty n\cdot x^{n-1}=\frac{1}{(1-x)^2}
$$
Умножаем на $x$:
$$
x\cdot \sum_{n=1}^\infty n\cdot x^{n-1}=\frac{x}{(1-x)^2}
$$
$$
\Downarrow
$$

Ответ:
$$
\sum_{n=1}^\infty n\cdot x^n=\frac{x}{(1-x)^2}
$$
\end{ex}

\begin{ex}\label{ex-21.4.3} Найти сумму степенного ряда
$$
  \sum_{n=1}^\infty (-1)^n \cdot \frac{x^n}{2n+1}
$$
Здесь вычисления начинаются с формулы \eqref{21.3.12}, которую мы перепишем с
переменной $t$ вместо $x$:
$$
\sum_{n=0}^\infty (-1)^n\cdot \frac{t^{2n+1}}{2n+1}=\arctg t
$$
Мы делим ее на $t$
$$
\frac{1}{t}\cdot \sum_{n=0}^\infty (-1)^n\cdot \frac{t^{2n+1}}{2n+1}=
\frac{\arctg t}{t}
$$
$$
\Downarrow
$$
$$
\sum_{n=0}^\infty (-1)^n\cdot \frac{t^{2n}}{2n+1}= \frac{\arctg t}{t}
$$
Потом делаем замену переменной $t^2=x$ и получается

Ответ:
$$
\sum_{n=0}^\infty (-1)^n\cdot \frac{x^n}{2n+1}= \frac{\arctg
\sqrt{x}}{\sqrt{x}}
$$
\end{ex}

\begin{ers} Найдите сумму ряда:

1. $\sum\limits_{n=1}^\infty n^2 \cdot x^n$;

2. $\sum\limits_{n=1}^\infty 2^n \cdot x^n$;

3. $\sum\limits_{n=1}^\infty \frac{x^{3n}}{n}$;

4. $\sum\limits_{n=1}^\infty \frac{n}{n+1}\cdot x^n$;

5. $\sum\limits_{n=1}^\infty \frac{x^n(1+x)}{n}$;

6. $\sum\limits_{n=1}^\infty n\cdot x^n(1+x)$;

7. $\sum\limits_{n=1}^\infty \frac{x^n(1+x^n)}{n}$;

8. $\sum\limits_{n=1}^\infty n\cdot x^n(1+x^n)$.
 \end{ers}

\end{multicols}\noindent\rule[10pt]{160mm}{0.1pt}

\section{Аналитические последовательности и производящие функции}

\subsection{Определение аналитической последовательности}

 \bit{
\item[$\bullet$] Числовая последовательность $\{a_n;\ n\in\Z_+\}$ называется
{\it аналитической}, если она удовлетворяет следующим эквивалентным условиям:
 \bit{
\item[(i)] степенной ряд $\sum_{n=0}^\infty a_n\cdot x^n$ имеет ненулевой
радиус сходимости:
 \beq\label{DEF:anal-posled-1}
 \exists r>0\qquad \text{ряд $\sum_{n=0}^\infty a_n\cdot r^n$
сходится}
 \eeq

\item[(ii)] для некоторого числа $A>0$ последовательность $A^n$ мажорирует
последовательность $|a_n|$:
 \beq\label{DEF:anal-posled-2}
\forall n\in\Z_+\qquad |a_n|\le A^n
 \eeq
 }\eit
Множество всех аналитических последовательностей обозначается символом
$\mathcal{A}$.

\item[$\bullet$] Радиус сходимости степенного ряда $\sum_{n=0}^\infty a_n\cdot
x^n$ называется {\it радиусом сходимости последовательности} $a$ и обозначается
$\rho(a)$:
$$
\rho(a)=\sup\left\{r>0:\quad \text{ряд $\sum_{n=0}^\infty a_n\cdot r^n$
сходится} \right\}
$$

 }\eit

\noindent\rule{160mm}{0.1pt}\begin{multicols}{2}

\begin{ex}\label{ex-21.1.5-sec} Последовательность
$$
a_n=n!
$$
не будет аналитической, потому что в силу примера \ref{ex-21.1.5}, порождаемый
ею ряд
$$
  \sum_{n=1}^\infty n!\cdot x^n
$$
имеет нулевой радиус сходимости. К такому же выводу можно прийти, используя
критерий \eqref{DEF:anal-posled-2}.
\end{ex}

\begin{ex}\label{ex-21.1.6-sec} Наоборот, последовательность
$$
a_n=\frac{1}{n!}
$$
будет аналитической, потому что, как мы убедились в примере \ref{ex-21.1.6},
порождаемый ею ряд
$$
  \sum_{n=1}^\infty \frac{x^n}{n!}
$$
имеет радиус сходимости $R=\infty$ (для нас главное, что $R\ne 0$). Критерий
\eqref{DEF:anal-posled-2} дает то же самое.

\end{ex}

\end{multicols}\noindent\rule[10pt]{160mm}{0.1pt}

\subsection{Алгебраические операции с аналитическими
последовательностями}

 \bit{
\item[$\bullet$] Если $a\in\mathcal{A}$, то есть $a=\{a_n;\ n\in\Z_+\}$ --
аналитическая последовательность, то {\it противоположная последовательность}
$-a$ определяется формулой
$$
(-a)_n:=-a_n
$$
Из \eqref{DEF:anal-posled-2} сразу следует, что $-a\in\mathcal{A}$.

\item[$\bullet$] Если же $a,b\in\mathcal{A}$, то есть $a=\{a_n;\ n\in\Z_+\}$ и
$b=\{b_n;\ n\in\Z_+\}$ -- аналитические последовательности, то их {\it сумма}
$a+b$ и {\it свертка} $a*b$  определяются формулами
$$
(a+b)_n:=a_n+b_n,\qquad (a*b)_n:=\sum_{k=0}^n a_k\cdot b_{n-k}
$$
Сейчас мы покажем, что $a+b\in\mathcal{A}$ и $a*b\in\mathcal{A}$.
 }\eit
\bpr В силу \eqref{DEF:anal-posled-2}, условия $a\in\mathcal{A}$ и
$b\in\mathcal{A}$ означают, что для некоторых $A>0$ и $B>0$ выполняются
неравенства
$$
|a_n|\le A^n,\qquad |b_n|\le B^n\qquad (n\in\Z_+)
$$
Отсюда следует, во-первых,
$$
|(a+b)_n|=|a_n+b_n|\le |a_n|+|b_n|\le A^n+B^n\le \sum_{k=0}^n C_n^k\cdot
A^k\cdot B^{n-k}=(A+B)^n
$$
то есть последовательность $|(a+b)_n|$ мажорируется последовательностью
$(A+B)^n$, и значит, по условию \eqref{DEF:anal-posled-2}, $a+b\in\mathcal{A}$.

И, во-вторых,
$$
|(a*b)_n|=\left|\sum_{k=0}^n a_k\cdot b_{n-k}\right|\le \sum_{k=0}^n |a_k|\cdot
|b_{n-k}|\le \sum_{k=0}^n A^k\cdot B^{n-k}\le \sum_{k=0}^n C_n^k\cdot A^k\cdot
B^{n-k}=(A+B)^n
$$
то есть последовательность $|(a*b)_n|$ мажорируется последовательностью
$(A+B)^n$, и опять в силу \eqref{DEF:anal-posled-2}, $a*b\in\mathcal{A}$. \epr

\bigskip

\centerline{\bf Свойства алгебраических операций над аналитическими
последовательностями:}

 \bit{\it
\item[$1^\circ$.] Операция сложения на множестве $\mathcal{A}$ удовлетворяет
следующим тождествам:
 \beq
a+b=b+a,\quad (a+b)+c=a+(b+c),\quad  a+0=a,\quad a+(-a)=0
 \eeq
где $0$ -- последовательность, состоящая из нулей:
$$
0_n:=0
$$

\item[$2^\circ$.] Операция свертки на множестве $\mathcal{A}$ удовлетворяет
следующим тождествам:
 \beq
a*b=b*a,\quad (a*b)*c=a*(b*c),\quad  a*1=a
 \eeq
где $1$ -- последовательность, у которой на нулевом месте стоит единица, а на
остальных местах нули:
 \beq\label{DEF:1_n}
1_n:=\begin{cases}1,& n=0\\ 0,& n>0 \end{cases}
 \eeq

\item[$3^\circ$.] Операции $+$ и $*$ связаны между собой тождеством
дистрибутивности:
 \beq
(a+b)*c=a*c+b*c
 \eeq
 }\eit

\bpr Свойство $1^\circ$ мы считаем очевидным и сразу переходим к $2^\circ$.
Коммутативность свертки:
$$
(a*b)_n=\sum_{k=0}^n a_k\cdot b_{n-k}=\left|\begin{matrix}n-k=i,\quad k=n-i\\
0\le k\le n \quad\Leftrightarrow\quad 0\le i\le
n\end{matrix}\right|=\sum_{i=0}^n a_{n-i}\cdot b_i=(b*a)_n
$$
Ассоциативность свертки:
 \begin{multline*}
((a*b)*c)_n=\sum_{k=0}^n (a*b)_k\cdot c_{n-k}=\sum_{k=0}^n \l\sum_{i=0}^k
a_i\cdot b_{k-i}\r\cdot c_{n-k}=\sum_{\scriptsize\begin{matrix}0\le k\le n\\
0\le i\le k\end{matrix}} a_i\cdot b_{k-i}\cdot c_{n-k}=\\=
\sum_{\scriptsize\begin{matrix}0\le i\le n\\
i\le k\le n\end{matrix}} a_i\cdot b_{k-i}\cdot c_{n-k}= \sum_{i=0}^n a_i\cdot\l
\sum_{k=i}^n b_{k-i}\cdot c_{n-k}\r=\left|\begin{matrix}k-i=j,\quad k=i+j\\
i\le k\le n \quad\Leftrightarrow\quad 0\le j\le n-i\end{matrix}\right|=\\=
\sum_{i=0}^n a_i\cdot\l \sum_{j=0}^{n-i} b_j\cdot c_{n-i-j}\r= \sum_{i=0}^n
a_i\cdot (b*c)_{n-i}=(a*(b*c))_n
 \end{multline*}
Единица относительно свертки:
$$
(a*1)_n=\sum_{k=0}^n a_k\cdot\kern-7pt
\underbrace{1_{n-k}}_{\scriptsize\begin{matrix}\text{\rotatebox{90}{$=$}}\\
\begin{cases}1,& k=n\\ 0,& k\ne n\end{cases}\end{matrix}}\kern-7pt=\sum_{k=n} a_k\cdot
1=a_n
$$
Остается свойство дистрибутивности $3^\circ$:
$$
((a+b)*c)_n=\sum_{k=0}^n (a+b)_k\cdot c_{n-k}=\sum_{k=0}^n (a_k+b_k)\cdot
c_{n-k}=\sum_{k=0}^n a_k\cdot c_{n-k}+\sum_{k=0}^n b_k\cdot
c_{n-k}=(a*c)_n+(b*c)_n
$$
\epr

\subsection{Производящая функция}

 \bit{
\item[$\bullet$] {\it Производящей функцией} или {\it генератриссой} $\Gen_a$
аналитической последовательности $a\in\mathcal{A}$ называется функция,
определенная в окрестности нуля $(-\rho(a),\rho(a))$ формулой
$$
\Gen_a(x):=\sum_{n=0}^\infty a_n\cdot x^n
$$
По теореме \ref{TH:diff+int-step-ryada} эта функция является гладкой на
интервале $(-\rho(a),\rho(a))$.
 }\eit

\bigskip

\centerline{\bf Свойства производящих функций:}

 \bit{\it
\item[$1^\circ$.] Производящая функция суммы равна сумме производящих функций,
 \beq\label{Gen(a+b)(x)=Gen(a)(x)+Gen(b)(x)}
\Gen_{a+b}(x)=\Gen_a(x)+\Gen_b(x),
 \eeq
причем $\rho(a+b)\ge\min\{\rho(a),\rho(b)\}$.

\item[$2^\circ$.] Производящая функция свертки равна произведению производящих
функций,
 \beq\label{Gen(a*b)(x)=Gen(a)(x)-cdot-Gen(b)(x)}
\Gen_{a*b}(x)=\Gen_a(x)\cdot\Gen_b(x),
 \eeq
причем $\rho(a*b)\ge\min\{\rho(a),\rho(b)\}$.

 }\eit

Для доказательства второго из этих свойств нам понадобится следующая

\blm Пусть $\sum_{k=0}^\infty u_k$ и $\sum_{l=0}^\infty v_l$ -- два абсолютно
сходящихся ряда. Положим
$$
w_n=\sum_{k=0}^n u_k\cdot v_{n-k}
$$
Тогда ряд $\sum_{n=0}^\infty w_n$ абсолютно сходится, и его сумма равна
произведению сумм рядов $\sum_{k=0}^\infty u_k$ и $\sum_{l=0}^\infty v_l$:
 \beq\label{vspom-rav-dlya-Gen(a*b)}
\sum_{n=0}^\infty w_n=\lim_{N\to\infty}\sum_{n=0}^N \sum_{k=0}^n u_k\cdot
v_{n-k}=\left(\lim_{N\to\infty}\sum_{k=0}^N
u_k\right)\cdot\left(\lim_{N\to\infty}\sum_{l=0}^N
v_l\right)=\left(\sum_{k=0}^\infty u_k\right)\cdot\left(\sum_{l=0}^\infty
v_l\right)
 \eeq
 \elm
\bpr Для всякого $N\in\N$ получаем:
 \begin{multline*}
\sum_{n=0}^N |w_n|\le \sum_{n=0}^N \left|\sum_{k=0}^n u_k\cdot
v_{n-k}\right|\le \sum_{n=0}^N \sum_{k=0}^n |u_k|\cdot
|v_{n-k}|=\underbrace{\sum_{k+l\le N} |u_k|\cdot |v_l|\le \sum_{\max\{k,l\}\le
N} |u_k|\cdot |v_l|}_{k+l\le N\quad\Longrightarrow\quad \max\{k,l\}\le N}
=\\=\sum_{k=0}^N\sum_{l=0}^N |u_k|\cdot|v_l|=\l\sum_{k=0}^N
|u_k|\r\cdot\l\sum_{l=0}^N|v_l|\r\le \l\sum_{k=0}^\infty
|u_k|\r\cdot\l\sum_{l=0}^\infty |v_l|\r<\infty
 \end{multline*}
Отсюда следует, что ряд $\sum_{n=0}^\infty w_n$ абсолютно сходится. С другой
стороны, снова для любого $N\in\N$ имеем:
 \begin{multline*}
\left|\sum_{n=0}^{2N} w_n-\l\sum_{k=0}^N u_k\r\cdot\l\sum_{l=0}^N v_l\r\right|=
\left|\sum_{n=0}^{2N}\sum_{k=0}^n u_k\cdot v_{n-k}-\sum_{k=0}^N \sum_{l=0}^N
u_k\cdot v_l\right|=\\= \underbrace{\Bigg|\sum_{k+l\le 2N} u_k\cdot v_l
-\sum_{\max\{k,l\}\le N} u_k\cdot v_l\Bigg|}_{k+l\le 2N\quad\Longleftarrow\quad
\max\{k,l\}\le N}=\Bigg|\sum_{\scriptsize \left\{\begin{matrix}k+l\le 2N\\
\max\{k,l\}>N \end{matrix}\right\}} u_k\cdot v_l\Bigg|\le
\sum_{\scriptsize \left\{\begin{matrix}k+l\le 2N\\
\max\{k,l\}>N \end{matrix}\right\}} |u_k|\cdot |v_l|\le\\ \le
\sum_{\scriptsize \left\{\begin{matrix} \max\{k,l\}\le 2N\\
\max\{k,l\}>N \end{matrix}\right\}} |u_k|\cdot |v_l|=\sum_{N<\max\{k,l\}\le 2N}
|u_k|\cdot |v_l|=\sum_{\scriptsize\begin{matrix}0\le k\le N\\ N<l\le
2N\end{matrix}} |u_k|\cdot |v_l|+\sum_{\scriptsize\begin{matrix}N< k\le 2N\\
0\le l\le 2N\end{matrix}} |u_k|\cdot |v_l|=\\= \l\sum_{0\le k\le N}
|u_k|\r\cdot \l \sum_{N<l\le 2N} |v_l|\r+\l\sum_{N< k\le 2N}
|u_k|\r\cdot\l\sum_{0\le l\le N} |v_l|\r=\\= \underbrace{\l\sum_{k=0}^\infty
|u_k|\r}_{\scriptsize\begin{matrix}\text{\rotatebox{-90}{$<$}}\\
\infty \end{matrix}}\cdot
\underbrace{\l \sum_{l>N} |v_l|\r}_{\scriptsize\begin{matrix}\downarrow\\
\phantom{,}0,\\ \text{при $N\to\infty$}\end{matrix}}+\underbrace{\l\sum_{k>N}
|u_k|\r}_{\scriptsize\begin{matrix}\downarrow\\
\phantom{,}0,\\ \text{при
$N\to\infty$}\end{matrix}}\cdot\underbrace{\l\sum_{l=0}^\infty |v_l|\r}_{\scriptsize\begin{matrix}\text{\rotatebox{-90}{$<$}}\\
\infty \end{matrix}}\underset{N\to\infty}{\longrightarrow} 0
 \end{multline*}
И поэтому
$$
\sum_{n=0}^\infty w_n=\lim_{N\to\infty}\sum_{n=0}^{2N}
w_n=\lim_{N\to\infty}\left[\l\sum_{k=0}^N u_k\r\cdot\l\sum_{l=0}^N
v_l\r\right]= \l\sum_{k=0}^\infty u_k\r\cdot\l\sum_{l=0}^\infty v_l\r
$$
\epr

\bpr[Доказательство свойств
\eqref{Gen(a+b)(x)=Gen(a)(x)+Gen(b)(x)}-\eqref{|Gen(a)(x)|-le-Gen(|a|)(|x|)}]
1. Если $|x|<\min\{\rho(a),\rho(b)\}$, то
 \begin{multline*}
\Gen_{a+b}(x)=\lim_{N\to\infty}\sum_{n\le N}(a+b)_n\cdot
x^n=\lim_{N\to\infty}\sum_{n\le N}(a_n+b_n)\cdot
x^n=\\=\lim_{N\to\infty}\l\sum_{n\le N}a_n\cdot x^n+\sum_{n\le N}b_n\cdot
x^n\r= \sum_{n=0}^\infty a_n\cdot x^n+\sum_{n=0}^\infty b_n\cdot x^n=
\Gen_a(x)+\Gen_b(x)
 \end{multline*}
Это верно для любого $x$ такого, что $|x|<\min\{\rho(a),\rho(b)\}$, поэтому
$\rho(a+b)\ge \min\{\rho(a),\rho(b)\}$.

2. Если $|x|<\min\{\rho(a),\rho(b)\}$, то
 \begin{multline*}
\Gen_{a*b}(x)=\lim_{N\to\infty}\sum_{n\le N}(a*b)_n\cdot
x^n=\lim_{N\to\infty}\sum_{n\le N}\l\sum_{k=0}^n a_k\cdot b_{n-k}\r\cdot
x^n=\\=\lim_{N\to\infty}\sum_{n\le N}\sum_{k=0}^n \big(a_k\cdot x^k\big)\cdot
\big(b_{n-k}\cdot x^{n-k}\big)=\eqref{vspom-rav-dlya-Gen(a*b)}=
\l\sum_{k=0}^\infty a_k\cdot x^k\r\cdot\l\sum_{l=0}^\infty b_l\cdot x^l\r=
\Gen_a(x)\cdot\Gen_b(x)
 \end{multline*}
Это верно для любого $x$ такого, что $|x|<\min\{\rho(a),\rho(b)\}$, поэтому
$\rho(a+b)\ge \min\{\rho(a),\rho(b)\}$.

\epr

\subsection{Порядок аналитической последовательности}

 \bit{
\item[$\bullet$] {\it Порядком} аналитической последовательности
$b\in\mathcal{A}$ называется число
$$
\omega(a)=\min\{n\in\Z_+:\quad a_n\ne 0\}
$$
 }\eit

\btm Справедливо неравенство:
 \beq\label{omega(a*b)-ge-omega(a)+omega(b)}
\omega(a*b)\ge\omega(a)+\omega(b)
 \eeq
\etm \bpr Справедлива логическая цепочка:
$$
(a*b)_n=\sum_{k=0}^n a_k\cdot b_{n-k}\ne 0
$$
$$
\Downarrow
$$
$$
\exists k\in\{0,...,n\}\quad a_k\cdot b_{n-k}\ne 0
$$
$$
\Downarrow
$$
$$
\exists k\in\{0,...,n\}\quad \begin{cases}a_k\ne 0\\
b_{n-k}\ne 0\end{cases}
$$
$$
\Downarrow
$$
$$
\exists k\in\{0,...,n\}\quad \begin{cases}k\ge\omega(a)\\
n-k\ge\omega(b)\end{cases}
$$
$$
\Downarrow
$$
$$
\exists k\in\{0,...,n\}\quad \begin{cases}k\ge\omega(a)\\
n\ge\omega(b)+k\ge\omega(b)+\omega(a)
\end{cases}
$$
$$
\Downarrow
$$
$$
n\ge\omega(b)+\omega(a)
$$
Из нее получаем:
$$
\omega(a*b)=\min\big\{n\in\Z_+:\ (a*b)_n\ne 0\big\}\ge\omega(b)+\omega(a)
$$
\epr

\subsection{Сравнение аналитических последовательностей}

 \bit{
\item[$\bullet$] Неравенство $a\le b$ для последовательностей
$a,b\in\mathcal{A}$ означает, что каждая компонента $a$ не превосходит
соответствующую компоненту $b$:
$$
a\le b\qquad\Longleftrightarrow\qquad \forall n\in\Z_+\quad a_n\le b_n
$$
 }\eit

\bprop Производящая функция сохраняет неравенства для неотрицательных
аргументов:
 \beq\label{a-le-b=>Gen(a)(r)-le-Gen(b)(r)}
a\le b\qquad\Longrightarrow\qquad \forall r\in(0,\rho(b))\quad \Gen_{a}(r)\le
\Gen_b(r)
 \eeq
 \eprop
\bpr $ \Gen_{a}(r)=\underbrace{\sum_{n=0}^\infty a_n\cdot r^n\le
\sum_{n=0}^\infty b_n\cdot r^n}_{a_n\le b_n} \le \Gen_b(r). $ \epr

\subsection{Модуль аналитической последовательности}

 \bit{
\item[$\bullet$] {\it Модулем} последовательности $a\in\mathcal{A}$ называется
последовательность $|a|\in\mathcal{A}$, определенная формулой:
 \beq\label{modul-analit-posl}
|a|_n:=|a_n|
 \eeq
Очевидно, последовательность $a=\{a_n;\ n\in\Z_+\}$ является аналитической
тогда и только тогда, когда ее модуль  $|a|=\{|a_n|;\ n\in\Z_+\}$ является
аналитической последовательностью:
 \beq\label{a-in-A<=>|a|-in-A}
a\in\mathcal{A}\qquad\Longleftrightarrow\quad |a|\in\mathcal{A}
 \eeq
 }\eit

\bigskip

\centerline{\bf Свойства модуля аналитической последовательности:}

 \bit{\it
\item[$1^\circ$.] Модуль суммы не превосходит суммы модулей:
 \beq\label{|a+b|-le-|a|+|b|-analit}
|a+b|\le |a|+|b|
 \eeq

\item[$2^\circ$.] Модуль свертки не превосходит свертки модулей:
 \beq\label{|a*b|-le-|a|*|b|}
|a*b|\le |a|*|b|
 \eeq

\item[$3^\circ$.] Радиус сходимости модуля совпадает с радиусом сходимости
исходной последовательности:
 \beq\label{rho(|a|)=rho(a)}
\rho(|a|)=\rho(a)
 \eeq
а производящие функции удовлетворяют неравенству:
 \beq\label{|Gen(a)(x)|-le-Gen(|a|)(|x|)}
|\Gen_{a}(x)|\le \Gen_{|a|}(|x|)
 \eeq

 }\eit

\bpr Первые два свойства мы считаем очевидными, и перейдем сразу к $3^\circ$. В
нем неравенство \eqref{|Gen(a)(x)|-le-Gen(|a|)(|x|)} тоже очевидно
$$
|\Gen_{a}(x)|=\left|\sum_{n=0}^\infty a_n\cdot x^n\right|\le \sum_{n=0}^\infty
|a_n|\cdot |x|^n=\Gen_{|a|}(|x|)
$$
А равенство \eqref{rho(|a|)=rho(a)} доказывается двукратным применением
логической цепочки Абеля (см. с.\pageref{Abel-log-chain}). С одной стороны,
$$
0<r<\rho(a)
$$
$$
\Downarrow
$$
$$
\text{числовой ряд $\sum_{n=0}^\infty a_n\cdot r^n$ сходится}
$$
$$
  \Downarrow
$$
$$
\text{общий член стремится к нулю: $a_n\cdot r^n\underset{n\to
\infty}{\longrightarrow} 0$}
$$
$$
  \Downarrow
$$
$$
\text{последовательность $a_n\cdot r^n$ ограничена:}\quad
\sup_{n\in\Z_+}|a_n\cdot r^n|=M<\infty
$$
$$
  \Downarrow
$$
$$
\forall \sigma\in(0,r)\quad\Longrightarrow\quad |a_n|\cdot\sigma^n=|a_n\cdot
r^n|\cdot\left|\frac{\sigma}{r}\right|^n\le M\cdot
\left|\frac{\sigma}{r}\right|^n\qquad \l\left|\frac{\sigma}{r}\right|<1\r
$$
$$
  \Downarrow
$$
$$
\text{для всякого $\sigma\in(0,r)$ числовой ряд $\sum_{n=0}^\infty |a_n|\cdot
\sigma^n$ сходится}
$$
$$
  \Downarrow
$$
$$
\text{радиус сходимости ряда $\sum_{n=0}^\infty |a_n|\cdot x^n$ не меньше $r$}
$$
$$
  \Downarrow
$$
$$
\kern150pt \rho(|a|)\ge r \quad\leftarrow\ \text{верно для любого
$r\in(0,\rho(a))$}
$$
$$
  \Downarrow
$$
$$
\rho(|a|)\ge \rho(a)
$$

А с другой стороны,
$$
0<r<\rho(|a|)
$$
$$
\Downarrow
$$
$$
\text{числовой ряд $\sum_{n=0}^\infty |a_n|\cdot r^n$ сходится}
$$
$$
  \Downarrow
$$
$$
\text{общий член стремится к нулю: $|a_n|\cdot r^n\underset{n\to
\infty}{\longrightarrow} 0$}
$$
$$
  \Downarrow
$$
$$
\text{последовательность $|a_n|\cdot r^n$ ограничена:}\quad
\sup_{n\in\Z_+}|a_n|\cdot r^n=M<\infty
$$
$$
  \Downarrow
$$
$$
\forall \sigma\in(0,r)\quad\Longrightarrow\quad |a_n\cdot\sigma^n|=|a_n\cdot
r^n|\cdot\left|\frac{\sigma}{r}\right|^n\le M\cdot
\left|\frac{\sigma}{r}\right|^n \qquad \l\left|\frac{\sigma}{r}\right|<1\r
$$
$$
  \Downarrow
$$
$$
\text{для всякого $\sigma\in(0,r)$ числовой ряд $\sum_{n=0}^\infty a_n\cdot
\sigma^n$ сходится}
$$
$$
  \Downarrow
$$
$$
\text{радиус сходимости ряда $\sum_{n=0}^\infty a_n\cdot x^n$ не меньше $r$}
$$
$$
  \Downarrow
$$
$$
\kern150pt \rho(a)\ge r \quad\leftarrow\ \text{верно для любого
$r\in(0,\rho(|a|))$}
$$
$$
  \Downarrow
$$
$$
\rho(a)\ge \rho(|a|)
$$

 \epr

\subsection{Степень аналитической последовательности}

 \bit{
\item[$\bullet$] {\it Степень $b^k$} аналитической последовательности
$b\in\mathcal{A}$ определяется индуктивно правилами
$$
b^0=1,\qquad b^{k+1}=b^k* b,\qquad k\in\Z_+
$$
(здесь 1 -- единичная последовательность, определенная равенством
\eqref{DEF:1_n}).

 }\eit

\bigskip

\centerline{\bf Свойства степени:}

 \bit{\it
\item[$1^\circ$.] Модуль степени не превосходит степени модуля:
 \beq\label{|b^k|-le-|b|^k}
|b^k|\le |b|^k
 \eeq

\item[$2^\circ$.] Производящая функция степени равна степени производящей
функции:
 \beq\label{Gen(b^k)(x)=(Gen(b)(x))^k}
\Gen_{b^k}(x)=\Big(\Gen_b(x)\Big)^k
 \eeq

\item[$3^\circ$.] Порядок и степень аналитической последовательности связаны
неравенством:
 \beq\label{poraydok-i-stepen}
\forall k\in\Z_+\qquad \omega(b^k)\ge k\cdot\omega(b)
 \eeq
 }\eit

\bpr Свойство $1^\circ$ следует из \eqref{|a*b|-le-|a|*|b|}, $2^\circ$ -- из
\eqref{Gen(a*b)(x)=Gen(a)(x)-cdot-Gen(b)(x)}, а $3^\circ$ -- из
\eqref{omega(a*b)-ge-omega(a)+omega(b)} . \epr

\subsection{Композиция аналитических последовательностей}

 \bit{
\item[$\bullet$] Пусть порядок аналитической последовательности
$b\in\mathcal{A}$ отличен от нуля, то есть не меньше единицы:
$$
\omega(b)\ge 1.
$$
Тогда по формуле \eqref{poraydok-i-stepen}, порядок ее степени не меньше
показателя степени:
$$
\omega(b^k)\ge k,\qquad k\in\Z_+.
$$
То есть
$$
\forall n<k\qquad (b^k)_n=0
$$
и это можно интерпретировать так, что при фиксированном номере $n\in\Z$ почти
все числа $(b^k)_n;\ k\in\Z_+$ равны нулю:
 \beq\label{(b^k)_n=0}
\forall n\in\Z_+\qquad \forall k>n\qquad (b^k)_n=0
 \eeq
Отсюда следует, что для любой последовательности $a\in\mathcal{A}$ и любого
$n\in\Z_+$ в ряде $\sum_{k=0}^\infty a_k\cdot (b^k)_n$ только первые $n$
элементов (с индексами $k=1,...,n$) могут быть отличны от нуля, поэтому он
сходится. Его сумма обозначается $(a\circ b)_n$:
 \beq\label{circ-analit-posl}
(a\circ b)_n=\sum_{k=0}^\infty a_k\cdot (b^k)_n=\sum_{k=0}^n a_k\cdot (b^k)_n
 \eeq
Эта формула определяет некую последовательность
$$
a\circ b=\Big\{(a\circ b)_n;\ n\in\Z_+\Big\},
$$
называемую {\it композицией последовательностей} $a$ и $b$. Более коротко ее
определение можно записать формулой:
 \beq
a\circ b=\sum_{k=0}^\infty a_k\cdot b^k=\sum_{k=0}^n a_k\cdot b^k.
 \eeq
 }\eit

\btm\label{TH:compos-analit-posl} Композиция аналитических последовательностей
тоже является аналитической последовательностью:
 \beq
a,b\in\mathcal{A}\quad \&\quad \omega(b)\ge 1\quad\Longrightarrow\quad a\circ
b\in\mathcal{A}
 \eeq
 \etm
\bpr Пусть $a,b\in\mathcal{A}$ и $\omega(b)\ge 1$. Выберем
$\e\in(0,\rho(|a|))$, то есть $\e>0$ такое, что
$$
\Gen_{|a|}(\e)=\sum_{n=0}^\infty |a_n|\cdot \e^m<\infty
$$
Условие $\omega(b)\ge 1$ означает, что $b_0=0$, поэтому производящая функция
последовательности $|b|$ имеет вид
$$
\Gen_{|b|}(x)=\sum_{m=0}^\infty |b_m|\cdot x^m=\sum_{m=1}^\infty |b_m|\cdot x^m
$$
Отсюда в свою очередь следует, что
$$
\Gen_{|b|}(0)=0
$$
и, поскольку $\Gen_{|b|}$ -- непрерывная функция, существует $\delta>0$ такое,
что
$$
|\Gen_{|b|}(\delta)|<\e
$$
Теперь получаем:
$$
\phantom{\scriptsize\eqref{|b^k|-le-|b|^k}}\forall k\in\Z_+ \qquad |b^k|\le
|b|^k\qquad {\scriptsize\eqref{|b^k|-le-|b|^k}}
$$
$$
\phantom{\scriptsize\eqref{a-le-b=>Gen(a)(r)-le-Gen(b)(r)}}
\qquad\Downarrow\qquad{\scriptsize\eqref{a-le-b=>Gen(a)(r)-le-Gen(b)(r)}}
$$
$$
\forall k\in\Z_+ \qquad \Gen_{|b^k|}(\delta)\le
\Big(\Gen_{|b|}(\delta)\Big)^k<\e^k
$$
$$
\Downarrow
$$
 \begin{multline*}
\forall K,N\in\Z_+,\qquad K\ge N\quad\Longrightarrow\\
\sum_{n=0}^N |(a\circ b)_n|\cdot \delta^n= \sum_{n=0}^N \left|\sum_{k=0}^K
a_k\cdot (b^k)_n \right|\cdot \delta^n\le \sum_{n=0}^N \sum_{k=0}^K |a_k|\cdot
|(b^k)_n|\cdot \delta^n=\sum_{k=0}^K \sum_{n=0}^N |a_k|\cdot |(b^k)_n|\cdot
\delta^n=\\=\sum_{k=0}^K |a_k|\cdot \sum_{n=0}^N |(b^k)_n|\cdot \delta^n\le
\sum_{k=0}^K |a_k|\cdot\underbrace{\sum_{n=0}^\infty |(b^k)_n|\cdot
\delta^n}_{\Gen_{|b^k|(\delta)}}=\sum_{k=0}^K |a_k|\cdot\Gen_{|b^k|(\delta)}<
\sum_{k=0}^K |a_k|\cdot\e^k\le \sum_{k=0}^\infty
|a_k|\cdot\e^k=\\=\Gen_{|a|}(\e)<\infty
 \end{multline*}
$$
\Downarrow
$$
$$
\Gen_{|a\circ b|}(\delta)=\sum_{n=0}^\infty|(a\circ b)_n|\cdot \delta^n\le
\Gen_{|a|}(\e)<\infty
$$
Это означает, что $|a\circ b|\in\mathcal{A}$, и по свойству
\eqref{a-in-A<=>|a|-in-A}, $a\circ b\in\mathcal{A}$. \epr

\bigskip

\centerline{\bf Свойства композиции:}

 \bit{\it

\item[$1^\circ$.] Модуль композиции аналитических последовательностей не
превосходит композиции модулей:
 \beq
|a\circ b|\le |a|\circ |b|
 \eeq

\item[$2^\circ$.] Выполняется следующее правило дистрибутивности:
 \beq\label{(a+b)-circ-c=a-circ-c+b-circ-c}
(a+b)\circ c=a\circ c+b\circ c
 \eeq

\item[$3^\circ$.] Производящая функция композиции аналитических
последовательностей равна композиции производящих функций в некоторой
окрестности нуля:
 \beq\label{Gen(a-circ-b)=Gen(a)-circ-Gen(b)}
\exists \delta>0\qquad \forall x\in(-\delta,\delta)\qquad \Gen_{a\circ
b}(x)=(\Gen_a\circ\Gen_b)(x)
 \eeq

 }\eit
\bpr Свойство $1^\circ$ доказывается цепочкой:
 \begin{multline*}
|a\circ b|_n=\eqref{modul-analit-posl}=|(a\circ
b)_n|=\eqref{circ-analit-posl}=\left|\sum_{k=0}^n a_k\cdot (b^k)_n\right|\le
\sum_{k=0}^n |a_k|\cdot |(b^k)_n|=\eqref{modul-analit-posl}=\\=\sum_{k=0}^n
|a|_k\cdot |b^k|_n\le\eqref{|b^k|-le-|b|^k}\le \sum_{k=0}^n |a|_k\cdot
(|b|^k)_n=\eqref{circ-analit-posl}=(|a|\circ |b|)_n
 \end{multline*}
Свойство $2^\circ$:
 \begin{multline*}
((a+b)\circ c)_n=\eqref{circ-analit-posl}=\sum_{k=0}^n (a+b)_k\cdot
(c^k)_n=\sum_{k=0}^n (a_k+b_k)\cdot (c^k)_n=\\=\sum_{k=0}^n a_k\cdot
(c^k)_n+\sum_{k=0}^n b_k\cdot (c^k)_n=\eqref{circ-analit-posl}=(a\circ
c)_n+(b\circ c)_n=(a\circ c+b\circ c)_n
 \end{multline*}
Свойство $3^\circ$ мы сначала докажем для случая, когда $a$ -- финитная
последовательность, то есть такая, у которой почти все элементы равны нулю:
 \beq\label{sup(k:a_k-ne-0)=K<infty}
\sup\{k\in\Z_+:\quad a_k\ne 0 \}=K<\infty
 \eeq
В этом случае для всякого $x\in(-\rho(b),\rho(b))$ мы получим:
 \begin{multline*}
\Gen_{a\circ b}(x)=\sum_{n=0}^\infty (a\circ b)_n\cdot x^n= \sum_{n=0}^\infty
\sum_{k=0}^\infty a_k\cdot (b^k)_n\cdot
x^n=\eqref{sup(k:a_k-ne-0)=K<infty}=\sum_{n=0}^\infty \sum_{k=0}^K a_k\cdot
(b^k)_n\cdot x^n=\\=\sum_{k=0}^K a_k\cdot \underbrace{\sum_{n=0}^\infty
(b^k)_n\cdot x^n}_{\Gen_{b^k}(x)}=\sum_{k=0}^K a_k\cdot
\Gen_{b^k}(x)=\eqref{Gen(b^k)(x)=(Gen(b)(x))^k}= \sum_{k=0}^K
a_k\cdot\Big(\Gen_{b}(x)\Big)^k=\\=\sum_{k=0}^\infty
a_k\cdot\Big(\Gen_{b}(x)\Big)^k=\Gen_a\Big(\Gen_{b}(x)\Big)
 \end{multline*}
После этого рассматривается общий случай. Для любой последовательности
$a\in\mathcal{A}$ и для всякого $N\in\Z_+$ обозначим через $P_Na$
последовательность, определенную формулой
$$
(P_Ka)_k=\begin{cases}a_k, & k\le K \\ 0, & k>K\end{cases}
$$
Это будет финитная последовательность, поэтому по уже доказанному,
 \beq\label{Gen(a-circ-b)=Gen(a)-circ-Gen(b)-K<infty}
\Gen_{P_Ka\circ b}(x)=\Gen_{P_Ka}\Big(\Gen_b(x)\Big),\qquad
x\in(-\rho(b),\rho(b))
 \eeq
Зафиксируем $\e>0$ такое, что $\Gen_{|a|}(\e)<\infty$, и подберем $\delta>0$
такое, что $|\Gen_{|b|}(\delta)|<\e$. Покажем, что при $x\in(-\delta,\delta)$ и
$K\to\infty$ равенство \eqref{Gen(a-circ-b)=Gen(a)-circ-Gen(b)-K<infty}
превращается в \eqref{Gen(a-circ-b)=Gen(a)-circ-Gen(b)}:
 \beq\label{Gen(a-circ-b)=Gen(a)-circ-Gen(b)-K->infty}
\underbrace{\Gen_{P_Ka\circ b}(x)}_{\scriptsize\begin{matrix}(K\to\infty)\ \downarrow \ \phantom{(K\to\infty)}\\
\Gen_{a\circ b}(x)
\end{matrix}}=\underbrace{\Gen_{P_Ka}\Big(\Gen_b(x)\Big)}_{\scriptsize\begin{matrix}\phantom{(K\to\infty)}\
\downarrow \ (K\to\infty)
\\
\Gen_{a\circ b}(x) \end{matrix}},\qquad x\in(-\delta,\delta)
 \eeq
-- это и будет доказательством для \eqref{Gen(a-circ-b)=Gen(a)-circ-Gen(b)} в
общем случае.

Предельный переход справа в \eqref{Gen(a-circ-b)=Gen(a)-circ-Gen(b)-K->infty}
очевиден:
$$
\Gen_{P_Ka}(y)=\sum_{k=0}^K a_k\cdot
y^k\underset{K\to\infty}{\longrightarrow}\sum_{k=0}^\infty a_k\cdot y^k=
\Gen_{a}(y),\qquad y\in (-\e,\e)
$$
$$
\Downarrow
$$
$$
\Gen_{P_Ka}\Big(\Gen_b(x)\Big)\underset{K\to\infty}{\longrightarrow}\Gen_a\Big(\Gen_b(x)\Big),\qquad
x\in (-\delta,\delta)
$$
Докажем левый предельный переход. Для этого заметим прежде всего, что
 \beq\label{a_k-(P_Ka)_k}
a_k-(P_Ka)_k=\begin{cases}0, & k\le K \\ a_k, & k>K\end{cases}
 \eeq
Тогда для всякого $x\in(-\delta,\delta)$ мы получим:
 \begin{multline*}
\left|\Gen_{a\circ b}(x)-\Gen_{P_Ka\circ
b}(x)\right|=\eqref{Gen(a+b)(x)=Gen(a)(x)+Gen(b)(x)}=\left|\Gen_{a\circ
b-P_Ka\circ
b}(x)\right|=\eqref{(a+b)-circ-c=a-circ-c+b-circ-c}=\left|\Gen_{(a-P_Ka)\circ
b}(x)\right|\le\eqref{|Gen(a)(x)|-le-Gen(|a|)(|x|)}\le\\ \le
\Gen_{|(a-P_Ka)\circ b|}(|x|)\le \Gen_{|a-P_Ka|\circ
|b|}(|x|)=\sum_{n=0}^\infty (|a-P_Ka|\circ|b|)_n\cdot |x|^n=\sum_{n=0}^\infty
\sum_{k=0}^\infty |a-P_Ka|_k\cdot (|b|^k)_n\cdot |x|^n=\\= \sum_{n=0}^\infty
\sum_{k=0}^\infty |a_k-(P_Ka)_k|\cdot (|b|^k)_n\cdot
|x|^n=\eqref{a_k-(P_Ka)_k}=\sum_{n=0}^\infty \sum_{k=K+1}^\infty |a_k|\cdot
(|b|^k)_n\cdot |x|^n=\\= \sum_{k=K+1}^\infty |a_k|\cdot \sum_{n=0}^\infty
(|b|^k)_n\cdot |x|^n=\sum_{k=K+1}^\infty |a_k|\cdot\Gen_{|b|^k}(|x|)\le
\sum_{k=K+1}^\infty |a_k|\cdot\Big(\Gen_{|b|}(|x|)\Big)^k\le\\ \le
\sum_{k=K+1}^\infty |a_k|\cdot\Big(\Gen_{|b|}(\delta)\Big)^k\le
\sum_{k=K+1}^\infty |a_k|\cdot\e^k\underset{K\to\infty}{\longrightarrow} 0,
 \end{multline*}
потому что $\sum_{k=0}^\infty |a_k|\cdot\e^k=\Gen_{|a|}(\e)<\infty$. Отсюда уже
следует нужное нам соотношение
$$
\Gen_{P_Ka\circ b}(x)\underset{K\to\infty}{\longrightarrow} \Gen_{a\circ
b}(x),\qquad x\in (-\delta,\delta)
$$
\epr

\section{Ряд Тейлора и аналитические функции}

\subsection{Формулы Тейлора}

\paragraph{Теорема об остатке.}

 \bit{
\item[$\bullet$] Функция $f(x)$ называется

\bit{

\item[---] {\it дифференцируемой порядка $n$} на интервале $(a,b)$, если на
$(a,b)$ она имеет производные до порядка $n$ включительно, то есть существует
конечная последовательность функций $f_0,f_1,...,f_m$ на $(a,b)$ таких, что
$$
f_0=f,
$$
и для всякого $k=0,...,m-1$ функция $f_k$ дифференцируема на $(a,b)$ и ее
производная равна $f_{k+1}$:
$$
f'_k=f_{k+1},\qquad k<m.
$$
 }\eit

 }\eit

\begin{tm}\label{TH:Taylor-Zorich} Пусть нам даны:
 \bit{
\item[1)] функция $f$, дифференцируемая порядка $n+1$ на отрезке $[a,b]$, и

\item[2)] функция $g$, непрерывная на отрезке $[a,b]$, и имеющая внутри этого
отрезка ненулевую производную:
$$
\forall t\in(a,b)\qquad g'(t)\ne 0
$$
 }\eit
Тогда найдутся:
 \bit{
\item[--] точка $\zeta\in(a;b)$, для которой будет справедливо равенство:
\begin{equation}\label{Taylor-Zorich+}
f(b)=\sum_{k=0}^n \frac{f^{(k)}(a)}{k!}\cdot (b-a)^k+
      \frac{g(b)-g(a)}{g'(\zeta)\cdot n!}\cdot f^{(n+1)}(\zeta)\cdot (b-\zeta)^n
    \end{equation}
\item[--] точка $\eta\in(a;b)$, для которой будет справедливо равенство:
\begin{equation}\label{Taylor-Zorich-}
f(a)=\sum_{k=0}^n \frac{f^{(k)}(b)}{k!}\cdot (a-b)^k+
      \frac{g(a)-g(b)}{g'(\eta)\cdot n!}\cdot f^{(n+1)}(\eta)\cdot (a-\eta)^n
    \end{equation}
  }\eit
\end{tm}

\brem Формулы \eqref{Taylor-Zorich+} и \eqref{Taylor-Zorich-} можно объединить
в одну, усложнив формулировку теоремы следующим образом: {\it пусть $x_0$ и $x$
обозначают концы отрезка $[a,b]$, причем возможны оба варианта, как
$\left\{\begin{matrix}x_0=a \\ x=b\end{matrix}\right\}$, так и
$\left\{\begin{matrix}x=a \\ x_0=b\end{matrix}\right\}$; тогда найдется точка
$\xi\in(a;b)$ такая, что}
\begin{equation}\label{Taylor-Zorich}
f(x)=\sum_{k=0}^n \frac{f^{(k)}(x_0)}{k!}\cdot (x-x_0)^k+
      \frac{g(x)-g(x_0)}{g'(\xi)\cdot n!}\cdot f^{(n+1)}(\xi)\cdot (x-\xi)^n
    \end{equation}
В такой формулировке случай $\left\{\begin{matrix}x_0=a \\
x=b\end{matrix}\right\}$ превращает формулу \eqref{Taylor-Zorich} в
\eqref{Taylor-Zorich+}, а при $\left\{\begin{matrix}x=a \\
x_0=b\end{matrix}\right\}$ формула \eqref{Taylor-Zorich} превращается в
\eqref{Taylor-Zorich-}.
 \erem

\bpr 1. Для доказательства \eqref{Taylor-Zorich+} рассмотрим вспомогательную
функцию
$$
F(t)=f(b)-\sum_{k=0}^n\frac{f^{(k)}(t)}{k!}\cdot (b-t)^k
$$
Применим к паре функций $F$ и $g$ теорему Коши об отношении приращений
\ref{Cauchy-II}: должна существовать точка $\zeta\in(a;b)$ такая, что
 $$
\frac{F(b)-F(a)}{g(b)-g(a)}=\frac{F'(\zeta)}{g'(\zeta)}
 $$
Или, иначе
 \beq\label{taylor-zorich-1}
F(b)-F(a)=\frac{g(b)-g(a)}{g'(\zeta)}\cdot F'(\zeta)
 \eeq
Выразим в этой формуле $F$ через $f$. Во-первых,
 $$
F(b)=f(b)-\sum_{k=0}^n\frac{f^{(k)}(b)}{k!}\cdot (b-b)^k=
f(b)-\sum_{k=0}^n\frac{f^{(k)}(b)}{k!}\cdot\kern-20pt\underbrace{
0^k}_{\scriptsize\begin{matrix} \| \\ \begin{cases} 0, & \text{при $k>0$} \\
 0, & \text{при $k=0$}\end{cases} \end{matrix}}\kern-20pt= f(b)-\underbrace{\frac{f^{(0)}(b)}{0!}
 }_{\scriptsize\begin{matrix} \| \\ f(b)\end{matrix}}=0
$$
Во-вторых,
$$
F(a)=f(a)-\sum_{k=0}^n\frac{f^{(k)}(a)}{k!}\cdot (b-a)^k
$$
И, в-третьих,
 \begin{multline*}
F'(t)=\frac{\d}{\d t}\left(f(b)-\sum_{k=0}^n\frac{f^{(k)}(t)}{k!}\cdot
(b-t)^k\right)=\frac{\d}{\d
t}\left(f(b)-f(t)-\sum_{k=1}^n\frac{f^{(k)}(t)}{k!}\cdot (b-t)^k\right)=\\=
0-f'(t)-\sum_{k=1}^n\frac{\d}{\d t}\left(\frac{f^{(k)}(t)}{k!}\cdot
(b-t)^k\right)=\\=-f'(t)-\sum_{k=1}^n\left(\frac{\d}{\d
t}\left(\frac{f^{(k)}(t)}{k!}\right)\cdot
(b-t)^k+\frac{f^{(k)}(t)}{k!}\cdot\frac{\d}{\d t}\left((b-t)^k\right)\right)=
\\= -f'(t)-\sum_{k=1}^n\left(\frac{f^{(k+1)}(t)}{k!}\cdot
(b-t)^k+\frac{f^{(k)}(t)}{k!}\cdot k\cdot (b-t)^{k-1}\cdot (-1)\right)=
\\= \underbrace{-\kern-10pt\underbrace{f'(t)}_{\scriptsize\begin{matrix}\| \\ \frac{f^{(1)}(t)}{0!}\cdot
(b-t)^0\end{matrix}}\kern-10pt-\sum_{k=1}^n \frac{f^{(k+1)}(t)}{k!}\cdot
(b-t)^k}_{\text{объединяем в одну сумму}}+\underbrace{\sum_{k=1}^n
\frac{f^{(k)}(t)}{(k-1)!}\cdot (b-t)^{k-1}}_{\text{заменяем $k-1$ на $i$}}=
\\= -\underbrace{\sum_{k=0}^n \frac{f^{(k+1)}(t)}{k!}\cdot
(b-t)^k}_{\text{отщепляем последнее слагаемое}}+\sum_{i=0}^{n-1}
\frac{f^{(i+1)}(t)}{i!}\cdot (b-t)^i=\\=-\frac{f^{(n+1)}(t)}{n!}\cdot (b-t)^n
\underbrace{-\sum_{k=0}^{n-1}\frac{f^{(k+1)}(t)}{k!}\cdot
(b-t)^k+\sum_{i=0}^{n-1} \frac{f^{(i+1)}(t)}{i!}\cdot
(b-t)^i}_{\text{сокращаем}}=-\frac{f^{(n+1)}(t)}{n!}\cdot (b-t)^n
 \end{multline*}
Подставляя $\zeta$ вместо $t$ получим:
 $$
F'(\zeta)=-\frac{f^{(n+1)}(\zeta)}{n!}\cdot (b-\zeta)^n
 $$
Теперь подставим полученные выражения в \eqref{taylor-zorich-1}:
 $$
0-\left(f(a)-\sum_{k=0}^n\frac{f^{(k)}(a)}{k!}\cdot (b-a)^k\right)=
\frac{g(b)-g(a)}{g'(\zeta)}\cdot\left(-\frac{f^{(n+1)}(\zeta)}{n!}\cdot
(b-\zeta)^n\right)
 $$
Убрав минус в обеих частях, получим равенство, эквивалентное
\eqref{Taylor-Zorich+}:
 $$
f(a)-\sum_{k=0}^n\frac{f^{(k)}(a)}{k!}\cdot (b-a)^k=
\frac{g(b)-g(a)}{g'(\zeta)}\cdot\frac{f^{(n+1)}(\zeta)}{n!}\cdot (b-\zeta)^n
 $$

2. Для доказательства \eqref{Taylor-Zorich-} нужно рассмотреть другую
вспомогательную функцию:
$$
F(t)=f(a)-\sum_{k=0}^n\frac{f^{(k)}(t)}{k!}\cdot (a-t)^k
$$
Тогда
$$
F(b)=f(a)-\sum_{k=0}^n\frac{f^{(k)}(b)}{k!}\cdot (a-b)^k
$$
$$
F(a)=f(a)-\sum_{k=0}^n\frac{f^{(k)}(a)}{k!}\cdot (a-a)^k=0
$$
$$
F'(t)=-\frac{f^{(n+1)}(t)}{n!}\cdot (a-t)^n
$$
И при подстановке в \eqref{taylor-zorich-1} получается равенство, эквивалентное
\eqref{Taylor-Zorich-}:
$$
f(a)-\sum_{k=0}^n\frac{f^{(k)}(b)}{k!}\cdot (a-b)^k-0=
\frac{g(b)-g(a)}{g'(\zeta)}\cdot\left(-\frac{f^{(n+1)}(\zeta)}{n!}\cdot
(a-\zeta)^n\right)
$$

 \epr

\paragraph{Формула Тейлора  с остаточным членом в форме Лагранжа.}

\btm[\bf Тейлор-Лагранж]\label{TH:Taylor-Lagrange} Пусть функция $f$ определена
и дифференцируема $n+1$ раз на интервале $(a,b)$. Тогда для любых дух точек
$x_0,x\in (a,b)$ найдется точка $\xi$, лежащая между ними (то есть
$\xi\in(x_0,x)$, если $x_0<x$, и $\xi\in(x,x_0)$, если $x<x_0$) такая, что
выполняется равенство
\smallskip
\begin{equation}\label{Taylor-Lagrange}\boxed{\quad
\begin{split}\phantom{\frac{\frac{\frac{\frac{1}{1}}{1}}{1}}{1}}
f(x)&=\sum_{k=0}^n \frac{f^{(k)}(x_0)}{k!}\cdot (x-x_0)^k+
      \frac{f^{(n+1)}(\xi)}{(n+1)!}\cdot (x-x_0)^{n+1}
    \phantom{\frac{1}{\frac{1}{\frac{1}{1}}}}
    \end{split}\quad}
    \end{equation}
 \smallskip
 \etm

 \bit{
\item[$\bullet$] Равенство \eqref{Taylor-Lagrange} называется  {\it формулой
Тейлора с остаточным членом в форме Лагранжа}\index{формула!Тейлора с
остаточным членом в форме!Лагранжа} {\it порядка} $n$ для функции $f$ в точке
$x_0$.
 }\eit
 \bpr
Этот факт следует из теоремы \ref{TH:Taylor-Zorich}, если в ней выбрать в
качестве $g$ функцию
$$
g(t)=(x-t)^{n+1}
$$
 \epr

\paragraph{Формула Тейлора  с остаточным членом в форме Коши.}

\btm[\bf Тейлора-Коши]\label{TH:Taylor-Cauchy} Пусть функция $f$ определена и
дифференцируема $n+1$ раз на интервале $(a,b)$. Тогда для любых дух точек
$x_0,x\in (a,b)$ найдется точка $\xi$, лежащая между ними (то есть
$\xi\in(x_0,x)$, если $x_0<x$, и $\xi\in(x,x_0)$, если $x<x_0$) такая, что
выполняется равенство
\smallskip
\begin{equation}\label{Taylor-Cauchy}\boxed{\quad
\begin{split}\phantom{\frac{\frac{\frac{\frac{1}{1}}{1}}{1}}{1}}
f(x)&=\sum_{k=0}^n \frac{f^{(k)}(x_0)}{k!}\cdot (x-x_0)^k+
      \frac{f^{(n+1)}(\xi)}{n!}\cdot (x-\xi)^n\cdot (x-x_0)
    \phantom{\frac{1}{\frac{1}{\frac{1}{1}}}}
    \end{split}\quad}
    \end{equation}
 \smallskip
 \etm
\bit{ \item[$\bullet$] Равенство \eqref{Taylor-Cauchy} называется {\it формулой
Тейлора с остаточным членом в форме Коши}\index{формула!Тейлора с остаточным
членом в форме!Коши} {\it порядка} $n$ для функции $f$ в точке $x_0$:
\smallskip
 \smallskip
}\eit
 \bpr
Это следует из теоремы \ref{TH:Taylor-Zorich}, если в ней выбрать в качестве
$g$ функцию
$$
g(t)=x-t
$$
 \epr

\subsection{Дифференциал функции и его вычисление}
\label{SEC-differentsial-i-diff-termy}

\paragraph{Дифференциал функции и запись формулы Тейлора с его помощью}

 \bit{
\item[$\bullet$] {\it Дифференциалом $0$-го порядка}\index{дифференциал!функции
одной переменной!$0$-го порядка} функции $f$ в точке $x$ называется функция
$p\mapsto\d f(x)[p]$ от нового переменного $p\in \R$, которая в каждой точке
$p\in \R$ определяется равенством:
 \begin{equation}\label{opredelenie-differentsiala-0-v-R^1}
\d^0 f(x)[p]=f(x)
 \end{equation}

\item[$\bullet$] {\it Дифференциалом (1-го порядка)}\index{дифференциал!функции
одной переменной!1-го порядка} функции $f$ в точке $x$ называется функция
$p\mapsto\d f(x)[p]$ от нового переменного $p\in \R$, которая в каждой точке
$p\in \R$ определяется равенством:
 \begin{equation}\label{opredelenie-differentsiala-1-v-R^1}
\d f(x)[p]=\lim_{t\to 0}\frac{f(x+tp)-f(x)}{t}=\left(\frac{\d}{\d
t}f(x+tp)\right)\Bigg|_{t=0}
 \end{equation}

\item[$\bullet$] {\it Дифференциалом $k$-го порядка}\index{дифференциал!функции
одной переменной!$k$-го порядка} функции $f$ в точке $x$ называется функция
$p\mapsto\d f(x)[p]$ от нового переменного $p\in \R$, которая в каждой точке
$p\in \R$ определяется равенством:
 \begin{equation}\label{opredelenie-differentsiala-k-v-R^1}
\d^k f(x)[p]=\left(\frac{\d^k}{\d t^k}f(x+tp)\right)\Bigg|_{t=0}
 \end{equation}
 }\eit

\btm Справедлива формула
$$
\d^k f(x)[p]=f^{(k)}(x)\cdot p^k
$$
из которой следует модификация формулы Тейлора-Лагранжа:
 \beq
f(x+p)=\sum_{k=0}^m\frac{1}{k!}\cdot\d^k f(x)[p]+\frac{1}{(m+1)!}\cdot\d^{m+1}
f(x+\theta\cdot p)[p],\qquad \theta\in[0,1]
 \eeq
\etm

\noindent\rule{160mm}{0.1pt}\begin{multicols}{2}

\paragraph{Дифференциальные выражения.}

\biter{

\item[$\bullet$] {\it Дифференциальным выражением} (или {\it дифференциальным многочленом}) называется обобщение понятий числового выражения (введенного на с. \pageref{opredelenie-chislovogo-terma}) и однородного дифференциального выражения степени 1 (введенного на с.\pageref{diff-vyrazh-step-1}), описываемое следующими правилами:
 \biter{

\item[0)] всякое числовое выражение ${\mathcal P}$ считается также дифференциальным выражением (такие выражения называются {\it однородными дифференциальными многочленами степени 0});

\item[1)] если ${\mathcal U}$ и ${\mathcal V}$ --- два дифференциальных выражения, то записи
$$
{\mathcal U}+{\mathcal V},\qquad {\mathcal U}-{\mathcal V},\qquad {\mathcal U}\cdot {\mathcal V}
$$
тоже считаются дифференциальными выражениями;

\item[2)] если ${\mathcal P}$ -- числовое выражение, а ${\mathcal U}$ -- дифференциальное выражение, то запись
$$
\frac{\mathcal U}{\mathcal P}
$$
тоже считается дифференциальным выражением;

\item[3)] если ${\mathcal U}$ --- дифференциальное выражение, то запись
$$
\d {\mathcal U},
$$
также считается дифференциальным выражением, называемым {\it дифференциалом} выражения ${\mathcal U}$;

 }\eiter
 }\eiter

\paragraph{Равенство дифференциальных выражений.}

\biter{

\item[$\bullet$] Два дифференциальных выражения ${\mathcal U}$ и ${\mathcal V}$ считаются {\it
равными на множестве $x\in X$}, и изображается это формулой
$$
{\mathcal U}\underset{x\in X}{=}{\mathcal V},
$$
или, в случае, если $X=\D(\mathcal U)\cap\D(\mathcal V)$, то формулой
$$
{\mathcal U}={\mathcal V},
$$
и, наконец, в случае, если $X=\D(\mathcal U)=\D(\mathcal V)$, то формулой
$$
{\mathcal U}\equiv{\mathcal V},
$$
если их можно отождествить с помощью следующих индуктивных правил:

 \biter{
\item[(a)] умножение дифференциальных выражений подчинено правилам (здесь ${\mathcal P}$ -- числовое выражение, а ${\mathcal U}$, ${\mathcal V}$ и ${\mathcal W}$ -- дифференциальные):
\kern-20pt \begin{align}
& {\mathcal U}\cdot {\mathcal V}\equiv{\mathcal V}\cdot {\mathcal U} \label{kommut-diff-vyrazh-8} \\
& \frac{1}{\mathcal P}\cdot {\mathcal U}\equiv\frac{\mathcal U}{\mathcal P}\equiv{\mathcal U}\cdot \frac{1}{\mathcal P} \label{drob-diff-vyrazh-8} \\
& {\mathcal U}\cdot \Big({\mathcal V}\cdot {\mathcal W}\Big)\equiv\Big({\mathcal U}\cdot {\mathcal V}\Big)\cdot {\mathcal W}
\label{assoc-diff-vyrazh-8} \\
& \Big({\mathcal U}+{\mathcal V}\Big)\cdot {\mathcal W}\equiv{\mathcal U}\cdot {\mathcal W}+ {\mathcal V}\cdot {\mathcal W} \label{distribut-diff-vyrazh-8} \\
& 1\cdot {\mathcal U}\equiv{\mathcal U} \label{edin-v-diff-vyrazh-8}
 \end{align}

\item[(b)] если $x$ -- переменная, ${\mathcal P}$ и ${\mathcal Q}$ -- числовые выражения, а ${\mathcal U}$ и ${\mathcal V}$ -- дифференциальные  выражения, то
 \begin{multline}\label{ravenstvo-diff-vyrazh}
{\mathcal P}+{\mathcal U}\d x\underset{x\in X}{=}{\mathcal Q}+{\mathcal V}\d x \quad\Longleftrightarrow\\ \Longleftrightarrow
\quad \begin{cases}{\mathcal P}\underset{x\in X}{=}{\mathcal Q} \\ {\mathcal U}\underset{x\in X}{=}{\mathcal V} \end{cases}
 \end{multline}

\item[(c)] для элементарных числовых выражений дифференциалы задаются таблицей дифференциалов на с.\pageref{tablitsa-differentsialov};

\item[(d)] дифференциал от суммы, разности, произведения и частного
дифференциальных выражений подчиняется правилам (здесь ${\mathcal Q}$ --  числовое, а ${\mathcal U}$ и ${\mathcal V}$ -- дифференциальные выражения):
\begin{align}
& \kern-20pt \d \Big({\mathcal U}\pm {\mathcal V}\Big)\equiv
\d {\mathcal U}\pm \d {\mathcal V}  \label{differentsial-ot-summy-8} \\
& \kern-20pt \d \Big({\mathcal U}\cdot {\mathcal V}\Big)\equiv  \Big(\d {\mathcal U}\Big)\cdot
{\mathcal V}+ {\mathcal U}\cdot \Big(\d {\mathcal V} \Big)
\label{differentsial-ot-proizvedeniya-8} \\
& \kern-20pt \d \l \frac{{\mathcal U}}{{\mathcal Q}}\r \equiv  \frac{\Big(\d {\mathcal
U}\Big)\cdot {\mathcal Q} -{\mathcal U}\cdot \Big(\d {\mathcal
Q}\Big)}{{\mathcal Q}^2} \label{differentsial-ot-drobi-8}
\end{align}

\item[(e)] если ${\mathcal P}$ --- одноместное числовое выражение от переменной
$x$, а ${\mathcal Q}$ --- одноместное числовое выражение от переменной $y$, то
для выражения, полученного подстановкой, выполняется равенство \eqref{differentsial-slozhnogo-terma}:
 \begin{align*}
& \kern-20pt \d \left({\mathcal Q}\Big|_{y={\mathcal P}}\right)\equiv  \left(\frac{\d}{\d
y}{\mathcal Q}\right)\Big|_{y={\mathcal P}}\cdot \d {\mathcal P}
 \end{align*}

\item[(f)] второй дифференциал от любой переменной $x$ удовлетворяет равенству:
 \beq\label{d-d-x}
\d(\d x)=0\cdot(\d x)^2
 \eeq

 }\eiter
 }\eiter

\paragraph{Нулевое дифференциальное выражение.}

\btm
Для всякой переменной $x$ справедливы равенства:
 \beq\label{0=0.dx=...=o.dx^k=...}
0\equiv 0\cdot\d x\equiv...\equiv  0\cdot\big(\d x\big)^k\equiv...
 \eeq
\etm
 \biter{
\item[$\bullet$] Дифференциальное выражение в формуле \eqref{0=0.dx=...=o.dx^k=...} называется {\it нулевым дифферен\-циальным выражением} от переменной $x$ и обозначается символом 0:
    $$
    0:=0\equiv 0\cdot\d x\equiv...\equiv  0\cdot\big(\d x\big)^k\equiv...
    $$
 }\eiter
\bpr
$$
0+\underbrace{0\cdot\d x}_{\scriptsize\begin{matrix}\text{переносим}\\ \text{вправо}\end{matrix}}\equiv \kern-10pt\overbrace{0}^{\scriptsize\begin{matrix}\text{переносим}\\ \text{влево}\end{matrix}}\kern-10pt+0\cdot\d x
$$
$$
\Updownarrow
$$
$$
\underbrace{0-0}_{\scriptsize\begin{matrix}\text{\rotatebox{90}{$\equiv$}}\\ 0\end{matrix}}\equiv
\underbrace{0\cdot\d x-0\cdot\d x}_{\scriptsize\begin{matrix}\text{\rotatebox{90}{$\equiv$}}\\ (0-0)\cdot\d x\\
\text{\rotatebox{90}{$\equiv$}}\\ 0\cdot\d x\end{matrix}}
$$
$$
\Updownarrow
$$
$$
0\equiv 0\cdot\d x
$$
Это первое тождество в \eqref{0=0.dx=...=o.dx^k=...}. Из него по индукции следуют все остальные. Например, умножая на $\d x$, получаем:
$$
0\cdot\d x\equiv 0\cdot\d x\cdot\d x\equiv 0\cdot\big(\d x\big)^2
$$
И так далее.
\epr

\bigskip

\centerline{\bf Свойства нулевого}
\centerline{\bf дифференциального  выражения:}

 \biter{\it

\item[$1^\circ.$] Умножение на нулевое дифференциальное выражение всегда дает нулевое дифференциальное выражение:
\beq\label{umnozh-na-nul-diff-vyrazh-8}
{\mathcal U}\cdot 0\equiv 0\equiv 0\cdot{\mathcal U}
\eeq

\item[$2^\circ.$] Прибавление нулевого дифференциального выражения не меняет дифференциальное выражение:
\beq\label{slozh-s-nul-diff-vyrazh-8}
{\mathcal U}+ 0\equiv{\mathcal U}\equiv 0+{\mathcal U}
\eeq

\item[$4^\circ.$] Дифференциал от произвольного параметра равен нулю:
\beq\label{differentsial-ot-parametra-8}
\d C\equiv 0
\eeq

\item[$5^\circ.$] Дифференциал от переменной наоборот, не может быть равен нулю:
\beq\label{differentsial-ot-peremennoi-8}
\d x\ne 0
\eeq

}\eiter
\bpr
\epr

\bcor
Всякий параметр можно выносить за знак
дифференциала:
 \beq\label{vynesenie-konstaty-za-differentsial}
\d\big( C\cdot {\mathcal P}\big)=C\cdot \d {\mathcal P}
 \eeq
\ecor

\paragraph{Разложение дифференциального выражения по степеням элементарного дифференциала.}

\bit{ \item[$\bullet$] Выражения вида
$$
(\d x)^2=\d x\cdot\d x,\quad (\d x)^3=(\d x)^2\cdot\d x,...
$$
называются {\it степенями элементарных дифференциалов} ($\d x$, $\d y$ и т.д.).
}\eit

Из \eqref{d-d-x} и \eqref{differentsial-ot-proizvedeniya} следует, что
дифференциал от любого такого выражения равен нулю:
$$
\d\Big((\d x)^k\Big)=0
$$
Отсюда, в свою очередь, для всякого дифференциального (или числового) выражения
${\mathcal D}$ получаем формулу:
 \beq\label{d-D-d-x}
\d\Big({\mathcal D}\cdot(\d x)^k\Big)=\Big(\d{\mathcal D}\Big)\cdot(\d x)^k
 \eeq

\btm[\bf о строении дифференциальных
выражений]\label{TH-o-stroenii-diff-termov} Всякое дифференциальное выражение
${\mathcal U}$ от переменной $x$ раскладывается по степеням элементарного
дифференциала $\d x$,
 \beq\label{razlozhenie-diff-terma}
{\mathcal U}=\sum_{k=0}^n {\mathcal P}_k\cdot \big(\d x\big)^k,
 \eeq
где коэффициенты ${\mathcal P}_k$ --- числовые выражения, причем такое
разложение единственно с точностью до равенства коэффициентов ${\mathcal P}_k$:
если ${\mathcal Q}_k$ --- какие-то другие числовые выражения с тем же
свойством,
$$
{\mathcal U}=\sum_{k=0}^n {\mathcal P}_k\cdot \big(\d x\big)^k=\sum_{k=0}^n
{\mathcal Q}_k\cdot \big(\d x\big)^k,
$$
то автоматически
$$
\forall k=0,...,n\qquad {\mathcal P}_k={\mathcal Q}_k
$$
При этом
\beq
 \D({\mathcal U})=\bigcap_{i=0}^k\D({\mathcal P}_i)
\eeq
 \etm

\paragraph{Связь между дифференциалом функции и дифференциалом выражения.}

 \biter{
\item[$\bullet$] Дифференциальное выражение ${\mathcal D}$ называется {\it
однородным}, если в его каноническом разложении \eqref{razlozhenie-diff-terma}
только одно слагаемое отлично от нуля, то есть если  ${\mathcal D}$ можно
представить в виде
 \beq\label{odnor-diff-term}
{\mathcal D}={\mathcal P}\cdot \big(\d x\big)^k
 \eeq
где ${\mathcal P}$ --- некоторое числовое выражение. Число $k\in\Z_+$ при этом
называется {\it степенью} однородного дифференциального выражения ${\mathcal
D}$.
 }\eiter

\btm[\bf о связи между дифференциалом функции и дифференциалом выражения] Если
стандартная функция $f$ задается выражением $\mathcal P$
$$
f(x)={\mathcal P}
$$
то ее дифференциал $\d^kf$ выражается через дифференциал выражения $\mathcal P$
по формуле:
 \beq
\d^kf(x)[p]=\Big(\d^k {\mathcal P}\Big)\Bigg|_{\d x=p}
 \eeq
\etm

Покажем, как описанные операции над дифференциальными выражениями позволяют упростить нахождение формул
Тейлора.

\bex В примере \ref{ex-f-taylor-dlya-sin} мы выписывали формулу Тейлора порядка
2 для функции
$$
f(x)=\sin x
$$
в точке $x=\frac{\pi}{4}$. С помощью исчисления выражений это делается так.
Сначала находим дифференциалы. Первый дифференциал определяющего эту функцию
выражения выписывается из таблицы дифференциалов:
$$
\d\sin x=\cos x\cdot \d x
$$
Второй получается применением правил вычисления:
 \begin{multline*}
\d^2\sin x=\d\left(\cos x\cdot \d x\right)=\eqref{d-D-d-x}=\\= \d\left(\cos
x\right)\cdot \d x= -\sin x\cdot\d x\cdot\d x=-\sin x\cdot(\d x)^2
 \end{multline*}
Третий:
 \begin{multline*}
\d^3\sin x=\d\left(-\sin x\cdot(\d x)^2\right)=\eqref{d-D-d-x}=\\=-\d\left(\sin
x\right)\cdot (\d x)^2=-\cos x\cdot\d x\cdot(\d x)^2=\\=-\cos x\cdot(\d x)^3
 \end{multline*}
Теперь можно выписать дифференциалы функции $f$:
 \begin{align*}
\d^0 f(x)[p]&=f(x)=\sin x\\
\d^1 f(x)[p]&=\cos x\cdot \d x\Big|_{\d
x=p}=\\&=\cos x\cdot p \\
\d^2 f(x)[p]&=-\sin x\cdot (\d x)^2\Big|_{\d x=p}=\\&=-\sin x\cdot
p^2\\
\d^3 f(x)[p]&=-\cos x\cdot (\d x)^3\Big|_{\d x=p}=\\&=-\cos x\cdot p^3
 \end{align*}
Формула Тейлора-Лагранжа (порядка 2) для этой функции имеет вид:
\begin{multline*}
\sin(x+p)=\sin x+\cos x\cdot p-\frac{1}{2}\sin x\cdot p^2-\\-\frac{1}{6}\cos
(x+\theta p)\cdot p^3,\quad \theta\in[0,1]
\end{multline*}
В точке $x=\frac{\pi}{4}$ получаем:
\begin{multline*}
\sin\left(\frac{\pi}{4}+p\right)=\frac{\sqrt{2}}{2}+\frac{\sqrt{2}}{2} \cdot
p-\frac{\sqrt{2}}{4}\cdot p^2-\\-\frac{1}{6}\cos \l\frac{\pi}{4}+\theta p\r\cdot p^3,\quad \theta\in[0,1]
\end{multline*}
\eex

\bex Найдем дифференциалы до 3 порядка включительно от функции
$$
f(x)=\arctg x
$$
Первый дифференциал определяющего ее выражения выписывается из таблицы
дифференциалов:
$$
\d\arctg x=\frac{1}{1+x^2}\cdot \d x
$$
Второй получается применением правил вычисления:
 \begin{multline*}
\d^2\arctg x=\d\left(\frac{1}{1+x^2}\cdot \d x\right)=\\=
\d\left(\frac{1}{1+x^2}\right)\cdot \d x=-\frac{2x}{(1+x^2)^2}\cdot \d x\cdot
\d x=\\=-\frac{2x}{(1+x^2)^2}\cdot (\d x)^2
 \end{multline*}
Третий:
 \begin{multline*}
\d^3\arctg x=\d\left(-\frac{2x}{(1+x^2)^2}\cdot (\d x)^2\right)=\\=
\d\left(-\frac{2x}{(1+x^2)^2}\right)\cdot (\d x)^2=\\=
\frac{-2(1+x^2)^2+2x\cdot 2(1+x^2)\cdot 2x}{(1+x^2)^4}\cdot \d x\cdot (\d
x)^2=\\= \frac{-2(1+x^2)+2x\cdot 2\cdot 2x}{(1+x^2)^3}\cdot (\d x)^3=\\=
\frac{-2-2x^2+8x^2}{(1+x^2)^3}\cdot (\d x)^3=\\=
2\cdot\frac{3x^2-1}{(1+x^2)^3}\cdot (\d x)^3
 \end{multline*}
Теперь можно выписать дифференциалы функции $f$:
 \begin{align*}
\d^0 f(x)[p]&=f(x)=\arctg x\\
\d^1 f(x)[p]&=\frac{1}{1+x^2}\cdot \d x\Big|_{\d
x=p}=\\&=\frac{1}{1+x^2}\cdot p \\
\d^2 f(x)[p]&=-\frac{2x}{(1+x^2)^2}\cdot (\d x)^2\Big|_{\d
x=p}=\\&=-\frac{2x}{(1+x^2)^2}\cdot p^2 \\
\d^3 f(x)[p]&=2\cdot\frac{3x^2-1}{(1+x^2)^3}\cdot (\d x)^3\Big|_{\d
x=p}=\\&=2\cdot\frac{3x^2-1}{(1+x^2)^3}\cdot p^3
 \end{align*}
Формула Тейлора-Лагранжа (порядка 2) для этой функции имеет вид:
\begin{multline*}
\arctg (x+p)=\arctg x+\frac{1}{1+x^2}\cdot
p-\frac{1}{2}\cdot\frac{2x}{(1+x^2)^2}\cdot
p^2+\\+\frac{1}{3}\cdot\frac{3(x+\theta p)^2-1}{(1+(x+\theta p)^2)^3}\cdot p^3
\end{multline*}
Можно выписать ее в конкретной точке, например, в $x=1$:
\begin{multline*}
\arctg (1+p)=\frac{\pi}{4}+\frac{1}{2}\cdot p-\frac{1}{4}\cdot
p^2+\\+\frac{1}{3}\cdot\frac{3(1+\theta p)^2-1}{(1+(1+\theta p)^2)^3}\cdot p^3
\end{multline*}
\eex

\end{multicols}\noindent\rule[10pt]{160mm}{0.1pt}

\subsection{Ряд Тейлора}

Формулы Тейлора
\eqref{Taylor-Lagrange}-\eqref{Taylor-Cauchy}-\eqref{Taylor-Peano} подсказывают
следующую конструкцию.

 \bit{
\item[$\bullet$] Если $f$ -- гладкая функция на интервале $(a,b)$, то любой
точке $x_0\in(a,b)$ можно поставить в соответствие степенной ряд
 \beq\label{ryad-Teilora}
\sum_{n=0}^\infty \frac{f^{(n)}(x_0)}{n!}\cdot (x-x_0)^n
 \eeq
Этот ряд называется {\it рядом Тейлора функции $f$ в точке $x_0$}.
 }\eit

 \bit{
\item[$\bullet$] В частном случае, когда $x_0=0$ ряд Тейлора называется {\it
рядом Маклорена}:
 \beq
\sum_{n=0}^\infty \frac{f^{(n)}(0)}{n!}\cdot x^n
 \eeq
 }\eit

От ряда Тейлора естественно ожидать, что он будет сходиться к порождающей его
функции $f$ хотя бы в некоторой окрестности $(x_0-\delta;x_0+\delta)$ точки
$x_0$, то есть что выполняется следующее тождество, называемое {\it разложением
Тейлора} функции $f$ в окрестности $(x_0-\delta;x_0+\delta)$:
 \beq\label{razlozhenie-Teilora}
f(x)=\sum_{n=0}^\infty \frac{f^{(n)}(0)}{n!}\cdot (x-x_0)^n,\quad x\in
(x_0-\delta;x_0+\delta)
 \eeq
В частном случае, когда $x_0=0$, оно называется {\it разложением Маклорена} (в
окрестности $(-\delta;\delta)$):
 \beq\label{razlozhenie-Mclaurent}
f(x)=\sum_{n=0}^\infty \frac{f^{(n)}(0)}{n!}\cdot x^n,\quad x\in
(-\delta;\delta)
 \eeq

Часто так и бывает (мы приведем примеры такого рода на
с.\pageref{EXS:razlozh-Maclaurent-1}). Но не всегда. В следующих примерах мы
увидим, что ряд Тейлора функции $f$ необязательно сходится, а если сходится, то
его сумма необязательно совпадает с функцией $f$.

\noindent\rule{160mm}{0.1pt}\begin{multicols}{2}

\paragraph{Контрпримеры.}

\bex{\bf Ряд Маклорена, расходящийся везде, кроме точки $0$.} Рассмотрим
последовательность
$$
a_n=(n!)^2
$$
По теореме Бореля \ref{TH:Borel} существует гладкая функция $f$ на $\R$ такая,
что
$$
\forall n\in\Z_+\qquad f^{(n)}(0)=a_n=(n!)^2
$$
Ее ряд Маклорена имеет вид
$$
\sum_{n=0}^\infty \frac{f^{(n)}(0)}{n!}\cdot x^n=\sum_{n=0}^\infty
\frac{(n!)^2}{n!}\cdot x^n=\sum_{n=0}^\infty n!\cdot x^n
$$
В примере \ref{ex-21.1.5} мы уже рассматривали такой степенной ряд, и поняли,
что он сходится только в точке $0$. \eex

\bex{\bf Ряд Маклорена сходящийся всюду, но не к порождающей его
функции.}\label{EX:Mclaurent-shoditsya-ne-k-porozhd-func} Рассмотрим функцию из
примера \ref{EX:gladkaya-func-0<f(x)<1}, то есть гладкую на $\R$ функцию $f$ со
свойствами:
$$
\begin{cases}
f(x)=0, & x\le 0 \\
f(x)>0,& x>0
\end{cases}
$$
Поскольку на левой полупрямой эта функция нулевая, в точке $0$ все ее
производные равны нулю, поэтому ряд Маклорена будет нулевой:
$$
\sum_{n=0}^\infty
\overbrace{\frac{f^{(n)}(0)}{n!}}^{\tiny\begin{matrix}0\\
\| \end{matrix}}\cdot x^n=\sum_{n=0}^\infty 0\cdot x^n
$$
Сумма такого ряда, понятное дело, тоже будет нулевой,
$$
\sum_{n=0}^\infty \frac{f^{(n)}(0)}{n!}\cdot x^n=0
$$
и поэтому она не может совпадать с $f$ ни на каком интервале
$(-\delta,+\delta)$: тождество
$$
\sum_{n=0}^\infty \frac{f^{(n)}(0)}{n!}\cdot x^n=0=f(x), \qquad
x\in(-\delta,+\delta),
$$
невозможно потому что при $x>0$ функция $f$ отлична от нуля.
 \eex
\end{multicols}\noindent\rule[10pt]{160mm}{0.1pt}

\paragraph{Всякий сходящийся степенной ряд есть ряд Тейлора своей суммы.}

Несмотря на отмеченное нами нерегулярное поведение ряда Тейлора, эта
конструкция оказывается центральным примером в теории степенных рядов из-за
следующей теоремы:

\btm\label{TH:step-ryad=ryad-teilora} Если функция $f$ является суммой
какого-то степенного ряда в некоторой окрестности его центра
 \beq
f(x)=\sum_{n=0}^\infty c_n\cdot (x-x_0)^n,\qquad x\in(x_0-\delta,x_0+\delta),
 \eeq
то этот степенной ряд является рядом Тейлора функции $f$ в точке $x_0$:
$$
c_n=\frac{f^{(n)}(x_0)}{n!}
$$
 \etm
\begin{proof} В силу следствия \ref{cor-21.2.5}, функция $f$ является гладкой на
интервале $(x_0-\delta, x_0+\delta)$, и по формуле \eqref{21.2.12},
\begin{multline*}
f^{(k)}(x)= \sum_{n=k}^\infty \frac{n!}{(n-k)!}\cdot c_n\cdot (x-x_0)^{n-k}=\\=
\frac{k!}{(k-k)!}\cdot c_k\cdot (x-x_0)^0+ \frac{(k+1)!}{(k+1-k)!}\cdot
c_{k+1}\cdot (x-x_0)^1+ \frac{(k+2)!}{(k+2-k)!}\cdot c_{k+2}\cdot
(x-x_0)^2+...=\\= k!\cdot c_k+ \frac{(k+1)!}{1!}\cdot c_{k+1}\cdot (x-x_0)+
\frac{(k+2)!}{2!}\cdot c_{k+2}\cdot (x-x_0)^2+..., \qquad x\in (x_0-R, x_0+R)
\end{multline*}
Подставив сюда $x=x_0$, получим
$$
f^{(k)}(x_0)= k!\cdot c_k+ \frac{(k+1)!}{1!}\cdot c_{k+1}\cdot 0+
\frac{(k+2)!}{2!}\cdot c_{k+2}\cdot 0+...=k!\cdot c_k
$$
То есть, $c_k=\frac{f^{(k)}(x_0)}{k!}$ \end{proof}

\brem Из этой теоремы следует, что если степенной ряд $\sum_{n=0}^\infty
c_n\cdot (x-x_0)^n$ сходится хотя бы в одной точке $x$, отличной от его центра
$x_0$, то автоматически он является рядом Тейлора для некоторой функции $f$, а
именно, для своей суммы $f(x)=\sum_{n=0}^\infty c_n\cdot (x-x_0)^n$ (которая
будет определена как минимум в окрестности точки $x_0$ радиуса
$\delta=|x-x_0|$). \erem

\paragraph{Сходимость ряда Тейлора.}

Обсудим теперь, когда все-таки ряд Тейлора сходится к порождающей его функции.
Удобное достаточное условие для этого выглядит так:

\bprop\label{TH:f^n<infty=>f-analit} Пусть функция $f$ определена в
$\delta$-окрестности некоторой точки $x_0$, бесконечно гладкая на
$(x_0-\delta;x_0+\delta)$, и имеет ограниченные в совокупности производные:
 \beq\label{sup_n-sup_x-|f^n(x)|<infty}
\sup_{n\in\Z_+}\sup_{x\in(x_0-\delta;x_0+\delta)} |f^{(n)}(x)|=M<\infty
 \eeq
Тогда на интервале $(x_0-\delta;x_0+\delta)$ функция $f$ удовлетворяет
тождеству:
 \beq
f(x)= \sum_{n=0}^\infty \frac{f^{(n)}(x_0)}{n!}\cdot (x-x_0)^n
 \eeq
 \eprop
 \bpr
Выпишем формулу Тейлора-Лагранжа \eqref{Taylor-Lagrange} для $f$: для любой
точки $x\in(x_0-\delta;x_0+\delta)$ найдется точка $\xi$, лежащая между $x_0$ и
$x$, такая что выполняется равенство
 $$
f(x)=\sum_{k=0}^n \frac{f^{(k)}(x_0)}{k!}\cdot (x-x_0)^k+
\frac{f^{(n+1)}(\xi)}{(n+1)!}\cdot (x-x_0)^{n+1}
 $$
Из него мы получаем:
$$
f(x)-\sum_{k=0}^n \frac{f^{(k)}(x_0)}{k!}\cdot (x-x_0)^k=
\frac{f^{(n+1)}(\xi)}{(n+1)!}\cdot (x-x_0)^{n+1}
$$
$$
\Downarrow
$$
$$
\left|f(x)-\sum_{k=0}^n \frac{f^{(k)}(x_0)}{k!}\cdot (x-x_0)^k\right|=
\left|\frac{f^{(n+1)}(\xi)}{(n+1)!}\cdot (x-x_0)^{n+1}\right|\le
M\cdot\frac{\delta^n}{(n+1)!}\underset{n\to\infty}{\longrightarrow}0
$$
$$
\Downarrow
$$
$$
\sum_{k=0}^n \frac{f^{(k)}(x_0)}{k!}\cdot
(x-x_0)^k\underset{n\to\infty}{\longrightarrow}f(x)
$$
 \epr

Теперь мы можем рассмотреть примеры.

\noindent\rule{160mm}{0.1pt}\begin{multicols}{2}

\paragraph{Стандартные разложения Маклорена.}

\bexs\label{EXS:razlozh-Maclaurent-1} Следующие разложения Маклорена
справедливы в интервале $|x|<1$:
 \begin{align}
&\frac{1}{1-x}=1+x+x^2+...+x^k+... =\notag \\ &\qquad =\sum_{n=0}^\infty x^n, \label{mclaurent-1/(1-x)}\\
&\frac{1}{1+x}=1-x+x^2+...+ (-1)^k\cdot x^k+...=\notag \\ &\qquad=
\sum_{n=0}^\infty (-1)^n \cdot x^n \label{mclaurent-1/(1+x)}\\
&\ln (1+x) =x-\frac{x^2}{2}+\frac{x^3}{3}+...+(-1)^k \cdot \frac{x^k}{k}+... =
\notag \\
&\qquad=\sum_{n=1}^\infty (-1)^{n-1}\cdot \frac{x^n}{n} \label{mclaurent-ln(1+x)}\\
&\frac{1}{1+x^2} =1-x^2+x^4+...+ (-1)^k\cdot x^{2k}+...=\notag \\
&\qquad= \sum_{n=0}^\infty (-1)^n \cdot x^{2n} \label{mclaurent-1/(1+x^2)}\\
&\arctg x=x-\frac{x^3}{3}+\frac{x^5}{5}-...+ (-1)^k\cdot\frac{x^{2k+1}}{2k+1}+\notag \\
&\qquad+...=\sum_{n=0}^\infty (-1)^n\cdot \frac{x^{2n+1}}{2n+1} \label{mclaurent-arctg(x)} \\
& (1+x)^\alpha =1+\alpha \cdot x+\alpha (\alpha-1) \cdot \frac{x^2}{2!}+... \notag \\
&\qquad ...+ \alpha (\alpha-1) ... (\alpha-k+1) \cdot \frac{x^k}{k!}+... =\notag \\
&\qquad= 1+\sum_{n=1}^\infty \alpha (\alpha-1) ... (\alpha-n+1) \cdot
\frac{x^n}{n!} \label{mclaurent-(1+x)^a}
 \end{align}
 \eexs
\bpr Тождества \eqref{mclaurent-1/(1-x)}-\eqref{mclaurent-arctg(x)} мы уже
доказали на с.\pageref{21.3.11}, поэтому их доказывать нет необходимости, нужно
только заметить, что в силу теоремы \ref{TH:step-ryad=ryad-teilora} они будут
разложениями Маклорена (поскольку центром ряда здесь будет точка 0).

Таким образом, интерес в этом списке представляет лишь формула
\eqref{mclaurent-(1+x)^a}. Для ее доказательства заметим, что функция
$f(x)=(1+x)^\alpha$ имеет производные
$$
f^{(n)}(x)=\alpha(\alpha-1)... (\alpha-n+1) \cdot (1+x)^{\alpha-n}
$$
Применим к этой функции теорему Тейлора-Коши \ref{TH:Taylor-Cauchy}: для любого
$x\in(-1,1)$ и любого $n\in\N$ найдется точка $\xi$, лежащая между $0$ и $x$
такая, что
 \begin{multline}\label{taylor-cauchy-(1+x)^a}
f(x)=\sum_{k=0}^n \frac{f^{(k)}(0)}{k!}\cdot x^k+\\+
\frac{f^{(n+1)}(\xi)}{n!}\cdot (x-\xi)^n\cdot x
 \end{multline}
Заметим две вещи: во-первых,
$$
1-|x|\le 1-|\xi|\le|1+\xi|\le1+|\xi|\le 1+|x|
$$
$$
\Downarrow
$$
 \beq\label{otsenka-|1+xi|^(alpha-1)}
|1+\xi|^{\alpha-1}\le \max\left\{(1-|x|)^{\alpha-1};(1+|x|)^{\alpha-1}\right\}
 \eeq
И, во-вторых,
 \begin{multline}\label{otsenka-|x-xi|/|1+xi|}
\frac{|x-\xi|}{|1+\xi|}\le
 \frac{|x|-|\xi|}{1-|\xi|}=
\frac{1-|\xi|-(1-|x|)}{1-|\xi|}=\\=1-\frac{1-|x|}{1-|\xi|} \le
1-\frac{1-|x|}{1-0}=|x|
  \end{multline}
Теперь остаток ряда Маклорена можно оценить так:
 \begin{multline*}
\left|f(x)-\sum_{k=0}^n \frac{f^{(k)}(0)}{k!}\cdot
x^k\right|=\eqref{taylor-cauchy-(1+x)^a}=\\=
\left|\frac{f^{(n+1)}(\xi)}{n!}\cdot (x-\xi)^n\cdot x\right|=\\=
\left|\frac{\alpha(\alpha-1)... (\alpha-n) \cdot
(1+\xi)^{\alpha-n-1}}{n!}\right|\cdot\\ \cdot |x-\xi|^n\cdot |x|=\\=
\left|\alpha\l 1-\frac{\alpha}{2}\r... \l 1-\frac{\alpha}{n}\r\right| \cdot\\
\cdot |1+\xi|^{\alpha-n-1}\cdot |x-\xi|^n\cdot |x|=\\=
\left|\alpha\l 1-\frac{\alpha}{2}\r... \l 1-\frac{\alpha}{n}\r\right| \cdot\\
\cdot |1+\xi|^{\alpha-1} \cdot \bigg(\frac{|x-\xi|}{|1+\xi|}\bigg)^n\cdot
|x|\le\\ \le \eqref{otsenka-|1+xi|^(alpha-1)},
\eqref{otsenka-|x-xi|/|1+xi|}\le\\
\le \left|\alpha\l 1-\frac{\alpha}{2}\r... \l 1-\frac{\alpha}{n}\r\right|
\cdot\\ \cdot \max\left\{(1-|x|)^{\alpha-1};(1+|x|)^{\alpha-1}\right\}\cdot
|x|^{n+1}
\end{multline*}
Обозначим последнюю величину $M_n$
 \begin{multline*}
M_n=\left|\alpha\l 1-\frac{\alpha}{2}\r... \l 1-\frac{\alpha}{n}\r\right| \cdot\\
\cdot \max\left\{(1-|x|)^{\alpha-1};(1+|x|)^{\alpha-1}\right\}\cdot |x|^{n+1}.
 \end{multline*}
Нам нужно показать, что $M_n$ стремится к нулю при $n\to\infty$. Зафиксируем
какое-нибудь число $q$ такое, что
$$
|x|<q<1
$$
(такое $q$ найдется, поскольку $|x|<1$). Теперь мы получаем:
$$
\left|1-\frac{\alpha}{n+1}\right|\underset{n\to\infty}{\longrightarrow} 1
$$
$$
\Downarrow
$$
$$
\frac{M_{n+1}}{M_n}=\left|1-\frac{\alpha}{n+1}\right|\cdot
|x|\underset{n\to\infty}{\longrightarrow} |x|<q
$$
$$
\Downarrow
$$
$$
\exists N\in\N\quad \forall n\ge N\quad \frac{M_{n+1}}{M_n}<q
$$
$$
\Downarrow
$$
$$
\forall k\in\Z_+\quad M_{N+k}<M_N\cdot q^k
$$
$$
\Downarrow
$$
$$
\sum_{k=0}^\infty M_{N+k}\le \sum_{k=0}^\infty M_N\cdot q^k<\infty
$$
$$
\Downarrow
$$
$$
M_n\underset{n\to\infty}{\longrightarrow} 0
$$

\epr

\bexs Следующие разложения Маклорена справедливы на всей числовой прямой
$x\in\R$:
 \begin{align}
& \label{21.3.8} e^x =1+x+\frac{x^2}{2!}+...+\frac{x^k}{k!}+... =\notag \\
&\qquad= \sum_{n=0}^\infty \frac{x^n}{n!}\\
& \label{21.3.5} \sin x=x-\frac{x^3}{3!}+\frac{x^5}{5!}-...+(-1)^n \cdot
\frac{x^{2k+1}}{(2k+1)!}+...=\notag \\
&\qquad= \sum_{n=0}^\infty (-1)^n\cdot\frac{x^{2n+1}}{(2n+1)!}\\
& \label{21.3.6} \cos x =1-\frac{x^2}{2!}+\frac{x^4}{4!}-...+(-1)^k
\cdot\frac{x^{2k}}{(2k)!}+...=\notag \\ &\qquad=\sum_{n=0}^\infty (-1)^n\cdot
\frac{x^{2n}}{(2n)!}
 \end{align}
\eexs

\bpr

1. Из формулы для производных функции $e^x$
$$
\frac{\d^n}{(\d x)^n}e^x=e^x
$$
видно, что она удовлетворяет условию \eqref{sup_n-sup_x-|f^n(x)|<infty} на
любом интервале $(-R,R)$, $R>0$:
 \begin{multline*}
\sup_{n\in\Z_+}\sup_{x\in(-R,R)} |f^{(n)}(x)|=\\=
\sup_{n\in\Z_+}\sup_{x\in(-R,R)} |e^x|=e^R<\infty
 \end{multline*}
В силу предложения \ref{TH:f^n<infty=>f-analit}, это означает, что функцию
$f(x)=e^x$ представить сходящимся на интервале $(-R,R)$ рядом \eqref{21.3.8}.
Поскольку это верно для любого $R>0$, мы получаем, что $f(x)=e^x$ формула
\eqref{21.3.8} верна при любом $x\in\R$.

2. Производные функции $\sin$
$$
\frac{\d^n}{(\d x)^n}\sin x=\begin{cases}
 \sin x,& n=4k, \quad k\in\Z_+ \\
 \cos x,& n=4k+1,\quad  k\in\Z_+ \\
 -\sin x,& n=4k+2, \quad k\in\Z_+ \\
 -\cos x,& n=4k+3, \quad k\in\Z_+
\end{cases}
$$
ограничены в совокупности на $\R$ (и как следствие, на любом интервале
$(-R,R)$):
 $$
\sup_{n\in\Z_+}\sup_{x\in(-R,R)} |\sin^{(n)}(x)|\le 1<\infty
 $$
В силу предложения \ref{TH:f^n<infty=>f-analit}, это означает, функцию $\sin$
можно представить сходящимся на интервале $(-R,R)$ рядом \eqref{21.3.5}.
Поскольку это верно для любого $R>0$, мы получаем, что формула \eqref{21.3.5}
верна при любом $x\in\R$.

3. Производные функции $\cos$
$$
\cos^{(n)}x=\begin{cases}
 \cos x,& n=4k, \quad k\in\Z_+ \\
 -\sin x,& n=4k+1,\quad  k\in\Z_+ \\
 -\cos x,& n=4k+2, \quad k\in\Z_+ \\
 \sin x,& n=4k+3, \quad k\in\Z_+
\end{cases}
$$
ограничены в совокупности на $\R$. Поэтому функция $\cos$ удовлетворяет условию
\eqref{sup_n-sup_x-|f^n(x)|<infty} на любом интервале $(-R,R)$, $R>0$:
 $$
\sup_{n\in\Z_+}\sup_{x\in(-R,R)} |\sin^{(n)}(x)|\le 1<\infty
 $$
В силу предложения \ref{TH:f^n<infty=>f-analit}, $\cos$ можно представить
сходящимся на интервале $(-R,R)$ рядом \eqref{21.3.6}. Поскольку это верно для
любого $R>0$, равенство \eqref{21.3.6} будет выполняться для любого $x\in\R$.

\epr

\paragraph{Следствия формулы Тейлора-Коши.}

\end{multicols}\noindent\rule[10pt]{160mm}{0.1pt}

\subsection{Аналитические функции}

 \bit{
\item[$\bullet$] Функция $f$ называется {\it  аналитической на интервале
$(a,b)$}\index{функция!аналитическая!на множестве}\label{DEF:analitich-func},
если она определена на $(a,b)$ и удовлетворяет следующим равносильным условиям:
 \bit{
\item[(i)] у любой точки $x_0\in (a,b)$ найдется окрестность $(x_0-\delta;
x_0+\delta)$, в которой функция $f$ представима в виде суммы степенного ряда:
 \beq\label{21.3.1}
f(x)= \sum_{n=0}^\infty c_n\cdot (x-x_0)^n, \qquad x\in (x_0-\delta;
x_0+\delta)
 \eeq
иными словами, в окрестности  $(x_0-\delta; x_0+\delta)$ функция $f$ должна
быть сдвигом производящей функции некоторой аналитической последовательности
$c_n$:
 \beq\label{21.3.2}
f(x)=\Gen_c(x-x_0), \qquad x\in (x_0-\delta; x_0+\delta)
 \eeq

\item[(ii)] у любой точки $x_0\in (a,b)$ найдется окрестность $(x_0-\delta;
x_0+\delta)$, в которой функция $f$ (является гладкой и) представима в виде
суммы своего ряда Тейлора:
 \beq
f(x)= \sum_{n=0}^\infty \frac{f^{(n)}(x_0)}{n!}\cdot (x-x_0)^n, \qquad x\in
(x_0-\delta; x_0+\delta) \label{21.3.3}
 \eeq
иными словами, в окрестности  $(x_0-\delta; x_0+\delta)$ функция $f$ должна
быть сдвигом производящей функции последовательности своих коэффициентов
Тейлора:
 \beq \label{21.3.2-*}
f(x)=\Gen_c(x-x_0),\qquad c_n=\frac{f^{(n)}(x_0)}{n!}, \qquad x\in (x_0-\delta;
x_0+\delta)
 \eeq
 }\eit

\item[$\bullet$] Функция $f$ называется {\it целой}\index{функция!целая}, если
она определена всюду на $\R$ и удовлетворяет следующим равносильным условиям:
 \bit{
\item[(i)$^*$] $f$ представима в виде суммы всюду сходящегося на $\R$
степенного ряда:
 \beq\label{tselaya-func}
f(x)=\sum_{n=0}^\infty c_n\cdot x^n,\qquad x\in\R
 \eeq
иными словами, функция $f$ всюду на $\R$ должна совпадать с производящей
функции некоторой аналитической последовательности $c_n$:
 \beq\label{tselaya-func-*}
f(x)=\Gen_c(x), \qquad x\in\R
 \eeq

\item[(ii)$^*$] $f$ является гладкой и всюду на $\R$ совпадает с суммой своего
ряда Маклорена:
 \beq\label{tselaya-func-2}
f(x)=\sum_{n=0}^\infty \frac{f^{(n)}(0)}{n!}\cdot x^n,\qquad x\in\R
 \eeq
иными словами, функция $f$ всюду на $\R$ должна совпадать с производящей
функции некоторой последовательности своих коэффициентов Маклорена:
 \beq\label{tselaya-func-**}
f(x)=\Gen_c(x),\qquad c_n=\frac{f^{(n)}(0)}{n!}, \qquad x\in\R
 \eeq

 }\eit
 }\eit

\bpr Равносильность условий (i) и (ii) (а также (i)$^*$ и (ii)$^*$) следует из
теоремы \ref{TH:step-ryad=ryad-teilora}.  \epr

\noindent\rule{160mm}{0.1pt}\begin{multicols}{2}

Приведем сразу два типичных примера.

\bex{\bf Аналитическая, но не целая функция.}\label{EX:1/x-analit} Покажем, что
функция
$$
f(x)=\frac{1}{x}
$$
является аналитической (на области своего определения $x\in\R\setminus\{0\}$).
Зафиксируем точку $x_0\ne 0$. Тогда при $|x-x_0|<|x_0|$ справедливо разложение
в степенной ряд:
 \begin{multline*}
f(x)=\frac{1}{x}=\frac{1}{x_0+(x-x_0)}=\frac{1}{x_0}\cdot\frac{1}{1+\underbrace{\frac{x-x_0}{x_0}}_{\scriptsize\begin{matrix}
\text{по модулю}\\ \text{меньше}\\ \text{единицы}\end{matrix}}}=\\=
\eqref{mclaurent-1/(1+x)}=\frac{1}{x_0}\cdot\sum_{n=0}^\infty(-1)^n\cdot\l\frac{x-x_0}{x_0}\r^n=\\=
\sum_{n=0}^\infty \frac{(-1)^n}{x_0^{n+1}}\cdot(x-x_0)^n
 \end{multline*}
Таким образом, функция $f$ аналитическая. Ее, однако, нельзя считать целой хотя
бы потому что она определена не на всей прямой $\R$. Более того, ее нельзя
продолжить до целой функции на $\R$, потому что по теореме \ref{tm-21.2.2}
целая функция должна быть непрерывна на $\R$, а $f$ не продолжается до
непрерывной функции на $\R$ из-за соотношения
$$
\lim_{x\to 0}f(x)=\infty.
$$
\eex

\begin{ex}\label{EX:ex-21.1.6-analit} В примере \ref{ex-21.1.6} мы убедились, что ряд
$\sum_{n=1}^\infty \frac{x^n}{n!}$ сходится всюду на $\R$. В соответствии с
нашим определением, это означает, что порождаемая им функция
$$
f(x)=\sum_{n=1}^\infty \frac{x^n}{n!}
$$
является целой.
\end{ex}

\end{multicols}\noindent\rule[10pt]{160mm}{0.1pt}

\paragraph{Операции над аналитическими функциями.}

В \ref{SEC-nepr-funct}\ref{SUBSEC:oper-nad-nepr-func} главы \ref{ch-cont-f(x)}
мы показывали, что сумма, разность, произведение, частное и композиция
непрерывных функций снова является непрерывной функцией. Те же самые свойства
для дифференцируемых функций доказывались нами в
\ref{SEC-proizv-i-diff-functsii}\ref{SEC-pravila-proizvodnoi} главы
\ref{ch-f'(x)}: сумма, разность, произведение, частное и композиция
дифференцируемых функций снова является дифференцируемой функцией. Оказывается,
что то же самое верно и для аналитических функций.

\begin{tm}[\bf о композиции аналитических функций]\label{tm-21.3.6}
Если $f$ -- аналитическая функция на множестве $U$, а $g$ -- аналитическая
функция на множестве $V\supseteq f(U)$, то их композиция $h(x)=g(f(x))$ --
аналитическая функция на множестве $U$.
\end{tm}
\bpr Зафиксируем точку $x_0\in U$ и положим
$$
F(x)=f(x)-f(x_0),\qquad G(y)=g(y+f(x_0))
$$
Композиция этих функций тоже равна $h$:
$$
G(F(x))=g(f(x)-f(x_0)+f(x_0))=g(f(x))=h(x),\qquad x\in U
$$
С другой стороны, эти функции тоже аналитические, поэтому
$$
F(x)=\Gen_{b}(x-x_0),\qquad G(y)=\Gen_a(y),\qquad x\in U_\delta(x_0),\ y\in
U_\e(0)
$$
для некоторых аналитических последовательностей $a$ и $b$. При этом
$$
b_0=\Gen_{b}(0)=F(x_0)=f(x_0)-f(x_0)=0,
$$
то есть порядок аналитической последовательности $b\in\mathcal{A}$ отличен от
нуля:
$$
\omega(b)\ge 1.
$$
Значит, определена композиция последовательностей $a\circ b$, которая по
теореме \ref{TH:compos-analit-posl} тоже является аналитической: $a\circ
b\in\mathcal{A}$. А производящие функции этих последовательностей связаны
формулой \eqref{Gen(a-circ-b)=Gen(a)-circ-Gen(b)}, из которой получаем:
 $$
h(x)=(G\circ F)(x)=G(F(x))=
\Gen_a\big(\Gen_b(x-x_0)\big)=\eqref{Gen(a-circ-b)=Gen(a)-circ-Gen(b)}=\Gen_{a\circ
b}(x-x_0),\qquad x\in U_\delta(x_0)
 $$
Таким образом, в некоторой окрестности $U_\delta(x_0)$ точки $x_0$ функция $h$
совпадает с производящей функцией некоторой аналитической последовательности
($a\circ b$). Поскольку это верно для произвольной точки $x_0\in U$, функция
$h$ является аналитической на множестве $U$.
 \epr

\begin{tm}[\bf об арифметических операциях с аналитическими
функциями]\label{tm-21.3.5} Если $f$ и $g$ -- аналитические функции на
множестве $U$, то следующие функции -- тоже аналитические на множестве $U$:
$$
f+g, \quad f-g, \quad C\cdot f, \quad f\cdot g, \quad
$$
($C$ -- произвольная константа). А функция
$$
\frac{g}{f}
$$
-- аналитическая на множестве $U\setminus \{x: \ f(x)=0\}$.
\end{tm}
\bpr В первых четырех случаях утверждение следует из формул
\eqref{Gen(a+b)(x)=Gen(a)(x)+Gen(b)(x)}-\eqref{Gen(a*b)(x)=Gen(a)(x)-cdot-Gen(b)(x)}:
для произвольной точки $x_0\in U$ выберем окрестность $(x_0-\delta,x_0+\delta)$
и аналитические последовательности $a$ и $b$ так, чтобы
$$
f(x)=\Gen_a(x-x_0),\qquad g(x)=\Gen_b(x-x_0),\qquad x\in(x_0-\delta,x_0+\delta)
$$
Тогда
$$
(f+g)(x)=f(x)+g(x)=\Gen_a(x-x_0)+\Gen_b(x-x_0)=\eqref{Gen(a+b)(x)=Gen(a)(x)+Gen(b)(x)}=\Gen_{a+b}(x-x_0),\qquad
x\in(x_0-\delta,x_0+\delta)
$$
То есть в некоторой окрестности $(x_0-\delta,x_0+\delta)$ точки $x_0$ функция
$f+g$ совпадает с производящей функцией некоторой  аналитической
последовательности ($a+b$). Поскольку это верно для произвольной точки $x_0\in
U$, функция $f+g$ является аналитической на множестве $U$. То же самое верно
для $f-g$, $C\cdot f$, $f\cdot g$ (в последнем случае используется формула
\eqref{Gen(a*b)(x)=Gen(a)(x)-cdot-Gen(b)(x)}).

Теперь для доказательства аналитичности $\frac{g}{f}$ достаточно доказать
аналитичность обратной функции $\frac{1}{f}$. Это делается с помощью теоремы
\ref{tm-21.3.6}: поскольку функция $G(y)=\frac{1}{y}$ аналитична в силу примера
\ref{EX:1/x-analit}, ее композиция с $f$
$$
G(f(x))=\frac{1}{f(x)}
$$
-- тоже аналитична.\epr

\noindent\rule{160mm}{0.1pt}\begin{multicols}{2}

\bex{\bf Всюду определенная аналитическая, но не целая функция.} Функция
$$
f(x)=\frac{1}{1+x^2}
$$
является аналитический по теореме \ref{tm-21.3.5}, как отношение двух
аналитических функций. Но она не является целой, потому что ее ряд Маклорена
$$
\frac{1}{1+x^2}=\eqref{mclaurent-1/(1+x)}=\sum_{n=0}^\infty (-1)^n\cdot x^{2n}
$$
сходится только при $|x|<1$. \eex

\bex{\bf Гладкая, но не ана\-ли\-ти\-ческая
функция.}\label{EX-gladkaya-ne-analit} Рассмотрим функцию из примера
\ref{EX:gladkaya-func-0<f(x)<1}, то есть гладкую на $\R$ функцию $f$ со
свойствами:
$$
\begin{cases}
f(x)=0, & x\le 0 \\
f(x)>0,& x>0
\end{cases}
$$
В примере \ref{EX:Mclaurent-shoditsya-ne-k-porozhd-func} мы отмечали, что ряд
Маклорена этой функции
$$
\sum_{n=0}^\infty
\overbrace{\frac{f^{(n)}(0)}{n!}}^{\tiny\begin{matrix}0\\
\| \end{matrix}}x^n=0,
$$
не сходится к ней ни в какой окрестности нуля $(-\delta,+\delta)$, потому что
при $x>0$ наша функция отлична от нуля. Отсюда следует, что $f$ не может быть
аналитической в точке $x=0$.
 \eex

\end{multicols}\noindent\rule[10pt]{160mm}{0.1pt}

\section{Приложения степенных рядов}\label{SEC:prilozh-step-ryadov}

\subsection{Доказательство избыточности Аксиомы степеней.}\label{PR:teor-o-step-otobr}

Читатель помнит, наверное, что степенное отображение
$$
(a,b)\mapsto a^b,
$$
вводилось нами в общем случае (то есть для необязательно целых показателей $b$)
в \ref{SUBSEC-elem-funktsii} главы \ref{ch-ELEM-FUNCTIONS} с помощью Аксиомы
степеней B1 на с.\pageref{TH-o-step-otobr}, избыточность которой мы пообещали
доказать позже. Теперь, наконец, мы можем выполнить данное тогда обещание.
Собственно доказательству этой аксиомы (на с.\pageref{PROOF:TH-o-step-otobr}),
мы предпошлем 16 вспомогательных утверждений (леммы
\ref{LM:exp}-\ref{LM:monot-a^x}). Все это время читателю предлагается делать
вид, что до сих пор он не глядел на текст в двух колонках и, как следствие не
знает ничего об отображении $(a,b)\mapsto a^b$, по крайней мере, для случая
нецелого $b$ (случай целых степеней можно считать известным, поскольку для него
все необходимые факты были аккуратно доказаны нами еще в главе \ref{ch-R&N} в
теореме о степенном отображении \ref{TH-o-step-otobr-malaya}). Только после
приведенного здесь доказательства Аксиомы степеней все сказанное нами ранее на
этот счет вступит в законную силу, и этой цели посвящен настоящий раздел.

\noindent\rule{160mm}{0.1pt}\begin{multicols}{2}

\paragraph{Функция $\exp$.}

\blm\label{LM:exp} Формула
 \begin{align}\label{DF:exp}
&\exp x=\sum_{n=0}^\infty \frac{x^n}{n!}
\end{align}
определяет функцию $\exp$ всюду на прямой $\R$.
 \elm
\bpr Радиус сходимости этого степенного ряда равен бесконечности:
 \begin{multline*}
R=\lim_{n\to\infty}\left|\frac{c_n}{c_{n+1}}\right|=
\lim_{n\to\infty}\left|\frac{(n+1)!}{n!}\right|=\\=
\lim_{n\to\infty}(n+1)=\infty
 \end{multline*}
поэтому ряд сходится всюду на $\R$. \epr

 \blm\label{LM:exp'} Функция $\exp$, определенная рядом \eqref{DF:exp},
имеет производную:
 \begin{align}\label{exp'x}
&\frac{\d}{\d x}\exp x=\exp x
\end{align}
 \elm
\bpr Продифференцируем ряд \eqref{DF:exp} с помощью теоремы \ref{tm-21.2.4}:
  \begin{multline*}
\frac{\d}{\d x}\exp x=\frac{\d}{\d x}\sum_{n=0}^\infty
\frac{x^n}{n!}=\sum_{n=0}^\infty \frac{\d}{\d x}\left(\frac{x^n}{n!}\right)=\\=
\sum_{n=1}^\infty \frac{x^{n-1}}{(n-1)!}=\sum_{k=0}^\infty \frac{x^k}{k!}=\exp
x
  \end{multline*}
\epr

 \blm\label{LM:exp(x+y)} Функция $\exp$, определенная рядом \eqref{DF:exp},
 удовлетворяет следующему равенству и двум тождествам:
 \begin{align}
& \exp 0=1, \label{exp(0)} \\
& \exp(-x)=\frac{1}{\exp x} \label{exp(-x)} \\
& \exp(x+y)=\exp x\cdot \exp y,\label{exp(x+y)}
 \end{align}
 \elm
\bpr 1. Первое равенство доказывается простым вычислением:
$$
\exp 0=\sum_{n=0}^\infty
\frac{0^n}{n!}=\underbrace{\frac{0^0}{0!}}_{1}+\underbrace{\frac{0^1}{1!}}_{0}+\underbrace{\frac{0^2}{2!}}_{0}+...=1
$$

2. Докажем затем последнее равенство. Из леммы \ref{LM:exp} следует, что
функция $\exp$ является целой. Поэтому по следствию \ref{cor-21.3.2}, она
является суммой своего ряда Тейлора в произвольной точке $a\in\R$:
 \beq\label{exp(a+y)}
\exp(a+y)=\sum_{n=0}^\infty\frac{\exp^{(n)}a}{n!}\cdot y^n
 \eeq
По лемме \ref{LM:exp'}, производная функции $\exp$ равна самой этой функции,
поэтому вторая, третья, и все остальные производные $\exp$ тоже совпадают с
$\exp$:
$$
\exp^{(n)}=\exp
$$
Подставив это в \eqref{exp(a+y)}, мы получим:
 \begin{multline*}
\exp(a+y)=\sum_{n=0}^\infty\frac{\exp^{(n)}a}{n!}\cdot y^n=
\sum_{n=0}^\infty\frac{\exp a}{n!}\cdot y^n=\\=\exp
a\cdot\sum_{n=0}^\infty\frac{1}{n!}\cdot y^n=\exp a\cdot\exp y
 \end{multline*}
От \eqref{exp(x+y)} это отличается только заменой $a$ на $x$.

3. Второе тождество оказывается следствием первого и третьего:
  $$
1=\exp 0=\exp(x-x)=\exp x\cdot\exp(-x)
 $$
 $$
\Downarrow
 $$
 $$
\exp(-x)=\frac{1}{\exp x}
 $$
 \epr

 \blm\label{LM:exp(x)>0} Функция $\exp$, определенная рядом \eqref{DF:exp}, всюду положительна, а в нуле равна единице:
 \beq\label{exp(x)>0}
 \exp x>0
 \eeq
 \elm
\bpr Предположим, что в какой-то точке $a$ функция $\exp$ равняется нулю:
$$
\exp a=0.
$$
Тогда:
  \begin{multline*}
1=\exp 0=\exp(a-a)=\eqref{exp(x+y)}=\\=\underbrace{\exp a}_{0}\cdot\exp(-a)=0,
  \end{multline*}
что, конечно, невозможно. Таким образом, $\exp$ -- непрерывная функция,
определенная всюду на прямой $\R$, нигде не обращающаяся в нуль, а в точке 0
равная единице. Такое возможно только если $\exp$ всюду положительна.
  \epr

 \blm\label{LM:exp(x)-vozrastaet}
Функция $\exp$, определенная рядом \eqref{DF:exp}, возрастает на $\R$:
 \beq\label{exp(x)-vozrastaet}
x<y\qquad\Longrightarrow\qquad \exp x<\exp y
 \eeq
 \elm
\bpr Производная этой функции всюду положительна
$$
\frac{\d}{\d x}\exp x=\eqref{exp'x}=\exp x>\eqref{exp(x)>0}>0,
$$
поэтому по теореме \ref{Roll-cons} о строгой монотонности, эта функция должна
возрастать. \epr

 \blm\label{LM:predely-exp(x)}
Функция $\exp$, определенная рядом \eqref{DF:exp}, имеет следующие пределы на
бесконечности:
 \begin{align}\label{exp(x)-vozrastaet}
&\lim_{x\to-\infty}\exp x=0 &&\lim_{x\to+\infty}\exp x=+\infty
 \end{align}
  \elm
\bpr Второй предел следует из очевидного неравенства:
$$
\forall x>0\qquad \exp x=1+x+\frac{x^2}{2}+...>1+x
$$
После того, как он доказан, первый предел получается заменой переменных:
  \begin{multline*}
\lim_{x\to-\infty}\exp x=\begin{pmatrix}y=-x \\ y\to+\infty
\end{pmatrix}=\lim_{y\to+\infty}\exp(-y)=\\=\eqref{exp(-x)}=\lim_{y\to+\infty}\frac{1}{\exp
y}=\left(\frac{1}{\lim\limits_{y\to+\infty}\exp
y}\right)=\\=\left(\frac{1}{+\infty}\right)=0
  \end{multline*}
 \epr

\paragraph{Функция $\ln$.}

 \blm\label{LM:exp-biektsiya}
Функция $\exp$, определенная рядом \eqref{DF:exp}, биективно отображает прямую
$\R$ на интервал $(0;+\infty)$
 \beq
\exp:\R\to(0;+\infty)
 \eeq
и поэтому правило
 \beq\label{DF:ln}
t=\ln x\qquad\Longleftrightarrow\qquad \exp t=x
 \eeq
корректно определяет функцию $\ln$, обратную к $\exp$:
 \beq
\ln:(0;+\infty)\to\R
 \eeq
  \elm
 \bpr
1. Прежде всего, неравенство \eqref{exp(x)>0} означает, что отображение $\exp$
действительно переводит $\R$ в $(0;+\infty)$.

2. Далее из условия возрастания \eqref{exp(x)-vozrastaet} следует, что
отображение $\exp$ инъективно.

3. Покажем, что $\exp$ сюръективно отображает $\R$ на $(0;+\infty)$. Пусть
$C\in(0;+\infty)$. Поскольку $\lim_{x\to-\infty}\exp x=0$ (первое равенство в
\eqref{exp(x)-vozrastaet}), найдется $x\in\R$ такое, что $\exp x<C$. С другой
стороны, $\lim_{x\to+\infty}\exp x=+\infty$ (второе равенство в
\eqref{exp(x)-vozrastaet}), поэтому найдется $y\in\R$ такое, что $\exp y>C$. Мы
получаем, что
$$
\exp x<C<\exp y
$$
причем $\exp$ -- непрерывная функция по теореме \ref{tm-21.2.2}. Значит, по
теореме Коши о среднем значении \ref{Cauchy-I}, найдется точка $c\in\R$ такая,
что
$$
\exp c=C
$$
Мы получили, что какое ни возьми $C\in(0;+\infty)$, для него найдется точка
$c\in\R$, в которой функция $\exp$ принимает значение $C$. Это и означает
сюръективность $\exp$ как отображения из $\R$ в $(0;+\infty)$.

4. Инъективность и сюръективность вместе означают биективность $\exp:\R\to
(0;+\infty)$.
 \epr

 \blm\label{LM:ln-monot} Функция $\ln$, определенная правилом \eqref{DF:ln},
возрастает:
 \beq\label{ln-monot}
0<x<y\qquad\Longrightarrow\qquad \ln x<\ln y
 \eeq
 \elm
 \bpr
Если бы оказалось, что $\ln x\ge \ln y$, то, в силу возрастания $\exp$, мы
получили бы
  $$
\ln x\ge \ln y
 $$
$$
\Downarrow
$$
$$
\underbrace{\exp(\ln x)}_{x}\ge\underbrace{\exp(\ln y)}_{y}
$$
$$
\Downarrow
$$
$$
x\ge y
$$
 \epr

 \blm\label{LM:ln(xy)} Функция $\ln$, определенная правилом \eqref{DF:ln},
 удовлетворяет следующему равенству и двум тождествам:
 \begin{align}
& \ln 1=0, \label{ln(1)} \\
& \ln\frac{1}{x}=-\ln x \label{ln(1/x)}\\
& \ln(x\cdot y)=\ln x+ \ln y, \label{ln(xy)}
 \end{align}
 \elm
\bpr Первое равенство следует непосредственно из определения \eqref{DF:ln}:
$$
\exp 0=\eqref{exp(0)}=1 \qquad\Longleftrightarrow\qquad 0=\ln 1
$$

Чтобы доказать второе, достаточно вычислить значения функции $\exp$ в обеих его
частях:
 \begin{multline*}
\exp\left(\ln\frac{1}{x}\right)=\frac{1}{x}=\frac{1}{\exp(\ln
x)}=\\=\eqref{exp(-x)}=\exp(-\ln x)
 \end{multline*}
Поскольку $\exp$ -- инъективное отображение, равенство
$\exp\left(\ln\frac{1}{x}\right)=\exp(-\ln x)$ означает, что аргументы у $\exp$
должны быть равны: $\ln\frac{1}{x}=-\ln x$.

Точно так же доказывается последнее равенство:
  \begin{multline*}
\exp\left(\ln(x\cdot y)\right)=x\cdot y=\exp(\ln x)\cdot\exp(\ln
y)=\\=\eqref{exp(x+y)}=\exp(\ln x+\ln y)
  \end{multline*}
Опять, поскольку $\exp$ -- инъективное отображение, равенство
$\exp\left(\ln(x\cdot y)\right)=\exp(\ln x+\ln y)$ возможно только если
аргументы у $\exp$ равны: $\ln(x\cdot y)=\ln x+ \ln y$. \epr

\paragraph{Определение степеней $a^b$ и доказательство Аксиомы степеней.}

Теперь мы можем определить степенное отображение $(a,b)\mapsto a^b$. В его
определении участвует отображение $\sgn_2:\frac{\Z}{2\N-1}\to\{-1;+1\}$,
которое мы когда-то давно задавали формулой \eqref{sgn_2}:
 \beq\label{DF:a^b}
a^b:=\begin{cases}\exp(b\cdot\ln a),& a>0 \\ 0,& a=0,\ b>0 \\ 1,& a=0,\ b=0
\\ \sgn_2b\cdot\exp(b\cdot\ln |a|),& a<0,\ b\in\frac{\Z}{2\N-1} \end{cases}
 \eeq

Напомним, что еще в главе \ref{ch-R&N} (формула \eqref{DEF:a-vee-b}) мы
условились символом $a\vee b$ обозначать операцию, совпадающую с операцией
взятия максимума двух чисел:
$$
a\vee b=\max\{a,b\}
$$
Такая запись позволяет сделать более наглядными (или более ``алгебраическими'')
выкладки, которыми мы займемся в этом пункте.

\blm При $a\ne 0$ и $b\in\frac{\Z}{2\N-1}$ справедливо тождество:
 \beq\label{a^b-pri-a-ne-0-b-in-Z/2N-1}
a^b=\Big(\sgn_2 b\vee \sgn a \Big)\cdot\exp(b\cdot\ln |a|)
 \eeq
\elm
 \bpr
Если $a>0$, то мы получаем:
  \begin{multline*}
a^b=\eqref{DF:a^b}=\exp(b\cdot\ln a)=\\=\underbrace{\Big(\sgn_2 b\vee
\overbrace{\sgn a}^{1} \Big)}_{1}\cdot\exp(b\cdot\ln \underbrace{|a|}_{a})
  \end{multline*}
Если же $a<0$, то
 \begin{multline*}
a^b=\eqref{DF:a^b}=\sgn_2b\cdot\exp(b\cdot\ln |a|)=\\=\underbrace{\Big(\sgn_2
b\vee\overbrace{\sgn a}^{-1} \Big)}_{\sgn_2b}\cdot\exp(b\cdot\ln |a|)
 \end{multline*}
 \epr

\blm\label{LM:pokaz-zakony} Отображение $(a,b)\mapsto a^b$, определенное
формулой \eqref{DF:a^b}, удовлетворяет показательным законам
\eqref{a^0}-\eqref{a^(x+y)}:
 \begin{align}
&a^0=1, \label{a^(0)-1} \\
& a^{-x}=\frac{1}{a^x}, \label{a^(-x)-1} \\
&a^{x+y}=a^x\cdot a^y \label{a^(x+y)-1}
 \end{align}
(каждое из этих равенств верно всякий раз, когда обе его части определены).
\elm
 \bpr
Здесь нужно рассмотреть несколько случаев.

1. Проверим сначала эти равенства в случае, когда $a>0$. Тогда, во-первых,
$$
a^0=\eqref{DF:a^b}=\exp(0\cdot\ln a)=\exp 0=\eqref{exp(0)}=1
$$
Во-вторых,
 \begin{multline*}
a^{-x}=\eqref{DF:a^b}=\exp(-x\cdot\ln
a)=\\=\eqref{exp(-x)}=\frac{1}{\exp(x\cdot\ln a)}=\eqref{DF:a^b}=\frac{1}{a^x}
 \end{multline*}
И, в-третьих,
 \begin{multline*}
a^{x+y}=\eqref{DF:a^b}=\exp((x+y)\cdot\ln a)=\\=\exp(x\cdot\ln a+y\cdot\ln
a)=\eqref{exp(x+y)}=\\=\exp(x\cdot\ln a)\cdot\exp(y\cdot\ln
a)=\eqref{DF:a^b}=a^x\cdot a^y
 \end{multline*}

2. После этого перейдем к случаю $a=0$. Тогда, во-первых,
$$
a^0=0^0=\eqref{DF:a^b}=1
$$
Во-вторых, в соответствии с определением \eqref{DF:a^b}, выражение $0^{-x}$
имеет смысл только при $x\le 0$, а  и $0^x$ -- только при $x\ge 0$. Значит,
одновременно то и другое имеет смысл только при $x=0$. Мы получаем, что
равенство $0^{-x}=\frac{1}{0^x}$ нужно проверять лишь для случая $x=0$, но
тогда оно становится следствием уже отмеченного равенства $0^0=1$:
$$
0^{-0}=0^0=1=\frac{1}{1}=\frac{1}{0^0}
$$
В-третьих, выражения $0^{x+y}$, $0^x$, $0^y$ имеют смысл при $x\ge 0$ и $y\ge
0$. Поэтому для доказательства равенства $0^{x+y}=0^x\cdot 0^y$ придется
рассмотреть дополнительно четыре случая:
 \bit{
\item[---] если $x=y=0$, то
 \begin{multline*}
0^{x+y}=0^{0+0}=0^0=1=\\=1\cdot 1=0^0\cdot 0^0=0^x\cdot 0^y
 \end{multline*}

\item[---] если $x=0$, $y>0$, то
 \begin{multline*}
0^{x+y}=0^{0+y}=0^y=\eqref{DF:a^b}=0=\\=1\cdot 0=\eqref{DF:a^b}=0^0\cdot
0^y=0^x\cdot 0^y
 \end{multline*}

\item[---] если $x>0$, $y=0$, то
 \begin{multline*}
0^{x+y}=0^{x+0}=0^x=\eqref{DF:a^b}=0=\\=0\cdot 1=\eqref{DF:a^b}=0^x\cdot
0^0=0^x\cdot 0^y
 \end{multline*}

\item[---] если $x>0$, $y>0$, то
 \begin{multline*}
 0^{x+y}=\eqref{DF:a^b}=0=\\=0\cdot 0=\eqref{DF:a^b}=0^x\cdot 0^y
 \end{multline*}
 }\eit

3. Остается случай $a<0$. Во-первых,
 \begin{multline*}
a^0=\eqref{DF:a^b}=\underbrace{\sgn_2 0}_{1}\cdot
\exp(0\cdot\ln|a|)=\\=\exp(0)=\eqref{exp(0)}=1
 \end{multline*}
Во-вторых, в соответствии с определением \eqref{DF:a^b}, выражения $a^{-x}$ и
$a^x$ имеют смысл для $x\in\frac{\Z}{2\N-1}$. В этом случае получается:
 \begin{multline*}
a^{-x}=\eqref{DF:a^b}=\sgn_2(-x)\cdot
\exp(-x\cdot\ln|a|)=\\=\eqref{sgn_2(-x)=sgn_2(x)},\eqref{exp(-x)}=\sgn_2 x\cdot
\frac{1}{\exp(x\cdot\ln|a|)}=\\=\eqref{sgn_2=1/sgn_2}= \frac{1}{\sgn_2
x\cdot\exp(x\cdot\ln|a|)}=\eqref{DF:a^b}=\frac{1}{a^x}
 \end{multline*}
В-третьих, выражения $a^{x+y}$, $a^x$, $a^y$ имеют смысл при
$x,y\in\frac{\Z}{2\N-1}$, и мы получаем:
 \begin{multline*}
a^{x+y}=\eqref{DF:a^b}=\\=\sgn_2(x+y)\cdot
\exp((x+y)\cdot\ln|a|)=\\=\sgn_2(x+y)\cdot
\exp(x\cdot\ln|a|+y\cdot\ln|a|)=\\=\eqref{sgn_2(x+y)=sgn_2(x)-sgn_2(y)},\eqref{exp(x+y)}=\\=\sgn_2
x\cdot\sgn_2
y\cdot\exp(x\cdot\ln|a|)\cdot\exp(y\cdot\ln|a|)=\\=\underbrace{\sgn_2
x\cdot\exp(x\cdot\ln|a|)}_{a^x}\cdot\underbrace{\sgn_2
y\cdot\exp(y\cdot\ln|a|)}_{a^y}=\\=\eqref{DF:a^b}=a^x\cdot a^y
 \end{multline*}
\epr

\blm\label{LM:step-zakony} Отображение $(a,b)\mapsto a^b$, определенное
формулой \eqref{DF:a^b}, удовлетворяет степенным законам
\eqref{1^b}-\eqref{(xy)^b}:
 \begin{align}
&1^b=1, && \left(\frac{1}{x}\right)^b=\frac{1}{x^b}, &&(x\cdot y)^b=x^b\cdot
y^b  \label{(xy)^b-1}
 \end{align}
(каждое из этих равенств верно всякий раз, когда обе его части определены).
\elm
 \bpr
 1. Первое равенство имеет смысл при любом $b\in\R$, и мы получаем:
 \begin{multline*}
1^b=\eqref{DF:a^b}=\exp(b\cdot\ln 1)=\eqref{ln(1)}=\\=\exp(b\cdot 0)=\exp
0=\eqref{exp(0)}=1
 \end{multline*}

2. Для доказательства второго тождества нужно рассмотреть три случая:
 \biter{
\item[---] если $x>0$, то $b\in\R$, и
 \begin{multline*}
\left(\frac{1}{x}\right)^b=\eqref{DF:a^b}=\exp\left(b\cdot\ln\frac{1}{x}\right)=\\=\eqref{ln(1/x)}=
\exp(-b\cdot\ln x)=\eqref{exp(-x)}=\\=\frac{1}{\exp(b\cdot\ln
x)}=\eqref{DF:a^b}=\frac{1}{x^b}
 \end{multline*}

\item[---] если $x=0$, то выражение $\frac{1}{x}$ не имеет смысл, поэтому
проверять нечего;

\item[---] если $x<0$, то $b\in\frac{\Z}{2\N-1}$, и
 \begin{multline*}
\left(\frac{1}{x}\right)^b=\eqref{DF:a^b}=\\=\sgn_2 b\cdot
\exp\left(b\cdot\ln\left|\frac{1}{x}\right|\right)=\\=\sgn_2 b\cdot
\exp\left(b\cdot\ln\frac{1}{|x|}\right)=\eqref{ln(1/x)}=\\=\sgn_2 b\cdot
\exp(-b\cdot\ln|x|)=\eqref{exp(-x)}=\\=\sgn_2 b\cdot\frac{1}{\exp(b\cdot\ln
|x|)}=\eqref{sgn_2=1/sgn_2}=\\=\frac{1}{\sgn_2 b\cdot\exp(b\cdot\ln
|x|)}=\eqref{DF:a^b}=\frac{1}{x^b}
 \end{multline*}
 }\eiter

3. Для доказательства третьего тождества приходится рассматривать 6 случаев.
Во-первых, нужно отдельно поглядеть, что получается, когда $x$ или $y$ равно
нулю:
 \biter{
\item[---] если $x=0$, то чтобы $x^b$ имело смысл, нужно чтобы $b\ge 0$, и при
этом становится неважно, каким будет $y$:
 \biter{
\item[1)] при $b=0$ получаем:
 \begin{multline*}
(x\cdot y)^b=(0\cdot y)^0=0^0=\eqref{DF:a^b}=1=\\=1\cdot
1=\eqref{DF:a^b}=0^0\cdot y^0=x^b\cdot y^b
 \end{multline*}

\item[2)] при $b>0$ получаем:
 \begin{multline*}
(x\cdot y)^b=(0\cdot y)^b=0^b=\eqref{DF:a^b}=0=\\=0\cdot
y^b=\eqref{DF:a^b}=0^b\cdot y^b=x^b\cdot y^b
 \end{multline*}
 }\eiter

\item[---] по тем же причинам при $y=0$ нужно чтобы $b\ge 0$, при этом
становится неважно каким будет $x$, и тождество также выполняется.

 }\eiter

После этого остается рассмотреть 4 случая, когда $x$ и $y$ принимают разные
знаки:

 \biter{
\item[---] если $x>0$, $y>0$, то $b\in\R$, и
 \begin{multline*}
(x\cdot y)^b =\eqref{DF:a^b}=\exp(b\cdot\ln(x\cdot y))=\\=\eqref{ln(xy)}=
\exp(b\cdot(\ln x+\ln y))=\\=\exp(b\cdot\ln x+b\cdot\ln
y)=\eqref{exp(x+y)}=\\=\exp(b\cdot\ln x)\cdot\exp(b\cdot\ln
y)=\eqref{DF:a^b}=x^b\cdot y^b
 \end{multline*}

\item[---] если $x<0$, $y>0$, то $b\in\frac{\Z}{2\N-1}$, и
 \begin{multline*}
(x\cdot y)^b =\eqref{DF:a^b}=\sgn_2 b\cdot\exp(b\cdot\ln|x\cdot y|)=\\= \sgn_2
b\cdot\exp(b\cdot\ln(|x|\cdot y))=\eqref{ln(xy)}=\\=\sgn_2
b\cdot\exp(b\cdot\ln|x|+b\cdot \ln y)=\eqref{exp(x+y)}=\\= \underbrace{\sgn_2
b\cdot\exp(b\cdot\ln|x|)}_{x^b}\cdot\underbrace{\exp(b\cdot \ln
y)}_{y^b}=x^b\cdot y^b
 \end{multline*}

\item[---] то же самое, с точностью до очевидных перестановок, получается при
$x>0$, $y<0$ (тогда опять $b\in\frac{\Z}{2\N-1}$);

\item[---] наконец, если $x<0$, $y<0$, то $b\in\frac{\Z}{2\N-1}$, и
 \begin{multline*}
(\overbrace{x\cdot y}^{\scriptsize\begin{matrix}0\\
\text{\rotatebox{90}{$>$}}\end{matrix}})^b
=\eqref{DF:a^b}=\exp(b\cdot\ln(x\cdot y))=\\= \overbrace{(\sgn_2 b)^2}^{1}\cdot
\exp(b\cdot\ln(\overbrace{|x|\cdot|y|}^{x\cdot y}))=\\=\eqref{ln(xy)}=(\sgn_2
b)^2\cdot\exp(b\cdot\ln|x|+b\cdot \ln |y|)=\\=\eqref{exp(x+y)}=\\=
\underbrace{\sgn_2 b\cdot\exp(b\cdot\ln|x|)}_{x^b}\cdot\underbrace{\sgn_2
b\cdot\exp(b\cdot \ln |y|)}_{y^b}=\\=x^b\cdot y^b
 \end{multline*}
 }\eiter
 \epr

\blm\label{LM:nakop-zakon} Отображение $(a,b)\mapsto a^b$, определенное
формулой \eqref{DF:a^b}, удовлетворяет накопительному закону
\eqref{nakop-zakon}:
 \begin{align}\label{nakop-zakon-1}
&(a^x)^y=a^{x\cdot y}&& (a\ge 0)
 \end{align}
в следующих случаях:
 \biter{
\item[---] при $a>0$ и $x,y\in\R$,

\item[---] при $a=0$ и $x\ge 0$, $y\ge 0$,

\item[---] при $a<0$ и $x,y\in\frac{\Z}{2\N-1}$.
 }\eiter
 \elm
 \bpr 1. Пусть $a>0$ и $x,y\in\R$. Тогда:
  \begin{multline*}
(a^x)^y=\eqref{DF:a^b}=\exp(y\cdot\ln
a^x)=\eqref{DF:a^b}=\\=\exp\Big(y\cdot\underbrace{\ln\big(\exp(x\cdot\ln
a)\big)}_{x\cdot\ln a}\Big)=\exp(y\cdot x\cdot\ln
a)=\\=\eqref{DF:a^b}=a^{x\cdot y}
 \end{multline*}

2. Пусть $a=0$ и $x\ge 0$, $y\ge 0$. Тогда:
 \biter{
\item[---] если $x=0$, то
 \begin{multline*}
(a^x)^y=(\underbrace{0^0}_{1})^y=1^y=\eqref{(xy)^b-1}=\\=1=\eqref{DF:a^b}=0^0=0^{0\cdot
y}=a^{x\cdot y}
 \end{multline*}

\item[---] если $x>0$ и $y=0$, то
$$
(a^x)^y=(\underbrace{0^x}_{0})^0=\underbrace{0^0}_{1}=0^{x\cdot 0}=a^{x\cdot y}
$$

\item[---] если $x>0$ и $y>0$, то
$$
(a^x)^y=(\underbrace{0^x}_{0})^y=\underbrace{0^y}_{0}=0^{x\cdot y}=a^{x\cdot y}
$$
 }\eiter

3. Пусть $a<0$ и $x,y\in\frac{\Z}{2\N-1}$. Тогда нужно воспользоваться
формулами \eqref{a^b-pri-a-ne-0-b-in-Z/2N-1} и
\eqref{sgn_2(xy)=max(sgn_2(x),sgn_2(y))}:
 \begin{multline*}
({\underbrace{a}_{\scriptsize\begin{matrix}\text{\rotatebox{-90}{$<$}}\\
0 \end{matrix}}}^x)^y=\eqref{DF:a^b}=\\=\Big(
{\underbrace{\sgn_2 x\cdot\exp(x\cdot\ln |a|)}_{\scriptsize\begin{matrix}\text{\rotatebox{-90}{$\ne$}}\\
0 \end{matrix}}}\Big)^y=\eqref{a^b-pri-a-ne-0-b-in-Z/2N-1}=\\=\Big(\sgn_2 y\
\vee \ \underbrace{\sgn\Big(\sgn_2 x\cdot\exp(x\cdot\ln |a|)\Big)}_{\sgn_2 x}
\Big)\cdot\\ \cdot \exp\Big(y\cdot\ln \underbrace{\Big|\sgn_2
x\cdot\exp(x\cdot\ln |a|)\Big|}_{\exp(x\cdot\ln |a|)}\Big)=\\=
\underbrace{\big(\sgn_2 y\ \vee\ \sgn_2 x
\big)}_{\scriptsize\begin{matrix}\text{\rotatebox{-90}{$=$}}\\
\eqref{sgn_2(xy)=max(sgn_2(x),sgn_2(y))}\\ \text{\rotatebox{-90}{$=$}}
\\ \sgn_2( y\cdot x)
\end{matrix}}\cdot\\ \cdot\exp\bigg(y\cdot
\underbrace{\ln\Big(\exp(x\cdot\ln |a|)\Big)}_{\scriptsize
\begin{matrix}\text{\rotatebox{-90}{$=$}}\\ x\cdot\ln |a|
\end{matrix}}\bigg)=\\=
\sgn_2( y\cdot x)\cdot\exp\big(y\cdot x\cdot\ln |a|\big)=a^{y\cdot x}=a^{x\cdot
y}
 \end{multline*}
 \epr

\blm\label{LM:x^b>0} Отображение $(a,b)\mapsto a^b$, определенное формулой
\eqref{DF:a^b}, удовлетворяет условию сохранения знака \eqref{x^n>0}:
 \beq\label{x^n>0-1}
a>0\quad\Rightarrow\quad a^n>0\quad (n\in\Z)
 \eeq
 \elm
\bpr
$$
a>0\quad\Rightarrow\quad a^n\overset{\eqref{DF:a^b}}{=}\exp(b\cdot\ln
a)\overset{\eqref{exp(x)>0}}{>}0
$$
\epr

\blm\label{LM:monot-x^b} Отображение $(a,b)\mapsto a^b$, определенное формулой
\eqref{DF:a^b}, удовлетворяет условиям монотонности
\eqref{monot-x^b-b>0}-\eqref{monot-x^b-b<0}:
 \beq\label{monot-x^b-b>0-1}
\Big(b>0\quad\&\quad 0<x<y\Big) \quad\Longrightarrow\quad 0<x^b<y^b
 \eeq
 \beq\label{monot-x^b-b<0-1}
\Big(b<0\quad\&\quad 0<x<y\Big) \quad\Longrightarrow\quad x^b>y^b>0
 \eeq
 \elm
 \bpr
При $b>0$ получаем:
$$
0<x<y
$$
$$
\Downarrow
$$
$$
\ln x<\ln y
$$
$$
\Downarrow
$$
$$
b\cdot\ln x<b\cdot\ln y
$$
$$
\Downarrow
$$
$$
\underbrace{\exp(b\cdot\ln x)}_{x^b}<\underbrace{\exp(b\cdot\ln y)}_{y^b}
$$
А при $b<0$:
$$
0<x<y
$$
$$
\Downarrow
$$
$$
\ln x<\ln y
$$
$$
\Downarrow
$$
$$
b\cdot\ln x>b\cdot\ln y
$$
$$
\Downarrow
$$
$$
\underbrace{\exp(b\cdot\ln x)}_{x^b}>\underbrace{\exp(b\cdot\ln y)}_{y^b}
$$
 \epr

\blm\label{LM:monot-a^x} Отображение $(a,b)\mapsto a^b$, определенное формулой
\eqref{DF:a^b}, удовлетворяет условиям монотонности
\eqref{monot-a^x-a>1}-\eqref{monot-a^x-0<a<1}:
 \begin{align}
&\Big(a>1\quad\&\quad x<y\Big) \quad\Longrightarrow\quad a^x<a^y
\label{monot-a^x-a>1-1}
 \\
&\Big(0<a<1\quad\&\quad x<y\Big) \quad\Longrightarrow\quad a^x>a^y
\label{monot-a^x-0<a<1-1}
 \end{align} \elm
 \bpr
Если $a>1$, то $\ln a>\ln 1=0$ (поскольку функция $\ln$ возрастает), поэтому
при умножении на $\ln a$ знак неравенства не меняется:
$$
x<y
$$
$$
\Downarrow
$$
$$
x\cdot \ln a<y\cdot \ln a
$$
$$
\Downarrow
$$
$$
\underbrace{\exp(x\cdot\ln a)}_{a^x}<\underbrace{\exp(y\cdot\ln a)}_{a^y}
$$

Если же $0<a<1$, то $\ln a<\ln 1=0$ (поскольку функция $\ln$ возрастает),
поэтому при умножении на $\ln a$ знак неравенства меняется:
$$
x<y
$$
$$
\Downarrow
$$
$$
x\cdot \ln a>y\cdot \ln a
$$
$$
\Downarrow
$$
$$
\underbrace{\exp(x\cdot\ln a)}_{a^x}>\underbrace{\exp(y\cdot\ln a)}_{a^y}
$$
 \epr

\bpr[Доказательство избыточности Аксиомы B1]\label{PROOF:TH-o-step-otobr}
Формула \eqref{DF:a^b} определяет отображение $(a,b)\mapsto a^b$ в точности для
тех ситуаций, что перечислены в условии $P_0$ Аксиомы степеней. По леммам
\ref{LM:pokaz-zakony} и \ref{LM:step-zakony}, это отображение удовлетворяет
показательным и степенным законам, то есть условию $P_1$. По лемме
\ref{LM:nakop-zakon}, оно удовлетворяет накопительному закону, то есть условию
$P_2$. По лемме \ref{LM:x^b>0}, закон сохранения знака \eqref{x^b>0} тоже
выполняется, а в силу замечания \ref{REM:0^b}, формула \eqref{0^b} выполняется
автоматически, поэтому ее доказывать не нужно. Наконец, по леммам
\ref{LM:monot-x^b} и \ref{LM:monot-a^x} оно удовлетворяет условиям монотонности
$P_3$. \epr

\end{multicols}\noindent\rule[10pt]{160mm}{0.1pt}

\subsection{Доказательство избыточности Аксиомы тригонометрии.}\label{PR:osn-teor-trig}

Наступил момент вспомнить еще об одном важном результате, доказательство
которого (как в случае с Аксиомой степеней B1 на с.\pageref{TH-o-step-otobr})
было отложено нами на будущее. Это Аксиома тригонометрии B2 на
с.\pageref{osn-teor-trigonometrii}, с помощью которой мы вводили функции синус
и косинус в главе \ref{ch-ELEM-FUNCTIONS}. Теперь мы можем доказать и ее. Мы
предпошлем доказательству 11 вспомогательных утверждений (леммы
\ref{LM:sin-cos}-\ref{LM:0<sin-x<x-2}). Здесь (как в случае с Аксиомой степеней
B1 на с.\pageref{TH-o-step-otobr}) мы опять предлагаем читателю сделать вид,
что до настоящего времени он не глядел на текст в двух колонках и не знает
ничего из того, что писалось о синусе и косинусе вплоть до момента, когда мы
собственно примемся за доказательство Аксиомы тригонометрии (на
с.\pageref{PROOF:osn-teor-trigonometrii}). Только в самом доказательстве мы
вспомним, что следствия из Аксиомы тригонометрии уже рассматривались ранее нами
в иллюстрациях (и из этого будет выводиться единственность функций $\sin$ и
$\cos$).

\noindent\rule{160mm}{0.1pt}\begin{multicols}{2}

\blm\label{LM:sin-cos} Формулы
 \begin{align}
&\sin x=\sum_{n=0}^\infty (-1)^n\cdot \frac{x^{2n+1}}{(2n+1)!} \label{DF:sin}\\
& \cos x=\sum_{n=0}^\infty (-1)^n\cdot \frac{x^{2n}}{(2n)!} \label{DF:cos}
\end{align}
определяют функции $\sin$ и $\cos$ всюду на прямой $\R$.
 \elm
 \bpr
 1. Первый ряд перепишем в виде
 \begin{multline*}
\sum_{n=0}^\infty (-1)^n\cdot \frac{x^{2n+1}}{(2n+1)!}=x\cdot\sum_{n=0}^\infty
(-1)^n\cdot \frac{(x^2)^n}{(2n+1)!}=\\=x\cdot\sum_{n=0}^\infty (-1)^n\cdot
\frac{y^n}{(2n+1)!}\Bigg|_{y=x^2}
 \end{multline*}
и вычислим радиус сходимости у последнего степенного ряда (от переменной $y$)
по формуле Даламбера:
 \begin{multline*}
R=\lim_{n\to\infty}\left|\frac{c_n}{c_{n+1}}\right|=
\lim_{n\to\infty}\left|\frac{(-1)^n\cdot(2n+3)!}{(-1)^{n+1}\cdot(2n+1)!}\right|=\\=
\lim_{n\to\infty}(2n+2)(2n+3)=\infty
 \end{multline*}
То есть ряд
$$
\sum_{n=0}^\infty (-1)^n\cdot \frac{y^n}{(2n+1)!}
$$
сходится при любом $y$. Значит при подстановке $y=x^2$ получающийся ряд
$$
\sum_{n=0}^\infty (-1)^n\cdot \frac{(x^2)^n}{(2n+1)!}
$$
тоже сходится независимо от того, какое значение принимает $x$. После умножения
на $x$ получающийся ряд
$$
\sum_{n=0}^\infty (-1)^n\cdot \frac{x^{2n+1}}{(2n+1)!}
$$
тоже должен сходиться в силу свойства $1^\circ$ на с.\pageref{18.2.1}, каким бы
ни был $x$. Все это означает, что формула \eqref{DF:sin} определяет функцию
$\sin$ всюду на прямой $\R$.

2. Второй ряд можно переписать так:
 \begin{multline*}
\sum_{n=0}^\infty (-1)^n\cdot \frac{x^{2n}}{(2n)!}=\sum_{n=0}^\infty
(-1)^n\cdot \frac{(x^2)^n}{(2n)!}=\\=\sum_{n=0}^\infty (-1)^n\cdot
\frac{y^n}{(2n)!}\Bigg|_{y=x^2}
 \end{multline*}
Радиус сходимости последнего степенного ряда будет таким:
 \begin{multline*}
R=\lim_{n\to\infty}\left|\frac{c_n}{c_{n+1}}\right|=
\lim_{n\to\infty}\left|\frac{(-1)^n\cdot(2n+2)!}{(-1)^{n+1}\cdot(2n)!}\right|=\\=
\lim_{n\to\infty}(2n+1)(2n+2)=\infty
 \end{multline*}
Значит ряд
$$
\sum_{n=0}^\infty (-1)^n\cdot \frac{y^n}{(2n)!}
$$
сходится при любом $y$. Как следствие, при подстановке $y=x^2$ получающийся ряд
$$
\sum_{n=0}^\infty (-1)^n\cdot \frac{(x^2)^n}{(2n)!}
$$
тоже должен сходиться, независимо от того, какое значение принимает $x$. Мы
получаем, что формула \eqref{DF:cos} определяет функцию $\cos$ всюду на прямой
$\R$.
 \epr

 \blm Функции $\sin$ и $\cos$, определенные рядами \eqref{DF:sin}-\eqref{DF:cos},
 удовлетворяют тождествам:
 \begin{align}\label{sin(-x)-cos(-x)}
&\sin(-x)=-\sin x && \cos(-x)=\cos x
\end{align}
 \elm
 \bpr
 \begin{multline*}
\sin(-x)=\eqref{DF:sin}=\sum_{n=0}^\infty (-1)^n\cdot
\frac{(-x)^{2n+1}}{(2n+1)!}=\\=\sum_{n=0}^\infty (-1)^n\cdot
\frac{-(x^{2n+1})}{(2n+1)!}=\\=-\sum_{n=0}^\infty (-1)^n\cdot
\frac{x^{2n+1}}{(2n+1)!}=\eqref{DF:sin}=-\sin x
 \end{multline*}
 \begin{multline*}
\cos(-x)=\eqref{DF:cos}=\sum_{n=0}^\infty (-1)^n\cdot
\frac{(-x)^{2n}}{(2n)!}=\\=\sum_{n=0}^\infty (-1)^n\cdot
\frac{x^{2n}}{(2n)!}=\eqref{DF:cos}=\cos x
 \end{multline*}
 \epr

 \blm Функции $\sin$ и $\cos$, определенные рядами \eqref{DF:sin}-\eqref{DF:cos},
 удовлетворяют равенствам:
 \begin{align}\label{sin(0)-cos(0)}
&\sin 0=0 && \cos 0=1
\end{align}
 \elm
 \bpr
 \begin{multline*}
\sin 0=\eqref{DF:sin}=\sum_{n=0}^\infty (-1)^n\cdot
\frac{0^{2n+1}}{(2n+1)!}=\\=\sum_{n=0}^\infty (-1)^n\cdot 0=0
 \end{multline*}
 \begin{multline*}
\cos 0=\eqref{DF:cos}=\sum_{n=0}^\infty (-1)^n\cdot
\frac{0^{2n}}{(2n)!}=\\=\sum_{n=0}^\infty (-1)^n\cdot
\frac{0^{2n}}{(2n)!}=\underbrace{(-1)^0\cdot \frac{0^{2\cdot 0}}{(2\cdot
0)!}}_{1}+\\+\sum_{n=1}^\infty \underbrace{(-1)^n\cdot
\frac{0^{2n}}{(2n)!}}_{0}=1
 \end{multline*}
 \epr

 \blm Функции $\sin$ и $\cos$, определенные рядами \eqref{DF:sin}-\eqref{DF:cos},
имеют следующие производные:
 \begin{align}\label{sin'x-cos'x}
&\frac{\d}{\d x}\sin x=\cos x && \frac{\d}{\d x}\cos x=-\sin x
\end{align}
 \elm
 \bpr
 \begin{multline*}
\frac{\d}{\d x}\sin x=\eqref{DF:sin}=\\= \frac{\d}{\d x}\sum_{n=0}^\infty
(-1)^n\cdot \frac{x^{2n+1}}{(2n+1)!}=\\=\sum_{n=0}^\infty
(-1)^n\cdot\frac{\d}{\d x}\frac{x^{2n+1}}{(2n+1)!}=\\=\sum_{n=0}^\infty
(-1)^n\cdot\frac{(2n+1)\cdot x^{2n}}{(2n+1)!}=\\=\sum_{n=0}^\infty
(-1)^n\cdot\frac{x^{2n}}{(2n)!}=\cos x
 \end{multline*}
 \begin{multline*}
\frac{\d}{\d x}\cos x=\eqref{DF:cos}=\\= \frac{\d}{\d x}\sum_{n=0}^\infty
(-1)^n\cdot \frac{x^{2n}}{(2n)!}=\\=\sum_{n=0}^\infty (-1)^n\cdot\frac{\d}{\d
x}\frac{x^{2n}}{(2n)!}=\\= (-1)^n\cdot 0+ \sum_{n=1}^\infty
(-1)^n\cdot\frac{(2n)\cdot x^{2n-1}}{(2n)!}=\\=\sum_{n=1}^\infty
(-1)^n\cdot\frac{x^{2n-1}}{(2n-1)!}={\scriptsize\begin{pmatrix}\text{замена:}\\
k=n-1
\\ n=k+1\end{pmatrix}}=\\=\sum_{k=1}^\infty
(-1)^{k+1}\cdot\frac{x^{2k+1}}{(2k+1)!}=\\= -\sum_{k=1}^\infty
(-1)^k\cdot\frac{x^{2k+1}}{(2k+1)!}=-\sin x
 \end{multline*}
 \epr

 \blm\label{LM:sin(x+y)-cos(x+y)} Функции $\sin$ и $\cos$, определенные рядами \eqref{DF:sin}-\eqref{DF:cos},
удовлетворяют тождествам:
 \begin{align}
&\sin(x+y)=\sin x\cdot\cos y+\cos x\cdot\sin y \label{sin(x+y)}\\
& \cos(x+y)=\cos x\cdot\cos y-\sin x\cdot\sin y \label{cos(x+y)}
 \end{align}
 \elm
 \bpr
Зафиксируем какое-нибудь число $a\in\R$ и рассмотрим три функции:
 \begin{align*}
& F(x)=\sin(x+a)-\sin x\cdot\cos a-\cos x\cdot\sin a \\
& G(x)=\cos(x+a)-\cos x\cdot\cos a+\sin x\cdot\sin a \\
& H(x)=F(x)^2+G(x)^2
 \end{align*}
Заметим две вещи: во-первых,
 $$
\begin{cases}F(0)=\sin a-\sin 0\cdot\cos a-\cos 0\cdot\sin a=0\\
G(0)=\cos a-\cos 0\cdot\cos a+\sin 0\cdot\sin a=0\end{cases}
$$
$$
\Downarrow
$$
$$
H(0)=F(0)^2+G(0)^2=0
$$
И, во-вторых,
$$
\begin{cases}F'(x)=\cos(x+a)-\cos x\cdot\cos a+\\
\qquad
+\sin x\cdot\sin a=G(x)\\
G'(x)=-\sin(x+a)+\sin x\cdot\cos a+\\ \qquad+\cos x\cdot\sin a=-F(x)\end{cases}
$$
$$
\Downarrow
$$
 \begin{multline*}
H'(x)=2\cdot F(x)\cdot F'(x)+2\cdot G(x)\cdot G'(x)=\\=2\cdot F(x)\cdot
G(x)-2\cdot G(x)\cdot F(x)=0
 \end{multline*}
Вместе это означает, что функция $H$ должна быть тождественно равна нулю, а
значит $F$ и $G$ тоже тождественно равны нулю:
 $$
H(x)=F(x)^2+G(x)^2\equiv 0
$$
$$
\Downarrow
$$
$$
\begin{cases} F(x)\equiv 0 \\  G(x)\equiv 0
\end{cases}
 $$
Заменяя $a$ на $y$, мы из последних тождеств получаем
\eqref{sin(x+y)}-\eqref{cos(x+y)}.
 \epr

 \blm\label{LM:sin^2-x+cos^2-x=1} Функции $\sin$ и $\cos$, определенные рядами \eqref{DF:sin}-\eqref{DF:cos},
удовлетворяют тождеству:
 \beq
\sin^2 x+\cos^2 x=1 \label{sin^2-x+cos^2-x=1}
 \eeq
 \elm
 \bpr
Подставив $y=-x$ в \eqref{cos(x+y)}, мы получим:
 \begin{multline*}
1=\eqref{sin(0)-cos(0)}=\cos 0=\cos(x-x)=\\=\cos x\cdot\cos(-x)-\sin
x\cdot\sin(-x)=\\=\eqref{sin(-x)-cos(-x)}=\cos^2x+\sin^2 x
 \end{multline*}
 \epr

\blm\label{LM:sin>0-na-(0,2)} Функция $\sin$, определенная рядом
\eqref{DF:sin}, строго положительна на интервале $(0,2)$:
$$
\forall x\in (0,2)\qquad \sin x>0
$$
 \elm
\bpr Перепишем ряд, определяющий $\sin$ таким образом:
 \begin{multline*}
\sin x=\sum_{n=0}^\infty (-1)^n\cdot
\frac{x^{2n+1}}{(2n+1)!}=\\=x-\frac{x^3}{3!}+\frac{x^5}{5!}-\frac{x^7}{7!}+...=\\=
x\cdot\left(1-\frac{x^2}{2\cdot
3}\right)-\frac{x^5}{5!}\cdot\left(1-\frac{x^2}{6\cdot
7}\right)+...=\\=\sum_{k=0}^\infty\frac{x^{4k+1}}{(4k+1)!}\cdot\left(1-\frac{x^2}{(4k+3)\cdot(4k+4)}\right)
 \end{multline*}
И заметим, что при $x\in[0,2]$ множитель в скобках в последнем ряде
положителен:
 $$
 x\in[0,2]
 $$
 $$
\Downarrow
 $$
 $$
 x^2\le 4< (4k+3)\cdot(4k+4) \qquad (k\in\Z_+)
 $$
 $$
\Downarrow
 $$
 $$
 \frac{x^2}{(4k+3)\cdot(4k+4)}<1 \qquad (k\in\Z_+)
 $$
 $$
 \Downarrow
 $$
 \beq\label{1-x^2/(4k+3)(4k+4)}
 1-\frac{x^2}{(4k+3)\cdot(4k+4)}>0\qquad (k\in\Z_+)
 \eeq
Отсюда следует, что при $x\in(0,2)$ в полученном ряде все слагаемые тоже
положительны, и значит
 $$
\sin
x=\sum_{k=0}^\infty\frac{x^{4k+1}}{(4k+1)!}\cdot\left(1-\frac{x^2}{(4k+3)\cdot(4k+4)}\right)>0
 $$
 \epr

 \blm Функция $\cos$, определенная рядом \eqref{DF:cos},
имеет на промежутке $(0,2)$ ровно один нуль:
 \beq\label{DEF:cos-a=0}
\exists !\ a\in (0,2)\qquad \cos a=0
 \eeq
 \elm
 \bpr
Перепишем ряд, определяющий $\cos$, иначе:
 \begin{multline*}
\cos x=\sum_{n=0}^\infty (-1)^n\cdot
\frac{x^{2n}}{(2n)!}=\\=1-\frac{x^2}{2!}+\frac{x^4}{4!}-\frac{x^6}{6!}+\frac{x^8}{8!}-...=\\=
1-\frac{x^2}{2!}\cdot\left(1-\frac{x^2}{3\cdot
4}\right)-\frac{x^6}{6!}\cdot\left(1-\frac{x^2}{7\cdot
8}\right)+...=\\=1-\sum_{k=0}^\infty\frac{x^{4k+2}}{(4k+2)!}\cdot\left(1-\frac{x^2}{(4k+3)\cdot(4k+4)}\right)
 \end{multline*}
В силу \eqref{1-x^2/(4k+3)(4k+4)}, при $x\in[0,2]$ множитель в скобках в
последнем ряде положителен, поэтому слагаемые в этом ряду неотрицательны.
Значит, $\cos$  на отрезке  $x\in[0,2]$ можно оценить сверху, оборвав последний
ряд после первого слагаемого:
 \begin{multline*}
\sum_{k=0}^\infty\frac{x^{4k+2}}{(4k+2)!}\cdot\left(1-\frac{x^2}{(4k+3)\cdot(4k+4)}\right)\ge\\
\ge
\sum_{k=0}^0\frac{x^{4k+2}}{(4k+2)!}\cdot\left(1-\frac{x^2}{(4k+3)\cdot(4k+4)}\right)=\\=
\frac{x^2}{2!}\cdot\left(1-\frac{x^2}{3\cdot
4}\right)
 \end{multline*}
$$
\Downarrow
$$
 \begin{multline*}
\cos
x=1-\\-\sum_{k=0}^\infty\frac{x^{4k+2}}{(4k+2)!}\cdot\left(1-\frac{x^2}{(4k+3)\cdot(4k+4)}\right)\le\\
\le 1-\frac{x^2}{2!}\cdot\left(1-\frac{x^2}{3\cdot 4}\right)
 \end{multline*}
В частности, при $x=2$ получаем:
 \begin{multline*}
\cos 2\le 1-\frac{2^2}{2!}\cdot\left(1-\frac{2^2}{3\cdot
4}\right)=\\=1-\frac{2}{2!}\cdot\left(1-\frac{1}{3}\right)=-\frac{1}{3}<0
 \end{multline*}
С другой стороны, в силу \eqref{sin(0)-cos(0)},
$$
\cos 0=1
$$
Мы получаем, что на отрезке $[0,2]$ функция $\cos$ меняет знак. Значит, по
теореме Коши о промежуточном значении \ref{Cauchy-I}, она обращается в нуль в
некоторой точке $a\in(0,2)$.

Остается показать, что это единственный нуль на $(0,2)$. Действительно, по
лемме \ref{LM:sin>0-na-(0,2)}, производная этой функции отрицательна на
интервале $(0,2)$:
$$
\frac{\d}{\d x}\cos x=-\sin x<0
$$
Значит, по теореме \ref{Roll-cons} $\cos$ монотонно убывает на интервале
$(0,2)$, и поэтому не может дважды обращаться в нуль на нем.
 \epr

\blm Точка $a\in(0,2)$, определенная условием \eqref{DEF:cos-a=0}, обладает
следующими свойствами:
 \begin{align}
& \cos a=0 && \sin a=1 \\
& \cos2a=-1 && \sin2a=0 \\
& \cos4a=1 && \sin4a=0 \label{cos-4a=1-sin-4a=0}
 \end{align}
 \elm
\bpr Первая строчка:
 $$
\cos a=0
$$
$$
\Downarrow
$$
$$
\sin^2a=1-\cos^2a=1
$$
$$
\Downarrow
$$
$$
\left[\begin{matrix}\sin a=1 & \\
\sin a=-1 & \text{\scriptsize$\leftarrow$ невозможно, в силу леммы
\ref{LM:sin>0-na-(0,2)}}
\end{matrix}\right]
$$
$$
\Downarrow
$$
$$
\sin a=1
$$
Вторая строчка:
$$
\cos2a=\cos^2a-\sin^2a=0-1=-1,
$$
$$
\sin2a=2\cdot\sin a\cdot\cos a=2\cdot 1\cdot 0=0
$$
Третья строчка:
$$
\cos4a=\cos^22a-\sin^22a=0-1=-1,
$$
$$
\sin2a=2\cdot\sin a\cdot\cos a=2\cdot 1\cdot 0=0
$$
\epr

\blm\label{LM:periodichnost-sin-cos} Число $4a$, где $a\in(0,2)$ определено
условием \eqref{DEF:cos-a=0}, является периодом для функций $\sin$ и $\cos$,
определенных рядами \eqref{DF:sin}-\eqref{DF:cos}:
 \begin{align*}
&\sin(x+4a)=\sin x, && \cos(x+4a)=\cos x
 \end{align*}
\elm
 \bpr Здесь просто применяется \eqref{cos-4a=1-sin-4a=0}:
$$
\sin(x+4a)=\sin x\cdot\underbrace{\cos 4a}_{1}+\cos
x\cdot\underbrace{\sin4a}_{0}=\sin x
$$
$$
\cos(x+4a)=\cos x\cdot\underbrace{\cos 4a}_{1}-\sin
x\cdot\underbrace{\sin4a}_{0}=\cos x
$$
 \epr

\blm\label{LM:0<sin-x<x-2} При $x\in(0;1)$ для функций $\sin$ и $\cos$,
определенных рядами \eqref{DF:sin}-\eqref{DF:cos} выполняется следующее тройное
неравенство:
 \beq\label{0<sin-x<x-2}
0<x\cos x<\sin x<x.
 \eeq
\elm
 \bpr Здесь всюду используется прием, который мы применяли при доказательстве
леммы \ref{LM:sin>0-na-(0,2)} -- нужно представить ряд, которым записывается
функция, в удобном для оценок виде.

1. Первое неравенство:
 \begin{multline*}
x\cdot\cos x=x\cdot\sum_{n=0}^\infty (-1)^n\cdot \frac{x^{2n}}{(2n)!}=\\=
x-\frac{x^3}{2!}+\frac{x^5}{4!}-\frac{x^7}{6!}+...=\\=
x\cdot\left(1-\frac{x^2}{2}\right)+\frac{x^5}{4!}\cdot\left(1-\frac{x^2}{5\cdot
6}\right)+...=\\=
x\cdot\underbrace{\left(1-\frac{x^2}{2}\right)}_{\scriptsize\begin{matrix}
\text{\rotatebox{90}{$<$}}\\ \phantom{,} 0, \\ \text{при
$x\in(0,1)$}\end{matrix}}+\\+\sum_{k=1}^\infty
\frac{x^{4k+1}}{(4k)!}\cdot\underbrace{\left(1-\frac{x^2}{(4k+1)\cdot(4k+2)}\right)}_{\scriptsize\begin{matrix}
\text{\rotatebox{90}{$<$}}\\ \phantom{,} 0, \\ \text{при
$x\in(0,1)$}\end{matrix}}>0
 \end{multline*}

2. Второе неравенство:
 \begin{multline*}
\sin x-x\cdot\cos x=\\=\sum_{n=0}^\infty (-1)^n\cdot
\frac{x^{2n+1}}{(2n+1)!}-x\cdot\sum_{n=0}^\infty (-1)^n\cdot
\frac{x^{2n}}{(2n)!}=\\=\sum_{n=0}^\infty (-1)^n\cdot
\frac{x^{2n+1}}{(2n+1)!}-\sum_{n=0}^\infty (-1)^n\cdot
\frac{x^{2n+1}}{(2n)!}=\\= \sum_{n=0}^\infty (-1)^n\cdot
\left(\frac{x^{2n+1}}{(2n+1)!}-\frac{x^{2n+1}}{(2n)!}\right)=\\=
 \sum_{n=0}^\infty (-1)^n\cdot\frac{x^{2n+1}}{(2n+1)!}\cdot
\Big(1-(2n+1)\Big)=\\= \sum_{n=0}^\infty
(-1)^{n+1}\cdot\frac{x^{2n+1}}{(2n+1)!}\cdot 2n=\\= 0+ \frac{x^3}{3!}\cdot 2-
\frac{x^5}{5!}\cdot 4+\frac{x^7}{7!}\cdot 6-\frac{x^9}{9!}\cdot 8+...=\\=
\frac{x^3}{3!}\cdot \left(2- \frac{x^2}{5}\right)+\frac{x^7}{7!}\cdot\left(
6-\frac{x^2}{9}\cdot
8\right)+...=\\=\sum_{k=0}^\infty\frac{x^{4k+3}}{(4k+3)!}\cdot\underbrace{\left((4k+2)-
\frac{x^2}{4k+5}\right)}_{\scriptsize\begin{matrix}
\text{\rotatebox{90}{$<$}}\\ \phantom{,} 0, \\ \text{при
$x\in(0,1)$}\end{matrix}}>0
 \end{multline*}

3. Третье неравенство:
 \begin{multline*}
x-\sin x=x-\sum_{n=0}^\infty (-1)^n\cdot \frac{x^{2n+1}}{(2n+1)!}=\\=
\frac{x^3}{3!}-\frac{x^5}{5!}+\frac{x^7}{7!}-\frac{x^9}{9!}+...=\\=
\frac{x^3}{3!}\cdot\left(1-\frac{x^2}{4\cdot
5}\right)+\frac{x^7}{7!}\cdot\left(1-\frac{x^2}{8\cdot 9}\right)+...=\\=
\sum_{k=0}^\infty
\frac{x^{4k+3}}{(4k+3)!}\cdot\underbrace{\left(1-\frac{x^2}{(4k+4)\cdot(4k+5)}\right)}_{\scriptsize\begin{matrix}
\text{\rotatebox{90}{$<$}}\\ \phantom{,} 0, \\ \text{при
$x\in(0,1)$}\end{matrix}}>0
 \end{multline*}
 \epr

\bpr[Доказательство избыточности Аксиомы B2]
\label{PROOF:osn-teor-trigonometrii}

1. Лемма \ref{LM:sin-cos} определяет функции $\sin$ и $\cos$ всюду на числовой
прямой $\R$, а по лемме \ref{LM:periodichnost-sin-cos} эти функции являются
периодическими. Это доказывает условие $T_0$ Аксиомы тригонометрии.

2. Условие $T_1$ доказывается леммами \ref{LM:sin(x+y)-cos(x+y)} и
\ref{LM:sin^2-x+cos^2-x=1}.

3. Условие $T_1$ доказывается леммой \ref{LM:0<sin-x<x-2}.

4. Остается заметить, что единственность функций $\sin$ и $\cos$,
удовлетворяющих условиям $T_0$-$T_1$ Аксиомы тригонометрии, уже доказана: в
иллюстративном тексте мы, последовательно выводя следствия из Аксиомы B2,
получили, что если какие-то функции $\sin$ и $\cos$ удовлетворяют условиям
$T_0$-$T_2$ Аксиомы B2, то они автоматически должны описываться рядами
\eqref{21.3.5} и \eqref{21.3.6}, или, что то же самое, рядами
\eqref{DF:sin}-\eqref{DF:cos}.\epr

\end{multicols}\noindent\rule[10pt]{160mm}{0.1pt}

\subsection{Значение числа $\pi$}

В заключение разговора о степенных рядах и аналитических функциях отметим еще
одно приложение этой теории -- формулу для нахождения числа $\pi$, которую мы
обещали читателю еще на с.\pageref{DEF:pi}, когда определяли это число.

\noindent\rule{160mm}{0.1pt}\begin{multicols}{2}

\bprop Число $\pi$ представимо в виде суммы числового ряда формулой
 \beq\label{ryad-dlya-pi}
\pi=\sum_{n=0}^\infty\frac{(-1)^n\cdot 4}{2n+1},
 \eeq
и для всякого $N\in\N$ удовлетворяет оценке
 \beq\label{otsenka-dlya-pi}
\left|\pi-\sum_{n=0}^N \frac{(-1)^n\cdot 4}{2n+1}\right|\le \frac{4}{2N+3}
 \eeq
или, что то же самое, оценке
 \begin{multline}\label{otsenka-dlya-pi-1}
\sum_{n=0}^N \frac{(-1)^n\cdot 4}{2n+1}-\frac{4}{2N+3}<\pi<\\< \sum_{n=0}^N
\frac{(-1)^n\cdot 4}{2n+1}+\frac{4}{2N+3}
 \end{multline}
 \eprop
 \bpr
Рассмотрим функциональный ряд
 \beq\label{ryad-dlya-vych-pi}
\sum_{n=0}^\infty(-1)^n\cdot \frac{x^{2n+1}}{2n+1},\qquad x\in[0,1]
 \eeq
Если обозначить $b_n(x)=\frac{x^{2n+1}}{2n+1}$, то при $x\in[0,1]$ это будет
неотрицательная невозрастающая последовательность
$$
b_0(x)\ge b_1(x)\ge b_2(x)\ge ... \ge 0
$$
стремящаяся к нулю равномерно на $[0,1]$:
$$
||b_n||_{[0,1]}= \frac{1}{2n+1}\underset{n\to \infty}{\longrightarrow} 0
$$
То есть ряд \eqref{ryad-dlya-vych-pi} удовлетворяет условиям теоремы
\ref{TH:priznak-Leibnitza-ravnom-shod} (признак Лейбница равномерной
сходимости), и значит, он равномерно сходится на отрезке $[0,1]$. Обозначим его
сумму буквой $f$:
$$
f(x)=\sum_{n=0}^\infty(-1)^n\cdot \frac{x^{2n+1}}{2n+1},\qquad x\in[0,1]
$$
Как сумма равномерно сходящегося ряда из непрерывных функций, функция $f$
должна быть непрерывна на отрезке $[0,1]$. С другой стороны, в силу тождества
\eqref{21.3.12}, функция $f$ совпадает с функцией $\arctg$ на полуинтервале
$[0,1)$:
$$
\forall x\in[0,1)\qquad \arctg x=\sum_{n=0}^\infty(-1)^n\cdot
\frac{x^{2n+1}}{2n+1}=f(x)
$$
При этом, как мы знаем, функция $\arctg$ непрерывна на $\R$ (и значит на
$[0,1]$) в силу предложения \ref{PROP:nepr-arcsin-arccos}. Отсюда следует, что
функции $f$ и $\arctg$ совпадают в точке 1:
$$
\arctg 1=\lim_{x\to 1-0}\arctg x=\lim_{x\to 1-0}f(x)=f(1)
$$
Теперь мы получаем:
 \begin{multline*}
\underbrace{\frac{\pi}{4}=\arctg 1}_{\scriptsize\begin{matrix}\text{по
определению}\\ \text{арктангенса
\eqref{DEF:arctg}}\end{matrix}}=f(1)=\\=\sum_{n=0}^\infty(-1)^n\cdot
\frac{x^{2n+1}}{2n+1}\Big|_{x=1}=\sum_{n=0}^\infty(-1)^n\cdot \frac{1}{2n+1}
 \end{multline*}
Умножая это на 4, мы получим формулу \eqref{ryad-dlya-pi}, а из оценки остатка
знакочередующегося ряда \eqref{otsenka-ostatka-Leibnitz} получаем оценку
\eqref{otsenka-dlya-pi}.
 \epr

\brem Как и в случае с числом Непера $e$, значение которого, как мы помним
оценивались с помощью неравенства \eqref{0<e-sum-1/k!<1/n!n-1}, значения $\pi$
также можно оценивать, используя теперь неравенство \eqref{otsenka-dlya-pi-1},
однако качественное различие между этими двумя случаями будет в том, что нужная
точность для $\pi$ достигается гораздо медленнее, чем для $e$. Например, для
вычисления первых двух знаков после запятой для $\pi$ нужно брать сумму 1000
слагаемых: во-первых,
$$
\begin{matrix}
\frac{4}{2\cdot 1000+3}=0,0019970...
 \\
 \Downarrow
 \\
0,001<\frac{4}{2\cdot 1000+3}<0,002
 \\
 \Downarrow
 \\
-0,002<-\frac{4}{2\cdot 1000+3}<-0,001
\end{matrix}
$$
во-вторых,
 $$
\begin{matrix}
\sum_{n=0}^{1000}\frac{(-1)^n\cdot 4}{2n+1}=3,14259165... \\
\Downarrow \\
3,142<\sum_{n=0}^{1000}\frac{(-1)^n\cdot 4}{2n+1}<3,143
\end{matrix}
 $$
и вместе это дает
 \begin{multline*}
3,140=3,142-0,002<\\<\sum_{n=0}^{1000}\frac{(-1)^n\cdot
4}{2n+1}-\frac{4}{2\cdot 1000+3}<\pi<\\<\sum_{n=0}^{1000}\frac{(-1)^n\cdot
4}{2n+1}-\frac{4}{2\cdot 1000+3}<\\<3,142+0,002=3,144
 \end{multline*}
$$
\Downarrow
$$
 \beq\label{otsenka-pi}
\kern32pt\boxed{\quad 3,14<\pi<3,15\quad}
 \eeq
$$
\Downarrow
$$
 \beq\label{otsenka-pi-1}
\kern32pt\boxed{\quad \pi=3,14...\quad}
 \eeq
\erem

\end{multicols}\noindent\rule[10pt]{160mm}{0.1pt}

\subsection{Числа Фибоначчи}\label{Fibonacci}

Пример приложения теории степенных рядов, не упомянуть который нельзя, --
способ нахождения явных формул для последовательностей, заданных рекуррентно.
Мы проиллюстрируем его на последовательности Фибоначчи, описанной нами в
примере \ref{EX:rekurr-posledov}.

\noindent\rule{160mm}{0.1pt}\begin{multicols}{2}

Напомним, что числа Фибоначчи определяются рекуррентными соотношениями
\eqref{DEF:Fibonacci}:
$$
x_1=x_2=1,\quad x_{n+2}=x_{n+1}+x_n, \quad n\in\N
$$

\bprop Числа Фибоначчи описываются явно формулой
 \beq\label{formula-dlya-Fibonacci}
x_n=\frac{1}{\sqrt{5}}\left[\l\frac{1+\sqrt{5}}{2}\r^n-\l\frac{1-\sqrt{5}}{2}\r^n\right]
 \eeq
 \bpr
Рассмотрим производящую функцию этой последовательности:
$$
F(s)=\sum_{n=1}^\infty x_n\cdot s^n
$$
(суммирование начинается с номера $n=1$; это можно объяснить, например, так: мы
будем считать, что $x_0=0$). Помножив $F(s)$ на многочлен $s+s^2$, мы получим
цепочку:
 \begin{multline*}
(s+s^2)\cdot F(s)=(s+s^2)\cdot \sum_{n=1}^\infty x_n\cdot s^n=
\\=\underbrace{\sum_{n=1}^\infty x_n\cdot s^{n+1}}_{\text{замена:
$n+1=k$}}+\underbrace{\sum_{n=1}^\infty x_n\cdot s^{n+2}}_{\text{замена:
$n+2=k$}}=
\\= \sum_{k=2}^\infty x_{k-1}\cdot s^k+\sum_{k=3}^\infty x_{k-2}\cdot
s^k=
\\= \underbrace{x_1}_{\scriptsize\begin{matrix}\text{\rotatebox{90}{$=$}}
\\ 1\end{matrix}}\cdot s^2+ \sum_{k=3}^\infty x_{k-1}\cdot s^k+\sum_{k=3}^\infty
x_{k-2}\cdot s^k=
\\= s^2+ \sum_{k=3}^\infty \underbrace{(x_{k-1}+x_{k-2})}_{x_k}\cdot s^k=
\\= -s+\underbrace{s+s^2}_{\scriptsize\begin{matrix} \text{\rotatebox{90}{$=$}}\\ x_1\cdot s^1+x_2\cdot s^2\end{matrix}}+ \sum_{k=3}^\infty x_k\cdot s^k=
\\= -s+\sum_{k=1}^\infty x_k\cdot s^k=-s+F(s)
 \end{multline*}
$$
\Downarrow
$$
$$
(s+s^2)\cdot F(s)=-s+F(s)
$$
$$
\Downarrow
$$
$$
s=(1-s-s^2)\cdot F(s)
$$
$$
\Downarrow
$$
 \beq\label{F(s)-fibonacci}
F(s)=\frac{s}{1-s-s^2}
 \eeq
Разложим знаменатель на множители:
$$
1-s-s^2=-(s-a)(s-b)
$$
где
$$
a=\frac{-1-\sqrt{5}}{2},\quad b=\frac{-1+\sqrt{5}}{2}
$$
Заметим попутно, что
 \beq\label{b-a-Fibonacci}
b-a=\frac{-1+\sqrt{5}}{2}+\frac{1+\sqrt{5}}{2}=\sqrt{5}
 \eeq
 \begin{multline}\label{1/a-Fibonacci}
\frac{1}{a}=-\frac{2}{1+\sqrt{5}}=-\frac{2\cdot (1-\sqrt{5})}{(1+\sqrt{5})\cdot
(1-\sqrt{5})}=
\\=-\frac{2\cdot (1-\sqrt{5})}{1-5}=\frac{1-\sqrt{5}}{2}
 \end{multline}
 \begin{multline}\label{1/b-Fibonacci}
\frac{1}{b}=-\frac{2}{1-\sqrt{5}}=-\frac{2\cdot (1+\sqrt{5})}{(1-\sqrt{5})\cdot
(1+\sqrt{5})}=
\\=-\frac{2\cdot (1+\sqrt{5})}{1-5}=\frac{1+\sqrt{5}}{2}
 \end{multline}
Теперь получаем:
 \begin{multline*}
F(s)=\frac{s}{1-s-s^2}=-s\cdot\underbrace{\frac{1}{(s-a)(s-b)}}_{\scriptsize\begin{matrix}\text{раскладываем}\\
\text{на простейшие дроби}\\ \text{по правилам на
с.\pageref{pravila-razlozh-na-prost-drobi}}
\end{matrix}}=
\\=-s\cdot \l\frac{\frac{1}{a-b}}{s-a}+\frac{\frac{1}{b-a}}{s-b}\r=
\frac{s}{b-a}\cdot\l\frac{1}{s-b}-\frac{1}{s-a}\r=
\\=\frac{s}{b-a}\cdot\l\frac{1}{a-s}-\frac{1}{b-s}\r=
\\=\frac{s}{b-a}\cdot\l\frac{1}{a}\cdot\frac{1}{1-\frac{s}{a}}-\frac{1}{b}\cdot\frac{1}{1-\frac{s}{a}}\r=\eqref{mclaurent-1/(1-x)}=
\\=\frac{s}{b-a}\cdot\l\frac{1}{a}\cdot\sum_{n=0}^\infty\l\frac{s}{a}\r^n-\frac{1}{b}\cdot\sum_{n=0}^\infty\l\frac{s}{b}\r^n\r=
\\=\frac{s}{b-a}\cdot\l\sum_{n=0}^\infty\frac{s^n}{a^{n+1}}-\sum_{n=0}^\infty\frac{s^n}{b^{n+1}}\r=
\\=\frac{1}{b-a}\cdot\sum_{n=0}^\infty\l\frac{1}{a^{n+1}}-\frac{1}{b^{n+1}}\r\cdot
s^{n+1}=\left|n+1=k\right|=
\\=\frac{1}{b-a}\cdot\sum_{k=1}^\infty\l\frac{1}{a^k}-\frac{1}{b^k}\r\cdot
s^k=
\\=\eqref{b-a-Fibonacci},\eqref{1/a-Fibonacci},\eqref{1/b-Fibonacci}=
\\=\frac{1}{\sqrt{5}}\cdot\sum_{k=1}^\infty\left[\l\frac{1+\sqrt{5}}{2}\r^k-\l\frac{1-\sqrt{5}}{2}\r^k\right]\cdot s^k
 \end{multline*}
Выделяя коэффициенты перед $s^k$, мы теперь получим
\eqref{formula-dlya-Fibonacci}.
 \epr
\eprop

\end{multicols}\noindent\rule[10pt]{160mm}{0.1pt}

\chapter{ТРИГОНОМЕТРИЧЕСКИЕ РЯДЫ}\label{CH-trig-series}

 \bit{
\item[$\bullet$] {\it Тригонометрическим рядом}\index{ряд!тригонометрический}
называется функциональный ряд вида
 \begin{equation}\label{trig-ryad}
  \frac{a_0}{2}+\sum_{n=1}^\infty \Big\{ a_n\cdot\cos\frac{\pi nx}{T}+b_n\cdot\sin\frac{\pi nx}{T} \Big\}
  \end{equation}
Число $T$ называется {\it полупериодом} этого  тригонометрического ряда, а
числа $a_n$ и $b_n$ -- его {\it коэффициентами}.

\item[$\bullet$] Частным случаем тригонометрического ряда будет ряд, в котором
коэффициенты $a_n$ и $b_n$ найдены по формулам
 \begin{equation}\label{koeff-Fourier}
\begin{split}
& a_0=\frac{1}{T}\int_{-T}^T f(x) \, \d x, \\
& a_n=\frac{1}{T}\int_{-T}^T f(x) \cos \frac{\pi nx}{T} \, \d x, \\
& b_n=\frac{1}{T}\int_{-T}^T f(x) \sin \frac{\pi nx}{T} \, \d x
\end{split}
 \end{equation}
где $f$ -- некоторая интегрируемая на отрезке $[-T,T]$ функция. Такой
тригонометрический ряд называется {\it рядом Фурье} функции $f$ с полупериодом
$T$. Его коэффициенты \eqref{koeff-Fourier} называются {\it коэффициентами
Фурье} функции $f$, а частичные суммы
 \beq\label{mnogochlen-Fourier}
S_N(x)=\frac{a_0}{2}+\sum_{n=1}^N \Big\{ a_n\cdot\cos\frac{\pi
nx}{T}+b_n\cdot\sin\frac{\pi nx}{T} \Big\}
 \eeq
-- {\it многочленами Фурье} функции $f$.
 }\eit
В этой главе нас будет интересовать вопрос, какой должна быть функция $f$,
чтобы быть суммой своего ряда Фурье:
 \beq\label{f=summa-ryada-Fourier}
f(x)=\frac{a_0}{2}+\sum_{n=1}^\infty \Big\{ a_n\cdot\cos\frac{\pi
nx}{T}+b_n\cdot\sin\frac{\pi nx}{T} \Big\}
 \eeq
Понятно, что первым необходимым условием для этого должна быть интегрируемость
функции $f$ на отрезке $[-T,T]$ (потому что иначе будет непонятно, что такое
ряд Фурье для нее, то есть как понимать интегралы \eqref{koeff-Fourier}). С
другой стороны, функции $x\mapsto\cos\frac{\pi nx}{T}$ и $x\mapsto\sin\frac{\pi
nx}{T}$ в тригонометрическом ряде \eqref{trig-ryad} имеют период $2T$:
$$
\cos\frac{\pi n(x+2T)}{T}=\cos\l\frac{\pi nx}{T}+2\pi n\r=\cos\frac{\pi
nx}{T},\qquad \sin\frac{\pi n(x+2T)}{T}=\sin\l\frac{\pi nx}{T}+2\pi
n\r=\sin\frac{\pi nx}{T}
$$
поэтому, если функция $f$ является суммой своего ряда Фурье
\eqref{f=summa-ryada-Fourier}, то она тоже должна иметь период $2T$:
$$
f(x+2T)=f(x)
$$
Отсюда становится понятно, что предметом изучения для нас должны быть
периодические функции $f$, интегрируемые на отрезке $[-T,T]$, где $T$ --
период. Всякая такая функция будет интегрируема не только на отрезке $[-T,T]$,
но и вообще на любом отрезке $[a,b]\subset\R$, поэтому нам будет удобно
употреблять для таких функций следующий термин:
 \bit{
\item[$\bullet$] Функция $f:\R\to\R$ называется {\it локально интегрируемой},
если она интегрируема на каждом отрезке $[a,b]\subset\R$.
 }\eit

\section{Сходимость ряда Фурье в среднем квадратичном}

Изучение рядов Фурье удобно проводить для случая, когда в качестве полупериода
$T$ выбрано число $\pi$, потому что тогда формулы \eqref{trig-ryad} и
\eqref{koeff-Fourier} упрощаются: тригонометрический ряд \eqref{trig-ryad}
принимает вид
 \begin{equation}\label{trig-ryad-pi}
  \frac{a_0}{2}+\sum_{n=1}^\infty \Big\{ a_n\cos nx+b_n\sin nx \Big\},
  \end{equation}
коэффициенты Фурье \eqref{koeff-Fourier} вычисляются по формулам:
 \begin{equation}\label{koeff-Fourier-pi}
a_0=\frac{1}{\pi}\int\limits_{-\pi}^\pi f(x) \, \d x, \quad
a_n=\frac{1}{\pi}\int\limits_{-\pi}^\pi f(x) \cos nx \, \d x, \quad
b_n=\frac{1}{\pi}\int\limits_{-\pi}^\pi f(x) \sin nx \, \d x
 \end{equation}
а многочлены Фурье записываются в виде
 \beq\label{mnogochlen-Fourier-pi}
S_N(x):= \frac{a_0}{2}+\sum_{n=1}^N \Big\{ a_n\cos nx+b_n\sin nx \Big\}
 \eeq
(здесь $a_n$ и $b_n$ -- коэффициенты Фурье \eqref{koeff-Fourier-pi}).

В соответствии с этим мы будем всюду далее до теоремы \ref{tm-22.2.10}
рассматривать локально интегрируемые $2\pi$-периодические функции (и их ряды
Фурье с полупериодом $\pi$). Множество всех таких функций удобно обозначить
каким-нибудь символом:
 \bit{
\item[$\bullet$] Символом ${\mathcal R}(\pi)$ мы обозначаем множество всех
локально интегрируемых функций с полупериодом $\pi$ на $\R$.
 }\eit

\subsection{Основная теорема}

 \bit{
\item[$\bullet$] Пусть $f\in{\mathcal R}(\pi)$, то есть $f$ -- локально
интегрируемая функция с полупериодом $\pi$ на $\R$. Поскольку $f$ интегрируема
на отрезке $[-\pi,\pi]$, по свойству $8^0$ на с.\pageref{14.6.5} ее квадрат
$f^2=|f|^2$ тоже интегрируем на отрезке $[-\pi,\pi]$. Значит, определено число
 \beq\label{norma-v-srednem-kvadrat}
\norm{f}=\sqrt{\frac{1}{\pi}\int_{-\pi}^{\pi} |f(x)|^2\ \d x}
 \eeq
Оно называется {\it среднеквадратичной нормой} функции $f$ в пространстве
${\mathcal R}(\pi)$. Слово ``норма'' употребляется применительно к этой
величине потому, что, как мы убедимся ниже на с.\pageref{norm_2-ge-0},
отображение $f\mapsto\norm{f}$ обладает свойствами $1^\circ$-$3^\circ$
равномерной нормы, упоминавшимися на с \pageref{svoistva-ravnom-normy} (правда,
без условия $\norm{f}=0$ $\Leftrightarrow$ $f=0$).
 }\eit

Весь этот параграф посвящен доказательству следующего фундаментального факта:

\begin{tm}[\bf о сходимости ряда Фурье в среднем квадратичном]\label{TH:Fourier-v-R_2}
Всякая локально интегрируемая $2\pi$-периодическая функция $f$ на прямой $\R$
единственным образом представима в виде суммы в среднем квадратичном
тригонометрического ряда с полупериодом $\pi$
 \beq\label{razlozhenie-Fourier-v-R_2}
f(x)=\frac{a_0}{2}+\sum_{n=1}^\infty \Big\{ a_n\cos nx+b_n\sin nx \Big\},
 \eeq
и таким рядом будет ряд Фурье с полупериодом $\pi$ функции $f$.
 \end{tm}
\brem Эти слова следует понимать, как выполнение двух утверждений:
 \bit{
\item[(i)] если $a_n$ и $b_n$ -- коэффициенты Фурье \eqref{koeff-Fourier-pi}
функции $f$, а $S_N$ -- ее многочлены Фурье \eqref{mnogochlen-Fourier-pi}, то
справедливо соотношение
 \begin{equation}\label{||f-S_N||_2->0}
\norm{f-S_N}=\sqrt{\frac{1}{\pi}\int_{-\pi}^{\pi} |f(x)-S_N(x)|^2\ \d
x}\underset{N\to\infty}{\longrightarrow}0
 \end{equation}
(именно так интерпретируется в теореме \ref{TH:Fourier-v-R_2} равенство
\eqref{razlozhenie-Fourier-v-R_2}), и

\item[(ii)] если $a_n$ и $b_n$ -- две последовательности чисел, а $S_N$ --
частичные суммы соответствующего им тригонометрического ряда
\eqref{trig-ryad-pi}, то из соотношения \eqref{||f-S_N||_2->0} следует, что
числа $a_n$ и $b_n$ являются коэффициентами Фурье \eqref{koeff-Fourier-pi}
функции $f$ (а $S_N$ -- ее многочленами Фурье \eqref{mnogochlen-Fourier-pi}).
 }\eit
 \erem

Неподготовленному читателю формулировка теоремы \ref{TH:Fourier-v-R_2},
несомненно, покажется обескураживающей, потому что из сказанного до сих пор
должно быть решительно непонятно, отчего разложение функции $f$ в ряд Фурье
\eqref{razlozhenie-Fourier-v-R_2} следует понимать в таком экзотическом смысле
(а не, например, как поточечную или равномерную сходимость, к которым мы
привыкли в главе \ref{CH-functional-sequen}). Объяснение заключается в том, что
теорема \ref{TH:Fourier-v-R_2} является одним из результатов весьма глубокой и
имеющей долгую историю математической дисциплины, с которой мы столкнулись
только сейчас, -- {\it гармонического анализа}. То, что сформировалось в этой
науке, как результат длительного осмысления и привыкания, встретилось нам сразу
в виде готового результата, и поэтому выглядит таким странным и неожиданным.
При этом, если в гармоническом анализе и есть результат, который можно признать
центральным, то им, конечно, будет теорема о разложении функции в ряд Фурье в
смысле среднего квадратичного (а не поточечно или в каком-нибудь другом
смысле), и именно поэтому мы начинаем изложение теории рядов Фурье с теоремы
\ref{TH:Fourier-v-R_2}. Оставшаяся часть этого параграфа посвящена ее
доказательству.

\paragraph{Скалярное произведение функций.}

Оригинальная идея, приведшая к теореме \ref{TH:Fourier-v-R_2}, заключается в
том, чтобы посмотреть на функции из ${\mathcal R}(\pi)$ с несколько неожиданной
точки зрения: как на векторы в пространстве со скалярным произведением.

 \bit{
\item[$\bullet$] Пусть $f,g\in{\mathcal R}(\pi)$, то есть $f$ и $g$ -- две
локально интегрируемые $2\pi$-периодические функции. Поскольку они интегрируемы
на отрезке $[-\pi,\pi]$, их произведение $f\cdot g$ тоже будет интегрируемо на
$[-\pi,\pi]$ в силу свойства $8^0$ на с.\pageref{integrir-proizvedeniya}.
Поэтому можно рассмотреть величину
 \beq
\langle f,g\rangle = \frac{1}{\pi}\int_{-\pi}^{\pi} f(x)\cdot g(x) \, \d x
\label{22.3.1}
 \eeq
Она называется {\it скалярным произведением}\index{скалярное
произведение!функций} функций $f$ и $g$. Перечислим некоторые
 }\eit

\bigskip

\centerline{\bf Свойства скалярного произведения}

\bigskip

 \bit{\it
\item[$1^0$.] Положительная полуопределенность:
$$
  \langle f,f\rangle\ge 0
$$
\item[$2^0$.] Симметричность:
$$
  \langle f,g\rangle=\langle g,f\rangle
$$
\item[$3^0$.] Билинейность:
$$
  \langle \alpha \cdot f+\beta \cdot g,h\rangle=\alpha \cdot \langle f,h\rangle+\beta \cdot \langle g,h\rangle
$$

 }\eit

\bigskip

\begin{proof} Эти свойства очевидны, например, третье следует из
линейности интеграла:
\begin{multline*}\langle \alpha \cdot f+\beta \cdot g,h\rangle=
\frac{1}{\pi}\int_{-\pi}^{\pi}\Big\{ \alpha \cdot f(x)+\beta \cdot
g(x)\Big\}\cdot h(x) \, \d x= \frac{1}{\pi}\int_{-\pi}^{\pi}\Big\{ \alpha \cdot
f(x)\cdot h(x)+\beta \cdot g(x)\cdot h(x) \Big\}\, \d x=\\= \alpha \cdot
\frac{1}{\pi}\int_{-\pi}^{\pi} f(x)\cdot h(x) \, \d x+ \beta \cdot
\frac{1}{\pi}\int_{-\pi}^{\pi} g(x)\cdot h(x)\, \d x= \alpha \cdot \langle
f,h\rangle+\beta \cdot \langle g,h\rangle
\end{multline*}\end{proof}

Еще одно свойство скалярного произведения -- формула, связывающая его со
среднеквадратичной нормой \eqref{norma-v-srednem-kvadrat}:
 \beq
\norm{f}=\sqrt{\frac{1}{\pi}\int_{-\pi}^{\pi} |f(x)|^2\ \d x}=\sqrt{\langle
f,f\rangle}
 \eeq

\btm[\bf неравенство Коши-Шварца для
функций]\label{TH-Cauchy-Schwarz-dlya-funcij} Для любых функций
$f,g\in{\mathcal R}(\pi)$ справедливо следующее неравенство, называемое
неравенством Коши-Шварца\index{неравенство!Коши-Шварца!для функций}:
 \begin{equation}\label{Cauchy-Schwarz-dlya-functsij}
 \big|\langle f,g\rangle\big|\le  \norm{f}\cdot\norm{g}
 \end{equation}
 \etm
 \begin{proof}
Зафиксируем $f$ и $g$ и рассмотрим многочлен второй степени:
 $$
 p(t)=\langle f+tg,f+tg\rangle=
 \langle f,f\rangle+2t\langle f,g\rangle+t^2\langle g,g\rangle=
 \norm{f}^2+2t\langle f,g\rangle+t^2\norm{g}^2
 $$
 Поскольку $\langle f+tg,f+tg\rangle\ge 0$, этот многочлен должен
 при любом значении $t$ быть неотрицателен:
 $$
 p(t)=\norm{f}^2+2t\langle f,g\rangle+t^2\norm{g}^2\ge 0, \quad \forall
 t\in\R
 $$
 Значит, он не может иметь двух вещественных корней, и
 поэтому его дискриминант должен быть неположительным:
 $$
 D=4\langle f,g\rangle^2-4\norm{f}^2\cdot\norm{g}^2\le 0
 $$
 $$
 \Downarrow
 $$
 $$
 \langle f,g\rangle^2\le \norm{f}^2\cdot\norm{g}^2
 $$
 $$
 \Downarrow
 $$
 $$
 \langle f,g\rangle\le \Big|\langle f,g\rangle\Big|\le \norm{f}\cdot\norm{g}
 $$
 \end{proof}

\bigskip

\centerline{\bf Свойства среднеквадратичной нормы:}

\bit{\it

\item[$1^\circ$.] Неотрицательность:
 \begin{equation}\label{norm_2-ge-0}
\norm{f}\ge 0.
 \end{equation}

\item[$2^\circ$.] Однородность:
 \begin{equation}\label{norm_2-lambda-x}
\norm{\lambda\cdot f}=|\lambda|\cdot\norm{f},\qquad \lambda\in \R
 \end{equation}

\item[$3^\circ$.] Неравенство треугольника:
 \begin{equation}\label{norm_2-triangle}
\norm{f+g}\le \norm{f}+\norm{g}.
 \end{equation}

\item[$4^\circ$.] Тождество параллелограмма:
 \begin{equation}\label{norm_2-tozhd-parallelogramma}
\frac{\norm{f+g}^2 +\norm{f-g}^2}{2}=\norm{f}^2+\norm{g}^2.
 \end{equation}
 }\eit

\bpr Первые два свойства очевидны, докажем неравенство треугольника
\eqref{norm_2-triangle}:
 $$
 \norm{f+g}\le\norm{f}+\norm{g}
 \quad\Longleftrightarrow\quad
\underbrace{\norm{f+g}^2}_{\scriptsize\begin{matrix}\|
\\ \langle f+g,f+g\rangle\\ \| \\ \langle f,f\rangle+2\langle f,g\rangle +\langle g,g\rangle
\\ \end{matrix}}\le\underbrace{\big(\norm{f}+\norm{g}\big)^2}_{\scriptsize\begin{matrix}\|
\\ \norm{f}^2+2\cdot\norm{f}\cdot\norm{g}+\norm{g}^2\\ \| \\
\langle f,f\rangle+2\cdot\norm{f}\cdot\norm{g}+\langle g,g\rangle \end{matrix}}
 \quad\Longleftrightarrow\quad
\underbrace{\langle f,g\rangle \le
\norm{f}\cdot\norm{g}}_{\scriptsize\begin{matrix}\uparrow \\ \text{неравенство}\\
\text{Коши-Шварца \eqref{Cauchy-Schwarz-dlya-functsij}} \end{matrix}}
 $$
Тождество параллелограмма \eqref{norm_2-tozhd-parallelogramma} доказывается
прямым вычислением:
$$
\frac{\norm{f+g}^2 +\norm{f-g}^2}{2}=\frac{\norm{f}^2+2\langle f,g\rangle
+\norm{g}^2 +\norm{f}^2-2\langle f,g\rangle
+\norm{g}^2}{2}=\frac{2\norm{f}^2+2\norm{g}^2}{2}=\norm{f}^2+\norm{g}^2
$$
\epr

\btm[\bf о непрерывности скалярного произведения]\label{TH:neprer-skal-proizv}
Если функции $f_n\in{\mathcal R}(\pi)$ стремятся к функции $f\in{\mathcal
R}(\pi)$ в смысле среднеквадратичной нормы,
$$
\norm{f-f_n}\underset{n\to\infty}{\longrightarrow}0,
$$
то для любой функции $g\in{\mathcal R}(\pi)$
$$
\langle f_n,g\rangle\underset{n\to\infty}{\longrightarrow}\langle f,g\rangle
$$
\etm \bpr
$$
\big|\langle f,g\rangle-\langle f_n,g\rangle\big|=\big|\langle
f-f_n,g\rangle\big|\le\eqref{Cauchy-Schwarz-dlya-functsij}\le\norm{f-f_n}\cdot\norm{g}
\underset{n\to\infty}{\longrightarrow}0
$$
\epr

\paragraph{Ортогональность тригонометрической системы и
связь скалярного произведения с коэффициентами Фурье.}

\bit{

\item[$\bullet$] Говорят, что функции $f$ и $g$ из ${\mathcal R}(\pi)$ {\it
ортогональны}\index{ортогональность!функций}, если их скалярное произведение
равно нулю:
$$
  \langle f,g\rangle=0
$$

\item[$\bullet$] Последовательность функций $f_n\in{\mathcal R}(\pi)$
называется {\it ортогональной системой}, если в ней любые две функции
ортогональны:
$$
\forall i\ne j\qquad  \langle f_i,f_j\rangle=0
$$

\item[$\bullet$]\label{ex-22.3.1} Обозначим символами $\cos_n$ и $\sin_n$
функции
 \beq \label{22.3.2}
\cos_n (x)=\cos nx, \qquad \sin_n (x)=\sin nx
 \eeq
и условимся символом $1$ обозначать не только число $1$, но и функцию,
тождественно равную единице:
 \beq\label{22.3.3}
  1(x)=1, \quad x\in \R
 \eeq
 Функции $1$, $\cos_n$, $\sin_n$ называются {\it тригонометрической системой}.
 }\eit

\blm Справедливы следующие формулы, означающие, что тригонометрическая система
функций ортогональна:
 \begin{align}
& \langle 1,1\rangle=2, && \langle 1,\cos_n\rangle=0, && \langle
1,\sin_n\rangle=0, \label{22.3.4} \\
& \langle \cos_n, \sin_m\rangle=0, && \langle \cos_n, \cos_m\rangle= \lll
\begin{array}{c} 0, \quad n\ne m \\ 1, \quad n=m \end{array}\rrr,
&& \langle \sin_n, \sin_m\rangle=\lll \begin{array}{c} 0, \quad n\ne m
\\ 1, \quad n=m
\end{array}\rrr \label{22.3.5}
 \end{align}
\elm
 \bpr
Это проверяется прямым вычислением с помощью тригонометрических формул, которые
мы выводили в \ref{SUBSEC-elem-funktsii} главы \ref{ch-ELEM-FUNCTIONS}.
Например, для произведения синусов при $n\ne m$ мы получим:
 \begin{multline*}
\langle \sin_n, \sin_m\rangle=\frac{1}{\pi}\int_{-\pi}^{\pi}\sin nx\cdot \sin
mx \, \d x=\eqref{sin(x)sin(y)}= \frac{1}{2\pi}\int_{-\pi}^{\pi}\lll \cos
(n-m)x-\cos (n+m)x \rrr \, \d x=\\= \frac{\sin (n-m)
x}{2\pi(n-m)}\Big|_{x=-\pi}^{x=\pi} -\frac{\sin (n+m)
x}{2\pi(n+m)}\Big|_{x=-\pi}^{x=\pi}=0
 \end{multline*}
А при $n=m$ получаем
$$
\langle \sin_n, \sin_n\rangle=\frac{1}{\pi}\int_{-\pi}^{\pi}\sin nx\cdot \sin
nx \, \d x=\eqref{cos-2x}=\frac{1}{2\pi}\int_{-\pi}^{\pi}\lll 1+\cos 2nx \rrr
\, \d x=\frac{2\pi}{2\pi}=1
$$
 \epr

\begin{lm}\label{tm-22.3.2} Коэффициенты ряда Фурье произвольной
функции $f\in{\mathcal R}(\pi)$ связаны со скалярным произведением следующими
формулами:
 \beq\label{22.3.6}
a_0=\Big\langle f,1\Big\rangle, \qquad a_n=\Big\langle f,\cos_n \Big\rangle,
\qquad b_n=\Big\langle f,\sin_n \Big\rangle
 \eeq
\end{lm}\begin{proof} Эти формулы проверяются простым вычислением.
Например,
$$
\langle f,\cos_n\rangle=\eqref{22.3.1}= \frac{1}{\pi}\int_{-\pi}^{\pi}
f(x)\cdot \cos_n (x) \, \d x= \frac{1}{\pi}\int_{-\pi}^{\pi} f(x)\cdot \cos nx
\, \d x=a_n
$$
\end{proof}

\blm Многочлены Фурье произвольной функции $f\in{\mathcal R}(\pi)$
 \beq\label{22.3.7}
S_N=\frac{a_0}{2}\cdot 1+ \sum_{k=1}^N \Big\{ a_k\cdot \cos_k+b_k\cdot \sin_k
\Big\}
 \eeq
имеют те же коэффициенты Фурье порядка не больше $N$, что и функция $f$:
 \beq\label{22.3.8}
 \begin{split}
& \Big\langle S_N,1\Big\rangle=a_0=\Big\langle f,1\Big\rangle, \\
& \Big\langle S_N,\cos_n \Big\rangle=a_n=\Big\langle f,\cos_n \Big\rangle,\qquad n\le N \\
& \Big\langle S_N,\sin_n \Big\rangle=b_n=\Big\langle f,\sin_n
\Big\rangle,\qquad n\le N
 \end{split}
 \eeq
\elm \bpr Действительно,
 \begin{multline*}
\langle S_N,1\rangle= \l \frac{a_0}{2}\cdot 1+\sum_{k=1}^N \Big\{
a_k\cos_k+b_k\sin_k \Big\} , 1 \r= \frac{a_0}{2}\cdot \Big\langle 1,1
\Big\rangle+ \sum_{k=1}^N \Big\{ a_k\cdot \Big\langle  \cos_k, 1\Big\rangle+
b_k\cdot \Big\langle \sin_k, 1\Big\rangle \Big\} =\\= \eqref{22.3.4}=
\frac{a_0}{2}\cdot 2 +\sum_{k=1}^N \Big\{ a_k\cdot 0+b_k\cdot 0 \Big\} =a_0
 \end{multline*}
Аналогично, если $n\le N$, то
 \begin{multline*}
\langle S_N,\cos_n\rangle= \l \frac{a_0}{2}\cdot 1+\sum_{k=1}^N \Big\{
a_k\cos_k+b_k\sin_k \Big\} , \cos_n \r=\\= \frac{a_0}{2}\cdot \Big\langle  1,
\cos_n \Big\rangle+ \sum_{k=1}^N \Big\{ a_k\cdot \Big\langle \cos_k,
\cos_n\Big\rangle+ b_k\cdot \Big\langle \sin_k, \cos_n\Big\rangle \Big\} =\\=
\eqref{22.3.4},\eqref{22.3.5}= \frac{a_0}{2}\cdot 0 +\sum_{k=1}^N \Big\{
a_k\cdot \lll
\begin{array}{cl}
0, & \text{если}\, k\ne n \\ 1, & \text{если}\, k= n
\end{array}\rrr +b_k\cdot 0 \Big\} =a_n
 \end{multline*}
Точно также доказывается последняя формула из \eqref{22.3.8}. \epr

\paragraph{Минимальное свойство многочленов Фурье.}

 \bit{
\item[$\bullet$] Если $f,g\in{\mathcal R}(\pi)$, то {\it расстоянием} от $f$ до
$g$ называется величина
 $$
\norm{f-g}
 $$
 }\eit

\blm Среди всевозможных тригонометрических многочленов $P$ степени $N$ самым
близким к функции $f\in{\mathcal R}(\pi)$ по среднеквадратичной норме является
ее многочлен Фурье $S_N$ степени $N$, расстояние до которого равно $\langle
f,f\rangle-\frac{a_0^2}{2}-\sum_{n=1}^N \Big\{ a_n^2+b_n^2\Big\}$:
 \beq\label{uklonenie-trig-mnogochl}
\min_{P}\norm{f-P}=\norm{f-S_N}=\langle f,f\rangle-\frac{a_0^2}{2}-\sum_{n=1}^N
\Big\{ a_n^2+b_n^2\Big\}
 \eeq
 ($a_n$ и $b_n$ -- коэффициенты Фурье \eqref{22.3.6}).
\elm \bpr Зафиксируем $f$ и обозначим
$$
P= \frac{\alpha_0}{2}+\sum_{n=1}^N \Big\{ \alpha_n\cos_n+\beta_n\sin_n \Big\},
\qquad S_N= \frac{a_0}{2}+\sum_{n=1}^N \Big\{ a_n\cos_n+b_n\sin_n \Big\}
$$
-- здесь числа $a_n$ и $b_n$  вычисляются по формулам \eqref{22.3.6}, поэтому
их мы считаем фиксированными, а числа $\alpha_n$ и $\beta_n$, наоборот, могут
меняться. Тогда:
 \begin{multline*}
\norm{f-P}^2=\langle f-P,f-P\rangle=\left\langle
f-\frac{\alpha_0}{2}-\sum_{n=1}^N \Big\{ \alpha_n\cos_n+\beta_n\sin_n \Big\},
f-\frac{\alpha_0}{2}-\sum_{n=1}^N \Big\{ \alpha_n\cos_n+\beta_n\sin_n
\Big\}\right\rangle=\\=\langle f,f\rangle- \alpha_0\cdot\langle
f,1\rangle-2\cdot\sum_{n=1}^N \Big\{ \alpha_n\cdot\langle
f,\cos_n\rangle+\beta_n\cdot\langle f,\sin_n\rangle
\Big\}+\frac{\alpha_0^2}{4}\cdot\langle 1,1\rangle+\sum_{n=1}^N \Big\{
\alpha_n^2\cdot\langle \cos_n,\cos_n\rangle+\beta_n^2\cdot\langle
\sin_n,\sin_n\rangle\Big\}=\\= \langle f,f\rangle- \alpha_0\cdot
a_0-2\cdot\sum_{n=1}^N \Big\{ \alpha_n\cdot a_n+\beta_n\cdot b_n
\Big\}+\frac{\alpha_0^2}{2}+\sum_{n=1}^N \Big\{ \alpha_n^2+\beta_n^2\Big\}=\\=
\langle f,f\rangle+\underbrace{\frac{(\alpha_0-a_0)^2}{2}+\sum_{n=1}^N \Big\{
(\alpha_n-a_n)^2+(\beta_n-b_n)^2
\Big\}}_{\scriptsize\begin{matrix}\text{достигает минимума при $\alpha_n=a_n$,
$\beta_n=b_n$}\\ \text{то есть при
$P=S_N$}\end{matrix}}-\frac{a_0^2}{2}-\sum_{n=1}^N \Big\{ a_n^2+b_n^2\Big\}
 \end{multline*}
 \epr

\paragraph{Полнота тригонометрической системы.}

\blm\label{LM:polnota-trig-sistemy} Для всякой функции $f\in{\mathcal R}(\pi)$
и любого $\e>0$ найдется тригонометрический многочлен $P$, расстояние от
которого до $f$ меньше $\e$:
 \beq\label{polnota-trig-sistemy}
\norm{f-P}<\e
 \eeq
\elm \bpr Обозначим
$$
M=\sup_{x\in[-\pi,\pi]}|f(x)|
$$
Воспользуемся сначала теоремой \ref{TH:o-priblizh-integrir-func-nepreryvnyni},
и подберем непрерывную функцию $g$ на $[-\pi,\pi]$ так, чтобы выполнялись
условия:
 $$
\int_{-\pi}^{\pi} |f(x)-g(x)|\ \d x<\frac{\pi\cdot\e^2}{8M},\qquad
\sup_{x\in[-\pi,\pi]}|g(x)|\le M,\qquad g(-\pi)=g(\pi)
 $$
Тогда
 \begin{multline*}
\norm{f-g}^2=\frac{1}{\pi}\int_{-\pi}^{\pi} |f(x)-g(x)|^2\ \d x=
\frac{1}{\pi}\int_{-\pi}^{\pi}
\underbrace{|f(x)-g(x)|}_{\scriptsize\begin{matrix}\IA\\ |f(x)|+|g(x)|\\ \IA \\
2M
\end{matrix}}\cdot|f(x)-g(x)|\ \d x\le\\ \le \frac{2M}{\pi}\int_{-\pi}^{\pi} |f(x)-g(x)|\ \d
x<\frac{2M}{\pi}\cdot\frac{\pi\cdot\e^2}{8M}=\frac{\e^2}{4}
 \end{multline*}
 $$
 \Downarrow
 $$
$$
 \norm{f-g}<\frac{\e}{2}
$$
Далее заметим, что поскольку функция $g$ непрерывна на $[-\pi,\pi]$ и принимает
одинаковые значения на концах этого отрезка, ее можно продолжить до непрерывной
$2\pi$-периодической функции на всю прямую $\R$. Тогда воспользовавшись
теоремой \ref{TH:Weierstrass-approx-trig-mnogochl} мы сможем подобрать
тригонометрический многочлен $P$ так, чтобы
$$
\sup_{x\in[-\pi,\pi]}|g(x)-P(x)|<\frac{\e}{2\sqrt{2}}
$$
Для него мы получим:
$$
\norm{g-P}^2=\frac{1}{\pi}\int_{-\pi}^{\pi}
\underbrace{|g(x)-P(x)|^2}_{\scriptsize\begin{matrix}\text{\rotatebox{90}{$>$}}
\\ \frac{\e^2}{8}
\end{matrix}}\ \d x\le \frac{1}{\pi}\int_{-\pi}^{\pi}\frac{\e^2}{8}\ \d x=\frac{\e^2}{4}
$$
 $$
 \Downarrow
 $$
 $$
\norm{g-P}<\frac{\e}{2}
 $$
Вместе получается:
$$
\norm{f-P}=\norm{(f-g)+(g-P)}\le \eqref{norm_2-triangle}\le
\norm{f-g}+\norm{g-P}<\frac{\e}{2}+\frac{\e}{2}=\e
$$
 \epr

\paragraph{Неравенство Бесселя.}

\begin{lm}\label{tm-Bessel}
Пусть $a_n$ и $b_n$ -- коэффициенты Фурье \eqref{22.3.6} функции $f\in{\mathcal
R}(\pi)$. Тогда числовой ряд
$$
  \frac{a_0^2}{2}+\sum_{n=1}^\infty \Big\{ a_n^2+b_n^2 \Big\}
$$
сходится, причем его сумма не превосходит величины $\langle f,f\rangle$:
 \beq
\frac{a_0^2}{2}+\sum_{n=1}^\infty \Big\{ a_n^2+b_n^2 \Big\}\le \langle
f,f\rangle=\frac{1}{\pi}\int_{-\pi}^{\pi} f(x)^2 \, \d x \label{22.3.9}
 \eeq
\end{lm}\begin{proof} 1. Докажем следующую формулу:
 \beq
\langle f,S_N\rangle=\langle S_N,S_N\rangle=\frac{a_0^2}{2}+\sum_{n=1}^N \Big\{
a_n^2+b_n^2 \Big\}\label{22.3.10}
 \eeq
Действительно, с одной стороны,
\begin{multline*}\langle f,S_N\rangle= \left\langle f, \frac{a_0}{2}\cdot 1+
\sum_{k=1}^N \Big\{ a_k\cdot \cos_k+b_k\cdot \sin_k \Big\}\right\rangle=\\=
\frac{a_0}{2}\cdot \Big\langle  f,1 \Big\rangle+ \sum_{k=1}^N \lll a_k\cdot
\Big\langle f,\cos_k\Big\rangle+ b_k\cdot \Big\langle  f,\sin_k\Big\rangle
\rrr= \eqref{22.3.6}= \frac{a_0^2}{2}+\sum_{n=1}^N \Big\{ a_n^2+b_n^2
\Big\}\end{multline*} А с другой,
\begin{multline*}\langle S_N,S_N\rangle= \left\langle S_N, \frac{a_0}{2}\cdot 1+
\sum_{k=1}^N \Big\{ a_k\cdot \cos_k+b_k\cdot \sin_k \Big\}\right\rangle=\\=
\frac{a_0}{2}\cdot \Big\langle  S_N,1 \Big\rangle+ \sum_{k=1}^N \lll a_k\cdot
\Big( S_N,\cos_k\Big\rangle+ b_k\cdot \Big\langle S_N,\sin_k\Big\rangle \rrr=
\eqref{22.3.8}= \frac{a_0^2}{2}+\sum_{n=1}^N \Big\{ a_n^2+b_n^2
\Big\}\end{multline*}

2. Обозначим
$$
  R_N=f-S_N
$$
и заметим, что $S_N$ и $R_N$ ортогональны:
 \beq
\Big\langle  S_N, R_N \Big\rangle=0 \label{22.3.11}
 \eeq
Действительно,
$$
\Big\langle  S_N, R_N \Big\rangle= \Big\langle  S_N, f- S_N \Big\rangle=
\Big\langle  S_N, f \Big\rangle- \Big\langle  S_N,S_N
\Big\rangle=\eqref{22.3.10}=0
$$

3. Теперь докажем формулу
 \beq
  \Big\langle  f,f \Big\rangle=\Big\langle  S_N, S_N \Big\rangle+ \Big\langle  R_N, R_N \Big\rangle
\label{22.3.12}
 \eeq
Действительно,
\begin{multline*}\Big\langle  f,f \Big\rangle= \Big\langle  S_N+R_N, S_N+R_N
\Big\rangle= \Big\langle  S_N+R_N, S_N \Big\rangle+\Big\langle S_N+R_N, R_N
\Big\rangle=\\=
 \Big\langle  S_N, S_N \Big\rangle+
 \underbrace{\Big\langle  R_N, S_N \Big\rangle+
 \Big\langle  S_N, R_N \Big\rangle}_{\tiny\begin{matrix}\|
 \\ 0 \\ \text{в силу \eqref{22.3.11}}
 \end{matrix}}
 +
 \Big\langle  R_N, R_N \Big\rangle=
 \Big\langle  S_N, S_N \Big\rangle+\Big\langle R_N, R_N
\Big\rangle \end{multline*}

4. Из \eqref{22.3.12} следует:
$$
\Big\langle  f,f \Big\rangle=\eqref{22.3.12}=\Big\langle  S_N, S_N \Big\rangle+
 \kern-20pt
 \underbrace{\Big\langle  R_N, R_N \Big\rangle}_{\scriptsize \begin{matrix}
 \phantom{\tiny \text{свойство $1^0$ }}\ \text{\rotatebox{-90}{$\ge$}}\ {\tiny \text{свойство $1^0$ }}\\ 0 \end{matrix}}
  \kern-20pt\ge
 \Big\langle  S_N, S_N \Big\rangle+0=\eqref{22.3.10}=
\frac{a_0^2}{2}+\sum_{n=1}^N \Big\{ a_n^2+b_n^2 \Big\}
$$
 то есть,
 \beq
\frac{a_0^2}{2}+\sum_{n=1}^N \Big\{ a_n^2+b_n^2 \Big\}\le \Big\langle  f,f
\Big\rangle \label{22.3.13}
 \eeq
Это верно при любом $N\in \mathbb{N}$, поэтому ряд
$\frac{a_0^2}{2}+\sum\limits_{n=1}^\infty \Big\{ a_n^2+b_n^2 \Big\}$ сходится.
При этом, автоматически выполняется \eqref{22.3.9}:
$$
\frac{a_0^2}{2}+\sum_{n=1}^\infty \Big\{ a_n^2+b_n^2 \Big\}= \lim_{N\to
\infty}\l \frac{a_0^2}{2}+\sum_{n=1}^N \Big\{ a_n^2+b_n^2 \Big\}\r \le
\eqref{22.3.13}\le \langle f,f\rangle=\frac{1}{\pi}\int_{-\pi}^{\pi} f(x)^2 \,
\d x
$$
\end{proof}

\paragraph{Равенство Парсеваля.}

\begin{lm}\label{tm-Parseval} Для всякой функции $f\in{\mathcal R}(\pi)$
выполняется следующее равенство, называемое равенством Парсеваля:
 \beq\label{Parseval}
\frac{a_0^2}{2}+\sum_{n=1}^\infty \Big\{ a_n^2+b_n^2 \Big\}= \langle
f,f\rangle=\frac{1}{\pi}\int_{-\pi}^{\pi} f(x)^2 \, \d x
 \eeq
(здесь $a_n$ и $b_n$ -- коэффициенты Фурье \eqref{22.3.6} функции $f$).
\end{lm}

\bpr По лемме \ref{tm-Bessel}, ряд в левой части этой формулы должен сходиться,
причем его сумма, обозначим ее $C$, не превосходит величины $\langle
f,f\rangle$:
$$
\frac{a_0^2}{2}+\sum_{n=1}^\infty \Big\{ a_n^2+b_n^2 \Big\}=C\le\langle
f,f\rangle
$$
Наша задача -- убедиться, что $C\ge\langle f,f\rangle$. Зафиксируем $\e>0$ и по
лемме \ref{LM:polnota-trig-sistemy} подберем тригонометрический многочлен $P$,
для которого выполняется \eqref{polnota-trig-sistemy}:
$$
\norm{f-P}<\e
$$
Пусть $N$ -- степень этого тригонометрического многочлена. Тогда:
$$
\langle f,f\rangle-\frac{a_0^2}{2}-\sum_{n=1}^N \Big\{ a_n^2+b_n^2
\Big\}\overset{\eqref{uklonenie-trig-mnogochl}}{=}\norm{f-S_N}^2
\overset{\eqref{uklonenie-trig-mnogochl}}{\le} \norm{f-P}^2<\e^2
$$
$$
\Downarrow
$$
$$
\langle f,f\rangle<\e^2+\frac{a_0^2}{2}+\sum_{n=1}^N \Big\{ a_n^2+b_n^2
\Big\}\le\e^2+\frac{a_0^2}{2}+\sum_{n=1}^\infty \Big\{ a_n^2+b_n^2
\Big\}=\e^2+C
$$
$$
\Downarrow
$$
$$
\langle f,f\rangle<\e^2+C
$$
Это верно для любого $\e>0$, значит $\langle f,f\rangle\le C$.
 \epr

\paragraph{Доказательство основной теоремы.}

\blm Для всякой функции $f\in{\mathcal R}(\pi)$ выполняется соотношение
\eqref{||f-S_N||_2->0},
$$
\norm{f-S_N}\underset{N\to\infty}{\longrightarrow}0,
$$
в котором $S_N$ -- многочлены Фурье \eqref{mnogochlen-Fourier-pi} функции $f$.
\elm
\bpr
$$
\norm{f-S_N}^2\overset{\eqref{uklonenie-trig-mnogochl}}{=} \langle
f,f\rangle-\bigg(\underbrace{\frac{a_0^2}{2}+\sum_{n=1}^N \Big\{ a_n^2+b_n^2
\Big\} }_{\scriptsize\begin{matrix}\phantom{\tiny \begin{matrix} N\\ \downarrow\\
\infty\end{matrix}}\ \downarrow\ {\tiny \begin{matrix} N\\ \downarrow\\
\infty\end{matrix}}\\ \frac{a_0^2}{2}+\sum\limits_{n=1}^\infty \Big\{
a_n^2+b_n^2 \Big\} \\ \phantom{\eqref{Parseval}}\ \text{\rotatebox{90}{$=$}}\
\eqref{Parseval} \\ \langle f,f\rangle
\end{matrix}}\bigg)\underset{N\to\infty}{\longrightarrow} 0
$$
\epr

\blm Пусть $a_n$ и $b_n$ -- две последовательности чисел и
$$
S_N=\frac{a_0}{2}+\sum_{n=1}^N \Big\{ a_n\cdot\cos_n+b_n\cdot\sin_n \Big\}
$$
-- последовательность частичных сумм соответствующего им тригонометрического
ряда с полупериодом $\pi$. Если функции $S_N$ сходятся в среднем квадратичном к
некоторой локально интегрируемой $2\pi$-периодической функции $f$,
$$
\norm{f-S_N}\underset{N\to\infty}{\longrightarrow}0,
$$
то числа $a_n$ и $b_n$ являются коэффициентами Фурье \eqref{koeff-Fourier-pi}
функции $f$. \elm \bpr Для всякого $k\in\N$
$$
\norm{f-S_N}\underset{N\to\infty}{\longrightarrow}0
\qquad\overset{\text{теорема
\ref{TH:neprer-skal-proizv}}}{\Longrightarrow}\qquad
\kern-50pt\underbrace{\langle
S_N,\cos_k\rangle}_{\scriptsize\begin{matrix}\text{\rotatebox{90}{$=$}}\\
\Big\langle \frac{a_0}{2}+\sum_{n=1}^N \Big\{ a_n\cdot\cos_n+b_n\cdot\sin_n
\Big\}, \cos_k\Big\rangle \\ \text{\rotatebox{90}{$=$}}\\ \phantom{,} a_k, \\
\text{при $N\ge k$}
\end{matrix}}\kern-50pt\underset{N\to\infty}{\longrightarrow}\langle
f,\cos_k\rangle \qquad\Longrightarrow\qquad a_k=\langle f,\cos_k\rangle
$$
и точно также, заменяя $\cos_k$ на $\sin_k$ и на $1$, мы получаем равенства
$$
b_k=\langle f,\sin_k\rangle,\qquad a_0=\langle f, 1\rangle.
$$
\epr

\section{Поточечная и равномерная сходимость ряда Фурье}

Теорема \ref{TH:Fourier-v-R_2} о разложении в ряд Фурье локально интегрируемой
$2\pi$-периодической функции $f$ ничего не говорит о том, будет ли этот ряд
сходиться к функции $f$ поточечно, то есть, можно ли утверждать, что для
всякого $x\in\R$ выполняется соотношение
$$
\frac{a_0}{2}+\sum_{n=1}^N\Big\{ a_n\cos nx+b_n\sin nx
\Big\}\underset{N\to\infty}{\longrightarrow} f(x)
$$
(где $a_n$ и $b_n$ -- коэффициенты Фурье \eqref{koeff-Fourier-pi} функции $f$).
Следующий элементарный пример показывает, что это не всегда верно.

\noindent\rule{160mm}{0.1pt}\begin{multicols}{2}

\bex Рассмотрим функцию
$$
f(x)=\left\{\frac{x}{2\pi}\right\}
$$
(взятие дробной части от $\frac{x}{2\pi}$). Она локально интегрируема и имеет
период $2\pi$. Вычислим ее коэффициенты Фурье:
 \begin{multline*}
a_0=\frac{1}{\pi}\int_{-\pi}^{\pi}\left\{\frac{x}{2\pi}\right\}\ \d
x=\frac{1}{\pi}\int_0^{2\pi}\left\{\frac{x}{2\pi}\right\}\ \d x=\\=
\frac{1}{\pi}\int_0^{2\pi}\frac{x}{2\pi}\ \d
x=\frac{x^2}{4\pi^2}\Big|_{x=0}^{x=2\pi}=1
 \end{multline*}
 \begin{multline*}
a_n=\frac{1}{\pi}\int_0^{2\pi}\left\{\frac{x}{2\pi}\right\}\cdot\cos nx\ \d
x=\\=\frac{1}{\pi}\int_0^{2\pi}\frac{x}{2\pi}\cdot\cos nx\ \d
x=\\=\frac{1}{2\pi^2n}\int_0^{2\pi} x\ \d\sin nx=\\= \frac{1}{2\pi^2n}\cdot
\underbrace{x\cdot\sin nx\Big|_0^{2\pi}}_{\scriptsize\begin{matrix}\text{\rotatebox{90}{$=$}}\\
0\end{matrix}}-\frac{1}{2\pi^2n}\int_0^{2\pi} \sin nx\ \d x=\\=
\frac{1}{2\pi^2n^2}\cdot \cos nx\Big|_0^{2\pi}=0
 \end{multline*}
 \begin{multline*}
b_n=\frac{1}{\pi}\int_0^{2\pi}\left\{\frac{x}{2\pi}\right\}\cdot\sin nx\ \d
x=\\=\frac{1}{\pi}\int_0^{2\pi}\frac{x}{2\pi}\cdot\sin nx\ \d
x=\\=-\frac{1}{2\pi^2n}\int_0^{2\pi n} x\ \d\cos nx=\\= -\frac{1}{2\pi^2n}\cdot
\underbrace{x\cdot\cos nx\Big|_0^{2\pi}}_{\scriptsize\begin{matrix}\text{\rotatebox{90}{$=$}}\\
2\pi\end{matrix}}+\frac{1}{2\pi^2n}\int_0^{2\pi} \cos nx\ \d x=\\=-\frac{1}{\pi
n}+\frac{1}{2\pi^2n^2}\cdot \underbrace{\sin nx\Big|_0^{2\pi}}_{\scriptsize\begin{matrix}\text{\rotatebox{90}{$=$}}\\
0\end{matrix}}=-\frac{1}{\pi n}
 \end{multline*}
Ряд Фурье получается такой:
$$
\frac{1}{2}-\sum_{n=1}^\infty\frac{1}{\pi n}\cdot\sin nx
$$
и в точке $x=0$ он сходится не к значению $f(0)$:
$$
\frac{1}{2}-\underbrace{\sum_{n=1}^N\frac{1}{\pi n}\cdot\sin
0}_{\scriptsize\begin{matrix}\text{\rotatebox{90}{$=$}}\\
0\end{matrix}}\underset{N\to\infty}{\longrightarrow}\frac{1}{2}\ne 0=f(0)
$$
 \eex

\end{multicols}\noindent\rule[10pt]{160mm}{0.1pt}

Уже из этого примера видно, что если мы хотим, чтобы ряд Фурье сходился к
функции $f$ поточечно, то нужно как-то сузить класс изучаемых функций, перейдя
от класса ${\mathcal R}(\pi)$ (локально интегрируемых функций с полупериодом
$\pi$) к каким-то его подклассам. В этом параграфе мы рассмотрим два таких
подкласса, для которых поточечная сходимость ряда Фурье сохраняется, хотя и в
разных смыслах: это классы кусочно-непрерывных и кусочно-гладких функций  с
полупериодом $\pi$. Попутно мы сформулируем для них теоремы о равномерной
сходимости ряда Фурье.

\subsection{Кусочно-непрерывные и кусочно-гладкие функции на $\R$}

 \bit{
\item[$\bullet$] Функция $f$ на прямой $\R$ называется
 \biter{
\item[---] {\it кусочно-непрерывной}\index{функция!кусочно-непрерывная!на
прямой} (на $\R$), если она кусочно непрерывна на любом
отрезке\footnote{Напомним, что определение кусочно-непрерывной функции на
отрезке было дано выше на с. \pageref{DEF:kusochno-nepr-func}.}
$[a,b]\subset\R$; это означает выполнение двух условий:
 \biter{
\item[(1)] на любом отрезке $[a,b]$ она непрерывна во всех точках, кроме,
возможно, конечного набора;

\item[(2)] в каждой точке разрыва $c\in \R$ функция $f$ имеет конечные левый и
правый пределы
 \beq
f(c-0)=\lim_{x\to c-0} f(x), \quad f(c+0)=\lim_{x\to c+0} f(x) \label{22.1.2}
 \eeq
 }\eiter

\item[---] {\it кусочно-гладкой}\index{функция!кусочно-гладкая!на прямой} (на
$\R$), если она кусочно гладкая на любом отрезке\footnote{Определение
кусочно-гладкой функции на отрезке приводилось выше на с.
\pageref{DEF:kusochno-gladkaya-func}.} $[a,b]\subset\R$; это означает
выполнение следующих условий:
 \biter{
\item[(1)] на любом отрезке $[a,b]$ она дифференцируема во всех точках, кроме,
возможно, конечного набора;

\item[(2)] в каждой точке $c\in\R$, где $f$ дифференцируема, ее производная
непрерывна,
 \beq\label{22.1.2-3}
f'(c)=\lim_{x\to c}f'(x)
 \eeq

\item[(3)] в каждой точке $c\in\R$, где $f$ не дифференцируема, она обладает
следующими свойствами:
 \biter{
\item[(i)] $f$ имеет конечные левый и правый пределы
 \beq\label{22.1.2-3}
f(c-0)=\lim_{x\to c-0} f(x), \quad f(c+0)=\lim_{x\to c+0} f(x)
 \eeq

\item[(ii)] $f$ имеет конечные левую и правую производные
 \beq\label{22.1.3-3}
f'_{-}(c)=\lim_{x\to c-0}\frac{f(x)-f(c-0)}{x-c}, \quad f'_{+}(c)=\lim_{x\to
c+0}\frac{f(x)-f(c+0)}{x-c}
 \eeq

\item[(iii)] при приближении аргумента $x$ к $c$ справа или слева, значение
$f'(x)$ стремится соответственно к $f'_{-}(c)$ или $f'_{+}(c)$:
 \beq\label{22.1.3}
f'_{-}(c)=\lim_{x\to c-0} f'(x), \quad f'_{+}(c)=\lim_{x\to c+0}f'(x)
 \eeq
 }\eiter
 }\eiter
 }\eiter
 }\eit
Из теоремы \ref{tm-14.3.5} следует

\btm Любая кусочно-непрерывная функция $f$ на $\R$ локально интегрируема. \etm

\noindent\rule{160mm}{0.1pt}\begin{multicols}{2}

\ber\label{ex-22.1.2} Функция
$$
  f(x)=\sgn \sin x
$$
имеющая график

%\pucture{0pt}{0pt}{ii-23.pcx}

\vglue100pt \noindent является $2\pi$-периодической и кусочно-гладкой.

\eer

\ber Функция
$$
  f(x)=\arcsin \sin x
$$
с графиком

%\pucture{0pt}{0pt}{ii-24.pcx}

\vglue100pt \noindent тоже $2\pi$-периодическая и кусочно-гладкая. \eer

\bex А функция
$$
  f(x)=\sqrt{|\sin x|}
$$
с графиком

%\pucture{0pt}{0pt}{ii-25.pcx}

\vglue100pt \noindent хотя и $2\pi$-периодическая, но не кусочно-гладкая,
потому что в точке $x\in\pi\Z$ не имеет односторонних производных. \eex

\end{multicols}\noindent\rule[10pt]{160mm}{0.1pt}

\bit{ \item[$\bullet$] Говорят, что кусочно-непрерывная функция $f$ на прямой
$\R$ удовлетворяет {\it тождеству Лебега}\index{тождество!Лебега}, если в
каждой точке $c\in \R$ ее значение $f(c)$ равно среднему арифметическому левого
и правого пределов в этой точке:
 \beq\label{22.1.4}
  f(x)=\frac{f(x-0)+f(x+0)}{2}, \qquad x\in \R
 \eeq
Очевидно, это равенство автоматически выполняется в точках непрерывности, так
как в этом случае $f(x-0)=f(x)=f(x+0)$, и значит
$$
\frac{f(x-0)+f(x+0)}{2}=\frac{f(x)+f(x)}{2}=f(x)
$$
Поэтому тождество \eqref{22.1.4} нужно проверять только в точках разрыва
функции $f$.
 }\eit

\noindent\rule{160mm}{0.1pt}\begin{multicols}{2}

\begin{ex}\label{ex-22.1.3} Функция
$$
f(x)=\begin{cases} 0, \quad x<0 \\ \frac{1}{2}, \quad x=0 \\ 1, \quad x>0
\end{cases}
$$
с графиком

%\pucture{0pt}{0pt}{ii-26.pcx}

\vglue100pt \noindent очевидно, удовлетворяет тождеству Лебега.
\end{ex}

\begin{ex}
А наоборот, функция
$$
f(x)=\begin{cases} 0, \quad x<0 \\ 1, \quad x\ge 0 \end{cases}
$$
с графиком

%\pucture{0pt}{0pt}{ii-27.pcx}

\vglue100pt \noindent не удовлетворяет тождеству Лебега.
\end{ex}

\end{multicols}\noindent\rule[10pt]{160mm}{0.1pt}

\subsection{Суммирование ряда Фурье обычным способом}

\begin{tm}[\bf о поточечной и равномерной сходимости ряда Фурье]\label{TH:potoch-Fourier}
Для всякой кусочно-гладкой $2\pi$-периодической функции $f$ на прямой $\R$ ее
ряд Фурье с полупериодом $\pi$ сходится в каждой точке $x\in\R$ к среднему
арифметическому ее левого и правого пределов в этой точке:
 \begin{equation}\label{potoch-Fourier}
\frac{f(x-0)+f(x+0)}{2}=\frac{a_0}{2}+\sum_{n=1}^\infty \Big\{ a_n\cos
nx+b_n\sin nx \Big\},\qquad x\in\R.
 \end{equation}
В частности, если $f$ удовлетворяет тождеству Лебега \eqref{22.1.4}, то она
является поточечной суммой своего ряда Фурье:
 \begin{equation}\label{22.1.5}
f(x)=\frac{a_0}{2}+\sum_{n=1}^\infty \Big\{ a_n\cos nx+b_n\sin nx \Big\},\qquad
x\in\R.
 \end{equation}
Кроме того,
 \bit{
\item[(i)] если $f$ -- непрерывная (и по-прежнему, кусочно-гладкая и
$2\pi$-периодическая) функция на $\R$, то ряд \eqref{22.1.5} сходится к ней на
прямой $\R$ равномерно:\footnote{Норма $\norm{\cdot}_E$ была определена
формулой \eqref{ravnomernaya-norma}.}
 \beq\label{ravnom-shodimost-ryada-Fourier}
\norm{f(x)-\frac{a_0}{2}-\sum_{n=1}^N \Big\{ a_n\cos nx+b_n\sin nx
\Big\}}_{x\in\R}\underset{N\to\infty}{\longrightarrow}0;
 \eeq

 \item[(ii)] если $f$ -- бесконечно гладкая (и по-прежнему,
$2\pi$-периодическая) функция на $\R$, то ряд \eqref{22.1.5} сходится к ней на
прямой $\R$ равномерно по производным:\footnote{Норма $\norm{\cdot}_E^{(m)}$
была определена формулой \eqref{ravnomernaya-norma-po-proizvodnym}.}
 \beq\label{ravnom-po-proizv-shodimost-ryada-Fourier}
\forall m\in\Z_+\qquad \norm{f(x)-\frac{a_0}{2}-\sum_{n=1}^N \Big\{ a_n\cos
nx+b_n\sin nx \Big\}}_{x\in\R}^{(m)}\underset{N\to\infty}{\longrightarrow}0.
 \eeq
 }\eit
 \end{tm}

Доказательство этого факта мы разобьем на несколько лемм, и посвятим ему
оставшуюся часть этого пункта.

\paragraph{Лемма Римана.}

\blm Если функция $W$ интегрируема на отрезке $[a,b]$, то
 \beq\label{Riemann:p->8}
\int_a^b W(x)\cdot\cos p\ x\ \d
x\underset{p\to+\infty}{\longrightarrow}0,\qquad \int_a^b W(x)\cdot \sin p\ x\
\d x\underset{p\to+\infty}{\longrightarrow}0
 \eeq
 \elm
\bpr Эти соотношения доказываются одинаково, поэтому мы докажем только второе
(в действительности, только оно используется при доказательстве теоремы
\ref{TH:potoch-Fourier}).

1. Прежде всего, нужно заметить, что если функция $W$ постоянна, то утверждение
будет очевидно:
 \begin{multline*}
\left|\int_a^b W(x)\cdot\sin p x\ \d x\right|=\bigg|\overbrace{\int_a^b
C\cdot\sin p x\ \d x}^{\scriptsize\begin{matrix} \frac{C}{p}\cdot\Big(\cos p\
a-\cos p\ b\Big)\\
\text{\rotatebox{90}{$=$}}\\
-\frac{C}{p}\cdot\cos px\Big|_{x=a}^{x=b} \\
\text{\rotatebox{90}{$=$}}
\end{matrix}}\bigg|=\frac{|C|}{p}\cdot\Big|\cos p\ a-\cos p\
b\Big|\le\\ \le \frac{|C|}{p}\cdot\Big(|\cos p\ a|+|\cos p\ b|\Big)\le
\frac{|C|}{p}\cdot 2 \underset{p\to+\infty}{\longrightarrow}0
 \end{multline*}

2. Из этого следует, что если функция $W$ кусочно-постоянна,
$$
W(x)=C_i,\qquad x\in(x_{i-1},x_i)
$$
то соотношение также верно:
$$
\int_a^b W(x)\cdot\sin p x\ \d x=\sum_{i=1}^k\int_{x_{i-1}}^{x_i} C_i\cdot\sin
p x\ \d x \underset{p\to+\infty}{\longrightarrow}0
$$

3. Теперь если $W$ -- произвольная интегрируемая функция, то для всякого $\e>0$
можно по теореме \ref{TH:o-priblizh-integrir-func-kusoch-postoyann} выбрать
кусочно-постоянную функцию $g$ на $[a,b]$ так, чтобы выполнялось
$$
\int_a^b |W(x)-g(x)|\ \d x<\frac{\e}{2}
$$
Тогда, по уже доказанному, $\int_a^b g(x)\cdot\cos p x\ \d x
\underset{p\to+\infty}{\longrightarrow}0$, и значит можно выбрать $P$ так,
чтобы
$$
\forall p>P\qquad \left|\int_a^b g(x)\cdot\cos p x\ \d x\right|<\frac{\e}{2},
$$
В результате, для всякого $p>P$ мы получим:
 \begin{multline*}
\left|\int_a^b W(x)\cdot\sin p x\ \d x\right|= \left|\int_a^b
\Big(W(x)-g(x)\Big)\cdot\sin p x\ \d x+\int_a^b g(x)\cdot\sin p x\ \d
x\right|\le\\ \le \underbrace{\left|\int_a^b \Big(W(x)-g(x)\Big)\cdot\sin p x\
\d x\right|}_{\scriptsize\begin{matrix} \IA \\
\int_a^b |W(x)-g(x)|\cdot |\sin px|\ \d x\\ \IA \\
\int_a^b |W(x)-g(x)|\ \d x\\
\text{\rotatebox{90}{$>$}}\\
\frac{\e}{2}
\end{matrix}}+\underbrace{\left|\int_a^b g(x)\cdot\sin p x\ \d
x\right|}_{\scriptsize\begin{matrix}
\text{\rotatebox{90}{$>$}}\\
\frac{\e}{2}\end{matrix}}<\frac{\e}{2}+\frac{\e}{2}=\e
 \end{multline*}
Это доказывает второе соотношение в \eqref{Riemann:p->8}.
 \epr

\paragraph{Ядро Дирихле.}

 \bit{
\item[$\bullet$] {\it Ядром Дирихле}\index{Дирихле!ядро} называется функция
 \beq\label{22.3.16}
D_N(t):=\frac{1}{2}+\sum_{n=1}^N \cos nt
 \eeq
 }\eit

\medskip

\centerline{\bf Свойства ядра Дирихле:}\label{svoistva-yadra-Dirichlet}

 \bit{\it
\item[$1^\circ$.] Ядро Дирихле $D_N$ является четной непрерывной
$2\pi$-периодической функцией на $\R$.
 $$
D_N(-t)=D_N(t),\qquad D_N(t+2\pi)=D_N(t)
 $$

\item[$2^\circ$.] Интеграл по периоду от ядра Дирихле равен $\pi$:
 \beq\label{22.3.17}
\int_{-\pi}^{\pi} D_N(t) \, \d t=2\int_0^{\pi} D_N(t) \, \d t=\pi
 \eeq

\item[$3^\circ$.] Ядро Дирихле удовлетворяет тождеству:
 \beq\label{22.3.18}
D_N(t)=\begin{cases}N+\frac{1}{2}, & t\in2\pi\Z\\
\frac{\sin \l N+\frac{1}{2}\r t}{2\sin \frac{t}{2}}, & t\notin 2\pi\Z
\end{cases}
 \eeq
 }\eit
\begin{proof} Свойство $1^\circ$ вытекает из определения $D_N$ формулой
\eqref{22.3.16}, а формула \eqref{22.3.17} доказывается прямым вычислением.
Тождество \eqref{22.3.18} для точек вида $t=2\pi k$ доказывается вычислением:
$$
D_N(2\pi k)=\frac{1}{2}+\sum_{n=1}^N \cos (2\pi nk)=\frac{1}{2}+\sum_{n=1}^N
1=\frac{1}{2}+N
$$
а для точек $t\ne2\pi k$ -- умножением на $2\sin \frac{t}{2}$:
 \begin{multline*}
\Big(\frac{1}{2}+\sum_{n=1}^N \cos nt \Big)\cdot 2\sin \frac{t}{2}= \sin
\frac{t}{2}+\sum_{n=1}^N 2\cos nt\cdot \sin \frac{t}{2}= \sin
\frac{t}{2}+\sum_{n=1}^N \l -\sin \l n-\frac{1}{2}\r t+\sin \l n+\frac{1}{2}\r
t \r=\\= \underbrace{\sin \frac{t}{2}}
\put(-6,-19){\put(-6.5,-6){\line(1,0){49.3}}\put(-9,-4){$\uparrow$}\put(40.3,-4){$\uparrow$}}
+ \bigg( \underbrace{-\sin \frac{t}{2}}+\underbrace{\sin\frac{3t}{2}}
\put(-7.5,-19){\put(-6.8,-6){\line(1,0){60}}\put(-9,-4){$\uparrow$}\put(51,-4){$\uparrow$}}
\bigg)+ \bigg(\underbrace{-\sin \frac{3t}{2}}+\underbrace{\sin \frac{5t}{2}}
\put(-7.5,-19){\put(-6.8,-6){\line(1,0){33}}\put(-9,-4){$\uparrow$}\put(24,-4){$\uparrow$}}
\bigg)+...\put(7,-19){\put(-6.8,-6){\line(1,0){55}}\put(-9,-4){$\uparrow$}\put(46,-4){$\uparrow$}}
+ \bigg( \underbrace{-\sin \Big( N-\frac{1}{2}\Big) t}+\sin
\Big(N+\frac{1}{2}\Big) t \bigg)=\\= \sin \l N+\frac{1}{2}\r t
 \end{multline*}
 \end{proof}

\paragraph{Интегралы Дирихле.}

\begin{lm}\label{tm-22.3.6} Для всякой $2\pi$-периодической локально интегрируемой функции
$f$ ее многочлен Фурье $S_N$ с полупериодом $\pi$ выражается через ядро Дирихле
$D_N$ формулами
 \begin{equation}\label{22.3.19}
S_N(x)=\frac{1}{\pi}\int_{-\pi}^{\pi} f(x+t)\cdot D_N(t) \, \d t, \qquad x\in
\R
 \end{equation}
 \begin{equation}\label{22.3.20}
S_N(x)= \frac{1}{\pi}\int_0^{\pi}\Big\{ f(x+t)+f(x-t) \Big\}\cdot D_N(t) \, \d
t, \qquad x\in \R
 \end{equation}
\end{lm}

 \bit{
\item[$\bullet$] Интегралы в правых частых формул \eqref{22.3.19} и
\eqref{22.3.20} называются {\it интегралами Дирихле}.
 }\eit

 \begin{proof} Сначала
\eqref{22.3.19}:
\begin{multline*}
S_N(x)= \frac{a_0}{2}+\sum_{k=1}^N \Big\{ a_k\cos kx+b_k\sin kx \Big\}=
\eqref{koeff-Fourier-pi}=\\= \frac{1}{2}\cdot \frac{1}{\pi}\int_{-\pi}^{\pi}
f(y) \, \d y+ \sum_{k=1}^N \lll \cos kx \frac{1}{\pi}\int_{-\pi}^{\pi} f(y)\cos
ky \, \d y +\sin kx \frac{1}{\pi}\int_{-\pi}^{\pi} f(y)\sin ky \, \d y \rrr=\\=
\frac{1}{\pi}\int_{-\pi}^{\pi} f(y) \lll \frac{1}{2} + \sum_{k=1}^N \cos kx
\cdot \cos ky +\sin kx \cdot \sin ky \rrr \, \d y=
\frac{1}{\pi}\int_{-\pi}^{\pi} f(y) \lll \frac{1}{2} + \sum_{k=1}^N \cos k(y-x)
\rrr\, \d y=\\= \frac{1}{\pi}\int_{-\pi}^{\pi} f(y)\cdot D(y-x) \, \d y= \left|
\begin{array}{c}
y-x=t, \, y=x+t, \, \, \d y=\, \d t
\\
y\in [-\pi,\pi] \Leftrightarrow t\in [-\pi-x,\pi-x]
\end{array}\right|= \frac{1}{\pi}\int_{-\pi-x}^{\pi-x} f(x+t)\cdot D(t) \, \d
t=\\= {\smsize\begin{pmatrix}\text{интеграл от периодической}\\
\text{функции по периоду равен}\\
\text{интегралу по ``сдвинутому периоду''}\end{pmatrix}}=
\frac{1}{\pi}\int_{-\pi}^{\pi} f(x+t)\cdot D_N(t) \, \d t
\end{multline*}
После этого \eqref{22.3.19} следует \eqref{22.3.20}:
\begin{multline*}
S_N(x)= \frac{1}{\pi}\int_{-\pi}^{\pi} f(x+t)\cdot D_N(t) \, \d t=
\frac{1}{\pi}\int_{-\pi}^0 f(x+t)\cdot D_N(t) \, \d t+
\frac{1}{\pi}\int_0^{\pi} f(x+t)\cdot D_N(t) \, \d t=\\= {\smsize\begin{pmatrix}\text{в первом интеграле}\\
\text{делаем замену}\\
t=-s
\end{pmatrix}}= -\frac{1}{\pi}\int_{\pi}^0 f(x-s)\cdot D_N(-s) \, \d s+
\frac{1}{\pi}\int_0^{\pi} f(x+t)\cdot D_N(t) \, \d t=\\= {\smsize\begin{pmatrix}\text{вспоминием, что $D_N$}\\
\text{-- четная функция:}\\
D_N(-s)=D_N(s)
\end{pmatrix}}= -\frac{1}{\pi}\int_{\pi}^0 f(x-s)\cdot D_N(s) \, \d s+
\frac{1}{\pi}\int_0^{\pi} f(x+t)\cdot D_N(t) \, \d t=\\= {\smsize\begin{pmatrix}\text{в первом интеграле}\\
\text{переворачиваем}\\
\text{пределы интегрирования}\end{pmatrix}}= \frac{1}{\pi}\int_0^{\pi}
f(x-s)\cdot D_N(s) \, \d s+ \frac{1}{\pi}\int_0^{\pi} f(x+t)\cdot
D_N(t) \, \d t=\\= {\smsize\begin{pmatrix}\text{в первом интеграле}\\
\text{делаем замену}\\
s=t
\end{pmatrix}}= \frac{1}{\pi}\int_0^{\pi} f(x-t)\cdot D_N(t) \, \d t+
\frac{1}{\pi}\int_0^{\pi} f(x+t)\cdot D_N(t) \, \d t=\\=
\frac{1}{\pi}\int_0^{\pi}\Big\{ f(x+t)+f(x-t) \Big\}\cdot D_N(t) \, \d t
\end{multline*}
\end{proof}

\paragraph{Поточечная сходимость ряда Фурье.} Мы можем
теперь доказать ту часть теоремы \ref{TH:potoch-Fourier}, где речь идет о
поточечной сходимости:

 \blm
Для всякой кусочно-гладкой $2\pi$-периодической функции $f$ на прямой $\R$ ее
ряд Фурье с полупериодом $\pi$ сходится в каждой точке $x\in\R$ к среднему
арифметическому ее левого и правого пределов в этой точке:
 $$
\frac{f(x-0)+f(x+0)}{2}=\frac{a_0}{2}+\sum_{n=1}^\infty \Big\{ a_n\cos
nx+b_n\sin nx \Big\},\qquad x\in\R.
 $$
\elm
\begin{proof}
Пусть $f$ -- кусочно-гладкая $2\pi$-периодическая функция. Рассмотрим функции
$$
A(t)= \lll
\begin{array}{c}\frac{f(x+t)-f(x+0)}{t},
\quad t\in (0,\pi]\\
f_{+}'(0), \quad t=0
\end{array}\rrr, \quad B(t)= \lll
\begin{array}{c}\frac{f(x-t)-f(x-0)}{-t},
\quad t\in (0,\pi]\\
f_{-}'(0), \quad t=0
\end{array}\rrr
$$
Они должны быть кусочно-непрерывны на отрезке $[0,\pi]$, потому что в точке
$t=0$ они непрерывны, а при $t\in (0,\pi]$ они получаются умножением
кусочно-непрерывных функций $g(t)=f(x+t)-f(x+0)$ и $h(t)=f(x-t)-f(x-0)$ на
непрерывную функцию $\frac{1}{t}$.

Кроме того, если рассмотреть функцию
$$
C(t)= \lll
\begin{array}{c}\frac{t}{2\sin \frac{t}{2}}\quad t\in (0,\pi]\\
1, \quad t=0
\end{array}\rrr
$$
то она будет непрерывна на отрезке $[0,\pi]$.

Из этого можно сделать вывод, что функция
$$
 W(t)=\Big\{ A(t) - B(t) \Big\}\cdot C(t)
$$
кусочно-непрерывна на отрезке $[0,\pi]$.

Теперь рассмотрим разность между $\frac{f(x+0)+f(x-0)}{2}$ и частичной суммой
$S_N(x)$ ряда Фурье:
\begin{multline*}
S_N(x)-\frac{f(x+0)+f(x-0)}{2}=\\= \underbrace{\frac{1}{\pi}\int_0^{\pi}\Big\{
f(x+t)+f(x-t) \Big\}\cdot D_N(t) \, \d t}_{\scriptsize\begin{matrix}
\phantom{\eqref{22.3.20}}\ \text{\rotatebox{90}{$=$}}\ \eqref{22.3.20}\\
S_N(x)
\end{matrix}}- \frac{f(x+0)+f(x-0)}{2}\cdot
\underbrace{\frac{2}{\pi}\int_0^{\pi} D_N(t) \, \d
t}_{\scriptsize\begin{matrix}
\phantom{\eqref{22.3.17}}\ \text{\rotatebox{90}{$=$}}\ \eqref{22.3.17}\\
1
\end{matrix}}=\\=
\frac{1}{\pi}\int_0^{\pi}\Big\{ f(x+t)+f(x-t) \Big\}\cdot D_N(t) \, \d t-
\frac{1}{\pi}\int_0^{\pi}\Big\{f(x+0)+f(x-0)\Big\}\cdot D_N(t) \, \d t=\\=
\frac{1}{\pi}\int_0^{\pi}\Big\{ f(x+t)+f(x-t) - f(x+0)-f(x-0)\Big\}\cdot D_N(t)
\, \d t=\\= \frac{1}{\pi}\int_0^{\pi}\lll \Big[ f(x+t)-f(x+0)\Big] + \Big[
f(x-t)-f(x-0)\Big] \rrr \cdot D_N(t) \, \d t =\\= \eqref{22.3.18}=
\frac{1}{\pi}\int_0^{\pi}\lll \Big[ f(x+t)-f(x+0)\Big] + \Big[
f(x-t)-f(x-0)\Big] \rrr \cdot \frac{\sin \l N+\frac{1}{2}\r t}{2\sin
\frac{t}{2}}\, \d t=\\= \frac{1}{\pi}\int_0^{\pi}\underbrace{\bigg\{
\underbrace{\frac{f(x+t)-f(x+0)}{t}}_{\scriptsize\begin{matrix}
\text{\rotatebox{90}{$=$}}\\ A(t)
\end{matrix}} - \underbrace{\frac{f(x-t)-f(x-0)}{-t}}_{\scriptsize\begin{matrix}
\text{\rotatebox{90}{$=$}}\\ B(t)
\end{matrix}}\bigg\} \cdot
\underbrace{\frac{t}{2\sin \frac{t}{2}}}_{\scriptsize\begin{matrix}
\text{\rotatebox{90}{$=$}}\\ C(t)
\end{matrix}}}_{\scriptsize\begin{matrix}
\text{\rotatebox{90}{$=$}}\\ W(t)
\end{matrix}}\cdot \sin \l N+\frac{1}{2}\r t \, \d
t=\\=  \frac{1}{\pi}\int_0^{\pi} W(t) \cdot \sin \l N+\frac{1}{2}\r t \, \d t
\underset{N\to\infty}{\longrightarrow} 0
 \end{multline*}
Мы получили нужное соотношение:
$$
S_N(x) \underset{N\to \infty}{\longrightarrow} \frac{f(x+0)+f(x-0)}{2}
$$
\end{proof}

\paragraph{Дифференцирование ряда Фурье.}

Если $f$ -- непрерывная кусочно-гладкая $2\pi$-периодическая функция на $\R$,
то ее производная $f'$ может быть определена не всюду на $\R$ (потому что в
некоторых точках на $\R$ функция $f$ может быть недифференцируема). Тем не
менее, как мы помним, по теореме \ref{TH:interg-ot-v'-dlya-kusoch-glad-v},
будут однозначно определены интегралы вида
$$
\int_{-\pi}^{\pi} g(x)\cdot f'(x)\ \d x=\int_{-\pi}^{\pi} g(x)\ \d f(x)
$$
где $g$ -- произвольная кусочно-непрерывная функция. В частности, поэтому для
$f'$ можно определить коэффициенты ряда Фурье:
 \beq\label{koeff-Fourier-dlya-proizvodnoi}
\alpha_0=\frac{1}{\pi}\int_{-\pi}^{\pi} f'(x)\ \d x,\qquad
\alpha_n=\frac{1}{\pi}\int_{-\pi}^{\pi} f'(x)\cdot\cos nx\ \d x,\qquad
\beta_n=\frac{1}{\pi}\int_{-\pi}^{\pi} f'(x)\cdot\sin nx\ \d x,
 \eeq
Таким образом, $f'$ тоже порождает какой-то ряд Фурье (с коэффициентами
$\alpha_n$, $\beta_n$). Оказывается, что ряд Фурье для $f'$ получается из ряда
Фурье для $f$ почленным дифференцированием. Вот как точно выглядит это
утверждение:

\blm\label{LM:koeff-Fourier-proizvodnoi} Пусть $f$ -- непрерывная
кусочно-гладкая $2\pi$-периодическая функция на $\R$. Тогда коэффициенты Фурье
\eqref{koeff-Fourier-dlya-proizvodnoi} производной $f'$ связаны с
коэффициентами Фурье \eqref{22.1.6} функции $f$ формулами
 \beq\label{koeff-Fourier-proizvodnoi}
\alpha_0=0,\qquad \alpha_n=n\cdot b_n,\qquad \beta_n=-n\cdot a_n, \qquad n\in\N
 \eeq
 \elm
\brem Объясним, почему это интерпретируется, как возможность почленно
дифференцировать ряд Фурье. Обозначим символами $S_N[f]$ и $S_N[f']$
соответственно частичные суммы рядов Фурье функций $f$ и $f'$
$$
S_N[f](x)=\frac{a_0}{2}+\sum_{n=1}^N\Big\{a_n\cos nx+b_n\sin nx\Big\}, \qquad
S_N[f'](x)=\frac{\alpha_0}{2}+\sum_{n=1}^N\Big\{\alpha_n\cos nx+\beta_n\sin
nx\Big\}
$$
Тогда из леммы \ref{LM:koeff-Fourier-proizvodnoi} следует, что $S_N[f']$
получается из $S_N[f]$ дифференцированием:
 \beq\label{S_N[f']=(S_N[f])'}
S_N[f']=\Big(S_N[f]\Big)'
 \eeq
Действительно,
 \begin{multline*}
\Big(S_N[f]\Big)'(x)=\left(\frac{a_0}{2}+\sum_{n=1}^N\Big\{a_n\cos nx+b_n\sin
nx\Big\} \right)'=\sum_{n=1}^N\Big\{-n\cdot a_n\sin nx+n\cdot b_n\cos
nx\Big\}=\\=\eqref{koeff-Fourier-proizvodnoi}=\underbrace{\frac{\alpha_0}{2}}_{\scriptsize\begin{matrix}
\text{\rotatebox{90}{$=$}}\\ 0\end{matrix}}+ \sum_{n=1}^N\Big\{\beta_n\sin
nx+\alpha_n\cos nx\Big\}=\frac{\alpha_0}{2}+\sum_{n=1}^N\Big\{\alpha_n\cos
nx+\beta_n\sin nx\Big\}=S_N[f'](x)
 \end{multline*}
 \erem
 \bpr Для $\alpha_0$ мы используем формулу
 \eqref{newton-leibnitz-dlya-kusochno-gladkih-funktsij-1}:
$$
\alpha_0=\frac{1}{\pi}\int_{-\pi}^{\pi} f'(x)\ \d
x=\frac{1}{\pi}\int_{-\pi}^{\pi} \d
f(x)=\eqref{newton-leibnitz-dlya-kusochno-gladkih-funktsij-1}=\frac{1}{\pi}\cdot
f(x)\Big|_{x=-\pi}^{x=\pi}
$$
А для $\alpha_n$ и $\beta_n$ это доказывается интегрированием по частям (с
применением теоремы \ref{TH:integrir-po-chastyam-dlya-nepr-kus-glad-funk}):
 \begin{multline*}
\alpha_n=\frac{1}{\pi}\int_{-\pi}^{\pi} f'(x)\cdot\cos nx\ \d x=
\frac{1}{\pi}\int_{-\pi}^{\pi} \cos nx\ \d
f(x)=\eqref{integrir-po-chastyam-dlya-nepr-kus-glad-funk}=\\=\underbrace{\frac{1}{\pi}\cdot\cos
nx\cdot
f(x)\Big|_{x=-\pi}^{x=\pi}}_{\scriptsize\begin{matrix}\text{\rotatebox{90}{$=$}}\\
\phantom{,}0,\\ \text{потому что $\cos_n$ и $f$ --}
\\ \text{$2\pi$-периодические}\\ \text{функции}\end{matrix}}-\frac{1}{\pi}\int_{-\pi}^{\pi} f(x)\ \d\cos
nx= n\cdot\frac{1}{\pi}\int_{-\pi}^{\pi} f(x)\cdot \sin nx\ \d x=n\cdot b_n
 \end{multline*}
и
 \begin{multline*}
\beta_n=\frac{1}{\pi}\int_{-\pi}^{\pi} f'(x)\cdot\sin nx\ \d x=
\frac{1}{\pi}\int_{-\pi}^{\pi} \sin nx\ \d
f(x)=\eqref{integrir-po-chastyam-dlya-nepr-kus-glad-funk}=\\=\underbrace{\frac{1}{\pi}\cdot\sin
nx\cdot
f(x)\Big|_{x=-\pi}^{x=\pi}}_{\scriptsize\begin{matrix}\text{\rotatebox{90}{$=$}}\\
\phantom{,}0,\\ \text{потому что $\sin_n$ и $f$ --}
\\ \text{$2\pi$-периодические}\\ \text{функции}\end{matrix}}-\frac{1}{\pi}\int_{-\pi}^{\pi} f(x)\
\d\sin nx=-n\cdot\frac{1}{\pi}\int_{-\pi}^{\pi} f(x)\cdot \cos nx\ \d x=-n\cdot
a_n
 \end{multline*}
\epr

\paragraph{Равномерная сходимость ряда Фурье непрерывной кусочно-гладкой функции.}

Теперь мы можем доказать часть (i) теоремы \ref{tm-22.1.4}, то есть следующее
утверждение:

\blm\label{LM:ravnom-shod-Fourier} Непрерывная кусочно-гладкая
$2\pi$-периодическая функция $f$ является равномерной суммой своего ряда Фурье
\eqref{22.1.5}.
 \elm

 \bpr
Заметим, что, поскольку функция $f'$ интегрируема на $[-\pi,\pi]$, по лемме
\ref{tm-Bessel}, ряд из квадратов ее коэффициентов Фурье должен сходиться:
$$
  \sum_{n=1}^\infty \Big\{ \alpha_n^2+\beta_n^2
  \Big\}<\infty
$$
Отсюда следует такая цепочка:
$$
 \underbrace{
 \begin{matrix}
0\le\l|\beta_n|-\frac{1}{n}\r^2=\beta_n^2-\frac{2|\beta_n|}{n}+\frac{1}{n^2}  \\
\Downarrow  \\
\frac{2|\beta_n|}{n}\le\beta_n^2+\frac{1}{n^2}\\
\Downarrow  \\
|a_n|=\frac{|\beta_n|}{n}\le\frac{1}{2}\cdot\left\{\beta_n^2+\frac{1}{n^2}\right\}
 \end{matrix}
\qquad\qquad
 \begin{matrix}
0\le\l|\alpha_n|-\frac{1}{n}\r^2=\alpha_n^2-\frac{2|\alpha_n|}{n}+\frac{1}{n^2} \\
 \Downarrow \\
\frac{2|\alpha_n|}{n}\le\alpha_n^2+\frac{1}{n^2}\\
\Downarrow  \\
|b_n|=\frac{|\alpha_n|}{n}\le\frac{1}{2}\cdot\left\{\alpha_n^2+\frac{1}{n^2}\right\}
 \end{matrix}
 }
$$
$$
\Downarrow
$$
$$
\sum_{n=1}^\infty \Big\{ |a_n|+|b_n| \Big\}\le
\frac{1}{2}\cdot\sum_{n=1}^\infty \Big\{ \alpha_n^2+\beta_n^2
\Big\}+\sum_{n=1}^\infty\frac{1}{n^2}<\infty
$$
$$
\Downarrow
$$
$$
\sum_{n=1}^\infty \norm{a_n\cos nx+b_n\sin nx}_{x\in\R}\le\sum_{n=1}^\infty
\Big\{\norm{a_n\cos nx}_{x\in\R}+\norm{b_n\sin nx}_{x\in\R}\Big\}=
\sum_{n=1}^\infty \Big\{ |a_n|+|b_n| \Big\}<\infty
$$
$$
\phantom{\text{\scriptsize теорема Вейерштрасса
\ref{tm-20.5.4}}}\quad\Downarrow\quad\text{\scriptsize теорема Вейерштрасса
\ref{tm-20.5.4}}
$$
$$
\text{ряд $\sum_{n=1}^\infty \Big\{ a_n\cos nx+b_n\sin nx \Big\}$ сходится
равномерно на $\R$}
$$
Мы получили, что ряд Фурье функции $f$ сходится равномерно к некоторой функции
$S$. Но с другой стороны, по лемме \ref{LM:potochechn-summa-Fourier} этот ряд
должен поточечно сходиться к функции $f$. Значит, $S=f$, и ряд Фурье этой
функции сходится к ней равномерно.
 \epr

\paragraph{Равномерная по производным сходимость ряда Фурье непрерывной кусочно-гладкой функции.}

Теперь мы можем доказать часть (ii) теоремы \ref{tm-22.1.4}, то есть следующее
утверждение:

 \blm Бесконечно гладкая $2\pi$-периодическая
функция $f$ является равномерной по производным суммой своего ряда Фурье
\eqref{22.1.5}.
 \elm
 \bpr
Для всякого $k\in\Z_+$ рассмотрим производную $f^{(k)}$ порядка $k$ функции
$f$, и пусть $S_N[f^{(k)}]$ обозначает многочлен Фурье функции $f^{(k)}$:
$$
S_N[f^{(k)}](x)=\frac{a_0^{(k)}}{2}+\sum_{n=1}^N\Big\{a_n^{(k)}\cos
nx+b_n^{(k)}\sin nx\Big\},
$$
$$
\alpha_0^{(k)}=\frac{1}{\pi}\int_{-\pi}^{\pi} f^{(k)}(x)\ \d x,\qquad
\alpha_n^{(k)}=\frac{1}{\pi}\int_{-\pi}^{\pi} f^{(k)}(x)\cdot\cos nx\ \d
x,\qquad \beta_n^{(k)}=\frac{1}{\pi}\int_{-\pi}^{\pi} f^{(k)}(x)\cdot\sin nx\
\d x.
$$
Формула \eqref{S_N[f']=(S_N[f])'}, применительно к функции $f^{(k)}$, означает,
что $S_N[f^{(k+1)}]$ получается из $S_N[f^{(k)}]$ дифференцированием:
$$
S_N[f^{(k+1)}]=\Big(S_N[f^{(k)}]\Big)'
$$
Отсюда по индукции можно заключить, что $S_N[f^{(k)}]$ есть просто производная
порядка $k$ из $S_N[f]$ (частичной суммы ряда Фурье исходной функции $f$):
 \beq\label{S_N[f^(k)]=(S_N[f])^(k)}
S_N[f^{(k)}]=\Big(S_N[f]\Big)^{(k)}
 \eeq
Теперь заметим, что каждая функция $f^{(k)}$ является гладкой и
$2\pi$-периодической, поэтому по лемме \ref{LM:ravnom-shod-Fourier} ее ряд
Фурье должен сходиться к ней равномерно на $\R$, то есть
$$
\norm{f^{(k)}-S_N[f^{(k)}]}_{\R}\underset{N\to\infty}{\longrightarrow}0
$$
В силу \eqref{S_N[f^(k)]=(S_N[f])^(k)} это можно интерпретировать, как
равномерное стремление к нулю производной порядка $k$ от разности $f-S_N[f]$:
$$
\norm{\Big(f-S_N[f]\Big)^{(k)}}_{\R}=\norm{f^{(k)}-\Big(S_N[f]\Big)^{(k)}}_{\R}=\eqref{S_N[f^(k)]=(S_N[f])^(k)}=
\norm{f^{(k)}-S_N[f^{(k)}]}_{\R}\underset{N\to\infty}{\longrightarrow}0
$$
Отсюда для всякого $m\in\Z_+$ получаем нужное нам соотношение:
$$
\norm{\phantom{\Big|} f-S_N[f] \phantom{\Big|}}_{\R}^{(k)}
=\eqref{ravnomernaya-norma-po-proizvodnym}= \sum_{k=0}^m
\frac{1}{k!}\cdot\norm{\Big(f-S_N[f]\Big)^{(k)}}_{\R}
\underset{N\to\infty}{\longrightarrow}0
$$
 \epr

\subsection{Суммирование ряда Фурье методом арифметических средних}

Если $2\pi$-периодическая функция $f$ не является кусочно-гладкой, а, скажем,
только кусочно-непре\-рывной и удовлетворяющей тождеству Лебега \eqref{22.1.4},
то из теоремы \ref{TH:potoch-Fourier} не следует, что она должна быть
поточечной суммой своего ряда Фурье. Это не случайно: как оказывается,
существуют такие кусочно-непрерывные, и даже непрерывные, функции, у которых
ряд Фурье не сходится в некоторых точках. Несмотря на это, кусочно-непрерывную
$2\pi$-периодическую функцию $f$, удовлетворяющую тождеству Лебега, всегда
можно считать поточечной суммой своего ряда Фурье, но в несколько неожиданном
смысле: для этого нужно $f$ понимать как предел не многочленов Фурье $S_N$, а
их арифметических средних. Это оригинальное наблюдение принадлежит венгерскому
математику Липоту Фейеру, и в этом пункте мы поговорим о его результатах.

 \bit{

\item[$\bullet$] Пусть  $S_N$ -- частичные суммы ряда Фурье с полупериодом
$\pi$ для функции $f$:
$$
S_N(x)= \frac{a_0}{2}+\sum_{n=1}^N \Big\{ a_n\cos nx+b_n\sin nx \Big\}
$$
(коэффициенты вычисляются по формулам \eqref{koeff-Fourier-pi}). Арифметические
средние этой последовательности, то есть функции
 \beq\label{DEF:Fejer}
\sigma_N(x)=\frac{1}{N}\sum_{n=0}^{N-1}S_n(x)
 \eeq
называются {\it многочленами Фейера} функции $f$.
 }\eit

\btm[\bf Фейер]\label{TH:Fejer} Для всякой кусочно-непрерывной
$2\pi$-периодической функции $f$ на прямой $\R$ ее многочлены Фейера $\sigma_N$
сходятся в каждой точке $x\in\R$ к среднему арифметическому ее левого и правого
пределов в этой точке:
 \begin{equation}\label{potoch-Fejer}
\frac{f(x-0)+f(x+0)}{2}=\lim_{N\to\infty}\sigma_N(x),\qquad x\in\R.
 \end{equation}
В частности, если $f$ удовлетворяет тождеству Лебега \eqref{22.1.4}, то она
является поточечным пределом своей последовательности многочленов Фейера:
 \begin{equation}\label{potoch-Fejer-Lebesgue}
f(x)=\lim_{N\to\infty}\sigma_N(x),\qquad x\in\R.
 \end{equation}
Кроме того,
 \bit{
\item[(i)] если $f$ -- непрерывная (и, по-прежнему, $2\pi$-периодическая)
функция на $\R$, то то многочлены Фейера сходятся к ней на прямой $\R$
равномерно:\footnote{Норма $\norm{\cdot}_E$ была определена формулой
\eqref{ravnomernaya-norma}.}
 \beq\label{Fejer->ravnom}
\norm{f(x)-\sigma_N(x)}_{x\in\R}\underset{N\to\infty}{\longrightarrow}0;
 \eeq

 \item[(ii)] если $f$ -- бесконечно гладкая (и по-прежнему,
$2\pi$-периодическая) функция на $\R$, то многочлены Фейера сходятся к ней на
прямой $\R$ равномерно по производным:\footnote{Норма $\norm{\cdot}_E^{(m)}$
была определена формулой \eqref{ravnomernaya-norma-po-proizvodnym}.}
 \beq\label{ravnom-po-proizv-shodimost-Fejer}
\forall m\in\Z_+\qquad
\norm{f(x)-\sigma_N(x)}_{x\in\R}^{(m)}\underset{N\to\infty}{\longrightarrow}0.
 \eeq
 }\eit
 \etm

Как и в случае с теоремой \ref{TH:potoch-Fourier}, доказательство этого факта
мы разобьем на несколько лемм.

\paragraph{Ядро Фейера.}

 \bit{
\item[$\bullet$] {\it Ядром Фейера}\index{Фейера!ядро} называется
арифметическое среднее ядер Дирихле
 \beq\label{yadro-Fejera}
\varPhi_N(t)=\frac{1}{N}\sum_{n=0}^{N-1}D_n(x)=\frac{1}{N}\sum_{n=0}^{N-1}
\l\frac{1}{2}+\sum_{k=1}^n \cos kt\r
 \eeq
 }\eit

Эта функция обладает свойствами, похожими на свойства ядра Дирихле, но главная
разница заключается в том, что ядро Фейера оказывается неотрицательной функцией
(что оказывается очень полезно для дальнейших выводов).

\bigskip

\centerline{\bf Свойства ядра Фейера:}

 \bit{\it
\item[$1^\circ$.] Ядро Фейера $\varPhi_N$ является непрерывной, четной,
$2\pi$-периодической и неотрицательной функцией на $\R$:
 $$
\varPhi_N(-t)=\varPhi_N(t),\qquad \varPhi_N(t+2\pi)=\varPhi_N(t),\qquad
\varPhi_N(t)\ge 0\qquad \Big(t\in\R\Big)
 $$

\item[$2^\circ$.] Интеграл по периоду от ядра Фейера равен $\pi$:
 \beq\label{int-ot-Fejer}
\int_{-\pi}^{\pi} \varPhi_N(t) \, \d t=2\int_0^{\pi} \varPhi_N(t) \, \d t=\pi
 \eeq

\item[$3^\circ$.] Ядро Фейера удовлетворяет тождеству:
 \beq\label{tozhd-dlya-Fejer}
\varPhi_N(t)=\begin{cases} \frac{N}{2}, & t\in 2\pi\Z
\\
\frac{\sin^2 \frac{N}{2} t}{2N\cdot\sin^2\frac{t}{2}}, & t\notin
2\pi\Z\end{cases}
 \eeq

\item[$4^\circ$.] Для всякого $\delta>0$ равномерная норма функции $\varPhi$ на
отрезке $[\delta,\pi]$ стремится к нулю при $N\to\infty$:
 \beq\label{otsenka-dlya-Fejer}
\norm{\varPhi_N(t)}_{t\in[\delta,\pi]}=\max_{t\in[\delta,\pi]}\varPhi_N(t)
\underset{N\to\infty}{\longrightarrow}0
 \eeq
 }\eit

\bpr Свойства $1^\circ$ и $2^\circ$ вытекают сразу же из свойств ядра Дирихле
на с.\pageref{svoistva-yadra-Dirichlet}, за исключением неотрицательности
$\varPhi_N$, которая, в свою очередь, является следствием тождества
\eqref{tozhd-dlya-Fejer}. Поэтому для доказательства $1^\circ$, $2^\circ$ и
$3^\circ$ достаточно проверить \eqref{tozhd-dlya-Fejer}. Для точек $t=2\pi k$
мы получаем:
$$
\varPhi_N(2\pi k)=\frac{1}{N}\sum_{n=0}^{N-1} \underbrace{D_n(2\pi
k)}_{\scriptsize\begin{matrix}\phantom{\eqref{22.3.18}}\
\text{\rotatebox{90}{$=$}}\ \eqref{22.3.18}\\
n+\frac{1}{2}\end{matrix}}=\frac{1}{N}\sum_{n=0}^{N-1} \l n+\frac{1}{2}\r=
\frac{1}{N}\underbrace{\sum_{n=1}^{N-1}
n}_{\scriptsize\begin{matrix}\phantom{\eqref{arifm-progressija}}\
\text{\rotatebox{90}{$=$}}\ \eqref{arifm-progressija}\\
\frac{(N-1)\cdot
N}{2}\end{matrix}}+\frac{N}{2N}=\frac{N-1}{2}+\frac{1}{2}=\frac{N}{2}
$$
А для $t\ne 2\pi k$ получаем:
 \begin{multline*}
\varPhi_N(t)\cdot 2N\sin^2\frac{t}{2}=\frac{1}{N}\sum_{n=0}^{N-1}
\underbrace{D_n(x)}_{\scriptsize\begin{matrix}\phantom{\eqref{22.3.18}}\
\text{\rotatebox{90}{$=$}}\ \eqref{22.3.18}\\
\frac{\sin \l n+\frac{1}{2}\r t}{2\sin \frac{t}{2}}\end{matrix}}\cdot
2N\sin^2\frac{t}{2}=2\sum_{n=0}^{N-1}\underbrace{\sin \l n+\frac{1}{2}\r t\cdot
\sin\frac{t}{2}}_{\scriptsize\begin{matrix}\phantom{\eqref{sin(x)sin(y)}}\
\text{\rotatebox{90}{$=$}}\ \eqref{sin(x)sin(y)}\\
\frac{1}{2}\Big(\cos nt-\cos(n+1)t\Big)
\end{matrix}}=\\=\sum_{n=0}^{N-1}\Big(\cos nt-\cos(n+1)t\Big)=
\Big(\underbrace{\cos 0}_{\scriptsize\begin{matrix}
\text{\rotatebox{90}{$=$}}\\ 1
\end{matrix}}-\underbrace{\cos t}
\put(-4.5,-11){\put(-6.8,-6){\line(1,0){46.4}}\put(-9,-4){$\uparrow$}\put(37,-4){$\uparrow$}}
\Big)+\Big(\underbrace{\cos t}-\underbrace{\cos 2t}
\put(-6.7,-11){\put(-6.9,-6){\line(1,0){31.5}}\put(-9,-4){$\uparrow$}\put(22,-4){$\uparrow$}}
\Big)+...
\put(7,-13.5){\put(-6.9,-6){\line(1,0){45.5}}\put(-9,-4){$\uparrow$}\put(36,-4){$\uparrow$}}
+\Big(\underbrace{\cos (N-1)t}-\cos Nt\Big)=\\=1-\cos
Nt=\eqref{cos-2x}=\sin^2\frac{N+1}{2}t
 \end{multline*}
Поделив это на $2N\sin^2\frac{t}{2}$, мы получим нижнюю строчку в
\eqref{tozhd-dlya-Fejer}.

После того, как \eqref{tozhd-dlya-Fejer} доказано, свойство
\eqref{otsenka-dlya-Fejer} становится его следствием:
$$
\norm{\varPhi_N(t)}_{t\in[\delta,\pi]}=\max_{t\in[\delta,\pi]}\Big|\varPhi_N(t)\Big|=
\max_{t\in[\delta,\pi]}\left|\frac{\sin^2 \frac{N}{2}
t}{2N\cdot\sin^2\frac{t}{2}}\right|\le \frac{\max_{t\in[\delta,\pi]}\sin^2
\frac{N}{2} t}{2N\cdot
\min_{t\in[\delta,\pi]}\sin^2\frac{t}{2}}=\frac{1}{2N\cdot
\sin^2\frac{\delta}{2}} \underset{N\to\infty}{\longrightarrow}0
$$

 \epr

\paragraph{Интеграл Фейера.}

\begin{lm}\label{LM:int-Fejer} Для всякой $2\pi$-периодической локально интегрируемой функции
$f$ ее многочлен Фейера $\sigma_N$ выражается через ядро Фейера $\varPhi_N$
формулами
 \begin{equation}\label{int-Fejer}
\sigma_N(x)=\frac{1}{\pi}\int_{-\pi}^{\pi} f(x+t)\cdot \varPhi_N(t) \, \d t,
\qquad x\in \R
 \end{equation}
 \begin{equation}\label{int-Fejer-2}
\sigma_N(x)= \frac{1}{\pi}\int_0^{\pi}\Big\{ f(x+t)+f(x-t) \Big\}\cdot
\varPhi_N(t) \, \d t, \qquad x\in \R
 \end{equation}
\end{lm}

 \bit{
\item[$\bullet$] Интегралы справа в формулах \eqref{int-Fejer} и
\eqref{int-Fejer-2} называются {\it интегралами Фейера}.
 }\eit

\bpr Это следует из формул для интеграла Дирихле \eqref{22.3.19} и
\eqref{22.3.20}: во-первых,
 \begin{multline*}
\sigma_N(x)=\frac{1}{N}\sum_{n=0}^{N-1}S_n(x)=
\frac{1}{N}\sum_{n=0}^{N-1}\overbrace{\frac{1}{\pi}\int_{-\pi}^{\pi}
f(x+t)\cdot D_n(t) \, \d t}^{\scriptsize\begin{matrix} S_n(x)
\\
\phantom{\eqref{22.3.19}}\
\text{\rotatebox{90}{$=$}}\ \eqref{22.3.19}
\end{matrix}}=\\=\frac{1}{\pi}\int_{-\pi}^{\pi} f(x+t)\cdot
\underbrace{\frac{1}{N}\sum_{n=0}^{N-1}
D_n(t)}_{\scriptsize\begin{matrix}\phantom{\eqref{yadro-Fejera}}\
\text{\rotatebox{90}{$=$}}\ \eqref{yadro-Fejera}\\
\varPhi_N(t)
\end{matrix}} \, \d t
=\frac{1}{\pi}\int_{-\pi}^{\pi} f(x+t)\cdot \varPhi_N(t) \, \d t,
 \end{multline*}
и, во-вторых,
 \begin{multline*}
\sigma_N(x)=\frac{1}{N}\sum_{n=0}^{N-1}S_n(x)=
\frac{1}{N}\sum_{n=0}^{N-1}\overbrace{\frac{1}{\pi}\int_0^{\pi}\Big\{
f(x+t)+f(x-t) \Big\}\cdot D_n(t) \, \d t}^{\scriptsize\begin{matrix} S_n(x)
\\
\phantom{\eqref{22.3.20}}\ \text{\rotatebox{90}{$=$}}\ \eqref{22.3.20}
\end{matrix}}=\\=\frac{1}{\pi}\int_0^{\pi}\Big\{
f(x+t)+f(x-t) \Big\}\cdot \underbrace{\frac{1}{N}\sum_{n=0}^{N-1}
D_n(t)}_{\scriptsize\begin{matrix}\phantom{\eqref{yadro-Fejera}}\
\text{\rotatebox{90}{$=$}}\ \eqref{yadro-Fejera}\\
\varPhi_N(t)
\end{matrix}} \, \d t
=\frac{1}{\pi}\int_0^{\pi}\Big\{ f(x+t)+f(x-t) \Big\}\cdot \varPhi_N(t) \, \d
t.
 \end{multline*}
 \epr

\paragraph{Поточечная сходимость многочленов Фейера.}
Мы можем теперь доказать ту часть теоремы \ref{TH:Fejer}, где речь идет о
поточечной сходимости:

 \blm
Для всякой кусочно-непрерывной $2\pi$-периодической функции $f$ на прямой $\R$
ее многочлены Фейера сходятся в каждой точке $x\in\R$ к среднему
арифметическому ее левого и правого пределов в этой точке:
 $$
\frac{f(x-0)+f(x+0)}{2}=\lim_{N\to\infty}\sigma_N(x),\qquad x\in\R.
 $$
\elm
\begin{proof} Поскольку $f$ кусочно-непрерывна и периодична, она ограничена на
$\R$. Обозначим буквой $M$ ограничивающую ее константу:
$$
\sup_{x\in\R}|f(x)|=M\qquad (M\in\R).
$$
Зафиксируем точку $x\in\R$. Из кусочной непрерывности $f$ следует также, что
существуют пределы
$$
f(x+0)=\lim_{t\to+0}f(x+t),\qquad f(x-0)=\lim_{t\to+0}f(x-t)
$$
Поэтому если зафиксировать $\e>0$, то найдется $\delta>0$ такое, что
$$
\forall t\in(0,\delta)\qquad |f(x+t)-f(x+0)|<\e\quad\&\quad |f(x-t)-f(x-0)|<\e
$$
Подберем такое $N_0\in\N$, чтобы для всех $N>N_0$ выполнялось неравенство
 \beq\label{max(t-in-delta-pi)varPhi_N(t)<e/2M}
\max_{t\in[\delta,\pi]}\varPhi_N(t)<\frac{\e}{4M}
 \eeq
(это всегда можно сделать, в силу \eqref{otsenka-dlya-Fejer}). Тогда для всех
$N>N_0$ мы получим:
\begin{multline*}
\int_0^{\pi} \Big| f(x+t)-f(x+0)\Big| \cdot \varPhi_N(t) \, \d t=\\=
\int_0^{\delta} \underbrace{\Big|
f(x+t)-f(x+0)\Big|}_{\scriptsize\begin{matrix}\IA \\ \e\end{matrix}} \cdot
\varPhi_N(t) \, \d t+\int_{\delta}^{\pi} \underbrace{\Big|
f(x+t)-f(x+0)\Big|}_{\scriptsize\begin{matrix}\IA \\
|f(x+t)|+|f(x+0)|
\\ \IA \\ 2M
\end{matrix}} \cdot \varPhi_N(t) \, \d t\le\\
\le \e \cdot  \underbrace{\int_0^{\delta} \varPhi_N(t) \, \d
t}_{\scriptsize\begin{matrix}\IA \\ \int\limits_0^{\pi} \varPhi_N(t) \, \d t \\
\phantom{\eqref{int-ot-Fejer}} \ \text{\rotatebox{90}{$=$}}\
\eqref{int-ot-Fejer}
\\ \frac{\pi}{2}
\end{matrix}}+ 2M \cdot
\underbrace{\int_{\delta}^{\pi} \varPhi_N(t) \, \d t}_{\scriptsize\begin{matrix}\IA
\\ (\pi-\delta)\cdot \max\limits_{t\in[\delta,\pi]}\varPhi_N(t) \\
\phantom{\eqref{max(t-in-delta-pi)varPhi_N(t)<e/2M}} \
\text{\rotatebox{90}{$>$}} \ \eqref{max(t-in-delta-pi)varPhi_N(t)<e/2M}
\\ (\pi-\delta)\cdot\frac{\e}{4M}\\ \IA \\ \frac{\pi\cdot\e}{4M}
\end{matrix}}\le \frac{\pi\cdot\e}{2}+\frac{\pi\cdot\e}{2}=\pi\cdot\e
 \end{multline*}
Таким образом,
 \beq\label{dok-Fejer-+}
\forall N>N_0\qquad \frac{1}{\pi}\cdot\int_0^{\pi} \Big| f(x+t)-f(x+0)\Big|
\cdot \varPhi_N(t) \, \d t<\e
 \eeq
Точно так же доказывается, что
 \beq\label{dok-Fejer--}
\forall N>N_0\qquad \frac{1}{\pi}\cdot\int_0^{\pi} \Big| f(x-t)-f(x-0)\Big|
\cdot \varPhi_N(t) \, \d t<\e
 \eeq

Теперь оценим разность между $\frac{f(x+0)+f(x-0)}{2}$ и многочленом Фейера:
\begin{multline*}
\bigg|\sigma_N(x)-\frac{f(x+0)+f(x-0)}{2}\bigg|=\\=
\bigg|\underbrace{\frac{1}{\pi}\int_0^{\pi}\Big\{ f(x+t)+f(x-t) \Big\}\cdot
\varPhi_N(t) \, \d t}_{\scriptsize\begin{matrix}
\phantom{\eqref{int-Fejer-2}}\ \text{\rotatebox{90}{$=$}}\ \eqref{int-Fejer-2}\\
\sigma_N(x)
\end{matrix}}- \frac{f(x+0)+f(x-0)}{2}\cdot
\underbrace{\frac{2}{\pi}\int_0^{\pi} D_N(t) \, \d
t}_{\scriptsize\begin{matrix}
\phantom{\eqref{int-ot-Fejer}}\ \text{\rotatebox{90}{$=$}}\ \eqref{int-ot-Fejer}\\
1
\end{matrix}}\bigg|=\\=
\bigg|\frac{1}{\pi}\int_0^{\pi}\Big\{ f(x+t)+f(x-t) \Big\}\cdot \varPhi_N(t) \,
\d t- \frac{1}{\pi}\int_0^{\pi}\Big\{f(x+0)+f(x-0)\Big\}\cdot \varPhi_N(t) \,
\d t\bigg|=\\= \bigg|\frac{1}{\pi}\int_0^{\pi}\Big\{ f(x+t)+f(x-t) -
f(x+0)-f(x-0)\Big\}\cdot \varPhi_N(t) \, \d t\bigg|=\\=
\frac{1}{\pi}\cdot\bigg|\int_0^{\pi} \Big[ f(x+t)-f(x+0)\Big] \cdot
\varPhi_N(t) \, \d t + \int_0^{\pi}\Big[ f(x-t)-f(x-0)\Big] \cdot \varPhi_N(t)
\, \d
t\bigg|\le\\
\le \underbrace{\frac{1}{\pi}\cdot\int_0^{\pi} \Big| f(x+t)-f(x+0)\Big| \cdot
\varPhi_N(t) \, \d t}_{\scriptsize\begin{matrix} \phantom{\eqref{dok-Fejer-+}}\
\text{\rotatebox{90}{$>$}}\ \eqref{dok-Fejer-+}
\\ \e\end{matrix}} +\underbrace{\frac{1}{\pi}\cdot \int_0^{\pi}\Big| f(x-t)-f(x-0)\Big|
\cdot \varPhi_N(t) \, \d t}_{\scriptsize\begin{matrix}
\phantom{\eqref{dok-Fejer--}}\ \text{\rotatebox{90}{$>$}}\ \eqref{dok-Fejer--}
\\ \e\end{matrix}}<2\e
 \end{multline*}
Мы получили, что для произвольной точки $x\in\R$ и любого $\e>0$ можно
подобрать такое $N_0\in\N$, что для любого $N>N_0$ выполняется соотношение:
 $$
\bigg|\sigma_N(x)-\frac{f(x+0)+f(x-0)}{2}\bigg|<2\e
 $$
Это эквивалентно тому, что нам нужно:
$$
\forall x\in\R\qquad \sigma_N(x) \underset{N\to \infty}{\longrightarrow}
\frac{f(x+0)+f(x-0)}{2}
$$
\end{proof}

\paragraph{Равномерная сходимость многочленов Фейера.}

\blm\label{LM:ravnom-shod-Fejer} Непрерывная $2\pi$-периодическая функция $f$
является равномерной суммой своей последовательности многочленов Фейера
\eqref{DEF:Fejer}.
 \elm
\bpr Пусть $f$ -- непрерывная $2\pi$-периодическая функция на $\R$. Поскольку
$f$ непрерывна на $[-\pi,\pi]$, должна быть конечной величина
 \beq\label{Fejer-M}
M=\sup_{x\in\R}|f(x)|=\sup_{x\in[-\pi,\pi]}|f(x)|
 \eeq
Зафиксируем число $\e>0$ и выберем такое $\delta>0$, чтобы
 \beq\label{Fejer-eps}
\forall x\in[-\pi,\pi]\qquad \forall t\in(-\delta,\delta)\qquad
\Big|f(x)-f(x+t)\Big|<\frac{\e}{3}
 \eeq
(поскольку $f$ непрерывна, и значит, равномерно непрерывна на $[-\pi,\pi]$, это
всегда можно сделать).

После этого найдем $N_0$ такое, что
 \beq\label{Fejer-e/3M}
\forall N>N_0\qquad \norm{\varPhi_N(t)}_{t\in[\delta,\pi]}<\frac{\e}{6M}
 \eeq
Тогда для всякого $N>N_0$ мы получим:
 \begin{multline*}
\Big|f(x)-\sigma_N(x)\Big|=\bigg|f(x)\cdot\underbrace{\frac{1}{\pi}\int_{-\pi}^{\pi}
\varPhi_N(t) \, \d t}_{\scriptsize\begin{matrix}\phantom{\eqref{int-ot-Fejer}}\
\text{\rotatebox{90}{$=$}}\ \eqref{int-ot-Fejer}\\
1
\end{matrix}} -\underbrace{\frac{1}{\pi}\int_{-\pi}^{\pi} f(x+t)\cdot \varPhi_N(t)
\, \d t}_{\scriptsize\begin{matrix}\phantom{\eqref{int-Fejer}}\
\text{\rotatebox{90}{$=$}}\ \eqref{int-Fejer}\\
\sigma_N(x)
\end{matrix}}\bigg|=\\=\bigg|\frac{1}{\pi}\int_{-\pi}^{\pi}\Big(f(x)-f(x+t)\Big)\cdot \varPhi_N(t)
\, \d t\bigg|\le \frac{1}{\pi}\int_{-\pi}^{\pi}\Big|f(x)-f(x+t)\Big|\cdot
\varPhi_N(t) \, \d t=\\=
\frac{1}{\pi}\int_{-\pi}^{-\delta}\underbrace{\Big|f(x)-f(x+t)\Big|}_{\scriptsize\begin{matrix}
\IA \\ |f(x)|+|f(x+t)|\\ \phantom{\eqref{Fejer-M}}\ \IA\ \eqref{Fejer-M}  \\ 2M
\end{matrix}}\cdot
\varPhi_N(t) \, \d
t+\frac{1}{\pi}\int_{-\delta}^{\delta}\underbrace{\Big|f(x)-f(x+t)\Big|}_{\scriptsize\begin{matrix}
\phantom{\eqref{Fejer-eps}}\ \text{\rotatebox{90}{$>$}}\ \eqref{Fejer-eps} \\
\frac{\e}{3}
\end{matrix}}\cdot \varPhi_N(t)
\, \d
t+\\+\frac{1}{\pi}\int_{\delta}^{\pi}\underbrace{\Big|f(x)-f(x+t)\Big|}_{\scriptsize\begin{matrix}
\IA \\ |f(x)|+|f(x+t)|\\ \phantom{\eqref{Fejer-M}}\ \IA\ \eqref{Fejer-M}  \\ 2M
\end{matrix}}\cdot \varPhi_N(t)
\, \d t\le\\ \le \underbrace{\frac{2M}{\pi}\int_{-\pi}^{-\delta} \varPhi_N(t)
\, \d t}
\put(-36,-22){\put(-6.9,-6){\line(1,0){178}}\put(-9,-4){$\uparrow$}\put(168.5,-4){$\uparrow$}}
 +\frac{\e}{3\pi}\int_{-\delta}^{\delta} \varPhi_N(t) \,
\d t+\underbrace{\frac{2M}{\pi}\int_{\delta}^{\pi} \varPhi_N(t) \, \d t}=\\=
\frac{\e}{3\pi}\cdot\underbrace{\int_{-\delta}^{\delta} \varPhi_N(t) \, \d
t}_{\scriptsize\begin{matrix} \IA \\ \int_{-\pi}^{\pi} \varPhi_N(t) \, \d t\\
\phantom{\eqref{int-ot-Fejer}}\ \text{\rotatebox{90}{$=$}}\
\eqref{int-ot-Fejer} \\ \pi
\end{matrix}}+\frac{4M}{\pi}\cdot\kern-10pt\underbrace{\int_{\delta}^{\pi} \varPhi_N(t) \, \d
t}_{\scriptsize\begin{matrix} \IA \\ (\pi-\delta)\cdot\norm{\varPhi_N(t)}_{t\in[\delta,\pi]}\\
\phantom{\eqref{Fejer-e/3M}}\ \text{\rotatebox{90}{$>$}}\ \eqref{Fejer-e/3M}
\\ (\pi-\delta)\cdot\frac{\e}{6M}\\ \text{\rotatebox{90}{$>$}} \\ \frac{\pi\e}{6M}
\end{matrix}}\kern-10pt<\frac{\e}{3\pi}\cdot\pi+\frac{4M}{\pi}\cdot\frac{\pi\e}{6M}=\frac{\e}{3}+\frac{2\e}{3}=\e
 \end{multline*}
Отсюда
$$
\forall N>N_0\qquad
\norm{f(x)-\sigma_N(x)}_{x\in\R}=\sup_{x\in\R}\Big|f(x)-\sigma_N(x)\Big|\le\e
$$
Это доказывает \eqref{Fejer->ravnom}. \epr

\paragraph{Равномерная по производным сходимость многочленов Фейера.}

\blm\label{LM:gladk-shod-Fejer} Бесконечно гладкая $2\pi$-периодическая функция
$f$ является равномерной по производным суммой своей последовательности
многочленов Фейера \eqref{DEF:Fejer}.
 \elm
\bpr Это можно вывести как следствие пункта (ii) теоремы
\ref{TH:potoch-Fourier}: из соотношения
$$
\forall m\in\Z_+\qquad
\norm{f-S_N}_{\R}^{(m)}\underset{N\to\infty}{\longrightarrow}0.
$$
следует, что при фиксированных $m\in\Z_+$ и $\e>0$ существует $N_0$ такое, что
при $N>N_0$ выполняется
$$
\norm{f-S_N}_{\R}^{(m)}<\e
$$
Отсюда следует, что при $N>N_0$ будет выполняться
 \begin{multline*}
\norm{f-\sigma_N}_{\R}^{(m)}=\norm{f-\frac{1}{N}\sum_{n=0}^{N-1}S_n}_{\R}^{(m)}=
\norm{\frac{1}{N}\sum_{n=0}^{N-1}f-\frac{1}{N}\sum_{n=0}^{N-1}S_n}_{\R}^{(m)}=\\=
\norm{\frac{1}{N}\sum_{n=0}^{N-1}(f-S_n)}_{\R}^{(m)} \le
\frac{1}{N}\sum_{n=0}^{N-1}\norm{(f-S_n)}_{\R}^{(m)}<\frac{1}{N}\sum_{n=0}^{N-1}\e=\e
 \end{multline*}
Поскольку $m\in\Z_+$ и $\e>0$ выбирались произвольно, получаем
$$
\forall m\in\Z_+\qquad
\norm{f-\sigma_N}_{\R}^{(m)}\underset{N\to\infty}{\longrightarrow}0.
$$
 \epr

\noindent\rule{160mm}{0.1pt}\begin{multicols}{2}

\section{Примеры}\label{SEC:EX-Fourier}

Нам остается проиллюстрировать результаты этой главы конкретными вычислениями.
Нас будет интересовать поточечная сходимость, поэтому приводимые здесь примеры
будут иллюстрациями к теореме \ref{TH:potoch-Fourier}.

\paragraph{Разложение Фурье произвольной $2\pi$-периодической функции.}

\begin{ex}\label{ex-22.2.1} Найдем разложение Фурье функции
$$
  f(x)=\sgn \sin x=
\begin{cases}
1, \quad x\in (2\pi n , \pi+2\pi n)\\
0, \quad x=\pi n \\
-1, \quad x\in (-\pi + 2\pi n , 2\pi n)
\end{cases}
$$
Для этого нужно просто вычислить коэффициенты Фурье:
 \begin{multline*}
a_0= \frac{1}{\pi}\int\limits_{-\pi}^\pi \sgn \sin x \, \d x=
\frac{1}{\pi}\int\limits_{-\pi}^0 (-1) \, \d x+\\+
\frac{1}{\pi}\int\limits_0^{\pi} 1 \, \d x=-\frac{\pi}{\pi} +\frac{\pi}{\pi} =0
 \end{multline*}
 \begin{multline*}
a_n= \frac{1}{\pi}\int\limits_{-\pi}^\pi \sgn \sin x\cdot \cos nx \, \d x=\\=
\frac{1}{\pi}\int\limits_{-\pi}^0 (-1)\cdot \cos nx \, \d x+
\frac{1}{\pi}\int\limits_0^{\pi} 1\cdot \cos nx \, \d x=\\=
-\frac{1}{\pi}\int\limits_{-\pi}^0 \cos nx \, \d x+
\frac{1}{\pi}\int\limits_0^{\pi}  \cos nx \, \d x=\\= -\frac{1}{\pi n}\sin nx
\Big|_{x=-\pi}^{x=0}+ \frac{1}{\pi n}\sin nx \Big|_{x=0}^{x=\pi}=0+0=0
 \end{multline*}
 \begin{multline*}
b_n= \frac{1}{\pi}\int\limits_{-\pi}^\pi \sgn \sin x\cdot \sin nx \, \d x=\\=
\frac{1}{\pi}\int\limits_{-\pi}^0 (-1)\cdot \sin nx \, \d x+
\frac{1}{\pi}\int\limits_0^{\pi} 1\cdot \sin nx \, \d x=\\=
-\frac{1}{\pi}\int\limits_{-\pi}^0 \sin nx \, \d x+
\frac{1}{\pi}\int\limits_0^{\pi}  \sin nx \, \d x=\\= \frac{1}{\pi n}\cos nx
\Big|_{x=-\pi}^{x=0} -\frac{1}{\pi n}\cos nx \Big|_{x=0}^{x=\pi}=\\=
\frac{1}{\pi n}\{ 1- \cos  (-\pi n) \} -\frac{1}{\pi n}\{ \cos (\pi n) -1
\}=\\= \frac{2}{\pi n}\{ 1- \cos (\pi n) \}= \frac{2}{\pi n}\{ 1- (-1)^n \}
 \end{multline*}
Вывод:
$$
\sgn \sin x=\sum_{n=1}^\infty  \frac{2}{\pi n}\{ 1- (-1)^n \}\cdot \sin nx
$$
\end{ex}

\begin{ex}\label{ex-22.2.2} Вычислим разложение $2\pi$-периодической
функции $f$, удовлетворяющей тождеству Лебега и заданной на интервале
$(-\pi,\pi)$ формулой
$$
f(x)=\begin{cases}
1, \quad x\in (0 , \pi)\\
\frac{1}{2}, \quad x=0 \\
0, \quad x\in (-\pi, 0)
\end{cases}
$$

%\pucture{0pt}{0pt}{ii-28.pcx}

\vglue100pt \noindent Вычисляем коэффициенты Фурье:
 \begin{multline*}
a_0= \frac{1}{\pi}\int\limits_{-\pi}^\pi f(x) \, \d x=\\=
\frac{1}{\pi}\int\limits_{-\pi}^0 0 \, \d x+ \frac{1}{\pi}\int\limits_0^{\pi} 1
\, \d x=\frac{\pi}{\pi}=1
 \end{multline*}
 \begin{multline*}
a_n= \frac{1}{\pi}\int\limits_{-\pi}^\pi f(x)\cdot \cos nx \, \d x=\\=
\frac{1}{\pi}\int\limits_{-\pi}^0 0\cdot \cos nx \, \d x+
\frac{1}{\pi}\int\limits_0^{\pi} 1\cdot \cos nx \, \d x=\\= \frac{1}{\pi n}\sin
nx \Big|_{x=0}^{x=\pi}=0
 \end{multline*}
 \begin{multline*}
b_n= \frac{1}{\pi}\int\limits_{-\pi}^\pi f(x)\cdot \sin nx \, \d x=\\=
\frac{1}{\pi}\int\limits_{-\pi}^0 0\cdot \sin nx \, \d x+
\frac{1}{\pi}\int\limits_0^{\pi} 1\cdot \sin nx \, \d x=\\= -\frac{1}{\pi
n}\cos nx \Big|_{x=0}^{x=\pi}= -\frac{1}{\pi n}\{ \cos (\pi n) -1 \}=\\=
\frac{1}{\pi n}\{ 1- (-1)^n \}
 \end{multline*}
Вывод:
$$
f(x)=\frac{1}{2}+ \sum_{n=1}^\infty  \frac{1}{\pi n}\{ 1- (-1)^n \}\cdot \sin
nx
$$
\end{ex}

\begin{ex}\label{ex-22.2.3} Вычислим разложение $2\pi$-перио\-дичес\-кой
функции $f$, удовлетворяющей тождеству Лебега и заданной на интервале
$(-\pi,\pi)$ формулой
$$
f(x)=\begin{cases}
1, \quad x\in \left( 0 , \frac{\pi}{2}\right) \\
\frac{1}{2}, \quad x\in \lll 0, \frac{\pi}{2}\rrr \\
0, \quad x\in (-\pi, 0)\cup \left(\frac{\pi}{2},\pi \right)
\end{cases}
$$

%\pucture{0pt}{0pt}{ii-29.pcx}

\vglue100pt \noindent Вычисляем коэффициенты Фурье:
 \begin{multline*}
a_0= \frac{1}{\pi}\int\limits_{-\pi}^\pi f(x) \, \d x=\\=
\frac{1}{\pi}\int\limits_{-\pi}^0 0 \, \d x+
\frac{1}{\pi}\int\limits_{\pi}^{\frac{\pi}{2}} 1 \, \d x+
\frac{1}{\pi}\int\limits_{\frac{\pi}{2}}^{\pi} 0 \, \d x=
\frac{\pi}{\pi}=\frac{1}{2}
 \end{multline*}
 \begin{multline*}
a_n= \frac{1}{\pi}\int\limits_{-\pi}^\pi f(x)\cdot \cos nx \, \d x=\\=
\frac{1}{\pi}\int\limits_{\frac{\pi}{2}}^{\pi} 1\cdot \cos nx \, \d x=
\frac{1}{\pi n}\sin nx \Big|_{x=\frac{\pi}{2}}^{x=\pi}= \frac{\sin \frac{\pi
n}{2}}{\pi n}
 \end{multline*}
 \begin{multline*}
b_n= \frac{1}{\pi}\int\limits_{-\pi}^\pi f(x)\cdot \sin nx \, \d x=\\=
\frac{1}{\pi}\int\limits_{\frac{\pi}{2}}^{\pi} 1\cdot \sin nx \, \d x=
-\frac{1}{\pi n}\cos nx \Big|_{x=\frac{\pi}{2}}^{x=\pi}=\\= -\frac{1}{\pi n}\{
\cos  (\pi n) - \cos \frac{\pi}{2}\}= \frac{1}{\pi n}\{ \cos \frac{\pi}{2}-
(-1)^n \}
 \end{multline*}
Вывод:
 \begin{multline*}
f(x)=\frac{1}{4}+ \sum_{n=1}^\infty  \Big\{ \frac{\sin \frac{\pi n}{2}}{\pi
n}\cdot \cos nx+\\+ \frac{1}{\pi n}\{ \cos \frac{\pi}{2}- (-1)^n \}\cdot \sin
nx \Big\}
 \end{multline*}
\end{ex}

\begin{ers} Разложите в ряд Фурье
$2\pi$-периодические функции, заданные на интервале $(-\pi,\pi)$ формулами:

1) $f(x)=e^x$

2) $f(x)=\begin{cases} x, \quad x\in \left( 0 , \pi \right) \\
                0, \quad x\in \left( -\pi, 0 \right) \end{cases}$

3) $f(x)=\begin{cases} 0, \quad x\in \left( 0 , \pi \right) \\
                x, \quad x\in \left( -\pi, 0 \right) \end{cases}$

4) $f(x)=\begin{cases} x, \quad x\in \left( 0 , \pi \right) \\
                \frac{\pi}{2}, \quad x=0 \\
                \pi, \quad x\in \left( -\pi, 0 \right) \end{cases}$
\end{ers}

\paragraph{Разложение Фурье четных и нечетных $2\pi$-периодических функций.}

\bprop\label{tm-22.2.5} Если $2\pi$-периодическая функция $f$ четна
 \beq\label{22.2.2}
  f(-x)=f(x)
 \eeq
то в ее ряде Фурье отсутствуют слагаемые с синусами:
 \beq\label{22.2.3}
f(x)=\frac{a_0}{2}+\sum_{n=1}^\infty a_n\cos nx
 \eeq
а коэффициенты $a_0, \,  a_n$ можно вычислять по формулам
 \beq\label{22.2.4}
 \begin{split}
& a_0=\frac{2}{\pi}\int\limits_0^\pi f(x) \, \d x, \\
& a_n=\frac{2}{\pi}\int\limits_0^\pi f(x) \cos nx \, \d x
 \end{split}
  \eeq
 \eprop
\begin{proof} Коэффициенты $b_n$ обнуляются:
 \begin{multline*}
b_n=\frac{1}{\pi}\int\limits_{-\pi}^\pi f(x) \sin nx \, \d x=\\=
\underbrace{\frac{1}{\pi}\int\limits_{-\pi}^0 f(x) \sin nx \, \d
x}_{\text{замена: $x=-t$}}+ \frac{1}{\pi}\int\limits_0^{\pi} f(x) \sin nx \, \d
x=\\= \frac{1}{\pi}\int\limits_{\pi}^0
\kern-17pt\underbrace{f(-t)}_{\scriptsize\begin{matrix}\phantom{\eqref{22.2.2}}\
\text{\rotatebox{90}{$=$}}\ \eqref{22.2.2}\\
f(t)\end{matrix}}\kern-17pt\sin (-nt) \, \d (-t)+
\frac{1}{\pi}\int\limits_0^{\pi} f(x) \sin nx \, \d x=\\=
\underbrace{\frac{1}{\pi}\int\limits_{\pi}^0 f(t) \sin nt \, \d
t}_{\scriptsize\begin{matrix} \text{\rotatebox{90}{$=$}}\\
-\frac{1}{\pi}\int\limits_0^{\pi} f(t) \sin nt \, \d t\end{matrix}}+
\frac{1}{\pi}\int\limits_0^{\pi} f(x) \sin nx \, \d x=0
 \end{multline*}
Докажем формулы \eqref{22.2.4}: во-первых,
 \begin{multline*}
a_0=\frac{1}{\pi}\int\limits_{-\pi}^\pi f(x) \, \d x=\\=
\underbrace{\frac{1}{\pi}\int\limits_{-\pi}^0 f(x) \, \d x}_{\text{замена:
$x=-t$}}+ \frac{1}{\pi}\int\limits_0^{\pi} f(x) \, \d x=\\=
\frac{1}{\pi}\int\limits_{\pi}^0
\kern-17pt\underbrace{f(-t)}_{\scriptsize\begin{matrix}\phantom{\eqref{22.2.2}}\
\text{\rotatebox{90}{$=$}}\ \eqref{22.2.2}\\
f(t)\end{matrix}}\kern-17pt \, \d (-t)+ \frac{1}{\pi}\int\limits_0^{\pi} f(x)
\, \d x=\\= \underbrace{-\frac{1}{\pi}\int\limits_{\pi}^0 f(t) \, \d t}_{\scriptsize\begin{matrix} \text{\rotatebox{90}{$=$}}\\
\frac{1}{\pi}\int\limits_0^{\pi} f(t) \, \d t\end{matrix}}+
\frac{1}{\pi}\int\limits_0^{\pi} f(x)  \, \d x=\frac{2}{\pi}\int\limits_0^{\pi}
f(x) \, \d x
\end{multline*} и, во-вторых,
\begin{multline*}
a_n=\frac{1}{\pi}\int\limits_{-\pi}^\pi f(x) \cos nx \, \d x=\\=
\underbrace{\frac{1}{\pi}\int\limits_{-\pi}^0 f(x) \cos nx \, \d
x}_{\text{замена: $x=-t$}}+ \frac{1}{\pi}\int\limits_0^{\pi} f(x) \cos nx \, \d
x=\\= \frac{1}{\pi}\int\limits_{\pi}^0
\kern-17pt\underbrace{f(-t)}_{\scriptsize\begin{matrix}\phantom{\eqref{22.2.2}}\
\text{\rotatebox{90}{$=$}}\ \eqref{22.2.2}\\
f(t)\end{matrix}}\kern-17pt \cos (-nt) \, \d (-t)+\\+
\frac{1}{\pi}\int\limits_0^{\pi} f(x) \cos nx \, \d x= \\=
\underbrace{-\frac{1}{\pi}\int\limits_{\pi}^0 f(t) \cos nt \, \d t}_{\scriptsize\begin{matrix} \text{\rotatebox{90}{$=$}}\\
\frac{1}{\pi}\int\limits_0^{\pi} f(t)\cos nt \, \d t\end{matrix}}+
\frac{1}{\pi}\int\limits_0^{\pi} f(x) \cos nx \, \d x=\\=
\frac{2}{\pi}\int\limits_0^{\pi} f(x) \cos nx \, \d x
\end{multline*}\end{proof}

\bprop\label{tm-22.2.6} Если $2\pi$-периодическая функция $f$ нечетна
 \beq\label{22.2.5}
  f(-x)=-f(x)
 \eeq
то в ее ряде Фурье отсутствуют свободное слагаемое и слагаемые с косинусами:
 \beq\label{22.2.6}
f(x)=\sum_{n=1}^\infty b_n\sin nx
 \eeq
а коэффициенты $b_n$ можно вычислять по формулам
 \beq \label{22.2.7}
b_n=\frac{2}{\pi}\int\limits_0^\pi f(x) \sin nx \, \d x
 \eeq
 \eprop
\begin{proof} Для коэффициента $a_0$ получим:
 \begin{multline*}
a_0=\frac{1}{\pi}\int\limits_{-\pi}^\pi f(x) \, \d x=\\=
\underbrace{\frac{1}{\pi}\int\limits_{-\pi}^0 f(x)  \, \d x}_{\text{замена:
$x=-t$}}+ \frac{1}{\pi}\int\limits_0^{\pi} f(x)  \, \d x=\\=
\frac{1}{\pi}\int\limits_{\pi}^0
\kern-17pt\underbrace{f(-t)}_{\scriptsize\begin{matrix}\phantom{\eqref{22.2.5}}\
\text{\rotatebox{90}{$=$}}\ \eqref{22.2.5}\\
-f(t)\end{matrix}}\kern-17pt \, \d (-t)+ \frac{1}{\pi}\int\limits_0^{\pi} f(x)
\, \d x=\\= \underbrace{\frac{1}{\pi}\int\limits_{\pi}^0 f(t) \, \d t}_{\scriptsize\begin{matrix} \text{\rotatebox{90}{$=$}}\\
-\frac{1}{\pi}\int\limits_0^{\pi} f(t) \, \d t\end{matrix}}+
\frac{1}{\pi}\int\limits_0^{\pi} f(x) \, \d x=0
 \end{multline*}
а для коэффициентов $a_n$,
 \begin{multline*}
a_n=\frac{1}{\pi}\int\limits_{-\pi}^\pi f(x) \cos nx \, \d x=\\=
\underbrace{\frac{1}{\pi}\int\limits_{-\pi}^0 f(x) \cos nx \, \d
x}_{\text{замена: $x=-t$}}+ \frac{1}{\pi}\int\limits_0^{\pi} f(x) \cos nx \, \d
x=\\= \frac{1}{\pi}\int\limits_{\pi}^0
\kern-17pt\underbrace{f(-t)}_{\scriptsize\begin{matrix}\phantom{\eqref{22.2.5}}\
\text{\rotatebox{90}{$=$}}\ \eqref{22.2.5}\\
-f(t)\end{matrix}}\kern-17pt \cos (-nt) \, \d (-t)+\\+
\frac{1}{\pi}\int\limits_0^{\pi} f(x) \cos nx \, \d x=\\=
\underbrace{\frac{1}{\pi}\int\limits_{\pi}^0 f(t) \cos nt \, \d t}_{\scriptsize\begin{matrix} \text{\rotatebox{90}{$=$}}\\
-\frac{1}{\pi}\int\limits_0^{\pi} f(t)\cos nt \, \d t\end{matrix}}+
\frac{1}{\pi}\int\limits_0^{\pi} f(x) \cos nx \, \d x=0
 \end{multline*}
Докажем формулу \eqref{22.2.7}:
 \begin{multline*}
b_n=\frac{1}{\pi}\int\limits_{-\pi}^\pi f(x) \sin nx \, \d x=\\=
\underbrace{\frac{1}{\pi}\int\limits_{-\pi}^0 f(x) \sin nx \, \d
x}_{\text{замена: $x=-t$}}+ \frac{1}{\pi}\int\limits_0^{\pi} f(x) \sin nx \, \d
x=\\= \frac{1}{\pi}\int\limits_{\pi}^0
\kern-17pt\underbrace{f(-t)}_{\scriptsize\begin{matrix}\phantom{\eqref{22.2.5}}\
\text{\rotatebox{90}{$=$}}\ \eqref{22.2.5}\\
-f(t)\end{matrix}}\kern-17pt \sin (-nt) \, \d (-t)+\\+
\frac{1}{\pi}\int\limits_0^{\pi} f(x) \sin nx \, \d x= \\=
\underbrace{-\frac{1}{\pi}\int\limits_{\pi}^0 f(t) \sin nt \, \d t}_{\scriptsize\begin{matrix} \text{\rotatebox{90}{$=$}}\\
\frac{1}{\pi}\int\limits_0^{\pi} f(t)\sin nt \, \d t\end{matrix}}+
\frac{1}{\pi}\int\limits_0^{\pi} f(x) \sin nx \, \d x=\\=
\frac{2}{\pi}\int\limits_0^{\pi} f(x) \sin nx \, \d x
\end{multline*}\end{proof}

\begin{ex}\label{ex-22.2.7} Найдем разложение четной $2\pi$-периодической
функции $f$, удовлетворяющей тождеству Лебега и заданной на интервале $(0,\pi)$
формулой
$$
f(x)=x
$$
Ясно, что график имеет вид

%\pucture{0pt}{0pt}{ii-30.pcx}

\vglue100pt \noindent Поскольку $f$ четна, нужно вычислить только $a_0$ и
$a_n$:
$$
a_0=\frac{2}{\pi}\int\limits_0^\pi x \, \d x= \frac{2}{\pi}\cdot
\frac{x^2}{2}\Big|_0^\pi=\pi
$$
 \begin{multline*}
a_n= \frac{2}{\pi}\int\limits_0^\pi x\cdot \cos nx \, \d x= \frac{2}{\pi
n}\int\limits_0^\pi x \, \d \sin nx =\\= \frac{2}{\pi n}\lll x\cdot \sin nx
\Big|_0^\pi- \int\limits_0^\pi \sin nx \, \d x \rrr=\\= \frac{2}{\pi n}\lll
0+\frac{\cos nx}{n}\Big|_0^\pi \rrr= \frac{2}{\pi n^2}\lll \cos \pi n -1
\rrr=\\= \frac{2}{\pi n^2}\lll (-1)^n - 1 \rrr
 \end{multline*}
Вывод:
$$
f(x)=\frac{\pi}{2}+ \sum_{n=1}^\infty \frac{2}{\pi n^2}\lll (-1)^n - 1 \rrr
\cdot \cos nx
$$
\end{ex}

\begin{ex}\label{ex-22.2.8} Найдем разложение нечетной $2\pi$-периодической
функции $f$, удовлетворяющей тождеству Лебега и заданной на интервале $(0,\pi)$
формулой
$$
f(x)=x
$$
График здесь имеет вид

%\pucture{0pt}{0pt}{ii-31.pcx}

\vglue100pt \noindent Поскольку $f$ нечетна, нужно вычислить только $b_n$:
 \begin{multline*}
b_n= \frac{2}{\pi}\int\limits_0^\pi x\cdot \sin nx \, \d x= -\frac{2}{\pi
n}\int\limits_0^\pi x \, \d \cos nx =\\= -\frac{2}{\pi n}\lll x\cdot \cos nx
\Big|_0^\pi- \int\limits_0^\pi \cos nx \, \d x \rrr=\\= -\frac{2}{\pi n}\lll
\pi\cdot \cos \pi n- \frac{\sin nx}{n}\Big|_0^\pi \rrr=\\= -\frac{2}{\pi n}\lll
\pi\cdot (-1^n)-0 \rrr= -\frac{2}{n}\cdot (-1)^n
 \end{multline*}
Вывод:
$$
f(x)= -\sum_{n=1}^\infty \frac{2(-1)^n}{n}\cdot \sin nx
$$
\end{ex}

\bex\label{EX:Frourier-dlya-x^2} Разложим в ряд Фурье четную
$2\pi$-периодическую функцию, заданную на интервале $(-\pi,\pi)$ формулой
$$
f(x)=x^2
$$
График имеет вид

%\pucture{0pt}{0pt}{ii-30.pcx}

\vglue100pt \noindent Поскольку $f$ четна, нужно вычислить только $a_0$ и
$a_n$:
$$
a_0=\frac{2}{\pi}\int\limits_0^\pi x^2 \, \d x= \frac{2}{\pi}\cdot
\frac{x^3}{3}\Big|_0^\pi=\frac{2\pi^2}{3}
$$
 \begin{multline*}
a_n= \frac{2}{\pi}\int_0^\pi x^2\cdot \cos nx \, \d x= \frac{2}{\pi
n}\int\limits_0^\pi x^2 \, \d \sin nx =\\= \frac{2}{\pi n}\cdot\Big\{
\underbrace{x^2\cdot\sin
nx\Big|_0^\pi}_{\scriptsize\begin{matrix}\text{\rotatebox{90}{$=$}}\\
0\end{matrix}}- 2\int_0^\pi x\cdot\sin nx \, \d x \Big\}=\\= \frac{4}{\pi
n^2}\cdot\int_0^\pi x\ \d \cos nx=\\= \frac{4}{\pi n^2}\cdot\Big\{ x\cdot\cos
nx\Big|_0^\pi- \int_0^\pi \cos nx\ \d x\Big\}=\\=\frac{4}{\pi n^2}\cdot \Big\{
\pi\cdot(-1)^n - \underbrace{\frac{1}{n}\cdot\sin nx\Big|_0^{\pi}}_{\scriptsize\begin{matrix}\text{\rotatebox{90}{$=$}}\\
0\end{matrix}}\Big\}=4\cdot\frac{(-1)^n}{n^2}
 \end{multline*}
Вывод:
 \beq\label{Fourier=>sum-1/n^2}
f(x)=\frac{\pi^2}{3}+ 4\cdot\sum_{n=1}^\infty \frac{(-1)^n}{n^2}\cdot \cos nx
 \eeq
 \eex

\brem\label{REM:sum-1/n^2} Из тождества \eqref{Fourier=>sum-1/n^2}, между
прочим, следует формула \eqref{sum-1/n^2=pi^2/6}, которую мы обещали доказать в
начале главы \ref{CH-number-series}: если положить $x=\pi$, то мы получим
равенство
$$
\pi^2=\frac{\pi^2}{3}+ 4\cdot\sum_{n=1}^\infty \frac{(-1)^n}{n^2}\cdot(-1)^n,
$$
которое после упрощений как раз дает \eqref{sum-1/n^2=pi^2/6}
$$
\frac{\pi^2}{6}=\sum_{n=1}^\infty \frac{1}{n^2}.
$$
Разумеется, это наблюдение не дает никакого способа нахождения сумм рядов
(например, рядов Дирихле из примера \ref{ex-18.3.3}), на который читатель мог
надеяться после данного нами обещания доказать формулу
\eqref{sum-1/n^2=pi^2/6}. Такого алгоритма математики не нашли, однако в
утешение читателю мы скажем, что в \ref{SEC:asymp-integr} следующей главы мы
опишем алгоритмы в некотором смысле приближенного нахождения сумм рядов (а
также точного нахождения конечных сумм) для некоторых классов рядов.
 \erem

\begin{ers} Разложите в ряд Фурье
$2\pi$-периодические функции, заданные на интервале $(0,\pi)$, рассмотрев
отдельно случай когда $f$ четна и нечетна:

1) $f(x)=\begin{cases} 1, \quad x\in \left( 0 , \frac{\pi}{2}\right) \\
                \frac{1}{2}, \quad x=\frac{\pi}{2}\\
                0 \quad x\in \left(\frac{\pi}{2}, \pi \right) \end{cases}$

2) $f(x)=\begin{cases} x, \quad x\in \left( 0 , \frac{\pi}{2}\right) \\
\frac{\pi}{2}, \quad x\in \left(\frac{\pi}{2}, \pi \right)
\end{cases}$

3) $f(x)=
\begin{cases}
x, \quad x\in \left( 0 , \frac{\pi}{2}\right) \\
\frac{\pi}{4}, x=\frac{\pi}{2}\\
0,\quad x\in \left(\frac{\pi}{2}, \pi \right)
\end{cases}$
\end{ers}

\paragraph{Ряды Фурье функций с произвольным периодом.}

\begin{tm}\label{tm-22.2.10}
Для всякой кусочно-гладкой функции $f$ с полупериодом $T$ на прямой $\R$, ее
ряд Фурье (с коэффициентами \eqref{koeff-Fourier}) сходится в каждой точке к
среднему арифметическому ее левого и правого пределов
 \begin{multline}
\frac{f(x-0)+f(x+0)}{2}=\\=\frac{a_0}{2}+\sum_{n=1}^\infty \lll a_n\cos
\frac{\pi nx}{T}+ b_n\sin \frac{\pi nx}{T}\rrr \label{22.2.8}
 \end{multline}
При этом
 \biter{
\item[(i)] если функция $f$ четная, то
 \begin{align*}
& \frac{f(x-0)+f(x+0)}{2}=\frac{a_0}{2}+\sum_{n=1}^\infty a_n\cos \frac{\pi nx}{T}, \\
& a_0=\frac{2}{T}\int\limits_0^T f(x) \, \d x, \\
& a_n=\frac{2}{T}\int\limits_0^T f(x) \cos \frac{\pi nx}{T}\, \d x
 \end{align*}
\item[(ii)] если $f$ нечетная, то
 \begin{align*}
& \frac{f(x-0)+f(x+0)}{2}=\sum_{n=1}^\infty b_n\sin \frac{\pi nx}{T}, \\
& b_n=\frac{2}{T}\int\limits_0^T f(x) \sin \frac{\pi nx}{T}\, \d x
 \end{align*}
 }\eiter
\end{tm}\begin{proof} Заменой переменной
$$
x=\frac{Ty}{\pi}
$$
функция $f$ превращается в функцию с полупериодом $\pi$. После этого
применяется теорема \ref{TH:potoch-Fourier} и предложения \ref{tm-22.2.5} и
\ref{tm-22.2.6}. \end{proof}

\begin{ers} Разложите в ряд Фурье
$2 T$-периодические функции, заданные на полупериоде $(0,T)$, рассмотрев
отдельно случай когда $f(x)$ четна и нечетна:
 \biter{
\item[1)] $f(x)=x, \quad T=4$

\item[2)] $f(x)=3-x, \quad T=3$

\item[3)] $f(x)=\begin{cases} {0, \quad x\in (0,1)}\\
                {\frac{1}{2}, \quad x=1}\\
                {1 \quad x\in (1,2)}
                \end{cases}, \quad T=2$

\item[4)] $f(x)=\begin{cases} 1, \quad x\in (0,1] \\
                2-x \quad x\in (1,2) \end{cases}, \quad T=2$

\item[5)] $f(x)=\begin{cases} 2-x, \quad x\in (0,1] \\
                1 \quad x\in (1,2) \end{cases}, \quad T=2$

\item[6)] $f(x)=\begin{cases} x, \quad x\in (0,1] \\
                1 \quad x\in (1,2) \end{cases}, \quad T=2$
 }\eiter
\end{ers}

\end{multicols}\noindent\rule[10pt]{160mm}{0.1pt}

\chapter{АСИМПТОТИЧЕСКИЕ МЕТОДЫ}\label{ch-o(f(x))}

Поведение функции при аргументе, стремящемся к некоторой величине, является
предметом изучения специального раздела математического анализа, называемого
{\it асимптотическими методами}. По идеологии этой науки, свойство функции $f$
иметь (или не иметь) предел при $x\to a$ -- только самое первое, поверхностное
наблюдение, за которым скрывается гораздо более сложная и многоплановая
картина, изучая которую можно извлечь много полезной информации для решения
задач, где важно, как ведет себя $f$ при $x\to a$.

\noindent\rule{160mm}{0.1pt}\begin{multicols}{2}

\bex Чтобы было понятно, о чем идет речь, вспомним снова о том, как мы
вычисляли пределы. После всех усилий, что мы потратили на изучение этого
вопроса в \ref{SEC:stand-func-v-analize}\ref{SUBSEC:vychisl-predelov} главы
\ref{ch-ELEM-FUNCTIONS}, а затем в \ref{SEC:Lopital} главы \ref{ch-f'(x)}, для
читателя может быть неожиданностью, что встречаются все же пределы, нахождение
которых известными нам сейчас методами оказывается затруднительным или, даже,
невозможным. В качестве упражнения можно попробовать найти следующий, простой
на первый взгляд, предел:
 \beq\label{LIM:reklama-asimp-metodov}
\lim_{x\to 0}\frac{\ln (1+\sin x)-\sin \ln (1+x)}{x^4}
 \eeq
Здесь, если пользоваться правилом Лопиталя, то его придется применять 4 раза,
причем возникающие при дифференцировании выражения столь громоздки, что у
обыкновенного человека терпение истощается уже после первых двух шагов. \eex

\end{multicols}\noindent\rule[10pt]{160mm}{0.1pt}

Асимптотические методы предлагают остроумный способ упростить подобные задачи,
абстрагировавшись от несущественных свойств функций, входящих в выражение под
пределом и оставив только то, что важно. При этом оказывается, что используемый
прием применим и в разных других ситуациях, например, при исследовании на
сходимость несобственных интегралов и рядов. В этой главе мы познакомим
читателя с основными понятиями этой науки и простейшими ее приложениями (а
конкретный предел \eqref{LIM:reklama-asimp-metodov} мы вычислим на с.
\pageref{EX:pokaz-predel-v-asimp-metodah}).

\section{Асимптотические отношения и формулы}

\subsection{Асимптотические отношения}

\paragraph{Асимптотическая эквивалентность функций (символ $\sim$)}

Пусть  $f_1$ и $f_2$ -- две функции, причем
$$
\frac{f_1(x)}{f_2(x)}\underset{x\to a}{\longrightarrow} 1
$$
Тогда говорят, что {\it функция $f_1$ эквивалентна функции $f_2$
при}\index{символ!$\sim$ (асимптотическая эквивалентность)} $x\to a$ и
записывают это так
$$
f_1(x) \underset{x\to a}{\sim} f_2(x)
$$

\noindent\rule{160mm}{0.1pt}\begin{multicols}{2}

В качестве примеров выпишем несколько важных эквивалентностей.
 \begin{align}
\sin x &\underset{x\to 0}{\sim} x \label{11.1.1}\\
\tg x &\underset{x\to 0}{\sim} x\label{11.1.2}\\
1-\cos x &\underset{x\to 0}{\sim}\frac{x^2}{2}\label{11.1.3}\\
\ln (1+x) &\underset{x\to 0}{\sim} x \label{11.1.4}\\
\log_a (1+x) &\underset{x\to 0}{\sim}\frac{x}{\ln a}\label{11.1.5}\\
e^x-1 &\underset{x\to 0}{\sim} x \label{11.1.6}\\
a^x-1 & \underset{x\to 0}{\sim} x\cdot \ln a\label{11.1.7}\\
(1+x)^\alpha-1 &\underset{x\to 0}{\sim}\alpha \cdot  x \label{11.1.8}\\
\arcsin x &\underset{x\to 0}{\sim} x \label{11.1.9}\\
\arctg x &\underset{x\to 0}{\sim} x \label{11.1.10}
 \end{align}
 \begin{proof} Все эти соотношения
дока\-зываются одним и тем же способом -- простым вычислением,
например, с помощью правила Лопиталя. В качестве примера мы
рассмотрим \eqref{11.1.7}: по определению символа $\sim$,
$$
\quad a^x-1 \underset{x\to 0}{\sim} x \ln a \quad \Leftrightarrow
\quad \frac{a^x-1}{x\ln a}\underset{x\to 0}{\longrightarrow} 1
$$
Последнее соотношение доказывается с по\-мощью правила Лопиталя:
 \begin{multline*}\lim_{x\to 0}\frac{a^x-1}{x\ln a}= (\text{Лопиталь})=\\= \lim_{x\to
0}\frac{a^x\ln a}{\ln a}= \lim_{x\to 0} a^x=1
 \end{multline*}\end{proof}\end{multicols}\noindent\rule[10pt]{160mm}{0.1pt}

\bigskip

\centerline{\bf Свойства асимптотической эквивалентности:}
 \bit{\it
\item[$1^\circ$.] Симметричность:
$$
f_1(x)\underset{x\to a}{\sim} f_2(x) \quad \Longrightarrow \quad
f_2(x)\underset{x\to a}{\sim} f_1(x)
$$

\item[$2^\circ$.]  Транзитивность:
$$
f_1(x)\underset{x\to a}{\sim} f_2(x) \quad \& \quad
f_2(x)\underset{x\to a}{\sim} f_3(x) \quad \Longrightarrow \quad
f_1(x)\underset{x\to a}{\sim} f_3(x)
$$

\item[$3^\circ$.] Мультипликативность:
$$
f_1(x)\underset{x\to a}{\sim} f_2(x) \quad \& \quad
g_1(x)\underset{x\to a}{\sim} g_2(x) \quad \Longrightarrow \quad
f_1(x)\cdot g_1(x)\underset{x\to a}{\sim} f_2(x)\cdot g_2(x) \quad \&
\quad \frac{f_1(x)}{g_1(x)}\underset{x\to
a}{\sim}\frac{f_2(x)}{g_2(x)}
$$
 }\eit

Понятие асимптотической эквивалентности полезно при вычислении
сложных пределов, которые слишком утомительно находить с помощью
правила Лопиталя. Для этого используется следующая

\begin{tm}[\bf об эквивалентной замене]
Если $f_1(x) \underset{x\to a}{\sim} f_2(x)$ и $g_1(x)
\underset{x\to a}{\sim} g_2(x)$, то
$$
\lim_{x\to a}\frac{f_1(x)}{g_1(x)}=\lim_{x\to
a}\frac{f_2(x)}{g_2(x)}\qquad \lim_{x\to a} f_1(x)\cdot
g_1(x)=\lim_{x\to a} f_2(x)\cdot g_2(x)
$$
\end{tm}\begin{proof}
$$
\quad \frac{\lim\limits_{x\to
a}\frac{f_1(x)}{g_1(x)}} {\lim\limits_{x\to
a}\frac{f_2(x)}{g_2(x)}}= \lim\limits_{x\to
a}\frac{\frac{f_1(x)}{g_1(x)}}{\frac{f_2(x)}{g_2(x)}}=
\lim\limits_{x\to
a}\frac{\frac{f_1(x)}{f_2(x)}}{\frac{g_1(x)}{g_2(x)}}=
\frac{1}{1}=1 \quad \Rightarrow \quad \lim\limits_{x\to
a}\frac{f_1(x)}{g_1(x)}=\lim\limits_{x\to a}\frac{f_2(x)}{g_2(x)}
$$
 и аналогично,
$$
\frac{\lim\limits_{x\to a} f_1(x)\cdot g_1(x)}
{\lim\limits_{x\to a}f_2(x)\cdot g_2(x)}= \lim\limits_{x\to
a}\frac{f_1(x)}{f_2(x)}\cdot \frac{g_1(x)}{g_2(x)}= 1\cdot 1=1
\quad \Rightarrow \quad \lim\limits_{x\to a} f_1(x)\cdot g_1(x)=
\lim\limits_{x\to a} f_2(x)\cdot g_2(x)
$$
\end{proof}

\noindent\rule{160mm}{0.1pt}\begin{multicols}{2}

\begin{ex} Приведем иллюстрацию:
 \begin{multline*}\lim_{x\to 0}\frac{\sqrt{1+x}-1}{2^x-1}= {\smsize
\left|\begin{matrix}\sqrt{1+x}-1 \underset{x\to 0}{\sim}\frac{x}{2}\\
2^x-1 \underset{x\to 0}{\sim} x\cdot \ln 2
\end{matrix}\right|}=\\=
\lim_{x\to 0}\frac{\frac{x}{2}}{x\cdot \ln 2}= \frac{1}{2\ln 2}
 \end{multline*}\end{ex}

\begin{ex}
 \begin{multline*}\lim_{x\to 0}\frac{1-\cos x}{\ln^2 (1+x)}= {\smsize \left|
\begin{matrix}
1-\cos x \underset{x\to 0}{\sim}\frac{x^2}{2}\\
\ln(1+x) \underset{x\to 0}{\sim} x
\end{matrix}\right|}=\\= \lim_{x\to 0}\frac{\frac{x^2}{2}}{x^2}= \frac{1}{2}
 \end{multline*}\end{ex}

\end{multicols}\noindent\rule[10pt]{160mm}{0.1pt}

В более сложных примерах следует пользоваться еще одной теоремой:

\begin{tm}[\bf о замене переменной под знаком эквивалентности]
\label{ch-var-sym} Если $f_1(y) \underset{y\to b}{\sim} f_2(y)$ и
$\alpha(x) \underset{x\to a}{\longrightarrow} b$, причем
$\alpha(x) \ne b$ при $x\ne a$, то $f_1(\alpha(x)) \underset{x\to
a}{\sim} f_2(\alpha(x))$.
\end{tm}\begin{proof}
$$\lim\limits_{x\to a}\frac{f_1(\alpha(x))}{f_2(\alpha(x))}=
{\smsize {\smsize\begin{pmatrix} y=\alpha(x) \\
y \underset{x\to a}{\longrightarrow} b\\
y\ne b \,\, \text{при}\,\, x\ne a
\end{pmatrix}}}=
\lim\limits_{y\to b}\frac{f_1(y)}{f_2(y)}=1 \quad \Rightarrow
\quad f_1(\alpha(x)) \underset{x\to a}{\sim} f_2(\alpha(x))
$$
\end{proof}

\noindent\rule{160mm}{0.1pt}\begin{multicols}{2}

\begin{ex} Рассмотрим предел посложнее:
$$
\lim_{x\to 0}\frac{\arctg \frac{7}{4} x}{e^{-2x}-1}
$$
Чтобы его вычислить, заметим, что из формулы \eqref{11.1.10}, $\arctg
y \underset{y\to 0}{\sim} y$, можно заменой $y=\frac{7}{4} x$
$(y\underset{x\to 0}{\longrightarrow} 0)$ получить формулу $\arctg
\frac{7}{4} x \underset{x\to 0}{\sim}\frac{7}{4} x$. Точно так же, из
формулы \eqref{11.1.6}, $e^z-1 \underset{z\to 0}{\sim} z$, заменой
$z=-2x$ $(z\underset{x\to 0}{\longrightarrow} 0)$ получается
$e^{-2x}-1 \underset{x\to 0}{\sim} -2x$. Отсюда
 \begin{multline*}\lim_{x\to 0}\frac{\arctg \frac{7}{4} x}{e^{-2x}-1}=\\={\smsize
\left|\begin{matrix}\arctg y \underset{y\to 0}{\sim} y
 \,\Longrightarrow\,
\arctg \frac{7}{4} x\underset{x\to 0}{\sim}\frac{7}{4} x
 \\
e^z-1 \underset{z\to 0}{\sim} z
 \,\Longrightarrow\,
e^{-2x}-1 \underset{x\to 0}{\sim} -2x
\end{matrix}\right|}=\\=
\lim_{x\to 0}\frac{\frac{7}{4} x}{-2x}= \lim_{x\to
0}\frac{\frac{7}{4}}{-2}= -\frac{7}{8}
 \end{multline*}\end{ex}

\begin{ex}
 \begin{multline*}\lim_{x\to 0}\frac{1-\cos 5x}{\tg x^2}={\smsize \left|
\begin{array}{c}
1-\cos 5x \underset{x\to 0}{\sim}\frac{(5x)^2}{2}\\
\tg x^2 \underset{x\to 0}{\sim} x^2
\end{array}\right|}=\\=\lim_{x\to 0}\frac{\frac{(5x)^2}{2}}{x^2}=\lim_{x\to
0}\frac{25x^2}{2x^2}=\frac{25}{2}
 \end{multline*}\end{ex}

\begin{ex}
 \begin{multline*}\lim_{x\to -\infty}\frac{\ln (1+3^x)}{\ln (1+2^x)}={\smsize \left|
\begin{array}{c}   \ln (1+3^x) \underset{x\to -\infty}{\sim} 3^x
\\        \ln (1+2^x) \underset{x\to -\infty}{\sim} 2^x
\end{array}\right|}=\\= \lim_{x\to -\infty}\frac{3^x}{2^x}= \lim_{x\to
-\infty}\left(\frac{3}{2}\right)^x=0
 \end{multline*}\end{ex}

\begin{ex}
 \begin{multline*}\lim_{x\to 0}\frac{\sin (\sqrt{1+x}-1)}{10^x-1}=\\= {\smsize \left|
\begin{array}{c}   \sin (\sqrt{1+x}-1) \underset{x\to 0}{\sim}\sqrt{1+x}-1
\underset{x\to 0}{\sim}\frac{x}{2}\\        10^x-1 \underset{x\to
0}{\sim} x \ln 10
\end{array}\right|}=\\= \lim_{x\to 0}\frac{\frac{x}{2}}{x \ln 10}=\frac{1}{2 \ln
10}
 \end{multline*}\end{ex}

\begin{ers} Найдите пределы

1. $\lim_{x\to 0}\frac{\sqrt{1+x\sin x}-1}{e^{x^2}-1}$

2. $\lim_{x\to 0}\frac{e^{\sin 5x}-1}{\ln (1+ 2x)}$

3. $\lim_{x\to 0}\frac{\arctg(e^{x^3}-1)}{10^{\arcsin^3 x}-1}$

4. $\lim_{x\to 0}\frac{1-\cos x}{\sqrt[3] {1+x^2}-1}$
\end{ers}

\end{multicols}\noindent\rule[10pt]{160mm}{0.1pt}

\paragraph{Асимптотическое сравнение функций (символы $\ll$ и $\bold{o}$) и шкала бесконечностей} \label{CH-symbol-o}

Пусть  $f$ и $g$ -- две функции, причем
$$
\frac{f(x)}{g(x)}\underset{x\to a}{\longrightarrow} 0
$$
Тогда говорят, что {\it функция $f$ бесконечно мала по сравнению с функцией $g$
при}\index{символ!$\ll$ (асимптотическое сравнение)}\index{символ!$\bold{o}$
(асимптотическое сравнение)} $x\to a$ и записывают это так
$$
f(x) \underset{x\to a}{\ll} g(x)
$$
или так
$$
f(x) = \underset{x\to a}{\bold{o}}\Big( g(x) \Big)
$$

\noindent\rule{160mm}{0.1pt}\begin{multicols}{2}

\begin{ex} Покажем, что
$$
x^2 \underset{x\to 0}{\ll} x
$$
или, в других обозначениях,
$$
x^2 = \underset{x\to 0}{\bold{o}}(x)
$$
Действительно, $\frac{x^2}{x}=x \underset{x\to 0}{\longrightarrow}
0.$
\end{ex}

\begin{ex} Убедимся, что
$$
x \underset{x\to \infty}{\ll} x^2
$$
или, в других обозначениях,
$$
x = \underset{x\to \infty}{\bold{o}}(x^2)
$$
Действительно, $\frac{x}{x^2}=\frac{1}{x}\underset{x\to
\infty}{\longrightarrow} 0.$
\end{ex}

\end{multicols}\noindent\rule[10pt]{160mm}{0.1pt}

\bigskip

\centerline{\bf Свойства асимптотического сравнения}
 \bit{\it
\item[$1^\circ$] Транзитивность:
$$
f_1(x)\underset{x\to a}{\ll} f_2(x) \quad \& \quad
f_2(x)\underset{x\to a}{\ll} f_3(x) \quad \Longrightarrow \quad
f_1(x)\underset{x\to a}{\ll} f_3(x)
$$
}\eit

\noindent\rule{160mm}{0.1pt}\begin{multicols}{2}

\begin{tm}[\bf о сравнении степенных\\ функций]
\label{x^a<<x^b} Справедливы соотношения

$$
\boxed{\quad\phantom{\Big|} x^\alpha \underset{(\e>0)}{\underset{x\to
\infty}{\ll}} x^{\alpha+\e}, \qquad
x^{\alpha+\e}\underset{(\e>0)}{\underset{x\to 0}{\ll}} x^\alpha
\phantom{\Big|}\quad}
$$
или, в других обозначениях, при $\alpha<\beta$,
$$
x^\alpha  = \underset{x\to \infty}{\bold{o}}\left( x^\beta \right),
\quad x^\beta  = \underset{x\to 0}{\bold{o}}\left( x^\alpha \right).
$$
\end{tm}\begin{proof}
$$
\frac{x^\alpha}{x^{\alpha+\e}}=\frac{1}{x^{\e}}\underset{x\to\infty}{\longrightarrow}
0
$$
и это означает, что $x^\alpha \underset{x\to \infty}{\ll}
x^{\alpha+\e}$. Наоборот,
$$
\frac{x^{\alpha+\e}}{x^\alpha}=x^\e \underset{x\to
0}{\longrightarrow} 0
$$
и это означает, что $x^{\alpha+\e}\underset{x\to 0}{\ll} x^\alpha$.
\end{proof}

\begin{tm}[\bf о сравнении показательных функций]
\label{a^x<<b^x} Справедливы соотношения
$$
\boxed{\quad\phantom{\Big|} a^x \underset{(\e>0)}{\underset{x\to
+\infty}{\ll}} (a+\e)^x, \qquad (a+\e)^x
\underset{(\e>0)}{\underset{x\to -\infty}{\ll}} a^x
\phantom{\Big|}\quad}
$$
или, в других обозначениях, при $0<a<b$,
$$
a^x = \underset{x\to +\infty}{\bold{o}}\left( b^x \right), \quad
 b^x= \underset{x\to -\infty}{\bold{o}}\left( a^x \right).
$$
\end{tm}\begin{proof}
$$
 \frac{a^x}{(a+\e)^x}=
 \Big(\underbrace{\frac{a}{(a+\e)}}_{\tiny\begin{matrix}\wedge \\ 1
 \end{matrix}}\Big)^x \underset{x\to +\infty}{\longrightarrow} 0
$$
То есть, $a^x \underset{x\to +\infty}{\ll} (a+\e)^x$. Наоборот,
$$
\frac{(a+\e)^x}{a^x}=\Big(\underbrace{\frac{a+\e}{a}}_{\tiny\begin{matrix}\vee
\\ 1\end{matrix}}\Big)^x
\underset{x\to -\infty}{\longrightarrow} \infty,
$$
и это означает, что $(a+\e)^x \underset{x\to -\infty}{\ll} a^x$.
\end{proof}

\begin{tm}[\bf о шкале бесконечностей]
\label{scale-infties} Справедливы соотношения
 \beq\label{shkala-besk}
\boxed{\quad\phantom{\Big|} 1 \underset{x\to +\infty}{\ll}\,
\underset{(a>1)}{\log_a x}\, \underset{x\to +\infty}{\ll}\,
\underset{(\alpha>0)}{x^\alpha}\, \underset{x\to +\infty}{\ll}\,
\underset{(a>1)}{a^x}\phantom{\Big|}\quad}
 \eeq
или, в других обозначениях,
 \begin{align*}
1 &=\underset{x\to +\infty}{\bold{o}}\left(\log_a x  \right), &&
\text{если}\,\, a>1, \\
\log_a x &= \underset{x\to +\infty}{\bold{o}}\left( x^\alpha
\right), && \text{если}\,\, a>1, \alpha>0\\
x^\alpha &= \underset{x\to +\infty}{\bold{o}}\left( a^x \right), &&
\text{если}\,\, \alpha>0, a>1
 \end{align*}\end{tm}\begin{proof} Первое соотношение очевидно, потому что
$$
\lim_{x\to +\infty}\frac{1}{\log_a x}=
\left(\frac{1}{\infty}\right)=0
$$
Докажем второе:
 \begin{multline*}\lim_{x\to +\infty}\frac{\log_a x}{x^\alpha}=
{\smsize{\smsize\begin{pmatrix}\text{правило}\\
\text{Лопиталя}\end{pmatrix}}}=\\= \lim_{x\to +\infty}\frac{(\log_a
x)'}{(x^\alpha)'}= \lim_{x\to +\infty}\frac{\frac{1}{x \cdot \ln
a}}{\alpha \cdot x^{\alpha-1}}=\\= \lim_{x\to +\infty}\frac{1}{\alpha
\ln a  \cdot x \cdot x^{\alpha-1}}=\frac{1}{\alpha \ln a}\cdot
\lim_{x\to +\infty}\frac{1}{x^\alpha}=\\={\smsize
{\smsize\begin{pmatrix}\text{учитываем, что}\\
\alpha>0
\end{pmatrix}}}= 0 \end{multline*}

А теперь третье:
 \begin{multline*}\lim_{x\to +\infty}\frac{x^\alpha}{a^x}={\smsize
{\smsize\begin{pmatrix}\text{правило}\\
\text{Лопиталя}\end{pmatrix}}}= \lim_{x\to
+\infty}\frac{(x^\alpha)'}{(a^x)'}=\\= \lim_{x\to
+\infty}\frac{\alpha x^{\alpha-1}}{a^x \ln a}={\smsize
{\smsize\begin{pmatrix}\text{если $\alpha-1>0$, то еще раз}\\
\text{применяем правило Лопиталя}\end{pmatrix}}}=\\= \lim_{x\to
+\infty}\frac{\alpha \cdot (\alpha-1)\cdot x^{\alpha-2}} {a^x (\ln
a)^2}=\\={\smsize {\smsize\begin{pmatrix}\text{если $\alpha-2>0$, то
еще раз}\\
\text{применяем правило Лопиталя}\end{pmatrix}}}=\\= \lim_{x\to
+\infty}\frac{\alpha \cdot (\alpha-1)\cdot (\alpha-2)\cdot
x^{\alpha-3}}{a^x (\ln a)^3}=...=\\={\smsize {\smsize\begin{pmatrix}\text{применяем правило Лопиталя}\\
\text{до тех пор, пока в числителе}\\
\text{не получится $x^{\alpha-n}$, где
$\alpha-n<0$}\end{pmatrix}}}=...=\\ \text{\smsize $=\lim_{x\to
+\infty}\frac{\alpha \cdot (\alpha-1)\cdot (\alpha-2)\cdot ... \cdot
(\alpha-n+1) \cdot x^{\alpha-n}}{a^x (\ln a)^n}=$}\\={\smsize
{\smsize\begin{pmatrix}\text{поскольку}\,\,
\alpha-n<0, \\
x^{\alpha-n}\underset{x\to +\infty}{\longrightarrow} 0
\end{pmatrix}}}= \left(\frac{0}{\infty}\right)= 0
\end{multline*}\end{proof}

\end{multicols}\noindent\rule[10pt]{160mm}{0.1pt}

\paragraph{Принцип выделения главного слагаемого}

Символы $\ll$ и $\bold{o}$, описанные в предыдущем параграфе,
бывают полезны в различных ситуациях при вычислении пределов.
Самое простое их применение описывается в следующей теореме.

\begin{tm}[\bf о выделении главного слагаемого]\label{th-11.3.1}
Пусть даны несколько функций
$$
A_0(x), A_1(x),...,A_n(x),
$$
причем при $x\to a$ функция $A_0(x)$ бесконечно велика по
сравнению с любой другой функцией из этого списка:
\begin{equation}
A_1(x)\underset{x\to a}{\ll} A_0(x), \quad A_2(x)\underset{x\to
a}{\ll} A_0(x), \quad ... , \quad A_n(x)\underset{x\to a}{\ll} A_0(x)
\label{11.3.1}\end{equation} Тогда сумма этих функций эквивалентна
$A_0(x)$ при $x\to a$:
\begin{equation}
A_0(x)+A_1(x)+A_2(x)+...+A_n(x) \underset{x\to a}{\sim} A_0(x)
\label{11.3.2}\end{equation}\end{tm}\begin{proof} Из \eqref{11.3.1}
следует, что
$$
\frac{A_1(x)}{A_0(x)}\underset{x\to a}{\longrightarrow} 0, \quad
\frac{A_2(x)}{A_0(x)}\underset{x\to a}{\longrightarrow} 0, \quad
... \quad \frac{A_n(x)}{A_0(x)}\underset{x\to a}{\longrightarrow}
0
$$
поэтому
$$
\frac{A_0(x)+A_1(x)+A_2(x)+...+A_n(x)}{A_0(x)}=
1+\frac{A_1(x)}{A_0(x)}+\frac{A_2(x)}{A_0(x)}+...+\frac{A_n(x)}{A_0(x)}
\underset{x\to a}{\longrightarrow} 1+0+0+...+0=1
$$
то есть, выполняется \eqref{11.3.2}. \end{proof}

\bit{ \item[$\bullet$] Функция $A_0(x)$, для которой выполняются соотношения
\eqref{11.3.1}, называется {\it главным слагаемым} в сумме
$A_0(x)+A_1(x)+A_2(x)+...+A_n(x)$.
 }\eit

Теорема \ref{th-11.3.1} применяется в следующей ситуации. Пусть
нам необходимо вычислить предел
$$
\lim_{x\to a}\frac{A_0(x)+A_1(x)+A_2(x)+...+A_n(x)}
{B_0(x)+B_1(x)+B_2(x)+...+B_n(x)}
$$
причем известно, что $A_0(x)$ является главным слагаемым в сумме
$A_0(x)+A_1(x)+A_2(x)+...+A_n(x)$, а $B_0(x)$ является главным
слагаемым в сумме $B_0(x)+B_1(x)+B_2(x)+...+B_n(x)$:
$$
A_1(x)\underset{x\to a}{\ll} A_0(x), \quad A_2(x)\underset{x\to
a}{\ll} A_0(x), \quad ... \quad A_n(x)\underset{x\to a}{\ll}
A_0(x)
$$
$$
B_1(x)\underset{x\to a}{\ll} B_0(x), \quad B_2(x)\underset{x\to
a}{\ll} B_0(x), \quad ... \quad B_n(x)\underset{x\to a}{\ll}
B_0(x)
$$
Тогда
$$
A_0(x)+A_1(x)+A_2(x)+...+A_n(x) \underset{x\to a}{\sim} A_0(x)
$$
$$
B_0(x)+B_1(x)+B_2(x)+...+B_n(x) \underset{x\to a}{\sim} B_0(x)
$$
и поэтому
$$
\lim_{x\to a}\frac{A_0(x)+A_1(x)+A_2(x)+...+A_n(x)}
{B_0(x)+B_1(x)+B_2(x)+...+B_n(x)}= \lim_{x\to
a}\frac{A_0(x)}{B_0(x)}
$$
Таким образом, справедлив следующий важный принцип.

\bigskip

\begin{tm}[\bf принцип выделения главного слагаемого]
\index{принцип!выделения главного слагаемого}\label{main+princ}
Вычисляя предел
$$
\lim_{x\to a}\frac{A_0(x)+A_1(x)+A_2(x)+...+A_n(x)}
{B_0(x)+B_1(x)+B_2(x)+...+B_n(x)},
$$
следует выделить в числителе и знаменателе главные слагаемые, а
все остальные -- ``второстепенные'' -- просто зачеркнуть:
 $$
\lim_{x\to a}\frac{A_0(x)+A_1(x)+A_2(x)+...+A_n(x)}
{B_0(x)+B_1(x)+B_2(x)+...+B_n(x)}= \lim_{x\to a}\frac{ A_0(x)+
A_1(x)\kern-15pt{\Big/\kern-8pt\Big\backslash}\kern10pt
+A_2(x)\kern-15pt{\Big/\kern-8pt\Big\backslash}\kern10pt+...
+A_n(x)\kern-15pt{\Big/\kern-8pt\Big\backslash}} {B_0(x)
+B_1(x)\kern-15pt{\Big/\kern-8pt\Big\backslash}\kern10pt
+B_2(x)\kern-15pt{\Big/\kern-8pt\Big\backslash}\kern10pt+...
+B_n(x)\kern-15pt{\Big/\kern-8pt\Big\backslash}}= \lim_{x\to
a}\frac{A_0(x)}{B_0(x)}
 $$
\end{tm}

\bigskip

Рассмотрим примеры.

\noindent\rule{160mm}{0.1pt}\begin{multicols}{2}

\begin{ex}
 \begin{multline*}
\lim_{x\to +\infty}\frac{x+\ln x}{x^2+\sin x}=
{\smsize\begin{pmatrix}\text{по теореме \ref{scale-infties}}\\
\text{(о шкале бесконечностей),}\\
\ln x\underset{x\to +\infty}{\ll} x, \\
\text{кроме того,}\\ \text{как легко проверить,}\\
\sin x\underset{x\to +\infty}{\ll} x^2
\end{pmatrix}}=\\= \lim_{x\to +\infty}\frac{x+ \ln
x\kern-15pt{\Big/\kern-8pt\Big\backslash}\kern10pt }{x^2+ \sin
x\kern-15pt{\Big/\kern-8pt\Big\backslash}\kern10pt }= \lim_{x\to
0}\frac{x}{x^2}=0 \end{multline*}\end{ex}

\begin{ex} В следующем примере снова используется шкала бесконечностей
(теорема \ref{scale-infties}):
 \begin{multline*}
\lim_{x\to +\infty}\frac{2\sqrt{x} +2\log_5 x +5 }{ \cos x +\sqrt{x}
-1 }=\\= \lim_{x\to +\infty}\frac{2\sqrt{x} +2\log_5
x\kern-15pt{\Big/\kern-8pt\Big\backslash}\kern10pt
+5\kern-7pt{\Big/\kern-8pt\Big\backslash}\kern10pt }{ \cos
x\kern-15pt{\Big/\kern-8pt\Big\backslash}\kern10pt +\sqrt{x}
-1\kern-7pt{\Big/\kern-8pt\Big\backslash}\kern10pt }= \lim_{x\to
+\infty}\frac{2\sqrt{x}}{\sqrt{x}}=2
 \end{multline*}
\end{ex}

\begin{ex}
Здесь используется теорема \ref{x^a<<x^b} (о сравнении степенных
функций): поскольку $x\to \infty$, главными слагаемыми будут
максимальные степени у $x$:
 \begin{multline*}
\lim_{x\to \infty}\frac{5x+3\sqrt{x}-2x^5}{x-3x^4+8\sqrt{x}}=\\=
\lim_{x\to \infty}\frac{
5x\kern-9pt{\Big/\kern-8pt\Big\backslash}\kern10pt
+3\sqrt{x}\kern-15pt{\Big/\kern-8pt\Big\backslash}\kern10pt -2x^5 }{
x\kern-7pt{\Big/\kern-8pt\Big\backslash}\kern10pt -3x^4
+8\sqrt{x}\kern-15pt{\Big/\kern-8pt\Big\backslash}\kern10pt }=
\lim_{x\to \infty}\frac{-2x^5}{-3x^4}=\\= \lim_{x\to
\infty}\frac{-2x}{-3}=\infty
 \end{multline*}
\end{ex}

\begin{ex}
Здесь дробь та же, что и в предыдущем примере, но, поскольку $x\to
0$, главными слагаемыми, по той же теореме \ref{x^a<<x^b}, будут
наоборот, минимальные степени у $x$:
 \begin{multline*}
\lim_{x\to 0}\frac{5x+3\sqrt{x}-2x^5}{x-3x^4+8\sqrt{x}}=\\=
\lim_{x\to 0}\frac{
5x\kern-9pt{\Big/\kern-8pt\Big\backslash}\kern10pt +3\sqrt{x}
-2x^5\kern-15pt{\Big/\kern-8pt\Big\backslash}\kern10pt }{
x\kern-7pt{\Big/\kern-8pt\Big\backslash}\kern10pt
-3x^4\kern-15pt{\Big/\kern-8pt\Big\backslash}\kern10pt +8\sqrt{x} }=
\lim_{x\to 0}\frac{3\sqrt{x}}{8\sqrt{x}}= \frac{3}{8}
 \end{multline*}
\end{ex}

\begin{ex}
Здесь используется теорема \ref{a^x<<b^x} (о сравнении
показательных функций): поскольку $x\to +\infty$, главными
слагаемыми будут степени с наибольшим основанием:
 \begin{multline*}
\lim_{x\to +\infty}\frac{2^{x-1}-3^{x+1}}{2^{x+1}+3^{x-1}}=
\lim_{x\to +\infty}\frac{ \frac{1}{2}\cdot 2^x -3\cdot 3^x}{ 2\cdot
2^x +\frac{1}{3}\cdot 3^x}=\\= \lim_{x\to +\infty}\frac{
\frac{1}{2}\cdot 2^x\kern-15pt{\Big/\kern-8pt\Big\backslash}\kern10pt
-3\cdot 3^x}{ 2\cdot
2^x\kern-15pt{\Big/\kern-8pt\Big\backslash}\kern10pt
+\frac{1}{3}\cdot 3^x}= \lim_{x\to +\infty}\frac{ -3\cdot 3^x}{
\frac{1}{3}\cdot 3^x}=-9
 \end{multline*}
\end{ex}

\begin{ex}
Здесь дробь та же, что и в предыдущем примере, но, поскольку $x\to
-\infty$, главными слагаемыми, по той же теореме \ref{a^x<<b^x},
будут наоборот, степени с наименьшим основанием:
 \begin{multline*}
\lim_{x\to -\infty}\frac{2^{x-1}-3^{x+1}}{2^{x+1}+3^{x-1}}=
\lim_{x\to -\infty}\frac{ \frac{1}{2}\cdot 2^x -3\cdot 3^x}{ 2\cdot
2^x +\frac{1}{3}\cdot 3^x}=\\= \lim_{x\to -\infty}\frac{
\frac{1}{2}\cdot 2^x -3\cdot
3^x\kern-15pt{\Big/\kern-8pt\Big\backslash}\kern10pt }{ 2\cdot 2^x
+\frac{1}{3}\cdot
3^x\kern-15pt{\Big/\kern-8pt\Big\backslash}\kern10pt }= \lim_{x\to
-\infty}\frac{ \frac{1}{2}\cdot 2^x }{ 2\cdot 2^x }=\frac{1}{4}
 \end{multline*}
\end{ex}

\begin{ex}
В следующем примере мы сначала выделяем главное слагаемое, а затем
применяем теорему \ref{ch-var-sym} об эквивалентной замене:
 \begin{multline*}
\lim_{x\to 0}\frac{5x^2+3x^{10}}{2x^5+7\sin x^2}= \lim_{x\to 0}\frac{
5x^2 +3x^{10}\kern-15pt{\Big/\kern-8pt\Big\backslash}\kern10pt }{
2x^5\kern-15pt{\Big/\kern-8pt\Big\backslash}\kern10pt +7\sin x^2
}=\\= \lim_{x\to 0}\frac{5x^2}{7\sin x^2}= \left(\sin x^2
\underset{x\to 0}{\sim} x^2 \right)= \lim_{x\to
0}\frac{5x^2}{7x^2}=\frac{5}{7}
 \end{multline*}
\end{ex}

\begin{ers}$\phantom{.}$
 \biter{
\item[1.] $\lim_{x\to 0}\frac{2^x+ x^5}{x^2+5^x}$

\item[2.] $\lim_{x\to +\infty}\frac{2^x+ x^5}{x^2+5^x}$

\item[3.] $\lim_{x\to 0}\frac{\sqrt[3] {x}+ \sqrt {x} } {\sin x+
\sqrt[5]{x^2}}$

\item[4.] $\lim_{x\to +\infty}\frac{\sqrt[3] {x}+\sqrt{x}} {\sin x+
\sqrt[5]{x^2}}$

\item[5.] $\lim_{x\to 0}\frac{x+ \sqrt x}{e^x-1}$

\item[6.] $\lim_{x\to +\infty}\frac{x+ \sqrt x}{e^x-1}$

\item[7.] $\lim_{x\to 0}\frac{2^x+3^x}{x^2+x^3}$

\item[8.] $\lim_{x\to +\infty}\frac{2^x+3^x}{x^2+x^3}$

\item[0.] $\lim_{x\to -\infty}\frac{2^x+3^x}{x^2+x^3}$

\item[10.] $\lim_{x\to +\infty}\frac{x+\sin x}{x+\cos x}$
 }\eiter
 \end{ers}

\end{multicols}\noindent\rule[10pt]{160mm}{0.1pt}

\paragraph{Эквивалентность модулей, степеней и логарифмов}

\begin{tm}
Пусть даны две эквивалентные функции при $x\to a$:
$$
f_1(x) \underset{x\to a}{\sim} f_2(x)
$$
Тогда
 \bit{
\item[(i)] модули этих функций тоже эквивалентны:
$$
|f_1(x) |\underset{x\to a}{\sim} |f_2(x)|
$$
\item[(ii)] степени этих функций (с произвольным показателем)
тоже эквивалентны:
$$
f_1(x)^\alpha \underset{x\to a}{\sim} f_2(x)^\alpha, \quad \alpha\in \R
$$
\item[(iii)] если вдобавок $f_2(x)$ стремится к бесконечности или к нулю
$$
f_2(x) \underset{x\to a}{\longrightarrow} +\infty \quad
(\text{или}\quad +0)
$$
то логарифмы этих функций (с произвольным основанием)
эквивалентны:
$$
\log_a f_1(x) \underset{x\to a}{\sim}\log_a f_2(x), \quad 0<a\ne 1
$$
 }\eit
\end{tm}
\begin{proof}
 \begin{align*}
 1.\quad & f_1(x) \underset{x\to a}{\sim} f_2(x) && \Rightarrow \quad
\frac{f_1(x)}{f_2(x)}\underset{x\to a}{\longrightarrow} 1 \quad \Rightarrow
\quad \frac{|f_1(x)|}{|f_2(x)|}=\left|\frac{f_1(x)}{f_2(x)}\right|
\underset{x\to a}{\longrightarrow} 1 \\ &&&  \Rightarrow \quad
|f_1(x)| \underset{x\to a}{\sim} |f_2(x)|  \quad\Rightarrow \\
2.\quad & f_1(x) \underset{x\to a}{\sim} f_2(x) && \Rightarrow \quad
\frac{f_1(x)}{f_2(x)}\underset{x\to a}{\longrightarrow} 1 \quad \Rightarrow
\quad
\frac{f_1(x)^\alpha}{f_2(x)^\alpha}=\left(\frac{f_1(x)}{f_2(x)}\right)^\alpha
\underset{x\to a}{\longrightarrow} 1 \quad\Rightarrow  \\ &&& \Rightarrow \quad
f_1(x)^\alpha \underset{x\to a}{\sim} f_2(x)^\alpha
\\
3.\quad &
 f_1(x) \underset{x\to a}{\sim} f_2(x) \underset{x\to a}{\longrightarrow}\begin{cases}+\infty\\ +0\end{cases}
 && \Rightarrow\quad
 \frac{\log_a f_1(x)}{\log_a f_2(x)}= \frac{\log_a \left(f_2(x)\cdot
 \frac{f_1(x)}{f_2(x)}\right)}{\log_a f_2(x)}=\\ & && = \frac{\log_a f_2(x)+ \log_a
 \frac{f_1(x)}{f_2(x)}}{\log_a f_2(x)}=
 1+
 \overbrace{
 \frac{1}{\log_a
 \kern-8pt\underbrace{f_2(x)}_{\tiny\begin{matrix}\downarrow \\ +\infty \; \text{или}\; +0 \end{matrix}}
 }
 }^{\tiny\begin{matrix} 0 \\ \uparrow \end{matrix}}
 \cdot
 \overbrace{
 \log_a
 \underset{\tiny\begin{matrix}\downarrow \\ 1 \end{matrix}}
 {
 \boxed{
 \frac{f_1(x)}{f_2(x)}
 }
 }
 }^{\tiny\begin{matrix} 0 \\ \uparrow  \end{matrix}}
  \underset{x\to a}{\longrightarrow} 1
  \end{align*}
 \end{proof}

Из этой теоремы следует, что при вычислении пределов с модулями, степенями или
логарифмами также можно выделять главные слагаемые, {\it но нужно только
следить, чтобы аргументы у логарифмов стремились к $+\infty$ или к $+0$}.
Например,

\begin{multline*}\lim_{x\to a}\frac{\left|A_0(x)+A_1(x)+A_2(x)+...+A_n(x)\right|}
{\left(B_0(x)+B_1(x)+B_2(x)+...+B_n(x)\right)^\beta}=\\= \lim_{x\to
a}\frac{\left| A_0(x)+
A_1(x)\kern-15pt{\Big/\kern-8pt\Big\backslash}\kern10pt
+A_2(x)\kern-15pt{\Big/\kern-8pt\Big\backslash}\kern10pt+...
+A_n(x)\kern-15pt{\Big/\kern-8pt\Big\backslash}\kern10pt \right|}
{\left(B_0(x)
+B_1(x)\kern-15pt{\Big/\kern-8pt\Big\backslash}\kern10pt
+B_2(x)\kern-15pt{\Big/\kern-8pt\Big\backslash}\kern10pt+...
+B_n(x)\kern-15pt{\Big/\kern-8pt\Big\backslash}\kern10pt
\right)^\beta}= \lim_{x\to
a}\frac{|A_0(x)|}{B_0(x)^\beta}\end{multline*} Или, если $B_0(x)
\underset{x\to a}{\longrightarrow} +\infty \quad (\text{или}\quad
+0)$,

\begin{multline*}\lim_{x\to
a}\frac{\left(A_0(x)+A_1(x)+A_2(x)+...+A_n(x)\right)^\alpha} {\ln
\left(B_0(x)+B_1(x)+B_2(x)+...+B_n(x)\right)}=\\= \lim_{x\to a}\frac{
\left( A_0(x)+
A_1(x)\kern-15pt{\Big/\kern-8pt\Big\backslash}\kern10pt
+A_2(x)\kern-15pt{\Big/\kern-8pt\Big\backslash}\kern10pt+...
+A_n(x)\kern-15pt{\Big/\kern-8pt\Big\backslash}\kern10pt\right)^\alpha}
{\ln \left(B_0(x)
+B_1(x)\kern-15pt{\Big/\kern-8pt\Big\backslash}\kern10pt
+B_2(x)\kern-15pt{\Big/\kern-8pt\Big\backslash}\kern10pt+...
+B_n(x)\kern-15pt{\Big/\kern-8pt\Big\backslash}\kern10pt\right)}=
\lim_{x\to a}\frac{A_0(x)^\alpha}{\ln B_0(x)}\end{multline*} И в
таком духе. Рассмотрим несколько примеров.

\noindent\rule{160mm}{0.1pt}\begin{multicols}{2}

\begin{ex}
 \begin{multline*}
\lim_{x\to +\infty}\frac{|\sqrt[3] {x}
+3x-2x^2|}{\sqrt{x+x^2+x^4}}=\\= \lim_{x\to +\infty}\frac{| \sqrt[3]
{x}\kern-10pt{\Big/\kern-8pt\Big\backslash}\kern10pt
+3x\kern-10pt{\Big/\kern-8pt\Big\backslash}\kern10pt -2x^2|}{\sqrt{
x\kern-7pt{\Big/\kern-8pt\Big\backslash}\kern10pt
+x^2\kern-10pt{\Big/\kern-8pt\Big\backslash}\kern10pt +x^4}}=\\=
\lim_{x\to +\infty}\frac{|-2x^2|}{\sqrt{x^4}}= \lim_{x\to
+\infty}\frac{2x^2}{x^2}=2
 \end{multline*}
\end{ex}

\begin{ex}
 \begin{multline*}
\lim_{x\to 0}\frac{|\sqrt[3] {x} +3x-2x^2|}{\sqrt{x+x^2+x^4}}=\\=
\lim_{x\to 0}\frac{| \sqrt[3] {x}
+3x\kern-10pt{\Big/\kern-8pt\Big\backslash}\kern10pt
-2x^2\kern-12pt{\Big/\kern-8pt\Big\backslash}\kern10pt|}{\sqrt{ x
+x^2\kern-10pt{\Big/\kern-8pt\Big\backslash}\kern10pt
+x^4\kern-10pt{\Big/\kern-8pt\Big\backslash}\kern10pt }}=\\=
\lim_{x\to 0}\frac{|\sqrt[3] {x}|}{\sqrt{x}}= \lim_{x\to 0}
x^{-\frac{1}{6}}=\infty
 \end{multline*}
\end{ex}

\begin{ex}
 \begin{multline*}
\lim_{x\to +\infty}\frac{\ln(x^2+x)}{\ln (x^5+x^3)}= \lim_{x\to
+\infty}\frac{\ln (x^2 +x\kern-7pt{\Big/\kern-8pt\Big\backslash}\,)
}{\ln (x^5 +x^3\kern-10pt{\Big/\kern-8pt\Big\backslash}\,) }=\\=
\lim_{x\to +\infty}\frac{\ln x^2}{\ln x^5}= \lim_{x\to
+\infty}\frac{2\ln x}{5\ln x}= \frac{2}{5}
 \end{multline*}
\end{ex}

\begin{ex}
 \begin{multline*}
\lim_{x\to +0}\frac{\ln(x^2+x)}{\ln (x^5+x^3)}= \lim_{x\to
+0}\frac{\ln ( x^2\kern-10pt{\Big/\kern-8pt\Big\backslash}\kern10pt
+x) }{\ln ( x^5\kern-10pt{\Big/\kern-8pt\Big\backslash}\kern10pt
+x^3) }=\\= \lim_{x\to +0}\frac{\ln x}{\ln x^3}= \lim_{x\to
+0}\frac{\ln x}{3\ln x}= \frac{1}{3}
 \end{multline*}
\end{ex}

\begin{ex}
 \begin{multline*}
\lim_{x\to +\infty}\frac{3\ln x+x-2}
{\ln\left|5+x^3-e^{2x}\right|}=\\= \lim_{x\to +\infty}\frac{ 3\ln
x\kern-12pt{\Big/\kern-8pt\Big\backslash}\kern10pt +x
-2\kern-7pt{\Big/\kern-8pt\Big\backslash}\kern10pt } {\ln\left|
5\kern-7pt{\Big/\kern-8pt\Big\backslash}\kern10pt
+x^3\kern-10pt{\Big/\kern-8pt\Big\backslash}\kern10pt
-e^{2x}\right|}=\\= \lim_{x\to
+\infty}\frac{x}{\ln\left|-e^{2x}\right|}= \lim_{x\to
+\infty}\frac{x}{2x}=\frac{1}{2}
 \end{multline*}
\end{ex}

\begin{ers}
Вычислите пределы, обращая внимание на значение, к которому
стремится аргумент $x$:

1. $\lim\limits_{x\to +\infty}\frac {\ln \left| x-2\ln
x-3\sqrt{x}\right|} {\sqrt[3] {\ln x+2x\sqrt{x}-x^2}}$

2. $\lim\limits_{x\to 0}\frac {\ln \left| x-2\ln x-3\sqrt{x}\right|}
{\sqrt[3] {\ln x+2x\sqrt{x}-x^2}}$

3. $\lim\limits_{x\to +\infty}\frac {\ln \left(
1+\sqrt{x}+\sqrt[3]{x}\right)} {\ln \left(
1+\sqrt[3]{x}+\sqrt[4]{x}\right)}$

4. $\lim\limits_{x\to 0}\frac {\ln \left(
1+\sqrt{x}+\sqrt[3]{x}\right)} {\ln \left(
1+\sqrt[3]{x}+\sqrt[4]{x}\right)}$

5. $\lim\limits_{x\to 0}\frac {\ln \left(\sqrt{x}+\sqrt[3]{x}\right)}
{\ln \left(\sqrt[3]{x}+\sqrt[4]{x}\right)}$

6. $\lim\limits_{x\to +\infty}\frac{\ln \left( 1+3^x \right)}{\ln
\left( 1+2^x \right)}$

7. $\lim\limits_{x\to 0}\frac{\ln \left( 1+3^x \right)}{\ln \left(
1+2^x \right)}$
 \end{ers}

\end{multicols}\noindent\rule[10pt]{160mm}{0.1pt}

\subsection{Асимптотические формулы}

В математическом анализе имеется специальный язык {\it
асимптотических формул}, позволяющий описывать различные связи
между функциями с помощью символов $\bold{o}$. Для его понимания
достаточно просто знать, что если в какой-то формуле встречается
символ $\bold{o}\Big( f(x) \Big)$, то это означает, что на этом
месте стоит некоторая (неизвестная) функция $\alpha (x)$, для
которой справедливо соотношение $\alpha (x)=\underset{x\to
a}{\bold{o}\Big( f(x) \Big)}$.

Как и все на свете, это лучше понимать на примерах.

\noindent\rule{160mm}{0.1pt}\begin{multicols}{2}

 \bex
Скажем, асимптотическая фор\-мула
$$
\sin x=x+\underset{x\to 0}{\bold{o}(x)}
$$
означает, что
$$
\sin x=x+\alpha (x), \qquad \text{где}\quad \alpha (x)=\underset{x\to
0}{\bold{o}(x)}
$$
Ее нетрудно доказать:
 \begin{multline*}
\sin x=x+\underset{x\to 0}{\bold{o}(x)}\quad \Leftrightarrow \quad
\sin x-x=\underset{x\to 0}{\bold{o}(x)}\quad \Leftrightarrow \\
\Leftrightarrow \quad \frac{\sin x-x}{x}\underset{x\to
0}{\longrightarrow} 0 \quad \Leftrightarrow \quad
\\
\quad \Leftrightarrow \quad \frac{\sin x}{x}-1\underset{x\to
0}{\longrightarrow} 0 \quad \Leftrightarrow \quad \frac{\sin
x}{x}\underset{x\to 0}{\longrightarrow} 1 \end{multline*} Последнее
соотношение верно, поскольку представляет собой первый замечательный
предел.
 \eex
 \bex
Или вот такая асимптотическая формула:
$$
\cos x=1-\frac{x^2}{2}+\underset{x\to 0}{\bold{o}(x^2)}
$$
Она означает, что
$$
\cos x=1-\frac{x^2}{2}+\alpha (x), \qquad  \text{где}\quad \alpha
(x)=\underset{x\to 0}{\bold{o}(x^2)}
$$
и доказать ее также нетрудно:
 \begin{multline*}
\cos x=1-\frac{x^2}{2}+\underset{x\to 0}{\bold{o}(x^2)}\quad
\Leftrightarrow \\ \Leftrightarrow \quad \cos
x-1+\frac{x^2}{2}=\underset{x\to 0}{\bold{o}(x^2)}\quad
\Leftrightarrow \\ \Leftrightarrow \quad \frac{\cos
x-1+\frac{x^2}{2}}{x^2}\underset{x\to 0}{\longrightarrow} 0 \quad
\Leftrightarrow \\ \Leftrightarrow\quad \frac{\cos
x-1}{x^2}+\frac{1}{2}\underset{x\to 0}{\longrightarrow} 0 \quad
\Leftrightarrow \quad \frac{\cos x-1}{x^2}\underset{x\to
0}{\longrightarrow} -\frac{1}{2}\end{multline*} Последнее соотношение
доказывается напрямую:
 \begin{multline*}
\lim_{x\to 0}\frac{\cos x-1}{x^2}={\smsize
\begin{pmatrix}\text{правило}\\ \text{Лопиталя}\end{pmatrix}}= \lim_{x\to
0}\frac{-\sin x}{2x}=\\= {\smsize
\begin{pmatrix}\text{снова правило}\\ \text{Лопиталя}\end{pmatrix}}
= \lim_{x\to 0}\frac{-\cos x}{2}=-\frac{1}{2}\end{multline*}
 \eex

Выпишем несколько важных асимптотических формул:
\begin{align}
&\sin x=x+\underset{x\to
0}{\bold{o}(x)}\label{11.5.1}
\\
& \tg x=x+\underset{x\to
0}{\bold{o}(x)}\label{11.5.2}\\
& \cos x=1-\frac{x^2}{2}+\underset{x\to
0}{\bold{o}(x^2)}\label{11.5.3}\\
& \ln (1+x)=x+\underset{x\to
0}{\bold{o}(x)}\label{11.5.4}\\
& \log_a (1+x)=\frac{x}{\ln a}+\underset{x\to
0}{\bold{o}(x)}\label{11.5.5}\\
& e^x=1+x+\underset{x\to
0}{\bold{o}(x)}\label{11.5.6}\\
&  a^x=1+x\cdot \ln a+\underset{x\to
0}{\bold{o}(x)}\label{11.5.7}\\
& (1+x)^\alpha=1+\alpha x+\underset{x\to
0}{\bold{o}(x)}\label{11.5.8}\\
& \arcsin x= x+\underset{x\to
0}{\bold{o}(x)}\label{11.5.9}\\
& \arctg x= x+\underset{x\to 0}{\bold{o}(x)}\label{11.5.10}\end{align}

\begin{proof} Все эти формулы легко выводятся из \eqref{11.1.1} -- \eqref{11.1.10}.
В качестве примера рассмотрим \eqref{11.5.7}:
 \begin{multline*}
 \overbrace{a^x=1+x\cdot \ln a+\underset{x\to 0}{\bold{o}(x)}}^{
 \footnotesize
 \begin{matrix}
 \text{формула \eqref{11.5.7}}
 \\
 \downarrow
 \end{matrix}
 }
 \quad \Leftrightarrow \\
\Leftrightarrow \quad a^x-1-x\cdot \ln a=\underset{x\to
0}{\bold{o}(x)}\quad \Leftrightarrow
\\
\Leftrightarrow \quad \frac{a^x-1-x\cdot \ln a}{x}\underset{x\to
0}{\longrightarrow} 0 \quad \Leftrightarrow
\\
\Leftrightarrow \quad \frac{a^x-1}{x}-\ln a\underset{x\to
0}{\longrightarrow} 0 \quad \Leftrightarrow
\\
\Leftrightarrow\quad \frac{a^x-1}{x}\underset{x\to
0}{\longrightarrow}\ln a \quad \Leftrightarrow
\\
\Leftrightarrow \quad \frac{a^x-1}{x\cdot \ln a}\underset{x\to
0}{\longrightarrow} 1 \quad \Leftrightarrow \quad
 \underbrace{a^x-1\underset{x\to 0}{\sim} x\cdot \ln a}_{
 \footnotesize
 \begin{matrix}
 \uparrow \\
 \text{формула \eqref{11.1.7}}
 \end{matrix}
 }
 \end{multline*}\end{proof}

\end{multicols}\noindent\rule[10pt]{160mm}{0.1pt}

\paragraph[Связь между символами $\sim$ и $\bold{o}$]{Связь между
асимптотической эквивалентностью \newline и асимптотическим сравнением}

Связи между формулами \eqref{11.5.1} -- \eqref{11.5.10} и
\eqref{11.1.1} -- \eqref{11.1.10} отражают важную общую зависимость
между символами $\sim$ и $\bold{o}$:

\begin{tm}[\bf о связи между символами $\sim$ и $\bold{o}$]\label{TH:svyaz-<<-i-o}
Следующие утверждения эквивалентны:
 \bit{
\item[(i)] $f(x) \underset{x\to a}{\sim} g(x)$;
\item[(ii)] $f(x) =g(x)+\underset{x\to a}{\bold{o}}(g(x))$;
\item[(ii)] $g(x) =f(x)+\underset{x\to a}{\bold{o}}(f(x))$.
 }\eit
\end{tm}\begin{proof} Докажем, что $(ii)\Leftrightarrow
(i)$:

\begin{multline*}
(ii) \quad \Leftrightarrow \quad f(x) =g(x)+\underset{x\to
a}{\bold{o}}(g(x)) \quad \Leftrightarrow \quad
f(x)-g(x)=\underset{x\to a}{\bold{o}}(g(x)) \quad \Leftrightarrow
\quad
\\
\quad \Leftrightarrow \quad \frac{f(x)-g(x)}{g(x)}\underset{x\to
a}{\longrightarrow} 0 \quad \Leftrightarrow \quad
\frac{f(x)}{g(x)}-1\underset{x\to a}{\longrightarrow} 0 \quad
\Leftrightarrow \quad \frac{f(x)}{g(x)}\underset{x\to
a}{\longrightarrow} 1 \quad \Leftrightarrow \quad
\\
\quad \Leftrightarrow \quad f(x) \underset{x\to a}{\sim} g(x) \quad
\Leftrightarrow \quad (i) \end{multline*} А теперь -- что
$(iii)\Leftrightarrow (i)$:

\begin{multline*}
(iii) \quad \Leftrightarrow \quad g(x) =f(x)+\underset{x\to
a}{\bold{o}}(f(x)) \quad \Leftrightarrow \quad
g(x)-f(x)=\underset{x\to a}{\bold{o}}(f(x)) \quad \Leftrightarrow
\quad
\\
\quad \Leftrightarrow \quad \frac{g(x)-f(x)}{f(x)}\underset{x\to
a}{\longrightarrow} 0 \quad \Leftrightarrow \quad
\frac{g(x)}{f(x)}-1\underset{x\to a}{\longrightarrow} 0 \quad
\Leftrightarrow \quad \frac{g(x)}{f(x)}\underset{x\to
a}{\longrightarrow} 1 \quad \Leftrightarrow \quad
\\
\quad \Leftrightarrow \quad \frac{f(x)}{g(x)}\underset{x\to
a}{\longrightarrow} 1 \quad \Leftrightarrow \quad f(x)
\underset{x\to a}{\sim} g(x) \quad \Leftrightarrow \quad (i) \quad
 \end{multline*}\end{proof}

\paragraph{Свойства символа $\bold{o}$ и упрощение асимптотических формул}

Перечислим некоторые

\bigskip

\centerline{\bf Общие свойства символа $\bold{o}$}
 \bit{
 \item[$1^\bold{o}$]
 $C\cdot \underset{x\to a}{\bold{o}}\Big( f(x)\Big)=
\underset{x\to a}{\bold{o}}\Big( C\cdot f(x)\Big)= \underset{x\to
a}{\bold{o}}\Big( f(x)\Big),\quad C\ne 0$;

 \item[$2^\bold{o}$]
$\underset{x\to a}{\bold{o}}\Big( f(x)\Big) \pm \underset{x\to
a}{\bold{o}}\Big( f(x)\Big) = \underset{x\to a}{\bold{o}}\Big(
f(x)\Big) $;

 \item[$3^\bold{o}$]
$f(x)\underset{x\to a}{\sim} g(x) \quad \Rightarrow \quad
\underset{x\to a}{\bold{o}}\Big( f(x)\Big)= \underset{x\to
a}{\bold{o}}\Big( g(x)\Big) $;

 \item[$4^\bold{o}$]
$\underset{x\to a}{\bold{o}}\Big( f(x)\Big) \cdot \underset{x\to
a}{\bold{o}}\Big( g(x)\Big) = \underset{x\to a}{\bold{o}}\Big(
f(x)\cdot g(x)\Big) = f(x)\cdot \underset{x\to a}{\bold{o}}\Big(
g(x)\Big) $;

 \item[$5^\bold{o}$]
$\underset{x\to a}{\bold{o}}\Big(\underset{x\to a}{\bold{o}}\Big(
f(x)\Big) \Big) = \underset{x\to a}{\bold{o}}\Big( f(x)\Big) $.

 \item[$6^\bold{o}$]
 $\underset{x\to a}{\bold{o}}(f(x))= f(x)\cdot \underset{x\to
a}{\bold{o}}(1)$;

 \item[$7^\bold{o}$]
$\underset{x\to a}{\bold{o}}(1) \underset{x\to a}{\longrightarrow}
0$;

 \item[$8^\bold{o}$]
$\Big( f(x)+\underset{x\to a}{\bold{o}}\Big( f(x)\Big) \Big) ^\mu=
f(x)^\mu+\underset{x\to a}{\bold{o}}\Big( f(x)^\mu \Big) , \quad
\mu\ne 0$. }\eit

Эти формулы надо понимать вот как. Если Вы видите где-нибудь,
например, выражение
$$
C\cdot \underset{x\to a}{\bold{o}}(f(x)),
$$
то, по формуле $1^{\bold{o}}$, Вы можете заменить его на
$$
\underset{x\to a}{\bold{o}}(C\cdot f(x)),
$$
или на
$$
\underset{x\to a}{\bold{o}}(f(x)).
$$
Или, скажем, выражение
$$
\underset{x\to a}{\bold{o}}\Big( f(x)\Big) + \underset{x\to
a}{\bold{o}}\Big( f(x)\Big)
$$
можно, пользуясь формулой $2^{\bold{o}}$, заменять на выражение
$$
\underset{x\to a}{\bold{o}}\Big( f(x)\Big) .
$$

\begin{proof}
1. С одной стороны, $\alpha(x)=C\cdot \underset{x\to a}{\bold{o}}(f(x)) \quad
\Leftrightarrow \quad \frac{\alpha(x)}{C}=\underset{x\to a}{\bold{o}}(f(x))
\quad \Leftrightarrow \quad \frac{\frac{\alpha(x)}{C}}{f(x)}\underset{x\to
a}{\longrightarrow} 0 \quad \Leftrightarrow \quad \frac{\alpha(x)}{C\cdot
f(x)}\underset{x\to a}{\longrightarrow} 0 \quad \Leftrightarrow \quad
\alpha(x)=\underset{x\to a}{\bold{o}}(C\cdot f(x))$. С другой стороны,
$\alpha(x)=\underset{x\to a}{\bold{o}}(C\cdot f(x)) \quad \Leftrightarrow \quad
\frac{\alpha(x)}{C\cdot f(x)}\underset{x\to a}{\longrightarrow} 0 \quad
\Leftrightarrow \quad \frac{\alpha(x)}{f(x)}\underset{x\to a}{\longrightarrow}
0 \quad \Leftrightarrow \quad \alpha(x)=\underset{x\to a}{\bold{o}}(f(x))$.

2. $\begin{cases}\alpha(x)=\underset{x\to a}{\bold{o}}\Big( f(x)\Big)
\\ \beta(x)=\underset{x\to a}{\bold{o}}\Big( f(x)\Big)\end{cases}$ $\Leftrightarrow$
$\begin{cases} \frac{\alpha(x)}{f(x)}\underset{x\to a}{\longrightarrow} 0
\\ \frac{\beta(x)}{f(x)}\underset{x\to a}{\longrightarrow} 0\end{cases}$ $\Rightarrow$
$\frac{\alpha(x) \pm \beta(x)}{f(x)}= \frac{\alpha(x)}{f(x)}\pm
\frac{\beta(x)}{f(x)}\underset{x\to a}{\longrightarrow} 0$ $\Leftrightarrow$
$\alpha(x)\pm \beta(x)=\underset{x\to a}{\bold{o}}\Big( f(x)\Big)$.

3. Пусть $f(x)\underset{x\to a}{\sim} g(x)$, то есть
$\frac{f(x)}{g(x)}\underset{x\to a}{\longrightarrow} 1$. Тогда
$\alpha(x)=\underset{x\to a}{\bold{o}}\Big( f(x) \Big)$ $\Leftrightarrow$
$\frac{\alpha(x)}{f(x)}\underset{x\to a}{\longrightarrow} 0$ $\Leftrightarrow$
$\frac{\alpha(x)}{g(x)}=\frac{\alpha(x)}{f(x)}\cdot
\frac{f(x)}{g(x)}\underset{x\to a}{\longrightarrow} 0\cdot 1=0$
$\Leftrightarrow$ $\alpha(x)=\underset{x\to a}{\bold{o}}\Big( g(x) \Big)$.

4. Во-первых, $\begin{cases}\alpha=\underset{x\to a}{\bold{o}}\Big( f(x)\Big)
\\ \beta=\underset{x\to a}{\bold{o}}\Big( g(x)\Big)\end{cases}$ $\Rightarrow$
$\begin{cases}\frac{\alpha(x)}{f(x)}\underset{x\to a}{\longrightarrow} 0 \\
\frac{\beta(x)}{g(x)} \underset{x\to a}{\longrightarrow} 0\end{cases}$
$\Rightarrow$ $\frac{\alpha(x)\cdot \beta(x)}{f(x)\cdot g(x)}=
 \frac{\alpha(x)}{f(x)}\cdot \frac{\beta(x)}{g(x)}
 \underset{x\to a}{\longrightarrow} 0$ $\Rightarrow$ $\alpha(x)\cdot \beta(x)= \underset{x\to a}{\bold{o}}\Big(
 f(x)\cdot g(x)\Big)$. И, во-вторых, $\alpha(x)=\underset{x\to a}{\bold{o}}\Big( f(x)\cdot g(x)\Big)
 \quad \Rightarrow \quad \frac{\alpha(x)}{f(x)\cdot g(x)}
 \underset{x\to a}{\longrightarrow} 0 \quad \Rightarrow \quad
 \frac{\frac{\alpha(x)}{f(x)}}{g(x)}\underset{x\to a}{\longrightarrow} 0
 \quad \Rightarrow \quad
 \frac{\alpha(x)}{f(x)}=\underset{x\to a}{\bold{o}}\Big(  g(x)\Big)
 \quad \Rightarrow \quad \alpha(x)=f(x)\cdot \underset{x\to
 a}{\bold{o}}\Big(  g(x)\Big)$.

5. Формула $\alpha(x)=\underset{x\to a}{\bold{o}}
 \Big(\underset{x\to a}{\bold{o}}\Big( f(x)\Big) \Big)$ означает, что
 $\alpha=\underset{x\to a}{\bold{o}}(\beta(x))$ для некоторого
 $\beta(x)=\underset{x\to a}{\bold{o}}\Big( f(x)\Big)$. Тогда: $\frac{\alpha(x)}{\beta(x)}\underset{x\to
 a}{\longrightarrow} 0 \quad\&\quad \frac{\beta(x)}{f(x)}
 \underset{x\to a}{\longrightarrow} 0 $ $\Rightarrow$ $\frac{\alpha(x)}{f(x)}= \frac{\alpha(x)}{\beta(x)}\cdot
 \frac{\beta(x)}{f(x)}\underset{x\to a}{\longrightarrow} 0$ $\Rightarrow$
 $\alpha(x)=\underset{x\to a}{\bold{o}}\Big(f(x)\Big)$.

6. $\alpha(x)=\underset{x\to a}{\bold{o}}(f(x)) \quad
 \Leftrightarrow \quad \frac{\alpha(x)}{f(x)}\underset{x\to a}{\longrightarrow} 0
  \quad \Leftrightarrow \quad
 \frac{\frac{\alpha(x)}{f(x)}}{1}\underset{x\to a}{\longrightarrow} 0
 \quad \Leftrightarrow \quad
 \frac{\alpha(x)}{f(x)}=\underset{x\to a}{\bold{o}}(1) \quad
 \Leftrightarrow \quad \alpha(x)=f(x)\cdot \underset{x\to a}{\bold{o}}(1)$.

7. $\alpha(x)=\underset{x\to a}{\bold{o}}(1)$ $\Leftrightarrow$
$\alpha(x)=\frac{\alpha(x)}{1}\underset{x\to a}{\longrightarrow} 0$.

8. $\alpha(x)= \Big( f(x)+\underset{x\to a}{\bold{o}}\Big(
 f(x)\Big) \Big) ^\mu$ $\Rightarrow$ $\alpha(x)^\frac{1}{\mu}=f(x)+\underset{x\to a}{\bold{o}}\Big(
 f(x)\Big)$ $\Rightarrow$ $\alpha(x)^\frac{1}{\mu}-f(x)=\underset{x\to a}{\bold{o}}\Big(
 f(x)\Big)$ $\Rightarrow$ $\frac{\alpha(x)^\frac{1}{\mu}-f(x)}{f(x)}\underset{x\to
 a}{\longrightarrow} 0$ $\Rightarrow$ $\frac{\alpha(x)^\frac{1}{\mu}}{f(x)}\underset{x\to
 a}{\longrightarrow} 1$ $\Rightarrow$ $\frac{\alpha(x)}{f(x)^\mu}\underset{x\to a}{\longrightarrow}
 1^\mu=1$ $\Rightarrow$ $\frac{\alpha(x)-f(x)^\mu}{f(x)^\mu}= \frac{\alpha(x)}{f(x)^\mu}-1
 \underset{x\to a}{\longrightarrow} 0$ $\Rightarrow$ $\alpha(x)-f(x)^\mu =\underset{x\to a}{\bold{o}}\Big( f(x)^\mu
 \Big)$ $\Rightarrow$ $\alpha(x)=f(x)^\mu+\underset{x\to a}{\bold{o}}\Big( f(x)^\mu
 \Big)$.
 \end{proof}

К свойствам $1^{\bold{o}} - 8^{\bold{o}}$ полезно добавить еще
три, относящихся специально к случаю, когда $x\to 0$

\bigskip

\centerline{\bf Свойства символа $\bold{o}$ при $x\to
0$}

 \bit{
\item[$1^{\bold{o}\bold{o}}$.] $\alpha<\beta \quad \Longrightarrow \quad
x^\beta=\underset{x\to 0}{\bold{o}}\Big( x^\alpha \Big)$;

\item[$2^{\bold{o}\bold{o}}$.] $\alpha<\beta \quad \Longrightarrow \quad
\underset{x\to 0}{\bold{o}}\Big( x^\alpha \Big)+x^\beta= \underset{x\to
0}{\bold{o}}\Big( x^\alpha \Big)$;

\item[$3^{\bold{o}\bold{o}}$.] $\alpha<\beta \quad \Longrightarrow \quad
\underset{x\to 0}{\bold{o}}\Big( x^\alpha \Big)+ \underset{x\to
0}{\bold{o}}\Big( x^\beta \Big)= \underset{x\to 0}{\bold{o}}\Big( x^\alpha
\Big)$;
 }\eit

\begin{proof} 1. Если $\alpha<\beta$, то по теореме \ref{x^a<<x^b},
$x^\beta\underset{x\to 0}{\ll} x^\alpha$, то есть
$x^\beta=\underset{x\to 0}{\bold{o}}\Big( x^\alpha \Big)$.

2. Отсюда следует, что если $\alpha<\beta$, то $ \underset{x\to
0}{\bold{o}}\Big( x^\alpha \Big)+x^\beta= \underset{x\to
0}{\bold{o}}\Big( x^\alpha \Big)+ \underset{x\to 0}{\bold{o}}\Big(
x^\alpha \Big)= \underset{x\to 0}{\bold{o}}\Big( x^\alpha \Big)$.

3. Если $\alpha<\beta$, то $ \underset{x\to 0}{\bold{o}}\Big(
x^\alpha \Big)+ \underset{x\to 0}{\bold{o}}\Big( x^\beta \Big)=
\underset{x\to 0}{\bold{o}}\Big( x^\alpha \Big)+
x^{\beta-\alpha}\cdot\underset{x\to 0}{\bold{o}}\Big( x^\alpha
\Big)= \underset{x\to 0}{\bold{o}}\Big( x^\alpha \Big)+
\underset{x\to 0}{\bold{o}}\Big( 1 \Big)\cdot \underset{x\to
0}{\bold{o}}\Big( x^\alpha \Big)= \underset{x\to 0}{\bold{o}}\Big(
x^\alpha \Big)+ \underset{x\to 0}{\bold{o}}\Big( 1 \cdot x^\alpha
\Big)= \underset{x\to 0}{\bold{o}}\Big( x^\alpha \Big)+
\underset{x\to 0}{\bold{o}}\Big( x^\alpha \Big)= \underset{x\to
0}{\bold{o}}\Big( x^\alpha \Big)$
 \end{proof}

Покажем, как с помощью свойств $1^{\bold{o}} - 8^{\bold{o}}$ и
$1^{\bold{o}\bold{o}} - 3^{\bold{o}\bold{o}}$ можно упрощать
асимптотические формулы.

\noindent\rule{160mm}{0.1pt}\begin{multicols}{2}

\begin{ex}

\begin{multline*}
x+x^2+\underset{x\to 0}{\bold{o}}( x^2)+ \underset{x\to
0}{\bold{o}}\Big( x+x^2+\underset{x\to 0}{\bold{o}}(
x^2)\Big)=\\={\smsize (\text{свойство}\, 3^{\bold{o}})}=
x+x^2+\underset{x\to 0}{\bold{o}}( x^2)+\\+ \underset{x\to
0}{\bold{o}}\Big( x+ x^2\kern-10pt{\Big/\kern-8pt\Big\backslash}+
\underset{x\to 0}{\bold{o}}(x^2)
\kern-18pt{\Big/\kern-8pt\Big\backslash}\kern10pt \Big)=\\=
x+x^2+\underset{x\to 0}{\bold{o}}( x^2)+ \underset{x\to
0}{\bold{o}}(x)= {\smsize (\text{свойство}\,
3^{\bold{o}\bold{o}})}=\\= x+x^2+ \underset{x\to 0}{\bold{o}}( x^2)
\kern-18pt{\Big/\kern-8pt\Big\backslash}\kern10pt + \underset{x\to
0}{\bold{o}}(x)= x+x^2+\underset{x\to 0}{\bold{o}}(x)=\\={\smsize
(\text{свойство}\, 2^{\bold{o}\bold{o}})}= x+ x^2
\kern-10pt{\Big/\kern-8pt\Big\backslash}+ \underset{x\to
0}{\bold{o}}(x)= x+\underset{x\to 0}{\bold{o}}(x)
\end{multline*}\end{ex}

\begin{ex}\begin{multline*}
\l x+x^2+\underset{x\to 0}{\bold{o}}( x^2)\r^2= \\
=\l x+x^2+\underset{x\to 0}{\bold{o}}( x^2)\r \cdot \l
x+x^2+\underset{x\to 0}{\bold{o}}( x^2)\r=\\
{\smsize\text{$=x^2+x^4+\l \underset{x\to 0}{\bold{o}}( x^2)\r^2+
2x^3+x\underset{x\to
0}{\bold{o}}( x^2)+ x^2\underset{x\to 0}{\bold{o}}( x^2)=$}}\\
=x^2+x^4+\underset{x\to 0}{\bold{o}}\Big( x^4\Big)+
2x^3+\underset{x\to 0}{\bold{o}}\Big( x^3\Big)+ \underset{x\to
0}{\bold{o}}\Big( x^4\Big)=\\={\smsize (\text{свойство}\,
3^{\bold{o}\bold{o}})}=\\
=x^2+x^4+\underset{x\to 0}{\bold{o}}\Big( x^4\Big)
\kern-20pt{\Big/\kern-8pt\Big\backslash}\kern10pt +
2x^3+\underset{x\to 0}{\bold{o}}\Big( x^3\Big) + \underset{x\to
0}{\bold{o}}\Big( x^4\Big)
\kern-20pt{\Big/\kern-8pt\Big\backslash}\kern10pt=\\
=x^2+x^4+2x^3+\underset{x\to 0}{\bold{o}}\Big( x^3\Big)={\smsize
(\text{свойство}\, 2^{\bold{o}\bold{o}})}=\\
=x^2+x^4\kern-10pt{\Big/\kern-8pt\Big\backslash} +2x^3+\underset{x\to
0}{\bold{o}}\Big( x^3\Big)= x^2+2x^3+\underset{x\to 0}{\bold{o}}\Big(
x^3\Big) \end{multline*}\end{ex}

\end{multicols}\noindent\rule[10pt]{160mm}{0.1pt}

Последний пример особенно важен, и поэтому будет полезно найти
более экономный способ упрощать подобные выражения. Представим
себе, что нам нужно упростить асимптотическое выражение
\begin{equation}\l A_0x^l+A_1x^{l+1}+...+A_kx^{l+k}+ \underset{x\to
0}{\bold{o}}\Big( x^{l+k}\Big)\r\cdot \l
B_0x^m+A_1x^{m+1}+...+A_nx^{m+n}+ \underset{x\to 0}{\bold{o}}\Big(
x^{m+n}\Big)\r \label{11.7.1}\end{equation} Для этого можно
пользоваться следующим алгоритмом.

\bigskip

\centerline{\bf Алгоритм экономного раскрытия скобок}\centerline{\bf
при упрощении асимптотическихвыражений:}

\bigskip

1. Мысленно раскрываем скобки, и находим сначала {\it минимальную
степень $M$ переменной $x$, которая может получиться под символом}
$\bold{o}$. В нашем случае, если $A_0\ne 0, \, B_0\ne 0$,
получается
$$
M=\min \{ l+m+n, l+k+m\}.
$$
Это будет означать, что после упрощения наше выражение должно
оканчиваться символом
$$
\underset{x\to 0}{\bold{o}}\Big( x^M\Big)
$$

2. По свойствам $2^{\bold{o}\bold{o}}$ и $3^{\bold{o}\bold{o}}$, все слагаемые
со степенями $i>M$ при упрощении ``съедаются'' слагаемым $\underset{x\to
0}{\bold{o}}\Big( x^M\Big)$:
$$
 x^i\kern-10pt{\Big/\kern-8pt\Big\backslash}
+ \underset{x\to 0}{\bold{o}}\Big( x^M\Big)= \underset{x\to
0}{\bold{o}}\Big( x^M\Big)
$$
Поэтому, раскрывая скобки, мы должны лишь найти коэффициенты перед степенями
$x^i, \, i\le M$. Это делается так же, как при обычном умножении многочленов.

\bigskip

Покажем, как этот алгоритм действует на конкретных примерах.

\noindent\rule{160mm}{0.1pt}\begin{multicols}{2}

\begin{ex} Рассмотрим еще раз выражение из примера 7.2:
 \begin{multline*}
\l x+x^2+\underset{x\to 0}{\bold{o}}( x^2)\r^2=\\= \l
x+x^2+\underset{x\to 0}{\bold{o}}( x^2)\r \cdot \l
x+x^2+\underset{x\to 0}{\bold{o}}( x^2)\r
 \end{multline*}

1. Сначала замечаем, что минимальная степень переменной $x$,
стоящая под символом $\bold{o}$, которая может получиться при
раскрытии скобок равна 3. Это получается, когда слагаемое $x$
умножается на слагаемое $\underset{x\to 0}{\bold{o}}( x^2)$.
Значит, наше упрощенное выражение должно оканчиваться на
$$
\underset{x\to 0}{\bold{o}}( x^3)
$$

2. Выписываем коэффициенты перед степенями $x^i, \, i\le 3$,
которые получаются при раскрытии скобок:
\begin{align*}
x^0 & :  0 & {\smsize \begin{pmatrix}\text{при раскрытии скобок};\\
\text{степень $x^0$ вообще}\\ \text{получиться не может}\end{pmatrix}}; \\
x^1 & : 0 &
{\smsize \begin{pmatrix}\text{степень $x^1$}\\
\text{тоже получиться не может}\end{pmatrix}}; \\
x^2 & : 1 & {\smsize \begin{pmatrix}\text{степень $x^2$ может
получиться,}\\ \text{когда $x$ умножается на $x$}\end{pmatrix}}; \\
x^3 & : 1+1=2 & {\smsize \begin{pmatrix}\text{степень $x^3$ может
получиться,}\\ \text{когда $x$ умножается на $x^2$,}\\ \text{и когда
наоборот}\\ \text{$x^2$ умножается на $x$}\end{pmatrix}}.
\end{align*}

3. Выписываем результат:
$$
\l x+x^2+\underset{x\to 0}{\bold{o}}( x^2)\r^2=
x^2+2x^3+\underset{x\to 0}{\bold{o}}\Big( x^3\Big)
$$
\end{ex}

\begin{ex} Упростим выражение, отличающееся от 7.2 степенью:
 \begin{multline*}
\l x+x^2+\underset{x\to 0}{\bold{o}}( x^2)\r^3=\\= \l
x+x^2+\underset{x\to 0}{\bold{o}}( x^2)\r \cdot \l
x+x^2+\underset{x\to 0}{\bold{o}}( x^2)\r \cdot\\ \cdot \l
x+x^2+\underset{x\to 0}{\bold{o}}( x^2)\r
 \end{multline*}

1. Сначала замечаем, что минимальная степень переменной $x$,
стоящая под символом $\bold{o}$, которая может получиться при
раскрытии скобок равна 4. Это получается, когда два слагаемых $x$
умножаются на слагаемое $\underset{x\to 0}{\bold{o}}( x^2)$.
Значит, наше упрощенное выражение должно оканчиваться на
$$
\underset{x\to 0}{\bold{o}}\Big( x^4 \Big)
$$

2. Выписываем коэффициенты перед степенями $x^i, \, i\le 4$,
которые получаются при раскрытии скобок:
\begin{align*}
x^0 & : 0 & {\smsize \begin{pmatrix}\text{при раскрытии скобок}\\
\text{степень $x^0$ вообще}\\ \text{получиться не может}\end{pmatrix}}; \\
x^1 & : 0  & {\smsize \begin{pmatrix}\text{степень $x^1$ тоже}\\
\text{получиться не может}\end{pmatrix}}; \\
x^2 & : 0  & {\smsize \begin{pmatrix}\text{степень $x^2$ тоже}\\
\text{получиться не может}\end{pmatrix}}; \\
x^3 & : 1  & {\smsize
\begin{pmatrix}\text{степень $x^3$ может получиться,}\\
\text{когда $x$ умножается}\\
\text{на $x$ и еще раз на $x$}\end{pmatrix}}; \\
x^4 & : 1+1+1=3  & {\smsize \begin{pmatrix}\text{степень $x^3$
может}\\ \text{получиться в трех случаях,}\\ \text{когда в разных
вариантах $x$}\\ \text{умножается на $x$ и на $x^2$}\end{pmatrix}}.
\end{align*}

3. Выписыаем результат:
$$
\l x+x^2+\underset{x\to 0}{\bold{o}}( x^2)\r^3=
x^3+3x^4+\underset{x\to 0}{\bold{o}}\Big( x^4\Big)
$$
\end{ex}

\begin{ex} Упростим выражение
$$
\l x^2-\frac{1}{2}x^3+\underset{x\to 0}{\bold{o}}( x^3)\r \cdot \l
2x^3+x^5+\underset{x\to 0}{\bold{o}}( x^6)\r
$$

1. Замечаем, что минимальная степень переменной $x$, стоящая под
символом $\bold{o}$, которая может получиться при раскрытии скобок
равна 6. Это получается, когда слагаемое $\underset{x\to
0}{\bold{o}}( x^3)$ умножается на слагаемое $2x^3$. Значит, наше
упрощенное выражение должно оканчиваться на
$$
\underset{x\to 0}{\bold{o}}\Big( x^6 \Big)
$$

2. Выписываем коэффициенты перед степенями $x^i, \, i\le 6$,
которые получаются при раскрытии скобок:
\begin{align*}
x^0 & : 0  & {\smsize \begin{pmatrix}
\text{при раскрытии скобок степень $x^0$}\\
\text{вообще получиться не может}\end{pmatrix}}; \\
x^1 & : 0  &
{\smsize (-''-);} \\
x^2 & : 0 &
{\smsize (-''-);} \\
x^3 & : 0 &
{\smsize (-''-);} \\
x^4 & : 0 &
{\smsize (-''-);} \\
x^5 & : 2  & {\smsize \begin{pmatrix}\text{степень $x^5$ может
получиться,}\\ \text{когда $x^2$ умножается на $2x^3$}\end{pmatrix}}; \\
x^6 & : -1  & {\smsize \begin{pmatrix}\text{степень $x^6$ может
получиться,}\\ \text{когда $-\frac{1}{2}x^3$ умножается на
$2x^3$}\end{pmatrix}}.
\end{align*}

3. Выписываем результат:
 \begin{multline*}
\l x^2-\frac{1}{2}x^3+\underset{x\to 0}{\bold{o}}( x^3)\r \cdot \l
2x^3+x^5+\underset{x\to 0}{\bold{o}}( x^6)\r=\\=
2x^5-x^6+\underset{x\to 0}{\bold{o}}( x^6)
 \end{multline*}
\end{ex}

\begin{ers} Упростите выражения:
 \biter{\smsize
\item[1.] $\l -x^2+\frac{2}{3} x^3+\underset{x\to 0}{\bold{o}}(
x^3)\r^2$

\item[2.] $\l -x^2+\frac{2}{3} x^3+\underset{x\to 0}{\bold{o}}(
x^3)\r^3$

\item[3.] $\l 3x+\frac{1}{5} x^3+\underset{x\to 0}{\bold{o}}(
x^4)\r^2$

\item[4.] $\l 3x+\frac{1}{5} x^3+\underset{x\to 0}{\bold{o}}(
x^4)\r^3$

\item[5.] $ \l 3x^2+\frac{1}{4}x^3+\underset{x\to 0}{\bold{o}}(
x^3)\r \cdot \l -2x+x^2+\underset{x\to 0}{\bold{o}}( x^2)\r $

\item[6.] $ \l 3x^2+\frac{1}{4}x^3+\underset{x\to 0}{\bold{o}}(
x^3)\r \cdot \l -2x+x^2+\underset{x\to 0}{\bold{o}}( x^3)\r $

\item[7.] $ \l 3x^2+\frac{1}{4}x^3+\underset{x\to 0}{\bold{o}}(
x^4)\r \cdot \l -2x+x^2+\underset{x\to 0}{\bold{o}}( x^3)\r $

\item[8.] $ \l 4x+3x^2-2x^3+x^4+\underset{x\to 0}{\bold{o}}( x^4)\r
\cdot \\ \qquad\cdot\l -x+x^2+x^3-x^4+\underset{x\to 0}{\bold{o}}(
x^4)\r $

\item[9.] $ \l 4x+3x^2-2x^3+x^4+\underset{x\to 0}{\bold{o}}( x^5)\r
\cdot\\ \qquad\cdot \l -x+x^2+x^3-x^4+\underset{x\to 0}{\bold{o}}(
x^4)\r $

\item[10.] $ \l 4x+3x^2-2x^3+x^4+\underset{x\to 0}{\bold{o}}( x^5)\r
\cdot$\break $\qquad\qquad\cdot \l -x+x^2+x^3-x^4+\underset{x\to
0}{\bold{o}}( x^5)\r $
 }\eiter
\end{ers}

\end{multicols}\noindent\rule[10pt]{160mm}{0.1pt}

\paragraph{Вычисление пределов с помощью асимптотических формул}

Объясним теперь, зачем нужны асимптотические формулы. Рассмотрим
предел
\begin{equation}\lim_{x\to 0}\frac{x^6\cdot \arctg^4 x+\sin^5\sin\tg x^2}
{x^2\cdot \sin^8 x+\ln \cos 3x^5}\label{11.8.1}\end{equation} Для его
вычисления не подходит теорема об эквивалентной замене
\ref{ch-var-sym}, или принцип выделения главного слагаемого
\ref{main+princ} (не говоря уже о правиле Лопиталя
\ref{Lopital_x->a_0/0}, которое пришлось бы здесь применять 10 раз).

Так вот, в этой ситуации как раз бывают полезны асимптотические
формулы. Вместе с доказанными в предыдущем параграфе свойствами
$1^0 - 11^0$, они позволяют вычислять такие ``сложные'' пределы.
Покажем как это делается, начиная с более простых задач.

\noindent\rule{160mm}{0.1pt}\begin{multicols}{2}

\begin{ex}
 \begin{multline*}
\lim_{x\to 0}\frac{\arcsin 2x-\sin x}{x+\ln (1+3x)}=
{\smsize\begin{pmatrix}\text{используем формулы}\\
\text{\eqref{11.5.1} и \eqref{11.5.9}}\end{pmatrix}}=\\
=\lim_{x\to 0}\frac{ 2x+\underset{x\to 0}{\bold{o}}(2x)- \left(
x+\underset{x\to 0}{\bold{o}}(x) \right)} {x+\l 3x+\underset{x\to
0}{\bold{o}}(3x)\r}=\\
=\lim_{x\to 0}\frac{ x+\underset{x\to 0}{\bold{o}}(x)- \underset{x\to
0}{\bold{o}}(x)} {4x+\underset{x\to 0}{\bold{o}}(x)}= \lim_{x\to
0}\frac{ x+\underset{x\to 0}{\bold{o}}(x)} {4x+\underset{x\to
0}{\bold{o}}(x)}=\\= \lim_{x\to 0}\frac{ 1+\underset{x\to
0}{\bold{o}}(1)} {4+\underset{x\to 0}{\bold{o}}(1)}=
\frac{1+0}{4+0}=\frac{1}{4}\end{multline*}\end{ex}

\begin{ex}
 \begin{multline*}
\lim_{x\to 0}\frac{\ln(\cos x)}{2^{x^2}-1}=
{\smsize\begin{pmatrix}\text{используем формулы}\\
\text{\eqref{11.5.3} и \eqref{11.5.6}}\end{pmatrix}}=\\
=\lim_{x\to 0}\frac{\ln\l 1-\frac{x^2}{2}+\underset{x\to
0}{\bold{o}}(x^2)\r}{x^2\cdot \ln 2 +\underset{x\to
0}{\bold{o}}(x^2)}=\\
={\smsize\begin{pmatrix}\text{поскольку}\;
-\frac{x^2}{2}+\underset{x\to 0}{\bold{o}}(x^2)\underset{x\to
0}{\longrightarrow} 0, \\
\text{можно сделать замену}\\
y=-\frac{x^2}{2}+\underset{x\to 0}{\bold{o}}(x^2),
\\
\text{и тогда получится, что}\; y\underset{x\to 0}{\longrightarrow} 0, \\
\text{и по формуле \eqref{11.5.4} мы будем иметь}\\
\ln\l1-\frac{x^2}{2}+\underset{x\to 0}{\bold{o}}(x^2)\r=\\=
\ln(1+y)=y+\underset{y\to 0}{\bold{o}}(y)=\\=
-\frac{x^2}{2}+\underset{x\to 0}{\bold{o}}(x^2)+ \underset{x\to
0}{\bold{o}}\left(-\frac{x^2}{2}+\underset{x\to
0}{\bold{o}}(x^2)\right)
\end{pmatrix}}=\\
=\lim_{x\to 0}\frac{ -\frac{x^2}{2}+\underset{x\to 0}{\bold{o}}(x^2)+
\underset{x\to 0}{\bold{o}}\left(-\frac{x^2}{2}+\underset{x\to
0}{\bold{o}}(x^2)\right) }{x^2\cdot \ln 2 +\underset{x\to
0}{\bold{o}}(x^2)}=\\
={\smsize\begin{pmatrix}\text{используем}\\
\text{свойство}\; 3^{\bold{o}\bold{o}}\end{pmatrix}}=\\
=\lim_{x\to 0}\frac{ -\frac{x^2}{2}+\underset{x\to 0}{\bold{o}}(x^2)+
\underset{x\to 0}{\bold{o}}\left(-\frac{x^2}{2}+ \underset{x\to
0}{\bold{o}}(x^2) \kern-20pt{\Big/\kern-8pt\Big\backslash}\kern10pt
\right)}{x^2\cdot \ln 2 +\underset{x\to 0}{\bold{o}}(x^2)}=\\=
\lim_{x\to 0}\frac{ -\frac{x^2}{2}+\underset{x\to 0}{\bold{o}}(x^2)+
\underset{x\to 0}{\bold{o}}\left(-\frac{x^2}{2}\right)} {x^2\cdot \ln
2 +\underset{x\to 0}{\bold{o}}(x^2)}=\\=
{\smsize\begin{pmatrix}\text{используем}\\ \text{свойство}\;
1^{\bold{o}}\end{pmatrix}}= \lim_{x\to 0}\frac{
-\frac{x^2}{2}+\underset{x\to 0}{\bold{o}}(x^2)+ \underset{x\to
0}{\bold{o}}\left( x^2 \right)} {x^2\cdot \ln 2 +\underset{x\to
0}{\bold{o}}(x^2)}=\\= {\smsize\begin{pmatrix}\text{снова}\\
\text{используем}\; 2^0
\end{pmatrix}}=\lim_{x\to 0}\frac{ -\frac{x^2}{2}+\underset{x\to
0}{\bold{o}}(x^2)} {x^2\cdot \ln 2 +\underset{x\to
0}{\bold{o}}(x^2)}=\\= {\smsize\begin{pmatrix}\text{после этого}\\
\text{свойство}\; 6^{\bold{o}}\end{pmatrix}}= \lim_{x\to
0}\frac{-\frac{x^2}{2}+x^2\underset{x\to 0}{\bold{o}}(1)} {x^2\cdot
\ln 2 +x^2\underset{x\to
0}{\bold{o}}(1)}=\\= {\smsize\begin{pmatrix}\text{делим все}\\
\text{на}\; x^2
\end{pmatrix}}=
\lim_{x\to 0}\frac{-\frac{1}{2}+\underset{x\to 0}{\bold{o}}(1)} {\ln
2 +\underset{x\to 0}{\bold{o}}(1)}=\\=
{\smsize\begin{pmatrix}\text{применяем}\\  \text{свойство}\,
7^{\bold{o}}\end{pmatrix}}= \frac{-\frac{1}{2}+0}{\ln 2 +0}=
-\frac{1}{2\ln 2}\end{multline*}\end{ex}

\begin{ex}
 \begin{multline*}
\lim_{x\to 0} (x+e^x)^\frac{1}{x}= \lim_{x\to 0} e^{\ln
\left((x+e^x)^\frac{1}{x}\right)}=\\= e^{\lim_{x\to
0}\frac{1}{x}\cdot \ln (x+e^x)}= e^{\lim_{x\to 0}\frac{1}{x}\cdot \ln
(x+1+x+ \underset{x\to 0}{\bold{o}}(x))}=\\= e^{\lim_{x\to
0}\frac{1}{x}\cdot \ln (1+2x+ \underset{x\to 0}{\bold{o}}(x))}=\\=
e^{\lim_{x\to 0}\frac{1}{x}\cdot \left(2x+\underset{x\to
0}{\bold{o}}(x) +\underset{x\to 0}{\bold{o}}\left( 2x+\underset{x\to
0}{\bold{o}}(x) \right)\right)}=\\= e^{\lim_{x\to 0}\frac{1}{x}\cdot
\left(2x+\underset{x\to 0}{\bold{o}}(x) +\underset{x\to
0}{\bold{o}}\left( 2x+\underset{x\to 0}{\bold{o}}(2x)
\right)\right)}=\\= e^{\lim_{x\to 0}\frac{1}{x}\cdot
\left(2x+\underset{x\to 0}{\bold{o}}(x) +\underset{x\to
0}{\bold{o}}\left( 2x \right)\right)}=\\= e^{\lim_{x\to
0}\frac{1}{x}\cdot \left(2x+\underset{x\to 0}{\bold{o}}(x)
+\underset{x\to 0}{\bold{o}}\left( x \right)\right)}=\\=
e^{\lim_{x\to 0}\frac{1}{x}\cdot x\cdot \left(2+ \underset{x\to
0}{\bold{o}}(1) +\underset{x\to 0}{\bold{o}}\left( 1
\right)\right)}=\\= e^{\lim_{x\to 0}\left(2+ \underset{x\to
0}{\bold{o}}(1) +\underset{x\to 0}{\bold{o}}\left( 1
\right)\right)}=e^{2+0+0}=e^2
\end{multline*}\end{ex}

\begin{ex}
Вернемся теперь к примеру \eqref{11.8.1}:
 {\smsize
 \begin{multline*}
\lim_{x\to 0}\frac{x^6\cdot \arctg^4 x+\sin^5\sin\tg x^2} {x^2\cdot
\sin^8 x+\ln \cos 3x^5}=\\
=\lim_{x\to 0}\left\{x^6\cdot \left(x+\underset{x\to
0}{\bold{o}}(x)\right)^4+ \left(\sin\tg x^2+\underset{x\to
0}{\bold{o}}(\sin\tg x^2)\right)^5\right\} :\\:\left\{x^2\cdot
\left(x+\underset{x\to 0}{\bold{o}}(x)\right)^8+ \ln
\left(1-\frac{(3x^5)^2}{2}+ \underset{x\to 0}{\bold{o}}\left(
(3x^5)^2 \right)\right)\right\}=\\
=\lim_{x\to 0}\Bigg\{x^6\cdot \left(x^4+\underset{x\to
0}{\bold{o}}(x^4)\right)+\\+ \Bigg(\tg x^2+\underset{x\to
0}{\bold{o}}(\tg x^2) +\underset{x\to 0}{\bold{o}}\left(\tg
x^2+\underset{x\to 0}{\bold{o}}(\tg x^2) \right)\Bigg)^5\Bigg\}: \\
:\Bigg\{x^2\cdot \left(x^8+\underset{x\to 0}{\bold{o}}(x^8)\right)+
\ln \left(1-\frac{9x^{10}}{2}+ \underset{x\to 0}{\bold{o}}\left(
9x^{10}\right)\right)\Bigg\}=\\= \lim_{x\to 0}\left\{
x^{10}+\underset{x\to 0}{\bold{o}}(x^{10})+ \left(\tg
x^2+\underset{x\to 0}{\bold{o}}(\tg x^2) \right)^5\right\} : \\ :
\Bigg\{x^{10}+\underset{x\to 0}{\bold{o}}(x^{10}) -\frac{9x^{10}}{2}+
\underset{x\to 0}{\bold{o}}\left( 9x^{10}\right) +\\+\underset{x\to
0}{\bold{o}}\left( -\frac{9x^{10}}{2}+ \underset{x\to
0}{\bold{o}}\left( 9x^{10}\right) \right) \Bigg\}=\\
=\lim_{x\to 0}\Bigg\{ x^{10}+\underset{x\to 0}{\bold{o}}(x^{10})+\\+
\Big( x^2+\underset{x\to 0}{\bold{o}}(x^2) +\underset{x\to
0}{\bold{o}}\left(x^2+\underset{x\to 0}{\bold{o}}(x^2) \right)
\Big)^5\Bigg\}:\\:\Bigg\{x^{10}+\underset{x\to 0}{\bold{o}}(x^{10})
-\frac{9x^{10}}{2}+\underset{x\to 0}{\bold{o}}\left(
x^{10}\right)\Bigg\}=\\= \lim_{x\to 0}\frac{x^{10}+\underset{x\to
0}{\bold{o}}(x^{10})+ \left( x^2+\underset{x\to
0}{\bold{o}}(x^2)\right)^5}{-\frac{7x^{10}}{2}+\underset{x\to
0}{\bold{o}}\left( x^{10}\right)}=\\= \lim_{x\to
0}\frac{x^{10}+\underset{x\to 0}{\bold{o}}(x^{10})+
x^{10}+\underset{x\to 0}{\bold{o}}(x^{10})}{-\frac{7x^{10}}{2}+
\underset{x\to 0}{\bold{o}}\left( x^{10}\right)} =\\= \lim_{x\to
0}\frac{2x^{10}+\underset{x\to 0}{\bold{o}}(x^{10})}
{-\frac{7x^{10}}{2}+\underset{x\to 0}{\bold{o}}\left( x^{10}\right)}=
\lim_{x\to 0}\frac{2+\underset{x\to 0}{\bold{o}}(1)}
{-\frac{7}{2}+\underset{x\to 0}{\bold{o}}\left( 1 \right)}=\\=
\frac{2+0}{-\frac{7}{2}+0}=-\frac{4}{7}\end{multline*}
 }\end{ex}

\begin{ex}
 \begin{multline*}
\lim_{x\to 0}\frac{4\sqrt[3] {x}+3\sqrt[4] {x}-2\sin x} {1-\cos
x+\sqrt[7]{x^2}}=\\= {\smsize\begin{pmatrix}\text{применяем формулы}\\
\text{\eqref{11.5.1} и \eqref{11.5.3}}\end{pmatrix}}=\\= \lim_{x\to
0}\frac{4 x^\frac{1}{3}+3x^\frac{1}{4}- 2(x+\underset{x\to
0}{\bold{o}}(x))} {\frac{1}{2}x^2+\underset{x\to
0}{\bold{o}}(x^2)+x^\frac{2}{7}}=\\= {\smsize\begin{pmatrix}
x^\frac{1}{3} = \underset{x\to 0}{\bold{o}}(x^\frac{1}{4}) \\
x = \underset{x\to 0}{\bold{o}}(x^\frac{1}{4}) \\
x^2=\underset{x\to 0}{\bold{o}}(x^\frac{2}{7})
\end{pmatrix}}=\\= \lim_{x\to 0}\Big\{4 \underset{x\to
0}{\bold{o}}(x^\frac{1}{4}) +3x^\frac{1}{4}-\\- 2\big(\underset{x\to
0}{\bold{o}}(x^\frac{1}{4})+ \underset{x\to
0}{\bold{o}}(\underset{x\to 0}{\bold{o}}(x^\frac{1}{4}))\big)\Big\}
:\\:\Big\{\frac{1}{2}\underset{x\to 0}{\bold{o}}(x^\frac{2}{7})
+\underset{x\to 0}{\bold{o}} (\underset{x\to
0}{\bold{o}}(x^\frac{2}{7}))+x^\frac{2}{7}\Big\}=\\=
{\smsize\begin{pmatrix}\text{применяем}\\
\text{свойство $6^{\bold{o}}\; \S 7$}\end{pmatrix}}=\\= \lim_{x\to
0}\frac{4 \underset{x\to 0}{\bold{o}}(x^\frac{1}{4}) +3x^\frac{1}{4}-
2(\underset{x\to 0}{\bold{o}}(x^\frac{1}{4})+ \underset{x\to
0}{\bold{o}}(x^\frac{1}{4}))} {\frac{1}{2}\underset{x\to
0}{\bold{o}}(x^\frac{2}{7}) +\underset{x\to 0}{\bold{o}}
(x^\frac{2}{7})+x^\frac{2}{7}}=\\= {\smsize\begin{pmatrix}\text{применяем}\\
\text{свойства $1^0$ и $2^0 \; \S 7$}\end{pmatrix}}=\\= \lim_{x\to
0}\frac{\underset{x\to 0}{\bold{o}}(x^\frac{1}{4}) +3x^\frac{1}{4}-
2(\underset{x\to 0}{\bold{o}}(x^\frac{1}{4}))} {\underset{x\to
0}{\bold{o}}(x^\frac{2}{7})+x^\frac{2}{7}}=\\= {\smsize\begin{pmatrix}\text{применяем}\\
\text{свойства $1^0$ и $3^0 \; \S 7$}\end{pmatrix}}=\\= \lim_{x\to
0}\frac{\underset{x\to 0}{\bold{o}}(x^\frac{1}{4})
+3x^\frac{1}{4}-\underset{x\to 0}{\bold{o}}(x^\frac{1}{4})}
{\underset{x\to 0}{\bold{o}}(x^\frac{2}{7})+x^\frac{2}{7}}=\\=
\lim_{x\to 0}\frac{\underset{x\to 0}{\bold{o}}(x^\frac{1}{4})
+3x^\frac{1}{4}} {\underset{x\to
0}{\bold{o}}(x^\frac{2}{7})+x^\frac{2}{7}}=\\= \lim_{x\to
0}\frac{x^\frac{1}{4}\cdot \underset{x\to 0}{\bold{o}}(1)
+3x^\frac{1}{4}} {x^\frac{2}{7}\cdot \underset{x\to
0}{\bold{o}}(1)+x^\frac{2}{7}}=\\= \lim_{x\to 0}\frac{\underset{x\to
0}{\bold{o}}(1)+3} {x^\frac{1}{28}\cdot \underset{x\to
0}{\bold{o}}(1)+x^\frac{1}{28}}=\\
={\smsize\begin{pmatrix}\text{по свойству $8^0 \; \S 7$}\\
\underset{x\to 0}{\bold{o}}(1) \underset{x\to 0}{\longrightarrow}
0\end{pmatrix}}= \left(\frac{0+3}{0\cdot 0+0}\right)=\infty
\end{multline*}\end{ex}

\begin{ex}
 \begin{multline*}
\lim_{x\to \infty}\frac{4\sqrt[3] {x}+3\sqrt[4] {x}-2\sin x} {1-\cos
x+\sqrt[7]{x^2}}=\\= \lim_{x\to
\infty}\frac{4x^\frac{1}{3}+3x^\frac{1}{4}-2\sin x} {1-\cos
x+x^\frac{2}{7}}=\\= {\smsize\begin{pmatrix}
x^\frac{1}{4}=\underset{x\to \infty}{\bold{o}}(x^\frac{1}{3}) \\
\sin x=\underset{x\to \infty}{\bold{o}}(x^\frac{1}{3}) \\
\cos x=\underset{x\to \infty}{\bold{o}}(x^\frac{2}{7})\\
1=\underset{x\to \infty}{\bold{o}}(x^\frac{2}{7})
\end{pmatrix}}=\\=
\lim_{x\to \infty}\frac{4x^\frac{1}{3}+ 3\underset{x\to
\infty}{\bold{o}}(x^\frac{1}{3})- 2\underset{x\to
\infty}{\bold{o}}(x^\frac{1}{3})} {\underset{x\to
\infty}{\bold{o}}(x^\frac{2}{7})- \underset{x\to
\infty}{\bold{o}}(x^\frac{2}{7})+x^\frac{2}{7}}=\\
=\lim_{x\to \infty}\frac{4x^\frac{1}{3}+ \underset{x\to
\infty}{\bold{o}}(x^\frac{1}{3})} {\underset{x\to
\infty}{\bold{o}}(x^\frac{2}{7})+x^\frac{2}{7}}=\\= \lim_{x\to
\infty}\frac{4x^\frac{1}{3}+ x^\frac{1}{3}\cdot \underset{x\to
\infty}{\bold{o}}(1)} {x^\frac{2}{7}\cdot \underset{x\to
\infty}{\bold{o}}(1)+x^\frac{2}{7}}=\\= \lim_{x\to
\infty}\frac{4x^\frac{1}{21}+ x^\frac{1}{21}\cdot \underset{x\to
\infty}{\bold{o}}(1)} {\underset{x\to \infty}{\bold{o}}(1)+1}=\\=
{\smsize\begin{pmatrix}\text{по свойству $8^0 \; \S 7$}\\
\underset{x\to \infty}{\bold{o}}(1) \underset{x\to
\infty}{\longrightarrow} 0
\end{pmatrix}}=\\= \lim_{x\to \infty} x^\frac{1}{21}\cdot \frac{4+
\underset{x\to \infty}{\bold{o}}(1)} {\underset{x\to
\infty}{\bold{o}}(1)+1}=\\= \left(\infty \cdot
\frac{4+0}{0+1}\right)= \infty \end{multline*}\end{ex}

\begin{ers} Вычислите пределы:
 \bit{
\item[1.] $\lim\limits_{x\to 0}\frac{1-\cos x} {\ln (1+\tg^2 x)}$
\item[2.] $\lim\limits_{x\to 0}\frac{\ln (1+x+x^2)+\arcsin 3x-5x^3}
{\sin 2x+\tg^2 x+(e^x-1)^5}$ \item[3.] $\lim\limits_{x\to
0}\left(\cos x+\arctg^2 x \right)^\frac{1}{\arcsin^2 x}$
 }\eit
\end{ers}

\end{multicols}\noindent\rule[10pt]{160mm}{0.1pt}

\section{Формулы Пеано и асимптотика}\label{SEC:formuly-Peano-i-asimptotika}

\subsection{Формулы Пеано}

\paragraph{Формула Тейлора-Пеано.}

Асимптотические формулы \eqref{11.5.1} -- \eqref{11.5.10}, оказывается,
являются частными случаями одной общей формулы, которая описывается следующей
теоремой:

\begin{tm}[\bf Тейлора-Пеано]\label{TH:Taylor-Peano}
Пусть функция $f$ определена и дифференцируема порядка $n$ на интервале
$(a,b)$. Тогда для всякой точки $x_0\in (a,b)$ справедлива следующая
асимптотическая формула:
\medskip
 \begin{equation}\label{Taylor-Peano}
 \boxed{\quad
\begin{split}\phantom{\frac{\frac{\frac{\frac{1}{1}}{1}}{1}}{1}}
f(x) &= \sum_{k=0}^n \frac{f^{(k)}(x_0)}{k!}\cdot (x-x_0)^k+ \underset{x\to
x_0}{\bold{o}}\Big( (x-x_0)^n\Big)
\phantom{\frac{1}{\frac{1}{\frac{1}{1}}}}\end{split}\quad
 }\end{equation}\medskip
\end{tm}

\bit{\item[$\bullet$] Эта формула называется {\it формулой Тейлора с остаточным
членом в форме Пеано} \index{формула!Тейлора с остаточным членом в форме!Пеано}
{\it порядка} $n$ для функции $f$ в точке $x_0$.}\eit

\begin{proof}[Доказательство теоремы \ref{TH:Taylor-Peano}]
Рассмотрим многочлен Тейлора
$$
T_n(x)=f(x_0)+ \frac{f'(x_0)}{1!}\cdot (x-x_0)+ \frac{f''(x_0)}{2!}\cdot
(x-x_0)^2 +...+ \frac{f^{(n)}(x_0)}{n!}\cdot (x-x_0)^n= \sum_{k=0}^n
\frac{f^{(k)}(x_0)}{k!}\cdot (x-x_0)^k
$$
и заметим, что
\begin{equation}
T_n(x_0)=f(x_0), \quad T_n'(x_0)=f'(x_0), \quad T_n''(x_0)=f''(x_0), \quad ...
\quad T_n^{(n)}(x_0)=f^{(n)}(x_0) \label{11.9.2}\end{equation} Действительно,

\begin{multline*}
T_n(x_0)=f(x_0)+ \frac{f'(x_0)}{1!}\cdot (x_0-x_0)+ \frac{f''(x_0)}{2!}\cdot
(x_0-x_0)^2 +...+\frac{f^{(n)}(x_0)}{n!}\cdot (x_0-x_0)^n=\\=
f(x_0)+0+0+...+0=f(x_0) \end{multline*} Затем,
$$
T_n'(x)=0+ \frac{f'(x_0)}{1!}\cdot 1+ \frac{f''(x_0)}{2!}\cdot 2(x-x_0) +...+
\frac{f^{(n)}(x_0)}{n!}\cdot n (x-x_0)^{n-1}
$$
и поэтому
$$
T_n'(x_0)= f'(x_0)+ \frac{f''(x_0)}{2!}\cdot 2(x_0-x_0) +...+
\frac{f^{(n)}(x_0)}{n!}\cdot n (x_0-x_0)^{n-1}=
\frac{f'(x_0)}{1!}+0+...+0=f'(x_0)
$$
И так далее. Для последней производной мы получим
$$
T_n^{(n)}(x)=0+0+...+ \frac{f^{(n)}(x_0)}{n!}\cdot n!=f^{(n)}(x_0)
$$
и поэтому
$$
T_n^{(n)}(x_0)=f^{(n)}(x_0)
$$
Формулы \eqref{11.9.2} позволяют вычислить следующий предел:
\begin{multline*}\lim_{x\to x_0}\frac{f(x)-T_n(x)}{(x-x_0)^n}=
{\smsize\begin{pmatrix}\text{по первой формуле \eqref{11.9.2}}\\
\text{получаем неопределенность}\; \frac{0}{0},\\
\text{поэтому можно применить}\\ \text{правило Лопиталя}\end{pmatrix}}=
\lim_{x\to x_0}\frac{f'(x)-T_n'(x)}{n(x-x_0)^{n-1}}=\\=
{\smsize\begin{pmatrix}\text{по второй формуле \eqref{11.9.2}}\\
\text{получаем неопределенность}\; \frac{0}{0},  \\
\text{поэтому опять применяем}\\ \text{правило Лопиталя}\end{pmatrix}}=
\lim_{x\to x_0}\frac{f''(x)-T_n''(x)}{n (n-1)(x-x_0)^{n-2}}=\\=...=
{\smsize\begin{pmatrix}\text{по формулам \eqref{11.9.2}}\\
\text{всегда будет неопределенность}\; \frac{0}{0},\\
\text{поэтому правило Лопиталя}\\
\text{применяется $n$ раз}\end{pmatrix}}= ...= \lim_{x\to
x_0}\frac{f^{(n)}(x)-T_n^{(n)}(x)}{n!}=\\= {\smsize\begin{pmatrix}\text{применяем последнюю}\\
\text{формулу \eqref{11.9.2}}\end{pmatrix}}= \frac{0}{n!}=0
\end{multline*}
 Итак, мы получили
$$
\frac{f(x)-T_n(x)}{(x-x_0)^n}\underset{x\to x_0}{\longrightarrow} 0
$$
то есть,
$$
f(x)-T_n(x)= \underset{x\to x_0}{\bold{o}}((x-x_0)^n)
$$
или
$$
f(x)=T_n(x)+ \underset{x\to x_0}{\bold{o}}((x-x_0)^n)
$$
а это как раз и есть формула \eqref{Taylor-Peano}. \end{proof}

\noindent\rule{160mm}{0.1pt}\begin{multicols}{2}

\begin{ex}\label{ex-f-taylor-dlya-sin} Выпишем формулу Тейлора для функции $f(x)=\sin x$
в точке $x_0=\frac{\pi}{4}$ с точностью $n=2$:
 \begin{multline*}
\sin x=f(x)= f(x_0)+ \frac{f'(x_0)}{1!}\cdot (x-x_0)+\\+
\frac{f''(x_0)}{2!}\cdot (x-x_0)^2+ \underset{x\to x_0}{\bold{o}}\Big(
(x-x_0)^2 \Big)=\\= \sin \frac{\pi}{4}+ \frac{\cos\frac{\pi}{4}}{1!}\cdot \l
x-\frac{\pi}{4}\r+ \frac{-\sin (\frac{\pi}{4})}{2!}\cdot \l
x-\frac{\pi}{4}\r^2+\\+ \underset{x\to \frac{\pi}{4}}{\bold{o}}\l \l
x-\frac{\pi}{4}\r^2 \r= \frac{\sqrt{2}}{2}+ \frac{\sqrt{2}}{2}\cdot
\left(x-\frac{\pi}{4}\right)-\\-\frac{\sqrt{2}}{4}\cdot
\left(x-\frac{\pi}{4}\right)^2+ \underset{x\to
\frac{\pi}{4}}{\bold{o}}\left(\left(x-\frac{\pi}{4}\right)^2\right)
 \end{multline*}
Если убрать выкладки, получается следующий результат:
 \begin{multline*}
\sin x= \frac{\sqrt{2}}{2}+ \frac{\sqrt{2}}{2}\cdot
\left(x-\frac{\pi}{4}\right) -\frac{\sqrt{2}}{4}\cdot
\left(x-\frac{\pi}{4}\right)^2+\\+ \underset{x\to
\frac{\pi}{4}}{\bold{o}}\left(\left(x-\frac{\pi}{4}\right)^2\right)
 \end{multline*}
Можно увеличить точность, например, до $n=3$, и тогда получится
 \begin{multline*}
\sin x=f(x)= f(x_0)+ \frac{f'(x_0)}{1!}\cdot (x-x_0)+\\+
\frac{f''(x_0)}{2!}\cdot (x-x_0)^2+ \frac{f'''(x_0)}{3!}\cdot (x-x_0)^3+\\+
\underset{x\to x_0}{\bold{o}}((x-x_0)^3)= \sin \frac{\pi}{4}+
\frac{\cos\frac{\pi}{4}}{1!}\cdot \l x-\frac{\pi}{4}\r+\\+ \frac{-\sin
(\frac{\pi}{4})}{2!}\cdot \l x-\frac{\pi}{4}\r^2+ \frac{-\cos
(\frac{\pi}{4})}{3!}\cdot \l x-\frac{\pi}{4}\r^3+\\+ \underset{x\to
\frac{\pi}{4}}{\bold{o}}\left(\left(x-\frac{\pi}{4}\right)^3\right)=
\frac{\sqrt{2}}{2}+ \frac{\sqrt{2}}{2}\cdot \left(x-\frac{\pi}{4}\right)-\\
-\frac{\sqrt{2}}{4}\cdot \left(x-\frac{\pi}{4}\right)^2
-\frac{\sqrt{2}}{12}\cdot \left(x-\frac{\pi}{4}\right)^3+ \underset{x\to
\frac{\pi}{4}}{\bold{o}}\left(\left(x-\frac{\pi}{4}\right)^3\right)
 \end{multline*}
а без выкладок это будет выглядеть так:
 \begin{multline*}
\sin x= \frac{\sqrt{2}}{2}+ \frac{\sqrt{2}}{2}\cdot
\left(x-\frac{\pi}{4}\right) -\frac{\sqrt{2}}{4}\cdot
\left(x-\frac{\pi}{4}\right)^2 -\\-\frac{\sqrt{2}}{12}\cdot
\left(x-\frac{\pi}{4}\right)^3+ \underset{x\to
\frac{\pi}{4}}{\bold{o}}\left(\left(x-\frac{\pi}{4}\right)^3\right)
 \end{multline*}
\end{ex}

\begin{ex} Выпишем формулу Тейлора для функции $f(x)=\ln x$
в точке $x_0=1$ с точностью $n=3$:
 \begin{multline*}
\ln x=f(x)= f(x_0)+ \frac{f'(x_0)}{1!}\cdot (x-x_0)+\\+
\frac{f''(x_0)}{2!}\cdot (x-x_0)^2+ \frac{f'''(x_0)}{3!}\cdot
(x-x_0)^3+\\+\underset{x\to x_0}{\bold{o}}((x-x_0)^3)= \ln 1+
\frac{\frac{1}{x}|_{x=1}}{1!}\cdot (x-1)+\\+
\frac{-\frac{1}{x^2}|_{x=1}}{2!}\cdot (x-1)^2+
\frac{\frac{2}{x^3}|_{x=1}}{3!}\cdot (x-1)^3+\\+ \underset{x\to
x_0}{\bold{o}}((x-1)^3)= (x-1) -\frac{1}{2}\cdot (x-1)^2+\\+ \frac{1}{3}\cdot
(x-1)^3+ \underset{x\to x_0}{\bold{o}}((x-1)^3)
\end{multline*} Если убрать выкладки, получается следующий
результат:
$$
\ln x= (x-1) -\frac{(x-1)^2}{2}+ \frac{(x-1)^3}{3}+ \underset{x\to
x_0}{\bold{o}}((x-1)^3)
$$
\end{ex}

\end{multicols}\noindent\rule[10pt]{160mm}{0.1pt}

\paragraph{Формулы Маклорена-Пеано.}

 \bit{
 \item[$\bullet$]
При $x_0=0$ формула Тейлора называется {\it формулой
Маклорена}\index{формула!Маклорена}:
 \begin{multline}\label{11.9.3}
f(x)=f(0)+ \frac{f'(0)}{1!}\cdot x+ \frac{f''(0)}{2!}\cdot x^2+ +...+
\frac{f^{(n)}(0)}{n!}\cdot x^n+ \underset{x\to 0}{\bold{o}}(x^n)=\\=
\sum_{k=0}^n \frac{f^{(k)}(0)}{k!}\cdot x^k+ \underset{x\to 0}{\bold{o}}(x^n)
 \end{multline}
 }\eit

Некоторые примеры формул Маклорена нам уже знакомы:

\bigskip

\centerline{\bf Формулы Маклорена с точностью $n=1$ или $n=2$.}

\bigskip

Нетрудно заметить, что {\it асимптотические формулы \eqref{11.5.1} --
\eqref{11.5.10} являются частными случаями формулы Маклорена}. Например,
формула \eqref{11.5.1} получается, если подставить $f(x)=\sin x$ и $n=1$.
 $$
\sin x=f(x)=f(0)+ \frac{f'(0)}{1!}\cdot x+ \underset{x\to
0}{\bold{o}}(x)=0+x+\underset{x\to 0}{\bold{o}}(x)
 $$
А формула \eqref{11.5.2} -- если $f(x)=\cos x$ и $n=2$.
 $$
\cos x=f(x)=f(0)+ \frac{f'(0)}{1!}\cdot x+ \frac{f''(0)}{2!}\cdot
x^2+\underset{x\to 0}{\bold{o}}(x^2)= 1+0-\frac{x^2}{2}+\underset{x\to
0}{\bold{o}}(x^2)
 $$

\bigskip

\centerline{\bf Формулы Маклорена с точностью $n=3$}

\bigskip

Для тех же функций, что и в \eqref{11.5.1} -- \eqref{11.5.10} нетрудно выписать
соответствующие формулы с более высокой точностью, например, $n=3$ или $n=4$:
 \begin{align}
 &\sin x =x-\frac{x^3}{6}+\underset{x\to 0}{\bold{o}}(x^3)
 \label{11.9.4} \\
 &\tg x =x+\frac{x^3}{3}+\underset{x\to 0}{\bold{o}}(x^3)
 \label{11.9.5}\\
 &\cos x  =1-\frac{x^2}{2}+\underset{x\to 0}{\bold{o}}(x^3)
 \label{11.9.6}\\
 &\ln (1+x)  =x-\frac{x^2}{2}+\frac{x^3}{3}+\underset{x\to
 0}{\bold{o}}(x^3) \label{11.9.7}\\
 &\log_a (1+x) =\frac{x}{\ln a}-\frac{x^2}{2\ln a}+\frac{x^3}{3\ln
 a}+\underset{x\to 0}{\bold{o}}(x^3) \label{11.9.8}\\
 &e^x =1+x+\frac{x^2}{2}+\frac{x^3}{6}+\underset{x\to 0}{\bold{o}}(x^3)
 \label{11.9.9}\\
 &a^x =1+x\cdot \ln a+\frac{x^2\cdot \ln^2 a}{2}+
 \frac{x^3\cdot \ln^3 a}{6}+
 \underset{x\to 0}{\bold{o}}(x^3)\; \label{11.9.10}\\
 &(1+x)^\alpha = 1+\alpha x+\frac{\alpha (\alpha-1)
 x^2}{2}+\frac{\alpha (\alpha-1) (\alpha-2) x^3}{6}+ \underset{x\to
 0}{\bold{o}}(x^3) \label{11.9.11}\\
 &\arcsin x =x+\frac{x^3}{6}+\underset{x\to 0}{\bold{o}}(x^3)
 \label{11.9.12}\\
 &\arctg x =x-\frac{x^3}{3}+\underset{x\to 0}{\bold{o}}(x^3)
 \label{11.9.13}
 \end{align}

\bigskip
\centerline{\bf Формулы Маклорена с произвольной точностью}
\bigskip

Некоторые из формул Маклорена \eqref{11.9.4} -- \eqref{11.9.13} можно обобщить
на случай произвольной точности:
 \begin{align}
 &\sin x=x-\frac{x^3}{3!}+\frac{x^5}{5!}-...+(-1)^n \cdot
\frac{x^{2k+1}}{(2k+1)!} +\underset{x\to
0}{\bold{o}}(x^{2k+1})=\nonumber\\
 &\kern180pt =\sum_{i=1}^{k} (-1)^i\cdot
\frac{x^{2i+1}}{(2i+1)!} +\underset{x\to
0}{\bold{o}}(x^{2k+1}) \label{11.9.14}\\
 &\cos x=1-\frac{x^2}{2!}+\frac{x^4}{4!}-...+(-1)^k
\cdot\frac{x^{2k}}{(2k)!}+\underset{x\to 0}{\bold{o}}(x^{2k})= \sum_{i=0}^{k}
(-1)^i\cdot \frac{x^{2i}}{(2i)!} +\underset{x\to 0}{\bold{o}}(x^{2k})
\label{11.9.15}\\
 &\ln (1+x)=x-\frac{x^2}{2}+\frac{x^3}{3}+...+(-1)^{k-1}\cdot
\frac{x^k}{k} +\underset{x\to 0}{\bold{o}}(x^k)=\sum_{i=1}^{k} (-1)^{i-1}\cdot
\frac{x^i}{i} +\underset{x\to 0}{\bold{o}}(x^k)
\label{11.9.16}\\
 &e^x=1+x+\frac{x^2}{2!}+...+\frac{x^k}{k!}+ \underset{x\to
0}{\bold{o}}(x^k)= \sum_{i=0}^{k}\frac{x^i}{i!} +\underset{x\to
0}{\bold{o}}(x^k) \label{11.9.17}\\
 &(1+x)^\alpha=1+\alpha \cdot x+\alpha (\alpha-1) \cdot
\frac{x^2}{2!}+...+ \alpha (\alpha-1) ... (\alpha-k+1) \cdot
\frac{x^k}{k!}+ \underset{x\to 0}{\bold{o}(x^k)}=\nonumber\\
 &\kern180pt =
 1+\sum_{i=1}^{k}\alpha (\alpha-1) ... (\alpha-i+1) \cdot \frac{x^i}{i!}
+\underset{x\to 0}{\bold{o}}(x^k)
 \label{11.9.18}\\
 &\frac{1}{1+x}=1-x+x^2+...+ (-1)^k\cdot x^k+ \underset{x\to
0}{\bold{o}(x^k)}= \sum_{i=0}^{k} (-1)^k \cdot x^i +\underset{x\to
0}{\bold{o}}(x^k)
\label{11.9.19}\\
 &\frac{1}{1-x}=1+x+x^2+...+x^k+
\underset{x\to 0}{\bold{o}(x^k)}= \sum_{i=0}^{k} x^i +\underset{x\to
0}{\bold{o}}(x^k)
\label{11.9.20}\\
 &\arctg
x=x-\frac{x^3}{3}+\frac{x^5}{5}-...+ (-1)^k\cdot \frac{x^{2k+1}}{2k+1}+
\underset{x\to 0}{\bold{o}(x^k)}= \sum_{i=0}^{k} (-1)^i\cdot
\frac{x^{2i+1}}{2i+1} +\underset{x\to 0}{\bold{o}}(x^k)
\label{11.9.21}\end{align}

\paragraph{Вычисление пределов в помощью формул Маклорена-Пеано.}

Формулы Маклорена позволяют вычислять пределы, которые невозможно вычислить с
помощью асимптотических формул \eqref{11.5.1} -- \eqref{11.5.10}. Правила,
которые следует иметь в виду в таких случаях, состоят в следующем.

\bigskip

\centerline{\bf Правила вычисления пределов с помощью формул Маклорена:}\nobreak
 \bit{\it
 \item[(A)] При использовании формул Маклорена может не получиться ответа,
если возникнет неопределенность.

\item[(B)] Если возникает неопределенность, то
это означает, что нужно увеличить точность формул Маклорена.

\item[(C)] Если
ответ получен, то дальнейшее увеличение точности формул Маклорена не может
привести к изменению ответа.
 }\eit

\bigskip

Чтобы понять смысл этих слов, рассмотрим несколько примеров.

\noindent\rule{160mm}{0.1pt}\begin{multicols}{2}

\begin{ex}
$$
\lim_{x\to 0}\frac{\sqrt{1+2\tg x}-e^x+x^2}{\arcsin x-\sin x}
$$
Объясним сначала, почему здесь нужны фор\-мулы Маклорена. Если этот предел
вычислять с помощью правила Лопиталя, то Вы увидите, что дифференцировать
придется 3 раза, а поскольку функции в числителе и знаменателе не слишком
просты, это будет довольно утомительно. Покажем, с другой стороны, что и
асимптотические формулы \eqref{11.5.1} -- \eqref{11.5.10} здесь тоже не
помогают:
 \begin{multline*}
\lim_{x\to 0}\frac{\sqrt{1+2\tg x}-e^x+x^2}{\arcsin x-\sin x}=\\
={\smsize\begin{pmatrix}\text {применяем формулы}\\
\eqref{11.5.8},\eqref{11.5.6},\eqref{11.5.9},\eqref{11.5.1}
\end{pmatrix}}=\\
{\smsize\text{$=\lim_{x\to 0}\frac {1+\frac{1}{2} 2\tg x+\underset{x\to
0}{\bold{o}}(2\tg x) -(1+x+\underset{x\to 0}{\bold{o}}(x))+x^2}
{x+\underset{x\to 0}{\bold{o}}(x)-(x+\underset{x\to
0}{\bold{o}}(x))}=$}}\\
=\lim_{x\to 0}\frac {\tg x+\underset{x\to 0}{\bold{o}}(2\tg x)
-x-\underset{x\to 0}{\bold{o}}(x)+x^2} {\underset{x\to
0}{\bold{o}}(x)-\underset{x\to 0}{\bold{o}}(x)}=\\=
{\smsize\begin{pmatrix}\text {применяем свойства}\\
1^0, 3^0
\end{pmatrix}}=\\
=\lim_{x\to 0}\frac {\tg x+\underset{x\to 0}{\bold{o}}(\tg x) -x-\underset{x\to
0}{\bold{o}}(x)+x^2} {\underset{x\to
0}{\bold{o}}(x)}=\\= {\smsize\begin{pmatrix}\text {применяем}\\
(5.2)
\end{pmatrix}}=\\
{\smsize\text{$=\lim_{x\to 0}\frac {x+\underset{x\to
0}{\bold{o}}(x)+\underset{x\to 0}{\bold{o}} (x+\underset{x\to 0}{\bold{o}}(x))
-x-\underset{x\to 0}{\bold{o}}(x)+x^2}
{\underset{x\to 0}{\bold{o}}(x)}=$}}\\
={\smsize\begin{pmatrix}\text {теперь}\\
4^0
\end{pmatrix}}=\\
=\lim_{x\to 0}\frac {x+\underset{x\to 0}{\bold{o}}(x)+\underset{x\to
0}{\bold{o}}(x) -x-\underset{x\to 0}{\bold{o}}(x)+x^2} {\underset{x\to
0}{\bold{o}}(x)}=\\=
{\smsize\begin{pmatrix}\text {теперь}\\
2^0, 3^0
\end{pmatrix}}= \lim_{x\to 0}\frac {\underset{x\to 0}{\bold{o}}(x)+x^2}
{\underset{x\to 0}{\bold{o}}(x)}= {\smsize\begin{pmatrix}\text {теперь}\\
1^0
\end{pmatrix}}=\\
=\lim_{x\to 0}\frac {x\cdot \underset{x\to 0}{\bold{o}}(1)+x^2} {x\cdot
\underset{x\to 0}{\bold{o}}(1)}= \lim_{x\to 0}\frac {\underset{x\to
0}{\bold{o}}(1)+x} {\underset{x\to
0}{\bold{o}}(1)}=\\= {\smsize\begin{pmatrix}\text {теперь}\\
8^0
\end{pmatrix}}= \left(\frac{0+0}{0}\right)=?
 \end{multline*}
Тем не менее, ситуация не безвыходна. Покажем как этот предел можно вычислить с
помощью формул Маклорена. Заметим только сразу, что \eqref{11.5.1} --
\eqref{11.5.10} уже представляют собой формулы Маклорена с точностью $n=1$, а
то, что с их помощью не удалось найти ответ есть иллюстрация сформулированного
нами выше правила (A). Если же их заменить на другие формулы с большей
точностью, то чудесным образом пример становится решаемым.

Действительно, попробуем применить формулы \eqref{11.9.4} -- \eqref{11.9.13}
(то есть формулы Маклорена с точностью $n=3$), и убедимся, что после этого
сработает правило (B). Для этого сначала выпишем отдельно корень в числителе:
 \begin{multline*}
\sqrt{1+2\tg x}= {\smsize\begin{pmatrix}\text{по формуле \eqref{11.9.11},}\\
\sqrt{1+y}=\\= 1+\frac{1}{2} y-\frac{y^2}{8}+ \frac{y^3}{16}+ \underset{y\to
0}{\bold{o}}(y^3)
\end{pmatrix}}=\\
=1+\tg x-\frac{\tg^2 x}{2}+ \frac{\tg^3 x}{2}+ \underset{x\to
0}{\bold{o}}(\tg^3 x)=\\= \left(\text{применяем
\eqref{11.9.5}}\right)=\\
=1+\left(x+\frac{x^3}{3}+\underset{x\to 0}{\bold{o}}(x^3)\right)
-\\
-\frac{\left(x+\frac{x^3}{3}+\underset{x\to 0}{\bold{o}}(x^3)\right)^2}{2}+
\frac{ \left(x+\frac{x^3}{3}+\underset{x\to 0}{\bold{o}}(x^3)\right)^3}{2}+\\+
\underset{x\to 0}{\bold{o}}\left((x+\frac{x^3}{3}+\underset{x\to
0}{\bold{o}}(x^3))^3\right)=\\
= ...= {\smsize\begin{pmatrix}\text{упрощаем}\\ \text{с помощью}\\
\text{свойств}\; 1^0 - 8^0 \; \S 7
\end{pmatrix}}=...=\\
=1+x-\frac{x^2}{2}+\frac{5x^3}{6}+ \underset{x\to 0}{\bold{o}} (x^3)
 \end{multline*}
Теперь можно вычислять предел:
 \begin{multline*}
\lim_{x\to 0}\frac{\sqrt{1+2\tg x}-e^x+x^2}{\arcsin x-\sin x}=\\
=\lim_{x\to 0}\frac{ 1+x-\frac{x^2}{2}+\frac{5x^3}{6}+ \underset{x\to
0}{\bold{o}} (x^3) -e^x+x^2}{\arcsin x-\sin x}=\\
={\smsize\begin{pmatrix}\text{применяем формулы}\\
\eqref{11.9.9}, \eqref{11.9.12}, \eqref{11.9.13}
\end{pmatrix}}=\\
=\lim_{x\to 0}\Bigg\{ 1+x-\frac{x^2}{2}+\frac{5x^3}{6}+ \underset{x\to
0}{\bold{o}} (x^3)
-\\
-\left(1+x+\frac{x^2}{2}+\frac{x^3}{6}+\underset{x\to 0}{\bold{o}}(x^3)\right)
+x^2\Bigg\}:
\\
:\Bigg\{x+\frac{x^3}{6}+\underset{x\to 0}{\bold{o}}(x^3)
-\left(x-\frac{x^3}{6}+\underset{x\to
0}{\bold{o}}(x^3)\right)\Bigg\}=\\
=\lim_{x\to 0}\frac{\frac{2x^3}{3}+\underset{x\to 0}{\bold{o}} (x^3)}
{\frac{x^3}{3}+\underset{x\to 0}{\bold{o}}(x^3)}= \lim_{x\to
0}\frac{\frac{2x^3}{3}+x^3\cdot \underset{x\to 0}{\bold{o}} (1)}
{\frac{x^3}{3}+x^3\cdot \underset{x\to 0}{\bold{o}}(1)}=\\
=\lim_{x\to 0}\frac{\frac{2}{3}+\underset{x\to 0}{\bold{o}} (1)}
{\frac{1}{3}+\underset{x\to 0}{\bold{o}}(1)}=
\frac{\frac{2}{3}+0}{\frac{1}{3}+0}=2
 \end{multline*}
Число 2 будет ответом.

Теперь полезно убедиться в справедливости правила (C). Проверим, что если
увеличить точность формул Маклорена, то ответ от этого не изменится. Это можно
сделать, например, выписав более подробное разложение экспоненты:
 \begin{multline*}
\lim_{x\to 0}\frac{\sqrt{1+2\tg x}-e^x+x^2}{\arcsin x-\sin x}=\\
=\lim_{x\to 0}\frac{ 1+x-\frac{x^2}{2}+\frac{5x^3}{6}+ \underset{x\to
0}{\bold{o}} (x^3) -e^x+x^2}{\arcsin x-\sin x}=\\=
{\smsize\begin{pmatrix}\text{применяем формулы}\\ \eqref{11.9.12}, \eqref{11.9.13} \\
\text{и \eqref{11.9.17} с}\; k=4
\end{pmatrix}}=\\
=\lim_{x\to 0}\Bigg\{ 1+x-\frac{x^2}{2}+\frac{5x^3}{6}+ \underset{x\to
0}{\bold{o}} (x^3)
-\\
-\left(1+x+\frac{x^2}{2!}+\frac{x^3}{3!}+\frac{x^4}{4!}+ \underset{x\to
0}{\bold{o}}(x^4)\right) +x^2\Bigg\}:
\\
:\Bigg\{x+\frac{x^3}{6}+\underset{x\to 0}{\bold{o}}(x^3)
-\left(x-\frac{x^3}{6}+\underset{x\to
0}{\bold{o}}(x^3)\right)\Bigg\}=\\
=\lim_{x\to 0}\frac{\frac{2x^3}{3}+\underset{x\to 0}{\bold{o}} (x^3)
+\frac{x^4}{24}+\underset{x\to 0}{\bold{o}}(x^4)}
{\frac{x^3}{3}+\underset{x\to 0}{\bold{o}}(x^3)}=\\
=\lim_{x\to 0}\frac{\frac{2x^3}{3}+x^3\cdot \underset{x\to 0}{\bold{o}} (1)
+\frac{x^4}{24}+x^4\cdot \underset{x\to 0}{\bold{o}}(1)}
{\frac{x^3}{3}+x^3\cdot \underset{x\to 0}{\bold{o}}(1)}=\\= \lim_{x\to
0}\frac{\frac{2}{3}+\underset{x\to 0}{\bold{o}}(1) +\frac{x}{24}+\underset{x\to
0}{\bold{o}}(1)}
{\frac{1}{3}+\underset{x\to 0}{\bold{o}}(1)}=\\
=\frac{\frac{2}{3}+0+0+0}{\frac{1}{3}+0}=2
 \end{multline*}
То же самое получится, если поменять формулу Маклорена на более точную в любом
другом месте (или в нескольких местах).
\end{ex}

\begin{ex}
$$
\lim_{x\to 0}\frac{e^{\arctg x}-\frac{1}{1-x}+\frac{x^2}{2}} {\ln
\frac{1+x}{1-x}-2x}
$$
Убедимся, что, как и в предыдущем случае, этот предел не вычисляется с помощью
формул \eqref{11.5.1} -- \eqref{11.5.10}, то есть с помощью формул Маклорена с
точностью $n=1$:
 \begin{multline*}
\lim_{x\to 0}\frac{e^{\arctg x}-\frac{1}{1-x}+\frac{x^2}{2}} {\ln
\frac{1+x}{1-x}-2x}=\\
{\smsize\text{$=\lim_{x\to 0}\frac{ 1+\arctg x+\underset{x\to
0}{\bold{o}}(\arctg x)
-(1-x)^{-1}+\frac{x^2}{2}} {\ln (1+x)- \ln (1-x)-2x}=$}}\\
=\lim_{x\to 0}\Bigg\{ 1+\left(x+\underset{x\to 0}{\bold{o}}(x)\right)
+\underset{x\to 0}{\bold{o}}\left(x+\underset{x\to
0}{\bold{o}}(x)\right)-\\
-\left(1-x+\underset{x\to
0}{\bold{o}}(x)\right) +\frac{x^2}{2}\Bigg\}: \\
:\Bigg\{\left(x+\underset{x\to 0}{\bold{o}}(x)\right)
-\left(-x+\underset{x\to 0}{\bold{o}}(x)\right) -2x\Bigg\}=\\
=\lim_{x\to 0}\frac{ \underset{x\to 0}{\bold{o}}(x) +\frac{x^2}{2}}
{\underset{x\to 0}{\bold{o}}(x)}= \lim_{x\to 0}\frac{ x\cdot \underset{x\to
0}{\bold{o}}(1) +\frac{x^2}{2}} {x\cdot \underset{x\to
0}{\bold{o}}(1)}=\\
=\lim_{x\to 0}\frac{ \underset{x\to 0}{\bold{o}}(1) +\frac{x}{2}}
{\underset{x\to 0}{\bold{o}}(1)}= \left(\frac{0+0}{0}\right)=?
 \end{multline*}
Вы можете проверить самостоятельно, что даже если увеличить точность формул
Маклорена до $n=2$, это все равно не поможет.

Но если воспользоваться формулами Мак\-лорена с точностью $n=3$ (то есть
формулами \eqref{11.9.4} -- \eqref{11.9.13}), то этот предел удается сосчитать.
Только, как и в предыдущем примере, перед вычислениями здесь следует выписать
отдельно асимптотическое представление некоторых фрагментов:
 \begin{multline*}\ln \frac{1+x}{1-x}= \ln (1+x) -\ln (1-x)=
\\=\left(\text{применяем \eqref{11.9.7}}\right)=\\=
x-\frac{x^2}{2}+\frac{x^3}{3}+\underset{x\to 0}{\bold{o}}(x^3)-\\- \left(
-x-\frac{x^2}{2}-\frac{x^3}{3}+\underset{x\to
0}{\bold{o}}(x^3) \right)=\\
=2x+\frac{2x^3}{3}+\underset{x\to 0}{\bold{o}}(x^3)
 \end{multline*}
Далее,
 \begin{multline*}
e^{\arctg x}= \left(\text{применяем \eqref{11.9.9}}\right)=\\
=1+\arctg x+\frac{\arctg^2 x}{2}+\frac{\arctg^3 x}{6}+\\
+\underset{x\to 0}{\bold{o}}(\arctg^3 x)=\left(\text{применяем
\eqref{11.9.13}}\right)=\\
=1+\left(x-\frac{x^3}{3}+\underset{x\to 0}{\bold{o}}(x^3)\right)
+\\
+\frac{ \left(x-\frac{x^3}{3}+\underset{x\to 0}{\bold{o}}(x^3)\right)^2}{2}+\\+
\frac{\left(x-\frac{x^3}{3}+\underset{x\to 0}{\bold{o}}(x^3)\right)^3}{6}
+\\+\underset{x\to 0}{\bold{o}}\left(\left(x-\frac{x^3}{3}+\underset{x\to
0}{\bold{o}}(x^3)\right)^3\right)=\\=
...= {\smsize\begin{pmatrix}\text{упрощаем}\\ \text{с помощью}\\
\text{свойств}\; 1^0 - 8^0 \; \S 7
\end{pmatrix}}=...=\\= 1+x+\frac{x^2}{2}-\frac{x^3}{6}+ \underset{x\to
0}{\bold{o}}(x^3)
 \end{multline*}
И, наконец,
 \begin{multline*}
\frac{1}{1-x}=(1-x)^{-1}= \left(\text{применяем
\eqref{11.9.11}}\right)=\\= 1+(-1)(-x)+\frac{(-1)(-2) (-x)^2}{2}+\\
+\frac{(-1)(-2)(-3) (-x)^3}{6}+ \underset{x\to 0}{\bold{o}}(x^3)=\\
=1+x+x^2+x^3+\underset{x\to 0}{\bold{o}}(x^3)
 \end{multline*}
Теперь вычисляем предел:
 \begin{multline*}
\lim_{x\to 0}\frac{e^{\arctg x}-\frac{1}{1-x}+\frac{x^2}{2}} {\ln
\frac{1+x}{1-x}-2x}=\\
=\lim_{x\to 0}\Bigg\{ \left(1+x+\frac{x^2}{2}-\frac{x^3}{6}+ \underset{x\to
0}{\bold{o}}(x^3) \right) -\\- \left( 1+x+x^2+x^3+\underset{x\to
0}{\bold{o}}(x^3)\right)
+\frac{x^2}{2}\Bigg\}: \\
:\Bigg\{\left( 2x+\frac{2x^3}{3}+\underset{x\to 0}{\bold{o}}(x^3)
\right)-2x\Bigg\}=\\
=\lim_{x\to 0}\frac{-\frac{7 x^3}{6}+\underset{x\to 0}{\bold{o}}(x^3)}
{\frac{2x^3}{3}+\underset{x\to
0}{\bold{o}}(x^3)}=\\
=\lim_{x\to 0}\frac{-\frac{7 x^3}{6}+ x^3\cdot \underset{x\to 0}{\bold{o}}(1)}
{\frac{2x^3}{3}+x^3\cdot \underset{x\to
0}{\bold{o}}(1)}=\\
=\lim_{x\to 0}\frac{-\frac{7}{6}+ \underset{x\to 0}{\bold{o}}(1)}
{\frac{2}{3}+\underset{x\to 0}{\bold{o}}(1)}=
-\frac{7}{4}\end{multline*}\end{ex}

\begin{ex}
$$
\lim_{x\to 0}\frac {\ln \left( x+\sqrt{1+x^2}\right)-x+\frac{x^3}{6}}
{x-\arcsin x}
$$
Здесь мы сразу применим формулы Маклорена с точностью $n=3$. Выпишем отдельно
асимп\-тотическое представление некоторых фрагментов:
 \begin{multline*}
\ln \left( x+\sqrt{1+x^2}\right)= \ln \left( x+(1+x^2)^\frac{1}{2}\right)=\\=
{\smsize\begin{pmatrix}\text{применяем}\\ \eqref{11.9.11}
\end{pmatrix}}=\ln \Bigg( x+ 1+\frac{1}{2} x+
\frac{\frac{1}{2}\left(\frac{1}{2}-1\right) x^2}{2}+\\
+ \frac{\frac{1}{2}\left(\frac{1}{2}-1\right) \left(\frac{1}{2}-2\right) x^3}
{6}+ \underset{x\to
0}{\bold{o}}(x^3)\Bigg)=\\
=\ln \Bigg( 1+\frac{3}{2} x+\frac{\frac{1}{2}\left(-\frac{1}{2}\right) x^2}{2}+
\frac{\frac{1}{2}\left(-\frac{1}{2}\right)
\left(-\frac{3}{2}\right) x^3} {6}+\\
+\underset{x\to 0}{\bold{o}}(x^3)\Bigg)=\\
=\ln \left( 1+\frac{3}{2} x-\frac{x^2}{8}+ \frac{3 x^3} {48}+ \underset{x\to
0}{\bold{o}}(x^3)\right)=\\= {\smsize\begin{pmatrix}\text{применяем}\\
\eqref{11.9.7}
\end{pmatrix}}=\\
=\left(\frac{3}{2} x-\frac{x^2}{8}+ \frac{3 x^3} {48}+ \underset{x\to
0}{\bold{o}}(x^3)\right) -\\
-\frac{ \left(\frac{3}{2} x-\frac{x^2}{8}+ \frac{3 x^3} {48}+ \underset{x\to
0}{\bold{o}}(x^3)\right)^2}{2}+\\+ \frac{ \left(\frac{3}{2} x-\frac{x^2}{8}+
\frac{3 x^3} {48}+ \underset{x\to 0}{\bold{o}}(x^3)\right)^3}
{3}+\\+\underset{x\to 0}{\bold{o}}\left(\left(\frac{3}{2} x-\frac{x^2}{8}+
\frac{3 x^3} {48}+ \underset{x\to 0}{\bold{o}}(x^3)\right)^3\right)=\\=...=
{\smsize\begin{pmatrix}\text{упрощаем}\\
\text{с помощью}\\
1^0 - 9^0 \; \S 7
\end{pmatrix}}=...= x-\frac{x^3}{6}+\underset{x\to 0}{\bold{o}}(x^3)
 \end{multline*}
Теперь вычисляем предел:
 \begin{multline*}
\lim_{x\to 0}\frac {\ln \left( x+\sqrt{1+x^2}\right)-x+\frac{x^3}{6}}
{x-\arcsin x}=
\\=
{\smsize\begin{pmatrix}\text{подставляем только что найденное}\\ \text{разложение}\ln \left( x+\sqrt{1+x^2}\right) \\
\text{а также формулу \eqref{11.9.12}}\end{pmatrix}}=\\= \lim_{x\to 0}\frac
{\left( x-\frac{x^3}{6}+\underset{x\to 0}{\bold{o}}(x^3) \right)
-x+\frac{x^3}{6}} {x-\left( x+\frac{x^3}{6}+\underset{x\to
0}{\bold{o}}(x^3)\right)}=\\= \lim_{x\to 0}\frac {\underset{x\to
0}{\bold{o}}(x^3)} {-\frac{x^3}{6}-\underset{x\to 0}{\bold{o}}(x^3)}=
\lim_{x\to 0}\frac {x^3\cdot \underset{x\to 0}{\bold{o}}(1)}
{-\frac{x^3}{6}-x^3\cdot \underset{x\to 0}{\bold{o}}(1)}=\\= \lim_{x\to 0}\frac
{\underset{x\to 0}{\bold{o}}(1)} {-\frac{1}{6}-\underset{x\to 0}{\bold{o}}(1)}=
\frac{0}{-\frac{1}{6}-0}=0
\end{multline*}\end{ex}

\begin{ex}\label{EX:pokaz-predel-v-asimp-metodah}
Вспомним теперь о пределе \eqref{LIM:reklama-asimp-metodov}, которым мы
начинали эту главу:
$$
\lim_{x\to 0}\frac{\ln (1+\sin x)-\sin \ln (1+x)}{x^4}
$$
Как уже говорилось, вычислять его с помощью правила Лопиталя слишком тяжело. С
другой стороны, применять теорему об эквивалентной замене или выделять главное
слагаемое тут тоже бесполезно. Попробуем применить асимптотические формулы
\eqref{11.5.1} -- \eqref{11.5.10} (то есть  формулы Маклорена с точностью
$n=1$):
 \begin{multline*}
\lim_{x\to 0}\frac{\ln (1+\sin x)-\sin \ln (1+x)}{x^4}=\\=
{\smsize\begin{pmatrix}\text {применяем формулы}\\
\eqref{11.5.4},\eqref{11.5.1},\eqref{11.5.10}
\end{pmatrix}}=\\
=\lim_{x\to 0}\Bigg\{\sin x+\underset{x\to 0}{\bold{o}}(\sin x) -\\-\left(\ln
(1+x)+\underset{x\to 0}{\bold{o}}(\ln
(1+x)) \right)\Bigg\}:{x^4}=\\
={\smsize\begin{pmatrix}\text {применяем формулы \eqref{11.5.4},\eqref{11.5.1}}\\
\text {и свойство}\; 1^{\bold{o}\bold{o}}\; \S 7
\end{pmatrix}}=\\
=\lim_{x\to 0}\Bigg\{ \left(x+\underset{x\to 0}{\bold{o}}(x)\right)
+\underset{x\to 0}{\bold{o}}\left(x+\underset{x\to 0}{\bold{o}}(x)\right)
-\\-\left(\left(x+\underset{x\to 0}{\bold{o}}(x)\right) +\underset{x\to
0}{\bold{o}}\left(x+\underset{x\to 0}{\bold{o}}(x)\right)
\right)\Bigg\}:{x^4}=\\
={\smsize\begin{pmatrix}\text {упрощаем с помощью}\\
\text {свойств}\; 1^0 - 9^0 \; \S 7
\end{pmatrix}}= \lim_{x\to 0}\frac{ \underset{x\to 0}{\bold{o}}(x) }
{x^4+\underset{x\to 0}{\bold{o}}(x^4)}=\\= \lim_{x\to 0}\frac{ x\cdot
\underset{x\to 0}{\bold{o}}(1) } {x^4}= \lim_{x\to 0}\frac{ \underset{x\to
0}{\bold{o}}(1) } {x^3}= \left(\frac{0}{0}\right)=?
 \end{multline*}
Вы видите, что ничего хорошего не получается. В соответствии с правилом (B),
это означает, что нужно увеличить точность формул Маклорена. Возьмем $n=3$ (то
есть воспользуемся формулами \eqref{11.9.4} -- \eqref{11.9.13}). Тогда мы
получим:
 \begin{multline*}
\ln (1+\sin x)= \sin x-\frac{\sin^2 x}{2}+\frac{\sin^3 x}{3}+\\
+\underset{x\to 0}{\bold{o}}\left(\sin^3 x \right)= \left(
x-\frac{x^3}{6}+\underset{x\to 0}{\bold{o}}(x^3)\right)
-\\
-\frac{1}{2}\left( x-\frac{x^3}{6}+\underset{x\to 0}{\bold{o}}(x^3)\right)^2
+\\+\frac{1}{3}\left(
x-\frac{x^3}{6}+\underset{x\to 0}{\bold{o}}(x^3)\right)^3 +\\
+\underset{x\to 0}{\bold{o}}\left(\left( x-\frac{x^3}{6}+\underset{x\to
0}{\bold{o}}(x^3)\right)^3
\right)=\\
=...= x-\frac{x^2}{2}+\frac{x^3}{6}+ \underset{x\to 0}{\bold{o}}(x^3)
 \end{multline*}
и
 \begin{multline*}
\sin \ln (1+x)= \ln (1+x)-\frac{\left(\ln (1+x)\right)^3}{6}+
\\+\underset{x\to 0}{\bold{o}}\left(\left(\ln (1+x)\right)^3 \right)=
\Bigg( x-\frac{x^2}{2}+\frac{x^3}{3}+\\+ \underset{x\to 0}{\bold{o}}\left( x^3
\right) \Bigg)- \frac{\left( x-\frac{x^2}{2}+\frac{x^3}{3}+ \underset{x\to
0}{\bold{o}}\left( x^3
\right) \right)^3}{6}+ \\
+\underset{x\to 0}{\bold{o}}\left(\left( x-\frac{x^2}{2}+\frac{x^3}{3}+
\underset{x\to 0}{\bold{o}}\left( x^3 \right) \right)^3 \right)=\\=...=
x-\frac{x^2}{2}+\frac{x^3}{6}+ \underset{x\to 0}{\bold{o}}(x^3)
 \end{multline*}
и поэтому
 \begin{multline*}
\lim_{x\to 0}\frac{\ln (1+\sin x)-\sin \ln (1+x)}{x^4}=\\
=\lim_{x\to 0}\Bigg\{\left( x-\frac{x^2}{2}+\frac{x^3}{6}+
\underset{x\to 0}{\bold{o}}(x^3) \right) -\\
-\left(x-\frac{x^2}{2}+\frac{x^3}{6}+ \underset{x\to 0}{\bold{o}}(x^3)
\right)\Bigg\}:{x^4}=\\= \lim_{x\to 0}\frac{\underset{x\to
0}{\bold{o}}(x^3)}{x^4}= \lim_{x\to
0}\frac{x^3\cdot \underset{x\to 0}{\bold{o}}(1)}{x^4}=\\
=\lim_{x\to 0}\frac{\underset{x\to 0}{\bold{o}}(1)}{x}=
\left(\frac{0}{0}\right)=?
 \end{multline*}
Опять получилась неопределенность. Она снова означает, что нужно увеличить
точность формул Маклорена (согласно правилу (B)). Возьмем $n=4$ (то есть
воспользуемся формулами \eqref{11.9.14} и \eqref{11.9.16}, положив $n=4$):
$$
\sin x=x-\frac{x^3}{3!}+\underset{x\to 0}{\bold{o}}(x^4)
$$
$$
\ln (1+x)= x-\frac{x^2}{2}+\frac{x^3}{3}-\frac{x^4}{4} +\underset{x\to
0}{\bold{o}}(x^4)
$$
Тогда получим
 \begin{multline*}
\ln (1+\sin x)= \sin x-\frac{\sin^2 x}{2}+\frac{\sin^3
x}{3}-\frac{\sin^4 x}{4}+\\
+\underset{x\to 0}{\bold{o}}\left(\sin^4 x \right)= \left(
x-\frac{x^3}{3!}+\underset{x\to 0}{\bold{o}}(x^4) \right) -\\-
\frac{1}{2}\left( x-\frac{x^3}{3!}+\underset{x\to 0}{\bold{o}}(x^4)
\right)^2 +\\
+\frac{1}{3}\left( x-\frac{x^3}{3!}+\underset{x\to 0}{\bold{o}}(x^4)
\right)^3 -\\
-\frac{1}{4}\left( x-\frac{x^3}{3!}+\underset{x\to 0}{\bold{o}}(x^4) \right)^4
+\\+ \underset{x\to 0}{\bold{o}}\left(\left( x-\frac{x^3}{3!}+\underset{x\to
0}{\bold{o}}(x^4) \right)^4
\right)=\\
=...=x-\frac{x^2}{2}+\frac{x^3}{6}-\frac{x^4}{12}+ \underset{x\to
0}{\bold{o}}(x^4)
 \end{multline*}
и
 \begin{multline*}
\sin \ln (1+x)= \ln (1+x)-\frac{\left(\ln (1+x)\right)^3}{6}+
\\+\underset{x\to 0}{\bold{o}}\left(\left(\ln (1+x)\right)^4
\right)=\\
=\left(x-\frac{x^2}{2}+\frac{x^3}{3}-\frac{x^4}{4} +\underset{x\to
0}{\bold{o}}(x^4)\right) -\\-
\frac{1}{6}\left(x-\frac{x^2}{2}+\frac{x^3}{3}-\frac{x^4}{4} +\underset{x\to
0}{\bold{o}}(x^4)\right)^3 +\\+ \underset{x\to
0}{\bold{o}}\left(\left(x-\frac{x^2}{2}+\frac{x^3}{3}-\frac{x^4}{4}
+\underset{x\to 0}{\bold{o}}(x^4)\right)^4  \right)=\\=...=
x-\frac{x^2}{2}+\frac{x^3}{6}+0+ \underset{x\to 0}{\bold{o}}(x^4)
 \end{multline*}
Теперь можно, наконец, вычислить наш предел:
 \begin{multline*}
\lim_{x\to 0}\frac{\ln (1+\sin x)-\sin \ln (1+x)}{x^4}=\\
=\lim_{x\to 0}\Bigg\{ \left( x-\frac{x^2}{2}+\frac{x^3}{6}-\frac{x^4}{12}+
\underset{x\to
0}{\bold{o}}(x^4) \right) -\\
-\left( x-\frac{x^2}{2}+\frac{x^3}{6}+0+ \underset{x\to 0}{\bold{o}}(x^4)
\right)\Bigg\}:{x^4}=\\= \lim_{x\to 0}\frac{ -\frac{x^4}{12}+ \underset{x\to
0}{\bold{o}}(x^4) }{x^4}= \lim_{x\to 0}\frac{ -\frac{x^4}{12}+ x^4\cdot
\underset{x\to 0}{\bold{o}}(1) }{x^4}=\\= \lim_{x\to 0}\left(-\frac{1}{12}+
\underset{x\to 0}{\bold{o}}(1)\right)=-\frac{1}{12}
 \end{multline*}\end{ex}

\begin{ex}
$$
\lim_{x\to 0}\left(  \cos (x e^x) -\ln (1-x) -x \right)^{\ctg x^3}
$$
Здесь мы не будем утомлять читателя предварительными прикидками, а сразу
применим формулы Маклорена с ``нужной'' точностью. При этом, в разных местах
точность мы будем брать разной, но важно понимать, согласно нашим правилам, что
если ты не знаешь какой должна быть точность, то можно брать наугад любую, а
затем уже, если получится неопределенность (а это самое страшное, что может
получиться) -- повышать порядок точности до тех пор, пока не получится
определенность. И, кроме того, следует иметь в виду, что если на каком-то шаге
ответ уже получен, а ты все равно продолжаешь вычисления, увеличивая порядок
точности, то ты не сможешь получить другой ответ (если, конечно не делаешь
арифметических ошибок).

Выпишем получающиеся выражения:
 \begin{multline*}
\cos (x e^x)=\cos \left( x\cdot \left( 1+x+\underset{x\to
0}{\bold{o}}(x) \right) \right)=\\
=\cos \left( x+x^2+\underset{x\to 0}{\bold{o}}(x^2) \right)=\\=
1-\frac{1}{2}\left(x+x^2+\underset{x\to 0}{\bold{o}}(x^2)\right)^2
+\\
+\underset{x\to 0}{\bold{o}}\left(\left(x+x^2+\underset{x\to
0}{\bold{o}}(x^2)\right)^3 \right)=\\=...= 1-\frac{1}{2}x^2-x^3 +\underset{x\to
0}{\bold{o}}(x^3)
 \end{multline*}
 \begin{multline*}
\cos (x e^x) -\ln (1-x) -x =\\
=\left( 1-\frac{1}{2}x^2-x^3 +\underset{x\to 0}{\bold{o}}(x^3)
\right)-\\
-\left( -x-\frac{x^2}{2}-\frac{x^3}{3}+\underset{x\to 0}{\bold{o}}(x^3)
\right)-x=\\= 1-\frac{2x^3}{3}+\underset{x\to 0}{\bold{o}}(x^3)
\end{multline*}
$$
\kern-60pt\ctg x^3=\frac{1}{\tg x^3}=\frac{1}{x^3+\underset{x\to
0}{\bold{o}}(x^3)}
$$
Теперь вычисляем предел:
 \begin{multline*}
\lim_{x\to 0}\left(  \cos (x e^x) -\ln (1-x) -x \right)^{\ctg x^3}=\\
=\lim_{x\to 0}\left( 1-\frac{2x^3}{3}+\underset{x\to 0}{\bold{o}}(x^3)
\right)^{\frac{1}{x^3+\underset{x\to 0}{\bold{o}}(x^3)}}=\\= \lim_{x\to 0}
e^{\ln \left(\left( 1-\frac{2x^3}{3}+\underset{x\to 0}{\bold{o}}(x^3)
\right)^{\frac{1}{x^3+\underset{x\to 0}{\bold{o}}(x^3)}}\right)} =\\=\lim_{x\to
0} e^{ \frac{\ln \left(1-\frac{2x^3}{3}+\underset{x\to
0}{\bold{o}}(x^3)\right)} {x^3+\underset{x\to 0}{\bold{o}}(x^3)} } =\\=
\lim_{x\to 0} e^{ \frac{ \left(-\frac{2x^3}{3}+\underset{x\to
0}{\bold{o}}(x^3)\right) +\underset{x\to
0}{\bold{o}}\left(-\frac{2x^3}{3}+\underset{x\to 0}{\bold{o}}(x^3)\right) }
{x^3+\underset{x\to 0}{\bold{o}}(x^3)}
}=\\
=\lim_{x\to 0} e^{ \frac{ -\frac{2x^3}{3}+\underset{x\to 0}{\bold{o}}(x^3) }
{x^3+\underset{x\to 0}{\bold{o}}(x^3)} }= \lim_{x\to 0} e^{ \frac{-x^3\cdot
\frac{2}{3}+x^3\cdot \underset{x\to 0}{\bold{o}}(1)} {x^3+x^3\cdot
\underset{x\to 0}{\bold{o}}(1)}}=\\= \lim_{x\to 0} e^{
\frac{-\frac{2}{3}+\underset{x\to 0}{\bold{o}}(1)} {1+\underset{x\to
0}{\bold{o}}(1)}}= e^{\frac{-\frac{2}{3}+0}{1+0}}=
e^{-\frac{2}{3}}\end{multline*}\end{ex}

\begin{ers}
Найдите пределы:

1. $\lim\limits_{x\to 0}  \frac{\ln(1+x)-x}{x^2}$

2. $\lim\limits_{x\to 0}  \frac{e^x-1-x}{x^2}$

3. $\lim\limits_{x\to 0}  \frac{\cos x-1-\frac{x^2}{2}}{x^4}$

4. $\lim\limits_{x\to 0}  \frac{\arctg x-\arcsin x}{x^2}$

5. $\lim\limits_{x\to 0}  \frac{\tg x- x}{\sin x -x}$

6. $\lim\limits_{x\to 0}  \frac{\arctg x-\arcsin x}{\tg x -\sin x}$

7. $\lim\limits_{x\to 0}  \frac{2\arcsin x-\arcsin 2x}{x^3}$

8. $\lim\limits_{x\to 0}\frac{\sqrt[5]{1+2x}-1}{\sqrt[4]{1+x}-\sqrt{1-x}}$

9. $\lim\limits_{x\to 0}  \frac{1+x\cos x-\sqrt{1+2x}}{\ln (1+x)-x}$

10. $\lim\limits_{x\to 0}  \frac{e^x-\sqrt{1+2x}}{\ln \cos x}$

11. $\lim\limits_{x\to 0}  \frac{3\cos x +\arcsin x-3\sqrt[3] {1+x}} {\ln
(1-x^2)}$

12. $\lim\limits_{x\to 0}\frac{\ln(1+x^3)-2\sin x+2x\cos x^2}{\arctg x^3}$

13. $\lim\limits_{x\to 0}\frac{x\sqrt{1+\sin x}-\frac{1}{2}\ln(1+x^2)-x}{\tg^3
x}$

14. $\lim\limits_{x\to 0}\frac{e^{\sin x}-\sqrt{1+x^2}-x\cos x}{\ln^3(1-x)}$

15. $\lim\limits_{x\to 0}\frac{e^{\sin x \ln \cos
x}-\sqrt[4]{1+4x}+x-\frac{3}{2}x^2}{x\sin x^2}$

16. $\lim\limits_{x\to 0}\left(\frac{2\tg x}{x+\sin x}\right)^\frac{1}{1-\cos
x}$

17. $\lim\limits_{x\to 0}\left(\frac{x e^x +1}{x \pi^x+1}\right)^\frac{1}{x^2}$

18. $\lim\limits_{x\to 0}\left(\cos (\sin x)+\frac{1}{2}\arcsin^2 x \right)
^\frac{1}{x^2 (\sqrt{1+2x} -1)}$

19. $\lim\limits_{x\to 2}\left(\sqrt{3-x}+\ln {x}{2}\right) ^\frac{1}{\sin^2
(x-2)}$

20. $\lim\limits_{x\to 0}\left(\frac{1}{\sin x \arctg x}-\frac{1}{\tg x \arcsin
x}\right)$

21. $\lim\limits_{x\to 0}\frac{e^{\sin x}-1-\sin(e^x-1)}{\arctg^4 x}$

22. $\lim\limits_{x\to 0}\frac{\sin x^3}{\ln \cos x-\cos \ln (1+x)+1}$

23. $\lim\limits_{x\to 0}\frac{\sin((1+x)^\alpha-1)-(1+\sin x)^\alpha+1}{\tg^4
x}$
 \end{ers}

\end{multicols}\noindent\rule[10pt]{160mm}{0.1pt}

\subsection{Асимптотика}\label{SUBSEC:asimptotika}

Формула Тейлора-Пеано \eqref{Taylor-Peano}, с которой мы начинали этот
параграф, является частным случаем следующей общей конструкции.

 Пусть $a$ -- какое-нибудь число или символ бесконечности.
 \bit{
\item[$\bullet$] Последовательность функций $\{\ph_n\}$ называется {\it
асимптотической последовательностью} при $x\to a$, если
 \beq\label{DEF:asymp-posl}
\forall k\in\Z_+ \qquad \ph_k(x)\underset{x\to a}{\gg} \ph_{k+1}(x)
 \eeq
или, более наглядно,
$$
\ph_0(x)\underset{x\to a}{\gg} \ph_1(x)\underset{x\to a}{\gg} \ph_2(x)
\underset{x\to a}{\gg}...\underset{x\to a}{\gg} \ph_k(x)\underset{x\to
a}{\gg}...
$$

\item[$\bullet$] {\it Асимптотическим разложением} или {\it асимптотикой}
порядка $n$ функции $f$ вдоль асимптотической последовательности $\{\ph_n\}$
при $x\to a$, называется всякая асимптотическая формула вида
 \beq\label{DEF:asymp-razlozh}
f(x)=\sum_{k=0}^n \lambda_k\cdot\ph_k(x)+\underset{x\to
a}{\bold{o}}\Big(\ph_n(x)\Big),\qquad \lambda_k\in\R
 \eeq
(если она верна). Из \eqref{DEF:asymp-posl} следует, что коэффициенты
$\lambda_k\in\R$ в формуле \eqref{DEF:asymp-razlozh} определяются однозначно
(поэтому не бывает двух разных асимптотических разложений одного порядка
относительно данной асимптотической последовательности).
 }\eit

\noindent\rule{160mm}{0.1pt}\begin{multicols}{2}

\paragraph{Степенная последовательность в нуле.}
Функции
$$
\ph_k(x)=x^k,\qquad k\in\Z_+
$$
образуют асимптотическую последовательность при $x\to 0$, потому что
$$
1=x^0\underset{x\to 0}{\gg} x^1\underset{x\to 0}{\gg} x^2\underset{x\to
0}{\gg}...\underset{x\to 0}{\gg} x^k\underset{x\to 0}{\gg}...
$$
Для функции $f$, гладкой в окрестности точки $0$ разложения вдоль такой
последовательности -- это просто формулы Маклорена \eqref{11.9.3}, о которых мы
уже говорили выше:
$$
f(x)= \sum_{k=0}^n \frac{f^{(k)}(0)}{k!}\cdot x^k+ \underset{x\to
0}{\bold{o}}(x^n)
$$

\paragraph{Степенная последовательность на бесконечности.}

Функции
$$
\ph_k(x)=\frac{1}{x^k},\qquad k\in\Z_+
$$
образуют асимптотическую последовательность при $x\to\infty$, потому что
$$
1=\frac{1}{x^0}\underset{x\to\infty}{\gg}
\frac{1}{x^1}\underset{x\to\infty}{\gg}
\frac{1}{x^2}\underset{x\to\infty}{\gg}...\underset{x\to\infty}{\gg}
\frac{1}{x^k}\underset{x\to\infty}{\gg}...
$$
Асимптотические разложения вдоль этой последовательности обычно получаются
заменой $x=\frac{1}{t}$ с последующим применением формул Маклорена. Объясним
это на примерах.

\bex\label{EX:asymp-sqrt(x^2+1)} Пусть нам нужно найти асимптотическое
разложение порядка 3 при $x\to+\infty$ вдоль степенной последовательности для
функции
$$
f(x)=\sqrt{x^2+1}
$$
Сделаем замену $x=\frac{1}{t}$:
 \begin{multline*}
\sqrt{x^2+1}=\left|\scriptsize\begin{matrix}x=\frac{1}{t} \\
t\to+0\end{matrix}\right|=\sqrt{\frac{1}{t^2}+1}=\frac{\sqrt{1+t^2}}{t}=\\=\eqref{11.9.18}=\\=
\frac{1+\frac{1}{2}\cdot t^2+\frac{1}{2}\cdot\l\frac{1}{2}-1\r \cdot
\frac{t^4}{2!}+\underset{t\to+0}{\bold{o}}(t^4)}{t}=\\=
\frac{1}{t}+\frac{1}{2}\cdot t-\frac{1}{4}\cdot
\frac{t^3}{2!}+\underset{t\to+0}{\bold{o}}(t^3)=
\left|\scriptsize\begin{matrix}t=\frac{1}{x} \\
x\to+\infty\end{matrix}\right|=\\=
x+\frac{1}{2x}-\frac{1}{8x^3}+\underset{x\to+\infty}{\bold{o}}\l\frac{1}{x^3}\r
 \end{multline*}
\eex

\bex Найдем асимптотическое разложение порядка 1 при $x\to+\infty$ вдоль
степенной последовательности для функции
$$
f(x)=e^{\sqrt{x^2+x}-x}
$$
Опять делаем замену $x=\frac{1}{t}$:
 \begin{multline*}
e^{\sqrt{x^2+x}-x}=\left|\scriptsize\begin{matrix}x=\frac{1}{t} \\
t\to+0\end{matrix}\right|=e^{\sqrt{\frac{1}{t^2}+\frac{1}{t}}-\frac{1}{t}}=
e^{\frac{\sqrt{1+t}-1}{t}}=\\=\eqref{11.9.18}=\\= e^{\frac{1+\frac{1}{2}\cdot
t+\frac{1}{2}\cdot\l\frac{1}{2}-1\r \cdot
\frac{t^2}{2!}+\underset{t\to+0}{\bold{o}}(t^2)-1}{t}}=
e^{\frac{1}{2}-\frac{t}{8}+\underset{t\to+0}{\bold{o}}(t)}=\\=
e^{\frac{1}{2}}\cdot
e^{-\frac{t}{8}+\underset{t\to+0}{\bold{o}}(t)}=\eqref{11.9.17}=\\=
\sqrt{e}\cdot {\scriptsize\text{$\l
1-\frac{t}{8}+\underset{t\to+0}{\bold{o}}(t)+\underset{t\to+0}{\bold{o}}\l-\frac{t}{8}+\underset{t\to+0}{\bold{o}}(t)\r\r$}}=\\=
\sqrt{e}\cdot \l1-\frac{t}{8}+\underset{t\to+0}{\bold{o}}(t)\r=
\sqrt{e}-\frac{t\sqrt{e}}{8}+\underset{t\to+0}{\bold{o}}(t)=\\=
\sqrt{e}-\frac{\sqrt{e}}{8x}+\underset{x\to+\infty}{\bold{o}}\l\frac{1}{x}\r
\end{multline*} \eex

\paragraph{Другие асимптотические последовательности.}

Не всегда данная функция $f$ имеет асимптотическое разложение вдоль данной
асимптотической последовательности $\{\ph_n\}$.

\bex\label{EX:nevozm-razlozh-sqrt-3-(x)} Формула
$$
\sqrt[3]{x}=\sum_{k=0}^n \lambda_k\cdot x^k+\underset{x\to
0}{\bold{o}}\Big(x^n\Big)
$$
не будет верна ни при каких $n>0$ и $\lambda_0,...,\lambda_n\in\R$.
 \eex
 \bpr
Это достаточно доказать для $n=1$, потому что любую формулу большего порядка
можно записать в виде
$$
\sqrt[3]{x}=\lambda_0+\lambda_1\cdot x+\underbrace{\sum_{k=2}^n \lambda_k\cdot
x^k+\underset{x\to 0}{\bold{o}}\Big(x^n\Big)}_{\underset{x\to
0}{\bold{o}}\Big(x\Big)}
$$
и если эта формула верна, то будет верна и формула первого порядка:
 \beq\label{nevozm-razlozh-sqrt-3-(x)}
\sqrt[3]{x}=\lambda_0+\lambda_1\cdot x+\underset{x\to 0}{\bold{o}}\Big(x\Big)
 \eeq
Чтобы убедиться, что эта формула неверна, нужно с самого начала заметить, что
из нее следует, что $\lambda_0$ должно быть нулевым, потому что при $x\to 0$
получаем:
$$
\underbrace{\sqrt[3]{x}}_{\scriptsize\begin{matrix}\downarrow\\
0\end{matrix}}=\lambda_0+\underbrace{\lambda_1\cdot x}_{\scriptsize\begin{matrix}\downarrow\\
0\end{matrix}}+\underbrace{\underset{x\to 0}{\bold{o}}\Big(x\Big)}_{\scriptsize\begin{matrix}\downarrow\\
0\end{matrix}}
$$
Подставим $\lambda_0=0$ в \eqref{nevozm-razlozh-sqrt-3-(x)}:
$$
\sqrt[3]{x}=\lambda_1\cdot x+\underset{x\to 0}{\bold{o}}\Big(x\Big)
$$
Поделив все на $x$, при $x\to 0$ мы получим соотношение
$$
\underbrace{\frac{1}{\sqrt[3]{x^2}}}_{\scriptsize\begin{matrix}\downarrow\\
\infty \end{matrix}}=\lambda_1+\underbrace{\underset{x\to 0}{\bold{o}}\Big(1\Big)}_{\scriptsize\begin{matrix}\downarrow\\
0\end{matrix}}
$$
которое не может быть верно ни при каком $\lambda_1$. \epr

Приведенный пример \ref{EX:nevozm-razlozh-sqrt-3-(x)} имеет целью подготовить
читателя к мысли, что {\it асимптотическую последовательность $\{\ph_n\}$,
вдоль которой мы собираемся раскладывать данную функцию $f$, бывает удобно
выбирать отдельно, в зависимости от вида функции $f$}.

\bex Асимптотическое разложение функции
$$
f(x)=\sqrt[3]{x+x^2}
$$
при $x\to 0$ бессмысленно искать вдоль стандартной степенной последовательности
$$
\ph_k(x)=x^k
$$
потому что здесь возникает тот же эффект, что в примере
\ref{EX:nevozm-razlozh-sqrt-3-(x)}: никакая формула вида
\eqref{DEF:asymp-razlozh} не будет верна, если $n>0$.

Однако если взять последовательность
$$
x^{\frac{1}{3}}\underset{x\to 0}{\gg} x^2\underset{x\to 0}{\gg}
x^{\frac{11}{3}}\underset{x\to 0}{\gg}...\underset{x\to 0}{\gg}
x^{\frac{5k+1}{3}}\underset{x\to 0}{\gg}...
$$
то этот эффект исчезает, и разложение находится очень просто:
 \begin{multline*}
\sqrt[3]{x+x^2}=x^{\frac{1}{3}}\cdot
(1+x^{\frac{5}{3}})^{\frac{1}{3}}=\eqref{11.9.18}=\\= x^{\frac{1}{3}}\cdot
\text{\scriptsize $\l 1+\frac{1}{3}\cdot
x^{\frac{5}{3}}+\frac{1}{3}\cdot\l\frac{1}{3}-1\r\cdot\frac{\l
x^{\frac{5}{3}}\r^2}{2!}+ \underset{x\to 0}{\bold{o}\l \l
x^{\frac{5}{3}}\r^2\r}\r$}=
\\= x^{\frac{1}{3}}\cdot \l 1+\frac{x^{\frac{5}{3}}}{3}
-\frac{x^{\frac{10}{3}}}{9}+\underset{x\to 0}{\bold{o}\l
x^{\frac{10}{3}}\r}\r=\\= x^{\frac{1}{3}}+\frac{x^2}{3}
-\frac{x^{\frac{11}{3}}}{9}+\underset{x\to 0}{\bold{o}\l x^{\frac{11}{3}}\r}
 \end{multline*}
(одновременно из этих вычислений становится понятно, откуда берется
последовательность $\ph_k(x)=x^{\frac{5k+1}{3}}$). \eex

\bex Попробуем найти разложение той же функции
$$
f(x)=\sqrt[3]{x+x^2}
$$
при $x\to\infty$. Как и в предыдущих случаях, когда ищется асимптотика на
бесконечности, начнем с замены $x=\frac{1}{t}$:
 \begin{multline*}
\sqrt[3]{x+x^2}=\left|\scriptsize\begin{matrix}x=\frac{1}{t} \\
t\to 0\end{matrix}\right|=\sqrt[3]{\frac{1}{t}+\frac{1}{t^2}}=
\frac{(1+t)^\frac{1}{3}}{t^{\frac{2}{3}}}=\\=\eqref{11.9.18}=
\frac{1+\frac{1}{3}\cdot t+\frac{1}{3}\cdot\l\frac{1}{3}-1\r\cdot
\frac{t^2}{2!}+ \underset{t\to 0}{\bold{o}\l t^2\r}}{t^{\frac{2}{3}}}=\\=
t^{-\frac{2}{3}}+\frac{t^{\frac{1}{3}}}{3}-\frac{t^{\frac{4}{3}}}{9}+
\underset{t\to 0}{\bold{o}\l t^{\frac{4}{3}}\r}=\left|\scriptsize\begin{matrix}t=\frac{1}{x} \\
x\to\infty\end{matrix}\right|=\\= x^{\frac{2}{3}}+\frac{1}{3\cdot
x^{\frac{1}{3}}}-\frac{1}{9\cdot x^{\frac{4}{3}}}+
\underset{x\to\infty}{\bold{o}\l \frac{1}{x^{\frac{4}{3}}}\r}
 \end{multline*}
Из вычислений видно, что разложение в этом случае удобно искать вдоль
асимптотической последовательности
$$
x^{\frac{2}{3}}\underset{x\to\infty}{\gg}
x^{-\frac{1}{3}}\underset{x\to\infty}{\gg}
x^{-\frac{4}{3}}\underset{x\to\infty}{\gg}...\underset{x\to\infty}{\gg}
x^{\frac{2}{3}-k}\underset{x\to\infty}{\gg}...
$$
 \eex

\end{multicols}\noindent\rule[10pt]{160mm}{0.1pt}

\section{Асимптотика интегралов и сумм с переменным пределом}\label{SEC:asymp-integr}

\subsection{Асимптотика интегралов}\label{SUBSEC:asymp-integr}

\paragraph{Асимптотическая эквивалентность интегралов}

\begin{tm}
[\bf признак эквивалентности несобственных интегралов]\label{tm-17.5.7}

Пусть функции $f$ и $g$ обладают следующими свойствами:
 \bit{
\item[a)] $f$ и $g$ непрерывны на полуинтервале $[a;c)$;

\item[b)] $f(x)\underset{x\to c-0}{\sim} g(x)$;

\item[c)] $g(x)>0, \quad \forall x\in (a;c)$.
 }\eit
\noindent Соответствующая зависимость между несобственными интегралами коротко
записывается формулой:
$$
\int_a^c f(x) \, \d x \underset{c}{\sim}\int_a^c g(x) \, \d x
$$
Тогда сходимость интеграла $\int_a^c f(x) \, \d x$ эквивалентна сходимости
интеграла $\int_a^c g(x) \, \d x$:
 \beq\label{asymp-priznak-dlya-int}
\int_a^c f(x) \, \d x \quad \text{сходится}\quad \Longleftrightarrow \quad
\int_a^c g(x) \, \d x \quad \text{сходится},
 \eeq
При этом,
 \bit{
\item[(i)] если эти интегралы сходятся, то справедливо следующее соотношение,
показывающее, что скорость их сходимости будет одинаковой:
 \beq\label{asymp-ekviv-int-0}
\int_T^c f(x) \, \d x\underset{T\to c-0}{\sim}\int_T^c g(x) \, \d x
 \eeq

\item[(ii)] если эти интегралы расходятся, то справедливо следующее
соотношение, показывающее, что скорость их расходимости будет одинаковой:
 \beq\label{asymp-ekviv-int}
\int_a^T f(x) \, \d x\underset{T\to c-0}{\sim}\int_a^T g(x) \, \d x
 \eeq
}\eit
\end{tm}

\noindent\rule{160mm}{0.1pt}\begin{multicols}{2}

\brem Если интегралы \eqref{asymp-priznak-dlya-int} расходятся, то формула
\eqref{asymp-ekviv-int-0} утрачивает смыл, и поэтому перестает быть верной. А
если интегралы \eqref{asymp-priznak-dlya-int} сходятся, то формула
\eqref{asymp-ekviv-int}, хотя и сохраняет смысл, но также перестает работать.
Например,
$$
\frac{x-1}{x^3}\underset{x\to+\infty}{\sim}\frac{1}{x^2}
$$
однако
 \begin{multline*}
\int_1^T\frac{x-1}{x^3}\ \d x=\int_1^T\l\frac{1}{x^2}-\frac{1}{x^3}\r\ \d x=\\=
\l-\frac{1}{x}+\frac{1}{2x^2}\r\Bigg|_1^T=-\frac{1}{T}+1+\frac{1}{2T^2}-\frac{1}{2}
\underset{T\to+\infty}{\sim}\\
\underset{T\to+\infty}{\sim}\frac{1}{2}\underset{T\to+\infty}{\not\sim}
1\underset{T\to+\infty}{\sim}1-\frac{1}{T}=-\frac{1}{x}\
\bigg|_1^T=\\=\int_1^T\frac{1}{x^2}\ \d x
 \end{multline*}
То есть
$$
\int_1^T\frac{x-1}{x^3}\ \d
x\underset{T\to+\infty}{\not\sim}\int_1^T\frac{1}{x^2}\ \d x
$$
Более того, можно заметить, что если одна из функций не превосходит другую на
интервале интегрирования, например,
$$
f(x)\le g(x),\qquad x\in[a,c)
$$
то соотношение \eqref{asymp-ekviv-int} для сходящихся интегралов возможно
только, если подынтегральные функции совпадают. Действительно, из неравенства
$g(x)>0$ следует, что интеграл от $g$ не может быть нулевым,
$$
\int_a^c g(x)\ \d x=\lim_{T\to c-0}\int_a^T g(x) \, \d x\ne 0
$$
Поэтому соотношение \eqref{asymp-ekviv-int}, если оно выполняется, означает
просто равенство пределов:
$$
\lim_{T\to c-0}\int_a^T f(x) \, \d x=\lim_{T\to c-0}\int_a^T g(x) \, \d x
$$
Отсюда
 \begin{multline*}
\lim_{T\to c-0}\int_a^T
\overbrace{\Big(g(x)-f(x)\Big)}^{\scriptsize\begin{matrix}0\\
\text{\rotatebox{90}{$\ge$}}
\end{matrix}} \,
\d x=\\=\lim_{T\to c-0}\int_a^T g(x) \, \d x-\lim_{T\to c-0}\int_a^T f(x) \, \d
x=0
 \end{multline*}
$$
\Downarrow
$$
$$
f(x)=g(x),\qquad x\in[a,c)
$$

 \erem

\end{multicols}\noindent\rule[10pt]{160mm}{0.1pt}

Для доказательства теоремы \ref{tm-17.5.7} нам понадобится следующая

\begin{lm}\label{lm-17.5.8} Пусть функция $f$ непрерывна на полуинтервале
$[a;c)$ и пусть $b\in [a;c)$. Тогда сходимость несобственного интеграла от $f$
на промежутке $[a;c)$ эквивалентна сходимости на промежутке $[b;c)$:
$$
\int_a^c f(x) \, \d x \quad \text{сходится}\quad \Longleftrightarrow \quad
\int_b^c f(x) \, \d x \quad \text{сходится}
$$
\end{lm}\begin{proof} Это следует из того, что интегралы
по промежуткам $[a;c)$ и $[b;c)$ отличаются друг от друга на константу
\begin{multline*}\int_a^c f(x) \, \d x=\lim_{t\to c}\int_a^t f(x) \, \d x=
\lim_{t\to c}\l \int_a^b f(x) \, \d x+\int_b^t f(x) \, \d x \r=\\= \int_a^b
f(x) \, \d x+\lim_{t\to c}\int_b^t f(x) \, \d x = \int_a^b f(x) \, \d
x+\int_b^c f(x) \, \d x
 \end{multline*}
 \end{proof}

\begin{proof}[Доказательство теоремы \ref{tm-17.5.7}].
1. Сначала докажем \eqref{asymp-priznak-dlya-int}. Условие (b) означает, что
 \beq\label{asymp-prizn-ekviv:f(x)/g(x)_>1}
 \frac{f(x)}{g(x)}\underset{x\to c-0}{\longrightarrow} 1
 \eeq
По определению предела Коши (с.\pageref{DEF:Cauchy-lim(x->a-0)}), существует
$b\in [a;c)$ такое что $\forall x\in [b;c) \quad \frac{f(x)}{g(x)}\in
(\frac{1}{2};\frac{3}{2})$, то есть
$$
\forall x\in [b;c) \quad \frac{1}{2}<\frac{f(x)}{g(x)}<\frac{3}{2}
$$
Поскольку $g(x)\ge 0$, можно умножить это двойное неравенство на $g(x)$, и при
этом знаки неравенства не изменятся:
$$
\forall x\in [b;c) \quad \frac{1}{2}g(x)<f(x)<\frac{3}{2}g(x)
$$
Перепишем это иначе:
 \beq
\forall x\in [b;c) \quad
\begin{cases}{g(x)<2f(x)}\\{f(x)<\frac{3}{2}g(x)}\end{cases}\label{17.5.6}
 \eeq
Отсюда получаются две цепочки следствий. Во-первых,
$$
  \text{интеграл}\,\, \int_a^c f(x) \, \d x  \,\,\text{сходится, то}
$$
$$
\Downarrow\put(20,0){\smsize (\text{применяем лемму \ref{lm-17.5.8}})}
$$
$$
  \text{интеграл}\,\, \int_b^c f(x) \, \d x  \,\,\text{сходится}
$$
$$
\Downarrow
$$
$$
  \text{интеграл}\,\, \int_b^c 2f(x) \, \d x  \,\,\text{сходится}
$$
$$
\Downarrow\put(20,0){\smsize $\begin{pmatrix} \text{вспоминаем первое
неравенство}\\ \text{в \eqref{17.5.6} и теорему \ref{tm-17.5.1}}\end{pmatrix}$}
$$
$$
  \text{интеграл}\,\, \int_b^c g(x) \, \d x  \,\,\text{сходится}
$$
$$
\Downarrow\put(20,0){\smsize (\text{применяем лемму \ref{lm-17.5.8}})}
$$
$$
  \text{интеграл}\,\, \int_a^c g(x) \, \d x  \,\,\text{сходится}
$$
И, во-вторых,
$$
  \text{интеграл}\,\, \int_a^c g(x) \, \d x  \,\,\text{сходится, то}
$$
$$
\Downarrow\put(20,0){\smsize (\text{применяем лемму \ref{lm-17.5.8}})}
$$
$$
  \text{интеграл}\,\, \int_b^c g(x) \, \d x  \,\,\text{сходится}
$$
$$
\Downarrow
$$
$$
  \text{интеграл}\,\, \int_b^c \frac{3}{2} g(x) \, \d x  \,\,\text{сходится}
$$
$$
\Downarrow\put(20,0){\smsize $\begin{pmatrix}\text{вспоминаем второе
неравенство}\\ \text{в \eqref{17.5.6} и теорему \ref{tm-17.5.1}}\end{pmatrix}$}
$$
$$
  \text{интеграл}\,\, \int_b^c f(x) \, \d x  \,\,\text{сходится}
$$
$$
\Downarrow\put(20,0){\smsize (\text{применяем лемму \ref{lm-17.5.8}})}
$$
$$
  \text{интеграл}\,\, \int_a^c f(x) \, \d x  \,\,\text{сходится}
$$
Мы доказали \eqref{asymp-priznak-dlya-int}.

2. Остается доказать \eqref{asymp-ekviv-int-0} и \eqref{asymp-ekviv-int}. Обе
эти формулы доказываются с помощью правила Лопиталя. Сначала предположим, что
оба интеграла сходятся, и обозначим
$$
F(T)=\int_T^c f(x) \, \d x=\int_a^c f(x) \, \d x-\int_a^T f(x) \, \d x,\qquad
G(T)=\int_T^c g(x) \, \d x=\int_a^c g(x) \, \d x-\int_a^T g(x) \, \d x
$$
Тогда
$$
F'(T)=-f(T),\qquad G'(T)=-g(T)\qquad \Big(T\in(a,c)\Big)
$$
и поэтому
$$
\lim_{T\to c-0}\frac{\int_T^c f(x) \, \d x}{\int_T^c g(x) \, \d
x}=\l\frac{0}{0}\r=\eqref{9.1.1}=\lim_{T\to c-0}\frac{-f(T)}{-g(T)}=\lim_{T\to
c-0}\frac{f(T)}{g(T)}=\eqref{asymp-prizn-ekviv:f(x)/g(x)_>1}=1
$$
Наоборот, предположим, что оба интеграла расходятся, и обозначим
$$
F(T)=\int_a^T f(x) \, \d x,\qquad G(T)=\int_a^T g(x) \, \d x
$$
Тогда
$$
F'(T)=f(T),\qquad G'(T)=g(T)\qquad \Big(T\in(a,c)\Big)
$$
и поэтому
$$
\lim_{T\to c-0}\frac{\int_a^T f(x) \, \d x}{\int_a^T g(x) \, \d
x}=\l\frac{\infty}{\infty}\r=\eqref{9.2.1}=\lim_{T\to
c-0}\frac{f(T)}{g(T)}=\eqref{asymp-prizn-ekviv:f(x)/g(x)_>1}=1
$$

\end{proof}

\noindent\rule{160mm}{0.1pt}\begin{multicols}{2}

\bex Чтобы понять, сходится ли интеграл
$$
\int_1^{+\infty}\frac{\sqrt{x^2+1}}{x^2} \ \d x
$$
заметим, что подынтегральная функция положительна, и подберем к ней
эквивалентную:
$$
\frac{\sqrt{x^2+1}}{x^2}\underset{x\to+\infty}{\sim} \frac{x}{x^2}=\frac{1}{x}
$$
Теперь по теореме \ref{tm-17.5.7} получаем, что сходимость нашего интеграла
эквивалентна сходимости более простого интеграла:
 \begin{multline*}
\int_1^{+\infty}\frac{\sqrt{x^2+1}}{x^2} \ \d x
\underset{x\to+\infty}{\sim}\int_1^{+\infty} \frac{1}{x} \ \d x=\\=\ln
x\Bigg|_1^{+\infty}=+\infty
 \end{multline*}
Поскольку второй интеграл расходится, мы получаем, что наш исходный интеграл
тоже расходится:
$$
\int_1^{+\infty}\frac{\sqrt{x^2+1}}{x^2} \ \d x=\infty
$$
Из формулы \eqref{asymp-ekviv-int} можно вывести скорость расходимости:
$$
\int_1^T\frac{\sqrt{x^2+1}}{x^2} \ \d x\underset{T\to+\infty}{\sim}
\int_1^T\frac{1}{x} \ \d x=\ln T
$$
\eex

\begin{ex}\label{ex-17.5.9} Рассмотрим интеграл
$$
  \int_0^1 \frac{\sin \sqrt{x}}{e^x-1}\, \d x
$$
Здесь особой точкой ($x=0$) будет нижняя граница интегрирования, но теорема
\ref{tm-17.5.7} применима и к этому случаю, в очевидным образом
переформулированном виде. Опять подынтегральная функция положительна. Подберем
к ней эквивалентную функцию при стремлении $x$ к особой точке:
$$
\frac{\sin \sqrt{x}}{e^x-1}\underset{x\to
0}{\sim}\frac{\sqrt{x}}{x}=\frac{1}{\sqrt{x}}
$$
Мы получаем по теореме \ref{tm-17.5.7}, что сходимость исходного интеграла
эквивалентна сходимости более простого интеграла
$$
\int_0^1 \frac{\sin \sqrt{x}}{e^x-1}\, \d x \underset{0}{\sim}\int_0^1
\frac{1}{\sqrt{x}}\, \d x=2\sqrt{x}\ \Big|_0^1
$$
Второй интеграл сходится, значит, наш исходный интеграл тоже сходится. Его
скорость сходимости описывается формулой \eqref{asymp-ekviv-int-0}:
$$
\int_0^T \frac{\sin \sqrt{x}}{e^x-1}\, \d x \underset{T\to +0}{\sim}\int_0^T
\frac{1}{\sqrt{x}}\, \d x=2\sqrt{T}.
$$
\end{ex}

\begin{ex}\label{ex-17.5.10} Те же самые рассуждения для интеграла
$$
  \int_0^1 \frac{1-\cos x}{\ln (1+x^3)}\, \d x
$$
запишем короче:
 \begin{multline*}
\int_0^1 \frac{1-\cos x}{\ln (1+x^3)}\, \d x \underset{0}{\sim}\int_0^1
\frac{\frac{x^2}{2}}{x^3}\, \d x=\\= \frac{1}{2}\int_0^1 \frac{1}{x}\, \d x=\ln
x\ \Big|_0^1
 \end{multline*}
Второй интеграл расходится, значит, наш исходный интеграл тоже расходится.
Скорость расходимости описывается формулой
$$
\int_T^1 \frac{1-\cos x}{\ln (1+x^3)}\, \d x \underset{T\to
+0}{\sim}\frac{1}{2}\int_T^1 \frac{1}{x}\, \d x=-\ln T
$$
\end{ex}

\begin{ex}\label{ex-17.5.11} Снова интеграл по бесконечному промежутку:
 \begin{multline*}
\int_1^{+\infty}\frac{x+\sin x}{\sqrt{1+x^6}}\, \d x \underset{x\to
+\infty}{\sim}\int_1^{+\infty}\frac{x}{\sqrt{x^6}}\, \d x=\\=
\int_1^{+\infty}\frac{x}{x^3}\, \d x= \int_1^{+\infty}\frac{1}{x^2}\, \d
x=-\frac{1}{x}\ \Big|_1^{+\infty}
 \end{multline*}
Последний интеграл сходится, значит, наш исходный интеграл тоже сходится.
Скорость сходимости:
$$
\int_T^{+\infty}\frac{x+\sin x}{\sqrt{1+x^6}}\, \d x \underset{T\to
+\infty}{\sim}\int_T^{+\infty}\frac{1}{x^2}\, \d x=\frac{1}{T}
$$

\end{ex}

\begin{ex}\label{ex-17.5.12} Еще один интеграл по бес\-ко\-неч\-ному промежутку:
 \begin{multline*}
\int_1^{+\infty}\frac{x+\sin x}{\sqrt{1+x^3}}\, \d x \underset{x\to
+\infty}{\sim}\int_1^{+\infty}\frac{x}{\sqrt{x^3}}\, \d x=\\=
\int_1^{+\infty}\frac{1}{\sqrt{x}}\, \d x=2\sqrt{x}\ \Big|_1^{+\infty}
 \end{multline*}
Последний интеграл расходится, значит, наш исходный интеграл тоже расходится.
Скорость расходимсоти:
 \begin{multline*}
\int_1^T\frac{x+\sin x}{\sqrt{1+x^3}}\, \d x \underset{T\to+\infty}{\sim}
\int_1^T\frac{1}{\sqrt{x}}\, \d x=\\=2\sqrt{x}\
\Big|_1^T=2\sqrt{T}-2\underset{T\to+\infty}{\sim}2\sqrt{T}
 \end{multline*}

\end{ex}

\begin{ers} Исследуйте на сходимость интегралы и оцените скорость сходимости и расходимости:
 \biter{
\item[1)] $\int_2^{+\infty}\frac{1}{1+\sqrt{x}}\, \d x$;

\item[2)] $\int_2^{+\infty}\frac{1}{1+\sqrt{x}+x^5}\, \d x$;

\item[3)] $\int_{-\infty}^{-2}\frac{1}{1+x^5}\, \d x$;

\item[4)] $\int_0^1 \frac{\arctg x}{x\sqrt[3]{x}}\, \d x$;

\item[5)] $\int_1^{+\infty}\frac{\arctg x}{x\sqrt[3]{x}}\, \d x$;

\item[6)] $\int_0^1 \frac{\arctg x}{\sqrt[3]{x}}\, \d x$;

\item[7)] $\int_1^{+\infty}\frac{\arctg x}{\sqrt[3]{x}}\, \d x$;

\item[8)] $\int_0^1 \frac{e^x-1}{x^3}\, \d x$;

\item[9)] $\int_{-1}^0 \frac{e^x-1}{x^3}\, \d x$;

\item[10)] $\int_0^1 \frac{1}{\sqrt{x}-1}\, \d x$;

\item[11)] $\int_0^1 \frac{1}{\sqrt{x-1}}\, \d x$;

\item[12)] $\int_0^1 \frac{1}{e^x-\cos x}\, \d x$;

\item[13)] $\int_0^1 \frac{\ln (1+\sqrt[3] {x^2})}{e^x-1}\, \d x$;

\item[14)] $\int_e^{+\infty}\frac{\ln (1+\frac{1}{x})}{x+\ln x}\, \d x$;

\item[15)] $\int_e^{+\infty}\frac{\ln (1+\frac{1}{x}+x)}{1+\ln x+\sin x}\, \d
x$.
 }\eiter
 \end{ers}

\end{multicols}\noindent\rule[10pt]{160mm}{0.1pt}

\paragraph{Асимптотическое сравнение интегралов}

\begin{tm}[\bf асимптотический признак сравнения]\label{tm-17.7.5}
Пусть функции $f$ и $g$ обладают следующими свойствами:
 \bit{
\item[a)] $f$ и $g$ непрерывны на полуинтервале $[a;c)$;

\item[b)] $f(x)\underset{x\to c-0}{\ll} g(x)$ (то есть $f(x)=\underset{x\to
c-0}{\bold{o}}\l g(x) \r$);

\item[c)] $g(x)>0, \quad \forall x\in (a;c)$.
 }\eit
 \bit{\rm
\item[$\bullet$] Возникающая зависимость между несобственными интегралами
коротко записывается формулой
$$
\int_a^c f(x) \, \d x \underset{c}{\ll}\int_a^c g(x) \, \d x.
$$
 }\eit
Тогда
 \bit{
\item[(i)] из сходимости большего интеграла $\int_a^c g(x) \, \d x$ следует
сходимость меньшего интеграла $\int_a^c f(x) \, \d x$:
 \beq\label{asymp-sravn-priznak-dlya-int-1}
\int_a^c f(x) \, \d x \quad \text{сходится}\quad \Longleftarrow \quad \int_a^c
g(x) \, \d x \quad \text{сходится},
 \eeq
причем скорости сходимости оцениваются соотношением
 \beq\label{asymp-sravn-int-0}
\int_T^c f(x) \, \d x \underset{T\to c-0}{\ll}\int_T^c g(x) \, \d x;
 \eeq

\item[(ii)] из расходимости меньшего интеграла $\int_a^c f(x) \, \d x$ следует
расходимость большего интеграла $\int_a^c g(x) \, \d x$:
 \beq\label{asymp-sravn-priznak-dlya-int-2}
\int_a^c f(x) \, \d x \quad \text{расходится}\quad \Longrightarrow \quad
\int_a^c g(x) \, \d x \quad \text{расходится}
 \eeq
причем скорости расходимости оцениваются соотношением
 \beq\label{asymp-sravn-int}
\int_a^T f(x) \, \d x \underset{T\to c-0}{\ll}\int_a^T g(x) \, \d x;
 \eeq
 }\eit
\end{tm}\begin{proof}
1. Забудем на время о формулах \eqref{asymp-sravn-int-0} и
\eqref{asymp-sravn-int}. Оставшиеся в условиях (i) и (ii) утверждения
\eqref{asymp-sravn-priznak-dlya-int-1} и \eqref{asymp-sravn-priznak-dlya-int-2}
эквивалентны, поэтому из них достаточно доказать какое-нибудь одно, например,
первое. Условие (b) означает, что
 \beq\label{asim-prizn-sravn:f(x)/g(x)->0}
 \frac{f(x)}{g(x)}\underset{x\to c-0}{\longrightarrow} 0.
 \eeq
По определению предела Коши (с.\pageref{DEF:Cauchy-lim(x->a-0)}), существует $b\in [a;c)$
такое что $\forall x\in [b;c) \quad \frac{f(x)}{g(x)}\in
(-\frac{1}{2};+\frac{1}{2})$, то есть
$$
\forall x\in [b;c) \quad -\frac{1}{2}<\frac{f(x)}{g(x)}<+\frac{1}{2}
$$
Поскольку $g(x)>0$, можно умножить это двойное неравенство на $g(x)$, и при
этом знаки неравенства не изменятся:
$$
\forall x\in [b;c) \quad -\frac{1}{2}\ g(x)<f(x)<\frac{1}{2}\ g(x)
$$
Перепишем это иначе:
 \beq
\forall x\in [b;c) \quad |f(x)|<\frac{1}{2}\ g(x) \label{17.7.1}
 \eeq
Отсюда получается следующая цепочка следствий:
$$
  \text{интеграл}\,\, \int_a^c g(x) \, \d x  \,\,\text{сходится}
$$
$$
\Downarrow\put(20,0){\smsize \text{лемма \ref{lm-17.5.8}}}
$$
$$
  \text{интеграл}\,\, \int_b^c g(x) \, \d x  \,\,\text{сходится}
$$
$$
\Downarrow
$$
$$
  \text{интеграл}\,\, \int_b^c \frac{1}{2}\ g(x) \, \d x  \,\,\text{сходится}
$$
$$
\Downarrow\put(20,0){\smsize \text{неравенство \eqref{17.7.1} и теорема
\ref{tm-17.5.1}}}
$$
$$
  \text{интеграл}\,\, \int_b^c |f(x)| \, \d x  \,\,\text{сходится}
$$
$$
\Downarrow\put(20,0){\smsize \text{теорема \ref{tm-17.7.1}}}
$$
$$
  \text{интеграл}\,\, \int_b^c f(x) \, \d x  \,\,\text{сходится}
$$

2. Теперь вспомним о формулах \eqref{asymp-sravn-int-0} и
\eqref{asymp-sravn-int}. Они доказывается тем же приемом (с помощью правила
Лопиталя), что и формулы \eqref{asymp-ekviv-int-0}-\eqref{asymp-ekviv-int}
выше. Сначала предположим, что оба интеграла сходятся, и обозначим
$$
F(T)=\int_T^c f(x) \, \d x=\int_a^c f(x) \, \d x-\int_a^T f(x) \, \d x,\qquad
G(T)=\int_T^c g(x) \, \d x=\int_a^c g(x) \, \d x-\int_a^T g(x) \, \d x
$$
Тогда
$$
F'(T)=-f(T),\qquad G'(T)=-g(T)\qquad \Big(T\in(a,c)\Big)
$$
и поэтому
$$
\lim_{T\to c-0}\frac{\int_T^c f(x) \, \d x}{\int_T^c g(x) \, \d
x}=\l\frac{0}{0}\r=\eqref{9.1.1}=\lim_{T\to c-0}\frac{-f(T)}{-g(T)}=\lim_{T\to
c-0}\frac{f(T)}{g(T)}=\eqref{asim-prizn-sravn:f(x)/g(x)->0}=0
$$
Наоборот, предположим, что оба интеграла расходятся, и обозначим
$$
F(T)=\int_a^T f(x) \, \d x,\qquad G(T)=\int_a^T g(x) \, \d x
$$
Тогда
$$
F'(T)=f(T),\qquad G'(T)=g(T)\qquad \Big(T\in(a,c)\Big)
$$
и поэтому
$$
\lim_{T\to c-0}\frac{\int_a^T f(x) \, \d x}{\int_a^T g(x) \, \d
x}=\l\frac{\infty}{\infty}\r=\eqref{9.2.1}=\lim_{T\to
c-0}\frac{f(T)}{g(T)}=\eqref{asim-prizn-sravn:f(x)/g(x)->0}=0
$$
\end{proof}

\noindent\rule{160mm}{0.1pt}\begin{multicols}{2}

\begin{ex}\label{ex-17.7.6} Чтобы исследовать на сходимость несобственный
интеграл
$$
  \int_0^{+\infty}\frac{x^{10}\sin x}{e^x}\, \d x
$$
заметим, что
$$
\frac{x^{10}\sin x}{e^x}\underset{x\to
+\infty}{\ll}\frac{1}{e^{\frac{x}{2}}}=e^{-\frac{x}{2}}
$$
Поэтому
$$
\int_0^{+\infty}\frac{x^{10}\sin x}{e^x}\, \d x
\underset{+\infty}{\ll}\int_0^{+\infty}e^{-\frac{x}{2}}  \, \d
x=-2e^{-\frac{x}{2}}\ \Big|_0^{+\infty}
$$
Поскольку больший интеграл (справа от $\ll$) сходится, меньший (то есть наш
исходный) интеграл тоже должен сходиться. Его скорость сходимости оценивается
формулой \eqref{asymp-sravn-int-0}:
$$
\int_T^{+\infty}\frac{x^{10}\sin x}{e^x}\, \d x
\underset{T\to+\infty}{\ll}\int_T^{+\infty}e^{-\frac{x}{2}}  \, \d
x=2e^{-\frac{T}{2}}
$$
или, если упростить,
$$
\int_T^{+\infty}\frac{x^{10}\sin x}{e^x}\, \d x
\underset{T\to+\infty}{\ll}e^{-\frac{T}{2}}
$$
Применяя оценку
$$
x^{10}\underset{x\to +\infty}{\ll} e^{\e x}\qquad (\e>0),
$$
можно доказать более точную формулу:
$$
\int_T^{+\infty}\frac{x^{10}\sin x}{e^x}\, \d x \underset{T\to+\infty}{\ll}
a^T\qquad (a<e).
$$
\end{ex}

\begin{ex}\label{ex-17.7.7} Чтобы исследовать на сходимость несобственный
интеграл
$$
  \int_1^{+\infty}\frac{e^x}{x^{10}}\, \d x
$$
заметим, что
$$
\frac{e^x}{x^{10}}\underset{x\to +\infty}{\gg} e^{\frac{x}{2}}
$$
Поэтому
$$
\int_1^{+\infty}\frac{e^x}{x^{10}}\, \d x
\underset{+\infty}{\gg}\int_1^{+\infty} e^{\frac{x}{2}}\, \d
x=2e^{\frac{x}{2}}\ \Big|_1^{+\infty}
$$
Поскольку меньший интеграл расходится, больший тоже должен расходиться.
Скорость расходимости оценивается формулой
 \begin{multline*}
\int_1^T\frac{e^x}{x^{10}}\, \d x \underset{T\to+\infty}{\gg}\int_1^T
e^{\frac{x}{2}}\, \d
x=\\=2e^{\frac{T}{2}}-2e^{\frac{1}{2}}\underset{T\to+\infty}{\sim}
2e^{\frac{T}{2}}
 \end{multline*}
и, если упростить,
$$
\int_1^T\frac{e^x}{x^{10}}\, \d x \underset{T\to+\infty}{\gg} e^{\frac{T}{2}}
$$
Применяя оценку
$$
x^{10}\underset{x\to +\infty}{\ll} e^{\e x}\qquad (\e>0),
$$
можно доказать формулу
$$
\int_1^T\frac{e^x}{x^{10}}\, \d x \underset{T\to+\infty}{\gg} a^T\qquad (a<e)
$$
\end{ex}

\begin{ex}\label{ex-17.7.8}
 \begin{multline*}
\int_0^1 \frac{\ln x}{\sqrt[3] {x}}\, \d x={\smsize \left|
\begin{array}{c}\frac{1}{\sqrt[3] {x}}=y, \quad x=\frac{1}{y^3}\\
x\to +0 \Leftrightarrow y\to +\infty
\\
x=1 \Leftrightarrow y=1
\end{array}\right|}=\\= \int_{+\infty}^1 y\ln \frac{1}{y^3}\, \d \frac{1}{y^3}=
\int_{+\infty}^1 \l -3y \ln y\r \frac{-3}{y^4}\, \d y=\\= 9\int_{+\infty}^1
\frac{\ln y}{y^3}\, \d y= -9\int_1^{+\infty}\frac{\ln y}{y^3}\, \d y
\underset{+\infty}{\ll} \\ \underset{+\infty}{\ll} -9\int_1^{+\infty}
\frac{y}{y^3}\, \d y= -9\int_1^{+\infty}\frac{1}{y^2}\, \d y=\frac{9}{y}\
\Big|_1^{+\infty}
 \end{multline*}
Последний интеграл сходится, значит меньший (исходный) интеграл тоже сходится.
Для скорость сходимости с помощью той же оценки получается формула:
 \begin{multline*}
\int_0^T \frac{\ln x}{\sqrt[3] {x}}\, \d x={\smsize \left|
\begin{array}{c}\frac{1}{\sqrt[3] {x}}=y, \quad x=\frac{1}{y^3}\\
x\to +0 \Leftrightarrow y\to +\infty
\\
x=T \Leftrightarrow y=\frac{1}{T}
\end{array}\right|}=\\= -9\int_{\frac{1}{T}}^{+\infty}\frac{\ln y}{y^3}\, \d y
\underset{T\to +0}{\ll} -9\int_{\frac{1}{T}}^{+\infty} \frac{y}{y^3}\, \d y=\\=
\frac{9}{y}\ \Big|_\frac{1}{T}^{+\infty}=9T
 \end{multline*}
Сокращая константу 9, получаем:
$$
\int_0^T \frac{\ln x}{\sqrt[3] {x}}\, \d x \underset{T\to +0}{\ll} T
$$
Оценкой
$$
\ln y\underset{y\to +\infty}{\ll} y^\e\qquad (\e>0)
$$
можно добиться формулы
$$
\int_0^T \frac{\ln x}{\sqrt[3] {x}}\, \d x \underset{T\to +0}{\ll}
T^{\alpha}\qquad (\alpha<2).
$$
\end{ex}

\begin{ex}\label{ex-17.7.9}
 \begin{multline*}\int_0^1 e^{\frac{1}{x}}\, \d x={\smsize \left|
\begin{array}{c}\frac{1}{x}=y, \quad x=\frac{1}{y}\\
x\to +0 \Leftrightarrow y\to +\infty \\
x=1 \Leftrightarrow y=1
\end{array}\right|}=\\= \int_{+\infty}^1 e^y \, \d \frac{1}{y}= \int_{+\infty}^1
\frac{-e^y}{y^2}\, \d y= \int_1^{+\infty}\frac{e^y}{y^2}\, \d y
\underset{+\infty}{\gg}\\
\underset{+\infty}{\gg}\int_1^{+\infty}\frac{y^2}{y^2}\, \d y= \int_1^{+\infty}
1 \, \d y=y\ \Big|_1^{+\infty}
 \end{multline*}
Последний интеграл расходится, значит больший (исходный) интеграл тоже
расходится. Для скорости расходимости получаем с помощью той же оценки формулу:
 \begin{multline*}
\int_T^1 e^{\frac{1}{x}}\, \d x={\smsize \left|
\begin{array}{c}\frac{1}{x}=y, \quad x=\frac{1}{y}\\
x=T \Leftrightarrow y=\frac{1}{T}\\
x=1 \Leftrightarrow y=1
\end{array}\right|}= \int_{\frac{1}{T}}^1 e^y \, \d \frac{1}{y}
\underset{T\to+0}{\gg}\\
\underset{T\to+0}{\gg}\int_1^{\frac{1}{T}}\frac{y^2}{y^2}\, \d y=y\
\Big|_1^{\frac{1}{T}}=\frac{1}{T}-1\underset{T\to+0}{\sim}\frac{1}{T}
 \end{multline*}
Применяя формулу
$$
e^y\underset{y\to +\infty}{\gg} y^\e\qquad (\e>0)
$$
можно добиться оценки
$$
\int_T^1 e^{\frac{1}{x}}\, \d x\underset{T\to+0}{\gg} \frac{1}{T^n}\qquad
(n\in\N)
$$

\end{ex}

\begin{ers} Исследуйте на сходимость несобственные интегралы и найдите оценки сходимости и расходимости:
 \biter{
\item[1)] $\int_{-1}^0 e^{\frac{1}{x}}\, \d x$;

\item[2)] $\int_0^1 \frac{e^{\frac{1}{x}}}{x^3}\, \d x$;

\item[3)] $\int_{-1}^0 \frac{e^{\frac{1}{x}}}{x^3}\, \d x$;

\item[4)] $\int_0^\frac{1}{2}\frac{\, \d x}{x^3 \ln x}$;

\item[5)] $\int_e^{+\infty}\frac{\ln(1+x+x^2)}{1+x+x^2}\, \d x$;

\item[6)] $\int_e^{+\infty}\frac{\ln(1+\frac{1}{x}+\frac{1}{\sqrt{x}})}
{\ln(1+x+\sqrt{x})}\, \d x$;

\item[7)] $\int_e^{+\infty}\frac{\ln(1+\frac{1}{\ln x}+\frac{1}{\sqrt{x}})}
{x+x^\frac{5}{4}+x^\frac{4}{3}}\, \d x$;
 }\eiter
\end{ers}

\end{multicols}\noindent\rule[10pt]{160mm}{0.1pt}

\paragraph{Интегрирование асимптотических формул}

Из теоремы \ref{tm-17.7.5} сразу следует

\btm\label{TH:integrirovanie-asymp-formul} Пусть функции $f$, $g$ и $h$
обладают следующими свойствами:
 \bit{
\item[a)] $f$, $g$ и $h$ непрерывны на полуинтервале $[a;c)$;

\item[b)] $f$, $g$ и $h$ удовлетворяют асимптотическому соотношению
$$
f(x)=g(x)+\underset{x\to c-0}{\bold{o}}\l h(x) \r
$$

\item[c)] $h(x)>0$, $\forall x\in (a;c)$.
 }\eit
Тогда
 \bit{
\item[(i)] если несобственные интегралы $\int_a^c f(x) \, \d x$, $\int_a^c g(x)
\, \d x$ и $\int_a^c h(x) \, \d x$ сходятся, то справедлива асимптотическая
формула
 \beq\label{int-asymp-ostatka}
\int_x^c f(t) \, \d t=\int_x^c g(t) \, \d t+\underset{x\to c-0}{\bold{o}}\l
\int_x^c f(t) \, \d t\r
 \eeq

\item[(ii)] если несобственный интеграл $\int_a^c h(x) \, \d x$ расходится, то
справедлива асимптотическая формула
 \beq\label{int-asymp-nachala}
\int_a^x f(t) \, \d t=\int_a^x g(t) \, \d t+\underset{x\to c-0}{\bold{o}}\l
\int_a^x f(t) \, \d t\r
 \eeq
 }\eit
\etm

\noindent\rule{160mm}{0.1pt}\begin{multicols}{2}

Из теоремы \ref{TH:integrirovanie-asymp-formul} тут же следует метод нахождения
асимптотики интегралов
$$
F(x)=\int_a^x f(t)\ \d t,\qquad F(x)=\int_x^c f(t)\ \d t
$$
-- для этого нужно найти асимптотику для подынтегральной функции, а затем
проинтегрировать ее. Покажем, как это делается на примерах.

\bex Найдем асимптотику интеграла
$$
F(x)=\int_1^x\sqrt{t^2+1}\ \d t,\qquad x\to+\infty
$$
Оговоримся сразу, что этот интеграл можно вычислить явно (мы даже предлагали
это читателю в упражнении \ref{ERS:int-po-chastyam}), и как следствие, его
асимптотику можно вывести по аналогии с примерами
\ref{SEC:formuly-Peano-i-asimptotika}\ref{SUBSEC:asimptotika} (вычислив $F(x)$,
и, после преобразований, разложив полученное выражение с помощью формул Пеано).
Однако здесь наша цель -- объяснить, как находится асимптотика подобных
интегралов без их явного вычисления. Мы выбрали этот пример только чтобы не
усложнять вычисления. Читателю же мы предлагаем самостоятельно найти
асимптотику так, как мы описали, а затем сравнить ответы, полученные разными
методами.

Вспомним, что в примере \ref{EX:asymp-sqrt(x^2+1)} мы уже искали асимптотику
подынтегрального выражения, и (с точностью до замены переменной) мы получили
такой ответ:
 \beq\label{symptotika-sqrt(t^2+1)-n=3}
\sqrt{t^2+1}=t+\frac{1}{2t}-\frac{1}{8t^3}+\underset{t\to+\infty}{\bold{o}}\l\frac{1}{t^3}\r
 \eeq
К этой формуле, однако, невозможно применить часть (ii) теоремы
\ref{TH:integrirovanie-asymp-formul}, как нам бы хотелось, потому что
получающийся интеграл под символом $\bold{o}$ сходится. Но если огрубить эту
формулу до формулы
$$
\sqrt{t^2+1}=t+\frac{1}{2t}+\underset{t\to+\infty}{\bold{o}}\l\frac{1}{t}\r
$$
то интеграл под $\bold{o}$ становится расходящимся, и теорема
\ref{TH:integrirovanie-asymp-formul} (ii) становится применима:
 \begin{multline*}
\int_1^x\sqrt{t^2+1}\ \d t=\\=\int_1^x \l t+\frac{1}{2t} \r\ \d t+
\underset{t\to+\infty}{\bold{o}}\l \int_1^x \frac{1}{t}\ \d t\r=\\= \l
\frac{t^2}{2}+\frac{1}{2}\cdot\ln t \r\Big|_{t=1}^{t=x}+
\underset{t\to+\infty}{\bold{o}}\l \ln t\Big|_{t=1}^{t=x} \r=\\=
\frac{x^2}{2}+\frac{1}{2}\cdot\ln x \underbrace{-\frac{1}{2}+
\underset{t\to+\infty}{\bold{o}}(\ln x)}_{\underset{t\to+\infty}{\bold{o}}(\ln
x)}=\\=\frac{x^2}{2}+\frac{1}{2}\cdot\ln x+\underset{t\to+\infty}{\bold{o}}(\ln
x)
 \end{multline*}
Мы получили асимптотическую формулу:
$$
\int_1^x\sqrt{t^2+1}\ \d t=\frac{x^2}{2}+\frac{1}{2}\cdot\ln
x+\underset{t\to+\infty}{\bold{o}}(\ln x)
$$

Можно было бы поступить иначе: вместо того, чтобы огрублять формулу
\eqref{symptotika-sqrt(t^2+1)-n=3} перед тем как ее интегрировать, можно
переписать ее так, чтобы была применима часть (i) теоремы
\ref{TH:integrirovanie-asymp-formul}:
 \beq\label{symptotika-sqrt(t^2+1)-n=3-new}
\sqrt{t^2+1}-t-\frac{1}{2t}=-\frac{1}{8t^3}+\underset{t\to+\infty}{\bold{o}}\l\frac{1}{t^3}\r
 \eeq
Здесь правая часть бесконечно мала по сравнению с функцией $\frac{1}{t^2}$,
интеграл от которой на $[1,+\infty)$ сходится. Поэтому по теореме
\ref{tm-17.7.5}, интеграл от левой части тоже должен сходиться:
$$
\int_1^\infty \l\sqrt{t^2+1}-t-\frac{1}{2t}\r\ \d t=A\in\R
$$
Теперь, применяя к \eqref{symptotika-sqrt(t^2+1)-n=3-new} теорему
\ref{TH:integrirovanie-asymp-formul} (i), получаем:
 \begin{multline*}
\int_x^\infty \l\sqrt{t^2+1}-t-\frac{1}{2t}\r\  \d t=\\=-\int_x^\infty
\frac{1}{8t^3}\ \d t+\underset{t\to+\infty}{\bold{o}}\l\int_x^\infty
\frac{1}{t^3}\ \d t\r=\\= \frac{1}{16 t^2}\Big|_{t=x}^{t=\infty}+
\underset{t\to+\infty}{\bold{o}}\l-\frac{1}{2 t^2}\Big|_{t=x}^{t=\infty}\r=\\=
-\frac{1}{16 x^2}+ \underset{t\to+\infty}{\bold{o}}\l \frac{1}{2 x^2}\r=\\=
-\frac{1}{16 x^2}+ \underset{t\to+\infty}{\bold{o}}\l \frac{1}{x^2}\r
 \end{multline*}
Отсюда следует:
 \begin{multline*}
\int_1^x \l\sqrt{t^2+1}-t-\frac{1}{2t}\r\  \d t=\\= \underbrace{\int_1^\infty
\l\sqrt{t^2+1}-t-\frac{1}{2t}\r\  \d t}_{A}-\\-\underbrace{\int_x^\infty
\l\sqrt{t^2+1}-t-\frac{1}{2t}\r\  \d t}_{-\frac{1}{16 x^2}+
\underset{t\to+\infty}{\bold{o}}\l \frac{1}{x^2}\r}=\\= A+\frac{1}{16 x^2}+
\underset{t\to+\infty}{\bold{o}}\l \frac{1}{x^2}\r
 \end{multline*}
$$
\Downarrow
$$
 \begin{multline*}
\int_1^x \sqrt{t^2+1}\  \d t=\underbrace{\int_1^x
\l\sqrt{t^2+1}-t-\frac{1}{2t}\r\ \d t}_{\frac{x^2}{2}+\frac{1}{2}\cdot\ln x
-\frac{1}{2}}+\\+ A+\frac{1}{16 x^2}+ \underset{t\to+\infty}{\bold{o}}\l
\frac{1}{x^2}\r=\\= \frac{x^2}{2}+\frac{1}{2}\cdot\ln x -\frac{1}{2}+
A+\frac{1}{16 x^2}+ \underset{t\to+\infty}{\bold{o}}\l \frac{1}{x^2}\r
 \end{multline*}
Мы получили асимптотическое разложение нашего интеграла порядка 2 вдоль
степенной последовательности на бесконечности:
 \begin{multline*}
\int_1^x \sqrt{t^2+1}\  \d t=\\= \frac{x^2}{2}+\frac{1}{2}\cdot\ln x
-\frac{1}{2}+ A+\frac{1}{16 x^2}+ \underset{t\to+\infty}{\bold{o}}\l
\frac{1}{x^2}\r,
 \end{multline*}
где
$$
A=\int_1^\infty \l\sqrt{t^2+1}-t-\frac{1}{2t}\r\ \d t.
$$

Заметим, наконец, что, выписывая более точную асимптотику для подынтегральной
функции (то есть заменяя формулу \eqref{symptotika-sqrt(t^2+1)-n=3} на более
точную), можно по этому алгоритму получать и асимптотику интеграла нужной
точности.
 \eex

\end{multicols}\noindent\rule[10pt]{160mm}{0.1pt}

\paragraph{Нахождение асимптотики интегрированием по частям}

Описанный в предыдущем пункте способ нахождения асимптотики интеграла не всегда
работает. Например, асимптотику интегралов
$$
\int_1^x \frac{\sin t}{t^2}\ \d t,\qquad  \int_1^x t^2\cdot e^t\ \d t\qquad
(x\to+\infty)
$$
невозможно получить, раскладывая подынтегральную функцию в асимптотическую
формулу, а затем интегрируя ее (потому что асимптотику для подынтегральной
функции здесь трудно придумать). В таких случаях полезно держать в голове
другой метод -- интегрирование по частям. Здесь мы рассмотрим несколько
примеров на эту тему.

\noindent\rule{160mm}{0.1pt}\begin{multicols}{2}

\bex\label{EX:asymp-int-sin(t)/t^2} Найдем асимптотику интеграла
 \beq\label{asymp-int-sin(t)/t^2}
F(x)=\int_1^x \frac{\sin t}{t^2}\ \d t
 \eeq
Прежде всего заметим, что существует конечный предел
 $$
A=\lim_{x\to\infty}F(x)=\int_1^\infty \frac{\sin t}{t^2}\ \d t
 $$
-- поскольку интеграл $\int_1^\infty \frac{\sin t}{t^2}\ \d t$ сходится
абсолютно:
 $$
\int_1^\infty \left|\frac{\sin t}{t^2}\right|\ \d t\le \int_1^\infty
\frac{1}{t^2}\ \d t=-\frac{1}{t}\Big|_1^\infty=1
 $$
Запомним это число $A$. Тогда, интегрируя по частям один раз, получим:
 \begin{multline*}
F(x)=\int_1^x \frac{\sin t}{t^2}\ \d
t=\kern-13pt\overset{\scriptsize\begin{matrix}\int\limits_1^\infty \frac{\sin
t}{t^2}\ \d t
\\
\text{\rotatebox{90}{$=$}}\end{matrix}}{A}\kern-13pt-\int_x^\infty \frac{\sin
t}{t^2}\ \d t=\\=A+\int_x^\infty \frac{\d\cos t}{t^2}=A+\frac{\cos
t}{t^2}\Big|_x^\infty-2\int_x^\infty \frac{\cos t}{t^3}\ \d t=\\=
A-\kern-5pt\underbrace{\frac{\cos
x}{x^2}}_{\underset{x\to+\infty}{\bold{o}}\l\frac{1}{x}\r}\kern-5pt
+2\underbrace{\int_x^\infty \kern-5pt\underbrace{\frac{\cos
t}{t^3}}_{\underset{t\to+\infty}{\bold{o}}\l\frac{1}{t^2}\r}\kern-5pt\ \d
t}_{\underset{x\to+\infty}{\bold{o}}\l\frac{1}{x}\r}=\\=A+\underset{x\to+\infty}{\bold{o}}\l\frac{1}{x}\r
 \end{multline*}
Это первое асимптотическое разложение. Если проинтегрировать по частям два
раза, получится более точная формула:
 \begin{multline*}
F(x)=\int_1^x \frac{\sin t}{t^2}\ \d t=...=\\= A-\frac{\cos
x}{x^2}-2\int_x^\infty \frac{\cos t}{t^3}\ \d t=\\= A-\frac{\cos
x}{x^2}-2\int_x^\infty \frac{\d\sin t}{t^3}=\\= A-\frac{\cos
x}{x^2}-\frac{2\sin t}{t^3}\Big|_x^\infty-6\int_x^\infty \frac{\sin t}{t^4}\ \d
t=\\= A-\frac{\cos x}{x^2}+\kern-5pt\underbrace{\frac{2\sin
x}{x^3}}_{\underset{x\to+\infty}{\bold{o}}\l\frac{1}{x^2}\r}\kern-5pt
-6\underbrace{\int_x^\infty \kern-5pt\underbrace{\frac{\sin
t}{t^4}}_{\underset{t\to+\infty}{\bold{o}}\l\frac{1}{t^3}\r}\kern-5pt\ \d
t}_{\underset{x\to+\infty}{\bold{o}}\l\frac{1}{x^2}\r}=\\= A-\frac{\cos
x}{x^2}+\underset{x\to+\infty}{\bold{o}}\l\frac{1}{x^2}\r
 \end{multline*}
И так далее. Увеличивая число интегрирований по частям, мы будем повышать
точность асимптотического разложения. \eex

\bex Найдем асимптотику интеграла
$$
F(x)=\int_1^x t^\alpha\cdot e^t\ \d t,\quad x\to+\infty
$$
Интегрируя по частям, получаем:
 \begin{multline}\label{int_1^x-t^alpha-cdot-e^t-d-t}
F(x)=\int_1^x t^\alpha\cdot e^t\ \d t=\int_1^x t^\alpha\ \d e^t=\\=
t^\alpha\cdot e^t\Big|_{t=1}^{t=x}- \int_1^x e^t\ \d t^\alpha=\\= x^\alpha\cdot
e^x-e- \alpha\cdot \int_1^x t^{\alpha-1}\cdot e^t\ \d t
 \end{multline}
Поглядим внимательно на последний интеграл:
$$
t^{\alpha-1}\cdot e^t\underset{t\to+\infty}{\ll}t^\alpha\cdot e^t
$$
$$
\phantom{\scriptsize\eqref{asymp-sravn-int}}\quad\Downarrow\quad{\scriptsize\eqref{asymp-sravn-int}}
$$
$$
\int_1^x t^{\alpha-1}\cdot e^t\ \d t\underset{x\to+\infty}{\ll} \int_1^x
t^\alpha\cdot e^t\ \d t=F(x)
$$
Мы получаем:
$$
F(x)=x^\alpha\cdot e^x\underbrace{-e+\underset{x\to+\infty}{\bold{o}}\l
F(x)\r}_{\underset{x\to+\infty}{\bold{o}}\l F(x)\r}
$$
$$
\Downarrow
$$
$$
F(x)=x^\alpha\cdot e^x+\underset{x\to+\infty}{\bold{o}}\l F(x)\r
$$
$$
\Downarrow
$$
$$
F(x)=\int_1^x t^\alpha\cdot e^t\ \d t\underset{x\to+\infty}{\sim}x^\alpha\cdot
e^x
$$
$$
\Downarrow
$$
 \beq\label{int_1^x-t^alpha-cdot-e^t-d-t-0}
\int_1^x t^\alpha\cdot e^t\ \d t=x^\alpha\cdot
e^x+\underset{x\to+\infty}{\bold{o}}\l x^\alpha\cdot e^x\r
 \eeq
Это первое асимптотическое разложение. Из него можно выводить асимптотики более
высокого порядка. Заменив в \eqref{int_1^x-t^alpha-cdot-e^t-d-t-0} $\alpha$ на
$\alpha-1$, мы получим цепочку:
 \beq\label{int_1^x-t^alpha-cdot-e^t-d-t-00}
\int_1^x t^{\alpha-1}\cdot e^t\ \d t=x^{\alpha-1}\cdot
e^x+\underset{x\to+\infty}{\bold{o}}\l x^{\alpha-1}\cdot e^x\r
 \eeq
$$
\Downarrow
$$
 \begin{multline*}
\int_1^x t^\alpha\cdot e^t\ \d t=\eqref{int_1^x-t^alpha-cdot-e^t-d-t}=\\=
x^\alpha\cdot e^x-e- \alpha\cdot \int_1^x t^{\alpha-1}\cdot e^t\ \d
t=\eqref{int_1^x-t^alpha-cdot-e^t-d-t-00}=\\=x^\alpha\cdot e^x-e- \alpha\cdot\l
x^{\alpha-1}\cdot e^x+\underset{x\to+\infty}{\bold{o}}\l x^{\alpha-1}\cdot
e^x\r\r=\\=x^\alpha\cdot e^x- \alpha\cdot x^{\alpha-1}\cdot
e^x\underbrace{-e+\underset{x\to+\infty}{\bold{o}}\l x^{\alpha-1}\cdot
e^x\r}_{\underset{x\to+\infty}{\bold{o}}\l x^{\alpha-1}\cdot e^x\r}
=\\=x^\alpha\cdot e^x- \alpha\cdot x^{\alpha-1}\cdot
e^x+\underset{x\to+\infty}{\bold{o}}\l x^{\alpha-1}\cdot e^x\r
 \end{multline*}
Получается разложение:
 \beq\label{int_1^x-t^alpha-cdot-e^t-d-t-1}
\int_1^x t^\alpha\cdot e^t\ \d t=x^\alpha\cdot e^x- \alpha\cdot
x^{\alpha-1}\cdot e^x+\underset{x\to+\infty}{\bold{o}}\l x^{\alpha-1}\cdot
e^x\r
 \eeq
Здесь также можно заменить  $\alpha$ на $\alpha-1$ и после этого подставить
полученную формулу в \eqref{int_1^x-t^alpha-cdot-e^t-d-t} -- тогда появится
асимптотика более высокого порядка, и так далее. Понятно, что асимптотическая
последовательность здесь берется такая:
 $$
x^\alpha\cdot e^x\underset{x\to+\infty}{\gg} x^{\alpha-1}\cdot
e^x\underset{x\to+\infty}{\gg}x^{\alpha-2}\cdot
e^x\underset{x\to+\infty}{\gg}...
 $$
 \eex

\end{multicols}\noindent\rule[10pt]{160mm}{0.1pt}

\subsection{Асимптотика сумм}\label{SEC:asymp-summ}

\paragraph{Асимптотическая эквивалентность рядов.}

\begin{tm}[\bf признак эквивалентности рядов]\label{tm-18.3.13}
Пусть
 \beq\label{ekviv-ryadov:a_n-sim-b_n}
a_n\ge 0,\qquad b_n\ge 0,\qquad a_n\underset{n\to \infty}{\sim} b_n
 \eeq
 \bit{\rm
\item[$\bullet$] Коротко эти три условия записываются формулой
$$
\sum_{n=1}^\infty a_n\sim \sum_{n=1}^\infty b_n
$$
 }\eit
Тогда сходимость ряда $\sum\limits_{n=1}^\infty a_n$ эквивалентна сходимости
ряда $\sum\limits_{n=1}^\infty b_n$:
$$
\sum\limits_{n=1}^\infty a_n \quad \text{сходится}\quad \Longleftrightarrow
\quad \sum\limits_{n=1}^\infty b_n \quad \text{сходится}  \quad
$$
При этом,
 \bit{
\item[(i)] если эти ряды сходятся, то справедливо следующее соотношение,
показывающее, что скорость их сходимости будет одинаковой:
 \beq\label{asymp-ekviv-ryadov-0}
\sum_{n=N}^\infty a_n\underset{N\to\infty}{\sim}\sum_{n=N}^\infty b_n
 \eeq

\item[(ii)] если эти ряды расходятся, то справедливо следующее соотношение,
показывающее, что скорость их расходимости будет одинаковой:
 \beq\label{asymp-ekviv-ryadov}
\sum_{n=1}^N a_n\underset{N\to\infty}{\sim}\sum_{n=1}^N b_n
 \eeq
}\eit
\end{tm}

Для доказательства нам понадобится следующая

\begin{lm}[\bf лемма об остатке]\label{lm-18.3.13}
Пусть $\sum_{n=1}^\infty a_n$ -- числовой ряд, и $M\in \mathbb{N}$ --
произвольное число. Тогда
$$
\text{ряд $\sum_{n=1}^\infty a_n$ сходится}\quad \Longleftrightarrow \quad
\text{сходится ряд $\sum_{n=M}^\infty a_n$}\l \text{называемый остатком ряда
$\sum_{n=1}^\infty a_n$}\r
$$
\end{lm}\begin{proof}
Обозначим через $S_N$ и $R_N$ частичные суммы самого ряда и его остатка:
$$
S_N=\sum_{n=1}^N a_n, \qquad R_N=\sum_{n=M}^N a_n
$$
Очевидно, при $N\ge M$
$$
S_N=\sum_{n=1}^N a_n=\sum_{n=1}^{M-1} a_n+\sum_{n=M}^N a_n=C+R_N,
$$
где $C=\sum_{n=1}^{M-1} a_n$ -- константа, не зависящая от $N$. Теперь получаем
следующую логическую цепочку:
$$
\text{ряд $\sum_{n=1}^\infty a_n$ сходится}
$$
$$
\Updownarrow
$$
$$
\text{существует конечный предел}\,\, \lim_{N\to \infty} S_N
$$
$$
\Updownarrow
$$
$$
\text{существует конечный предел}\,\, \lim_{N\to \infty} R_N
$$
$$
\Updownarrow
$$
$$
\text{остаток $\sum_{n=M}^\infty a_n$ сходится}\quad
$$
 \end{proof}

\begin{proof}[Доказательство теоремы \ref{tm-18.3.13}] Зафиксируем сначала какое-нибудь
$\e\in(0,1)$ и заметим такую логическую цепочку:
$$
a_n\underset{n\to \infty}{\sim} b_n
$$
$$
\Downarrow
$$
$$
\frac{a_n}{b_n}\underset{n\to \infty}{\longrightarrow} 1
$$
$$
\Downarrow
$$
$$
\exists L \quad \forall n\ge L \quad \frac{a_n}{b_n}\in \left[ 1-\e;1+\e\right]
$$
$$
\Downarrow
$$
$$
\exists L \quad \forall n\ge L \quad 1-\e\le \frac{a_n}{b_n}\le 1+\e
$$
$$
\Downarrow
$$
$$
\exists L \quad \forall n\ge L \quad (1-\e)\cdot b_n\le a_n\le (1+\e)\cdot b_n
$$
$$
\Downarrow
$$
 \beq
\sum_{n=L}^\infty (1-\e)\cdot b_n\le \sum_{n=L}^\infty a_n\le \sum_{n=L}^\infty
(1+\e)\cdot b_n \label{18.3.2}
 \eeq
(последнее двойное неравенство понимается как почленное сравнение рядов,
определенное формулой \eqref{DEF:sravnenie-ryadov}).

1. Покажем теперь, что сходимость ряда $\sum\limits_{n=1}^\infty a_n$
эквивалентна сходимости ряда $\sum\limits_{n=1}^\infty b_n$. С одной стороны,
справедлива такая цепочка:
$$
\sum\limits_{n=1}^\infty a_n \quad \text{сходится}
$$
$$
\Downarrow\put(20,0){\smsize (\text{применяем лемму об остатке
\ref{lm-18.3.13}})}
$$
$$
\sum\limits_{n=L}^\infty a_n \quad \text{сходится}
$$
$$
\Downarrow\put(20,0){\smsize \text{$\begin{pmatrix}\text{применяем неравенство
\eqref{18.3.2}}\\ \text{и признак сравнения \ref{tm-18.3.7}}\end{pmatrix}$}}
$$
$$
\sum\limits_{n=L}^\infty (1-\e)\cdot b_n \quad \text{сходится}
$$
$$
\Downarrow\put(20,0){\smsize (\text{применяем свойство $1^0, \, \S 2$})}
$$
$$
\sum\limits_{n=L}^\infty b_n \quad \text{сходится}
$$
$$
\Downarrow\put(20,0){\smsize (\text{применяем лемму об остатке
\ref{lm-18.3.13}})}
$$
$$
\sum\limits_{n=1}^\infty b_n \quad \text{сходится}
$$
А, с другой стороны, справедлива цепочка в обратном направлении:
$$
\sum\limits_{n=1}^\infty b_n \quad \text{сходится}
$$
$$
\Downarrow\put(20,0){\smsize (\text{применяем лемму об остатке
\ref{lm-18.3.13}})}
$$
$$
\sum\limits_{n=L}^\infty b_n \quad \text{сходится}
$$
$$
\Downarrow\put(20,0){\smsize (\text{применяем свойство $1^0, \, \S 2$})}
$$
$$
\sum\limits_{n=L}^\infty (1+\e)\cdot b_n \quad \text{сходится}
$$
$$
\Downarrow\put(20,0){\smsize \text{$\begin{pmatrix}\text{применяем неравенство
\eqref{18.3.2}}\\ \text{и признак сравнения \ref{tm-18.3.7}}\end{pmatrix}$}}
$$
$$
\sum\limits_{n=L}^\infty a_n \quad \text{сходится}
$$
$$
\Downarrow\put(20,0){\smsize (\text{применяем лемму об остатке
\ref{lm-18.3.13}})}
$$
$$
\sum\limits_{n=1}^\infty a_n \quad \text{сходится}
$$

2. Теперь предположим, что оба ряда сходятся и докажем формулу
\eqref{asymp-ekviv-ryadov-0}. Из \eqref{18.3.2} получаем:
$$
\exists L\in\N\quad \forall N\ge L\quad (1-\e)\cdot\sum_{n=N}^\infty  b_n\le
\sum_{n=N}^\infty a_n\le (1+\e)\cdot\sum_{n=N}^\infty  b_n\qquad
\l\text{\scriptsize здесь уже $\sum_{n=N}^\infty a_n$ и $\sum_{n=N}^\infty b_n$
-- числа}\r
$$
$$
\Downarrow
$$
$$
\exists L\in\N\quad \forall N\ge L\qquad 1-\e\le\frac{\sum\limits_{n=N}^\infty
a_n}{\sum\limits_{n=N}^\infty b_n}\le1+\e\qquad \Big(\text{\scriptsize здесь
$\e\in(0;1)$ с самого начала выбиралось произвольным}\Big)
$$
$$
\Downarrow
$$
$$
\forall \e\in(0,1)\quad \exists L\in\N\quad \forall N\ge L\qquad
1-\e\le\frac{\sum\limits_{n=N}^\infty a_n}{\sum\limits_{n=N}^\infty b_n}\le1+\e
$$
$$
\Downarrow
$$
$$
\frac{\sum\limits_{n=N}^\infty a_n}{\sum\limits_{n=N}^\infty
b_n}\underset{N\to\infty}{\longrightarrow} 1
$$
то есть справедливо \eqref{asymp-ekviv-ryadov-0}.

3. Нам остается доказать \eqref{asymp-ekviv-ryadov} в предположении, что оба
ряда расходятся. Здесь применяется теорема Штольца \ref{Stoltz}. Из условия
\eqref{ekviv-ryadov:a_n-sim-b_n} следует, что начиная с некоторого номера $M$
все числа $b_n$ ненулевые, и поэтому положительны:
 \beq\label{ekviv-ryadov:b_n>0}
\forall n\ge M\qquad b_n>0
 \eeq
Зафиксируем это $M$ и обозначим
$$
A_N=\sum_{n=M}^N a_n,\qquad B_N=\sum_{n=M}^N b_n.
$$
Из \eqref{ekviv-ryadov:b_n>0} следует, что последовательность $B_N$ ($N\ge M$)
строго возрастает:
$$
B_M<B_{M+1}<...<B_N<...
$$
С другой стороны, она стремится к бесконечности, поскольку мы предполагаем ряд
$\sum_{n=1}^\infty b_n$ расходящимся:
$$
B_N=\sum_{n=M}^N b_n\underset{N\to\infty}{\longrightarrow}\infty
$$
Поэтому по теореме \ref{Stoltz},
$$
\lim_{N\to\infty}\frac{A_N}{B_N}=\lim_{N\to\infty}\frac{A_N-A_{N-1}}{B_N-B_{N-1}}=\lim_{N\to\infty}\frac{a_N}{b_N}=1
$$
Если теперь обозначить
$$
A=\sum_{n=1}^{M-1} a_n,\qquad B=\sum_{n=1}^{M-1} b_n,
$$
то мы получим:
$$
\lim_{N\to\infty}\frac{\sum\limits_{n=1}^N a_n}{\sum\limits_{n=1}^N
b_n}=\lim_{N\to\infty}\frac{\sum\limits_{n=1}^{M-1} a_n+\sum\limits_{n=M}^N
a_n}{\sum\limits_{n=1}^{M-1} b_n+\sum\limits_{n=M}^N
b_n}=\lim_{N\to\infty}\frac{A+A_N}{B+B_N}=\lim_{N\to\infty}\frac{
\boxed{\frac{A}{B_N}}\put(-5,20){\vector(1,2){5}\put(2,14){$0$}}+
\boxed{\frac{A_N}{B_N}}\put(-5,20){\vector(1,2){5}\put(2,14){$1$}}}{
\boxed{\frac{B}{B_N}}\put(-5,-18){\vector(1,-2){5}\put(2,-18){$0$}}+1}=1
$$
То есть справедливо \eqref{asymp-ekviv-ryadov}.
\end{proof}

\noindent\rule{160mm}{0.1pt}\begin{multicols}{2}

\begin{ex}\label{ex-18.3.15}
$$
\sum_{n=1}^\infty \frac{1}{n^2+n-1}\sim \underbrace{\sum_{n=1}^\infty
\frac{1}{n^2}}_{\text{сходится}}
$$
Вывод: ряд $\sum\limits_{n=1}^\infty \frac{1}{n^2+n-1}$ сходится.

В соответствии с теоремой \ref{tm-18.3.13}, можно дать асимптотическую оценку
скорости сходимости этого ряда:
$$
\sum_{n=1}^N \frac{1}{n^2+n-1}\underset{N\to\infty}{\sim} \sum_{n=1}^N
\frac{1}{n^2},
$$
однако это будет мало информативно, потому что до теорем
\ref{TH:ASYMP-COR-int-priznak-Cauchy} и \ref{TH:form-Euler} мы не сможем
выражать асимптотику рядов через стандартные функции.
\end{ex}

\begin{ex}\label{ex-18.3.16}
$$
\sum_{n=1}^\infty \frac{1+\sqrt{n}}{\sqrt{n^3+2n}}\sim \sum_{n=1}^\infty
\frac{\sqrt{n}}{\sqrt{n^3}}=\underbrace{\sum_{n=1}^\infty
\frac{1}{n}}_{\text{расходится}}
$$
Вывод: ряд $\sum\limits_{n=1}^\infty \frac{1+\sqrt{n}}{\sqrt{n^3+2n}}$
расходится.
\end{ex}

\begin{ex}\label{ex-18.3.17}
$$
\sum_{n=1}^\infty \l \sin \frac{1}{n}\r^\alpha \sim
\underbrace{\sum_{n=1}^\infty
\frac{1}{n^\alpha}}_{\scriptsize\begin{matrix}\text{сходится}\\ \text{при
$\alpha>1$}\end{matrix}}
$$
Вывод: ряд $\sum\limits_{n=1}^\infty \l \sin \frac{1}{n}\r^\alpha$ сходится при
$\alpha>1$.
\end{ex}

\begin{ex}\label{ex-18.3.18}
$$
\sum_{n=1}^\infty  \frac{\ln \l 1+\frac{1}{n}\r}{n^\alpha}\sim
\underbrace{\sum_{n=1}^\infty
\frac{1}{n^{\alpha+1}}}_{\scriptsize\begin{matrix}\text{сходится}\\ \text{при
$\alpha+1>1$}\\ \text{то есть при $\alpha>0$}\end{matrix}}
$$
Вывод: ряд $\sum\limits_{n=1}^\infty \frac{\ln \l 1+\frac{1}{n}\r}{n^\alpha}$
сходится при $\alpha>0$.
\end{ex}

\begin{ers} Исследуйте на сходимость ряды:
 \biter{
\item[1)] $\sum\limits_{n=1}^\infty \sqrt{n}\cdot \l \frac{1+n^2}{1+n^3}\r^2$;

\item[2)] $\sum\limits_{n=1}^\infty \frac{1}{\sqrt{n(n+1)(n+2)}}$;

\item[3)] $\sum\limits_{n=1}^\infty n\cdot \l \sqrt{n+1}-\sqrt{n-1}\r$;

\item[4)] $\sum\limits_{n=1}^\infty \frac{1}{\sqrt[3]{n}}\cdot \ln
\frac{n+1}{n-1}$;

\item[5)] $\sum\limits_{n=1}^\infty \frac{\ln \l n+\sqrt{n}+\ln n+1 \r} {3n-\ln
n +\arctg{n}-1}$;

\item[6)] $\sum\limits_{n=1}^\infty \frac{\ln \l 1+\frac{1}{n}+\frac{1}{n^2}\r}
{\ln \l n-n\sqrt{n}+\ln n+3 \r}$;

\item[7)] $\sum\limits_{n=1}^\infty n\cdot \l e^{\frac{1}{n}}-1 \r^\alpha$;

\item[8)] $\sum\limits_{n=1}^\infty \frac{1+n} {\l \frac{1}{n}-\sin
\frac{1}{n}\r^\alpha}$;

\item[9)] $\sum\limits_{n=1}^\infty \lll \l 1+ \frac{1}{n}\r^{\frac{3}{2}}-1
\rrr^\alpha$;

\item[10)] $\sum\limits_{n=1}^\infty \frac{\l \arctg
\frac{1}{n}\r^\alpha}{1-\cos \frac{1}{n}}$.
 }\eiter \end{ers}
\end{multicols}\noindent\rule[10pt]{160mm}{0.1pt}

\paragraph{Асимптотическое сравнение рядов.}

\begin{tm}[\bf асимптотический признак сравнения рядов]\label{tm-18.5.9}
Пусть последовательности $\{ a_n \}$ и $\{ b_n \}$ обладают следующими
свойствами:
 \bit{
\item[(1)] $a_n\underset{n\to \infty}{\ll} b_n$  (то есть $a_n=\underset{n\to
\infty}{\bold{o}}\l b_n \r$);

\item[(2)] $b_n\ge 0, \quad \forall n\in \mathbb{N}$.
 }\eit
 \bit{\rm
\item[$\bullet$] Коротко эти два условия записываются формулой:
$$
\sum_{n=1}^\infty a_n \ll \sum_{n=1}^\infty b_n
$$
 }\eit
Тогда
 \bit{
\item[(i)] из сходимости большего ряда $\sum\limits_{n=1}^\infty b_n$ следует
сходимость меньшего ряда $\sum\limits_{n=1}^\infty a_n$:
 \beq\label{sravn-ryadov:sum-a_n-shod<=sum-b_n-shod}
\sum\limits_{n=1}^\infty a_n \quad \text{сходится}\quad \Longleftarrow \quad
\sum\limits_{n=1}^\infty b_n \quad \text{сходится}
 \eeq
причем скорости сходимости оцениваются соотношением
 \beq\label{asymp-sravn-ryadov-0}
\sum_{n=N}^\infty a_n\underset{N\to\infty}{\ll}\sum_{n=N}^\infty b_n;
 \eeq

\item[(ii)] из расходимости меньшего ряда $\sum\limits_{n=1}^\infty a_n$
следует расходимость большего ряда $\sum\limits_{n=1}^\infty b_n$:
 \beq\label{sravn-ryadov:sum-a_n-rashod=>sum-b_n-rashod}
\sum\limits_{n=1}^\infty a_n \quad \text{расходится}\quad \Longrightarrow \quad
\sum\limits_{n=1}^\infty b_n \quad \text{расходится}
 \eeq
причем скорости расходимости оцениваются соотношением
 \beq\label{asymp-sravn-ryadov-2}
\sum_{n=1}^N a_n\underset{N\to\infty}{\ll}\sum_{n=1}^N b_n;
 \eeq
 }\eit
\end{tm}\begin{proof}
Зафиксируем сначала какое-нибудь $\e>0$ и заметим такую логическую цепочку:
$$
a_n\underset{n\to \infty}{\ll} b_n
$$
$$
\Downarrow
$$
$$
\frac{a_n}{b_n}\underset{n\to \infty}{\longrightarrow} 0
$$
$$
\Downarrow
$$
$$
\frac{|a_n|}{b_n}\underset{n\to \infty}{\longrightarrow} 0
$$
$$
\Downarrow
$$
$$
\exists L \quad \forall n\ge L \quad \frac{|a_n|}{b_n}<\e
$$
$$
\Downarrow
$$
$$
\exists L \quad \forall n\ge L \quad |a_n|<\e\cdot b_n
$$
$$
\Downarrow
$$
 \beq\label{18.3.2-0}
\sum_{n=L}^\infty |a_n|< \sum_{n=L}^\infty \e\cdot b_n
 \eeq
(последнее неравенство понимается как почленное сравнение рядов).

1. Докажем \eqref{sravn-ryadov:sum-a_n-shod<=sum-b_n-shod} и
\eqref{sravn-ryadov:sum-a_n-rashod=>sum-b_n-rashod}. Замечаем такую цепочку:
$$
  \text{ряд}\,\, \sum_{n=1}^\infty b_n \,\,\text{сходится}
$$
$$
\Downarrow\put(20,0){\smsize (\text{применяем лемму об остатке
\ref{lm-18.3.13}})}
$$
$$
  \text{ряд}\,\, \sum_{n=L}^\infty b_n  \,\,\text{сходится}
$$
$$
\Downarrow
$$
$$
  \text{ряд}\,\, \sum_{n=L}^\infty \e\cdot b_n  \,\,\text{сходится}
$$
$$
\Downarrow\put(20,0){\smsize \text{$\begin{pmatrix}\text{вспоминаем неравенство
\eqref{18.3.2-0}} \\ \text{и признак сравнения рядов
\ref{tm-18.3.7}}\end{pmatrix}$}}
$$
$$
  \text{ряд}\,\, \sum_{n=L}^\infty |a_n|  \,\,\text{сходится}
$$
$$
\Downarrow\put(20,0){\smsize (\text{применяем лемму об остатке
\ref{lm-18.3.13}})}
$$
$$
  \text{ряд}\,\, \sum_{n=1}^\infty |a_n|  \,\,\text{сходится}
$$
$$
\Downarrow
$$
$$
  \text{ряд}\,\, \sum_{n=1}^\infty a_n  \,\,\text{сходится}
$$
Это доказывает утверждение \eqref{sravn-ryadov:sum-a_n-shod<=sum-b_n-shod}. Его
можно переформулировать так:
$$
\text{\it невозможна ситуация, когда}\, \sum\limits_{n=1}^\infty a_n \,
\text{\it расходится. а}\, \sum\limits_{n=1}^\infty b_n \, \text{\it сходится}
$$
и отсюда следует, что если $\sum\limits_{n=1}^\infty a_n$ расходится, то
$\sum\limits_{n=1}^\infty b_n$ тоже должен расходиться. То есть выполняется и
\eqref{sravn-ryadov:sum-a_n-rashod=>sum-b_n-rashod}.

2. Теперь предположим, что оба ряда сходятся и докажем формулу
\eqref{asymp-sravn-ryadov-0}. Из \eqref{18.3.2-0} получаем:
$$
\exists L\in\N\quad \forall N\ge L\quad \sum_{n=N}^\infty |a_n|<
\e\cdot\sum_{n=N}^\infty  b_n\qquad \l\text{\scriptsize здесь уже
$\sum_{n=N}^\infty a_n$ и $\sum_{n=N}^\infty b_n$ -- числа}\r
$$
$$
\Downarrow
$$
$$
\exists L\in\N\quad \forall N\ge L\qquad \frac{\sum\limits_{n=N}^\infty
|a_n|}{\sum\limits_{n=N}^\infty b_n}< \e\qquad \Big(\text{\scriptsize здесь
$\e>0$ с самого начала выбиралось произвольным}\Big)
$$
$$
\Downarrow
$$
$$
\forall \e>0\quad \exists L\in\N\quad \forall N\ge L\qquad
\frac{\sum\limits_{n=N}^\infty |a_n|}{\sum\limits_{n=N}^\infty b_n}<\e
$$
$$
\Downarrow
$$
$$
\frac{\sum\limits_{n=N}^\infty |a_n|}{\sum\limits_{n=N}^\infty
b_n}\underset{N\to\infty}{\longrightarrow} 0
$$
$$
\Downarrow
$$
$$
\frac{\sum\limits_{n=N}^\infty a_n}{\sum\limits_{n=N}^\infty
b_n}\underset{N\to\infty}{\longrightarrow} 0
$$
то есть справедливо \eqref{asymp-sravn-ryadov-0}.

3. Нам остается доказать \eqref{asymp-sravn-ryadov-2} в предположении, что оба
ряда расходятся. Здесь применяется теорема Штольца \ref{Stoltz}. Из условия
$a_n\underset{n\to \infty}{\ll} b_n$ следует, что начиная с некоторого номера
$M$ все числа $b_n$ ненулевые, и поэтому положительны:
 \beq\label{sravn-ryadov:b_n>0}
\forall n\ge M\qquad b_n>0
 \eeq
Зафиксируем это $M$ и обозначим
$$
A_N=\sum_{n=M}^N a_n,\qquad B_N=\sum_{n=M}^N b_n.
$$
Из \eqref{sravn-ryadov:b_n>0} следует, что последовательность $B_N$ ($N\ge M$)
строго возрастает:
$$
B_M<B_{M+1}<...<B_N<...
$$
С другой стороны, она стремится к бесконечности, поскольку мы предполагаем ряд
$\sum_{n=1}^\infty b_n$ расходящимся:
$$
B_N=\sum_{n=M}^N b_n\underset{N\to\infty}{\longrightarrow}\infty
$$
Поэтому по теореме \ref{Stoltz},
$$
\lim_{N\to\infty}\frac{A_N}{B_N}=\lim_{N\to\infty}\frac{A_N-A_{N-1}}{B_N-B_{N-1}}=\lim_{N\to\infty}\frac{a_N}{b_N}=0
$$
Если теперь обозначить
$$
A=\sum_{n=1}^{M-1} a_n,\qquad B=\sum_{n=1}^{M-1} b_n,
$$
то мы получим:
$$
\lim_{N\to\infty}\frac{\sum_{n=1}^N a_n}{\sum_{n=1}^N
b_n}=\lim_{N\to\infty}\frac{\sum_{n=1}^{M-1} a_n+\sum_{n=M}^N
a_n}{\sum_{n=1}^{M-1} b_n+\sum_{n=M}^N
b_n}=\lim_{N\to\infty}\frac{A+A_N}{B+B_N}=\lim_{N\to\infty}\frac{
\boxed{\frac{A}{B_N}}\put(-5,20){\vector(1,2){5}\put(2,14){$0$}}+
\boxed{\frac{A_N}{B_N}}\put(-5,20){\vector(1,2){5}\put(2,14){$0$}}}{
\boxed{\frac{B}{B_N}}\put(-5,-18){\vector(1,-2){5}\put(2,-18){$0$}}+1}=0
$$
То есть справедливо \eqref{asymp-sravn-ryadov-2}.
\end{proof}

\noindent\rule{160mm}{0.1pt}\begin{multicols}{2}

Асимптотический признак сравнения рядов удобно применять, переписав и дополнив
для дискретного аргумента шкалу бесконечностей \eqref{shkala-besk} следующим
образом:
 \beq\label{shakala-besk-s-diskr-arg}
 \boxed{
\begin{split}
 \quad
1 \underset{n\to \infty}{\ll}\quad \underset{(a>1)}{\log_a n}\quad
\underset{n\to \infty}{\ll}\quad \underset{(\alpha>0)}{n^\alpha}\quad
\underset{n\to \infty}{\ll}\quad \\
\underset{n\to \infty}{\ll}\quad \underset{(a>1)}{a^n}\underset{n\to
\infty}{\ll}\quad n! \underset{n\to \infty}{\ll}\quad n^n \quad \end{split} }
 \eeq

\begin{ex}\label{ex-18.5.10}
$$
\sum_{n=1}^\infty \frac{\ln n}{n^2}\ll \sum_{n=1}^\infty
\frac{n^{\frac{1}{2}}}{n^2}=\kern-16pt \underbrace{\sum_{n=1}^\infty
\frac{1}{n^{\frac{3}{2}}}}_{\smsize
\begin{matrix}\text{сходится}, \\
\text{в силу примера \ref{ex-18.3.3}}\end{matrix}}
$$
Больший ряд сходится, значит и меньший ряд сходится.

Вывод: ряд $\sum\limits_{n=1}^\infty \frac{\ln n}{n^2}$ сходится.

В соответствии с теоремой \ref{tm-18.5.9}, можно дать асимптотическую оценку
скорости сходимости этого ряда:
$$
\sum_{n=1}^N \frac{\ln n}{n^2}\underset{N\to\infty}{\ll} \sum_{n=1}^N
\frac{1}{n^{\frac{3}{2}}},
$$
но как и в примерах, иллюстрировавших теорему \ref{tm-18.3.13}, это мало что
объяснит, поскольку, как мы уже говорили, до теорем
\ref{TH:ASYMP-COR-int-priznak-Cauchy} и \ref{TH:form-Euler} мы не сможем
выражать асимптотику рядов через стандартные функции.

\end{ex}

\begin{ex}\label{ex-18.5.11}
$$
\sum_{n=2}^\infty \frac{1}{\ln n}\gg \kern-16pt\underbrace{\sum_{n=2}^\infty
\frac{1}{n}}_{\smsize\begin{matrix}\text{расходится,}\\
\text{в силу примера \ref{ex-18.3.2}}\end{matrix}}
$$
Меньший ряд расходится, значит и больший ряд расходится.

Вывод: ряд $\sum\limits_{n=2}^\infty \frac{1}{\ln n}$ расходится.
\end{ex}

\begin{ers}  Исследуйте на сходимость ряды:
 \begin{multicols}{2}
1) $\sum\limits_{n=2}^\infty \frac{n^5}{2^n}$;

2) $\sum\limits_{n=2}^\infty (-1)^n \frac{5^n}{n!}$;

3) $\sum\limits_{n=2}^\infty \sin n \frac{\ln n}{n\sqrt{n}}$;

4) $\sum\limits_{n=2}^\infty \frac{e^n}{n^{100}}$;

5) $\sum\limits_{n=2}^\infty \frac{5^n}{\sqrt{n!}}$.
\end{multicols}\end{ers}

\end{multicols}\noindent\rule[10pt]{160mm}{0.1pt}

\paragraph{Связь с несобственными интегралами и формула суммирования Эйлера.}

Связь между рядами и несобственными интегралами, успевшая к настоящему времени
проявить себя в интегральном признаке сходимости (теорема \ref{tm-18.3.1}),
позволяет в некоторых случаях выразить асимптотику частичных сумм ряда в
стандартных функциях. Первый результат на эту тему -- очевидное следствие
теоремы \ref{tm-18.3.1}:

\btm\label{TH:ASYMP-COR-int-priznak-Cauchy} Пусть $f$ -- неотрицательная и
невозрастающая функция на полуинтервале $[1;+\infty)$. Тогда
 \beq\label{ASYMP-COR:int-priznak-Cauchy}
\sum_{n=1}^Nf(n)=\int_1^N f(x)\ \d x+C+\underset{N\to\infty}{\bold{o}}(1)
 \eeq
где $C$ -- постоянная Эйлера ряда $\sum_{n=1}^\infty f(n)$. Если вдобавок ряд
$\sum_{n=1}^\infty f(n)$ расходится, то
 \beq\label{ASYMP-COR:int-priznak-Cauchy-infty}
\sum_{n=1}^Nf(n)\underset{N\to\infty}{\sim}\int_1^N f(x)\ \d x
 \eeq
\etm \bpr Формула \eqref{ASYMP-COR:int-priznak-Cauchy} -- просто по-другому
переписанное соотношение \eqref{int-priznak-Cauchy}:
$$
\sum_{n=1}^Nf(n)-\int_1^N f(x)\ \d x\underset{N\to\infty}{\longrightarrow} C
$$
Если ряд $\sum_{n=1}^\infty f(n)$ расходится, то есть
$$
\sum_{n=1}^N f(n)\underset{N\to\infty}{\longrightarrow}\infty,
$$
то, мы получим:
$$
\sum_{n=1}^Nf(n)=\int_1^N f(x)\ \d x+
\kern-24pt\overbrace{C}^{\scriptsize\begin{matrix}
\underset{N\to\infty}{\bold{o}}\l \sum\limits_{n=1}^N f(n)\r\\
\text{\rotatebox{90}{$=$}} \end{matrix}}\kern-24pt +
\kern-18pt\underbrace{\underset{N\to\infty}{\bold{o}}(1)}_{\scriptsize\begin{matrix}\text{\rotatebox{90}{$=$}}\\
\underset{N\to\infty}{\bold{o}}\l \sum\limits_{n=1}^N f(n)\r\end{matrix}}
\kern-18pt =\int_1^N f(x)\ \d x+\underset{N\to\infty}{\bold{o}}\l \sum_{n=1}^N
f(n)\r
$$
а это эквивалентно формуле \eqref{ASYMP-COR:int-priznak-Cauchy} по теореме
\ref{TH:svyaz-<<-i-o}. \epr

\noindent\rule{160mm}{0.1pt}\begin{multicols}{2}

\bex\label{asimptotika-garm-ryada} Частичные суммы гармонического ряда
$$
\sum_{n=1}^\infty\frac{1}{n}
$$
по теореме \ref{TH:ASYMP-COR-int-priznak-Cauchy} удовлетворяют формуле
 \beq\label{harmonic-asymp}
\sum_{n=1}^N\frac{1}{n}=\underbrace{\ln
N}_{\scriptsize\begin{matrix}\text{\rotatebox{90}{$=$}}\\
\int\limits_1^N\frac{1}{x}\ \d x
\end{matrix}}+C+\underset{N\to\infty}{\bold{o}}(1)
 \eeq
где $C$ -- постоянная Эйлера \eqref{harmonic-Euler}. Мы отмечали в примере
\ref{ex-18.3.2}, и это видно из формулы \eqref{harmonic-asymp}, что этот ряд
расходится. Формула \eqref{ASYMP-COR:int-priznak-Cauchy-infty} для него
принимает вид:
$$
\sum_{n=1}^N\frac{1}{n}\underset{N\to\infty}{\sim}\ln N.
$$
 \eex

\bex Частичные суммы ряда Дирихле
$$
\sum_{n=1}^\infty\frac{1}{n^\alpha},
$$
который мы рассматривали в примере \ref{ex-18.3.3}, при $0<\alpha\ne 1$ должны
по теореме \ref{TH:ASYMP-COR-int-priznak-Cauchy} удовлетворять формуле
 \beq\label{Dirichlet-asymp}
\sum_{n=1}^N\frac{1}{n^\alpha}=\underbrace{\frac{N^{1-\alpha}-1}{1-\alpha}}_{\scriptsize
 \begin{matrix} \text{\rotatebox{90}{$=$}}\\
 \int\limits_1^N \frac{1}{x^{\alpha}}\, \d x \end{matrix}}+C_\alpha+\underset{N\to\infty}{\bold{o}}(1)
 \eeq
где $C_\alpha$ -- постоянные Эйлера \eqref{Dirichlet-Euler}. При $0<\alpha<1$
этот ряд расходится, поэтому выполняется соотношение
\eqref{ASYMP-COR:int-priznak-Cauchy-infty}
$$
\sum_{n=1}^N\frac{1}{n^\alpha}\underset{N\to\infty}{\sim}\frac{N^{1-\alpha}-1}{1-\alpha}.
$$
Его можно немного упростить, воспользовавшись теоремой \ref{TH:svyaz-<<-i-o},
$$
\frac{N^{1-\alpha}-1}{1-\alpha}=
\underbrace{\frac{N^{1-\alpha}}{1-\alpha}}_{\scriptsize\begin{matrix}\downarrow\\
\infty\end{matrix}}+\underbrace{\frac{1}{1-\alpha}}_{\scriptsize\begin{matrix}
\text{\rotatebox{90}{$=$}}
\\
\underset{N\to\infty}{\bold{o}}(1)\end{matrix}}
\underset{N\to\infty}{\sim}\frac{N^{1-\alpha}}{1-\alpha},
$$
и в результате мы получим
$$
\sum_{n=1}^N\frac{1}{n^\alpha}\underset{N\to\infty}{\sim}\frac{N^{1-\alpha}}{1-\alpha}\qquad
(0<\alpha<1).
$$
\eex

\end{multicols}\noindent\rule[10pt]{160mm}{0.1pt}

\btm\label{TH:form-Euler} Для любой гладкой функции $f$ на отрезке $[1,N]$,
$N\in\N$, справедлива формула:
 \beq\label{Euler-form}
\sum_{n=1}^N f(n)=f(1)+\int_1^N f(x)\ \d x+\int_1^N\{x\}\cdot f'(x)\ \d x
 \eeq
(здесь $\{x\}$ -- дробная часть $x$). \etm
 \bit{
\item[$\bullet$] Эта формула называется {\it формулой суммирования Эйлера}.
 }\eit

\brem Понятно, что выбор единицы в качестве нижнего предела суммирования здесь
несущественен -- сдвигом аргумента из \eqref{Euler-form} выводится, например,
что для любой гладкой функции $f$ на отрезке $[0,N]$ справедлива формула
$$
\sum_{n=0}^N f(n)=f(0)+\int_0^N f(x)\ \d x+\int_0^N\{x\}\cdot f'(x)\ \d x
$$
\erem

\bpr Рассмотрим функцию
$$
F(x)=\sum_{n=1}^{[x]} f(n)+\{x\}\cdot f(x)
$$
(здесь $[x]$ -- целая часть $x$).

1. Заметим, что функция $F$ кусочно-гладкая. Действительно, при $x\in (k,k+1)$,
$k\in\N$, выполняется равенство $[x]=k$, поэтому $\{x\}=x-[x]=x-k$, и значит
 \beq\label{dok-form-Euler}
F(x)=\sum_{n=1}^{[x]} f(n)+\{x\}\cdot f(x)=\sum_{n=1}^{k} f(n)+(x-k)\cdot f(x)
 \eeq
Отсюда видно, что на интервале $(k,k+1)$ функция $F$ гладкая и имеет конечные
односторонние производные на концах.

2. Заметим далее (и это неожиданно), что функция $F$ непрерывна. Мы уже
показали, что на каждом интервале $(k,k+1)$, $k\in\N$, она будет гладкой,
поэтому ее непрерывность нужно проверить только в целых точках. Зафиксируем
какую-нибудь точку $k\in\N$ и вычислим в ней значение $F$ и значения ее
односторонних пределов:
 \begin{align*}
 F(k)&=\sum_{n=1}^{[k]}
f(n)+\underbrace{\{k\}}_{\scriptsize\begin{matrix}\text{\rotatebox{90}{$=$}}\\
0\end{matrix}}\cdot f(k)=\sum_{n=1}^k f(n) \\
 \lim_{x\to k-0}F(x)&=\lim_{x\to k-0}\Big\{\sum_{n=1}^{[x]} f(n)+\{x\}\cdot
f(x)\Big\}=\sum_{n=1}^{k-1} f(n)+1\cdot f(k)=\sum_{n=1}^{k} f(n) \\
 \lim_{x\to k+0}F(x)&=\lim_{x\to k+0}\Big\{\sum_{n=1}^{[x]} f(n)+\{x\}\cdot
f(x)\Big\}=\sum_{n=1}^{k} f(n)+0\cdot f(k)=\sum_{n=1}^{k} f(n)
 \end{align*}
Видно, что эти три величины совпадают между собой, значит $F$ непрерывна в $k$.

3. Заметим, что из формулы \eqref{dok-form-Euler} следует, что на каждом
интервале $(k,k+1)$, $k\in\N$, производная функции $F$ имеет вид:
$$
F'(x)=\frac{\d}{\d x}\Big(\sum_{n=1}^k f(n)+(x-k)\cdot f(x)\Big)= \frac{\d}{\d
x}\Big((x-k)\cdot f(x)\Big)=1\cdot f(x)+(x-k)\cdot f'(x)=f(x)+\{x\}\cdot f'(x)
$$
Теперь по формуле перемещения непрерывной кусочно-гладкой функции
\eqref{peremeshenie-nepreryvnoi-kusochno-gladkoi-funktsii} получаем:
$$
\underbrace{F(N)}_{\scriptsize\begin{matrix}\text{\rotatebox{90}{$=$}}
\\ \sum\limits_{n=1}^N  f(n)+\underbrace{\{N\}}_{\tiny\begin{matrix}\text{\rotatebox{90}{$=$}}
\\ 0\end{matrix}}\cdot f(N)\\
 \text{\rotatebox{90}{$=$}}
\\ \sum\limits_{n=1}^N  f(n)
\end{matrix}}\kern-25pt-\kern-25pt\overbrace{F(1)}^{\scriptsize\begin{matrix}f(1)\\
\text{\rotatebox{90}{$=$}}\\
\sum\limits_{n=1}^1  f(n)+\overbrace{\{1\}}^{\tiny\begin{matrix}0\\
\text{\rotatebox{90}{$=$}}
\end{matrix}}\cdot f(1)\\
\text{\rotatebox{90}{$=$}}
\end{matrix}}\kern-25pt=\int_1^NF'(x)\ \d x=\int_1^N f(x)\ \d x+\int_1^N\{x\}\cdot f'(x)\ \d x
$$
$$
\Downarrow
$$
$$
\sum\limits_{n=1}^N  f(n)-f(1)=\int_1^N f(x)\ \d x+\int_1^N\{x\}\cdot f'(x)\ \d
x
$$
Остается перенести $f(1)$  в правую часть, и получится формула
\eqref{Euler-form}.
 \epr

Формула суммирования Эйлера сводит нахождение асимптотики суммы к асимптотике
интеграла. Чтобы объяснить, как ее применять, нам понадобится некое новое
понятие:

 \bit{
\item[$\bullet$]  {\it Функциями Эйлера} называются функции, определенные
индуктивным правилом
 \begin{align}
& E_0(x)=\{x\} \label{DEF:Euler-0}\\
& E_k(x)=\int_0^x E_{k-1}(t)\ \d t-x\cdot \int_0^1 E_{k-1}(t)\ \d t,
\label{DEF:Euler-k+1}
 \end{align}
Интегралы
 \beq\label{DEF:chisla-Euler}
\e_k=\int_0^1 E_k(t)\ \d t
 \eeq
называются {\it числами Эйлера}, и с их помощью формула \eqref{DEF:Euler-k+1}
записывается в виде:
 \beq\label{Euler-n+1-chisla-Euler}
E_k(x)=\int_0^x E_{k-1}(t)\ \d t-\e_{k-1}\cdot x
 \eeq
 }\eit

\bprop Каждая функция Эйлера $E_k$
 \bit{
\item[---] кусочно-гладка, а если $k>0$, то и непрерывна,

\item[---] периодична с периодом 1,
$$
E_k(x+1)=E_k(x),\qquad x\in\R
$$

\item[---] ограничена на $\R$,

\item[---] равна нулю в целых точках:
 \beq\label{Bernoulii=0-na-Z}
E_k(n)=0,\qquad n\in\Z
 \eeq

\item[---] имеет производную в нецелых точках:
 \beq\label{Bernoulii=0-na-Z}
E_k'(x)=E_{k-1}(x)-\e_{k-1},\qquad x\notin\Z
 \eeq

 }\eit
 \eprop
\bpr Кусочная гладкость и непрерывность при $k>0$ следует из теоремы
\ref{TH:stroenie-nepr-kusoch-glad-func}. Периодичность доказывается по
индукции: для $E_0(x)=\{x\}$ это верно, а если $E_{k-1}(x+1)=E_{k-1}(x)$, то
 \begin{multline*}
E_k(x+1)-E_k(x)=\int_0^{x+1} E_{k-1}(t)\ \d t-\e_{k-1}\cdot (x+1)-\l \int_0^x
E_{k-1}(t)\ \d t-\e_{k-1}\cdot x\r=\\=\int_x^{x+1} E_{k-1}(t)\ \d
t-\e_{k-1}=\int_0^1 E_{k-1}(t)\ \d t-\e_{k-1}=0.
 \end{multline*}
После этого ограниченность становится следствием непрерывности и периодичности.
А равенство нулю в целых точках следует из равенства нулю в нуле:
$$
E_k(0)=\int_0^0 E_{k-1}(t)\ \d t-\e_{k-1}\cdot 0=0.
$$
  \epr

Для приложений будет полезно вычислить несколько первых чисел Эйлера и функций
Эйлера на периоде $[0,1)$.

\bprop При $x\in[0,1)$ справедливы формулы:
 \begin{align}
& E_0(x)=x,  && \e_0=\frac{1}{2} \label{Euler-0}\\
& E_1(x)=\frac{x^2}{2}-\frac{x}{2},  && \e_1=-\frac{1}{12} \label{Euler-1}\\
& E_2(x)=\frac{x^3}{6}-\frac{x^2}{4}+\frac{x}{12},  && \e_2=0 \label{Euler-2}
 \end{align}
 \eprop
 \bpr
При $x\in[0,1)$ получаем цепочку:
$$
E_0(x)=\{x\}=x
$$
$$
\Downarrow
$$
$$
\e_0=\int_0^1 E_0(t)\ \d t=\int_0^1 t\ \d t=\frac{1}{2}
$$
$$
\Downarrow
$$
$$
E_1(x)=\int_0^x E_0(t)\ \d t-\e_0\cdot x=\int_0^x t\ \d t-\frac{x}{2}=
\frac{x^2}{2}-\frac{x}{2}
$$
$$
\Downarrow
$$
$$
\e_1=\int_0^1 E_1(t)\ \d t=\int_0^1 \l \frac{t^2}{2}-\frac{t}{2} \r\ \d t= \l
\frac{t^3}{6}-\frac{t^2}{4} \r\bigg|_0^1=\frac{1}{6}-\frac{1}{4}=-\frac{1}{12}
$$
$$
\Downarrow
$$
$$
E_2(x)=\int_0^x E_1(t)\ \d t-\e_1\cdot x=\int_0^x \l \frac{t^2}{2}-\frac{t}{2}
\r\ \d t+\frac{x}{12}= \frac{x^3}{6}-\frac{x^2}{4}+\frac{x}{12}
$$
$$
\Downarrow
$$
$$
\e_2=\int_0^1 E_2(t)\ \d t=\int_0^1 \l \frac{t^3}{6}-\frac{t^2}{4}+\frac{t}{12}
\r\ \d t= \l \frac{t^4}{24}-\frac{t^3}{12}+\frac{t^2}{24}
\r\bigg|_0^1=\frac{1}{24}-\frac{1}{12}+\frac{1}{24}=0
$$
\epr

Теперь перейдем к асимптотикам сумм.

\noindent\rule{160mm}{0.1pt}\begin{multicols}{2}

\bex В примере \ref{asimptotika-garm-ryada} мы нашли асимптотическую формулу
для частичной суммы гармонического ряда:
$$
\sum_{n=1}^N\frac{1}{n}=\ln N+C+\underset{N\to\infty}{\bold{o}}(1)\qquad\qquad
{\scriptsize\eqref{harmonic-asymp}}
$$
($C$ -- постоянная Эйлера \eqref{harmonic-Euler}).

Покажем теперь, как с помощью формулы Эйлера \eqref{Euler-form} можно эту
асимптотику уточнять. Взяв $f(x)=\frac{1}{x}$, мы по формуле Эйлера получим:
 \begin{multline}\label{asimp-sum-1/n-equiv-int}
\sum_{n=1}^N\frac{1}{n}=1+\int_1^N \frac{\d x}{x}-\int_1^N \frac{\{x\}}{x^2}\
\d x=\\=1+\ln N-\int_1^N \frac{\{x\}}{x^2}\ \d x.
 \end{multline}
Таким образом, нахождение асимптотики суммы свелось к нахождению асимптотики
интеграла
 \beq\label{asymp-int-(t)/t^2}
I(N)=\int_1^N\frac{\{t\}}{t^2}\ \d t=\int_1^N\frac{E_0(t)}{t^2}\ \d t
 \eeq
Он отличается от изучавшегося уже нами интеграла \eqref{asymp-int-sin(t)/t^2}
$$
\int_1^x\frac{\sin t}{t^2}\ \d t
$$
тем, что синус $\sin t$ в подынтегральном выражении заменен дробной частью
$\{t\}$, а верхний предел стал дискретным. Метод нахождения асимптотики для
\eqref{asymp-int-(t)/t^2} представляет собой модификацию метода для
\eqref{asymp-int-sin(t)/t^2}.

Коротко идею можно сформулировать так: {\it здесь тоже нужно интегрировать по
частям, но сначала нужно добиваться, чтобы в дифференциале стояла функция
Эйлера с подходящим индексом}. А это достигается добавлением и вычитанием чисел
Эйлера в периодической части подынтегрального выражения.

Продемонстрируем это на интеграле \eqref{asymp-int-(t)/t^2}. Прежде всего
замечаем, что существует конечный предел
 $$
A=\lim_{N\to\infty}I(N)=\int_1^\infty \frac{E_0(t)}{t^2}\ \d t
 $$
-- поскольку интеграл $\int_1^\infty \frac{E_0(t)}{t^2}\ \d t$ сходится
абсолютно:
 $$
\int_1^\infty \left|\frac{E_0(t)}{t^2}\right|\ \d t\le \int_1^\infty
\frac{1}{t^2}\ \d t=-\frac{1}{t}\Big|_1^\infty=1
 $$
Запомним это число $A$. Тогда:
 \begin{multline*}
I(N)=\int_1^N \frac{E_0(t)}{t^2}\ \d t=\\= \underbrace{\int\limits_1^\infty
\frac{E_0(t)}{t^2}\ \d t}_{\scriptsize\begin{matrix}
\text{\rotatebox{90}{$=$}}\\ A\end{matrix}}-\int_N^\infty \frac{E_0(t)}{t^2}\
\d t=\\
=A-\int_N^\infty \frac{\e_0+E_0(t)-\e_0}{t^2}\ \d
t=\\
=A-\kern-20pt\underbrace{\e_0}_{\scriptsize\begin{matrix} \eqref{Euler-0}\
\text{\rotatebox{90}{$=$}}\ \phantom{\eqref{Euler-0}}
\\
\frac{1}{2}\end{matrix}}\kern-20pt\int_N^\infty \frac{1}{t^2}\ \d t-
\int_N^\infty \frac{E_0(t)-\e_0}{t^2}\ \d t=\\= A+\frac{1}{2t}\ \Big|_N^\infty-
\int_N^\infty \frac{E_1'(t)}{t^2}\ \d t=\\= A-\frac{1}{2N}- \int_N^\infty
\frac{1}{t^2}\ \d
E_1(t)=\eqref{integrir-po-chastyam-dlya-nepr-kus-glad-funk}=\\= A-\frac{1}{2N}-
\frac{E_1(t)}{t^2}\ \Big|_N^\infty-2\int_N^\infty \frac{E_1(t)}{t^3}\ \d t =\\=
A-\frac{1}{2N}+\kern-14pt\underbrace{\frac{E_1(N)}{N^2}}_{\scriptsize\begin{matrix}
\eqref{Bernoulii=0-na-Z}\ \text{\rotatebox{90}{$=$}}\ \phantom{\eqref{Bernoulii=0-na-Z}}\\
0\end{matrix}}\kern-14pt- 2\underbrace{\int_N^\infty
\underbrace{\frac{E_1(t)}{t^3}}_{\underset{t\to\infty}{\bold{o}}\l\frac{1}{t^2}\r}\
\d t}_{\underset{N\to\infty}{\bold{o}}\l\frac{1}{N}\r} =\\=
A-\frac{1}{2N}+\underset{N\to\infty}{\bold{o}}\l\frac{1}{N}\r
 \end{multline*}
Подставим это в \eqref{asimp-sum-1/n-equiv-int}:
 \begin{multline*}
\sum_{n=1}^N\frac{1}{n}=1+\ln N-\int_1^N \frac{\{x\}}{x^2}\ \d x=\\= 1+\ln
N-A+\frac{1}{2N}+\underset{N\to\infty}{\bold{o}}\l\frac{1}{N}\r.
 \end{multline*}
Сравнив с \eqref{harmonic-asymp}, мы видим теперь, что постоянная $A$ связана с
постоянной Эйлера $C$ равенством
$$
1-A=C
$$
и поэтому полученную асимптотику можно переписать так:
 \beq\label{harmonic-asymp-1}
\sum_{n=1}^N\frac{1}{n}=\ln
N+C+\frac{1}{2N}+\underset{N\to\infty}{\bold{o}}\l\frac{1}{N}\r.
 \eeq
Это уточняет асимптотическую формулу \eqref{harmonic-asymp}.

Ее можно и дальше уточнять. Для этого нужно заметить, что, вычисляя асимптотику
для $I(N)$, мы попутно вывели формулу
$$
I(N)=A-\frac{1}{2N}-2\int_N^\infty \frac{E_1(t)}{t^3}\ \d t
$$
Подставим ее в \eqref{asimp-sum-1/n-equiv-int}:
 \begin{multline}\label{asimp-sum-1/n-equiv-int-1}
\sum_{n=1}^N\frac{1}{n}=1+\ln N-\int_1^N \frac{\{x\}}{x^2}\ \d x=\\= 1+\ln N-
A+\frac{1}{2N}+2\int_N^\infty \frac{E_1(t)}{t^3}\ \d t=\\=\ln
N+C+\frac{1}{2N}+2\int_N^\infty \frac{E_1(t)}{t^3}\ \d t
 \end{multline}
Отдельно найдем асимптотику последнего интеграла:
 \begin{multline*}
I_2(N)= 2\int_N^\infty \frac{E_1(t)}{t^3}\ \d t =\\= 2\int_N^\infty
\frac{\e_1+E_1(t)+\e_1}{t^3}\ \d t =\\=
2\kern-20pt\underbrace{\e_1}_{\scriptsize\begin{matrix} \eqref{Euler-1}\
\text{\rotatebox{90}{$=$}}\ \phantom{\eqref{Euler-1}}
\\
-\frac{1}{12}\end{matrix}}\kern-20pt\int_N^\infty \frac{1}{t^3}\ \d
t+2\int_N^\infty \frac{E_1(t)-\e_1}{t^3}\ \d t =\\= \frac{2}{2\cdot 12}\cdot
\frac{1}{t^2}\Big|_N^\infty+2\int_N^\infty \frac{\d E_2(t)}{t^3} =\\=
-\frac{1}{12 N^2}+\frac{2E_2(t)}{t^3}\Big|_N^\infty+6\int_N^\infty
\frac{E_2(t)}{t^4}\ \d t =\\= -\frac{1}{12
N^2}-\kern-12pt\underbrace{\frac{2E_2(N)}{N^3}}_{\scriptsize\begin{matrix}
\eqref{Bernoulii=0-na-Z}\ \text{\rotatebox{90}{$=$}}\ \phantom{\eqref{Bernoulii=0-na-Z}}\\
0\end{matrix}}\kern-12pt+6\underbrace{\int_N^\infty
\underbrace{\frac{E_2(t)}{t^4}}_{\underset{t\to\infty}{\bold{o}}\l\frac{1}{t^3}\r}\
\d t}_{\underset{N\to\infty}{\bold{o}}\l\frac{1}{N^2}\r} =\\= -\frac{1}{12
N^2}+\underset{N\to\infty}{\bold{o}}\l\frac{1}{N^2}\r
 \end{multline*}
Подставив это в \eqref{asimp-sum-1/n-equiv-int-1}, получим:
 \begin{multline*}
\sum_{n=1}^N\frac{1}{n}=\ln N+C+\frac{1}{2N}+2\int_N^\infty \frac{E_1(t)}{t^3}\
\d t=\\=\ln N+C+\frac{1}{2N}-\frac{1}{12
N^2}+\underset{N\to\infty}{\bold{o}}\l\frac{1}{N^2}\r
 \end{multline*}
Мы получили уточнение для формулы \eqref{harmonic-asymp-1}:
 \beq\label{harmonic-asymp-2}
\sum_{n=1}^N\frac{1}{n}=\ln N+C+\frac{1}{2N}+-\frac{1}{12
N^2}+\underset{N\to\infty}{\bold{o}}\l\frac{1}{N^2}\r.
 \eeq

И так далее. Применяя функции Эйлера нужное число раз, можно добиться любой
точности асимптотического разложения.
 \eex

\bex Найдем асимптотику суммы
$$
\sum_{n=1}^N\ln n
$$
По формуле Эйлера \eqref{Euler-form} получаем:
 \begin{multline*}
\sum_{n=1}^N\ln n =\overbrace{\ln 1}^{\scriptsize\begin{matrix}0\\
\text{\rotatebox{90}{$=$}}\end{matrix}}+\kern-10pt\overbrace{\int_1^N\ln x\ \d
x}^{\scriptsize\begin{matrix} N\cdot\ln N-(N-1)
\\
\text{\rotatebox{90}{$=$}}
\\
x\cdot\ln x\Big|_1^N-\int\limits_1^Nx\ \d\ln x \\
\text{\rotatebox{90}{$=$}}\end{matrix}}\kern-10pt+\int_1^N\frac{\{x\}}{x}\ \d
x=\\=  N\cdot\ln N-N+1+\int_1^N\frac{\{x\}}{x}\ \d x =\\= N\cdot\ln
N-N+1+\int_1^N\frac{E_0(x)}{x}\ \d x=\\=
 N\cdot\ln
N-N+1+\int_1^N\frac{\e_0+E_0(x)-\e_0}{x}\ \d x=\\= N\cdot\ln
N-N+1+\\+\kern-20pt\underbrace{\e_0}_{\scriptsize\begin{matrix}
\eqref{Euler-0}\ \text{\rotatebox{90}{$=$}}\ \phantom{\eqref{Euler-0}}
\\
\frac{1}{2}\end{matrix}}\kern-20pt\int_1^N\frac{\d
x}{x}+\int_1^N\frac{E_0(x)-\e_0}{x}\ \d x=\\= N\cdot\ln N-N+1+\frac{1}{2}\ln
N+\int_1^N\frac{\d E_1(x)}{x}\ \d x=\\= N\cdot\ln N-N+1+\frac{1}{2}\ln
N+\underbrace{\frac{E_1(x)}{x}\Big|_1^N}_{\scriptsize\begin{matrix}
\text{\rotatebox{90}{$=$}}\\
\frac{E_1(N)}{N}-\frac{E_1(1)}{1}\\
\eqref{Bernoulii=0-na-Z}\ \text{\rotatebox{90}{$=$}}\ \phantom{\eqref{Bernoulii=0-na-Z}}\\
0\end{matrix}} +\\+\int_1^N\frac{E_1(x)}{x^2}\ \d x= N\cdot\ln
N-N+1+\frac{1}{2}\ln N+\\+\int_1^N\frac{E_1(x)}{x^2}\ \d x
 \end{multline*}
Теперь заметим, что поскольку функция $E_1$ ограничена на $\R$, сходится
абсолютно интеграл $\int_1^\infty\frac{E_1(x)}{x^2}\ \d x$:
$$
\int_1^\infty\left|\frac{E_1(x)}{x^2}\right|\ \d x\le
\int_1^\infty\frac{M}{x^2}\ \d x<\infty
$$
Обозначим буквой $A$ его значение:
$$
A=\int_1^\infty\frac{E_1(x)}{x^2}\ \d x
$$
Тогда
 \begin{multline*}
\sum_{n=1}^N\ln n = N\cdot\ln N-N+1+\frac{1}{2}\ln
N+\\+\int_1^N\frac{E_1(x)}{x^2}\ \d x= N\cdot\ln N-N+1+\frac{1}{2}\ln N+\\+
\underbrace{\int_1^\infty\frac{E_1(x)}{x^2}\ \d x}_{\scriptsize\begin{matrix}
\text{\rotatebox{90}{$=$}}\\ A
\end{matrix}}-\underbrace{\int_N^\infty\underbrace{\frac{E_1(x)}{x^2}}_{\scriptsize\begin{matrix}
\text{\rotatebox{90}{$=$}}\\ \underset{x\to\infty}{\bold{o}}\l\frac{1}{x}\r
\end{matrix}}\ \d x}_{\scriptsize\begin{matrix}
\text{\rotatebox{90}{$=$}}\\ \underset{N\to\infty}{\bold{o}}(1)
\end{matrix}}
 \end{multline*}
И в качестве ответа получаем:
$$
\sum_{n=1}^N\ln n = N\cdot\ln N-N+\frac{1}{2}\ln
N+1-A+\kern-5pt\underset{N\to\infty}{\bold{o}}(1)
$$
($A$ определено выше). Понятно, что эту асимптотику можно уточнять.
   \eex

\ber Уточните асимптотическую формулу \eqref{Dirichlet-asymp} для частичных
сумм ряда Дирихле:
$$
\sum_{n=1}^\infty\frac{1}{n^\alpha},\qquad 0<\alpha\ne 1
$$
\eer

\paragraph{Точное вычисление сумм.}

Неожиданное следствие формулы Эйлера состоит в том, что, если под знаком суммы
стоит многочлен, то мы можем в точности сосчитать частичную сумму ряда.

\bex Выведем формулу для суммы
$$
\sum_{n=1}^N n^2
$$
По формуле Эйлера получаем:
 \begin{multline*}
\sum_{n=1}^N n^2=1+\int_1^N x^2\ \d x+2\int_1^N \{x\}\cdot x\ \d x=\\=
1+\frac{N^3-1}{3}+2\int_1^N \Big(\e_0+E_0(x)-\e_0\Big)\cdot x\ \d x=\\=
1+\frac{N^3-1}{3}+2\kern-20pt\underbrace{\e_0}_{\scriptsize\begin{matrix}
\eqref{Euler-0}\ \text{\rotatebox{90}{$=$}}\ \phantom{\eqref{Euler-0}}
\\
\frac{1}{2}\end{matrix}}\kern-20pt\int_1^N  x\ \d x+\\+2\int_1^N
\Big(E_0(x)-\e_0\Big)\cdot x\ \d
x=\\
=1+\frac{N^3-1}{3}+\frac{N^2-1}{2}+2\int_1^N x\ \d E_1(x)=\\
=1+\frac{N^3-1}{3}+\frac{N^2-1}{2}+\\+\underbrace{2x\cdot
E_1(x)\Big|_1^N}_{\scriptsize\begin{matrix}
\text{\rotatebox{90}{$=$}}\\
2N\cdot E_1(N)-2\cdot E_1(1)\\
\eqref{Bernoulii=0-na-Z}\ \text{\rotatebox{90}{$=$}}\ \phantom{\eqref{Bernoulii=0-na-Z}}\\
0\end{matrix}}- 2\underbrace{\int_1^N E_1(x)\ \d x}_{\scriptsize\begin{matrix}
\text{\rotatebox{90}{$=$}}\\
(N-1)\int_0^1 E_1(x)\ \d
x\\
\text{\rotatebox{90}{$=$}}\\
(N-1)\cdot\e_1 \\
\eqref{Euler-1}\ \text{\rotatebox{90}{$=$}}\ \phantom{\eqref{Euler-1}}\\
-\frac{N-1}{12}
\end{matrix}} =\\=1+\frac{N^3-1}{3}+\frac{N^2-1}{2}
+\frac{N-1}{6}=\\= \frac{N^3}{3}+\frac{N^2}{2}+\frac{N}{6}
 \end{multline*}
Вывод:
 \beq
\sum_{n=1}^N n^2=\frac{N^3}{3}+\frac{N^2}{2}+\frac{N}{6}
 \eeq
 \eex

\ber Найдите формулы для сумм:
 \biter{
\item[1)] $\sum_{n=1}^N n^3$,

\item[2)] $\sum_{n=1}^N n^4$.
 }\eiter
 \eer

\section{Формула Стирлинга}

Закончить разговор об асимптотических методах будет поучительно знаменитой {\it
формулой Стирлинга}, описывающей асимптотическое поведение факториала:
 \beq\label{Stirling}
n!\underset{n\to\infty}{\sim}\frac{n^n\cdot\sqrt{2\pi n}}{e^n}
 \eeq
Для ее доказательства нам понадобится еще одно соотношение, называемое {\it
формулой Валлиса}:
 \beq\label{Wallis}
\frac{1}{2k+1}\cdot\l\frac{(2k)!!}{(2k-1)!!}\r^2\underset{k\to\infty}{\longrightarrow}\frac{\pi}{2}
 \eeq
\bpr[Доказательство формулы Валлиса] Здесь применяется формула
\eqref{int_0^(pi/2)-sin^n-cos^n}. Обозначим
$$
x_k=\l\frac{(2k)!!}{(2k-1)!!}\r^2
$$
Тогда:
$$
\sin^{2k+1}x<\sin^{2k}x<\sin^{2k-1}x,\qquad x\in\left[0,\frac{\pi}{2}\right]
$$
$$
\Downarrow
$$
 \begin{multline*}
\int_0^\frac{\pi}{2}\sin^{2k+1}x\ \d x<\int_0^\frac{\pi}{2}\sin^{2k}x\ \d x<\\<
\int_0^\frac{\pi}{2}\sin^{2k-1}x\ \d x
 \end{multline*}
$$
\phantom{\scriptsize \eqref{int_0^(pi/2)-sin^n-cos^n}}\ \Downarrow\
{\scriptsize \eqref{int_0^(pi/2)-sin^n-cos^n}}
$$
$$
\frac{(2k)!!}{(2k+1)!!}<\frac{(2k-1)!!}{(2k)!!}\cdot\frac{\pi}{2}<\frac{(2k-2)!!}{(2k-1)!!}
$$
$$
\Downarrow
$$
$$
\frac{1}{2k+1}\cdot\underbrace{\l\frac{(2k)!!}{(2k-1)!!}\r^2}_{x_k}<\frac{\pi}{2}<\frac{1}{2k}\cdot
\underbrace{\l\frac{(2k)!!}{(2k-1)!!}\r^2}_{x_k}
$$
$$
\Downarrow
$$
$$
\frac{x_k}{2k+1}<\frac{\pi}{2}<\frac{x_k}{2k}
$$
$$
\Downarrow
$$
$$
\frac{2k}{2k+1}\cdot\frac{\pi}{2}<\frac{2k}{2k+1}\cdot\underbrace{\frac{x_k}{2k}}_{\scriptsize\begin{matrix}
\text{\rotatebox{90}{$<$}}\\ \frac{\pi}{2} \end{matrix}}
=\frac{x_k}{2k+1}<\frac{\pi}{2}
$$
$$
\Downarrow
$$
$$
\underbrace{\frac{2k}{2k+1}}_{\scriptsize\begin{matrix} \downarrow \\ 1
\end{matrix}}\cdot\frac{\pi}{2}<\frac{x_k}{2k+1}<\frac{\pi}{2}
$$
$$
\Downarrow
$$
$$
\frac{x_k}{2k+1}\underset{k\to\infty}{\longrightarrow}\frac{\pi}{2}
$$
\epr

\bpr[Доказательство формулы Стирлинга] Обо\-значим
$$
x_n=\frac{n!\cdot e^n}{n^{n+\frac{1}{2}}}
$$
Из маклореновского разложения для логарифма \eqref{mclaurent-ln(1+x)} получаем
такую цепочку:
$$
\ln(1+x)=\sum_{m=1}^\infty (-1)^{m-1}\cdot\frac{x^m}{m},\qquad x\in(0,1).
$$
$$
\Downarrow
$$
 \begin{multline*}
\ln\frac{1+x}{1-x}=\ln(1+x)-\ln(1-x)=\\= \sum_{m=1}^\infty
(-1)^{m-1}\cdot\frac{x^m}{m}-\sum_{m=1}^\infty
\l-\frac{x^m}{m}\r=\\=\underbrace{2x}_{\scriptsize\begin{matrix}\text{\rotatebox{90}{$<$}}\\
0\end{matrix}}\cdot
\sum_{k=0}^\infty \underbrace{\frac{x^{2k}}{2k+1}}_{\scriptsize\begin{matrix}\text{\rotatebox{90}{$<$}}\\
0\end{matrix}}>2x\cdot \underbrace{1}_{k=0}=2x
 \end{multline*}
$$
\phantom{\text{\scriptsize подстановка: $x=\frac{1}{2n+1}$}}
\quad\Downarrow\quad \text{\scriptsize подстановка: $x=\frac{1}{2n+1}$}
$$
$$
\ln\l 1+\frac{1}{n}\r
=\ln\frac{1+\frac{1}{2n+1}}{1-\frac{1}{2n+1}}>\frac{2}{2n+1}=\frac{1}{n+\frac{1}{2}}
$$
$$
\Downarrow
$$
$$
\ln\l 1+\frac{1}{n}\r^{n+\frac{1}{2}}=\l n+\frac{1}{2}\r\cdot \ln\l
1+\frac{1}{n}\r>1
$$
$$
\Downarrow
$$
$$
\l 1+\frac{1}{n}\r^{n+\frac{1}{2}}>e
$$
$$
\Downarrow
$$
 \begin{multline*}
\frac{x_n}{x_{n+1}}=\frac{\frac{n!\cdot
e^n}{n^{n+\frac{1}{2}}}}{\frac{(n+1)!\cdot
e^{n+1}}{(n+1)^{n+1+\frac{1}{2}}}}=\frac{(n+1)^{n+1+\frac{1}{2}}}{(n+1)\cdot
n^{n+\frac{1}{2}}\cdot
e}=\\=\frac{1}{e}\cdot\frac{(n+1)^{n+\frac{1}{2}}}{n^{n+\frac{1}{2}}}=\frac{1}{e}\cdot\l
1+\frac{1}{n}\r^{n+\frac{1}{2}}>1
 \end{multline*}
$$
\Downarrow
$$
$$
x_n>x_{n+1}
$$
Мы получили, что $\{x_n\}$ -- убывающая и ограниченная снизу нулем
последовательность. Значит, у нее имеется предел:
$$
x_n=\frac{n!\cdot e^n}{n^{n+\frac{1}{2}}}\underset{n\to\infty}{\longrightarrow}
a
$$
Это соотношение можно понимать как эквивалентность
 \beq\label{stirling-prelim}
n!\underset{n\to\infty}{\sim} a\cdot \frac{n^{n+\frac{1}{2}}}{e^n},
 \eeq
в которой $a\in\R$ -- пока неизвестное число. Для этого и нужна доказанная выше
формула Валлиса \eqref{Wallis}:
 \begin{multline*}
\frac{\pi}{2}\underset{\infty\gets k}{\longleftarrow}
\frac{1}{2k+1}\cdot\l\frac{(2k)!!}{(2k-1)!!}\r^2=\eqref{nechetnyi-dvoinoi-faktorial}=\\=
\frac{1}{2k+1}\cdot\l\frac{(2k)!!\cdot(2k-2)!!}{(2k-1)!}\r^2=\eqref{chetnyi-dvoinoi-faktorial}=\\=
\frac{1}{2k+1}\cdot\l\frac{2^k\cdot k!\cdot 2^{k-1}\cdot
(k-1)!}{(2k-1)!}\r^2=\\= \frac{1}{2k+1}\cdot\l\frac{2^k\cdot k!\cdot 2^k\cdot
k!}{2k\cdot (2k-1)!}\r^2=\\=
\frac{2^{4k}}{2k+1}\cdot\frac{(k!)^4}{\Big((2k)!\Big)^2}\overset{\eqref{stirling-prelim}}{\underset{k\to\infty}{\sim}}
\\ \underset{k\to\infty}{\sim}
\frac{2^{4k}}{2k}\cdot\frac{\l a\cdot \frac{k^{k+\frac{1}{2}}}{e^k} \r^4}{\l
a\cdot \frac{(2k)^{2k+\frac{1}{2}}}{e^{2k}}\r^2}=
\frac{2^{4k}}{2k}\cdot\frac{a^4\cdot \frac{k^{4k+2}}{e^{4k}}}{a^2\cdot
\frac{(2k)^{4k+1}}{e^{4k}}}=\\= \frac{2^{4k}}{2k}\cdot\frac{a^2\cdot
k^{4k+2}}{(2k)^{4k+1}}= \frac{2^{4k}}{2k}\cdot\frac{a^2\cdot
k^{4k+2}}{2^{4k+1}\cdot k^{4k+1}}=\frac{a^2}{4}
 \end{multline*}
$$
\Downarrow
$$
$$
2\pi=a^2
$$
$$
\Downarrow
$$
$$
\sqrt{2\pi}=a
$$
Это превращает \eqref{stirling-prelim} в \eqref{Stirling}. \epr

\end{multicols}\noindent\rule[10pt]{160mm}{0.1pt}

\end{document}